\documentclass[graybox,envcountchap]{svmono}
\usepackage{helvet}
\usepackage{courier}
\usepackage{mathptmx}
\usepackage{amssymb,amscd,bm}  
\usepackage{mathtools}
\usepackage{tikz-cd}
\usepackage{pictexwd,dcpic}

\usepackage{parskip}

\usepackage{graphicx}
\usepackage{epsf}
\usepackage{epsfig}
\usepackage{bbm}
\usepackage[colorlinks=true]{hyperref}
\usepackage{marginnote}
\usepackage{mathrsfs}
\usepackage{pdfsync}
\usepackage{setspace}

\usepackage{nomencl}
\usepackage{bbding}
\usepackage{wasysym}
\usepackage{type1cm}         

\usepackage{makeidx}         
\usepackage{multicol}        
\usepackage[bottom]{footmisc}

\providecommand{\norm}[1]{\lVert#1\rVert}
\providecommand{\abs}[1]{\lvert#1\rvert}

\newcommand{\ssc}{\text{sc}}

\renewcommand{\epsilon}{\varepsilon}

\newcommand{\what}{\widehat}
\newcommand{\wh}{\widehat}

\newcommand{\wt}{\widetilde}
\newcommand{\wtilde}{\widetilde}

\newcommand{\gr}{\text{graph}}

\newcommand{\ov}{\overline}

\newcommand{\Fred}{\operatorname{Fred}}

\newcommand{\cl}{\operatorname{cl}}

\newcommand{\id}{\operatorname{id}}
\newcommand{\ind}{\operatorname{ind}}

\newcommand{\real}{\operatorname{Re}}

\newcommand{\supp}{\operatorname{supp}}

\newcommand{\R}{{\mathbb R}}

\renewcommand{\qedsymbol}{$\square$}

\newcommand{\mo}{\mathcal O}
\providecommand{\ker}[1]{$\text{ker}\ {#1}$}

\providecommand{\coker}[1]{\text{coker}\ {#1}}

\newcommand{\N}{{\mathbb N}}
\newcommand{\Q}{{\mathbb Q}}
\newcommand{\pr}{\text{pr}}

\newcommand{\bls}{\boldsymbol}
 \newcommand{\sg}{S_{\Theta}^{\text{g}}}
\renewcommand{\epsilon}{\varepsilon}
\newcommand{\morp}{\text{mor}}
\newcommand{\sba}{S_{\Theta}^{\text{b}}}
\newcommand{\co}{\circ}

\newenvironment{Myitemize}{%
\begin{itemize}}{\end{itemize}}

\providecommand{\ker}[1]{$\text{ker}\ {#1}$}

\newcommand{\bssc}{\text{Sc}}
\newcommand{\bx}{{\bf X}}

\makeindex

\smartqed

\author{Helmut Hofer\\ Krzysztof Wysocki\\  Eduard Zehnder}

\title{Polyfold and Fredholm Theory}
\begin{document}
\maketitle

\frontmatter

%
%
%

\begin{dedication}
Two lives were tragically lost while writing this book.
My coauthor and friend Kris Wysocki passed away on February 16, 2016 and my son David S. Hofer on March 1, 2016.\linebreak
\noindent I dedicate this book to them.
\vspace{\baselineskip}
\begin{flushright}\noindent
 \hfill {\it Helmut Hofer}\\
\end{flushright}
\phantom{X}\par

\phantom{X}\par

\phantom{X}\par

\noindent I dedicate this book to my friend and coauthor Kris Wysocki who lost his long fight against cancer on February 16, 2016, and to his wife and our friend, Beata Wysocka.\linebreak
We all miss Kris.

\vspace{\baselineskip}
\begin{flushright}\noindent
 \hfill {\it Eduard Zehnder}\\
\end{flushright}

\end{dedication}
%
%

\preface
Polyfold theory is a systematic and efficient approach providing a language and a large body of results for dealing with larger classes of nonlinear elliptic equations involving compactness
(bubbling-off), smoothness (varying domains and targets) and (geometric) transversality issues.
The abstract theory, proposed in the series of papers \cite{HWZ2,HWZ3,HWZ3.5},\  
generalizes  differential geometry and  nonlinear Fredholm theory to a class of spaces which are much more general than (Banach) manifolds. These spaces may have (locally) varying dimensions and are described locally by retracts on scale Banach spaces, replacing the open sets of Banach spaces in the familiar local description of manifolds. The theory also allows to equip certain kinds of categories with smooth structures.
Current applications include the study of moduli spaces of pseudoholomorphic curves in symplectic geometry, f.e. the construction of Symplectic Field Theory (SFT),
see \cite{HWZ5,H2,FH2,FH3}.\par

The theory grew out of the attempts to find an appropriate framework to develop  SFT,  a general theory of symplectic invariants outlined in \cite{EGH}.
The SFT constructs invariants of symplectic cobordisms by analyzing the structure of solutions of elliptic partial differential equations  of Cauchy-Riemann type, from Riemann surfaces into compact symplectic manifolds.  The partial differential equations are defined on varying domains and map into varying targets. The occurring singular limits, like bubbling off phenomena, give rise to serious compactness and transversality problems. Yet, these phenomena are needed to derive the underlying rich algebraic structure of SFT. Therefore 
they have to be `embraced' and accurately accounted for. The analytical difficulties in dealing with SFT become apparent
in the series of papers \cite{H1,HWZ-Convex,HWZ-tight,HWZ7.1} culminating in the SFT-compactness paper \cite{BEHWZ}.\par

This book is written as a reference volume for polyfold theory and the accompanying Fredholm theory.  We do not give any applications, since they are developed elsewhere.
For example in \cite{FH2} the theory is used to develop methods to equip specific groupoidal categories with smooth structures and   to apply these ideas to construct the polyfolds of  SFT. Using the polyfold setup for SFT developed in \cite{FH2} the transversality 
from the present book is applied in \cite{FH3} to construct SFT. The publication \cite{HWZ5} illustrates the polyfold theory by its application to Gromov-Witten theory and makes  transparent how the theory deals with the occurring issues. Another application is given
in \cite{Suhr} to the Weinstein conjecture. The paper \cite{HWZ8.7} illustrates part of  the theory by examples relevant for applications. 

A reader who  understands the inner-workings of the machinery should feel comfortable to use more stream-lined versions, as for example developed in   \cite{FH2}, which exploit
the full strength of the material described here.
 It is clear that the polyfold theory, which grew out of the  structurally very rich SFT-example, also applies to other classes or families of nonlinear elliptic partial differential equations as well. \par

In Part I of this book, based on \cite{HWZ2,HWZ3},  we develop  the Fredholm theory in a class
of spaces called M-polyfolds.  These spaces can be viewed as a generalization of the notion of a manifold (finite or infinite dimensional).
In  Part II M-polyfolds, based on \cite{HWZ3.5},    will be generalized  to a class of spaces called ep-groupoids, which can be viewed as the polyfold generalization
of \'etale proper Lie groupoids. This generalization is useful in dealing with problems having local symmetries, which is needed
 in more  advanced applications.  
 In  Part III we develop, following  \cite{HWZ3,HWZ3.5,HWZ7,HWZ7err}, the Fredholm theory in ep-groupoids. Transversality and symmetry are antagonistic concepts, since perturbing a problem to achieve 
 transversality might not be possible when we also want to preserve the symmetries.  Although it  seems hopeless  in many cases to achieve transversality, there 
 is  a version of transversality theory over the rational numbers, which is based on a perturbation theory using set-valued perturbations.  Such an idea was introduced
 by Fukaya and Ono, \cite{FO}, in  their  Kuranishi framework for the Arnold conjectures. The use of set-valued operators in nonlinear analysis is older, and a description can be found in \cite{Zeidler}. 
 It seems that the latter was never used to develop a transversality theory in a context of symmetries. We merge these ideas into a nice formalism which was in a very special case
 suggested in \cite{CRS}.
 In some sense  the transversality theory in a symmetric setting replaces locally the equivariant problem by a symmetric weighted family of problems. \par

 A by-product of our considerations is a Fredholm theory,  generalizing what in classical terms would be 
 a Fredholm theory on (Banach) orbifolds with boundary with corners.
 In summary, Parts {I}-{III} describe a wide array of nonlinear functionalanalytic tools to study perturbation and transversality questions
for a large class of so-called Fredholm sections and Fredholm section functors in ep-groupoids.\par

The whole theory can be generalized even a step further to equip certain categories (groupoidal categories) with smooth structures and to develop a theory of Fredholm functors.  This is carried out in Part IV. 
An outline of such ideas in the case of Gromov-Witten theory is given in \cite{H2}. Another example is the construction of  SFT  which    most elegantly is described within such a framework. These ideas are carried out in \cite{FH2,FH3}.\par

The series of videos \cite{H-Video}, \cite{A_Video}, \cite{Wehrheim_video},  illustrate the motivation behind the polyfold theory. We also recommend \cite{FFGW} for the intuitive ideas involved in the polyfold theory.\par

Other approaches to the type of problems considered in  applications were put forward in \cite{FO,FOOO}, \cite{LiuT,LT},
\cite{McW,McW2,McW3,McW4},  \cite{Yang, Yang2}, and \cite{Pardon}.
Particularly the carefully written papers by McDuff and Wehrheim are a very good introduction to the finite-dimensional 
Kuranishi-type approaches. The work of Yang is concerned with showing the relationship between the polyfold and Kuranish type 
approaches.

\vspace{\baselineskip}
\begin{flushright}\noindent
Princeton and Z\"urich,\hfill {\it Helmut Hofer}\\
July 2017\hfill {\it Eduard Zehnder}\\
\end{flushright}

%
%

\extrachap{Acknowledgements}

The first author was  partially supported by the NSF grants DMS-0603957 and DMS-1104470, the second author  was partially supported by  the NSF grant DMS-0906280.  The second and the third author would like to thank the Institute for Advanced Study (IAS) in Princeton for the support and hospitality, the second author  would like to thank the Forschungsinstitut f{\"u}r Mathematik (FIM) in Zurich for the support and hospitality.
The authors would like to thank  Joel Fish, Dusa McDuff,  and Katrin Wehrheim for many useful and enlightening discussions.
The first author also would like to thank  Michael Jemison, Jake Solomon,   Dingyu Yang and  the participants of the workshop on polyfolds at Pajaro Dunes in August 2012 for their valuable feedback.

\tableofcontents

\mainmatter

%
%
%

\begin{partbacktext}
\part{Basic Theory in M-Polyfolds}
\noindent In Part I we introduce the notion of an sc-Banach space.  It can be viewed as a Banach space $E$ with an additional structure  provided by a compact scale of Banach spaces. 
Although such gadgets are known from interpolation theory the new viewpoint interprets them as some kind of smooth structure on $E$.  
Then we introduce a new notion of differentiability and the smooth maps in this context are called sc-smooth maps. There are many more sc-smooth maps
than in the usual context of classical differentiability, but still the maps have enough structure to recognize boundaries with corners, unlike homeomorphisms
which only recognize boundaries.  Also the chain rule holds
$$
T(g\circ f)=Tg\circ Tf,
$$
where $Tg$ and $Tf$ are well-defined tangent maps.
 Of particular interest are the sc-smooth maps $r:U\rightarrow U$ which satisfy $r\circ =r$.  They are called sc-smooth retractions.
 
If a subset $O$ of the sc-Banach space $E$  can be written as $O=r(U)$ for some sc-smooth retraction it turns out that the definition $TO:=Tr(TU)$ 
is independent of the sc-smooth retraction onto $O$ and defines what is called the tangent space.  At this point we can take the sc-smooth retracts $O$ as the local models for an extended differential geometry.
One should note that sc-smooth retractions also exist in the classical context, but their images are always submanifolds, so that taking them as local model gives nothing new.
In our context the sc-smooth retracts can be very general. For example there exist finite-dimensional ones with local varying dimensions, which however cannot be embedded smoothly 
into a  finite-dimensional vector spaces.
The new sc-smooth spaces which are build on the new local models are studied in detail.  They can be viewed as a generalization of the standard smooth manifolds.
We shall study sc-smooth subspaces and sc-smooth bundles.  It is clear that taking any classical text book on differentiable manifolds many of the concepts 
can be generalized, and that in addition new concepts arise.  

Since the notion of differentiability is rather weak smooth maps in general do not satisfy an implicit function theorem. However, we are able to identify a class
of smooth maps for which it holds. These will be the so called sc-Fredholm sections. We study these in detail, derive a perturbation and transversality theory
and develop a theory of orientation.  With all these ingredients we obtain the `Fredholm Package' which provides for a Sard-Smale type theory, but in the context 
of much more general spaces. The notion of smoothness is so general that in applications bubbling-off phenomena can be described within this abstract framework.  The basic core ideas go back to \cite{HWZ2,HWZ3}, and 
we refer the reader to \cite{HWZ8.7} and \cite{FFGW} for illustrations of the ideas in concrete contexts.

\end{partbacktext}

\chapter{Sc-Calculus}
The basic concept is the  sc-structure on a Banach space.

\section{Sc-Structures and Differentiability}
We begin with the  linear sc-theory. 
\begin{definition}\label{sc-structure}
A {\bf sc-structure}  \index{D- Sc-Structure} (or scale structure) on a  Banach space $E$ consists of 
a decreasing sequence $(E_m)_{m\geq 0}$ of Banach spaces,   
$$
E=E_0\supset E_1\supset E_2\supset \ldots  ,
$$
such that  the following two conditions are satisfied, 
\begin{itemize}
\item[(1)]\ The inclusion operators $E_{m+1}\rightarrow E_m$ are  compact.
\item[(2)]\ The intersection $E_\infty:=\bigcap_{i\geq 0} E_i$\index{$E_i,\ E_{\infty}$} is dense in every $E_m$.
\end{itemize}
\qed
\end{definition}

Sc-structures, where ``sc'' is short for scale, are known from  linear interpolation theory, see \cite{Tr}. However, our  interpretation is that of a smooth structure.    It should be noted that scales have been used in geometric settings
  before; see for example the work by D. Ebin,  \cite{Ebin},  and
  H. Omori, \cite{Omori}\footnote{The first author thanks T. Mrowka for pointing out these two references.}.
The latter of these is somewhat closer to the viewpoint in our paper.
  However Omori does not use \emph{compact} scales, and this turns out
  to be a crucial condition for our applications. 
Without the compactness assumption the maps we shall be considering do not satisfy the chain rule.
 Moreover, one does not obtain the
  new local models (the sc-smooth retracts)  for a generalized differential geometry in which 
  bubbling-off and trajectory-breaking is a smooth phenomenon.

In the following, a sc-Banach space or a sc-smooth Banach space stands for a Banach space equipped with a sc-structure.
A finite-dimensional Banach space $E$  has precisely one sc-structure, namely the constant structure $E_m=E$ for all $m\geq0$.  If $E$ is an infinite-dimensional Banach space, the constant structure is not a sc-structure, because it fails property (1).  We shall see that the sc-structure leads to interesting new phenomena in nonlinear analysis.

Points in $E_\infty$ are called {\bf smooth points}, points in $E_m$ are called points of {\bf regularity $m$}\index{Points of regularity $m$}. A subset $A$ of a sc-Banach space $E$ inherits a filtration  $(A_m)_{m\geq 0}$ defined by $A_m=A\cap E_m$. It is, of course, possible that $A_\infty=\emptyset $.   We adopt the convention that 
$A^k$ stands for the set $A_k$ equipped with the induced filtration  
$$(A^k_m)_{m\geq 0}=(A_{k+m})_{m\geq 0}.$$
The direct sum $E\oplus F$ of sc-Banach spaces is a sc-Banach space, whose sc-smooth structure is defined by $(E\oplus F)_m:=E_m\oplus F_m$ for  all $m\geq 0$.

\begin{example}\label{sc-example1}
An example of a sc-Banach space, which is relevant in applications, is as follows.  We choose a strictly increasing sequence $(\delta_m)_{m\geq 0}$  of real numbers starting with $\delta_0=0$.  
We consider the Banach spaces  (in fact Hilbert spaces)  $E=L^2(\R\times S^1)$ and  $E_m=H^{m,\delta_m}(\R\times S^1)$, where  the space  $H^{m,\delta_m}(\R\times S^1)$ consists of those  elements in $E$  having
 weak  partial derivatives up to order $m$ which, if  weighted by $e^{\delta_m |s|}$, belong to $E$.  Using Sobolev's compact embedding theorem for bounded domains and the assumption that the sequence $(\delta_m)$ is strictly increasing, one sees that the sequence $(E_m)_{m\geq 0}$ defines a  sc-structure on  $E$. 
 \qed
\end{example}

\begin{definition}  A linear operator $T\colon 
E\rightarrow F$ between sc-Banach spaces is called a {\bf sc-operator},  \index{D- Sc-operator}  if  $T(E_m)\subset F_m$ and the induced operators $T\colon 
E_m\rightarrow F_m$  are continuous for all $m\geq 0$.  A linear {\bf sc-isomorphism}\index{D- Sc-isomorphism}  is a bijective sc-operator whose  inverse is also a sc-operator. 
\qed
\end{definition}

A special class of linear sc-operators  are $\ssc^+$-operators,  defined as follows.
\begin{definition}\label{sc_plus_operators}
 A sc-operator $S\colon 
E\to F$ between sc-Banach spaces  is called a {\bf $\ssc^{\pmb{+}}$-operator}\index{D- Sc$^+$-operator},  if 
 $S(E_m)\subset F_{m+1}$  and  $S\colon 
E\rightarrow F^1$ is a sc-operator.
\qed
\end{definition}
In view of Definition \ref{sc-structure}, the inclusion operator $F_{m+1}\to F_m$ is compact, implying that a $\ssc^+$-operator $S\colon 
E\to F$ is a sc-compact operator  in the sense that  $S\colon 
E_m\to F_m$ is compact for every $m\geq 0$. Hence, given a sc-operator $T\colon 
E\rightarrow F$ and a sc$^+$-operator $S\colon 
E\rightarrow F$, the 
 operator $T+S$ can be viewed, on every level, as a  perturbation of $T$ by the compact operator $S$.

\begin{definition}\label{sc-subspace}
A subspace   $F$ of a sc-Banach space  $E$ is called a {\bf sc-subspace}\index{D- Sc-subspace}, provided $F$ is closed
and the sequence $(F_m)_{m\geq 0}$ given by  $F_m=F\cap E_m$ defines  a sc-structure on $F$.
\qed
\end{definition} 

A sc-subspace $F$ of the sc-Banach space $E$ has a sc-complement provided there exists an algebraic complement $G$ of $F$ so that $G_m=E_m\cap G$
defines an sc-structure on $G$, such  that on every level $m$ we have a topological direct sum 
$$E_m=F_m\oplus G_m.$$
Such a splitting $E=F\oplus G$\index{$E\oplus F$} is called a {\bf sc-splitting}\index{Sc-splitting}. 

The following result from \cite{HWZ2}, Proposition 2.7, 
will be used frequently.
\begin{proposition}\label{prop1}\index{P- Finite-dimensional sc-subspaces}
Let $E$ be a sc-Banach space. A finite-dimensional subspace $F$ of $E$  is a sc-subspace  if and only if $F$ belongs to $E_\infty$.
A finite-dimensional sc-subspace always  has a sc-complement.
\end{proposition}
\begin{proof}
If $F\subset E_\infty$, then $F_m:=F\cap E_m=F$ and $F$ is equipped with the constant sc-structure, so that $F$ is a sc-subspace of $E$. Conversely, 
if $F$ is a sc-subspace, then, by definition,  $F_\infty:=\bigcap_{m\geq 0}F_m=F\cap E_\infty$ is dense in $F$.  Consequently, $F=F_\infty\subset E_\infty$, since $F$ and $F_\infty$ are finite dimensional.  

Next,  let $e_1,\ldots , e_k$ be a basis for a finite-dimensional sc-subspace $F$. In view of the above discussion, $e_i\in E_\infty$. By the  Hahn-Banach theorem, the dual basis  can be extended to linear functionals $\lambda_1,\ldots ,\lambda_k$, which are continuous on $E$ and hence on $E_m$ for every $m$. The map $P\colon 
E\to E$,  defined by $P(x)=\sum_{1\leq i\leq k}\lambda_i(x)e_i$,  has its image in $F\subset E_\infty$. It induces a continuous map from $E_m$ to $E_m$,  and since $P\circ P=P$, it is a sc-projection. Introduce the closed subspace $G:=(\mathbbm{1}-P)(E)$ and let $G_m=G\cap E_m$. Then $E=F\oplus G$ and $G_{m+1}\subset G_m$.  The set $G_\infty:=\bigcap_{m\geq 0}G_m=G\cap E_\infty$ is also dense in $G_m$. Indeed, if $g\in G_m$, we find a sequence $f_n\in E_\infty$, such that $f_n\to g$ in $E_m$. Setting $g_n:=(\mathbbm{1}-P)(f_n)\in G_\infty$ we have $g_n\to (\mathbbm{1}-P)(g)=g$, as claimed.  Consequently, $(G_m)_{m\geq 0}$ defines a sc-structure on $G$ and the proof of the proposition is complete.
 \qed \end{proof}
Next we describe  the quotient construction in the sc-framework. 

\begin{proposition}\label{thm_quotient}\index{P- Sc-quotients}
Assume that $E$ is a sc-Banach space and $A\subset E$ a sc-subspace. Then $E/A$ equipped with the 
filtration $E_m/A_m$ is a sc-Banach space. Note that $E_m/A_m= \{ (x+ A)\cap E_m \,  \vert \, x\in E\}$.
\end{proposition}
\begin{proof}
By definition of a sc-subspace,  the filtration on $A$ is given by $A_m=A\cap E_m$ and $A_m\subset E_m$ is a closed subspace.
Hence 
$$
E_m/A_m = E_m/(A\cap E_m) = \{ (x + A)\cap E_m\, \vert \, x\in E_m\}.
$$
We identify an element $x+A_m$ in $E_m/A_m$ with the  element $x+A$ of  $E/A$, so that algebraically we can view
$E_m/A_m\subset E/A$.
The inclusion $E_{m+1}\rightarrow E_m$ is compact,  implying that the quotient map  $E_{m+1}\rightarrow E_m/A_m$ is compact.  We claim that the inclusion map 
$E_{m+1}/A_{m+1}\rightarrow E_m/A_m$ is compact. In order to show this, we take a sequence $(x_k+A_{m+1})\subset E_{m+1}/A_{m+1}$ satisfying 
$$
\norm{x_k+A_{m+1}}_{m+1}:=\inf \{ \abs{x+a|_{m+1}}\, \vert \, a\in A_{m+1}\}\leq 1, 
$$
and choose a sequence $(a_{k})$ in $A_{m+1}$ such that  $\abs{x_k+a_k}_{m+1}\leq 2$. 
We may assume, after taking a subsequence,  that $x_k+a_k\rightarrow x\in E_{m}$. 
The image of the sequence $(x_k+A_{m+1})$
under the inclusion $E_{m+1}/A_{m+1}\rightarrow E_m/A_m$ is the sequence $(x_k+a_k+A_m)$. Then
\begin{equation*}
\begin{split}
&\norm{(x_k+a_k+A_{m})-(x+A_m)}_m\\
&\quad =\norm{(x_k+a_k-x)+A_m}_m\\
&\quad =\inf\{\abs{x_k+a_k-x)+a}_m\, \vert \,  a\in A_m\}\\
&\quad\leq \abs{x_k+a_k-x}_m\to  0,
\end{split}
\end{equation*}
showing  that  the inclusion  $E_{m+1}/A_{m+1}\rightarrow E_m/A_m$ is compact. Finally let us show that $\bigcap_{j\geq 0} (E_j/A_j)$ is dense in every $E_m/A_m$.
Let us first note that $(E/A)_\infty:=\bigcap_{j\geq 0} (E_j/A_j)$ consists of all elements of the form $x+A_\infty$ with $x\in E_\infty$.
Since $E_\infty$ is dense in $E_m$,  the image under the continuous quotient map $E_m\rightarrow E_m/A_m$ is dense. This completes the proof of 
Proposition \ref{thm_quotient} 
\qed \end{proof}

A distinguished class of sc-operators is the class of sc-Fredholm operators.
\begin{definition}
A sc-operator $T\colon 
E\rightarrow F$  between sc-Banach spaces is called  {\bf sc-Fredholm}\index{D- Sc-Fredholm operator},  provided there exist sc-splittings $E=K\oplus X$ and $F=C\oplus Y$, having the following properties.
\begin{itemize}
\item[(1)]\ $K$ is the kernel of $T$ and is finite-dimensional.
\item[(2)]\  $C$ is finite-dimensional and  $Y$ is the image of $T$.
\item[(3)]\ $T\colon 
X\rightarrow Y$ is a sc-isomorphism.
\end{itemize}
\qed
\end{definition}
In view of Proposition \ref{prop1}  the definition implies that the kernel $K$ consists of smooth points and $T(X_m)=Y_m$ for all $m\geq 0$.  In particular,  we have the topological direct sums
$$
E_m=K\oplus X_m\quad  \text{and}\quad \ F=C\oplus T(E_m).
$$
The index of a sc-Fredholm operator $T$,  denoted by $\ind(T)$, is  as usual defined  by
$$
\ind (T)=\dim(K)-\dim(C).\index{$\ind (T)$}
$$
Sc-Fredholm operators have the following regularizing property.
\begin{proposition}\label{regular}\index{P- Regularizing property}
A sc-Fredholm operator  $T\colon 
E\rightarrow F$ is regularizing, i.e. if $e\in E$ satisfies $T(e)\in F_m$,  then $e\in E_m$.
\end{proposition}
\begin{proof}
By assumption,  $T(e)\in F_m=C\oplus T(X_m)$, so that $T(e)=T(x)+c$ for some $x\in X_m$ and $c\in C$.
Then $T(e-x)=c$ and $e-x\in E_0$. Since $T(E_0)\oplus C=F_0$ it follows that $c=0$ and $e-x\in  K\subset E_\infty$, implying that $x$ and $e$ are on the same level $m$. 
\qed \end{proof}

The following stability result will be crucial later on. The proof, reproduced in Appendix \ref{A1.0},  is taken from \cite{HWZ2},  Proposition 2.1.\begin{proposition}[Compact Perturbation]\label{prop1.21}\index{P- sc$^+$-perturbations}
Let $E$ and $F$ be sc-Banach spaces. If $T\colon 
E\rightarrow F$ is a
sc-Fredholm operator and $S\colon 
E\rightarrow F$ a sc$^+$-operator, then
$T+S$ is also a  sc-Fredholm operator.
\qed
\end{proposition}

The next concept will be crucial later on for the definition of boundaries.

\begin{definition}\label{partial_quadrant}
A {\bf partial quadrant} \index{D- Partial Quadrant}  
in a sc-Banach space $E$  is  a closed convex subset $C$ of $E$, 
such that  there exists a sc-isomorphism $T\colon 
E\rightarrow {\mathbb R}^n\oplus W$ satisfying   $T(C)=[0,\infty)^n\oplus W$.
\qed
\end{definition}

We now consider tuples $(U,C,E)$\index{$(U,C,E)$}, in which  $U$ is a relatively open subset of the partial quadrant $C$ in the sc-Banach space $E$.  
\begin{definition}\label{sc_continuous}
If  $(U, C, E)$ and $(U', C', E')$ are two such  tuples, then a  map $f\colon 
U\rightarrow U'$ 
is called {\bf  $\ssc^0$}  (or of {\bf class $\ssc^0$}, or {\bf sc-continuous}),  provided $f(U_m)\subset V_m$ and the induced maps $f\colon 
U_m\rightarrow V_m$ are  continuous for all $m\geq 0$. 
\qed
\end{definition}

\begin{definition}\label{sc-tangent}\index{D- Tangent tuple} 
The tangent $T(U,C,E)$ of the tuple  $(U,C,E)$ is defined as the tuple 
$$
T(U,C,E)=(TU,TC,TE)\index{$T(U,C,E)$}
 $$where
$$
TU=U^1\oplus E,\quad  TC=C^1\oplus E, \quad  \text{and}\quad  TE=E^1\oplus E.
$$
\qed
\end{definition}
Note that $T(U,C,E)$ is a tuple  consisting again  of a relatively open subset $TU$ in the  partial quadrant $TC$ in  the sc-Banach space $TE$.

Sc-differentiability is a  new notion of differentiability in sc-Banach spaces, which is considerably weaker than
the familiar  notion of Fr\'echet differentiability. The new notion of differentiability is the following.

\begin{definition}\label{scx}\index{D- Sc-differentiability}
We consider  two tuples $(U,C,E)$ and $(U',C',E')$  and a map $f\colon U\to U'$. The map $f$ is called $\pmb{\ssc^{1}}$ (or of {\bf class} $\pmb{\ssc^{1}}$)  provided
the following conditions are satisfied.
\begin{itemize}
\item[(1)]\ The map $f$ is $\ssc^0$.
\item[(2)]\ For every $x\in U_1$ there exists a bounded linear operator $Df(x)\colon E_0\to F_0$ such  that
for $h\in E_1$ satisfying  $x+h\in U_1$,
$$
\lim_{\abs{h}_1\rightarrow 0} \frac{ \abs{f(x+h)-f(x)-Df(x)h}_0}{\abs{h}_1}=0.
$$
\item[(3)]\ The map $Tf\colon TU\to TU'$,  defined by 
$$Tf(x,h)=(f(x),Df(x)h), \quad \text{$x\in U^1$ and $h\in E$}, $$
 is a  $\ssc^0$-map.
The map $Tf\colon TU\to TU'$ is called the {\bf tangent map} of $f$.
\end{itemize}
\qed
\end{definition}

\begin{remark}\index{R- Important remark about (non-)continuity of $Df(x)$} In general,  the map $U_1\to   \mathscr{L}(E_0,F_0)$, defined by $x\to  Df(x)$,  will not(!) be continuous, if the space of bounded linear operators is equipped with the operator norm. However,  if we equip it with the compact open topology it will be continuous. The $\ssc^1$-maps between finite dimensional Banach spaces are the familiar $C^1$-maps.
\qed
\end{remark}

Proceeding inductively,  we define  what it means for  the map $f$  to be  $\ssc^k$ or $\ssc^\infty$. Namely, a $\ssc^0$--map $f$ is said to be a $\ssc^2$--map, if  it is $\ssc^1$ and if  its tangent map $Tf\colon 
TU\to TV$ is  $\ssc^1$. By Definition \ref{scx} and Definition \ref{sc-tangent},  the  tangent map of $Tf$,
$$T^2f\colon 
=T(Tf)\colon 
T^2(U)=T(TU)\to T^2(V)=T(TV),$$
is of class $\ssc^0$.
If the tangent map $T^2f$  is $\ssc^1$, then $f$ is said to be $\ssc^3$,  and so on. The map $f$ is 
{\bf $\ssc^{\pmb{\infty}}$} or {\bf $\ssc$-smooth}, if it is $\ssc^k$ for all $k\geq 0$.
\begin{remark}\index{R- Possible generalizations of domains} \label{rem116}
The above consideration can be generalized. Instead of taking a partial quadrant $C$ one might take a closed convex partial cone $P$ with nonempty interior. This means $P\subset E$ is a closed subset, so that for a real number $\lambda\geq 0$ it holds $\lambda\cdot P\subset P$. Moreover $P$ is convex and has a nonempty interior. If then $U$ is an open subset of $P$ we can define the notion of being sc$^1$ as in the previous definition. The tangent map $Tf$ is then defined on $TP=P^1\oplus E$. We note that $TP$ is a closed cone with nonempty interior. In this context one should be able to deal with 
sc-Fedholm problems, on M-polyfolds with `polytopal' boundaries and even more general situations.
Many of the results in this book can be carried over to this generality as well, but one should check carefully which arguments carry over.
\qed
\end{remark}

\begin{remark}\index{R- Other notions of differentiability: ssc-smoothness}\label{rem117}
Let us mention that there is, of course, also a stronger notion of differentiability which resembles the classical notion.
Namely consider $f:U\rightarrow U'$, where $(U,C,E)$ and $(U',C',D')$  are as in Definition \ref{scx}. We say that $f$ is  ssc$^{1}$
provided id is sc$^0$ and for every $m$ the map $f:U_m\rightarrow U_m'$ is $C^1$, i.e. once continuously differentiable.
Similarly we can define being of class ssc$^k$ and ssc$^\infty$.  The theory associated to this notion of differentiability runs completely parallel 
to the study of classically smooth maps between Banach spaces and there are no surprises. This theory is, of course,  useful in applications as well. To work out the sc-results 
presented later on  in the ssc-framework is left to the reader. Note that {\bf ssc} stands for {\bf strong sc}.
\qed
\end{remark}

\section{Properties of Sc-Differentiability}

In this section we  shall discuss the relationship between the classical smoothness in the Fr\'echet sense and the sc-smoothness. The proofs of the following results can be found  in  \cite{HWZ8.7}.

\begin{proposition}[{\bf Proposition 2.1, \cite{HWZ8.7}}] \label{x1}\index{P- Sc-differentiability}
Let $U$ be a relatively open subset of a partial  quadrant in  a sc-Banach space $E$ and let $F$ be another $\ssc$-Banach space. 
Then a  $\ssc^0$-map $f\colon 
U\to F$ is of class $\ssc^1$  if and only if  the following conditions hold true.
\begin{itemize}\label{sc-100}
\item[(\em{1})]\ For every $m\geq 1$,  the induced map
$f\colon 
U_m\to  F_{m-1}$
is of class $C^1$. In particular,  the derivative  
$$
df\colon 
 U_m\to  \mathscr{L}(E_m,F_{m-1}),\quad x\mapsto df(x)
$$
is a  continuous map.
\item[(\em{2})]\ For every  $m\geq 1$ and every $x\in U_m$,  the bounded  linear operator
$df(x)\colon 
 E_m\to F_{m-1}$ has an extension to a bounded  linear operator $Df(x)\colon 
E_{m-1}\to F_{m-1}$. In addition, the map
\begin{equation*}
U_m\oplus  E_{m-1}\to  F_{m-1}, \quad  (x,h)\mapsto  Df(x)h
\end{equation*}
is continuous.
\end{itemize}
\qed
\end{proposition}

In particular, if $x\in U_\infty$ is  a smooth point in $U$  and $f\colon 
U\to F$ is a $\ssc^1$-map, then the linearization 
$$Df(x)\colon E\to F$$
is a sc-operator.

A  consequence of Proposition \ref{x1} is the following result about lifting the indices.
\begin{proposition}[{\bf Proposition 2.2, \cite{HWZ8.7}}] \label{sc_up}\index{P- Sc-differentiability under lifts}
Let  $U$  and $V$ be  relatively open subsets of partial quadrants in sc-Banach spaces, and let $f\colon 
U\rightarrow V$ be $\ssc^k$. 
Then $f\colon 
U^1\rightarrow V^1$ is also $\ssc^k$.
\qed
\end{proposition}

\begin{proposition}[{\bf Proposition 2.3, \cite{HWZ8.7}}] \label{lower}\index{P- Sc-smoothness versus $C^k$}
Let $U$ and $V$ be relatively open subsets of partial quadrants in sc-Banach spaces.
If $f\colon 
U\to  V$ is $\ssc^k$, then for every $m\geq 0$,  the map $f\colon 
U_{m+k}\to  V_m$ is of class $C^k$.  Moreover,  $f\colon 
U_{m+l}\to  V_m$ is of class $C^l$ for every $0\leq l\leq k$.
\qed
\end{proposition}

The next result is very useful in proving that a given map between sc-Banach spaces is sc-smooth.
\begin{proposition}[{\bf Proposition 2.4, \cite{HWZ8.7}}] \label{ABC-x}
Let $U$ be a relatively open subset of a partial  quadrant in  a sc-Banach space $E$ and let $F$ be another $\ssc$-Banach space. 
Assume that for every $m\geq 0$ and $0\leq l\leq k$,  the map $f\colon 
U\rightarrow V$  induces  a map
$$
f\colon 
U_{m+l}\rightarrow  F_m,
$$
which is of  class $C^{l+1}$. Then $f$ is $\ssc^{k+1}.$
\qed
\end{proposition}

In the case that the target  space $F= \R^N$, Proposition \ref{ABC-x} takes the following form.
\begin{corollary}[{\bf Corollary 2.5, \cite{HWZ8.7}}] \label{ABC-y}
Let $U$ be a  relatively open subset of a partial quadrant in a sc-Banach space and  $f\colon 
U\to  \R^N$. If 
for some $k$ and all $0\leq l \leq k$ the map $f\colon 
U_l \to  \R^N$ belongs to  $C^{l +1}$, then $f$ is $\ssc^{k+1}$.
\qed
\end{corollary}

\section{The Chain Rule and Boundary Recognition}\label{subsection_boundary_recognition}
The cornerstone of the sc-calculus, the chain rule, holds true.
\begin{theorem}[{\bf Chain Rule}]\label{sccomp}\index{T- Chain rule}
Let $U\subset C\subset E$ and $V\subset D\subset F$ and $W\subset Q\subset G$ be relatively open subsets
of partial quadrants in sc-Banach spaces and let $f\colon 
U\rightarrow V$ and $g\colon 
V\rightarrow W$ 
be $\ssc^1$-maps.
Then the composition $g\circ f\colon 
U\rightarrow W$ is also $\ssc^1$ and
$$
T(g\circ f) =(Tg)\circ (Tf).
$$
\qed
\end{theorem}
The result is  proved in \cite{HWZ2} as  Theorem 2.16 for open sets $U$, $V$ and $W$.
We give the adaptation to the somewhat more general setting here. The proof, which is given in Appendix \ref{A1.1},  is very close 
to the one in \cite{HWZ2} and does not require any new ideas.

\begin{remark}\index{R- Remark on chain rule}
The result is somewhat surprising, since differentiability is only guaranteed
under the loss of one level $U_1\rightarrow F_0$, so that one might expect
 for a composition a loss of two levels. However, one is saved
by the compactness of the embeddings $E_{m+1}\rightarrow E_m$.
\qed
\end{remark}

The next basic result concerns  the boundary recognition. Let $C$ be a partial quadrant in the  sc-Banach space  $E$. 
We choose a linear $\ssc$-isomorphism $T\colon 
E\to  \R^k\oplus W$ satisfying $T(C)=[0,\infty)^k\oplus W$. 
If $x\in C$,  then 
$$T(x)=(a_1, \ldots, a_k, w)\in [0,\infty )^k\oplus W,$$
and we define the integer $d_C(x)$ by 
$$
d_C(x)\colon 
=\#\{i\in \{1,\ldots, k\}\vert \, a_i=0\}.$$
\begin{definition}\label{degeneracy_index_1}\index{D- Degeneracy index}
The map $d_C\colon 
C\rightarrow {\mathbb N}$ is called the {\bf degeneracy index}.
\qed
\end{definition}

Points $x\in C$ satisfying $d_C(x)=0$\index{$d_C$} are interior points of $C$, the  points satisfying  $d_C(x)=1$ are honest boundary points, and the  points with $d_C(x)\geq 2$ are corner points. The size of the index gives the complexity of the corner.

It is not difficult to see that this definition is independent of  the choice of a sc-linear isomorphism $T$.
\begin{lemma}\label{corner_recognition_linear}
The  map  $d_C$  does not depend on the choice of a linear sc-isomorphism $T$.
\qed
\end{lemma}
The proof is given in Appendix \ref{A1.2}.

\begin{theorem}\label{ppp22}\index{T- Corner recognition}
Given the tuples $(U, C, E)$ and $(U', C', E')$ and a germ  $f\colon 
(U,x)\rightarrow (U',x')$ of a local sc-diffeomorphism satisfying $x'=f(x)$, then
$$d_C(x)=d_{C'}(f(x)).$$
\qed
\end{theorem}

The proof of a more general result is given in \cite{HWZ2}, Theorem 1.19. 
We  note that sc-smooth diffeomorphisms recognize  boundary points and corners, whereas
homeomorphisms recognize only  boundaries, but no corners.

\section{Appendix}
\subsection{Proof of the sc-Fredholm Stability Result}\label{A1.0}
\begin{proof}[Proposition \ref{prop1.21}]
Since $S\colon 
E_m\rightarrow F_m$ is compact for every level,  the operator 
$T+S\colon 
E_m\rightarrow F_m$ is Fredholm for every $m$. Denoting  by $K_m$ the
kernel of this operator,  we claim that $K_{m}=K_{m+1}$ for every $m\geq 0$. Clearly, $K_{m+1}\subseteq K_{m}$. If  $x\in K_m$,  then  $Tx=-Sx\in F_{m+1}$ and hence, by Proposition \ref{regular},
$x\in E_{m+1}$. Thus $x\in K_{m+1}$, implying  $K_m\subseteq K_{m+1}$, so that  $K_m=K_{m+1}$. Set $K=K_0$.  By
Proposition \ref{prop1}, $K$ splits the sc-space $E$,  since it is a
finite dimensional subset of $E_{\infty}$. Therefore, we have the
sc-splitting $E=K\oplus X$  for a suitable sc-subspace  $X$ of $E$.

Next define $Y=(T+S)(E)=(T+S)(X)$ and $Y_m\colon 
=Y\cap F_m$ for $m\geq 0$.  Since $T+S\colon 
E\to F$ is Fredholm,  its  image $Y$ is closed.  We claim, that the set  $Y_\infty$, defined by $Y_\infty:=\bigcap_{m\geq 0}Y_m=Y\cap F_\infty$,  is dense in  $Y_m$ for every $m\geq 0$. Indeed, if $y\in Y_m$, then $y=(T+S)(e)$ for some $e\in E$. Since  by Proposition \ref{regular}, $T$ is regularizing, we conclude that $e\in E_m$. Since $E_\infty$ is dense in $E_m$,  there exists a sequence $(e_n)\subset E_\infty$ such that $e_n\to e$ in $E_m$. Then, setting $y_n:=(T+S)(e_n)\in Y_\infty$, we conclude that $y_n\to (T+S)(e)=y$, which shows that $Y_\infty$ is dense in $Y_m$.  Finally, we consider the projection map $p\colon 
F\to  F/Y$. Since $p$  is a continuous surjection onto a finite-dimensional space and $F_\infty$ is dense in $F$, it follows that $\ov{p(F_\infty)}=F/Y$.  Hence, we can choose  a basis in $F/Y$ whose representatives belong to $F_\infty$.  Denoting by  $C$ the span of these representatives, we obtain  $F_m=Y_m\oplus C$.   Consequently, we have a sc-splitting $F=Y\oplus C$ and the proof of Proposition \ref{prop1.21} is complete.
\qed \end{proof}

\subsection{Proof of the Chain Rule}\label{A1.1}
\begin{proof}[Theorem \ref{sccomp}]
The maps 
 $f\colon 
U_1\to F$ and $g\colon 
V_1\to G$  are of class $C^1$ in view of Proposition \ref{lower}. Moreover, $Dg(f(x))\circ Df(x)\in \mathscr{L}(E, G)$, if
$x\in U_1$. Fix $x\in U_1$ and 
 $h\in E_1$ sufficiently small satisfying $x+h\in U_1$ and  $f(x+h)\in V_1$.
Then, using  the postulated
properties of $f$ and $g$,
\begin{equation*}
\begin{split}
&g(f(x+h))-g(f(x))-Dg(f(x))\circ Df(x)h\\
&=\int_{0}^{1} Dg(tf(x+h)+(1-t)f(x))\ [f(x+h)-f(x)-Df(x)h] dt\\
&\phantom{=}+\int_{0}^{1}\bigl( [Dg(tf(x+h)+(1-t)f(x))-Dg(f(x))] \circ Df(x)h \bigr)dt.
\end{split}
\end{equation*}
The first integral after dividing by $\abs{h}_1$ takes the form
\begin{equation}\label{integral1}
\int_{0}^{1} Dg(tf(x+h)+(1-t)f(x)) \frac{f(x+h)-f(x)-Df(x)h}{\abs{h}_1}\  dt.
\end{equation}
For  $h\in E_1$ such that  $x+h\in U_1$ and $\abs{h}_1$ is small, we have that $tf(x+h)+(1-t)f(x)\in V_1$. Moreover, the maps
$
[0,1]\rightarrow F_1$ defined by $t\mapsto tf(x+h)+(1-t)f(x)
$
are continuous and converge in $C^0([0,1],F_1)$
to the constant map
$t\mapsto  f(x)$ as  $\abs{h}_1\rightarrow 0$.
Since $f$ is  of class $\ssc^1$, the quotient
 $$
a(h)\colon 
=\frac{ f(x+h)-f(x)-Df(x)h}{\abs{h}_1}$$
converges  to $0$ in $F_0$  as $\abs{h}_1\to 0$. Therefore, by the continuity assumption (3) in Definition \ref{scx},
 the map $$
(t,h)\mapsto  Dg(tf(x+h)+(1-t)f(x)) [a(h)],
$$
as a map from $[0,1]\times E_1$ into $G_{0}$, converges to $0$ as
$\abs{h}_1\to 0$  uniformly in $t$.
Consequently,   the expression in \eqref{integral1}  converges to $0$ in $G_0$ as $\abs{h}_1\to 0$, if  $x+h\in U_1$.
Next we 
consider the integral
\begin{equation}\label{integral2}
\int_{0}^{1} \bigl[  Dg(tf(x+h)+(1-t)f(x))-Dg(f(x))\bigr] \circ Df(x)
\frac{h}{\abs{h}_{1}}\ dt.
\end{equation}
By Definition \ref{sc-structure},  the inclusion operator $E_1\to E_0$ is compact,  so that 
the set of all $\frac{h}{\abs{h}_{1}}\in E_1$ has  a compact closure in $E_{0}$.  Therefore, since $Df(x)\in \mathscr{L}(E_{0}, F_{0})$ is a continuous map  by Definition \ref{scx},  the closure of
the set of all
$$
Df(x) \frac{h}{ \abs{h}_{1}}
$$
is compact  in $F_{0}$. Consequently, again by  Definition \ref{scx}, every sequence $(h_n)\subset E_1$, satisfying $x+h_n\in U_1$ and $\abs{h_n}_1\to 0$, possesses a subsequence having the property that the integrand of the integral in \eqref{integral2} converges to $0$ in $G_{0}$ uniformly in $t$. Hence, the integral \eqref{integral2} also  converges to $0$ in $G_0$, as $\abs{h}_1\to 0$ and  $x+h\in U_1$.  We
have proved  that
$$
\frac{\abs{g(f(x+h))-g(f(x))-Dg(f(x))\circ Df(x)h}_{0}}{\abs{h}_{1}}\rightarrow 0
$$
as $\abs{h}_1\to 0$  and  $x+h\in U_1$.   Consequently,  condition (2) of Definition \ref{scx} is satisfied for the composition $g\circ f$ with the linear operator
$$D(g\circ f)(x)=Dg(f(x))\circ Df(x)\in  \mathscr{L}(E_0, G_0),
$$
where $x\in U_1$. We conclude that the tangent map
$T(g\circ f)\colon 
TU\to TG$,
$$(x, h)\mapsto (\ g\circ f(x), D(g\circ f)(x)h\ )$$
is $\ssc$-continuous and, moreover,
$T(g\circ f)=Tg\circ Tf$. The proof of Theorem  \ref{sccomp} is complete.
\qed \end{proof}

\subsection{Proof Lemma \ref{corner_recognition_linear} }\label{A1.2}
\begin{proof}[Lemma \ref{corner_recognition_linear}]
We take a second sc-isomorphism  $T'\colon 
E\to \R^{k'}\oplus W'$,   satisfying $T'(C)=[0,\infty)^{k'}\oplus W'$ and first show  that $k=k'$. 
We consider the composition $S=T'\circ T^{-1}\colon 
\R^k\oplus W\to \R^{k'}\oplus W'$, where  $S$ is a sc-isomorphism, mapping 
$[0,\infty )^k\oplus W\to [0,\infty )^{k'}\oplus W'$. We claim that  $S(\{0\}^k\oplus W)=\{0\}^{k'}\oplus W'$.  Indeed,  suppose $S(0, w)=(a, w')$ for some $(a, w')\in [0,\infty)^{k'}\oplus W'$.
Then,  we conclude 
$$S(0, tw)= tS(0, w)=t(a, w')=(ta, tw')\in  [0,\infty)^{k'}\oplus W'$$
for all $t\in \R$.  Consequently, $a=0$. Hence, $S(\{0\}^k\oplus W)\subset \{0\}^{k'}\oplus W'$ and since $S$ is an isomorphim, we obtain equality of these sets.  
The set  $\{0\}^k\oplus W\equiv W$ has codimension $k$ in $\R^k\oplus W$, and since $S$ is an isomorphism  we have 
$$S(\R^k\oplus W)=S(\R^k\oplus \{0\})\oplus S(\{0\}^k\oplus W)=S(\R^k\oplus \{0\})\oplus ( \{0\}^{k'}\oplus W'),$$
so that the codimension of $ \{0\}^{k'}\oplus W'$ in $\R^{k'}\oplus W'$ is equal to $k$. On the other hand, this codimension is equal to $k'$. Hence $k'=k$, as claimed.
To prove our result, it now suffices  to show that if $x=(a, w)\in [0,\infty )^k\oplus W$ and $y=S(x)=(a', w')$, then 
\begin{equation}\label{independence_0}
\# I(x)=\# I' (y)
\end{equation}
where $I(x)$ is the set of indices $1\leq i\leq k$ for which $x_i=0$. The set $I'(y)$ is defined similarly. With $I(x)$ we associate the subspace $E_x$ of $\R^k\oplus W$, 
\begin{align*}
E_x&=\{(b, v)\in R^k\oplus W\vert \,\text{$b_i=0$ for all $i\in I(x)$}\}.
\end{align*}
The set $E'_{y}$ is defined analogously. Observe that $\# I(x)$ is equal to the codimension of $E_x$ in $\R^k\oplus W$.  In order to prove \eqref{independence_0}, it is enough to show that 
$S(E_x)\subset E'_y$. If this is the case, then,  since $S$ is an isomorphism, we obtain $S(E_x)=E'_y$,  which shows that the codimensions of $E_x$ and $E'_y$ in $\R^k\oplus W$ and $\R^{k'}\oplus W'$, respectively, are the same. Now, to see that $S(E_x)\subset E'_y$, we take $u=(b, v)\in E_x$ and note that $x+tu\in [0,\infty)\oplus W$ for $\abs{t}$ small. Then, 
$S(x+tu)=S(x)+tS(u)=y+tS(u).$ If $i\in I'(y)$, then $(S(x+tu))_i=y_i+t(S(u))_i=t(S(u))_i\in [0,\infty)^k\oplus W'$, which is only possible if  $(S(u))_i=0$.  Hence, $S(E_x)\subset E'_y$ and the proof is complete.
\qed \end{proof}

\chapter{Retracts}
The crucial  concept of this book is that of a sc-smooth retraction, which we  are going to introduce next.
\section{Retractions and Retracts}\label{section2.1}
In this subsection we introduce several basic  concepts.
\begin{definition} \label{sc-smooth_retraction}\index{D- Sc-smooth retraction}
We consider a tuple $(U, C, E)$ in which $U$ is a relatively open subset of the partial quadrant $C$ in the sc-Banach space $E$. A sc-smooth map $r\colon U\to U$ is called 
a {\bf sc-smooth retraction}  on $U$, if 
$$
r\circ r=r.
$$
\qed
\end{definition}

As a side remark, we observe that  if $U$ is an open subset of a Banach space $E$ and $r\colon U\rightarrow U$ is a $C^\infty$-map (in the classical sense) satisfying $r\circ r=r$, then the image of $r$ is a smooth submanifold of the Banach space $E$, as the following proposition shows.

\begin{proposition}\label{smooth_retract}\index{Cartan's last theorem}
Let $r\colon 
U\to U$ be a $C^\infty$-retraction defined on an open subset $U$ of a Banach space $E$. Then $O:=r(U)$ is a $C^\infty$-submanifold of $E$. More precisely, for  every point $x\in O$,  there exist an open neighborhood  $V$ of $x$ and  an open neighborhood $W$  of $0$ in $E$, a splitting $E =R \oplus N$, and a smooth diffeomorphism 
$\psi \colon 
 V\to  W$  satisfying  $\psi (0)=x$ and 
$$\psi (O\cap V ) = R \cap W.$$
\qed
\end{proposition}
An elegant proof is due to Henri Cartan in \cite{H.Cartan}\footnote{The first author would like to thank E. Ghys for pointing out this reference during his visit at IAS in 2012. It seems that this is H. Cartan's last mathematical paper.}.

In sharp contrast to the conclusion of this proposition from classical differential geometry,  there are $\ssc^\infty$-retractions which have,  for example, locally varying finite dimensions, as we shall see later on. 
\begin{definition}\label{sc_smooth_retract}\index{D- Smooth retract}
A tuple $(O,C,E)$ is called a {\bf sc-smooth retract}, if there exists a relatively open subset $U$ of the partial quadrant $C$ and a  sc-smooth retraction $r\colon 
U\to U$ satisfying 
$$r(U)=O.$$
\qed
\end{definition}

In case $x\in U^1$ we deduce  from $r\circ r=r$ by the chain rule, 
$ Dr (r(x))\circ Dr(x)=Dr(x)$. Hence if $x\in O^1$, then $r(x)=x$ so that  $Dr (x)\circ Dr(x)=Dr(x)$ and $Dr(x)\colon E\to E$ is a projection.

\begin{proposition}\label{independent_of_retractions}\index{P- Intrinsic geometry of  retracts}
Let  $(O, C, E)$ be  a sc-smooth  retract and assume  that $r\colon 
U\to U$ and $r'\colon 
U'\to U'$ are two sc-smooth retractions defined on relatively open subsets $U$ and $U'$ of $C$  and satisfying $r(U)=r'(U')=O$. Then
$$Tr (TU)=Tr'(TU').$$
\end{proposition}
\begin{proof}
If $y\in U$, then there exists $y'\in U'$, so that $r(y)=r'(y')$.  Consequently, 
$r'\circ r(y)=r'\circ r'(y')=r'(y')=r(y)$, and hence $r'\circ r=r$.  Similarly, one sees that $r\circ r'=r'$.
If  $(x, h)\in Tr (TU)$,  then $(x, h)=Tr(y, k)$ for a pair  $(y, k)\in TU$.
Moreover,  $x\in O_1\subset U'_1$, so that $(x, h)\in TU'$.  From $r'\circ r=r$ it follows,  using the chain rule,  that
$$Tr'(x, h)=Tr' \circ Tr(y,k)=T(r'\circ r)(y, k)=Tr(y, k)=(x, h),$$
implying $Tr (TU)\subset Tr'(TU')$.  Similarly,  one shows that $Tr' (TU')\subset Tr(TU)$, and the proof of the proposition is complete.
\qed \end{proof}

Proposition \ref{independent_of_retractions} allows us  to  define the tangent of a sc-smooth retract  $(O,C,E)$\index{$(O,C,E)$} as follows.
\begin{definition}\index{D- Tangent of a retract}
The {\bf tangent of the sc-smooth retract}  $(O,C,E)$  is defined as the triple
$$
T(O,C,E)=(TO,TC,TE),\index{$T(O,C,E)$}
$$
in which $TC=C^1\oplus E$ is the tangent of the partial quadrant $C$, $TE=E^1\oplus E$,  and  $TO:= 
Tr(TU)$, where  $r\colon U\to U$ is any sc-smooth retraction onto $O$.
\qed
\end{definition}
We recall that, explicitly, 
$$
TO=Tr(TU)=Tr(U^1\oplus E)=\{(r(x), Dr(x)h)\, \vert \, x\in U^1, h\in E\}.
$$
Starting with a relatively open subset $U$ of a partial quadrant $C$ in $E$,  the tangent $TC$ of $C$  is a partial quadrant in $TE$, and $TU$ is a relatively open subset of $TC$.

In the following we shall quite often just write $O$ instead of $(O,C,E)$
and  $TO$ instead of $(TO,TC,TE)$. However, we would like to point out that the ``reference'' $(C,E)$ is important, since 
it is possible that $O$ is a sc-smooth retract with respect to some nontrivial $C$, but would not be a sc-smooth retract for $C=E$.

\begin{proposition}\label{preparation_retraction}\index{P- Tangent map}
Let $(O,C,E)$ and $(O',C',E')$  be sc-smooth retracts and let $f\colon 
O\rightarrow O'$  be a map. If  $r\colon 
U\to U$ and $s\colon 
V\to V$  are sc-smooth retractions onto $O$ for the triple $(O, C, E)$, then the following holds.
\begin{itemize}
\item[\em(1)]\  If $f\circ r\colon 
U\rightarrow E'$ is sc-smooth, then $f\circ s\colon 
V\rightarrow E'$ is sc-smooth, and vice versa.
\item[\em(2)]\ If  $f\circ r$ is sc-smooth, then  
$T(f\circ r)|TO=T(f\circ s)|TO$.
\item[\em(3)]\ The  tangent map $T(f\circ r)\vert TO$ maps $TO$ into $TO'$.
\end{itemize}
\end{proposition}
\begin{proof}
We assume that $f\circ r\colon 
U\to  E'$ is $\ssc^\infty$. Since  $s\colon 
V\to  U\cap V$ is $\ssc^\infty$, the chain rule implies that the composition $f\circ r\circ s\colon 
V\to  F$ is $\ssc^\infty$.
Using  the identity  $f\circ r\circ s=f\circ s$, we conclude that $f\circ s$ is $\ssc^\infty$. Interchanging the role of $r$ and $s$,  the first part of the lemma is proved.
If $(x, h)\in TO$, then $(x, h)=Ts (x, h)$ and using the identity $f\circ r\circ s=f\circ s$ and the chain rule, we conclude
\begin{equation*}
T(f\circ r)(x,h)=T(f\circ r)(Ts)(x,h)=T(f\circ r\circ s)(x,h)=T(f\circ s)(x,h)
\end{equation*}
Now we take any  sc--smooth retraction $\rho\colon 
U'\to  U'$ defined on a  relatively open subset $U$ of the  partial quadrant $C'$ in $E'$  satisfying $\rho(U')=O'$. Then $\rho\circ f= f$ so that $\rho\circ f\circ r=f\circ r$. Application of  the chain rule  yields the identity
$$
T(f\circ r)(x,h)=T(\rho\circ f\circ r)(x,h)=T\rho\circ T(f\circ r)(x,h)
$$
for all $(x, h)\in Tr(TU).$
Consequently, $T(f\circ r)\vert Tr(TU)\colon 
TO\to TO'$ and this map is  independent of the choice of a sc-smooth retraction onto $O$.
\qed \end{proof}

In view of  Proposition  \ref{preparation_retraction},  we can define the sc-smoothness of a map between sc-smooth retracts, as well as its tangent as follows.
\begin{definition}\label{tangent_retract}\index{D- Tangent map}
Let  $f\colon 
O\rightarrow O'$  be a map between sc-smooth retracts $(O,C,E)$ and $(O',C',E')$, and let  $r\colon 
U\to U$ be a sc-smooth retraction for   $(O, C, E)$. Then {\bf the  map $f$ 
is sc-smooth},  if  the  composition
$$
U\rightarrow E', \quad x\mapsto f\circ r(x)
$$
is sc-smooth.  In this case the {\bf tangent map} $Tf\colon 
TO\rightarrow TO'$ is defined by
$$
Tf\colon 
=T(f\circ r)\vert TO.\index{$Tf$}
$$
\qed
\end{definition}
Proposition  \ref{preparation_retraction} shows that   the definition does not depend on the choice of  the sc-smooth retraction $r$, as long as it retracts onto  $O$. With Definition \ref{tangent_retract}, the chain rule for sc-smooth maps between sc-smooth retracts  follows from Theorem \ref{sccomp}.
\begin{theorem}
Let $(O,C,E)$, $(O',C',E')$ and $(O'',C',E'')$ be two sc-smooth retracts and let $f\colon 
O\rightarrow O'$ and $g\colon 
O'\rightarrow O''$ be sc-smooth. Then $g\circ f\colon 
O\rightarrow O''$ is sc-smooth and $T(g\circ f)=(Tg)\circ (Tf)$.
\qed
\end{theorem}

If $f\colon O\to O'$ is a sc-smooth map between sc-smooth retracts as in Definition \ref{tangent_retract}, we abbreviate the linearization of $f$ at the point $o\in O^1=r(U^1)$ (on level $1$) by 
$$
Tf(o)=T(f\circ r)(o)\vert T_oO,
$$
so that 
$$Tf(o)\colon T_oO\to T_{f(o)}O'$$
is a continuous linear  map between the tangent spaces 
$T_oO=Dr(o)E$ and $T_{f(o)}O'=Dr'(f(o))E'$.

At this point, we have defined a new {\bf category ${\mathcal R}$}\index{${\mathcal R}$} whose  objects  are  sc-smooth
retracts $(O,C,E)$ and whose morphisms are sc-smooth maps 
$$
f\colon 
(O,C,E)\rightarrow (O',C',E')
$$
 between sc-smooth retracts.
The map $f$ is only defined between $O$ and $O'$, but the other data $(O,C,E)$ are  needed to define the differential geometric properties of $O$.

We also have a  well-defined functor, namely the {\bf tangent functor} \index{Tangent functor}
$$
T\colon 
{\mathcal R}\rightarrow {\mathcal R}.\index{$T\colon 
{\mathcal R}\rightarrow {\mathcal R}$}
$$
It maps  the object $(O,C,E)$ into its tangent $(TO,TC,TE)$ and the morphism
$f\colon 
(O,C,E)\rightarrow (O',C',E')$ into its tangent map $Tf$,
\begin{gather*}
Tf\colon 
(TO,TC,TE)\rightarrow (TO',TC',TE')\\
Tf(x,h)=(f(x),Df(x)h)\quad  \text{for all $(x,h)\in TO$.}
\end{gather*}
 The chain rule guarantees the functorial property 
 $$T(g\circ f)=(Tg)\circ (Tf).$$
 Clearly, this is enough to build a differential geometry  whose  local models
are sc-smooth retracts. The details will be carried out  in the next subsection. Let us note that $T(O,C,E)$ has evidently
more structure than $(O,C,E)$, for example $TO\rightarrow O^1$ seems to have some kind of bundle structure.
We shall discuss this briefly and 
 introduce another category of retracts. 
 
Consider a tuple $(U, C, E)$, where $U$ is a relatively open subset of the partial quadrant $C$ in the sc-Banach space $E$, and let $F$ be another sc-Banach space.
\begin{definition}\index{D- Bundle retraction}\index{D- Bundle model}
Let $p\colon 
U\oplus F\rightarrow U$ be a trivial sc-bundle  defined by the projection $p\colon 
U\oplus F\rightarrow U$ onto $U$, and let $R\colon 
U\oplus F\rightarrow U\oplus F$ be  a sc-smooth map of the form $R(u,v)=(r(u),\rho(u,v))$, satisfying $R\circ R=R$, and where $\rho(u,v)$ is linear in $v$. We call $R$ a {\bf sc-smooth  bundle retraction}  (covering the sc-smooth retraction $r$)
and $B=R(U\oplus F)$ the associated sc-smooth bundle retract, and denote by $p\colon 
B\rightarrow O$ the induced projection onto $O=r(U)$.   We shall refer to $p:B\rightarrow O$ as a {\bf sc-smooth bundle model} or {\bf sc-smooth bundle retract} and sometimes write $(B,C\oplus F,E\oplus F)$.
\qed
\end{definition}
Given two sc-smooth bundle retracts $p\colon 
B\rightarrow O$ and $p'\colon 
B'\rightarrow O'$,  a sc-smooth map $\Phi\colon 
B\rightarrow B'$
of the form $\Phi(u,v)=(a(u),\phi(u,v))$, where $\phi(u,v)$ is linear in $v$ and $p'\circ \Phi=p$ is called a {\bf sc-smooth (local) bundle map}.\index{D- Bundle map}

By $\mathcal{BR}$\index{$\mathcal{BR}$} we denote the category whose objects are sc-smooth bundle retracts and whose  morphisms are the sc-smooth bundle maps.
There is  a natural forgetful functor $\mathcal{BR}\rightarrow \mathcal{R}$, which on objects associates with the
bundle $B\xrightarrow{p}O$, the total space  $B$ and which views a bundle map $\Phi$ just as a sc-smooth map.

If $r\colon 
U\rightarrow U$ is a sc-smooth  retraction, then its  tangent map $Tr\colon 
TU\rightarrow TU$ is a sc-smooth bundle retraction. Moreover, the tangent map is a sc-smooth bundle map. Consequently,  the tangent functor can be  viewed as the functor 
 $$
 T\colon 
\mathcal{R }\rightarrow \mathcal{BR}.\index{$ T\colon 
\mathcal{R }\rightarrow \mathcal{BR}$}
 $$
The functor $T$  associates with  the object $(O,C,E)$\index{$(O,C,E)$} the triple $(TO,TC,TE)$\index{$T(O,C,E)$}, in which we  view $TE$  as the bundle $TE\rightarrow E^1$,
$TC$ as the bundle $TC\rightarrow C^1$,  and $TO$ as the bundle $TO\rightarrow O^1$. With a  sc-smooth map $f\colon 
(O,C,E)\rightarrow (O',C',E')$, the functor $T$  associates the 
sc-smooth bundle map $Tf$.

There is another class of retractions, called strong bundle retractions, which will be introduced in the later parts of Section \ref{SEC2.2}. 

\begin{remark}\index{R- Possible generalizations of domains}
Many of the aspects of the theory developed later on can be generalized as follows.
We take as local models $(O,C,E)$,
where $C$ is a closed convex set in the sc-Banach space $E$ with nonempty interior,  and $O\subset C$, so that there exists 
a relatively open subset $U\subset C$ with $O\subset U\subset C$ and an sc-smooth map $r:U\rightarrow U$ with $r\circ r=r$ and $O=r(U)$.
The fact that the interior of $C$ is open allows the necessary generalization of  the notion of sc-differentiability.
A particular case, which is of interest, is concerned 
with a situation described up to sc-isomorphism as follows. We assume that $E={\mathbb R}^n\oplus W$ and 
$C= D\oplus W$, where $D$ is the closed convex hull of a finite number of points in ${\mathbb R}^n$, so that $\text{int}(D)\neq \emptyset$.
We leave such generalizations to the reader.\qed
\end{remark}

\section{M-Polyfolds and Sub-M-Polyfolds}\label{SEC2.2}
We start with the following observation about sc-retractions and sc-retracts.

\begin{proposition}\label{new_retract_1}\index{P- Restrictions of retracts}
Let $(O, C, E)$ be a sc-smooth retract.
\begin{itemize}
\item[\em(1)]\ If $O'$ is an open subset of $O$, then $(O',C,E)$ is a sc-smooth retract.
\item[\em(2)]\ Let  $V$ be  an open subset of $O$ and $s\colon 
V\rightarrow V$  a sc-smooth map, satisfying  $s\circ s=s$.  If $O'=s(V)$, then $(O',C,E)$ is a sc-smooth retract.
\end{itemize}
\end{proposition}
\begin{proof}
(1)\, 
By assumption, there exists a sc-smooth retraction $r\colon 
U\to U$ defined on a relatively open subset $U$ of  the partial quadrant $C\subset E$, whose image is $O=r(U)$.  Since the map $r\colon 
U\rightarrow U$ is continuous,  the set $U':=r^{-1}(O')$ is an open subset of $U$ and, therefore, a relatively open subset of $C$.  Clearly,  $O'\subset U'$ 
and the restriction  $r'=r\vert U'$ defines a sc-smooth retraction $r'\colon 
U'\rightarrow U'$ onto  $r'(U')=O'$.
Consequently, the triple  $(O',C,E)$ is a sc-smooth retract, as claimed.\par

\noindent (2)\, The triple $(V, C, E)$ is, in view of (1),  a sc-smooth retract. Hence, there exists a sc-smooth retraction $r\colon U\to U$ onto $r(U)=V$, where $U\subset C\subset E$ 
We define the map $\rho\colon U\to U$ by 
$$\rho=s\circ r.$$
\par

Then $\rho$ is sc-smooth and $\rho\circ \rho=(s\circ r)\circ (s\circ r).$ If $x\in U$, then $r(x)\in V$, hence $s(r(x))\in V$ and consequently, $r(s(r(x))=s(r(x))$, so that 
$$r\circ \rho=\rho.$$
From $s\circ s=s$, we conclude that  $\rho\circ \rho=\rho$. Hence $\rho$ is a sc-smooth retraction onto the subset 
$$\rho(U)=s\circ r(U)=s(V)=O'.$$
We see that $(O', C, E)$ is a sc-smooth retract, as claimed.
\qed \end{proof}

In the following every sc-smooth map $r\colon 
O\rightarrow O$, defined on a sc-smooth retract  $O$ and satisfying $r\circ r=r$ will also be called a {\bf sc-smooth retraction}.
\begin{definition}\label{sc-charts}\index{D- Sc-charts}
Let $X$ be a topological space and $x\in X$.
{\bf A chart around $x$}  is a tuple $(V,\phi,(O,C,E))$,\index{$(V,\phi,(O,C,E))$} in which  $V\subset X$ is an open neighborhood of  $x$ in $X$, 
and $\phi\colon 
V\rightarrow O$ is a homeomorphism. Moreover, $(O,C,E)$ is a sc-smooth retract.
Two charts $(V,\phi,(O,C,E))$ and $(V',\psi,(O',C',E'))$ are called {\bf sc-smoothly compatible},  if  the transition maps 
$$\psi\circ\phi^{-1}\colon 
\phi(V\cap V')\rightarrow \psi(V\cap V')\quad \text{and}\quad \phi\circ \psi^{-1}\colon 
\psi(V\cap V')\rightarrow \phi(V\cap V')$$
are sc-smooth maps (in the sense of Definition \ref{tangent_retract}).
\qed
\end{definition}
Let us observe that $\phi(V\cap V')$ is an open subset of $O$ and so, in view of part (1) of Proposition  \ref{new_retract_1}, the tuple $(\phi(V\cap V'), C, E)$ is a sc-smooth retract.  So,  the above  transition maps are defined on  sc-smooth retracts.
\begin{definition}\label{sc_atlas}\index{D- Sc-smooth atlas}
{\bf A sc-smooth atlas} on the topological space $X$ consists
of a set of charts 
$$
(V,\varphi,(O,C,E)),
$$ 
such  that any two of them are sc-smoothly compatible and the open sets $V$ cover $X$.
Two {\bf sc-smooth atlases} on $X$  are said to be {\bf equivalent}, if their  union is again a  sc-smooth atlas.
\qed
\end{definition}

\begin{definition} \index{D- M-polyfold}
{\bf A M-polyfold} $X$  is a Hausdorff paracompact  topological
space equipped with an equivalence class of sc-smooth atlases.
\qed
\end{definition}

Analogous to the smoothness of maps between manifolds we shall define sc-smoothness of maps between M-polyfolds.

\begin{definition}
A  map $f\colon 
X\to Y$ between two M-polyfolds is called {\bf sc-smooth} \index{D- Sc-smooth map}
if its local coordinate representations are sc-smooth. In detail, this requires the following. If $f(x)=y$ and $(V,\varphi, (O, C, E))$ is a chart around $x$ belonging to the atlas of $X$  and 
$(V',\varphi', (O', C', E'))$ is a chart around $y$ belonging to the atlas of $Y$, so that $f(V)\subset V'$, then the map 
$$\varphi'\circ f\circ \varphi^{-1}\colon 
O\to O'$$
is a sc-smooth map between the sc-smooth retracts in the sense of Definition \ref{tangent_retract}.
\qed
\end{definition}
The definition does not depend on the choice of sc-smoothly compatible charts.

We recall that a Hausdorff topological space is paracompact, provided  every open cover of $X$ has an open locally finite refinement.  
It is a well-known fact,  that, given an open cover $\mathscr{U}={(U_i)}_{i\in I}$,  one can find a refinement $\mathscr{V}=(V_i)_{i\in I}$ satisfying  $V_i\subset U_i$  for all $i\in I$ (some  of the sets $V_i$ might be  empty).

Given  a  M-polyfold $X$, we say that a point $x\in X$  is on the {\bf level $m$}\index{Level $m$ points}, if  there exists a chart $(V, \varphi, (O, C, E))$ around $x$, such that $\varphi(x)\in O_m$. Of course, this definition is independent of  the choice of a chart around $x$. We denote the collection of all points on the level $m$ by $X_m$.   The topology on $X_m$ is defined as follows.  We abbreviate by  $\mathscr{B}$ the collection of  all sets  $\varphi^{-1}(W)$, where $W$ is an open subset of $O_m:=O\cap E_m$ in the chart $(V, \varphi, (O, C, E))$ on $X$.  Then $\mathscr{B}$ is a basis for a topology on $X_m$.   With this topology, the set $V_m:=V\cap X_m$  is an open subset of $X_m$ and  $\varphi\colon 
V_m\to O_m$ is a homeomorphism, so that  the tuple $(V_m, \varphi, (O_m, C_m, E_m))$ is a chart on $X_m$. 
Any two such  charts are sc-smoothly compatible and  the collection of such charts is an atlas on $X_m$.
We shall see below that the topology on $X_m$ is Hausdorff and paracompact. Hence the above sc-smoothly compatible
charts define a M-polyfold structure on $X_m$. We shall denote $X_m$ with this M-polyfold structure by $X^m$\index{$X^m$}
\index{Raising the index of $X$} and say
that it is obtained from $X$ by {\bf raising the index by $m$}. 
Also note that the M-polyfold $X$ inherits from the charts the  filtration 
$$
X=X_0\supset X_1\supset \cdots \supset X_\infty=\bigcap_{i\geq 0}X_i.
$$

\begin{lemma}\label{inculsion_continuous}
The inclusion map $i:X_{m+1}\to X_m$ is continuous for all $m\geq 0$.
\end{lemma}
\begin{proof} 
In view of the definition of topologies on $X_{m}$ and $X_{m+1}$, it suffices to show, that if $(V, \varphi, (O, C, E))$ is a chart on $X$ and $W$ is an open subset of $O_m$, then $\varphi^{-1}(W)\cap X_{m+1}$ is open in $X_{m+1}$. If $x\in \varphi^{-1}(W)\cap X_{m+1}$, then $\varphi (x)\in W\cap O_{m+1}\subset W\cap E_{m+1}$. Since $W$ is open in $O_{m}$, there exists an open set  $W'$ in $E_m$, so that $W=W'\cap O$. Hence, $\varphi (x)\in  (W'\cap E_{m+1})\cap O=W''\cap O$, where $W''=W'\cap E_{m+1}$ is an open subset of $E_{m+1}$. Hence $W''\cap O$ is an open subset of $O_{m+1}$ and $ \varphi^{-1}(W)\cap X_{m+1}=\varphi^{-1}(W''\cap O)$,  proving our claim.
\qed \end{proof}

\begin{theorem}\label{X_m_paracompact}\index{T- Metrizability of $X_m$}
Let $X$ be a M-polyfold. For every $m\geq 0$, the space $X_m$ is metrizable and, in particular, paracompact. In addition, the space $X_\infty$ is metrizable.
\qed
\end{theorem}
The proof is  postponed to  Appendix \ref{A2.2}.

In order to define the tangent $TX$ of the M-polyfold $X$, we start with its  local description in a chart.
\begin{definition}[{\bf Tangent space $T_xO$}]\index{D- Tangent space}
Let $(O, C, E)$ be a retract and $T(O,C,E)=(TO,TC,TE)$ its tangent,  so that $p\colon 
TO\rightarrow O^1$ is the tangent bundle over $O$.
The {\bf tangent space} $T_xO$ at a point $x\in O^1$ is the pre-image $p^{-1}(x)$, which is a Banach space. Note that only in the case that $x$ is  a smooth point, the tangent space  $T_xO$ has a natural sc-structure.
\qed
\end{definition} 
In the case that $x\in O_1$,  the tangent space $T_xO$\index{$T_xO$} is the image of the projection $Dr(x)\colon 
E\rightarrow E$, where $r\colon 
U\rightarrow U$ is any sc-smooth retraction
associated with  $(O,C,E)$ satisfying $r(U)=O$.

Next we consider tuples  $(x, V, \varphi,(O,C,E),h)$, in which $x\in X_1$ is a point in the M-polyfold $X$ on  level $1$ and  $(V,\varphi,(O,C,E))$ is a chart around the point $x$. Moreover,  $h\in T_{\varphi (x)}O$. 
Two such tuples, 
$$
(x, V, \varphi,(O,C,E),h)\quad \text{and} \quad (x',V', \psi, (O',C',E'),h'),
$$
are called equivalent, if  
$$x=x'\quad  \text{and}\quad T(\varphi\circ\psi^{-1})(\psi (x))h'=h.$$ 
\begin{definition}\label{tangent_space}\index{D- Tangent space}
The {\bf tangent space $TX$  of $X$}  as a set, is the collection of all equivalence classes $[(x,V, \varphi,(O,C,E),h)]$.  
\qed
\end{definition}

If $x\in X_1$ is fixed, the tangent space $T_xX$ of the M-polyfold $X$ at the point $x \in X$ is the subset of equivalence classes
$$T_xX=\{[(x, \varphi, V, (O, C, E), h)]\ \vert \,  h\in T_{\varphi (x)}O\}.$$
It has the structure of a vector space defined by 
\begin{eqnarray*}
&\lambda\cdot [(x, \varphi, V, (O, C, E), h_1)]+
\mu\cdot [(x, \varphi, V, (O, C, E), h_2)]&\\
&=
[(x, \varphi, V, (O, C, E), (\lambda h_1+\mu h_2)],&
\end{eqnarray*}
 where $\lambda,\mu\in \R$ and $h_1,h_2\in T_{\varphi(x)}O$.
Clearly,
\begin{equation}\label{equation_tangent_bundle}
TX=\bigcup_{x\in X_1}\{x\}\times T_xX.
\end{equation}

To  define a topology on the tangent space $TX$,  we first fix a chart $(V,\varphi,  (O, C, E))$ on $X$ and associate with it a subset $TV$ of $TX$,  defined as   
$$TV:
=\{[(x, \varphi, V, (O, C, E), h)]\ \vert \, x\in V,\, h\in T_{\varphi (x)}O\}.$$

We introduce the {\bf tangent map}  $T\varphi\colon TV\to TO$\index{Tangent map}  by 
$$T\varphi ([(x, \varphi, V, (O, C, E), h)])=(\varphi (x), h),$$
where $h\in T_{\varphi(x)}O$. If  $x\in X_1$ is fixed, then the map
$$T\varphi\colon 
T_xV\to T_{\varphi (x)}O$$
is a linear isomorphism.
Therefore, $T_xX$, the tangent space at $x \in X_1$ inherits the Banach space structure from $T_{\varphi(x)} O \subset E = E_0$. If $x \in X_\infty$, the tangent space $T_xX$ is a sc-Banach space, because the tangent space at the smooth point $\varphi(x) \in O_\infty$ is a sc-Banach space.

If $W$ is an an open subset of $TO$,  we define  the subset $\wt{W}\subset TV$ by $\wt{W}:=(T\varphi )^{-1}(W).$  
We denote by  $\mathscr{B}$ the collection of all such sets  $\wt{W}$ obtained by taking all the charts  $(\varphi, V, (O, C, E))$  of the atlas and all open subsets $W$ of the corresponding tangents $TO$. 

\begin{proposition}\label{op}\mbox{}\index{P- Properties of $TX$}
The following holds.
\begin{itemize}
\item[\em(1)]\ The collection $\mathscr{B}$ defines a basis for a Hausdorff  topology on $TX$.  
\item[\em(2)]\ The projection $p\colon 
TX\to X^1$ is a continuous and an open map.
\item[\em(3)]\ With the topology defined by $\mathscr{B}$,  the tangent space $TX$ of the M-polyfold $X$  is metrizable and hence, in particular,  paracompact.
\end{itemize}
\qed
\end{proposition}
The proof is postponed to  Appendix \ref{A2.11}.
In view of our definition of the topology on $TX$, the map 
$T\varphi\colon 
 TV\to TO$, associated with the chart $(\varphi, V, (O, C, E))$,  is a homeomorphism. Moreover, given two such maps 
$$T\varphi \colon 
TV\to TO\quad \text{and}\quad T\varphi'\colon 
TV'\to TO', $$ the composition $T\varphi' \circ (T\varphi)^{-1}\colon 
 T\varphi (TV\cap TV')\to T\varphi' (TV\cap TV')$ is explicitly of the form 
\begin{equation}\label{transition_1}
T\varphi' \circ (T\varphi)^{-1}(a, h)=(\varphi'\circ \varphi^{-1}(a), 
D(\varphi'\circ \varphi^{-1})(a)\cdot h)=T(\varphi'\circ\varphi^{-1})(a,h).
\end{equation}
Since the transition map  $\varphi' \circ  \varphi^{-1}$ between sc-retracts is sc-smooth, the composition in \eqref{transition_1} is also sc-smooth. 
In addition, $(TO, TC, TE)$ is a sc-smooth retract. Consequently, the tuples $(TV, T\varphi,  (TO, TC, TE))$  define a sc-smooth atlas on $TX$. Since, as we have proved above, $TX$ is paracompact,  the tangent space $TX$ of the M-polyfold $X$ is also  a M-polyfold.  The projection map $p\colon 
TX\to X^1$  is locally built on the bundle retractions $TO\rightarrow O^1$,  and the transition maps of the charts
are sc-smooth bundle maps. Therefore, 
$$p\colon 
TX\rightarrow X^1$$ 
is a sc-smooth M-polyfold bundle.

\begin{figure}[htb]
\centering
\def\svgwidth{30ex}
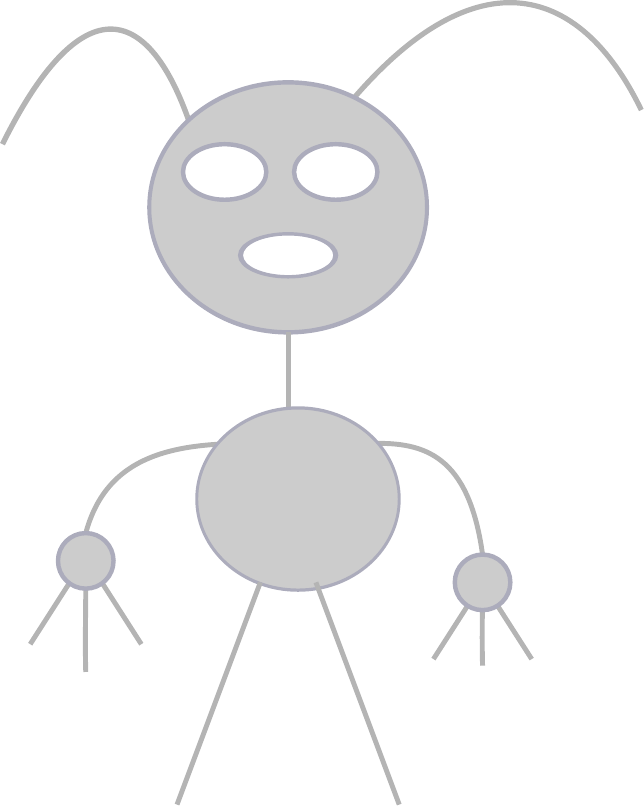
\caption{This chap has a  M-polyfold structure for which every point is smooth. }\label{fig:pict1}
\end{figure}

The  M-polyfold is the notion of a smooth manifold in our extended universe.  If 
$X$ is a  M-polyfold,  which consists entirely of smooth points, then it  has a tangent space at every point. 
There are finite-dimensional examples. For example,  the  chap depicted in Figure \ref{fig:pict1} has a M-polyfold structure,
for which $X=X_\infty$. It illustrates, in particular, that M-polyfolds allow to describe in a smooth way geometric objects having locally varying dimensions.
 For details in the construction of the chap and further illustrations, we refer to \cite{HWZ8.7}, Section 1, in particular, Example 1.22.

Next we introduce the notion of a sub-M-polyfold.

\begin{definition}\label{def_sc_smooth_sub_M_polyfold}\index{D- Sub-M-polyfold}
Let $X$ be a M-polyfold and let $A$ be a subset of $X$. The subset  $A$ is called a {\bf sub-M-polyfold} of $X$, 
if every $a\in A$ possesses  an open neighborhood $V$ and a sc-smooth retraction
$r\colon 
V\rightarrow V$, such that 
$$r(V)=A\cap V.$$
\qed
\end{definition}

\begin{proposition}\label{sc_structure_sub_M_polyfold}\index{P- Properties of sub-M-polyfolds}
A sub-M-polyfold $A$ of a  M-polyfold $X$ has, in a natural way, the structure of a  M-polyfold,
for which  the following holds. 
\begin{itemize}
\item[{\em (1)}]\ The inclusion map $i\colon A\rightarrow X$ is sc-smooth and a homeomorphism onto its image.
\item[{\em (2)}]\ For every $a\in A$ and every sc-smooth  retraction $r\colon V\rightarrow V$ satisfying  $r(V)=A\cap V$ and $a\in V$,
the map $i^{-1}\circ r\colon V\rightarrow A$ is sc-smooth.
\item[{\em (3)}]\ The tangent space $T_aA$ for a smooth $a\in A$ has a sc-complement in $T_aX$.
\item[{\em (4)}]\ If $a$ is a smooth point and $s\colon 
W\rightarrow W
$ is a sc-smooth retraction satisfying  $s(W)=W\cap A$ and $a\in W$, then the induced map
$W\rightarrow A$ is sc-smooth and $Ts(a)T_aX=T_aA$.
\end{itemize}
\end{proposition}
\begin{proof}
We  first define a sc-smooth atlas for $A$. We choose a point $a\in A$ and let $( \varphi, V, (O, C, E))$ be a chart of the M-polyfold $X$ around the point $a$.  By definition of a sub-M-polyfold,  there exists an open neighborhood $U$ of $a$ in $X$ and a sc-smooth retraction $r\colon 
U\to U$ satisfying  $r(U)=U\cap A$.  The set $W\subset X$, defined by  $W:=r^{-1}(U\cap V)\cap (U\cap V)$, is  open in $X$ and satisfies  $r(W)\subset W$ and $r(W)=W\cap A$. Hence, $\varphi (W)$ is an open subset of $O$, so that, in view of of part (1) 
Proposition  \ref{new_retract_1}, the tuple $(\varphi (W), C, E)$ is a sc-retract.  We may therefore  assume without loss of generality that $U=V=W$.  We define the sc-smooth map   $\rho\colon O\to O$ by 
$$\rho=\varphi\circ r\circ \varphi^{-1}.$$
From $r\circ r=r$ we deduce  $\rho\circ \rho=\rho$, so that $\rho$ is a sc-smooth retraction onto $\rho(O)=O'$. By the statement (2) in  Proposition  \ref{new_retract_1},  the triple $(O', C, E)$ is a sc-retract. Therefore there exists a relatively open subset $U'$ of the partial quadrant $C$ in $E$ and a sc-smooth retraction $r'\colon 
U'\to U'$ onto  $r'(U')=O'$.  Restricting the map $\varphi$ to $W\cap A$,  we set $\psi:=\varphi\vert W\cap A$ and compute,
$$\psi (W\cap A)=\varphi (W\cap A)=\varphi \circ r(W)=\varphi\circ r\circ \varphi^{-1}(O)=\rho(O)=O'.$$
Consequently,  
$$
\psi\colon 
 W\cap A\to O'$$
  is a homeomorphism and the triple  $( \psi, W\cap A, , (O', C, E))$ is a chart on $A$. 

In order to consider the chart transformation we take a second compatible chart $( \varphi', V', (\wh{O}, C', E'))$ of the M-polyfold $X$ around the point $a\in A$ and use it two construct the second chart $( \psi', W'\cap A, , (O'', C', E'))$ of $A$. 
We shall show that the second chart is compatible with the already constructed chart  $( \psi, W\cap A, , (O', C, E))$.  
The domain $\psi ((W\cap A)\cap (W'\cap A))$ of  the transition map 
\begin{equation}\label{transition_map_0}
\psi'\circ \psi^{-1}\colon 
\psi ((W\cap A)\cap (W'\cap A))\to \psi' ((W'\cap A)\cap (W\cap A))
\end{equation}
 is an open subset of $O'$, so that, in view of  (2) of Proposition \eqref{new_retract_1}, there exists a relatively open subset $U''\subset C$ and a sc-smooth retraction $s''\colon 
U''\to U''$ onto  $s''(U'')=\psi ((W\cap A)\cap (W'\cap A))$.  By construction,
 $$(\psi'\circ \psi^{-1})\circ s''=(\phi'\circ \phi^{-1})\circ s''.$$
 The chart transformation $\varphi'\circ \varphi^{-1}\colon 
O\to \wh{O}$ is sc-smooth so that by the chain rule 
the right-hand side is sc-smooth. Therefore, also the left-hand  is a sc-smooth map. So, in view of Definition \ref{tangent_retract} of a sc-smooth map between retracts, the transition map $\psi'\circ \psi^{-1}$ is sc-smooth.
  
We have shown that  the collection of charts $(\psi, W\cap A, (\psi (W\cap A), C, E))$  defines a sc-smooth atlas for $A$.  The sc-smooth structure on $A$  is defined by its equivalence class.

In order to prove the statement (2) in 
Proposition \ref{sc_structure_sub_M_polyfold} we use the above local coordinates, assuming, as above, that $U=V=W$. 
The inclusion map
$i\colon 
A\rightarrow X$ is, in the local coordinates, the inclusion 
$j\colon 
O'=r'(U')\to O$ which  is a  sc-smooth map since $r'\colon 
U'\to U'$ is a  sc-retraction. Conversely, the relations $\varphi\circ r\circ \varphi^{-1}(O)=O'$ and $U=\varphi^{-1}(O)$ show that the map
$i^{-1}\circ r\colon 
U\to A$ is sc-smooth because the retraction $r\colon 
U\to U$ is, by assumption, sc-smooth. This  proves  statement (2) of the proposition.

In order to prove the statement (3) we work in local coordinates,  and assume that  $X$ is given by the triple $(O,C,E)$ in which $O=r(U)$ and $r\colon 
U\to U$ is a retraction of the relatively open subset $U$ of $C$ in $E$. Then  $A$ is a subset of $O$ having the property that every point $a\in A$ possesses  an open neighborhood $V$ in $O$ and a sc-smooth retraction $s\colon 
V\rightarrow V$ onto $s(V)=A\cap V$. We now assume that $a\in A$ is  a smooth point and introduce  the map $t=s\circ r\colon 
U\to U$. Then 
$t\circ t=t$ and $t$ is a sc-smooth retraction onto the set $V\cap A$. Hence the tangent space $T_aA$ is defined by 
$$T_aA=Dt(a)E=Ds(a)\circ Dr(a)E=Ds(a)T_aO.$$
From $r\circ t=r\circ s\circ r=s\circ r=t$ we conclude 
$Dr(a)\circ Dt(a)E=Dt(a)E$ and hence $Dt(a)E\subset T_aO.$
Therefore, 
\begin{equation*}
\begin{split}
T_aO=Dr(a)E&=Ds(a)\circ Dr (a)+(I-Ds(a))Dr(a)E\\
&=
(T_aA)\oplus (I-Ds(a))(T_aO).
\end{split}
\end{equation*}
This proves the statement (3).  Using the same arguments,  the statement (4) follows and the proof of Proposition \ref{sc_structure_sub_M_polyfold} is complete.
\qed \end{proof}

\section{Degeneracy Index and Boundary Geometry}
Abbreviating ${\mathbb N}=\{0,1,2,3,..\}$\index{${\mathbb N}$}
  we shall introduce on the  M-polyfold $X$ the map $d_X\colon 
X\rightarrow {\mathbb N}$, called degeneracy index,  as follows.
We first take a smooth chart  $(V,\phi,(O,C,E))$ around the point $x$ 
and define the integer 
$$
d(x, V,\phi,(O,C,E))= d_C(\phi(x)),\index{$d(x, V,\phi,(O,C,E))$}
$$
where $d_C$ is  the index defined in 
Section \ref{subsection_boundary_recognition}.
In  other words, we record how many vanishing coordinates the image point $\phi(x)$ has in the partial quadrant $C$.
\begin{definition}\label{M_polyfold_degeneracy _index}\index{D- Degeneracy index $d_X$}\index{$d_X$}
The {\bf degeneracy} $d_X(x)$ at the point $x\in X$ is the minimum of all numbers 
$d(x, V, \phi,(O,C,E))$, where $(V,\phi,(O,C,E))$ varies over all smooth charts around  the point $x$.
The {\bf degeneracy index}  of the M-polyfold $X$ is the map $d_X\colon 
X\rightarrow {\mathbb N}$.
\qed
\end{definition}

The next lemma is evident.
\begin{lemma}\index{L- Local property of $d_X$}
Every point $x\in X$ possesses an open neighborhood $U(x)$, such  that $d_X(y)\leq d_X(x)$ for all $y\in U(x)$. 
\qed
\end{lemma}

From the definitions one deduces immediately the following result.
\begin{proposition}\label{newprop2.24}\index{P- Diffeomorphism invariance of $d_X$}
If $X$ and $Y$ are M-polyfolds and if $ f\colon (U,x)\rightarrow (V,f(x))$ is a germ of sc-diffeomorphisms
around the points $x\in X$ and $f(x)\in Y$,  then $$d_X(x)=d_Y(f(x)).$$
\qed
\end{proposition}

The index $d_X$ quantifies to which extend a point $x$ has to be seen as a boundary point. A more degenerate point has a higher index.

\begin{definition}\label{boundary_M_polyfold}\index{D- Boundary of an M-polyfold}
The subset  $\partial X=\{x\in X\ |\ d_X(x)\geq 1\}$ of $X$  is called the {\bf boundary of $X$}.
A  M-polyfold $X$ for which $d_X\equiv 0$ is called
a  {\bf M-polyfold without boundary}. \qed
\end{definition}

The relationship between $d_A$ and $d_X$  where  $A$ is a sub-M-polyfold of $X$ is described in the following lemma.
\begin{lemma}\label{k^0}\index{L- $d_A$ versus $d_X$}
If $X$ is  a M-polyfold and $A\subset X$ a sub-M-polyfold of $X$,   then
$$
d_A(a)\leq d_X(a)
$$
for all $a\in A$.
\end{lemma}
\begin{proof}
We take a point $a\in A$ and choose a chart $(\varphi, V, (O, C, E))$ around the point $a$, belonging to the atlas for $X$ and satisfying 
$$d_X(a)=d(a, \varphi , V, (O, C, E)).$$
The integer on the right-hand side remains unchanged, if we take a smaller domain, still containing $a$,  and replace $O$ by its image. Then, arguing as in 
Proposition \ref{sc_structure_sub_M_polyfold}, we may assume, that the  corresponding chart of the atlas  for $A$ is $(\psi, V', (O', C, E))$, where we have abbreviated $V'=V\cap A$, $\psi=\varphi\vert {V'}$, and $O'=\psi (V')$.  Then, 
$$d(a, \psi, V', (O', C, E))=d(a, \varphi, V, (O, C, E))=d_X(a),$$
and hence, taking the minimum on the left-hand side,
$d_A(a)\leq d_X(a),$ as claimed in the lemma.
\qed \end{proof}

If $E={\mathbb R}^k\oplus W$ is the sc-Banach space and $C=[0,\infty)^k\oplus W$ the partial quadrant in $E$,  we   define the linear subspace $E_i$ of $E$ by
$$
E_i=\{(a_1,\ldots, a_k,w)\in {\mathbb R}^k\oplus W\ \vert \ a_i=0\}.\index{$E_i$}
$$
With a subset $I\subset \{1,\ldots ,k\}$,  we associate the subspace 
$$
E_I=\bigcap_{i\in I} E_i.\index{$E_I$}
$$
In particular, $E_{\emptyset}=E$, $E_{\{i\}}=E_i$, and $E_{\{1,\ldots , k\}}=\{0\}^k\oplus W\equiv W$.  If $x\in C$, we denote by $I(x)$ the  set 
of indices $ i\in \{1,\ldots ,k\}$ for which $x\in E_i.$ We 
abbreviate
$$
E_x:=E_{I(x)}.
$$
Associated with  $E_i$ we have the closed half space $H_i$,  consisting of all elements $(a,w)$ in $ {\mathbb R}^k\oplus W$, satisfying  $a_i\geq 0$, 
$$
H_i=\{(a_1,\ldots, a_k,w)\in E\ |\ a_i\geq 0\}.\index{$H_i$}
$$
If $x\in C$, we define  the partial cone $C_x$ in $E$ as 
$$
C_x=C_{I(x)}:=\bigcap_{i\in I(x)}H_i.\index{$C_x$}
$$
We observe that $E_x\subset C_x$.

As an illustration we take  the standard quadrant $C\subset {\mathbb R}^2$ consisting of all $(x,y)$ with $x,y\geq 0$.
Then $C_{(0,0)}=C$, $C_{(1,0)}=\{(x,y)\ |\ y\geq 0\}$, $C_{(0,1)}=\{(x,y)\ |\ x\geq 0\}$ and $C_{(1,1)}={\mathbb R}^2$.
Moreover $E_{(0,0)}=\{0\}$, $E_{(1,0)}={\mathbb R}\oplus \{0\}$,  $E_{(0,1)}=\{0\}\oplus {\mathbb R}$, and $E_{(1,1)}={\mathbb R}^2$.
One should view  $C_{(x,y)}$ as a partial quadrant in the tangent space $T_{(x,y)}C= {\mathbb R}^2$ and $E_{(x,y)}\subset C_{(x,y)}$
as the maximal linear subspace.

In the following we shall put some additional structure on a sc-smooth retract $(O,C,E)$, which turns out to be useful.
We call a subset $C$ of a Banach space (or sc-Banach space) 
a {\bf cone} \index{Cone}provided it is closed, convex, and satisfies ${\mathbb R}^+C=C$ and $C\cap (-C)=\{0\}$.
If all properties except the last one hold,  we call $C$  a {\bf  partial cone}.\index{Partial cone}
In the following considerations $E={\mathbb R}^k\oplus W$ and $C=[0,\infty)^k\oplus W$.
The aim of the following  is to extract some information from the geometry of a retract $(O,C,E)$.
\begin{definition}\label{reduced_cone_tangent}\index{D- Reduced tangent space}\index{D- Partial cone}
Let $(O,C,E)$ be a sc-smooth retract and let $x\in O_\infty$  be a smooth point in $O$. The {\bf partial cone} $C_xO$  at  $x$ is defined as the following subset of the tangent space at $x$, 
$$
C_xO:=T_xO\cap C_x,\index{$C_xO$}
$$
where $C_x=\bigcap_{i\in I(x)}H_i$.
The {\bf reduced tangent space}  $T^{\textrm{R}}_xO$ is defined as the following subset of the tangent space at $x$, 
$$
T_x^{\textrm{R}} O =T_xO\cap E_x.\index{$T_x^{\textrm{R}}O$}
$$
\qed
\end{definition}
\begin{remark}\index{R- Remark on $C_xO$}
We have the inclusions
$$
T^R_xO\subset C_xO\subset T_xO, 
$$
and naively one might expect that $C_xO$ is a partial quadrant in $T_xO$. However, this is in general not the case.
\qed
\end{remark}
The reduced tangent space and the partial cone are  characterized in the next lemma.
\begin{lemma}\label{characterization_reduced_tangent}\index{L- Reduced tangent}
Let $(O,C,E)$  be as described before  and let  $x\in O_\infty$ be a smooth point in $O$. Then the following holds, where the  $\varepsilon$ occurring in the definitions may depend on $\alpha$.
\begin{itemize}
\item[\em(1)]\ $T^{\textrm{R}}_xO=\textrm{cl}( \{ \text{$\dot{\alpha}(0)\vert \, \alpha\colon 
(-\varepsilon, \varepsilon)\to O$ is sc-smooth and $ \alpha(0)=x$}\} )$.
\item[\em(2)]\ $C_xO=\text{cl}( \{\text{$\dot{\alpha}(0)\vert \,  \alpha\colon 
[0, \varepsilon)\rightarrow O$ is sc-smooth and $ \alpha(0)=x$}\})$.
\item[\em(3)]\ $T_xO =C_xO-C_xO$.
\end{itemize}
Here $\dot{\alpha}(0)=\frac{d}{dt}\alpha (t)_{\vert t=0}$ stands for the derivative of the sc-smooth path $\alpha$ in the parameter $t$ varying in $(-\varepsilon, \varepsilon)$ resp. in $[0, \varepsilon)$.
\end{lemma}
\begin{proof}  We assume that $U\subset C$ is relatively open and $r:U\rightarrow U$ is a sc-smooth retraction with $O=r(U)$.
Let us denote for $x\in {\mathbb R}^k\oplus W$
by $x_1\ldots,x_k$ its  coordinates in ${\mathbb R}^k$.\par
 
(1)\,  We first  introduce the set  
\begin{equation}\label{set_gamma}
\Gamma= \{ \text{$\dot{\alpha}(0)\vert \, \alpha\colon 
(-\varepsilon, \varepsilon)\to O$ is sc-smooth and $ \alpha(0)=x$}\},
\end{equation}
and observe, that $\alpha ((-\varepsilon, \varepsilon))$ is contained in $ O_\infty$, and   $\dot{\alpha}(0)\in E_\infty$.
Since the closure of $(T^{\textrm{R}}_xO)_\infty$ is  equal to $T^{\textrm{R}}_xO$, it is enough to prove that  $(T^{\textrm{R}}_xO)_\infty\subset \Gamma$ and $\Gamma\subset (T^{\textrm{R}}_xO)_\infty$.  As for the  first inclusion,  we take $v\in (T^{\textrm{R}}_xO)_\infty$. Hence  $x+tv\in U$ for $\abs{t}$ small. This follows from the fact that for $i\in I(x)$,  we have $a_i=v_i=0$ and for $i\not \in I(x)$, we have $x_i>0$. Since $x$ and $v$ are smooth points, $x+tv\in U_\infty$. Then $\alpha (t)=r(x+tv)$ is defined for $\abs{t}$ small, takes values in $O_\infty$,  and $\alpha (0)=x$.  By the chain rule, 
$$\dot{\alpha}(0)=Dr(x)v=v,$$
since $v\in T_xO=\text{image of $Dr(x)$}$. Hence, $v\in \Gamma$ and $(T^{\textrm{R}}_xO)_\infty\subset \Gamma$, as claimed. Conversely,  if $\alpha\colon 
(-\varepsilon, \varepsilon)\to O$ is a  sc-smooth path  satisfying  $\alpha (0)=x$,   then $r(\alpha (t))=\alpha (t)$, so that, applying the chain rule, we find 
$$Dr(x)\dot{\alpha}(0)=\dot{\alpha} (0).$$
This shows that  $\dot{\alpha} (0)$ is a smooth point belonging to $T_xO$. If $i\in I(x)$, then $\alpha_i(0)=x_i=0$ and since 
$\alpha_i(t)\geq 0$ for all $t\in (-\varepsilon, \varepsilon)$, we conclude  that $\dot{\alpha}_i(0)=0$, so that 
$\dot{\alpha}(0)\in E_x$.  
Hence, 
$\dot{\alpha}(0)\in T_xO\cap E_x=T_x^{\textrm{R}}O$  
and since 
$\dot{\alpha}(0)$ 
is a smooth point, 
$\dot{\alpha}(0)\in (T_x^{\textrm{R}}O)_\infty$.  
Therefore,
 $\Gamma \subset (T_x^{\textrm{R}}O)_\infty$,  and the proof of (1) is complete.\par
 
(2)\, The proof is along the same line, except that considering a sc-smooth path $\alpha\colon 
[0,\infty)\to O,$ we conclude that $Dr(x)\dot{\alpha}(0)\in C_x$ and since $Dr(x)\dot{\alpha}(0)=\dot{\alpha}(0)$, we find that $C_xO$ is a subset of the right-hand side of (2). Conversely, we take a smooth point $v\in C_xO$ and consider the path $\alpha (t)=r(x+tv)$ defined for $t\geq 0$ small. Then, $\dot{\alpha}(0)=Dr(x)v$ belongs to the right-hand side of (2) and since $C_xO$ is closed, the result follows.\par

(3)\, Clearly, $C_xO-C_xO\subset T_xO$. Conversely, let $h\in T_xO$. Then $h=Dr(x)k$, where $k=(a, w)\in E.$  
If $i\in I(x)$, we set $$a_i^\pm=\frac{\abs{a_i}\pm a_i}{2}.$$ 

Then, $a_i^\pm\geq 0$ and $a_i=a^+_i-a^-_i$. Now we define elements $k^\pm\in E$ as follows.  First, $k^+=(b, w),$ where $b_i=a_i$ if $i\not \in I(x)$ and $b_i=a^+_i$ if $i\in I(x)$. The element $k^-$ is defined as $k^-=(c, 0),$ where $c_i=0$ if $i\not \in I(x)$ and $c_i=a_i^-$ if $i\in I(x)$. Then $k=k^+-k^-$ and if $h^\pm =Dr(x)k^\pm$, then, by (2), we have $h^\pm\in  C_xO$ and 
$h=h^+-h^-\in C_xO-C_xO$. The proof of (3) and hence the proof of Lemma \ref{characterization_reduced_tangent}  is complete.

\qed \end{proof}

From the characterization of $T_xO$ and $C_xO$ in 
Lemma \ref{characterization_reduced_tangent} we  deduce immediately the next  proposition,
where $(O,C,E)$ is as described before and similarly $(O',C',E')$ with $E'={\mathbb R}^{k'}\oplus W'$ and $C'=[0,\infty)^{k'}\oplus W'$ 
\begin{proposition}\label{reduced_tangent_under_sc_diff}
Let $(O,C,E)$ and $(O',C',E')$ be sc-smooth retracts as just described,   and let $x\in O_\infty$. If $f\colon 
(O,x)\rightarrow (O',f(x))$
is  a germ of a sc-diffeomorphism mapping $x\in O_\infty$ onto $f(x)=y\in O'_\infty$, then 
$$
Tf(x)T^{\textrm{R}}_xO= T^{\textrm{R}}_{y}O'\quad  \text{and}\quad Tf(x)C_xO= C_yO'.
$$
\qed
\end{proposition}
We note that in view of the previous proposition the definition of the spaces  $T^R_xO$ and $C_xO$ is natural
and they are respected by the tangent of a  local sc-diffeomorphism.  Hence we can define
for an arbitrary  sc-retract $(O,C,E)$ and a point $x\in O_\infty$ the {\bf reduced tangent}\index{Reduced tangent} $T^R_xO$ and the {\bf partial cone}\index{Partial cone} $C_xO$.
In case  there exists a linear sc-isomorphism  $S:T_xO\rightarrow {\mathbb R}^n\oplus W$ satisfying $S(C_x)=[0,\infty)^n\oplus W$ we shall call $C_x$ a {\bf partial quadrant}\index{Partial quadrant}.
\begin{remark}\index{R- Tameness and partial quadrants}
Later on we shall introduce the notion of a tame M-polyfold. In view of the following Theorem \ref{hofer} and the later Proposition \ref{tame_equality} the tameness
implies that for this particular case the partial cones are always partial quadrants.
\qed
\end{remark}

A sc-smooth retract $O$ associated with a triple 
$(O, C, E)$ is a M-polyfold and, recalling Definition \ref{M_polyfold_degeneracy _index}, its degeneracy index $d_O(x)$ at the point $x$ is the integer 
$$d_O(x)=\min d_{C'}(\varphi (x))$$
where the minimum is taken over all germs of sc-diffeomorphisms 
$$
\varphi\colon (O, x)\to (O', \varphi (x))
$$
 into sc-smooth retracts $O'$ associated with $(O', C', E')$. The integer $d_{C'}(\varphi (x))$ is introduced in 
Section \ref{subsection_boundary_recognition}.

\begin{theorem}\label{hofer}\index{T- Basic properties of $d_O$}
Let  $(O,C,E)$ be a smooth retract and let $d_O$ be the degeneracy index of $O$. 
If $x\in O$ is a smooth point, we have the 
inequality
$$
\dim (T_xO/T^R_xO)\leq  d_O(x).
$$
Moreover, if
$\dim(T_xO/T^R_xO)=d_O(x)$, then 
$C_xO$ is a partial quadrant in $T_xO$. 
\end{theorem}
\begin{proof}
Let $x$ be a smooth point of the retract $O$. 
By 
Proposition \ref{reduced_tangent_under_sc_diff} 
the dimension of $T_xO/T^R_xO$ is preserved 
under germs of sc-diffeomorphisms.  Hence, in view of the definition of $d_O(x)$, we may assume, without loss of generality,  that
$$d_O(x)=d_C(x)\equiv d$$ 
Moreover, without loss of generality, we may  assume that $E=\R^k\oplus W$, 
$C=[0,\infty)^k\oplus W$
and $x=(0,\ldots ,0,x_{d+1},\ldots , x_k, w)$, where $x_i>0$ for $d+1\leq i\leq k$ and 
$w\in W$. We recall that if $v=(a, b, w)\in T_x^RO\subset \R^d\oplus R^{k-d}\oplus W$, then $a=0$,  and if $v=(a, b, w)\in C_xO$, then $a_i\geq 0$ for $1\leq i\leq d$.

In order to prove the first statement, we choose smooth vectors $v^1,\ldots v^l$ in $T_xO$ such that $(v^j+T_x^RO)_{1\leq j\leq l}$ are linearly independent in the vector space $T_xO/T_x^RO$. 
Representing $v^j=(a^j,b^j, w^j)\in \R^d\oplus R^{k-d}\oplus W$, we claim that the vectors $(a^j)_{1\leq j\leq l}$ are linearly independent in $\R^d$. 
Indeed, assuming  that  
$\sum_{j=1}^l\lambda_ja^j=0$, we have  
$$\sum_{j=1}^l\lambda_jv^j=\bigl(0, \sum_{j=1}^l\lambda_jb^j, \sum_{j=1}^l\lambda_jw^j\bigr)\in T^R_xO,$$
hence 
$$\sum_{j=1}^l\lambda_j\bigl(v^j+T_x^RO\bigr)=\bigl(\sum_{j=1}^l\lambda_jv^j\bigr)+T_x^RO=T_x^RO.$$
Since $(v^j+T_x^RO)_{1\leq j\leq l}$ 
are linearly independent  in $T_xO/T_x^RO$,  we conclude that $\lambda_1=\ldots =\lambda_l=0$, proving our claim. This implies that the vectors $a^1,\ldots ,a^l$ are linearly independent in $\R^d$. Therefore, $l\leq d$ and hence $\dim(T_xO/T^R_xO)\leq  d=d_O(x)$, proving the first statement of the theorem.

In order to prove the second statement, 
we assume that $\dim(T_xO/T^R_xO)=d_O(x)$. 
If now $(v^j+T^R_xO)_{1\leq j\leq d}$ is a basis of $T_xO/T^R_xO$, then  representing  $v^j=(a^j,b^j, w^j)\in T_xO\subset   \R^d\oplus \R^{k-d}\oplus W$ and arguing as above,   the vectors $a^j$ for $1\leq j\leq d$ form a basis of  $\R^d$. 
Consequently, the map 
$\Phi\colon 
T_xO/T^R_xO\to \R^d$, defined by 
$$\Phi (v+T^R_xO)=\Phi ((a, b, w)+T^R_xO)=a,$$
is a linear  isomorphism. Moreover, if $v=(a, b, w)\in C_xO$  so that $a_j\geq 0$ for $1\leq j\leq d$, then 
$$\Phi (v+T^R_xO)=\Phi ((a, b, w)+T^R_xO)\in [0,\infty)^d.$$

Denoting by $e^j$ for $1\leq j\leq d$ the standard basis of $\R^d$, we introduce $\Phi^{-1}(e^j)=\wh{v}^j+T^R_xO$. By definition, the vectors $v^j$ are of the form 
$\wh{v}^j=(e^j, \wh{b}^j, \wh{w}^j)$ and  are linearly independent in $\R^d\oplus \R^{k-d}\oplus W$.   
If now $v=(a, b, w)\in T_xO$, we have the decomposition 
\begin{equation*}
\begin{split}
v=(a, b, w)&=\sum_{j=1}^da_j\wh{v}^j+\bigl(v-\sum_{j=1}^da_j\wh{v}^j\bigr)
\end{split}
\end{equation*}
where  $a=(a_1,\ldots ,a_d)\in \R^d$. Since the second term  on the right-hand side belongs to $T_x^RO$,  we  have the following decompositions of the tangent space $T_xO$ and of $C_xO$,
\begin{align*}
T_xO=\R \wh{v}^1\oplus \ldots \oplus \R \wh{v}^d\oplus T_x^RO\quad \text{and}\quad 
C_xO=\R^+\wh{v}^1\oplus \ldots \oplus \R^+\wh{v}^d\oplus T_x^RO.
\end{align*}
Therefore, the map $T\colon T_xO\to \R^d\oplus T_xO$, defined  by 
$$T(\lambda_1\wh{v}^1,\ldots ,\lambda_d\wh{v}^d, w)=(\lambda_1, \ldots ,\lambda_d, w)$$
is a sc-linear isomorphism
satisfying $T(C_xO)=[0,\infty )^d\oplus T_x^RO$.
 Hence $C_xO$ is a partial quadrant in $T_xO$ and the second statement of Theorem \ref{hofer} is proved.  
\qed \end{proof}

The  partial quadrant $C$ in $E$ is the  image of the sc-smooth retraction $r=\mathbbm{1}_{C}\colon C\to C$. In particular, $C$ is a M-polyfold which we denote by $X_C$. 
In Section \ref{subsection_boundary_recognition}, we have defined the map $d_C\colon 
C\to \N$. Above, we have defined the degeneracy index $d_{X_C}\colon 
X_C\to \N$ of the M-polyfold $X_C$.   By definition, 
$d_{X_C}\leq d_C$ and we shall prove that $d_{X_C}=d_C$.  We may assume  without loss of generality
that $C=[0,\infty)^k\oplus W$ and $E={\mathbb R}^k\oplus W$.  For a smooth point $x\in C$ we have $T_xC=E$ and $T_x^RC=E_x$,  implying 
$\dim(T_xC/T_x^RC)=d_C(x)$. This, of course, also holds for any partial quadrant $C\subset E$. Hence, 
$$
d_C(x) = \dim(T_xX_C/T_x^RX_C)
$$
for a smooth point $x\in C$. From this we deduce  the following corollary of Theorem \ref{hofer}
\begin{corollary}\label{equality_of_d}\index{C- Computation of $d_{X_C}$}
Let $C$ be a partial quadrant in a sc-Banach space. Considering $C$ as a M-polyfold,  denoted by $X_C$,  we have the equality 
$$
d_C=d_{X_C}.
$$
\end{corollary}
\begin{proof}
Let us first take  a smooth point $x\in C$. In view of Theorem \ref{hofer}
$$
\dim(T_x X_C/T_x^R X_C)\leq d_{X_C}(x).
$$
Since $d_C(x)=\dim(T_x X_C/T_x^R X_C)$ as we have just seen,  it follows that $d_C(x)\leq d_{X_C}(x)$. By definition of  $d_{X_C}(x)$,  we always have the 
inequality $d_{X_C}(x)\leq d_C(x)$. Consequently,  
$$
d_{X_C}(x)=d_C(x)\quad  \text{if $x\in C_\infty$.}
$$
If $x\in C$ is arbitrary, we take  a sequence of smooth points $x_k\in C$ converging to $x$ and satisfying $d_C(x_k)=d_C(x)$. Hence, 
$$
d_C(x)=d_C(x_k) =d_{X_C}(x_k).
$$
In view of the definition of $d_{X_C}$ we  find a sc-diffeomorphism $f\colon 
U(x)\rightarrow O'$, where $U(x)\subset C$ is relatively open,  and $(O',C',E')$ is a sc-smooth retract, so that
$$
d_{X_C}(x)=d_{C'}(f(x)).
$$
Then $f(x_k)\rightarrow f(x)$ and trivially $d_{C'}(f(x_k))\leq d_{C'}(f(x))$ for large $k$. Hence, for large $k$, 
\begin{equation*}
\begin{split}
d_C(x)&=d_C(x_k)=d_{X_C}(x_k)\\
&=d_{X_{C'}}(f(x_k))\leq  d_{C'}(f(x_k))\leq d_{C'}(f(x)) =d_{X_C}(x).
\end{split}
\end{equation*}
Since $d_{X_C}(x)\leq d_C(x)$,  we  conclude 
$d_{X_C}(x)=d_C(x)$ and the proof of 
Corollary \ref{equality_of_d} is complete.
\qed \end{proof}
From now on we do not have to distinguish between  the index $d_C$  defined for partial quadrants and the degeneracy index $d_{X_C}$, 
where we view $C$ as a  M-polyfold.

\section{Tame M-polyfolds}\label {subsec_tame_m_polyfolds}
In order to define spaces whose boundaries have more structure, we introduce the notion of a tame M-polyfold
and of tame retractions and tame retracts. We start with some basic geometry.

Let $C\subset E$ be a partial quadrant in a sc-Banach space $E$. We begin with the particular case
$E={\mathbb R}^k\oplus W$ and $C=[0,\infty)^k\oplus W$.  We recall  the linear sc-subspace $E_i$ of codimension $1$ defined as  
\begin{equation}\label{k^4}
E_i=\{(a_1,\ldots, a_k,w)\in {\mathbb R}^k\oplus W\ \vert \,  a_i=0\}.
\end{equation}
Associated with  $E_i$, there is  the closed half space $H_i$  consisting of all elements $(a,w)$ in $ {\mathbb R}^k\oplus W$ satisfying  $a_i\geq 0$, 
\begin{equation}\label{k^5}
H_i=\{(a_1,\ldots, a_k,w)\in E\ |\ a_i\geq 0\}.
\end{equation}

In the general case of  a partial quadrant $C$ in $E$,  we can describe the above definitions in a more  intrinsic way as follows.
 We  consider the set $\{e\in C\,\vert \, d_C(e)=1\}$ of boundary points. This set has exactly $k$ connected components, which we denote by $A_1,\ldots ,A_k$. Each component $A_i$ lies in the smallest subspace $f_i$ of $E$ containing $A_i$. We call $f_i$ an {\bf extended face} and denote by ${\mathcal F}$ the set of all extended faces. The set ${\mathcal F}$ contains exactly $k$ extended faces. 
 Given an extended face $f\in {\mathcal F}$,  we denote  by $H_f$ the closed half subspace of $E$ which contains $C$.  
 In the special case $C=[0,\infty)^k\oplus W\subset \R^k\oplus W$, 
 the  extended  faces  $f_i$ are the subspaces $E_i$, and the half spaces $H_{f_i}$ are the half subspaces $H_i$. 
 
 If  $e\in C$,  we introduce the set of all extended faces containing $e$ by 
 $$ {\mathcal  F}(e)=\{f\in {\mathcal F} \vert \, e\in f\}.$$
Clearly, 
 $$d_C(e)=\# {\mathcal F}(e).$$
 \begin{definition}\label{new_def_2.33}
The {\bf partial quadrant $C_e$}\index{D- Partial quadrant $C_e$}  associated  with  $e\in C$ is defined as 
$$
C_e:= \bigcap_{f\in {\mathcal F}(e)} H_f.\index{$C_e$}
$$
The  {\bf minimal linear subspace}\index{D- Minimal linear subspace} associated with  $e\in C$ is defined by
$$
E_e:=\bigcap_{f\in {\mathcal F}(e)} f.\index{$E_e$}
$$
\qed
\end{definition}
Clearly, the following inclusions hold,
$$
E_e\subset C_e\subset E.
$$
For an interior point $x\in C$, i.e. a point satisfying $d_C(x)=0$, we set $E_e=C_e=E$.  The codimension of $E_e$  in $E$ is precisely $d_C(e)$.
The maximal value $d_{C_e}$ attains  is $d_C(e)$.

Next we introduce a special class of sc-smooth retracts. 
\begin{definition}[{\bf Tame sc-retraction}] \label{tame_retarctions}\index{D- Tame sc-smooth retraction}
Let $r\colon U\to U$ be a sc-smooth retraction defined on a relatively open subset $U$ of a partial quadrant $C$ in the sc-Banach space $E$.
The sc-smooth retraction  $r$ is  called a {\bf tame sc-retraction}, if  the following two conditions are satisfied.
\begin{itemize}
\item[(1)]\ \text{ $d_C(r(x))=d_C(x)$ for all  $x$ in $U$.}
\item[(2)]\ At every smooth point $x$ in $O:=r(U)$, there exists a sc-subspace $A\subset E$, such that $E=T_xO\oplus A$ and $A\subset E_x$. 
\end{itemize}
A sc-smooth retract $(O,C,E)$ is called a {\bf tame sc-smooth retract},  if $O$ is the image of a sc-smooth  tame
retraction.
\qed
\end{definition}

Let $x$ be a smooth point in the tame sc-retract $O$ and let $A\subset E_x$ be a sc-complement of the tangent space $T_xO$ as guaranteed by condition (2) in Definition \ref{tame_retarctions}, so that 
\begin{equation}\label{eq_T_oplus_A}
E=T_xO\oplus A.
\end{equation}
We recall that $T_xO=Dr(x)E$ and hence 
\begin{equation}\label{eq_T_oplus_A_1}
\begin{split}
E&=Dr(x)E+(\mathbbm{1}-Dr(x))E\\
&=T_xO+(\mathbbm{1}-Dr(x))E.
\end{split}
\end{equation}
Applying  the projection $\mathbbm{1}-Dr(x)$ to the equation \eqref{eq_T_oplus_A} and using that $(\mathbbm{1}-Dr(x))T_xO=0$, we obtain 

\begin{equation}\label{eq_T_oplus_A_2}
(\mathbbm{1}-Dr(x))E=(\mathbbm{1}-Dr(x))A.
\end{equation}

We claim that $(\mathbbm{1}-Dr(x))E\subset E_x$. In order to prove this claim we recall that $A\subset E_x$, so that, in view of \eqref{eq_T_oplus_A_2}, 
it is sufficient to prove that $Dr(x)E_x\subset E_x$.  We may assume that $E=\R^k\oplus W$ and $x=(0,\ldots ,0,x_{d+1},\ldots,x_k,w)$ with 
$x_{d+1},\ldots,x_k>0$. We choose a smooth point $y\in E_x$ so that $y=(0,\ldots,0,y_{d+1},\ldots,y_k,v)$. If $\abs{\tau}$ is small, then 
$x_\tau:=x+\tau y$ belongs to  the partial quadrant $C$ and has 
 the first $d$ coordinates vanishing. By condition (1) in Definition \ref{tame_retarctions} of a tame sc-retraction $r$, we have  $d_C(r(x+\tau y))=d_C(x+\tau y).$  Hence 
 the first 
 $d$ coordinates of $r(x+\tau y)$ 
 vanish, and from
 $$
 \dfrac{d}{d\tau}r(x+\tau y)\bigl\lvert_{\tau=0}=Dr(x)y
 $$
 we conclude that the first $d$ coordinates of $Dr(x)y$ 
 vanish.  The same is true if $y$ is on level $0$ in $E_x$. So, $Dr(x)E_x\subset E_x$ and hence, in view of \eqref{eq_T_oplus_A_2}, we have verified that $(\mathbbm{1}-Dr(x))E\subset E_x$, as claimed.
In view of  \eqref{eq_T_oplus_A_1}, we can therefore always assume, without loss of generality,  that in condition (2) of Definition \ref{tame_retarctions} the complement $A$ of $T_xO$ is equal to $A=(\mathbbm{1}-Dr(x))E$. 
Summarizing the discussion we have established the following proposition. 
\begin{proposition}\label{IAS-x}\index{P- Properties of tame retractions}
Let $U\subset C\subset E$ be a relatively open subset of a partial quadrant in a sc-Banach space and let $r\colon 
U\rightarrow U$  be 
a sc-smooth tame retraction. Then,  for every smooth point $x\in O=r(U)$,  the sc-subspace $(\mathbbm{1}-Dr(x))E$ is a subspace of $E_x$, so that in condition (2) of Definition \ref{tame_retarctions} we can take $A=(\mathbbm{1}-Dr(x))E$.
\qed
\end{proposition}

A particular example of sc-smooth retractions occurring in applications are splicings.
\begin{definition}\index{D- Splicings} 
A  sc-smooth  {\bf splicing} consists of the following data. A relatively open neighborhood  $V$ of $0$ in $[0,\infty)^k\times {\mathbb R}^{n-k}$
and a map which associates to $a\in V$ a sc-projection operator $\pi_a:E\rightarrow E$, where $E$ is an sc-Banach space so that 
 the map $V\times E\rightarrow E:(a,e)\rightarrow \pi_a(e)$ is sc-smooth. 
 \qed
 \end{definition}
 \begin{lemma}\index{L- Tameness of splicings}
 Let $V\ni a\rightarrow \pi_a$ be an sc-smooth splicing. Then 
 defining $U=V\times E$ and 
$r:U\rightarrow U$ by
$$
r(a,e)=(a,\pi_a(e))
$$
we obtain a sc-smooth tame retraction. 
\end{lemma}
\begin{proof}It is clear that with $D=[0,\infty)^k\times {\mathbb R}^{n-k}\times E$ the set $U$ is relatively open in $D$ 
and it holds that  
$d_D(a,e) = d_D(r(a,e))$.   Suppose $(a,e)$ is a smooth point in $O=r(U)$ implying that $\pi_a(e)=e$.
Then $T_{(a,e)}O$  consists of all $(b,k)$ with 
$$
k=\pi_{a}(k) +\sum_{i=1}^n (\partial (\pi_a(e))/\partial a_i)\cdot b_i.
$$
This in particular implies that 
\begin{eqnarray}\label{EQN209}
\sum_{i=1}^n (\partial (\pi_a(e))/\partial a_i)\cdot b_i\in  R(Id-\pi_a)=\ker(\pi_a).
\end{eqnarray}
Denote by $F$ the collection of all $(0,\ell)$ with $\pi_a(\ell)=0$ and observe that $F\subset \{0\}^n\times E\subset T^R_{(a,e)}(D)$.
If $(b,f)\in F\cap T_{(a,e)}O$ it follows that $b=0$ and $\pi_a(f)=f$ and $\pi_a(f)=0$, implying $(b,f)=(0,0)$.
On the other hand if $(b,f)\in T_{(a,e)}(U)=T_{(a,e)}(D)$, we can write,  abbreviating $\tau:= \sum_{i=1}^n (\partial (\pi_a(e))/\partial a_i)\cdot b_i$
$$
(b,f) = (b,  \pi_a(f) + \tau) + (0, (I-\pi_a)(f) - \tau).
$$
It holds $\pi_a((I-\pi_a)(f) - \tau)=-\pi_a(\tau)=0$ in view of (\ref{EQN209}), so that $(0, (I-\pi_a)(f) - \tau)\in F$. Further
$$
(Id-\pi_a)( \pi_a(f) + \tau)=\tau 
$$
showing that $ (b,  \pi_a(f) + \tau)\in T_{(a,e)}O$. Hence we have the sc-direct sum $T_{(a,e)}U= T_{(a,e)}O \oplus F$, with $F\subset T_{(a,e)}^R(U)$ implying that
$r$ is a tame sc-smooth retraction.
\qed \end{proof}
Next  discuss the tame sc-smooth retracts in more detail. 
\begin{proposition}\label{tame_equality}\index{P- Basic equality for tame retracts}
Let $(O,C,E)$ be a tame sc-smooth retract,  and let $x\in O_\infty$ be a smooth point of $O$.
Then, $T_x^{\textrm{R}}O$ is a sc-Banach space of codimension $d_O(x)$ in $T_xO$, so that 
\begin{equation}\label{opp}
\dim(T_xO/T^R_xO)=d_O(x). 
\end{equation}
In particular,   $C_xO$ is a partial quadrant in the tangent space $T_xO$.
In addition,  for every point $x\in O$, we have the equality
\begin{equation}\label{equality_d_O_and_d_C_1}
d_O(x)=d_C(x).
\end{equation}
\end{proposition}
\begin{proof}
We assume that $O=r(U)$ is the retract in $U\subset C\subset E$ associated to  the tame sc-retraction $r\colon 
U\to U$. Moreover, we may assume that $C=[0,\infty )^k\oplus W\subset E=\R^k\oplus W$. 

Then the  condition (2) of Definition \ref{tame_retarctions}  says that,  for every $x\in O_\infty$ there exists a sc-subspace $A$ of $E$, satisfying  $E=\text{im}\ Dr(x)\oplus A =T_xO\oplus A$ and $A\subset E_x$, where $E_x=E_{I(x)}$.  
Therefore, 
\begin{equation}\label{direct_sum_1}
E_x=(T_xO\oplus A)\cap E_x= (T_xO\cap E_x)\oplus (A\cap E_x)=T_x^{\textrm{R}}O\oplus A, 
\end{equation}
since  $A\subset E_x$. The space $E_x$ has codimension $d_C(x)$ in $E$, so that 
\begin{equation}\label{direct_sum_2}
d_C(x)=\dim\bigl(E/E_x \bigr)= \dim \bigl(T_xO\oplus A /T_x^{\textrm{R}}O\oplus A \bigr)=  \dim \bigl(T_xO /T_x^{\textrm{R}}O \bigr).
\end{equation}
Employing Theorem \ref{hofer},  we conclude from 
\eqref{direct_sum_2} that
$d_C(x)\leq d_O(x) $ for all $x\in O_\infty$.
By definition,  $d_O(x)\leq d_C(x)$ for all $x\in O$,  and hence, 
\begin{equation}\label{eq_d_C=d_O}
d_C(x)=d_O(x)\quad \text{if  $x\in O_\infty$.}
\end{equation}
Consequently, $\dim T_xO/T_x^RO=d_O(x)$. Employing 
Theorem \ref{hofer} once more, $C_xO$ is a partial quadrant in the tangent space $T_xO$.

It remains to prove that $d_O(x)=d_C(x)$ for all (not necessarily smooth) points $x\in O$. If $x=(a, w)\in O$,  we take a sequence $(w_j)\subset W_\infty$  converging to $w$  in $W$. 
Then the sequence $x_j\in O$, defined by $x_j=r(a, w_j)$, consists of smooth points and converges to $x$ in $E$.
Since, by assumption, the  retraction $r$ is tame, 
$$d_C(x_j)=d_C(r(a, w_j))=d_C(a, w_j)=d_C(a, w)=d_C(x)$$
and, using  \eqref{eq_d_C=d_O}, we conclude 
that $d_O(x_j)=d_C(x_j)=d_C(x).$
By definition of $d_O$, we find an open neighborhood $V'\subset O$ around $x$, such that 
$d_O(y)\leq d_O(x)$ for all $y\in V'$. Consequently, for large $j$,
$$d_O(x)\geq \lim_{j\to \infty}d_O(x_j)=d_C(x).$$
In view of $d_O(x)\leq d_C(x)$, we conclude 
$d_O(x)=d_C(x)$, and the proof of Proposition \ref{tame_equality} is complete.
\qed \end{proof}

Among all M-polyfolds, there is a distinguished class of M-polyfolds which are modeled on tame retracts.
These turn out to have an interesting and useful boundary geometry.

\begin{definition}\label{def_tame_m-polyfold}\index{D- Tame M-polyfold}
A  {\bf tame M-polyfold} $X$ is a M-polyfold which possesses an equivalent sc-smooth atlas whose charts are all modeled on {\bf tame} sc-smooth retracts.
\qed
\end{definition}
By Proposition \ref{reduced_tangent_under_sc_diff}, the following concepts for any M-polyfold  are well-defined and independent of the choice of the charts.

\begin{definition}\label{def_partial_cone_reduced_tangent}\index{D- Reduced tangent space}
Let $X$ be a  M-polyfold. For a smooth point $x\in X$,  the {\bf reduced tangent space $T^R_xX$}\index{$T^R_xX$}  is, by definition, the sc-subspace of the tangent space
$T_xX$   which, by definition,   is the preimage of $T_o^RO$ under any  chart $\psi\colon 
(V,x)\rightarrow (O,o)$, so that 
$$ 
T\psi(x) (T_x^RX) = T_o^RO = T_oO \cap E_o.
$$
 For a smooth point $x\in X$,  the {\bf partial cone} $C_xX$\index{Partial cone}\index{$C_xX$} is the closed convex subset of $T_xX$ which under a sc-smooth chart $\psi$ as above is mapped onto $C_oO$, i.e.
$$
T\psi(x)(C_xX) =  T_oO \cap C_o.
$$
\qed
\end{definition}

\begin{remark}\index{R- Tameness and partial quadrants}\label{tamecone}
From Proposition \ref{tame_equality} we conclude, if  
the M-polyfold $X$ is tame and if $x\in X$ is a smooth point, that $C_xX$ is a partial quadrant in the tangent space $T_xO$, and  we have the identity
$$
d_X(x)=\dim(T_xX/T_x^RX) = d_{C_xX}(0_x),
$$
where $0_x$ is the zero vector in $T_xX$.  
\qed
\end{remark}

We are going to show that the boundary of a tame M-polyfold has an additional structure. 
\begin{definition} \label{DEF248}\index{D- Face of an M-polyfold}\index{D- Face}
Let $X$ be a tame M-polyfold.
A {\bf face} $F$ of $X$  is the closure of a connected component in the subset $\{x\in X\ |\ d(x)=1\}$.
The  M-polyfold $X$ is called {\bf face-structured},  if every point $x\in X$ lies in exactly $d_X(x)$ many faces.
\qed
\end{definition}
Before we study faces in more detail we have a look at the local situation.
\begin{proposition}\label{PROP248} \index{P- Properties of tame $r$}
Assume that $E={\mathbb R}^k\oplus W$,  $C=[0,\infty)^k\oplus W$, $U\subset C$
is a relatively open subset  and $r:U\rightarrow U$ is a tame sc-smooth
retraction.  Suppose $U$ contains a point of the form $(0,w)$. Then the following holds.
\begin{itemize}
\item[{\em(1)}]\ The set $O=r(U)$ contains a smooth point of the form $(0,v)$.
\item[{\em(2)}]\ For a point $(0,v)\in O$  denote for $i=1,...,k$ by
$U_i$ the connected component in $U\cap E_i$ containing $(0,v)$. 
Then $r(U_i)\subset U_i$.
\item[{\em(3)}]\ For every $i=1,..,k$ the set $O\cap E_i$ contains points $y$  with $d_C(y)=1$
and these points form an open and dense subset of $O\cap E_i$.
\end{itemize}
\end{proposition}
\begin{proof}
(1) Since $d_C(0,w)=k$ we deduce that $d_C(r(0,w))=k$. This means $r(0,w)$ has the form $(0,v)$. Replacing
$w$ by a nearby smooth point we deduce that $O$ contains a smooth point of the form $(0,v)$.\par

(2)  Since $r(0,v)=(0,v)$, taking a nearby smooth $v'$ to $v$ the point $r(0,v')$ is smooth and close to $(0,v)$.
Since $d_C(r(0,v'))=d_C(0,v')=k$ it follows that $r(0,v')$ has the form $(0,w)$.
We also note that $(0,v)$ and $(0,w)$ belong both to $U_i$. Hence we may assume without loss of generality
that the given $(0,v)$ was already smooth.
By assumption $T_{(0,v)}O$ has an sc-complement contained 
in $\{0\}\oplus W$. This implies that we find smooth vectors of the form $e_1=(1,0,...,0,w_1),...,e_k=(0,...,1,w_k)$ in $T_{(0,v)}O$.
It holds
$$
Dr(0,v)(e_i)=e_i.
$$
Define $\wh{e}_i=\sum_{j\neq i}e_i$ which takes the form $(1,1,..,1,0,1,1,...1,q)$ for a smooth $q$
and $0$ occurs at the $i$-th coordinate. We note that for small $t\geq 0$ the point $(0,v)+t\wh{e}_i$ belongs to
$U\cap E_i$ and we can consider the path 
$$
t\rightarrow\Phi(t):= r((0,v)+t\wh{e}_i).
$$
We claim that for $t>0$ and small,  only the $i$-coordinate of $\Phi(t)$ vanishes.
We know that for $t>0$ small $d_C(\Phi(t))=d_C((0,v)+t\wh{e}_i)=1$. Consequently
precisely one of the first $k$ coordinates vanishes. We compute
\begin{eqnarray*}
\Phi(t)-(0,v)&=&\Phi(t)-\Phi(0)\\
&=&\int_0^t \left(Dr((0,v)+\tau \wh{e}_i)\wh{e}_i\right) d\tau \\
&=&t\cdot \wh{e}_i + \int_0^t \left((Dr((0,v)+\tau\wh{e}_i)-Dr(0,v))( \wh{e}_i) \right)d\tau\\
&+&  t \wh{e}_i  + t\delta(t),
\end{eqnarray*}
where $|\delta(t)|_0\rightarrow 0$. This implies that  the vanishing coordinate  must be the $i$-th one. 
This discussion shows that given $(0,v)$ in $O$ we can connect $(0,v)$ to a smooth $(0,v')$
via a path of the form $(0,v_t)$. Moreover, we have shown that there exists a smooth point
$(a,q)\in O$ with $a_i=0$ and $d_C(a,q)=1$ which moreover belongs to $U_i$. Given any other
point $(a',q')\in U_i$ we can connect it with $(a,q)$ by a continuous path $\gamma(t)\in U_i$
with $d_C(\gamma(t))=1$. Since $d_C(r(\gamma(t)))=1$ and the $i$-coordinate of $\gamma(0)$ vanishes
it has to vanish for all $t\in [0,1]$. This implies by a continuity argument that $r(U_i)\subset U_i$. \par

(3) Let $y\in O\cap E_i$. We can write $y=(a,w)$, where $a_i=0$. Moving $(a,w)$ slightly
to some $(b,v)$ with $v$ smooth and $b_i=0$, $b_j>0$ for $j\neq i$ we find arbitrarily close
to $(a,w)$ an element $r(b,v)$ of the form described.  Note that $d_C(r(b,v))=1$.
This shows that the desired elements are dense. The openness statement is trivial.
\qed \end{proof}

\begin{lemma}\label{ert}\index{L- Properties of tame $r$}
Let $(O,C,E)$ be a sc-smooth, tame retract and let ${\mathcal F}$ be the collection of extended faces,
associated with the partial quadrant  $C$.  Then the faces $F$  in $O$ are the connected components of the sets $f\cap O$, where 
$f\in {\mathcal F}$ are extended faces. The connected components of $f\cap O$ and $f'\cap O$,
containing a point $x \in O$, are equal if and  only if  $f=f'$.
\end{lemma}
\begin{proof}
By definition of $O$, the retract $O$ itself  is a tame M-polyfold. We denote by $U\subset  C\subset E$
the  relatively open subset in the  partial quadrant $C$, on which a tame sc-smooth retraction satisfies 
$r(U)=O$. Let $\wh{F}$ be a connected component
of the subset $\{x\in O\, \vert \, d_O(x)=1\}$. By Proposition \ref{tame_equality}, $d_O=d_C$ so that this set is  the same as a connected component
of $\{x\in O\, \vert \, d_C(x)=1\}$. If we look at  the isomorphic case $E={\mathbb R}^k\oplus W$, we see immediately, that there
exists an index  $i$, such  that $\wh{F}$
must lie in the subset of $C$, consisting of points $(a,w)$, for which $a_i=0$. We conclude that there exists
an extended face $f\in {\mathcal F}$, such that $\wh{F}\subset f$. Therefore, $\wh{F}$ is contained in a connected subset of $O\cap f$, 
and the closure of $\wh{F}$ in $O$ lies in $O\cap f$. This shows that a face of $O$ is contained in a connected component of some $O\cap f$.\par

Next we consider  a connected component $Q$ of $O\cap f$, where  $f$ is an extended face of $C$. Let $e\in Q$. By assumption, $e\in f$ and we can take a (suitable) vector $h\in f$ so small, that $e+th\in U\cap f$ for $t\in (0,1]$,
and such that  for $t\in (0,1]$ the points $e+th$ do not belong to any extended face other than $f$.
This follows from Proposition \ref{PROP248} and can be constructed explicitly, using the model $E={\mathbb R}^k\oplus W$ for $E$. Since the retraction $r$ is tame,  we have $d_{C}(r(e+th))=1$ for $t\in (0,1]$, and 
for $t\in (0,1]$ the points $r(e+th)$ belong
to the same connected component of $\{x\in O\ |\ d_O(x)=1\}$. Using  this argument repeatedly,
 we can connect any two points in $Q$ by a continuous path $\gamma\colon 
[0,1]\rightarrow Q$,
such that  the points $y\in \gamma(0,1)$ satisfy  $d_C(y)=1$. This implies that $Q$ is contained in the closure of a connected component $\wh{F}$ in $\{x\in X\ |\ d_O(x)=1\}$.\\
The discussion so far characterizes the faces of $O$ as the connected components of the sets $O\cap f$, where $f$ is some extended face.\par

The assertion that  connected components of $f\cap O$ and $f'\cap O$,
containing a point $x \in O$, are equal if and  only if  $f=f'$ is trivial. 
\qed \end{proof}

\begin{corollary}\label{cor_2.42}\index{C- Characterization of $d_O$ for tame $O$}
If  $(O,C,E)$ is a tame sc-smooth retract, then every point $x\in O$ lies in precisely $d_O(x)$ many faces.
 \end{corollary}
 \begin{proof}
 In view of Proposition \ref{tame_equality}, $d_O(x)=d_C(x)$, so that $x$ belongs to precisely $d:=d_C(x)$ extended faces $f\in {\mathcal F}(x)$, let us say, it belongs to the faces $f_1,\ldots ,f_d$ in the sc-Banach space $E=\R^k\oplus W$, defined as the subspaces $f_i=\{(a_1,\ldots, a_k, w)\vert \, a_i=0\}$, $1\leq i\leq d$. The point $x$ is represented  by $x=(0,\ldots ,0, a_{d+1}, \ldots, a_k, w)$, where 
 $a_j\neq 0$ for $d+1\leq j\leq k$. Let $r\colon 
U\to U$ be the tame  retraction onto $O=r(U)$. Then,  the paths $\gamma_i(t)$ in $O$ for $1\leq i\leq d$, starting at $\gamma_i(0)=x$ are  defined by $$
 \gamma_i(t) = r(t,\ldots ,t, a_i=0, t,\ldots, t, a_{d+1},\ldots ,a_k,w).
 $$
Because the retraction $r$ is tame we have, $d_C(\gamma_i(t))=1$ if $t>0$ and  the points $\gamma_i(t)$ belong to $f_i\cap O$, but not to any other face  $f_j\cap O$ with $j\neq i$. Hence, the connected components  of the sets $O\cap f_j$ containing $x$ for $1\leq j\leq d$ are all different. Of course, there cannot be additional components containing $x$.
 \qed \end{proof}
Corollary \ref{cor_2.42} implies immediately the following result.
\begin{proposition}\index{P- Number of local faces}
If  $X$ is  a tame M-polyfold, then every point $x\in X$ has an open neighborhood
$V$, so that $y\in V$ lies in precisely $d_X(y)$ many faces of $V$ and  $d_V(y)=d_X(y)$ for all $y\in V$.
In particular, globally, a point $x\in X$ lies in at most $d_X(x)$-many faces.
\qed
\end{proposition}

The following technical result turns out to be useful later on.  
 \begin{proposition}[Properties of faces]\label{FACE_X}\index{P- Properties of faces}
Let $(O,C,E)$ be a tame sc-smooth retract associated with the tame sc-smooth retraction $r\colon 
U\to U$ onto $O=r(U)$, and let $F$ be a face of $O$.  Then, there exists an open neighborhood $V^\ast$ of $F$ in $U$ and a sc-smooth retraction $s\colon 
V^\ast\rightarrow V^\ast$ onto  $s(V^\ast)=F$. Moreover, defining $V\subset O$
by $V=O\cap V^\ast$, the restriction $s\colon 
V\rightarrow V$ is a sc-smooth retraction onto  $s(V)=F$, so that $F$ is a sc-smooth  sub-M-polyfold of $O$.  Further, $F$ is tame, i.e. it admits a compatible sc-smooth atlas consisting of  tame local models.
In addition, 
$$
d_F(x)=d_O(x)-1\quad  \text{for all $x\in F$}.
$$
\end{proposition}
\begin{proof}
We may assume that $E={\mathbb R}^k\oplus W$ and $C=[0,\infty)^k\oplus W$. Let $U\subset C$ be a relatively open subset of the partial quadrant $C$ and let
$r\colon 
U\rightarrow U$ be a tame sc-smooth retraction onto $r(U)=O$. In view of Lemma \ref{ert}
we may assume, that our face $F$ is a connected component of $O\cap f_1$, where $f_1$ is the extended face consisting of those
points whose first coordinate vanishes. Consequently,  
 for every  $(a,w)\in F$ we find an open  neighborhood
$U_{(a,w)}\subset U$ in $C$, such  that, if $(b, v)\in U_{(a, w)}$, then $r(0,b_2,\ldots ,b_k,v)$  has the first coordinate vanishing.
Taking the union of these neighborhoods, we find an open neighborhood $U^\ast$ of $F$ which is contained in $U$,
such that for every $(b,v)\in U^\ast$,  the point $(0,b_2,\ldots ,b_k,v)$ belongs to $ U$ and $r(0,b_2,\ldots ,b_k,v)$ has the first coordinate vanishing
and belongs to $F$.
Hence, we can define the sc-smooth map
$$
s\colon 
U^\ast\rightarrow U\quad \text{by}\quad s(b,v):=r(0,b_2,\ldots,b_k,v),
$$
which has its image in $F$, so that,  in view of $F\subset U$,  we may assume that
$$
s\colon 
U^\ast\rightarrow U^\ast.
$$
It follows that the face $F$ is a sc-smooth  sub-M-polyfold of $O$. 

We know that $F$ is a connected component of $f_1\cap O$.  We find an open neighborhood $V$
of $O$ in $f_1\cap C\subset f_1$, so that $r(V)=F$. Clearly,  $s=r\vert V\colon 
V\rightarrow V$ preserves $d_{C\cap f_1}$, since it preserves $d_C$ and $d_C=d_{C\cap f_1}+1$ on $C\cap f_1$.

If $x\in F$ is a smooth point, then $T_xF\subset T_xO$ is of codimension $1$ and hence has a complement
$P$ of dimension $1$, so that 
$$
T_xO=T_xF\oplus P.
$$
We know from Definition \ref{tame_retarctions}  of a tame retract that $T_xO$ has a sc-complement $A\subset E_x$ in $E$. Hence, 
$$
E=T_xO\oplus A = T_xF\oplus P\oplus A.
$$

In view of $P\cap f_1 = \{0\}$ and $T_xF \cap f_1 = T_xF$, we deduce
$$
f_1 = E \cap f_1 = (T_xF\oplus P \oplus A) \cap f_1 = T_xF \oplus (f_1 \cap A).
$$
Noting that $(f_1 \cap A) \subset (f_1 \cap E_x) = f_1$, the sc-smooth retraction $s \colon 
 V \rightarrow V$ onto $s(V) = F$ is tame. Since the last statement of Proposition \ref{FACE_X} is obvious, the  proof  of Proposition \ref{FACE_X} is complete.
\qed \end{proof}
In order to formulate another result along the same lines we first recall  that, given any point $x$ in a tame 
M-polyfold $X$,  we find an open neighborhood $U=U(x)$ such  that $U\cap X$ has precisely $d_X(x)$ many faces containing $x$. 

The local discussion can be summarized as follows.
\begin{proposition}\label{FACE_XXXX}
Given a tame M-polyfold $X$ every point $x$ has $d_X(x)$ many local faces containing $x$. The local faces are tame M-polyfolds
with $d_F(x)=d_X(x)-1$.  The intersection
of local faces $F_\sigma$ associated to a subset $\sigma$ of ${\mathcal F}_x$ is a tame M-polyfold of codimension $\#\sigma$
and $d_{F_\sigma}(x)=d_X(x)-\#\sigma$.
\qed
\end{proposition}

Finally we  sum up the discussion about faces in the following theorem.

\begin{theorem}\label{FACE_XX}\index{T- Faces as sub-M-polyfolds}
The interior of a face $F$ in a tame M-polyfold $X$ is a sc-smooth sub-M-polyfold of $X$.  If $X$ is face-structured, then every face $F$ is a sub-M-polyfold  and the induced M-polyfold structure is tame.
The inclusion map $i\colon F\rightarrow X$ satisfies $d_X(i(x)) = d_F(x)+1$. 
Every point $x\in X$ has an open neighborhood $U$ such that every point $y\in U$
lies in precisely $d_X(y)$-many faces of $U$. If $X$ is face-structured, then every point $x \in X$ lies in precisely 
$d_X(x)$-many global faces of $X$.
\qed
\end{theorem}

Later on we shall frequently use the following proposition about degeneracy index in fibered products.

\begin{proposition}\label{fibered-x}\index{P- Fibered products}
Let $f\colon  X\rightarrow Z$ be a local sc-diffeomorphism  and $g\colon 
Y\rightarrow Z$ a sc-smooth map.
Then the fibered product  $X{_{f}\times_g}Y$,  defined by
$$
X{_{f}\times_g}Y=\{(x,y)\in X\times Y \vert \, f(x)=g(y)\}, 
$$
is a sub M-polyfold of $X\times Y$. If $Y$ is tame, then also $X{_{f}\times_g}Y$ is tame and 
$$
d_{X{_{f}\times_g}Y}(x,y)=d_Y(y).
$$
If both $X$ and $Y$ are tame and both $f$ and $g$ are local sc-diffeomorphism, then 
$$
d_{X{_{f}\times_g}Y}(x,y)=\frac{1}{2} d_{X\times Y}(x,y)=\frac{1}{2}[d_X(x)+d_Y(y)] = d_X(x)=d_Y(y).
$$
\end{proposition}
\begin{proof}
Clearly the product $X\times Y$ is an M-polyfold,  and if both $X$ and $Y$ are tame, then  also $X\times Y$ is tame. 
To see that $X{_{f}\times_g}Y$ is a sub-M-polyfold of $X\times Y$, we take a point $(x,y)\in X{_{f}\times_g}Y$ and  fix  open neighborhoods $U$ of $x$ in $X$ and $W$ of  $f(x)$ in $Z$ so that $f\vert U\colon 
 U\rightarrow W$ is a sc-diffeomorphism. Next we choose an open neighborhood $V$ of $y$ in $Z$ such that $g(V)\subset W$. We define the map 
$R\colon 
U\times V\rightarrow U\times V$ by
$$
R(x',y') = ((f\vert U)^{-1}\circ g (y'),y'), 
$$
It is sc-smooth and satisfies $R\circ R =R$ and $R(U\times V) = (U\times V)\cap (X{_{f}\times_g}Y).$
Consequently, $X{_{f}\times_g}Y$ is a sub M-polyfold of $X\times Y$.   Next we assume, in addition,  that $Y$ is tame.  With the  point $(x,y)\in X{_{f}\times_g}Y$ and open sets $U$ and $V$ as above, we assume that $(V, \varphi', (O', C', E'))$  is a chart around $x$ such that $ (O', C', E')$ is a tame sc-smooth retract.
Then we define a map  $\Phi\colon 
(U\times V)\cap X{_{f}\times_g}Y\to O'$  by setting 
$$\Phi (x', y')=\varphi' (y').$$
The map $\Phi$ is injective since given a pair $(x', y')\in  X{_{f}\times_g}Y$, we have $x'=(f\vert U)^{-1}\circ g (y')$ and hence  $\Phi$ is homeomorphism onto $O'$.
So, the tuple 
$$
((U\times V)\cap X{_{f}\times_g}Y, \Phi, (O', C', E'))
$$
 is a chart on $X{_{f}\times_g}Y$ and any two such charts are sc-smooth compatible.  Consequently, 
 $ X{_{f}\times_g}Y$ is tame. 
To prove the formula for the degeneracy index, we 
consider the projection $\pi_2\colon 
 X{_{f}\times_g}Y\to X$, $(x, y)\to y$, onto the second  component.   In view of the above discussion, 
$\pi_2$ is a local sc-diffeomorphism and consequently preserves the degeneracy index.  
Hence  $d_{X{_{f}\times_g}Y}(x,y)=d_Y(y)$. 

If both $f$ and $g$ are sc-diffeomorphisms and $X$ and $Y$ are both tame, then it follows from the previous case that  $d_{X{_{f}\times_g}Y}(x,y)=d_X(x)$, so that 
$$d_{X{_{f}\times_g}Y}(x,y)=d_X(x)=d_Y(y).$$
Finally, since $X\times Y$ is tame, then  $d_{X\times Y}(x,y)=d_X(x)+d_Y(y)$, so that 
$$d_{X\times Y}(x,y)=d_X(x)+d_Y(y)=2d_X(x)=2d_Y(y)=2d_{X{_{f}\times_g}Y}(x,y),$$
for every $(x, y)\in X{_{f}\times_g}Y$.  The proof of Proposition \ref{fibered-x} is complete.
\qed \end{proof} 

\section{Strong Bundles}\label{section2.5_sb}
As a preparation for the study of sc-Fredholm sections in the next chapter we shall introduce in this section the notion of a strong bundle over a M-polyfold. 

As usual we shall first describe the new notion in local charts 
of a strong bundle and consider $U\subset C\subset E$, where $U$ is a relatively open subset of the partial quadrant $C$ in the sc-Banach space $E$. If $F$ is another sc-Banach space, we define the non- symmetric product
 $$
U\triangleleft F\index{$U\triangleleft F$}
$$
 as follows. As a set, $U\triangleleft F$ is the product $U\times F$, but it possesses  a double filtration
$(m,k)$ for $0\leq k\leq m+1$, defined by
$$
(U\triangleleft F)_{m,k}:=U_m\oplus F_k.
$$
We view $U\triangleleft F\rightarrow U$ as a trivial bundle. We define for $i=0,1$ 
the sc-manifolds $(U\triangleleft F)[i]$\index{$(U\triangleleft F)[i]$} by their filtrations
$$
((U\triangleleft F)[i])_{m}:= U_{m}\oplus F_{m+i},\quad m\geq 0.
$$
\begin{definition}\index{D- Strong bundle map}
\noindent A {\bf strong bundle map}
$\Phi\colon  U\triangleleft F\rightarrow U'\triangleleft F'$
is a map which preserves the double filtration and is of the form
$\Phi(u,h)=(\varphi(u),\Gamma(u,h))$,
where the map $\Gamma$ is linear in $h$. In addition,  for $i=0,1$,  the maps
$$
\Phi\colon 
(U\triangleleft F)[i]\rightarrow (U'\triangleleft F')[i]
$$
are  sc-smooth.\par

\noindent  A {\bf  strong bundle isomorphism}\index{D- Strong bundle isomorhism} is an  invertible strong bundle map whose inverse is also a strong bundle map.\par

\noindent  A {\bf  strong bundle retraction}  \index{strong bundle retraction} is a strong  bundle map
$R\colon 
U\triangleleft F\rightarrow U\triangleleft F$
satisfying, in addition,  $R\circ R=R$.  The map  $R$ has the form
$$
R(u,h)=(r(u),\Gamma(u,h)),
$$
where $r\colon 
U\rightarrow U$ is a sc-smooth retraction.  \par

\noindent A {\bf tame strong bundle retraction}\index{D- Tame strong bundle retraction}  is one for which the retraction $r$ is tame.  
\qed
\end{definition}
The condition $R\circ R=R$ of the retraction $R$ requires that 
$$\bigl( r(r(u)), \Gamma (r(u), \Gamma (u, h))\bigr)=
\bigl( r(u), \Gamma (u, h)\bigr).$$
Hence, if $r(u)=u$, then $\Gamma \bigl(u, \Gamma (u, h)\bigr)=\Gamma (u, h)$, and the bounded linear operator $h\mapsto \Gamma (u, h)\colon  F\to F$ is a projection. If $r(u)=u$ is, in addition, a smooth point, then the projection is a sc-operator.
\par

We continue to denote by $U\subset C\subset E$ a relatively open subset of the partial quadrant $C$ in the sc-Banach space $E$ and let $F$ be another sc-Banach space.
\begin{definition}\label{def_loc_strong_b_retract}\index{D- Strong local bundle}
A {\bf local  strong bundle retract}, denoted by 
$$
(K,C\triangleleft F, E\triangleleft F),\index{$(K,C\triangleleft F, E\triangleleft F)$}
$$
consists of a subset $K\subset C\triangleleft F$, which is the image,  
$K=R(U\triangleleft F)$, of a strong bundle retraction
$R\colon 
U\triangleleft F\rightarrow U\triangleleft F$
of the form $R(u, h)=(r(u),\Gamma (u, h)).$
Here, $r\colon 
U\to U$ is a sc-smooth retraction onto $r(U)=O$.
The local strong bundle retract $(K,C\triangleleft F, E\triangleleft F)$ will sometimes be abbreviated by 
$$
p\colon 
K\to O\index{$p\colon K\to O$}
$$
where $p$ is the map induced by the projection onto the first factor. If the strong bundle retraction $R$  is tame, the  local strong bundle retract is called a {\bf tame local strong bundle retract}. \index{D- Tame local strong bundle retract}
\qed
\end{definition}

The retract  $K$ inherits the {\bf double filtration}\index{Double filtration}  $K_{m,k}$\index{$K_{m,k}$} for $m\geq 0$ and $0\leq k\leq m+1$, defined by 
\begin{equation*}
\begin{split}
K_{m,k}&=K\cap (U_m \oplus F_k)\\
&=\{(u, h)\in U_m \oplus F_k\vert \, R(u, h)=(u, h)\}\\
&=\{(u, h)\in O_m \oplus F_k\vert \, \Gamma (u, h)=h\}.
\end{split}
\end{equation*}
The associated spaces $K[i]$, defined by 
$$
K[i]=K_{0,i},\quad i=0,1,\index{$K[i]$}
$$
are equipped with the filtrations
$$
K[i]_m=K_{m,m+i}
\quad \text{for all $m\geq 0$.}$$
The projection maps 
$p\colon  K[i]\to O$
are sc-smooth for $i=0,1$. In fact, $p:K[i]\rightarrow O$,  for $i=0,1$, is a  sc-smooth bundle.

\begin{definition}\label{def_section_loc_strong_bundle}\index{D- Section}
A {\bf section} of the local strong bundle retract $p\colon 
K\to O$ is a map $f\colon 
O\to K$ satisfying $p\circ f=\mathbbm{1}_O$.
The section $f$ is called {\bf sc-smooth}\index{Sc-smooth section}, if $f$ is a section of the bundle 
$$p\colon 
K(0)\to O.$$
The section $f$ is called  {\bf $\ssc^{\pmb{+}}$-smooth}, \index{D- Sc$^+$-section} if $f$ is a sc-smooth section of the bundle 
$$p\colon 
K(1)\to O.$$
\qed
\end{definition}
A section $f\colon 
O\to K$ is of the form $f(x)=(x, {\bf f}(x))\in O\times F$ and the map ${\bf f}\colon 
O\to F$ is called {\bf principal part of the section}\index{Principal part of a section}. We shall usually denote the principal part with the same letter as the section.

At this point, we can introduce the {\bf category $\mathcal{SBR}$}.\index{$\mathcal{SBR}$} Its  objects are the local,  strong bundle retractions $(K,C\triangleleft F,E\triangleleft F)$. The morphisms of the category   are maps
$\Phi\colon 
K\rightarrow K'$ between local strong bundle retracts, which are linear in the fibers and preserve the double filtrations. Moreover,  the induced maps 
$\Phi[i]\colon 
K[i]\rightarrow K'[i]$ are sc-smooth for $i=0,1$.
We recall that, by definition, {\bf the map $\Phi[i]\colon 
K[i]\rightarrow K'[i]$ between retracts is sc-smooth}, if the composition with the retraction,
$$(\Phi \circ R) \colon 
(U\triangleleft F)[i]\to (E'\triangleleft F')[i],$$
is sc-smooth.
There are  two forgetful functors 
$$
\text{forget}[i]\colon 
\mathcal{SBR}\rightarrow \mathcal{BR},
$$
into the category of  $\mathcal{BR}$ of sc-smooth bundle retractions introduced in Section \ref{section2.1}. They are defined by 
\begin{equation*}
\text{forget}[i](K)=K[i]
\end{equation*}
on the objects $K$ of the category, and 

\begin{equation*}
\text{forget}[i](\Phi)=\Phi[i]
\end{equation*}
on the morphisms between the objects.

We are in a position to introduce the notion of a strong bundle over a M-polyfold $X$. 
We consider a continuous surjective map
$$P\colon 
Y\to X$$
from the paracompact Hausdorff space $Y$ onto the M-polyfold $X$.  We assume, for every $x\in X$,  that  the fiber $P^{-1}(x)=Y_x$
has the structure of a Banachable  space.

\begin{definition}\label{def_strong_bundle_chart}\index{D- Strong bundle chart}
A {\bf strong bundle chart} for the bundle $P\colon 
Y\to X$ is a tuple
$$
(\Phi, P^{-1}(V), (K, C\triangleleft F,E\triangleleft F)),
$$
in which $(K, C\triangleleft F,E\triangleleft F)$ is a strong bundle retract, say $p:K\rightarrow O$. In addition, $V\subset X$ is an open subset of $X$, homeomorphic to the retract $O$ by a  homeomorphism $\varphi\colon 
V\to O$. In addition, 
$$\Phi\colon 
P^{-1}(V)\to K$$
is a homeomorphism from $P^{-1}(V)\subset Y$ onto the retract $K$, covering the homeomorphism $\varphi\colon 
V\to O$, so that the diagram 
\begin{equation*}
\begin{CD}
P^{-1}(V)@>\Phi>>K\\
@VPVV @VVpV \\
V@>\varphi>>O
\end{CD}
\end{equation*}
commutes. The map $\Phi$ has the property that,  in the fibers over $x\in V$,  the map $\Phi\colon 
P^{-1}(x)\to p^{-1}(\varphi (x))$ is a 
bounded linear isomorphism between Banach spaces.\par

\noindent {\bf Two strong bundle charts} $\Phi\colon 
P^{-1}(V)\to K$ and 
$\Phi'\colon 
P^{-1}(V')\to K'$,  satisfying $V\cap V'\neq \emptyset$ are 
{\bf compatible},\index{Compatibility of strong bundle charts} if the transition maps
$$
\Phi'\circ \Phi^{-1}[i]\colon 
\Phi (P^{-1}(V\cap V'))[i]\to 
\Phi' (P^{-1}(V\cap V'))[i]
$$
are sc-smooth diffeomorphisms for $i=0,1.$
\qed
\end{definition}
As usual, one now proceeds to define a strong  bundle atlas, consisting of compatible strong bundle charts covering $Y$ and the equivalence between two such atlases.
\begin{definition}\label{def_strong_bundle} \index{D- Strong bundle}
The continuous surjection 
$$
P\colon 
Y\to X
$$
from the paracompact Hausdorff space $Y$ onto the M-Polyfold $X$, equipped with an equivalence class of strong bundle atlases is called a {\bf strong bundle over the M-polyfold $X$}.
\qed
\end{definition}
Induced by the strong bundle charts, the M-polyfold $Y$ is equipped with a natural double filtration into subsets $Y_{m,k}$, $m\geq 0$ and $0\leq k\leq m+1$. Therefore, we can distinguish the underlying M-Polyfolds $Y[i]$ for $i = 0,1$ with the filtrations
$$
Y[i]_m = Y_{m,m+i}  ,  m \geq 0.
$$
The projections
$$
P[i] \colon 
 Y[i] \rightarrow X
$$
are sc-smooth maps between M-Polyfolds.
Correspondingly, we distinguish two types of sections of the strong bundle $P \colon 
 Y \rightarrow X$.

\begin{definition}\label{def_sc_inft_sections}
\noindent A section of the strong bundle $P \colon 
 Y \rightarrow X$ is a map $f \colon 
 X \rightarrow Y,$ satisfying $P\circ f = {\mathbbm 1}_X.$\\
 The section $f$ is called a {\bf sc-smooth section}\index{D- Sc-smooth section}, if $f$ is an sc-smooth section of the bundle
$$
P[0] \colon 
 Y[0] \rightarrow X.
$$
The section $f$ is called a ${\bf sc^+}$-{\bf section} of $P \colon 
 Y \rightarrow X$, \index{D- Sc$^+$-section} if $f$ is an sc-smooth section of the bundle
$$ 
P[1] \colon 
 Y[1] \rightarrow X.
$$
\qed
\end{definition}
If we say $f$ is an sc-section of $P:Y\rightarrow X$ we mean that it is an sc-smooth section of $Y[0]\rightarrow X$.
An sc$^+$-section or sc$^+$-smooth section of $P:Y\rightarrow X$ is an sc-smooth section of $Y[1]\rightarrow X$.

\begin{definition}[{\bf Pull-back bundle}]
Let $P\colon Y\to X$ be a strong bundle over the M-polyfold $X$ and let  $f\colon Z\to X$  be an sc-smooth map  from the M-polyfold $Z$ into $X$. The pull-back bundle of $f$,
$$P_f\colon f^\ast Y\to Z,$$ 
is defined by the set 
$f^\ast Y=\{(z, y)\in Z\times Y\, \vert \, P(y)=f(z)\}$
and the projection $P_f (z, y)=z$, so that with the projection 
$P'\colon f^\ast Y\to Y$, defined by $P'(z, y)=y$, the diagram 
\begin{equation*}
\begin{CD}
f^\ast Y@>P'>>Y\\
@VP_fVV @VVPV \\
Z@>f>>X
\end{CD}
\end{equation*}
commutes.
\end{definition}

As already shown in \cite{HWZ2}, Proposition 4.11, the strong bundle  structure of $P:Y\rightarrow X$ induces a natural strong bundle structure of the pull-back bundle $P_f$.
\begin{proposition}\label{pull_back_strong_bundle}\index{P- Pull-back of strong bundles}
The pull-back bundle $P_f\colon f^\ast Y\to Z$ 
carries a natural structure of a strong M-polyfold bundle over the  M-polyfold $Z$.
\qed
\end{proposition}
The easy proof is left to the reader.  
\section{Appendix}

\subsection{Proof of Proposition  \ref{smooth_retract}}\label{A2.1}
Let us recall the statement.\par

\noindent {\bf Proposition.} Let $r\colon 
U\to U$ be a $C^\infty$-retraction defined on an open subset $U$ of a Banach space $E$. Then $O:=r(U)$ is a $C^\infty$-submanifold of $E$. More precisely, for  every point $x\in O$,  there exist an open neighborhood  $V$ of $x$ and  an open neighborhood $W$  of $0$ in $E$, a splitting $E =R \oplus N$, and a smooth diffeomorphism 
$\psi \colon 
 V\to  W$  satisfying  $\psi (0)=x$ and 
$$\psi (O\cap V ) = R \cap W.$$
\begin{proof} 
The $C^\infty$-retraction $r$ satisfies $r\circ r=r$ and hence,  by the chain rule,  
$$Dr(r(x))\circ Dr(x)=Dr (x)\quad \text{for every $x\in U$}.$$
Therefore, if $r(x)=x$ 
 \begin{equation}\label{c_infty_retract_eq0}
 Dr (x)\circ Dr (x)=Dr(x)
  \end{equation}
and hence the bounded linear operator  $Dr(x)\in \mathscr{L}(E, E)$ is a projection at every point $x\in O$.

Now we take a point $x\in O=r(U)$ and,  for simplicity,  assume that $x=0$.  In view of \eqref{c_infty_retract_eq0}, the Banach space $E$ splits into $E=R\oplus N$, where 
\begin{align*}
R&=\text{range}\ Dr(0)=\text{ker}\ (\mathbbm{1}-Dr (0))\\
N&=\text{ker}\ Dr(0)=\text{range}\ (\mathbbm{1}-Dr (0)).
\end{align*}
According to the splitting $E=R\oplus N$, we use the equivalent norm $\abs{(a, b)}=\max\{\norm{a}, \norm{b}\}$. By $B(\varepsilon)$ we denote an open ball of radius $\varepsilon$ (with respect to the norm $\abs{\cdot }$) centered at the origin.

Now we introduce  the map $f\colon 
 (R\oplus  N)\cap U\to E$, defined   by 
$$f(a, b)=r(a)+(\mathbbm{1}-Dr(a))b.$$
At $(a, b)=(0, 0)$, we have $f(0, 0)=r(0)=0$ and 
$$Df(0, 0)[h, k]=h+k$$
for all $(h, k)\in R\oplus N.$ Consequently, in view of the inverse function theorem, $f$ is a local $C^\infty$-diffeomorphism, and we assume without loss of generality that $f$ is a diffeomorphism on $U$. 

We claim that there exists a positive number $\delta$,  such that, if 
 $f(a, b)\in O=r(U)$ for $(a, b)\in B(\delta)$, then $b=0$.   The proposition then follows by setting  $W=B(\delta)$, $V=f(W)$,  and defining  the map $\psi\colon 
W\to V$  by $\psi=(f\vert W)^{-1}$.

It remains to  prove the  claim that $b=0$.  Since  $r$ is smooth, we find a constant $\varepsilon>0$, such that  
 $B(2\varepsilon)\subset U$ and 
 \begin{equation}\label{c_infty_retract_eq1}
 \abs{Dr(x)-Dr(0)}<\frac{1}{3}\quad \text{for all $\abs{x}<\varepsilon$}.
 \end{equation}
Moreover,  since $r(0)=0$, there exists a constant $0<\delta <\frac{3}{4}\varepsilon$, such that 
  \begin{equation}\label{c_infty_retract_eq2}
 \abs{r(x)}<\varepsilon \quad \text{for all $\abs{x}<\delta$}.
 \end{equation}
If $x\in B(\delta)$, then for $ \abs{h}<\varepsilon$, 
 \begin{equation}\label{c_infty_retract_eq3}
  r(r(x)+h)=r(r(x))+Dr(r(x))h+o(h)=r(x)+Dr(r(x))h+o(h),
 \end{equation}
where $o(h)$  is an $E$-valued map satisfying $\frac{o(h)}{\abs{h}}\to 0$, as $\abs{h}\to 0$.
Taking a smaller $\varepsilon$, we may assume that 
 \begin{equation}\label{c_infty_retract_eq3a}
 \abs{o(h)}\leq \frac{\abs{h}}{2}\quad  \text{ for all $\abs{h}<\varepsilon. $}
 \end{equation}
Let  $x=(a, b)\in B(\delta)$  and $f(a, b)\in O$. This means that  
 \begin{equation}\label{c_infty_retract_eq4}
r\bigl(r(a)+(\mathbbm{1}-Dr(a))b\bigr)=r(a)+(\mathbbm{1}-Dr(a))b.
 \end{equation}
By \eqref{c_infty_retract_eq1} and the fact that, $Dr(0)b=0$ for  $b\in N$, the norm of  $(\mathbbm{1}-Dr(a))b$ for $\abs{a}<\delta$ can be estimated as 
\begin{equation}\label{c_infty_retract_eq5}
\abs{b-Dr(a))b}=\abs{b-\bigl(Dr(a))-Dr(0)\bigr)b} \leq \frac{4}{3}\abs{b}.
\end{equation}
Inserting $x=a$ and $h=(\mathbbm{1}-Dr(a))b$ into  \eqref{c_infty_retract_eq3} and using  the identity  
\eqref{c_infty_retract_eq4}, we obtain
\begin{equation}\label{c_infinity_retraction_eq7} 
Dr(r(a))h=h+o(h).
\end{equation}
The left-hand side is equal to 
\begin{equation*}\label{c_infinity_retraction_eq8} 
Dr(r(a))h=Dr(r(a))b-Dr(r(a))Dr(a)b=Dr(r(a))b-Dr(a)b,
\end{equation*}
and the right-hand side is equal to $h+o(h)=b-Dr(a))b+o(h)$,
so that 
\eqref{c_infinity_retraction_eq7}  takes the form 
\begin{equation}\label{c_infty_retract_eq10}
Dr(r(a))b=b+o(h).
\end{equation} 
Arguing by contradiction, we assume  that $b\neq 0$. Since   $\abs{r(a)}< \varepsilon$,  the norm of the left-hand side is, in view of  \eqref{c_infty_retract_eq1},  bounded from above by 
\begin{equation}\label{c_infty_retract_eq11}
\abs{Dr(r(a))b}=\abs{Dr(r(a))b-Dr(0)b}<\frac{1}{3}\abs{b}.
\end{equation} 
On the other hand, using \eqref{c_infty_retract_eq3a} and \eqref{c_infty_retract_eq5}, we estimate  the norm of the right-hand side  of \eqref{c_infty_retract_eq10} from below as  
\begin{equation}\label{c_infty_retract_eq12}
\abs{b+o(h)}\geq \abs{b}-\frac{1}{2}\abs{h} \geq   \abs{b}-\frac{2}{3}\abs{b}=\frac{1}{3}\abs{b}. 
\end{equation} 
Consequently, $\frac{1}{3}\abs{b}>\frac{1}{3}\abs{b}$ 
which is absurd. Therefore,  $b=0$ and the proof is complete.
\qed \end{proof}

\subsection{Proof of Theorem \ref{X_m_paracompact}}\label{A2.2}
We recall the statement of the result which we would like to prove.\par

\noindent{\bf Theorem.} 
Let $X$ be a M-polyfold. For every $m\geq 0$, the space $X_m$ is metrizable and, in particular, paracompact. In addition, the space $X_\infty$ is metrizable.
\par

We  make use of the following theorem of metrizability due to Yu. M. Smirnov \cite{Smirnov}.
\begin{theorem}[{\bf {Yu. M. Smirnov}}]\label{Smirnov}\index{T- Smirnov's metrizability theorem}
Let $X$ be a space that is paracompact Hausdorff, and assume that  every point has an open  neighborhood in $X$ that is metrizable. Then $X$ is metrizable.\end{theorem}

 Recall that a regular topological space is one where a closed subset $A$ and a point $p\not\in A$ can be separated by suitable open neighborhoods $U$ and $V$  of $A$ and $p$, respectively, i.e. $U\cap V=\emptyset$. 
Moreover, a regular Hausdorff space is what is called sometimes a $T_3$-space.  We also note that in \cite{Michael} the definition 
of paracompactness includes the property of being Hausdorff. For us paracompactness just means that for every open covering there exists 
a locally finite subcovering.  With these conventions the following quoted result is an obvious reformulation.
\begin{lemma}[\cite{Michael}, Lemma 1] \label{equivalent_definitions_paracompact}\index{L- Characterization of paracompactness}
A  regular Hausdorff space $X$ is paracompact, if and only if every open cover of $X$ has a locally finite refinement consisting of closed sets.
\qed
\end{lemma}

We shall use the above lemma in the proof of the following result.

\begin{proposition}\label{union_of_paracompact}\index{P- Paracompactness}
Let $Y$ be  a regular Hausdorff space, and let $(Y_i)_{i\in I}$ be a locally  finite family of closed subspaces  of $Y$, so that  $Y=\bigcup_{i\in I}Y_i$. If every subspace $Y_i$ is paracompact, then 
$Y$ is  paracompact.
\end{proposition}
\begin{proof}
Given an open cover  $\mathscr{U}=(U_j)_{j\in J}$ of $Y$,  it is, in view of 
Lemma \ref{equivalent_definitions_paracompact}, enough to show, that there exists a closed locally finite refinement 
$\mathscr{C}=(C_j)_{j\in J}$ of $\mathscr{U}$.  In order to prove this, we  consider for every $i$  the cover  $\mathscr{Y}_i=(Y_i\cap U_j)_{j\in J}$ of $Y_i$, consisting of open sets in $Y_i$.  Since $Y_i$ is paracompact, Lemma \ref{equivalent_definitions_paracompact} implies that there exists a closed locally finite refinement  $(C_j^i)_{j\in J}$ of  $\mathscr{Y}_i$.  The sets $C_j^i$ are closed in $Y$ and by  definition of refinement satisfy 
\begin{equation}\label{eq1}
C_j^i\subset Y_i\cap U_j\quad \text{and}\quad \bigcup_{j\in J}C^i_j=Y_i.
\end{equation}  
Now, for every $j\in J$,  we define the set
$$C_j:=\bigcup_{i\in I}C_j^i$$
and claim that the family $\mathscr{C}=(C_j)_{j\in J}$ is a closed locally finite refinement of $\mathscr{U}$.  We start with  showing that  $C_j$ is closed. Let 
 $x\in X\setminus C_j=\bigcap_{i\in I}(X\setminus C^i_j)$.   Since $(Y_i)_{i\in I}$ is locally finite,  we find  an open neighborhood  $V(x)$ of $x$ in $Y$, which  intersects $Y_i$ for at most finitely many indices $i$, say for $i$ belonging to the  finite subset $I'\subset I$. 
If $i\not \in I'$, then, by \eqref{eq1},  $V(x)\subset  Y\setminus Y_i\subset Y\setminus C_j^i$,  so that $V(x)\subset \bigcap_{i\not \in  I'}(Y\setminus Y_i).$  
Then the set $U(x)=V(x)\cap  \bigcap_{l\in I'}(Y\setminus C_j^i)$ is an open neighborhood of $x$ in $Y$ contained in $Y\setminus C_j$. This shows that  $C_j$ is closed as claimed.
That  $\mathscr{C}$ is a cover of $X$ and a refinement of $\mathscr{U}$ follows from  \eqref{eq1}, 
\begin{gather*}
\bigcup_{j\in J}C_j=\bigcup_{j\in J}\bigcup_{i\in I}C_j^i=\bigcup_{i\in I}\bigcup_{j\in J}C_j^i=\bigcup_{i\in I}Y_i=X\\
\intertext{and}
V_j=\bigcup_{i\in I}C_j^i\subset \bigcup_{i\in I}Y_i\cap U_j=U_j. 
\end{gather*}
It remains to show that $\mathscr{C}$ is locally finite.  Local finiteness of $(Y_i)_{i\in I}$ implies that a given point $x\in Y$ belongs to finitely many $Y_i$'s, say $x\in Y_i$, if and only if $i$ belongs to the finite subset  $I'\subset I$.  Moreover, there exists an open neighborhood $V(x)$ of $x$ in $Y$ intersecting at most finitely many $Y_i$'s. 
 Since $Y$ is regular, we can replace  $V(x)$  by  a smaller  open set intersecting only those $Y_i$'s whose indices $i$ belong to $ I'$.  Also,  $(C_j^i)_{j\in J}$ is a locally finite family in $Y_i$, so that  there exists an open neighborhood $V_i(x)$ of $x$ in $Y_i$ intersecting $C_j^i$  for, at most, finitely many indices $j$,  which belong to the finite subset $J_i\subset J$. Each $V_i(x)$ is of the form 
$V_i(x)=U_i(x)\cap Y_i$ for some open neighborhood $U_i(x)$ of $x$ in $Y$. Then $U(x)=V(x)\cap \bigcap_{i\in I'}U_i(x)$ is an open neighborhood of $x$ in $Y$  intersecting 
$Y_i$ only for $i\in I'$ and the set $U(x)\cap Y_i$ intersects  $C_j^i$ only for  $j\in J_i$.   This implies that  $U(x)$ can have a nonempty intersection with $C_j$ only for $j\in \bigcup_{i\in I'}J_i$.  Hence,  the family $\mathscr{C}=(C_j)_{j\in J}$ is locally finite,  and the proof is complete.  

\qed \end{proof}
Now we are ready to prove the theorem.
\begin{proof}[Theorem \ref{X_m_paracompact}]
We start with $m=0$ and take an atlas 
$$
(U^j, \varphi^j, (O^j, C^j, E^j))_{j\in J}.
$$  
Since $\varphi^j\colon 
U^j\to O^j$
 is a homeomorphism and $O^j$ is a metric space, $U^j$ is metrizable. 
Hence, $X$ is locally metrizable. Since by assumption, $X$ is Hausdorff and paracompact,  the Smirnov metrizability theorem implies that $X$ is metrizable. 

In order to prove that $X_m$ is metrizable for $m\geq 1$,  we fix $m\geq 1$  and consider the topological  space $X_m$.  Since by assumption, $X$ is Hausdorff,  given two distinct points $x$ and $x'$ in $X_m$, there exist  two disjoint open neighborhoods $U$ and $U'$ of $x$ and $x'$ in $X$.  Hence, the sets $U_m$ and $U'_m$ are disjoint open neighborhoods of $x$ and $x'$ in $X_m$, so that also $X_m$ is a  Hausdorff space. 
Moreover, the maps $\varphi\colon 
U_m\to O^j_m$ are homeomorphisms and since $O^j_m$ is a metric space, $U^j_m$ are metrizable. So,  $X_m$ is locally metrizable. To prove that $X_m$ is metrizable, it suffices to show, that, in view of the Smirnov metrizability theorem, $X_m$ is paracompact.  Using  the paracompactness  of $X$ and 
Lemma \ref{equivalent_definitions_paracompact}, we find a closed locally finite refinement $\mathscr{C}=(C_j)_{j\in J}$. In particular, $C_j\subset U_j$, so that $C^j_m\subset U^j_m$. Since $U^j_m$ is metrizable, it is also paracompact, so that $C^j_m$ is paracompact as a closed subspace of $U^j_m$.   Hence, $\mathscr{C}_m=(C^j_m)_{j\in J}$ is a locally  finite family of closed subsets of $X_m$ and each $C^j_m$ is paracompact. Thus,  by Proposition \ref{union_of_paracompact}, the space $X_m$ is paracompact, and since it is Hausdorff and locally metrizable, it is metrizable. 
Finally, choosing a metric $d_m$ defining the topology on $X_m$ , we set
$$
d(x,y)=\sum_{m=0}^{\infty} \frac{1}{2^m}\cdot  \frac{d_m(x,y)}{1+d_m(x,y)}.
$$
The metric $d$ defines the topology on $X_\infty$.
\qed \end{proof}

\subsection{Proof of Proposition \ref{op}}\label{A2.11}
We recall the statement of the proposition.\par

\noindent{\bf Proposition.}
The following holds.
\begin{itemize}
\item[(1)]\ The collection $\mathscr{B}$ defines a basis for a Hausdorff  topology on $TX$.  
\item[(2)]\ The projection $p\colon 
TX\to X^1$ is a continuous and an open map.
\item[(3)]\ With the topology defined by $\mathscr{B}$,  the tangent space $TX$ of the M-polyfold $X$  is metrizable and hence, in particular,  paracompact.
\end{itemize}
\begin{proof}
(1)\, In order to prove that $\mathscr{B}$ defines a basis for a topology on $TX$, we take two  sets   $\wt{W}_1$ and $\wt{W}_2$ in  $\mathscr{B}$ and assume that $\alpha=[(x, \varphi, V, (O, C, E), h)]\in \wt{W}_1\cap \wt{W}_2$. We claim that  there exists a set $\wt{W}\in \mathscr{B}$, satisfying  $\alpha\in \wt{W}\subset \wt{W}_1\cap \wt{W}_2.$ By definition,  $\wt{W}_i=(T\varphi_i)^{-1}(W_i)$, where $(\varphi_i, V_i, (O_i, C_i, E_i))$ is a chart and $W_i$ is an open subset of $TO_i$, containing 
$$
\alpha=[(x_i, \varphi_i, V_i, (O_i, C_i, E_i), h_i)]\quad \text{for $i=1,2$}.
$$
This means that 
\begin{equation}\label{basis_top_1}
x_i=x\quad \text{and}\quad h_i=T(\varphi_i \circ \varphi^{-1})(\varphi (x))h\quad \text{for $i=1,2$},
\end{equation}
and, moreover, $(\varphi_i (x_i), h_i)=(\varphi_i (x), h_i)\in W_i$. We define  
$$
W_i'= [ T(\varphi_i \circ \varphi^{-1})(\varphi (x))]^{-1}(W_i),
$$
 and observe that the $W_i$'s  are open subsets of $TO_i$. By \eqref{basis_top_1},  $(\varphi (x), h)\in W_1'\cap W_2'$ and if  $W=W'_1\cap W_2'$ and $\wt{W}=(T\varphi )^{-1}(W)$, then 
$$[(x, \varphi, V, (O, C, E), h)]\in \wt{W}\subset \wt{W}_1\cap \wt{W}_2.$$
Consequently, $\mathscr{B}$ defines a topology on $TX$. To prove that this topology is Hausdorff, we take two distinct elements 
$$
\alpha=[(x, \varphi, V, (O, C, E), h)]\quad \text{and}\quad \alpha'=[(x', \varphi', V', (O', C', E'), h')]
$$
in the tangent space $TX$.  Since $\alpha\neq \alpha'$, either $x\neq x'$ or if $x=x'$, then $h'\neq T(\varphi'\circ \varphi)^{-1}(\varphi (x))h.$ In the first case $x\neq x'$ we may replace $V$ and $V'$ by smaller open neighborhoods of $x$ and $x'$, so that $V\cap V'=\emptyset$, and then replace $(O, C, E)$, resp. $(O, C', E')$,  by the retracts $(\varphi (V), C, E)$, resp. $(\varphi'(V'), C', E')$. If $W$ (resp. $W'$)  is an open neighborhood of  $(\varphi (x), h)$ (resp. $(\varphi'(x'), h')$)  in $TO$ (resp. $TO'$), then $\wt{W}=(T\varphi)^{-1}(W)$ (resp. $\wt{W}'=(T\varphi')^{-1}(W')$ is an open neighborhood  
of $\alpha$ (resp. $\alpha'$) in $TX$ and  $\wt{W}\cap \wt{W'}=\emptyset$.  In the second case $x=x'$ and $h'\neq T(\varphi'\circ \varphi)^{-1}(\varphi (x))h$. We choose an open neighborhood 
 $W$ of $(\varphi' (x), h')$ in $TO'$ and an open neighborhood $W$ of $(\varphi (x), h)$ in $TO$, so that $W'\cap T(\varphi')^{-1}(W)=\emptyset$. 
  Then $\alpha\in \wt{W}=(T\varphi)^{-1}(W)$ and $\alpha'\in \wt{W'}= (T\varphi')^{-1}(W')$. Moreover, both sets are open  and their intersection is empty. Consequently, the topology defined by $\mathscr{B}$ is Hausdorff.\par
  
\noindent (2)\,   We start by proving that  $p\colon 
TX\to X^1$ is an open map. It suffices to show that $p(\wt{W})$ is open in $X^1$ for every element
 $\wt{W}\in\mathscr{B}.$
Let $(\varphi, V, (O, C, E))$ be a chart on $X$ and let $T\varphi\colon 
TV\to TO$ be the associated map defined above and introduce $\wt{W}=(T\varphi )^{-1}(W)$ for the open subset $W$ of $TO$. 
We denote by $\pi\colon 
TO\to O_1$ the projection onto the first factor. This map is continuous and open. Moreover,  
$$p(\wt{W})=\varphi^{-1}\circ \pi\circ (T\varphi)(\wt{W}).$$
Since, by construction, the map $T\varphi\colon 
TV\to TO$ is open and $\varphi\colon 
V\to O$ is a homeomorphism, the composition on the right hand side is an open subset of $X^1$. Hence, $p(\wt{W})$ is open in $X^1$ as claimed.

To show that the projection map  $p\colon TX\to X^1$ is continuous, it suffices to show, that given a chart $(V, \varphi, (O, C, E))$ on $X$ and an open subset $U$ of $X^1$, satisfying $U\subset V_1$, the preimage $p^{-1}(U)$ is open. For such a chart and open set $U$ we have 
$$p^{-1}(U)=(T\varphi )^{-1}\bigl((\varphi (U)\times E)\cap TO\bigr).$$
Since   $(\varphi (U)\times E)\cap TO$ is open in $TO$, the set on the left-hand side belongs to $\mathscr{B}$. Hence, the set $p^{-1}(U)$ is open, and the projection $p$ is continuous as claimed. \par

\noindent (3)\,  
We start with an atlas $\mathscr{V}=(V^j, \varphi^j, (O^j, C^j, E^j))_{j\in J}$, such  that the family $\mathscr{V}=(V^j)_{j\in J}$ of domains is an open,  locally finite cover of $X$.  The associated maps 
$T\varphi^j\colon 
TV^j\to TO^j$ are homeomorphisms, and since $TO^j$ is metrizable, the same holds for the open sets $TV^j$.  Hence, $TX$ is locally metrizable.  To show that $TX$ is metrizable, it remains  to show that $TX$ is paracompact.  By Theorem \ref{equivalent_definitions_paracompact}, there exists a closed, locally finite refinement $\mathscr{C}=(C^j)_{j\in J}$ of $\mathscr{V}$, so that $\mathscr{C}_1=(C^j_1)_{j\in J}$ is a closed, locally finite refinement of $\mathscr{V}_1=(V^j_1)_{j\in J}$.  The sets $K^j=TO\vert \varphi(C^j_1):=\bigcup_{x\in \varphi (C^j_1)}T_xO$ are closed in $TO$, so that 
the sets $\wt{K^j}:=(T\varphi^j)^{-1}(K^j)$ are closed subsets of $TV$. In particular, each $\wt{K^j}$ is paracompact as a closed subset of metrizable space.  Also, the family $(\wt{K^j})_{j\in J}$ is locally finite. Indeed, let $\alpha=[(x, V, \varphi, (O, C, E), h)]\in TX$. Then there exists an open neighborhood $U(x)$ of $x$ in $X_1$ intersecting at most finitely many $C^j_1$'s, say with indices $j$, belonging to a finite subset $J'\subset J$. Moreover, since $X_1$ is regular, $U(x)$ can be taken so small that also $x\in C^j_1$ for $j\in J'$. Then, setting $W(x):=TO\vert \varphi (U (x)):=\bigcup_{y\in \varphi (U(x))}T_yO$, the set $W(x)$ is an open subset of $TO$ and $\wt{W(x)}:=(T\varphi)^{-1}(W(x))$ intersects only those $\wt{K^j}$ whose indices $j$ belong to $J'$. 
Now applying  Proposition \ref{union_of_paracompact}, we conclude, that the tangent space $TX$ is paracompact and hence metrizable, in view of the Smirnov metrizability  theorem. This finishes the proof of the proposition.
\qed \end{proof}

\chapter{Basic Sc-Fredholm Theory}
In this chapter  we start with  the Fredholm theory in the sc-framework. Since sc-maps are more flabby than $C^\infty$-maps, we do not have an implicit function theorem for all sc-smooth maps.  However, for a restricted class, which occurs in the applications 
of the theory,  such a theorem is available.

\section[Sc-Fredholm Sections]{Sc-Fredholm Sections and Some of the Main Results}
The section is devoted to the basic notions and the description of the results leading  to  implicit function theorems. Our overall goal is the notion of a sc-Fredholm section  of a 
strong bundle $P\colon Y\rightarrow X$, as defined in
 Definition \ref{def_strong_bundle}. 
Sometimes we need to require that the underlying M-polyfold  $X$ is a tame M-polyfold as defined in  
Definition \ref{def_tame_m-polyfold} in order to have good versions of the implicit function theorem near the boundary.
 The section will end with some useful implicit function theorems.
The more sophisticated perturbation and transversality results are described in a later section.

We start by introducing various types of germs in the sc-context. As usual we denote  by $E$ be a sc-Banach space and by $C\subset E$ a partial  quadrant of $E$.  The sc-Banach space $E$ is equipped with the filtration 
$$E_0=E\supset E_1\supset\cdots \supset E_\infty =\bigcap_{m\geq 0}E_m.$$

\begin{definition}\label{def_germ_nbgh}\index{ Sc-germ of neighborhoods}
A  {\bf sc-germ of neighborhoods around $0\in C$}, 
denoted by ${\mathcal O}(C,0)$,  is 
a decreasing sequence 
$$
U_0\supset U_1\supset U_2\supset \cdots 
$$
where $U_m$ is a relatively open neighborhoods of $0$ in 
$C\cap E_m$. \index{D- Tangent of ${\mathcal O}(C,0)$}\par

\noindent The {\bf tangent of ${\mathcal O}(C,0)$},  denoted by $T{\mathcal O}(C,0)$, is  the decreasing sequence \index{${\mathcal O}(C,0)$}\index{$T{\mathcal O}(C,0)$}
$$
U_1\oplus E_0\supset U_2\oplus E_1\supset\cdots \, 
$$
\qed
\end{definition}
A special example of a sc-germ of neighborhoods around $0\in C\subset E$ is  a relatively open neighborhood $U$ of $C$ containing $0$ which is equipped with the induced sc-structure defined by the filtration $U_m=U\cap E_m$, $m\geq 0$.  If $E$ is infinite dimensional, the sets $U_m$ is this example 
are not bounded in $E_m$, since the inclusions $E_m\to E_0$ are compact operators. 
Another example is the decreasing sequence $U_m=U
\cap B_{1/m}^{E_m}(0)$, where $B_{1/m}^{E_m}(0)$ is the open ball in $E_m$ centered at $0$ and radius $1/m$, presents a sc-germ of neighborhoods. Here the sets $U_m$ are bounded in $E_m$ for $m>0$. We point out that the size of the sets  $U_m$ in Definition \ref{def_germ_nbgh} does not matter. In the applications the size of $U_m$ quite often decreases rapidly.

\begin{definition}\label{germy}\index{D- Sc-mapping germ}
A  {\bf $\ssc^{\pmb{0}}$-germ} $f\colon {\mathcal O}(C,0)\rightarrow F$\index{$f\colon {\mathcal O}(C,0)\rightarrow F$} into the sc-Banach space $F$, is a continuous map
$f\colon U_0\rightarrow F$ such that $f(U_m)\subset F_m$ and
$f\colon U_m\rightarrow F_m$ is continuous.
A  {\bf $\ssc^{\pmb{1}}$-germ} $f\colon {\mathcal O}(C,0)\rightarrow F$ is a  $\ssc^0$-germ which is of class  $\ssc^1$ in the sense, that there exists 
for every $x\in U_1$  a bounded linear operator $Df(x)\in L(E_0,F_0)$ such that for $h\in U_1$ with $x+h\in U_1$, 
$$
\lim_{\abs{h}_1\rightarrow 0} \frac{\abs{f(x+h)-f(x)-Df(x)h}_0 }{\abs{h}_1}=0.
$$
Moreover,  $Tf\colon U_1\oplus E_0\rightarrow TF$,  defined by $Tf(x,h)=(f(x),Df(x)h)$,  satisfies $Tf(U_{m+1}\oplus E_m)\subset F_{m+1}\oplus F_m$ and
$$
Tf\colon T{\mathcal O}(C,0)\rightarrow TF\index{$Tf\colon T{\mathcal O}(C,0)\rightarrow TF$}
$$
is a $\ssc^0$-germ. We say $f$ is a $\ssc^2$-germ provided $Tf$ is  $\ssc^1$, etc.  If  the germ $f$ is a $\ssc^k$-germ for every $k$ we call it a sc-smooth germ.  
If we write $f\colon {\mathcal O}(O,0)\rightarrow (F,0)$ we indicate
that $f(0)=0$.
\qed
\end{definition}
We shall be mostly interested in $\ssc^\infty$-germs $f\colon {\mathcal O}(C,0)\rightarrow F$. \index{Sc$^\infty$-germs}

From Definition \ref{def_strong_bundle_chart} we recall the strong bundle chart $(\Phi, P^{-1}(V), K, U\triangleleft F)$ of a strong bundle $P\colon Y\to X$ over the M-polyfold $X$, illustrated by the commutative diagram
\begin{equation*}
\begin{CD}
P^{-1}(V)@>\Phi>>K\\
@VPVV @VVpV \\
V@>\varphi>>O .
\end{CD}\,
\end{equation*}
In the diagram, $V\subset X$ is an open set and the maps
$\Phi$ and $\varphi$  are homeomorphisms. Moreover, 
$K=R(U\triangleleft F)$ is the image of the strong bundle retraction $R$, and $O=r(U)$ is the image of the sc-smooth retraction $r\colon U\to U$ of the relatively open set $U$ of the the partial quadrant $C$ in the sc-Banach space $E$. 

\begin{definition}\index{D- Sc-section germ}
A {\bf sc-smooth section germ} $(f,x_0)$ of the strong bundle $P\colon Y\to X$ is a continuous section $f\colon V\rightarrow P^{-1}(V)$ defined on  a (sufficiently small) open neighborhood $V$ of the smooth point $x_0$,  for which the following holds.
There exists a strong bundle chart 
$(\Phi,  P^{-1}(V), (K, C\triangleleft F,E\triangleleft F))$  satisfying $\varphi (x_0)=0\in O$ in which the principal part $\wh{{\bm{f}}}$ of the local continuous section $\wh{f}=\Phi \circ f\circ \varphi^{-1}\colon O\to K$ has the property that the composition 
$$\wh{{\bm{f}}}\circ r\colon {\mathcal O}(C, 0)\to F$$
is a sc-smooth germ as defined in 
Definition \ref{germy}.
\qed
\end{definition}
To recall, the section $\wh{f}\colon  O\to K\subset U\triangleleft F$, as continuous section,  is  of the form $\wh{f}(p)=(p, \wh{\bm{f}}(p))$ for $p\in O$,  and $\wh{\bm{f}}\colon O\to F$ is called its {\bf principal part}.  The smoothness properties are concerned with the point $x_0$ and encapsulated in the behavior
of $\wh{\bm{f}}\circ r $ defined on a sc-germ of neighborhoods.  By abuse of the notation we shall often use the same letter for the principal part as for the section.

In the next step we introduce the useful  notion of a  filling of a  sc-smooth section germ $(f,0)$  of  a 
strong local bundle $K\rightarrow O$ near the given smooth point  $0$.  The notion of  a filling is a new concept specific to the world of retracts. In all known applications it deals successfully with bubbling-off phenomena and similar singular phenomena.

\begin{definition}[{\bf Filling}]\label{x-filling}\index{D- Filling}
We consider a 
strong local bundle $K\to O$. We recall that $K=R(U\triangleleft F)$  where  $U\subset C\subset E$ is a relatively open neighborhood  of $0$ in the partial quadrant $C$ of the  sc-Banach space $E$ and  $F$ is a sc-Banach space. Moreover,  $R$ is a strong bundle retraction of the form 
$$R(u, h)=(r(u), \rho(u)(h)), $$
covering the tame retraction $r\colon U\to U$ onto $O=r(U),$
and $\rho (u)\colon F\to F$ is a bounded linear operator.  
We  assume that $r(0)=0$. 
A sc-smooth section germ 
$(f, 0)$ of the bundle $K\to O$ possesses a  {\bf filling}\index{D- Filling of a sc-smooth section germ}
if there exists a  sc-smooth section germ $(g, 0)$ of the bundle $U\triangleleft F\to U$ extending $f$ and having  the following properties.
\begin{itemize}
\item[(1)]\ $f(x)=g(x)$  for $x\in O$ close to $0$.
\item[(2)]\ If $g(y)=\rho (r(y))g(y)$ for a point $y\in U$ near $0$, then $y\in O$.
\item[(3)]\ The linearization of the map
$
y\mapsto  [\mathbbm{1} -\rho(r(y))]\cdot g(y)
$
at the point $0$, restricted to $\ker(Dr(0))$, defines a topological linear  isomorphism
$$
\ker(Dr(0))\rightarrow \ker(\rho (0)).
$$
\end{itemize}
\qed
\end{definition}
The crucial property of a filler is the fact that the 
solution sets $\{y\in O\, \vert \, f(y)=0\}$ and $\{y\in U\, \vert \, g(y)=0\}$ coincide near $y=0$. 
Indeed, if $y\in U$ is a solution of the filled section $g$ so that $g(y)=0$, then it follows from (2) that $y\in O$ and from (1) that $f(y)=0$. The section $g$ is, however, much  easier to analyze than  the section $f$,  whose domain of definition has a rather complicated structure. It turns out that in the applications these extensions $g$ are surprisingly easy to  detect. In the Gromov-Witten theory and the SFT they  seem almost canonical.

 The condition (3) plays a role in the comparison of the linearizations $Df(0)$ and $Dg(0)$, assuming that $f(0)=0=g(0)$, as we are going to explain next.
 
 It  follows from the definition of a retract that
$\rho (r(y))\circ \rho (r(y))=\rho (r(y)).$ Hence, since $y=0\in O$ we have $r(0)=0$ and $\rho(0)\circ \rho(0)=\rho (0)$ so that $\rho (0)$ is a linear sc-projection in $F$ and we obtain the sc-splitting
$$
F=\rho (0)F\oplus (\mathbbm{1} -\rho (0))F.
$$
Similarly, it follows from $r(r(y))=r(y)$ for $y\in U$ that $Dr (0)\circ Dr (0)=Dr(0)$ so that $Dr (0)$ is a linear sc-projection in $E$ which gives rise to the sc-splitting
$$
\alpha \oplus \beta \in  E=Dr (0)E\oplus (\mathbbm{1} -Dr (0))E.
$$
We recall that the linearization $Tf(0)$  of the section $f\colon O\to K$ at $y=0=r(0)$ is defined as the restriction of the derivative $D(f\circ r)(0)$ of the map $f\circ r\colon U\to F$ to $T_0O$. From $\rho (r(y))f(r(y))=f(r(y))$ for $y\in U$  close to $0$,  we obtain, using $f(0)=0$,  by linearization 
at $y=0$ the relation $\rho (0)Tf(0)=Tf(0)$. From $g(r(y))=f(r(y))$ for $y\in U$ near $0$ we deduce $Tg(0)=Tf(0)$ on $T_0O$. 
 From the identity
 \begin{gather*}
(\mathbbm{1} -\rho (r(y))g(r(y))=0\quad \text{for $y\in U$ near $0$},
\end{gather*}
we deduce, using $g(0)=0$, the relation
$(\mathbbm{1} -\rho (0))Dg (0)\circ Dr(0)=0$. 
Hence  the matrix representation of $D g (0)\colon E\to F$ with respect to the above splittings of $E$ and $F$ looks as follows,$$
Dg(0)\begin{bmatrix}\alpha \\ \beta
\end{bmatrix}=
\begin{bmatrix}Tf (0)&\rho (0)Dg(0)({\mathbbm 1}-Dr(0))\\0&(\mathbbm{1} -\rho (0)) Dg (0)({\mathbbm 1}-Dr(0))\end{bmatrix}\cdot
\begin{bmatrix}\alpha\\ \beta
\end{bmatrix}.
$$
In view of property (3), the linear map $\beta \mapsto (\mathbbm{1} -\rho (0))\circ Dg (0)({\mathbbm 1}-Dr(0))\beta$ from $(\mathbbm{1} -Dr (0))E$ to $(\mathbbm{1} -D\rho(0))F$ is an isomorphism of Banach spaces. Therefore, 
$$
\text{kernel}\ Dg (0)=(\text{kernel}\ Tf (0))\oplus \{0\}.$$
Moreover the filler has the following additional properties.

\begin{proposition}[{\bf Filler}]\label{filler_new_1}\mbox{}\index{P- Properties of a filler} Assume $f$ has the filling $g$ and $f(0)=0$.
\begin{itemize}
\item[{\em (1)}]\ The operator $Tf(0)\colon Dr(0)E\to \rho (0)F$ is surjective if and only if the operator $Dg(0)\colon E\to F$ is surjective.
\item[{\em (2)}]\ $Tf(0)$ is a Fredholm operator (in the classical sense) if and only if $Dg(0)$ is a Fredholm operator and $\ind Tf(0)=\ind Dg(0)$.
\end{itemize}
\end{proposition}
\begin{proof}
(2)\, To simplify the notation we abbreviate the above matrix representing $Dg(0)$ by 
$$
Dg(0)=\begin{bmatrix}A&B\\0&C\end{bmatrix}
$$
and abbreviate the above splittings by $E=E_0\oplus E_1$ and $F=F_0\oplus F_1$.
The operators in the matrix are bounded between corresponding Banach spaces and $C\colon E_1\to F_1$ is an isomorphism of Banach spaces. Therefore, if $B=0$, the operator $A=Df(0)\colon E_0\to F_0$ is Fredholm if and only if the operator 
$$\begin{bmatrix}
A&0\\0&C\end{bmatrix}\colon E\to F
$$
is Fredholm in which case their indices agree. The statement (2) now follows from the composition formula
$$
\begin{bmatrix}
{\mathbbm 1}&BC^{-1}\\
0&{\mathbbm 1}
\end{bmatrix}
\begin{bmatrix}
A&0\\
0&C
\end{bmatrix}=
\begin{bmatrix}
A&B\\
0&C
\end{bmatrix}
$$
since the first factor is an isomorphism from $F$ to $F$, and hence has index equal to $0$,  and the Fredholm indices of a composition are additive. \\[0.5ex]
(1)\, The statement (1) is an immediate consequence of our assumption that $C$ is an isomorphism.
\qed \end{proof}

To sum up the role of a filler, instead of studying the solution set of the section $f\colon O\to K$ we can as well study the solution set of the filled section $g\colon U \to U\triangleleft F$, which is defined on the relatively open set $U$ of the partial quadrant $C$ in the sc-space $E$ and which, therefore,  is easier to analyze.
\begin{definition}[{\bf Filled version}]\label{filled_version_def}\index{D- Filled version}
If $f$ is a sc-smooth section of the 
strong bundle $P\colon Y\to X$ and $x_0\in X$ a smooth point, we 
 say that the {\bf section germ $(f, x_0)$ has a filling}, if there exists a strong bundle chart  as defined in 
 Definition \ref{def_strong_bundle_chart}, 
$$
\Phi\colon \Phi^{-1}(V)\to  K\quad  \text{covering $\varphi\colon (V,x_0)\mapsto  (O,0)$,}
$$
 where $K\rightarrow O$ is a 
  strong local bundle with $0\in O\subset U$,  such that the section germ $\Phi\circ  f\circ \varphi^{-1}\colon O\to K\subset (U\triangleleft F)$ has a filling 
 $g\colon U\to U\triangleleft F$ near $0$. 
 We shall refer to the section germ $(g, 0)$ as a 
{\bf filled version of $(f, x_0)$}.
\qed
\end{definition}

The next concept is that of a   basic germ.
\begin{definition}[{\bf Basic germ}]\label{BG-00x}\index{D- Basic germ}
Let $W$ be a sc-Banach space and $C=[0,\infty )^k\oplus \R^{n-k}\oplus W$ a partial quadrant. Then 
a {\bf basic germ} is a sc-smooth germ 
$$f\colon {\mathcal O}(C, 0)\rightarrow \R^N\oplus W,$$
satisfying $f(0)=0$ and having the following property.
If  $P\colon  \R^N\oplus W\to  W$ denotes the projection, 
the germ $P\circ f\colon {\mathcal O}(C,0)\rightarrow (W,0)$ has the form
$$
P\circ f(a,w)=w-B(a,w) 
$$
for $(a,w)\in ([0,\infty )^k\oplus \R^{n-k})\oplus W$ where $B(0,0)=0$. 
Moreover, for every $\varepsilon>0$ and  every integer $m\geq 0$, the estimate
$$
\abs{B(a,w)-B(a,w')}_m\leq \varepsilon\cdot \abs{w-w'}_m
$$
holds, if $(a, w)$ and $(a, w')$ are close enough to $(0,0)$ on level $m$.
\qed
\end{definition}

\begin{remark}\index{R- On the notion of basic class}
The notion of basic class was introduced in \cite{HWZ3} where, however, we did not require 
$f(0)=0$. Instead we required that $P\circ (f-f(0))$ has  a form as described in Definition \ref{oi}. The later developments convinced us that it is more convenient to require that $f(0)=0$. 
\qed
\end{remark}

\begin{lemma}\label{new_Lemma3.9}\index{L- Derivative of a basic germ}
Let $B\colon [0,\infty)^k\oplus \R^{n-k}\oplus W\to W$ be a sc-smooth germ around $0$ satisfying the properties described in Definition \ref{BG-00x}. Then 
for every $\varepsilon>0$ and $m\geq 0$,
\begin{equation}\label{est_second_der_B}
\abs{D_2B(a, w)\zeta}_m\leq \varepsilon\abs{\zeta}_m
\end{equation}
for all $\zeta\in W_m$, if $(a, w)\in E_{m+1}$ is close enough to $(0, 0)$ in $E_m$. In particular, 
$$D_2B(0, 0)=0.$$
\end{lemma}
\begin{proof}
Since $W_{m+1}\subset W_m$ is dense, it is sufficient to verify the estimate for $\zeta\in W_{m+1}$ satisfying $\abs{\zeta}_{m+1}=1$. For such a $\zeta$ we know  from the definition of the linearization,  recalling Proposition \ref{sc_up},  that $B(a, w)-B(a, w+t\zeta)-D_2B(a, w)(t\zeta)=o(t)$, where $w+t\zeta\in C$,  $\abs{o(t)/t}_m\to 0$ as $t\to 0$. Therefore, 
\begin{equation*}
\begin{split}
\abs{D_2B(a, w)\zeta}_m&=\abs{\dfrac{1}{t}D_2B(a, w)(t\zeta)}_m\\
&\leq \dfrac{1}{\abs{t}}\abs{B(a, w)-B(a, w+t\zeta)}_m+ \dfrac{1}{\abs{t}}\abs{o(t)}_m.
\end{split}
\end{equation*}
The first  term on the right hand side is estimated by $\varepsilon\abs{\zeta}_m$ if $(a, w)$ and $(a, w+t\zeta)$ are sufficiently small in $E_m$. Therefore, the estimate \eqref{est_second_der_B} follows as $t\to 0$.
\qed \end{proof}

We will see that basic germs have special properties as already the following application of Lemma \ref{new_Lemma3.9} demonstrates.
\begin{proposition}\label{Newprop_3.9}\index{P- Fredholm property of basic germs}
Let $f\colon [0,\infty)^k\oplus \R^{n-k}\oplus W=E \to \R^N\oplus W$ be a sc-smooth germ around $0$ of the form $f=h+s$ where $h$ is a basic  germ and $s$ a $\ssc^+$-germ. Then 
$$Df(0)\colon \R^n\oplus W\to \R^N\oplus W$$
is a sc-Fredholm operator and its index is equal to 
$$\ind Df(0)=n-N.$$
Moreover, for every $m\geq 0$,
$$Df(a, w)\colon \R^n\oplus W_m\to \R^N\oplus W_m$$ is a Fredholm operator having index $n-N$,  if $(a, w)\in E_{m+1}$ is sufficiently small in $E_m$.
\end{proposition}

\begin{proof}

With the sc-projection $P\colon \R^N\oplus W\to W$, the linearization of $f$ at the smooth point $0$, 
$Df(0)=P\circ Df(0)+({\mathbbm 1}-P)\circ Df(0)$, is explicitly given by the formula
\begin{equation*}
\begin{split}
Df(0)(\delta a, \delta w)&=\delta w-D_2B(0)\delta w-D_1B(0)\delta a\\
&\phantom{=}+({\mathbbm 1}-P)\circ Df(0)(\delta a, \delta w) +Ds(0)(\delta a, \delta w).
\end{split}
\end{equation*}
By Lemma \ref{new_Lemma3.9}, $D_2B(0)=0$. Therefore, the operator $Df(0)$ is a $\ssc^+$-perturbation of the operator 
\begin{equation}\label{new_equation_sc-pert}
\R^n\oplus W\to \R^N\oplus W, \quad (\delta a, \delta w)
\mapsto (0,\delta w).
\end{equation}
The operator \eqref{new_equation_sc-pert} is a sc-Fredholm operator whose 
kernel is equal to $\R^n$ and whose cokernel is $\R^N$, so that its Fredholm index is equal to $n-N$. 
Since $Df(0)$ is a $\ssc^+$-perturbation of a sc-Fredholm operator, it is also a   sc-Fredholm operator by Proposition \ref{prop1.21}. Because $\ssc^+$-operators are compact, if considered on the same level, the Fredholm index is unchanged and so, 
$\ind Df(0)=n-N$. 

The second statement follows from the fact that the linear operator 
$$
Df(a, w)\colon \R^n\oplus W_m\to \R^N\oplus W_m
$$
 is a compact perturbation of the operator 
\begin{equation}\label{new_equation_fredholm_pert}
(\delta a, \delta w)
\mapsto (0,({\mathbbm 1}-D_2B(a, w))\delta w).
\end{equation}
Choosing $0<\varepsilon <1$ in Lemma \ref{new_Lemma3.9}, the operator 
${\mathbbm 1}-D_2B(a, w)\colon W_m\to W_m$ is an isomorphism of Banach spaces if $(a, w)\in E_{m+1}$ is sufficiently small in $E_m$. Hence the operator \eqref{new_equation_fredholm_pert} is a Fredholm operator of index $n-N$ and the proposition follows.

\qed \end{proof}

Finally,  we are in a position to introduce the  sc-Fredholm germs.

\begin{definition}[{\bf sc-Fredholm germ}]\label{oi}\index{D- Sc-Fredholm germ}
Let $f$ be a sc-smooth section of the strong bundle $P\colon Y\to X$ over the 
M-polyfold $X$, and let $x_0\in X$ be a smooth point.  Then $(f,x_0)$ is a {\bf sc-Fredholm germ} provided it possesses  a filled version 
$$
(g,0)\colon {\mathcal O}(U,0)\to U\triangleleft F
$$
 according to Definition \ref{filled_version_def} and having the following property. 
There exists a local $\ssc^+$-section $s\colon U\to U\triangleleft F$ satisfying $s(0)=g(0)$ such that  the germ 
$(g-s,0)\colon {\mathcal O}(U,0)\to U\triangleleft (F,0)$
is conjugated to a basic germ. 
The last condition requires the existence of a strong bundle isomorphism  defined near $0\in U$ covering the sc-diffeomorphism 
$\psi$ such that the push-forward section 
$$
\Psi \circ (g-s)\circ \psi^{-1}
\colon {\mathcal O}(U',0) \to (U'\triangleleft F',(0,0))
$$
is a basic germ (near $0$).
\qed
\end{definition}
From Proposition \ref{Newprop_3.9} we deduce that the linearization $D(g-s)(0)$ at the point $0$ is a sc-Fredholm operator. Consequently,  $Dg(0)$ is a sc-Fredholm operator by Proposition \ref{prop1.21}, and so, the tangent map 
$Tf(x_0)\colon T_{x_0}X\to Y_{f(x_0)}Y$ is a linear sc-Fredholm operator having the same index as $Dg(0)$, namely $\ind Tf(x_0)=n-N$, in view of the properties of a filler in Proposition \ref{filler_new_1}.

\begin{remark}\index{R- Comments on the notion of sc-Fredholm germ}
The above definition of a sc-Fredholm germ looks  complicated; one first has to find a filled version, which then, after some correction by a $\ssc^+$-section, is conjugated to a basic germ.
It turns out that the definition of a sc-Fredholm germ is extremely practicable in the applications we have in mind. 
By experience one may say that the fillings, which are usually only needed
near data describing bubbling-off situations seem almost natural, i.e. ``if one sees one example one has seen them all''. 
Examples of fillings  in the  Gromov-Witten, SFT and Floer Theory can be found in  \cite{HWZ5,FH2}.
The subtraction of a suitable $\ssc^+$-section is in applications essentially the removing of lower order terms of a nonlinear differential operator
and therefore allows simplifications of the expressions before one tries to conjugate them to a basic germ.
 One also has to keep in mind that the sc-Fredholm theory is designed to cope with spaces whose tangent spaces have locally varying dimensions, on which, on the analytical side, one studies systems of partial differential equations on varying domains into varying codomains.  Later on 
we shall give criteria (which in practice are easy to check)  to verify  that a section is conjugated to a basic germ.
In \cite{FH2} a pre-Fredholm theory has been developed which allows to build complicated sc-Fredholm sections from simpler pieces.
\qed
\end{remark}

Sc-Fredholm germs possess a useful local compactness property.
\begin{theorem}[{\bf Local Compactness for sc-Fredholm Germs}] \label{compact-x}\index{T- Local compactness}
Let $f$ be a sc-smooth section of the 
strong bundle $P\colon Y\rightarrow X$, and $x_0\in X$  a smooth point. We assume that $(f,x_0)$ is a sc-Fredholm germ satisfying $f(x_0)=0$. Then there exist a nested sequence
of open neighborhood ${\mathcal O}(i)$ of $x_0$ in $X_0$,  for $i\geq 0$, 
$$
{\mathcal O}(0)\supset  {\mathcal O}(1)\supset \cdots  \supset {\mathcal O}(i)\supset   {\mathcal O}(i+1) \supset \cdots  \, ,
$$
such  that for every $i$ the $X_0$-closure 
 $\cl_{X_0}(\{x\in {\mathcal O}(i)\,  \vert \,  f(x)=0\})$ is a compact subset of $ X_i$.
 \qed
\end{theorem}
\begin{remark}\index{R- On $O(i)$}
We emphasize that the ${\mathcal O}(i)$ are open neighborhoods in $X$, i.e. on level $0$.
\qed
\end{remark}

Theorem \ref{compact-x} is an immediate consequence of 
Theorem \ref{save}, which will be introduced later,
and has the following corollary.
\begin{corollary}\index{C- Local compactness and regularity}
Let $f$ be a sc-smooth section of the 
strong bundle $P\colon Y\rightarrow X$, and $x_0\in X$  a smooth point. We assume that $(f,x_0)$ is a sc-Fredholm germ satisfying $f(x_0)=0$. 
If $(x_k)\subset X$ is a sequence satisfying   $f(x_k)=0$ and $x_k\rightarrow x$ in $X_0$,  then it follows, for every  given any $m\geq 0$, that 
$x_k\in X_m$ for $k$ large and $x_k\rightarrow x$ in $X_m$.
\qed
\end{corollary}
Next we give the basic definition of an sc-Fredholm section.

\begin{definition}[{\bf Sc-Fredholm section}]\index{D- Sc-Fredholm section}\label{DEFX3116}
A section  $f$ of the 
strong bundle $P\colon Y\rightarrow X$ over the M-polyfod $X$ is called  {\bf sc-Fredholm section}, if it has 
the following three properties.
\begin{itemize}
\item[(1)]\ $f$ is sc-smooth.
\item[(2)]\ $f$ is regularizing, i.e.,  if $x\in X_m$ and  $f(x)\in Y_{m,m+1}$, then $x\in X_{m+1}$.
\index{D- Regularizing section}
\item[(3)]\ The germ $(f,x)$ is a sc-Fredholm germ at  every smooth point $x\in X$.
\end{itemize}
\qed
\end{definition}

\begin{remark}\index{R- On the sc-implicit funtion theorem}
The  implicit function theorem,  introduced later on,  is applicable to sc-Fredholm sections and will lead to the following local  result near a smooth interior point $x_0\in X$ ( i.e.,  $d_X(x_0)=0$). Assuming that  $f(x_0)=0$  the linearization $f'(x_0)\colon T_{x_0}X\rightarrow Y_{x_0}$ 
is a sc-Fredholm operator. Moreover, if  $f'(x_0)$ is surjective, then the solution set
$\{x\in X\, \vert \, f(x)=0\}$ near $x_0$ has the structure of a finite dimensional smooth manifold (in the classical sense) whose dimension agrees with the Fredholm index.  Its smooth structure is in a canonical way  induced 
from the M-polyfold structure of $X$.
In case that $x$ is a boundary point, so that  $d_X(x)\geq 1$, we can only expect the solution set to be reasonable provided 
the kernel of $f'(x)$ lies in good position to the boundary of $X$ and the boundary $\partial X$ is sufficiently well-behaved.
In order that $\partial X$ is regular enough we shall require 
that $X$ is a tame M-polyfold,  so that we can ask  $\ker(f'(x))$ to be in good position to the partial quadrant $C_xX$ in $T_xX$, 
a  notion which we shall introduce later on.
\qed
\end{remark}

 If $P\colon Y\to X$ is a 
 strong bundle,  we denote by $\Gamma(P)$\index{$\Gamma(P)$} the vector space of sc-smooth sections; by  $\text{Fred}(P)$\index{$\text{Fred}(P)$} we denote the subset of $\Gamma(P)$ consisting of sc-Fredholm sections. Finally, by  
$\Gamma^+(P)$\index{$\Gamma^+(P)$} we denote  
the vector space of $\ssc^+$-sections as introduced in Definition 
\ref{def_sc_inft_sections}.
The following stability property   of a sc-Fredholm section will be crucial for the transversality theory.

\begin{theorem}[{\bf Stability under $\ssc^{\pmb{+}}$-perturbations}]\label{stabxx}  \index{T- Stability of sc-Fredholm sections}
Let $P\colon Y\rightarrow X$ be a  strong bundle over the 
M-polyfold $X$.
If  $f\in \Fred (P)$ and $s\in \Gamma^+(P)$, then $f+s\in \Fred (P)$.
\qed
\end{theorem}

In order to prove the theorem we need  two lemmata for local strong bundles. We recall  the local strong bundle retract  $(K, C\triangleleft F, E\triangleleft F)$ from  Definition \ref{def_loc_strong_b_retract}, consisting of the retract $K=R(U\triangleleft F)$, where $R\colon U\triangleleft F\to U\triangleleft F$ is a strong bundle retraction of the form 
$$R(u, h)=(r(u), \rho (u)h), $$
in which $r\colon U\to U$ is a smooth retraction onto $O=r(U)\subset  U$. We shall denote the principal parts  of the sections  of the 
bundles $U\triangleleft F\to U$ and $K\to O$  by bold letters.

\begin{lemma}\label{rio}
A $\ssc^+$-section $s\colon O\to K$ (as defined in Definition \ref{def_section_loc_strong_bundle}) possesses an extension to a $\ssc^+$-section 
$\wt{s}\colon U\to U\triangleleft F$, 
$\wt{s}(u)=(u, \wt{\bm{s}}(u))$, 
having the following properties.
\begin{itemize}
\item[{\em (1)}]\ $\wt{s}(u)=s(u)$ if $u\in O$.
\item[{\em (2)}]\  $R(r(u),\wt{\bm{s}}(u))=s(r(u))$ if $u\in U$.
\end{itemize} 
\end{lemma}
\begin{proof}
If $s\colon O\to K$  is given by $s(u)=(u, {\bm{s}}(u))$, $u\in O$, we define the section $\wt{s}\colon U\to U\triangleleft F$ by 
$$\wt{s}(u)=(u, \wt{\bm{s}}(u))=\bigl(u, {\bm{s}}(r(u))\bigr),\quad u\in U.$$
Clearly, $\wt{s}$ is a $\ssc^+$-section  of the bundle 
$U\triangleleft F\to U$ and we claim that its restriction to $O$ agrees with the section $s$. Indeed, if $u\in O$, 
then $r(u)=u$, implying $\wt{s}(u)=\wt{s}(r(u))=\bigl( r(u), {\bm{s}}(r\circ r(u))\bigr)=\bigl( r(u), {\bm{s}}(r(u))\bigr)=
\bigl(u, {\bm{s}}(u)\bigr)=s(u)$ as claimed.  Moreover, using that $s(u)=R(s(u))$ if $u\in O$,
$$R\bigl(r(u), \wt{\bm{s}}(u)\bigr)
=R\bigl(r(u), {\bm{s}}(r(u))\bigr)=R\bigl({ s}(r(u))\bigr)=s(r(u))$$
for $u\in U$.
\qed \end{proof}
\begin{lemma}\label{filler_extension_s}
Let $f\colon O\to K$ be a sc-smooth section of the (previous) local strong bundle retract, and let $s\colon O\to K$ be a $\ssc^+$-section. 
If $f$ possesses the filler $g\colon U\to U\triangleleft F$, then $f+s$ has the filler $g+\wt{s}\colon U\to U\triangleleft F$, where $\wt{s}$ is the extension of $s$ constructed in the previous lemma.
\end{lemma}
\begin{proof}
We have to verify that the section $g+\wt{s}$ meets the three conditions in Definition \ref{x-filling}. The properties  (1) and (2) for $g+\wt{s}$ follow immediately from the properties (1) and (2)   for the filler $g$ and the properties (1) and (2) for the section $\wt{s}$ in Lemma \ref{rio}. In order to verify property (3) of a filler we have to linearize the map 
$$u\mapsto [{\mathbbm 1}-\rho (r(u))](g(u)+\wt{s}(u))$$
at the point $u=0$. Since 
$\bigl({\mathbbm 1}-\rho (r(u))\bigr) \wt{s}(u)=
\bigl({\mathbbm 1}-\rho (r(u))\bigr) s(r(u))=0$ by property (2) of Lemma \ref{rio}, the linearization agrees with the linearization of the map $\bigl({\mathbbm 1}-\rho (r(u))\bigr) g(u)$ which satisfies the required property (3), since $g$ is a filler. The proof of Lemma \ref{filler_extension_s} is finished.
\qed \end{proof}

\begin{proof}[Theorem \ref{stabxx}]
Let $f$ be a sc-Fredholm section of the 
strong bundle $P\colon Y\to X$ and let $s\colon X\to Y$ be a $\ssc^+$-section of $P$. Then $f+s$ is a sc-smooth section which is also regularizing. It remains to verify that $(f+s, x)$ is a sc-Fredholm germ for every smooth point $x\in X$. By definition of being sc-Fredholm,
$(f, x)$ is a sc-Fredholm germ at the smooth point $x$. Therefore, there exists an open neighborhood $V$ of $x$ and a strong bundle chart 
$(V, P^{-1}(V), (K, C\triangleleft F,E\triangleleft F))$ as defined in Definition
 \ref{def_strong_bundle_chart} and satisfying $\varphi (x)=0\in O$, such that the local representation $\wt{f}=\Phi_\ast (f)=\Phi\circ f\circ \varphi^{-1}\colon O\to K$ of the section $f$ possesses a filled version $g\colon U\to 
U\triangleleft F$ around $0$,  which after subtraction of a suitable $\ssc^+$-section,  is conjugated to a basic germ around $0$. Define $t=\Phi_\ast (s)$. Then $t\colon O\to K$ is a $\ssc^+$-section. By Lemma \ref{rio}  
there is a particular $\ssc^+$-section $\wt{t}\colon U\to U\triangleleft F$ extending $t$. By Lemma \ref{filler_extension_s}, the section $g+\wt{t}\colon U\to U\triangleleft F$ is a filling of $\wt{f}+t$.
In view of the sc-Fredholm germ property, there exists a $\ssc^+$-section $t'$ satisfying $t'(0)=g(0)$ and such that $g-t'$ is conjugated to a basic germ. Now taking the $\ssc^+$-section $\wt{t}+t'\colon U\to U\triangleleft F$, we have $(g+\wt{t})(0)=(\wt{t}+t)(0)$. Moreover, $(g+\wt{t})-(\wt{t}+t')=g-t'$ which, as we already know, is conjugated to a basic germ. To sum up, we have verified that $(f+s, x)$ is a sc-Fredholm germ. This holds true for every smooth point $x\in X$. Consequently, the section $f+s$ is a Fredholm section and the proof of Theorem \ref{stabxx} is complete.
\qed \end{proof}

In order to formulate  a parametrized version of Proposition \ref{stabxx} 
we assume that $P\colon Y\rightarrow X$ is a  strong bundle and $f$ a sc-Fredholm section. The sc-smooth projection 
$$
\pi\colon \R^n\times X\rightarrow X,\quad  (r,x)\mapsto  x,
$$
pulls back the bundle $P$ to the  strong bundle   $\pi^\ast(P)\colon \pi^\ast Y\rightarrow \R^n\times X$.
The section $\wt{f}$ of $\pi^\ast(P)$,  defined by
$$
\wt{f}(r,x)=((r,x),f(x)), 
$$
 is a sc-Fredholm section as is readily verified.  
If $s_1,\ldots ,s_n$ are $\ssc^+$-sections of $P$, then $\wt{s}(r,x)\colon =\bigl( (r,x),\sum_{i=1}^n r_i\cdot s_i(x)\bigr)$ is a $\ssc^+$-section
 of  the pull back bundle $\pi^\ast(P)$ and,  by the stability 
 Theorem \ref{stabxx}, the section 
 $$
 (r,x)\mapsto \wt{f}(r,x)+\wt{s}(r,x)
 $$
 is a sc-Fredholm section of $\pi^\ast(P)$. Hence we have proved the following stability result.
 \begin{theorem}[{\bf Parameterized Perturbations}]\label{corro}\index{T- Parameterized perturbations}
 Let $P\colon Y\rightarrow X$ be a 
 strong bundle and $f$ a sc-Fredholm section. If  $s_1,\ldots ,s_n\in\Gamma^+(P)$, 
then the map
 $$
 \R^n\times X\rightarrow Y,\quad (r,x)\mapsto f(x)+\sum_{i=1}^n r_i\cdot s_i(x)
 $$
 defines a sc-Fredholm section of the bundle $\pi^\ast(P)\colon \pi^\ast Y\rightarrow \R^n\times X$.
 \qed
 \end{theorem}
This theorem and refined versions of the theorem play a role  in the perturbation and transversality theory. As already pointed out,  
the distinguished class of sc-Fredholm sections allows  to apply an implicit function theorem in the usual sense. 
We first formulate the  implicit function theorem at an interior point. 
\begin{theorem}[{\bf Implicit Function Theorem}]\label{implicit-x}\index{T- Implicit function theorem}
Assume that $P\colon Y\rightarrow X$ is a strong bundle over the  M-polyfold $X$ satisfying $d_X\equiv 0$,  and $f$ a sc-Fredholm section.
Suppose that $x_0\in X$ is a smooth point in $X$,  such  that $f(x_0)=0$. Then the linearization  $f'(x_0)\colon T_{x_0}X\rightarrow Y_{x_0}$  is a sc-Fredholm operator. If $f'(x_0)$  surjective,  
then there exists an open neighborhood $U$ of $x_0\in X$ such  that the solution set $S(f,U)=\{x\in U\, \vert \,  f(x)=0\}$ in $U$ has in a natural way the structure
of a smooth finite dimensional manifold whose dimension agrees with the Fredholm index.  In addition,  $U$ can be chosen  in such a way that the linearization $f'(y)\colon T_yX\rightarrow Y_y$ for $y\in S(f,U)$ is surjective  and  $\ker(f'(y))=T_yS(f,U)$ is the tangent space.
\qed
\end{theorem}
Theorem \ref{implicit-x} is an immediate consequence of Theorem \ref{IMPLICIT0} in Section \ref{ssec3.4}.
As we shall see, the smooth manifold structure on $S(f,U)$ is induced from the M-polyfold structure of $X$.
\begin{remark}\index{R- On the definition of sc-Fredholm section}\label{remark.3_20}
We should point out that in the original proof of Theorem \ref{IMPLICIT0} in \cite{HWZ3} (Theorem 4.6 and Proposition 4.7) the sc-Fredholm section is defined slightly differently, namely as follows. In \cite{HWZ3}
a sc-smooth section $f$ of the strong bundle $Y\to X$ is called sc-Fredholm, if it possesses around all smooth points of $X$ a filled version 
 $(g,0)$ such that $g-g(0)$ near $0$ is conjugated to a basic germ. It is, in this case, not true that $f+s$ is sc-Fredholm for $\ssc^+$-section $s$. However, by a nontrivial theorem in \cite{HWZ3} (Theorem 3.9), increasing the level, the section 
$(f+s)^1$ of the the strong bundle $Y^1\rightarrow X^1$ is sc-Fredholm. Although not harmful in practice,  this looks unsatisfactory.

This is why we have introduced the new definition (Definition \ref{oi}) of sc-Fredholm, where we require for the filled version $(g, 0)$ that there exists a local $\ssc^+$-section $s\colon U\to U\triangleleft F$ such that $s(0)=g(0)$ and $g-s$ is locally conjugated to a basic germ.
If now $t$ is a $\ssc^+$-section, then the section $f+t$ of the bundle $Y\to X$ is automatically sc-Fredholm in view of Theorem \ref{stabxx}. 

The difficulty of the nontrivial theorem is now hidden in the proof of the implicit function theorem, which has to incorporate
the arguments of the nontrivial theorem. With the new definition, even if we want to study $f$ only, we have
only normal forms for the perturbed expression, which might be unrelated to our problem.  However, writing $f=(f-s) + s$, we know how $f-s$ looks like, and we know that $s$ is a compact perturbation of $f-s$.  We  combine  these facts to gain sufficient information about $f^1$ to determine, in view the regularizing property of sc-Fredholm sections, the solution set of $f$.
\qed
\end{remark}

Next we discuss some consequences of Theorem \ref{implicit-x}.  Considering a sc-Fredholm section $f$ of the 
strong bundle $Y\rightarrow X$ over the  M-polyfold $X$ having no boundary (i.e., $d_X\equiv 0$),
we assume, in addition,  that the M-polyfold $X$ admits a sc-smooth partition of unity.
Then there exists, for two given smooth points $x\in X$ and $e\in Y_x$, a $\ssc^+$-section $s$ supported near $x$ and satisfying $s(x)=e$.  For the easy proof we refer to \cite{HWZ3}. As we shall see later it suffices to assume the existence of sc-smooth bump functions instead of sc-smooth partitions of  unity.
If $f(x_0)=0$ and $f'(x_0)$ is not surjective, we find finitely many smooth elements $e_1,\ldots ,e_k\in Y_{x_0}$ satisfying $R(f'(x_0))\oplus \R e_1\oplus \cdots \oplus \R e_k=Y_{x_0}$. Taking  $\ssc^\infty$-sections $s_i$ satisfying $s_i(x_0)=e_i$, we define the map $\wt{f}\colon \R^k\oplus X\to Y$  by 
$$\wt{f}(r, x)=f(x)+\sum_{i=1}^k r_is_i(x).$$
The linearization of $\wt{f}$ at  the distinguished point $(0,x_0)\in  \R^k\oplus X$ is the continuous linear map 
$$
\wt{f}'(0, x_0)(h,u)= f'(x_0)u +\sum_{i=1}^k h_i e_i,
$$
which is surjective. By Theorem \ref{corro}  and Theorem \ref{IMPLICIT0} we find an open neighborhood $U$ of $(0,x_0)\in \R^k\oplus X$
such  that the solution set of $S(\wt{f}, U)=\{(r,x)\in U\, \vert \,  f(x)+\sum_{i=1}^k r_i\cdot s_i(x)=0\}$ is a smooth finite-dimensional manifold.
The trivial bundle $S(\wt{f},U)\times \R^k\rightarrow S(\wt{f},U)$ has the canonical section $(r,y)\mapsto r$. 
The zero set of this section is precisely the unperturbed solution set of $f(y)=0$ for $y\in U$.
If the solution set of $f(y)=0$, $y\in X$, is compact we can carry out the previous construction globally, which gives rise to  global finite-dimensional reduction. This will be discussed later on.

In order to formulate the boundary version of the implicit function theorem we start with some preparation.

\begin{definition}[{\bf In good position}] \label{mission1}\index{D- Good position}\index{D- Good complement}
Let $C\subset E$ be a partial quadrant in the sc-Banach space $E$ and $N\subset E$ be a finite-dimensional sc-subspace of $E$.  The subspace $N$ is in  {\bf good position to the partial quadrant $C$}, if the interior of $N\cap C$ in $N$  is non-empty, and if $N$ possesses a  sc-complement $P$, so that $E=N\oplus P$, having the following property. There exists  $\varepsilon>0$, such that 
for pairs $(n, p)\in N\times P$ satisfying $\abs{p}_{0}\leq \varepsilon\abs{n}_{0}$ the statements 
$n\in C$ and $n+p\in C$ are equivalent. We call such a sc-complement $P$ a {\bf good complement}.
\qed
\end{definition}

The choice of the right complement $P$ is important. One cannot take a random sc-complement of $N$, in general,  as Lemma 
\ref{cone1} demonstrates.

In view of Proposition \ref{prop1} the finite dimensional subspace $N$ in Definition \ref{mission1}, possessing the sc-complement $P$, is necessarily a smooth subspace. A finite dimensional subspace $N$ which is not necessarily smooth is called in good position to the partial quadrant $C$ in the sc-Banach space $E$ if there exists a (merely) topological complement $P$ in $E$ satisfying the requirements of Definition \ref{mission1} for some $\varepsilon>0$.

\begin{figure}[htb]
\centering
\def\svgwidth{65ex}
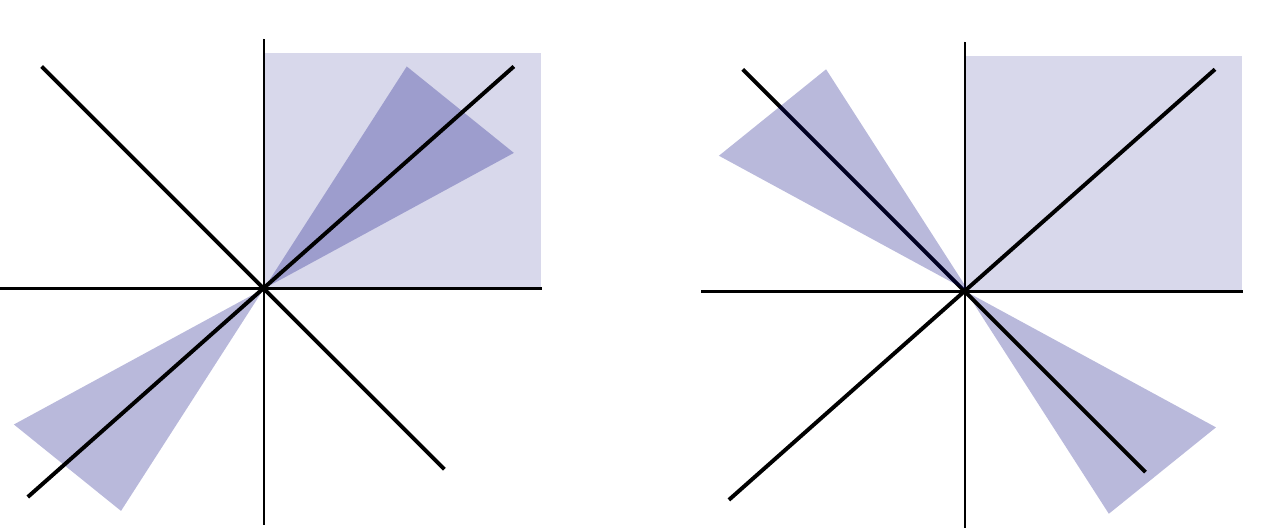
\caption{In figure (a) $N$ is in good position to $C$ while in figure (b) $N$ is not in good position to $C$.}\label{fig:pict3}
\end{figure}

The following result  is taken from \cite{HWZ3}, Proposition 6.1. Its proof is reproduced in Appendix \ref{pretzel-A}.
\begin{proposition}\label{pretzel}\index{P- Good position and partial quadrants}
If  $N$ is a finite-dimensional sc-subspace in good position to the partial quadrant  $C$ in $E$, then
$N\cap C$ is a partial quadrant in $N$.
\qed
\end{proposition}

The boundary version of the implicit function theorem is formulated in the next theorem. The proof is again given later.

\begin{theorem}[{\bf Implicit Function Theorem: Boundary Case}]\label{bound}\index{T- Implicit function theorem {II}}
We assume that $P\colon Y\rightarrow X$ is a strong bundle over the tame M-polyfold $X$, and $f$ is a sc-Fredholm section.
Suppose that $x\in X$ satisfies $f(x)=0$ and the following two properties.
\begin{itemize}
\item[{\em (a)}]\ The linearisation $f'(x)\colon T_xX\rightarrow Y_x$ is surjective.
\item[{\em (b)}]\ The kernel $N$ of $f'(x)$ is in good position to the boundary of $X$, i.e. $N$ is in good position to the partial quadrant $C_xX$ in the tangent space $T_xX$.
\end{itemize}
Then there exists an open neighborhood $U$ of $x$ such  that the following holds.
\begin{itemize}
\item[{\em (1)}]\ The local solution set $S(f,U):=\{y\in U\ |\ f(y)=0\}$, which consists of smooth points, is a tame sub-M-polyfold of $X$. 
\item[{\em (2)}]\ The tame sub-M-polyfold $S(f,u)$ admits a uniquely determined sc-smooth equivalent structure as a smooth manifold with boundary with corners. 
\item[{\em (3)}]\ For every $z\in S(f,U)$ the linearization $f'(z):T_zX\rightarrow Y_z$ is surjective and $\ker(f'(z))=T_zS(f,U)$.
\end{itemize}
\qed
\end{theorem}
Theorem \ref{bound} is a consequence of Theorem \ref{IMPLICIT0} in Section \ref{ssec3.4}.
In the proof of Theorem \ref{bound} we shall describe the manifold structure on the solution space in  detail.
Here we just indicate  how it looks like.  Since $S(f,U)\subset X_\infty$ and $X$  is a tame M-polyfold, we can take for a point $y\in S(f,U)$
a sc-diffeomorphism $\Psi\colon U(y)\rightarrow O=O(0)$, where $(O,C,E)$ is a tame retract. Then if $t\colon V\rightarrow V$ satisfies
$O=t(V)$ we have the splitting $E=T_0O\oplus (({\mathbbm 1}-Dt(0))E)$, with $Y=({\mathbbm 1}-Dt(0))E$ contained in $T_0^RC$. 
Also the proof will show that $T_0O$ is finite-dimensional. Let $p=Dt(0)$ be the projection onto $T_0O$.
Then it will be shown that near $0$ the projection $p\colon {\mathcal U}'(0)\rightarrow {\mathcal V}(0)$ is a sc-diffeomorphism,
where ${\mathcal U}'$ is an open neighborhood of $0$ in $O$ and ${\mathcal V}$ is an open neighborhood of $0$ in $C_0O$.
Then, for ${\mathcal U}=\Psi^{-1}({\mathcal U}')$ the map
$$
{\mathcal U}\rightarrow {\mathcal V}, \quad  y\mapsto  p\circ \Psi(y)
$$
is a sc-diffeomorphism and its mage lies in an relatively open neighborhood of $0$ in the partial quadrant $C_0O$ in $T_0O$.
The associated transition maps for any two such sc-diffeomorphisms are (trivially) sc-diffeomorphisms between relatively open subsets in partial quadrants 
of finite-dimensional vector spaces. Hence they are classically smooth. This shows that the system of such sc-diffeomorphisms
defines a smooth atlas for the structure of a manifold with boundary with corners and by construction this structure is compatible with the existing 
M-polyfold structure on $S(f,U)$.

As a corollary of Theorem \ref{bound} we shall obtain the following result.
\begin{corollary}\index{C- Global implicit function theorem}
Let $P\colon Y\rightarrow X$ be  a strong bundle over the tame M-polyfold $X$ and $f$ be a sc-Fredholm section.
Suppose that for every $x\in X$ satisfying  $f(x)=0$ the linearisation $f'(x)\colon T_xX\rightarrow Y_x$ is surjective 
and the kernel of $f'(x)$ is in good position to the boundary of $X$ (the latter being an empty condition if $d_X(x)=0$). 
Then the solution set $M:=\{x\in X\ |\ f(x)=0\}\subset X$
is a  sub-M-polyfold of $X$ for which the induced structure is tame, and which, moreover,  admits  a sc-smoothly  equivalent structure as a smooth manifold with boundary with corners.
\qed
\end{corollary}

The remaining subsections are devoted to the proofs of the above  results. Since the Fredholm theory is one of the main parts 
of the polyfold theory and draws heavily on the possibilities offered  in the sc-smooth theory we shall carry out  the constructions in great detail.

\section{Subsets with Tangent Structure}
The solution sets of sc-Fredholm sections will come with a certain structure, which in the generic case will induce a natural smooth manifold 
on the solution set. This subsection studies this structure.  Recall the definition of a smooth finite-dimensional subspace $N$ in good position to the partial quadrant $C$ (Definition \ref{mission1}). For such a subspace, $N\cap C$ is a partial quadrant in $N$ (Proposition \ref{pretzel}).

\begin{definition}[{\bf $n$-dimensional tangent germ property}]\label{toast}\index{D- Tangent germ property}
We consider a tame M-polyfold $X$ and a subset $M\subset X$ of  $X$. The subset $M$ has the {\bf $\bm{n}$-dimensional tangent germ property}
provided the following holds.
\begin{itemize}
\item[(a)]\  $M\subset X_\infty$.
\item[(b)]\  Every point $x\in M$  lies in an open neighborhood $U\subset X$ of $x$ such that there exists a sc-smooth chart
$\varphi\colon (U,x)\rightarrow (O,0)$ onto a tame  retract  $(O,C,E)$. Moreover, there exists  a $n$-dimensional smooth subspace $N\subset E$ in good position to the partial quadrant $C$, which possesses a good complement $Y$ so that $E=N\oplus Y$.  In addition,  there exists a relatively open neighborhood $V$ of $0$ in the partial quadrant $N\cap C$  and a continuous map
$\delta\colon V\rightarrow Y$ having the following properties.
\begin{itemize}
\item[(1)]\ $\varphi(M\cap U) =\{v+\delta(v)\, \vert \, v\in V\}\subset N\oplus Y$.
\item[(2)]\ $\delta\colon {\mathcal O}(N\cap C,0)\rightarrow (Y,0)$ is a $\ssc^\infty$-germ satisfying  $\delta (0)=0$ and $D\delta(0)=0$.
\end{itemize}
\end{itemize}
\qed
 \end{definition}

\begin{figure}[htb]
\begin{centering}
\def\svgwidth{50ex}
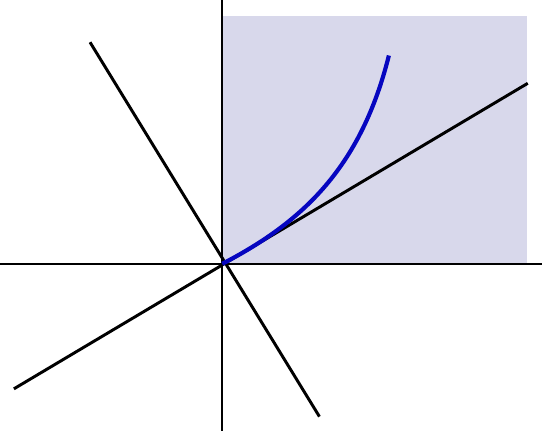
\caption{}\label{fig:pict4}
\end{centering}
\end{figure}

Recalling  the definition of a $\ssc^\infty$-germ (Definition \ref{germy}) we note that here $V$ is a relatively open neighborhood of $0$ in the partial quadrant $N\cap C$ where $N$ is a smooth  finite-dimensional space.
There exists a nested sequence $(V_m)$ of relatively open neighborhoods of $0$ in $N\cap C$, say $V=V_0\supset V_1\supset V_2\supset \ldots$,  such  that
$\delta(V_m)\subset Y_m$ and $\delta\colon V_m\rightarrow Y_m$ is continuous.  Denoting this sequence of neighborhoods by ${\mathcal O}(N\cap C,0)$, its tangent is the nested sequence 
$TV_1\supset TV_2\supset TV_2 \ldots $ denoted by $T{\mathcal O}(N\cap C,0)$.  If $x\in V_1$, then the map $D\delta(x)$ is defined, and since  $\delta$ is  a $\ssc^\infty$-germ, the tangent map 
 $T\delta\colon T{\mathcal O}(N\cap C,0)\rightarrow TY$  is again of class $\ssc^1$. Iteratively it follows that $T\delta$ is a $\ssc^\infty$-germ.

\begin{proposition}\index{P- Invariance of tangent germ property}
For a pair $(X,M)$ in which  $X$ is a tame M-polyfold and $M$ a subset of $X$,  the property that $M$ has the n-dimensional tangent germ property, 
is a sc-diffeomorphism invariant. More precisely, if $(X',M')$ is a second pair in which  $X'$ is a tame M-polyfold and $M'$ a subset of $X'$ and if $\psi\colon X\rightarrow X'$ is a sc-diffeomorphism satisfying $\psi(M)=M'$,  then $M'$ has the n-dimensional tangent germ property if and only if
$M$  has the n-dimensional tangent germ property.
\end{proposition}
\begin{proof}
We show that if $M\subset X$ has the n-dimensional tangent germ property, then $M'\subset X'$ has this property  too.
Since $M\subset X_\infty$ and $\psi$ is sc-smooth,  we see that $M'=\psi(M)\subset X'_\infty$. Let $m'\in M'$ and  choose a point  $m\in M$ satisfying $\psi(m)=m'$.
By assumption there exists a sc-smooth chart $\varphi\colon (U,m)\rightarrow (O,0)$, where $(O,C,E)$ is a tame retract. By assumption there exists
a smooth n-dimensional linear subspace $N$ in good position to $C$ with sc-complement $Y$ and a continuous map
$\delta\colon V\rightarrow Y$, where $V$ is a relatively open neighborhood of $0$ in $C_N: =C\cap N$, which satisfies
\begin{itemize}
\item  $\varphi(M\cap U)=\{v+\delta(v)\ |\ v\in V\}\subset N\oplus Y$.
\item  $\delta\colon {\mathcal O}(C_N,0)\rightarrow (Y,0)$ is a sc$^\infty$-germ satisfying $\delta (0)=0$ and $D\delta(0)=0$.
\end{itemize}
To deduce the corresponding construction for $(X',M')$ we define the open neighborhood $U'\subset X'$ of $m'=\psi (m)$ by $U'=\psi (U)$ and the sc-smooth chart by $\varphi'=\varphi\circ \psi^{-1}\colon (U', m')\to (O, 0)$. Then  
$$
\varphi'(U'\cap M') = \varphi(U\cap M) =\{v+\delta(v)\, \vert \, v\in V\}
$$
and the lemma follows.
\qed \end{proof}

The important aspect of the $n$-dimensional tangent germ property of a subset is the following result. 

\begin{theorem}\label{HKL}\index{T- Characterization of subsets with tangent germ property}
If  $X$ is a tame M-polyfold and $M\subset X$ a subset possessing the $n$-dimensional tangent germ property, then the following holds.
\begin{itemize}
\item[{\em (1)}]\ $M$ is a  sub-M-polyfold of $X$ whose  induced M-polyfold structure is tame.  Moreover,  the induced M-polyfold structure on M is sc-smoothly equivalent to a smooth structure  of a manifold with boundary with corners on $M$. 
\item[{\em (2)}]\ If  $x\in M$ is given, we denote by $U$,  $\varphi$, $N$, $V \subset N$, and $\delta\colon V\rightarrow Y$ the data described in condition (2) of Definition \ref{toast}.
Denoting  by $\pi\colon N\oplus Y\rightarrow N$ the sc-projection, the map $U\cap M\rightarrow V$, given by $y\mapsto  \pi\circ \varphi(y)$,  defines a smooth chart  on $M$ around the point $x$. 
\end{itemize}
\end{theorem}

\begin{proof}

We choose a point $x\in M$ and  find an open neighborhood $U\subset X$ of $x$ and a M-polyfold chart
$\varphi\colon (U,x)\rightarrow (O,0)$ onto the retract $(O, C,E)$  so that the set $M\cap U$ is represented as 
$$
\varphi(M\cap U)=\{v+\delta(v)\, \vert \,   v\in V\}\subset N\oplus Y.
$$
The map $\delta\colon V\to Y$ possesses all the properties listed in condition (2) of Definition \ref{toast}.  The map $\varphi\colon M\cap U\to V$  is of the form $\varphi (y)=v(y)+\delta (v(y))$.
With the sc-projection $\pi\colon N\oplus Y\to N$ onto $N$, the map 
\begin{equation}\label{edi}
\pi\circ \varphi\colon M\cap U\to V,\quad \pi\circ \varphi (y)=v(y)
\end{equation}
is continuous  and bijective onto $V$. It is the restriction of the sc-smooth map $\pi\circ \varphi\colon U\to N$, which maps the point $x\in M$ onto $0\in N$. Then the  inverse $\gamma$  of \eqref{edi},
$$
\gamma\colon V\rightarrow M\cap U,\quad  \gamma(v)=\varphi^{-1}(v+\delta(v))
$$
has its image in $X_\infty$ and,  as a map into any $X_m$,  has arbitrarily high regularity if only $v$ is close enough to $0$, depending on $m$. 

Next we shall show that the map $\gamma$ is $\ssc^\infty$ on all of $V$ and not only at the point $0\in V$.

To this aim we choose a $v_0\in V$ and put  $x_0=\varphi^{-1}(v_0+\delta(v_0))$. By construction,  $x_0\in M$,  and by our assumption there is a sc-smooth chart $\psi\colon (U',x_0)\to (O', 0)$ satisfying 
$$\psi(M\cap U')=\{w+\tau (w)\, \vert\, w\in V'\}$$
where $V'$ is a relatively open neighborhood of $0$ in the smooth $n$-dimensional subspace $N'\subset E'$ possessing the tangent germ property. The map $\tau\colon V'\to Y'$ possesses the properties listed in condition (2) of Definition \ref{toast}. In particular, $\tau$ is a $\ssc^\infty$-germ at the point $w=0$.

For $v\in V$ near $v_0$ and $w\in V'$ near $0$ we consider the equation
\begin{equation}\label{edi1}
v=\pi\circ\varphi\circ \psi^{-1}(w+\tau(w)).
\end{equation}

If $v=v_0$, we have the solution $w=0$. Near $0$ the map $\tau$ possesses arbitrary high classical differentiability into any level. Linearizing the right-hand side of the equation at the point $0$, and recalling that $D\tau (0)=0$, we obtain the linear isomorphism 
$$h\mapsto \pi\circ T(\varphi\circ \psi^{-1})(0)h$$
from $N'$ onto $N$.  By the classically implicit function theorem we obtain a germ $v\mapsto w(v)$ for $v$ close to $v_0$ satisfying  $w(v_0)=0$ and solving the equation \eqref{edi1}. The germ has arbitrary high classical differentiability once we are close enough to $v_0$. Now consider the map 
\begin{equation}\label{edi2}
v\mapsto  \varphi\circ \psi^{-1}(w(v)+\tau(w(v))
\end{equation}
for $v$ near $v_0$. Since $v\mapsto w(v)$ has arbitrarily high differentiability at $v_0$ and $\tau$ is $\ssc^\infty$-germ near  $w=0$, we see that the map \eqref{edi2} has, into  any given level,  arbitrarily high differentiability for $v$ near $v_0$.  Consequently, the map is a $\ssc^\infty$-germ near $w=0$. The image of the map  lies in the infinity level. Applying the sc-projection  ${\mathbbm 1}-\pi\colon N\oplus Y\to Y$, we obtain the identity 
$$
\delta(v) =({\mathbbm 1}-\pi)\circ \varphi\circ \psi^{-1}(w(v)+\tau(w(v))
$$
which implies that $\delta$ is a $\ssc^\infty$-germ near $v_0$.
Since $v_0$ is arbitrary in $V$ we see that
$v\mapsto \delta(v)$
is a $\ssc^\infty$-germ around every $v_0\in V$ as we wanted to show. 

Moreover, we conclude that the map
$$
V\rightarrow X,\quad  v\rightarrow \varphi^{-1}(v+\delta(v))
$$
is an injective sc-smooth map whose image is equal to  $M\cap U$.

Next we shall verify that the set $M$ is a sub-M-polyfold of $X$ according to Definition \ref{sc_structure_sub_M_polyfold}.  By construction, we have, so far, at every point $x\in M$ an open neighborhood $U=U(x)\subset X$  and a sc-smooth chart $\varphi\colon U\rightarrow O$,
satisfying  $\varphi(x)=0$, where $(O,C,E)$ is a tame sc-smooth retract. Moreover,  recalling the sc-splitting 
 $
 E=N\oplus Y
 $
 there is  a relatively open neighborhood $V$ of $0$ in the partial quadrant $N\cap C$ of $N$ and a sc-smooth map
 $$
\text{$ \delta \colon V\rightarrow Y$ satisfying $\delta (0)=0$ and
 $D\delta (0)=0$}, 
 $$
 such that  $\varphi(M\cap U) =\{v+\varphi(v)\,  \vert \, v\in V\}$. The map $V\rightarrow U$, 
 $$
 v\mapsto \varphi^{-1}(v+\delta(v)), 
 $$
 is sc-smooth and injective.  The subset $\Sigma\subset C$, defined by 
 $$
 \Sigma=\{v+y\in C\, \vert \, v\in V,\ y\in Y\},
 $$
is relatively open in $C$ and contains $0$. Since $O$ is a tame retract,  there exist a relatively open subset $W$ of $C$ and a tame sc-smooth retraction $r\colon W\rightarrow  W$ onto $O=r(W)$. Consequently, in view of 
 $$
 r\colon r^{-1}(\Sigma\cap O)\rightarrow r^{-1}(\Sigma\cap O),
 $$
the subset $\Sigma\cap O$ is also a sc-smooth retract.

By construction, $v+\delta (v)\in O$ and also $v+\delta (v)\in \Sigma$ and we define the map $t\colon \Sigma\cap O\rightarrow \Sigma\cap O$  by
 $$
 t(v+w)=v+\delta(v).
 $$
The map $t$ is sc-smooth 
 and satisfies $t\circ t=t$, so that  $t$ is a sc-smooth retraction defined on a relatively  open neighborhood of $0$ in $C$ 
 and $t(\Sigma\cap O) = \varphi(M\cap U).$ Therefore, the composition  $s=\varphi^{-1}\circ t\circ\varphi$ defines a sc-smooth retraction
\begin{equation}\label{edi3}
 s\colon U\rightarrow U
 \end{equation}
onto $M\cap U=s(U)$, proving that the subset $M$ is a sc-smooth sub-M-polyfold of $X$.

 The map $u\mapsto C\cap N$, $u\mapsto \pi\circ \varphi (s(u))$ is sc-smooth. Therefore, the map 
 $$M\cap U\to V,\quad m\mapsto \pi\circ \varphi (m)$$
 is a sc-smooth M-polyfold chart on $M$ for the induced M-polyfold structure. The image of the chart is the local model 
 $(V, N\cap C, N)$  so that the transition maps are classically smooth maps and define on $M$  the structure of a manifold with boundary with  corners. The proof of Theorem \ref{HKL} is complete.
 \qed \end{proof}

 As an aside we mention that, in general,  we can not find a local retraction $s$ in \eqref{edi3} which is  tame, as the 
example $X=[0,\infty)^2$ and $M=\{(x,x)\,  \vert , \ x\geq 0\}$ shows.
 
The strength  of the theorem stems  from the fact that in our sc-Fredholm theory the machinery produces subsets $M\subset X$, which have the $n$-dimensional tangent germ property. 
\section{Contraction  Germs}
The notion of a contraction germ is a slight modification of a basic germ. These germs are convenient for the proof of the implicit function theorem (Theorem \ref{newthm5.4}),  which is the main result of this section. It turns out that the local geometry of sc-Fredholm  germs are intimately related to contraction germs. In the generic case they are used to prove that the zero set of a sc-Fredholm section 
must have the n-dimensional tangent germ property. It follows that  the zero set  is  in a natural  way a smooth manifold with boundary with corners.

In the following we abbreviate by $\wt{C}$ the partial quadrant $\wt{C}=[0,\infty)^{k}\oplus \R^{n-k}$ in $\R^n$ so that 
$C=\wt{C}\oplus W$ is a partial quadrant in the sc-Banach space $E=\R^n\oplus W$.

\begin{definition}[$\ssc^0$-{\bf contraction germ} ] \label{BG}\index{D- Contraction germs}
A $\ssc^0$-germ 
$$
f\colon {\mathcal O}(C,0)\rightarrow ( W,0)
$$
  is called  a 
$\ssc^0$-{\bf contraction germ} if the following holds.
The germ $ f$
has the form
$$
 f(a,w)=w-B(a,w)
$$
for $(a,w)$ close to $(0,0)\in C$. Moreover, for every $\varepsilon>0$ and $m\geq 0$,  the estimate
$$
\abs{B(a,w)-B(a,w')}_m\leq \varepsilon\cdot \abs{w-w'}_m
$$
holds for all $(a,w),(a,w')$ on level $m$ sufficiently close to $(0, 0)$, depending on $\varepsilon$ and $m$.
\qed
\end{definition}

More precisely, the  $\ssc^0$-contraction germ requires for given $\varepsilon>0$,  that we can choose a perhaps smaller germ ${\mathcal O}(C,0)$ of neighborhoods 
$U_0\supset U_1\supset U_2\supset U_3\supset \ldots $ of the point $(0, 0)$ in $[0,\infty)^k\oplus \R^{n-k}\oplus W$ such that 
$
\abs{B(a,u)-B(a,v)}_m\leq \varepsilon \abs{u-v}_m
$
holds if $(a, u), (a, v)\in U_m$.

Starting on level $0$, 
the  parametrized version of Banach  fixed point theorem together with  $B(0,0)=0$, guarantee the existence of relatively open and connected neighborhood $V=V_0$ of $0$ in $[0,\infty)^k\oplus \R^{n-k}$  and  
a uniquely determined continuous map $\delta\colon V\rightarrow W_0$ satisfying  $\delta(0)=0$ and solving the equation
$$
\delta(a)=B(a,\delta(a))\quad \text{for all $a\in V$}.
$$
Going to level $1$ we find, again using the fixed point theorem, an 
 open neighborhood $V_1\subset V_0$ of $0$  and a continuous map $\delta_1\colon V_1\to W_1$ satisfying $\delta (0)=0$ and solving 
 the equation on level $1$. 
From the uniqueness of the solutions of the Banach fixed point problem we conclude that $\delta_1 =\delta\vert V_1$. Continuing  this way, we obtain a decreasing sequence 
of relatively open neighborhoods of $0$ in $[0,\infty)^{k}\times \R^{n-k}$, 
$$
V=V_0\supset V_1\supset V_2\supset \ldots 
$$
such that the continuous solution $\delta\colon V\to W$  satisfies 
$\delta(0)=0$ and $\delta (V_m)\subset W_m$ and $\delta\colon V_m\to W_m$ is continuous. 
In  other words,  we obtain  a $\ssc^0$-solution germ $\delta\colon {\mathcal O}([0,\infty)^k\oplus \R^{n-k},0)\rightarrow (W,0)$.

Summarizing  the discussion we have proved the  following theorem from \cite{HWZ3}, Theorem 2.2.

\begin{theorem}[{\bf Existence}]\label{thm5.2}\index{T- Local solution germ}
A $\ssc^0$-contraction germ $f\colon {\mathcal O}(\wt{C}\oplus W,0)\rightarrow (W,0)$ admits a uniquely determined 
sc$^0$-solution germ 
$$
\delta\colon {\mathcal O}(\wt{C},0)\rightarrow (W,0)
$$
solving 
$$
f\circ \gr (\delta)=0.
$$
Here $\gr (\delta)$ is the associated graph germ $a\mapsto(a,\delta(a))$.
\qed
\end{theorem}

Our next aim is the regularity of the unique continuous solution germ $\delta$ of the equation $f(v, \delta (v))=0$ guaranteed by Theorem  \ref{thm5.2}, and we are going to prove that the solution germ $\delta$ is of class $\ssc^k$ if the given germ $f$ is of class $\ssc^k$. 
 By a somewhat tricky induction it turns out that we actually  only have to know that if  $f$ is  $\ssc^1$,  then the solution germ $\delta$ is 
 $\ssc^1$ as well. Here we shall make use of the following regularity result from \cite{HWZ3}, Theorem 2.3, which is the hard part of the regularity theory.

\begin{theorem}\label{thm5.3}\index{T- Regularity of solution germ}
If the $\ssc^0$-contraction germ $f\colon \mo(\wt{C}\oplus  W,0)\rightarrow (
W,0)$ is of class $\ssc^1$, then the solution germ $\delta\colon \mo(\wt{C},
0)\to (W, 0)$ in Theorem \ref{thm5.2} is also of class  $\ssc^1$. 
\qed
\end{theorem}

Theorem \ref{thm5.3} shows that a $\ssc^0$-contraction germ $f$
of class $\ssc^1$ has a solution germ $\delta$ satisfying $f(v, \delta (v))=0$ which is 
also of class $\ssc^1$. We shall use this to verify by induction that $\delta$ is of class $\ssc^k$ if $f$ is of class $\ssc^k$.
We start with the following lemma. 

\begin{lemma}\label{newlemma5.4}\index{L- Higher regularity of solution germ}
Let $f\colon {\mathcal O}(\wt{C}\oplus W,0)\rightarrow (W,0)$ be 
a $\ssc^0$-contraction germ of class $\ssc^k$ where  $k\geq 1$.
Moreover, we assume that the solution germ $\delta$ is of class $\ssc^j$.
(By Theorem \ref{thm5.3}, $\delta$ is  at least of class $\ssc^1$.)
We define the germ $f^{(1)}$ by
\begin{equation*}
f^{(1)}\colon {\mathcal O}(T\wt{C}\oplus TW,0)\rightarrow TW,
\end{equation*}
\begin{equation}\label{germequation}
\begin{split}
f^{(1)}(v,b,u,w)&=\left(u-B(v,u),w-DB(v,\delta(v))\left(b,w\right)\right)\\
&=(u,w)-B^{(1)}(v,b,u,w),
\end{split}
\end{equation}
where the last line defines the map $B^{(1)}$. Then $f^{(1)}$ is an
$\ssc^0$-contraction germ and of class $\ssc^{\min\{k-1,j\}}$.
\end{lemma}
\begin{proof}
For $v$ small, the map $B^{(1)}$ has the contraction property with
respect to $(u, w)$.
 Indeed,  on the $m$-level of $(TW)_m=W_{m+1}\oplus W_m$, i.e.,
 for $(u, w)\in W_{m+1}\oplus W_m$,  we  estimate for given $\varepsilon>0$ and $v$ sufficiently small,
  \begin{equation*}
\begin{split}
&|B^{(1)}(v,b,u',w')-B^{(1)}(v,b,u,w)|_m\\
&\phantom{===}=|B(v,u')-B(v,u) |_{m+1}\\
&\phantom{=====}+
|DB(v,\delta(v))(b,w')-DB(v,\delta(v))(b,w)| _{m}\\
&\phantom{===}\leq \varepsilon |u'-u|_{m+1} + |D_2B(v,\delta(v))[w'-w] |_m,
\end{split}
\end{equation*}
which, using the estimate $\norm{D_2B(v, \delta (v))}_m\leq \varepsilon$ for the operator norm from Lemma \ref{new_Lemma3.9},  is estimated by
\begin{equation*}
\phantom{===}\leq \varepsilon \cdot \bigl(|
u'-u|_{m+1}+|w'-w|_m\bigr)= \varepsilon \cdot | (u', w')-(u,
w)|_m.
\end{equation*}
Consequently,  the  germ $f^{(1)}$ is an  $\ssc^0$-contraction germ.
If  now $f$ is of class $\ssc^k$ and $\delta$ of class $\ssc^j$,
then the germ $f^{(1)}$ is of class $\ssc^{\min\{k-1,j\}}$, as one
verifies by comparing the tangent map $Tf$ with the map $f^{(1)}$
and using the fact that the solution $\delta$ is of class $\ssc^j$.
By Theorem \ref{thm5.2}, the solution germ $\delta^{(1)}$ of
$f^{(1)}$ is at least of class $\ssc^0$. It solves the equation
\begin{equation}\label{chap5eq9}
f^{(1)}(v,b, \delta^{(1)}(v, b))=0.
\end{equation}
But also the tangent germ $T\delta$,  defined by $T\delta
(v,b)=(\delta (v), D\delta (v)b)$,  is a solution of \eqref{chap5eq9}.
From the uniqueness we conclude  that $\delta^{(1)}=T\delta $.
\qed \end{proof}

To prove higher regularity we will also make use of  the next lemma.
\begin{lemma}\label{newcontrlem}
Assume we are given a $\ssc^0$-contraction germ  $f$ of class
$\text{sc}^k$ and a solution germ $\delta$ of class $\ssc^j$
with $j\leq k$. Then  there exists a $\ssc^0$-contraction germ
$f^{(j)}$ of class $\ssc^{\min\{k-j, 1\}}$ having
$\delta^{(j)}:=T^j\delta$ as the solution  germ.
\end{lemma}
\begin{proof}
We prove the lemma by induction with respect to $j$. If $j=0$ and
$f$ is a $\ssc^0$-contraction germ of class $\ssc^k$, $k\geq 0$,
then we  set $f^{(0)}=f$ and $\delta^{(0)}=\delta$.  Hence the
result holds true if $j=0$. Assuming the result has been proved for
$j$,  we show it is true for $j+1$. Since $j+1\geq 1$ and $k\geq
j+1$,   the map $f^{(1)}$,  defined by \eqref{germequation},  is of
class $\ssc^{\min \{k-1, j+1\}}$ in view of Lemma \ref{newlemma5.4}.
Moreover, the solution germ $\delta^{(1)}=T\delta$  satisfies
$$
f^{(1)}\circ \gr (\delta^{(1)})=0,
$$
and is of class $\ssc^j$.  Since $\min \{k-1, j+1\}\geq j$, by  the
induction hypothesis there exists  a map
${(f^{(1)})}^{(j)}=:f^{(j+1)}$ of regularity class
$\min\{\min\{k-1,j+1\}-j,1\}=\min\{k-(j+1),1\}$ so that
$$
f^{(j+1)}\circ\gr ({(\delta^{(1)})}^{(j)})=0.
$$
Setting  $\delta^{(j+1)}={(\delta^{(1)})}^{(j)}=T^{j}(T\delta
)=T^{j+1}\delta$,   Lemma \ref{newcontrlem} follows.
\qed \end{proof}

The main result of this section is the following germ-implicit functions theorem.

\begin{theorem}[{\bf Germ-Implicit Function Theorem}] \label{newthm5.4}\index{T- Germ implicit function theorem}
If 
$$
f\colon \mo(\wt{C}\oplus {W},0)\to  ({ W},0)
$$
 is  an
$\text{sc}^0$-contraction germ which is, in addition, of class
$\ssc^k$, then the solution germ
$$\delta\colon  \mo(\wt{C},0)\rightarrow ({ W},0)$$
satisfying
$$f(v,\delta (v))=0$$
is also of class $\ssc^k$.
\qed
\end{theorem}

From Theorem \ref{newthm5.4}, using  Proposition \ref{sc_up} and \ref{lower}
and Proposition \ref{save}  we deduce the following properties of the section germ $\delta$ under the additional assumptions, that  $f$ is a sc-smooth germ.

\begin{corollary}\label{new_cor_3.33}\index{C- Classical smoothness properties of solution germs}
If the $\ssc^0$-contraction germ $f$ is a sc-smooth germ, there exists for every $m\geq 0$ and $k\geq 0$ a relatively open neighborhoods $V_{m, k}$ of $0$ in $\wt{C}$ such that 
\begin{itemize}
\item[{\em (1)}\ ]\ $\delta (V_{m, k})\subset W_m$.
\item[{\em (2)}\ ]\ $\delta\colon  V_{m, k}\to W_m$ is of class $C^k$.
\end{itemize}
In particular, the solution germ $\delta$ is sc-smooth at the smooth point $0$.
\qed
\end{corollary}

Theorem \ref{newthm5.4} will be one of the 
building blocks for all future versions of implicit function
theorems,  as well as for the transversality theory.

\begin{proof}[Theorem \ref{newthm5.4}]
Arguing by contradiction we assume that the solution germ $\delta$ is
of class $\text{sc}^j$ for some $j<k$ but not of class $\text{sc}^{j+1}$. In view of Lemma \ref{newcontrlem}, there exists an
$\text{sc}^0$-contraction germ $f^{(j)}$ of class $\text{sc}^{\min
\{k-j, 1\}}$ such  that  $\delta^{(j)}=T^j\delta$ satisfies
$$
f^{(j)}\circ \gr (\delta^{(j)})=0.
$$
Since also $k-j\geq 1$, it follows that $f^{(j)}$ is at least of
class $\text{sc}^1$. Consequently,  the solution germ $\delta^{(j)}$
is at least of class $\ssc^{1}$. Since
$\delta^{(j)}=T^j\delta$, we conclude that $\delta$ is at least of
class $\text{sc}^{j+1}$ contradicting our assumption. The proof of
Theorem \ref{newthm5.4} is complete.
\qed \end{proof}

\begin{remark}\index{R- On the properties of solution germs}\label{hofer-rem}
For  later use we reformulate Corollary \ref{new_cor_3.33} in quantitative terms. If $f$ is a $\ssc^0$-contraction germ which, in addition, is a sc-smooth germ satisfying $f(0)=0$. then the solution germ $\delta$ possesses the following properties of existence, uniqueness, and regularity. 

There exist monotone decreasing sequences  $(\varepsilon_i)$ for $i\geq 0$ and $(\tau_i)$ for $i\geq 0$ such that 
\begin{itemize}
\item[(1)]\ $\delta\colon \{a\in [0,\infty)^k\oplus \R^{n-k}\, \vert \, |a|_0\leq \varepsilon_0\}\rightarrow
\{w\in W\ |\ |w|_0\leq \tau_0\}$
is a continuous solution of $f(a, \delta (a))=0$ satisfying $\delta (0)=0$. 
 \item[(2)]\ If the solution  $f(a,w)=0$ satisfies  $|a|_0\leq \varepsilon_0$ and $|w|_0\leq \tau_0$, then $w=\delta(a)$.
\item[(3)]\ If  $|a|_0\leq \varepsilon_i$, then  $\delta(a)\in W_i$ and $|\delta(a)|_i\leq \tau_i$ for every  $i\geq 0$.
\item[(4)]\ The germ  $\delta\colon \{a\in [0,\infty)^k\oplus \R^{n-k} \, \vert \,  |a|_0\leq \varepsilon_i\}\rightarrow W_i$ is of class $C^i$,  for every $i\geq 0$.
\end{itemize}
\end{remark}

\section{Stability  of Basic Germs}

All the maps considered in the section are sc-smooth maps.  
Let us recall (from Definition \ref{BG-00x}) the notion of  a basic germ.
\begin{definition}[{\bf The basic class $\mathfrak{C}_{basic}$}]\index{D- Basic class}\index{$\mathfrak{C}_{basic}$} Let $W$ be a sc-Banach space.  
A {\bf basic germ} $f\colon {\mathcal O}([0,\infty)^k\oplus\R^{n-k}\oplus W,0)\rightarrow (\R^N\oplus W,0)$
is a sc-smooth germ having the property that the germ $P\circ f$ is a $\ssc^0$-contraction  germ, where 
$P\colon \R^N\oplus W\rightarrow W$ is the sc-projection. 
We denote the class of all basic germs by $\mathfrak{C}_{basic}$. 
\end{definition}
In view of Definition \ref{oi}, the basic germs are the local models for the germs of sc-Fredholm sections.

 \begin{theorem}[{\bf Weak Stability of Basic Germs}]\label{arbarello}\index{T- Weak stability of basic germs}
We consider a basic germ 
$$
f\colon {\mathcal O}(([0,\infty)^k\oplus \R^{n-k})\oplus W,0)\rightarrow (\R^N\oplus W,0), 
$$
which we can view as the principal part of a sc-smooth section of the obvious strong bundle. We assume that
$s$ is the principal part of a $\ssc^+$-section of the same bundle satisfying $s(0)=0$.  Then there exists a strong bundle isomorphism
$$
\Phi\colon U\triangleleft (\R^N\oplus W)\rightarrow U'\triangleleft (\R^{N'}\oplus W'),
$$
where $U$ is an open neighborhood of $0$ in $[0,\infty)^k\oplus \R^{n-k}\oplus W$, and $U'$ is an open neighborhood
of $0$ in $[0,\infty)^k\oplus \R^{n'-k}\oplus W'$, 
covering the sc-diffeomorphism $\varphi\colon (U,0)\rightarrow (V,0)$, 
so that ${(\Phi_\ast(f+s))}^1$ is a basic germ. 
Here ${(\Phi_\ast(f+s))}^1$  is the germ $\Phi_\ast(f+s)\colon {\mathcal O}(V^1,0)\to {\mathcal O}((\R^{N'}\oplus W')^1,0)$, where the levels are raised by $1$.
\end{theorem}

Recalling  the Fredholm index of a basic germ in Proposition \ref{Newprop_3.9}, we conclude that $n-N=n'-N'$, because the Fredholm index is invariant under strong bundle isomorphisms. The integer $k$ is the  degeneracy index $k=d_C(0)$ of the point $0$ which is, in view of Proposition \ref{newprop2.24} and 
Corollary \ref{equality_of_d}, preserved under the sc-diffeomorphism $\varphi$ satisfying $\varphi (0)=0$.
Although Theorem \ref{arbarello}  was not explicitly formulated  in \cite{HWZ3}, it follows from the proof of Theorem 3.9  in \cite{HWZ3}.

\begin{proof}
Denoting by  $P\colon \R^N\oplus W\rightarrow W$ the sc-projection, the  composition $P\circ f$ is, by definition, of  the form
$$
 P\circ f(a,w)=w-B(a,w), 
$$
and has the property that for every $\varepsilon>0$ the estimate $\abs{B(a,w)-B(a,w')}_m\leq \varepsilon\cdot |w-w'|_m$ holds,   if $(a, w)$ and $w'$ are sufficiently small on level $m$.

Linearizing the  $\ssc^+$-section  $s$ with respect to the variable $w\in W$ at the point $0$, we introduce the sc-operator 
$$A:=P\circ D_2s(0)\colon W\to W.$$
Since $s$ is a
$\ssc^+$-section and $0$ is smooth point, the operator $A\colon W\rightarrow W$ is a $\ssc^+$-operator. Therefore,  the operator ${\mathbbm 1}+A\colon  W\to W$ is a  $\ssc^+$-perturbation of the identity and hence a sc-Fredholm operator by Proposition \ref{prop1.21}.  Because $A$ is level wise compact, 
the index $\ind ({\mathbbm 1}+A)$ is equal to $0$. 
The associated sc-decompositions of the sc-Banach space $W$ are the following,
$$
{\mathbbm 1}+A\colon W=C\oplus X\to W=R\oplus Z,
$$
where $C=\text{ker}({\mathbbm 1}+A)$ and $R=\text{range}\ ({\mathbbm 1}+A)$ and
$\text{dim}\ C=\text{dim}\ Z<\infty$.

Since $s$ is a $\ssc^+$-section, we conclude from Proposition \ref{lower} that the restriction  $s\colon U_m\to \R^N\oplus W_m$ is of class $C^1$, for every 
$m\geq 1$.  
From the identity 
$P\circ s(a, w)=
P\circ D_2 s(0)w+(P\circ s(a, w)-P\circ D_2 s(0)w)$,
one  deduces the following representation for 
$P\circ s$, on every level $m\geq 1$, 
$$
\text{$P\circ s(a,w)=Aw+S(a, w)$\quad  and \quad 
 $D_2S(0,0)=0$}.$$
Therefore,  $S$ is, with respect to the second
variable $w$, a arbitrary small contraction on every level $m\geq 1$, if $a$ and $w$ are sufficiently small  depending 
on the level $m$ and the contraction constant.  
We can make the arguments which follow only on the levels $m\geq 1$. This explains  the reason for the index  raise by $1$ in the theorem.

We can write
\begin{equation*}
\begin{split}
P\circ (f+s)(a, w)&=w-B(a, w)+Aw+S(a, w)\\
&=({\mathbbm 1}+A)w-[ B(a, w)-S(a, w) ]\\
&=({\mathbbm 1}+A)w-\ov{B}(a, w),
\end{split}
\end{equation*}
where we have abbreviated
$$\ov{B}(a, w)=B(a, w)-S(a, w).$$
By assumption,  the map $B$ belongs to the $\ssc^0$-contraction germ and hence the map $\ov{B}$ is a contraction in the second
variable on every level $m\geq 1$ with arbitrary small contraction
constant $\varepsilon>0$ if $a$ and $w$ are sufficiently small
depending on the level $m$ and the contraction constant $\varepsilon$.  Introducing  the canonical projections by
\begin{align*}
P_1&\colon W=C\oplus X\to X\\
P_2&\colon W=R\oplus Z\to R,
\end{align*}
we abbreviate 
\begin{equation*}
\begin{split}
\varphi (a, w)&:=P_2\circ P\circ (f+s)(a,w)\\
&=P_2  [({\mathbbm 1}+A)w-\ov{B}(a, w)]\\
&=P_2[ ({\mathbbm 1}+A)P_1 w- \ov{B}(a, w)].
\end{split}
\end{equation*}
We have used the relation  $({\mathbbm 1}+A)({\mathbbm 1}-P_1)=0$. The operator
$L:=({\mathbbm 1}+A)\vert X\colon X\to R$ is a sc-isomorphism. In view of
$L^{-1}\circ P_2\circ ({\mathbbm 1}+A)P_1w=P_1w$, we obtain the formula
$$
L^{-1}\circ \varphi (a, w)=P_1w-L^{-1}\circ P_2\circ \ov{B}(a, w).$$

\noindent Writing $w=({\mathbbm 1}-P_1)w\oplus P_1w$,  we  shall consider $(a,
({\mathbbm 1}-P_1)w)$ as our  new finite parameter,  and
correspondingly  define the map $\wh{B}$ by
$$
\wh{B}((a, (1-P_1)w),
P_1w)=L^{-1}\circ P_2\circ  \ov{B}(a, ({\mathbbm 1}-P_1)w+P_1w).
$$
Since $\ov{B}(a, w)$ is a contraction in the second variable on every level $m\geq 1$ with arbitrary small contraction constant if $a$ and $w$ are sufficiently small depending on the level $m$ and the contraction constant,  the right hand side of
$$
L^{-1}\circ \varphi (a, ({\mathbbm 1}-P_1)w+P_1w)=P_1w-\wh{B}(a, ({\mathbbm 1}-P_1)w, P_1w)
$$
possesses the required  contraction normal
form with respect to the variable $P_1w$ on all levels $m\geq 1$,
again if $a$ and $w$ are small enough depending on $m$ and the contraction constant.

 It remains
to prove that the above normal form is the result of an admissible
coordinate transformation of the perturbed section $f+s$. Choosing  a
linear isomorphism $\tau \colon Z\to C$, we  define the fiber
transformation $\Psi\colon \R^N\oplus W\to \R^N\oplus X\oplus C$ by
\begin{equation*}
\begin{split}
\Psi (\delta a\oplus \delta w):= \delta a \oplus L^{-1}\circ P_2 \cdot \delta w \oplus  \tau  \circ ({\mathbbm 1}-P_2)\cdot \delta w.
\end{split}
\end{equation*}
We shall view $\Psi$ as a strong bundle map covering the
sc-diffeomorphism $\psi\colon V\oplus W\to V\oplus C\oplus X$ defined by $\psi (a,
w)=(a, (1-P_1)w, P_1w)$ where $V=[0,\infty )^k\oplus \R^{n-k}.$ With
the canonical projection
\begin{gather*}
\ov{P}\colon  (\R^N\oplus C)\oplus X\to X\\
\ov{P}(a\oplus ({\mathbbm 1}-P_1)w\oplus P_1w)=P_1w,
\end{gather*}
and the relation  $\ov{P}\circ \Psi \circ ({\mathbbm 1}-P)=0$,  we obtain the
desired formula
\begin{gather*}
\ov{P}\circ \Psi \circ(f+s)\circ \psi^{-1} (a, ({\mathbbm 1}-P_1)w, P_1w)\\
=P_1w-\wh{B}(a, ({\mathbbm 1}-P_1)w, P_1w).
\end{gather*}
The proof of Theorem \ref{arbarello}  is complete.

\qed \end{proof}

The theorem has the following corollary, where we use the standard notations, denoting, as usual,  by $C$  a partial quadrant in a sc-Banach space $E$.
We also use a second sc-Banach space $F$.

\begin{corollary}\label{op-perp}\index{C- Weak stability of basic germs}
We assume that the sc-germ  $g\colon {\mathcal O}(C,0)\rightarrow (F,0)$ is equivalent by a strong bundle isomorphism  $\Phi$ to the basic germ $\Phi_\ast g$, and assume that  $s\colon {\mathcal O}(C,0)\rightarrow (F,0)$ is a $\ssc^+$-germ. Then there exists a strong bundle isomorphism $\Psi$ such  that 
 $((\Psi\circ \Phi)_\ast(g+s))^1$ is a basic germ.
\end{corollary}

\begin{proof}
By assumption there exist open neighborhoods $U$ of $0$ in $C$ and $V$ of $0$ in $[0,\infty)^k\oplus \R^{n-k}\oplus W$
and a sc-diffeomorphism $\varphi\colon (U,0)\rightarrow (V,0)$ which is covered by a strong bundle isomorphism
$\Phi\colon U\triangleleft F\rightarrow V\triangleleft(\R^N\oplus W)$ such  that $\Phi_\ast g$ is a basic germ $h$.
Then $t=\Phi_\ast s$ defines a $\ssc^+$-section satisfying  $t(0)=0$. Clearly $\Phi_\ast(g+s)=h+t$,  and applying Theorem \ref{arbarello}, we find a  second strong bundle isomorphism $\Psi$ such  that $(\Psi_\ast(h+t))^1$ is a basic germ.  Taking the composition $\Gamma =\Psi\circ \Phi$, we conclude that 
$\Gamma_\ast (g+s)^1$ is a basic germ.  This completes the proof of Corollary \ref{op-perp}.
\qed \end{proof}

In order to illustrate the corollary, we now consider the sc-smooth germ  
$$
h\colon {\mathcal O}(C,0)\rightarrow F
$$
 for which we know that there exists a $\ssc^+$-germ $s$
satisfying  $s(0)=h(0)$,  and assume that the germ $h-s$ around $0$  is equivalent to the  basic germ $g=\Phi_\ast(h-s)$.
We observe that  $h-h(0)=(h-s) +(s-h(0))$, where  $s-h(0)$ is a $\ssc^+$-section. Then $t=\Phi_\ast(s-h(0))$ is a $\ssc^+$-section and $g+t$ is a perturbation by a $\ssc^+$-section of a basic germ. By the previous corollary we find 
a strong bundle coordinate change such  that $(\Psi_\ast(g-s))^1$ is a basic germ, or in  other words, 
$((\Psi\circ\Phi)_\ast( h-h(0)))^1$ is a basic germ.

Note that for the implicit function theorem it does not matter whether  we work with $f$, or $f^1$, ore even $f^{(501)}$.
It  matters that our coordinate change is compatible with the original sc-structure. 

We also  point out  that a strong bundle coordinate change 
for $h^1$ is not the same as  a strong bundle coordinate change for $h$ followed by a subsequent raise of the index.

\section{Geometry of Basic Germs}
In this section we shall study in detail sc-smooth germs 
\begin{equation}\label{poi1}
f\colon {\mathcal O}(([0,\infty)^k\oplus \R^{n-k})\oplus W,0)\rightarrow (\R^N\oplus W,0)
\end{equation}
around $0$
of the form
\begin{equation}\label{pi2}
f=h+s
\end{equation}
where $h$ is a basic germ and $s$ is a $\ssc^+$-germ satisfying $s(0)=0$. 
We already know from Corollary \ref{Newprop_3.9} that 
$$
Df(0)\colon \R^n\oplus W\to \R^N\oplus W
$$
 is a sc-Fredholm operator of index $\ind Df(0)=n-N$.\par

In the following we abbreviate $E=\R^n\oplus W$,  $C=([0,\infty)^k\oplus \R^{n-k})\oplus W$, and  $F=\R^N\oplus W$  and by  $P\colon \R^N\oplus W\to W$ the sc-projection.

\begin{theorem}[{\bf Local Regularity and Compactness}]\label{save}\index{T- Local regularity and compactness}
Let $U$ be a relatively open neighborhood of $0$ in $C$.  We  assume that $f\colon U\rightarrow F$ is a sc-smooth map satisfying 
$f(0)=0$ and of the form $f=h+s$ where $h$ is a basic germ and $s$ is a  $\ssc^+$-germ  satisfying $s(0)=0$. We denote by $S=\{(a,w)\in U \, \vert \, f(a,w)=0\}$ the solution set of $f$ in $U$.  Then there exists 
a nested sequence 
$$U\supset {\mathcal O}(0)\supset {\mathcal O}(1)\supset {\mathcal O}(2)\supset \ldots $$
of relatively open neighborhoods of $0$ in $C$ on level $0$ such that  for every $m\geq 0$, the closure of $S\cap {\mathcal O}(m)$ in $C\cap E_0=C_0$ is contained in $C\cap E_m=C_m$, i.e., 
$$
\cl_{C_0}(S\cap {\mathcal O}(m))\subset C_m.
$$
\end{theorem}

\begin{remark} \index{R- Comment on solution regularity}
Theorem \ref{save} says, in particular, that 
${\mathcal O}(m)\cap S\subset C_m$ for all $m\geq 0$.  Therefore,  the regularity of solutions $(a, w)$ of the equation $f(a, w)=0$ is the higher, the closer to $0$  they are on the level $0$.  Moreover, the solution set on level $m$ sufficiently close to $0$ on level $0$ has a closure on level $0$, which still belongs to level $m$. 
Moreover, the solution set on level $m$, sufficiently close to $0$ on level $0$, has a closure on level $0$, which belongs to level $m$.
\qed
\end{remark}
\begin{proof}[Proof of Theorem \ref{save}]
We construct the sets ${\mathcal O}(m)$ inductively  by showing that there exists a decreasing sequence $(\tau_m)_{m\geq 0}$ of positive numbers such that the sets 
$$
{\mathcal O}(m)=\{(a, w)\in C\, \vert \, \text{$\abs{a}_0<\tau_m$ and $\abs{w}_0<\tau_m$}\}
$$
have the desired properties. We begin with the construction of ${\mathcal O}(0)$. 

By definition of a basic germ,  the composition $P\circ h$ is of the form 
$$P\circ h(a, w)=w-B(a, w),$$ 
where $B(0, 0)=0$ and $B$ is a contraction in $w$ locally near $(0, 0)$. 
Moreover, $s$ is a $\ssc^+$-germ satisfying $s(0)=0$.

We choose $\tau_0'>0$ such that the closed set $\{(a, w)\in C\, \vert \, \abs{a}_0\leq  \tau_0',\abs{w}_0\leq \tau_0'\}$ in  $E_0$ is contained in $U$  and such that, in addition, 
\begin{itemize}
\item[$(0_1)$\, ]\ \ \ $\abs{B(a, w)-B(a, w')}_0\leq \dfrac{1}{4}\abs{w-w'}_0$.
\end{itemize}
for all $a$, $w$, and $w'\in W_0$ satisfying $\abs{a}_0\leq \tau_0'$, 
$\abs{w}_0\leq \tau_0'$, and $\abs{w'}_0\leq \tau_0'$.
Using $B(0)=0$ and $s(0)=0$,  we can choose $0<\tau_0<\tau_0'$ such that 
\begin{itemize}
\item[$(0_2)$\, ]\ \ \ $\abs{B(a, 0)}_0\leq \dfrac{1}{4}\tau_0'$ for all 
$\abs{a}_0\leq \tau_0.$
\item[$(0_3)$\, ]\ \ \ $\abs{P\circ s(a, w)}_0\leq \dfrac{1}{4}\tau_0'$ for all 
$\abs{a}_0\leq \tau_0$ and $\abs{w}_0\leq \tau_0$.
\end{itemize}
For these choices of the constants $\tau_0'$ and $\tau_0$, we introduce the closed set 
$$\Sigma_0=\{(a, z)\in C\, \vert \, \abs{a}_0\leq \tau_0,\, \abs{z}_0\leq \tau_0'/4\}, $$ and denote by $\ov{B}_0(\tau_0')\subset W_0$ the  closed ball in $W_0$ centered  at $0$ and having radius $\tau_0'$.
We define the map 
$F_0\colon \Sigma_0\times \ov{B}_0(\tau_0')\to  W_0$  by 
$$F_0(a, z, w)=B(a, w)-z.$$
If $(a, z)\in \Sigma_0$ and $w, w'\in \ov{B}_0(\tau_0')$, we estimate  using $(0_1)$ and $(0_2)$, 
\begin{equation*}
\begin{split}
\abs{F_0(a, z, w)}_0&=\abs{B(a, w)-z}_0\leq \abs{B(a, w)-B(a, 0)}_0+\abs{B(a, 0)}_0+\abs{z}_0\\
&\leq \dfrac{1}{4}\tau_0'+\dfrac{1}{4}\tau_0'+\dfrac{1}{4}\tau_0'=\dfrac{3}{4}\tau_0'<\tau_0',
\end{split}
\end{equation*}
and
\begin{equation*}
\abs{F_0(a, z, w)-F_0(a, z, w')}_0\leq \dfrac{1}{4}\abs{w-w'}_0.
\end{equation*}
Hence $F_0(a, z,\cdot )\colon \ov{B}_0(\tau_0')\to \ov{B}_0(\tau_0')$ is a contraction,  uniform in  $(a, z)\in \Sigma_0$. 
Therefore, by the parametrized version of Banach's  fixed point theorem there exists a unique continuous function 
$\delta_0\colon \Sigma_0\to \ov{B}_0(\tau_0')$ solving the 
equation 
$$\delta_0(a, z)=B(a, \delta_0 (a, z))-z$$ 
for all  $(a, z)\in \Sigma_0$.
Now we define the open neighborhood ${\mathcal O}(0)$ by 
$${\mathcal O}(0)=\{(a, w)\in C\, \vert \, \text{$\abs{a}_0<\tau_0$ and $\abs{w}_0<\tau_0$}\}.$$
Clearly, the set ${\mathcal O}(0)$ satisfies $\text{cl}_{C_0}(S\cap {\mathcal O}(0))\subset C_0$. 
We observe that if $(a, w)\in {\mathcal O}(0)$, then $\abs{P\circ s(a, w)}_0\leq {\tau_0'}/4$ by $(0_3)$ so that $\delta_0(a, P\circ s(a, w))$ is defined. If, in addition, $f(a, w)=0$, then  $P\circ f(a, w)=0$  so that $w=B(a, w)-P\circ s(a, w)$ and we claim that  
\begin{equation}\label{new_number_1}
w=\delta_0(a, P\circ s(a, w))\quad \text{for all $(a, w)\in {\mathcal O}(0)$.}
\end{equation} 
Indeed, since $\delta_0(a, P\circ s(a, w))=B(a, \delta_0(a, P\circ s(a, w)))-P\circ s(a, w)$, we estimate, using $(0_1)$,  
\begin{equation*}
\begin{split}
\abs{w-\delta_0(a, P\circ s(a, w))}_0&=\abs{B(a, w)-B(a, \delta_0(a, P\circ s(a, w)))}_0\\
&\leq 
\dfrac{1}{4}\abs{w-\delta_0(a, P\circ s(a, w))}_0,
\end{split}
\end{equation*}
implying $w=\delta_0(a, P\circ s(a, w))$ as claimed.


We next construct the set ${\mathcal O}(1)\subset {\mathcal O}(0)$.  Since the embedding $W_1\to W_0$ is continuous, there is a constant $c_1>0$ such that 
$\abs{\cdot }_0\leq c_1\abs{\cdot}_1$. With the constant $\tau_0$ defined above, we choose $0<\tau_1'<\min \{\tau_0, \tau_0/c_1\}$ such that the following holds. 
The set $\{(a, w)\in C\, \vert \,\text{$\abs{a}_0\leq \tau_1'$ and $\abs{w}_1\leq \tau_1'$}\}$ is contained in $U$,  and 
\begin{itemize}
\item[$(1_1)$\, ]\ \ \ $\abs{B(a, w)-B(a, w')}_1\leq \dfrac{1}{4}\abs{w-w'}_1$.
\end{itemize}
for all $a$, $w, w'\in W_1$ satisfying $\abs{a}_0\leq \tau_1'$, 
$\abs{w}_1\leq \tau_1'$, and $\abs{w'}_1\leq \tau_1'$.
We choose $0<\tau_1<\tau_1'$ such that 
\begin{itemize}
\item[$(1_2)$\, ]\ \ \ $\abs{B(a, 0)}_1\leq \dfrac{1}{4}\tau_1'$ for all 
$\abs{a}_0\leq \tau_1$
\item[$(1_3)$\, ]\  \ \ $\abs{P\circ s(a, w)}_1\leq \dfrac{1}{4}\tau_1'$ for all 
$\abs{a}_0\leq \tau_1$ and $\abs{w}_0\leq \tau_1$.
\end{itemize}
Proceeding as in the construction of ${\mathcal O}(0)$, we introduce the closed set $\Sigma_1$ in $E_1$ by 
$$\Sigma_1=\{(a, z)\in C_1\,\vert \, \abs{a}_0\leq  \tau_1,\, \abs{z}_1\leq \tau_1'/4\}, $$
and abbreviate by $\ov{B}_1(\tau_1')$ the closed ball in $W_1$ having its center at $0$ and radius $\tau_1'$.  We define the map 
$F_1\colon \Sigma_1\times \ov{B}_1(\tau_1')\to W_1$ by 
$F_1(a, z, w)=B(a, w)-z$. By  $(1_1)$ and $(1_3)$, the map 
$F_1\colon \Sigma_1\times \ov{B}_1(\tau_1')\to \ov{B}_1(\tau_1')$ is a contraction, uniform in $(a, z)\in \Sigma_1$. 
Again using the Banach fixed point theorem, we find a unique continuous map
$\delta_1\colon \Sigma_1\to \ov{B}_1(\tau_1')$  solving the equation $\delta_1(a, z)=B(a, \delta_1(a, z))-z$ for all $(a, z)\in \Sigma_1$.  
Now we define the open neighborhood ${\mathcal O}(1)$ as 
$${\mathcal O}(1)=\{(a, w)\in C_0\, \vert \, \abs{a}_0<\tau_1,\, \abs{w}_0< \tau_1\}.$$
By our definition of $\tau_1$ we have $\tau_1\leq \tau_1'<\tau_0$ so that  
${\mathcal O}(1)\subset {\mathcal O}(0)$.

We next claim that 
\begin{equation}\label{delta_0_equals_delta_1}
\delta_0(a, P\circ s(a, w))=\delta_1(a, P\circ s(a, w))\quad \text{for all $(a, w)\in {\mathcal O}(1)$}.
\end{equation}
To verify the claim,  we note that if $(a, w)\in {\mathcal O}(1)$, then, by $(1_3)$,  $\abs{P\circ s(a, w)}_1\leq \tau'_1/4$. Hence $\delta_1(a,  P\circ s(a, w))$ is defined  and its norm  satisfies $\abs{\delta_1(a, P\circ s(a, w))}_1\leq \tau_1'$ because $\delta_1$ takes its values in the ball $\ov{B}_1(\tau_1')$. 
This implies, recalling that $\abs{\cdot}_0\leq c_1\abs{\cdot}_1$ and 
$\tau_1'\leq \tau_0/c_1$,  the estimate
$$\abs{\delta_1(a, P\circ s(a, w))}_0\leq c_1\abs{\delta_1(a, P\circ s(a, w))}_1\leq c_1\tau_1'\leq \tau_0\leq \tau'_0.$$
Therefore, by construction, the map 
$(a, w)\mapsto \delta_1(a, P\circ s(a, w))$ solves the equation $\delta_1(a, P\circ s(a, w))=B(a, \delta_1(a, P\circ s(a, w))-P\circ s(a, w)$ for all $(a, w)\in {\mathcal O}(1)$.  On the other hand, it follows from ${\mathcal O}(1)\subset {\mathcal O}(0)$ and $(0_3)$ that $\abs{P\circ s (a, w)}_0\leq \tau_0'/4$ and hence the map $(a, w)\mapsto \delta_0(a, P\circ s(a, w))$ solves, by construction, the same equation $\delta_0(a, P\circ s(a, w))=B(a, \delta_0(a, P\circ s(a, w))-P\circ s(a, w)$ for all $(a, w)\in {\mathcal O}(1)$. The claim \eqref{delta_0_equals_delta_1} now follows from the uniqueness of the Banach fixed point theorem on the level $0$.

If $(a, w)\in {\mathcal O}(1)$ satisfies, in addition, $f(a, w)=0$, we deduce from \eqref{new_number_1} and 
\eqref{delta_0_equals_delta_1} that 
\begin{equation}\label{w_equals_delta_1}
w=\delta_1(a, P\circ s(a, w))\in W_1\quad \text{for all $(a, w)\in {\mathcal O}(1)$}.
\end{equation}

In order to verify the desired  property of ${\mathcal O}(1)$ we fix 
$(a, w)\in \text{cl}_0(S\cap {\mathcal O}(1))$. Then there exists a sequence  
$(a_n, w_n)\in S\cap {\mathcal O}(1)$ such that $(a_n, w_n)\to (a, w)$ on level $0$. 
From \eqref{w_equals_delta_1} it follows that  $w_n=\delta_1(a_n, P\circ s(a_n, w_n))\in W_1$  for all $n$.   Since $s$ is $\ssc^+$,  we know that $P\circ s(a_n, w_n)\to P\circ s(a, w)$ on level $1$. From the continuity of $\delta_1$ we conclude the convergence 
$w_n=\delta_1(a_n, P\circ s(a_n, w_n))\to \delta_1(a, P\circ s(a, w))=w$ on level $1$.  Consequently, $(a, w)\in C_1$ as desired.

The induction step is now clear and the further details are left to the reader.

\qed \end{proof}

The previous result has a useful corollary.

\begin{corollary}\label{corex1}\index{C- Stability of surjectivity}
We assume that $U$ is a relatively open neighborhood of $0$ in a partial quadrant $C=[0,\infty)^n\oplus W$ in a sc-Banach space $E=\R^n\oplus W$ and let $F=\R^N\oplus W$.  Let $f\colon U\to F$ be a sc-smooth map satisfying $f(0)=0$ and  admitting the decomposition
$f=h+s$ where $h\in \mathfrak{C}_{basic}$ and $s$  is a $\ssc^+$-map satisfying $s(0)=0$. We assume, in addition, that 
$Df(0,0)$ is surjective.  Then there exists a 
 relatively open neighborhood $U'\subset U$ on  level $0$ such that the following holds.
\begin{itemize} 
\item[{\em (1)}]\ If $(a,w)\in U'$ satisfies  $f(a,w)=0$, then $(a,w)$ is on level $1$.
\item[{\em (2)}]\ If $(a,w)\in U'$ and $f(a, w)=0$, then $Df(a,w)\colon E\rightarrow F$ is a surjective Fredholm operator of index $n-N$.
\end{itemize}
\end{corollary}

\begin{proof}
In view of Theorem \ref{save}  we know that if $f(a,w)=0$ and $(a,w)$ is sufficiently close to $(0,0)$ on level $0$, then $(a, w)\in E_1$. Hence the linearization $Df(a,w)$ is a well-defined as a bounded linear operator from $E_0$ to $F_0$. By Proposition \ref{Newprop_3.9}, the linearization $Df(0,0)\colon E\to F$ is a  Fredholm operator whose index is equal to $\ind Df(0)=n-N$. 

By assumption, the Fredholm operator $Df(0,0)\colon E\to F$ is surjective. Denoting by $K$ its kernel, we have the splitting $E=K\oplus N$ and conclude that the restriction $Df(0)\vert N\colon N\to F$ is an isomorphism of 
Banach spaces and hence a Fredholm operator of index $0$. To see this, we observe that 
$$
Df(a, w)\vert N=Df(0, 0)\vert N+(Df(a, w)\vert N-Df(0, 0)\vert N), 
$$
and the  second term is a admissible perturbation which does not affect the Fredholm character nor the index. Indeed, 
\begin{equation*}
\begin{split}
\bigl( Df(a, w)-Df(0, 0)\bigr)(\alpha, \zeta)=D_2B(a, w)\zeta+C(a, w)(\alpha, \zeta),
\end{split}
\end{equation*}
where $C(a, w)\colon E\to F$ is a compact operator, and 
$\norm{D_2B(a, w)}\leq \varepsilon$ for every $\varepsilon>0$ if $(a, w)$ sufficiently small in $E_0$ depending on $\varepsilon$,  in view of 
Lemma \ref{new_Lemma3.9}.

It follows that $Df(a, w)\vert N\colon N\to F$ is a Fredholm operator of index $0$, if $(a, w)\in E_1$ is sufficiently small in $E_0$.
\par

We finally show that $Df(a, w)\vert N\colon N\to F$ is a surjective operator if $(a, w)\in E_1$ solves, in addition, $f(a, w)=0$.
Since the index is equal to $0$, the kernel of $Df(a, w)\vert N$ has the same dimension as the cokernel of $Df(a, w)\vert N$ in $F$. Hence we have to prove that the kernel of $Df(a, w)\vert N$ is equal to $\{0\}$, if $(a, w)\in E$ is close to $(0, 0)$ in $E_0$ and solves $f(a, w)=0$.

Arguing by contradiction we assume that  there exist a sequence $(a_k, w_k)\in C\cap E_1$ satisfying $f(a_k, w_k)=0$ and $\abs{(a_k, w_k)}_0\to 0$ as $k\to \infty$. Moreover, there exists a sequence 
$(\alpha_k, \zeta_k)\in (\R^n\oplus W)\cap N$  satisfying $\abs{(\alpha_k, \zeta_k)}_0=1$ and $Df(a_k, w_k)(\alpha_k, \zeta_k)=0$.  Consequently, 
\begin{equation}\label{eq_surj_1}
 \zeta_k-D_2B(a_k,w_k) \zeta_k = D_1B(a_k,w_k)\alpha_k - PDs(a_k,w_k)(\alpha_k,  \zeta_k).
\end{equation}
Without loss of generality we  may assume that $\alpha_k\rightarrow \alpha$.  In view of the proof of the previous theorem, for large values of $k$, the sequence $(a_k, w_k)$ is bounded on level $2$. Consequently, since the embedding $W_2\to W_1$ is compact, we may assume that $(a_k, w_k)\to (0, 0)$ on level $1$. Therefore, 
$D_1B(a_k,w_k)\alpha_k \rightarrow D_1B(0,0)\alpha$ in $E_0$. 
In addition, since $s$ is a $\ssc^+$-operator, the map $E_1\oplus E_0\to F_1$,  defined  by $(x, h)\mapsto PDs(x)h$,  is continuous. Hence,  there exists $\rho>0$ such that $\abs{PDs(x)h}_1\leq 1$ if $\abs{x}_1\leq \rho$ and $\abs{h}_0\leq \rho$. This implies that there is a constant $c>0$ such that $\abs{PDs(x)h}_1\leq c$ for all $\abs{x}_1\leq \rho$ and $\abs{h}_0\leq 1$.  From this estimate, we conclude that the sequence 
$PDs(a_k,w_k)(\alpha_k, \zeta_k)$ is bounded in $W_1$. Since the embedding $W_1\to W_0$ is compact, we may assume that the sequence $PDs(a_k,w_k)(\alpha_k, \zeta_k)$ converges to some point $w_0$  in $W_0$. Denoting by $z_k$ the right-hand side of \eqref{eq_surj_1}, we have proved that $z_k\to z_0=D_1B(0,0)\alpha-w_0$ in $W_0$.  Choosing $0<\varepsilon<1$ in Lemma \ref{new_Lemma3.9}, the  operators ${\mathbbm 1}-D_2B(a_k, w_k)$ have a bounded inverse for large $k$ and we obtain from \eqref{eq_surj_1} that 
\begin{equation*}
\zeta_k=\bigl({\mathbbm 1}-D_2B(a_k, w_k)\bigr)^{-1}z_k=\sum_{l\geq 0}\bigl(D_2B(a_k, w_k)\bigr)^{l}z_k
\end{equation*}
for large $k$. We claim that the sequence $\zeta_k$ converges to $z_0$ in $W_0$. Indeed, take  $\rho>0$. From $\norm{D_2B(a_k, w_k)}\leq \varepsilon$ for $k$ large and the fact that the sequence  $(z_k)$ is bounded in $W_0$, it follows that there exists $l_0$ such that 
$$\abs{\sum_{l\geq l_0}\bigl(D_2B(a_k, w_k)\bigr)^{l}z_k}_0\leq \rho/2.$$
Hence 
\begin{equation*}
\begin{split}
\abs{\zeta_k-z_0}_0&\leq \abs{z_k-z_0}_0+\sum_{l=1}^{l_0}\abs{\bigl(D_2B(a_k, w_k)\bigr)^{l}z_k}_0+\abs{\sum_{l>l_0}\bigl(D_2B(a_k, w_k)\bigr)^{l}z_k}_0\\
&\leq \abs{z_k-z_0}_0+\sum_{l=1}^{l_0}\abs{\bigl(D_2B(a_k, w_k)\bigr)^{l}z_k}_0+\rho/2.
\end{split}
\end{equation*}
From  $\abs{z_k-z_0}_0\to 0$ and $D_2B(a_k, w_k)z_k\to D_2B(0, 0)z_0=0$, it follows that 
$$
\limsup_{k\to \infty}\abs{\zeta_k-z_0}_0\leq \rho/2
$$
 and,  since $\rho$ was arbitrary,  that $\zeta_k\to z_0$ in $W_0$.
We have proved that $(\alpha_k, \zeta_k)\to (\alpha, z_0)$ in $E_0$.  Hence $(\alpha, z_0)\in N$ and 
$\abs{(\alpha, z_0)}_0=1$. On the other hand,
\begin{equation*}
0=\lim_{k\to \infty}Df(a_k, w_k)(\alpha_k, \zeta_k)=Df(0)(\alpha, z_0).
\end{equation*}
Since $Df(0, 0)\vert N$ is an isomorphism, $(\alpha, z_0)=(0, 0)$, in  contradiction to $\abs{(\alpha, z_0)}_0=1$. The proof of 
Corollary \ref{corex1} is complete.
\qed \end{proof}

In the following theorem we denote, as usual, by $C$ the  partial quadrant $[0,\infty )^k\oplus \R^{n-k}\oplus W$ in the  sc-Banach space $E=\R^n\oplus W$ and by $U$ a relatively open neighborhood of $0$ in $C$. Moreover, $F$ is another sc-Banach space of the form $F=\R^N\oplus W$. 

We consider a sc-smooth germ $f\colon U\to F$ satisfying $f(0)=0$ of the form $f=h+s$ where $h$ is a basic germ and $s$ is a $\ssc^+$-germ satisfying $s(0)=0$.
By Theorem \ref{arbarello} there exists a strong bundle isomorphism 
$$\Phi\colon U\triangleleft (\R^N\oplus W)\to U'\triangleleft (\R^{N'}\oplus W').
$$
where $U'$ is a relatively open neighborhood of $0$ in the partial quadrant $C'=[0,\infty)^k\oplus \R^{n'-k}\oplus W'$, covering the sc-diffeomorphism 
$\varphi\colon (U, 0)\to (U',0)$ such that the section $g=\Phi\circ f\circ \varphi^{-1}$ has the property that $g^1=(\Phi_\ast (h+s))^1\colon (U')^1\to (\R^{N'}\oplus W')^1$ is a basic germ. Clearly, $\varphi \bigl(\{x\in U\, \vert \, f(x)=0\}\bigr) =\{x'\in U'\, \vert \, g(x')=0\}.$

We abbreviate $F'=\R^{N'}\oplus W'$ and denote by $P'$ the sc-projection $P'\colon \R^{N'}\oplus W'\to W'$.
By definition of a basic germ, the composition $P'\circ g^1$ is a $\ssc^0$-contraction germ which is sc-smooth. Therefore, in view of Remark \ref{hofer-rem} applied to the levels $m\geq 1$,
there are monotone decreasing sequences $(\varepsilon_i')$ for $i\geq 1$ and 
$(\tau_i')$ for $i\geq i$ such that, abbreviating by $\ov{B}_i(\tau_i')$ the closed ball $W'_i$ of center $0$ and radius $\tau_i'$, and by $U'_i$ the neighborhood $U_i'=\{a\in [0,\infty )^k\oplus \R^{n'-k}\, \vert \, \abs{a}_0\leq \tau_i'\}$, the following statements (1)-(4) hold. 
\begin{itemize}
\item[(1)]\ There exists a unique continuous map 
$$
\delta\colon U'_1\rightarrow \ov{B}_1(\tau_1') 
$$
  satisfying  $P'\circ g(a,\delta(a))=0$ and $\delta(0)=0$.
\item[(2)]\ If $(a, w)\in U_1'\oplus \ov{B}_1(\tau_1')$ solves the equation $P'\circ g(a,w)=0$, then $w=\delta(a)$.
\item[(3)]\ If $a\in U_i'$, then $\delta (a)\in W_{i}'$ and $|\delta(a)|_{i}\leq \tau_i'$ for every 
$i\geq 1$.
\item[(4)]\ $\delta\colon U_i'\to W_i'$ is of class $C^{i-1}$ for all $i\geq 1$.
\end{itemize}

\begin{lemma}\label{new_lemmat_level_0}
There exists $\varepsilon_0'$ and $\tau_0'$ having the following properties. If $(a', w')$ solves the equation $g(a', w')=0$ and satisfies $\abs{a'}_0\leq \varepsilon'_0$ and $\abs{w'}_0\leq \tau_0'$, then $w'=\delta (a')$.
Moreover, if $\abs{a'}_0\leq \varepsilon'_0$, then $\abs{\delta (a')}_0\leq \tau_0'$.
\end{lemma}

\begin{proof}

The solutions of $f(a)=0$ and $g(x')=0$ are related by the sc-diffeomorphism 
$\varphi$ via $x'=\varphi (x)$. Applying the proof of Theorem \ref{save} to the sc-germ $f=h+s$ we find,  for every $\sigma>0$,  a constant $\tau>0$ such that if $x\in C$ is a solution of $f(x)=0$ satisfying $\abs{x}_1<\tau$, then $\abs{x}_1<\sigma$. Using the continuity of $\varphi$ on the level $1$, we choose now  $\sigma>0$ such that if $\abs{x}_1<\sigma$, then $x'=\varphi (x)=(a', w')$ satisfies $\abs{a'}_0\leq \varepsilon_1'$ and $\abs{w'}_1\leq \tau_1'$. Then we conclude from $\abs{x}_0<\tau$ that $\abs{x}_1<\sigma$. Using the continuity of $\varphi^{-1}$ on level $0$ we next choose the desired constants $\varepsilon_0'$ and $\tau_0'$ such that the estimates $\abs{a'}_0\leq \varepsilon_0'$ and $\abs{w'}_0\leq \tau_0'$ imply that $x=\varphi^{-1}(a, w)$ satisfies $\abs{x}_0<\tau$. Assuming now that  $x'=(a', w')$ is a solution of $g(a', w')=0$ satisfying $\abs{a'}_0<\varepsilon_0'$ and $\abs{w'}_0< \tau_0'$, we conclude that $\abs{\varphi^{-1}(x')}_0<\tau$. It follows that 
$\abs{\varphi^{-1}(x')}_1<\sigma$, which implies $\abs{a'}_0<\varepsilon_1'$ and 
$\abs{w'}_1< \tau_1'$. From property (2) we conclude that $w'=\delta (a')$ proving the first statement of the lemma.

Using that the embedding $W_1\to W_0$ is continuous, we find a constant $c>0$ such that 
$\abs{\delta (a')}_0\leq c\abs{\delta (a')}_1$. Taking $\varepsilon_0'$ and $\tau_0'$ smaller we can achieve that 
 $\abs{\delta (a')}_0\leq \tau_0'$ if $\abs{a'}_0\leq \varepsilon_0'$ and the proof of the lemma is complete.
 
 \qed \end{proof}
 
\begin{theorem}[{\bf Local Germ-Solvability I}]\label{LGS}\index{T- Local germ-solvability {I}}
Assume that  $f\colon U\rightarrow F$ is  a sc-smooth  germ  satisfying $f(0)=0$ and of the form $f=h+s$,  where $h$ is a basic germ 
and $s$ is a $\ssc^+$-section satisfying $s(0)=0$. We assume 
that the linearization $Df(0)\colon E\to F$ is surjective and the  kernel $K=\ker Df(0)$ is in good position to the partial quadrant $C$. 
Let 
$$\Phi\colon U\triangleleft (\R^N\oplus W)\to U'\triangleleft (\R^{N'}\oplus W')$$
be the strong bundle isomorphism covering the sc-diffeomorphism $\varphi\colon (U,0)\rightarrow (U',0)$ guaranteed by Theorem \ref{arbarello}. 
Here $U'$ is a relatively open neighborhood of $0$ 
in the  partial quadrant $C'=[0,\infty )^k\oplus \R^{n'-k}\oplus W'$ sitting in   the sc-Banach space $E'=\R^{n'}\oplus W'$.   

Then, denoting by $g=\Phi_\ast (f)$ the push-forward section, the following holds.
The kernel  $K'=\ker Dg(0)=T\varphi(0)K$  is in good position to the partial quadrant $C'$ and there is a good complement $Y'$  of $K'$ in $E'=K'\oplus Y'$,  and a $C^1$-map $\tau\colon V\rightarrow Y_1'$, defined on the  relatively  open neighborhood $V$ 
of $0$ in $K'\cap C'$  such that 
\begin{itemize}
\item[{\em (1)}]\ $\tau\colon {\mathcal O}(K'\cap C',0)\rightarrow (Y',0)$ is a sc-smooth germ.
\item[{\em (2)}]\ $\tau (0)=0$ and $D\tau(0)=0$.
\item[{\em (3)}]\ After perhaps suitably shrinking $U$ it holds
$$\varphi(\{x\in U\ |\ f(x)=0\}) =
\{x'\in U'\, \vert \, g(x')=0\}=\{y+\tau(y)\,  \vert \,  y\in V\}.
$$
\end{itemize}

\end{theorem}

\begin{proof}
The linearization $Dg(0)\colon \R^{n'}\oplus W'\to \R^{N'}\oplus W'$ is a surjective map because, by assumption, $Df(0)\colon \R^{n}\oplus W\to \R^{N}\oplus W$ is surjective. Moreover, $\ker Dg(0)=T\varphi (0)\ker Df(0)$. In view of Proposition \ref{newprop2.24} and Corollary \ref{equality_of_d} we  have $T\varphi (0)C=C'$.
We also  recall that (above Lemma \ref{new_lemmat_level_0}) we have introduced the sc-projection $P'\colon \R^{N'}\oplus W'\to W'$, the neighborhood $U'_i=\{a\in [0,\infty)^k\oplus \R^{n'-k}\, \vert \, \abs{a}_0\leq \tau'_i\}$ of $0$ in $\R^{n'}$, and the map $\delta\colon U_i'\to W_i'$ satisfying $P'\circ g(a, \delta (a))=0$ and $\delta (a)=0$.

We introduce the map $H\colon U_1'\to \R^{N'}$ by 
$$H(a)=({\mathbbm 1}-P')\circ g(a, \delta (a)).$$
It satisfies $H(0)=0$ since $g(0)=0$.  In the following we need the map $H$ to be of class $C^1$. In view of property (4), the map $H$ restricted to $U_i'$ is of class $C^{i-1}$. Since, by Proposition \ref{lower}, the sc-smooth map $g$ induces a $C^1$-map $g\colon E_{m+1}\to E_m$ for every $m\geq 0$, we shall restrict the domain of $H$ to the set $U'_3$. In order to prove the theorem we shall first relate the solution set $\{H=0\}$ to the solution set $\{g=0\}$, and start with the following lemma.

\begin{lemma}\label{new_lemma_relation}
The linear map $\alpha \mapsto (\alpha , D\delta (0)\alpha)$ from $\R^{n'}$ into $E'$ induces a linear isomorphism 
$$\wt{K}:=\ker DH(0)\to \ker Dg(0)=:K'.$$
\end{lemma}

\begin{proof}[Proof of the lemma]
We first claim that  
$$\ker Dg(0)=\{(\alpha ,D\delta (0)\alpha)\, \vert \, \alpha \in \ker DH(0)\}.$$  Indeed, if $(\alpha, \zeta)\in \ker Dg(0)$, then  $P'\circ Dg(0)(\alpha, \zeta)=0$ and  $({\mathbbm 1}-P')\circ Dg(0)(\alpha, \zeta)=0$. By assumption, $P'\circ g(a, w)=w-B(a, w)$. Hence,   differentiating at the point $0$ in the direction of $(\alpha, \zeta)$ and  recalling from Lemma \ref{new_Lemma3.9} that $D_2B(0)=0$, we obtain 
$$0=P'\circ Dg(0)(\alpha, \zeta)= \zeta-D_1B(0)\alpha-D_2B(0)\zeta= \zeta-D_1B(0)\alpha.$$
 On the other hand, differentiating  the identity $P'\circ g(a, \delta (a))=0$ at $a=0$ and evaluating the derivative at $\alpha$, we find,  in view of previous equation,  
$$0=P'\circ g(0)(\alpha, D\delta (0)\alpha)=D\delta (0)\alpha-D_1B(0)\alpha.$$
Consequently,  $\zeta=D\delta (0)\alpha$ as claimed.

From  
$$DH(0)\alpha =({\mathbbm 1}-P')\circ Dg(0)(\alpha , D\delta (0)\alpha )$$
 it follows, if $\alpha \in \ker DH(0)$, in view of 
 $P'\circ Dg(0)(\alpha ,D\delta (0)\alpha)=0$ for all $\alpha$, that 
 $(\alpha, D\delta (0)\alpha)\in \ker Dg(0)$.
 
Conversely, if $(\alpha, \zeta)\in \ker Dg(0)$, then $0=P'\circ Dg(0)(\alpha, \zeta)=
\zeta-D_1B(0)\alpha$ and $({\mathbbm 1}-P')\circ Dg(0)(\alpha, \zeta)=0$. Since 
$P'\circ Dg(0)(\alpha, D\delta (0)\alpha)=0$, we conclude that $\zeta=D\delta (0)\alpha$ and hence $({\mathbbm 1}-P')\circ Dg(0)(\alpha, D\delta (0)\alpha)=DH(0)\alpha=0$, so that $\alpha \in \ker DH(0)$. The proof of 
Lemma \ref{new_lemma_relation} is complete.

\qed \end{proof}

By assumption, the kernel  $K=\ker Df(0)$ is in good position to the partial quadrant $C$. We recall that this requires  that 
$K\cap C$ has a nonempty interior in $K$ and there exists a complement of $K$, denoted by $K^\perp$, so that $K\oplus K^\perp=E$, having the following property. There exists 
$\varepsilon>0$ such that 
if $k+k^\perp\in K\oplus K^\perp$  satisfies  $\abs{k^\perp}_0\leq \varepsilon\abs{k}_0$, then 
\begin{equation}\label{eq_good_position}
\text{$k+k^\perp\in C$ \quad if and only if \quad $k\in C$}.
\end{equation}

The kernel $K'=\ker Dg(0)=T\varphi (0)K$ has a complement in $E'$,  denoted by $(K')^\perp=T\varphi (0)(K^\perp),$ so that 
$K'\oplus (K')^\perp =E'$.
\begin{lemma} 
The complement $(K')^\perp$ is a good complement of $K'$ in $E'$. 
\end{lemma}
\begin{proof}
Since $T\varphi (0)\colon \R^n\oplus W\to \R^{n'}\oplus W'$ is a topological isomorphism, $T\varphi (0)C=C'$ and $K'=T\varphi (0)K$, it follows that $K'\cap C'$ has a nonempty interior in $K'$. Next we choose 
$\varepsilon'>0$ satisfying $\varepsilon'\norm{T\varphi (0)}\norm{(T\varphi (0))^{-1}}\leq \varepsilon$ and take 
$x'\in K'$ and $y'\in (K')^\perp$ satisfying $\abs{y'}_0\leq \varepsilon'\abs{x'}_0$. 
If  $x=(T\varphi (0))^{-1}x'$ and $y=(T\varphi (0))^{-1}y'$, then
\begin{equation*}
\begin{split}
\abs{y}_0&
=\abs{(T\varphi (0))^{-1}y'}_0\leq \norm{(T\varphi (0))^{-1}}\abs{y'}_0\leq \varepsilon' \norm{(T\varphi (0))^{-1}}\abs{x'}_0\\
&=\varepsilon' \norm{(T\varphi (0))^{-1}}\abs{T\varphi (0)}_0\leq 
\varepsilon' \norm{(T\varphi (0))^{-1}}\norm{T\varphi (0)}\abs{x}_0\leq \varepsilon\abs{x}_0.
\end{split}
\end{equation*}
By \eqref{eq_good_position}, $x+y\in C$ if and only if $x\in C$ and since  $T\varphi (0)C=C'$,  we conclude that $x'+y'\in C'$ if and only if $x'\in C'$.  
\qed \end{proof}
The next lemma is proved in Appendix \ref{pretzel-B}.
\begin{lemma}\label{new_lemma_Z}
The kernel  $\wt{K}:=\ker DH(0)\subset \R^{n'}$ of the linearization $DH(0)$ is in good position to the partial quadrant $\wt{C}=[0,\infty)^{k}\oplus \R^{n'-k}$ in $\R^{n'}$. Moreover, there exists a good complement $Z$ of $\wt{K}$ in $\R^{n'}$, hence 
$\wt{K}\oplus Z =\R^{n'}$, having the property that $Z\oplus W'$ is a good complement of $K'=\ker Dg(0)$ in $E'$, 
$$E'=K'\oplus ( Z\oplus W').$$
\qed
\end{lemma}
 In view of  Lemma \ref{new_lemma_Z},  Theorem \ref{help-you} in Appendix \ref{implicit_finite_partial_quadrants} is applicable to the $C^1$-map $H\colon U_2'\to \R^{N'}$,  where 
$U_2'=\{a\in [0,\infty )^k\oplus \R^{n'-k}\, \vert \, \abs{a}_0<\tau_2'\}\subset \R^{n'}$ is introduced before the statement of Theorem \ref{LGS}. According to Theorem 
\ref{help-you}  there exists a relatively open neighborhood $\wt{V}$ of $0$ in $\wt{K}\cap \wt{C}$ and a $C^1$-map 
$$\sigma\colon \wt{V}\to Z$$
satisfying $\sigma(0)=0$ and $D\sigma(0)=0$ and solving the equation  
$$
\text{$H(a+\sigma(a))=0$ for all $a\in \wt{V}$.}
$$   
The situation is illustrated in the following Figure \ref{fig:pict5}.
\begin{figure}[htb]
\begin{centering}
\def\svgwidth{60ex}
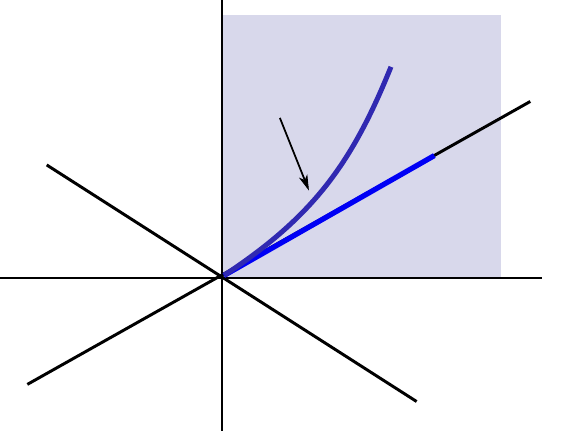
\caption{}\label{fig:pict5}
\end{centering}
\end{figure}
Since $H$ is the more regular the closer we are  to $0$,  the map  $\tau$ has the same property.   Recalling the solution germ 
$\delta\colon U_1'\subset \R^{n'}\to \ov{B}_1' \subset W_1'$ of 
$g(a, \delta (a))=0$, we define the map
$\gamma\colon \wt{V}\rightarrow Z\oplus W'_1$ by 
$$\gamma (a)=\sigma(a)+\delta(a+\sigma(a)),\quad v\in \wt{V}.
$$
By the properties of $\sigma$ and $\delta$ 
there exists a nested sequence $\wt{V}_1\supset \wt{V}_2\supset \ldots $ of relatively open neighborhoods of $0$ in $\wt{K}\cap \wt{C}$ such that the restriction
$$
\gamma \vert \wt{V}_i\colon \wt{V}_i\to Z\oplus W_i'
$$
is of class $C^i$,  for every $i\geq 1$.  Recalling from Lemma 
\ref{new_lemma_relation} that $K'=\ker Dg(0)=\{(v, D\delta (0)v)\, \vert \, v\in \wt{K}\}$, we denote by $\pi\colon K'\to\wt{K}$
the projection
$$\pi\colon (v, D\delta (0)v)\mapsto v.$$
The map $\pi$  is a sc-isomorphism. 

We next introduce the map $\tau$ from a relatively open neighborhood of $0$ in $K'\cap C'$ into $Z\oplus W'$ by 
$$\tau (y)=-D\delta (0)\circ \pi (y)+\gamma (\pi (y)),\quad y\in K'\cap C'.$$
It satisfies $\tau (0)=0$. Differentiating $\tau$  at $y=0$, and using  $D\sigma(0)=0$, we obtain
$$
D\tau(0)\eta =-D\delta(0)\circ \pi ( \eta)+D\delta(0)\circ \pi ( \eta)=0
$$
for all $\eta\in K'$ and hence $D\tau (0)=0$. 
Recalling that $y=(v, D\delta (0)v)$ where $v\in \ker DH(0)$, and $\pi (y)=v$, we  compute 
\begin{equation*}
\begin{split}
y+ \tau (y)&=(v +D\delta(0)v) +(-D\delta(0)v+\gamma (v))\\
&=v+\gamma (v)=(v+\sigma (v))+\delta (v+\sigma (v)),
\end{split}
\end{equation*}
where $v+\sigma (v)\in \R^{n'}$. Hence, by definition of the solution germ $\delta$ we conclude from Lemma \ref{new_lemmat_level_0} that 
$$g(y+\tau (y))=0$$
for all $y\in K'\cap C'$ near $0$ on level $0$. By construction, the map $\tau$ is a sc-smooth germ satisfying $\tau (0)=0$ and $D\tau (0)=0$ and there exists a nested sequence of relatively open subsets $\wt{O}_i=\pi^{-1} (V_i)$ of $0$ in $K'$, where 
$O_i=\{a\in [0,\infty)^k\oplus\R^{n'-k}\ |\ |a|_0<\varepsilon_i'\}$
such that the restrictions satisfy 
$$\tau\vert \wt{O}_i\in C^{i}(\wt{O}_i,Z\oplus W_{i}'),\quad \text{for $ i\geq 1$}.$$
Moreover, if $g(y+w)=0$ for $y\in K'\cap C'$ and $w\in Z\oplus W'$ sufficiently small on level $0$, then $w=\tau (y)$.
To sum up, there exists a relatively open subset $U$ of $0$ in $C$ diffeomorphic to  the relatively open subset $U'=\varphi (U)$ of $0$ in $C'$, and a relatively open subset $V$ of $0$ in $K'\cap C'$ such that  
$$
\varphi \bigl(\{x\in U\  |\ f(x)=0\}\bigr) =\{x'\in U'\, \vert \, g(x')=0\}=\{y+\tau(y)\, \vert \, y\in V\}.
$$
The proof of Theorem \ref{LGS} is finished.

\qed \end{proof}

\begin{corollary}\label{LGS2}\index{C- Germ-solvability {II}}
Let  $f\colon U\rightarrow F$ be  a sc-smooth  germ  satisfying $f(0)=0$ and of the form $f=h+s$,  where $h$ is a basic germ 
and $s$ is a $\ssc^+$-section satisfying $s(0)=0$. We assume 
that the linearization $Df(0)\colon E\to F$ is surjective and the  kernel $K=\ker Df(0)$ is in good position to the partial quadrant 
$C\subset E$. 

Then there exists a good complement $Y$ of $K$ in $E$, so that $K\oplus Y=E$, and  there exists a $C^1$-map 
$\sigma\colon V\rightarrow Y_1$ defined on relatively open neighborhood $V$ of $0$ in $K\cap C$, having the following properties.  
\begin{itemize}
\item[{\em (1)}]\ $\sigma\colon {\mathcal O}(K\cap C,0)\rightarrow (Y,0)$ is a sc-smooth germ.
\item[{\em (2)}]\ $\sigma (0)=0$ and $D\sigma(0)=0$.
\item[{\em (3)}]\ 
$\{x\in U\, \vert \, f(x)=0\} =\{v+\sigma(v)\, \vert \, v\in V\}$ after perhaps replacing $U$ by a smaller open neighborhood.
\end{itemize}

\end{corollary}
\begin{proof}

In view of  Theorem \ref{LGS}, the sc-diffeomorphism $\varphi\colon (U, 0)\to (U', 0)$ between  relatively open subsets $U$ and $U'$ of the partial quadrants $C=[0,\infty)^k\oplus \R^{n-k}\oplus W$ and 
$C'=[0,\infty)^k\oplus \R^{n'-k}\oplus W'$ satisfies 
$T\varphi (0)K=K'=\ker Dg(0)$ and $K'$ is in good position to $C'$.   If $Y'$ is the  good sc-complement in Theorem \ref{LGS} of $K'$ in $E'$, we define $Y=T\varphi(0)^{-1}(Y')$. Then $Y$  is the desired good complemement of $K$ in $E$ with respect to the partial quadrant $C$, so that $E=K\oplus Y$. Let $\pi\colon K\oplus Y\to K$ be the sc-projection.

Recalling  the map $\tau\colon V\subset K'\cap C'\to Y_1'$ from Theorem \ref{LGS},  we define the map 
$\psi\colon V\subset K'\cap C'\rightarrow K$  by 
$$\psi (v)=\pi\circ\varphi^{-1}(v+\tau(v)).$$
The map $\psi$ satisfies  $\psi (0)=0$ and its derivative $D\psi (0)\colon K'\to K$  is equal to 
$D\psi (0)=T\varphi(0)^{-1}\vert K'$, hence it is an isomorphism. 

We claim that the map $\psi$ preserves the degeneracy index, that is 
\begin{equation}\label{index_preserving}
d_{K'\cap C'}(v)=d_{K\cap C}(\psi (v))
\end{equation}
for $v\in V$ close to $0$.  To see this we first assume that $v\in V\subset K'\cap C'$ belongs to $K'\cap (\R^{n'-k}\oplus W')$ where $\R^{n'-k}\oplus W'$ is identified with $\{0\}^k\oplus \R^{n'-k}\oplus W'.$ Hence $v$ is of the form $v=(0, w)\in \{0\}^k\oplus \R^{n'-k}\oplus W'$. 

If $K'$ is one-dimensional, then, in view of Lemma \ref{ll1},  $v=t(a, w')\in \R^{k}\oplus \R^{n'-k}\oplus W'$ where $a=(a_1,\ldots ,a_k)$ satisfies $a_j>0$. Consequently, $t=0$, so that $v=0$. 
Since $\psi (0)=0$,  we have $d_{C'\cap K'}(0)=d_{C\cap K}(\psi (0))$.

If $\dim K'\geq 2$, then,  denoting by $\wt{K'}$ the algebraic complement of $K\cap (\R^{n'-k}\oplus W')$ in $K'$,  we may assume, after linear change of coordinates, that the following holds. 
\begin{itemize}
\item[(a)]\ The finite dimensional subspace $\wt{K}'$ is spanned by the vectors $e'_j=(a'_j, b'_j,  w'_j)$ for $1\leq j\leq m=\dim \wt{K}'$ in which $a_j'$ are the vectors of the standard basis of $\R^{l}$, $b_j'=(b_{j, m+1}',\ldots, b_{j,k}')$ with $b_{j,i}'>0$, and $w_j'\in \R^{n'-k}\oplus W'.$
\item[(b)]\ The good complement $Y'$ of $K'$ in $E'=\R^{n'}\oplus W'$ is contained in $\{0\}^l\oplus \R^{k-l}\oplus \R^{n'-k}\oplus W'$.
\end{itemize}
By adding a basis $e_j'=(0, w_j')\in \{0\}^k\oplus \R^{n'-k}\oplus W'$,  $m+1\leq j\leq \dim K'$,  of $K'\cap (\R^{n'-k}\oplus W')$ we obtain a basis of $K'$.

We have a similar statement for  the subspace $K$ with the algebraic complement $\wt{K}$ of $K\cap (\R^{n-k}\oplus W)$ in $K$, $e'_j$ replaced by $e_j$, and the good complement $Y$ of $K$ contained in $\{0\}^m\oplus \R^{k-m}\oplus \R^{n-k}\oplus W$.

Now if $v=(0, w_1')\in \{0\}^k\oplus \R^{n'-k}\oplus W'$, then, in view of (a), $d_{K'\cap C'}(v)=m=\dim \wt{K}'$.  
Recalling $K'$ is in good position to the partial quadrant $C'$ and $Y'$ is a good complement of $K'$ in $E'$. there exists a constant $\gamma$  such that if $n\in K'$ and $y\in Y'$, then $\abs{y}_0\leq \gamma\abs{n}_0$ implies that $n+y\in C$ if and only if $n\in C$.  It follows from $\tau (0)=0$ and $D\tau (0)=0$ that 
$\abs{\tau (v)}_0\leq \gamma\abs{v}_0$ and,  since $v\in C$,  we conclude that $v\pm \tau (v)\in C$. Since 
for $v$ as above $v_j=0$ for all $1\leq j\leq k$, we find that also $\tau_j (v)=0$ for all $1\leq j\leq k$. Hence 
$v+\tau (v)=(0, w_2')\in \{0\}^k\oplus \R^{n'-k}\oplus W'$.  Because   the map $\varphi$ is a sc-diffeomorphism, $\varphi^{-1}(v+\tau (v))=\varphi^{-1}((0, w_2'))=(0, w)\in  \{0\}^k\oplus \R^{n-k}\oplus W$.

We shall show that $\pi ((0, w))\in K\cap C$ and that 
$d_{K\cap C}(\pi ((0, w))=m$. We decompose $(0, w)$ according to the direct sum  $E=K\oplus Y$ as 
$(0, w)=k+y$ where $k\in K$ and $y\in Y$, and write 
$k=(\alpha_1, \ldots , \alpha_k, w_1)\in \R^k\oplus \R^{n-k}\oplus W$, $y=(y_1,\ldots , y_k, w_2)\in Y$. Then,   using 
$Y\subset \{0\}^m\oplus \R^{k-m}\oplus \R^{n-k}\oplus W$, we find that  $y_1=\ldots y_m=0$. This implies that 
also $\alpha_1=\ldots =\alpha_m=0$. Consequently, $\pi ((0, w))=k=(0, w_1)\in K\cap C$ and 
$d_{K\cap C}(\pi ((0, w))=m$. We have proved the identity 
\eqref{index_preserving} in the case $v\in K'\cap  (\R^{n'-k}\oplus W')$.

In the case $v=(a, w)\in K\cap C\setminus K\cap  (\R^{n-k}\oplus W)$, we use  Lemma \ref{big-pretzel_1a} and  compute  
$$
d_{C\cap K}(\psi(v))=d_C(\varphi^{-1}(v+\tau(v)))=d_{C'}(v+\tau(v))=d_{C'\cap K'}(v).
$$
Hence $\psi\colon V\rightarrow C$ is a $C^1$-map  satisfying 
$\psi (0)=0$, $D\psi (0)(K'\cap C')=K\cap C$, and preserving  the degeneracy index. Consequently, 
applying  the inverse function theorem for partial quadrants in $\R^n$,  Theorem \ref{QIFT},  to the map $D\psi (0)^{-1}\circ \psi$,  we find relatively open neighborhoods $V_0\subset V$ and $V_1\subset C$ of $0$ such that 
$$
\psi \colon V_0\rightarrow V_1
$$
is a $C^1$-diffeomorphism, which has higher and higher differentiability closer and closer to $0$. Considering the map
$$
w\mapsto  \varphi^{-1}(\psi^{-1}(w)+\tau\circ\psi^{-1}(w)), \quad w\in V_1,
$$
we obtain   $\pi\circ  \varphi^{-1}(\psi^{-1}(w)+\tau\circ\psi ^{-1}(w)) =w$,  and define the map $\sigma\colon V_1\to Y$ by 
$$
\sigma(w)=\varphi^{-1}\bigl(\psi^{-1}(w)+\tau\circ\psi^{-1}(w)\bigr)-w =({\mathbbm 1}-\pi)\varphi^{-1}\bigl(\psi^{-1}(w)+\tau\circ\psi^{-1}(w)\bigr).
$$
Finally,   we can now choose  open neighborhoods $V\subset K\cap C$ of $0$ and $U$ of $0$ in $C$ appropriately 
and see that the map $\sigma$ has the desired properties (1)-(3) of Corollary \ref{LGS2}

\qed \end{proof}

\mbox{}
\section{Implicit Function Theorems}\label{ssec3.4}

Let us  recall from Definition \ref{x-filling}
 the concept  of a filled version.
\begin{definition}[{\bf Filling}]\index{D- Filling}
We consider a  strong local bundle $K\to O$, where $K=R(U\triangleleft F)$, and  let  the set $U\subset C\subset E$ be  a relatively open neighborhood  of $0$ in the partial quadrant $C$ of the  sc-Banach space $E$.  Here  $F$ is a sc-Banach space and $R$ is a strong bundle retraction of the form 
$$R(u, h)=(r(u), \rho(u)(h))$$
covering the tame retraction $r\colon U\to U$ onto $O=r(U)$. We  assume that $r(0)=0$.

A sc-smooth section germ 
$(f, 0)$ of the
strong bundle $K\to O$ possesses a  {\bf filling}
if there exists a sc-smooth section germ $(g, 0)$ of the bundle $U\triangleleft F\to U$ satisfying the following properties.
\begin{itemize}
\item[(1)]\ $f(x)=g(x)$  for $x\in O$ close to $0$
\item[(2)]\ If $g(y)=\rho (r(y))g(y)$ for a point $y\in U$ near $0$, then $y\in O$.
\item[(3)]\ The linearisation of the map
$
y\mapsto  [{\mathbbm 1} -\rho(r(y))]\cdot g(y)
$
at the point $0$, restricted to $\ker(Dr(0))$, defines a topological linear  isomorphism
$$
\ker(Dr(0))\rightarrow \ker(\rho (0)).
$$
\end{itemize}
\qed
\end{definition}

\begin{remark}
By replacing $O$ by a smaller set, we may assume in (1) that $f(x)=g(x)$ for all $x\in O$ and in (2) that 
$g(y)=\rho(r(y))g(y)$ for $y\in U$ implies that $y\in O$.
\qed
\end{remark}

We also  recall the notions  sc-Fredholm germ and sc-Fredholm section.

\begin{definition}[{\bf Sc-Fredholm germ}]\index{D- Sc-Fredholm germ}
Let $f$ be a sc-smooth section of the 
strong bundle 
$P\colon Y\to X$ and let $x\in X$ be a smooth point. 
Then $(f, x)$ is a {\bf sc-Fredholm germ}, if there exists a strong bundle chart around $x$ (as defined in 
Definition \ref{def_strong_bundle_chart})
$$\text{$\Phi\colon \Phi^{-1}(V)\to K$ covering $\varphi \colon (V, x)\to (O, 0)$,}$$
where $K\to O$ is a 
strong local bundle containing $0\in O\subset U$, such that the local section germ $\Phi\circ f\circ \varphi^{-1}\colon O\to K$ has a filling $g\colon U\to U\triangleleft F$ near $0$ which possesses the following additional property. There exists a local $\ssc^+$-section $s\colon U\to U\triangleleft F$ satisfying $s(0)=g(0)$ such that $g-s$ is conjugated near $0$ to a basic germ.

\end{definition}
\begin{definition}[{\bf sc-Fredholm section}]\label{def_again_fedholm_sect}
A sc-smooth section $f$ of the 
strong bundle $P\colon Y\to X$ is a {\bf sc-Fredholm section}, if 
\begin{itemize}
\item[(1)]\ $f$ is regularizing, i.e., if $x\in X_m$ and $f(x)\in Y_{m, m+1}$, then $x\in X_{m+1}$.
\item[(2)]\ The germ $(f, x)$ is a sc-Fredholm germ at every smooth point $x\in X$.
\end{itemize}
\qed
\end{definition}
From Theorem \ref{stabxx} we know that sc-Fredholm sections are stable under  
$\ssc^+$-perturbations; if $f$ is a sc-Fredholm section, then $f+s$ is also a sc-Fredholm section, for every $\ssc^+$-section $s$.
\begin{proposition}\label{prop3_52}
Let $x\in X$ be a solution of $f(x)=0$, where  $f$ is a sc-Fredholm section of the
strong bundle $P\colon Y\to X$.  Then $Tf(x)\colon T_xX\to Y_x$  is a sc-Fredholm operator.
\end{proposition}
\begin{proof}
Since $(f, x)$ is a sc-Fredholm germ, there exists a
strong bundle chart $\Phi\colon \Phi^{-1}(V)\to K$ around $x$ such that the section $\wt{f}=\Phi\circ f\circ \varphi^{-1}\colon O\to K$ has a filling $g\colon U\to U\triangleleft F$ for which there exists a $\ssc^+$-section $s\colon U\to U\triangleleft F$ satisfying $s(0)=g(0)$ such that $g-s$  is conjugated near $0$ to a basic germ. Hence, by Proposition \ref{Newprop_3.9}, the linearization $D(g-s)(0)\colon E\to F$ is a sc-Fredholm operator. Since $Ds(0)$ is a $\ssc^+$-operator,  the linearization $Dg(0)$ is sc-Fredholm operator,  in view of Proposition \ref{prop1.21}. From Proposition \ref{filler_new_1} (2) about fillers it follows that $T\wt{f}(0)\colon T_0O\to K_0$ is a sc-Fredholm operator of index 
$\ind (T\wt{f}(0))=\ind (Dg(0))=\ind (D(g-s)(0))$,  and the proposition follows.
\qed \end{proof}
We now  focus on the solution set $\{f=0\}$ of a sc-Fredholm section near a point $x\in X$ of  $f(x)=0$. In the case of a boundary,  i.e., $d_X(x)\geq 1$,  we require that $P:Y\rightarrow X$ is a strong bundle over a tame M-polyfold and that $\ker f'(x)$ is in good position to the partial quadrant $C_xX$
in the tangent space $T_xX$, see 
Definition  \ref{reduced_cone_tangent}, Theorem \ref{hofer}, and Remark \ref{tamecone}. From this it follows that in the case that $X$ 
is tame the partial cones $C_x$ are, in fact, partial quadrants.
\begin{definition}[{\bf Good position of a sc-Fredholm  germ}]\label{new_good_position_def}\index{D- Good position}
A sc-Fred\-holm germ $(f, x)$ of a tame strong bundle
$Y\rightarrow X$  satisfying $f(x)=0$ is in {\bf good position},   if 
\begin{itemize}
\item[(1)]\ $f'(x)\colon T_xX\rightarrow Y_x$ is surjective.
\item[(2)]\  If $d_X(x)\geq 1$,  then  $\ker (f'(x))\subset T_xX$ is in good position to the partial quadrant  $C_xX$ in the tangent space $T_xX$.
\end{itemize}
\qed
\end{definition}
\begin{remark}\index{R- Remark on good position}
Good position at points $x\in X\setminus\partial X$ satisfying $f(x)=0$ just means that $f'(x):T_xX\rightarrow Y_x$ is surjective
and does not require the tameness of $X$, which is an assumption about the boundary geometry of $X$.
In case $x\in\partial X$ we require in addition that $\ker(f'(x))$  is in good position to the partial quadrant $C_x$.
It would be sufficient to assume that $X$ is only tame on a neighborhood of the points $x\in \partial X$ satisfying $f(x)=0$.
We also should remark that it seems possible to have implicit function theorems when $X$ is not tame, but that would require
a detailed analysis of the boundary structure of $X$ and the local geometry of $(f,x)$.
\end{remark}
The fundamental implicit function theorem is as follows.

\begin{theorem}\label{IMPLICIT0}\index{T- Implicit function theorem}
Let $f$ be a sc-Fredholm section  of a tame strong bundle
$Y\rightarrow X$. If $f(x)=0$, and if the sc-Fredholm germ $(f, x)$ is in good position, then there exists an open neighborhood $V$ of $x\in X$ such that  
the solution set $S=\{y\in V\, \vert \, f(y)=0\}$ in $V$ has the following properties.
\begin{itemize}
\item[{\em(1)}]\ At  every point $y\in S$, the sc-Fredholm germ  $(f,y)$ is in good position.
\item[{\em(2)}]\ $S$ is a sub-M-polyfold of $X$ and  the induced  M-polyfold 
structure is equivalent to a smooth manifold structure with boundary with corners. 
\end{itemize}
\end{theorem}
The requirement that $(f,x)$ is in good position, where $f(x)=0$,  means in particular that $f'(x)$ is surjective.
We note that Theorem \ref{implicit-x} and Theorem \ref{bound} are immediate consequences of Theorem  \ref{IMPLICIT0}.

\begin{proof}[{\bf Proof of Theorem \ref{IMPLICIT0}}]
The sc-smooth sc-Fredhom section $f$   of the tame strong bundle $P\colon Y\to X$ over the M-polyfold $X$ is regularizing so that the solutions $y\in X$ of $f(y)=0$ are smooth points. Moreover, at every smooth point $y\in X$, the sc-germ $(f, y)$ is a sc-Fredholm germ.

We now focus  on a neighborhood of the solution $x$ of $f(x)=0$. By assumption, the sc-Fredholm germ $(f, x)$ is in good position according to Definition \ref{new_good_position_def}, so that the kernel $N:=\ker f'(x)$ is in good position to the partial quadrant
 $C_xX\subset T_xX$, and the linearization $f'(x)\colon T_xX\to Y_x$ is surjective. In view of Proposition \ref{prop3_52}, the linear operator $f'(x)$ is a Fredholm operator. Its index is equal to $\ind (f'(x))=\dim N.$

By definition of a sc-Fredholm germ, there exists a strong bundle chart  of the tame strong bundle $P\colon Y\to X$
\begin{equation*}
\begin{CD}
P^{-1}(V)@>\Phi>>K\\
@VPVV @VVpV \\
V@>\varphi>>O 
\end{CD}\,
\end{equation*}
covering the sc-diffeomorphism $\varphi\colon V\to O$, which is defined on the open neighborhood $V\subset X$ of the given point $x\in X$ and satisfies  $\varphi (x)=0$. The retract $K=R(U\triangleleft F)$ is the image of the strong bundle retraction $R\colon U\triangleleft F\to U\triangleleft F$ of the form 
$R(u, h)=(r(u), \rho (u)h)$  and which covers the tame sc-smooth retraction $r\colon U\to U$ onto $O=r(U)$. As usual, the set $U\subset C$ is a relatively  open subset of the partial quadrant $C=[0,\infty)^k\oplus \R^{n-k}\oplus W$ in the sc-Banach space $E=\R^n\oplus  W$.
Still by definition of a sc-Fredholm germ, the push-forward section
$$\wt{f}=\Phi_\ast (f)\colon O\to K$$
possesses a filled version $g\colon U\to U\triangleleft F$. It has, in particular, the property that 
\begin{equation}\label{g=wtilde_f}
g(u)=\wt{f}(u)\quad \text{if $u\in O$}.
\end{equation}
In addition, there exists a strong bundle isomorphism 
\begin{equation*}
\begin{CD}
U\triangleleft F@>\Psi>>U'\triangleleft F'\\
@VVV @VVV \\
U@>\psi>>U'
\end{CD}\,
\end{equation*}
where $F'=\R^{N}\oplus W$. It covers the sc-diffeomorphism $\psi\colon U\to U'$ satisfying  $\psi (0)=0$. The set $U'$ is a relatively open subset of the partial quadrant $C'=T\psi (0)C$ in the sc-space $E'=\R^N\oplus W$. In addition,  by definition of a sc-Fredholm germ, there exists a $\ssc^+$-section $s\colon U\to U\triangleleft F$ satisfying $s(0)=g(0)$ such that the push-forward section 
$$\Psi_\ast (g-s)=h\colon U'\to F'$$
is a basic germ according to Definition \ref{BG-00x}. Therefore, 
\begin{equation}\label{new_equ_50}
\Psi_\ast (g)=h+t,
\end{equation}
where $t=\Psi_\ast (s)$ is a $\ssc^+$-section $U'\to F'$.

In view of Proposition \ref{Newprop_3.9}, the linearization $D(h+t)(0)\colon E'\to F'$ is a sc-Fredholm operator, so that, by Proposition \ref{prop1.21}, the operator $Dg(0)\colon E\to F$ is also a sc-Fredholm operator.

From the postulated surjectivity of the linearization $f'(x)\colon T_xX\to Y_x$ we deduce that 
$T\wt{f}(0)\colon T_0O\to K_0$ is surjective. Hence, by Proposition \ref{filler_new_1}, the operator $Dg(0)$ is surjective. Since $\ker f'(x)$ is in good position to the partial quadrant $C_xX\subset T_xX$, and since $T\varphi (x)(C_xX)=C_0O=T_0O\cap C$, the kernel $\ker T\wt{f}(0)$ is in good position to the partial quadrant  $C_0O$. The retract $O=r(U)$ is, by assumption, tame. Hence,  the tangent space $T_0O$ has, in view Proposition \ref{IAS-x}, the sc-complement $\ker Dr(0)$ in $E$ so that 
$$E=T_0O\oplus \ker Dr(0),$$
and $\ker Dr(0)\subset E_x$ at the point $x=0$ ($E_x$ is defined in Definition \ref{new_def_2.33}).
If $Z\subset T_0O$ is a good complement of $\ker T\wt{f}(0)$ in $T_0O=\ker T\wt{f}(0)\oplus Z$, the space $Z\oplus \ker Dr(0)$ is a good complement of 
$\ker T\wt{f}(0)$ in $E$, so that 
$$E=\ker T\wt{f}(0)\oplus (Z\oplus \ker Dr(0)).$$

Since, by the properties of the filler,  
$$\ker T\wt{f}(0)\oplus \{0\}=\ker Dg(0),$$
we conclude that $\ker Dg(0)$ is in good position to $C$ in $E$. Therefore, $N'=\ker D(h+t)(0)$ is in good position to $C'=T\psi (0)C$ in $E'$. Let now $Y'\subset E'$ be the  good complement of $N'$ in $E'=N'\oplus Y'$  from Theorem \ref{LGS2}. 

Then we can apply Theorem \ref{LGS2} 
about the local germ solvability and find an open neighborhood $V'\subset N'\cap C'$ of $0$,  and a map 
$$\sigma\colon V'\to (Y')_1$$ 
possessing the following properties.
\begin{itemize}
\item[(1)]\ $\sigma\colon V'\to (Y')_1$ is of class $C^1$ and satisfies $\sigma (0)=0$ and $D\sigma (0)=0$.
\item[(2)]\ $\sigma\colon {\mathcal O}(N'\cap C', 0)\to (Y', 0)$ is a sc-smooth germ.
\item[(3)]\ 
\begin{align*}
\psi\circ \varphi (&\{y\in V\, \vert \, f(y)=0\}&& \text{}\\
&=\psi (\{u\in O\, \vert \, \wt{f}(u)=0\})&& \text{}\\
&=\psi (\{u\in O\, \vert \, g(u)=0\})&& \text{by \eqref{g=wtilde_f}}\\
&=\{u'\in U'\, \vert \, (h+t)(u')=0\}&& \text{by \eqref{new_equ_50}}\\
&=\{v+\sigma (v)\, \vert \, v\in V'\}&& \text{by Corollary \ref{LGS2}.}
\end{align*}
\end{itemize}
Moreover,
\begin{itemize}
\item[(4)]\ For every $y\in V$ satisfying $f(y)=0$, so that  $\psi\circ \varphi (y)=v+\sigma (v)$, the kernel $\ker f'(y)$ is in good position to the partial quadrant  $C_yX\subset T_yX$, and $f'(y)\colon T_yX\to Y_y$ is surjective.
\end{itemize}
 Property (4) follows from Corollary \ref{corex1} and from the following lemma, whose proof is postponed to Appendix \ref{geometric_preparation}.

\begin{lemma}\label{good_pos}
Let $C\subset E$ be  a partial quadrant in the  sc-Banach space $E$ and $N\subset E$ a finite-dimensional smooth subspace in good position to $C$ and let $Y$ be a good complement of $N$ in $E$, so that $E=N\oplus Y$.
We assume that $V\subset N\cap C$ is a relatively open neighborhood
of $0$ and $\tau\colon V\rightarrow Y_1$ a map of class $C^1$ satisfying $\tau(0)=0$ and $D\tau(0)=0$.

Then there exists a relatively open neighborhood $V'\subset V$ of $0$ such  that the following holds.
\begin{itemize}
\item[{\em (1)}]\ $v+\tau(v)\in C_1$ for $v\in V'$.
\item[{\em (2)}]\ For every $v\in V'$,  the linear subspace $N_v=\{n+D\tau(v)n\, \vert \, n\in N\}$ has the Banach space $Y=Y_0$ as a topological complement.
\item[{\em (3)}]\ For every $v\in V'$,  there exists a  constant $\gamma_v>0$ such  that  if  $n\in N_v$ and $y\in Y$ satisfy 
$\abs{y}_0\leq \gamma_v\cdot \abs{n}_0$,  the statements $n\in C_z$ and $n+y\in C_z$ are equivalent, where $z=v+\tau (v)$. 
\end{itemize}

\end{lemma}

So far we are confronted with the following situation. The open neighborhood $V\subset X$ of the smooth point $x$ is a M-polyfold and we denote the solution set of the sc-Fredholm section $f$ of the tame strong bundle $P^{-1}(V)\to V$ by $S=\{y\in V\, \vert \, f(y)=0\}$. It consists of smooth points. For every $y\in S$, the germ $(f, y)$ is a sc-Fredholm germ and the linearization $f'(y)\colon T_yV\to Y_y$ is a surjective Fredholm operator whose kernel $\ker f'(y)$ is in good position to the partial cone $C_yV\subset T_yV$, so that, proceeding as above we can construct a map $\sigma$ satisfying the above properties (1)-(3). Consequently, abbreviating $d=\dim (\ker f'(x))$, the solution set $S\subset V$ possesses the $d$-dimensional tangent germ property according to Definition \ref{toast}. Therefore we can apply Theorem \ref{HKL} to conclude that the solution set 
$S=\{y\in V\, \vert \, f(y)=0\}$ is a sub-M-polyfold, whose induced M-polyfold structure is equivalent to the structure of a smooth manifold with boundary with corners.

This completes the proof of Theorem \ref{IMPLICIT0}.
\qed \end{proof}
Finally,  we note two immediate consequences of Theorem \ref{IMPLICIT0}.
{\begin{theorem}[{\bf Global implicit function theorem {I}}]\label{io-xx}\index{T- Global implicit function theorem}
If $P\colon Y\rightarrow X$ is a strong bundle over an  M-polyfold $X$ satisfying  $\partial X=\emptyset$,  and if $f$ a sc-Fredholm section having the property that at  every point $x$ in the solution set $S=\{y\in X\ |\ f(y)=0\}$,  the linearization $f'(x)\colon T_xX\rightarrow Y_x$ is surjective. Then $S$ is a sub-M-polyfold
of $X$ and the induced M-polyfold structure on $S$ is equivalent to the structure of a smooth manifold  without boundary.
\end{theorem}

In a later section we shall study the question how to perturb a sc-Fredholm section to guarantee the properties required 
in the hypotheses of Theorem \ref{io-xx} and the following  boundary version.
\begin{theorem}[{\bf Global implicit function theorem {II}}]\label{io-xxx}\index{T- $\partial$-Global implicit function theorem}
Let $P\colon Y\rightarrow X$ be a strong bundle over the tame M-polyfold $X$,   and $f$ a sc-Fredholm section having the property that
at  every point $x$ in the solution set $S=\{y\in X\ |\ f(y)=0\}$,  the linearization $f'(x)\colon T_xX\rightarrow Y_x$ is surjective and the kernel $\ker(f'(x))$ is in good position
to the partial cone $C_xX\subset T_xX$.  Then $S$ is a sub-M-polyfold
of $X$ and the induced M-polyfold structure on $S$ is equivalent to the structure of a smooth manifold with boundary with corners.
\end{theorem}

\section{Conjugation to a Basic Germ}

A useful criterion to decide  in practice whether  a filled version is conjugated to a basic germ is given in Theorem \ref{basic_germ_criterion} below.
We would like to point out a similar result,  due to K. Wehrheim,  in  \cite{Wehr}.
 The following criterion was 
introduced in \cite{HWZ5} and employed to show that the nonlinear Cauchy-Riemann operator occurring in the Gromov-Witten theory
defines a sc-Fredholm section.
\begin{theorem}[{\bf Basic Germ Criterion}]\label{basic_germ_criterion}\index{T- Basic germ criterion}
Let  $U$ be  a relatively open neighborhood of $0$ in the partial quadrant $C$ of the sc-Banach space $E$, and let 
${\bm{f}}\colon U\rightarrow F$ be a sc-smooth map into the sc-Banach space $F$, which satisfies the following conditions.
\begin{itemize}
\item[{\em (1)}]\ At every smooth point $x\in U$ the linerarization $D{\bm{f}}(x)\colon E\rightarrow F$
is a sc-Fredholm operator and the index does not depend on $x$.
\item[{\em (2)}]\ There is a sc-splitting $E=B\oplus X$ in which $B$ is a finite-dimensional subspace of $E$ containing the kernel
of $D{\bm{f}}(0)$, and $X\subset C$, such  that the following holds for $b\in B\cap U$ small enough.  If $(b_j)\subset B\cap U$
is a sequence converging to $b$ and if $(\eta_j)\subset X$ is a sequence which is bounded on level $m$ and satisfying
$$
D{\bm{f}}(b_j)\eta_j=y_j+z_j,
$$
where $y_j\rightarrow 0$ in $F_m$ and where the sequence $(z_j)$ is bounded in $F_{m+1}$, then the sequence $(\eta_j)$ possesses a convergent
subsequence in $X_m$.
\item[{\em (3)}]\ For every  $m\geq 0$ and $\varepsilon>0$, the estimate 
$$
\abs{[D_2{\bm{f}}(b,0)-D_2{\bm{f}}(b,x)]h}_m\leq \varepsilon\cdot \abs{h}_m
$$
holds for all $h\in E_m$, and for  $b\in B$ sufficiently small,  and for $x\in X_{m+1}$ sufficiently small on level $m$, i.e.
we require $x$ to be on level $m+1$, but $|x|_m$ small only for level $m$.
\end{itemize}
Then the section ${\bf f-s}$ is conjugated near $0$ to a basic germ, where ${\bm{s}}$ is the constant $\ssc^+$-section ${\bm{s}}(x)={\bm{f}}(0)$.
\end{theorem}

\begin{proof}
Abbreviating  $Y:=D{\bm{f}}(0)X\subset F$, the restriction 
$$
D{\bm{f}}(0)\vert X\colon X\rightarrow Y
$$
is an injective and surjective sc-operator. Since $D{\bm{f}}(0)\colon E\to F$ is a sc-Fredholm operator, we find a finite-dimensional sc-complement $A$ of $Y$ in $F$, such that 
$$
F=Y\oplus A.
$$
Denoting  by $P\colon Y\oplus A\rightarrow Y$ the sc-projection, we consider the family $b\mapsto L(b)$ of bounded linear operators, defined by  
$$L(b)=P\circ D{\bm{f}}(b)\vert X\colon X\to Y.$$
It is not assumed  that the operators $L(b)$ depend continuously on $b$. Since $B$ is a finite-dimensional
sc-smooth space,  the map 
$$
(B\cap U)\oplus X\rightarrow Y\colon (b,x)\mapsto  P\circ D{\bm{f}}(b)x
$$
is sc-smooth. If we raise the index by one, then the map  
$$
(B\cap U)\oplus X^1\rightarrow Y^1\colon (b,x)\mapsto  P\circ D{\bm{f}}(b)x
$$
is also sc-smooth by Proposition \ref{sc_up}.

\begin{lemma}\label{trivial_kernel}\index{L- Basic germ criterion {I}}
There exists a relatively open neighborhood $O$ of $0$ in $B\cap C$ such  that  the composition  
$$
P\circ D{\bm{f}}(b)\colon X\rightarrow Y,
$$
has a trivial kernel for every $b\in O$.
\end{lemma}
\begin{proof}
Assuming  that such an open set $O$ does not exist,  we find a sequence $b_j\in B\cap C$ satisfying  $b_j\rightarrow 0$ and a sequence $(h_j)\subset X_0$ satisfying $|h_j|_0=1$ such  that
$P\circ D{\bm{f}}(b_j)h_j=0$. Then $D{\bm{f}}(b_j)h_j = z_j$ is a bounded sequence in $A$,  in fact on every level since $A$ is a smooth subspace, and we consider  the level $1$ for the moment.
From property (2) we deduce that $(h_j)$ has a convergent subsequence in $X_0$. So,  without loss of generality,  
we may assume $h_j\rightarrow h$ in $X_0$ and  $|h|_0=1$. Hence
$$
P\circ D{\bm{f}}(0)h=0, 
$$
in contradiction to the injectivity of the map $P\circ D{\bm{f}}(0)\vert X$. 
The lemma is proved.
\qed \end{proof}

From the property (1) and the fact that $P$ is a sc-Fredholm operator,  we conclude that $P\circ D{\bm{f}}(x)$ for smooth $x$ are all sc-Fredholm operators having  the same index. In particular, if $b\in O$, then $P\circ D{\bm{f}}(b)\colon X\rightarrow Y$ is an  injective sc-Fredholm operators of index $0$ in view of Lemma \ref{trivial_kernel},  and 
hence  a  sc-isomorphism. Next we sharpen this result.

\begin{lemma}\index{L- Basic germ criterion {II}}
We take a relatively  open neighborhood $\wt{O}$ of $0$ in $C\cap B$ whose compact closure is contained in $O$. Then
for every level $m$ there exists a number $c_m>0$ such  that,  for every $b\in \wt{O}$,  we have the estimate
$$
|P\circ D{\bm{f}}(b)h|_m\geq c_m\cdot |h|_m\quad  \text{for all $h\in X_m$.}
$$
\end{lemma}

\begin{proof}

Arguing indirectly we find a level $m$ for which there is no such constant $c_m$. Hence there are sequences
$(b_j)\subset \wt{O}$ and $(h_j)\subset X_m$ satisfying $|h_j|_m=1$ and  
$\abs{P\circ D{\bm{f}}(b_j)h_j}_m\rightarrow 0$.
After perhaps taking a subsequence we may assume that $b_j\rightarrow b$ in $O$. 

From 
$$
D{\bm{f}}(b_j)h_j =P\circ D{\bm{f}}(b_j)h_j=y_j\to 0\quad \text{in $F_m$},
$$
we conclude,  in view of the property (2),  for a subsequence,  that $h_j\to h$ in $X_m$. By continuity, $P\circ D{\bm{f}}(b)h=0$ and $\abs{h}_m=1$, in contradiction to the fact, that $P\circ D{\bm{f}}(b)\colon X\to Y$ is a 
sc-isomorphism for $b\in O$.

\qed \end{proof}

So far we have verified that the family $b\mapsto L(b)$ meets the assumptions 
of the following lemma,  taken from \cite{HWZ8.7}, Proposition 4.8.

\begin{lemma}\label{family_maps}\index{L- Basic germ criterion {III}}
We assume that $V$ is a relatively open subset in the partial quadrant of a finite-dimensional vector space $G$.
We suppose further that $E$ and $F$ are sc-Banach spaces and consider a family of linear operators $v\rightarrow L(v)$ 
having  the following properties.
\begin{itemize}
\item[{\em (a)}]\ For every $v\in V$, the linear operator  $L(v)\colon E\rightarrow F$ is a sc-isomorphism.
\item[{\em (b)}]\ The map
$$
V\oplus E\rightarrow F,\quad (v,h)\mapsto L(v)h
$$
is sc-smooth.
\item[{\em (c)}]\ For every level $m$ there exists a constant $c_m>0$ such that for $v\in V$ and all $h\in E_m$
$$
|L(v)h|_m\geq c_m\cdot |h|_m.
$$
\end{itemize}
Then the well-defined map
$$
V\oplus F\rightarrow E\colon (v,k)\mapsto  L(v)^{-1}(k)
$$
is sc-smooth.
\qed
\end{lemma}

Let us emphasize that it is  not assumed  that the operators $v\rightarrow [L(v)\colon E_m\rightarrow F_m]$ depend continuously as operators 
on $v$.

In view of Lemma \ref{family_maps}, the map 
$B\oplus Y\to Y$,
$$(b, y)\mapsto \bigl(P\circ Df(b)\vert X\bigr)^{-1}y,$$
is sc-smooth.

We may assume that the finite-dimensional space  $B$ is equal to $\R^n$ and that $E=\R^n\oplus X$ and $C=[0,\infty)^k\oplus \R^{n-k}\oplus X$ is the partial quadrant in $E$. Hence $B\cap C=[0,\infty)^k\oplus \R^{n-k}$. Moreover, we may identify the finite-dimensional subspace $A$ of $F$ with $\R^N=A$. Replacing, if necessary, the  relatively open neighborhood $U \subset C$ of $0$ by a smaller one
we may assume, in addition,  that $(b,x)\in U$ implies that $b\in \wt{O}$.

 We now define a strong bundle map
$$
\Phi\colon U\triangleleft (\R^N\oplus Y)\rightarrow U\triangleleft (\R^N\oplus Y)
$$
covering the identity $U\to U$ by
$$
\Phi((b,x),(c,y))= ((b,x), (c,[P\circ D{\bm{f}}(b)|X]^{-1}(y))).
$$
We define  the sc-smooth germ $({\bm{h}},0)$ by ${\bm{h}}(b,x)={\bm{f}}(b,x)-{\bm{f}}(0,0)$, where  $(b,x)\rightarrow {\bm{f}}(0,0)$ is a constant $\ssc^+$-section.
We shall show that the  push-forward germ ${\bm{k}}=\Phi_\ast({\bm{h}})$ is a basic germ. Using $D{\bm{h}}=D{\bm{f}}$,  we compute 
$$
{\bm{k}}(b,x) = (({\mathbbm 1}-P){\bm{h}}(b,x), [P\circ D{\bm{h}}(b)|X]^{-1}P{\bm{h}}(b,x))\in \R^N\oplus X.
$$
The germ ${\bm{k}}$ is  a sc-smooth germ 
$$ 
{\mathcal O}([0,\infty)^k\oplus \R^{n-k}\oplus X,0)\rightarrow (\R^N\oplus X,0).
$$
Denoting by $Q\colon \R^N\oplus X\rightarrow X$ the sc-projection, we shall verify that $Q{\bm{k}}$ is a $\ssc^0$-contraction germ.  We define the sc-smooth  germ $H\colon \R^n\oplus X\to X$ by 
$$
H(b,x) = x- Q{\bm{k}}(b,x)=x - [P\circ D{\bm{h}}(b)|X]^{-1}P{\bm{h}}(b,x).
$$
By construction,  $H(0,0)=0$.  The family $L(b)\colon X\to Y$ of bounded linear operators, defined by 
$$
L(b):=P\circ D{\bm{f}}(b)\vert X,\quad b\in B,
$$
satisfies the assumptions of Lemma \ref{family_maps}. Recalling now the condition (3) of Theorem \ref{basic_germ_criterion}, we choose $m\geq 0$ and $\varepsilon>0$ and accordingly take $b\in B$ small and $x, x'\in X_{m+1}$ small on level $m$. Then using the estimates in condition (3) and in Lemma  \ref{family_maps}, we estimate, recalling that $P\circ D_2h(b, x)=P\circ D_2f(b, x)$, 

\begin{equation*}
\begin{split}
&|H(b,x)-H(b,x')|_m=\abs{L(b)^{-1} [L(b)(x-x') - P{\bm{h}}(b,x) +P{\bm{h}}(b,x')}_m\\
&\quad \leq c_m^{-1}\cdot \abs{P \circ \bigl[D_2{\bm{f}}(b, 0)(x-x')  - {\bm{h}}(b,x) +{\bm{h}}(b,x')] }_m\\
&\quad =c_m^{-1}\cdot \abs{ 
\int_0^1 
P\circ  \bigl[D_2{\bm{f}}(b, 0)- D_2{\bm{f}}(b,tx +(1-t)x')\bigr](x-x')dt}_m\\
&\quad \leq  c_m^{-1}\cdot d_m\cdot\varepsilon\cdot |x-x'|_m.
\end{split}
\end{equation*}
The map $(b, x)\mapsto H(b, x)$ is sc-smooth. Therefore, using the density of $X_{m+1}$ in $X_m$, we conclude,  
for every $m\geq 0$ and $\varepsilon>0$,  that the estimate 
$$\abs{H(b, x)-H(b, x')}_m\leq \varepsilon\abs{x-x'}_m$$
holds, if $b\in B$ is small enough and $x, x'\in X_m$ are sufficiently small.  Having verified that the push-forward germ ${\bm{k}}=\Phi_\ast ({\bm{h}})$ is a basic germ, the proof of Theorem \ref{basic_germ_criterion} is complete.

\qed \end{proof}
\section{Appendix}

\subsection{Proof of Proposition \ref{pretzel}}\label{pretzel-A} 
We need to show the following:\par

\noindent{\bf Proposition} If  $N$ is a finite-dimensional sc-subspace in good position to the partial quadrant  $C$ in $E$, then
$N\cap C$ is a partial quadrant in $N$.\par

As a preparation for the proof we recall some tools and results from the appendix in \cite{HWZ3} and begin with the geometry of closed convex cones and quadrants in finite dimensions.
A  {\bf closed convex cone} $P$ in a finite-dimensional vector space
 $N$ is a closed convex
subset satisfying  $P\cap (-P)=\{0\}$ and $\R^+\cdot P=P$.
 An
{\bf extreme ray} in a closed convex cone $P$ is a subset $R$ of the
form
$$
R=\R^+\cdot x, 
$$
where  $x\in P\setminus \{0\}$, having the property that  if   $y\in P$ and  $x-y\in P$, then $y\in R$. If the cone $P$ has a nonempty interior, then  it generates the vector space  $N$, that is $N=P-P$.

A {\bf quadrant} $C$ in a vector space $N$ of dimension $n$ is a closed convex cone such that there exists a linear isomorphism $T\colon N\rightarrow \R^n$ mapping $C$ onto $[0,\infty)^n$. We observe that a quadrant
in $N$ has precisely $\dim(N)$ many extreme rays.  

The following  version of the Krein-Milman theorem
is well-known, see exercise 30 on page 72 in \cite{Schaefer}.  A proof can be found in the appendix of \cite{HWZ3}, Lemma 6.3.
\begin{lemma}\label{kreinmilman}
A closed convex cone $P$ in a finite-dimensional vector space $N$ is
the  closed convex hull of its extreme rays.
\end{lemma}
 A closed convex cone $P$ is called {\bf finitely generated} provided $P$ has
finitely many extreme rays. If this is the case, then  $P$ is the convex hull of
its finitely many extreme rays.

 For example,  if $C$ is a partial
quadrant in the sc-Banach space $E$ and $N\subset E$ is a finite-dimensional subspace of
$E$ such  that $C\cap N$ is a closed convex cone, then $C\cap N$ is
finitely generated.

\begin{lemma}\label{nomer}
Let $N$ be a finite-dimensional vector space and $P\subset N$ a
closed convex cone having a nonempty interior. Then $P$ is a quadrant if
and only if it has $\dim(N)$-many extreme rays.
\end{lemma}
The proof is given in  \cite{HWZ3}, Lemma  6.4. 

We consider the sc-Banach space $E =\R^n\oplus W$  containing the partial quadrant $C =[0,\infty)^n\oplus W$. A point  $a\in C=[0,\infty )^n\oplus W\subset \R^n\oplus W$, has  the representation
$a=(a_1,\ldots ,a_n, a_{\infty})$,  where $(a_1,\ldots ,a_n)\in [0,\infty )^n$ and $a_{\infty}\in W$.
By $\sigma_a$ we shall denote the collection of all  indices $i\in \{1,\ldots ,n\}$ for which
$a_i=0$ and denote by $\sigma_a^c$  the complementary set of indices in $\{1,\ldots ,n\}$. Correspondingly, we introduce the following subspaces in $\R^n$,
\begin{align*}
\R^{\sigma_a}&=\{x\in \R^n\ \vert \ \text{$x_j=0$ for all $j\not \in \sigma_a$}\}\\
\R^{\sigma^c_a}&=\{x\in \R^n\ \vert \ \text{$x_j=0$ for all $j\not \in \sigma^c_a$}\}.
\end{align*}
The next lemma and its proof is taken from the appendix in \cite{HWZ3}. The hypothesis that $C\cap N$ is a closed convex cone is crucial.

\begin{lemma}\label{roxy}
Let $N\subset E_{\infty}$ be  a finite-dimensional smooth subspace of  $E=\R^n\oplus W$ such  that $N\cap C$ is a closed convex cone. If $a\in N\cap C$ is a generator of  an extreme ray $R=\R^+\cdot x$ in $N\cap C$, then
$$
\dim(N)-1\leq \#\sigma_a.
$$
If,  in addition, $N$ is in good position to $C$, then
$$
\dim(N)-1= \#\sigma_a.
$$
\end{lemma}
\begin{proof}
We assume $R=\R^+\cdot a$ is an extreme ray in $C\cap N$ and abbreviate $\sigma=\sigma_a$ and its complememnt in $\{1, \ldots ,n\}$ by  $\sigma^c$.  Then $R\subset N\cap C\cap (\R^{\sigma^c}\oplus W)$. Let $y\in C\cap N\cap (\R^{\sigma^c}\oplus
W)$ be a nonzero element. Since $a_i>0$ for all $i\in\sigma^c$,  there exists $\lambda>0$ so that $\lambda a-y\in N\cap C\cap
(\R^{\sigma^c}\oplus W)\subset N\cap C$. We conclude $y\in R$ because  $R$ is an extreme ray. Given any element $z\in
N\cap(\R^{\sigma^c}\oplus W)$ we find $\lambda>0$ so that
$\lambda a+z\in N\cap C\cap(\R^{\sigma^c}\oplus W) $ and infer,
by the previous argument,  that $\lambda a+z\in R$. This implies that $z\in \R \cdot a$. Hence
\begin{equation}\label{pol}
 \dim(N\cap(\R^{\sigma^c}\oplus W))=1.
\end{equation}
The projection $P\colon \R^n\oplus W=\R^{\sigma}\oplus (\R^{\sigma^c}\oplus W) \to \R^{\sigma}$ induces a linear map
\begin{equation}\label{ohx}
P\colon N\rightarrow \R^{\sigma}
\end{equation}
which by \eqref{pol} has  an one-dimensional kernel. Therefore,
$$
\#\sigma=\dim(\R^{\sigma})\geq \dim R(P)=\dim N-\dim \ker P=\dim(N)-1.
$$
Next assume $N$ is in good position to $C$. Hence there exist a constant $c>0$ and a sc-complement $N^{\perp}$ such that $N\oplus N^{\perp}=\R^n\oplus W$ and if $(n, m)\in N\oplus N^{\perp}$ satisfies $|m|_0\leq c |n|_0$, then $n+m\in C$ if and only if $n\in C$.  We claim that  $N^{\perp}\subset \R^{\sigma^c}\oplus W$. Indeed, let $m$ be any element of $N^{\perp}$.  Multiplying $m$ by a real number we may assume  $|m|_0\leq c |a|_0$ .  Then $a+m\in C$ since $a\in C$. This implies that
$m_i\geq 0$ for  all indices $i\in\sigma_a$. Replacing $m$ by $-m$,  we conclude
$m_i=0$ for all $i\in\sigma_a$.  So $N^{\perp}\subset  \R^{\sigma^c}\oplus W$ as claimed.
Take  $k\in \R^{\sigma_a}$
and write $(k, 0)=n+m\in N\oplus N^{\perp}$.  From  $N^{\perp}\subset {\mathbb
R}^{\sigma^c}\oplus W$, we conclude $P(n)=k$. Hence the map $P$  in \eqref{ohx}  is surjective  and the desired result follows.
\qed \end{proof}

Having studied the geometry of closed convex cones and partial quadrants in finite dimensions we shall next study finite dimensional subspaces $N$ in good position to a partial quadrant $C$ in a sc-Banach space.  In this case,  $N\cap C$ can be a partial quadrant rather than a quadrant, which requires, some additional arguments.

We assume that $N$ is a smooth  finite-dimensional subspace of
$E=\R^n\oplus W$ which is in good position to the partial quadrant $C=[0,\infty)^n\oplus
W$. Thus, by definition, there is  a sc-complement,  denoted by  $N^{\perp}$,  of $N$ in $E$ and a constant $c>0$ such 
that if  $(n,m)\in N\oplus N^{\perp}$ satisfies  $|m|_0\leq c |n|_0$, then  the statements
$n\in C$ and $n+m\in C$ are equivalent.
We introduce the subset $\Sigma$ if $\{1,\ldots,n\}$ by 
$$\Sigma=\bigcup_{a\in C\cap N, a\neq 0}\sigma_{a}\subset \{1,\ldots ,n\}.$$
and denote by $\Sigma^c$ the complement $\{1,\ldots ,n\}\setminus \Sigma$. The associated subspaces of $\R^n$ are defined by
$\R^{\Sigma}=\{x\in \R^n\ \vert \ \text{$ x_j=0$ for $j\not \in \Sigma$}\}$ and
$\R^{{\Sigma}^c}=\{x\in \R^n\ \vert \ \text{$ x_j=0$ for $j\not \in \Sigma^c$}\}$.

\begin{lemma}\label{cone1}
$N^{\perp}\subset {\mathbb
R}^{\Sigma^c}\oplus W$.
\end{lemma}
\begin{proof}
Take   $m\in N^{\perp}$. We  have to show that $m_i=0$ for all $i\in\Sigma$.
So fix an index  $i\in\Sigma$ and let $a$ be a nonzero element of  $C\cap N$ such that
$i\in\sigma_a$. Multiplying $a$ by a suitable positive number we may
assume $|m|_0\leq c |a|_0$.  Since $a\in C$,   we infer that $a+m\in C$. This implies
that $a_i+m_i\geq 0$. By definition of $\sigma_a$, we have $a_i=0$  implying $m_i\geq 0$.
Replacing $m$ by $-m$ we find $m_i=0$. Hence $N^{\perp}\subset {\mathbb
R}^{\Sigma^c}\oplus W$ as claimed.
\qed \end{proof}

Identifying $W$ with $\{0\}\oplus W$,  we take   an algebraic complement $\wt{N}$ of $N\cap W$ in $N$ so that
\begin{equation}\label{hoferx}
N=\wt{N}\oplus (N\cap W)\quad \text{and}\quad\quad E= \wt{N}\oplus(N\cap W)\oplus N^{\perp}.
\end{equation}
Let us note that the projection $\pi\colon \R^n\oplus W\rightarrow  \R^n$ restricted to $\tilde{N}$ is injective,  so that
\begin{equation}\label{hoferx0}
\dim(\pi(\tilde{N}))=\dim(\tilde{N}).
\end{equation}
\begin{lemma}\label{ntilde}
If the subspace $N$ of $E$ is in good position to the partial quadrant  $C$, then
$\wt{N}$ is also  in good position to $C$ and the  subspace $\wt{N}^{\perp}:=(N\cap W)\oplus N^{\perp}$ is a good complement of $\wt{N}$ in $E$.
\end{lemma}
\begin{proof}
We define $\abs{x}:=|x|_0$.
Since $N$ is in good position to the quadrant $C$ in $E$, there exist a constant $c>0$ and a sc-complement $N^{\perp}$ of $N$ in $E$ such that if $(n, m)\in N\oplus N^{\perp}$ satisfies $\abs{m}\leq c\abs{n}$, then
 the statements $n\in C$ and $n+m\in C$ are equivalent.  Since $E$ is a Banach space and $N$ is a finite dimensional subspace of $E$, there is a constant $c_1>0$ such that
 $\abs{n+m}\geq c_1[ \abs{n}+\abs{m}]$ for all $(n, m)\in N\oplus N^{\perp}$.  To prove  that $\wt{N}$ is in good position to $C$,  we shall show that $\wt{N}^{\perp}:=(N\cap W)\oplus N^{\perp}$ is a good complement of $\wt{N}$ in $E$. Let $(\wt{n}, \wt{m})\in \wt{N}\oplus \wt{N}^{\perp}=E$ and assume that  $\abs{\wt{m}}\leq c_1c\abs{\wt{n}}$. Write $\wt{m}=n_1+n_2\in (N\cap W)\oplus N^{\perp}$. Since $c_1[\abs{n_1}+\abs{n_2}]\leq \abs{n_1+n_2}=\abs{\wt{m}}\leq c_1c\abs{\wt{n}}$, we get $\abs{n_2}\leq c\abs{\wt{n}}$.  Note that $\wt{n}+\wt{m}=\wt{n}+n_1+n_2\in C$ if and only if $\wt{n}+n_2\in C$ since $n_1\in  \{0\}\oplus W$.  Since  $\abs{n_2}\leq c\abs{\wt{n}}$, this is equivalent to $\wt{n}\in C$. It remains to show that $\wt{N}\cap C$ has a nonempty interior. By assumption $N\cap C$ has nonempty interior. Hence there is a point $n\in N\cap C$ and $r>0$ such that the ball $B^{N}_r(n)$ in $N$ is contained in $N\cap C$. Write $n=\wt{n}+w$ where $\wt{n}\in \wt{N}$ and $w\in N\cap W$. Since $n\in C$ and  $w\in W$, we conclude that $\wt{n}\in C$. Hence $\wt{n}\in \wt{N}\cap C$.
Take $\nu \in B^{\wt{N}}_r(\wt{n})$,  the open ball in $\wt{N}$ centered at
  $\wt{n}$ and of radius $r>0$.  We want to prove that $\nu \in C$. Since $C=[0,\infty )^n\oplus W\subset \R^n\oplus W$, we have to prove for $\nu= (\nu', \nu'')\in \R^n\oplus W$ that $\nu'\in [0,\infty )^n$.  We estimate  $\abs{(\nu +w)-n}=\abs{(\nu +w)-(\wt{n}+w)}=\abs{\nu -\wt{n}}<r$ so that $\nu +w\in B^{N}_r(n)$ and hence  $\nu +w\in N\cap C$ . Having identified $W$ with $\{0\}\oplus W$, we have $w=(0, w'')\in \R^n\oplus W$.  Consequently, $\nu+w=(\nu', \nu''+w'')\in N\cap C$ implies $\nu' \in [0,\infty )^n$.  Since also $\nu \in \wt{N}$, one concludes  that $\nu \in \wt{N}\cap C$ and that $\wt{n}$ belongs to the interior of $\wt{N}\cap C$ in $\wt{N}$.  The proof of Lemma \ref{ntilde} is complete.
 \qed \end{proof}
 Since by  Lemma \ref{cone1},  $N^{\perp}\subset \R^{\Sigma^c}\oplus W$,  the  good complement   $\wt{N}^{\perp}=(N\cap W)\oplus N^{\perp}$ satisfies
$$
\wt{N}^{\perp}\subset \R^{\Sigma^c}\oplus W.
$$
We claim that $C\cap \wt{N}$ is a closed convex cone . 
It suffices to verify that if $a\in C\cap \wt{N}$ and $-a\in C\cap \wt{N}$, then $a=0$. We write 
$a=(a', a_\infty)$ where $a'\in \R^n$ and $a_\infty\in W$. Then $a', -a'\in [0,\infty)^n$ implies that $a'=0$ so that $a=(0, a_\infty)\in \{0\}\oplus W$. Since $a\in \wt{N}$ and $\wt{N}$ is an algebraic complement of $N\cap W$ in $N$, we conclude that $a=0$.  
Moreover,  by  Lemma  \eqref{ntilde},   $\wt{N}$ is in good position to $C$. 
Hence,  recalling  Lemma \ref{roxy}, we have proved  the following result.

\begin{lemma}
The intersection $C\cap\wt{N}$ is a closed convex cone in
$\wt{N}$. If  $a$ is a generator of an extreme ray in $C\cap \wt{N}$, then  
\begin{equation}\label{Nsigma}
\dim \wt{N}-1 =\#{\sigma_a}= d_C(a).
\end{equation}
\end{lemma}

The position of $\wt{N}$ with
respect to $\R^{\Sigma^c}\oplus W$ is described in the next lemma.
\begin{lemma}\label{ll1}
Either $\wt{N}\cap (\R^{\Sigma^c}\oplus W)=\{0\}$ or
$\wt{N}\subset \R^{\Sigma^c}\oplus W$.  In the second case
$\dim \wt{N}=1$ and $\Sigma=\emptyset$.
\end{lemma}

\begin{proof}
Assume that  $\wt{N}\cap (\R^{\Sigma^c}\oplus
W)\neq\{0\}$. Take  a nonzero point $x\in \wt{N}\cap (\R^{\Sigma^c}\oplus W)$. We know that $\wt{N}\cap C$ has a
nonempty interior in $\wt{N}$ and is therefore generated as the convex hull of
its extreme rays by Lemma \ref{kreinmilman}.  Let $a\in C\cap \wt{N}$ be a generator of an
extreme ray $R$. Then $a_i>0$ for all $i\in\Sigma^c$ and  hence $\lambda a +x\in C\cap \wt{N}$ for large  $\lambda >0$. Taking another large  number  $\mu>0$, we get   $\mu a- (\lambda a+x)\in C\cap \wt{N}$. Since $R=\R^+\cdot a$ is an extreme ray, we conclude $\lambda a+x\in \R^+\cdot a$ so that $x\in \R\cdot a$. Consequently,
there is only one extreme ray in $\wt{N}\cap C$, namely
 $R=\R^+\cdot a$ with  $a\in \R^{\Sigma^c}\oplus W$.  Since $\wt{N}\cap C$ has a nonempty interior in $\wt{N}$, we conclude that $\dim \wt{N}=1$  . Hence $\wt{N}=\R\cdot a$ and
 $\wt{N}\subset  \R^{\Sigma^c}\oplus W$. From equation \eqref{Nsigma} we also conclude that $a_i>0$ for all $1\leq i\leq n$.  This in turn implies that $\Sigma=\emptyset$ since  $a\in \R^{\Sigma^c}\oplus W$. The proof of Lemma \ref{ll1} is complete.
\qed \end{proof}

We finally come to the proof of Proposition  \ref{pretzel}.

\begin{proof}[Proposition \ref{pretzel}] 
We consider, according to Lemma \ref{ll1},  two cases.
Starting with the first case  we assume that $\wt{N}\cap (\R^{\Sigma^c}\oplus
 W)=\{0\}$. The projection $P\colon \wt{N}\oplus \wt{N}^{\perp}=\R^{\Sigma}\oplus (\R^{\Sigma^c}\oplus W)\to \R^{\Sigma}$  induces the linear map
 \begin{equation}\label{opl}
 P\colon \wt{N}\rightarrow \R^{\Sigma}.
 \end{equation}
Take $k\in\R^{\Sigma}$ and write
 $(k,0)=n+m\in\wt{N}\oplus\wt{N}^{\perp}$. Since $\wt{N}^{\perp}\subset
\R^{\Sigma^c}\oplus W$,  we conclude that
 $$
 P(n+m)=P(n)=k
 $$
so that $P$ is surjective. If $n\in \wt{N}$ and $P(n)=0$, then $n\in  \wt{N}\cap ({\mathbb
 R}^{\Sigma^c}\oplus W)=\{0\}$ by assumption. Hence the map in
 \eqref{opl} is a bijection.  By  Lemma \ref{nomer},
$C\cap\wt{N}$ is a quadrant in $\wt{N}$. We shall show that $P$ maps the quadrant $C\cap \wt{N}$ onto the standard quadrant $Q^{\Sigma}=[0,\infty )^{\Sigma}$ in $\R^{\Sigma}$.
Let $a$ be a nonzero element in $C\cap \wt{N}$ generating  an extreme ray $R=\R^+\cdot a$. Then, by Lemma \ref{roxy},
 $$
 \dim \wt{N}-1=\sharp\sigma_a, 
 $$
 and since $\sharp \Sigma=\text{dim}\ \wt{N}$ there is
exactly one index $i\in\Sigma$ for which
$a_i>0$. Further,  $a_i>0$ for all $i\in\Sigma^c$ by definition of $\Sigma$.  This implies that there can be
at most $\dim(\wt{N})$-many extreme rays. Indeed, if $a$ and $a'$
generate extreme rays and
$a_i,a_i'>0$  for some $i\in \Sigma$, then $a_k=a_k'=0$ for all
$k\in\Sigma\setminus\{i\}$.  Hence,  from  $a_j>0$ for all $j\in \Sigma^c$,   we conclude $\lambda
a-a'\in C$ for large $\lambda>0$. Therefore,
$a'\in \R^+a$ implying that $a$ and $a'$ generate the same
extreme ray. As a consequence,  $\wt{N}\cap C$ has  precisely $\dim \wt{N}$-many extreme rays because $\wt{N}\cap C$ has a nonempty interior in view of Lemma \ref{ntilde}.  Hence  the map $P$ in \eqref{opl}  induces an isomorphism
$$
(\wt{N}, \wt{N} \cap C)\rightarrow (\R^\Sigma,Q^\Sigma).
$$
This implies that $(N,C\cap N)$ is isomorphic to
$\bigl( {\R}^{\dim(N)}, [0,\infty)^{\sharp\Sigma}\oplus {\R}^{\dim(N)-\sharp\Sigma}\ \bigr).$

In the second case we assume  that $\wt{N}\subset \R^{\Sigma^c}\oplus W$.  From Lemma  \ref{ll1},  $\Sigma=\emptyset$ and $\wt{N}=\R\cdot a$ for an element $a\in
C\cap\wt{N}$  satisfying  $a_i>0$ for all $1\leq i\leq n$. Hence
$(\wt{N},\wt{N}\cap C)$ is isomorphic to
$( {\R}, {\R}^+)$ and therefore $(N,N\cap C)$ is isomorphic to
$( {\R}, \R^+)$ since in this case $N=\wt{N}$. The proof of Proposition \ref{pretzel} is complete.

\qed \end{proof}

We would like to add two results, which makes use of the previous discussion.

\begin{proposition}\label{big-pretzel}
Let $N$ be a sc-smooth finite-dimensional subspace of $E=\R^n\oplus W$ in good position to the partial quadrant $C=[0,\infty)^n\oplus W$.
If $x\in N\cap C$, then 
\begin{itemize}
\item[{\em (1)}]\ $d_{N\cap C}(x)= d_C(x)$  if  $x\not \in N\cap W$.
\item[{\em (2)}]\ $d_{N\cap C}(x)= \dim(N)-\dim(N\cap W)$ if $x\in N\cap W$. 
\end{itemize}
Here we identify $W$ with $\{0\}^n\oplus W$.
\end{proposition}

\begin{proof}
We make use of the previous notations and distinguish the two case $\dim N=1$ and $\dim N>1$.

If $\dim N=1$, then $N$ is spanned by a 
vector $e=(a_1,\ldots ,a_n, w)\in \R^n\oplus W$ in which $a_j>0$ for all $1\leq j\leq n$. This  implies that $(N, N\cap C)$ is isomorphic to $(\R, [0,\infty))$. 
If $x\in N\cap C$, then $x=te$ for some $t\geq 0$. If, in addition, 
$x\not \in N\cap W$, then $t>0$ and hence $d_{N\cap C}(x)=0=d_C(x)$. If  however, $x\in N\cap W$, then $t=0$, implying that $d_{N \cap C}(x)=1=\dim N.$

Now we assume that $\dim N>1$ and that $\wt{N}$ is the algebraic complement of $N\cap W$ in $N$ so that 
$N=\wt{N}\oplus (N\cap W)$. We may assume, 
after a linear change of coordinates, that $\wt{N}$ is represented as follows.
\begin{itemize}
\item[$(\ast)$]\ The linear subspace $\wt{N}$ is spanned by the vectors $e^j\in \R^n\oplus W$,  for $1\leq j\leq m=\dim \wt{N}$ of the form 
$e^j=(a^j, b^j, w^j)\in \R^m\oplus \R^{n-m}\oplus W$ where $a^j$  is  the standard basis vector in $\R^m$, $b^j=(b^j_1, \ldots ,b^j_{n-m})\in \R^{n-m}$ satisfies  $b^j_i>0$ for all $1\leq i\leq n-m$, and $w^j\in W$. 
\end{itemize}

If $e^j=(0,w^j)\in \{0\}^n\oplus W$ for $m+1\leq j\leq k$,$k=\dim N$,  is a basis of the finite-dimensional subspace $N\cap W$, then 
the vectors $e^1,\ldots ,e^k$ form a basis of $N=\wt{N}\oplus (N\cap W)$.  It follows that 
$(N, N\cap C)$ is isomorphic to $(\R^k, [0,\infty )^m\oplus \R^{k-m})$.

Now we assume that $x\in N\cap C$. 
 If $x\in N\cap C$, then 
\begin{equation}\label{algebraic_eq}
x=\sum_{j=1}^m\lambda_je^j+\sum_{j=m+1}^k\lambda_j e^j\in \wt{N}\oplus(N\cap W)
\end{equation} 
with $\lambda_j\geq 0$ for all $1\leq j\leq m$ and $d_{N\cap C}(x)=\#\{1\leq j\leq m\ \vert \ \lambda_j=0\}$. If, in addition, $x\in N\cap W$, then $\lambda_j=0$ for all $1\leq j\leq m$, and we conclude that 
$d_{N\cap C}(x)=m=\dim N-\dim (N\cap W)$.

If $x\in ( N\cap C)\setminus (N\cap W)$, then there is at least one index $j_0$ in $1\leq j_0\leq m$  for which $\lambda_{j_0}>0$. This implies, in view of $(\ast)$  and \eqref{algebraic_eq} that 
$$x=\sum_{j=1}^me^j=(x_1, \ldots, x_n,w)\in \R^n\oplus W$$
satisfies $x_s>0$ for all $m+1\leq s\leq n$ since $b^j_i>0$ for  all $1\leq i\leq n-m$ and $m+1\leq j\leq n$. Using that $a^j$ in $(\ast)$ are vectors of the standard basis in $\R^m$,  it follows  for $1\leq s\leq m$
that $x_s=0$ if and only if $\lambda_s=0.$ This shows that $d_{N\cap C}(x)=d_C(x)$ if $x\not \in N\cap W$ and completes the proof of Proposition \ref{big-pretzel}. 

\qed \end{proof}

\begin{lemma}\label{big-pretzel_1a}
Let $N$ be a sc-smooth finite-dimensional subspace of $E=\R^n\oplus W$ in good position to the partial quadrant $C=[0,\infty)^n\oplus W$ whose good complement in $E$ is $Y$.  We assume that $\tau\colon N\cap C\to Y_1$ is  a $C^1$-map satisfying $\tau (0)=0$ and $D\tau (0)=0$.  Then 
$$d_{N\cap C}(v)=d_{C}(v+\tau (v))$$
for all $v\in N\cap C\setminus N\cap W$ close to $0$.
\end{lemma}

\begin{proof}

Since $N$ is in good position to $C$ and $Y$ is a good complement of $N$ is $E$, there exists a constant 
exists a constant  $\gamma>0$ such that if $n\in N$ and $y\in Y$ satisfy $\abs{y}_0\leq \gamma \abs{n}_0$, then 
$n\in C$ if and only if $n+y\in C$. It follows, in view of $\tau (0)=0$ and $D\tau (0)=0$, that $\abs{\tau (v)}_0\leq \gamma\abs{v}_0$ for $v\in N\cap C$ close to $0$. 

Now proceeding as in the proof of Proposition \ref{big-pretzel}, we 
distinguish the two cases $\dim N=1$ and $\dim N\geq 2$.

In the first case, $N$ is spanned by a vector $e=(a, w)\in \R^n\oplus W$ in which $a=(a_1,\ldots ,a_n)$ satisfies $a_j>0$ for all $1\leq j\leq n$. Take  $v\in N\cap C\setminus N\cap W$, then $v=t(a, w)$ for $t>0$, so that $d_{N\cap C}(v)=0$. 
We  already know that $v_j+\tau_j(v)\geq 0$ for all $1\leq j\leq n$. So, to prove the lemma, it suffices to show that 
$v_j+\tau_j(v)>0$ for $1\leq j\leq n$ and $v\in N\cap C\setminus N\cap W$ close to $0$.  Arguing by contradiction we assume that there exists a sequence $t_k\to 0$ such  that if $v^k=t_ka$, then $v^k_i+\tau_i(v^k)=0$ for some index $1\leq i\leq n$. Then 
$$0\neq \dfrac{a}{\abs{a}_0}=\dfrac{v_i^k}{\abs{v^k}_0}=-\dfrac{\tau_i(v^k)}{\abs{v^k}_0},$$
which leads to   a contradiction since the right-hand side converges to $0$ in view of $\tau (0)=0$ and $D\tau (0)=0$. Consequently, $d_C(v+\tau (v))=0$ for all $v\in  N\cap C\setminus N\cap W$ which are close to $0$.\par

In the case $\dim N\geq 2$, we use the notation introduced in the proof of Proposition \ref{big-pretzel}. If $v=(v_1,\ldots ,v_k, w)\in 
 N\cap C\setminus N\cap W$, then, in view of (1) in the proof of Proposition \ref{big-pretzel}, $d_{N\cap C}(v)=\#\{1\leq i\leq k\, \vert \, v_i=0\}$ and  there exists at least one index $1\leq i\leq k$ for which  $v_i>0$. This implies,  in particular,  that 
 $v_j>0$ for all $k+1\leq j\leq n$.
 Moreover,  in view of Lemma \ref{cone1}, the good complement $Y$ of $N$ is contained in $\{0\}^m\oplus \R^{n-m}\oplus W$. From this we conclude that $v_j+\tau_j (v)=v_j$ for all $1\leq j\leq m$. Hence in order to finish the proof it suffices to show that $v_j+\tau_j(v)>0$ for all $m+1\leq j\leq n$. In order to verify  this we introduce the sc-linear map $T\colon N\to \R^k$ defined by 
 $T(e^l)=\ov{e}^l$ for $1\leq l\leq k:=\dim N$ where  $\ov{e}^1, \ldots ,\ov{e}^k$ of $\R^k$. Then $T(N\cap C)=C':=[0,\infty )^m\oplus \R^{k-m}$ and  $d_{N\cap C}(v)=d_{C'}(T(v)).$

We consider  the map $g\colon C'\to \R^n\oplus W$,  defined by 
$$g(v')=T^{-1}(v')+\tau (T^{-1}(v'))$$ 
for $v'=(v_1', \ldots ,v'_k)\in C'$ close to $0$. Since $v+\tau (v)\in C$, $g(v')\in C$ for $v'\in C'$ close to $0$. 
We prove our claim by showing that $g_j(v')>0$ for $m+1\leq j\leq n$  and nonzero $v'\in C'$ close to $0$.
Arguing by contradiction we assume that there exists a point $v'\in C'$  different from  $0$ at which $g_j(v')=0$ for some $m+1\leq j\leq n$. 
We have  $v_i'>0$ for some $1\leq i\leq m$ so that  $v'+t\ov{e}^i\in C'$ for $\abs{t}$ small.  Then we compute, 
$$\dfrac{d}{dt}g_j(v'+t\ov{e}^i)\vert_{t=0}=b^i_j+D\tau_j(v)T^{-1}\ov{e}^i=b^i_j+D\tau_j(v)e^i.$$
Since the map $\tau$ is of class $C^1$ and $D\tau (0)=0$ and $b^i_j>0$ we conclude that the derivative is positive, and hence the function $t\mapsto g_j(v'+t\ov{e}^i)$ is strictly increasing for $\abs{t}$ small.  By assumption $g_j(v')=0$, hence  $g_j(v'+t\ov{e}^i)<0$ for $t<0$ small, contradicting  $g_j(v'+t\ov{e}^i)\geq 0$ for $\abs{t}$ small.

Consequently, $v_j+\tau_j(v)>0$ for all $m+1\leq j\leq n$ and all $v\in N\cap C\setminus N\cap W$ sufficiently close to $0$ and since $v_j+\tau_j(v)=v_j$ for $1\leq j\leq m$, we conclude $d_{N\cap C}(v)=d_C(v+\tau (v))$ 
for $v\in N\cap C\setminus N\cap W$ close to $0$. The proof of Lemma \ref{big-pretzel_1a} is complete. 

\qed \end{proof}

\subsection{Proof of Lemma \ref{new_lemma_Z}}\label{pretzel-B}
We recall the statement for the convenience of the reader

\begin{lemma}\label{new_lemma_Z_1}
The kernel  $\wt{K}:=\ker DH(0)\subset \R^{n}$ of the linearization $DH(0)$ is in good position to the partial quadrant $\wt{C}=[0,\infty)^{k}\oplus \R^{n'-k}$ in $\R^{n'}$. Moreover, there exists a good complement $Z$ of $\wt{K}$ in $\R^{n'}$, so that 
$\wt{K}\oplus Z =\R^{n'}$, having the property that $Z\oplus W'$ is a good complement of $K'=\ker Dg(0)$ in $E'$, 
$$E'=K'\oplus ( Z\oplus W').
$$
\end{lemma}
\begin{proof}
The space $E'=\R^{n'}\oplus W'$ is equipped with the norm $\abs{\cdot }_0$. We use the  equivalent norm defined by $\abs{(a, w)}=\max\{ \abs{a}_0,\abs{w}_0\}$ for $(a, w)\in \R^{n'}\oplus W'$. The kernel $K'=\ker Dg(0)\subset \R^{n'}\oplus W'$ is in good position to the partial quadrant $C'=[0,\infty )^k\oplus \R^{n'-k}\oplus W'$ and $Y'$ is a good complement of $K'$ in $E'$, so that $K'\oplus Y'=E'$.

We recall that by the definition of good position, there exists $\varepsilon>0$ such that,  
if  $k\in K'$ and $y\in Y'$ satisfy 
\begin{equation}\label{again_good_1}
\abs{y}\leq \varepsilon \abs{k},
\end{equation}
then 
\begin{equation}\label{again_good_2}
\text{$k+y\in C'$\quad  if and only if\quad  $k\in C'$.}
\end{equation}

Let $Y_0$ be an algebraic complement of $Y'\cap W'$ in  $Y'$, where we have  identified  $W'$ with $\{0\}^{n'}\oplus W'$, so that $Y'=Y_0\oplus (Y'\cap W').$

We claim that the projection $P'\colon \R^{n'}\oplus W'\to \R^{n'}$, restricted to $Y_0$,  is injection. Indeed, assume that $x=(a, w)\in Y_0\subset \R^{n'}\oplus W'$ satisfies that  $P'(x)=0$. Then $a=0$ and $x=(0, w)$. Since $x=(0, w)\in Y'\cap W'$ and $Y_0\cap (Y'\cap W')=\{0\}$, we conclude that $x=0$ so  that indeed $P'\vert Y_0$ is an injection.

Introducing the subspace  $Z'=P'(Y_0)$, we claim that $\R^{n'}=\wt{K}+Z'$.  To verify the claim, we take $a\in \R^{n'}$ and let $x=(a, 0)$. Since $\R^{n'}\oplus W'=K'\oplus Y'$, there are unique elements $k\in K'$ and $y\in Y'$ such that $x=k+y$. In view of Lemma 3.44,
the element $k\in K'$ is of the form $k=(\alpha, D\delta (0)\alpha)$ for a unique $\alpha \in \wt{K}$. The element $y\in Y'$ is of the form $y=(b, w)\in \R^{n'}\oplus W'$. Hence 
$x=(a, 0)=k+y=(\alpha , D\delta (0)\alpha)+(b, w)=(\alpha+b, D\delta (0)\alpha+w)$, showing $a=\alpha+b$. On the other hand, since $Y'=Y_0\oplus (Y'\cap W')$,  we have  $(b, w)=(b, w_1)+(0, w_2)$ for $w_1, w_2\in W'$ such that $(b, w_1)\in Y'$ and $(0, w_2)\in Y'\cap W'.$ Hence $b=P'(b, w_1)\in Z'$ and, therefore, $a=\alpha +b\in \wt{K}+Z'$ as claimed.

From this it follows that  $\dim Z'\geq n'-\dim \wt{K}$ and we  choose a subspace $Z$ of $Z'$ of dimension  $\dim Z=n'-\dim \wt{K}$, so that $\R^{n'}=\wt{K}\oplus Z$.  Recalling that  the projection $P'\colon Y_0\to Z'$ is an isomorphism, we have $(P')^{-1}(Z)\subset Y'$.

We shall prove $\wt{K}$ is in good position to the partial quadrant $\wt{C}=[0,\infty )^k\oplus R^{n'-k}$ in $\R^{n'}$ and $Z$ is a good complement of $\wt{K}$, so that 
$\wt{K}\oplus Z=\R^{n'}$.  In view of  the fact that the map $(P')^{-1}\colon Z'\to Y_0$ is an isomorphism and $Z\subset Z'$, there exists a constant $A$ such that 
\begin{equation}\label{again_good_3}
\text{$\abs{(P')^{-1}(z)}\leq A\abs{z}$\quad  for all $z\in Z$.}
\end{equation}
Moreover, by Lemma \ref{new_lemma_relation}, the map $L\colon K'\to \wt{K}$, defined by $L(\alpha , D\delta (0)\alpha)=\alpha$,  is an isomorphism, and hence there exists a constant $B$ such that 
\begin{equation}\label{again_good_4}
\text{$\abs{\alpha}=\abs{L(\alpha ,D\delta  (0)\alpha)}\leq B\abs{(\alpha ,D\delta (0)\alpha)}$ \quad  for all $\alpha\in \wt{K}$.}
\end{equation}

We choose $\varepsilon'>0$ such that $\varepsilon' \cdot A\cdot B<\varepsilon$, where 
$\varepsilon>0$ is the constant from the condition \eqref{again_good_1}. 
We claim that if  $z\in Z$ and $\alpha\in \wt{K}$ satisfy 
\begin{equation}\label{again_good_5}
\abs{z}\leq \varepsilon' \abs{\alpha}, 
\end{equation} 
then 
\begin{equation}\label{again_good_6}
\text{$z+\alpha \in \wt{C}$\quad  if and only if \quad $\alpha \in \wt{C}.$}
\end{equation}
Since $z\in Z\subset Z'$, there exists a unique $w\in W'$ such that $(P')^{-1}(z)=(z, w)\in Y_0$. 
Then we estimate,  using \eqref{again_good_3}, \eqref{again_good_5}, and then \eqref{again_good_4}, 
$$\abs{(z, w)}=\abs{(P')^{-1}(z)}\leq A\abs{z}\leq \varepsilon' A\cdot B\abs{(\alpha ,D\delta  (0)\alpha)}<\varepsilon  \abs{(\alpha ,D\delta (0)\alpha)}.$$
Thus, $(z, w)\in Y_0\subset Y'$ and $(\alpha, D\delta (0)\alpha )\in K'$ satisfy the estimate \eqref{again_good_1} and it follows that 
\begin{equation}\label{again_good_7}
\text{$(z, w)+(\alpha, D\delta (0)\alpha )\in C$ \quad if and only if\quad  $(\alpha, D\delta (0)\alpha)\in C'.$}
\end{equation}
Since $(z, w)+(\alpha, D\delta (0)\alpha)=(z+\alpha ,w+D\delta (0)\alpha)\in \R^{n'}\oplus W'$, 
\eqref{again_good_7} implies that $z+\alpha\in \wt{C}$ if and only if $\alpha \in \wt{C}$, proving our claim \eqref{again_good_6} and we see that $Z$ is a good complement of $\wt{K}$ in $\R^{n'}$.

Next we set 
$$\ov{Y}=Z\oplus W'$$
and claim that $K'\oplus \ov{Y}=\R^{n'}\oplus W'$ and that $\ov{Y}$ is a good complement of $K'$ in $E'=\R^{n'}\oplus W'.$ We first verify that $\R^{n'}\oplus W'=K'+\ov{Y}$. Take $(a, w)\in \R^{n'}\oplus W'$.  Then, since 
$\R^{n'}\oplus W'=K'\oplus Y'$, there are unique elements $(\alpha, D\delta (0)\alpha)\in K'$, where $\alpha \in \wt{K}$,   and $(b, w_1)\in Y'$ such that 
\begin{equation}\label{again_good_8}
(a, w)=(\alpha , D\delta (0)\alpha)+(b, w_1)\in K'\oplus Y'.
\end{equation}
Hence $a=\alpha +b.$ Since $\R^{n'}=\wt{K}\oplus Z$, we may decompose $b$ as $b=\alpha_1+b_1$ with $\alpha_1\in \wt{K}$ and $b_1\in Z$  and then \eqref{again_good_8} can be written as 
\begin{equation*}
\begin{split}
(a, w)&=(\alpha , D\delta (0)\alpha)+(b, w_1)=(\alpha , D\delta (0)\alpha)+(\alpha_1+b_1, w_1)\\
&=(\alpha , D\delta (0)\alpha)+(\alpha_1, D\delta (0)\alpha_1)+(b_1, w_1-D\delta (0)\alpha_1)\\
&=(\alpha+\alpha_1, D\delta (0)(\alpha+\alpha_1))+(b_1, w_1-D\delta (0)\alpha_1)
\end{split}
\end{equation*}
where  $(\alpha+\alpha_1, D\delta (0)(\alpha+\alpha_1))\in K'$ (since $\alpha+\alpha_1\in \wt{K}$) and $(b_1, w-D\delta (0)\alpha_1)\in Z\oplus W'=\ov{Y}$ (since $b_1\in Z$). Hence $\R^{n'}\oplus W'=K'+\ov{Y}$ as claimed.

If $(\alpha ,D\delta (0)\alpha)\in K'\cap \ov{Y}=K'\cap (Z\oplus W')$, then $\alpha \in \wt{K}\cap Z$, so that $\alpha=0$ and $K'\cap \ov{Y}=\{0\}$ and we  have proved that  $Z\oplus \ov{Y}=\R^{n'}\oplus W'$.

Finally, we shall show that $\ov{Y}=Z\oplus W'$ is a good complement of $K'$ in $\R^{n'}\oplus W'$.  Recalling the isomorphism $L\colon K'\to \wt{K}$, $L(\alpha, D\delta (0)\alpha)=\alpha$, there exists a constant $M$ such that $\abs{L^{-1}(\alpha)}\leq M\abs{\alpha}$, i.e.,  
$$\text{$\abs{(\alpha ,D\delta (0)\alpha)}\leq M\abs{\alpha}$\quad  for all $\alpha\in \wt{K}$.}$$ 

We choose $\varepsilon''>0$ such that $\varepsilon'' M<\varepsilon'$ where $\varepsilon'$ is defined in \eqref{again_good_5}. We shall show that if 
$(\alpha, D\delta (0)\alpha)\in K'$ and $(a, w)\in Z$
satisfy the estimate
\begin{equation}\label{again_good_9}
\abs{(a, w)}\leq \varepsilon'' \abs{(\alpha, D\delta (0)\alpha)},
\end{equation}
then 
\begin{equation}\label{again_good_10}
\text{$(a, w)+(\alpha, D\delta (0)\alpha)\in C'$\quad if and only if $(\alpha, D\delta (0)\alpha)\in C'$.}
\end{equation}
Since $\abs{(a, w)}=\max\{\abs{a}_0, \abs{w}_0\}\geq \abs{a}$, we conclude from  \eqref{again_good_10} and 
\eqref{again_good_9} the estimate
$$\abs{a}\leq \abs{(a, w)}\leq \varepsilon'' \abs{(\alpha, D\delta (0)\alpha)}\leq \varepsilon'' M\abs{\alpha}.$$
In view of  \eqref{again_good_5} and 
\eqref{again_good_6}, 
$$\text{$a+\alpha \in \wt{C}$\quad if and only if $\alpha \in \wt{C}$}.$$
This is equivalent to 
$$\text{$(a+\alpha, w+D\delta (0)\alpha)\in C'$ \quad if and only if $(\alpha, D\delta (0)\alpha)\in C'.$}$$
Hence \eqref{again_good_10} holds and the proof that $\ov{Y}=Z\oplus W'$ is a good complement of $K'$ is complete. The proof of Lemma \ref{new_lemma_Z} is finished.

\qed \end{proof}

\subsection{Proof of Lemma \ref{good_pos} }\label{geometric_preparation}
We recall the statement for the convenience of the reader.
\begin{lemma}\label{good_pos_1}
Let $C\subset E$ be  a partial quadrant in the  sc-Banach space $E$ and $N\subset E$ a finite-dimensional smooth subspace in good position to $C$ and let $Y$ be a good complement of $N$ in $E$, so that $E=N\oplus Y$.
We assume that $V\subset N\cap C$ is a relatively open neighborhood
of $0$ and $\tau\colon V\rightarrow Y_1$ a map of class $C^1$ satisfying $\tau(0)=0$ and $D\tau(0)=0$.

Then there exists a relatively open neighborhood $V'\subset V$ of $0$ such that the following holds.
\begin{itemize}
\item[{\em (1)}]\ $v+\tau(v)\in C_1$ for $v\in V'$.
\item[{\em (2)}]\ For every $v\in V'$,  the Banach space $Y=Y_0$ is a topological complement of the linear subspace $N_v=\{n+D\tau(v)n\, \vert \, n\in N\}$.
\item[{\em (3)}]\ For every $v\in V'$,  there exists a  constant $\gamma_v>0$ such  that  if  $n\in N_v$ and $y\in Y$ satisfy 
$\abs{y}_0\leq \gamma_v\cdot \abs{n}_0$,  then  $n\in C_x$ if and only if  $n+y\in C_x$ are equivalent, where $x=v+\tau (v)$. 
\end{itemize}
\end{lemma}

\begin{proof}[{\bf Proof of Lemma \ref{good_pos}}]
We shall use the notations of Proposition 
\ref{big-pretzel}
Since  the smooth finite dimensional subspace $N$ of $E$ is in good position to $C$,  and $Y$ is its good complement so that  $E=N\oplus Y$,  there exists a constant $\gamma>0$ such that if $n\in N$ and $y\in Y$ satisfy
\begin{subequations}
\begin{gather}
\abs{y}_0\leq \gamma \abs{n}_0\label{eq_lemma3.55_1}\\
\text{then}\ \ 
\text{$y+n\in C$\quad if and only if\quad  $n\in C$.} \label{eq_lemma3.55_2}
\end{gather}
\end{subequations}
\noindent(1)\,  The $C^1$-map $\tau\colon V\to Y_1$ satisfies $\tau (0)=0$ and $D\tau (0)=0$, so that  
$$\lim_{\abs{v}\to 0}\dfrac{\abs{\tau (v)}_0}{\abs{v}_0}=0.
$$
Consequently, there exists a relatively open neighborhood $V'\subset V$ of $0$ in $N$ such that $\abs{\tau (v)}_0\leq \gamma \abs{v}_0$ for all $v\in V'$.  Since $V'\subset V\subset C$, we conclude from \eqref{eq_lemma3.55_2} that 
$v+\tau (v)\in C$ for all $v\in V'$.  In addition, since  $V'$ consists of smooth points and $\tau (v)\in Y_1$,  we have $v+\tau (v)\in C_1$ for all $v\in V'$.\par

\noindent  (2)\,  By assumption,  $E=N\oplus Y$ and $D\tau (v)n\in Y$ for all $n\in N$. If   $x=n+y\in N\oplus Y\in E$, then $y=(n+D\tau (v)n)+(y-D\tau (v)n)$, showing that $E=N_v +Y$.  If $(n+D\tau (v)n)+y=(n'+D\tau (v)n')+y'$ for some $n, n'\in N$ and $y, y'\in Y$, then $n-n'=(y'-y)+D\tau (v)(n'-n)\in N\cap Y$. Hence $n=n'$ and $y=y'$, showing that $N_v\cap Y=\{0\}$ and hence $E=N_v\oplus Y$. That the subspace $Y$ is a topological complement of $N_v$ follows from the fact that $N_v$ is a finite dimensional subspace of $E$. \par

\noindent (3)\,  
Without loss of generality we may assume $E=\R^n\oplus W$ and $C=[0,\infty)^n\oplus W$. 
We recall that for  $x\in C$, the set $C_x$ is defined by   
$$
C_x=\{(a, w)\in \R^n\oplus W\ \vert \ \text{$a_i\geq 0$ for all $1\leq i\leq n$ for which $x_i=0$}\},
$$ 
and if $x_i>0$ for all $\leq i\leq n$, we set $C_x=E$. Clearly, if $x=(0, w)\in \R^n\oplus W$, then $C_x=C$.\par

If $v=0$, then $x=v+\tau (v)=0$, so that $C_x=C$ and the statement of the corollary follows from the fact that $N_v=N$ and $N$ is in good position to $C$ having the good complement $Y$.

In order to prove the statement for $v\neq 0$ we distinguish thew  two cases $\dim N=1$ and $\dim N>1$. In 
the first case, the proof of Lemma \ref{big-pretzel_1a} shows that $v_j+\tau_j (v)>0$ for all $1\leq j\leq n$.
This implies that $C_x=E$ and that, in view of (2),  $Y$ is a good complement of $N_v$ with respect to $C_x$.\par

 In the second case $\dim N>1$, and we denote by $\wt{N}$ the  algebraic complement of $N\cap W$ in $N$,  where $W$ is identified with $W=\{0\}^n\oplus W$,  so that $N=\wt{N}\oplus (N\cap W)$.   We may assume,  after a linear change of coordinates, that the following holds. 
\begin{itemize}
\item[(a)]\ The linear subspace $\wt{N}$ is spanned by the vectors $e^j$ for $1\leq j\leq m=\dim \wt{N}$ of the form 
$e^j=(a^j, b^j, w^j)$ where $a^j$ are the vectors of the standard basis in $\R^m$, $b^j=(b^j_1, \ldots ,b^j_{n-m})$ satisfy $b^j_i>0$ for all $1\leq i\leq n-m$, and $w^j\in W$. 
\item[(b)]\ The good complement $Y$ of $N$  in $E$  is contained in $\{0\}^m\oplus \R^{n-m}\oplus W.$
\end{itemize}
Denoting by   $e^j=(0, w^j)\in \{0\}^n\cap W$, $m+1\leq j\leq k=\dim N$, a basis of  $N\cap W$, the vectors $e^1,\ldots ,e^k$  form a basis of $N$.   
We choose $0<\varepsilon<\gamma/2$ where $\gamma$ is the constant from \eqref{eq_lemma3.55_1} and define the open neighborhood $V''$ of $0$ in $N\cap C$ consisting  of points $v\in V'$ satisfying $\abs{D\tau (v)n}_0<\varepsilon\abs{n}_0$ for all $n\in N$.

We consider points  $v\in V'$ belonging to $N\cap W$ where we identify $W$ with $\{0\}^n\oplus W$.  Hence  $v$ is  
of the form $v=(0, w)\in \{0\}^n\oplus W$.  Then $\abs{\tau(v)}_0\leq \gamma \abs{v}_0$ and,  since $v\in C$,  we conclude that $v+\tau (v)$ and $v-\tau (v)$ belong to the partial quadrant $C$. This implies that $\tau_i(v)=0$ for all $1\leq i\leq n$ and  $C_{x}=C$ where $x=v+\tau (v)=(0, w')$. We choose a positive constant $\gamma_v$ such that $\gamma_v(1+\varepsilon)+\varepsilon <\gamma$.  Abbreviating $h=n+D\tau (v)n\in N_v$ for some $n\in N$ and taking $y\in Y$, we assume that 
\begin{equation}\label{eq_lemma3.55_3}
\abs{y}_0\leq \gamma_0\abs{h}_0=\gamma_v \abs{n+D\tau (v)n}_0.
\end{equation}
 
 We claim that $y+h\in C_x$ if and only if $h\in C_x$. Suppose that $y+h=y+n+D\tau (v)n\in C_x=C$. Then we estimate using \eqref{eq_lemma3.55_3}, 
\begin{equation}\label{eeqq1}
\begin{split}
\abs{y+D\tau (v)n}_0&\leq \abs{y}_0+\abs{D\tau (v)n}_0\leq \gamma_v\abs{h+D\tau (v)n}+\abs{D\tau (v)n}_0\\
&\leq (\gamma_v+\varepsilon\gamma_v+\varepsilon)\abs{n}_0\leq \gamma \abs{n}_0.
\end{split}
\end{equation}

Since $y+D\tau (v)n\in Y$ and by assumption $y+h=y+D\tau (v)n+n\in C_x=C$, we deduce from \eqref{eeqq1}  and 
\eqref{eq_lemma3.55_2} that $n\in C$. 
Then the estimate $\abs{D\tau (v)n}_0\leq \varepsilon \abs{n}_0<\gamma\cdot \abs{n}_0$ and again \eqref{eq_lemma3.55_2}, imply that $h=n+D\tau (v)n\in C$, as claimed.\\
Conversely, the  assumption  that $h=n+D\tau (v)n\in C$ implies, in view of  $\abs{D\tau (v)h}_0\leq \abs{h}_0$ and \eqref{eq_lemma3.55_2}, that $n\in C$. Then using  \eqref{eeqq1} we find hat $y+D\tau (v)n +h=y+h\in C$. We have proved that if $h\in N_v$ and $y\in Y$ satisfy $\abs{y}_0\leq \gamma_v\abs{h}_0$, then $y+h\in C_x$ if and only if $h\in C_x$.\par

Finally we consider points $v=(v_1,\ldots ,v_n, w)\in N\cap C\setminus N\cap W$. Then $v_i>0$ for some $1\leq i\leq m$ and $v_j>0$ for all $m+1\leq j\leq k$. Moreover, in view of the proof of the statement (1), 
$\abs{\tau (u)}_0\leq \gamma\abs{v}_0$ which implies that $v+\tau (v)\in C$ if $v$ is close to $0$.
It follows from the proof of Lemma \ref{big-pretzel_1a} that  $v_j+\tau_j (v)>0$ for all $m+1\leq j\leq n$. Abbreviating  $x=v+\tau (v)$ and denoting by $\Lambda_x$ the set of indices $1\leq j\leq n$ for which $x_j>0$,  we conclude 
$$C_x=\{(a, w)\in \R^n\oplus W\, \vert \, \text{$a_j\geq 0$ for all $1\leq j\leq m$ satisfying $j\not \in \Lambda_x$}\}.$$
Recall that $N_v=\{n+D\tau (v)n\, \vert \, n\in N\}$ and $E=N_v\oplus Y$ in view of the statement (2).  Now the inclusion $Y\subset  \{0\}^m\oplus \R^{n-m}\oplus W$ and the definition of $C_x$ above show that, 
if $h=n+D\tau (v)n\in N_v$ and  $y\in Y$, then $y+h\in C_x$ if and only if $h\in C_x$.  This completes the proof of Lemma \ref{good_pos}.

\qed \end{proof}

\subsection{Diffeomorphisms Between Partial 
Quadrants }\label{finite_partial_inverse}

In the sections \ref{finite_partial_inverse} and \ref{implicit_finite_partial_quadrants} 
we shall prove an inverse function theorem and an implicit function theorem for smooth maps defined on partial quadrants in $\R^n$. We recall that the closed subset $C\subset \R^n$ is a partial quadrant, if there exists an isomorphism $L$ of $\R^n$ mapping $C$ onto $L(C)=[0,\infty)^k\oplus \R^{n-k}$ for some $k$.
We begin  with a definition.
\begin{definition}[{\bf Class $C^1$}] Let $C$ be a partial quadrant in $\R^n$ and $f\colon U\rightarrow \R^m$ a map defined on a relatively open subset $U\subset C$. The map  $f$ is said to be of {\bf class  $C^1$}  if,  for every $x\in U$,  there exists a bounded linear map $Df(x)\colon \R^n\rightarrow \R^m$
such that
$$
\lim_{x+h\in U,h\to  0} \dfrac{ \abs{f(x+h)-f(x)-Df(x)h}}{\abs{h}}=0,
$$
and,  in addition,  the map 
$$Df\colon U\to {\mathscr L}(\R^n, \R^m), \quad x\mapsto Df(x),$$
is continuous.
\end{definition}
We note that the tangent 
$TU=U\oplus {\mathbb R}^n$ is a relatively open subset of the partial quadrant $C\oplus \R^n$ and the map 
$Tf\colon TU\rightarrow T\R^m$,  defined by 
$$Tf(x,h)=\bigl(f(x), Df(x)h ),$$
is continuous. If the maps $f$ and $Tf$ are of class $C^1$, then  we say that $f$ is of class $C^2$. Inductively, the map $f$ is $C^k$ if the maps $f$ and $Tf$ are of class $C^{k-1}$.  Finally, $f$  is called smooth 
(or $C^\infty$) if  $f$ is of class $C^k$ for all $k\geq 1$.

We consider the following situation.  Assume  $C$ is  a partial quadrant in $\R^n$ and $U\subset C$  a relatively open neighborhood of $0$ in $C$. Let $f\colon U\to C$ be  a smooth map such that $f(0)=0$ and $Df(0)\colon \R^n\to \R^n$ an isomorphism. 
Even if $Df(0)C=C$ one cannot guarantee that the image $f(O)$  of an open neighborhood $O$ of $0$ in $C$  is an open neighborhood of $0$ in $C$ as Figure \ref{fig:pict6}  illustrates.

\begin{figure}[htb]
\begin{centering}
\def\svgwidth{69ex}
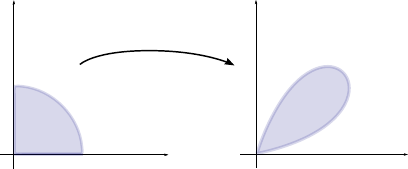
\caption{The boundary tangents at $(0,0)$ are the $x-$ and $y$-axis, but the domain in the right picture is not an open neighborhood of $(0,0)$.}
\label{fig:pict6}
\end{centering}
\end{figure}


Our aim in this section is to derive an inverse and implicit function theorem in this context.

\begin{theorem}[{\bf Quadrant Inverse Function Theorem}]\label{QIFT}
We assume that $C$ is a partial quadrant in ${\mathbb R}^n$ and $U\subset C$ a relatively open neighborhood of $0$, and consider  a smooth map  $f\colon U\rightarrow \R^n$ satisfying $f(U)\subset C$ and having the following properties:
\begin{itemize}
\item[{\em (1)}]\  $f(0)=0$ and $Df(0)C=C$.
\item[{\em (2)}]\  $d_C(x)=d_C(f(x))$ for every $x\in U$.
\end{itemize}
Then there exist two relatively open neighborhoods $U'$ and $V'$ of $0$ in $C$  such that $U'\subset U$ and the map 
$$f\colon U'\rightarrow V'$$ is a diffeomorphism.
\end{theorem}

From $\R^n=C-C$ and $Df(0)C=C$,  it follows that $Df(0)\colon \R^n\to \R^n$ is an isomorphism.  Therefore, it suffices to study the map $g\colon U\to C$, defined by 
$$g(x)=Df(0)^{-1}f(x).$$
The map $g$ has  the same properties (1) and (2), but has the simplifying feature that 
$Dg(0)={\mathbbm1}$. 
Theorem \ref{QIFT} is a consequence of the following proposition.

\begin{proposition}\label{hongkong}
We assume that $C$ is a partial quadrant in ${\mathbb R}^n$ and $U\subset C$ a relatively open convex neighborhood of $0$, and let $f\colon U\rightarrow C$ be a smooth map having the following properties:
\begin{itemize}
\item[{\em (1)}]\ $f(0)=0$.
\item[{\em (2)}]\ $\abs{Df(x)-{\mathbbm 1}}\leq 1/2$ for all $x\in U$.
\item[{\em (3)}]\  $d_C(x)=d_C(f(x))$ for every $x\in U$.
\end{itemize}
Then $V:=f(U)$ is a relatively open neighborhood of $0$ in $C$,  and the map $f\colon U\rightarrow V$ is a diffeomorphism, that is, $f$ and $f^{-1}$ are smooth
in the sense defined above.
\end{proposition}

The proof of Proposition \ref{hongkong} follows from three lemmata proved below.

\begin{lemma}\label{new_lemma_3.73}
With  the assumptions  of Proposition \ref{hongkong} , 
the image  $V=f(U)$ is a relatively open subset of  $C$ and $f\colon U\rightarrow V$ is an open map. Moreover, 
\begin{equation}\label{new_equ_56}
\dfrac{1}{2}\abs{x-x'}\leq \abs{f(x)-f(x')}\leq \dfrac{3}{2} \abs{x-x'}, \quad \text{$x, x'\in U$}.
\end{equation}
\end{lemma}

\begin{proof}

Abbreviating $g(x)=x-f(x)$ and $G(x)= {\mathbbm 1}-Df(x)$, and using  the  convexity of the set $U$, we obtain the identity   
$$g(x')=g(x)+\biggl( \int_0^1G\bigl(\tau x'+(1-\tau )x\bigr)\ d\tau\biggr)\cdot (x'-x)$$
for $x, x'\in U$, from which the desired estimate 
\eqref{new_equ_56} follows, in view of  property (2).
We conclude, in particular,  that the map  $f\colon U\to V$ is injective.

In order to show that the image $V:=f(U)$ is a relatively open subset of $C$, we take a point $y'\in V$. Then $f(x')=y'$ for a unique $x'\in U$ and since 
$U$ is a relatively open subset of $C$, there exists $r>0$ such that $B(x', 3r)\cap C\subset U$. Here  we denoted by  $B(x', 3r)$  the open ball in $\R^n$ centered at $x'$ and having radius $3r$. We claim that $B(y', r)\cap C\subset V$.  In order to prove our claim, we first show that $\{y\in B(x', r)\cap C\, \vert \, d_C(y)=0\}\subset V$. To see this, 
we  take  a point $x_0\in B(x', 3r)\cap C$ such that $d_C(x_0)=0$ and let $y_0=f(x_0)$. By property (3), $d_C(y_0)=0$. Next  we take an arbitrary point $y_1\in B(y', r)\cap C$,  also satisfying $d_C(y_1)=0$,  and consider the points $y_\tau=(1-\tau )y_0+\tau y_1$ for $0\leq \tau \leq 1$.  Abbreviating  $\Sigma=\{\tau \in [0,1]\, \vert \, y_\tau\in V\}$ and $\tau^*=\sup \Sigma$, we assume that $\tau^*<1$.   
We note that $\Sigma$ is non-empty since $0\in \Sigma$ and $d_C(y_{\tau^*})=0$. Then we  choose a sequence $(\tau_n)\subset \Sigma$ such that $\tau_n\to \tau'$ and the  corresponding sequence of points $(y_{\tau_n})$  belonging  to $B(y', r)\cap C$ satisfying $y_{\tau_n}\to y_{\tau^*}$. Since $y_{\tau_n}\in V$, we find points $x_n\in U$ such that $f(x_n)=y_{\tau_n}$ and $d_C(x_n)=d_C(y_{\tau_n})=0$.
By \eqref{new_equ_56},  
$$\abs{x_n-x'}\leq 2\abs{f(x_n)-f(x')}=2\abs{y_{\tau_n}-y'}<2r$$
and 
$$\abs{x_n-x_m}\leq 2\abs{y_{\tau_n}-y_{\tau_m}},$$
which show that $(x_n)$ is a Cauchy sequence belonging to $B(x', 2r).$
Hence $(x_n)$ converges to some point $x^*$  belonging to $\ov{B}(x', 2r).$ By continuity of $f$,  $f(x^*)=y_{\tau^*}$, and the point $x^*$ belongs to the interior of $C$ since, by property (3),  $d_C(x^*)=d_C(y_{\tau^*})=0$.  Now  the classical  inverse function theorem implies that there are  open neighborhoods $W$ and $W'$ of $x^*$ and $y_{\tau^*}=f(x^*)$ both contained in the interior of $C$ such that the map $f\colon W\to W'$ has a continuous inverse $f^{-1}\colon W'\to W$. In particular, if  $\tau>0$ is small, then $y_{\tau^*+\tau}\in V$,  contradicting $\tau^*<1$. Summing up, we have proved our claim that 
$$\{y\in B(y', r)\cap C\, \vert \, d_C(y)=0\}\subset V.$$
Next we take any point $y\in B(y', r)\cap C$ satisfying $d_C(y)\geq 1$. We find sequences $(y_n)\subset B(y', r)\cap C$ and 
$(x_n)=(f^{-1}(y_n))$ satisfying  $d_C(y_n)=d_C(x_n)=0$ and $\abs{y_n-y}\to 0$. 
From  \eqref{new_equ_56}, we  have 
$\abs{x_n-x'}\leq 2\abs{y_n-y'}<2r$ and $\abs{x_n-x_m}\leq 2\abs{y_n-y_m}$
from which we deduce the convergence of the sequence $(x_n)$  to some point $x\in B(x', 3 r)\cap C.$  Hence  $f(x)=y$, implying that $y\in V$. This shows that  $ B(y', r)\cap C\subset V$ and that $V$ is a relatively open subset of $C$.

Finally, to see that $f\colon U\to V$  maps relatively open subsets onto relatively open subsets,  we take an open subset $U'$ of $U$. Since $U'$ can be written as a union of relatively open convex subsets of $C$,   employing the previous arguments we conclude that $V'=f(U')$ is an open subset of $V$. The proof of
 Lemma \ref{new_lemma_3.73} is complete.

\qed \end{proof}

From the above lemma we obtain immediately the following corollary.
 \begin{corollary}\label{f_homeomrphism}
With the assumptions  of Proposition \ref{hongkong},  the map $f\colon U\to V$ is a smooth homeomorphism between the relatively open subsets $U$ and $V$ of $C$.
\qed
 \end{corollary}

We would like to point out that we cannot employ the usual implicit function theorem (it only applies at interior points). Therefore,  we must provide an argument  for the smoothness of the inverse map  $f^{-1}$. One could try to avoid work,  by first showing that $f$ can be extended to a smooth map defined on an open neighborhood
 of $0$ in ${\mathbb R}^n$. However,  our notion of smoothness does not stipulate that $f$ is near every point the restriction
 of a smooth map defined on an open subset of ${\mathbb R}^n$. This would complicate 
 the construction of a smooth extension.

 \begin{lemma}
With  the assumptions of Proposition \ref{hongkong},  the iterated tangent map $T^kf:T^kU\rightarrow T^kV$ is 
 a smooth homeomorphism.
 \end{lemma}
 
 \begin{proof}

The set $T^kC$ is a partial quadrant in $T^k\R^n$ and the sets $T^kU$ and $T^kV$ are relatively open subsets of $T^kC$. 
Then the iterated tangent map $T^kf\colon T^kU\rightarrow T^kV$ is smooth since the map $f\colon U\to V$ is smooth.  Hence we only have to show that $T^kf$ has an inverse 
$(T^kf)^{-1}\colon T^kV\to T^kU$ which is continuous. 
 In order to prove this we proceed by induction starting with $k=1$.  
 We introduce the map $\Phi_1\colon TV\to TU$, defined by 
 $$\Phi_1(y,l)=\bigl( f^{-1}(y),[ Df(f^{-1}(y)) ]^{-1}l).$$
 The map $\Phi_1$ is continuous since,  by Corollary \ref{f_homeomrphism},  the map $f\colon U\to V$ is a homeomorphism.
Moreover, 
 $$
 \Phi_1\circ Tf (x,h) = \bigl(f^{-1}(f(x)),[Df(f^{-1}(f(x))) ]^{-1}\circ Df(x)h\bigr)=(x,h), 
 $$
and similarly $Tf\circ \Phi_1=\mathbbm{1}$.  Hence $\Phi_1$ is an inverse of the tangent map $Tf$.  This  together with the continuity of $Tf$ and $\Phi_1$ show that $Tf$ is a homeomorphism, as claimed.

Now we assume that result holds for $k\geq 1$ and we show that 
$$
T^{k+1}f\colon T^{k+1}U\to T^{k+1}V
$$
 is a homeomorphism.  Recalling  that $T^{k+1}U=T^kU\oplus T^k\R^n$, we define the map 
\begin{gather*}
\Phi_{k+1}\colon T^kV\oplus T^k\R^n \to T^kU\oplus T^k\R^n, \\
\Phi_{k+1}(y, l)=\bigl((T^kf)^{-1}(y), [D(T^kf)((T^kf)^{-1}(y))]^{-1}l\bigr).
\end{gather*}
By the inductive assumption, the map $T^kf\colon T^kV\to T^kU$ is a homeomorphism and hence $\Phi_{k+1}$ is well-defined and continuous. 
Moreover, 
\begin{equation*}
\begin{split}
\Phi_{k+1}&\circ T^{k+1}f (x, h)=\Phi_{k+1}\bigl(T^kf(x), (D(T^kf))(x)h\bigr)\\
&=\bigl((T^kf)^{-1}\bigl(T^kf(x)\bigl),  [D(T^kf)((T^kf)^{-1}(T^kf(x)   ))]^{-1}(D(T^kf))(x)h\bigr)\\
&=\bigl(x,  [D(T^kf)(x)]^{-1}(D(T^kf))(x)h\bigr)=\bigl(x, h),
\end{split}
\end{equation*} 
and similarly $ T^{k+1}f\circ \Phi_{k+1}=\mathbbm{1}$. Hence the map $\Phi_{k+1}$ is a continuous inverse of the iterated tangent map $T^{k+1}f$. This completes the inductive step and the proof of the lemma.

\qed \end{proof}

\begin{lemma}\label{addd}
Let  $U$ be  a relatively open subset of a partial quadrant $C\subset {\mathbb R}^n$ and $f\colon U\rightarrow C$ a smooth map satisfying the assumptions of Theorem \ref{hongkong}.  We assume that the image $V=f(U)$ is a relatively open subset of $C$ and $f\colon U\rightarrow V$ a homeomorphism such that $Df(x)$ is invertible at every point $x\in U$. 
Then the inverse map $g=f^{-1}$ is of class $C^1$
 and its derivative at the point $y\in V$ is given by
 $$
 Dg(y) =[Df(g(y))]^{-1}.
 $$
 \end{lemma}
 
\begin{proof}

By assumption,  the map $f\colon U\to V$ is a homeomorphism. This implies that $g\colon V\to U$ is also a homeomorphism and the map 
$V\to {\mathscr L}(\R^n)$, defined by 
$$y\mapsto [Df(g(y))]^{-1},$$
is continuous.  Next we take  $y_0\in V$ and $h\in \R^n$ satisfying  $y_0+h\in V$. 
Then there are unique points  $x_0$ and $x_0+\delta\in U$ such that $f(x_0)=y_0$ and $f(x_0+\delta )=y_0+h.$
Recalling the abbreviation $g=f^{-1}$ we note that 
$g(y_0)=x_0$, $g(y_0+h)-g(y_0)=\delta$ and $f(x_0+\delta)-f(x_0)=h$ and compute,
\begin{equation*}
\begin{split}
&\dfrac{1}{\abs{h}}
\abs{g(y_0+h)-g(y_0)- Df(x_0)^{-1}h}\\
&\quad= \dfrac{1}{\abs{h}}\abs{\delta - Df(x_0)^{-1}(f(x_0+\delta )-f(x_0))}\\
&\quad =\dfrac{\abs{\delta}}{\abs{h}}\cdot 
\dfrac{1}{\abs{\delta}}\cdot  \abs{Df(x_0)^{-1}\bigl[ f(x_0+\delta)-f(x_0)-Df(x_0)\delta\bigr] }\\
\end{split}
\end{equation*}

From \eqref{new_equ_56} in Lemma \ref{new_lemma_3.73} we obtain 
$\frac{1}{2}\abs{\delta}\leq \abs{h}\leq \frac{3}{2}\abs{\delta}$ so that $\abs{h}\to 0$ if and only if $\abs{\delta}\to 0$. 
The limits vanish since $f$ is differentiable at $x_0$. The proof of  the lemma is finished.

\qed \end{proof}
 
Now we are in the position to prove Proposition \ref{hongkong}.

\begin{proof}[Proof of Proposition \ref{hongkong}]
Under the hypotheses of the proposition the previous discussion shows for every $k$,  
that $T^kf:T^kU\rightarrow T^kV$ satisfies the hypotheses of Lemma \ref{addd} for a suitable choice 
of data, i.e.  taking $T^kC$ as partial quadrant in $T^k{\mathbb R}^n$. Hence we conclude that $(T^kf)^{-1}$ is $C^1$,
which precisely means that $T^kg$ is $C^1$.  In  other words,  $g$ is of class $C^{1+k}$. Since $k$ is arbitrary we conclude that $g$ is $C^\infty$.
\qed \end{proof}

\subsection{An Implicit Function Theorem in Partial Quadrants}\label{implicit_finite_partial_quadrants}

We shall prove a version for the classical implicit function theorem for  maps defined on  partial quadrants.

\begin{theorem}\label{help-you}

We assume that $U$ is relatively open neighborhood of $0$ in the  partial quadrant $C=[0,\infty)^k\oplus {\mathbb R}^{n-k}$ in $\R^n$, and consider a map 
$f\colon U\rightarrow {\mathbb R}^N$ of class $C^j$, $j\geq 1$,  satisfying $f(0)=0$. Moreover, we  assume that  
$Df(0)\colon \R^n\to \R^N $ is surjective  and the kernel $K:=\ker Df(0)$ is in  good position to $C$, and let  
$Y$ be a good  complement of $K$ in $\R^n$, so that $\R^n=K\oplus Y$. 

Then  there exist an open neighborhood $V$ of $0$ in the partial quadrant $K\cap C$, and a map  
$\tau\colon V\rightarrow Y$ of class $C^j$,  and  positive constants $\varepsilon, \sigma$ having the following properties.
\begin{itemize}
\item[{\em (1)}]\ $\tau(0)=0$, $D\tau(0)=0$, and $k+\tau(k)\in U$ if $k\in V$.
\item[{\em (2)}]\ If $k\in V$  satisfies  $\abs{k}\leq \varepsilon$, then $\abs{\tau(k)}\leq \sigma$.
\item[{\em (3)}]\ If $f(k+y)=0$ for $k\in V$ satisfying 
$\abs{k}\leq \varepsilon$  and $y\in Y$ 
satisfying $\abs{y}\leq\sigma$, then $y=\tau (k)$.
\end{itemize}

\end{theorem}

The proof will follow from several lemmata, where we shall us the notations 
$$x=(k, y)=k+y\in K\oplus Y$$
for $k\in K$ and $y\in Y$  interchangeably.

The restriction $Df(0)\vert_Y\colon Y\rightarrow {\mathbb R}^N$ is an isomorphism and we abbreviate its inverse by   $L=\bigl[ Df(0)\vert_Y\bigr]^{-1}\colon Y\to \R^N$. Instead of studying the  solutions $f(y)=0$ we can as well study  the solutions of $\wt{f}(x)=0$, where $\wt{f}$ is the composition $\wt{f}=L\circ f$.It satisfies $D\wt{f}(0)k=0$ if $k\in K$ and 
$D\wt{f}(0)y=y$ if $y\in Y$ .  In abuse of notation we shall in the following  denote  the composition  $\wt{f}=L\circ f$ by the old letter $f$ again. 
We  may therefore assume that  
$$
f\colon U\cap (K\oplus Y)\rightarrow Y
$$
has the form 
\begin{equation}\label{new_equ_57}
f(k, y)=y-B(k, y),
\end{equation}
where 
\begin{equation}\label{oipu}
B(0,0)=0\quad  \text{and}\quad DB(0,0)=0.
\end{equation}

By $B_K(a)$ , we denote the open ball in $K$ of radius $a$ centered at $0$. Similarly, $B_Y(b)$ is  the open ball in $Y$ of radius $b$ centered  at $0$.

\begin{lemma}\label{LEMMA1}
There exist constants $a>0$ and $b>0$ such  that 
\begin{itemize}
\item[{\em (1)}]\ $\bigl( B_K(a)\oplus B_Y(b)\bigr)\cap C\subset U$.
\item[{\em (2)}]\ $\abs{B(k,y)-B(k,y')}\leq \dfrac{1}{2}\abs{y-y'}$ for 
all $k\in B_K(a)\cap C$ and $y,y'\in B_Y(b)\cap C$. 
\end{itemize}
Moreover, if $f(k,y)=f(k,y')$ for  $(k, y)$ and $(k,y')\in \bigl( B_K(a)\oplus B_Y(b)\bigr)\cap C$, then $y=y'$. 
\end{lemma}

\begin{proof}

It is clear that there exists constants $a, b>0$ such that 
(1) holds. 
To verify (2) we estimate for 
$(k, y)$ and $(k, y')\in \bigl( B_K(a)\oplus B_Y(b)\bigr)\cap C$, 
\begin{equation}\label{est_contraction_B}
\abs{B(k, y)-B(k, y')}\leq\biggl[ \int_0^1\abs{D_2B(k, sy+(1-s)y')}ds\biggr](y-y').
\end{equation}

In view of $DB(0, 0)=0$ we find smaller $a, b$ such that $ \abs{D_2B(k, y'')}\leq 1/2$ for all $(k, y'')\in \bigl( B_K(a)\oplus B_Y(b)\bigr)\cap C$.  So, 
\begin{equation*}
\abs{B(k, y)-B(k, y')}\leq \dfrac{1}{2}\abs{y-y'}
\end{equation*}
as claimed in (2).

If $f(k,y)=f(k,y')$, for two point $(k,y), (k,y')\in \bigl( B_K(a)\oplus B_Y(b)\bigr)\cap C$, 
then, $y-B(k, y)=y'-B(k, y')$ and using (2),  
$$\abs{y-y'}=\abs{B(k,y)-B(k,y')}\leq \dfrac{1}{2}\abs{y-y'}$$
which implies  $y=y'$. This completes the proof of Lemma \ref{LEMMA1}.

\qed \end{proof}

To continue with the proof of the theorem, we recall that $K$ is in  good position to $C$ and $Y$ is a good complement of $K$ in $\R^n$. Therefore,  there exists $\varepsilon>0$ such that for $k\in K$ and $y\in Y$ satisfying 
$\abs{y}\leq \varepsilon\abs{k}$  the statements
$k+y\in C$ and $k\in C$ are equivalent.

\begin{lemma}\label{LEMMA2}
Replacing $a$ by a smaller number, while keeping $b$,  we have
\begin{equation*}
\abs{B(k,y)}\leq \varepsilon \abs{k}\quad \text{ for all $k\in B_K(a)\cap C $ and $\abs{y}\leq \varepsilon \abs{k}$.}
\end{equation*}
\end{lemma}
\begin{proof}
First we replace $a$ by a  perhaps smaller number such that $\varepsilon a<b$. If $k\in B_K(a)\cap C $ and $\abs{y}\leq \varepsilon \abs{k}$, then $k+y\in C$  and $\abs{k+y}\leq (1+\varepsilon )\abs{k}<(1+\varepsilon )a$. Introducing 
$$c(a):= (1+\varepsilon)\cdot \max_{x\in B(1+\varepsilon )a)} \abs{DB(x)},$$
 we  observe that $c(a)\to 0$ as $a\to 0$,  in view of $DB(0, 0)=0$. Then we estimate for $k\in B_K(a)\cap C $  and $y\in Y$ satisfying $\abs{y}\leq \varepsilon \abs{k}$, 
 using that $B(0, 0)=0$, 
\begin{equation}\label{new_equ_70}
\begin{split}
\abs{B(k,y)}& \leq \biggl[ \int_0^1 \abs{DB(sk,sy)}\ ds\biggr] (k+y)\\
&\leq \max_{x\in B(1+\varepsilon )a)} \abs{DB(x)}(1+\varepsilon )\abs{k}=c(a)\abs{k}.
\end{split}
\end{equation}
The assertion follows if we choose $a$ so small that 
$c(a)\leq \varepsilon$.
\qed \end{proof}

Now we take $a>0$ so small that $c(a)\leq \varepsilon/2$ and keep the original $b>0$.
For every $k\in B_K(a)\cap C$, we set 
$$X_k:=\ov{B}_Y(\varepsilon \abs{k})\quad\text{and}\quad X=\bigcup_{k\in B_K(a)\cap C}X_k. $$
Clearly, $X\subset U$. By  Lemma \ref{LEMMA2},  
$B(k,\cdot )\colon X_k\to X_k$, and by 
 Lemma \ref{LEMMA1}, the map $B(k, \cdot )$ is a contraction. Therefore, it  has a unique fixed point $\tau (k)\in X_k$ satisfying 
$$\tau (k)=B(k, \tau (k)).$$
We claim that if   
$y\in B_Y(b)\cap C$ and $B(k, y)=y$ for some $k\in B_K(a)\cap C$, then 
$y=\tau (k)$. 

Indeed, since $B(k, \tau (k))=\tau (k)$, it follows that 
$f(k, y)=f(k,\tau (k))$,  and Lemma \ref{LEMMA1} shows that 
$y=\tau (k)$, as claimed.

\begin{lemma}
The map $\tau\colon B_K(a)\cap C \to B_Y(b)$  is continuous
and satisfies $\abs{\tau(k)}\leq \varepsilon  \abs{k}/2$.
\end{lemma}
\begin{proof}
In view of \eqref{new_equ_70}, the estimate 
$\abs{\tau (k)}\leq \varepsilon \abs{k}/2$ follows from 
$$
\abs{\tau (k)}=\abs{B(k, \tau (k)}\leq c(a)\abs{k}
$$
 and from the choice $c(a)\leq \varepsilon/2$.
 In order to prove continuity of $\tau$, we fix a point $k_0\in  B_K(a)\cap C$ and  use  Lemma \ref{LEMMA1}
 to estimate, 
\begin{equation*}
\begin{split}
\abs{\tau (k)-\tau (k_0)}&=\abs{B(k, \tau (k))-B(k_0, \tau (k_0))}\\
&\leq \abs{B(k, \tau (k))-B(k, \tau (k_0))}+\abs{B(k, \tau (k_0))-B(k_0, \tau (k_0))}\\
&\leq \dfrac{1}{2}\abs{\tau (k)-\tau (k_0)} +\abs{B(k, \tau (k_0))-B(k_0, \tau (k_0))}.
\end{split}
\end{equation*}
This implies, 
\begin{equation}\label{continuity_tau}
\begin{split}
\abs{\tau (k)-\tau (k_0)}&\leq 2\abs{B(k, \tau (k_0))-B(k_0, \tau (k_0))}\\
&\leq 2\biggl[ \int_0^1\abs{D_1B(sk+(1-s)k_0, \tau (k_0))}\ ds\biggr] \abs{k-k_0}\\
&\leq \varepsilon \abs{k-k_0},
\end{split}
\end{equation}
where we have used our choice of $a$ that $c(a)= (1+\varepsilon)\cdot \max_{x\in B(1+\varepsilon )a)} \abs{DB(x)}\leq \varepsilon/2$. This finishes the proof of the continuity of the map $\tau$.
\qed \end{proof}


\begin{lemma}
If $f\colon U\to \R^N$ is of class $C^j$, $j\geq 1$, then the map $\tau\colon B_K(a)\cap C \to B_Y(b)$ is also of class $C^j$.
\end{lemma}

\begin{proof}

We fix the  point $k\in B_K(a)\cap C$ and assume that  
$k+\delta k\in B_K(a)\cap C.$
Then we abbreviate 
\begin{align*}
A(\delta k)&:=\int_0^1\bigl(D_2B(k+\delta k, s\tau (k+\delta k)+(1-s)\tau (k))-D_2B(k, \tau (k))\bigr)\ ds\\
R(\delta k)&:=B(k+\delta k, \tau (k))-
B(k, \tau (k))-D_1B(k, \tau (k))\delta k.
\end{align*}
Using the fixed point property $\tau (k)=B(k,\tau (k))$ and 
$\tau (k+\delta k)=B(k+\delta k,\tau (k+\delta k))$ we compute,
\begin{equation*}
\begin{split}
&\tau (k+\delta k)-\tau (k)-D_1B(k,\tau (k))\delta k\\
&\quad =B(k+\delta k, \tau (k+\delta k))-B(k, \tau (k))-D_1B(k,\tau (k))\delta k\\
&\quad=\bigl[ B(k+\delta k, \tau (k))-B(k, \tau (k))-D_1B(k_0,\tau (k_0))\bigr] \\
&\quad\phantom{=}+
\bigl[ B(k+\delta k, \tau (k+\delta k)) -B(k+\delta k, \tau (k))\bigr] \\
&\quad=R(\delta k)+A(\delta k)\bigl[ \tau (k+\delta k)-\tau (k)\bigr]+D_2B(k, \tau (k))\bigl[ \tau (k+\delta k)-\tau (k)\bigr].
\end{split}
\end{equation*}
Therefore,
\begin{equation*}
\begin{split}
 \bigl[{\mathbbm 1}-D_2B(k, \tau (k))\bigr]&(\tau (k+\delta k)-\tau (k))-D_1B(k, \tau (k))\delta k\\
&=R(\delta k)+A(\delta k)(\tau (k+\delta k)-\tau (k)).
\end{split}
\end{equation*}
Now, $R(\delta k)/\abs{\delta k}\to 0$ and $A(\delta k)\to 0$ as $\abs{\delta k}\to 0$. Moreover, in view of  
\eqref{continuity_tau}, 
$\abs{\tau(k+\delta k)-\tau (k)}/\abs{\delta k}\leq \varepsilon \abs{\delta k}/\abs{\delta k}\leq \varepsilon$. 
Consequently,
$$
\lim_{\abs{\delta k}\to 0}\dfrac{1}{\abs{\delta k}}\abs{\tau (k+\delta k)-\tau (k)- \bigl[{\mathbbm 1}-D_2B(k, \tau (k))\bigr]^{-1}D_1B(k, \tau (k))\delta k}=0.
$$
This shows that the map $\tau$ is differentiable at the point $k$ and its  derivative $D\tau (k)\in {\mathscr L}(K, Y)$ is given by the familiar formula
$$D\tau (k)\delta k= \bigl[{\mathbbm 1}-D_2B(k, \tau (k_0))\bigr]^{-1} D_1B(k,\tau (k))\delta k.$$
The formula shows that $\tau \mapsto D\tau (k)$ is continuous 
since $\tau $ and $DB$ are continuous maps. Consequently,  the map  $\tau$ is of class $C^1$.

We also note that, differentiating $\tau (k)=B(k, \tau (k))$, we obtain 
\begin{equation*}
\begin{split}
D\tau(k)\delta k& = D_2B(k,\tau(k))D\tau(k)\delta k + D_1B(k,\tau(k))\delta k\\
&=DB(k,\tau(k))(\delta k, D\tau(k)\delta k). 
\end{split}
\end{equation*}
If now $f$ is of class $C^2$, we define the $C^1$-map 
$$
f^{(1)}\colon V\oplus K\oplus Y\oplus Y\to Y\oplus Y
$$
by 
\begin{equation*}
\begin{split}
f^{(1)}(k, h, y,\eta)&=\bigl( y-B(k, y), \eta -DB(k, y)(h, \eta)\bigr)\\
&=(y,\eta)-\bigl( B(k, y), DB(k,y)(h, \eta)\bigr)
\end{split}
\end{equation*}
and find by the same reasoning as before, using that  $D_2B(k, y)$ for small $(k,y)$ is contractive, that the associated family of fixed points
$$(k,h)\mapsto  (\tau(k),D\tau(k)h)$$
is of class $C^1$ near $(0,0)$. In particular, the map 
$\tau$ is of class $C^2$ near $k=0$. Proceeding by induction we conclude that $\tau$ is of class $C^j$ 
on a possibly smaller neighborhood $V$ of $0$ in $K$.

\qed \end{proof}
The proof of Theorem \ref{help-you} is finished.

\chapter{Manifolds and Strong Retracts}\label{chap4}
\label{section_tame_manifolds}

The previous chapter showed that  a solution set of a sc-Fredholm  is a sub-M-polyfold whose induced polyfold structure is equivalent to the structure of a finite dimensional smooth manifold with boundary with corners. In this chapter we shall study these objects in more details. 

\section{Characterization}\label{subsect_characterization}

A smooth manifold with boundary with corners is a paracompact Hausdorff   space $M$ which admits an atlas of smooth compatible quadrant charts $(V, \varphi, (U, C,\R^n))$, where $\varphi\colon V\to U$ is a homeomorphism 
from an open subset $V\subset M$ onto a relatively open set $U$ of the partial quadrant $C$ of $\R^n$.
A M-polyfold $X$ is a paracompact Hausdorff space which, in addition, is equipped with a sc-structure. This prompts the following definition.

\begin{definition}\index{D- Equivalent manifold structure}
A M-polyfold $X$ has a {\bf compatible smooth manifold structure with boundary with corners}, if  
it admits an atlas ${\mathcal A}$ consisting of sc-smoothly  compatible  charts 
$(V,\varphi, (U, C, \R^n))$, where $\varphi\colon V\to U$ is a sc-diffeomorphism from the open set $V\subset X$ onto the relatively open set $U\subset C$ of the partial quadrant 
$C$ in $\R^n$ (or in a finite-dimensional vector space $E$). 
\qed
\end{definition}

The transition maps between the  charts are sc-diffeomorphisms between relatively open subsets of partial quadrants in finite-dimensional vector spaces,  
and therefore are classically smooth maps, i.e., of class $C^\infty$. Consequently, the atlas ${\mathcal A}$ defines the structure of a smooth manifold with boundary with corners.
We also note that from the definition it follows immediately that a compatible smooth manifold structure with boundary with corners, if it exists, is unique.

The aim of Section \ref{subsect_characterization}
is the proof of the following characterization.
\begin{theorem}[{\bf Characterization}]
\label{XXX--}\index{T- Characterization of smooth manifolds}
For a tame M-polyfold $X$ the following statements are equivalent.
\begin{itemize}
\item[{\em (1)}]\  $X$
has a compatible smooth manifold structure with boundary with corners.
\item[{\em (2)}]\ $X_0=X_\infty$ and the identity  map ${\mathbbm 1}^1_0\colon X\rightarrow X^1$, $x\mapsto x$ is sc-smooth. Moreover, every point $x\in X$ is contained in a tame sc-smooth polyfold chart $(V, \psi ,(O, C, E))$ satisfying $\psi (x)=0\in O$, and the tangent space $T_0O$ is finite dimensional and in good position to the partial quadrant $C$ in the sc-Banach space $E$.
\end{itemize}\qed
\end{theorem}

\begin{remark}\label{r:r}
If $X_0=X_\infty$, the identity map 
${\mathbbm 1}^1_0\colon X\rightarrow X^1$ has as its inverse the identity map ${\mathbbm 1}^0_1\colon X^1\rightarrow X$. The map 
${\mathbbm 1}^0_1$ is 
always sc-smooth. 
Hence, if  the map ${\mathbbm 1}^1_0\colon X\to X^1$ is sc-smooth, it is a sc-diffeomorphism. 
There is a  subtle point in (2). If $X_0=X_\infty$,  then $X_m=X_0$ as sets for all $m\geq 0$.  However,  it is possible that $X^1\neq X$ as sc-smooth spaces
despite the fact that ${(X^1)}_m = X_{m+1}=X_i$ as sets. In  other words, it  is possible that  as M-polyfolds $X^1\neq X$  
even if the underlying sets are the same.
Indeed, recalling the definition
of sc-differentiability, one realizes that  the sc-smoothness of ${\mathbbm 1}\colon X\rightarrow X$ 
is  completely different from the 
sc-smoothness of ${\mathbbm 1}\colon X\rightarrow X^1$.
This is the reason we use the notation ${\mathbbm 1}_0^1$ 
for the identity map $X\to X^1$.
In contrast, if $X$ is  a finite-dimensional vector space equipped with the constant sc-structure, then $X=X^1$ as sc-spaces.
\qed
\end{remark}

The proof of Theorem \ref{XXX--} requires some preparations and we start with a definition where we denote as usual with $U\subset C\subset E$ a relatively open set $U$ in the partial quadrant $C$ of the sc-Banach space $E$.

\begin{definition}[{\bf Sc$^+$-retraction}] \label{sc_plus_retraction}\index{D- Sc$^+$-retraction}
 A sc-smooth retraction $r\colon U\to U$  is called a 
 {\bf $\ssc^+$-retraction}, if 
$r(U_m)\subset U_{m+1}$ for all $m\geq 0$ and if  $r\colon U\to U^1$  is sc-smooth.
Similarly, if $V\subset X$ is an open subset of a  M-polyfold $X$, we call the sc-smooth retraction $r\colon V\rightarrow V$ a {\bf $\ssc^+$-retraction},  if $r\colon V\rightarrow V^1$ is sc-smooth.
\qed
\end{definition}

\begin{lemma}\label{L-L}\index{L- Sc$^+$-retractions}
Let $(O,C,E)$ be  a sc-smooth retract and suppose that there exists a relatively open subset  $U\subset C$  and a 
$\ssc^+$-retraction
$t\colon U\rightarrow U$ onto  $t(U)=O$. Then every   sc-smooth retraction $s\colon V\rightarrow V$ of a relatively open subset $V\subset C$ satisfying $s(V)=O$
is a $\ssc^+$-retraction. 
\end{lemma}

\begin{proof}
By assumption, the sc-smooth map  $t:U\rightarrow U$  satisfies $t\circ t=t$ and $t(U)=O$. In addition, $t\colon U\to U^1$ is  sc-smooth. 
If now $s\colon V\rightarrow V$ is a sc-smooth retraction onto $s(V)=0$, then $t(s(v))=s(v)$ and hence 
$$
s = t\circ s.
$$
In view of the properties of $t$, the composition 
$t\circ s\colon V\to V^1$ is sc-smooth and hence 
$s\colon V\rightarrow V^1$ is sc-smooth, as claimed. 
 \qed \end{proof}

\begin{definition}[{\bf Sc$^+$-retract}] \index{D- Sc$^+$-retract}
The sc-smooth retract $(O,C,E)$ is a 
{\bf $\ssc^+$-retract}, if there exists a relatively open subset 
$U\subset C$ and a $\ssc^+$-retraction $t\colon U\rightarrow U$ onto  $t(U)=O$.
\qed
\end{definition}

In view of the Lemma \ref{L-L} the choice of $t$ is irrelevant, being a $\ssc^+$-retract is an intrinsic property of the sc-smooth retract $(O, C, E)$.

\begin{lemma}\label{new_lemma_4.8}\index{L- Image of sc$^+$-retraction}
If $r\colon U\to U$ is a $\ssc^+$-retraction onto the sc-smooth retract $(O, C, E)$, then 
\begin{itemize}
\item[{\em (1)}]\ $r(U)=O=O_\infty$, i.e., consists of smooth points.
\item[{\em (2)}]\ At every point $o\in O$, the tangent space $T_oO=Dr(o)E$ is finite-dimensional.
\end{itemize}
\end{lemma}
\begin{proof}
(1)\, 
 If $x\in U$, then, by definition, $r(x)\in U_1$ and using 
$r\circ r=r$, $r(x)=r(r(x))\in U_2$. Continuing this way we conclude that $r(x)\in \bigcap_{m\geq 0} U_m=U_\infty$. Hence $r(U)=O\cap U_\infty=O_\infty$.\par

\noindent (2)\, At the point $o\in O$, the linearization 
$Dr(o)\colon E\to E$ is, in view of (1),  well-defined, and it is a 
$\ssc^+$-operator, since $r$ is a $\ssc^+$-map. Consequently, $Dr(o)\colon E\to E$ is a compact operator between every level. Therefore, the image of the projection $Dr(o)$, namely $Dr(o)E=T_oO$ must be finite-dimensional. This proves the lemma.
 \qed \end{proof}

The $\ssc^+$-retracts are characterized by the following proposition.
\begin{proposition}\label{frank}\index{P- Characterization of sc$^+$-retractions}
A sc-smooth retract $(O,C,E)$ is a $\ssc^+$-retract  if and only if  $O_0=O_\infty$ and ${\mathbbm 1}_0^1\colon O\rightarrow O^1$ is sc-smooth.
\end{proposition}

\begin{proof}

If $(O,C,E)$ is a $\ssc^+$-retract, there exists a $\ssc^+$-retraction $r\colon U\to U$ onto  $r(U)=O$. 
The retract has the induced M-polyfold structure  defined  by the scaling $O_m=r(U_m)$ and $O^1$ inherits this structure, so that 
$(O^1)_m=O_{m+1}$, for all $m\geq 0$.
In view of Lemma \ref{new_lemma_4.8}, $O=O_\infty$. Hence the identity map ${\mathbbm 1}_0^1\colon O\to O^1$ is well-defined. According to Definition \ref{tangent_retract}, the map  
${\mathbbm 1}_0^1$ is sc-smooth, provided  the composition 
${\mathbbm 1}_0^1\circ r \colon U\to E^1$ is sc-smooth, which is the case because  $r$ is a $\ssc^+$-map.

In order to prove the opposite direction we assume that $O=O_\infty$ and ${\mathbbm 1}_0^1\colon O\to O^1$ is sc-smooth. Let $r\colon U\to U$ be a sc-smooth retraction onto $O=r(U)$. Then the composition 
$$
U\xrightarrow{r} O\xrightarrow{{\mathbbm 1}_0^1} O^1\xrightarrow{\textrm{inclusion}} U^1
$$
is sc-smooth. The composition agrees with the map $r\colon U\to U^1$ and the lemma is proved.
 \qed \end{proof}

The proof of Theorem \ref{XXX--} will make use of the following technical result.

\begin{proposition}
\label{FFF}\index{T- Tangent representation}

Let $(O,C,E)$  be a $\ssc^+$-retract and $t\colon U\to U$ a $\ssc^+$-retraction of  the relatively open subset $U\subset C$ onto $O=t(U)$. We assume that $0\in O$ and the the tangent space $T_0O$ (which  by Lemma \ref{new_lemma_4.8} is finite-dimensional) is in good position to the partial quadrant $C$.   We denote by $Y$ a good complement of $T_0O$, so that $E=T_0O\oplus Y,$  and denote by 
$$
p\colon E=T_0O\oplus Y\rightarrow T_0O
$$ 
the associated sc-projection. 
Then there exist  open neighborhoods ${\mathcal U}$ of $0$ in  $O$ and ${\mathcal V}$ of $0$ in $T_0O\cap C$ such  that
$$
p\colon {\mathcal U}\rightarrow {\mathcal V}
$$
is a sc-diffeomorphism.

\end{proposition}
 
\begin{remark}We 
note that $({\mathcal V},T_0O\cap C,T_0O)$ is a local M-polyfold model, because $T_0O\cap C$ is a  partial quadrant in $T_0O$, in view of Proposition \ref{pretzel}. 
The proof of Proposition \ref{FFF} is based on the implicit function theorem for sc-Fredholm sections proved in the previous chapter.
\qed
\end{remark}
\begin{proof}[Proof of Proposition \ref{FFF}]

We recall that $U$ is a  relatively open neighborhood of $0$ in  the  partial quadrant  $C$ of the sc-Banach space $E$.
Moreover,  $t\colon U\rightarrow U$ is  a $\ssc^+$-retraction onto $O=t(U)$ which contains $0$. In view of 
Lemma \ref{new_lemma_4.8}, the tangent space $T_0O=Dt(0)E$ is a smooth subspace of finite dimensions. By assumption, $T_0O$ lies in good position to $C$. Accordingly, there exists a good sc-complement $Y$ of $T_0O$ in $E$ so that 
$$
E=T_0O\oplus Y.
$$
 It has the property that there 
exists  $\varepsilon_0>0$  such that for $a+y\in T_0O\oplus Y$ satisfying $\abs{y}_0\leq \varepsilon \abs{a}_0$ the  statements $a\in C$    
and $a+y\in C$ are equivalent.

We shall use in the following the notation $u=a+y\in T_0O\oplus Y$ and denote by $p$ the  sc-smooth projection 
$$
p\colon E=T_0O\oplus Y\rightarrow T_0O.
 $$
Now we consider the local strong bundle 
$$
\pi\colon U\triangleleft Y\rightarrow U,
$$
and  the sc-smooth section $f\colon U\to U\triangleleft Y$ , $f(u)  = (u,{\bm{f}}(u))$, whose principal part ${\bm{f}}\colon U\to Y$
is defined  by
$$
 {\bm{f}}(u) = ({\mathbbm 1}-p)\bigl( u-t(u)\bigr).
$$
We observe that if $u\in O$, then $u=t(u)$ and hence ${\bm{f}}(u)=0.$ We shall show later on that there are no other solutions $u\in U$ of ${\bm{f}}(u)=0$ locally near $u=0$.

\begin{lemma}\label{new_lemma4.11}
The section $f$ is sc-Fredholm.
\end{lemma}
\begin{proof}
In order to verify that the sc-smooth section $f$ is regularizing, we take $u=a+y\in U_m$ and assume that 
${\bm{f}}(u)=y-({\mathbbm 1}-p)\circ t (a+y)\in Y_{m+1}.$
Since $t$ is $\ssc^+$-smooth, we conclude  that $y\in Y_{m+1}$ and since $a$ is smooth, that  $a+y\in U_{m+1}$, showing that $f$ is indeed regularizing. 

We continue by studying the principal part ${\bm{f}}$ of $f$.
The derivative of ${\bm{f}}$ at the smooth point $u\in U$ has the form
\begin{equation*}
\begin{split}
D{\bm{f}}(u)h &=({\mathbbm 1}-p)(h - Dt(u)h)\\
& = h -(p(h)-({\mathbbm 1}-p)Dt(u)h)
\end{split}
\end{equation*}
for $h\in E$. Since  $Dt(u)$  and $p$ are  $\ssc^+$-operators, the operator $D{\bm{f}}(u)$ is a perturbation of the identity operator by a $\ssc^+$-operator, and therefore a sc-Fredholm operator, in view of Proposition 
\ref{prop1.21}.

Next we have to verify that at a given smooth point $u\in U$, a filled version of ${\bm{f}}$, after modification by a $\ssc^+$-section ${\bm{s}}$, is conjugated to a basic germ. However, in the case at hand, we already 
work on relatively open subsets of a partial quadrant, so that  a filling is not needed. We merely have to find a suitable $\ssc^+$-section to obtain a section  which is conjugated to a basic germ.

In general, the good complement $Y$ in the decomposition 
$E=T_0O\oplus Y$ is not a subset of $C$. 
However, we claim that there exists  a finite-dimensional sc-Banach space $B$ and a sc-Banach space $W$ contained in $C$,  such  that $Y=B\oplus W$,  which leads to the sc-decomposition
$$
E= (T_0O\oplus B)\oplus W.
$$
Indeed, if $E=\R^n\oplus F$ and the partial quadrant $C\subset E$ is of  the form $C =[0, \infty)^n\oplus F$ where $F$ is a sc-Banach space, then, identifying the sc-Banach space $F$ with $\{0\}^n\oplus F$, we let $B$ to be an  algebraic complement of $Y\cap F$ in $Y$ and $W=Y\cap F$, so that $Y=B\oplus  W$. Clearly, $W\subset C$ and  $B$ is finite-dimensional. 
 In the general case, there exists a sc-isomorphism $L\colon E\to \R^n\oplus F$ mapping $C$ onto 
the partial quadrant $C'=[0,\infty)^n\oplus F$.  Then the subspace $L(T_0O)$ is in good position to the partial quadrant $C'$ and $L(Y)$ is a good complement of $L(T_0O)$ in $\R^n\oplus F$. 
By the above argument, we find a finite-dimensional subspace $B'$ and a sc-Banach space $W'$ contained in $C'$, so that $L(T_0O)=B'\oplus W'$.  With $B=L^{-1}B'$ and $W=L^{-1}W'$,  we conclude that $B$ is finite-dimensional, $W$ is contained in $C$ and $Y=B\oplus W$, as claimed.
Since $W\subset C$,
$$
C=\bigl(C\cap (T_0O\oplus B)\bigr)\oplus W.
$$
The finite-dimensional space $T_0O\oplus B$ is in good position to $C$. Therefore, the set $C\cap (T_0O\oplus B)$ is a partial quadrant in $T_0O\oplus B$, in view of Proposition \ref{pretzel}. We shall represent an element $u\in U$ as 
$$u=a+b+w\in T_0O\oplus B\oplus W,$$
where $b+w=y\in Y=B\oplus W.$
Furthermore, we denote by  $P$ the projection 
$$P\colon Y=B\oplus W\to W.$$
In accordance with the notation in the definition of a basic germ, we view $Y$ as $\R^N\oplus W$, where 
$N=\dim(B)$,  and $P$  as the sc-projection $P\colon \R^N\oplus W\to W$. 

If $u\in C$ is a smooth point near $0$ we denote by $C_u$ the associated partial quadrant introduced in Definition \ref{new_def_2.33}. The partial quadrant $C_u$ contains $C$, so that 
$$C_u=(C_u\cap (T_0O\oplus B))\oplus W.$$
Fixing the smooth point $u\in C$, the map $v\mapsto u+v$ for $v$ near $0$ in $C_u$, maps an open neighborhood of $0$ in $C_u$ to an open neighborhood of $u$ in $C$. We now study the principal part ${\bm{f}}\colon U\rightarrow Y$ near the fixed smooth point $u$, make  the change of coordinate $v\mapsto  u+v$, and define the section 
${\bm{g}}(v)$ by 
$${\bm{g}}(v):={\bm{f}}(u+v).$$
In addition, we define a $\ssc^+$-section ${\bm{s}}$ near 
near $0$ in $C_u$,  by 
$${\bm{s}}(v) ={\bm{f}}(u) +( {\mathbbm 1}-p)[  t(u)-t(u+v)]. $$
It has the property that  at $v=0$, 
$$
{\bm{g}}(0)-{\bm{s}}(0)=0.
$$
In the decomposition $v=(a+b)+w$ we interpret $a+b\in C_u\cap (T_0O\oplus B)$ as a finite-dimensional parameter and consider 
the germ 
$$(a+b,w)\mapsto  ({\bm{g}}-{\bm{s}})(a+b+w)\in Y.$$
Then 
\begin{equation*}
\begin{split}
 P({\bm{g}}-{\bm{s}})(a+b+w)=P({\mathbbm 1}-p)(a+b+w) =P(b+w)=w.
\end{split}
\end{equation*}
Identifying an open neighborhood of $0$ in $C_u\cap (T_0O\oplus B)$ with an open neighborhood of $0$ in some partial quadrant $[0,\infty)^k\oplus {\mathbb R}^{n-k}$, we conclude that ${\bm{g}}-{\bm{s}}$ is a basic germ, whose contraction part happens to be  identically zero. We have used the local  strong bundle  isomorphism $(u+v,b+w)\mapsto (v,b+w)$.

Consequently,  $(f,u)$ is a sc-Fredholm germ and the proof of lemma \ref{new_lemma4.11}
is  complete.
 \qed \end{proof}

The sc-Fredholm section $f$ vanishes at the point $u=0$, so that ${\bm{f}}(0)=0$. The  linearization of its principal part is  is given by 
$$
D{\bm{f}}(0)(h) = ({\mathbbm 1}-p)h,\quad h\in E.
$$
Therefore, $D{\bm{f}}(0)$ has the kernel $\ker(D{\bm{f}}(0))=T_0O=Dt(0)E$,  which, by assumption,  is in good position to the partial quadrant $C$. Moreover, its image is $ ({\mathbbm 1}-p)E=Y$, so that $D{\bm{f}}(0)$ is surjective.

We are in position to apply the implicit function theorem, Corollary  \ref{LGS2},  and conclude that the local solution set of ${\bm{f}}$ near $u=0$ is represented by 
$$\{u\in U_0\, \vert \, {\bm{f}}(u)=0\}=\{u=a+\delta (a)\in T_0O\oplus Y\, \vert \, a\in V\}\subset U.$$
Here  $U_0\subset U$ is an open neighborhood of $0$ in $U$ and $V$ is an open neighborhood of $0$ in the partial quadrant $T_0O\cap C$ of $T_0O$. The map $\delta\colon V\to Y$ is a sc-smooth map satisfying $\delta (0)=0$ and $D\delta (0)=0$.

That all sufficiently small solutions of ${\bm{f}}(u)=0$ are of the form $a+\delta (a)$, is, in our special case, easily verified directly.

\begin{lemma}\label{new_lemma4.12}
Let $f(a+y)=0$ and $\abs{a}_0<\varepsilon$, $\abs{y}_0<\varepsilon$. If $\varepsilon$ is sufficiently small, then 
$y=\delta (a)$.
\end{lemma}

\begin{proof}
From $f(a+y)=0$ and $f(a+\delta (a))=0$, we obtain 
$y=({\mathbbm 1}-p)t(a+y)$ and 
$
\delta (a)=
({\mathbbm 1}-p)t(a+\delta (a))
$. 
Consequently,  
$$
y-\delta (a)=({\mathbbm 1}-p)[t(a+y)-t(a+\delta (a))].
$$
As $f$ is regularizing, the solutions are smooth points. Since $t$ is a $\ssc^+$-map, the map 
$u\mapsto ({\mathbbm 1}-p)t(u)$ is of class $C^1$ between 
the $1$-levels, and the derivative vanishes at $u=0$, 
$({\mathbbm 1}-p)Dt(0)=0$. Therefore, there exists 
$\varepsilon_1>0$ such that 
\begin{equation}\label{new_eq1_em4.12}
\norm{({\mathbbm 1}-p)Dt(u)}_{L(E_1,E_1)}\leq 1/2
\end{equation}
if $\abs{u}_1<\varepsilon_1$. 
Consequently, on level $1$,
\begin{equation}\label{new_eq2_em4.12}
y-\delta (a)=\biggl(\int_0^1({\mathbbm 1}-p)Dt(u(\tau))\ d\tau \biggr)\cdot (y-\delta (a)),
\end{equation}
where $u(\tau )=a+\tau y+(1-\tau )\delta (a)$.
Since $t$ is a $\ssc^+$-map satisfying $t(0)=0$, we conclude from 
$y=({\mathbbm 1}-p)t(a+y)$ that $\abs{y}_1$ is small if $a$ and $y$ are small on level $0$. Moreover, the map $\delta$ is sc-smooth and satisfies $\delta (0)=0$ and hence 
$\abs{\delta (a)}_1$ is small, if $a$ is small on level $0$. Summing up, 
$\abs{u(\tau)}_1<\varepsilon_1$ if $\varepsilon$ is sufficiently small,  and we conclude from \eqref{new_eq1_em4.12} and 
\eqref{new_eq2_em4.12} the estimate 
$$
\abs{y-\delta (a)}_1\leq \dfrac{1}{2}\abs{y-\delta (a)}_1,
$$
so that indeed $y=\delta (a)$ is $\varepsilon$ is sufficiently small, as claimed in Lemma \ref{new_lemma4.12}.
 \qed \end{proof}

We finally verify that the solutions $f(a+\delta (a))=0$ belong to $O$ if $a$ is small on level $0$.

\begin{lemma}\label{new_lemma_4.13}
Let $f(a+\delta (a))=0$ and $\abs{a}_0<\varepsilon.$
If $\varepsilon>0$ is sufficiently small, then 
$a+\delta (a)\in O$.
\end{lemma}
\begin{proof}
We have to confirm that $a+\delta (a)=t(a+\delta (a))$. From 
$f(a+\delta (a))=0$ we conclude that $\delta (a)=
t(a+\delta (a))-p\circ t(a+\delta (a))$, so that our aim is to prove that 
$$a=p\circ t(a+\delta (a)).$$
Applying the retraction $t$ to the identity
\begin{equation}\label{new_eq3_lem4.13}
p\circ t(a+\delta (a))+\delta (a)=t(a+\delta (a)),
\end{equation}
and using $t\circ t=t$,  we obtain
$$
t\bigl(p\circ t(a+\delta (a))+\delta (a)\bigr)=t(a+\delta (a)).
$$
Hence, abbreviating 
$$a_1:=p\circ t(a+\delta (a)),$$
we arrive at the identity
$$t(a_1+\delta (a))-t(a+\delta (a))=0.$$
Going to level $1$, abbreviating $a(\tau)=\tau a_1+(1-\tau)a+\delta (a)$ and observing that 
$Dt(0)(a_1-a)=a_1-a$,  we estimate
\begin{equation*}
\begin{split}
0&=\abs{t(a_1+\delta (a))-t(a+\delta (a))}_1\\
&=\left|\int_0^1Dt(a(\tau ))\ d\tau (a_1-a)\right|_1\\
& \geq \abs{a_1-a}-\biggl(\int_0^1\norm{Dt(a(\tau ))-Dt(0)}_{L(E_1, E_1)}\ dt \biggr)\cdot \abs{a_1-a}_1
\end{split}
\end{equation*}
Since $\norm{Dt(u)}_{L(E_1, E_1)}$ is continuous in $u$ on level $1$, there exists $\varepsilon_1>0$ such that 
$$\norm{Dt(v))-Dt(0)}_{L(E_1, E_1)}\leq 1/2$$
if $\abs{v}_1\leq \varepsilon_1$. Arguing as in the previous lemma, $\abs{a(\tau)}_1\leq \varepsilon_1$ if $\abs{a}_0<\varepsilon$ and $\varepsilon$ is sufficiently small. Consequently,
$$
0\geq \abs{a_1-a}_1-\dfrac{1}{2}\abs{a_1-a}_1=\dfrac{1}{2}\abs{a_1-a}_1.
$$
Therefore, $a_1=a$ and hence $p\circ t(a+\delta (a))=a$ if 
$\varepsilon>0$ is sufficiently small, and Lemma \ref{new_lemma_4.13} is proved.
 \qed \end{proof}
In order to complete the proof of Proposition \ref{FFF}, we set 
${\mathcal V}=\{a\in T_oO\cap C\, \vert \, \abs{a}_0<\varepsilon\}.$ 
The map $p$, satisfying 
$$p (a+\delta (a))=a,$$
is a sc-diffeomorphism from the open set 
${\mathcal U}=\{ a+\delta (a)\, \vert \, a\in {\mathcal V} \}$ in $O$, to the open set ${\mathcal V}$ having the sc-smooth map $a\mapsto a+\delta (a)$ as its inverse. This completes the proof of Proposition \ref{FFF}.
 \qed \end{proof}

\begin{proof}[Theorem \ref{XXX--}]
We assume that (1) holds true: the M-polyfold $X$  has a compatible smooth manifold structure with boundary with corners. 
Correspondingly there is an atlas of sc-smoothly compatible partial quadrants charts $(V, \varphi, (U, C, \R^n))$ where $\varphi\colon V\to U$ is a sc-diffeomorphism from the open subset $V\subset X$ onto the relatively open subset $U\subset C$ of the partial quadrant $C$ in $\R^n$. We take in $\R^n$ the unique constant sc-structure. 

Then the identity map ${\mathbbm 1}\colon U\to U$ is a tame sc-smooth retraction of $(U, C, \R^n)$.

 Clearly, $U=U_\infty$ and ${\mathbbm 1}_0^1\colon U\to U^1$ is sc-smooth. Consequently, these charts define a tame M-polyfold structure on $X$ for which  $X_0=X_\infty$ and
 ${\mathbbm 1}_0^1\colon X\to X^1$ is sc-smooth and the statement (2) follows.
 
Conversely, if (2) holds for the tame M-polyfold $X$, then $X_0=X_\infty$ and ${\mathbbm 1}_0^1\colon X\to X^1$ is sc-smooth. If $x\in X$, then we find a tame sc-smooth polyfold chart $(V,\psi, (O, C, E))$ such that the sc-diffeomorphism $\psi\colon V\to O$ satisfies $\psi (x)=0$. 
It follows that $O_0=O_\infty$ and 
${\mathbbm 1}_0^1\colon O\to O^1$ is sc-smooth. By Proposition \ref{frank} the sc-retract $(O, C, E)$ is a $\ssc^+$-retract and hence there is a $\ssc^+$-retraction $t\colon U\to U$ of a relatively open sbset $U$ in the partial quadrant $C$ satisfying $O=t(U)$. The tangent space $T_0O$ is a finite-dimensional smooth subspace of the sc-Banach space $E$ and, by assumption,  in good position to $C$.

Using Proposition \ref{FFF} we find an open neighborhood ${\mathcal O}'$ of $0$ in $O$ and an open neighborhood ${\mathcal V}$ of $0$ in $T_0O\cap C$ and a sc-diffeomorphsim $p\colon {\mathcal O}'\to {\mathcal V}$ onto the finite-dimensional polyfold model $({\mathcal V}, T_0O\cap C, T_0O)$. Taking the open set $V_0=\psi^{-1}({\mathcal O}')\subset V$, the composition 
$$
V_0\xrightarrow{\psi}{\mathcal O'}\xrightarrow{p}{\mathcal V}
$$
is a sc-diffeomorphism onto the finite-dimensional model $({\mathcal V}, T_0O\cap C, T_0O)$. Carrying out this construction around every point $x\in X$, we obtain a smoothly compatible atlas of tame charts and the statement (1) follows. This completes  the proof of Theorem \ref{XXX--}.
 \qed \end{proof}

\section{Smooth Finite Dimensional  Submanifolds}

We start with the definition.

\begin{definition}\label{Def2.40}\index{D- Smooth submanifold}
A {\bf smooth finite-dimensional submanifold} of the M-polyfold $X$  is a subset $A\subset X$ having the property that every point $a\in A$ possesses 
 an open  neighborhood $V\subset X$ and a  $\ssc^+$-retraction $s\colon V\to V$ onto 
 $$
 s(V)=A\cap V.
 $$
 \qed
\end{definition}

In contrast to the sc-smooth sub-M-polyfold of an M-polyfold in Definition 
\ref{def_sc_smooth_sub_M_polyfold}, the retracts in the above definition are  $\ssc^+$-retracts.
We also point out, that a smooth finite dimensional submanifold is not defined as a sub-M-polyfold whose induced structure admits an equivalent structure of a smooth manifold with boundary with corners.
\begin{remark} \index{R- On the notion of  finite-dimensional manifold}
The subset $A\subset X$ in 
Definition \ref{Def2.40}  is called a smooth finite-dimensional manifold for the following reasons. The retraction $s\colon V\to V$, satisfying $s\circ s=s$, is a $\ssc^+$-map. Therefore, the image $s(V)\subset X_\infty$ consists of smooth points, so that every point $a\in s(V)$ possesses a tangent space $T_aA$ defined by $T_aA=Ts(a)(T_aX)$. From the $\ssc^+$-smoothness of the map $s$, it follows that the projection $Ts (a)\colon T_aX\to T_aX$ is a $\ssc^+$-operator  and 
therefore level-wise compact. 
Hence its image, $T_aA$ is finite-dimensional.  Since $A$,  as we shall see,  is in particular a sub-M-polyfold it has an inherited M-polyfold structure
so that its degeneracy index is defined. In general for $a\in A$ it holds that $d_A(a)\leq d_X(a)$. One can show that $A$ near points $a\in A$ with $d_A(a)=0$
has a natural smooth manifold structure. Near points with $d_A(a)=1$ it has a natural structure as a smooth manifold with boundary. If $d_A(a)\geq 2$ the boundary structure
can be quite complex.  However, under the  additional assumption that the tangent space $T_aA$ is in good position to the partial quadrant $C_aX$ in $T_aX$ the subset $A$
has near $a$ in a natural way the structure of a smooth manifold with boundary with corners.  We shall see that $A$ inherits the structure of a M$^+$-polyfold. These differ from
usual classical finite-dimensional manifolds by allowing more general boundary behavior, but if $d_A\leq 1$ they are the same as smooth manifolds with or without boundary.
\qed
\end{remark}

Since an sc$^+$-retraction is an sc-smooth retraction it follows that  a smooth finite-dimensional submanifold $A$ inherits 
from its ambient space $X$ a natural M-polyfold structure. As 
such the degeneracy index $d_A(a)$ is well-defined,  and the boundary $\partial A$ is defined as the subset $\{a\in A\, \vert \, d_A(a)\geq 1\}$.  However, being locally an sc$^+$-retraction implies additional properties.

Recall that a subset $A$ of  a topological space $X$ is called {\bf locally closed}\index{Locally closed} if every point $a\in A$ possesses an open neighborhood $V(a)\subset X$ having the property that a point $b\in V(a)$ belongs to $A$, if $U\cap A\neq \emptyset$ for all open neighborhoods $U$ of $b$. 

\begin{proposition}\label{new_prop4.15}\index{P- Basic properties of submanifolds}
A smooth finite-dimensional submanifold $A\subset X$  of the M-polyfold $X$ has the
following properties.
\begin{itemize}
\item[{\em (1)}]\ $A\subset X_\infty$,  and  $A$  inherits the M-polyfold structure induced from $X$.  In particular,  $A$  possesses a tangent space at every point in $A$, and 
the degeneracy index $d_A$ is defined on $A$.
\item[{\em (2)}]\ $A$ is locally closed in $X$.
\item[{\em (3)}]\ $\partial A=\{a\in A\, \vert \, d_A(x)\geq 1\}\subset \partial X$. 
\item[{\em (4)}]\ ${\mathbbm 1}_0^1\colon A\to A^1$ is sc-smooth.
\end{itemize}
\end{proposition}
\begin{proof}\mbox{}  (1)\, A $\ssc^+$-retraction is, in particular, a sc-retraction and hence $A$ inherits the M-polyfold structure from $X$, in  view of Proposition \ref{sc_structure_sub_M_polyfold}. Since a $\ssc^+$-retraction $s$ has its image in $X_\infty$,  we conclude that $A\subset X_\infty$,  implying that every point in $a\in A$ has 
the  tangent space $T_aA=Ts(s)(T_aX)$. Since $A$ is a M-polyfold,  the degeneracy index $d_A$ is defined on $A$.\\[0.5ex]
(2)\, In order to prove (2) we choose a point   $a\in A$ and an open neighborhood $V\subset X$ of $a$ such that a  suitable 
$\ssc^+$-retraction $s\colon V\to V$  retracts onto $s(V)=A\cap V.$
If  $b\in V$ lies in the closure of $A$, then there exists a sequence $(a_k)\subset A$
satisfying $a_k\rightarrow b$. Since $V$ is open $b\in V$, we conclude that  
$a_k\in A\cap V$ for large $k$.  Hence $a_k=s(a_k)$ and,  using that $b\in V$,  we conclude that $s(b)=\lim_k s(a_k)=b,$ implying that $b\in A\cap V$.\\[0.5ex]
(3)\,  By definition,  $\partial A=\{a\in A\,  \vert \, d_A(x)\geq 1\}$.
We assume that $a\in A$ satisfies $d_A(a)\geq 1$ and show that $d_X(a)\geq 1$. If $d_X(a)=0$,  we find an open neighborhood
$V\subset X$  of $a$ which is sc-diffeomorphic to a retract $(O,E,E)$. This implies that there exists an open neighborhood 
$U\subset A$ of $a$ which is sc-diffeomorphic to a retract $(O',E,E)$, so that $d_A(a)=0$.  \\[0.5ex]
(4)\, The postulated $\ssc^+$-retraction is the identity on the retract $A$, and $A=A_\infty$, so that ${\mathbbm 1}_0^1\colon A\to A^1$ is sc-smooth.

The proof of 
Proposition \ref {new_prop4.15}
is complete.
 \qed \end{proof}

\begin{question} It is an open question,  whether  there 
is a difference between a smooth finite-dimensional submanifold $A$
and a sub-M-polyfold  satisfying $A=A_\infty$ and $\dim T_aA$ being finite and locally constant.  In fact, the question is if under the latter conditions ${\mathbbm 1}_0^1\colon A\rightarrow A^1$ is sc-smooth.  We conjecture that this is not always the case. Rather we would expect,
for example,  if the dimension of the tangent space is equal to $1$, the set $A$ might be something like a branched one-dimensional manifold
as defined in \cite{Mc}. It would be interesting to see if such type of examples  can be constructed. For example,  
consider the subset $T:=\{(x,0)\in {\mathbb R}^2\, \vert\, -1<x<1\}\cup\{(0,y)\in {\mathbb R}^2\, \vert \, y\in [0,1)\}$ with the induced topology.
Is it possible to find a sc-smooth retract $(O,E,E)$ where $E$ is a sc-Banach space so that $T$ is homeomorphic to $O$
and $O=O_\infty$?  From our previous discussion this is impossible if we require $O$ to be a $\ssc^+$-retract. 
\qed
\end{question}

So far not much can be said about the structure of the smooth finite-dimensional submanifold $A$ at the boundary without additional assumptions. In order to formulate such an assumption, we consider a M-polyfold $X$ which is required to be tame. Thus, at every smooth point $a\in X$, the cone $C_aX$ is a partial quadrant in the tangent space $T_aX$, in view of Proposition \ref{tame_equality}, and we can introduce the following definition.

\begin{definition}[{\bf Good position at  $a\in A$}] \label{DEF424}\index{D- Good position at $a\in A$}
A smooth finite-dimensional submanifold $A$ of the tame M-polyfold $X$ is {\bf in good position at the point $a\in A$}, if the finite-dimensional linear subspace $T_aA\subset T_aX$ is in good position to the partial quadrant $C_aX$ in $T_aX$.
\end{definition}

We also need the next definition.

\begin{definition}[{\bf Tame submanifold}] \index{D- Tame submanifold}
A smooth finite-dimensional submanifold $A\subset X$ of the M-polyfold $X$ is called 
{\bf tame}, if, equipped with its induced M-polyfold structure, the M-polyfold $A$ is tame.
\end{definition}

The main result of the section is as follows.

\begin{theorem}\label{thm-basic}\index{T- Characterization of tame submanifolds}
Let $X$ be a  tame M-polyfold and  $A\subset X$ a smooth finite-dimensional submanifold of $X$. If $A$ is at every point $a\in A$ in good position, then the induced M-polyfold structure  on $A$ is equivalent to the structure of a smooth manifold with boundary with corners. In particular, the M-polyfold $A$ is tame.  
\end{theorem}

\begin{proof}

The result will be based on Proposition \ref{XXX--}. We focus on a point $a\in A$. Then $a$ is a smooth point, and we find a sc-smooth tame chart $\varphi\colon (V, a)\mapsto (O, 0)$, where $\varphi$ is a sc-diffeomorphism from the open neighborhood 
$V\subset X$ of $a$ onto the retract $O$ satisfying $\varphi (a)=0$. The retract $(O, C, E)$ is  a tame local model. By assumption, $T_aA$ is in good position to the partial quadrant $C_aX$ in $T_aX$ and there is a good complement $Y'\subset T_aX$, so that $T_aX=T_aA\oplus Y'$. Therefore, the tangent space $T_0O=T\varphi (a)(T_aX)$ has the sc-splitting 
$$T_0O=N\oplus Y,$$
in which $N=T\varphi(a)(T_aA)$ and $Y=T\varphi (a)Y'$. Since $(O, C, E)$ is tame, the tangent space $T_0O$ has a sc-complement $Z$ contained in $C$, in view of Proposition \ref{IAS-x}. Hence 
$$
E=T_0O\oplus Z= N\oplus Y\oplus Z.
$$

\begin{lemma}\label{new_lemma4.20}
The finite-dimensional subspace $N=T\varphi(a)(T_aA)$ is in good position to $C$ and $Y\oplus Z$ is a good complement in $E$.
\end{lemma}

\begin{proof}[{\bf Proof of Lemma \ref{new_lemma4.20}}] 
By assumption, $N$ is in good position to $C_0=T_0O\cap C=T\varphi (a)(C_aX)\subset T_0O$, with the good complement $Y$. Hence there exists $\gamma>0$ such that  for $(n, y)\in N\oplus Y$
satisfying $\abs{y}_0\leq \gamma \abs{n}_0$ we have $n\in C_0$ if and only if $n+y\in C_0$. Take the norm $\abs{(y, z)}_0=\abs{y}_0+\abs{z}_0$ on $Y\oplus Z$ and consider $(n, y, z)$ satisfying $\abs{(y, z)}_0\leq \gamma \abs{n}_0$, and note that $n\in C$ if and only if $n\in C_0$. If $n\in C_0$, we conclude from 
$\abs{(y, z)}_0\leq \gamma \abs{n}_0$  that $n+y\in C_0$ which implies $n+y\in C$. Since $z\in C$, we conclude that $n+y+z\in C$. Conversely, we assume that $n+y+z\in C$ satisfies $\abs{(y, z)}_0\leq \gamma \abs{n}_0$. Then it follows from $z\in Z\subset C$ that $n+y\in C$ and hence $n+y\in C_0=T_0O\cap C$. From $\abs{y}_0\leq \gamma \abs{n}_0$ we deduce that 
$n\in C_0$ and hence $n\in C$. Having verified that $N$ is in good position to $C$ and $Y\oplus Z$ is a good complement, the proof of lemma is complete.
 \qed \end{proof}

Continuing with the proof of Theorem \ref{thm-basic}, we recall that $(O, C, E)$ is  a tame retract. Hence there exists a relatively open  subset $U\subset C$ and a  sc-smooth tame retraction $r\colon U\to U$ onto $O=r(U)$. Since $A\subset X$ is a finite-dimensional smooth submanifold of $X$, we find an open neighborhood $V\subset X$ of $a$ and a $\ssc^+$-retraction $s\colon V\to V$ onto $s(V)=A\cap V$. Taking $V$ and $U$ small, we may assume that $V=\varphi^{-1}(O)$. 
Then, defining  $t\colon U\rightarrow U$ by
$$
t(u)=\varphi\circ s\circ \varphi^{-1}\circ r (u),
$$
we compute, using $r\circ \varphi =\varphi$ and $s\circ s=s$, 
\begin{equation*}
\begin{split}
t\circ t& = \varphi\circ s\circ\varphi^{-1}\circ r\circ  \varphi\circ s\circ \varphi^{-1}\circ r\\
&=\varphi\circ s\circ\varphi^{-1}\circ  \varphi\circ s\circ \varphi^{-1}\circ r\\
&= \varphi\circ s\circ s \circ \varphi^{-1}\circ r\\
&= \varphi\circ s\circ \varphi^{-1}\circ r\\
&=t.
\end{split}
\end{equation*}
We see that $t$ is 
a retraction. Moreover, it is a $\ssc^+$-retraction since $s$ is a 
$\ssc^+$-retraction. It retracts onto $t(U)=\varphi (A\cap V)$ and 
$Dt(0)E=T\varphi (a)Ts(a)(T_aX)=T\varphi (a)T_aA=N$. 
Therefore,  $Q:=t(U)\subset O$ is a $\ssc^+$-retract and $T_0Q=N$ 
is in good position to $C$, by Lemma \ref{new_lemma4.20}.

Now we can apply Proposition  \ref{FFF}  and conclude that 
that there are open neighborhoods $V_1$ of $0$ in $Q$ and $V_0$ of $0$ in $C\cap N$ such that $p=Dt(0)$ is a sc-diffeomorphism 
$$p\colon V_1\to V_0.$$
In view of Proposition \ref{pretzel},  $C\cap N$ is a partial quadrant in $N$, so that $(V_0, C\cap N, N)$ is a local model which is tame since $V_0\subset C\cap N$ is open.

Defining the open neighborhood $V(a)=\varphi^{-1}(V_0)\cap A$ of $a$ in $A$,  the map 
$$p\circ \varphi \colon V(a)\to V_0$$
is a sc-diffeomorphism defining a sc-smooth tame chart on the M-polyfold $A$. The collection of all these charts defines a sc-smoothly  compatible structure of a smooth manifold with boundary with corners. Moreover, the M-polyfold $A$ is tame. 
The proof of Theorem \ref{thm-basic} is complete.
 \end{proof}

\begin{remark}\index{R- Good position}
If $A$ is in good position at $a$, then there exists and open neighborhood $V(a)\subset X$, so that $A$ is in good position for all $b\in V(a)\cap A$. 
In order to see this we observe that Theorem \ref{FFF} defines a chart from the sole knowledge 
that we are in good position at $0$. Since the tangent spaces  of $A$  move only slowly, 
Lemma \ref{good_pos} implies that we are in good position at the nearby points as well.
\end{remark}

\begin{definition}[{\bf General position}] \index{D- Submanifold in general position}
Let $A\subset X$ be a smooth finite-dimensional submanifold of the 
M-polyfold $X$. Then $A$ is in {\bf general position} at the point $a\in A$, if 
 $T_aA$ has in $T_aX$ a sc-complement 
 contained in the reduced tangent space $T^R_aX$ as defined in Definition
 \ref{def_partial_cone_reduced_tangent}.
 \qed
\end{definition}

The following  local result describes an $A\subset X$ in general position.

\begin{theorem}[{\bf Local structure of smooth submanifolds}]\index{T- Local structure of submanifolds}\label{local-str}
Let $A$ be a smooth submanifold of the tame M-polyfold $X$. We assume that $A$ is in general position at the point $x\in A$, assuming that $T_xA$ has in $T_xX$ a sc-complement contained in $T_x^RX$. Then there exists an open neighborhood $U\subset X$  on which the following holds.
\begin{itemize}
\item[{\em (1)}]\ There are precisely $d=d_X(x)$  many local faces ${\mathcal F}_1,\ldots,{\mathcal F}_d$ contained in $U$. 
\item[{\em (2)}]\ $d_X(a)=d_A(a)$ for $a\in A\cap U$.
\item[{\em (3)}]\ If $a\in A\cap U$, then the tangent space $T_aA$ has in $T_aX$ a sc-complement contained in $T^R_aX$.
\item[{\em (4)}]\ If $\sigma\subset \{1, \ldots ,d\}$, then the intersection ${\mathcal F}_\sigma:=\bigcap_{i\in \sigma}{\mathcal F}_i$ of local faces is a tame M-polyfold and $A\cap {\mathcal F}_\sigma$ is a tame smooth submanifold of ${\mathcal F}_\sigma$. At every point $a\in (A\cap U)\cap {\mathcal F}_\sigma$, the tangent space $T_a(A\cap {\mathcal F}_\sigma)$ has in $T_a{\mathcal F}_\sigma$ a sc-complement contained in  $T_a^R{\mathcal F}_\sigma$.
 \end{itemize}
\end{theorem}

\begin{proof}
The result is local and,  going into a tame M-polyfold chart of $X$ around the point $x$,  we may assume that $X=O$ and $x=0\in O$, where $O$ belongs to the tame retract 
$(O, C, E)$ in which $C=[0,\infty)^n\oplus W$, $n=d_O(0)$, is a partial  quadrant in the sc-Banach space $E={\mathbb R}^n\oplus W$.
Moreover, $A\subset O$ is a smooth submanifold of $O$. We recall that the faces of $C$ are the subsets of $E$, defined by 
$$
F_i=\{(a_1,\ldots ,a_n,w)\in E\, \vert \,\text{$a_i=0$ and $a_j\geq 0$ for $j\neq  i$ and $w\in W$}\}.$$
For a subset $\sigma\subset \{1,\ldots ,n\}$ we introduce $F_\sigma=\bigcap_{i\in \sigma}F_i$. Then the M-polyfold $O$ has the faces ${\mathcal F}_i=O\cap F_i$ and we abbreviate the intersection by ${\mathcal F}_\sigma:=O\cap F_\sigma$.

By Definition \ref{reduced_cone_tangent}, $T_0^RO=T_0O\cap (\{0\}\oplus W)$, in our local coordinates.
By assumption of the theorem, the tangent space $T_0A\subset T_0O$ has in $T_0O$ a sc-complement $\wt{A}\subset T_0^RO$, so that $T_0O=T_0A\oplus \wt{A}$. In view of the definition of tame, the tangent space $T_0O$ has in $E$ a sc-complement $Z\subset \{0\}\oplus W$. Hence $E=T_0O\oplus Z=T_oA\oplus \wt{A}\oplus Z$ and it follows that $T_0A$ has in $E$ a sc-complement $Y$ contained in $W=\{0\}\oplus W$,
$$
E=T_0A\oplus Y,\quad Y\subset W.
$$
We denote by $N\subset T_aA$ a sc-complement of $T_0A\cap W$ in $T_0A$,
$$
T_0A=N\oplus (T_0A\cap W).
$$
From $E\cap W=W=(T_0A\cap W)\oplus (Y\cap W)=(T_0A\cap W)\oplus Y$, we obtain 
$$E=T_0A\oplus Y=N\oplus (T_0A\cap W)\oplus Y=N\oplus W.$$
Using $Y\subset W\subset C$, we deduce 
$$
C=(C\cap N)\oplus W$$
and
$$
C=(T_0A\cap C)\oplus Y.
$$
Therefore, the projection $\R^n\oplus W\to \R^n$ induces an isomorphism from $N$ to $\R^n$ and from $C\cap N$ onto $[0,\infty )^n$. Moreover, the tangent space $T_0A$ is in good position to $C$, and $Y$ is a good complement.
Hence $T_0A\cap C$ is a partial quadrant of $T_0A$ by Proposition \ref{pretzel}.

As proved in Proposition \ref{FFF} the smooth submanifold $A\subset O$ is represented local as the graph 
$$A=\{v+\delta (v)\, \vert \, v\in V\}$$
of a sc-smooth map $\delta\colon V\to Y$ defined on an open neighborhood $V$ of $0$ in the partial quadrant $T_0A\cap C$ of $T_0A$ satisfying $\delta (0)=0$ and $D\delta (0)=0$. 
The projection 
$$p\colon A\to (T_0A\cap C)\cap V,$$
defined by $p(v+\delta (v))=v$, is a sc-diffeomorphism from the M-polyfold $A$ and 
we conclude, by Proposition \ref{newprop2.24}, that 
$$d_A(a)=d_{T_0A\cap C}(p(a))$$
for all $a\in A\cap U$.
The points in $T_0A\cap C$  are represented by $(n, m)\in (N\cap C)\oplus (T_0A\cap W)$ so that the points $a\in A$ are represented by 
$$a=((n, m), \delta (n, m)).$$

Introducing 
$$
(T_0A)_\sigma:=T_0A\cap T_0{\mathcal F}_\sigma=\{(n, m)\in T_0A\, \vert \, 
\text{$n_i=0$ for $i \in \sigma$}\},
$$
the diffeomorphism $p\colon A\to  V\subset T_0A\cap C$, maps the intersection of faces $A\cap {\mathcal F}_\sigma$ to $(T_0A)_\sigma \cap C$. In particular, if $p(a)\in (T_0A)_\sigma\cap C$, then $a\in A\cap {\mathcal F}_\sigma$. This implies that 
$$d_{T_0A\cap C}(p(a))=d_C(a).$$
In view of Proposition \ref{tame_equality}, $d_C(a)=d_O(a)$, and we have verified for $a\in A\cap U$ that $d_A(a)=d_O(a)$ as claimed in the statement (2) of the theorem.

In order to prove the statement (3) we choose $a\in A\cap U$. Then $a=v+\delta (v)$ with $v\in T_0A\cap C$. The tangent space $T_aA$ is the image of the linear map $\alpha\colon h\mapsto h+D\delta (v)h$, $h\in T_0A$. If $a=0$, then $T_0A+W=E$ and we conclude, by means of the particular form of the map $\alpha$, that also $T_aA+W=E$ for all $a\in A\cap U$. Intersecting  with $T_aO$ and using $T_aA\subset T_aO$ leads to 
\begin{equation*}
\begin{split}
T_aO&=T_aA+(W\cap T_aO)\\
&=T_aA+T^R_aO,
\end{split}
\end{equation*}
implying that $T_aA$  has in $T_aO$ a sc-complement which is contained in $T_a^RO$. This proves the statement (3) of the theorem.

As for the last statement we observe that, in view of the arguments in the proof of Proposition \ref{FACE_XXXX}, the sets ${\mathcal F}_\sigma$ are tame M-polyfolds. Moreover, $A\cap {\mathcal F}_\sigma$ are smooth submanifolds of ${\mathcal F}_\sigma$. The previous arguments, but now applied to $A\cap {\mathcal F}_\sigma$ conclude the proof of Theorem \ref{local-str}.
 \qed \end{proof}

We conclude this section by defining M$^+$-polyfolds. If $A\subset X$ is a finite-dimen\-sional submanifold then it will naturally inherit such a structure.
We leave the details to the reader.
If $X$ is a M-polyfold and $A\subset X$ a finite-dimensional submanifold it inherits the structure of a M-polyfold. However, note that the natural local models
are $(O,C,E)$ where $O$ is an sc$^+$-retract. This prompts the following definition.
\begin{definition}\index{D- M$^+$-polyfold}
Let $Z$ be a metrizable topological space. 
\begin{itemize}
\item[(1)]\ A M$^+$-polyfold structure for $Z$ is given by an equivalence class of sc$^+$-atlases.
A sc$^+$-atlas consists of sc-smoothly compatible charts where the local models are sc$^+$-retracts.
\item[(2)]\ A M$^+$-polyfold structure is called {\bf tame} provided there exists an sc$^+$-atlas so that 
the underlying local models are tame.
\end{itemize}
\end{definition}
As a consequence of Theorem \ref{XXX--} and the previous discussions  we have the following result.
\begin{proposition}
Let $Z$ be a M$^+$-polyfold. Then the following holds. 
\begin{itemize}
\item[{\em(1)}]\ If $d_Z\equiv 0$ the topological space $Z$ has a unique  compatible smooth manifold structure (without boundary)  which is sc-smoothly compatible with the M$^+$-polyfold structure.
\item[{\em(2)}]\ If $d_Z\leq 1$ the topological space has a  unique  smooth structure as manifold with boundary, which sc-smoothly compatible with the M$^+$-polyfold structure.
\item[{\em(3)}]\ If the M$^+$-polyfold structure is tame there is a unique smooth structure as manifold with boundary with corners on $Z$ which is sc-smoothly compatible with the tame M$^+$-polyfold structure.
\end{itemize}
\end{proposition}
\begin{proof}
This is just a reformulation of earlier results and details are left to the reader.
 \qed \end{proof}

\section{Families and an Application of Sard's Theorem}
If $X$ is a tame M-polyfold we introduce the family
$$Z:=\R^m\oplus X=\{(\lambda , x)\, \vert \, \lambda \in \R^m, x\in X\}.$$
The family $\R^m\oplus X$ has a natural tame M-polyfold structure defined by the product of the corresponding charts and we denote by 
$$P\colon \R^m\oplus X\to \R^m,\quad P(\lambda ,x)=\lambda$$
the sc-smooth projection.
Using Theorem \ref{local-str}
and Sard's theorem we are going to establish the following result.

\begin{theorem}\label{SARD}\index{T- Fibered families and Sard}
Let $A\subset Z$ be a smooth submanifold of the tame M-polyfold $Z$. We assume that the induced 
M-polyfold structure is tame and the closure of $A$ is compact. Denoting by $p:=P\vert A\colon A\to \R^m$ the restriction of $P$ to $A$ we assume, in addition, that there exists $\varepsilon>0$ such that the following holds.
\begin{itemize}
\item[{\em(1)}]\ For every $\lambda\in \R^m$ satisfying $\abs{\lambda}\leq \varepsilon$, the set $p^{-1}(\lambda )\subset A$ is compact and non-empty.
\item[{\em(2)}]\ For $z\in A$ there exists a sc-complement $T_zA$ in $T_zZ$ which is contained in $T_z^RZ$.
\end{itemize}
Then there exists a set of full measure $\Sigma\subset \{\lambda\in {\mathbb R}^m\, \vert \, \abs{\lambda}<\varepsilon\}$ of regular values of $p$, having full measure in $B^m(\varepsilon)$ such that for $\lambda\in\Sigma$,  the set 
$$A_\lambda:=\{x\in X \, \vert \,  (\lambda,x)\in A\}$$ is a smooth compact  manifold with boundary with corners, having the additional property that the tangent space $T_xA_{\lambda}$ at $x\in A_\lambda$ has in $T_xX$ a sc-complement contained in $T^R_xX$. 
\end{theorem}

\begin{remark}\index{R- Comment on Proposition \ref{tame_equality}}
By Proposition  \ref{tame_equality}, 
the codimension of $T^R_xX$ is equal to $d_X(x)$. Therefore, if $\lambda\in \Sigma$, then the point $x\in A_\lambda$ can only belong to a corner if $d_X(x)\leq \dim (A_\lambda)$. For example, zero-dimensional manifolds $A_\lambda$ have to lie in $X\setminus \partial X$, one-dimensional $A_\lambda$ can only hit the $d_X=1$ part of the boundary,  etc.
\end{remark}

\begin{proof}

Fixing a point $z\in A$ there exists, in view of Theorem \ref{local-str} an open neighborhood $U(z)\subset Z$ of $z$ in $Z$ such that there are precisely $d=d_X(z)$ local faces ${\mathcal F}_1,\ldots ,{\mathcal F}_d$ in $Z\cap U$. Moreover, for all subsets $\sigma\subset \{1, 2,\ldots ,d\}$, the intersection of $A$ with ${\mathcal F}_\sigma=\bigcap_{i\in \sigma}{\mathcal F}_i$ is a tame smooth submanifold in ${\mathcal F}_\sigma$ of dimension $\dim A-\#\sigma$. Moreover, the tangent space $T_z(A\cap {\mathcal F}_\sigma)$ has a sc-complement in $T_z{\mathcal F}_\sigma$ which is contained in $T_z^R{\mathcal F}_\sigma$.

We cover the set $\{z\in A\, \vert \, \abs{p(z)}\leq \varepsilon\}$ with the finitely many such neighborhoods $U(z_1),\ldots , U(z_k)$ and first
 study the  geometry of the problem in one of these neighborhoods $U(z_i)$ which, for simplicity of notation, we denote by $U(z)$.

If $\sigma\subset \{1,\ldots ,d\}$
and ${\mathcal F}_\sigma$ the  associated 
intersection of local faces, we denote by 
${\mathcal F}^{\circ}_\sigma$ the interior of ${\mathcal F}_\sigma$, i.e., the set ${\mathcal F}_\sigma$ with its boundary removed. The union of all 
${\mathcal F}^{\circ}_\sigma$ is equal to $U(z)$. Here we use the convention that for the empty subset of 
$\{1,\ldots ,d\}$, the empty intersection is equal to $U(z)$ from which the boundary is removed. We note that ${\mathcal F}_{\{1,\ldots ,d\}}^\circ$ is a M-polyfold without boundary. By Theorem \ref{local-str}
the intersection $A\cap {\mathcal F}^\circ_{\sigma}$ is a smooth manifold without boundary having dimension $\dim (A)-\#\sigma$.

With a subset $\sigma$ we associate the smooth projection 
$$p_\sigma\colon A\cap {\mathcal F}_\sigma^\circ\to \R^m,\quad p_\sigma (\lambda , x)=\lambda.$$
Using Sard's theorem we find a subset $\Sigma_\sigma\subset B^m(\varepsilon)$ of regular values of $p_{\sigma}$ having full measure. The intersection 
$$\Sigma=\bigcap_{\sigma \subset \{1,\ldots ,d\}}\Sigma_\sigma$$
has full measure in $B^m(\varepsilon)$ and consists of regular values for all the maps in (1).

Now we fix $\lambda \in \Sigma$. Then the preimage $p^{-1}(\lambda )\subset A$ has the form $\{\lambda \}\times A_\lambda$ and, by construction,
$$p_\sigma^{-1} (\lambda)=(\{\lambda\}\times A_\lambda)\cap {\mathcal F}_\sigma^\circ ,$$ 
which is a smooth submanifold of $A\cap {\mathcal F}_\sigma^\circ$. Moreover, the tangent space $T_{z'}(A\cap {\mathcal F}_\sigma^\circ )$ at the point $z'=(\lambda, x)\in A\cap {\mathcal F}_\sigma^\circ$ is equal to 
$$T_{z'}(A\cap {\mathcal F}_\sigma^\circ )=T_{z'}\bigl((\{\lambda \}\times A_\lambda)\cap {\mathcal F}_\sigma^\circ\bigr)\oplus \xi_{z'}$$
where the linearized projection 
$$Tp_\sigma (z')\colon T_{z'}(A\cap {\mathcal F}_\sigma^\circ )\to T_{p(z')}\R^m$$
maps $\xi_{z'}$  isomorphically onto $\R^m$.

From our discussion we conclude, in view of the assumption (2) of the theorem,  for $z'=(\lambda, x)\in U(z)$ that the tangent space $T_{z'}(\{\lambda \}\times A_\lambda)$ has a sc-complement in $T^R_{z'}X$. Consequently, $A_\lambda\cap U(z)$ is a smooth submanifold of $X$ in general position.

The argument above applies  to every $z_1,\ldots ,z_k$. Hence for every $z_i$ there exists a subset $\Sigma_{z_i}\subset B^m (\varepsilon)$ of full measure consisting of regular values. The subset $\wh{\Sigma}$, 
$$\wh{\Sigma}=\bigcap_{i=1}^k\Sigma_{z_i}\subset B^m(\varepsilon ),$$
has full measure and it follows, for every $\lambda \in \wh{\Sigma}$, that $A_\lambda\subset X$ is a smooth compact submanifold of the M-polyfold $X$ which, moreover, is in general position.
\qed \end{proof}

\section{Sc-Differential Forms}\label{subs_sc_differential}
Assume $X$ is a M-polyfold. Then $TX\rightarrow X^1$ has the structure of an sc-smooth bundle and 
it is an easy exercise that the associated Whitney sum of $k$-copies $TX\oplus..\oplus TX\rightarrow X^1$ is an sc-smooth bundle.
Following  \cite{HWZ7} we start with the definition.
\begin{definition}
Let $X$ be a M-polyfold and $TX\rightarrow X^1$ its tangent bundle. A sc-differential $k$-form  $\omega$ is a sc-smooth map
$$\omega:\bigoplus_k TX\rightarrow \R$$
 which is linear in each argument and skew-symmetric. 
The vector space of sc-differential k-forms on $X$ is denoted by $\Omega^k(X)$. \index{$\Omega^k(X)$}
\end{definition}

By means of the inclusion maps $X^i\rightarrow X$ the 
sc-differential $k$-form $\omega$ on $X$ pulls back to a sc-differential $k$-form on $X^i$, defining this way the directed system
$$
\Omega^k(X)\to \ldots \to  \Omega^k(X^i)\rightarrow\Omega^{k}(X^{i+1}\to \ldots 
$$
We denote the direct limit of the system by $\Omega_\infty^\ast (X)$\index{$\Omega_\infty^\ast (X)$}
and introduce the set $\Omega_\infty^\ast (X)$ of sc-differential form on $X_\infty$ by 
$$
\Omega_\infty^\ast (X)=\bigoplus_k \Omega_\infty^k (X).
$$

Next we define the exterior differential. For this we use the Lie bracket of vector fields which has to be generalized to our context.
As shown in \cite{HWZ7} Proposition 4.4,  the following holds.
\begin{proposition}\index{P- Lie bracket}
Let $X$ be a M-polyfold and given two sc-smooth vector fields $A$ and $B$ on $X$ one can define the Lie-bracket by the usual formula which defines
a sc-smooth vector field $[A,B]$ on $X^1$, that is,  
$[A, B]$ is a section of the tangent bundle $T(X^1)\to X^2.$\index{$[A,B]$}
\end{proposition}

In order to define the exterior derivative 
$$d:\Omega^k(X^{i+1})\to \Omega^{k+1}(X^{i}),$$
we take a sc-differential $k$-form and $(k+1)$ sc-smooth vector fields $A_0,A_1,\ldots, A_k$ on $X$ and define $(k+1)$-form $d\omega$ on $X$ by the following familiar formula
\begin{equation*}
\begin{split}
d\omega(A_0,A_1,\ldots, A_k)&=\sum_{i=0}^k(-1)^iD(\omega(A_0,\ldots,\what{A}_i,\ldots, A_k)\cdot A_i \\
&\phantom{=}+ \sum_{i<j} (-1)^{i+j}\omega([A_i,A_j],A_0,\ldots,\what{A}_i,\ldots,\what{A}_j,\ldots,A_k).
\end{split}
\end{equation*}

The exterior derivative $d$ commutes with the inclusion map $X^{i+1}\to X^i$  occurring in the directed system, and consequently induces a map
$$d:\Omega_\infty^\ast (X)\rightarrow\Omega_\infty^\ast (X)$$ 
having the property  $d^2=0$. The pair $(\Omega_\infty^\ast (X), d)$\index{$(\Omega_\infty^\ast (X), d)$} is a graded differential algebra which we call the de Rham complex of the M-polyfold $X$.

\begin{definition}\label{def_de_Rham_cohomology}
The sc-de Rham cohomology of the M-polyfold $X$ is defined as $H^\ast_{sc}(X):=\text{ker}(d)/\text{im}(d)$.\index{$H^\ast_{sc}(X)$}
\end{definition}

There is also a relative version.
If $X$ is a tame M-polyfold the inclusion map $\partial X\rightarrow X$ restricted to local faces is sc-smooth. Local faces are naturally tame M-polyfolds and the same is true for 
the intersection of local faces. Therefore it makes sense to talk about differential forms on $\partial X$ and $\partial X^i$. We define  the differential algebra $\Omega^\ast_\infty(X,\partial X)$ by
$$
\Omega_\infty^\ast(X,\partial X):=\Omega_\infty^\ast(X)\oplus \Omega_\infty^{\ast-1}(\partial X)\index{$\Omega_\infty^\ast(X,\partial X)$}
$$
with differential
$$
d(\omega,\tau)=(d\omega,j^\ast\omega-d\tau)
$$
where $j:\partial X\rightarrow X$ is the inclusion. One easily verifies that $d\circ d=0$ and we denote the associated cohomology by $H^{\ast}_{dR}(X,\partial X)$.

Clearly, a sc-differential $k$-form $\omega\in \Omega^\ast_{\infty}(X)$ induces a classical smooth  differential form 
on a smooth finite-dimensional submanifold $N$ of the M-polyfold $X$.  The following version of Stokes'  theorem holds true.

\begin{theorem}\index{T- Stokes}
Let $X$ be a M-polyfold and let $N$ be an oriented smooth n-dimensional compact tame submanifold of $X$ whose boundary 
$\partial N$, a  union of smooth faces $F$, is equipped with the induced orientation. If $\omega$ is  a sc-differential $(n-1)$-form  on $X$, then 
$$
\int_M d\omega  = \sum_{F} \int_F \omega.
$$
\qed
\end{theorem}
Here, the submanifold $N$ is not assumed to  be face structured.

The sc-smooth map $f:X\to Y$ between two M-polyfolds induces 
the map  $f^\ast:\Omega^\ast_\infty(Y)\rightarrow\Omega^\ast_\infty(X)$ in the usual way. There is also a version of the Poincar\'e Lemma, formulated and proved in \cite{HWZ7}.
The general the theory of sc-differential forms on M-polyfolds can be worked out as for the classical smooth manifolds. We leave it to the reader to carry  out the details.

\chapter{Fredholm Package for M-Polyfolds}\label{sec_package}

Chapter \ref{sec_package}
is devoted to compactness properties 
of sc-Fredholm sections, to their perturbation theory, and to the transversality theory.

\section{Auxiliary Norms}

Recalling Section \ref{section2.5_sb}, we consider the strong bundle 
$$
P\colon Y\rightarrow X
$$
over the 
M-polyfold $X$.  The subset $Y_{0,1}$\index{$Y_{0,1}$ } of biregularity $(0, 1)$\index{Biregularity $(0, 1)$} is a topological space and
the map  $P\colon Y_{0,1}\rightarrow X$ is continuous. 
The fibers $Y_x:=P^{-1}(x)$ have the structure of Banach spaces.

We  first introduce the notion
of an auxiliary norm. This concept allows us to quantify the size of admissible perturbations in the transversality and perturbation theory.
We point out  that our definition is more general than the earlier one  given in \cite{HWZ3}.

\begin{definition}
An {\bf auxiliary norm} \index{D- Auxiliary norm} is a continuous map $N\colon Y_{0,1}\rightarrow{\mathbb R}$, which has the following properties.
\begin{itemize}
\item[(1)]\ The  restriction of $N$ to a fiber is a complete norm.  
\item[(2)]\ If $(w_k)$ is a sequence in $Y_{0,1}$ such that $P(w_k)\rightarrow x$ in $X$
and $N(w_k)\rightarrow 0$, then $w_k\rightarrow 0_x$ in $Y_{0,1}$.
\end{itemize}
\qed
\end{definition}

Any two auxiliary norms are locally compatible according to the later Proposition \ref{peter},  which is an immediate corollary of  the following local comparison result.
\begin{lemma}\label{erde}\index{L- Comparison of  auxiliary norms}
Let $P\colon Y\rightarrow X$ be a strong bundle over the M-polyfold $X$,  and  let   $N\colon Y_{0,1}\rightarrow{\mathbb R}$ be  a continuous map which fiber-wise is a complete norm.
Then the following statements are equivalent.
\begin{itemize}
\item[{\em (1)}]\  $N$ is an auxiliary norm.
\item[{\em (2)}]\ For every $x\in X$  there exists,  for  a suitable open neighborhood $V$ of $x$,  a strong bundle isomorphism 
$\Phi\colon Y\vert V\rightarrow K$ whose underlying sc-diffeomorphism $\varphi\colon V\rightarrow O$ maps $x$ to a point $o\in O$,  and constants $0<c <C$ such that for all $w\in P^{-1}(V)$, 
$$
c\cdot N(w)\leq \abs{h}_1\leq C\cdot N(w),
$$
where $\Phi(w)=(p,h)$. Here $K\rightarrow O$ is a local strong bundle model.
\end{itemize}
\end{lemma}
\begin{proof}

Assume that (1) holds.   We fix  $x\in X$ and choose for a suitable open neighborhood $V=V(x)\subset X$ a  strong bundle isomorphism $\Phi\colon Y\vert V\rightarrow K$, 
where $K\rightarrow O$ is a local strong bundle and  $K\subset U\triangleleft F\subset C\triangleleft F\subset E\triangleleft F$ is the retract  $K=R(V\triangleleft F)$.  Define $q\in O$ by  $(q,0)=\Phi(0_x)$.
An open neighborhood of $0_x\in  Y_{0,1}$ consists of all $w\in Y_{0,1}$ for which the set  of all  $(p,h)=\Phi(w)$ belongs to some open neighborhood of $(q,0)$ in $K_{0,1}$.
In view of the continuity of the function $N\colon Y_{0,1}\to \R$ there exists $\varepsilon>0$ such that $N(y)<1$ for all $y = \Phi^{-1}(p, h)$ satisfying $\abs{p-q}_0+\abs{h}_1<\varepsilon.$
In particular,  $\abs{p-q}_0 <\varepsilon/2$ and $\abs{h}_1=\varepsilon/2$ imply  $N(w)\leq 1$.
Using that $N(\lambda w) =\abs{\lambda} N(w)$ and similarly for the norm $\abs{\cdot }_1$,  we infer  for $y$ close enough to $x$,  
$P(w)= y$, and $\Phi (w)=(p, h)$, that 
$$
N(w)\leq \frac{2}{\varepsilon}\cdot \abs{h}_1.
$$

On the other hand,  assume there is no constant $c>0$ such that 
$$
c\cdot \abs{h}_1\leq N(w)
$$
for $(p,h)=\Phi(w)$ and $p$ close to $q$.  Then we find  sequences $y_k\rightarrow x$ and $(w_k)$ satisfying  $P(w_k)=y_k$ and $\abs{h_k}_1=1$
such that $N(w_k)\rightarrow 0$. Since $N$ is an auxiliary norm,  we conclude that $w_k\rightarrow 0_x$ in $Y_{0,1}$ 
which implies the convergence $\Phi(w_k)=(q_k,h_k)\rightarrow (q,0)$, contradicting  $\abs{h_k}_1=1$.  At this point we have proved that,  given an auxiliary norm $N$,  there exist for every $x\in X$ constants $0<c<C<\infty$ depending on $x$
such  that 
$$
c\cdot N(w)\leq \abs{h}_1\leq C\cdot N(w)
$$
for all $w\in Y_{0,1}$  for which  $P(w)$ is close to $x$ and $(p,h)=\Phi(w)$,  for a suitable strong bundle isomorphism $\Phi$ to a local strong bundle model. Hence (1) implies (2).

The other direction of the proof is obvious: since $N$ is continuous and fiber-wise a complete norm
we see that  the property (1) in  the definition of an auxiliary norm holds.  The estimate  in the statement of the lemma  implies that also property (2) holds.

\end{proof}

As a corollary of the lemma we immediately obtain the following proposition.

\begin{proposition}\index{P- Local equivalence of auxiliary norms}\label{peter}
Let $P\colon Y\rightarrow X$ be a strong bundle. Then there exists an auxiliary norm $N$.
Given two auxiliary norms $N_0$ and $N_1$, then there exists a continuous function $f\colon X\rightarrow (0,\infty)$ such that for all $h\in P^{-1}(x)\subset Y_{0,1}$,  
$$
f(x)\cdot N_0(h)\leq N_1(h)\leq \frac{1}{f(x)}\cdot N_0(h).
$$
\end{proposition}
\begin{proof}
The existence follows from a (continuous) partition of unity argument using the paracompactness of $X$, pulling back by strong bundle maps the standard norm $\abs{\cdot }_1$
to the fibers of the strong bundles. The local compatibility implies the existence of $f$.
\end{proof}

For later considerations we also introduce the notion of a reflexive auxiliary norm, which leads to particularly useful compactness considerations for sc-Fredholm sections.
Consider $P:Y\rightarrow X$, which is a strong bundle over the M-polyfold $X$. 
\begin{definition}\index{D- Reflexive $(0,1)$-fibers} \index{D- Reflexive $1$-fibers}
We say that $P:Y\rightarrow X$ has {\bf reflexive $(0,1)$-fibers} (or more sloppily: {\bf reflexive $1$-fibers}) provided there exists a strong bundle atlas so that for every strong bundle chart
the local model $K\rightarrow O$, has the property that $K\subset C\triangleleft F$, where $F_1$ is a reflexive Banach space.
\qed
\end{definition}

In a strong bundle $P:Y\rightarrow X$ with reflexive $(0,1)$-fibers we can introduce a particular kind of convergence called mixed convergence for a sequence $(y_k)\subset Y_{0,1}$.
It is a  mixture of strong convergence in the base on level $0$ and a weak convergence in the fiber on level $1$.  

\begin{lemma}
Let $P:Y\rightarrow X$ be a strong bundle over a M-polyfold with reflexive $(0,1)$-fibers. Let  $(y_k)\subset Y_{0,1}$ be  a sequence for which $x_k:=P(y_k)$ converges to some $x\in X$ on level $0$.
Given two  a strong bundle charts around $x$, say $\Psi :Y|U(x)\rightarrow K$ and $\Psi': Y|U'(x)\rightarrow K'$, where $K\subset C\triangleleft F$ with $F_1$ being a reflexive Banach space,
and similarly $K'\subset C'\triangleleft F'$, write $(a_k,h_k)=\Psi(y_k)$ and $(a_k',h_k')=\Psi'(y_k)$.  Then $(h_k)\subset F_1$ converges weakly if and only if $(h_k')\subset F_1'$ converges weakly.
\end{lemma}

\begin{proof}
If $\Phi:K\rightarrow K'$
is a strong bundle isomorphism (arising as a transition map) and $\Phi(a,h)=(\phi(a),A(a,h))$, we consider a sequence $(a_k,h_k)\in K$ satisfying $a_k\rightarrow a$ in $E_0$ and 
$h_k\rightharpoonup h$ in $F_1$. The latter implies that $h_k\rightarrow h$ in $F_0$.
Then $\Phi(a_k,h_k)=:(a'_k,h'_k)$ satisfies $a'_k\rightarrow b$ in $E_0'$ and $h'_k\rightarrow  h'$ in $F_0'$.
For large $k$,  the operator norms of continuous linear operators $A(a_k,\cdot )\colon F_1\rightarrow F_1'$ are uniformly bounded.  Therefore,  the sequence 
$(h'_k)=(A(a_k,h_k))$ is bounded in $F_1'$.
From  $h_k\rightarrow h$ in $E_0$ we conclude the convergence  $(a'_k,h'_k)=\Phi(a_k,h_k)\rightarrow \Phi(a,h)=:(a',h')$ in $E_0'\oplus F_0'$. The boundedness of $(h'_k)$ in $F_1'$ and the convergence of $(a'_k,h'_k)$ on level $0$ 
implies the weak convergence $h'_k\rightharpoonup h'$ in $F_1'$. 
\end{proof}
In view of this lemma we can define mixed convergence for sequences in $Y_{0,1}$, where $P:Y\rightarrow X$ has reflexive $(0,1)$-fibers as follows.
 \begin{definition}[{\bf Mixed convergence}]\index{D- Mixed convergence}
 Let $P:Y\rightarrow X$ be a strong bundle over the M-polyfold $X$ having reflexive $(0,1)$-fibers. 
 A sequence $(y_k)\subset Y$ of bi-regularity $(0,1)$, i.e. $y_k\in Y_{0,1}$, is said to be {\bf mixed convergent} to
 an element $y\in Y_{0,1}$ provided $P(y_k)\rightarrow P(y)=:x$ in $X_0$, and there exists a strong bundle chart
 $\Psi:Y|U(x)\rightarrow K\subset E\triangleleft F$ ($F_1$ being reflexive),  
 such  that for  $\Psi(y_k)=(a_k,h_k)$ and $\Psi(y)=(a,h)$  the sequence $h_k$ converges weakly to $h$ in $F_1$,  denoted by $h_k\rightharpoonup h$ in $F_1$.
 We shall write
 $$
 y_k\stackrel{m}{\longrightarrow} y
 $$
if  $y_k$ is mixed convergent to $y$.
\qed
\end{definition}
The next definition introduces the important concept of a reflexive auxiliary norm.
\begin{definition}\index{D- Reflexive auxiliary norm}
Let $P:Y\rightarrow X$ be a a strong bundle over a M-polyfold  with reflexive $1$-fibers. A {\bf reflexive auxiliary norm}
$N:Y_{0,1}\rightarrow [0,\infty)$ is a continuous map possessing  the following properties.
\begin{itemize}
\item[(1)]\ For every $x\in X$ the function $N$,  restricted to the fiber of $Y_{0,1}$ over $x$,  is a complete norm.
\item[(2)]\ If $(y_k)$ is a sequence in $Y_{0,1}$  satisfying  $P(y_k)\rightarrow x$ in $X$ and $N(y_k)\rightarrow 0$,
then $y_k\rightarrow 0_x$ in $Y_{0,1}$.
\item[(3)]\ If $(y_k)\subset Y_{0,1}$ is mixed convergent to $y\in W_{0,1}$, then 
$$
N(y)\leq \text{liminf}_{k\rightarrow 0} N(y_k).
$$
\end{itemize}
\qed
\end{definition}

 A version of Proposition \ref{peter} holds for reflexive auxiliary norms and we leave the proof to the reader.
 In the context of strong bundles over ep-groupoids a similar result is given in Theorem \ref{EXTTT}.
 \begin{proposition}
 Let $P:Y\rightarrow X$ be a strong bundle over a M-polyfold with reflexive $(0,1)$-fibers.
 Then there exists a reflexive auxiliary norm $N$.  Given an auxiliary norm $N_0$ and a reflexive auxiliary norm
 $N_1$ there exists a continuous map $f:X\rightarrow (0,\infty)$ satisfying  for all $y\in Y_{0,1}$ with $x=P(y)$
 $$
 f(x)\cdot N_0(y) \leq N_1(y)\leq \frac{1}{f(x)}\cdot N_0(y).
 $$
 \end{proposition}

\section{Compactness Results}

There are several different kinds of compactness requirements on a sc-smooth section which are useful in practice. It will turn out that they are all equivalent for sc-Fredholm sections. We note that compactness is a notion on the $0$-level of $X$.

\begin{definition}
Let $f$ be a sc-smooth section of a strong bundle $P\colon Y\rightarrow X$. 
\begin{itemize}
\item[(1)]\ We say that 
$f$ has a {\bf compact solution set} \index{D- Compact solution set} if  $f^{-1}(0)=\{x\in X\ \vert \, f(x)=0\}$ is compact in $X$ (on level $0$).
\item[(2)]\ The section  $f$ is called {\bf weakly proper} \index{D- Weakly proper}
if  it has a compact solution set and if for every auxiliary norm $N$ there exists an open neighborhood
$U$ of $f^{-1}(0)$ such  that for every $\ssc^+$-section $s$ having support in $U$  and satisfying $N(s(y))\leq 1$ for all $y$,  the solution set 
$$
\{x\in X\,  \vert \, f(x)=s(x)\}
$$
is compact in $X$. 
\item[(3)]\ The section  $f$ is called {\bf  proper} \index{D- Proper} if $f$ has a compact solution set and if for  every auxiliary norm $N$ there exists an open neighborhood $U$ of
$f^{-1}(0)$ such that the closure in $X_0$ of the set $\{x\in U\,  \vert \,  N(f(x))\leq 1\}$ is compact.
\end{itemize}
\qed
\end{definition}

In point (3) we adopt the convention  that if $f(x)$ does not have bi-regularity $(0,1)$, then $N(f(x))=\infty$.
Obviously  proper implies weakly proper, which in turn implies compactness.  
$$
\text{\bf proper}\, \Longrightarrow\, \text{\bf weakly proper}\ \  \Longrightarrow\, \text{\bf compact solution set}.
$$

In general, for a sc-smooth section $f$,  these notions are not equivalent.  
The basic result that all previous compactness notions coincide for a sc-Fredholm section is given by the following theorem.

\begin{theorem}\label{x-cc}\index{T- Equivalence of compactness notions}
Assume that $P\colon Y\rightarrow X$ is a strong M-polyfold bundle over the M-polyfold $X$  and $f$ a sc-Fredholm section.
If $f$ has a compact solution set, then $f$ is  proper. In particular,  for a sc-Fredholm section the properties of being  proper,  or being weakly proper, or having a compact solution set are equivalent.
\qed
\end{theorem}

\begin{proof} 

We denote by  $N$ the auxiliary norm on a strong M-polyfold bundle  $P$.  Fixing a solution $x\in X$ of $f(x)=0$, we have  to find an open neighborhood $U\subset X$ of $x$, such  that the closure 
$\cl_{X}(\{y\in U\, \vert \,  N(f(y))\leq 1\})$ is compact. Since $f$ is regularizing,  the point $x$ is smooth. 
There is no loss of generality in assuming  that we work
with a filled version $g$ of $f$ for which we have a $\ssc^+$-section $s$ such  that $g-s$ is conjugated to a basic germ.
 Hence, without loss of generality,  we may  assume that  
 we work in local coordinates and $f=h+t$  where $h\colon {\mathcal O}(C,0)\to \R^N\oplus W$ is a basic germ and $t\colon {\mathcal O}(C,0)\to \R^N\oplus W$ a $\ssc^+$-germ satisfying $t(0)=0$. Here $C$ is the partial quadrant $C=[0,\infty)^k\oplus {\mathbb R}^{n-k}\oplus W$ in the sc-Banach space $E=\R^n\oplus W$. Then it suffices to find in local coordinates an open neighborhood $U$ of $0$ in the partial quadrant $C$ such that 
the closure (on level $0$) of the set 
\begin{equation}\label{eq_est_0_5}
\{(a, w)\in U\; \vert \; \abs{(h+t)(a, w)}_1\leq c\}
\end{equation}
is compact. 

The section  $t$ is a $\ssc^+$-section satisfying  $t(0)=0$. Therefore, we find  a constant $\tau'>0$ such that 

\begin{equation}\label{eq_est_1_5}
\text{$\abs{t(a, w)}_1\leq c$\quad for $(a, w)\in C$ satisfying   $\abs{a}_0<\tau'$, $\abs{w}_0<\tau'.$ }
\end{equation}
We denote by  $P\colon \R^N\oplus W\to W$ the sc-projection.   By assumption, $h$ is a basic germ and hence $P\circ h$ is of the form 
$$P\circ h(a, w)=w-B(a, w)\quad \text{for $(a, w)\in C$ near $0$.}$$
Here  $B$ is a sc-contraction germ. Moreover, $({\mathbbm 1}-P)h$ takes values in $\R^N$,  so that its range consists of smooth point. 
In view of the sc-contraction property of $B$ and since any two norms on $\R^N$ are equivalent, replacing $\tau'>0$ by a smaller number,  we may assume that the estimates 

\begin{equation}\label{eq_est_2_5}
\abs{B(a, w)-B(a, w')}_0\leq \dfrac{1}{4}\abs{w-w'}_0
\end{equation}
and
\begin{equation}\label{eq_est_2_5b}
\abs{({\mathbbm 1}-P)h(a, w)}_1\leq c
\end{equation}
are satisfied 
for all $(a, w), (a, w')\in C$ such that $\abs{a}_0\leq \tau'$, $\abs{w}_0\leq \tau'$,  and $\abs{w'}_0\leq \tau'$.
Since $B(0,0)=0$, we find $0<\tau\leq \tau'$ such that 
\begin{equation}\label{eq_est_3_5}
\abs{B(a, 0)}_0\leq \tau'/8\quad \text{for all $a\in \R^n$ such that  $\abs{a}_0\leq \tau$.}
\end{equation}
We introduce the closed set 
$$\Sigma=\{(a, z)\in \R^n\oplus W\, \vert \, \abs{a}_0\leq \tau,\ \abs{z}_0\leq \tau'/2\}$$
and denote by $\ov{B}(\tau')$ the closed ball in $Y_0$ centered at $0$ and having radius $\tau'$.
Then we  define the map $F\colon \Sigma \times \ov{B}(\tau')\to Y_0$ by 
$$F(a, z, w)=B(a, w)+z.$$
If $(a, z, w)\in \Sigma \times \ov{B}(\tau')$, then
\begin{equation*}
\begin{split}
\abs{F(a, z, w)}_0&\leq \abs{B(a, w)-B(a, 0)}_0+\abs{B(a, 0)}_0+\abs{z}_0\\
&\leq \dfrac{1}{4}\abs{w}_0+\abs{B(a, 0)}_0+\abs{z}_0\leq \tau'/4+\tau'/8+\tau'/2=3\tau'/4<\tau',
\end{split}
\end{equation*}
and, if  $(a, z)\in \Sigma$ and $w, w'\in \ov{B}(\tau')$, then 
$$
\abs{F(a, z, w)-F(a, z, w')}_0=\abs{B(a, w)-B(a, w')}_0\leq \dfrac{1}{4}\abs{w-w'}_0.
$$
We see that $F$ is a parametrized contraction of $\ov{B}(\tau')$, uniform in $(a, z)\in \Sigma$.
 In view of the parametrized Banach fixed point theorem, there exists a unique continuous map $\delta \colon \Sigma \to \ov{B}(\tau')$ satisfying 
 $F(a, z, \delta (a, z))=\delta (a, z)$ for every $(a, z)\in \Sigma$. 
Thus 
 $$\delta (a, z)=B(a, \delta (a, z))+z,\quad \text{for all $(a, z)\in \Sigma$}.$$
In particular, if $(a,z, w)\in \Sigma\times \ov{B}(\tau')$ and $z=w-B(a, w)$, then $w=\delta (a, z).$

Now we define the open neighborhood $U$ of $0$ in $C$ by 
 \begin{equation}\label{eq_est_4_5}
 U=\{(a, w)\, \vert \, \abs{a}_0<\tau,\, \abs{w}_0<\tau'/4\}.
 \end{equation}
Assume that $(a, w)\in U$ and  that $z'=h(a, w)$ belongs to $\R^n\oplus Y_1$ and satisfies $\abs{z'}_1\leq c$.  Then  
\begin{equation*}
\begin{split}
z'=h(a, w)+t(a, w)&=P\circ h(a, w)+\bigl(({\mathbbm 1}-P)\circ h+t\bigr) (a, w)\\
&=w-B(a, w)+\bigl(({\mathbbm 1}-P)\circ h+t\bigr)(a, w)
\end{split}
\end{equation*}
and 
\begin{equation}\label{eq_est_5}
w-B(a, w)=z\quad \text{where $z=z'-\bigl(({\mathbbm 1}-P)\circ h+t\bigr) (a, w)$}.
\end{equation}
In view of the estimates \eqref{eq_est_1_5} and \eqref{eq_est_2_5b}, 
\begin{equation}\label{eq_est_5_5}
\abs{z}_1\leq 3c. 
\end{equation}
The  norm of $z$ on level $0$,  can be estimated as 
\begin{equation}\label{eq_est_6_5}
\begin{split}
\abs{z}_0&=\abs{w-B(a, w)}_0\leq \abs{w}_0+\abs{B(a, w)}_0\\
&\leq \abs{w}_0+\abs{B(a, w)-B(a, 0)}_0+\abs{B(a, 0)}_0\\
&\leq \abs{w}_0+\dfrac{1}{4}\abs{w}_0+\abs{B(a, 0)}_0\\
&\leq \tau'/4+\tau'/16+\tau'/8=7\tau'/16<\tau'/2. 
\end{split}
\end{equation}
We conclude that  $(a, z,w )\in \Sigma\times \ov{B}(\tau')$ and $w=\delta (a, z).$ 

At this point we can verify the claim that the closure on level $0$ of the set defined in 
\eqref{eq_est_0_5} is compact. With a sequence $(a_n, w_n)\in \{(a, w)\in U\, \vert \, \abs{(h+t)(a, w)}_1<c\}$, we  consider the   corresponding sequence $z_n=w_n-B(a_n, w_n)$ defined by \eqref {eq_est_5}. Then the estimates  \eqref{eq_est_5_5} and \eqref{eq_est_6_5} give 
$$ \abs{z_n}_1\leq 3c \quad  \text{and}\quad  \abs{z_n}_0\leq \tau'/2, $$
which implies that $w_n=\delta (a_n, z_n)$. Since the embedding $E_1\to E_0$ is compact and $a_n\in \R^n$, we conclude, after taking a subsequence,  that $(a_n, z_n)\to (a, z)$ in $E_0=\R^n\oplus W_0$. From the continuity of the map  $\delta$, we deduce the convergence  
$w_n=\delta (a_n, z_n)\to w:=\delta (a, z)$ in $W_0$. Therefore, 
the sequence $(a_n, w_n)$ converges to $(a, w)$ in $E_0$. Since, by assumption, the solution set of $f$ is compact, the proof of the properness of the Fredholm section is complete.

\end{proof}

There are several other useful considerations. The first is the following  consequence of  the local Theorem \ref{save}.
\begin{theorem}[{\bf Local Compactness}]\label{save-1}\index{T- Local compactness}
Assume that $P\colon Y\rightarrow X$ is strong bundle over the M-polyfold $X$ and $f$ a sc-Fredholm section. Then for a given solution  $x\in X$ of  $f(x)=0$ there exists a nested sequence of open neighborhoods $U(i)$ of $x$ on level zero,
say
$$
x\in U(i+1)\subset U(i)\subset X_0,\quad i\geq 0,
$$
such  that, for all $i\geq 0$,  $\cl_{X_0}(\{y\in U(i)\ |\ f(y)=0\})$ is a compact subset of $X_i$. 
\qed
\end{theorem}

The next result shows that compactness, a notion on the $0$-level, also  implies compactness on higher levels.
\begin{theorem}\index{T- Compactness properties}
Assume that $P\colon Y\rightarrow X$ is a strong bundle over the  M-polyfold $X$ and $f$ is a sc-Fredholm section with compact solution set
 $$
S=\{x\in X\, \vert \, f(x)=0\}
$$
 in $X_0$.
Then $S$  is a compact subset of $X_\infty$.
\end{theorem}
\begin{proof}
By assumption $S$ is compact in $X_0$. Since $f$ is regularizing,  $S\subset X_\infty$.
As was previously shown $X_\infty$ is a metric space. Hence we can argue with sequences.
Take a sequence $(x_k)\subset S$. We have to show that it has a convergent subsequence in $X_\infty$.
After perhaps taking a subsequence we may assume that $x_k\rightarrow x\in S$ in $X_0$. 
From  Theorem \ref{save-1} we conclude that $x_k\rightarrow x$ on every level $i$. 
This implies the convergence  $x_k\rightarrow x$ in $X_\infty$.
\end{proof}

We recall that a sc-smooth section $f$ of the strong bundle $Y\rightarrow X$ defines,  by raising the index,  a sc-smooth section
$f^i$ of $Y^i\rightarrow X^i$. In view of the  previous theorem we conclude that a sc-Fredholm section with compact solution set
produces a sc-Fredholm section $f^i$ with compact solution set. Note that it is a priori clear that $f^i$ is sc-Fredholm. Hence we obtain the following corollary.
\begin{corollary}
Let $f$ be a sc-Fredholm section of the strong bundle $P\colon Y\rightarrow X$ over the M-polyfold $X$ and suppose that the solution set $f^{-1}(0)$ is compact.
Then $f^i$ is a sc-Fredholm section of $P^i\colon Y^i\to X^i$ with compact solution set.
\qed
\end{corollary}

Next we consider another compactness property which is very useful in applications and involves reflexive auxiliary norms.
Given an auxiliary norm $N:Y_{0,1}\rightarrow [0,\infty)$ we define $N:Y\rightarrow [0,+\infty]$ by setting $N(y)=+\infty $ if $y\in Y\setminus Y_{0,1}$.
\begin{definition}\index{D- Reflexive local compactness property}\label{DEFX526}
Assume that $P:Y\rightarrow X$ is a strong bundle over a M-polyfold with reflexive $(0,1)$-fiber and let $f$ be an sc-smooth section of $P$,
which in addition is regularizing.
We say that $f$ has the {\bf reflexive local compactness property}  provided for every reflexive auxiliary norm $N$ 
and every point $x\in X$ there exists an open neighborhood $U(x)$ in $X$, i.e. on level $0$, so that 
$\cl_X\left(  \{q\in U(x)\ |\ N(f(x))\leq 1\}\right) $ is a compact subset of $X$. 
\qed
\end{definition}

\begin{proposition}
Assume that $P:Y\rightarrow X$ is a strong bundle over a M-polyfold with reflexive $(0,1)$-fiber and let $f$ be an sc-smooth section of $P$,
which in addition is regularizing and has the reflexive local compactness property.
\begin{itemize}
\item[{\em(1)}]\ If $x\in X_0\setminus X_1$ there exists for every reflexive auxiliary norm $N$ an open neighborhood $V(x)$ in $X_0$ 
such that $N(f(q))>1$ for all $q\in V(x)$.
\item[{\em(2)}]\ Given a reflexive auxiliary norm $N$ and a point $x\in X$ with $N(f(x))>1$ there exists an open neighborhood $V(x)$ such that for every
$q\in V(x)$ it holds that $N(f(q))>1$.
\end{itemize}
\end{proposition}
\begin{proof}
\noindent (1) is a special case of (2) since $N(f(x))=+\infty$.\par

\noindent (2) Assume that $N(f(x))\in (1,\infty]$. 
Arguing indirectly we find $(x_k)$ with $N(f(x_k))\leq 1$ and $x_k\rightarrow x$. 
Define $h_k=f(x_k)$ and note that $(h_k)\subset Y_{0,1}$. Since $N(h_k)\leq 1$ and $P(h_k)\rightarrow x$
we may assume after perhaps taking a subsequence that $h_k\xrightarrow{m} h$ for some $h\in Y_{0,1}$ with $P(h)=x$.
Since $N$ is a reflexive auxiliary norm we must have
$$
N(h)\leq \text{liminf}_{k\rightarrow\infty} N(h_k)\leq 1.
$$
Moreover $h_k\rightarrow h$ in $Y_{0,0}$ implying that $f(x)=h$. Hence $N(f(x))\leq 1$ contradicting $N(f(x))>1$.
\end{proof}

Under the assumptions of the proposition assume that $U$ is an open subset of $X$ so that the closure of $\{x\in U\ |\ N(f(x))\leq 1\}$
is compact. If $x_k\in U$ with $N(f(x_k))\leq 1$ and $x_k\rightarrow x$ on level $0$ it follows that $h_k=f(x_k)$ is mixed convergent 
to $h=f(x)$. Hence $x\in \cl_X(U)$ and $N(f(x))\leq 1$. Hence $\cl_X(\{x\in X\ |\ N(f(x))\leq 1 \})\subset \{x\in\cl_X(U)\ |\ N(f(x))\leq 1\}$.

A very useful result in applications is given by the following theorem, where we refer to \cite{FH2} for an application to SFT.

\begin{theorem}\label{THM528}\index{T- Extension of controlling $(U_\partial, N)$}
Assume that $P:Y\rightarrow X$ is a strong bundle over a tame M-polyfold with reflexive $(0,1)$-fiber.  Suppose that
$f$ is an sc-Fredholm section with compact solution set having the reflexive local compactness property.  Denote by  $N$ a reflexive auxiliary norm and 
with  $S_\partial:=\{x\in \partial X\ |\ f(x)=0\}$ let  $U_\partial \subset \partial X$ be an open neighborhood of $S_\partial$ in $\partial X$ 
so that the closure of the set $\{x\in U_\partial\ |\ N(f(x))\leq 1\}$ is compact.  Then there exists an open neighborhood 
$U$ of $S=\{x\in X\ |\ f(x)=0\}$ having the following properties.
\begin{itemize}
\item[{\em(1)}]\ $U\cap \partial X= U_\partial$.
\item[{\em(2)}]\ The closure of $\{x\in U\ |\ N(f(x))\leq 1\}$ is compact.
\end{itemize}
\end{theorem}
\begin{remark}\index{R- Controlling compactness}
The result says that under suitable assumptions if $(U_\partial, N)$ controls the compactness on the boundary, one 
can construct $(U,N)$ controlling compactness with $U\cap \partial X=U_\partial$.
\end{remark}
\begin{proof}
Pick a point $x_0\in U_\partial$.  Assume first that $N(f(x_0))>1$. Then we find an open neighborhood $U(x_0)$ in $X$ such that
$N(f(x))>1$ for all $x\in U(x_0)$. By perhaps taking a smaller such neighborhood we may assume that $\partial X\cap U(x_0)\subset U_\partial$.
Next assume that $N(f(x_0))\leq 1$. We find an open neighborhood $U(x_0)$ such that $\{x\in U(x_0)\ |\ N(f(x_0))\leq 1\}$ has compact closure.
At this point we have constructed for every $x_0\in U_\partial$ an open neighborhood $U(x_0)$ with one of the two  specified properties.
\begin{itemize}
\item $N(f(x))>1$ for $x\in U(x_0)$ and $U(x_0)\cap \partial X \subset U_\partial$.
\item  $\cl_X(\{x\in U(x_0)\ |\ N(f(x))\leq 1\})$ is compact and $U(x_0)\cap \partial X \subset U_\partial$.
\end{itemize}

The  closure of $\{x\in U_\partial \ |\ N(f(x))\leq 1\}$ is compact and consists of points $q$ with $N(f(q))\leq 1$. For each such point
we have a constructed an open neighborhood $U(q)$ in $X$ such that the closure of $\{p\in U(q)\ |\ N(f(p))\leq 1\}$ is compact.
We find finitely many $\bar{x}_1,...,\bar{x}_\ell\in \partial X$ 
satisfying 
$$
\cl_{X}(\{x\in U_\partial \ |\ N(f(x))\leq 1\})\subset U(\bar{x}_1)\cup..\cup U(\bar{x}_\ell).
$$
Hence 
$$
\cl_{X}(\{x\in U_\partial \ |\ N(f(x))\leq 1\})\subset \{x\in U(\bar{x}_1)\cup..\cup U(\bar{x}_\ell)\ |\ N(f(x))\leq 1\}.
$$
We also have that $(U(\bar{x}_1)\cup...\cup U(\bar{x}_\ell))\cap \partial X\subset U_\partial$.  If $x\in U_\partial \setminus (U(\bar{x}_1)\cup..\cup U(\bar{x}_\ell))$
we must have $N(f(x))>1$. Consider the union 
$$
\wt{U}:= U(\bar{x}_1)\cup..\cup U(\bar{x}_\ell) \cup \bigcup_{x\in U_\partial \setminus (U(\bar{x}_1)\cup..\cup U(\bar{x}_\ell))} U(x).
$$
This is an open subset of $X$ satisfying $\wt{U}\cap \partial X=U_\partial$.   Moreover
the closure of the set $\{x\in \wt{U}\ |\  N(f(x))\leq 1\}$ is compact. Indeed, take a sequence $(x_k)\subset \wt{U}$ satisfying $N(f(x_k))\leq 1$.
Then we must have  $(x_k)\subset U(\bar{x}_1)\cup..\cup U(\bar{x}_\ell)$ and since the closure of each of the finitely many sets $\{x\in U(\bar{x}_i)\ |\ N(f(x))\leq 1\}$
is compact, we see that $(x_k)$ has a convergent subsequence. The set $\{x\in X\setminus \wt{U}\ |\ f(x)=0\}$ is compact 
and we find finitely many $\wt{U}(\wt{x}_i)\subset X\setminus \partial X$, say $i=1,...,e$, covering this compact set so that in addition the closure of each set $\{x\in \wt{U}(\wt{x}_i)\ |\ N(f(x))\leq 1\}$ is compact. 
Finally define 
$$
U = \wt{U}\cup \wt{U}(\wt{x}_1)\cup...\cup \wt{U}(\wt{x}_e).
$$
Then $U\cap \partial X=U_\partial$ and $\cl_X(\{x\in U\ |\ N(f(x))\leq 1\})$ is compact. The proof is complete.
\end{proof}
\begin{remark}\index{R Compactness control}
The theorem is crucial in many of the inductive perturbation arguments occurring in Floer-type 
problems, for example SFT.
\qed
\end{remark}
\section{Perturbation Theory and Transversality}
For the perturbation theory we need to assume that $X$ is tame. This is needed to control the boundary behavior.
Hence let  $P\colon Y\rightarrow X$ be a strong bundle over a tame  M-polyfold.
 We shall study  for a given sc-Fredholm section 
 $f$ the perturbed section $f+s$, where $s\in\Gamma^+(P)$.

We need a supply of $\ssc^+$-sections,  which (in a complete abstract situation) we can only guarantee if we have enough sc-smooth bump functions
(this is a weaker requirement than having sc-smooth partitions of unity. See Appendix \ref{POU} for a more detailed discussion.)
This is,  for example,  the case if $X$ is built on sc-Hilbert spaces.  

\begin{definition}
The  M-polyfold $X$ {\bf admits sc-smooth bump functions} if for every $x\in X$ and every  open neighborhood $U(x)$
there exists a sc-smooth map $f\colon X\rightarrow {\mathbb R}$ satisfying $f\neq 0$ and $\supp(f)\subset U(x)$.
\qed
\end{definition}

In Section  \ref{POU} (Proposition \ref{prop-x5.36}) the following useful statement is proved.

\begin{proposition}
If the M-polyfold $X$ admits sc-smooth bump functions then for every $x\in X$ and every open neighborhood $U(x)$ there exists
a sc-smooth function $f\colon X\rightarrow [0,1]$ satisfying  $f(x)=1$ and $\supp(f)\subset U(x)$. One can even achieve that $f(y)=1$ for all $y$ close to $x$.
\qed
\end{proposition}

We start with an existence result.

\begin{lemma}[{\bf Existence of $\ssc^+$-sections}]
\index{L- Sc$^+$-sections I}
We assume that $P\colon Y\rightarrow X$ is a strong bundle over the tame M-polyfold $X$ which admits sc-smooth bump functions.
Let $N$ be an auxiliary norm for $P$.
Then for every smooth point $x\in X$,  every smooth point $e$ in the fiber $Y_x=P^{-1}(x)$,  and every $\varepsilon>0$ there exists,  for  a given open neighborhood
$U$ of $x$,  a $\ssc^+$-section $s\in \Gamma^+(P)$ satisfying  $s(x)=e$, $\supp (s)\subset U$,  and $N(s(x))<N(e)+\varepsilon$.
\end{lemma}

The proof is an adaption of the proof in \cite{HWZ3}.
\begin{proof}
Since $X$ is metrizable,  it is a normal space. We choose  an open neighborhood $Q$ of $x$ so that $W\vert Q$ is strong bundle isomorphic
to the tame local bundle $p\colon K\rightarrow O$ mapping $x$ into $0\in O$. We find an open neighborhood $Q'$  of $x$  whose closure in $X$
is contained in $Q$.  If we construct $s\in\Gamma^+(W\vert Q)$ with support in $Q'$ we can extend it by $0$ to $X$. 
It suffices to construct a section $s$ suitably in local coordinates.
Hence we work in $K\rightarrow O$ 
and choose a smooth point $e\in p^{-1}(0)$. The open set $Q'$ corresponds to an open neighborhood $O'$ of $0$ contained in $O$.
We write $e=(0,\wt{e})$.
Let $N\colon K\rightarrow {\mathbb R}$
be the auxiliary norm. Using the local strong bundle retraction $R$, we define the $\ssc^+$-section $t\colon O\to K$ by $t(y)=R(y,\wt{e})$, so that $t(0)=(0, \wt{e})=e$.  
If $\varepsilon>0$ there exists $\delta>0$ such that $N(t(y))<N(t(0))+\varepsilon=N(e)+\varepsilon$ for $y\in O$ and $\abs{y}_0<\delta$. Using Proposition \ref{prop-x5.36}, we find a sc-smooth function $\beta\colon O\to [0,1]$ satisfying $\beta (0)=1$ and having support in $O'\cap \{y\in O\, \vert \, \abs{y}_0<\varepsilon\}$. The $\ssc^+$-sections $s(y)=\beta (y)t(y)$ of $K\to O$ has the required properties.
\end{proof}

Next we discuss a  perturbation and transversality result in the case that the M-polyfold does not have a boundary.
Our usual notation will be $P\colon Y\rightarrow X$ for the strong bundle. In case we have an auxiliary norm $N$
and an open subset $U$ of $X$ we denote by $\Gamma^{+,1}_U(P)$\index{Allowable sc$^+$-sections, $\Gamma^{+,1}_U(P)$} the space of all $s\in\Gamma^+(P)$ satisfying  $\supp(s)\subset U$ and $N(s(x))\leq 1$ for all $x\in X$. In our applications  $U$ is the open neighborhood
of the compact solution set of a sc-Fredholm section $f$.  We shall refer to $\Gamma^{+,1}_U(P)$ as the space of {\bf allowable sc$^+$-sections}.
\index{Allowable sc$^+$-sections}
The space  $\Gamma^{+,1}_U(P)$\index{$\Gamma^{+,1}_U(P)$} becomes a metric space with respect to the uniform distance
defined by 
$$
\rho(s,s')=\sup_{x\in X}\{N(s(x)-s'(x))\,  \vert \, x\in X\}.
$$
The metric space  $(\Gamma^{+,1}_U(P),\rho)$ is not complete.

\begin{remark}
We note, however, that if a section  $s$ belongs to the completion of ${\Gamma}^{+,1}_U(P)$,  then the solution set of $f(x)+s(x)=0$ is still compact provided $f$ has a compact solution set and $U$ is an open neighborhood adapted
to the auxiliary norm $N$,  in the sense that $N(f(x))\leq 1$ for $x\in U$,  has the properties stipulated in the  properness result.
\qed
\end{remark}

\begin{theorem}[{\bf Perturbation: interior case}]\label{p:=}\index{T- Perturbation: interior case}
Let $P\colon Y\rightarrow X$ be  a strong bundle over the   M-polyfold $X$ satisfying $\partial X=\emptyset$, and assume that $X$ admits sc-smooth bump functions.
Let $f$ be a sc-Fredholm section with a compact solution set and $N$ an auxiliary norm. Then there exists an open neighborhood
$U$ of the solution set $S=\{x\in X\, \vert \, f(x)=0\}$ such  that for sections $s\in\Gamma^+(P)$ having the property that  
$\supp(s)\subset U$ and  $N(f(x))\leq 1$ for all $x\in X$,  the solution set
$S_{f+s}=\{x\in X\,  \vert \, f(x)+s(x)=0\}$ is compact. Moreover,   there exists a dense subset ${\mathcal O}
$ of the metric space $(\Gamma^{+,1}_U(P),\rho)$ such  that for every $s\in {\mathcal O}$ the solution set $S_{f+s}$ has the property that
$(f+s)'(x)\colon T_xX\rightarrow Y_x$ is surjective for all $x\in S_{f+s}$, i.e. for all $x\in S_{f+s}$ the germ $(f,x)$ is in good position.
In particular,  $S_{f+s}$ is a sub-M-polyfold whose induced structure is equivalent to the structure of a compact smooth manifold
without  boundary.
\qed
\end{theorem}

One can follow the proof of Theorem 5.21 in \cite{HWZ3}. In \cite{HWZ3}  we still worked with splicing cores.
For the convenience of the reader we therefore sketch the proof, the details can be filled in using the arguments  from  \cite{HWZ3}.
\begin{remark}
In \cite{HWZ3} we assume that the fibers of the strong bundle are separable sc-Hilbert spaces. This is in fact not needed
due to an improved  treatment of the compactness in the present text. Also,  originally an auxiliary norm
had to satisfy more properties involving weak convergence, which,  again due to the improved compactness results,  is not needed.
The strategies of the proofs in this more general context are the same.
\qed
\end{remark}

\begin{definition}\label{Def537}\index{D- Transversal to the zero-section}
The  sc-Fredholm section $g$ of the tame strong M-polyfold bundle $P\colon Y\to X$ is called {\bf transversal to the zero-section} if,  at every point $x$ satisfying  $g(x)=0$, 
the linearization $g'(x)\colon T_xX\rightarrow Y_x$ is surjective.
\qed
\end{definition}

\begin{proof}[Theorem \ref{p:=}]
We  choose  $s_0'$ in $\Gamma^{+,1}_{U}(P)$ and for given $\varepsilon>0$ we find $0<\delta<1$ such  that
$$
\abs{s_0'(x)-\delta s_0'(x)}_1\leq \varepsilon/2 \quad  \text{for all $x\in X$}.
$$
Define $s_0=\delta s_0'$. Consider the solution set $S_{f+s_0}$ which we know is compact. 
For every $x\in S_{f+s_0}\subset X_\infty$ we find finitely  many allowable $\ssc^+$-sections $s_1^x,\ldots, s^{k_x}_x$ $\Gamma^{+,1}_U(P)$ such  that the range of 
$(f+s_0)'(x)$ together with the sections $s_j^x$ span $Y_x$.
Then, abbreviating $\lambda=(\lambda_1,\ldots ,\lambda_{k_x})$,  the map
$$
{\mathbb R}^{k_x}\oplus X\rightarrow W, \quad  (\lambda,y)\rightarrow f(y)+s_0(y)+\sum_{j=1}^{k_x} \lambda_j \cdot s^j_x(y)
$$
is a sc-Fredholm section of the obvious pull-back bundle,  in view of the theorem about parameterized perturbations,
Theorem \ref{corro}. We also note that the linearization at the point $(0,x)$ is surjective. In view of  the interior case 
of the implicit function theorem, Theorem \ref{implicit-x}, there exists an open neighborhood $U(x)\subset X$ of $x$
such  that for every $(0,y)$ with $y\in U(x)\cap S_{f+s_0}$ the linearization of the above section is surjective.
We can carry out the above construction for every $x\in S_{f+s_0}$, and  obtain an open covering $(U(x))_{x\in S_{f+s_0}}$
of $S_{f+s_0}$. Hence we find a finite open cover $(U(x_i))_{i=1,\ldots ,p}$.  For every $i$ we have sections 
$s^{j}_{x_i}$, $1\leq j\leq k_{x_i}$.  For simplicity of notation, we list all of  them as $t_1,\ldots,t_m$. Then,  by construction, the section
$$
{\mathbb R}^m\oplus X\rightarrow W, \quad (\lambda,y)\mapsto f(y)+s_0(y)+\sum_{j=1}^m \lambda_j\cdot t_j(y)
$$
is sc-Fredholm,  and its  linearization at every point $(0,x)$ with $x\in S_{f+s_0}$  is surjective. There is a number $\delta_0>0$ such  that
$\abs{\lambda}<\delta_0$ implies $N(\sum_{j=1}^m \lambda_j\cdot t_j(y))<\varepsilon/2$.  By the implicit function theorem
we find an open neighborhood ${\mathcal U}\subset \{\lambda \in \R^m\, \vert \, \abs{\lambda}<\delta_0 \}\oplus X$ of $\{0\}\times S_{f+s_0}$ such that the section 
$$
F\colon {\mathcal U}\rightarrow Y,\quad  F(\lambda,y)=f(y)+s(\lambda,y):=f(y)+s_0(y)+\sum_{j=1}^m \lambda_j\cdot t_j(y)
$$
has a surjective linearization at every solution 
$(\lambda, y)\in {\mathcal U}$ of the equation $F(\lambda,y)=0$. Moreover, $\wt{S}:=F^{-1}(0)$ is a smooth manifold containing $\{0\}\oplus S_{f+s_0}$. 
Taking a regular value $\lambda_0$ for the  projection
$$
\wt{S}\rightarrow {\mathbb R}^m, \quad (\lambda,y)\mapsto  \lambda,
$$ 
it is easily verified  that
$f+s_0+s(\lambda_0, \cdot )$ is transversal to the zero section, see Theorem 5.21 in \cite{HWZ3}. By construction,  $N(s(\lambda_0,y))\leq \varepsilon/2$
and $N(s_0(y)+s(\lambda_0,y))\leq 1-\varepsilon/2$ for all $y$. 
This completes the proof of Theorem \ref{p:=}. 
\end{proof}

\begin{remark}\index{R- On homotopies}\label{rem_homotopy}
In practice we need to homotope from one sc-Fredholm operator
to the other. For example assume that $f_0$ and $f_1$ are sc-Fredholm sections for $P\colon Y\rightarrow X$, 
both transversal to the zero section and having compact solution sets.
 Suppose further that $f_t$, $t\in [0,1]$,  is an interpolating arc satisfying the following.
 First of all, the section 
 $$
 [0,1]\times X\rightarrow Y,\quad (t,x)\rightarrow f_t(x)
 $$
 is sc-Fredholm and has a compact solution set. Now we can use the above construction for a given auxiliary norm
 to find an open neighborhood $U$  of the solution set and a small perturbation $s$ supported in $U$ satisfying 
 $s(t,\cdot )=0$  for $t$ close to $t=0,1$ (we already have transversality at the boundaries) and such  that $(t,x)\mapsto  f_t(x)+s(t,x)$ is transversal to the 
 zero-section. Then the solution set is a compact smooth cobordism between the originally given solution sets  $S_{f_i}$ for $i=0,1$.
Here  we have to deal with the boundary situation which, in this 
special case,  is trivial. The reader will be able to carry out this construction in more detail
 once we have finished  our general discussion of the boundary case.
 \qed
 \end{remark}
 
 The next result shows that under a generic perturbation we are able to bring the solution set into a general position to the boundary
 and achieve the transversality to the zero-section. The solution space is then a smooth manifold with boundary with corners.
 
 \begin{definition}[{\bf General position}]\index{D- General position of $(f,x)$}\label{DEF_539}
Let $P\colon  Y\rightarrow X$ be a strong bundle over the tame  M-polyfold $X$ and let $f$ be a  sc-Fredholm section. We say that $(f,x)$, where $f(x)=0$, is in {\bf general position}
if  $f'(x)\colon T_xX\rightarrow Y_x$ is surjective and the kernel $\ker(f'(x))$ has a sc-complement contained in the reduced tangent space $T^R_xX$.
\qed
\end{definition}

The associated result is the following.
\begin{theorem} [{\bf Perturbation: general position}]\label{thm_pert_and_trans}\index{T- Perturbation: general position}
Assume that we are given a strong bundle $P\colon  Y \rightarrow X$  over the  tame  M-polyfold $X$ which admits sc-smooth bump functions.
Let $f$ be a sc-Fredholm section with compact solution set and $N$ an auxiliary norm. Then there exists an open neighborhood $U$ of the solution set 
$S=\{x\in X\,  \vert \, f(x)=0\}$ such  that for a section $s\in\Gamma^+(P)$ satisfying  $\supp (s)\subset U$ and $N(s(x))\leq 1$ for all $x\in X$,  the solution set
$S_{f+s}=\{x\in X\,  \vert \,  f(x)+s(x)=0\}$ is compact. Moreover,   there exists a dense subset ${\mathcal O}
$ of the metric space $(\Gamma^{+,1}_U(P),\rho)$ such  that,  for every $s\in {\mathcal O}$,  the solution set $S_{f+s}$ has the property 
that for every $x\in S_{f+s}$,  the pair $(f+s,x)$ is in general position.
In particular,  $S_{f+s}$ is a sub-M-polyfold whose induced structure is  equivalent to the structure of a compact smooth manifold with boundary with corners and  $d_{S_{f+s}}(x)=d_X(x)$ for all $x\in S_{f+s}$. 
\end{theorem}

We follow the ideas of the proof of Theorem 5.22 in \cite{HWZ3}.
\begin{proof}
By assumption,  $f$ has a compact solution set $S_f$. Given an auxiliary norm  we can find an open neighborhood $U$ of $S_f$
such  that the solution set $S_{f+s}$ is compact for every section $s\in \Gamma^{+,1}_U(P)$. 

In order to prove the result we will choose  $s_0\in \Gamma^{+,1}_U(P)$ and perturb nearby by introducing
suitable $\ssc^+$-sections. 
We note that  a good approximation of $s_0\in\Gamma^{+,1}_U(P)$ is $\delta s_0$ where  $\delta<1$ is close to $1$ so that we have to find a small perturbation of the latter.
If we take $s_0$ satisfying  $N(s_0(y))\leq 1-\varepsilon$ for all $y\in X$, 
there is no loss of generality assuming that $s_0=0$ by replacing $f+s_0$ by $f$ and $N$ by $cN$ for some large $c$. 
 The general strategy already appears in the proof of Theorem \ref{p:=}.
Here the only complication is that we would like to achieve additional  properties of the perturbed problem. This requires a more sophisticated set-up. 
However, for the reader familiar with finite-dimensional transversality questions, there are no surprises in the proof.

In the next step we choose enough sections
in $\Gamma^{+,1}_U(P)$, say $s_1,\ldots, s_m$ such  that,  near every $(0,x)$ with $x\in S_f$,  the section 
$$
F(\lambda,y)=f(y)+\sum_{i=1}^m\lambda_i\cdot s_i(y)
$$
has suitable properties. Namely,  we require the following properties:
\begin{itemize}
\item[(i)]\ $F'(0,x)\colon  {\mathbb R}^m\oplus T_xX\rightarrow Y_x$ is surjective.
\item[(ii)]\ $\ker(F'(0,x))$ is transversal to ${\mathbb R}^m\oplus  T^R_xX\subset {\mathbb R}^m\oplus T_xX$.
\end{itemize}
The strategy of the proof is the same as the strategy in  the proof of Theorem \ref{p:=}.  We fix a point $(0,x)$ with $x\in S_f$ and observe that if we have $\ssc^+$-sections
so that the properties (i)-(ii) hold at this specific $(0,x)$, then adding more sections,  the properties (i)-(ii) will still hold. 
Furthermore,   if for a section the properties (i)-(ii) hold at the specific $(0,x)$ , then they  will also hold 
at $(0,y)$ for  $y\in S_f$ close to $x$,  say for $y\in U(x)\cap S_f$. As a consequence we only have to find
the desired sections for a specific $x$ and then,  noting that the collection of neighborhoods $(U(x))$ is an open cover of $S_f$, 
we can choose  finitely many  points $x_1,\ldots ,x_p$ such  that the neighborhoods $U(x_1),\ldots , U(x_p)$ cover $S_f$. The collection of sections associated to these finitely many points 
then possesses  the  desired properties. Therefore,  it is enough to give the argument at a general point $(0,x)$ for  $x\in S_f$.
The way to  achieve  property (i) at $(0,x)$ is  as in the proof of Theorem \ref{p:=}. We take enough $\ssc^+$-sections to obtain the surjectivity.
We take a linear subspace $L$ complementing $T_x^RX$ in $T_xX$ and add sections,  which at $x$ span the image $f'(x)L$.  At this point the combined system of sections already satisfies (i) and (ii) and,  
taking the finite union of all these sections,  the desired properties at $(0,x)$ hold. 
By the previous discussion this completes the construction.

Since $S_f$ is compact  and since,  by construction, the section 
$(F,(x,0))$ is in general position at every  point $(0,x)$, we can apply the implicit function theorem to the section
$$
F(\lambda,y)=f(y)+\sum_{i=1}^m \lambda_i\cdot s_i(y).
$$
We deduce  the existence of $\varepsilon>0$ such  that the set 
$$
\wt{S}=\{(\lambda,y)\in {\mathbb R}^m\oplus X\,  \vert \, \text{$\abs{\lambda}<2 \varepsilon$ and $y\in X$}\}
$$
is a smooth manifold with boundary with corners containing $\{0\}\times S_f$. In addition, the  properties (i)-(ii) hold
for all $(\lambda,y)\in \wt{S}$ and not only for the points $(0,x)\in \{0\}\times S_f$.

To be precise,  $\wt{S}\subset {\mathbb R}^m\oplus X$ is a smooth submanifold with boundary with corners 
so that for every $z=(\lambda,x)\in \wt{S}$ the tangent space $T_z\wt{S}$ has a sc-complement in $ T^R_z({\mathbb R}^m\oplus X)$.  If  $p\colon  \wt{S}\rightarrow {\mathbb R}^m$
is  the projection,  the set $p^{-1}(\{\lambda\, \vert \ \abs{\lambda}\leq \varepsilon\})$ is compact. Hence we can apply 
Theorem \ref{SARD} and find,  for a subset $\Sigma$ of $B^m_\varepsilon$ of full measure,   that
the subset $S_\lambda=\{x\in X\,  \vert \,  (\lambda,x)\in \wt{S}\}$ of $X$ is a smooth submanifold with boundary with corners 
having the property  that every point is in general position, so that   for every  point $x\in S_\lambda$ and every parameter $\lambda\in\Sigma$, the tangent space $T_zS_\lambda$ has a sc-complement contained in $T^R_xX$.

\end{proof}

The third result is concerned with a relative perturbation, which vanishes at the boundary in case we already know that at the boundary we are in a good position. If we have a sc-Fredholm germ $(f,x)$ then a {\bf good position}\index{Good position of $(f,x)$} requires $f'(x)$ to be surjective and the
kernel to be in good position to the partial quadrant $C_xX$.

Such a germ would be in {\bf general  position}\index{General position of $(f,x)$} if we require in addition to the  surjectivity of $f'(x)$, that the kernel has a sc-complement in $T^R_xX$.  Clearly general position implies good position.
\begin{remark}\index{R- On general position}
In SFT  we have a lot of algebraic structure combining a possibly 
infinite family of Fredholm problems. In this case perturbations should respect the algebraic structure and genericity in these cases
might mean genericity within the algebraic constraints. In some of these cases general position is not achievable, but one can  still achieve a good
position. The perturbations occurring in such a  context  are very often constructed inductively,
so that at each step the problem is already in good position at the boundary, but has to be extended to a generic problem.
\qed
\end{remark}
The following theorem is a sample result along these lines.

\begin{theorem}[{\bf Perturbation and transversality: good position}]\index{T- Perturbation: good position}\label{THMX5.3.12}
Assume that $P\colon  Y\rightarrow X$ is a strong bundle 
over the the tame  M-polyfold $X$ which admits sc-smooth bump functions.
Let $f$ be a sc-Fredholm section with compact solution set and $N$ an auxiliary norm. Further,  assume 
that,  for every $x\in \partial X$ solving  $f(x)=0$, the pair $(f,x)$ is in good position. Then there exists an open neighborhood
$U$ of the solution set $S=\{x\in X\, \vert \,  f(x)=0\}$ so that,  for every  section $s\in\Gamma^{+,1}_U(P)$,  the solution set
$S_{f+s}=\{x\in X\ |\ f(x)+s(x)=0\}$ is compact. Moreover,  there exists an arbitrarily small section $s\in\Gamma^{+,1}_U(P)$ satisfying $s(x)=0$ near $\partial X$
such  that $f+s$ is transversal to the zero-section and for every $x\in S_{f+s}$ the pair $(f,x)$ is in good position.
In particular,  $S_{f+s}$ is a M-subpolyfold whose  induced structure is equivalent to the structure of a compact smooth manifold with  boundary with corners.
\end{theorem}

\begin{proof}

By our previous compactness considerations there exists,  for a given auxiliary norm $N$,  an open neighborhood $U$ of $S_f$
so that for $s\in \Gamma^{+,1}_U(P)$ the solution set $S_{f+s}$ is compact.
By the usual recipe already used in the previous proofs we can find finitely many sections $s_1,\ldots ,s_m$ in $\Gamma^{+,1}_U(P)$ 
which are vanishing near $\partial X$ so that for every $x\in S_f$ the image $R(f'(x))$ and the $s_i(x)$ span $Y_x$.
Of course, in the present construction we are allowed to have sections which vanish near $\partial X$, since by assumption
for $x\in S_f\cap \partial X$ we are already in good position (which in fact implies that $S_f$ is already a manifold with boundary with corners
near $\partial X$). Then we consider as before the section 
$$
(\lambda,x)\mapsto   f(x)+\sum_{i=1}^m \lambda_i\cdot s_i(x),
$$
and,  for a generic value of $\lambda$, which we can take as small as we wish,  we conclude that the associated section $s_\lambda=
\sum_{i=1}^m \lambda_i\cdot s_i$ has the desired properties.

\end{proof}

The next result deals with a homotopy $t\mapsto f_t$ of sc-Fredholm sections during which also the bundle  changes. 
\begin{definition}[{\bf Generalized compact homotopy}]\label{gen_comp_hom}
Consider two sc-Fredholm sections $f_i$ of tame strong bundles $P_i\colon  Y_i \rightarrow X_i$ having  compact solution sets.
We shall refer to $(f_i,P_i)$ as two {\bf compact sc-Fredholm problem}\index{D- Compact sc-Fredholm problem}. 
Then a {\bf generalized compact homotopy between the two compact  sc-Fredholm problems}\index{D- Compact homotopy}  consists of a tame strong bundle $P\colon  Y\rightarrow X$
and a sc-Fredholm section $f$, where $X$ comes with a sc-smooth surjective map $t\colon  X\rightarrow [0,1]$, so that the preimages
$X_t$ are tame M-polyfolds and $f_t=f|X_t$ is a sc-Fredholm section of the bundle $Y\vert X_t$. Moreover,  $f$ has a compact solution set
and $(f\vert X_i,P\vert Y_i)=(f_i,P_i)$ for $i=0,1$.
\qed
\end{definition}

\begin{remark}
Instead of requiring $(f\vert X_i,P\vert Y_i)=(f_i,P_i)$ for $i=0,1$ one should better require that the problems 
are isomorphic and make this part of the data. But in applications the isomorphisms are mostly clear, so that we allow ourselves
to be somewhat sloppy.
\qed
\end{remark}

\begin{theorem}[{\bf Morse-type structure}]\label{MORSE-type}\index{T- Fredholm homotopy, Morse-type}
We assume that all occurring M-polyfolds admit sc-smooth bump functions. 
Let $f$ be a sc-Fredholm section of the tame strong bundle $P\colon Y\to X$ which is a generalized compact homotopy between the compact Fredholm  problems $(f_i, P_i)$ for $i=0, 1$, as in Definition \ref{gen_comp_hom}.
We assume that $P_i\colon  Y_i\rightarrow X_i$ are  strong bundles over M-polyfolds $X_i$ having no boundaries,  and that  the Fredholm sections $f_i$
are already generic in the sense that for all $x\in S_{f_i}$ the germ $(f_i,x)$ is in general position.
We assume further that $\partial X=X_0\sqcup  X_1$.  Let $N$ be an auxiliary norm on $P$. Then there exists an open neighborhood $U$ 
of the solution set $S_f$ in $X$ such  that,  for all sections $s\in \Gamma^{+,1}_U(P)$,  the solution set $S_{f+s}$ is compact. Moreover,  there exists an arbitrarily small section 
$s_0\in\Gamma^{+,1}_U(P)$ which vanishes near $\partial X$,  possessing the following properties.
\begin{itemize}
\item[{\em (1)}]\ For every $x\in S_{f+s_0}$,  the germ $(f+s_0,x)$ is in general position.
\item[{\em (2)}]\ The smooth function $t\colon  S_{f+s_0}\rightarrow [0,1]$ has only Morse-type critical points.
\end{itemize}
\end{theorem}
\begin{proof}
Applying the  previous discussions we can achieve property  (1) for a suitable section $s_0$. The idea then is to perturb $f+s_0$ further
to achieve also property (2). Note that (1) is still true after a small perturbation.  Hence it suffices to assume that $f$ already has the property that
$(f,x)$ is in good position for all $x$ solving  $f(x)=0$. The solution set $S=\{x\in X\, \vert \, f(x)=0\}$ is a compact manifold
with smooth boundary components. Moreover,  $d(t\vert S)$ has no critical points near $\partial S$ in view of  the assumption
that $(f_i,X_i)$ are  already in general position. If we take a finite number of $\ssc^+$-sections $s_1, \ldots, s_m$ of $P\colon Y\to X$, which vanish near $\partial X$ and are supported near $S$ (depending on the auxiliary  norm $N$),
then the solution set
$\wt{S}=\{(\lambda,x)\,\vert \,  f(x)+\sum_{i=1}^m\lambda_i\cdot s_i(x)=0,\ \abs{\lambda}<\varepsilon\}$ is a smooth manifold for $\varepsilon$ is small enough. Using the classical standard implicit function theorem, we find a family of smooth embeddings $\Phi_\lambda\colon  S\rightarrow \wt{S}$ having the property that $\Phi_0(x)=(0,x)$ and $\Phi_\lambda(S)=\{\lambda\}\times S_\lambda$, where
$S_\lambda=\{y\in X\, \vert \,  f(x)+\sum_{i=1}^m\lambda_i\cdot s_i(y)=0\}$. 

Our aim is to construct the above sections 
$s_i$ in such a way
that, in addition, the map
$$
(\lambda,y)\mapsto d(t\circ\Phi_\lambda(y))\in T_y^\ast S
$$
is transversal at $\{0\}\times S$ to the zero section in cotangent bundle $T^\ast S$.
Then, after having achieved this, the 
parameterized version of Sard's theorem 
will guarantee values of the parameter 
$\lambda$ arbitrarily close to $0$, for which the smooth section
$$
S\ni y\mapsto  d(t\circ\Phi_\lambda(y))\in T^\ast S
$$
is transversal to the zero-section and hence the function $y\mapsto  t\circ\Phi_\lambda(y)$ will be  a Morse function on $S$.
Since $\Phi_\lambda\colon  S\rightarrow S_\lambda$ is a diffeomorphism,  we conclude  that $t\vert S_\lambda$ is a Morse-function.
Having constructed this way the desired section $s_0$, the proof of the theorem will then be complete.

It  remains to construct the desired family of $\ssc^+$-sections $s_1, \ldots ,s_m$.
We fix a critical point $x\in S$ of the function $t\vert S\to [0,1]$, hence $d(t\vert S)(x)=0$. Then $T_xS\subset \ker (dt(x))$. Since $x$ is a smooth point,  
we find a one-dimensional smooth linear sc-subspace $Z\subset T_xX$ such that
$$
T_xX = Z\oplus \ker(dt(x)).
$$
The proof of the following trivial observation is left to the reader.
\begin{lemma}\label{lemma5.29}
For every element $\tau\in T_x^\ast S$ there exists a  uniquely determined sc-operator 
$$
b_\tau \colon  T_xS\to  Z
$$
satisfying 
$$
dt(x)\circ b_\tau = \tau.
$$
If $a\colon  T_xS\rightarrow \ker(dt(x))$ is a sc-operator, then
$$
dt(x)\circ (b_\tau + a) =\tau.
$$
\qed
\end{lemma}

Slightly more difficult is the next  lemma.
\begin{lemma}\label{lemma5.30}
Assume that $\tau \in T^\ast_xS$ is given.  Then there exists a $\ssc^+$-section $s$ with sufficiently small support around $x$ having the following properties.
\begin{itemize}
\item[{\em (1)}]\ $s(x)=0$.
\item[{\em (2)}]\ $f'(x)\circ b_\tau +  (s'(x)\vert T_xS)=0.$
\end{itemize} 
We note that if $a\colon  T_xS\rightarrow T_xS$,  then also property (2) holds with $b_\tau $ replaced by $b_\tau+a$,  in view of $T_xS = \ker(f'(x))$.
\end{lemma}

\begin{proof}
Since $f$ is a sc-Fredholm section and 
$f(x)=0$, 
the linearization 
$$
f'(x)\colon  T_xX\rightarrow Y_x
$$
 is surjective and the kernel is equal to  $T_xS$.
Then $f'(x)\circ b_\tau \colon   T_xS\rightarrow Y_x$ is a sc-operator, i.e.,  the image of any vector in $T_xS$ belongs to $Y_\infty$. 
If this operator is the zero operator we can take $s=0$. Otherwise the operator has a one-dimensional image spanned
by some  smooth point $e\in Y_x$. We are done if we can construct a $\ssc^+$-section $s$ satisfying  $s(x)=0$,  and having support close to $x$, and  
 $s'(x)\vert T_xS\colon  T_xS\rightarrow Y_x$  has a one-dimensional image spanned by $e$,  and
$\ker(s'(0)\vert T_xS)= \ker(f'(x)\circ b_\tau)$. Then a suitable multiple of $s$ does the job. 

Denote by $K\subset T_xS$ the kernel of $f'(x)\circ b_\tau$  and by $L$ a complement of $K$ in $T_xS$. Then $L$ is one-dimensional.
We work now in local coordinates in order to construct $s$. We may assume that $x=0$ and represent $S$ near $0$ as a graph over the tangent space $T_xS$, say $q\mapsto  q+\delta(q)$ with $\delta(0)=0$ and $D\delta(0)=0$. Here $\delta\colon  {\mathcal O}(T_xS,0)\rightarrow V$, where $V$ is a sc-complement of $T_xS$ in the sc-Banach space $E$. The points in $E$ in a neighborhood of $0$ are of  the form
$$
q+\delta(q) + v,
$$
where  $v\in V$. We note that $q+\delta(q)\in O$, where $O$ is the local model for $X$ near $0$.
We split $T_xS=K\oplus L$ and correspondingly write $q=k+l$. Then we can represent  the points in a neighborhood of $0\in E$ in the form
$$
k+l +\delta(k+l) +v.
$$
Choosing  a linear isomorphism $j\colon  L\rightarrow {\mathbb R}$, we define the section $\wt{s}$ for $(k,l,v)$ small by
$$
k+l+\delta(k+l)+v\mapsto  R(k+l +\delta(k+l) +v,\beta(k+l+\delta(k+l)+v)j(l)e),
$$
where $\beta$ has support around $0$ (small) and $\beta$  takes the value $1$ near $0$. The section $s$ is then the restriction of $\wt{s}$ to
$O$. If we restrict $s$ near $0$ to $S$ we obtain 
$$
s(k+l+\delta(k+l)) = R(k+l+\delta(k+l),j(l) e).
$$
Hence  $s(0)=0$,  and the linearization of $s$ at $0$ restricted to $T_xS$ is given by
$$
s'(0)(\delta k+\delta l) = j(\delta l) e.
$$
This implies that $s'(0)\vert T_xS$ and $f'(0)\circ b_\tau$ have the same kernel and their image is spanned by $e$. Therefore,  $s$,  multiplied
by a suitable scalar,  has the desired properties and the proof of Lemma \ref{lemma5.30} is complete.
\end{proof}

Continuing with the proof of Theorem \ref{MORSE-type}  we focus as before on the critical point $x\in S$ of $t\vert S$, which satisfies   $d(t\vert S)(x)=0$. Associated with  a basis  $\tau_1,\ldots  ,\tau_m$ of $T_x^\ast S$ the previous lemma produces the sections $s_1, \ldots ,s_m$.
Consider the solution set $\wt{S}$ of solutions  $(\lambda,y)$ of $f(y)+\sum_{i=1}^m\lambda_i\cdot s_i(y)=0$.   Since $f'(y)$ is onto for all $y\in S$,  the solution set 
$S_\lambda=\{y\, \vert \,  f(y) +\sum_{i=1}^m \lambda_i\cdot s_i(y)=0\}$ is a compact manifold (with boundary) diffeomorphic to 
$S$ if $\lambda$ small. Moreover,  $\wt{S}$ fibers over a neighborhood of zero via the map $(\lambda,y)\mapsto  \lambda$. 

The smooth map
\begin{equation}\label{eqp}
(\lambda,y)\mapsto  d(t\vert S_\lambda)(y)\in T_yS_\lambda.
\end{equation}
is a smooth section of the bundle over $\wt{S}$ whose  fiber at $(\lambda,y)$ is equal to $T_y^\ast M_\lambda$.

We  show that the linearization of \eqref{eqp} at $(0,x),$ which is a map
$$
{\mathbb R}^d\oplus T_xS\rightarrow T_x^\ast S, 
$$
is surjective. Near $(0,x)\in \wt{S}$ we can parameterize $\wt{S}$,  using the implicit function theorem,  in the form
$$
(\lambda,y)\mapsto  (\lambda,\Phi_\lambda(y))
$$
where  $\Phi_0(y)=y$, and $\frac{\partial\Phi}{\partial\lambda_i}(0,x)=0$.
Using \eqref{eqp}  and the map $\Phi$ we obtain,  after a coordinate change on the base
for  $\lambda$ small and $z\in S$ near $x$, the map 
$$
(\lambda,z)\mapsto  d(t\circ \Phi_\lambda)(z) = dt(\Phi_\lambda(z))T\Phi_\lambda(z),
$$
where $d$ acts only on the $S$-part. By construction,  the section vanishes at $(0,x)$.
Recall that,  by construction,  $\Phi_\lambda(x)=x$. Hence for fixed $\delta \lambda$ the map 
$$
z\mapsto  \sum_{i=1}^m \delta\lambda_i \cdot \frac{\partial\Phi}{\partial\lambda_i}(0,z)
$$
is a vector field defined near  $x\in S$ which vanishes at $x$. Therefore, it has a well-defined linearization at $x$.
The linearization ${\mathbb R}^d\oplus T_xS\rightarrow T_x^\ast S$ at $(0,x)$ is  computed to be the mapping  
\begin{equation}\label{new_eq5_68}
\begin{split}
(\delta \lambda,\delta z)\mapsto &(dt\vert S)'(x)\delta z + dt(x)\biggl(\sum_{i=1}^m \delta\lambda_i\biggl( \frac{\partial\Phi}{\partial\lambda_i}(0,\cdot )\biggr)'(x)\biggr)\\
&=(dt\vert S)'(x)\delta z + \sum_{i=1}^m \delta\lambda_i \cdot dt(x)\circ \biggl(\frac{\partial\Phi}{\partial\lambda_i}(0,\cdot )\biggr)'(x).
\end{split}
\end{equation}
The derivative  $(dt\vert S)'(x)$ determines  the Hessian of the map $t\vert S$ at the point $x\in S$. The argument is complete if we can show that
\begin{equation}\label{question1}
\biggl(\frac{\partial\Phi}{\partial\lambda_i}(0,\cdot )\biggr)'(x) =b_{\tau_i} +a_i,
\end{equation}
where the image of $a_i$ belongs to  $T_xS$.
By the  previous discussion,   $dt(x)\circ b_{\tau_i}=\tau_i$ so that the map \eqref{new_eq5_68} can be rewritten, using $dt(x)\circ a_i=0$, as
$$
(\delta \lambda,\delta z)\mapsto   (dt\vert S)'(x)\delta z + \sum_{i=1}^m \delta\lambda_i \cdot \tau_i,
$$
which then proves our assertion. So, let us show that the identity \eqref{question1} holds. We first linearize the equation
$$
f(\Phi_\lambda(z))+\sum_{i=1}^d \lambda_i s_i(\Phi_\lambda(z))=0
$$
with respect to $\lambda$ at $\lambda=0$, which  gives 
$$
Tf(z)\biggl(\sum_{i=1}^m \delta\lambda_i\cdot \frac{\partial\Phi}{\partial\lambda_i}(0,z)\biggr) +\sum_{i=1}^d \delta\lambda_i\cdot s_i(z)=0.
$$
Next we linearize with respect to $z$ at $z=x$, leading  to 
$$
\sum_{i=1}^m \delta\lambda_i\cdot \biggl(f'(x)\circ \biggl(\frac{\partial\Phi}{\partial\lambda_i}\biggr)'(0,x) + s'_i(x)\vert T_xS\biggr)=0
$$
for all $i=1,\ldots  ,m$. Hence $f'(x)\circ \bigl(\frac{\partial\Phi}{\partial\lambda_i}\bigr)'(0,x) + s'_i(x)\vert T_xS=0$. This implies that 
$\bigl(\frac{\partial\Phi}{\partial\lambda_i}\bigr)'(0,x) =b_{\tau_i} +a_i$, where the image of $a_i$ lies  in the kernel of $f'(x)$, i.e., in  $T_xS$. 
At this point we have proved that the linearization of \eqref{eqp} at $(0,x)$ is surjective. Since the section is smooth,  there exists an open neighborhood $U(x)$ of $x$  in $S$ so that,  if at $(0,y)$ we have $d(t\vert S)(y)=0$, then the linearization of
$(\lambda,z)\rightarrow d(t\vert S_\lambda)(z)$ at $(0,y)$ is surjective.

We can now apply the previous discussion to all points $x$ solving $d(t\vert S)(x)=0$ and,  using the compactness,  we find finitely many such points $x_1,\ldots  ,x_k$
so that the union of all the $U(x_i)$ covers the critical points of $t\vert S$. For every $i$ we have $\ssc^+$-sections $s_1^i,\ldots  ,s_{m_i}^i$ possessing  the desired properties.  In order to simplify the notation we denote the union of these sections by $s_1,\ldots  ,s_d$. Then we consider the solutions of 
$$
f(y)+\sum_{i=1}^d \lambda_i\cdot s_i(y)=0.
$$
Again we denote the solution set by $\wt{S}$. It fibers over an open neighborhood of $0$ in ${\mathbb R}^d$. By construction, 
the smooth map
$$
(\lambda,z)\mapsto  d(t\vert S_\lambda)(z)\in T_zS_\lambda
$$
has,  at every point $(0,y)$ satisfying $d(t\vert S)(y)=0$,  a linearization
$$
{\mathbb R}^d\oplus T_yS\rightarrow T_y^\ast S
$$
which is surjective. Now we take a regular value $\lambda$ (small) for the projection $\wt{M}\rightarrow {\mathbb R}^d$
and find,   by the parameterized version of Sard's theorem, that
$$
d(t\vert S_\lambda)
$$
is indeed a Morse-function. The proof of  Theorem \ref{MORSE-type} is complete.
\end{proof}

\section{Remark on  Extensions of Sc\texorpdfstring{$^+$}{piy}-Sections}
In this section we state a general result concerning the extension of sc$^+$-sections defined over the boundary of a tame M-polyfold
to the whole M-polyfold. The theorem, given without proof, is a special case of a result given in the context of ep-groupoids, i.e.
in the case of local symmetries. More general results are stated and proved in Section \ref{SECR12.3} and Section \ref{SECR14.1}.

Start with a tame M-polyfold $X$ and recall from Proposition \ref{FACE_XXXX} that every point $x\in X$ belongs to precisely
$d_X(x)$-many local faces. These local faces  at a point $x\in X$, as well as their intersections are sub-M-polyfolds of $X$ which contain $x$.
Assume that $P:W\rightarrow X$ is a strong bundle.
\begin{definition}\index{D- Sc$^+$-section over $\partial X$}
A section $s$ of $W|\partial X\rightarrow \partial X$ is called a sc$^+$-section over $\partial X$, provided
for every $x\in \partial X$ the restrictions of $s$ to the local faces are sc$^+$-sections.
\end{definition}
The basic extension result in our situation, which is a special case of Theorem \ref{p-main-p}, is given by the following theorem.
\begin{theorem}[$\partial$-Extension Theorem]\index{T- Extension theorem for sc$^+$-sections}
Assume  $P:W\rightarrow X$ is   a strong bundle over a tame  M-polyfold $X$  admitting  sc-smooth partitions of unity. 
Let $s$ be a sc$^+$-section of $W|\partial X\rightarrow \partial X$ and let 
$N:W\rightarrow [0,\infty]$ be an auxiliary norm.  If $\wt{U}$ is an open neighborhood in $X$
of $\supp(s)$ and $f:X\rightarrow [0,\infty)$ a continuous map  supported in $\wt{U}$, satisfying
$N(s(x))<f(x)$ for all $x\in \supp(s)$.  Then there exists  a sc$^+$-section $\overline{s}$ of $P$
having  the following properties:
\begin{itemize}
\item[\em{(1)}]\ $N(\overline{s}(x))\leq f(x)$ for all $x\in X$.
\item[\em{(2)}]\ $\supp(\overline{s})\subset \wt{U}$.
\item[\em{(3)}]\ $\overline{s}\vert \partial X=s$.
\end{itemize}
\end{theorem}
\begin{proof}
Here is a sketch of the proof.  Using the tameness of $X$ we can employ in suitable local coordinates
the local extension result Proposition \ref{hucky}. Then these local extensions can be glued together
using a sc-smooth partition of unity.
\qed
\end{proof}

\section{Notes on Partitions of Unity and Bump Functions}\label{POU}

An efficient tool for globalizing local construction in M-polyfolds are sc-smooth partitions of unity.

\begin{definition}
A M-polyfold $X$ {\bf admits sc-smooth partitions of unity} if  for every open covering of $X$ there exists a subordinate sc-smooth partition of unity.
\qed
\end{definition}

So far we did not  need sc-smooth partitions of unity for our constructions on M-polyfolds. We would need them, for example, for the construction of sc-connections  M-polyfolds. We point out that a (classically) smooth partition of unity on a Banach space, which is  equipped with a sc-structure,  induces a sc-smooth partition of unity, in view of  Corollary \ref{ABC-y}. However, many Banach spaces do not admit smooth partition of unity subordinate to a given open cover.
Our discussion in this section is based on the survey article \cite{Fry} by Fry  and McManusi on smooth bump functions on Banach spaces. The article contains many interesting open questions.

For many constructions one does not need sc-smooth partitions of unity, but only sc-smooth bump functions.
\begin{definition}
A M-polyfold $X$ admits {\bf admits sc-smooth bump functions} if  for every point $x\in X$
and every open neighborhood $U(x)$ there exists a sc-smooth function $f\colon  X\rightarrow {\mathbb R}$ which is not identically zero and has support in $U(x)$.
\qed
\end{definition}

For example, if $E$ is a Hilbert space with a scalar product $\langle \cdot, \cdot \rangle$ and associated norm $\norm{\cdot }$, the map $x\mapsto \norm{x}^2$, $x\in E$, is smooth. Choosing a smooth function $\beta\colon \R\to \R$ of  compact support and satisfying $\beta (0)=1$, the function $f(x)=\beta (\norm{x}^2)$, $x\in E$, is a smooth, non-vanishing function of bounded support. Therefore, if $E$ is equipped with a sc-structure, the map is a sc-smooth function on the Hilbert space in view of Corollary \ref{ABC-y}. Using the Hilbert structure to rescale and translate, we see that the Hilbert space $E$ admits sc-smooth bump functions.

Currently it is an open problem whether the existence of sc-smooth partition of unity on M-polyfolds $X$ is equivalent to the existence of sc-smooth bump functions. The problem is related to an unsolved classical question in 
Banach spaces: is, in every Banach space, the existence of a single smooth bump function (we can use the Banach space structure to rescale and translate) equivalent 
to the existence of smooth partitions of unity subordinate to given  open  covers? 
For the discussion on this problem we refer to \cite{Fry}.

Sc-smooth bump functions can be used for the construction of functions having special properties, as the following example shows.

\begin{proposition}\label{prop-x5.36}
We assume that the M-polyfold $X$ admits sc-smooth bump functions. Then for every point $x\in X$ and every open neighborhood
$U(x)$ there exists a sc-smooth function $f\colon  X\rightarrow [0,1]$ with $f(x)=1$ and support in $U(x)$. In addition we can choose  $f$ in such a way 
that $f(y)=1$ for all $y$ near $x$.
\end{proposition}
\begin{proof}
By  assumption,  there exists  a sc-smooth bump function $g$ having support in $U(x)$   and satisfying $g(0)=1$. In order to achieve that the image is contained in $[0,1]$ we choose a smooth map $\sigma\colon  {\mathbb R}\rightarrow [0,1]$
satisfying $\sigma(s)=0$ for $s\leq 1$ and $\sigma(s)=1$ for $s\geq 1$ and define  the sc-smooth function $f$ by $f=\sigma\circ g$. If, in addition, we wish $f$ to be constant near $x$,  we take $f(y)=\sigma(\delta\cdot g(y))$ for  $\delta>1$. 
\end{proof}

The survey paper \cite{Fry} discusses, in particular, bump functions on Banach spaces, which are classically differentiable.

\begin{definition} [{\bf $C^k$-bump function}]  A Banach space $E$ 
{\bf admits a $C^k$- bump function}\index{D- $C^k$- bump function},  if there exists a $C^k$-function 
$f\colon  E\rightarrow {\mathbb R}$,  not identically zero and having bounded support.
Here $k\in \{0, 1, 2, \ldots\}\cup \{\infty\}$.
\qed
\end{definition}

The existence of $C^k$-bump functions on $L_p$ spaces is a consequence of the following result due to Bonic and Frampton, see Theorem 1 in \cite{Fry}.
\begin{theorem}[Bonic and Frampton]\label{Bonic_Frampton}
For the $L_p$-spaces the following holds for the usual $L^p$-norm $\norm{\cdot }$.
Let $p\geq 1$ and let $\norm{\cdot}$ be the usual norm on $L^p$.  

\begin{itemize}
\item[{\em (1)}]\ If $p$ is an even integer,  then $\norm{\cdot }^p$ is of class $C^\infty$.
\item[{\em (2)}]\ If i$p$ s an odd integer, then $\norm{\cdot }^p$ is of class $C^{p-1}$.
\item[{\em (3)}]\ If $p\geq 1$ is not an integer, then $\norm{\cdot }^p$ is of class $C^{[p]}$, where $[p]$ is the integer part of $p$.
\end{itemize}
\qed
\end{theorem}
We deduce immediately for the Sobolev spaces $W^{k,p}(\Omega)$ that  the usual norms 
$$
\norm{u}_{W^{k,p}}^p=\sum_{|\alpha|\leq k} \norm{D^\alpha u}^p_{L_p}.
$$
have  the same differentiability as $\norm{\cdot}_{L_p}$.

Taking a non-vanishing smooth function $\beta\colon \R\to \R$ of compact support, 
the map $f(x)=\beta (\norm{x}^p)$, is a bump function of $L_p$, whose smoothness depends on $p$ as indicated in Theorem \ref{Bonic_Frampton}.

If 
$\Omega$ is a bounded domain in ${\mathbb R}^n$ and $E=W^{1,4}(\Omega)$ is equipped with the sc-structure $E_m=(W^{1+m,4}(\Omega)$, $m\geq 0$, we deduce from (1) in Theorem \ref{Bonic_Frampton} that $E$ admits sc-smooth bump functions.
In contrast, if  $E=W^{1,3/2}(\Omega)$ is equipped with the  sc-structure  $E_m=W^{1+m,3/2}(\Omega)$, then the straightforward bump function constructed by using (3) of Theorem \ref{Bonic_Frampton} would only be of class $\ssc^1$. Does there exist a sc-smooth bump function on $E$, i.e., sc-smooth and of bounded support in $E$?

The existence of sc-smooth bump functions is a local property.
\begin{definition}
A local M-polyfold model $(O,C,E)$ has the {\bf sc-smooth bump function property} if,   for every $x\in O$ and every open neighborhood
$U(x)\subset O$ satisfying  $\cl_E(U(x))\subset O$,  there exists a  sc-smooth function $f\colon  O\rightarrow {\mathbb R}$ satisfying  $f\neq 0$ and $\supp(f)\subset U(x)$.
\qed
\end{definition}

Clearly,  the following holds.
\begin{theorem}
A M-polyfold $X$ admits sc-smooth bump functions if and only if it admits a sc-smooth atlas whose  the local models have the sc-smooth bump function property.
\qed
\end{theorem}

The following class of spaces have the sc-bump function property.
\begin{proposition}\label{P548}
Assume that $(O,C,E)$ is a local M-polyfold model in which the $0$-level $E_0$ of the sc-Bananch space $E$ is a Hilbert space. 
Then $(O,C,E)$ has the sc-bump function property.
\qed
\end{proposition}

\begin{proof}
Let $\langle \cdot ,\cdot \rangle$ be the inner product and $\norm{\cdot}$ the associated norm of $E_0$. We choose a smooth function $\beta\colon \R\to \R$ of compact support and satisfying $\beta (0)=1$. Then the function $f(x)=\beta (\norm{x}^2)$, $x\in E$, defines, in view of Corollary \ref{ABC-y}, a sc-smooth function on $E$. Using scaling, translating, and composition with the sc-smooth retraction onto $O$ the proposition follows.
\end{proof}

Next we study the question of the existence of sc-smooth partitions of unity.
By definition, a M-polyfold $X$ is paracompact and, therefore, there exist  continuous partitions of unity. Hence  it is not surprising that the existence of sc-smooth partitions of unity is connected to local properties of $X$,  namely to   the local  approximability 
of continuous functions  by sc-smooth functions.

\begin{definition}\label{approximation_property}
A local M-polyfold model $(O,C,E)$ has the {\bf sc-smooth approximation property} provided the following holds.
Given $(f,V,\varepsilon)$,  where $V$ is an open subset of $O$ such that $\cl_C(V)\subset O$,  $f\colon  O\rightarrow [0,1]$ is a continuous function with support contained  $V$,  and $\varepsilon>0$,
there exists a sc-smooth map $g\colon  O\rightarrow [0,1]$ supported  in  $V$  and satisfying $\abs{f(x)-g(x)}\leq \varepsilon$ for all $x\in O$.
\qed
\end{definition}

\begin{theorem}\label{partition_approximation}
The following two statements are equivalent.
\begin{itemize}
\item[{\em (1)}]\ A M-polyfold $X$ admits sc-smooth partitions of unity subordinate to any given open cover.
\item[{\em (2)}]\ The M-polyfold $X$ admits an atlas consisting of local models having the sc-smooth approximation property.
\end{itemize}
\end{theorem}

\begin{proof} 
Let us first show that (1) implies (2). For the M-polyfold $X$,  we take an atlas of M-polyfold charts $\phi\colon U\to O$. We shall show that the local models $(O, C, E)$ posses  the sc-smooth approximation property. Let $(f,V,\varepsilon)$ be as in Definition \ref{approximation_property}. 
The sc-smooth function $f\circ \phi$ is defined on $\phi^{-1}(O)\subset X$ and we extend it by $0$ to all of $X$ and obtain a continuous function 
$g\colon X\to [0,1]$ whose  support is contained in the open set $W=\phi^{-1}(V)$, satisfying $\cl_X(W)\subset U$. 
Define the open subset $\wt{W}$ of $X$ by
$$
\wt{W}=\{x\in W\, \vert \, g(x)>\varepsilon/4\}.
$$
Then $\cl_X(\wt{W})\subset W$.
For every $x\in \cl_X(\wt{W})$,  there exists open neighborhood $U_x$ of $x$ such that
$\cl_X(U_x)\subset W$ and $\abs{g(x)-g(y)} <\varepsilon/2$ for all $y\in U_x$. Take the  open cover of $X$ consisting of $U_0=\{x\in X\, \vert \, g(x)<\varepsilon/2\}$ and $(U_x)$, $x\in \cl_X(\wt{W})$. By assumption, there exists a subordinate sc-smooth partition of unity consisting of $\beta_0$ with the support in $U_0$ and $(\beta_x)$ supported in $(U_x)$,  $x\in \cl_X(\wt{W})$.  Then we  define the function $\wh{g}\colon  X\rightarrow [0,1]$ by 
$$
\wh{g}(y)=\sum_{x\in \cl_X(\wt{W})} \beta_x(y)g(x).
$$
Since the collection of supports of $\beta_0$ and $(\beta_x)$ is locally finite, the sum is locally finite and, therefore, $\wh{g}$ is sc-smooth.

We claim that $\abs{g(x)-\wh{g}(x)}<\varepsilon$ for all $x\in X$. In order to show this we first show that $\supp (\wh{g})\subset W$. We assume that $y\in \supp (\wh{g})$ and let 
$(y_k)$ be a sequence  satisfying $\wh{g}(y_k)>0$ and $y_k\rightarrow y$. Since $\wh{g}(y_k)>0$, every open neighborhood of $y$ intersects $\supp (\beta_x)$ for some $x\in \cl_X(W)$. Since  the collection of supports $(\supp (\beta_x))$ is locally finite, there exists an open neighborhood $Q=Q(y)$ of $y$ in $X$ and finitely many points $x_1,\ldots ,x_m$ such that only the supports of the functions $\beta_{x_i}$ intersect  $Q$.
Hence
$$
0<\wh {g}(y_k) =\sum_{i=1}^m \beta_{x_i}(y_k)g(x_i), 
$$
showing  that  $y_k\subset \bigcup_{i=1}^m \cl_X(U_{x_i})\subset W$ for large $k\geq 1$.  Consequently,  $y\in W$ and hence $\supp (\wh{g})\subset W$.

Next, given $z\in U$,  there exists an  open neighborhood $Q=Q(z)$ of $z$ in $X$  such that 
$Q$ intersects only finitely many supports of functions belonging to the partition of unity.
If $Q$ intersects the support of $\beta_0$ and none of the supports of the functions $\beta_x$, then $Q\subset U_0$ and $\wh{g}(y)=0$ for all $y\in Q$. Consequently, 
$$\abs{g(y)-\wh{g}(y)}=\abs{g(y)}<\varepsilon/2<\varepsilon$$
for all $y\in Q$. If $Q$ intersects supports of the functions $\beta_x$, there are finitely many points $x_1,\ldots, x_l$ such that only the supports of the functions $\beta_{x_i}$ intersect  $Q$.
Hence 
$$
\wh{g}(y)=\sum_{i=1}^l \beta_{x_i}(y)g(x_i)\quad  \text{for all $y\in Q$}.
$$
Since $\abs{g(x)-g(y)}<\varepsilon/2$ if $y\in U_x$ and $g(y)<\varepsilon/2$ if $y\in U_0$,  we conclude that 
\begin{equation*}
\begin{split}
|g(y)-\wh{g}(y)|
&=\abs{\beta_0(x)g(y)+\sum_{i=1}^l\beta_{x_i}(y)g(y)-\sum_{i=1}^l \beta_{x_i}(y)g(x_i)}\\
&< \varepsilon/2 +\sum_{i=1}^l \beta_{x_i}(z)\abs{g(y)-g(x_i)}\\
&\leq \varepsilon/2 +(\varepsilon/2)\sum_{i=1}^l \beta_{x_i}(y)< \varepsilon.
\end{split}
\end{equation*}
for all $y\in Q$. 
Consequently, $\abs{g(y)-\wh{g}(y)}< \varepsilon$ for all $y\in U$ and 
$\abs{f(x)-\wh{g}\circ\phi^{-1}(x)}< \varepsilon$ for $x\in O$. This completes the proof that (1) implies (2).

Next we show that (2) implies (1).  For the M-polyfold $X$ we take an atlas of M-polyfold charts $\phi_\tau\colon U_\tau \to O_\tau$, $\tau \in T$, where the local models $(O_\tau, C_\tau, E_\tau)_{\tau \in T}$ posses  the sc-smooth approximation property. We assume that ${(V_\lambda)}_{\lambda\in\Lambda}$ is an open cover of $X$.
Then there exists a refinement $(W_\lambda)_{\lambda \in \Lambda}$ (some sets may be empty) with the same index set $\Lambda$, which is locally finite so that $W_\lambda\subset V_\lambda$ and  for every $\lambda\in \Lambda$ there exists an index $\tau(\lambda)\in T$ such that $\cl_X(W_\lambda)\subset  U_{\tau(\lambda)}$. 
Since $X$ is metrizable, we find open sets 
$Q_\lambda$ such  that
$$
 Q_\lambda\subset \cl_X(Q_\lambda)\subset W_\lambda
$$
and $(Q_\lambda)_{\lambda \in \Lambda}$ is a locally finite open cover of $X$. Using again  the metrizability  of $X$, we find continuous functions  $f_\lambda\colon  X\rightarrow [0,1]$ satisfying
$$
f\vert Q_\lambda\equiv 1\quad \text{and} \quad  \supp(f_\lambda)\subset W_\lambda.
 $$
 Let $\varepsilon=1/2$, $V_\lambda'=\phi_{\tau(\lambda)}(W_\lambda)$,   and $f_\lambda'=f_\lambda\circ\phi^{-1}_{\tau(\lambda)}$.   In view of the hypothesis (2), for the triple $(f_\lambda', V_\lambda', 1/2)$ there exists  a sc-smooth function $g_\lambda\colon  O_{\tau(\lambda)}\rightarrow [0,1]$ having support in $V_\lambda'$ and satisfying 
 $$
\abs{ f_\lambda'(x)-g_\lambda(x)}< 1/2\quad  \text{for all $x\in O_{\tau(\lambda)}$}.
$$
Going back to $X$ and extending $g_\lambda\circ \phi_{\tau (\lambda)}^{-1}$ onto $X$ by $0$ outside of $W_\lambda$,  we obtain the sc-smooth functions $\wh{g}_\lambda\colon  X\rightarrow [0,1]$ satisfying
$\wh{g}_{\lambda}(x)>0$ for $x\in Q_\lambda$. Then we define 
$\gamma_\lambda\colon  X\rightarrow [0,1]$ by
$$
\gamma_\lambda(x)=\frac{\wh{g}_\lambda(x)}{\sum_{\lambda\in\Lambda} \wh{g}_{\lambda}(x)}.
$$
The family of functions  $(\gamma_\lambda)_{\lambda \in \Lambda}$ is the desired sc-smooth partition of unity subordinate to the given open cover $(V_\lambda)$ of $X$. 

\end{proof}

An immediate consequence of Theorem \ref{partition_approximation} is the following result.

\begin{proposition}
Assume that the sc-Banach space admits smooth partitions of unity. Then a local model $(O, C, E)$ has the sc-smooth approximation property.\qed
\end{proposition}

The result below, due to Tor\'unczyk (see \cite{Fry}, Theorem 30),  gives a complete characterization of Banach spaces admitting $C^k$-partitions of unity.
This criterion reduces the question to a problem  in the geometry of Banach spaces. Though this criterion is not easy to apply it serves as one of the main tools in the investigation of the question, see \cite{Fry}. In order to formulate the theorem we need a definition.

\begin{definition}
If  $\Gamma$ is a set, we denote by $c_0(\Gamma)$ the Banach space of functions $f\colon  \Gamma\rightarrow {\mathbb R}$ having the property that for every $\varepsilon>0$
the number of $\gamma\in\Gamma$ with $\abs{f(\gamma)}>\varepsilon$ is finite.  The vector space operations are obvious and the norm is defined by
$$
\abs{f}_{c_0}=\text{max}_{\gamma\in\Gamma} \abs{f(\gamma)}.
$$
A  homeomorphic embedding $h\colon  E\rightarrow c_0(\Gamma)$ is {\bf coordinate-wise $C^k$}, \index{Coordinate-wise $C^k$-embedding}
if for every $\gamma\in\Gamma$ the map $E\rightarrow {\mathbb R}$,  $e\rightarrow (h(e))({\gamma})$ is $C^k$.
\qed
\end{definition}

If $\Gamma=\N$, then 
$c_0({\mathbb N})$ is the usual space $c_0$ of sequences converging to $0$.
\begin{theorem}[Tor\'unczyk's Theorem]
A Banach space $E$ admits a $C^k$-partition of unity if and only if there exists a set $\Gamma$ and a coordinate-wise $C^k$ homeomorphic embedding
of  $E$ into $c_0(\Gamma)$. \qed
\end{theorem}

An important class of Banach spaces are those which are are weakly compactly generated. They have good smoothness properties and  will provide
us with examples of sc-Banach spaces admitting sc-smooth partitions of unity.

\begin{definition}
A Banach space $E$ is called  {\bf weakly compactly generated} (WCG) if there exists a weakly compact set $K$ in $E$ such  that the closure of the span of $K$ is the whole space, 
$$
E=\cl_E(\text{span}(K)).
$$
\qed
\end{definition}

There are  two  useful examples of WCG Banach spaces.
\begin{proposition}
Let $E$ be a Banach space.
\begin{itemize}
\item[{\em (1)}]\ If $E$ is reflexive, then $E$ is WCG.
\item[{\em (2)}]\ If $E$ is separable, then $E$ is WCG.
\end{itemize}
\end{proposition}
\begin{proof}
In case that $E$ is reflexive it is known that the closed unit ball ${B}$ is  compact in the weak topology. Clearly,  $E=\text{span}(B)$.
If $E$ is separable, we  take a dense sequence ${(x_n)}_{n\geq  }$ in the unit ball and 
define $K=\{0\}\cup\{\frac{1}{n}x_n\, \vert \,  n\geq1\}$. Then $K$ is compact and,  in particular,  weakly compact. 
\qed
\end{proof}
The usefulness of WCG spaces lies in the following result from \cite{GTWZ}, see also \cite{Fry},  Theorem 31. 
\begin{theorem}[\cite{GTWZ}]\label{wcg_partition}
If the  WCG-space $E$ admits a $C^k$-bump function, then it also admits  $C^k$-partitions of unity.\qed
\end{theorem}


\begin{corollary}\label{C5417}
Let $(O,C,E)$ be a local model, where $E_0$ is a Hilbert space. Then $(O,C,E)$ has the sc-smooth approximation property.
\end{corollary}
\begin{proof}
A Hilbert space is reflexive and hence a WCG-space. 
We have already seen that a Hilbert space equipped with a sc-structure admits sc-smooth bump functions and conclude from 
Theorem \ref{wcg_partition} that it admits sc-smooth partition of unity, and consequently has the smooth approximation property, in view of Theorem 
\ref{partition_approximation}.
\end{proof}

\chapter{Orientations}\label{ORIENTXX}

In this chapter we introduce the notion of a linearization of a sc-Fredholm section and  discuss orientations and invariants associated to proper sc-Fredholm sections. We refer the reader to \cite{DK} for some of the basic ideas around determinants of linear Fredholm operators,
and to \cite{FH} for applications of the more classical ideas to problems arising in symplectic geometry, i.e. linear Cauchy-Riemann type operators. In polyfold theory, the central issue is that the occurring linear Fredholm operators, which are linearizations of nonlinear sections,
 do in general not depend as operators continuously on the points where the linearization was taken. On the other hand there is some weak continuity
property which allows to introduce orientation bundles. However, it is necessary to develop some new ideas.

\section{Linearizations of Sc-Fredholm Sections}

Let $P\colon Y\rightarrow X$ be a strong bundle over the  M-polyfold $X$ and $f$ a sc-smooth section of $P$.
If $x$ is a smooth point in $X$ and $f(x)=0$, there exists a well-defined {\bf linearization}\index{Linearization} 
$$f'(x)\colon T_xX\rightarrow Y_x$$
which is a sc-operator.  In order to recall the definition, we identify, generalizing 
a classical fact of vector  bundles, the tangent space $T_{0_x}Y$ at the element $0_x$ with the sc-Banach space $T_xX\oplus Y_x$ where $Y_x=P^{-1}(x)$ is the fiber over $x$. Denoting by $P_x\colon T_xX\oplus Y_x\to Y_x$ the sc-projection, the linearization of $f$ at the point $x$ is the following operator,
$$f'(x):=P_x\circ Tf(x)\colon T_xX\to Y_x.$$

As in the case of vector bundles there is, in general, no intrinsic notion of a  linearization of the section $f$ at the smooth point $x$ if $f(x)\neq 0$. However, dealing with a strong bundle we can profit from the additional structure. We simply take a local $\ssc^+$-section $s$ defined near $x$ and satisfying $s(x)=f(x)$, so that the linearization 
$$(f-s)'(x)\colon T_xX\to Y_x$$
is well-defined. To find such a $\ssc^+$-section we take a strong bundle chart around $x$. Denoting the sections in the local charts by the same letters, we let $R$ be the local strong bundle retraction associated with the local strong bundle. It satisfies $R (x, f(x))=f(x)$ and we define the desired section $s$ by $s(y)=R(y, f(x))$.  Since $f(x)$ is a smooth point, $s$ is a $\ssc^+$-section satisfying $s(x)=f(x)$ at the distinguished point $x$, as desired.
If $t$ is another $\ssc^+$-section satisfying $t(x)=f(x)=s(x)$, then 
$$
(f-s)'(x)=(f-t)'(x)+(t-s)'(x)
$$
and the linearization $(t-s)'(x)$ is a $\ssc^+$-operator. It follows from Proposition \ref{prop1.21}, that $(f-s)'(x)$ is a sc-Fredholm operator if and only if $(f-t)'(x)$ is a sc-Fredholm operator, in which case their Fredholm indices agree because a $\ssc^+$-operator is level  wise a compact operator.

Let now $f$ be a sc-Fredholm section of the bundle $P$ and $x$ a smooth point in $X$. Then there exists  a local $\ssc^+$-section $s$ satisfying $s(x)=f(x)$ and, moreover, $(f-s)'(x)$ is a sc-Fredholm operator.

\begin{definition}\index{D- Space of linearizations}\index{$\text{Lin}(f,x)$}\label{LINUXX}
If $f$ is a sc-Fredholm section $f$ of the strong bundle $P\colon Y\rightarrow X$ and $x$ a smooth point in $X$, then the {\bf space of linearizations} of $f$ at $x$ is the set of sc-operators from $T_xX$ to $Y_x$ defined as 
$$
\text{Lin}(f,x)=\{(f-s)'(x)+a\, \vert \,  \text{$a\colon T_xX\to Y_x$ is a $\ssc^+$-operator}\}.
$$
\qed
\end{definition}

The operators in $\text{Lin}(f,x)$\index{$\text{Lin}(f,x)$}
all differ by  linear $\ssc^+$-operators and are all sc-Fredholm operators having the same Fredholm index. Hence the space 
$\text{Lin}(f,x)$ is a convex subset of sc-Fredholm operators $T_xX\rightarrow Y_x$. 
Given an sc-Fredholm section of $P:Y\rightarrow X$ taking the linearizations at smooth points defines over $X_\infty$ some kind of `bundle' $\text{Lin}(f)$
where the fibers are the $\text{Lin}(f,x)$
$$
\text{Lin}(f)\rightarrow X_\infty.
$$
We shall use this structure later on for defining orientations.
For the moment it allows us to define
the {\bf  index of the sc-Fredholm germ $(f,x)$} \index{Index of sc-Fredholm germ}
by
$$
\text{ind}(f,x):=\dim \ker \bigl((f-s)'(x)\bigr)-\dim\bigl(Y_x/(\text{Im}(f-s)'(x))\bigr).\index{$\text{ind}(f,x)$}
$$
We shall show that this index is locally constant. The proof has to cope with the difficulty caused by the fact that, in general, the linearizations do not depend continuously as operators on the smooth point $x$.

Recall that a M-polyfold
is locally path connected, and that,  moreover,  any two smooth points in the same connected component can be connected by a $\ssc^+$-smooth path $\phi\colon [0,1]\rightarrow X$.
\begin{proposition}[Stability of $\text{ind}(f,x)$]\index{P- Stability of $\text{ind}(f,x)$}\label{sst}
Let $P\colon Y\rightarrow X$ be a strong bundle over the tame M-polyfold $X$ and $f$  a sc-Fredholm section. If 
$x_0$ and $x_1$ are smooth points in $X$ connected by a sc-smooth path $\phi\colon [0,1]\rightarrow X$, then
$$
\text{ind}(f,x_0)=\text{ind}(f,x_1).
$$
\end{proposition}
\begin{remark}\index{R- On tameness}
The tameness assumption is in all likelihood not needed. However, it allows to apply a trick.
\end{remark}
\begin{proof}
We shall show that the map $t\mapsto  \text{ind}(f,\phi(t))$ is locally constant. The difficulty is that the linearizations, even if picked sc-smoothly will,  in general,  not depend as operators continuously on $t$. On top of it we have possibly varying  dimensions of the spaces so that we need to change the filled version at every point. However, one can prove the result with a trick, which will also be used 
in dealing with orientation questions later on. We consider the tame M-polyfold $[0,1]\times X$ and consider the graph of the path $\phi$. We fix $t_0\in [0,1]$ and choose  a locally defined $\ssc^+$-section $s(t, x)$ satisfying $s(t,\phi(t))=f(\phi(t))$ for  $(t,x)\in [0,1]\times X$ near  $(t_0,\phi(t_0))$. (We do not need a sc-smooth partition of unity. If we had one  available then we could define such a section $s$ which satisfies  $s(t,\phi(t))=f(\phi(t))$ for all $t\in [0,1]$.) We choose  finitely many smooth points $e_1,\ldots ,e_m$
such  that the image of $(f-s(t_0,\cdot ))'(\phi(t_0))$ together with the $e_i$ span $Y_{\phi(t_0)}$.  Next we take a smooth finite-dimensional
linear subspace $L$ of $T_{\phi(t_0)}X$ which has a sc-complement in $T^R_{\phi(t_0)}X$. Then the image
of $L$ under $(f-s(t_0,\cdot ))'(\phi(t_0))$ is a smooth finite-dimensional subspace of $Y_{\phi(t_0)}$ of dimension $r$, say.
We choose  smooth vectors $p_1,\ldots ,p_r$ spanning this space. Next we take $m+r$ many $sc^+$-sections depending on $(t,x)$ (locally defined)
so that at $(t_0,\phi(t_0))$ they take the different values $e_1,\ldots ,e_m$ and $p_1,\ldots ,p_r$. Now the section
$$
F(\lambda,t,x)= f(x)-s(t,x)+\sum_{i=1}^{m+p} \lambda_i\cdot s_i(t,x)
$$
is defined near $(0,t_0,\phi(t_0))$ and takes values in $Y$. The linearization at $(0,t_0,\phi(t_0))$ with respect to the first and second variable
is surjective  and the kernel of the linearization has a complement contained in 
$$
T_{(0,t_0,\phi(t_0))}^R({\mathbb R}^{m+p}\oplus [0,1]\oplus X)={\mathbb R}^{m+p}\oplus T^R_{t_0}[0,1]\oplus T^R_{\phi(t_0)}X.
$$
Hence $\ker(F'(0,t_0,\phi(t_0)))$ is in good position to the boundary. Employing  the implicit function theorem
for the boundary case, we  obtain a solution manifold $S$ of $F=0$ containing $(0,t,\phi(t))$ for $t\in [0,1]$ close to $t_0$.
Moreover,  if  $(\lambda,t,x(t))\in S$, then  $\ker(F'(\lambda,t,\phi(t)))=T_{(\lambda,t,\phi(t))}S$ and
$F'(\lambda,t,\phi(t))$ is surjective. Hence 
$$
t\mapsto  \dim(T_{(0,t,\phi(t))}S)
$$
is locally constant for $t\in [0,1]$ near $t_0$. By construction, 
\begin{equation*}
\begin{split}
\text{ind}(f,\phi(t_0))+m+p&=\text{ind}(F,(0,t_0,\phi(t_0)))\\
&=\text{ind}(F,(0,t,\phi(t)))=\text{ind}(f,\phi(t))+m+p.
\end{split}
\end{equation*}
Therefore, $\ind (f,\phi (t_0))=\ind (f,\phi (t))$ for all $t$ near $t_0$.
\qed \end{proof}

\section{Linear Algebra and Conventions}\label{sect_conventions}

In he following  we  are concerned with the orientation  which is  crucial in our applications.
We follow  to a large extent the  appendix in \cite{HWZ5}. The ideas of the previous proof are also useful
in dealing with orientation questions. There we did not use sc-smooth partitions of unity.  {\bf In the following, however,}
{\bf we shall assume the existence of sc-smooth partitions of unity} to simplify the presentation, but the proofs
could be modified arguing as in the index stability theorem.

We begin with standard facts about determinants and wedge products.
Basically all the constructions are natural, but  usually depend on conventions, which have  to be stated apriori.
Since different authors  use different conventions,  their natural isomorphisms can be different.  To avoid these difficulties we state our conventions carefully. 
We also would like to point out that A. Zinger has written a paper dealing with these type of issues, \cite{Zinger}.
He also describes  some of the mistakes occurring in the literature as well as deviating conventions by different authors.
Since the algebraic treatment of SFT (one of the important applications of the current theory)
relies on the orientations
of the moduli spaces and the underlying conventions we give a comprehensive treatment of orientation questions.
  Using  the notation introduced by Zinger in  \cite{Zinger}, we define 
$$
\lambda(E):=\Lambda^{max} E\quad   \text{and}\quad   \lambda^\ast(E):=(\lambda(E))^\ast,\index{$\lambda(E)$}\index{$\lambda^\ast(E)$}
$$
where $(\lambda (E))^\ast$ is the dual of the vector space $\lambda(E)$.
A linear map $\Phi\colon E \rightarrow F$ between finite-dimensional vector spaces of the same dimension  induces the linear map
$$
\lambda(\Phi)\colon \lambda(E)\rightarrow\lambda (F), \index{$\lambda(\Phi)$}
$$
defined by $\lambda(\Phi)(a_1\wedge\ldots \wedge a_n):= \Phi(a_1)\wedge\ldots \wedge \Phi(a_n).
$
The  map $\lambda(\Phi)$ is nontrivial if and only if $\Phi$ is an isomorphism. The dual map $\Phi^\ast\colon F^\ast\to E^\ast$ of $\Phi\colon E \rightarrow F$ induces the map 
$$\lambda (\Phi^\ast)\colon \lambda (F^\ast )\to \lambda (E^\ast).$$
Moreover, we denote by 
$$\lambda^\ast  (\Phi)\colon \lambda^\ast  (F )\to \lambda^\ast  (E)$$
the dual of the map $\lambda(\Phi)\colon \lambda(E)\rightarrow\lambda(F)$.
The composition of  the two maps 
$$
E\xrightarrow{\Phi}F\xrightarrow{\Psi} G
$$
between vector spaces of the  same dimension satisfies 
$$
\lambda(\Psi\circ\Phi) =\lambda(\Psi)\circ \lambda(\Phi).
$$
There are different canonical isomorphisms
$$
\lambda ( E^\ast) \rightarrow \lambda^\ast( E)
$$
depending on different conventions. Our convention is the following.

\begin{definition}
If $E$ is a finite-dimensional real vector space and $E^\ast$ its dual,  the {\bf natural isomorphism}
$$
\iota\colon \lambda(E^\ast)\rightarrow\lambda^\ast(E)\index{$\iota\colon \lambda(E^\ast)\rightarrow\lambda^\ast(E)$}
$$
is defined by 
$$\iota( e_1^\ast\wedge\ldots \wedge e_n^\ast)(a_1\wedge \ldots \wedge a_n)=\det (e_i^\ast(a_j)),
$$
where $n=\dim(E)$. If $n=0$, we set 
$\lambda( E^\ast) =\lambda(\{0\}^\ast)=\R$.  
\qed
\end{definition}
From this definition we deduce  for a basis $e_1,\ldots ,e_n$ of $E$ and its dual basis $e_1^\ast,\ldots ,e_n^\ast$  the formula
$$
\iota(e^\ast_1\wedge\ldots \wedge e^\ast_n)=(e_1\wedge\ldots \wedge e_n)^\ast,
$$
where the dual vector $v^\ast$ of  a vector $v\neq 0$ in a one-dimensional vector space is determined by $v^\ast(v)=1$. 

Indeed, 
$$
(e_1\wedge\ldots \wedge e_n)^\ast (e_1\wedge\ldots \wedge e_n)=1= \det((e_i^\ast(e_j))=\iota(e_1^\ast\wedge\ldots \wedge e_n^\ast)(e_1\wedge\ldots \wedge e_n).
$$
The definition of $\iota$ is compatible with the previous definition of induced maps.

\begin{proposition}\index{P- Naturality of $\iota$}
If  $\Phi\colon E\rightarrow F$ is an  isomorphism between two finite-dimensional vector spaces and $\Phi^\ast\colon F^\ast\rightarrow E^\ast$ is its dual, then the following diagram is commutative, 
$$\begin{CD}
\lambda(F^\ast)@>\lambda(\Phi^\ast)>> \lambda(E^\ast)\\
@VV\iota V   @VV\iota V\\
\lambda^\ast(F)@>\lambda^\ast(\Phi)>> \lambda^\ast(E)
\end{CD}$$
\end{proposition}
\begin{proof}
Let $f_1,\ldots ,f_n$ be a basis of $F$ and $f_1^\ast,\ldots ,f^\ast_n$ its dual basis of $F^\ast$. Then we define the basis $e_1,\ldots ,e_n$  of $E$ by  $\Phi(e_i)=f_i$.
Its dual basis in $E^\ast$ is given by 
$e_1^\ast=f_1^\ast \circ \Phi,\ldots , 
e_n^\ast=f_n^\ast \circ \Phi$ and we compute,
\begin{equation*}
\begin{split}
&(\lambda^\ast  (\Phi)\circ\iota(f_1^\ast\wedge\ldots \wedge f^\ast_n))(e_1\wedge\ldots \wedge e_n)\\
&\quad= (\lambda (\Phi^\ast)((f_1\wedge\ldots \wedge f_n)^\ast))(e_1\wedge\ldots \wedge e_n)\\
&\quad = (f_1\wedge\ldots \wedge f_n)^\ast\circ \lambda(\Phi)(e_1\wedge\ldots \wedge e_n)\\
&\quad=(f_1\wedge\ldots \wedge f_n)^\ast(f_1\wedge\ldots \wedge f_n)\\
&\quad= 1.
\end{split}
\end{equation*}
Similarly,
\begin{equation*}
\begin{split}
&( \iota\circ\lambda(\Phi^\ast)(f_1^\ast\wedge\ldots \wedge f_n^\ast))(e_1\wedge\ldots \wedge e_n)\\
&\quad = (\iota( f_1^\ast\circ\Phi\wedge\ldots \wedge f_n^\ast\circ\Phi))(e_1\wedge\ldots \wedge e_n)\\
&\quad =(\iota(e_1^\ast\wedge\ldots \wedge e_n^\ast))(e_1\wedge\ldots \wedge e_n)\\
&\quad = {(e_1\wedge\ldots \wedge e_n)}^\ast(e_1\wedge\ldots \wedge e_n)\\
&\quad =1.
\end{split}
\end{equation*}
Hence
$
\lambda^\ast(\Phi)\circ \iota=\iota\circ \lambda(\Phi^\ast)
$
and the commutativity of the diagram is proved. 
\qed \end{proof}

Next we consider the  exact sequence ${\bm{E}}$ of finite-dimensional linear vector spaces, \index{${\bm{E}}$, exact sequence}
$$
{\bm{E}}:\quad  0\rightarrow A\xrightarrow{\alpha} B\xrightarrow{\beta} C\xrightarrow{\gamma} D\rightarrow 0.
$$
We recall that  the sequence is exact at $B$, for example, if  $\text{im}(\alpha)=\ker (\beta)$. 
  

We deal with the exact sequence as follows, and  take  a complement $Z\subset B$ of $\alpha (A)$ so that $B=\alpha (A)\oplus Z$, and 
a complement $V\subset C$ of $\beta (B)$ so that  $C=\beta (B)\oplus V$. Then the exact sequence ${\bm{E}}$ becomes 
$$
{\bm{E}}:\quad  0\rightarrow A\xrightarrow{\alpha} \alpha (A)\oplus Z \xrightarrow{\beta} \beta (B)\oplus V\xrightarrow{\gamma} D\rightarrow 0.$$
Here the first  nontrivial map is $a\mapsto (\alpha (a), 0)$, the second is $(b, z)\mapsto 
(\beta (z), 0)$, and the third is $(c,v)\mapsto \gamma (v)$. The maps $\alpha\colon A\to \alpha (A)$, $\beta \colon Z\to \beta (Z)$, and $\gamma \colon V\to D$ are isomorphisms. 
From the exact sequence ${\bm{E}}$ we are going  to construct several natural isomorphisms,  fixing again some  conventions.
The first natural isomorphism is the isomorphism
$$
\Phi_{\bm{E}}\colon \lambda(A)\otimes\lambda^\ast(D)\rightarrow \lambda(B)\otimes\lambda^\ast(C)\index{$\Phi_{\bm{E}}$}
$$
constructed as follows.
We abbreviate $n=\dim(A)$, $m=\dim(B)$, $k=\dim(C)$, and $l=\dim(D)$.  
$\Phi_{\bm{E}}$ maps $0$ to $0$. 
Next we take a nonzero vector
$$
h:=(a_1\wedge\ldots \wedge a_n)\otimes (d_1\wedge\ldots \wedge d_l)^\ast\in \lambda (A)\otimes \lambda^\ast (D),
$$
where  $a_1,\ldots ,a_n$ is a basis of  $A$ and $d_1,\ldots ,d_l$ is a basis for $D$.  Then we define the basis $b_1, \ldots ,b_n$ of $\alpha(A)$ by $b_i=\alpha (a_i)$ and the basis $c_1,\ldots ,c_l$ of $V$ by $\gamma (c_i)=d_i$, $i=1,\ldots ,l$.  Now we  choose a basis $b_1',\ldots ,b_{m-n}'$ of $Z$ and define the basis $c_1',\ldots ,c'_{m-n}$ of $\beta (B)\subset C$ by $c_i'=\beta (b_i')$.
Finally, we define $\Phi_{\bm{E}}(h)\in\lambda(B)\otimes\lambda^\ast(C)$ by 
\begin{equation}\label{eq_Phi_E}
\begin{split}
&\Phi_{\bm{E}}((a_1\wedge \ldots \wedge a_n)\otimes (d_1\wedge \ldots \wedge d_l)^\ast)\\
& = (\alpha(a_1)\wedge\ldots \wedge\alpha(a_n)\wedge b_1'\wedge\ldots \wedge b_{m-n}')\otimes (c_1\wedge\ldots \wedge c_l\wedge c_1'\wedge\ldots \wedge c_{m-n}')^\ast.
\end{split}
\end{equation}
The {\bf two conventions}  here are  that 
$b_1',\ldots ,b_{m-n}'$ are listed after the
$\alpha (a_1),\ldots ,\alpha (a_n)$ and
$c_1',\ldots ,c_{m-n}'$ are listed after $c_1,\ldots ,c_l$.
Apart from these two conventions how to list the vectors, the resulting definition does not depend on the choices involved.
\begin{lemma}\label{o6.6}\index{L- Well-definedness of $\Phi_{\bm{E}}$}
With the above two conventions, the definition of $\Phi_{\bm{E}}$ does not depend on the choices.
Hence $\Phi_{\bm{E}}$ is a natural isomorphism. 
\qed
\end{lemma}
The  proof is carried out  in Appendix \ref{oo6.6}.

Later we need to compare specific exact sequences. 
Given two exact sequences ${\bm{E}}$ and ${\bm{E}}'$ we assume that we have the commutative diagram
\begin{eqnarray}\label{EQ6240}
\begin{CD}
0 @>>> A@>\alpha>> B @>\beta >> C @>\gamma >> D @>>> 0\\
@.  @V \mathsf{A} VV  @V \mathsf{B} VV @V \mathsf{C}VV @V \mathsf{D} VV @.\\
0 @>>> A'@>\alpha'>> B' @>\beta' >> C' @>\gamma' >> D' @>>> 0,
\end{CD}
\end{eqnarray}
where the vertical maps are isomorphisms.
The horizontal arrows define linear isomorphisms $\lambda(A)\otimes\lambda^\ast(D)\rightarrow \lambda(B)\otimes\lambda^\ast(C)$ and 
$\lambda(A')\otimes\lambda^\ast(D')\rightarrow \lambda(B')\otimes\lambda^\ast(C')$, respectively.  The vertical arrows induce isomorphisms
which complete the obvious commutative diagram
\begin{eqnarray}\label{EQ62401}
\begin{CD}
\lambda(A)\otimes\lambda^\ast(D)@>>> \lambda(B)\otimes\lambda^\ast(C)\\
@VVV @VVV\\
\lambda(A')\otimes\lambda^\ast(D')@>>> \lambda(B')\otimes\lambda^\ast(C').
\end{CD}
\end{eqnarray}

Let $a_1,...,a_n$ be a a basis of $A$ and $d_1,...,d_\ell$ a basis of $D$. Then the top horizontal arrow is obtained by
extending $b_1=\alpha(a_1),..,b_n=\alpha(a_n)$ via a choice of $b_1',..,b_{m-n}'$ to a basis of $B$. We pick $c_1,..,c_\ell$
so that $\gamma(c_i)=d_i$. Then the $c_i$ are linearly independent. 
Then we consider
the linearly  independent vectors $c_1'=\beta(b_1'),...,c_{m-n}'=\beta(b_{m-n}')$. 
The horizontal arrow maps $(a_1\wedge..\wedge a_n)\otimes (d_1\wedge..\wedge d_\ell)^\ast$ to 
$ (\alpha(a_1)\wedge\ldots \wedge\alpha(a_n)\wedge b_1'\wedge\ldots \wedge b_{m-n}')\otimes (c_1\wedge\ldots \wedge c_l\wedge c_1'\wedge\ldots \wedge c_{m-n}')^\ast$.
Using that the vertical arrows are isomorphisms we can map
\begin{eqnarray*}
& a_1,..,a_n \rightarrow \mathsf{A}(a_1),..,\mathsf{A}(a_n)&\\
& \alpha(a_1),..,\alpha(a_n), b_1',..,b_{m-n}'  \rightarrow \mathsf{B}(\alpha(b_1)),..\mathsf{B}(\alpha_n),\mathsf{B}(b_1'),..,\mathsf{B}(b_{m-n}')&\\
&c_1,..,c_\ell,c_1',..,c_{m-n}'\rightarrow \mathsf{C}(c_1),..,\mathsf{C}(c_\ell),\mathsf{C}(c_1'),..\mathsf{C}(c_{m-n}')&\\
&d_1,..,d_\ell\rightarrow \mathsf{D}(d_1),..,\mathsf{D}(d_\ell).&
\end{eqnarray*}
We map for example  $h=(a_1\wedge..\wedge a_n)\otimes(d_1\wedge..\wedge d_\ell)^\ast$ to 
$ (\mathsf{A}(a_1)\wedge..\wedge\mathsf{A}(a_n))\otimes(\mathsf{D}(d_1)\wedge..\wedge\mathsf{D}(d_\ell))^\ast$ for the first vertical arrow.
By the previous lemma the horizontal arrows do not depend on the choices involved. From this it follows
also that the above definition does not depend on the choice. We leave additional arguments to the reader and summarize our findings as follows.
\begin{lemma}\label{LEM624x}
A commutative diagram (\ref{EQ6240}) with vertical maps being isomorphisms induces a natural (depending on our conventions) commutative diagram 
(\ref{EQ62401}).
\qed
\end{lemma}

Associated with the exact sequence ${\bm{E}}$ there exists also a second natural isomorphism\index{$\Psi_{\bm{E}}$}
\begin{equation}\label{secon_iso}
\Psi_{\bm{E}}\colon  \lambda(C)\otimes\lambda(A)\otimes\lambda^\ast(D)\rightarrow \lambda(B)
\end{equation}
constructed as follows.
We first map  
$(c_1\wedge\ldots \wedge c_k)\otimes (a_1\wedge\ldots \wedge a_n)\otimes(d_1\wedge\ldots \wedge d_l)^\ast$ into
$(c_1\wedge\ldots \wedge c_k)\otimes\Phi_{\bm{E}}((a_1\wedge\ldots \wedge a_n)\otimes(d_1\wedge\ldots \wedge d_l)^\ast)$ which belongs to  $\lambda(C)\otimes\lambda(B)\otimes\lambda^\ast(C)$
and then we compose this map  with the  isomorphism
$$
\bar{\iota}\colon \lambda(C)\otimes\lambda(B)\otimes\lambda^\ast(C)\to \lambda (B),
 $$
defined by $v\otimes b\otimes v^\ast\rightarrow b.$ Since $\Phi_{\bm{E}}$ is well-defined, so is $\Psi_{\bm{E}}$. For convenience we present a more explicit formula for  $ \Psi_{\bm{E}}$, using the notations of $\Phi_{\bm{E}}$.

\begin{proposition}\label{oger}\index{P- Well-definedness of $\Psi_{\bm{E}}$}
We fix the basis  $a_1,\ldots ,a_n$ of $A$ and 
the basis $d_1,\ldots ,d_l$ of  $D$ and abbreviate their wedge products by $a$ and by $d$. In the complement of $\beta(B)\subset C$ we have the basis  $c_1,\ldots ,c_l$ defined by $\gamma(c_i)=d_i$.
The vectors  $b_1,\ldots ,b_n\in B$ defined by
$b_i=\alpha(a_i)$ are a basis of  of $\alpha (A)\subset B$, and in the complement of $\alpha(A)$ in $B$ we choose the basis  a $b_1',\ldots ,b_{m-n}'$ and define the basis $c_1',\ldots ,c_{m-n}'$ of $\beta (B)\subset C$ by $c_i'=\beta (b_i')$. Abbreviating 
$c=c_1'\wedge\ldots \wedge c_l'\wedge c_1\wedge\ldots \wedge c_{m-n}$, we obtain the formula
$$
\Psi_{\bm{E}}(c\otimes a\otimes d^\ast)= b_1\wedge\ldots \wedge b_n\wedge b_1'\wedge\ldots \wedge b_{m-n}'\in \lambda (B).
$$
\end{proposition}
      
\begin{proof}
We already know from Lemma \ref{o6.6} that $\Psi_{\bm{E}}$ is well-defined since $\Phi_{\bm{E}}$ is well-defined. By construction, 
$$
\Phi_{\bm{E}}(a\otimes d^\ast) = (b_1\wedge\ldots \wedge b_n\wedge b_1'\wedge\ldots \wedge b_{m-n}')\otimes (c_1'\wedge\ldots \wedge c'_l\wedge c_1\wedge\ldots \wedge c_{m-n})^\ast.
$$
Here $c_i=\beta(b_i)$ and $\gamma(c_i') =d_i$. Then,  abbreviating the wedge product of the vectors $b_1,\ldots b_n$, 
$b_1', \ldots, b_{n-m}'$ by $b$, we obtain 
\begin{equation*}
\Psi_{\bm{E}}(c\otimes a\otimes d^\ast)
=\bar{\iota}(c\otimes b\otimes c^\ast)
=b.
\end{equation*}
\qed \end{proof}
                             
                          

\section{The Determinant of a Fredholm Operator}

We shall use some of the classical theory of determinants which can be found, for example,  in \cite{Zinger}.
A key fact is that many constructions in this area are naturally isomorphic, which too often is a source 
of sign errors in orientation questions.  It is important to note
that there are usually many natural isomorphisms. Specific ones   depend on a priori conventions, which very often differ from author to author
without being specified.
\begin{definition}\label{def_determinant_1}\index{D- Determinant}\index{$\text{det}(T)$}
The {\bf determinant} $\text{det}(T)$ of a bounded linear Fredholm operator $T\colon E\to F$ between real Banach spaces is the $1$-dimensional real vector space defined by
$$
\text{det}(T)= \lambda( \ker(T))\otimes \lambda^\ast( \text{coker(T)}).
$$
\qed
\end{definition}
An alternative definition used by some authors is  
$$
\text{det} (T)=\lambda(\ker(T))\otimes\lambda(\text{coker}(T)^\ast).
$$
The two definitions are naturally isomorphic given a convention how to identify $\lambda^\ast(A)$ and $\lambda(A^\ast)$. 
\begin{definition}
An {\bf orientation}\index{D- Orientation of $T$} of an Fredholm operator $T\colon E\rightarrow F$ is an orientation of the real line $\det(T)$.
\qed
\end{definition}

We begin this subsection by deriving exact sequences associated to Fredholm operators.  
\begin{definition}
Let $T\colon E\rightarrow F$ be a Fredholm operator between two Banach spaces. A {\bf good left-projection}\index{D- Good left-projection} for $T$ is a bounded projection $P\colon F\rightarrow F$ having the following two properties.
\begin{itemize}
\item[(1)]\ $\dim(F/R(P))<\infty$.
\item[(2)]\ $R(P\circ T)=R(P)$.
\end{itemize}
By $\Pi_T$\index{$\Pi_T$} we denote the collection of all good left-projections for $T$.
\qed
\end{definition}
In view of (1) the projection $P$ satisfies $\dim \text{coker} (P)=\dim \ker (P)<\infty$. Hence $P$ is a Fredholm operator of index $0$. Since $T$ is Fredholm, the composition $P\circ T$ is Fredholm and $\text{ind}(P\circ T)=\text{ind} (P)+\text{ind} (T)=\text{ind} (T).$
The definition of being a good left-projection has a certain amount of stability build in, as can be seen from the next lemma.
\begin{lemma}
Let  $T_0:E\rightarrow F$ be a Fredholm operator between Banach spaces and $P_0\in\Pi_{T_0}$. Then there exists $\varepsilon>0$ so that 
for every $T\in {\mathcal L}(E,F)$, $P\in {\mathcal L}(F)$ with $P\circ P= P$,  $\norm{T-T_0}<\varepsilon$ and $\norm{P-P_0}<\varepsilon$  it holds that $T$ is Fredholm,
$P\in \Pi_T$, and $\text{ind}(T)=\text{ind}(T_0)$.
\end{lemma}
\begin{proof}
Define $H=(I-P_0)F$ and denote by ${\mathcal P}(F)$ the spaces of bounded linear projections $P:F\rightarrow F$. We define the continuous map
$$
{\mathcal P}(F)\times {\mathcal L}(E,F)\rightarrow {\mathcal L}(E\times H,F):(P,T)\rightarrow \Phi_{(P,T)}
$$
by
$$
\Phi_{(P,T)}(x,h)= P\circ T(x) +(I-P)(h).
$$
We note that $\Phi_{(P_0,T_0)}$ is a surjectve Fredholm operator. Also the subset of Fredholm operators ${\mathcal F}(E\times H,F)$ is open in ${\mathcal L}(E\times H,F)$.
We find $\varepsilon>0$ so that the following holds for $\norm{P-P_0}<\varepsilon$, $\norm{T-T_0}<\varepsilon$.
\begin{itemize}
\item[(1)]\ $\Phi_{(P,T)}\in {\mathcal F}(E\times H,F)$.
\item[(2)]\ $\Phi_{(P,T)}$ is surjective.
\item[(3)]\ $(I-P):H\rightarrow (I-P)(F)$ is surjective.
\item[(4)]\ The dimension of $F/P(F)$ equals that of $F/H$. 
\end{itemize}
From this we immediately conclude that for $\norm{P-P_0}<\varepsilon$ and $\norm{T-T_0}<\varepsilon$
$T:E\rightarrow F$ is Fredholm, $P\circ T:E\rightarrow P(F)$ is surjective, and $R(P)$ has finite co-dimension.
This means that $P\in \Pi_T$.
\qed \end{proof}
\begin{remark}\index{R- On $\det(T)$}
The previous lemma has some important consequence which will be used in later constructions.
If $T_0$ is Fredholm and $P\in \Pi_{T_0}$ then a $T$ near $T_0$ is Fredholm and $P\in \Pi_{T_0}$.
It is an easy exercise that $\ker(PT)$ varies continuously as $T$ (near $T_0$) varies and we obtain
an honest local line bundle. This trivially holds for the cokernel $F/R(PT)=F/R(P)$.
As we shall see there is a natural isomorphism between $\det(T)$ and $\det(PT)$ and this fact will be crucial, since it allows
to define a structure of a linear bundle on the union of all $\det(T)$ as $T$ varies.
\qed
\end{remark}
The set $\Pi_T$ of projections  possesses a {\bf partial ordering}\index{Partial ordering on $\Pi_T$} $\leq $ defined by
$$
\text{$P\leq Q$ \quad if and only if \quad $P=PQ=QP$}.
$$
It is in general not true that for given $P,Q\in\Pi_T$ there exists $L\in\Pi_T$ with $L\leq P$ and $L\leq Q$.
\begin{remark}\label{REM635}
In the context of Hilbert spaces one usually takes $\Pi_T^{orth}$, which is the subset of $\Pi_T$ consisting of orthogonal projections.
In this case there exists for $P,Q\in \Pi_T^{orth}$ a $L\in \Pi_T^{orth}$ with $L\leq Q$ and $L\leq P$. This fact simplifies later constructions in the Hilbert case space. 
This particular case, which is classical and well-known,  was discussed in \cite{HWZ5}.
\qed
\end{remark}
Given two projections $P,Q:F\rightarrow F$ having images with finite co-dimensions, the space $R(P)\cap R(Q)$ has finite co-dimension as well.
\begin{lemma}\label{lem634}
Let $Q:F\rightarrow F$ and $P:F\rightarrow F$ be  two continuous projections which have images of finite co-dimension and let $L:F\rightarrow F$ be a continuous projection  satisfying $R(L)=R(P)\cap R(Q)$. Then the bounded linear operators $LP$ and $LQ$ are projections onto $R(L)$ satisfying
 $LQ\leq Q$ and $LP\leq P$.
 \end{lemma}
 \begin{proof}
For $y\in F$ we compute since $R(L)\subset R(P)$  
$$
LPLP(y)= LLP(y)=LP(y)
$$
and similarly for $LQ$. This shows that $LP$ and $LQ$ are projections. Since $R(L)=R(P)\cap R(Q)$ it follows that $R(LQ)=R(LP)=R(L)$.
Finally we compute that
$$
P(LP) = LP = L(PP)=(LP)P
$$
which precisely means that $LP\leq P$. Similarly it follows that $LQ\leq Q$.
  \qed \end{proof}

\begin{definition}
Assume  $H\subset F$ is a closed linear subspace with finite co-dimension. 
We denote by ${\mathcal P}_H$\index{${\mathcal P}_H$}
the subset of the space of continuous linear operators ${\mathcal L}(F)$, consisting of projections onto $H$.
\qed
\end{definition}
Continuous projections $P:F\rightarrow F$
with image $H$ are in 1-1 correspondence with closed subspaces $K$ of $F$ satisfying
$$
\dim(K) = \dim(F/H)\ \ \text{and}\ \ \ K\cap H=\{0\}.
$$
Given an ordered pair $(K,K')$ of such linear subspaces, there is an associated linear map $\Gamma_{(K,K')}:K\rightarrow H$ uniquely determined by
\begin{eqnarray}\label{eq635}
K'=\{k+\Gamma_{(K,K')}(k)\ |\ k\in K\}.
\end{eqnarray}
We also note the following trivial result.
\begin{lemma}
Let $P,Q\in {\mathcal P}_H$.  Then 
$$
Q(x) = P(x) + Q((I-P)(x))
$$
for all $x\in F$. 
\end{lemma}
\begin{proof}
By assumption $R(P)=R(Q)=H$. We obtain two   topological direct sum decompositions 
$F=H\oplus (I-P)(F)$ and $F=H\oplus(I-Q)(F)$.  Given $x\in F$ we can write it as 
$$
x= P(x) + (I-P)x = Q(x) +(I-Q)x
$$
and applying $Q$ we find that 
$$
Q(x) =QP(x) +Q(I-P)(x)= P(x)+Q(I-P)(x).
$$
\qed \end{proof}

At the heart of the following constructions are suitable exact sequences.
Associated with  the projection $P\in\Pi_T$ there is the exact sequence\index{${\bm{E}}_{(T,P)}$, exact sequence}
$$
{\bm{E}}_{(T,P)}:\quad 0\rightarrow\ker(T)\xrightarrow{j_T^P} \ker(PT) \xrightarrow{\Phi^P_T} F/R(P)\xrightarrow{\pi_T^P}\text{coker}(T)\rightarrow 0,
$$
where $j_T^P$\index{${\bm{E}}_{(T,P)}$} is the inclusion map, $\pi_T^P$ is  defined by 
$$
\pi_T^P(f+R(P))= (I-P)f+R(T), \ \ f\in F,
$$
  and 
$$
\Phi^P_T(x)=T(x)+R(P),\ \  x\in \ker (PT).
$$
\begin{proposition}\index{P- Exactness of ${\bm{E}}_{(T,P)}$}
The sequence ${\bm{E}}_{(T,P)}$ is exact.
\end{proposition}

\begin{proof}
The inclusion $j_T^P$ is injective and $\Phi^P_T\circ j_T^P=0$. If $\Phi^P_T(x)=0$,  then $T(x)\in R(P)$ implying $(I-P)T(x)=0$.
Since $x\in \ker(PT)$ we conclude $T(x)=0$. This proves exactness at $\ker(PT)$. Assuming  $x\in \ker(PT)$,  hence 
$PT(x)=0$, we compute
\begin{equation*}
\begin{split}
\pi_T^P\circ \Phi^P_T(x)&=\pi_T^P(T(x)+R(P))\\
&=\pi_T^P(T(x)+R(PT))\\
&=(I-P)T(x)+R(T)\\
&=T(x)+R(T)\\
&=R(T),
\end{split}
\end{equation*}
i.e. the composition $\pi_T^P\circ \Phi^P_T$ vanishes. If $\pi_T^P(y+R(P))=0$, hence $(I-P)y\in R(T)$, there exists 
$x\in E$ solving  $T(x)=(I-P)y$. Consequently, 
$$
\Phi_T^P(x) = T(x) + R(P) = (I-P)y+R(P)= y+R(P),
$$
proving the  exactness at $F/R(P)$. Finally we show the surjectivity of the map $\pi_T^P$.
Given $f+R(T)\in \text{coker}(T)$, we  choose  $x\in E$ satisfying  $PT(x)=Pf$ and  compute,
$$
\pi_T^P(f-T(x)+R(P))= (I-P)(f-T(x)) +R(T)= f +R(T),
$$
which finishes the proof of the exactness.
\qed \end{proof}

We discuss the ramifications of the previous discussion.
Starting with the  Fredholm operator $T\colon E\rightarrow F$ and the projection $P\in\Pi_T$ we obtain from  the exact sequence
$$
{\bm{E}}_{(T,P)}:\quad 0\rightarrow\ker(T)\xrightarrow{j_T^P} \ker(PT) \xrightarrow{\Phi^P_T} F/R(P)\xrightarrow{\pi_T^P}\text{coker}(T)\rightarrow 0,
$$
recalling $R(P)=R(PT)$, \index{${\bm{E}}_{(T,P)}$, exact sequence} and employing Lemma \ref{o6.6}
the natural isomorphism
$$
\Phi_{{\bm{E}}_{(T,P)}}\colon\lambda(\ker(T))\otimes\lambda^\ast(\text{coker}(T))\rightarrow \lambda(\ker(PT))\otimes \lambda^\ast(\text{coker}(PT)).
$$
By  definition of the determinant,  this means that
$$
\Phi_{{\bm{E}}_{(T,P)}}\colon  \det(T)\rightarrow \det(PT)
$$
is an isomorphism.
We rename this isomorphism for the further discussion and,  setting $\Phi_{{\bm{E}}_{(T,P)}}=\gamma_T^P$, we have the isomorphism
\begin{eqnarray}\label{EQ6311}
\gamma_T^P\colon \det(T)\rightarrow\det(PT).\index{$\gamma^P_T$}
\end{eqnarray}
It is important to note that the righthand side $\det(PT)$ has a local continuity property when $T$ varies. 
For further constructions it will be important to study (\ref{EQ6311}) if we change $P$ to a $Q$ satisfying $Q\leq P$
and ultimately to $Q$ not order related to $P$.

}

{
We consider the Fredholm operator $T\colon E\rightarrow F$ and the projections $P,Q\in \Pi_T$ satisfying $P\leq Q$, so that  $P=PQ=QP$. Then, as we have already seen, the composition  $QT$ is Fredholm. Moreover, 
$P\in \Pi_{QT}$. Therefore, the associated exact sequences ${\bm{E}}_{(T,P)}$, ${\bm{E}}_{(T,Q)}$, and
${\bm{E}}_{(QT,P)}$ produce the following isomorphisms,
\begin{align*}
\gamma^P_T\colon&\det(T)\rightarrow\det(PT),\\ \gamma^Q_T\colon&\det(T)\rightarrow\det(QT),\\ \gamma^P_{QT}\colon&\det(QT)\rightarrow\det(PT).
\end{align*}
\index{$\gamma^P_T:\det(T)\rightarrow\det(PT)$}\index{$\gamma^Q_T:\det(T)\rightarrow\det(QT)$}
\index{$\gamma^P_{QT}:\det(QT)\rightarrow\det(PT)$}
The crucial observation is the following result, whose proof is postponed to Appendix \ref{oooo6.6}.
\begin{proposition}\label{ooo6.6}
If  $T\colon E\rightarrow F$ is a Fredholm operator and $Q,P\in\Pi_T$ satisfy $P\leq Q$, 
then $QT\colon E\rightarrow F$ is a Fredholm operator and $P\in\Pi_{QT}$. Moreover, 
\begin{eqnarray*}
\gamma^P_{QT}\circ \gamma^Q_T=\gamma^P_T.
\end{eqnarray*}
\qed
\end{proposition}
From Proposition \ref{ooo6.6} we conclude for the three 
projections 
$P, R,S\in \Pi_T$ satisfying $P\leq R$ and $P\leq S$, the relation 
$\gamma^P_{TR}\circ \gamma_T^R=\gamma_T^P=\gamma_{ST}^P\circ \gamma_T^S$, so that the 
diagram
\begin{equation*}
\begin{CD}
\det(T)@> \gamma_T^R >> \det(RT)\\
@VV\gamma_T^S V   @VV\gamma_{RT}^P V\\
\det(ST) @>\gamma_{ST}^P>>  \det(PT)
\end{CD}
\end{equation*}
is commutative, implying the following corollary.
\begin{corollary}\label{KPD}
With the assumptions of Proposition \ref{ooo6.6} the map 
$$
\det(ST)\rightarrow \det(RT), \quad  h\mapsto {(\gamma_{RT}^P)}^{-1}\circ\gamma_{ST}^P(h)
$$
is independent of the choice of $P$ as long as $P\leq S$ and $P\leq R$. Moreover,  this map is equal to the isomorphism 
\begin{eqnarray}\label{EQ6313}
\det(ST)\rightarrow \det(RT),\quad h\mapsto \gamma^R_T\circ {(\gamma_T^S)}^{-1}(h).
\end{eqnarray}
\qed
\end{corollary}
The last statement (\ref{EQ6313})  may be viewed as some kind of transition map for
the isomorphisms  $\det(T)\rightarrow \det(PT)$ and $\det(T)\rightarrow \det(QT)$ provided 
 $P,Q\in \pi_T$ and there exists $S\leq P$ and $S\leq Q$. The important fact is that, as already pointed out,
 $\det(PT)$ and $\det(QT)$ change locally continuously in $T$.
Unfortunately it is  not always possible to find for given $P,Q\in \Pi_T$ a projection $S\in \Pi_T$ with $S\leq P$ and $S\leq Q$.  Therefore we 
need a result comparing the constructions for $P$ and $Q$, both in $\Pi_T$, respectively.  
Recalling Remark \ref{REM635} we would like to point out, that in the Hilbert space  setting, taking only
orthogonal projections the previous constructions are enough to define the general determinant bundle.

For the Banach space case we need to understand the relationship when we use two non-compatible 
projections. In view of Lemma \ref{lem634} given $Q,P\in \Pi_T$ we can find $L\in\Pi_T$
such that $LP, LQ\in \Pi_T\cap {\mathcal P}_{R(Q)\cap R(P)}$, $LP\leq P$ and $LQ\leq Q$. 
This gives us the natural maps
\begin{eqnarray}\label{HKL99}
\det(T)\xrightarrow{\gamma^P_T}\det(PT) \ \ \text{and}\ \ \det(T)\xrightarrow{\gamma^Q_T}\det(QT)
\end{eqnarray}
and we need in principle to understand the continuity property of ${(\gamma^Q_T)}\circ{(\gamma^P_T)}^{-1}$.
In order to simplify this task, one uses the fact that 
the natural maps 
$$
\det(PT)\rightarrow\det(LPT)\ \ \text{and}\ \ \det(QT)\rightarrow \det(LQT)
$$
vary (near $T$) continuously. In view of Proposition \ref{ooo6.6},  it suffices therefore to establish
instead of the (local) continuous dependence on $T$ of the transition map associated to  (\ref{HKL99}),  the continuous dependence
of the  transition maps associated to
$$
\det(T)\xrightarrow{\gamma^P_T}\det(LPT) \ \ \text{and}\ \ \det(T)\xrightarrow{\gamma^Q_T}\det(LQT).
$$
The advantage here is that $LQ,LP\in \Pi_{T}\cap {\mathcal P}_{R(Q)\cap R(P)}$. 
Consider  the following diagram which has exact horizontal rows
$$
\begin{CD}
0@>>> \ker(T)@> j^{LP}_T >> \ker(LPT) @> \Phi^{LP}_T >> F/H @> \pi^{LP}_T >> \text{coker}(T) @>>> 0\\
@.         @|                        @V \bm{?}VV                    @V \bm{??}     VV           @| @.\\
0@>>> \ker(T)@> j^{LQ}_T >> \ker(LQT) @> \Phi^{LQ}_T >> F/H @> \pi^{LQ}_T >> \text{coker}(T) @>>> 0.
\end{CD}
$$
We would like to show that we can fill in for $\bm{?}$ and $\bm{??}$ so that the diagram is commutative 
and so that the `fill-in' maps introduce (locally when $T$ varies) continuous maps between the kernel bundles 
associated to $\ker(LPT)$ and $\ker(LQT)$ and similarly for the constant bundles associated to $F/H$. 
Having the diagram filled in by the maps $\bm{?}$ and $\bm{??}$ will give us in view of Lemma \ref{LEM624x}
$$
\begin{CD}
\det(T)@> \gamma^{LP}_T >> \det(LPT)\\
@|    @V\gamma_{{\bm{?,??}}} VV\\
\det(T) @>\gamma^{LQ}_T >>\det(LQT).
\end{CD}
$$
Of course, since $ \gamma_{\bm{?,??}} = \gamma^{LQ}_T\circ {(\gamma^{LP}_{T})}^{-1}$, the proposed procedure
would show, that the transition map is continuously depending when $T$ is locally changed.
The  desired result follows from the next proposition, where we change somewhat the notation to simplify the representation.
It is concerned with the question of filling in the previous (big) diagram.
\begin{proposition}\label{PROP6310}
Let $T:E\rightarrow F$ be a linear Fredholm operator and consider  for $P,Q\in \Pi_T\cap {\mathcal P}_H$ the associated exact sequences
${\bm{E}}_{(T,P)}$ and ${\bm{E}}_{(T,Q)}$.  Then there is a construction for linear maps  $A_T$ and $B_T$,  which are isomorphisms and make the following diagram
commutative
$$
\begin{CD}
0@>>> \ker(T)@> j^P_T >> \ker(PT) @> \Phi^P_T >> F/H @> \pi^P_T >> \text{coker}(T) @>>> 0\\
@.         @|                        @V A_T    VV                    @V B_T     VV           @| @.\\
0@>>> \ker(T)@> j^Q_T >> \ker(QT) @> \Phi^Q_T >> F/H @> \pi^Q_T >> \text{coker}(T) @>>> 0.
\end{CD}
$$
Varying $T$ to a $T'$ near $T$ so that $P,Q\in \Pi_{T'}$ the constructions of $A_T$ and $B_T$ depend continuously on $T$.
\end{proposition}
\begin{proof}
Since $\text{ind}(PT)=\text{ind}(QT)=\text{ind}(T)$ and $R(PT)=R(QT)= H$ it follows that $\dim(\ker(PT))=\dim(\ker(QT))$.
We can take a closed linear subspace $X$ of $E$ which at the same time is a topological linear complement of $\ker(QT)$ and $\ker(PT)$.
Hence
$$
E=X\oplus\ker(QT)=X\oplus \ker(PT).
$$
For  given $e\in \ker(PT)$
denote by $x(e)\in X$ the unique solution in $X$ of 
\begin{eqnarray}\label{EQ234}
QT(x)= QT(e).
\end{eqnarray}
Then $QT(e-x(e))=0$ and we define
$$
A_T\colon\ker(PT)\rightarrow \ker(QT)\colon A_T(e)=e-x(e).
$$
Assume that $e\in \ker(PT)\cap \ker(QT)$ (This intersection  contains $\ker(T)$.). 
Then $QT(x(e))=QT(e)=PT(e)=0$ implying that $x(e)=0$ and therefore
$$
A_T|(\ker(PT)\cap \ker(QT))=Id.
$$
The map  $A_T$ is an isomorphism since $A_T(e)=0$ implies $e\in X$, which together with $e\in \ker(PT)$
implies $e\in X\cap\ker(PT)=\{0\}$.
We also note that the equation (\ref{EQ234}) is uniquely solvable for $T'$ near $T$, in addition for such closeby $T'$ 
we also have $Q,P\in \Pi_{T'}$, and the solution depends continuously on $T'$.

In order  to define the map $B_T:F/H\rightarrow F/H$. 
We consider the bounded linear map
$$
\Phi\colon X\times (I-Q)F\rightarrow F:(x,\wt{f})\rightarrow T(x)+\wt{f}.
$$
This map is an isomorphism. Indeed, if $T(x)+\wt{f}=0$ we find that $QT(x)=0$ implying $x=0$, which implies $\wt{f}=0$.
Given $f\in F$ we find a unique $x\in X$ solving $QT(x)=Qf$. We then pick $\wt{f}:= (I-Q)f-(I-Q)T(x)$ and compute
$$
T(x)+\wt{f} = QT(x) +(I-Q)T(x) +(I-Q)f-(I-Q)T(x)=Qf+(I-Q)f =f.
$$
Now we solve for $f$ with $(I-P)f=f$
$$
T(x)+\wt{f}=f.
$$
The solution is denoted by $(x(f),\wt{f}(f))$ and we  define
$$
B_T\colon F/H\rightarrow F/H\colon B_T(f+H) = \wt{f}+H.
$$
Next we show that
$$
B_T\circ \Phi^P_T= \Phi^Q_T\circ A_T.
$$
For $e\in \ker(PT)$ we compute, using that $(I-P)T(e)=T(e)$ with $T(x(e))+\wt{f}(e)=T(e)$ and observing that $QT(x(e))=QT(e)$ (the defining equation for $A_T$)
\begin{eqnarray*}
B_T\circ \Phi^P_T(e) &=&B_T( T(e)+H)\\
&=& \wt{f}(e)+H\\
&=& T(e-x(e)) + H\\
&=& T(A_T(e))+H\\
&=&\Phi^Q_T(A_T(e)).
\end{eqnarray*}
Further we calculate where $(I-P)f=f$ and $T(x) +\wt{f} =f$
\begin{eqnarray*}
\pi^Q_T\circ B_T(f+H) &=& \pi^Q_T(\wt{f}+H)\\
&=& \wt{f} +R(T)\\
&=& \wt{f}+T(x)+R(T)\\
&=& f+R(T)\\
&=& \pi^P_T(f+H).
\end{eqnarray*}
In order to calculate $B_T((I-P)T(e)+H) $ for $e\in \ker(PT)$, we need to solve for $x\in X$
$$
T(x)+\wt{f}= T(e),
$$
implying $QT(x) =QT(e)$ and consequently  $x=x(e)$. Hence $\wt{f}= T(e-x(e))=T(A_T(e))$ from which we  infer that
$$
\wt{f}+H = T(A_T(e))+H =\Phi^Q_T\circ A_T(e).
$$
We  calculate
\begin{eqnarray*}
\pi^Q_T\circ B(q)&=& \pi^Q_T(\wt{f}(q)+H)\\
&=&(I-Q)\wt{f}(q) + R(T)\\
&=& \wt{f}(q)+R(T)\\
&=& f-T(x(f))+R(T)\\
&=& f+R(T)\\
&=&\pi^P_T(f+H)\\
&=&\pi^P_T(q).
\end{eqnarray*}
Let $e\in \ker(PT)$ and solve $T(x)+\wt{f} = (I-P)T(e)$, so that $B^T_X((I-P)T(e) +H) =\wt{f}+H$.
\begin{eqnarray*}
&& B_X^T\circ \Phi_T^P(e)\\
&=& B_X^T (T(e)+H)\\
&=& \wt{f}+H.
\end{eqnarray*}
Moreover,
\begin{eqnarray*}
&& \Phi_T^Q\circ A(e)\\
&=&\Phi_T^Q(e-x(e))\\
&=& T(e-x(e))+H.
\end{eqnarray*}
Since the diagram has been shown to be commutative, and $A$ is an isomorphism, the five-lemma implies that $B$ is an isomorphism.

By construction $A_T$ is obtained for $e\in\ker(PT)$ by solving $QT(x)=QT(e)$ which is the same as solving
$T(x)+\wt{f} =T(e)$, and $A_T$ maps $e$ to $e-x$.  The map $B_T$ is obtained by solving for $f$ with $(I-P)f=f$
the equation $T(x)+\wt{f}=f$ and mapping $f+H$ to $\wt{f}+H$.  The important fact is that
$$
X\times (I-Q)F\rightarrow F\colon (x,\wt{f})\rightarrow T(x)+\wt{f}
$$
is a linear isomorphism.  If we keep $X$ and perturb $T$ to $T'$  the equation is still uniquely solvable 
and $\ker(PT')\oplus X=\ker(QT')\oplus X=E$. This shows that the maps 
$B_{T'}$ and $A_{T'}$ for small perturbation are continuously depending on $T'$.
\qed \end{proof}
We already discussed the following informally.
In view of Lemma \ref{lem634}
we can take a projection $L$ onto $R(P)\cap R(Q)$ and  the previous discussion applies to $LP\leq P$ and $LQ\leq Q$.
Then we only need to study the relationship between the constructions for  $LP$ and $LQ$. Both projections belong to 
${\mathcal P}_H$, where $H=R(L)=R(LP)=R(LQ)$. In view of Proposition \ref{PROP6310} we can relate the constructions
for $LP$ and $LQ$.
Hence, by  the previous discussions we have natural  isomorphisms 
$$
\det(T)\rightarrow\det(PT)\rightarrow\det(LPT)\ \ \text{and} \ \ \det(T)\rightarrow\det(QT)\rightarrow\det(LQT)
$$
and from Proposition \ref{PROP6310} 
$$
\begin{CD}
0@>>> \ker(T)@> j^{LP}_T >> \ker(LPT) @> \Phi^{LP}_T >> F/H @> \pi^{LP}_T >> \text{coker}(T) @>>> 0\\
@.         @|                        @V A_T    VV                    @V B_T    VV           @| @.\\
0@>>> \ker(T)@> j^{LQ}_T >> \ker(LQT) @> \Phi^{LQ}_T >> F/H @> \pi^{LQ}_T >> \text{coker}(T) @>>> 0
\end{CD}
$$
from which 
we obtain the commutative  diagram
$$
\begin{CD}
\det(T) @>>> \det(LPT)\\
@| @V\gamma_T VV\\
\det(T) @>>> \det(LQT)
\end{CD}
$$
Here $\gamma_T$ is obtained from $A_T$ and $B_T$ and therefore shows continuous dependence on $T$.

Postponing  the construction of determinant bundles there is an important stabilization construction,  which is the linearized version
of a construction occurring in perturbation theory,  and which has to be understood from the point of view of orientations of determinants. 
We start with a Fredholm operator $T\colon E\rightarrow F$ and assume that $\phi\colon {\mathbb R}^n\rightarrow F$ is a linear map
such  that the map $T_\phi\colon  E\oplus {\mathbb R}^n\rightarrow F$, defined by 
$$
T_\phi(e,r)=T(e)+\phi(r),\index{$T_{\phi}$}
$$ 
is surjective. Writing ${\mathbb R}^n$ after 
$E$ is convenient and goes hand in hand with the conventions in  the definition of the isomorphism $\Phi_{\bm{E}}$ associated with the  exact sequence ${\bm{E}}$.  We would like to introduce 
a convention relating the orientations of  $\det(T)$ and $\det(T_\phi)$,  knowing that  ${\mathbb R}^n$ possesses  the   preferred orientation as 
direct sum of $n$-many copies of ${\mathbb R}$,  each of which is  oriented by $[1]$.
We introduce  the exact sequence
$$
{\bm{E}}:\quad  0\rightarrow \ker(T)\xrightarrow{j} \ker(T_\phi)\xrightarrow{\pi} {\mathbb R}^n\xrightarrow{[\phi]} F/R(T)\rightarrow 0,
$$
in which  $j(e)=(e,0)$,  $\pi(e,r)= r$, and  $ [\phi](r)=\phi(r)+R(T)$. 
\begin{lemma}
The sequence ${\bm{E}}$ is exact.
\end{lemma}
\begin{proof}
Clearly,  $\pi\circ j=0$ and $j$ is injective. If $\pi(e,r)=0$, then $T(e)=0$ and $j(e)=(0,e)=(e,r)$. If $(e,r)\in\ker(T_\phi)$, then $T(e)+\phi(r)=0$
implying that $\phi(r)\in R(T)$. Hence $[\phi(r)]=R(T)$. If $[\phi](r)=0$,  we have $\phi(r)=T(-e)$ for some $e\in E$ and hence $(e,r)\in\ker(T_\phi)$.
Moreover,  $\pi(e,r)=r$. The last map is surjective since,   by  assumption,  $T_\phi$ is surjective. 
\qed \end{proof}
From  the  exact sequence ${\bm{E}}$ we deduce the previously constructed natural isomorphism
$$
\Phi_{\bm{E}}\colon \det(T)\rightarrow \lambda(\ker(T_\phi))\otimes \lambda^\ast({\mathbb R}^n).
$$
In the applications the operator $T$ is oriented and the auxiliary constructions to obtain transversality  yield,  on the linearized level,  operators of the type $T_\phi$.
The above isomorphism can be used to relate the orientations of  $\det(T)$ and $\det(T_\phi)$. However, this requires  an additional convention. For this it suffices 
to fix an isomorphism $\lambda^\ast({\mathbb R}^n)\rightarrow {\mathbb R}^\ast =\lambda^\ast(\{0\})$. This gives rise to an isomorphism
$$
\lambda(\ker(T_\phi))\otimes \lambda^\ast({\mathbb R}^n)\rightarrow \lambda(\ker(T_\phi))\otimes {\mathbb R}^\ast=\det(T_\phi),
$$
which then yields an isomorphism $\det(T)\rightarrow \det(T_\phi)$. The isomorphism we choose is defined by 
$$
\psi\colon \lambda^\ast({\mathbb R}^n)\rightarrow \lambda^\ast(\{0\})={\mathbb R}^\ast,\quad  (e_1\wedge\ldots \wedge e_n)^\ast\mapsto 1^\ast,
$$
where $e_1,\ldots ,e_n$ is the standard basis of ${\mathbb R}^n$. It then follows for every basis $c_1,\ldots ,c_n$ of  ${\mathbb R}^n$ that 
\begin{equation}\label{iso--}
\psi \bigl((c_1\wedge\ldots \wedge c_n)^\ast\bigr)=\frac{1}{\det([c_1,\ldots ,c_n])}1^\ast,
\end{equation}
where $\det$ denotes the determinant of a $n\times n$ matrix of the column vectors.

Summing up the previous discussion and using the definition of the isomorphism $\Phi_{\bm{E}}$ associated to the exact sequence ${\bm{E}}$ we can summarize the findings 
in  the following proposition.
\begin{proposition}
With the chosen isomorphism $\phi\colon \lambda^\ast({\mathbb R}^n)\rightarrow\lambda^\ast(\{0\})={\mathbb R}^\ast$, $(e_1\wedge\ldots \wedge e_n)^\ast\rightarrow 1^\ast$, we define the  isomorphism \index{$\iota_\phi\colon \det(T)\rightarrow \det(T_\phi)$}
$$
\iota_\phi\colon \det(T)\rightarrow \det(T_\phi)
$$
by the following formula. We choose  a basis $a_1,\ldots ,a_k$ of $\ker(T)$ and a basis $\bar{d}_1=\phi(\bar{c}_1)+R(T),\ldots ,\bar{d}_m=\phi(\bar{c}_m)+R(T)$ of  $F/R(T)$.
Let $h=(a_1\wedge\ldots \wedge a_k)\otimes (\bar{d}_1\wedge\ldots \wedge \bar{d}_m)^\ast$.  
Extend $(a_1,0),\ldots ,(a_k,0)$ to a basis of  $\ker(T_\phi)$ by adding
$\bar{b}_1=(b_1,r_1),\ldots ,\bar{b}_l=(b_l,r_l)$,  and notice that $r_1=\pi(b_1,r_1),\ldots ,r_l=\pi(b_l,r_l)$ together with $\bar{c}_1,\ldots ,\bar{c}_m$ form a basis of ${\mathbb R}^n$.
Then the isomorphism $\iota_\phi$ is defined by the formula
$$
\iota_\phi(h) = \frac{1}{ \det(\bar{c}_1,\ldots , \bar{c}_m, r_1\,\ldots ,r_l)}((a_1,0)\wedge\ldots \wedge (a_k,0)\wedge \bar{b}_1\wedge\ldots \wedge \bar{b}_l)\otimes1^\ast.
$$
\end{proposition}
This follows immediately from the definition of the isomorphism $\Phi_{\bm{E}}$ associated with  the exact sequence ${\bm{E}}$. Our choice of isomorphism $\lambda^\ast({\mathbb R}^n)\rightarrow {\mathbb R}^\ast$ is, of course, dual
to a uniquely determined choice of isomorphism ${\mathbb R}\rightarrow\lambda({\mathbb R}^n)$,  defined by the mapping,  which maps $1$ to the wedge of the standard basis.
Its  inverse is the map 
\begin{equation}\label{iso--1}
c_1\wedge\ldots \wedge c_n \rightarrow \det([c_1,\ldots ,c_n]).
\end{equation}
Here again $\det$ is the determinant of a $n\times n$ matrix.  Having chosen  the isomorphisms (\ref{iso--}) and (\ref{iso--1}) we deduce the following natural identifications
for every finite dimensional vector space $V$,
$$
V\otimes \lambda({\mathbb R}^n)\rightarrow V\quad  v\otimes (e_1\wedge\ldots \wedge e_n)\mapsto  v\otimes 1\rightarrow v, 
$$
and
$$
V\otimes\lambda^\ast({\mathbb R}^n)\rightarrow V:v\otimes (e_1\wedge\ldots \wedge e_n)^\ast\rightarrow v\otimes 1^\ast\rightarrow v.
$$
We need another convention.
\begin{definition}\index{D- Natural isomorphism for $\det(T\oplus S)$}\index{$\det(T)\otimes \det(S)\rightarrow \det(T\oplus S)$, natural isomorphism}
 Assume that $T\colon E\rightarrow F$ and $S\colon E'\rightarrow F'$ are Fredholm operators. Then the direct sum $T\oplus S\colon E\oplus E'\rightarrow F\oplus F'$ is a Fredholm operator, and we define the {\bf natural isomorphism} 
$$
\det(T)\otimes \det(S)\rightarrow \det(T\oplus S)
$$
as follows. We choose  $h=(a_1\wedge\ldots \wedge a_n)\otimes (d_1\wedge\ldots \wedge d_l)^\ast\in \det (T)$ and  $h'=(a_1'\wedge\ldots \wedge a_{n'}')\otimes (d_1'\wedge\ldots \wedge d_{l'}')^\ast\in \det (S)$.
Then we map $h\otimes h'$ to the vector $g\in \det (T\oplus S)$,  defined by
$$
g=(a_1\wedge\ldots \wedge a_n\wedge a_1'\wedge\ldots \wedge a_{n'}')\otimes(d_1\wedge\ldots \wedge d_l\wedge d_1'\wedge\ldots \wedge d_{l'}')^\ast.
$$
\qed
\end{definition}

\section{Classical Local Determinant Bundles}
We continue  with some local constructions in the neighborhood of  a Fredholm operator $T\colon E\rightarrow F$ between two Banach spaces.
\begin{lemma}
Given a Fredholm operator $T\colon E\rightarrow F$ and a projection  $P\in \Pi_T$, we take a topological complement $Y$ of $\ker(PT)$ in $E$.
Then there exists $\varepsilon>0$ such  that for every $S\in {\mathcal L}(E,F)$ satisfying $\norm{S-T}<\varepsilon$,  the following holds.
\begin{itemize}
\item[{\em (1)}]\ $P\in\Pi_S$.
\item[{\em (2)}]\ $PS\vert Y\colon Y\rightarrow R(P)$ is a topological isomorphism.
\end{itemize}
\end{lemma} 
\begin{proof}
We estimate $\norm{PT\vert Y-PS\vert Y}\leq \norm{P}\cdot\norm{S-T}$. Since $PT\vert Y\colon Y\rightarrow R(P)$ is a topological linear isomorphism, 
it follows from the openness of invertible linear operators $Y\rightarrow R(P)$ that (2) holds for a suitable $\varepsilon$.
Since $PS\colon Y\rightarrow R(P)$ is a topological linear isomorphism,  we conclude that $R(PS)=R(P)$,  which implies
$P\in\Pi_S$.
\qed \end{proof}

We assume the hypotheses of  the previous lemma, split $E=\ker(PT)\oplus Y$,  and  correspondingly write $e=k+y$.
For $S$ satisfying  $\norm{S-T}<\varepsilon$ we consider the equation
$PS(k+y)=0$ which can be rewritten as
$$
PS(y)=-PS(k).
$$
Hence $y=y(k,S))=-(PS\vert Y)^{-1}(PS(k))$. The map $(k,S)\mapsto  y(k,S)$ is continuous. 
We equipe  the topological space
$$
\bigcup_{\norm{S-T}<\varepsilon} \{S\}\times \ker(PS)\subset {\mathcal L}(E,F)\oplus E
$$
with the induced topology. 
The projection
$$
\pi\colon  \bigcup_{\norm{S-T}<\varepsilon} \{S\}\times \ker(PS)\rightarrow  \bigcup_{\norm{S-T}<\varepsilon} \{S\}
 $$
 is the restriction of the continuous projection ${\mathcal L}(E,F)\oplus E\rightarrow {\mathcal L}(E,F)$ and therefore continuous. The fibers of $\pi$ are finite-dimensional vector spaces of the same dimension.  
\begin{lemma}\label{reddat}
The topological  space $\bigcup_{\norm{S-T}<\varepsilon} \{S\}\times \ker(PS)$ together with the  continuous projection $\pi$ has the structure of a trivial  vector bundle.
\end{lemma}
\begin{proof}
The continuous and bijective map
$$
\{S\, \vert \,  \norm{S-T}<\varepsilon\}\times\ker(PT)\rightarrow \bigcup_{\norm{S-T}<\varepsilon} \{S\}\times \ker(PS) 
$$
is defined by
$$
(S,k)\mapsto (S,k+y(k,S)).
$$
Its inverse is the restriction of the continuous map
$$
\{S\,  \vert \, \norm{S-T}<\varepsilon\}\times E\rightarrow\{S\,  \vert \,  \norm{S-T}<\varepsilon\}\times \ker(PT),\quad (S,k+y)\mapsto (S,k), 
$$
so that our map is indeed a topological bundle trivialization. 
\qed \end{proof}

Now we have the trivial bundle 
$$
\pi\colon \bigcup_{ \norm{S-T}<\varepsilon} \{S\}\times\ker(PS)\rightarrow \{S\,  \vert \, \norm{S-T}<\varepsilon\},
$$
and the  product bundle 
$$
\pi_0\colon \{S\,  \vert \, \norm{S-T}<\varepsilon\}\times(F/R(P))\rightarrow \{S\,  \vert \,  \norm{S-T}<\varepsilon\}.
$$
The obvious viewpoint about these bundles is the following. Given the family of Fredholm operators $S\mapsto  PS$ defined for $S$ satisfying  $\norm{S-T}<\varepsilon$, 
the kernel dimension and cokernel dimension is constant. Therefore,  we obtain the topological kernel bundle as well as the cokernel bundle.
We need the following standard lemma whose proof is left to the reader.

\begin{lemma}\label{lemma6.21}
Given finite-dimensional topological vector bundles $E$ and $F$ over the topological space $X$,  the bundles $\lambda(E)\rightarrow X$, $\lambda^\ast(F)\rightarrow X$
as well as $E\oplus F\rightarrow X$, and $E\otimes F\rightarrow X$ have the structure of topological vector bundles in a natural way.
\qed
\end{lemma}

{
From Lemma \ref{lemma6.21} we deduce  immediately the following result.
\begin{proposition}\label{erty-o}
Let $T:E\rightarrow F$ be a Fredholm operator and $P\in \Pi_T$. Then the topological space 
$$
\bigcup_{\norm{S-T}<\varepsilon} \{S\}\times \det(PS), 
$$
together with the projection onto the set $\{S\, \vert \, \norm{S-T}<\varepsilon\}$ 
has in a natural way the structure of topological line bundle.
\qed
\end{proposition}
}

{We introduce,  for the Fredholm operator $T$ and the projection $P\in\Pi_T$,  the abbreviation
$$
\text{DET}(T,P,\varepsilon)=\bigcup_{\norm{S-T}<\varepsilon} \{S\}\times \det(PS),
$$
where $\varepsilon>0$ is guaranteed by Lemma \ref{reddat}. We shall call $\text{DET}(T,P,\varepsilon)$ the {\bf local
determinant bundle} \index{Local determinant bundle} associated with  the Fredholm operator  $T$, the projection $P\in\Pi_T$,  and $\varepsilon$.
}

{
For the projection  $Q\in\Pi_T$ satisfying  $Q\leq P$ we abbreviate
$$
\text{DET}(T,Q,\varepsilon)=\bigcup_{\norm{S-T}<\varepsilon} \{S\}\times \det(QS)
$$
which,  again in a natural way,  is a topological line bundle. 
\begin{lemma}\label{kimportant}
The algebraic isomorphism
$$
\what{\gamma}_{T,P,Q,\varepsilon}:\text{DET}(T,P,\varepsilon)\rightarrow \text{DET}(T,Q,\varepsilon):(S,h)\rightarrow (S,\gamma^Q_{PS}(h))
$$ 
is a topological line bundle isomorphism.
\end{lemma}
\begin{proof}
This is trivial and follows from an inspection of the maps $\gamma_{PS}^Q$.  A sketch of the proof goes as follows.
We start with the exact sequence
$$
0\rightarrow \ker(PS)\rightarrow \ker(QS)\rightarrow F/R(Q)\rightarrow F/R(P)\rightarrow 0.
$$
Here $\norm{S-T}<\varepsilon$. Since $Q,P\in\Pi_T$ and $Q\leq P$,  the kernels of $PS$ and $QS$ define topological bundles if 
we vary the Fredholm operators $S$. This is  trivially true for the cokernel bundles associated to $PS$ and $QS$,  which are honest product bundles,  in view of
$R(PS)=R(P)$ and $R(QS)=R(Q)$. We can take linearly independent continuous sections which span  the kernels of the $PS$. Similarly,  for the cokernels we can take constant sections. Now going through the construction of the maps $\gamma_{PS}^Q$,  we see that we can extend the kernel sections for the $PS$ to
a family of continuous sections which are  point-wise linearly independent and  span the kernels of the $QS$. Proceeding  the same way  with the cokernels, we obtain at the end a continuous family of point-wise linearly independent sections spanning $\ker(QS)$ and $F/R(Q)$. Taking the appropriate wedges
and tensor products we see  that a continuous section of the first line bundle is mapped to a continuous section of the second.
The argument can  be reversed to verify  the continuity of the inverse map.
\qed \end{proof}
}

{
Now we introduce  the bundle
$$
\text{DET}_{(E,F)}=\bigcup_{S\in{\mathcal F}(E,F)} \{S\}\times \det (S)\index{$\text{DET}_{(E,F)}$}
$$
over the space of all Fredholm operators.
At first sight this set seems not to have a lot of structure globally (though near certain $S$ it has some due to the previous discussion). The problem
is that  in the definition of $\det(S)$ the `ingredients' 
$\ker(S)$ and $F/R(S)$ associated to  $S$ might have varying local dimensions. However, in view of the discussion in the previous subsection, we shall see, that we can equip $\text{DET}_{(E,F)}$  with the structure
of a topological line bundle over the space of Fredholm operators ${\mathcal F}(E,F)$. This structure will turn out to be natural, i.e.,  the structure  does not depend on the choices  involved in  its construction. Here, it is important to follow the conventions already introduced.
}

{
We have the projection map
$$
\text{DET}_{(E,F)}\rightarrow {\mathcal F}(E,F),\quad (S,h)\mapsto  S.
$$
The base ${\mathcal F}(E,F)$ has a topology, but the total space $\text{DET}_{(E,F)}$ at this point has not.
Given $T\in{\mathcal F}$ we can invoke Proposition \ref{erty-o} and find $\varepsilon>0$ and $P\in\Pi_T$
such  that 
$$
\gamma_{(T,P,\varepsilon)}\colon \text{DET}_{(E,F)}\vert \{S\, \vert \, \norm{S-T}<\varepsilon\}\rightarrow \text{DET}(T,P,\varepsilon),\quad (S,h)\mapsto  (S,\gamma^P_S(h))
$$
is a bijection which is fibers-wise linear and covers the identity on the base.
}

{
\begin{definition}
Denote by ${\mathcal B}$ the collection of all subsets $B$ of $\text{DET}_{(E,F)}$ having the property  that there exists $\gamma_{(T,P,\varepsilon)}$ for which 
$\gamma_{(T,P,\varepsilon)}(B)$ is open in $\text{DET}(T,P,\varepsilon)$. 
\qed
\end{definition} 
}

{
The fundamental topological observation is  the following proposition.
\begin{proposition}
The collection of sets ${\mathcal B}$ is a basis for a topology on $\text{DET}_{(E,F)}$.
\qed
\end{proposition}
}

Postponing the proof for the moment, we denote the topology associated with ${\mathcal B}$ by ${\mathcal T}$,  and  obtain the following result.
\begin{proposition}\index{P- Line bundle structure of $\text{DET}_{(E,F)}$}
The set $\text{DET}_{(E,F)}$ with its linear structure on the fibers, equipped with the topology ${\mathcal T}$, 
is a topological line bundle. The maps
$$
\gamma_{(T,P,\varepsilon)}\colon \text{DET}_{(E,F)}|\{S\ |\ \norm{S-T}<\varepsilon\}\rightarrow \text{DET}(T,P,\varepsilon),\quad (S,h)\rightarrow (S,\gamma^P_S(h))
$$
are topological bundle isomorphisms. 
\qed
\end{proposition}
In order to prove the previous two propositions we recall  the maps $\gamma_{(T,P,\varepsilon)}$ and $\gamma_{(T',P',\varepsilon')}$. We are interested in the transition map
$\gamma_{(T',P',\varepsilon')}\circ \gamma_{(T,P,\varepsilon)}^{-1}$. The domain of this map  consists of all pairs $(S,h)$ satisfying  $\norm{T-S}<\varepsilon$ and $\norm{S-T'}<\varepsilon'$ and $h\in \det(PS)$. The transition map preserves the base point and maps $h$ to an element in $\det(P'S)$. This transition map is obviously an algebraic isomorphism between
two topological line bundles, but as we have discussed in the previous subsection these transition maps are continuous.
\begin{lemma}
The transition map $\gamma_{(T',P',\varepsilon')}\circ \gamma_{(T,P,\varepsilon)}^{-1}$ is a topological isomorphism between topological line bundles.
\qed
\end{lemma}
The two previous propositions follows from this lemma, which summarizes previous assertions,
i.e. the discussion around Proposition \ref{ooo6.6} and Proposition \ref{PROP6310}.

Let us finally remark that it is known that the line bundle $\text{DET}_{(E,F)}$ is non-orientable. However,  over certain subsets it is orientable and that will be used in the discussion of sc-Fredholm sections.

\section{Local Orientation Propagation}\label{SECTX65}

Viewed locally, a classical Fredholm section $f$ is a smooth map $f\colon U\to F$ from an open subset $U$ of some Banach space $E$
into  a Banach space $F$. The derivatives $Df(x)$,  $x\in U$,  form a  family of Fredholm operators depending continuously as operators on $x\in U$.
The determinant line bundle of this family of Fredholm operators is a continuous (even smooth)  line bundle over $U$. If $U$ is contractible, this line bundle possesses precisely two possible orientations because a line has precisely two possible orientations.  Indeed, in this case, the orientation of one single line in the bundle determines 
a natural orientation of all the other lines. We might view 
this procedure as a local continuous propagation of an orientation.

This method, however, is not applicable in the sc-Fredholm setting because the linearizations do not, in general,  depend continuously on the base point. 
 There are even additional difficulties. For example, we are confronted with nontrivial bundles where the dimensions 
can locally jump. Nevertheless there is enough structure 
which allows to define a propagation of the orientation 
and this is the core of the orientation theory in the following section.

We consider  the sc-Fredholm section $f$ of the strong M-polyfold bundle $P\colon Y\rightarrow X$ over the tame $M$-polyfold $X$, choose  for a smooth point $x\in X$
a locally defined $\ssc^+$-section $s$ satisfying  $s(x)=f(x)$,  and take the linearization $(f-s)'(x)\colon T_xX\rightarrow Y_x$.  The linearization $(f-s)'(x)$ belongs to
the space of linearizations 
$\text{Lin}(f,x)$ which is a subset of linear Fredholm 
 operators $T_xX\to Y_x$ possessing the induced metric defined by the norm of the space ${\mathscr L}(T_xX, Y_x)$ of bounded operators.
Moreover, $\text{Lin}(f,x)$ is a convex subset and therefore contractible. Introducing 
$$
\text{DET}(f,x):=\bigcup_{L\in \text{Lin}(f,x)}\{L\}\times \text{det}(L)\index{$\text{DET}(f,x)$}
$$
over the convex set $\text{Lin}(f,x)$,  the previous 
discussion shows that $\text{DET}(f,x)$ is naturally a topological line bundle.
\begin{proposition}\index{P- Line bundle structure for $\text{DET}(f,x)$}
$\text{DET}(f,x)$ has in a natural way the structure of a topological line bundle over $\text{Lin}(f,x)$ and consequently has two possible orientations since the base space
is contractible.
\qed
\end{proposition}

\begin{definition}[{\bf Orientation}]\index{D- Orientations of $\text{DET}(f,x)$}
Let $f$ be a sc-Fredholm section of the strong M-polyfold bundle $P\colon Y\to X$ over the tame M-polyfold $X$ and $x\in X$ a smooth point. Then an {\bf orientation for the pair $(f,x)$},  denoted by 
$\mathfrak{o}_{(f,x)}$,  is a choice
of one of the two possible orientations of  
$\text{DET}(f,x)$.
\qed
\end{definition}

Let us denote by ${\mathcal O}_f$ the orientation space associated to $f$,  which is the collection of all pairs $(x,\mathfrak{o})$, in which $x\in X_\infty$ and $\mathfrak{o}$ is an orientation of $\text{DET}(f,x)$.
We consider a category whose  objects are the sets $Or_x^f:=\{\mathfrak{o}_x,-\mathfrak{o}_x\}$ of the two possible orientations of $\text{DET}(f,x)$. The reader should not be confused by our notation.
There is no distinguished class $\mathfrak{o}_x$. We only know there are two possible orientations, namely $\mathfrak{o}_x$ and $-\mathfrak{o}_x$.
For any two smooth points $x$ and $y$ 
in the same path component of $X$ the associated morphism set consists of two isomorphisms. The first one maps $\pm\mathfrak{o}_x\rightarrow \pm\mathfrak{o}_y$
and the second one $\pm\mathfrak{o}_x\rightarrow \mp\mathfrak{o}_y$.

\begin{proposition}\label{umbra}\index{P- Basic construction for orientation {I}}
If $f$ is a  sc-Fredholm section of the strong M-polyfold bundle $P\colon Y\rightarrow X$,  then every smooth point $x$ possesses  two open neighborhoods $U=U(x)$ 
and $U'=U'(x)$ having the following properties.
\begin{itemize}
\item[{\em (1)}]\  $U'\subset U$ is sc-smoothly contractible in $U$ to the point $x$.
\item[{\em (2)}]\ The solution set 
$\{x\in \cl_X (U)\, \vert \, f(x)=0\}$ is compact.
\item[{\em (3)}]\ Given a sc-smooth path $\phi\colon [0,1]\rightarrow U$,  there exists a $\ssc^+$-section $s\colon U\times [0,1]\rightarrow Y$
having  the property that $s(\phi(t),t)=f(\phi(t))$. 
\item[{\em (4)}]\ If  $t_0\in [0,1]$, $y_0\in U_\infty$,  and $e_0$ is  a smooth point in $Y_{y_0}$, then  there exists a $\ssc^+$-section $q\colon U\times [0,1]\rightarrow Y$ satisfying $q(y_0,t_0)=e_0$.
\end{itemize}
\end{proposition}

\begin{proof}
Taking local strong bundle coordinates we may assume that the bundle is $P\colon K\rightarrow O$ and $x=0$.  Here $K=R(V\triangleleft F)$ is the strong bundle retract covering the sc-retraction $r:V\rightarrow V$ satisfying  $r(V)=O$.
As usual,  $V\subset C$ is an open neighborhood of $0$ in the partial quadrant $C$ of the sc-Banach space $E$. 
Take a smaller open neighborhood $V'$ of $0\in C$ which is convex and contained in $V$.
Then we can define 
$$
\wh{r}\colon [0,1]\times V'\rightarrow O\colon \wh{r}(t,x) =r(tx).
$$
We abbreviate $O'=r(V')$ and note that $t\rightarrow \wh{r}(t,V')$ retracts $O'$ at $t=1$ to $x$ at $t=0$ while staying during the deformation
in $O$. We can always replace $V'$ by a smaller convex open neighborhood of $0$ and keep the property just described.

The  sc-Fredholm section $f$ possesses the local compactness property. We can therefore choose  an open  neighborhood $V$ of $0\in C$
such  that $U=r(V)$ satisfies  $\cl_C(U)\subset O$ and the restriction $f\vert \cl_C(U)$ has a compact solution set.
At this point we have constructed an open neighborhood $U$ of $x$ which has the first two properties.
Shrinking $U$ (and $U'$) further will keep these properties.  We  have to show that an additional shrinking will guarantee (3) and (4).

The section $f$ of $K\rightarrow O$ has the form $f(y)=(y,{\bf f}(y))$ where $R(f(y))=f(y)$.
The map  ${\bf f}$ is the  principle part of the section $f$ and is sc-smooth as a map from 
$O$ to $F$.
If $\phi\colon [0,1]\rightarrow U$ is a sc-smooth path, $t\in [0,1]$ and $y\in U$, then 
$(y,{\bf f}(\phi(t)))\in O\triangleleft F\subset V\triangleleft F$  and we define the map
$s\colon U\times [0,1]\rightarrow K$ by 
$$
s(y,t)=R(y,{\bf f}(\phi(t))).
$$
The map  $s$ is a $\ssc^+$-section of the pull-back of the bundle $K\rightarrow O$ by the map $U\times [0,1]\rightarrow O$, $(y, t)\mapsto y$.
By construction,  $s(\phi(t),t)=R(\phi(t),{\bf f}(\phi(t)))=R(f(\phi(t)))=f(\phi(t))$.  This proves (3).

The statement (4) was already proved in some variation  in the sections about Fredholm theory.  Formulated in the local
coordinates, the statement (4) assumes that  the  smooth point $e_0=(y_0,{\bf e_0})\in K$ and 
$t_0\in 0,1]$ are given.
The required  section $q\colon O\times [0,1]\to K$ can be then defined as the section 
$$
q(y,t)=R(y,{\bf e_0}).
$$
It satisfies  $q(y_0,t_0)=R(y_0,{\bf e_0})=(y_0,{\bf e_0})=e_0$, as desired.
\qed \end{proof}
Now we are in the position to define a propagation mechanism. Since it will involve several choices we have to make sure that the end result is independent of the choices involved. 
The set up is as follows. We have a strong bundle $P\colon Y\rightarrow X$ over the M-polyfold $X$ and an sc-Fredholm section $f$ of $P$. Around a smooth point $x\in X$ we choose  the  open neighborhoods $U=U(x)$, $U'(x)$,  so that the statements (1)-(4) of Proposition \ref{umbra} hold.

There is no loss of generality assuming that $X=U$. If 
 $\phi\colon [0,1]\rightarrow X$ is a sc-smooth path, we employ Proposition \ref{umbra}
and choose  a $\ssc^+$-section 
$s:X\times [0,1]\rightarrow W$ satisfying  $s(\phi(t),t)=f(\phi(t))$.

Adding  finitely many $\ssc^+$-sections $s_1,\ldots ,s_k$ defined on $X\times [0,1]$, we obtain the 
sc-Fredholm section
$$
F\colon  [0,1]\times X\times {\mathbb R}^k\rightarrow W,\quad F(t,y,\lambda)=f(y)-s(y,t)+\sum_{i=1}^k\lambda_i\cdot s_i(y,t)
$$
having the property that  for $\lambda=0$ the points $(t,\phi(t),0)$ are solutions of  
$$
F(t,\phi(t),0)=0.
$$
 For fixed $t\in [0,1]$ we introduce the  sc-Fredholm section
$$
F_t\colon X\times {\mathbb R}^k\rightarrow Y,\quad (y,\lambda)\rightarrow F(t,y,\lambda).
$$
The discussion, so far,  is true for all choices $s_1,\ldots ,s_k$. Using the results from the transversality theory we can choose  
$s_1,\ldots ,s_k$ such  that, in addition, $F_t$ has a good boundary behavior.

\begin{lemma}  \label{lem_property_star}
There exist finitely many $\ssc^+$-sections $s_1,\ldots ,s_k$ of the bundle  $Y\rightarrow [0,1]\times X$
such  that the sc-smooth Fredholm section $F$ of the bundle $Y\rightarrow [0,1]\times {\mathbb R}^k\times X$,  defined by
$$
F(t,y,\lambda) = f(y)-s(y,t)+\sum_{i=1}^k \lambda_i\cdot s_i(y,t),
$$
has the following property (P).
\begin{itemize}
\item[{\em (P)}]\ For every fixed $t\in [0,1]$,  the section $F_{t}$ is at the point $(\phi(t),0)$  in general position to the boundary of $X\times {\mathbb R}^k$.
\end{itemize}
\qed
\end{lemma}
Property (P) automatically implies that $F$ is in general position to the boundary of $[0,1]\times X\times {\mathbb R}^k$
at all points $(t,\phi(t),0)$ for $t\in [0,1]$. As a consequence of the implicit function theorem for the  boundary case we obtain the following result.
\begin{lemma}
The solution set 
$$
S=\{ (t, y, \lambda)\, \vert \,  F(t,y,\lambda)=0\}
$$
 for $(t,y,0)$ near $(t,\phi(t),0)$ is a smooth manifold with boundary with corners and
 $$
 \{(t,\phi(t),0)\, \vert \,  t\in [0,1]\}\subset S.
 $$
  Moreover,  it follows from the property (P) that the projection 
$$
\pi\colon S\rightarrow [0,1], \quad (t,y,\lambda)\mapsto  t
$$
is a submersion. The set  $S_t$ defined by  $\{t\}\times S_t:=\pi^{-1}(t)$  is a manifold with boundary with corners contained
in $ X\times {\mathbb R}^k$.
\qed
\end{lemma}
We do not claim that $S_t$ is compact, but the manifold $S_t$  lies in such a way in $X\times {\mathbb R}^k$ that its intersection with $(\partial X)\times {\mathbb R}^k$
carves out a boundary with corners on the manifold $S_t$.
The tangent space $T_{(\phi (t), 0)}S_t$ at the point $(\phi (t), 0)\in S_t$ agrees with the kernel $\ker (F_t'(\phi (t), 0))$ and, abbreviating  
$$
L_t:=T_{(\phi(t),0)}S_t=\ker (F_t'(\phi (t), 0))\quad \text{for $t\in [0,1]$},$$
we introduce the bundle  $L$ of tangent spaces along the path $\phi (t)\in S$, for $t\in [0,1]$, by 
$$
L=\bigcup_{t\in [0,1]} \{t\}\times L_t.
$$
The bundle $L$ is a smooth vector bundle over $[0,1]$. By $\lambda (L)$ we denote the line bundle associated with $L$,
$$\lambda (L)=\bigcup_{t\in [0,1]} \{t\}\times \lambda (L_t).$$

An orientation of the line $\lambda (L_t)$ determines by continuation an orientation of all the other lines. In particular, an orientation of $\lambda (L_0)$ at $t=0$ determines an orientation of $\lambda (L_1)$ at $t=1$, and we shall relate these orientations to the orientations of
$\text{DET}(f,\phi(0))$ and $\text{DET}(f,\phi(1))$.
To this aim we introduce, for every fixed $t\in [0,1]$, the exact sequence ${\bm{E}}_t$ defined by 
\begin{equation}\label{eexact}
\begin{gathered}
{\bm{E}}_t:\quad 0\rightarrow \ker( (f-s(\cdot,t ))'(\phi(t)) )\xrightarrow{j} \ker ({F}_t'(\phi(t),0))\xrightarrow{p}\\
\xrightarrow{p} {\mathbb R}^k\xrightarrow{c} Y_{\phi(t)}/R((f-s(\cdot,t ))'(\phi(t)))\rightarrow 0
\end{gathered}
\end{equation}
in which  $j$ is the inclusion map, $p$ the projection onto the ${\mathbb R}^k$-factor, and the  map $c$ is defined
by 
$$
c(\lambda) =\left(\sum_{i=1}^k \lambda_is_i(t,\phi(t))\right)+ R((f-s(\cdot,t ))'(\phi(t))).
$$

\begin{lemma}\label{lem_exact_seq_2}
The sequence \eqref{eexact} is exact.
\end{lemma}
\begin{proof}
The inclusion map $j$ is injective and $p\circ j=0$. From  $p(h,\lambda)=0$ it follows that 
$\lambda=0$
so that $h\in \ker((f-s(\cdot,t ))'(\phi(t)))$. If $(h,\lambda)\in \ker(F_t'(\phi(t),0))$,  then $\sum \lambda_is_i(\phi(t),t)$
belongs to the image of $(f-s(\cdot,t ))'(\phi(t))$ which implies $c\circ p=0$. It is also immediate that an
element $\lambda\in {\mathbb R}^k$ satisfying  $c(\lambda)=0$ implies that $\sum \lambda_is_i(\phi(t),t)$ belongs to the image
of $(f-s(\cdot,t ))'(\phi(t))$. This allows us to construct an element $(h,\lambda)\in \ker(F_t'(\phi(t),0))$ satisfying  $p(h,\lambda,)=\lambda$. 
Finally, it follows from the property $(P)$ in Lemma \ref{lem_property_star}  that the map $c$ is  surjective. 
The proof of Lemma 
\ref{lem_exact_seq_2} is complete.
\qed \end{proof}

From the exact sequence ${\bm{E}}_t$ we deduce the natural isomorphism $\Phi_{{\bm{E}}_t}$ introduced in Section \ref{sect_conventions}, 
\begin{equation*}
\begin{split}
\Phi_{{\bm{E}}_t}&\colon \lambda \bigl(\ker (f-s(\cdot , t))'(\phi (t))\bigr)\otimes \lambda^\ast (\coker  (f-s(\cdot , t))'(\phi (t))\bigr)\\
&\to \lambda \bigl(\ker F_t'(\phi (t), 0)\otimes \lambda^\ast(\R^k),
\end{split}
\end{equation*}
and obtain, in view of the the definition of the determinant, the isomorphism
$$
\Phi_{{\bm{E}}_t}\colon \det ((f-s(\cdot ,t))'(\phi (t)))\to \lambda (L_t)\otimes \lambda^\ast(\R^k).
$$
Now we assume that we have chosen the orientation 
$\mathfrak{o}_0$ of $\text{DET}(f,\phi(0))$  at $t=0$. It induces an orientation of $\det ((f-s(\cdot ,0))'(\phi (0))).$
The isomorphism $\Phi_{{\bm{E}}_0}$ determines an orientation of $\lambda (L_0)\otimes \lambda^\ast (\R^k)$. By continuation we extend this orientation to an orientation of  $\lambda (L_1)\otimes \lambda^\ast (\R^k).$ Then the isomorphism $(\Phi_{{\bm{E}}_1})^{-1}$ at $t=1$ gives us an orientation of $\det ((f-s(\cdot ,1))'(\varphi (1)))$ and consequently an orientation $\mathfrak{o}_1$ of 
$\text{DET}(f,\phi(1))$  at $t=1$, denoted by 
$$\mathfrak{o}_1=\phi_\ast \mathfrak{o}_0,$$
and called the {\bf push forward orientation}.
A priori our procedure might depend on the choices of the sections $s(y, t)$ and $s_1\ldots ,s_k$ and we shall prove next that it actually does not.
\begin{proposition}\label{local-xyz}\index{P- Basic construction for orientation {II}}
We assume that $f$ is the sc-Fredholm section of the strong M-polyfold bundle $P\colon Y\rightarrow X$ over the tame $X$, and $U=U(x)$, $U'=U'(x)\subset U$ are  open neighborhood  around a smooth point $x$, for which the conclusions of Proposition \ref{umbra} hold. Let $\phi\colon [0,1]\rightarrow U'$ be a  sc-smooth path. Then the construction of the map
$\mathfrak{o}\rightarrow\phi_\ast\mathfrak{o}$,
associating with  an orientation of $\text{DET}(f,\phi(0))$ an orientation $\phi_\ast\mathfrak{o}$  of $ \text{DET}(f,\phi(1)) $ does not depend on the choices involved as long as the hypotheses of Lemma \ref{lem_property_star} hold. 
\qed
\end{proposition}

The proof of the proposition follows from two lemmata.
\begin{lemma}\label{opas}
Under the assumption of the proposition we assume that 
the sc-Fredholm section 
$$
F(t,y,\lambda)=f(y)-s(y,t)+\sum_{i=1}^k\lambda_i\cdot s_i(y,t),
$$ 
satisfies the property  (P) from Lemma \ref{lem_property_star}.
We view $F$ as section of the bundle $Y\rightarrow [0,1]\times X\times {\mathbb R}^k$. Adding more $\ssc^+$-sections we 
also introduce the second sc-Fredholm section 
$$
\ov{F}(t,y,\lambda,\mu)=F(t,y,\lambda)+\sum_{j=1}^l\mu_j\cdot \bar{s}_j(y,t),
$$
viewed as a section of the bundle 
$Y\rightarrow [0,1]\times X\times {\mathbb R}^{k+l}$. 
Then both sections define the same propagation of the orientation along the sc-smooth path $\phi$.
\end{lemma}
\begin{proof}
The key is to view $F$ as the section
$$
\wh{F}(\tau, t,y,\lambda, \mu)= f(y)-s(y,t)+\sum_{i=1}^k\lambda_i\cdot s_i(y,t)
$$
of the strong M-polyfold bundle $Y\rightarrow [0,1]\times [0,1]\times X\times {\mathbb R}^{k+l}$.
There is a sc-Fredholm section 
$$\wt{F}\colon
[0,1]\times[0,1]\times X\times{\mathbb R}^{k+l}\rightarrow Y, 
$$
defined 
by
$$
\wt{F}(\tau,t,y,\lambda,\mu)=\wh{F}(\tau, t,y,\lambda,\mu)+\sum_{j=1}^l \mu_j\cdot \bar{s}_j(y,t).
$$
The induced map $X\times{\mathbb R}^{k+l}\rightarrow Y$ for fixed $(\tau,t)\in [0,1]\times [0,1]$ satisfies the property (P).
We obtain a solution manifold $\tilde{S}$ with boundary with corners  defined near the points $(\tau,t,\phi(t),0)$,which fibers over $[0,1]\times [0,1]$.
Let $\tilde{\pi}\colon \tilde{S}\rightarrow [0,1]\times [0,1]$ and consider its tangent map.
Then take the vector bundle $\tilde{L}\rightarrow [0,1]\times [0,1]$ whose fiber over $(\tau,t)$ is the preimage of the zero section of $T([0,1]\times [0,1])$ under $T\tilde{\pi}$.
For $\tau=0$ we obtain the bundle $\tilde{L}^0\rightarrow \{0\}\times [0,1]=[0,1]$. This is of the form
$$
L\oplus ([0,1]\times {\mathbb R}^l)\rightarrow [0,1],
$$
where $L$ is the bundle associated to $F$ which we used originally to define the propagation along $\phi$. 
Applying the construction now using $\tilde{L}^0$ one verifies by a simple computation that it defines the same 
propagation.

Hence it remains to verify that every bundle $\wt{L}^\tau$, $\tau\in [0,1]$,  defines  the same propagation. 
First of all we note that as $\tau$ varies,  the bundle varies continuously. Moreover, we relate the end points 
to $\det((f-s(\cdot,0))'(\phi(0)))$ and $\det((f-s(\cdot,1))'(\phi(1)))$. The data in the occurring exact sequence
relating  the latter with  the corresponding orientations of $\tilde{L}^\tau$,  vary continuously.  This implies that the propagation definition 
does not depend on $\tau$.
\qed \end{proof}

Next we show that the choice of the $\ssc^+$-section $s$ satisfying  $s(\phi(t),t)=f(\phi(t))$ does not affect the definition of the propagation.
\begin{lemma}\label{opas1}
We assume that the $\ssc^+$-sections $s^i$ for $i=0,1$ satisfy 
$$
s^i(\phi(t),t)=f(\phi(t))
$$
 and consider the homotopy $s^\tau=\tau \cdot s^1+(1-\tau)\cdot s^0$.
We choose  additional sections $s^\tau_1,\ldots ,s^\tau_k$  so that the property (P) holds together with $s^\tau$  for every $\tau$.
Then the  propagation of the orientation along the sc-smooth path $\phi$ using each of the collections $s^\tau,s_1^\tau,\ldots ,s_k^\tau$ is the same.
 \end{lemma}
 
 \begin{proof}
 The proof is similar as the lemma above. We obtain a vector bundle $\wtilde{L}\rightarrow [0,1]\times [0,1]$
 with parameters $(\tau,t)$ in the base. Then we can define for fixed $\tau$ a vector bundle over $t\in[0,1]$, say $\wt{L}^\tau$.
 Clearly $\wt{L}^\tau$ varies continuously in $\tau$.  At  $t=0,1$ the linearizations
 $(f-s^\tau(\cdot,0))'(\phi(0))$ and $(f-s^\tau(\cdot,1))'(\phi(1))$ vary continuously as operators in $\text{DET}(f,\phi(0))$ and $\text{DET}(f,\phi(1))$, respectively.
 This implies that the bundle $\tilde{L}^\tau$ induces,  independently of $\tau\in [0,1]$,  the same propagation of the orientation  along the path $\phi$.
\qed \end{proof}
Now we are in the position to finish the proof of Proposition \ref{local-xyz}.
We start with two  collections $s^i,s_1^i,\ldots ,s_{k^i}^i$ to define the propagation and add for $i=0$ and $i=1$
additional $\ssc^+$-sections such that the following holds.
\begin{itemize}
\item[(1)]\ For both situations $i=0$ and $i=1$ we have the same number of sections.
\item[(2)]\ For  a convex homotopy parametrized by $\tau\in [0,1]$, the property (P) holds for each fixed $\tau$.
\end{itemize}
By Lemma \ref{opas} adding section does not change the propagation. We can apply this for $i=0,1$.
Then we can use Lemma \ref{opas1} to show that these two propagations (after adding sections) are the same.

Next assume we have the same hypotheses and $U=U(x)$, $U'=U'(x)$  have the initially stated properties.
If the two sc-smooth paths $\phi^0,\phi^1\colon [0,1]\rightarrow U'$ have the same starting and end points,  we homotope sc-smoothly from one to the other with end points fixed,  using that $U' $ is sc-smoothly contractible in $U$.  We 
denote this sc-smooth homotopy by
$$
\Phi\colon [0,1]\times [0,1]\rightarrow U
$$
where $\phi^i=\Phi(i,0)$ for $i=0,1$. Then we  construct the $\ssc^+$-section $s$ satisfying $s(\Phi(\tau,t),\tau,t)=f(\Phi(\tau,t))$ and define the sc-Fredholm section  
$$
(\tau,t,y,\lambda)\rightarrow f(y)-s(y,\tau,t)+\sum_{i=1}^k\lambda_i\cdot s_i(y,\tau,t)
$$
possesing the obvious properties. 
For every $\tau$ we obtain a vector bundle $\wt{L}^\tau\rightarrow [0,1]$. Since all the data change continuously in these
bundles and the operators for $t=0,1$ in $\text{DET}(\ast)$  change continuously in $\tau$,  we see that every $\wt{L}^\tau$ defines the same 
propagation of the orientation along paths $\phi^\tau=\Phi(\tau,\cdot )$. Hence we have proved the following statement.

\begin{theorem}\label{thm_6.35}\index{T- Propagation of orientation}
We assume that $P\colon Y\rightarrow X$ is a strong bundle over the tame M-polyfold and  $x$ a smooth point in $X$.
Then there exists an open  neighborhood $\wt{U}=\wt{U}(x)$ such  that for every sc-smooth path $\phi\colon [0,1]\rightarrow \wt{U}$
there exists a well-defined propagation 
$$
\mathfrak{o}\rightarrow \phi_\ast\mathfrak{o}
$$
of an orientation $\mathfrak{o}$ of $\text{DET}(f,\phi(0))$. Moreover,  if $\phi^0$ and $\phi^1\colon [0,1]\rightarrow \wt{U}$ are two sc-smooth paths from the  same starting points
to the same end point, then the  propagation along $\phi^0$ and $\phi^1$ is the same, i.e., 
$$
\phi^0_\ast\mathfrak{o}=\phi^1_\ast\mathfrak{o}.
$$
\qed
\end{theorem}
Having this theorem there are precisely two possible ways to orient the family $\text{DET}(f,y)$, $y\in \wt{U}_\infty$,  so that 
these orientations are related by propagation along paths. Namely,  we fix $x$ and take for $y\in \wt{U}_\infty$ a sc-smooth path
$\phi\colon [0,1]\rightarrow \wt{U}$, starting at $x$ and ending at $y$.  Fixing  an orientation $\mathfrak{o}_x$ of $\text{DET}(f,x)$, 
we define $\mathfrak{o}_y=\phi_\ast\mathfrak{o}_x$. The definition is independent of the choice of the path $\phi$.
At this point we have a map which associates to a smooth point $y\in \wt{U}$, i.e. $y\in \wt{U}_\infty$,  an orientation $\mathfrak{o}_y$ of $\text{DET}(f,y)$. It follows from Theorem \ref{thm_6.35} that if  $\psi$ is a sc-smooth path connecting $y_1$ with $y_2$, then $\psi_\ast\mathfrak{o}_{y_1}=\mathfrak{o}_{y_2}$.

\begin{definition}\label{SDEF6511}.
Given a strong bundle $P\colon Y\rightarrow X$ over the tame M-polyfold $X$ and a sc-Fredholm section $f$, we call an open path connected 
neighborhood $\wt{U}(x)$ of a smooth point $x$ on which the local propagation construction can be carried out an {\bf orientable neighborhood}\index{D- Orientable neighborhood} of $(f,x)$. Any two of the possible orientations of $\text{DET}(f,y)$ for $y\in \wt{U}_\infty$, which have the propagation property,  is called  a  {\bf continuous orientation}\index{D- Continuous orientations}.
\qed
\end{definition}
In view of Theorem \ref{thm_6.35}, given an sc-Fredholm section $f$ of a strong bundle $P:Y\rightarrow X$ over a tame $X$, for every  $(f,x)$ with $x\in X_\infty$,
there exists an orientable open neighborhood $\wt{U}(x)$.
A suitable orientable neighborhood  $\wt{U}(x)$ has precisely two continuous orientations. 
Now we globalize the local propagation of orientation constructions to a more global procedure.
Let $X$ be an M-polyfold and $\phi\colon [0,1]\rightarrow X$ a sc-smooth map. We denote by $[\phi]$ its sc-smooth homotopy class with end points fixed. Associated with  $[\phi]$ we have the source $s([\phi])$ which is the starting point $\phi(0)$ and the target $t([\phi])$ which is the end point $\phi(1)$.
Given two sc-smooth homotopy classes  $[\phi]$ and $[\psi]$ satisfying $t([\phi])=s([\psi])$,  we  define the composition $[\psi]\ast[\phi]$ as the class of $\gamma$ defined as follows.
Take a smooth map $\beta\colon [0,1]\rightarrow [0,2]$ satisfying  that $\beta([0,1/2])=[0,1]$, $\beta(0)=0$ and $\beta(s)=1$ for $s$ near $1/2$. Moreover, 
$\beta([1/2,1])=[1,2]$ and $\beta(1)=2$. Then we define the sc-smooth path $\gamma$ by $\gamma(t)=\phi(\beta(s))$ for $s\in [0,1/2]$ and $\gamma(t)=\psi(\beta(s))$ for $s\in [1/2,1]$.
This way we obtain a category ${\mathcal P}_X$ whose objects are the smooth points in $X$ and whose morphisms $x\rightarrow y$ are the homotopy  classes $[\phi]$ with source $x$ and target $y$.
The main result of this section is  described  by the following theorem which is a consequence of the local constructions.

\begin{theorem}\index{T- Main orientation result}
Let $P\colon Y\rightarrow X$ be a strong bundle over the tame M-polyfold $X$ and $f$ a sc-Fredholm section.
Then there exists a uniquely determined functor 
$$
{\mathcal P}_X\xrightarrow{\Gamma_f}{\mathcal O}_f,\index{${\mathcal P}_X\xrightarrow{\Gamma_f}{\mathcal O}_f$}
$$
which associates
with  the smooth point $x$ the set $Or_x^f$,  and with the sc-smooth homotopy class  $[\phi]\colon x\rightarrow y$ a morphism $\Gamma_f([\phi])\colon Or_x^f\rightarrow Or_y^f$,  which has the following property.\\[1ex]
$({\bf \ast})$\quad 
Given any smooth point $x$ in $X$  and  an orientable open neighborhood $\wt{U}=\wt{U}(x)$
the maps 
$$
\phi_\ast\colon Or_{\phi(0)}^f\rightarrow Or_{\phi(1)}^f\ \ \text{and}\ \ \Gamma^f([\phi])\colon Or_{\phi(0)}^f\rightarrow Or_{\phi(1)}^f
$$
coincide for a sc-smooth path $\phi\colon [0,1]\rightarrow \wt{U}$.
\end{theorem}
\begin{proof}
We consider the  sc-smooth path $\phi:[0,1]\rightarrow X$, fix $t_0\in [0,1]$, and choose an orientable neighborhood $\wt{U}(\phi(t_0))$.
Then we find an open interval $I(t_0)$ so that an orientation
$\mathfrak{o}_t$ for $(f,\phi(t))$, where $t\in I(t_0)\cap [0,1]$, determines an orientation for all $s\in I(t_0)\cap [0,1]$.  
Now using the compactness of $[0,1]$ and a finite covering, we can transport a given orientation $\mathfrak{o}_0$ of $\text{DET}(f,\phi(0))$ to an orientation 
$\mathfrak{o}_1$ of $\text{DET}(f,\phi(1))$. The map $\pm\mathfrak{o}_0\rightarrow \pm\mathfrak{o}_1$ is denoted by $\phi_\ast$. By a covering argument one verifies that the map $\phi_\ast$ only depends on the homotopy class $[\phi]$ for fixed end-points. It is obvious that $\Gamma_f$ is unique.
\qed \end{proof}
Finally we can give two equivalent definitions of orientability of a sc-Fredholm section.
\begin{definition}\index{D- Orientable sc-Fredholm section {I}}
A sc-Fredholm section  $f$  of the strong bundle $P\colon Y\rightarrow X$ over the tame M-polyfold $X$ is called  {\bf orientable}  provided for every pair of smooth points $x$ and $y$
in the same path component,  the morphism $\Gamma_f([\phi])$ does not depend on the choice of the homotopy class $[\phi]\colon x\rightarrow y$. 
\qed
\end{definition}
An equivalent version is the following.
\begin{definition} \index{D- Orientable sc-Fredholm section {II}}
Let $f$ be a sc-Fredholm section of the strong bundle $P\colon Y\rightarrow X$ over the tame M-polyfold $X$. Then $f$ is  called {\bf orientable}  if there exists 
a map which associates with a point  $y\in X_\infty$ an orientation $\mathfrak{o}_y$ of $\text{DET}(f,y)$ having the following property. For every smooth point $x$ there exists
an orientable neighborhood $U(x)$ for which  $U_\infty(x)\ni y\mapsto \mathfrak{o}_y$ is one of the two continuous orientations.
\qed
\end{definition}
Finally we  define an orientation of the sc-Fredholm section as follows.
\begin{definition}[{\bf Orientation of sc-Fredholm sections}]\index{D- Orientation of an sc-Fredholm section}
Let $P\colon Y\rightarrow X$ be a strong bundle over the tame M-polyfold $X$ and assume that $f$ is a  sc-Fredholm section.
An orientation $\mathfrak{o}^f$ for $f$ is a map which associates with every  smooth  point $y\in X$ an orientation $\mathfrak{o}_y^f$ of $\text{DET}(f,y)$
such  that,  for every smooth point $x$ and orientable neighborhood $U=U(x)$,  the restriction $\mathfrak{o}^f\vert U_\infty$ is one of the two possible continuous orientations.
\qed
\end{definition}
We close this subsection with a general result.
\begin{theorem}\index{T- Orientation of solution set}
We assume that $P\colon Y\rightarrow X$ is a strong bundle over the tame M-polyfold $X$. Let $f$ be a proper sc-Fredholm section
such  that for every $x\in X$ solving $f(x)=0$ the pair $(f,x)$ is in 
general position to the boundary $\partial X$.
Suppose $\mathfrak{o}^f$ is an orientation for $f$. Then the solution set $S=f^{-1}(0)$ is a smooth compact manifold with boundary with corners
possessing  a natural orientation.
\end{theorem}
\begin{proof}
We already know from the transversality discussion that $S$ is a compact manifold with boundary with corners.
If  $x\in S$, then $\ker(f'(x))=T_xS$ and because  $f'(x)$ is surjective,  every tangent space $T_xS$ is oriented by $\mathfrak{o}_x^f$.
Since $\mathfrak{o}^f$ has the local continuation property,  one verifies that the differential geometric local prolongation of the orientation on $S$ is
the same as the Fredholm one.
\qed \end{proof}

If $f$ is not  generic to start with,  we have to take a small perturbation $s$ supported near the zero set of $f$, so that  the solution set $S_{f+s} $ has a natural orientation coming from $\mathfrak{o}^f$ since $\det((f+s)'(x))$ has an orientation coming from $\text{DET}(f,x)$ for $x\in S_{f+s}$.

\begin{remark} \index{R- On orientations for sc-Fredholm sections}
If $(f,\mathfrak{o})$ is an oriented Fredholm section of the strong bundle $P\colon Y\rightarrow X$ and $s$ is a $\ssc^+$-section of $Y\rightarrow X\times [0,1]$, then the sc-Fredholm section 
$\wh{f}$,  defined on  $X\times[0,1]$ by $\wh{f}(x,t)=f(x)+s(x,t)\in Y$, 
has a natural orientation associated with  
$\mathfrak{o}$ and the standard orientation of $[0,1]$. 
We choose  a local $\ssc^+$-section $t$ near the smooth point $x$  satisfying  $t(x)=f(x)$ and consider the orientation 
$\mathfrak{o}_x$ of  $\det((f-t)'(x))$, where $L=(f-t)'(x)\colon T_xX\rightarrow Y_x$. Taking the vector 
$h=(a_1\wedge\ldots \wedge a_n)\otimes {(b_1\wedge\ldots \wedge b_l)}^\ast$ determining the orientation of $L$, we orient 
$$
\wt{L}\colon T_{(x,t)}(X\times [0,1])\rightarrow Y_x,\quad (a,b)\rightarrow La
$$
by the vector 
$\wt{h}=(a_1\wedge\ldots \wedge a_n\wedge e)\otimes (b_1\wedge\ldots \wedge b_l)^\ast$,
where the vector $e$ defines the standard orientation of $[0,1].$
The linear map $\wt{L}$ belongs to the linearization space of $\wh{f}$ at $(x,t)$, so that the latter obtains an orientation.
\end{remark}

\section{Invariants}
We consider the oriented and proper sc-Fredholm section $(f, \mathfrak{o})$ of the strong bundle $P\colon Y\to X$ over the tame M-polyfold $X$ admitting sc-smooth bump functions. The solution set $S=\{x\in X\, \vert \, f(x)=0\}$ is compact. If $N$ is an auxiliary norm on $P$ we know from the perturbation and transversality result, Theorem \ref{thm_pert_and_trans}, that there exists an open neighborhood $U\subset X$ of $S$ such that for every $\ssc^+$-section $s$ of $P$ which is supported in $U$ and satisfies $N(s(x))\leq 1$ for all $x\in X$, the solution set $S_{f+s}=\{x\in X\, \vert \, f(x)+s(x)=0\}$ is compact. Moreover, there exist distinguished such 
$\ssc^+$-sections having the additional property that,  for all $x\in S_{f+s}$, the  pair $(f+s, x)$ is in general position to the boundary $\partial X$. Its solution set is a compact manifold with boundary with corners and, as we have seen above, possesses a natural orientation induced from the orientation $\mathfrak{o}$ of the sc-Fredholm section $f$.
In the following theorem we shall call these distinguished $\ssc^+$-sections $s$ {\bf admissible for the pair $(N, U)$.}

In Section \ref{subs_sc_differential} we have introduced the differential algebra 
$\Omega^\ast_\infty(X,\partial X)=\Omega_\infty^\ast(X)\oplus \Omega_\infty^{\ast-1}(\partial X).$ Its differential $d(w, \tau)=(d\omega, j^\ast w-d\tau)$, in which $j\colon \partial X\to X$ is the inclusion map, satisfies $d\circ d=0$ and we denote the associated cohomology by 
$H^\ast_{\textrm{dR}}(X, \partial X)$.

\begin{theorem}\index{T- Invariants for oriented proper sc-Fredholm sections}
Let $P\colon Y\rightarrow X$ be a strong bundle over the tame M-polyfold $X$. We assume that $X$ admits sc-smooth bump functions.
Then there exists a well-defined map which associates with  a proper and oriented sc-Fredholm section $(f, \mathfrak{o})$ of the bundle  $P$ a linear map 
$$
\Psi_{(f,\mathfrak{o})}\colon H^\ast_{\textrm{dR}}(X,\partial X)\rightarrow {\mathbb R}
$$
having the following properties.
\begin{itemize}
\item[{\em (1)}]\ If $N$ is an auxiliary norm on $P$ and $s$ a $\ssc^+$-section of $P$ which is admissible for the pair $(N, U)$, then the solution set  $S_{f+s}=\{x\in X\, \vert \, f(x)+s(x)=0\}$ is an oriented and compact manifold with boundary with corners, and 
$$
\Psi_{(f,\mathfrak{o})}([\omega,\tau]) = \int_{S_{f+s}} \omega -\int_{\partial S_{f+s}}\tau,
$$ 
holds for every cohomology class $[\omega,\tau]$ in $H^\ast_{dR}(X,\partial X)$. The integrals on the right-hand side are defined to be zero, if the dimensions of the forms and manifolds do not agree.  
\item[{\em (2)}]\ For  a proper and oriented homotopy between two oriented sc-Fredholm sections $(f_0,\mathfrak{o}_0)$ and $(f_1,\mathfrak{o}_1)$ of the bundle $P$  we have the identity
$$
\Psi_{(f_0,\mathfrak{o}_0)}=\Psi_{(f_1,\mathfrak{o}_1)}.
$$
\end{itemize}
\end{theorem}
\begin{proof}
(1)\,  The properness of the sc-Fredholm section $f$ implies the compactness of its solution set $S_f=\{x\in X\,\vert \, f(x)=0\}$. If $N$  is an auxiliary norm $N$, then there exist, by Theorem 
\ref{thm_pert_and_trans}, an open neighborhood $U\subset X$ of $S_f$ and a $\ssc^+$-section of $s_0$ of $P$ which is admissible for the pair $(N, U)$ such that the solution set $S_0=S_{f+s_0}=\{x\in X\, \vert \, f(x)+s_0(x)=0\}$ is, in addition, an oriented compact manifold with boundary with corners of dimension $n=\dim S_{f+s_0}=\ind (f'(x))$ for $x\in S_{f+s_0}$. 
The orientation $(f, \mathfrak{o})$ induces an orientation on $S_{f+s}$ and its local faces.  We define the map $\Psi_{(f,\mathfrak{o})}$ by 
\begin{equation}\label{int_eq_1}
\Psi_{(f,\mathfrak{o})}([\omega,\tau]):=\int_{S_{f+s_0}}\omega -\int_{\partial S_{f+s_0}}\tau
\end{equation}
for a cohomology class $[\omega, \tau]$ in $H^\ast_{\textrm{dR}}(X,\partial X)$. The integrals on the right-hand side are defined to be zero 
if the dimensions of the forms and manifolds do not agree.

In order to verify that the 
definition \eqref{int_eq_1} does not depend on the choice of admissible section $s_0$,  we take a second $\ssc^+$-section $s_1$ admissible for the pair $(N, U)$ so that the solution set $S_1=S_{f+s_1}=\{x\in X\, \vert \, f(x)+s_1(x)=0\}$ is a compact oriented manifold with boundary with corners. As explained in Remark \ref{rem_homotopy}, 
there exists a proper $\ssc^+$- homotopy $s_t$ connecting  $s_0$ with  $s_1$
such the sc-Fredholm section  $F$, defined on $X\times [0,1]$ by $F(x,t)=f(x)+s_t(x)$, is in 
general position to the boundary of $ X\times [0,1]$ for every  $(x,t)$ in the solution set   $S_F=\{(x, t)\in X\times [0,1]\, \vert \, f(x)+s_t(x)=0\}$, which is an 
oriented compact manifold with boundary with corners, whose orientation is induced from the orientation  $\mathfrak{o}$ of $f$ and the standard orientation of $[0,1]$.  We recall that  $[\omega,\tau]$ satisfies $d\omega=0$ and $j^\ast\omega=d\tau$ on $\partial X$ where $j\colon \partial X\to X$ is the inclusion map. Extending $\omega$ to the whole space $X\times [0,1]$, we introduce the form $\ov{\omega}$ on $X\times [0,1]$ by 
$$\ov{\omega}=p_1^\ast\omega\wedge p_2^\ast dt,$$
where $p_1\colon X\times [0,1]\to X$ and $p_2\colon X\times [0,1]\to [0,1]$ are the sc-projection maps. It follows that $d\ov{\omega}=0$ and $\ov{\omega}\vert \partial X\times [0,1]=d\ov{\tau}$, where $d\ov{\tau}=\ov{p}_1^\ast (d\tau)\wedge \ov{p}_2^\ast dt.$ The maps $\ov{p}_i$ for $i=1, 2$ are the restrictions of $p_i$ to $\partial X\times [0,1]$. 
By Stokes theorem,
$$\int_{S_F}d\ov{\omega}=\int_{\partial S_F}\ov{\omega}=0.$$
We decompose the boundary $\partial S_F$ of the manifold $S_F$ into 
$$\partial S_F=S_0\cup S_1\cup \wt{S},$$
where $S_0\subset \partial X\times \{0\}$, $S_1\subset\partial X\times \{1\}$, and $\wt{S}\subset \partial X\times (0, 1)$, as illustrated in Figure \ref{pict7}.

\begin{figure}[htb]
\begin{centering}
\def\svgwidth{70ex}
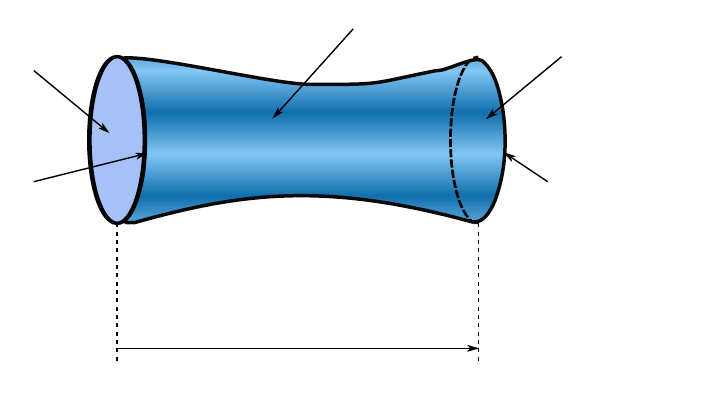
\caption{}\label{pict7}
\end{centering}
\end{figure}

Hence,
$$
0=\int_{\partial S_F}\ov{\omega}=\int_{\partial S_F}d\ov{\tau},
$$
and we compute the two integrals, taking the orientations into account,
$$0=\int_{\partial S_F}\ov{\omega}=\int_{S_0}\omega-\int_{S_1}\omega-\int_{\wt{S}}\ov{\omega}.$$
Integration of $d\ov{\tau}$ over $\partial S_F$ using Stokes theorem gives 
\begin{equation*}
\begin{split}
0&=\int_{\partial S_F}d\ov{\tau}=\int_{S_0}d\tau-
\int_{S_1}d\tau+\int_{\wt{S}}\ov{\omega}\\
&=\int_{\partial S_0}\tau-\int_{\partial S_1}\tau
+\int_{\wt{S}}\ov{\omega}.
\end{split}
\end{equation*}
Comparing the integrals, we 
obtain 
$$\int_{S_0}\omega-\int_{\partial S_0}\tau=\int_{S_1}\omega-\int_{\partial S_1}\tau,$$
which shows that the definition \eqref{int_eq_1} does indeed not depend on the choice of the distinguished section $s_0$.

The formula in (2) is verified by the same homotopy argument as in (1).
\qed \end{proof}

In the special situation $\partial X=\emptyset$ of no boundary, the map $\Psi_{(f,\mathfrak{o})}\colon H^n_{\textrm{dR}}(X)\rightarrow {\mathbb R}$ is defined by 
$$
\Psi_{(f,\mathfrak{o})}([\omega])=\int_{S_{f+s}}\omega
$$ 
for $[\omega]\in H^n_{\textrm{dR}}(X)$,  if $n$ is equal to the index of the sc-Fredholm section $f$, and zero otherwise.

\section{Appendix}
\subsection{Proof of Lemma \ref{o6.6}}\label{oo6.6}
The goal is to prove that given an exact sequence
$$
{\bm{E}}\colon\ \   0\rightarrow A\xrightarrow{\alpha} B\xrightarrow{\beta} C\xrightarrow{\gamma} D\rightarrow 0
$$
there is a well-defined natural isomorphism $\Phi_{\bm{E}} :\lambda(A)\otimes\lambda^\ast(D)\rightarrow \lambda(B)\otimes \lambda^\ast(C)$.
\begin{proof}[Lemma \ref{o6.6}]
Recalling the exact sequence 
$$
{\bm{E}}:\quad  0\rightarrow A\xrightarrow{\alpha} B\xrightarrow{\beta} C\xrightarrow{\gamma} D\rightarrow 0,
$$
we choose linear subspaces $V\subset B$ and $W\subset C$ such  that
$$
B=\alpha(A)\oplus V\quad   \text{and}\quad    C=\beta(B) \oplus W
$$
and denote this choice of subspaces by $(V, W)$. The maps 
$$\beta_V\colon V\to \beta (B)\quad \text{and}\quad \gamma_W\colon W\to D$$
are isomorphisms and  $\beta(B)=\beta(V)$ by exactness. Fixing the basis 
$a_1,\ldots ,a_n$ of $A$ and  $d_1,\ldots ,d_l$ of $D$, we choose any basis $b_1,\ldots ,b_{m-n}$ of $V\subset B$ and define the basis 
$c_1,\ldots ,c_{m-n}$ of $\beta (B)\subset C$ by  $c_i=\beta_V(b_i)$, and  define the basis  $c_1',\ldots ,c_l'$ of $W$ by $\gamma_W(c_i')=d_i$.

In order to show that  the vector
$$
(\alpha(a_1)\wedge\ldots \wedge\alpha(a_n)\wedge b_1\wedge\ldots \wedge b_{m-n})\otimes ( c_1'\wedge\ldots \wedge c_l'\wedge c_1\wedge\ldots \wedge c_{m-n})^\ast
$$
does not depend on the choice of the basis $b_1,\ldots ,b_{m-n}$ of $V$,  we choose a second basis
 $\ov{b}_1,\ldots ,\ov{b}_{m-n}$ of $V$, so that there is a linear isomorphism
$\sigma\colon V\rightarrow V$ mapping one basis into the other by $\ov{b}=\sigma (b)$.
Hence 
$$\ov{b}_1\wedge \ldots \wedge \ov{b}_{m-n}=\det (\sigma )\cdot b_1\wedge \ldots \wedge b_{m-n},$$
where $\det(\sigma)$ is the usual determinant of the  linear map $\sigma$. 
The associated basis $\ov{c}_1,\ldots ,\ov{c}_{m-n}$ of $\beta (B)$  is then defined by 
$\ov{c}_i =\beta_V(\ov{b}_i)$. The two basis 
$c_i$ and $\ov{c}_i$ of $\beta (B)$ are therefore related by the isomorphism
$$
 \beta_V\circ \sigma\circ \beta_V^{-1}\colon \beta(B)\rightarrow\beta(B)
 $$
and consequently,
$$
\ov{c}_1\wedge\ldots \wedge \ov{c}_{m-n} =\det( \beta_V\circ \sigma\circ \beta_V^{-1})c_1\wedge\ldots \wedge c_{m-n}=\det(\sigma)c_1\wedge\ldots \wedge c_{m-n}.
$$
From this we obtain
\begin{equation*}
\begin{split}
&(\alpha(a_1)\wedge\ldots \wedge\alpha(a_n)\wedge\bar{b}_1\wedge\ldots \wedge \ov{b}_{m-n})\otimes ( c_1'\wedge\ldots \wedge c_l'\wedge \ov{c}_1\wedge\ldots \wedge \ov{c}_{m-n})^\ast\\
&\quad \quad =\det (\sigma) (\alpha(a_1)\wedge\ldots \wedge\alpha(a_n)\wedge b_1\wedge\ldots \wedge b_{m-n})\\
&\phantom{\quad \quad =}\otimes \dfrac{1}{\det (\sigma)}( c_1'\wedge\ldots \wedge c_l'\wedge c_1\wedge\ldots \wedge c_{m-n})^\ast,
\end{split}
\end{equation*}
so that the basis change has indeed no influence.

Next we replace the basis  $a_1,\ldots ,a_n$ of $A$  by the basis $\ov{a}_1,\ldots ,\ov{a}_n$ by  means of the linear isomorphism  $\sigma\colon A\rightarrow A$ and replace  the basis $d_1,\ldots ,d_l$ of $D$ 
by the basis $\bar{d}_1,\ldots ,\bar{d}_l$  by means of the isomorphism $\varepsilon\colon D\rightarrow D$. Then 
$$(\ov{a}_1\wedge \ldots \wedge \ov{a}_n)\wedge (\ov{d}_1\wedge \ldots \wedge \ov{d}_l)^\ast=
\dfrac{\det (\tau)}{\det (\varepsilon)}\cdot 
(a_1\wedge \ldots \wedge a_n)\wedge (d_1\wedge \ldots \wedge d_l)^\ast,$$
and the vector 
$(a_1\wedge \ldots \wedge a_n)\wedge (d_1\wedge \ldots \wedge d_l)^\ast$ is independent of the choices of the basis if and only if $\det (\tau)=\det (\varepsilon)$. So far, the definition of $\Phi_{\bm{E}}$ could  only depend on the choice 
 of the complements $(V,W)$.

So, we  assume that $(V,W)$ and $(V',W')$ are two choices. We assume that $a_1,\ldots ,a_n$ and $d_1,\ldots ,d_l$ are fixed bases for $A$ and $D$, respectively.
Now fixing  any basis $b_1,\ldots ,b_{m-n}$ of  $V$, we find vectors $q_1,\ldots ,q_{m-n}$ of  $\alpha (A)$ such  that $\ov{b}_1,\ldots ,\ov{b}_{m-n}$ with $\ov{b}_i=b_i+q_i$ is a basis of  $V'$. We define $c_i\in W$ by  $\gamma_W(c_i)=d_i$. Then we choose  $p_i\in\beta(C)$ so that $\ov{c}_i=c_i +p_i\in W'$
and $\gamma_{W'}(\ov{c}_i)=d_i$. Now we use the fact that for  $(V,W)$ and $(V',W')$ other choices do not matter.
Using that $q_j$ is a linear combination of 
$\alpha (a_1),\ldots ,\alpha (a_n)$ and $p_j$ a linear combination of $c_1',\ldots ,c_l'$, it follows that 
\begin{equation*}
\begin{split}
&(\alpha(a_1)\wedge\ldots \wedge\alpha(a_n)\wedge b_1\wedge\ldots \wedge b_{m-n})\otimes ( c_1'\wedge\ldots \wedge c_l'\wedge c_1\wedge\ldots \wedge c_{m-n})^\ast\\
&\quad =(\alpha(a_1)\wedge\ldots \wedge\alpha(a_n)\wedge\ov{b}_1\wedge\ldots \wedge \ov{b}_{m-n})\otimes ( c_1'\wedge\ldots \wedge c_l'\wedge \ov{c}_1\wedge\ldots \wedge \ov{c}_{m-n})^\ast,
\end{split}
\end{equation*}
which completes the proof that $\Phi_{\bm{E}}$ is well-defined.
\qed \end{proof}

\subsection{Proof of Proposition \ref{ooo6.6}}\label{oooo6.6}
\begin{proof}[Proof of Proposition \ref{ooo6.6}]
Let $P, Q\in\Pi_T$, $P\leq Q$,  and recall the exact sequence \begin{equation*}
{\bm{E}}_{(T,Q)}\colon\quad 0\rightarrow\ker(T)\xrightarrow{j^Q_T}\ker(QT)\xrightarrow{\Phi^Q_T} F/R(Q)\xrightarrow{\pi^Q_T}  \text{coker}(T)\rightarrow 0.
\end{equation*}
By definition of the determinant,  $h\in \det(T)$ is of  the form
$$
h=(a_1\wedge,\ldots ,\wedge a_n)\otimes ((d_1+R(T))\wedge\ldots \wedge (d_l+R(T)))^\ast,
$$
where $a_1,\ldots ,a_n$ is a basis of  $\ker(T)$ and $d_1+R(T),\ldots , d_l+R(T)$ is a basis of $F/R(T)$ satisfying $d_i\in R(I-Q)$.
The latter condition can be achieved since, by definition of $\Pi_T$,  $R(QT)=R(Q)$.
We choose  the  vectors $b_1,\ldots ,b_m$ in $\ker(QT)$
so that $a_1,\ldots ,a_n,b_1,\ldots ,b_m$ form a basis of  $\ker(QT)$ and take the linearly independent vectors  $T(b_1)+R(Q),\ldots ,T(b_m)+R(Q)$ in $F/R(Q)$.
By definition of the isomorphism $\gamma^Q_T\colon \det T\to \det (QT)$, 
\begin{equation*}
\begin{split}
&\gamma^Q_T(h)\\
&\quad =(a_1\wedge\ldots \wedge a_n\wedge b_1\wedge\ldots \wedge b_m)\otimes\\
&\quad ((d_1+R(Q))\wedge\ldots \wedge (d_l+R(Q))\wedge(T(b_1)+R(Q))\wedge\ldots \wedge (T(b_m)+R(Q)))^\ast .
\end{split}
\end{equation*}
We note that 
\begin{equation} 
T(a_1)=\ldots =T(a_n)=0\ \text{and}\ \ QT(b_1)=\ldots =QT(b_m)=0.
\end{equation}
Next we consider the exact sequence 
\begin{equation*}
0\rightarrow\ker(QT)\xrightarrow{j_{QT}^P}\ker(PT)\xrightarrow{\Phi^P_{QT}} F/R(P)\xrightarrow{\pi_{QT}^P}\text{coker}(QT)\rightarrow 0.
\end{equation*}
In order to  compute $\gamma^P_{QT}(\gamma^Q_T(h))$ we take as basis of  $\ker(QT)$ the vectors
$$
a_1,\ldots ,a_n,b_1,\ldots ,b_m,
$$
  and as basis of $F/R(Q)$ the vectors
   $$
   d_1+R(Q),\ldots ,d_l +R(Q),T(b_1)+R(Q),\ldots ,T(b_m)+R(Q).
   $$
Since  $\ker(QT)\subset\ker(PT)$,  we choose the vectors  $\wt{b}_1,\ldots ,\wt{b}_k$  in $\ker (PT)$ so that
$a_1,\ldots ,a_n,b_1,\ldots ,b_m,\wt{b}_1,\ldots ,\wt{b}_k$ form  a basis of  $\ker(PT)$.  We define 
$$
\alpha^P=a_1\wedge\ldots \wedge a_n\wedge b_1\wedge \ldots \wedge b_m\wedge \wt{b}_1\wedge\ldots \wedge \wt{b}_k,
$$
and note  that $\pi^P_{QT}(d_i+R(P))=(I-P)d_i +R(Q)=d_i+R(Q)$.  The vectors 
\begin{gather*}
d_1+R(P),\ldots ,d_l+R(P),\quad T(b_1)+R(P),\ldots ,T(b_m)+R(P),\\
QT(\wt{b}_1)+R(P),\ldots ,QT(\wt{b}_k)+R(P)
\end{gather*}
are a basis 
of $F/R(P)$ and we abbreviate  their  wedge product by $\beta^P$, so that 
$$
\gamma^P_{QT}(\gamma^Q_T(h)) = \alpha^P\otimes {(\beta^P)}^\ast.
$$
Recall that $d_i=(I-Q)d_i$ and $(I-Q)T(b_i)=T(b_i)$ for $i=1,\ldots ,m$. Abbreviating 
$[d]=d+R(P)$, the projection  $Q\colon F\rightarrow F$ induces the  projection
$\wt{Q}\colon F/R(P)\rightarrow F/R(P)$,  defined by $[d]\mapsto [Qd]$. Then $(I-\wt{Q})[d]=[(I-Q)d]$,  and therefore   
$$
F/R(P)= \ker(\wt{Q})\oplus \ker(I-\wt{Q}).
$$
The vectors 
 \begin{gather*}
 [d_1],\ldots ,[d_l],  [T(b_1)],\ldots ,[T(b_m)]\in \ker(\wt{Q}),
[QT(\wt{b}_1)],\ldots ,[QT(\wt{b}_k)]\in\ker(I-\wt{Q})
\end{gather*}
form  a basis for $F/R(P)$. By the  standard properties of the wedge product, 
\begin{equation*}
\begin{split}
\beta^P&=[d_1]\wedge\ldots \wedge[d_l]\wedge [T(b_1)]\wedge\ldots \wedge[T(b_m)]\wedge[QT(\wt{b}_1)]\wedge\ldots \wedge [QT(\wt{b}_k)]\\
&=[d_1]\wedge\ldots \wedge[d_l]\wedge [T(b_1)]\wedge\ldots \wedge[T(b_m)]\wedge[T(\wt{b}_1)]\wedge\ldots \wedge [T(\wt{b}_k)].
\end{split}
\end{equation*}
Hence we arrive for the composition 
$\gamma^P_{QT}\circ \gamma^Q_T$ at the formula
\begin{eqnarray*}
&&\gamma^P_{QT}(\gamma^Q_T(h))\\
& =& \alpha^P\otimes ([d_1]\wedge\ldots \wedge[d_l]\wedge [T(b_1)]\wedge\ldots \wedge[T(b_m)]\wedge[T(\wt{b}_1)]\wedge\ldots \wedge [T(\wt{b}_k)])^\ast.
\end{eqnarray*}
Next we compute $\gamma_T^P(h)$. In order to do so we start with the basis $a_1,\ldots ,a_n$ of  $\ker(T)$
and extend it to a basis of  $\ker(PT)$  by choosing   
 $b_1,\ldots ,b_m,\wt{b}_1,\ldots ,\wt{b}_k$. The wedge of all these vectors  is $\alpha^q$.
For $F/R(T)$ we have the basis 
$d_1+R(T),\ldots ,d_l +R(T)$.   Recall that $(I-Q)d_i=d_i$ so that $d_i\in \ker(\wt{Q})$. Further, 
\begin{equation*}
\begin{split}
\pi^P_T(d_i +R(P))&=(I-P) d_i+R(T)\\
&=(I-P)(I-Q)d_i+R(T)\\
&=(I-Q)d_i+R(T)\\
&=d_i+R(T).
\end{split}
\end{equation*}
Then we take  the basis of $F/R(P)$ formed  (with the previous convention) by the vectors
$$
[d_1],\ldots ,[d_l], [T(b_1)],\ldots ,[T(b_m)],[T(\wt{b}_1)],\ldots ,[T(\wt{b}_k)],
$$
and note that their  wedge product is equal to $\beta^P$. By definition,
$$
\gamma^P_T(h)=\alpha^P\otimes {(\beta^P)}^\ast,
$$
which agrees with  
$\gamma^P_{QT}\circ \gamma^Q_T(h)$.  The proof  of Proposition \ref{ooo6.6} complete.
\qed \end{proof}

%
%
%

\begin{partbacktext}
\part{Ep-Groupoids}
\noindent The basic notion in Part II is that of an ep-groupoid.  It generalizes the classical notion of an \'etale proper Lie groupoid
to the polyfold setup. From the view point of Fredholm theory it  allows to efficiently describe problems with local symmetries. 
Different ep-groupoids can be used to describe the same problem and it useful to understand the relationship between various descriptions.
This is discussed in the context of equivalences of categories and localization and leads to the introduction of the notion of a generalized
isomorphism.  An important topic is the study of branched ep$^+$-subgroupoids 
which naturally arise as solution spaces of sc-Fredholm section functors. Another topic important in applications is concerned with proper covering functors.
The discussion ends with the study of the various notions and their behavior under generalized isomorphisms and the introduction of a new kind of space called a
polyfold. In some sense it  generalizes the notion of an orbifold with boundary and corners to the world of a generalized differential geometry based on sc-smooth retracts. This basic ideas in the part  are  based on the paper \cite{HWZ3.5}

\end{partbacktext}

\chapter{Ep-Groupoids}
An ep-groupoid is the generalization of an \'etale proper Lie groupoid to the sc-smooth world.
Philosophically it should be viewed as a generalization of an atlas (in the theory of manifolds) in two different ways.
It generalizes smoothness to sc-smoothness and  introduces geometric structures beyond manifolds.  The development of the present part of the theory is quite parallel to parts of the theory Lie groupoids, see \cite{Mj}.
Many of the ideas about Lie groupoids and their importance go back to Haefliger, \cite{Haefliger1,Haefliger1.5, Haefliger2,Haefliger3}.

\section{Ep-Groupoids and Basic Properties}\label{sec_1_1}

We start with the notion of a groupoid following the expositions in  \cite{Mj}-\cite{MM} and recall that a {\bf small category}\index{Small category} is a category whose object class as well as the morphism class are sets.
 \begin{definition}[{\bf Groupoid}] \label{groupoid}\index{D- Groupoid}
A  groupoid $X=(X,\bm{X})$ is a small category consisting of the set $X$ of objects, denoted by $x$, and the set $\bm{X}$ of morphisms (arrows) which are assumed to be invertible, and the following five structure maps $(s, t, m, u, i)$.
\begin{itemize}
\item[(1)]\ {\bf Source and target maps $s$ and $t$}.\index{D- Target map}\index{D- Source map}
The source and the target maps $s, t\colon \bm{X}\to X$ assign to every morphism $\phi\colon x\to  y$ in  $\bm{X}$ its source $s(\phi )=x$ and its target $t(\phi )=y$. 
 \item[(2)]\ {\bf Multiplication map $m$}.\index{D- Multiplication map}
The associative multiplication (composition) map 
$$m\colon \bm{X}{_{s}\times_t}\bm{X}\to \bm{X},\quad m(\phi, \psi)=\phi\circ \psi$$
is defined on the fibered product 
$$ \bm{X}{_{s}\times_t}\bm{X}=\{(\phi, \psi)\in \bm{X}\times \bm{X}\, \vert \, s(\phi)=t(\psi)\}$$
so that if $\psi\colon x\to z$ and $\phi\colon z\to y$, then $\phi\circ \psi\colon x\to y$.
\item[(3)]\ {\bf Unit map $u$}.\, \index{D- $1$-map}
For every object $x\in X$ there exists the unit morphism $1_x\colon x\to x$ in $\bm{X}$, which is a $2$-sided unit for the composition, that is $\phi \circ 1_x=\phi$ and $1_x\circ \psi=\psi$ for all morphisms $\phi, \psi\in \bm{X}$ satisfying $s(\phi )=x=t(\psi)$. These unit morphisms together define the unit map 
$$u\colon X\to \bm{X},\quad u(x)=1_x.$$
\item[(4)]\ {\bf Inverse map $\iota$}.\index{D- Inverse map}
For every morphism $\phi\colon x\to y$ in $\bm{X}$ there exists the inverse morphism $\phi^{-1}\colon y\to x$ which is a $2$-sided inverse (isomorphism), that is $\phi\circ \phi^{-1}=1_y$ and 
$\phi^{-1}\circ \phi=1_x$. These inverses together define the inverse map 
$$\iota\colon \bm{X}\to \bm{X}, \quad \iota(\phi )=\phi^{-1}.$$
\end{itemize}
\qed
\end{definition}
Rather than writing $(X,\bm{X})$ for a groupoid we shall often write simply $X$.
With 
$\text{mor}(x, y)\subset \bm{X}$
we shall abbreviate the set of all morphisms $\phi\in \bm{X}$ satisfying $s(\phi)=x$ and $t(\phi )=y$.
The {\bf orbit space}\index{Orbit space} of the groupoid $(X, \bm{X})$,  denoted by 
$$
\abs{X}=X/\sim,
$$
is the quotient of the set of objects by the equivalence relation defined by $x\sim y$ if and only if there exists a morphism $\phi\colon x\to y$ between the two objects. The equivalence class will be denoted by 
$$\abs{x}:=\{y\in X\, \vert \, y\sim x\},$$
and we have the quotient map
$$\pi\colon X\to \abs{X},\quad \pi (x)=\abs{x}.$$
\begin{definition}[{\bf Isotropy group $G_x$}] \label{isotropy_group}
For fixed $x\in X$ we denote by $G_x$ the isotropy group (stabilizer group) of $x$\index{D- Isotropy group} defined by 
$$G_x=\{\phi \in \bm{X}\, \vert \, \phi\colon x\to x\}.$$
\qed
\end{definition}

The crucial concept of  ep-groupoids defined next can be viewed as  M-polyfold version of  \'etale and proper Lie groupoids discussed  f.e. in \cite{Mj} and \cite{MM}. The ideas in the classical case go back to A. Haefliger, \cite{Haefliger1,Haefliger1.5, Haefliger2,Haefliger3}.

\begin{definition}[{\bf Ep-groupoid}]\label{ep-groupoid_def}
An {\bf ep-groupoid}\index{D- Ep-groupoid} $X=(X,\bm{X})$ consists of a groupoid equipped  with   M-polyfold structures on the object set and the morphism set $\bm{X}$ having  the following properties.
\begin{itemize}
\item[(1)]\ ({\bf etal\'e}).\, \index{D- \'Etale}
  The source and target maps $s,t\colon \bm{X}\rightarrow X$ are  local sc-diffeomorphisms.
\item[(2)]\ The unit map $u$ and  the inverse  map $\iota$ are sc-smooth.
\item[(3)]\ ({\bf Properness}).\, \index{D- Ep-groupoid, properness} Every point $x\in X$ possesses an open neighborhood $V(x)$ of $x$ in $X$ such  that the map
$$t\colon s^{-1}\bigl(\cl_X(V(x)) \bigr)\rightarrow X$$
 is proper.
\end{itemize}
Since $s$ is a  local sc-diffeomorphism, the 
fibered product $\bm{X}{_{s}\times_t}\bm{X}$ is in a natural way a M-polyfold in view of Lemma 2.8 in \cite{HWZ3.5} and we require
\begin{itemize}
\item[(4)]\ The multiplication map $m\colon \bm{X}{_{s}\times_t}\bm{X}\rightarrow \bm{X}$
is sc-smooth.
\end{itemize}
In particular, all the structure maps $(s, t, m, u, i)$ are sc-smooth maps. 
\qed
\end{definition}

Since $s, t\colon \bm{X}\to X$ are local sc-diffeomorphism, every morphism $\phi\colon x\to y$ of an ep-groupoid has an extension to a local sc-diffeomorphism  $\wh{\phi}$ satisfying $\wh{\phi}(x)=y$ which is defined as follows. There exist open neighborhoods $\bm{U}(\phi)\subset \bm{X}$ and $U(x), U(y)\subset X$ such that $s\colon \bm{U}(\phi )\to U(x)$ and $t\colon \bm{U}(\phi)\to U(y)$ are sc-diffeomorphisms,  and the local sc-diffeomorphism 
$\wh{\phi}$\index{$\wh{\phi}$} is defined by 
$$
\wh{\phi}=t\circ (s\vert \bm{U}(\phi))^{-1}\co U(x)\to U(y).
$$
We also note that in the case that $\phi$ is smooth, i.e. $\phi\in \bm{X}_\infty$, there is an associated
sc-isomorphism 
$$
T\phi:T_{s(\phi)}X\rightarrow T_{t(\phi)}
$$
defined by $T\phi:= T\wh{\phi}(x)$\index{$T\phi$}. More generally we have $T\wh{\phi}:TU(s(\phi))\rightarrow TU(t(\phi))$. We note that $T\phi$ can be defined for $\phi$ on level at least $1$ and will give 
a bounded linear isomorphism on level $0$. Only in the case
that it is smooth $T\phi$ will be an sc-operator.
\begin{definition}\index{D- Tangent of $\phi$}
Associated to a smooth morphism $\phi$ we have the so-called {\bf associated local sc-diffeomorphism}
$\wh{\phi}$, as well as its {\bf tangent} $T\phi$.
\qed
\end{definition}

We point out that the above properness assumption in the definition of an ep-groupoid is different from the properness  notion in the classical Lie-groupoids which requires that the map $(s, t)\colon \bm{X}\to X\times X$ is proper. The latter  requirement is too weak for some of our  constructions.

The M-polyfolds $X$ and $\bm{X}$ are equipped with a filtration,
$$X=X_0\supset X_1\supset X_2\supset\ldots \supset X_\infty=\bigcap_{k\geq 0}X_k,$$
and a similar filtration on the morphism set $\bx$.
The set $X_\infty$ is dense in every space $X_k$. If $x$ belongs to the M-polyfold $X$, we denote by $\text{ml}(x)\in \N\cup \{\infty\}$\index{$\text{ml}(x)$}
the largest non-negative integer $m$ or $\infty$ such that $x\in X_m$. Here  ${\mathbb N}=\{0,1,2,3,..\}$.
Since in an ep-groupoid $(X, \bm{X})$ the source and the target maps are local sc-diffeomorphisms and therefore preserve by definition the levels, we conclude that 
$$\text{ml}(x)=\text{ml}(y)=\text{ml}(\phi)$$
for every morphism $\phi\colon x\to y$ in $\bm{X}.$
Consequently, the above filtration on the M-polyfold $X$ induces 
the filtration 
$$\abs{X}=\abs{X_0}\supset \abs{X_1}\supset\ldots \supset  \abs{X_\infty}=\bigcap_{k\geq 0}\abs{X_k}$$
on the orbit space $\abs{X}=X/\sim$.

For the ep-groupoid $(X, \bm{X})$ the degeneracy maps $d_X\colon X\to \N$  and $d_{\bm{X}}\colon \bm{X}\to \N$ are both defined. By Proposition \ref{newprop2.24}, these indices are invariant under local sc-diffeomorphisms and we deduce 
that 
$$d_X(s(\phi ))=d_{\bm{X}}(\phi )=d_{X}(t(\phi))$$
for every morphism $\phi \in \bm{X}$. In particular,
$$d_{X}(x)=d_X(x')\quad \text{if\quad $x\sim x'$},$$
so that  we obtain the non-degeneracy index on the orbit space $\abs{X}$ as follows.

\begin{definition}[{\bf Induced index $d_{\abs{X}}$}]
\index{D- Induced degeneracy index} 
The {\bf induced degeneracy index} $d_{\abs{X}}\colon \abs{X}\to \N$ is defined by 
$$d_{\abs{X}}(\abs{x})=d_X(x).$$
\qed
\end{definition}

\begin{definition}[{\bf Tame ep-groupoid}]\index{D- Tame ep-groupoid}
The ep-groupoid $(X, \bm{X})$ is called {\bf tame} if the object space $X$ is a tame M-polyfold as introduced in Definition \ref{def_tame_m-polyfold}.
\qed
\end{definition}

In a tame ep-groupoid $(X, \bm{X})$ not only $X$ is a tame M-polyfold but also $\bm{X}$ is a tame M-polyfold. This follows from the facts that $\bm{X}$ is locally sc-diffeomorphic by $s$ to $X$ and  the tameness condition is invariant under local sc-diffeomorphisms.

\begin{definition}[{\bf Category ${\mathcal E}{\mathcal P}$}]
\index{D- Category of ep-groupoids}
By ${\mathcal E}{\mathcal P}$\index{${\mathcal E}{\mathcal P}$} 
we denote the category whose objects are the ep-groupoids and whose morphisms $f\colon (X, \bm{X})\rightarrow (Y, \bm{Y})$
are the functors which are sc-smooth maps between the object spaces as well as the morphism spaces.
\qed
\end{definition}

Actually, it is sufficient to require that the functor $f$ is sc-smooth between the object spaces $X\to Y$ because the sc-smoothness of $f\colon \bm{X}\to \bm{Y}$ follows again from the properties of the source map $s\colon \bm{X}\to X$ as follows. 
Given $\phi\in \bm{X}$ and $f(\phi)=\psi\in \bm{Y}$ there exist open neighborhoods $\bm{U}(\phi)$, $\bm{U}(\psi)$, $U(s(\phi))$, and $U(s(\psi))$ satisfying 
$f(U(s(\phi)))\subset U(s(\psi))$ and $f(\bm{U}(\phi))\subset \bm{U}(\psi)$. Moreover,  $s\colon \bm{U}(\phi)\rightarrow U(s(\phi))$ and $s\colon \bm{U}(\psi)\rightarrow U(s(\psi))$ are sc-diffeomorphisms. Therefore, the claim follows from the formula  
$$
f(\gamma)= (s\vert \bm{U}(\psi))^{-1}(f(s(\gamma)))
$$
for all $\gamma\in \bm{U}(\phi).$

We equip the orbit space  $\abs{X}$ 
with the {\bf quotient topology}\index{Quotient topology}, i.e. the finest topology for which the projection  $\pi\colon X\rightarrow \abs{X}$, $x\rightarrow \abs{x}$, is  continuous. A subset $V$ of $\abs{X}$ is open if and only if the set of all objects $x\in X$ which belong to a class in $V$ is open, i.e.,  $\pi^{-1}(V)$ is open in $X$.

\begin{proposition}[{\bf Openness of $\pi$}]\label{q-open}\index{P- Openness of $\pi$}
Let $(X,\bm{X})$ be an ep-groupoid and assume that the orbit space $\abs{X}$  is equipped with the quotient topology. Then
the projection map $\pi\colon X\rightarrow \abs{X}$ is open. 
\end{proposition}
\begin{proof}
Let $U$ be an open subset of $X$. Then $\pi(U)\subset \abs{X}$ is open if and only if $\pi^{-1}(\pi(U))$ is open in $X$.
A point $y\in X$ belongs to $\pi^{-1}(\pi(U))$ if and only if there exist $x\in U$ and $\phi\in \bm{X}$
satisfying  $s(\phi)=x$ and $t(\phi)=y$. From the property (1) of an ep-groupoid we conclude  that there exist open neighborhoods
$\bm{V}(\phi)$, $V(x)$ and $V(y)$ such  that
$$
s\colon \bm{V}(\phi)\rightarrow V(x)\quad  \text{and}\quad t\colon \bm{V}(\phi)\rightarrow V(y)
$$
are sc-diffeomorphisms. We can choose  an open neighborhood $O(y)\subset V(y)$ such that $s\circ (t|V(\phi))^{-1}(O(y))$ is an open neighborhood of $x$ contained in $U$. This shows that $O(y)\subset \pi^{-1}(\pi(U))$, which implies, since $y$ was arbitrary in $\pi^{-1}(\pi(U))$, that the set $\pi^{-1}(\pi(U))$  is open in $X$.
\qed \end{proof}

\begin{proposition}[{\bf Countable neighborhood basis}]\index{P- Countable neighborhood basis}
Let $X$ be an ep-groupoid and let $\abs{X}$ be equipped with the quotient topology ${\mathcal T}$. Then every point $z\in \abs{X}$
has a countable basis of open neighborhoods $(V_k)$, which can be chosen  in a monotone way, 
$$
z\in V_k\subset V_{k-1}\subset \ldots \subset V_1 \quad \text{for all $k$.}
$$
\end{proposition}
\begin{proof}
For  $z\in \abs{X}$ we choose  $x\in X$ satisfying  $\pi(x)=z$. Since $X$ is metrizable, in view of 
Theorem \ref{X_m_paracompact}   we can fix a metric
and take the monotone sequence $U_k=B^X_{1/k}(x)$ of open balls. Then the monotone sequence $V_k=\pi(U_k)$ of open neighborhoods is the desired countable neighborhood basis of $z$.
\qed \end{proof}

\begin{proposition}[{\bf Local closure relation}]\label{Prop7.1.9}\index{P- Local closure relation}
Let $X$ be an ep-groupoid. Then for every $z\in \abs{X}$ and every representative $x\in X$ of $z$ there exists an open neighborhood
$U$ of $x$ in $X$ such  that for every subset $A\subset U$ the equality
$$
\pi(\cl_X(A))= \cl_{\abs{X}}(\pi(A))
$$
holds.
\end{proposition}
\begin{proof}
We take an open neighborhood $U$ of $x$ representing $z$ so that $\pi(x)=z$,  which has the property that the map
 $t\colon  s^{-1}(\cl_X(U))\rightarrow X$ is proper.
If $A$ is  a subset of $U$, then $\cl_X(A)\subset \cl_X(U)$ and obviously $\pi(\cl_X(A))\subset \cl_{|X|}(\pi(A))$. To verify  the reversed inclusion we let  $z'\in \cl_{|X|}(\pi(A))$. In this case every open neighborhood of $z'$ intersects $\pi(A)$
and we  have to show that $z'\in \pi(\cl_X(A))$. We choose  a representative $y\in X$ for $z'$ and take a monotone countable neighborhood basis $(V_k)$ for $z'$ constructed from such a basis $(O_k)$ for $y$, i.e. $V_k=\pi(O_k)$. By assumption,  we find  sequences
$(x_k)\subset A$ and $(y_k)\subset X$, $y_k\in  O_k$,  so that there exists a sequence  $(\phi_k)\subset \bm{X}$ satisfying $s(\phi_k)=x_k$ and $t(\phi_k)=y_k$. Clearly, 
$$
t(\phi_k)=y_k\rightarrow y\quad  \text{in $X$}.
$$
Using that $A\subset U$, the properness assumption implies,  after perhaps taking a subsequence,   the convergence 
$\phi_k\rightarrow \phi$ in $\bm{X}$. Hence $x_k=s(\phi_k)\rightarrow s(\phi)=:a'\in\cl_X(A)$ and $y_k=t(\phi_k)\rightarrow t(\phi)$.
Since $y_k\rightarrow y$,  we conclude  that $\phi\colon a'\rightarrow y$,  proving that $z'=\pi(y)=\pi(a')\in \pi(\cl_X(A))$. This completes the proof.
\qed \end{proof}

 We should point out that we refer to the ``arrows'' $\phi\colon x\rightarrow y$
usually as to morphisms, but note that by the invertibility assumption they are all isomorphisms. If we want to stress the latter we call the arrow $\phi$  the {\bf isomorphism} $\phi$  between $x$ and $y$.

\begin{lemma}\label{finite_x}
Assume that $x, x'\in X$ and that $U(x')$ is an open neighborhood of $x'$ such that the map $t\colon s^{-1}(\cl_X (U(x'))\to \bx$ is proper. Then there exist only finitely many morphisms  $\psi$ satisfying $s(\psi )=x$ and $t(\psi)\in U(x')$.
\end{lemma}
\begin{proof}
Arguing by contradiction there exists an infinite sequence of distinct morphisms $\psi_n : s(\psi_n)\to x$ where $s(\psi_n)\in U(x').$
The properness of the map 
$$
t: s^{-1}(\cl_X (U(x'))\to \bx
$$
implies  the convergence of a subsequence  of $(\psi_n)$ to the morphism 
$\psi : s(\psi)\to x$ where $s(\psi)\in s^{-1}(\cl_X (U(x'))$. The target map $t$ is a local sc-diffeomorphism and therefore we find open neighborhoods $V(x)$ and $\bm{V}(\psi)$ such that $t: \bm{V}(\psi)\to V(y)$ is a sc-diffeomorphism. Consequently, for large $n$ and $m\neq n$, $\psi_n$ and $\psi_m\in \bm{V}(\psi)$ and $t(\psi_n)=t(\psi_m)=x$ imply that $\psi_n=\psi_m$. This contradicts the assumption $\psi_n\neq \psi_m$ if $n\neq m$ and proves the lemma.
\qed \end{proof} 
As a consequence we have the following proposition.

\begin{proposition}[{\bf Finiteness of Isotropy Group}]\index{P- Finiteness of $G_x$}
Let $X$ be an ep-groupoid and $x$ an object in $X$. Then the isotropy group $G_x$ is finite.
\qed 
\end{proposition}

If $x$ and $y$ are isomorphic the associated isotropy groups are isomorphic as well.  In general,  there is no canonical identification of $G_x$ 
with  $G_y$. Any identification depends on an initial choice of isomorphism $\phi\colon x\rightarrow y$ and is given by
$$
G_x\rightarrow G_y,\quad  g\mapsto \phi\circ g\circ \phi^{-1}.
$$
A basic result  from  \cite{HWZ3.5}, Theorem 2.3, describes the following local structure of an ep-groupoid.
\begin{theorem}[{\bf Natural representation of $G_x$}] \label{x-local-x}\index{T- Basic structure of ep-groupoids}
Let $X$ be an ep-groupoid and $x$ an object with the isotropy  group $G=G_x$.  Then for every open neighborhood $V$ of $x$ we find an open neighborhood $U\subset V$ of $x$,  a group homomorphism $\Phi\colon G\rightarrow \text{Diff}_{\textrm{sc}}(U)$, $g\mapsto  \Phi (g)$,  and a sc-smooth map
$\Gamma\colon G\times U\rightarrow \bm{X}$\index{$\Gamma$} such that the following holds.
\begin{itemize}
\item[{\em (1)}]\ $\Gamma(g,x)=g$.
\item[{\em (2)}]\ $s(\Gamma(g,y))=y$ and $t(\Gamma(g,y))=\Phi (g)(y)$ for all $y\in U$ and $g\in G$.
\item[{\em (3)}]\ If $h\colon y\rightarrow z$ is a morphism between points in $U$, then there exists a unique element $g\in G$ satisfying $\Gamma(g,y)=h$, i.e., 
$$
\Gamma\colon G\times U\rightarrow \{\phi\in \bm{X}\, \vert \,  \text{$s(\phi)$ and $t(\phi)\in U$}\}
$$
is a bijection.
\end{itemize}
\qed
\end{theorem}
We shall refer to $(\Phi,\Gamma)$ as the data for the {\bf natural representation}\index{Natural representation} of $G_x$.
The proof  is a  slight modification of the proof in \cite{HWZ3.5} and will be given for the convenience of the reader in Appendix \ref{k-natural}.

In the following we shall often abbreviate the group action of the isotropy group $G=G_{x}$ on $U$ by 
$$g\ast y=\Phi (g)(y)\quad \text{if $g\in G$ and $y\in U$}.
$$
We now discuss informally some of the consequences of Theorem \ref{x-local-x}. To fill in the details the reader might consult the proof in  Appendix \ref{k-natural}. If $x\in X$, and if $U(x)=U\subset X$ is the open neighborhood 
of $x$ guaranteed by  Theorem \ref{x-local-x}, we introduce 
the set of morphisms
$$
\bm{U}(x)=\{\phi\in \bm{X}\, \vert \, \text{$s(\phi)$ and $t(\phi)\in U(x)$}\}.\index{$\bm{U}(x)$} 
$$
The pair  $(U(x),\bm{U}(x))$ is the full subcategory of $(X, \bm{X})$ associated  with $U(x)$. By the properness property of an ep-groupoid we can choose $U(x)$ so small that $\bm{U}(x)$  fits into an arbitrarily small open neighborhood of the finite group $G_x\subset \bm{X}$. 
Using that $s$ is a local sc-diffeomorphism, we choose the connected open neighborhood $U(x)\subset X$ so small that $\bm{U}(x)$ splits into $\# G_x$-many (disjoint) connected components
enumerated  by the group elements  $g\in G_x$ so that 
$$
\bm{U}(x)=\sqcup_{g\in G_x} \bm{U}(g)\index{$\bm{U}(g)$}
$$
and, moreover, the source and target maps 
$$
s,t\colon \bm{U}(g)\rightarrow U(x)
$$ 
are sc-diffeomorphisms.  We introduce the category 
$$
G_x\ltimes U(x):=(U(x), G_x\times U(x)),
$$
called the {\bf translation groupoid}\index{Translation groupoid} associated with $x$ and $U(x)$, 
and defined as follows. We abbreviate
$$U=U(x)\quad \text{and}\quad G=G_x.$$
The object set of the category $G\ltimes U=(U, G\times U)$ \index{$G\ltimes U$}  consists of the points in $U$, and the morphism set consists of the pairs $(g, y)\in G\times U$.  The space $G\times U$ has a natural M-polyfold structure as the finite disjoint  union of copies of the M-polyfold $U$.
The source and the target maps $s, t\colon G\times U\to U$ are defined as 
$$
s(g,y)=y\quad \text{and}\quad t(g,y)=\Phi (g)(y).
$$
The composition of the morphisms 
$(g,y)$ and $(h,z)\in G\times U$ satisfying $y=s(g,y)=t(h,z)=\Phi (h)(y)$  is defined by 
$$
(g,y)\circ (h,z)=(g\circ h,z)\in G\times U.
$$
Using the proof of Theorem \ref{x-local-x} one verifies that the category $G\ltimes U$ is an ep-groupoid. The natural functor 
$$
\bm{\Gamma}\colon G\ltimes U=(U, G\times U)\rightarrow  (X, \bm{X})
$$
between the two categories  maps, on object level,  $y\in U$ to $y\in X$ and, on the morphism level, the  morphism $(g,y)\in G\times U$ into the unique morphism   $\Gamma(g,y)\in \bm{X}$.  By 
property (3) of Theorem \ref{x-local-x}  the functor $\Gamma$ is {\bf full}\index{Full functor} and {\bf faithful}\index{Faithful functor} (fully faithful). 
{\bf Full} requires, that for two objects
$y$ and $z$ in $G\ltimes U$,  the functor $\bm{\Gamma}$ induces a surjection between the morphisms sets of morphisms $y\to z$,
$$
\text{mor}_{G\ltimes U}(y,z)\rightarrow \text{mor}_X(y,z).
$$
{\bf Faithful} means that these maps are injective. 
 The image of the functor $\bm{\Gamma}$ is the full subcategory $(U(x), \bm{U}(x))$ of 
 $(X,\bm{X})$. Hence $\bm{\Gamma}$ identifies the translation groupoid $G_x\ltimes U(x)$ with the category 
$(U(x),\bm{U}(x))$.

Fixing  $x\in X$, we choose as before a sufficiently small connected open neighborhood $U(x)$ in which we have the natural $G_x$-action $\Phi$, and the associated set $\bm{U}(x)$ of morphisms, which is the disjoint union of the connected open neighborhoods 
$\bm{U}(g)\subset \bm{X}$ for  $g\in G_x$. 
By property (3) of  Theorem \ref{x-local-x}, for a given morphism $\phi\in \bm{U}(x)$ there exists a unique $h\in G_x$ satisfying $\phi\in \bm{U}(h)$, and
$$
\phi=\Gamma (h, s(\phi))\quad \text{and}\quad t(\phi)=\Phi (h)(s(\phi)).
$$
Abbreviating  the group action of $G_x$ by 
$$
g\ast z:=\Phi (g)(z),\quad g\in G_x, z\in U(x),
$$
we define the map 
$$
\wh{\Phi}\colon G_x\to \text{Diff}_{\textrm{sc}}(\bm{U}(x))
$$
by 
$$
\wh{\Phi}(g)(\phi )=\Gamma \bigl(g\circ h\circ g^{-1},g\ast (s(\phi))\bigr)\in \bm{U}(g\circ h\circ g^{-1})
$$
if $\phi \in \bm{U}(h)$.  By property (2) of Theorem \ref{x-local-x}, 
$s(\wh{\Phi}(g)(\phi ))=g\ast s(\phi)$ and 
$$
t(\wh{\Phi}(g)(\phi ))=(g\circ h\circ g^{-1})\ast (g\ast s(\phi))=g\ast (h\ast(s(\phi)))=g\ast t(\phi).
$$
If $g'\in G_x$, we compute
\begin{equation*}
\begin{split}
\wh{\Phi}(g')\bigl(\wh{\Phi}(g)(\phi)\bigr)
&=\wh{\Phi}(g')(\Gamma (g\circ h\circ g^{-1},g\ast 
s(\phi)))\\
&=\Gamma \bigl(g'\circ (g\circ h\circ g^{-1})\circ (g')^{-1}, g'\ast(g\ast s(\phi)) \bigr)\\
&=\wh{\Phi}(g'\circ g)(\phi).
\end{split}
\end{equation*}

For $1_x\in G_x$ we obtain $\wh{\Phi}(1_x)(\phi)=\phi$
for all $\phi \in \bm{U}(x)$.   Therefore, the map $\wh{\Phi}\colon G_x\to \text{Diff}_{\textrm{sc}}(\bm{U}(x))$ 
is a group homomorphism satisfying 
$$\wh{\Phi}(g)(\bm{U}(h))=\bm{U}(g\circ h\circ g^{-1})$$ 
for all $g,h\in G_x$. Moreover, 
$$
s(\wh{\Phi}(g)(\phi))=\Phi (g)(s(\phi))\ \ \text{and}\ \ t(\wh{\Phi}(g)(\phi))= \Phi (g)(t(\phi)).
$$

To sum up the discussion, we can choose a small connected neighborhood $U(x)\subset X$ such that $\bm{U}(x)$  is the disjoint union of the connected open subsets
$\bm{U}(g)\subset \bm{X}$ for  $g\in G_x$,   and for $g,h\in G_x$ the sc-diffeomorphism
$\wh{\Phi}(g):\bm{U}(x)\rightarrow \bm{U}(x)$ induces the commutative diagrams
of sc-diffeomorphisms 
$$
\begin{array}{cc}
\begin{CD}
\bm{U}(h) @> s >> U(x)\\
@ V\wh{\Phi}(g)|\bm{U}(h) VV    @V \Phi(g) VV\\
\bm{U}(g\circ h\circ g^{-1}) @> s>> U(x)
\end{CD}
&\ \ \ \ \ \ 
\begin{CD}
\bm{U}(h) @> t >> U(x)\\
@ V\wh{\Phi}(g)|\bm{U}(h) VV    @V \Phi(g) VV\\
\bm{U}(g\circ h\circ g^{-1}) @> t>> U(x)
\end{CD}
\end{array}
$$
The following proposition extends the theorem about natural actions of isotropy groups.

\begin{proposition}\label{delaney}\index{P- Natural representation of $Y$}
 Let $(X, \bm{X})$ be an ep-groupoid and $Y=\{y_1, \ldots ,y_k\}$ a finite set of distinct objects in $X$ which are mutually isomorphic. 
 Denote by $Y=(Y, \bm{Y})$ the full subcategory associated with $Y$ whose  morphism set is $\bm{Y}$.  Then,  given a neighborhood $V$
 of $Y$ in $X$,  there exists a functor $f$, which associates with  the objects $y$ in $Y$  mutually disjoint connected open neighborhoods $U(y)$ contained in $V$, and with  a morphism $\psi\colon y\rightarrow y'$ in $\bm{Y}$  a sc-diffeomorphism
 $f (\psi)\colon U(y)\rightarrow U(y')$. The set of morphisms  $\bm{U}=\{\phi\in \bm{X}\, \vert \,  \text{$s(\phi)$ and $t(\phi)\in \bigcup_{y\in Y} U(y)$}\}$
decomposes into union of connected  components $\bm{U}(\psi)$, $\psi\in \bm{Y}$, and the following holds true.
 \begin{itemize}
 \item[{\em (1)}]\ The isotropy group $G_y$ for $y\in Y$ acts via the natural action on $U(y)$.
 \item[{\em (2)}]\ For $\psi\in \bm{Y}$ the maps $s:\bm{U}(\psi)\rightarrow U(y)$ and $t:\bm{U}(\psi)\rightarrow U(y)$ are sc-diffeo\-mor\-phisms.
 \item[{\em (3)}]\ If  $\psi\colon y\rightarrow y'$ and  $z\in U(y)$, then  
 $$
f(\psi) (z) =t\circ (s|\bm{U}(\psi))^{-1}(z).
 $$
 \end{itemize}
 \end{proposition}
 \begin{proof}
 We prove the result by induction on the number of objects  of the set $Y$. We begin with the case $Y=\{y_1\}$. Morphisms $y_1\to y_1$ are elements of the isotropy group $G_{y_1}$. Taking a $G_{y_1}$-invariant open connected neighborhood $U(y_1)$ of $y_1$ contained in $V$, every morphism $\phi\in \bm{U} (y_1, y_1)$ is of the form $\phi=\Gamma( g, s(\phi))$ for a unique $g\in G_x$. The sets $\bm{U} (g):=\{\varphi =\Gamma (g, z)\, \vert \, z\in U(y_1)\}$ are mutually disjoint  connected open neighborhoods of $g$ and 
 $$
\bm{U}=\bm{U}(y_1)= \bm{U} (y_1, y_1)=\bigcup_{g\in G_{y_1}}\bm{U} (g).
 $$
Moreover, the source and target maps 
$$
s, t: \bm{U} (g)\to U(y_1)
$$
 are sc-diffeomorphims for every $g\in G_{y_1}$.   This is an immediate consequence of Theorem \ref{x-local-x}.
 
Considering  the case of $Y=\{y_1, y_2\}$, we fix  a morphism  $\gamma_{1}: y_1\to y_2$. We take open neighborhoods $U'(y_1)$ and $U'(y_2)$, both contained in $V$ and invariant with respect to the actions of the isotropy groups $G_{y_1}$ and $G_{y_2}$. Using that the source and target maps $s, t$ are local sc-diffeomorphisms, we find a $G_{y_1}$-invariant open neighborhood $V(y_1)\subset U'(y_1)$,  and open neighborhood $V(y_2)\subset  U'(y_2)$, and an open neighborhood $\bm{V}(\gamma_1)$ such that $s: \bm{V}(\gamma_1)\to V(y_1)$ and $t: \bm{V}(\gamma_1)\to V(y_2)$ are sc-diffeomorphisms. Hence  the  composition $\wh{\gamma}_1=t\circ (s\vert \bm{V}(\gamma_1))^{-1}: V(y_1)\to V(y_2)$ 
is also  a sc-diffeomorphism.  Choosing a $G_{y_2}$-invariant open neighborhood 
$W(y_2)$ contained in $V(y_2)$, we find a $G_{y_1}$-invariant connected open neighborhood $U(y_1)$ contained in $V(y_1)$ such that $\wh{\gamma_1}(U(y_1))\subset W(y_2)$.  Denoting  $U(y_2)=\wh{\gamma_1}(U(y_1))$, we claim that $U(y_2)$ is a $G_{y_2}$-invariant neighborhood of $y_2$.  To see this, let $z\in U(y_2)$ and $z'=g\ast z$ where $g\in G_{y_2}$. 
Then $z'\in W(y_2)$ and,  by means of sc-diffeomorphism $\wh{\gamma}_1$, we find points $u\in U(y_1)$, $u'\in V(y_1)$,  and morphisms
$\phi: u\to z$ and $\phi: u'\to z'$. This implies that there is a morphism $u\to u'$ of the form $\Gamma (g', u)$ for some $g'\in G_{y_1}$.
Since $U(y_1)$ is $G_{y_1}$-invariant and $u\in U(y_1)$,  we conclude that $u'\in U(y_1)$, so that   $z'=\wh{\gamma}_1(u')\in U(y_2)$. This shows that $U(y_2)$ is $G_{y_2}$-invariant. It is also  which implies that $z'\in U(y_2)$ since $\wh{\gamma}_1$ maps $U(y_1)$ onto $U(y_2)$.
The source and target maps $s: \bm{U} (\gamma_1)\to U(y_1)$ and $t: \bm{U} (\gamma_1)\to U(y_2)$ are sc-diffeomorphisms, and hence $\wh{\gamma}_1=t\circ (s\vert \bm{U} (\gamma_1))^{-1}: U(y_1)\to U(y_2)$ is also a sc-diffeomorphism. 
   
Any other morphism from $y_1\to y_2$ is of the form $\gamma_g=g\circ \gamma_1$ for some  $g\in G_{y_2}$.  Fixing $g\in G_{y_1}$, we consider the map   $\Phi_g: \bm{U} (\gamma_1)\to \bx$, defined by
$$\Phi_g(\phi)=\Gamma (g, t(\phi))\circ \phi.$$
Then  $\Phi_g$ maps  $\gamma_1$ onto $\gamma_g$ and is a  sc-diffeomorphism onto its image $\bm{U} (\gamma_g):=\Psi_g(\bm{U} (\gamma_1))$.  The images  $\bm{U} (\psi_g)$ are connected and pairwise disjoint neighborhoods of $\gamma_g$
since for any $\phi$ morphism satisfying $s(\phi)\in U(y_1)$ and $t(\phi)\in U(y_2)$, there $\phi=\Gamma (g, s(\phi))\circ \psi$ for uniquely determined  $g\in G_{y_2}$ and $\psi\in \bm{U} (\gamma_1)$. Therefore, the set of morphisms $\bm{U} (y_1, y_2)$ is the union of connected components $\bm{U} (\gamma_g)$ for all $g\in G_{y_2}$.  Also  the source and target maps $s: \bm{U} (\gamma_g)\to U(y_1)$ and $t: \bm{U} (\gamma_g)\to U(y_2)$ are sc-diffeomorhisms.

If $\phi: y_2\to y_1$, then the morphism $\phi^{-1}: y_1\to y_2$ is one of the morphisms $\gamma_g$ for some $g\in G_{y_2}$. Then we define 
$\bm{U} (\phi)=i(\bm{U} (\phi^{-1})),$ where $i:\bx \to \bx$ is the inverse map $i(\phi)=\phi^{-1}$. The source and target maps $s: \bm{U} (\phi)\to U(y_2)$ and $t: \bm{U} (\phi)\to U(y_1)$ are given by 
$$s\vert  \bm{U} (U(\phi))=(t\vert \bm{U} (\phi^{-1}))\circ i\quad \text{and}\quad 
t\vert  \bm{U}(\phi)\bm{U}=(s\vert \bm{U} (\phi^{-1}))\circ i,$$
and so they are sc-diffeomorphisms. This finishes the case of $Y=\{y_1,y_2\}$. 

Let $k\geq 3$ and assume that the statement holds for any set $Y'$ of $k$ distinct objects which are mutually isomorphic. Considering  the set $Y=\{y_1,\ldots, y_{k+1}\}$ of  distinct mutually isomorphic points  $y_i$,  let $V$ be an open neighborhood of $Y$. Let  $Y'=\{y_1,\ldots, y_k\}$ and $\bm{Y}'$ be the set of all morphisms between points in $Y$.  By assumption,  every point $y_i\in Y'$ possesses  a $G_{y_i}$-invariant connected and open neighborhood $U'(y_i)$ contained in $V$ and every morphism $\phi: y_i\to y_j$ possesses  a connected open neighborhood $\bm{U}' (\phi)$ for which the source and target maps  $s: \bm{U}' (\phi)\to U'(y_i)$ and $t: \bm{U}' (\phi)\to U'(y_j)$ are sc-diffeomorphisms. The composition 
$\wh{\gamma}_j: U'(y_i)\to U'(y_j)$ is a sc-diffeomorphism and 
the set of morphisms from $\bigcup_{y\in Y'}U(y)$ and $\bigcup_{y'\in Y'}U'(y')$ is the union of connected components $\bm{U}' (\phi)$ for all  $\phi\in \bm{Y}'$. 

We fix a $G_{y_{k+1}}$-invariant neighborhood $U'(y_{k+1})$ of $y_{k+1}$ contained in $V$ and for every $1\leq j\leq k$ a morphism $\delta_j: y_{k+1}\to  y_j$.  Similarly as in the case $k=2$, we find a $G_{y_{k+1}}$-invariant connected open neighborhood $U(y_{k+1})\subset U'(y_{k+1})$, $G_{y_j}$-invariant connected open neighborhoods $U(y_j)\subset U'(y_j)$, and connected open neighborhoods $\bm{U} (\delta_j)$ such that 
$s: \bm{U} (\delta_j)\to U(y_{k+1})$, $y: \bm{U} (\delta_j)\to U(y_{j})$, and 
$\wh{\delta}_j: U(y_{k+1})\to U(y_j)$ are sc-diffeomorphisms.

Now, as in the case $k=1$, one constructs neighborhoods $\bm{U} (\phi)$ for all other morphisms $\phi$  satisfying $s(\phi)=y_{k+1}$ and $t(\phi)\in Y'$ as well as their inverses. 

Still we have to adjust neighborhoods $\bm{U} '(\phi)$ for morphisms $\phi: y_i\to y_j$ for $i, j\leq k$.  By assumption,  if $\gamma: y_i\to y_j$, then there are neighborhoods $U'(y_i)$, $U'(y_j)$, and  $\bm{U}' (\gamma)$ such that 
$s: \bm{U} '(\phi)\to U'(y_i)$, $t: \bm{U} '(\phi)\to U'(y_j)$, and
$\wh{\gamma}=t\circ (s\vert \bm{U}' (\gamma))^{-1}: U'(y_i)\to U'(y_j)$  are  sc-diffeomorphisms. 

We define $\bm{U} (\gamma):=(s\vert \bm{U}' (\gamma))^{-1}(U(y_i))$ and denote $W(y_j)=t(\bm{U} (\gamma))$. The set $W(y_j)$ is  $G_{y_j}$-invariant and we claim that $W(y_j)=U(y_j)$.  To verify that $U(y_j)\subset W(y_j)$,  we take $z\in U(y_j)$. By means of the sc-diffeomorphisms $\wh{\delta}_i^{-1}: U'(y_i)\to U(y_{k+1})$ and $\wh{\delta}_j: U(y_{k+1})\to U(y_j)$ we find a point $u\in U(y_i)$ and morphisms $\phi\in \bm{U} (\delta_i^{-1})$, $\phi'\in \bm{U} (\delta_j)$ such that 
$\phi'\circ \phi: u\to z$. By means of the sc-diffeomorphism $\wh{\gamma}: U'(y_i)\to U(y_j)$, we find a morphism $\psi: u\to z'$ where $z'=t(\psi)\in W(y_j)$ since $u\in U(y_i).$  Hence there is a morphism  $z'\to z$  which, because $U'(y_j)$ is $G_{y_j}$-invariant, is of the form  
$\Gamma (g, z'): z'\to z$ for some $g\in G_{y_i}$. But $z'\in W(y_j)$ and $W(y_j)$ is $G_{y_j}$-invariant, and we conclude that $z\in W(y_j)$, as claimed. Similar arguments show that $W(y_j)\subset U(y_j)$, and we have proved that $t(\bm{U} (\gamma))=U(y_j)$ for every $1\leq j\leq k$.

The constructed neighborhoods $U(y)$,  for $y\in Y$, and the neighborhoods $\bm{U} (\phi)$ of morphisms between points in $Y$, and the associated source and target maps have the desired properties. This concludes the proof of the inductive step and the proof of the proposition.
\qed \end{proof}

\begin{definition}\index{D- Natural representation of $Y$}
Let $X$ be an ep-groupoid and $Y$ the full subcategory associated with  a finite set of isomorphic objects.
A functor $f$ as constructed in Proposition \ref{delaney} with its stated properties is called 
 the {\bf natural representation of $Y$}  on neighborhoods contained in $V$ or,  alternatively,  on a small neighborhood of $Y$.
 \qed
 \end{definition}

\begin{proposition}\label{go-go}\index{P- Existence of local uniformizers}
Let $X$ be an ep-groupoid and $z\in \abs{X}$. We fix  $x\in X$  satisfying  $z=\abs{x}$ and denote the isotropy group of $x$ by $G_x$ and 
take an open neighborhood $U(x)$ of $x$ on which we have the natural $G_x$-action $\Phi$, which we abbreviate as $g\ast y:=\Phi (g)(y)$ for $g\in G_x$ and $y\in U(x)$.
Then
the map $\tau:{_{G_x}\backslash} U(x)\rightarrow \abs{U(x)}=\pi (U(x)),\, G_x\ast y\rightarrow \abs{y}$ is a homeomorphism. In  other words,  the obvious functor
$$
G_x\ltimes U(x)\rightarrow X,
$$
which,  on the object level,  is the inclusion map and,  on the morphism level,  the map
$$
(g,y)\rightarrow \Gamma(g,y),
$$
is a fully faithful functor, which on the orbit spaces defines a homeomorphism between the orbit spaces $|U(x)|_{G_x\ltimes U(x)}={_{G_x}\backslash} U(x)$
and $|U(x)|=\pi(U(x))$.
\end{proposition}
\begin{proof}
Fixing a point $x\in X$, we take an open neighborhood  $U(x)$ of $x$ as guaranteed by Theorem \ref{x-local-x}\  
 and equip ${_{G_x}\backslash} U(x)$ with the quotient topology.
Then the map 
\begin{equation}\label{rok}
\tau:{_{G_x}\backslash} U(x)\rightarrow |U(x)|,\quad G_x\ast y\rightarrow \abs{y}
\end{equation}
is well-defined and  surjective by construction. It is also injective. Indeed, if  $\abs{y}=\abs{y'}$ for $y,y'\in U(x)$, then there exists a morphism $\phi:y\rightarrow y'$ in $\bm{X}$.
By point (3) of Theorem \ref{x-local-x} there exists $g\in G_x$ satisfying $\Gamma(g,y)=\phi$ and  $s(\phi)=y$ and,  moreover,  
$$
y'=t(\phi)=t(\Gamma(g,y))=\Phi (g)(y)=g\ast y.
$$
Hence $y'\in G_x\ast y$ and  our map $\tau$ is a bijection.

 The target $|U(x)|$ is open in $|X|$ by  Proposition \ref{q-open}.
The map 
$$
\pi\vert U(x)\colon U(x)\rightarrow |U(x)|
$$
 is continuous and the preimage of a point $\pi (y)$ is precisely the orbit $G_x\ast y$. Denoting the quotient map
$U(x)\rightarrow {_{G_x}\backslash} U(x)$ by $\sigma$,  we have the relationship
$$
\tau\circ\sigma=\pi\vert U(x).
$$
We note that $\sigma$ is open and continuous. Continuity is implied by the definition of the quotient topology
and the openness follows since $G_x$ acts by sc-diffeomorphisms.  If $O\subset |U(x)|$ is open we find that
$(\pi|U(x))^{-1}(O)= \sigma^{-1}(\tau^{-1}(O))$ is open, which,  by the definition of the quotient topology on 
${_{G_x}\backslash} U(x)$,  precisely means that $\tau^{-1}(O)$ is open. Hence $\tau$ is continuous.
If $Q\subset {_{G_x}\backslash} U(x)$ is open, then by definition,  $\sigma^{-1}(Q)$ is open. 
From 
$$
(\pi\vert U(x))^{-1}(\tau(Q))=\sigma^{-1}(Q)
$$
we conclude,  by definition of the quotient topology on $|X|$,  that $\tau(Q)$ is open. In summary $\tau$ is a homeomorphism.
 This completes the proof of Proposition \ref{go-go}.
\qed \end{proof}
\begin{definition}[{\bf Local uniformizer}]\index{D- Local uniformizer}
Let $X$ be an ep-groupoid and $x$ an object in $X$. A {\bf local uniformizer around $x$} is a  fully faithful
embedding
$$
\Psi_x\co G_x\ltimes U(x)\rightarrow X,
$$
where $U(x)$ is an open neighborhood of $x$ equipped with the natural $G_x$-action, such that $\Psi_x$ induces a homeomorphism
between $|U(x)|_{G_x\ltimes U(x)}$ and  the  open neighborhood $\abs{U(x)}=\pi (U(x))$ of $\pi(x)$ in $\abs{X}$.
\qed
\end{definition}
There exists a local uniformizer around every object $x$ as  Proposition \ref{go-go} shows. We also note that a 
a local uniformizer is injective on objects as well as morphisms.
Next we introduce the notion of an ep-subgroupoid. 
\begin{definition}\label{Xep-subgroupoid}
Let $X$ be an ep-groupoid. An {\bf ep-subgroupoid}\index{D- Ep-subgroupoid} $A$ of $X$ consists of the full subcategory associated
with  a subset $A$ of $X$ which has the following properties.
\begin{itemize}
\item[(1)]\ $A$  is a sub-M-polyfold of the object space $X$.
\item[(2)]\ $A$ is saturated, i.e. $\pi^{-1}(\pi(A))=A$.
\end{itemize}
\qed
\end{definition}
Let us show that the full-subcategory $A$ has in a natural way the structure of an ep-groupoid.
From the definition it follows, as proved in Proposition \ref{sc_structure_sub_M_polyfold}, that the object set $A$ has an induced M-polyfold structure 
for which the following holds.
\begin{itemize}
\item[(1)]\ The inclusion map $i:A\rightarrow X$ is sc-smooth and a homeomorphism onto its image. 
\item[(2)]\ For every $a\in A$ and every sc-smooth retraction $r\colon V\rightarrow V$ of an open neighborhood $V$ of $a$ in $X$ satisfying 
$r(V)=A\cap V$,  the map $i^{-1}\circ r:V\rightarrow A$ is sc-smooth.
\item[(3)]\ At a  smooth point $a\in A$ the tangent space $T_aA$ has a sc-complement in $T_aX$.
\item[(4)]\ If $W\subset X$ is open and  equipped with a sc-smooth retraction 
$r:W\rightarrow W$ satisfying $r(W)=A\cap W$, then the induced map $r:W\rightarrow A$ is sc-smooth (see (2)). Moreover, 
at a smooth point $a\in A$ the tangent space $T_aA$  is equal to $(Tr)(a)T_aX$.
\end{itemize}
In a next step we shall show that the set of morphisms 
$$
\bm{A}=\{\phi\in \bm{X}\, \vert \, \text{$s(\phi)$ and $t(\phi)\in A$}\}
$$
 has in  a natural way the structure of a M-polyfold.
Choosing  $\phi\in \bm{A}$, we find open neighborhoods $\bm{U}(\phi)\subset \bm{X}$ and $U(a)\subset X$, where $a=s(\phi)$, such that
$s\colon \bm{U}(\phi)\rightarrow U(a)$ is a sc-diffeomorphism. By taking a possibly smaller neighborhood $U(a)$ and adjusting $\bm{U}(\phi)$, we may assume that
we have a sc-smooth retraction $r: U(a)\rightarrow U(a)$ satisfying  $r(U(a))=U(a)\cap A$.  We  introduce  the sc-smooth retraction
$$
\bm{r}\colon \bm{U}(\phi)\rightarrow \bm{U}(\phi), \quad \bm{r}(\psi)=(s\vert \bm{U}(\phi))^{-1}\circ r \circ s(\psi).
$$
By construction $s(\bm{r}(\psi))= r(s(\psi))\in A$ and $t(\bm{r}(\psi))\in U(a)$, so that
$\bm{r}(\psi)\in \bm{U}(a)$.
Moreover, since $A$ is saturated we infer also that $\bm{r}(\psi)\in \bm{A}$, implying
$$
\bm{r}(\bm{U}(\phi))=\bm{U}(\phi)\cap\bm{A}.
$$
This shows that $\bm{A}$ is a sub-M-polyfold of $\bm{X}$ and the properties for $A$ listed above  also hold for $\bm{A}$.
From the construction of the sub-M-polyfold structures on $A$ and $\bm{A}$ it follows that the maps $s,t: \bm{A}\rightarrow A$ are sc-smooth.
 Moreover, if  $\phi\in A$ and $a=s(\phi)$, the sc-diffeomorphism $s: \bm{U}(\phi)\rightarrow U(a)$ between open subsets 
of $\bm{X}$ and $X$, respectively, maps the intersection $ \bm{U}(\phi)\cap \bm{A}$ bijectively onto $U(a)\cap A$. This shows that the source map $s: \bm{A}\rightarrow A$ is a  local sc-diffeomorphisms. The same holds for the target map $t$. The unit map, the inversion map, and the multiplication map are restrictions of sc-smooth maps
and therefore sc-smooth. The properness assumption also holds true. Indeed, if $a\in A$, we take an open neighborhood $U(a)$ in $X$ equipped with the natural $G_a$-action, such  that
$t: s^{-1}(\cl_X(U(a)))\rightarrow X$ is proper. Let us show that the open neighborhood $U(a)\cap A$ of $a$ in $A$ has the desired properties. Since $A$ is saturated, 
the set $U(a)\cap A$ is invariant under $G_a$-action.  In order to verify that  the map 
$$
t: s^{-1}(\cl_A(U(a)\cap A))\rightarrow A
$$
is proper we take a  sequence $(\phi_k)\subset \bm{A}$ satisfying $t(\phi_k)\rightarrow c$ in $A$ and $s(\phi_k)\in \cl_A(U(a)\cap A)$.  Then  $s(\phi_k)\in \cl_X(U(a))$ and hence we may assume, without loss of generality, the convergence  $\phi_k\rightarrow \phi$ in $\bm{X}$. It follows that  $t(\phi)=c$ and $s(\phi)\in \cl_X(U(a))$.
Since $A$ is saturated, we conclude from $t(\phi)\in A$ that also $s(\phi)\in A$. In particular,  $\phi_k\rightarrow \phi$ in $\bm{A}$. From  $s(\phi_k)\in \cl_A(U(a)\cap A)$ it follows that $s(\phi)\in \cl_A(U(a)\cap A)$.
Summarizing the argument we have proved the following result.
\begin{proposition}\index{P- Ep-subgroupoids}
Let $X$ be an ep-groupoid and $A$ a saturated subset of the object space $X$ such  that $A$ is also a sub-M-polyfold. Then the full subcategory associated with  $A$,  also denoted by
$A$,  has in a natural way the structure of an ep-groupoid.
In particular,  $A$ and $\bm{A}=\{\phi\in \bm{X}\ |\ s(\phi),t(\phi)\in A\}$ are sub-M-polyfolds of $X$ and $\bm{X}$, respectively and all the structure maps are the restrictions
of the structure maps of $X$ and $\bm{X}$, respectively.
\qed
\end{proposition}

\section{Effective and Reduced Ep-Groupoids}\label{subsec-Effective_Reduced}
In the previous section we have seen that the local morphism structure of an ep-groupoid around an object $x$ is given
by a group homomorphism $G_x\rightarrow\text{Diff}_{\textrm{sc}}(U(x))$.  We can use these local structures
to distinguish a class of ep-groupoids which is useful  for the applications. 
\begin{definition}[{\bf Effective ep-groupoid}]\index{D- Effective ep-groupoid}
An ep-groupoid $X$ is said to be {\bf effective}\index{D- Effective} provided for every $x\in X$ and $g\in G_x\setminus\{1_x\}$, given an open neighborhood 
$\bm{U}(g)\subset \bm{X}$, there exists  a morphism $\psi\in \bm{U}(g)$ satisfying $s(\psi)\neq t(\psi)$. An object $x\in X$ is said to be {\bf regular}\index{D- Regular point} provided $G_{x}=\{1_x\}$.
\qed
\end{definition}
If $\phi\in \bm{X}$ and $s(\phi)$ is a regular object, the same has to be true for $t(\phi)$. \\[1ex]
We define the {\bf regular part}\index{Regular part of an ep-groupoid} $X_{\textrm{reg}}$ of an ep-groupoid as the subset
$$
X_{\textrm{reg}}=\{x\in X\, \vert \, G_x=\{1_x\}\}
$$ 
of the object space $X$.
\begin{lemma}[{\bf Openness of regular set}]\label{ooopen}\index{L- Openness of regular set}
For an ep-groupoid $X$ the regular part $X_{\textrm{reg}} $ is open.
\end{lemma}
\begin{proof}
If  $x\in X_{\text{reg}} $, then $G_x=\{1_x\}$ and we find an open neighborhood $U(x)$ equipped with the natural $G_x$-action $\Phi$. Every morphism $\phi\colon a\rightarrow b$
between  $a,b\in U(x)$ has the form $\phi=\Gamma(1_x,a)$ so that $s(\phi)=a$ and $b=\Phi (1_x)(a)=a$. In particular, there is exactly one morphism
$a\rightarrow a$ and therefore  $\phi=1_a$.  This proves that the isotropy group $G_a$ for every $a\in U(x)$ is equal to $G_a=\{1_a\}$ and therefore $U(x)\subset X_{\textrm{reg}}$.
\qed \end{proof}
Examples show that the regular set can be empty.

\begin{lemma}\label{klsonja}\index{L- Characterization of effectiveness}
For  an ep-groupoid $X$ the following statements are equivalent.
\begin{itemize}
\item[{\em (1)}]\ $X$ is effective.
\item[{\em (2)}]\ For every object $x\in X$ and every open neighborhood $V(x)$  there exists an open neighborhood $U(x)\subset V(x)$ equipped with the natural $G_x$-action such that the 
group homomorphism $\Phi \colon G_x\rightarrow \text{Diff}_{\textrm{sc}}(U(x))$ is injective.
\end{itemize}
\end{lemma}
\begin{proof}
We assume that (1) holds and choose an open neighborhood  $U(x)\subset V(x)$ of $x$ equipped with the natural $G_x$-action  $\Phi$. Every morphism $\phi\in \bm{X}$ satisfying $s(\phi),t(\phi)\in U(x)$ is of the form 
$$
\phi=\Gamma(g,s(\phi))
$$
for a uniquely determined $g\in G_x$, and
$$
t(\phi) =\Phi (g)(s(\phi)).
$$
Arguing by contradiction we  assume that there exists $h\in G_x\setminus\{1_x\}$ for which  $\Phi (h)=\mathbbm{1}$.
Since $X$ is effective,  we find a converging sequence $(\phi_k)\subset \bm{X}$ satisfying  $\phi_k\rightarrow h$ and $s(\phi_k)\neq t(\phi_k)$.
For $k$ large we must have $s(\phi_k),t(\phi_k)\in U(x)$ and hence there exist  uniquely determined elements $(h_k)\subset G_x$ such that 
$$
\phi_k=\Gamma(h_k,s(\phi_k))\quad  \text{and}\quad  t(\phi_k)=\Phi (h_k)(s(\phi_k))\neq s(\phi_k).
$$
Clearly,   $h_k\neq h$, since otherwise $s(\phi_k)=t(\phi_k)$. Since $G_x$ is finite, after perhaps taking a subsequence,  we may assume for all $k$ that
$h_k=g\neq h$. From $\Phi (g)(s(\phi_k))\neq s(\phi_k)$ we deduce that
\begin{equation}\label{nancy}
\Phi (g)\neq \mathbbm{1}.
\end{equation}
Now we can pass to the limit and find, in view of $s(h)=t(h)=x$, that 
$$
h=\Gamma(g,x)\quad  \text{and}\quad \Phi (g)(x)=x.
$$
However, it is one of the properties of $\Gamma$ that $\Gamma(h,x)=h$ for all $h\in G_x$, which implies $h=g$, contradicting $h\neq g$. We have proved that (1) implies (2).

Next we assume that (2) holds. Hence for given $x\in X$ we  find an open neighborhood $U(x)$ equipped with the natural $G_x$-action such  that $\Phi\colon G_x\rightarrow \text{Diff}_{\textrm{sc}}(U(x))$ 
is injective. Assume that $g\in G_x\setminus\{1_x\}$ is given. Then $\Phi (g)$ is not the identity. Hence we find $a\in U(x)$ satisfying 
$\Phi (g)(a)\neq a$ and define the morphism $\phi=\Gamma(g,a)$. By construction,  $s(\phi)=a\neq \Phi (g)(a)=t(\phi)$.
By taking smaller and smaller neighborhoods $U(x)$ and adapting the corresponding restrictions of $\Gamma$,  we find a sequence $(\phi_k)$  satisfying $s(\phi_k)\neq t(\phi_k)$ such that $s(\phi_k),t(\phi_k)\rightarrow x$ and $\phi_k\rightarrow g$.  Hence (1) is verified and the proof of Lemma \ref{klsonja} is complete.
\qed \end{proof}

As the following proposition shows, effective ep-groupoids have an open and dense subset in the object space $X$ consisting of points  whose 
isotropy groups are trivial.
\begin{proposition}\label{prop_density}\index{P- Density result}
If  $X$ is an effective ep-groupoid, then the set $X_{reg}=\{x\in X\, \vert \, G_x=\{1_x\}\}$ is open and dense in $X$.
\end{proposition}
\begin{proof}
In view of Lemma \ref{ooopen} the set $X_{\textrm{reg}}$ is open in $X$.
If $x_0\in X$, we shall show that every  open neighborhood of $x_0$ contains a point with trivial isotropy.
Given an open neighborhood we find a smaller open neighborhood $U(x_0)$ equipped with  the natural action of $G_{x_0}$.
We shall write $g\ast x$ instead of $\Phi (g)(x)$.
 We may assume that $G_{x_0}\neq \{1_{x_0}\}$
since otherwise we are already done. Take a point $x\in U(x_0)$ whose isotropy group $G_x$ satisfies
\begin{equation}\label{min_eq}
\#G_x = \min \{\# G_y\, \vert \, y\in U(x_0)\}.
\end{equation}
If $\# G_x=1$ we are done. So let us assume this is not the case and derive a contradiction. 

In view of Lemma \ref{klsonja}, 
there exists an open neighborhood $U(x)\subset U(x_0)$ equipped with the $G_x$-action and such that the homomorphism $\Phi: G_x\to \text{Diff}_{\textrm{sc}}(U(x))$ is an injection. Hence, given $g\in G_x\setminus \{1_x\}$, there exists $a\in U(x)$ for which $\Phi (g)(a)=b\neq a$. 
The isotropy group $G_a$ at the point $a$ is isomorphic to the subgroup $H
=\{h\in G_x\,\vert \, \Phi (h)(a)=a\}$ of $G_x$. 
From \eqref{min_eq} we conclude that 
$$\#G_x\leq \#G_a=\#H\leq \#G_x,$$
so that $\#G_x=\#H$ and hence $G_x=H$. This shows that 
$g\in H$, so that $\Phi (g)a=a$ contradicting $\Phi (g)a=b\neq a$. 

This contradiction shows that $\#G_x=1$ and that $X_{\textrm{reg}}$  is dense.  This completes the proof  Proposition \ref{prop_density}.
\qed \end{proof}

We recall from Theorem \ref{x-local-x} that every point $x$ in the object set of an ep-groupoid $(X, \bm{X})$ possesses 
a suitable open neighborhood $U(x)$, which can be chosen arbitrarily small, on which the isotropy group $G_x$ acts by the homomorphism $\Phi\colon G_x\to \text{Diff}_{\textrm{sc}}(U(x)).$

We say that  $g\in G_x$ {\bf acts trivially}\index{Trivially acting element in $G_x$} if there exists an open neighborhood $U(x)$ equipped with the $G_x$-action 
$\Phi: G_x\to \textrm{Diff}_{\textrm{sc}}(U(x)$ such that  $g$ belongs to the kernel of the homomorphism $\Phi$, i.e., if $\Phi (g)=\mathbbm{1}.$
Clearly,  $1_x$ is always acting trivially, but there might be other elements $g\in G_x$ acting trivially.  Since we are dealing with a finite group $G_x$, we can choose $U(x)$ so small that the kernel
of $\Phi$ is equal to the subgroup of $G_x$ acting trivially on $U(x)$.

We denote by $T_x\subset G_x$ the subgroup of $G_x$ acting trivially on $U(x)$, which,  as the kernel of a group homomorphism,  is a normal subgroup of $G_x$.

\begin{definition}\label{DEF725}\index{D- Non-effective part}\index{D- Reduced morphisms}
The {\bf non-effective part $\bm{X}_T\subset \bm{X}$}
is the subset of morphisms defined as 
$$
\bm{X}_T=\{\phi \in T_x\subset G_x\, \vert \, x\in X\}.
$$
The subset $\bm{X}_{\textrm{R}}$ of {\bf reduced morphisms of $X$} is defined as follows. Two morphisms $\phi$ and 
$\phi'\in \bm{X}$ are called equivalent, $\phi\sim \phi'$,  if the following two conditions are satisfied.
\begin{itemize}
\item[(1)]\ $s(\phi)=s(\phi')=x$.
\item[(2)]\ $\phi=\phi'\circ g$ for some $g\in T_x$.
\end{itemize}
The set $\bm{X}_{\textrm{R}}$ is the collection of all equivalence classes $[\phi]$ of morphisms in $\bm{X}$.
\qed
\end{definition}

\begin{proposition}\label{reduced_m_pol}
\index{P- Reduced ep-groupoid}
The reduced morphism set $\bm{X}_{\textrm{R}}$ has a natural M-polyfold structure.
\end{proposition}
\begin{proof}
We equip $\bm{X}_{\textrm{R}}$ with the quotient topology so that 
the quotient map $\pi: \bm{X}\rightarrow \bm{X}_{\textrm{R}}$ is continuous. For every morphism $\phi_0\in \bm{X}$ we find open neighborhoods
$\bm{U}(\phi_0)\subset \bm{X}$ and $U(s(\phi_0))\subset X$ for which  $s\colon \bm{U}(\phi_0)\rightarrow U(s(\phi_0))$ is a sc-diffeomorphism.  Then the map 
$$\pi': \bm{U} (\phi_0)\to \bx_{\textrm{R}},\quad \phi\mapsto [\phi]$$
is continuous.  If  $[\psi]=[\psi']$ for $\psi$ and $\psi'\in \bm{U}(\phi_0)$, then $s(\psi)=s(\psi')$ implies that  $\psi=\psi'$. Hence the map is injective and therefore a bijection onto the image. 

Next we  shall show that the image $\bm{U} ([\phi_0]):=\pi'(\bm{U} (\phi_0))$ is open in $\bx_{\textrm{R}}$. In order to verify this we have to show that the preimage  $\pi^{-1}(\bm{U} ([\phi_0]))$ is open in $\bx$. To see this let  $\psi\in \pi^{-1}(\bm{U} ([\phi_0]))$. Then there is 
$\phi\in \bm{U} (\phi_0)$ such that $\phi\sim \psi$. This means that $s(\phi)=s(\psi)$ and there exists $g\in T_{s(\phi)}=T_{s(\psi)}$ such that $\phi=\psi\circ g.$ 
We have to show that there exists an open neighborhood $\bm{U} (\psi)$ of $\psi$ such that every morphism $\psi'\in \bm{U} (\psi)$ is equivalent to some morphism $\phi'$ belonging to $\bm{U} (\phi_0)$.
We denote by  $U(s(\phi))$ 
an open neighborhood of $s(\phi)$ contained in $U(s(\phi_0))$ and such that $T_{s(\phi)}=T_{s(\psi)}$ is equal to the kernel of the homomorphism $\Phi: G_{s(\psi)}\to \textrm{Diff}_{\textrm{sc}}(U(s(\psi)))$.  
Choosing  an open neighborhood $\bm{U} (\psi)$ of $\psi$ in $\bx$ such that the source map $s: \bm{U} (\psi)\to U(s(\psi))$ is a sc-diffeomorphism, we define  
the map  $\delta: \bm{U} (\psi)\to \bx$ by 
$$\delta (\psi')=(s\vert \bm{U} (\psi))^{-1}(s(\psi'))\circ \Gamma (g, s(\psi')).$$
If  $\psi'=\psi$, then 
$$\delta (\psi)=(s\vert \bm{U} (\psi))^{-1}(s(\psi))\circ \Gamma (g, s(\psi))=\psi\circ g=\phi.$$
Since the map $\delta$ is continuous and $\delta (\psi)=\phi$, there exists an open neighborhood $\bm{U}'(\psi)\subset \bm{U} (\psi)$ of $\psi$ such that $\delta (\bm{U}'(\psi))\subset \bm{U} (\phi_0)$.  Every morphism $\psi'$ in $\bm{U}'(\psi)$ is of the form $\psi'=(s\vert \bm{U} (\psi))^{-1}(s(\psi))$, so that 
$$\delta (\psi')=\psi'\circ \Gamma (g, s(\psi')).$$
Since 
the element $g\in T_{s(\psi)}$, it follows that the morphism 
$\Gamma (g, s(\psi'))$ belongs to $T_{s(\psi')}$, which implies that every morphism $\psi'\in \bm{U}'(\psi)$ is equivalent to the morphism $\delta (\psi')$ belonging to $\bm{U} (\phi_0)$.  This shows that $\bm{U}' (\psi)\subset \pi^{-1}(\bm{U} ([\phi_0])$, implying that  the preimage $\pi^{-1}(\bm{U} ([\phi_0])$ is open in $\bx$. Hence the set $\bm{U} ([\phi_0])$ is open in $\bx_{\textrm{R}}$, as claimed.

The same argument also shows that the map
$$
\bm{\Psi}_{\phi_0}: \bm{U}(\phi_0)\rightarrow \bm{U}([\phi_0]),\quad \phi\mapsto  [\phi]
$$
is an open map. We have shown that the map $\bm{\Psi}_{\phi_0}$ is a homeomorphism between open neighborhoods.  

Next we consider two such maps 
$$\bm{\Psi}_{\phi_0}:\bm{U} (\phi_0)\rightarrow \bm{U} ([\phi_0])\quad \text{and}\quad \bm{ \Psi}_{\psi_0}:\bm{U} (\psi_0)\rightarrow \bm{U} ([\psi_0]),$$
where $U([\phi_0])\cap \bm{U}([\psi_0])\neq\emptyset$, and we claim that the transition map 
$$\bm{\Psi}_{\psi_0}^{-1}\circ \bm{\Psi}_{\phi_0}: \bm{ \Psi}_{\phi_0}^{-1}(\bm{U} ([\phi_0])\cap \bm{U} ([\psi_0]))\to 
\bm{ \Psi}_{\psi_0}^{-1}(\bm{U} ([\phi_0])\cap \bm{U} ([\psi_0]))$$
is sc-smooth. To verify this claim we take a morphism $\phi\in 
\bm{ \Psi}_{\phi_0}^{-1}(\bm{U} ([\phi_0])\cap \bm{U} ([\psi_0]))$ and assume that 
$\bm{ \Psi}_{\psi_0}^{-1}(\phi)=\psi\in \bm{ \Psi}_{\psi_0}^{-1}(\bm{U} ([\phi_0])\cap \bm{U} ([\psi_0]))$. Hence 
$[\phi]=[\psi]\in \bm{U} ([\phi_0])\cap \bm{U}([\psi_0])$. Then $s(\phi)=s(\psi)$ and there exists $g\in T_{s(\phi)}=T_{s(\psi)}$ for which 
$$\psi=\phi\circ g.$$
Choosing a sufficiently small open neighborhood $U(s(\phi))\subset U(s(\phi_0))$ such that $T_x$ agrees with the kernel 
of the $G_{s(\phi)}$-action $\Phi: 
G_{s(\phi)}\to \textrm{Diff}_{\textrm{sc}}(U(s(\phi)))$, we find an open neighborhood $\bm{U} (\phi)\subset \bm{ \Psi}_{\phi_0}( \bm{U} ([\phi_0])\cap \bm{U} ([\psi_0]))$ for which the source map $s: \bm{U} (\phi)\to U(s(\phi))$ is a sc-diffeomorphims. Then every morphism $\phi'$ in $\bm{U} (\phi)$ is of the form $\phi'=(s\vert \bm{U} (\phi))^{-1}(s(\phi'))$. We define the map 
$$\delta (\phi')=(s\vert \bm{U} (\phi))^{-1}(s(\phi'))\circ \Gamma (g, s(\phi'))$$
for $\phi'\in \bm{U} (\phi)).$ 
Since  
$\delta (\phi)=\psi$ and the map $\delta$ is sc-smooth, in particular continuous, there exists an open neighborhood $\bm{U}'(\phi)\subset \bm{U} (\phi)$ such that $\delta (\bm{U}'(\phi)\subset \bm{U} (\psi_0)$. 
The morphism $ \Gamma (g, s(\phi'))$ belongs to $T_{s(\phi')}$ showing that 
$[\phi']=[\delta (\phi')]$ and $\delta (\phi')\in \bm{ \Psi}_{\psi_0}^{-1}(\bm{U} ([\phi_0])\cap \bm{U} ([\psi_0]))$ for all $\phi'\in \bm{U} '(\phi)$.  We have proved that the transition map $\bm{ \Psi}_{\psi_0}^{-1}\circ \bm{ \Psi}_{\phi_0}=\delta$ is sc-smooth,  as claimed.

Summing up,  the maps 
$$
\bm{ \Psi}_{\phi_0}: U(\phi_0)\rightarrow \bm{U} ([\phi_0])
$$
 are homeomorphisms between open subsets of the M-polyfold $\bm{X}$
and open subsets of $\bm{X}_{\textrm{R}}$, and the transition maps are sc-smooth. This defines the natural M-polyfold structure on $\bm{X}_{\textrm{R}}$. 
\qed \end{proof}
Now we define the {\bf  reduced ep-groupoid}\index{reduced ep-groupoid} $X_{\textrm{R}}$ associated with the ep-groupoid $X$. 
On the object level we put $X_{\textrm{R}}=X$ and
as morphism M-polyfold we take $\bm{X}_{\textrm{R}}$. The source and target maps $s, t: \bm{X}_{\textrm{R}} \to X$  
are defined by 
$$
s([\phi]):=s_{X}(\phi)\quad  \text{and}\quad t([\phi])=t_{X}(\phi).
$$
They  are well-defined since two morphisms in an equivalence class have the same source and the target. 
It is clear that $s$ and $t$ are surjective.
For  the local sc-diffeomorphisms 
$\bm{ \Psi}_{\phi}: \bm{U}(\phi)\rightarrow \bm{U} ([\phi])$  we have the identities 
$$
s(\bm{ \Psi}_{\phi}(\phi'))=s_{X}(\phi')\quad  \text{and}\quad t(\bm{ \Psi}_{\phi}(\phi'))=t_{X}(\phi')
$$
which show that the source and target maps are smoothness of the natural maps can be reduced by our local charts to the corresponding statements for $X$. The properness can be proved similarly.
We summarize these arguments in the following proposition.
\begin{proposition}\index{P- Effectiveness of $X_{\textrm{R}}$}
Every ep-groupoid $X$ possesses an associated reduced ep-grou\-poid $X_{\textrm{R}}$. Moreover,  $X_{\textrm{R}}$ is effective.
\end{proposition}
\begin{proof}
The first part is already proved.  The effectiveness of $X_{\textrm{R}}$ is seen as follows.
Fixing the object $x\in X$, we denote by $\wt{G}_{x}$ the isotropy group of the object $x$ in the ep-groupoid $X_{\textrm{R}}$ and observe that $[1_x]=T_x$. If $\tau\in \wt{G}_x\setminus [1_x]$, then $\tau=[g]$ for some $g\in G_x\setminus T_x$. Let $\bm{V}(\tau)$ be an open neighborhood of $\tau$ in $\bm{X}_{\textrm{R}}$ and let $U(x)=U(s(g))$ be an  open neighborhood of $x$ in $X$ which is equipped with the $G_x$-action $\Phi$ such that 
$T_x=\ker \Phi$. The neighborhood $U(x)$ can be taken as small as we wish. Then we choose an open neighborhood $\bm{U} (g)$ such that the source map $s: \bm{U} (g)\to U(x)$ is a sc-diffeomorphim.   It follows from the proof of Proposition \ref{reduced_m_pol} that the $\bm{ \Psi}_{g}: \bm{U} (g)\to \bm{U} ([g])=\bm{U} (\tau)$ is a homeomorphism. We take a neighborhood $U(x)$ so small that $\bm{U} ([g])\subset \bm{V}(\tau)$. 

Since $g\in G_x\setminus T_x$, there exists $a\in U(x)$ for which $\Phi (g)(a)=b\neq a$. Then for  the morphism $\phi=\Gamma (g, a)$ connecting $a$ with $b$,   we conclude that 
$s([\phi])=s_{X}(\phi)=a\neq b=t_{X}(\phi)=t([\phi]).$
Hence the ep-groupoid $X_{\textrm{R}}$ is effective, as claimed.
\qed \end{proof}

If we start with an ep-groupoid $X$ and take the associated reduced ep-groupoid $X_R$ we recall that the underlying object spaces
are the same, but the morphism spaces differ. If $x\in X$ and $G_x$ is the associated isotropy group, then the isotropy group in the reduced 
ep-groupoid is $G_x/T_x$, where $T_x$ is defined before Definition \ref{DEF725}.
\begin{definition}\label{DEF728}
For and object $x$ in an ep-groupoid $X$ we shall call $G_x/T_x$ the effective isotropy group of $x$ and denote it by $G_x^{\text{eff}}$.
\end{definition}

\section{Topological Properties of Ep-Groupoids}\label{section1.3_top_prop}
In this subsection we shall study the topological properties of ep-groupoids. The following result  summarizes the topological properties.
\begin{theorem}[{\bf Basic Topological Properties of $|X|$}]\label{topos-x}\index{T- Topological properties}\label{ATHOME}
Let $(X,\bm{X})$ be an ep-groupoid. Then  the orbit space $|X|$ is a locally metrizable, regular, Hausdorff topological space.
If in addition $|X|$ is paracompact then it also is metrizable, which is a consequence of the Nagata/Smirnov metrization theorem.
\qed
\end{theorem}
The theorem follows from several lemmata.
\begin{lemma}[{\bf Local metrizability}]\label{local_metrizability_x}\index{L- Local metrizability of orbit spaces}
Let $(X,\bm{X})$ be an ep-groupoid. Then every point $z\in |X|$ has an open neighborhood $V=V(z)$ on which the topology is metrizable, i.e. the orbit space 
$|X|$ is locally metrizable. In fact one can  take $|U(x)|$ where $U(x)$ is an open neighborhood equipped with the natural $G_x$-action.
\end{lemma}
\begin{proof}
We choose  a representative $x\in X$ of $z$ so that  $z=\abs{x}$ and  take an open neighborhood $U=U(x)$ of the form guaranteed by Theorem \ref{x-local-x}.
In particular,  we have a $G_x$-action on $U$ and $V(z):=|U|$ is open by Proposition \ref{q-open}.
The subset $U$ of the M-polyfold $X$ is metrizable and since it is invariant under the $G_x$-action, which acts by sc-diffeomorphisms,
we can define, using a metric $d$ on $U$, a $G_x$-invariant metric $\rho$ on $U$ by
$$
\rho(y,y')=\text{max}_{g\in G_x} d(g\ast y,g\ast y').
$$
The metric $\rho$ on $U$ is equivalent to the original one in the sense that it induces the same topology, and is,  in addition, invariant under the action of $G_x$,
$$
\rho(h\ast y,h\ast y') = \rho (y,y')\quad  \text{for $h\in G_x$.}
$$
On ${_{G_x}\backslash}U$ we obtain the  induced metric $\wh{\rho}$,  defined by
$$
\wh{\rho}([y],[y'])= \text{min}\{\rho(g\ast y,h\ast y')\, \vert \,  g,h\in G_x\}, 
$$
which defines the quotient topology. By Proposition \ref{go-go} the map ${_{G_x}\backslash} U\rightarrow V(z)$ is a homeomorphism.
Hence the push-forward of the metric $\what{\rho}$ to $|U|=V(z)$ shows that $V(z)\subset |X|$ is metrizable.
\qed \end{proof}
The next result follows from the previous one.
\begin{lemma}
Let $(X,\bm{X})$ be an ep-groupoid. Then for the quotient topology on $|X|$ every point $z\in \abs{X}$ has a countable neighborhood basis, consisting of open sets,
$$
{\mathcal V}_z =\{V_k(z)\, \vert \, k=1,2,3,\ldots \}
$$
satisfying  $z\in V_{k+1}(z)\subset V_k(z)$.
\end{lemma}
\begin{proof} In view of  Lemma \ref{local_metrizability_x}
every point in $|X|$ has a metrizable neighborhood.
\qed \end{proof}

That $|X|$ is always Hausdorff is mainly a consequence of the properness property (3)  of an ep-groupoid as the next lemma shows.

\begin{lemma}[{\bf Hausdorffness of orbit space}]\label{lemma1.37}\index{L- Hausdorffness of orbit spaces}
The orbit space $|X|$ of an ep-groupoid  $(X,\bm{X})$ is Hausdorff. 
\end{lemma}
\begin{proof}
We fix two different points $z,z'$ in $|X|$ and choose representatives $x,x'\in X$ satisfying  $\pi(x)=z$ and $\pi(x')=z'$. We can take countable monotonic neighborhood
bases $(U_k(x))$ and $(U_k(x'))$ in $X$ and obtain  for $z,z'$ the sets $V_k(z)=|U_k(x)|$ and $V_k(z')=|U_k(x')|$ in $\abs{X}$. It suffices to show that  $V_k(z)\cap V_k(z')=\emptyset$ for $k$ large. 
Arguing indirectly we find $z_k\in V_k(z)\cap V_k(z')$. This implies the existence of sequences   $x_k\in U_k(x)$ and $x_k'\in U_k(x')$ satisfying  $\pi(x_k)=\pi(x_k')=z_k$. We conclude that  $x_k\rightarrow x$ and $x_k'\rightarrow x'$ in $X$ since $(U_k(x))_{k\geq 1}$ and $(U_k(x')_{k\geq 1}$ are  monotone neighborhood bases.
We take a sequence $(\phi_k)\subset \bm{X}$ satisfying  $s(\phi_k)=x_k$ and $t(\phi_k)=x_k'$. Using the properness condition (3), we conclude after taking a suitable 
subsequence, that there exists a morphism $\phi_0\in \bm{X}$ with $s(\phi_0)=x$ and $t(\phi_0)=x'$ implying $z=z'$ and contradicting $z\neq z'$.
\qed \end{proof}
We recall that a topological space is {\bf regular}\index{Regular topological space} if  for a nonempty closed subset $A$
and a point $z$ not belonging to $A$ there exist open sets $U$ and $V$ such  that $z\in U$, $A\subset V$ and $U\cap V=\emptyset$.
A slight strengthening of the previous result shows that the orbit space $|X|$ of an ep-groupoid is a regular topological space.
Since we already know that $|X|$ is Hausdorff this implies that $|X|$ is a regular Hausdorff space.
\begin{lemma}[{\bf Regularity}]\index{L- Regularity of orbit spaces}
Let $(X,\bm{X})$ be an ep-groupoid.  Then $|X|$ is a regular topological space.
\end{lemma}
\begin{proof}
Assuming  we are given a closed subset $A$ of $|X|$ and a point $z\in |X|\setminus A$,  we choose a representative 
$x\in X$ satisfying $\pi(x)=z$.  Since $|X|\setminus A$ is open and contains $z$,  we find by continuity of $\pi$ an open neighborhood
$O$ of $x$ in $X$  such that  $\pi(O)\subset |X|\setminus A$. Since $X$ is metrizable and therefore normal,  we find an open neighborhood $O'$ of $x$ in $X$ such that 
$$
x\in O'\subset \cl_X(O')\subset O.
$$
By the properness assumption we may assume that $O$ was chosen  in such a way that
$$
t\colon  s^{-1}(\cl_X(O))\rightarrow X
$$
is proper.  We  choose a point $y\in \pi^{-1}(A)$ in  the preimage of $A$ and show that there exists an open neighborhood
$V(y)$ such  that there exists no morphism $\phi$ satisfying $s(\phi)\in O'$ and $t(\phi)\in V(y)$.
Since $X$ is a metrizable space,  we find otherwise a sequence $(\phi_k)$ of morphisms satisfying 
$s(\phi_k)\in O'$ and $t(\phi_k)\rightarrow y$. By the properness assumption we may assume 
after perhaps taking a subsequence that $\phi_k\rightarrow \phi$. Then
$s(\phi)\in \cl(O')$ and $t(\phi)=y$. Since $\cl(O')\subset O$, this implies the contradiction $\pi(y)\in \pi(O)\cap A$.
Hence we find for every $y\in\pi^{-1}(A)$ an open neighborhood $V(y)$ for which  $\pi(V(y))\cap \pi(O')=\emptyset$. The set 

$$
\wh{V}=\bigcup_{y\in \pi^{-1}(A)} V(y)
$$
is an open subset of $X$ and therefore $V:=\pi(\wh{V})$
is an open subset of $|X|$ containing $A$. By construction,  $V\cap \pi(O')=\emptyset$. 
\qed \end{proof}

In view of the above lemmata, the proof of Theorem \ref{topos-x} is complete.
\qed

As a consequence of Urysohn's metrizability theorem we can conclude,  under an additional  condition,  that the orbit space  $|X|$ is (globally) metrizable.
\begin{corollary}[{\bf Metrizability criterion}]\index{C- Metrizability criterion}
Let $X$ be an ep-groupoid whose induced  topology ${\mathcal T}$ on $|X|$ is second countable, i.e. ${\mathcal T}$ has a countable basis.
Then $|X|$ is metrizable.
\end{corollary}
\begin{proof} Urysohn's metrizability result states that a second countable, regular Hausdorff space is metrizable.
As we have shown the orbit space $|X|$ of an ep-groupoid is always a regular Hausdorff space. The additional assumption that 
$|X|$ is second countable gives via Urysohn's theorem the desired result.
\qed \end{proof}
Continuing with a general ep-groupoid $X$  we observe that the induced topology on a given subset  of $|X|$ is again locally metrizable, regular and Hausdorff.
\begin{definition}\index{D- Component compact subset}
If  $X$ is an ep-groupoid, we call a subset $K$ of the orbit space $|X|$ {\bf component compact} provided its intersection with every connected component of $|X|$ 
is compact. 
\qed
\end{definition}
Theorem \ref{ATHOME} has a second corollary, which employs the  Nagata-Smirnov metrizability theorem.

\begin{corollary}\label{NHOOD}\index{C- Metrizability criterion}
Let $(X,\bm{X})$ be an ep-groupoid and $K\subset |X|$ a  component compact subset. Then there exists an open neighborhood $U$ of $K$ for which 
$\cl_{|X|}(U)$ is paracompact. Since $\cl_{|X|}(U)$ is a regular Hausdorff space which is locally metrizable, it follows, in particular,  that $\cl_{|X|}(U)$ is metrizable 
and therefore $U$ is metrizable.
\end{corollary}
\begin{proof}
In order to prove the corollary we need the Nagata-Smirnov metrizability theorem.  We work in a connected component of $K$
and may assume without loss of generality that $K$ is compact. Every point $z\in K$ has an open neighborhood
$U(z)$ whose  closure in $|X|$ is metrizable. In particular,  $\cl_{|X|}(U(z))$ is paracompact. Using the compactness
of $K$ we find finitely many points $z_1,\ldots ,z_k$  in $K$ such that $U(z_1),\ldots ,U(z_k)$ cover $K$. Then the union 
$$
\ov{U}=\bigcup_{i=1}^k \cl_{|X|}(U(z_i))
$$
is paracompact. Moreover,  as a subset of $|X|$ its induced topology  is Hausdorff, regular and locally metrizable. 
By the Nagata Smirnov Theorem $\ov{U}$ is metrizable and the same holds then for the set $U=\bigcup_{i=1}^k U(z_i)$.
 \qed \end{proof}
 
 There are many new features in the sc-theory which do not occur in the classical theory. For example, 
we can raise the index of a M-polyfold $X$ to obtain the M-polyfold $X^1$. Doing this, the question arises whether for an ep-groupoid $X$ we obtain an ep-groupoid $X^1$. The answer is yes, and this is 
the next proposition.

\begin{proposition}[{\bf Index Lifting}]\label{raising-1}\index{P- Index lifting}
If $X=(X, \bm{X})$ is an ep-groupoid, then  the same is true for $X^1=(X^1, \bm{X}^1)$. If $X$ is tame so is $X^1$.
\end{proposition}
\begin{proof}
We note that $X^1=(X^1, \bm{X}^1)$ satisfies trivially  properties (1), (2) and (4) of Definition \ref{ep-groupoid_def}. 
We shall show that  the property (3) holds.   We take a point $x\in X_1$.  Then $x\in X$ and there exists an open neighborhood $V(x)$ of $x$ in $X$ such that 
\begin{equation}\label{ep_groupoid_eq_0}
t: s^{-1}(\textrm{cl}_X (V(x)))\to X
\end{equation}
is proper. Since the inclusion $i:X_1\to X$ is continuous, the set $W(x):=V(x)\cap X_1$ is open in $X_1$ and we also have $\textrm{cl}_{X_1}(W(x)) \subset \textrm{cl}_X (V(x))$, where $\textrm{cl}_{X_1}$ stands for the closure of a set in $X_1$. We claim that the map
\begin{equation}\label{ep_groupoid_eq_1}
t: s^{-1}(\textrm{cl}_{X_1}(W(x))\to X_1
\end{equation}
is proper.  It suffices to show that every sequence $(\phi_k)\subset s^{-1}(\cl_{X_1}(W(x))$ such that  $(t(\phi_k))$ belongs  to a compact subset of $X_1$, has a convergent subsequence in $\bm{X}_1$. Using the properness on level $0$ and the pre-compactness of $(t(\phi_k))$ on level $1$
we may assume after taking a suitable subsequence the following properties of the sequence $(\phi_k)$.
\begin{itemize}
\item[(a)]\ $\phi_k\rightarrow \phi$ in  $\bm{X}_0$.
\item[(b)]\ $s(\phi_k)\rightarrow s(\phi)\in \cl_X(V(x))$  and $t(\phi_k)\rightarrow t(\phi)$ in $X_0$.
\item[(c)]\ $t(\phi_k)\rightarrow t(\phi)$ in $X_1$.
\end{itemize}
From the fact that $t$ is a local sc-diffeomorphism and using (c) we deduce that 
\begin{itemize}
\item[(d)]\ $\phi\in \bm{X}_1$.
\end{itemize}
We  find open neighborhoods $\bm{U}(\phi)$ and $U(t(\phi))$ in $\bm{X}_0$ and $X_0$, respectively, such  that
$$
t\colon \bm{U}(\phi)\rightarrow U(t(\phi))
$$
is a sc-diffeomorphism. For large $k$, the sequences $\phi_k$ and $t(\phi_k)$ belong to these neighborhoods.
So far,  by construction, 
$$
\phi_k=(t\vert \bm{U}(\phi))^{-1}(t(\phi_k))
$$
for large $k$. In view of (d),  $t(\phi_k)\rightarrow t(\phi)$ in $X_1$ and since $(t\vert \bm{U}(\phi))^{-1}$ is $\ssc^0$,  the convergence 
$$
\phi_k=(t\vert \bm{U}(\phi))^{-1}(t(\phi_k))\rightarrow (t\vert \bm{U}(\phi))^{-1}(t(\phi))=\phi\quad \text{in $\bm{ X}_1$}
$$
follows. 
This proves the properness on level $1$.
That the tameness of $X$ implies the tameness of $X^1$,  is trivial. The proof of Proposition \ref{raising-1} is complete.
\qed \end{proof}

\section{Paracompact Orbit Spaces}
We have shown that the orbit space $|X|$ of an ep-groupoid $(X, \bm{X})$ is a locally metrizable, regular  Hausdorff topological space.
This  means that single points are closed subsets. Moreover,  a point $z\in \abs{X}$ and a closed subset $A$ of $|X|$ not containing $z$
can be separated by open subsets containing them. In practice many of the interesting ep-groupoids have orbit spaces whose  topologies  are  paracompact.
In view of the Nagata-Smirnov metrizability theorem such an orbit space will be metrizable. This has, of course, pleasant features when one has to carry out constructions.

\begin{definition}
The ep-groupoid  $(X,\bm{X})$ is an {\bf ep-groupoid with paracompact orbit space}\index{D- Ep-groupoid with paracompact orbit space}
if, in addition to the properties (1)--(4) of Definition \ref{ep-groupoid_def}, it also possesses the following property  (5).
\begin{itemize}
\item[(5)]\ The orbit space $|X|$ is paracompact.
\end{itemize}
\qed
\end{definition}
We recall that a topological space is {\bf paracompact}\index{paracompact} if  every open covering 
has a refinement which is a locally finite open covering.  We do not require a paracompact space to be Hausdorff. It is presently not known whether the property  (5) is a consequence of the properties (1)--(4).  In case
$X$ has a second countable topology we have already seen that (1)-(4) imply (5).
\begin{remark}
In our paper \cite{HWZ3.5} we always required the M-polyfold to be second countable.
We later realized that this  requirement is unnecessarily restrictive and that a better requirement is (5).
\qed
\end{remark}

\begin{proposition}[{\bf Paracompactness under lifts}]\index{P- Paracompactness under lifts}
Let $(X,\bm{X})$ be an ep-groupoid with paracompact orbit space. Then also $(X^1,\bm{X}^1)$ is an ep-groupoid with paracompact orbit space.
\end{proposition}
\begin{proof}
In view of Proposition \ref{raising-1} we only have to show that $|X^1|$ is paracompact and we already know that $|X|$ and $|X^1|$
are regular Hausdorff spaces.
We  recall from Lemma \ref{equivalent_definitions_paracompact}  that a regular Hausdorff space $Z$ is a  paracompact  space if and only if every open covering of $Z$ 
has a locally finite refinement of closed sets. 

For every $x\in X$ we choose an open neighborhood $U(x)$ of $x$ in $X$ equipped with the natural $G_x$-action. 
Since $\abs{X}$ is Hausdorff and paracompact,  the open cover $(\abs{U(x)})_{x\in X}$ of $\abs{X}$ has a refinement $(A_x)_{x\in X}$ having the following properties,
\begin{itemize}
\item[(1)]\  $(A_x)$ is a locally finite collection of closed sets.
\item[(2)]\ $A_x\subset \abs{U(x)}$.
\item[(3)]\ $\bigcup_{x\in X} A_x =\abs{X}$.
\end{itemize}
We note that some of the sets $A_x$ may be empty. We abbreviate  $A_{x,1}=A_x\cap |X_1|$. The set $A_{x,1}$ is the preimage of $A_x$ under the continuous map $|X_1|\rightarrow |X|$
and therefore closed. The collection $(A_{x,1})_{x\in X}$ inherits the property that $(A_x)$ is locally finite and,  by construction, 
$$
|X_1|=|X_1|\cap |X|= |X_1|\bigcap\left( \bigcup_{x\in X} A_x\right)=\bigcup_{x\in X} A_{x,1}.
$$
We shall show that each set $A_{x,1}$ is paracompact. Then the Hausdorff orbit space $\abs{X_1}$ as the union of a locally finite collection of closed paracompact sets is paracompact.   
Since $U(x)$ has the natural $G_x$-action, the same is true for $U(x)_1\subset X_1$. As in the proof of Lemma \ref{local_metrizability_x}
one can take any metric on $U(x)_1$ inducing the topology, and then construct a $G_x$-invariant metric, which then passes to the quotient
$|U(x)|_1=|U(x)_1|$. By construction,  $A_{x,1}\subset \abs{U(x)}_1$ so that $A_{x,1}$  inherits an induced metric and therefore is paracompact.
\qed \end{proof}

In the following a {\bf saturated subset}\index{Saturated subset} $U$ of an ep-groupoid is a subset  of the object space which contains with the object  $x\in U$ also all the objects  isomorphic  to $x$,  i.e., 
$U=\pi^{-1}(\pi(U))$, where $\pi\colon X\rightarrow |X|$ is the projection onto the orbit space.

The importance of paracompactness comes from the existence of partitions of unity which, for example,  is used in the next theorem.

\begin{theorem}[{\bf Continuous Partitions of Unity}]\label{c_partition_of_unity}\index{T- Continuous partitions of unity}
For an ep-grou\-poid $X$ with paracompact orbit space, there is  for every open cover ${(O_\lambda)}_{\lambda\in\Lambda}$ of $|X|$
a subordinate partition of unity ${(\sigma_\lambda)}_{\lambda\in\Lambda}$ (some of the $\sigma_\lambda\equiv 0$).
In particular,  every covering $(U_\lambda)$ of $X$ by saturated open sets possesses  a subordinated continuous partition of unity
$(\beta_\lambda)$,  where every function $\beta_\lambda $ takes the same values on isomorphic objects.
\end{theorem}
\begin{proof}
Due to paracompactness
the first part is an immediate consequence of a partition of unity. If  $(U_\lambda)$ is an open covering of $X$ by saturated open sets, 
then $O_\lambda=\pi(U_\lambda)$ is an open cover of  $|X|$. We choose  a subordinated partition of unity $(\sigma_\lambda)$ and 
define $\beta_\lambda=\sigma_\lambda\circ \pi$, where $\pi\colon X\rightarrow |X|$ is the quotient map.
Then  the support of $\beta_\lambda$ lies in $\pi^{-1}(O_\lambda)=U_\lambda$.  Taking  a point $x\in X$, we obtain $z=\pi(x)\in \abs{X}$ and  find an open neighborhood
$V=V(z)$ so that only finitely many functions $\sigma_\lambda$ have a support intersecting $V$. This means that only finitely many $\beta_\lambda $
have a support intersecting the open neighborhood $U=\pi^{-1}(V)$ of $x$.
\qed \end{proof}

For many differential geometric constructions we need the  sc-smooth partitions of unity or at least sc-smooth bump functions as already discussed in Appendix \ref{POU}.

\begin{definition}[{\bf Sc-Smooth Partitions of Unity}]\index{D- Sc-smooth partitions of unity for ep-groupoids}
An {\bf ep-groupoid $\bm{X}$ with sc-smooth partitions of unity} is an ep-groupoid $X$ having the property  that for every saturated open covering ${(U_\lambda)}_{\lambda\in\Lambda}$ 
of the object M-polyfold $X$ there exists
an subordinate sc-smooth partition of unity ${(\beta_\lambda)}_{\lambda\in\Lambda}$ such that 
\begin{itemize}
\item[(1)]\ For every point $x\in X$ there exists a saturated open neighborhood intersecting only finitely many supports $\text{supp}(\beta_\lambda)$.
\item[(2)]\ For every $\lambda\in\Lambda$ the map $\beta_\lambda$ has the property that $\beta_\lambda(z)=\beta_\lambda(y)$ if $y$ and $z$ are isomorphic objects.
\end{itemize}
\qed
\end{definition}

\begin{remark}\index{R- On sc-smoothness partitions of unity}
Unfortunately this topic,  apart from the case of sc-Hilbert spaces, seems to be generally a zoo, as is the question of classically smooth partitions of unity
on Banach spaces. In the case of sc-smoothness it could potentially be easier to show their existence on concrete scales.\hfill \qed
\end{remark}

The basic result about sc-smooth partitions of unity for ep-groupoids is given by the following theorem.
\begin{theorem}\label{pert}\index{T- Partitions of unity}
Let $X$ be an ep-groupoid with paracompact orbit space $|X|$. We assume that the object M-polyfold $X$ admits sc-smooth partitions of unity. Then, for every  open cover of $X$ by saturated open subsets ${(U_\lambda)}_{\lambda\in\Lambda}$,  there exists a subordinate
sc-smooth partition of unity of the ep-groupoid $X$.
\qed
\end{theorem}
The proof of Theorem \ref{pert} is given in Appendix \ref{SCPART}.

We end this subsection with a characterization of closed subsets of the orbit space $|X|$ assuming that $|X|$ is paracompact.
\begin{definition}\label{proper*}\index{D- Properness property}
 We call a subset $Q$ of $X$ {\bf proper} provided $t: s^{-1}(Q)\rightarrow X$ is proper.
 We shall say that an open subset $U$ of $X$ has the {\bf properness property} provided $\cl_X(U)$ is proper. 
 \qed
\end{definition}
\begin{proposition}\label{proper-closed}
If $Q\subset X$ is proper, then it is a closed subset of $X$.
\end{proposition}
\begin{proof}
For a point $x\in \text{cl}(Q)$ we take a sequence $(x_n)\subset Q$ converging to $x$ and claim that $x\in Q$. Let  $\phi_n$ be the morphisms $1_{x_n}: x_n\to x_n$. The set $K=\{x_n\, \vert\, n\geq 1 \}\cup \{x\}$ is a compact subset of $X$. By the properness assumption, the preimage $(t\vert s^{-1}(Q))^{-1}(K)$  of $K$ under the map $t: s^{-1}(Q)\rightarrow X$ is a compact subset of $s^{-1}(Q)$.  Since $s(\phi_n)=x_n\in Q$ and $t(\phi_n)=x_n$, the sequence $(\phi_n)$ belongs to the compact subset $(t\vert s^{-1}(Q))^{-1}(K)$  of $s^{-1}(Q)$. Without loss of generality we may assume that $\phi_n\to \phi$ for some morphism $\phi\in (t\vert s^{-1}(Q))^{-1}(K) \subset s^{-1}(Q)$. In particular, 
$s(\phi )\in Q$. Since $s(\phi_n)=x_n$ and $x_n\to x$, we conclude after taking a limit that $s(\phi )=x\in Q$, as claimed.
\qed \end{proof}

\begin{proposition}\index{P- Closedness property}\label{closedness}
Let $X$ be an ep-groupoid.
\begin{itemize}
\item[{\em (1)}]\ If $Q\subset X$  is a proper subset, then $\pi(Q)$ is closed in $|X|$.
\item[{\em (2)}]\ If $|X|$ is paracompact, then a subset  $A$ of $|X|$ is closed if  and only if there exists a proper subset  $Q\subset X$ satisfying  $\pi (Q)=A$.
\end{itemize}
\end{proposition}
\begin{proof}
We first prove (1) and assume that $Q$ is proper subset of $X$. 
We abbreviate  $A:=\pi(Q)$. We take $z\in \text{cl}(A)$ and let $y\in X$ be a representative of $z$ so that $\pi (y)=z$. 
Taking a monotone neighborhood basis $(U_k(y))_{k\geq 1}$ we abbreviate $V_k(z)=\pi (U_k(y))\subset \abs{X}$. The open sets $V_k(z)$ form a monotone neighborhood basis of $z$ in $\abs{X}$. Since $z\in \text{cl}(A)$, for every $k\geq 1$, there exists $z_k\in A\cap V_k(z)$. 
In particular, $z_k\to z$. Moreover, there are $x_k\in Q$  and $y_k\in U_k(y)$ satisfying $\pi (x_k)=\pi (y_k)=z_k$. 
Hence there is a sequence of morphisms $(\phi_k)$ in $\bm{X}$ such that $s(\phi_k)=x_k$ and $t(\phi_k)=y_k$. Since  $y_k\to y$, the set $K=\{y_k\, \vert \, k\geq 1\}\cup \{y\}$ is compact in $X$. By assumption, 
$(t\vert s^{-1}(Q))^{-1}(K)$ is a compact subset of $s^{-1}(Q)$.  Then, since $(\phi_k)\subset  (t\vert s^{-1}(Q))^{-1}(K)$, we may assume that $\phi_k\to \phi$ for some morphism $\phi\in s^{-1}(Q)$. This implies that $x=s(\phi)\in Q$. From the convergence $z_k\to z$ and $z_k=\pi (s(\phi_k))\to \pi (s(\phi ))=\pi (x)$, we conclude that $\pi (x)=z\in A$ proving (1).

In order to prove (2) we assume that $|X|$ is paracompact. In view of  (1) it suffices to prove  that a closed subset  $A\subset |X|$ can be written
as $\pi(Q)=A$ for some proper subset $Q$ of $X$.  For every $z\in A$ we choose a representative  $x_z\in X$ satisfying  $\pi(x_z)=z$ and  an open neighborhood
$U(x_z)$ for which we have the natural $G_{x_z}$-action so that $Q_z:=\cl_X(U(x_z))$ is proper and therefore $\pi(Q_z)$ is closed. 
Then the open sets $\pi(U(x_z))$ are an open covering of  $A$. Since $A$ is a closed subset of a paracompact space it is paracompact.
By  Lemma \ref{equivalent_definitions_paracompact}  we  find a locally finite covering of $A$ by closed subsets,  say $A_i\subset A$, $i\in I$, so that 
$A_i\subset \pi(Q_{z_i})$ for a suitable map $i\mapsto z_i$. Define the  closed subset $P_i\subset \cl_X(U_{x_{z_i}})$ by
$$
P_i = (\pi | \cl_X(U_{x_{z_i}}))^{-1}(A_i).
$$
Finally,  introducing  $Q\subset X$ by
$$
Q=\bigcup_{i\in I} P_i,
$$
then,  by construction,  $\pi(Q)=A$ and we  only need to show that $Q$ is proper.  Take a sequence of morphisms $(\phi_k)$ satisfying 
$s(\phi_k)\in Q$ and $t(\phi_k)\rightarrow y$. Since $A$ is closed we see that $z=\pi(y)\in A$.  Since the covering $(A_i)$ is locally finite,
we find an open neighborhood $V(z)\subset |X|$ which  intersects only finitely many $A_i$, say $A_1,\ldots ,A_l$.
For $k$ large, $\pi(t(\phi_k))\in V(z)$. This implies that $\pi(s(\phi_k))\in A_1\cup\ldots \cup A_l$ for large $k$.
Hence, for large $k$, by construction $s(\phi_k)\in P_1\cup\ldots \cup P_l$. After perhaps taking a subsequence we may assume that
$s(\phi_k)\in P_1$. At this point we know that $s(\phi_k)\in P_1$ and $t(\phi_k)\rightarrow y$. By construction,  $P_1$ is proper and hence 
we find a convergent subsequence. This shows that $Q\subset X$ is proper and concludes the proof of Proposition \ref{closedness}.
\qed \end{proof}

\section{Appendix}
\subsection{The Natural Representation}\label{k-natural}
We recall the statement for the convenience of the reader
\begin{theorem}\label{LocalStr1thm}
Given an ep-groupoid $X$, an object  $x_0\in X$, and an open neighborhood $V\subset X$ of $x_0$. Then there exist  an open
neighborhood $U\subset V$ of $x_0$,  a group homomorphism
$$
\Phi  : {G}_{x_0}\rightarrow \text{Diff}_{sc}(U), \quad g\mapsto
\Phi (g)=t_g\circ s_g^{-1},
$$
and an  sc-smooth map
$$
\Gamma : G_{x_0}\times U\rightarrow \bm{X}
$$
having the following properties:
\begin{Myitemize}
\item $\Gamma(g,x_0)=g$.
\item $s(\Gamma(g,y))=y$ and $t(\Gamma(g,y))=\Phi (g)(y)$ for all $y\in U$ and $g\in G_{x_0}$.
\item If $h : y\rightarrow z$  is a morphism  connecting  two objects $y,z\in U$,  then there exists a
unique $g\in G_{x_0}$ such that  $\Gamma(g,y)=h$.
\end{Myitemize}
\end{theorem}
In particular, every morphism between points in $U$ belongs to the
image of the map $\Gamma$.  We call  the group homomorphism $\Phi : G_x\rightarrow
\text{Diff}_{sc}(U)$ a {\bf natural representation} of the
stabilizer  group   $G_{x_0}$.
\begin{proof}

For
every $g\in G_{x_0}$ we choose two  contractible open neighborhoods
$N^t_g$ and $N^s_g\subset \bm{X}$ on which the target and source
maps $t$ and $s$ are sc-diffeomorphisms onto some open neighborhood
$U_0\subset X$ of $x_0$. Since  the isotropy group $G_{x_0}$ is
finite we can assume that the open sets $N^t_g\cup N^s_g$ for $g\in
G_{x_0}$ are disjoint and define the disjoint open neighborhoods
$N_g\subset \bm{X}$ of $g$ by
$$N_g:=N_g^t\cap N_g^s,\quad g\in G_{x_0}.$$
We  abbreviate the restrictions of the source and target maps by
$$s_g:=s\vert N_g\quad \text{and}\quad t_g:=t\vert N_g.$$

\begin{lemma}\label{LocalStr1lem}
With the choices made above there exists an open neighborhood
$U_1\subset U_0$ of $x_0$ so that every morphism $h\in \bm{X}$ satisfying
$s(h)$ and $t(h)\in U_1$ belongs to $N_g$ for some  $g\in G_{x_0}$.
\end{lemma}
\begin{proof}
Arguing indirectly we find a sequence $h_k\in \bm{X}$ with $h_k\not
\in N_g$ for all $g\in G_{x_0}$ and satisfying $s(h_k), t(h_k)\to
x_0$ as $k\to \infty$. By the properness assumption of ep-polyfolds
there is a convergent subsequence $h_{k_l}\to h\in \bm{X}$.
Necessarily $h\in G_{x_0}$ and hence $h\in N_g$ for some $g\in
G_{x_0}$. This contradiction implies the lemma.
\qed \end{proof}

\begin{lemma}\label{LocalStr2lem}
If  $U_1$ is  the open neighborhood of $x_0$ guaranteed by Lemma
\ref{LocalStr1lem}, then there exists an open neighborhood $U_2\subset
U_1$ of $x_0$ so that the open neighborhood $U$ of $x_0$,  which is
defined as the union
$$U:=\bigcup_{g\in G_{x_0}}t_g\circ s_g^{-1}(U_2),$$
is contained in $U_1$ and invariant
 under all the maps $t_g\circ s_g^{-1}$ for $g\in G_{x_0}$.
\end{lemma}
\begin{proof}
We choose an open neighborhood $U_2\subset U_1$ of $x_0$  so small
that the union $U$ and also $t_g\circ s_g^{-1}(U)$ are contained in
$U_1$ for all $g\in G_{x_0}$. Consider the map $t_g\circ
s^{-1}_g : U\to X$ and choose $x\in U$. Then we can represent it as
$x=t_h\circ s_h^{-1}(u)$ for some $h\in G_{x_0}$ and some $u\in
U_2$. Now, $v:=t_g\circ s_g^{-1}\circ t_h\circ s_h^{-1}(u)$ belongs
to $U_1$ and the formula implies the existence of a morphism $u\to
v$ in $\bm{X}$. By Lemma \ref{LocalStr1lem} the morphism has
necessarily the form $v=t_{g'}\circ s^{-1}_{g'}(u)$ for some $g'\in
G_{x_0}$. Since $u\in U_2$ it follows that $v=t_g\circ
s_g^{-1}(x)\in U$ implying the desired invariance of $U$.
\qed \end{proof}

In view of Lemma \ref{LocalStr2lem} we can associate  with every $g\in
G_{x_0}$ the sc-diffeomorphism
$$\Phi (g):=t_g\circ s_g^{-1} : U\to U$$
of the open neighborhood $U$ of $x_0$, and obtain the mapping
$$\Phi  :  G_{x_0}\to \text{Diff}_{sc}(U), \quad g\mapsto \Phi (g).
$$
 Since the neighborhoods $N_g\subset \bm{X}$ of $g$ are
disjoint and since the structure maps are continuous we conclude
that $\Phi $ is a homomorphism of groups, in the following called
the {\bf natural representation of the stabilizer group $G_{x_0}$}
by sc-diffeomorphisms of the open neighborhood $U\subset X$ of
$x_0$. Then we define the map
$$
\text{$\Gamma :  G_{x_0}\times U\rightarrow \bm{X}$ \quad by \quad $\Gamma(g,y)= s_g^{-1}(y).$}
$$
Summing up the consequences of Lemma \ref{LocalStr1lem} and Lemma
\ref{LocalStr2lem} we have proved Theorem \ref{LocalStr1thm}.
\qed \end{proof}

\subsection{Sc-Smooth Partitions of Unity}\label{SCPART}
In this part of the appendix we establish  the existence of sc-smooth partitions of unity on certain ep-groupoids.
We recall Theorem \ref{pert} for the convenience of the reader.
\begin{theorem}
Let $X$ be an ep-groupoid with paracompact orbit space $|X|$. We assume that the object M-polyfold $X$ admits sc-smooth partitions of unity. Then, for every  open cover of $X$ by saturated open subsets ${(U_\lambda)}_{\lambda\in\Lambda}$,  there exists a subordinate
sc-smooth partition of unity of the ep-groupoid $X$.
\end{theorem}

We start the proof with a lemma. 

\begin{lemma}\label{lemma1.126}
We assume that $X$ is an ep-groupoid with  a paracompact  orbit space $|X|$ and the object M-polyfold $X$ admits sc-smooth partitions of unity. 
Suppose further that $x\in X$ and $U(x)$ and $V(x)$ are open neighborhoods of $x$ satisfying  $U(x)\subset \cl_X(U(x))\subset V(x)$ and such that $V(x)$ admits the natural $G_x$-action  and $U(x)$ is invariant.  

If $\beta: X\rightarrow [0,1]$ is a continuous function with support in $U(x)$ and invariant  under the $G_x$-action on $V(x)$, 
then there exists a sc-smooth function  $f:X\rightarrow [0,\infty)$ with support in $V(x)$, which is,  in addition , invariant under the $G_x$-action on $V(x)$ 
and satisfies $\beta(y)\leq f(y)$ for all $y\in X$.
\end{lemma}
\begin{proof}
For every $z\in U(x)$ we take an open neighborhood $W(z)$ with 
$$
\cl_X(W(z))\subset U(x)
$$
 and a number $a_z>0$ so that
$$
\beta(y)< a_z\ \text{for}\ y\in W(z).
$$
For every $z\in V(x)\setminus U(x)$ we can find an open neighborhood $W(z)$ so that $W(z)\subset V(x)$ and $\beta|W(z)=0$. In this case define $a_z=0$.
Since $V(x)$ is a subset of the metrizable object M-polyfold $X$ it is paracompact. Using that $(W(z))$ is an open cover of $V(x)$
we find a subordinate sc-smooth partition of unity $\sigma_z$, $z\in V(x)$ and define the function $\gamma: V(x)\to [0,\infty)$ by 
$$
\gamma (y) =\sum_{z\in V(x)} \sigma_z(y) a_z.
$$
By construction $\gamma \geq \beta$ on $V(x)$. 
Hence $\gamma$ is a  sc-smooth function with support in $\cl_X(U(x))\subset V(x)$. We  average with respect to the $G_x$-action and obtain
the function $f$, defined by 
$$
f(y) =\frac{1}{|G_x|}\sum_{g\in G_x} \gamma (g\ast y).
$$
By construction,  $\beta (g\ast y) \leq \gamma (g\ast y)$ which implies $\beta \leq f$.
\qed \end{proof}

\begin{proof}[Theorem \ref{pert}]
We consider the open covering ${(|U_\lambda|)}_{\lambda\in\Lambda}$ of $|X|$. In a first step we take,  for every $x\in X$,  an open neighborhood $V(x)\subset X$ satisfying the following properties:
\begin{itemize}
\item[(1)]\  There exists an open neighborhood $W(x)$ equipped with the  natural $G_x$-action such that 
$V(x)\subset \cl_X(V(x))\subset W(x)$  and $V(x)$ is invariant.
\item[(2)]\ $t:s^{-1}(\cl_X(V(x))\rightarrow X$ is proper.
\item[(3)]\ For every $x\in X$ there exists a $\lambda$ with $V(x)\subset U_\lambda$.
\item[(4)]\  $( (\pi^{-1}(\abs{V(x)})\setminus V(x))\cap W(x)=\emptyset$.
\end{itemize}
The family $(|V(x)|)_{x\in X}$ is an open covering of $|X|$ refining $(|U_\lambda|)$ and we find a locally finite refinement $(V_i')_{i\in I}$. In particular,  there is a map
$i\rightarrow x_i$ such  that $V_i'\subset \abs{V(x_i)}$. We take,  for every $i\in I$,  an open subset $V_i''$  of $V_i'$ satisfying 
$$
V_{i}''\subset cl_{|X|}(V_i'')\subset V_i'
$$
and, in addition, 
$$
|X| =\bigcup_{i\in I} V_i'' .
$$
For every $i\in I$, $V_i''\subset V_i'\subset |V(x_i)|$,  and we introduce $U_i''= V(x_i)\cap \pi^{-1}(V_i'')$ and   let $U_i$ be its saturation.
Then 
$$
U_i =\pi^{-1}(V_i'').
$$
The open set $U_i''$ is invariant under the $G_{x_i}$ action on $V(x_i)$. Next we take a continuous partition of unity $(\beta_i)$ subordinate
to the locally finite covering $(V_i'')$.  The function   $(\beta_i\circ \pi)|V(x_i)$ has support in $U_i''\subset V(x_i)$. 
Using the existence of smooth partitions of unity for the M-polyfold $X$ we can approximate $ (\beta_i\circ \pi)|V(x_i)$ using Lemma \ref{lemma1.126} by a sc-smooth function $f_i$
supported in $V(x_i)$ satisfying $f_i\geq (\beta_i\circ \pi)|V(x_i)$. We average using the $G_{x_i}$-action,  extend to the saturation of $V(x_i)$,  and then extend by $0$ outside. This gives a function $h_i$ compatible with morphisms.
The sum $h$ of all $h_i$  is locally finite and always bigger than $0$. Consequently,  $(h_i/h)$ is the desired partition of unity for  the ep-groupoid $X$.
\qed \end{proof}

\begin{corollary}\index{C- Partitions of unity}
Let $X$ be an ep-groupoid with paracompact orbit space $|X|$
We assume that the object M-polyfold admits an atlas 
consisting of charts whose local models 
$(O,C,E)$  have the property  that the $E_0$ are Hilbert spaces.  
Then every  open covering of $X$ by saturated open subsets ${(U_\lambda)}_{\lambda\in\Lambda}$ possesses a subordinate
sc-smooth partition of unity of the ep-groupoid $X$.
\end{corollary}
\begin{proof}
We only have to note that the assumption on the atlas implies that the object M-polyfold $X$ admits sc-smooth partitions of unity. The proof can be found in Appendix \ref{POU}, see specifically Proposition \ref{P548},
Theorem \ref{partition_approximation}, and Corollary \ref{C5417}.
Hence, the corollary  follows from Theorem \ref{pert}.
\qed \end{proof}

\chapter{Bundles and Covering Functors}
After discussing the tangent of an ep-groupoid, we shall define strong bundles over ep-groupoids.
Finally we discuss the notion of a proper covering functor between ep-groupoids and strong bundles,
respectively.
\section{The Tangent of an Ep-Groupoid}\label{section_Tangent_Ep-Groupoid}
The {\bf tangent} $T(X, \bm{X})$ of an ep-groupoid $(X,  \bm{X})$ is defined as the pair of 
tangent spaces 
\begin{equation}\label{tangentep_groupoid_eq1}
T(X, \bm{X})=(TX, T\bm{X}).
\end{equation}
Our aim is to equip the pair with the structure of an ep-groupoid, whose object set is the tangent space $TX$ and whose morphism set is the tangent space $T\bm{X}$. The main results are Theorem \ref{main_KK} and Theorem \ref{Walter_xx}.

In order to formulate the theorems we need some preparation. The object M-polyfold is going to be the tangent space $TX$. That $TX$ is a M-polyfold we know from 
Section \ref{section2.1}. The tangent space $T\bm{X}$, again a M-polyfold, is viewed as a the set of morphisms defined as follows. The tangent vector $(\phi, h)\in T_\phi\bm{X}$ is viewed as a morphism 
\begin{equation}\label{earl}
(\phi,h)\colon Ts(\phi)(h)\rightarrow Tt(\phi)(h).
\end{equation}
between the two objects in $TX$.
If $\phi\colon x\to y$ is a morphism in $\bm{X}$, then $Ts(\phi )(h)\in T_xX$ and $Tt(\phi )(h)\in T_yX$.

In view of \eqref{earl} we define the source and 
the target maps $s_{TX}, t_{TX}: T\bm{X}\to TX$ of $(TX, T\bm{X})$ by 
\begin{equation}\label{source-target_TX}
s_{TX}(\phi,h)=Ts(\phi)(h)\quad  \text{and}\quad  t_{TX}(\phi,h)=Tt(\phi)(h).
\end{equation}
In order to introduce the multiplication $(\psi, k)\circ (\phi, h)$ of two morphisms in $T\bm{X}$
we have to assume that 
\begin{equation}\label{source-target_TX_1}
s_{TX}(\psi, k)=Ts(\psi)(k)=Tt(\phi)(h)=t_{TX}(\phi, h)
\end{equation}
which implies that  $s(\psi)=t(\phi)$. Using that the source map $s: \bm{X}\to X$ is a local sc-diffeomorphism, we define the multiplication of the two morphisms 
$(\phi, h)$ and $(\psi, k)$ in $T\bm{X}$ by  

\begin{equation}\label{kaufmann}
(\psi,k)\circ (\phi,h) = (\psi\circ\phi, \bigl(Ts(\psi\circ\phi)\bigr)^{-1}(Ts(\phi)(h)))\in T_{\phi\circ \psi}\bm{X},
\end{equation}
where $k=Ts(\psi)^{-1}\bigl(Tt(\phi)(h)\bigr).$  The notation $(\phi, h)\in T_\phi\bm{X}$ for the tangent vector contains a redundancy.
If $\bm{p}: T\bm{X}\to \bm{X}^1$ is the projection of the tangent bundle, then $\bm{p}(h)=\phi$ if $h\in T\bm{X}$. Therefore, the definition \eqref{kaufmann} for the composition can be expressed shorter as 
\begin{equation}\label{composition_1}
k\circ h = \bigl(Ts(\bm{p}(k)\circ \bm{p} (h))\bigr)^{-1}\bigl(Ts(\bm{p} (h))(h)\bigr)
\end{equation}
and $t_{TX}(h)=s_{TX}(k).$  Depending on the situation we shall or shall not use the short  notation for the composition.

The unit morphism $1_a\in T\bm{X}$ associated with the object $a\in T_xX$ is the morphism 
\begin{equation}\label{unit_TX_1}
1_a=(1_x, h_a)\in T\bm{X},
\end{equation}
where $h_a\in T_{1_x}\bm{X}$ satisfies 
$Ts(1_x)(h_a)=a=Tt(1_x)(h_a).$
Consequently, the unit map $u_{TX}: TX\to T\bm{X}$ is given by 
\begin{equation}\label{unit_TX_2}
u_{TX}(a)=\bigl(1_x, Ts(1_x)^{-1}(a)\bigr)
\end{equation}
for $a\in T_xX$.

If $\iota_X: \bm{X}\to \bm{X}$ is the inverse map of the groupoid $(X, \bm{X})$,  we define  the inverse map $\iota_{TX}: T\bm{X}\rightarrow T\bm{X}$ by
\begin{equation}\label{inversion_TX}
\iota_{TX}(\phi,k) = \bigl(\phi^{-1}, T\iota_X(\phi)(k)\bigr).
\end{equation}

\begin{theorem}\label{main_KK}\index{T- Tangent of ep-groupoids}
The tangent $T(X,\bm{X})=(TX,T\bm{X})$ of the  ep-groupoid $(X,\bm{X})$,  together  with the structure maps introduced above, is an  ep-groupoid. 

If the ep-groupoid $(X, \bm{X})$ is tame, then also the ep-groupoid $T(X, \bm{X})$ is tame. 
Moreover, the two projections $p:TX\rightarrow X^1$ and $\bm{p}:T\bm{X}\rightarrow \bm{X}^1$ together define a sc-smooth functor
$(p,\bm{p}): T(X,\bm{X})\rightarrow (X,\bm{X})^1$.
\end{theorem}
The  theorem will follow from a series of lemmata.
By definition,  the multiplication map of the 
ep-groupoid $(X, \bm{X})$,
$$m: \bm{X}{_{s}\times_t}\bm{X}\to \bm{X},\quad m(\psi, \phi)=\psi\circ \phi,$$
where $s(\psi)=t(\phi)$, is a sc-smooth map and we can pass to its tangent map 
$Tm$. To derive a formula for $Tm$ we recall that  $\bm{X}{_{s}\times_t}\bm{X}$ is a sub-M-polyfold of $\bm{X}\times \bm{X}$ so that the tangent space $T(\bm{X}{_{s}\times_t}\bm{X})$ can be identified
with a sub-M-polyfold of $T\bm{X}\times T\bm{X}$ by the identification
$$T(\bm{X}{_{s}\times_t}\bm{X})
=T\bm{X}{_{Ts}\times_{Tt}}T\bm{X},$$
so that  
$Ts(\psi )(k)=Tt(\phi)(h)$ for 
$(\psi, k), (\phi, h)\in T(\bm{X}{_{s}\times_t}\bm{X}).$

\begin{lemma}\label{multiplication}
If $m$ is the multiplication  of isomorphisms in $\bm{X}$, then the multiplication of morphisms 
in $T\bm{X}$, defined in \eqref{kaufmann}, 
can be written as 
$$(\psi, k)\circ (\phi, h)=Tm(\psi, \phi)(k, h)=
\bigl(\psi\circ \phi, Ts(\psi\circ\phi))^{-1}\circ Ts(\phi)(h) \bigr),
$$
keeping in mind that $k=Ts(\psi)^{-1}(Tt(\phi)(h)).$
\end{lemma}
\begin{proof}
In order to derive a formula for $Tm(\psi,\phi)(k, h)$ for $k\in T_\psi \bm{X}$ and $h\in T_\phi \bm{X}$, we vary $\phi'$ in a neighborhood $\bm{U} (\phi)\subset \bm{X}$ and $\psi'$ in a neighborhood $\bm{U} (\psi)\subset \bm{X}$ where $\phi'$ and $\psi'$ are restricted by $s(\psi')=t(\phi')$. Therefore, $\psi'$ is the function $F$ of $\phi'$ defined  by 
$$\psi'=F (\phi')=\bigl(s\vert \bm{U} (\psi )\bigr)^{-1}(t(\phi')).$$
Thus the multiplication becomes
$$m(\psi', \phi')=m\bigl(\bigl(s\vert \bm{U} (\psi )\bigr)^{-1}(t(\phi')), \phi'\bigr).$$
Applying the (restricted) source map leads to the identity 
$$\bigl(s\vert \bm{U} (\phi \circ \psi)\bigr)\circ m\bigl(\bigl(s\vert \bm{U} (\psi )\bigr)^{-1}(t(\phi')), \phi'\bigr)=s(\phi').$$

Differentiating  the identity at  $\phi\in \bm{X}_1$ in the direction $h\in T_\phi\bm{X}$,  we obtain
$$
Ts(\psi\circ\phi)\bigl((Tm(T(s\vert \bm{U} (\psi))^{-1}(Tt(\phi)(h)\bigr),h)=Ts(\phi)(h).
$$
Hence, again using that $s$ is a local sc-diffeomorphism,
$$
Tm\bigl((T(s|U(\psi))^{-1}(Tt(\phi)(h)\bigr),h)= (Ts(\psi\circ\phi))^{-1}\circ Ts(\phi)(h).
$$
Recalling  $Ts(\psi)(k)=Tt(\phi)(h)$ we finally obtain
$$
Tm(k,h)=(Ts(\psi\circ\phi))^{-1}\circ Ts(\phi)(h).
$$
and the lemma follows.
\qed \end{proof}

We abbreviate the morphisms in $T\bm{X}$ as $\Phi=(\phi,h)$ and $\Psi=(\psi, k)$.
\begin{lemma}\label{lemma1.54}
The definitions of $s_{TX}$ and $t_{TX}$ are compatible with the composition, i.e.,
$$
s_{TX}(\Psi\circ\Phi)=s_{TX}(\Phi)\ \text{and}\ t_{TX}(\Psi\circ\Phi)=t_{TX}(\Psi).
$$
\end{lemma}
\begin{proof}
Using the definitions \eqref{source-target_TX_1} and \eqref{kaufmann}, we compute
\begin{equation*}
\begin{split}
s_{TX}((\psi,k)\circ (\phi,h))=
&=Ts(\psi\circ\phi)(T(s(\psi\circ\phi))^{-1}(Ts(\phi)(h)))\\
&=Ts(\phi)(h)\\
&=s_{TX}(\phi,h).
\end{split}
\end{equation*}
We have verified that $s_{TX}(\Psi\circ \Phi)=s_{TX}(\phi).$ In order to verify the formula for $t_{TX}$ we 
fix suitable open neighborhoods $\bm{U} (\psi\circ\phi)$ and $U(x)$ and $U(y)$ with $x=s(\phi)$ and $y=t(\psi)$
so that $t: \bm{U}(\psi\circ\phi)\rightarrow U(y)$ and $: \bm{U}(\psi\circ \phi)\rightarrow U(x)$ are sc-diffeomorphisms.
We define the map $\alpha$ by $\alpha(b)=t\circ (s\vert \bm{U} (\psi\circ\phi))^{-1}(b)$.  Then taking the tangent of $\alpha$ at $s(\phi)$ gives
$$
T\alpha(s(\phi)) = Tt(\psi\circ\phi)\circ (Ts(\psi\circ\phi))^{-1}(s(\phi)).
$$
Now we consider for $\phi'$ near $\phi$ the map
$$
\phi'\mapsto  \alpha\circ s(\phi').
$$
In view of $s(\psi\circ \phi )=s(\phi)$ and 
$t(\psi\circ \phi )=t(\psi)$ we have the identity
\begin{equation*}
\begin{split}
\alpha(s(\phi')) &=t\circ (s\vert \bm{U} (\psi\circ\phi))^{-1}(s(\phi'))\\
&= t\circ (s\vert \bm{U} (\psi))^{-1}(t(\phi')).
\end{split}
\end{equation*}
Therefore,  taking the derivative at $\phi'=\phi$, we obtain
$$
T(\alpha\circ s)(\phi)(h) = Tt(\psi)(T((s|U(\psi))^{-1}(t(\phi))(Tt(\phi)h)).
$$
Using  the definition \eqref{source-target_TX} 
and recalling $Ts(\psi)(k)=Tt)\phi )(h)$, we compute,
\begin{equation*}
\begin{split}
t_{TX}((\psi,k)\circ (\phi,h))
&=Tt(\psi\circ\phi)((Ts(\psi\circ\phi))^{-1}(Ts(\phi)(h)))\\
&=T\alpha(s(\phi))((Ts(\phi)(h))\\
&=T(\alpha\circ s)(\phi)(h)\\
&=Tt(\psi)(T(s\vert \bm{U} (\psi)))^{-1}Tt(\phi)(h)\\
&=Tt(\psi)(T(s\vert \bm{U} (\psi)))^{-1}Ts(\psi)(k)\\
&=Tt(\psi)(k)\\
&=t_{TX}(\psi,k),
\end{split}
\end{equation*}
as claimed in Lemma \ref{lemma1.54}.
\qed \end{proof}
\begin{lemma}
The composition of morphisms in  $T\bm{X}$  is associative. 
\end{lemma}
\begin{proof}
Assuming that the corresponding sources and targets match we compute,
\begin{equation*}
\begin{split}
&\phantom{=;;}(\sigma,l)\circ ((\psi,k)\circ(\phi,h))\\
&=(\sigma,l)\circ (\psi\circ\phi,(Ts(\psi\circ\phi))^{-1}Ts(\phi)(h))\\
&=(\sigma\circ\psi\circ\phi, (Ts(\sigma\circ\psi\circ\phi))^{-1}Ts(\psi\circ \phi)(Ts(\psi\circ\phi))^{-1}Ts(\phi)(h))\\
&=(\sigma\circ\psi\circ\phi, (Ts(\sigma\circ\psi\circ\phi))^{-1}Ts(\phi)(h))\\
&=(\sigma\circ\psi,(Ts(\sigma\circ\psi))^{-1}Ts(\psi)(k))\circ(\phi,h)\\
&=((\sigma,l)\circ(\psi,k))\circ(\phi,h).
\end{split}
\end{equation*}
\qed \end{proof}

\begin{lemma}\label{lemma1_58}
The $1$-morphism associated with the object $q\in T_xX$, defined as 
$$
1_q=(1_x, Ts(1_x)^{-1}q)\in T_{1_x}\bm{X},
$$
is a $2$-sided identity for the composition.
\end{lemma}
\begin{proof}
We first  observe that  $1_q$ is a $1$-morphisms. Indeed, 
$$
s_{TX}(1_q)=Ts(1_x)(Ts(1_x)^{-1}q)=q\quad  \text{and}\quad  t_{TX}(1_q)=Tt(1_x)((Ts(1_x))^{-1}(q))=q.
$$
In the last equation we have used that 
$t(1_b)=b=s(1_b)$ for all $b\in X$, and therefore, since $t$ and $s\colon \bm{X}\to X$ are local sc-diffeomorphisms, we have the identity $t\circ (s\vert \bm{U} (1_x))^{-1}(b)=b$ for all $b$ near $x$, and by the chain rule
$Tt(s^{-1}(b))\circ Ts(s^{-1}(b))^{-1}=1$, so that $Tt(1_x)\circ Ts(1_x)^{-1}=1$ at $b=x$.

Next we assume that the morphism
 $(\phi,k)$ satisfies  $s_{TX}(\phi,k)=t_{TX}(1_q)$ which is equivalent to  $Ts(\phi)(k)=Tt(1_x)\circ Ts(1_x)^{-1}(q)=q$. Then 
\begin{equation*}
\begin{split}
(\phi, k)\circ 1_q&=(\phi,k)\circ(1_x,Ts(1_x)^{-1}q)\\
&=(\phi\circ 1_x, Ts(\phi)^{-1}Ts(1_x)Ts(1_x)^{-1}q)\\
&=(\phi, Ts(\phi)^{-1}(q))=(\phi,k).
\end{split}
\end{equation*}
If $(\phi,k)$ satisfies $t(\phi)=x$ and $t_{TX}(\phi,k)=s_{TX}(1_q)$,  which is equivalent to  $Tt(\phi)(k)=q$, we obtain
\begin{equation*}
1_q\circ (\phi,k)=(\phi, Ts(\phi)^{-1}(Ts(\phi)(k)))=(\phi,k),
\end{equation*}
as claimed. 
\qed \end{proof}
\begin{lemma}\label{lemma1_59}
The inverse morphism $\iota_{TX}$  of the morphism $(\phi, k)\in T\bm{X}$, as defined by 
$$
\iota_{TX}(\phi,k) = (\phi^{-1}, T\iota_X(\phi)(k)),
$$
is a $2$-sided inverse for the composition.
\end{lemma}
\begin{proof}
We first verify that $\iota_{TX}(\phi,k)$ is an 
inverse morphism of $(\phi, k)$. Indeed, from the identities $s(\phi^{-1})=t(\phi)$ and 
$t(\phi^{-1})=s(\phi)$ for $\phi\in \bm{X}$ 
we obtain by the chain rule $Ts(\phi^{-1})\circ T\iota_X (\phi)=Tt(\phi)$ and 
$Tt(\phi^{-1})\circ T\iota_X (\phi)=Ts(\phi)$. 
Consequently, 
$$
s_{TX}(\iota_{TX}(\phi,k))=Ts(\phi^{-1})(T\iota_{X}(\phi)(k))=Tt(\phi)(k)=t_{TX}(\phi,k)
$$
and
$$
t_{TX}(\iota_{TX}(\phi,k))=Tt(\phi^{-1})(T\iota(\phi)(k))=Ts(\phi)(k)=s_{TX}(\phi,k).
$$
For a given morphism $(\phi,k)\in T_\phi \bm{X}$  satisfying  $x=s(\phi)$ we compute, 
\begin{equation*}
\begin{split}
\iota_{TX}(\phi,k)\circ (\phi,k)
&=(\phi^{-1},T\iota_X(\phi)k)\circ(\phi,k)\\
&=(1_x,(Ts(1_x))^{-1}Ts(\phi)k)\\
&= 1_{Ts(\phi)k}.
\end{split}
\end{equation*}
The composition from the right is verified the same way and we have proved that  $\iota_{TX}$ is the $2$-sided inverse for the composition.
\qed \end{proof}

Next we prove the differential geometric properties for the category $TX$ consisting of the object set $TX$ and the morphism set $T\bm{X}$.

\begin{lemma}
For the category with object M-polyfold  $TX$  and morphism M-polyfold $T\bm{X}$ the source and target maps $s_{TX}$ and $t_{TX}$
are surjective local sc-diffeomorphisms and all structure maps are sc-smooth.
\end{lemma}
\begin{proof}
By construction,  $s_{TX}$ and $t_{TX}$ are  local sc-diffeomorphisms. The $1$-map $TX\rightarrow T\bm{X}$ associates with  $q\in T_xX$
the pair $1_q=(1_x,(Ts(1_x))^{-1}(q))$, which sc-smoothly depends on $q$. The inversion map $T\bm{X}\rightarrow T\bm{X}$
is the tangent map of  the sc-smooth inversion map $\bm{X}\rightarrow \bm{X}:\phi\rightarrow \phi^{-1}$ and therefore a sc-diffeomorphism.
Due to the properties of $s_{TX}$ and $t_{TX}$ the fibered product $T\bm{X}{_{s_{TX}}\times_{t_{TX}}}T\bm{X}$ is a M-polyfold. The multiplication map
is given by
$$
T\bm{X}{_{s_{TX}}\times_{t_{TX}}}T\bm{X}\rightarrow T\bm{X}: ((\psi,k),(\phi,h))\mapsto (\psi\circ \phi, (Ts(\psi\circ\phi))^{-1} Ts(\phi)(h)).
$$
which is sc-smooth by inspection. One can also view this  as follows. The tangent space of $\bm{X}{_{s}\times_t}\bm{X}$
can naturally be identified with $T\bm{X}{_{s_{TX}}\times_{t_{TX}}}T\bm{X}$ and our multiplication map $m_{TX}$ is precisely $Tm$, where
$m$ is the multiplication for the original ep-groupoid which  is sc-smooth. This completes the proof of the lemma.
\qed \end{proof}
In order to verify  that $(TX,T\bm{X})$ is an ep-groupoid we still need to confirm  the properness property.
\begin{lemma}
If  $(X,\bm{X})$ be an ep-groupoid, then $(TX,T\bm{X})$ has the properness property. 
\end{lemma}
\begin{proof}
Let $x\in X_1$ and take an open neighborhood $V(x)$ in $X$ such  that 
\begin{equation}\label{pol_proper_x}
t:s^{-1}(\cl_X(V(x)))\rightarrow X\ \text{is proper. }
\end{equation}
In order to show that the map 
\begin{equation}\label{pol_proper}
t_{TX}:s^{-1}(\cl_{TX}(TV))\rightarrow TX
\end{equation}
is proper, we may assume that we have a sequence of morphisms $((\phi_k,h_k))\subset T\bm{X}$ satisfying  $s_{TX}(\phi_k,h_k)\in \cl_{TX}(TV)$ and $(t(\phi_k,h_k))$
lying in a compact subset $K\subset TX$.  Due to the latter, after perhaps taking a subsequence, we may assume  the convergence 
$$
t(\phi_k,h_k)=Ts(\phi_k)(h_k)\rightarrow q\in T_xX.
$$
We note that 
\begin{equation}\label{pol_proper_t}
t(\phi_k)=:x_k\rightarrow x\quad \text{in\ $  X_1$},
\end{equation}
and,    moreover,  
\begin{equation}
t(\phi_k)\rightarrow x\quad \text{in\  $  X_0$.}
\end{equation}
Since $s(\phi_k)\in \cl_X(V)$,  we deduce from the properness in (\ref{pol_proper_x}) that
\begin{equation*}
\phi_k\rightarrow \phi\quad  \text{in\ ${\bf  X}$.}
\end{equation*}
Because $t$ is a sc-diffeomorphism between suitable open neighborhoods of $\phi$ and $x$, the convergence of $(x_k)$ in \eqref{pol_proper_t}
implies the stronger convergence
$$
\phi_k\rightarrow \phi\ \text{in}\ \bm{X}_1.
$$
At this point we know about our sequence $((\phi_k,h_k))$ that $\phi_k\rightarrow \phi$ in $\bm{X}_1$ and $Ts(\phi_k)(h_k)\rightarrow q\in T_xX$ in $TX$.
We take  open neighborhoods $\bm{U}(\phi)\subset \bm{X}$ and $U(x)\subset X$ such  that $s:\bm{U}(\phi)\rightarrow U(x)$ is a sc-diffeomorphism.
Then $Ts\colon T\bm{U}(\phi)\rightarrow TU(x)$ is a sc-diffeomorphism. Since for large $k$ we have $(\phi_k,h_k)\in T\bm{U}(\phi)$ and  the sequence $Ts(\phi_k,h_k)=q_k$ converges,
the same is  true for $(\phi_k,h_k)$. Hence $(\phi_k,h_k)\rightarrow (\phi,h)$ in $T\bm{X}$.  This completes the proof of the properness.
\qed \end{proof}

So far we have proved that the small category $(TX, T\bm{X})$ is an ep-groupoid, 
if $(X, \bm{X})$ is an ep-groupoid. 

\begin{lemma}\label{lemma1_62}
If the ep-groupoid $(X, \bm{X})$ is tame, so is the ep-groupoid $(TX, T\bm{X})$. 
\end{lemma}
\begin{proof}
If $\varphi\colon U\to O$, $O=(O, C, E)$, is a tame chart of the object $M$-polyfold $X$ of $(X, \bm{X})$ and $O=r(U)\subset C$ is a tame retract, then also $Tr\colon TU\to TU$ is a tame retraction. It follows that if $(X, \bm{X})$ has a tame object space, then also $(TX, T\bm{X})$ has a tame object space and hence is, by definition, a tame ep-groupoid. 

\qed \end{proof}

\begin{proof}[Proof of Theorem \ref{main_KK}]
The category $(TX, T\bm{X})$ is in  view of Lemmata \ref{multiplication} - \ref{lemma1_62} 
   a (tame) ep-groupoid, if $(X, \bm{X})$ is a (tame) ep-groupoid. The last statement about the projections being functors just says that if $$
   (\phi,h)\colon Ts(\phi)(h)\rightarrow Tt(\phi)(h),
   $$
then $\bm{p}(\phi,h)=\phi:s(\phi)\rightarrow t(\phi)$ is a morphism compatible with composition.
This  is obvious. The proof of Theorem \ref{main_KK} is complete.
\qed \end{proof}
There are equivalent ways to  equip  the tangent  $(TX, T\bm{X})$ of the ep-groupoid 
with a ep-groupoid structures.
We discuss another possibility in the following remark.
\begin{remark}\index{R- $TX$ as ep-groupoid}
The above definition of $TX$ as an ep-groupoid makes precise in which way we  view the elements in
$T\bm{X}$ as morphisms. There are equivalent ways to turn $TX$ into an ep-groupoid. Recall that a morphism
$\phi:x\rightarrow y$ for $x\in X_1$ (equivalently $y\in X_1$) has a well-defined associated linear isomorphism
$T\phi:T_xX\rightarrow T_yY$. Here $T\phi$ is defined as the tangent map  of $t\circ (s|U(\phi))^{-1}$ at $x$.
One interprets the fact that if $T\phi$ maps $h\in T_xX$ to $k\in T_yY$ as having a morphism $(T\phi,h)$ with
$s(T\phi,h)=h$ and $t(T\phi,h)=T\phi(h)$. The composition  $(T\psi,h)\circ (T\phi,l)$ should be defined,  
provided $h=T\phi(l)$, as follows, 
$$
(T\psi,h)\circ (T\phi,l)=(T\psi\circ T\phi,l)=(T(\psi\circ\phi), l).
$$
We note that this new definition and the previous one are equivalent  by means of the mapping
$$
G: (\phi,h)\rightarrow (T\phi,Ts(h)).
$$
Indeed, $G$ is compatible with the functorial structures. For example, if  $Tt(\phi)(k)=Ts(\psi)(h)$, 
\begin{equation*}
\begin{split}
G((\psi,h)\circ(\phi,k))
&=(T(\psi\circ\phi),Ts((Ts(\psi\circ\phi))^{-1}(Ts(\phi)(k))))\\
&=(T(\psi\circ\phi),Ts(\phi)(k))\\
&=(T\psi\circ T\phi,Ts(\phi)(k))\\
&=(T\psi,T\phi(Ts(\phi)(k)))\circ (T\phi,Ts(\phi)(k))\\
&=(T\psi,Tt(\phi)(k))\circ(T\phi,Ts(\phi)(k))\\
&=(T\psi,Ts(\psi)(h))\circ(T\phi,Ts(\phi)(k))\\
&=G(\psi,h)\circ G(\phi,k).
\end{split}
\end{equation*}

In the latter description the object space is $TX$ as before , but the morphisms are the pairs $(T\phi,h)$, where $\phi\in \bm{X}_1$
and  $h\in T_{s(\phi)}X$.  The map 
$$
G: T\bm{X}\rightarrow \{(T\phi,a)\ |\ \phi\in \bm{X}_1,\ a\in T_{s(\phi)}X\}, 
$$
is a bijection and we can equip the right-hand side with a M-polyfold structure such that 
$G$ becomes a sc-diffeomorphism. Of course we can identify the pairs $(T\phi,a)$
with the triplets  $(a,T\phi,b)$, where $b=T\phi(a)$. This gives another description of the morphism set.
Some further remarks are in order.  For  $\phi\in \bm{X}_1$ we abbreviate  $x=s(\phi)$ and $y=t(\phi)$ in $X_1$.
The linear isomorphism
$$
Tt(\phi)\circ Ts(\phi)^{-1}:T_xX\rightarrow T_yX.
$$
is  the linearization of $t\circ (s|U(\phi))^{-1}$ at $\phi$ and therefore,  by definition,  equal to $T\phi$. 
Using the definitions of $s_{TX}$ and $t_{TX}$,  we find that
$$
t_{TX}\circ (s_{TX}|T_xX)^{-1} = T\phi \colon T_xX\rightarrow T_yX
$$
 is a linear isomorphism. 
\end{remark}\qed

At this point we have constructed the tangent for the ep-groupoids, but not yet for the sc-smooth functors between them. 

\begin{theorem}\label{theorem1.63}\index{T- Tangent of functors}
If $(\Phi, \bm{\Phi})\colon (X, \bm{X})\rightarrow (Y, \bm{Y})$ is  a sc-smooth functor between two ep-group\-oids, then the tangent map 
$$
(T\Phi,T\bm{\Phi})\colon T(X,\bm{X})\rightarrow T(Y,\bm{Y})
$$
is a sc-smooth functor between the corresponding tangents of the ep-group\-oids.
\end{theorem}
\begin{proof}
Clearly, since $\Phi:X\rightarrow Y$ and $\Phi:\bm{X}\rightarrow \bm{Y}$ are sc-smooth maps, it follows
that the tangent maps $T\Phi:TX\rightarrow TY$ and $T\Phi:T\bm{X}\rightarrow T\bm{Y}$ are sc-smooth.  Let us consider the compatibility with 
the composition we have defined for the elements in $T\bm{X}$ and $T\bm{Y}$, respectively.
The functoriality of $\Phi$ requires that 
$$
\Phi\times\Phi:\bm{X}\times \bm{X}\rightarrow \bm{Y}\times \bm{Y}
$$
preserves the fibered products on which the multiplications  are  defined. Consequently,  
if  $(\psi,\phi)\in \bm{X}{_{s}\times_t}\bm{X}$,  then 
$$
\Phi(m_X(\psi,\phi))=m_Y(\Phi(\psi),\Phi(\phi)).
$$
Differentiating with respect to the pair  $(\psi,\phi)\in \bm{X}{_{s}\times_t}\bm{X}$ in the direction $(k,h)\in T\bm{X}{_{Ts}\times_{Tt}}T\bm{X}$
we obtain
$$
T\Phi\circ Tm_X(k,h)= Tm_Y(T\Phi(k),T\Phi(h)).
$$
Therefore, by  Lemma \ref{multiplication},
$$
T\Phi(\psi\circ\phi)(k\circ h) = (T\Phi(\psi)(k))\circ (T\Phi(\phi)(h)),
$$
which is the desired result. To verify  that $T\Phi$ preserves the units,  is left to the reader.

For the  morphisms  $\alpha$ we have the identities
$$
s(\Phi(\alpha))=\Phi(s(\alpha))\ \text{and}\ \ t(\Phi(\alpha))=\Phi(t(\alpha)).
$$
We differentiate at $\alpha_0\in \bm{X}_1$ in the direction $p\in T_{\alpha_0}\bm{X}$ and obtain, setting  $x_0=s(\alpha_0)$ and $x_0'=t(\alpha_0)$, 
\begin{align*}
s_{TY}(T\Phi(\alpha_0)(p))&=T\Phi(x_0)(s_{TX}(p))\\
t_{TY}(T\Phi(\alpha_0)(p))&=T\Phi(x_0')(t_{TX}(p)).
\end{align*}
In other words,  if $h\in T_{x_0}X$ and $k\in T_{x_0'}X$ are objects with a morphism
$(\phi,p)$ between them, i.e. $p\in T_\phi\bm{X}$  satisfying $Ts_X(p)=h$ and $Tt_X(p)=k$, then
$T\Phi(\phi)(p)$ is the corresponding morphism $T\Phi(x_0)(h)\rightarrow T\Phi(x_0')(k)$.
\qed \end{proof}

Recalling  the {\bf category 
$\mathcal{EP}$}  whose class of objects 
are the ep-groupoids and whose morphisms are the sc-smooth functors between the ep-groupoids, we can sum up the previous discussion in the following theorem.

\begin{theorem}[{\bf Functoriality property of $T$ on ep-groupoids}]\label{Walter_xx}\index{T- Functoriality of $T$ on ep-groupoids}
The tangent functor $T$ gives rise to a functor 
$$T\colon \mathcal{EP}\to \mathcal{EP}$$
associating with  an ep-groupoid $X$ its tangent $TX$ and with  a smooth functor $\Phi:X\rightarrow Y$ the tangent functor $T\Phi$.  In particular,  $T\Phi$ preserves the composition and the other structure maps. Furthermore,  $T$ preserves tameness; if $X$ is a tame ep-groupoid so is $TX$.
\qed
\end{theorem}

Finally we shall show that  the tangent functors are compatible with natural transformations.
\begin{definition}[{\bf Natural transformation, natural equivalence}]\label{natural_transf_def}\index{D- Natural transformation}\index{Natural equivalence}
If  
$$
F\colon  X \to Y\ \ \text{and}\ \  G\colon  X \to Y
$$
 are two sc-smooth functors  between the ep-
groupoids $X$ and $Y$, then  $F$  is  called {\bf naturally equivalent} to $G$ and denoted by $F\simeq G$,  if there exists a sc-smooth map
$\tau \colon X\to \bm{Y}$  which associates with every object $x\in X$   a morphism 
$$\tau (x)\colon F(x)\to G(x)$$ 
in $\bm{Y}$ such that 
$$\tau (x')\circ \bm{F}(h)=\bm{G}(h)\circ \tau (x)
$$
for every morphism $h\colon x\to x'$ in $\bm{X}$, i.e.,  
the diagram of morphisms 
$$
\begin{CD}
F(x)@>\tau(x)>> G(x)\\
@VV \bm{F}(h )V  @VV \bm{G}(h)V\\
F(x')@>\tau(x')>>G(x')
\end{CD}
$$
commutes for every morphism $h\colon x\to x'$ in $\bm{X}$.

The sc-smooth map $\tau \colon X\to \bm{Y}$ is called a {\bf natural transformation} between the two sc-functors $F, G\colon  X \to Y$ and will be referred to by the formula 
$$\tau: F\to G.$$
\qed
\end{definition}

The relation $F\simeq G$  between sc-smooth functors $F, G\colon  X \to Y$ of ep-groupoids is an equivalence relation.

\begin{proposition}\label{169}\index{P- Tangent of a natural transformation}
If $\tau : \Phi\rightarrow \Psi$  is a natural transformation  between the two sc-smooth functors 
$$\Phi,\Psi\colon X= (X, \bm{X})\rightarrow Y=(Y, \bm{Y})$$  between ep-groupoids, then 
the tangent  $T\tau\colon TX\rightarrow T\bm{Y}$
is a natural transformation 
$T\tau : T\Phi\rightarrow T\Psi$ between the two tangent functors 
$$
T\Phi, T\Psi\colon (TX, T\bm{X})\to (TY, T\bm{Y}).
$$ 
\end{proposition}
\begin{proof}
If $x\in X$, then 
$\tau (x)\colon \Phi(x)\rightarrow\Psi(x)$ is a morphism in $\bm{Y}$. We note that $s(\tau(x))=\Phi(x)$ and $t(\tau(x))=\Psi(x)$.
Differentiating $\Phi(x)=s(\tau(x))$ and $\Psi(x)=t(\tau(x))$ at $x=x_0\in X_1$ yields,  with $\psi_0=\tau(x_0)$, 
$$
T\Phi(x_0)(h) =Ts(\psi_0)(T\tau(x_0)(h))\ \text{and}\ \ T\Psi(x_0)(h)=Tt(\psi_0)(T\tau(x_0(h)).
$$
Consequently,  $s_{TY}(T\tau(x_0)(h))= T\Phi(x_0)(h)$ and $t_{TY}(T\tau(x_0)(h))=T\Psi(x_0)(h)$.  By definition of the morphisms
in $T\bm{Y}$ we have the morphism 
$$
T\tau(x_0)(h):T\Phi(x_0)(h)\rightarrow T\Psi(x_0)(h).
$$
Since $\tau$ is a natural transformation,  the diagram
$$
\begin{CD}
\Phi(x)@>\tau(x)>> \Psi(x)\\
@VV\Phi(\alpha )V  @VV\Psi(\alpha)V\\
\Phi(x')@>\tau(x')>>\Psi(x')
\end{CD}
$$
commutes for all morphisms $\alpha\colon x\to x'$. 
We  differentiate the relationship defined  by the diagram as follows. We fix $x_0, x_0'\in X_1$ and the morphism $\alpha_0\colon x_0\rightarrow x_0'$.
The latter then also belongs to level $1$. For the morphisms $\alpha$ near $\alpha_0$ on level $0$ we obtain the commuting diagrams
$$
\begin{CD}
\Phi(s(\alpha))@>\tau(s(\alpha))>> \Psi(s(\alpha))\\
@VV\Phi(\alpha )V  @VV\Psi(\alpha)V\\
\Phi(t(\alpha))@>\tau(t(\alpha))>>\Psi(t(\alpha)).
\end{CD}
$$
We differentiate at $\alpha=\alpha_0$ satisfying $s(\alpha_0)=x_0$ and $t(\alpha_0)=x_0'$ and obtain,  for every  $p\in T_{\alpha_0}\bm{X}$, $h=Ts_X(\alpha_0)p$, and $h'=Tt_X(\alpha_0)p$, the commutative diagram
$$
\begin{CD}
T\Phi(x_0)(h)@>T\tau(x_0)h>>T \Psi(x_0)h\\
@VVT\Phi(\alpha_0 )pV  @VVT\Psi(\alpha_0)p V\\
T\Phi(x_0')h' @>T\tau(x_0)h>>T\Psi(x_0')h'.
\end{CD}
$$

Recalling  that  $T\Phi(\alpha_0)p$ is the morphism $T\Phi(x_0)h\rightarrow T\Phi(x_0')h'$ in $T\bm{Y}$ and similarly for $\Psi$, the diagram
above  shows that the tangent map $T\tau: TX\rightarrow T\bm{Y}$, $(\phi,h)\mapsto  T\tau(x_0)(h)$ is a natural transformation $T\tau : T\Phi \to T\Psi$.
This completes the proof of Proposition \ref{169}.
\qed \end{proof}

If $X$ is an ep-groupoid with paracompact orbit space $|X|$  we have seen that $|X|$ is metrizable.  The tangent $TX$ of an ep-groupoid is also an ep-groupoid
as we have seen in Theorem \ref{main_KK}. As we shall show now the paracompactness of $|X|$ is inherited by $|TX|$.
\begin{theorem}\label{paraXCV} \index{T- Inherited paracompactness of $\abs{TX}$}
If $X$ is an ep-groupoid with paracompact orbit space $|X|$, then the tangent ep-groupoid $TX$  has a parcompact orbit space $|TX|$ as well.
\end{theorem} 
\begin{proof}
By assumption $|X|$ is paracompact and the same holds for the ep-groupoid $X^1$, i.e. $|X^1|$ is paracompact.
For every point $z\in |X^1|$ take a representative $x_z\in X^1$ and pick an open neighborhood $U(x_z)$ in $X$ with the natural $G_{x_z}$-action  so that 
$TU(x_z)$ is as a bundle over $U(x_z)^1$ sc-isomorphic to some local bundle model $TO\rightarrow O^1$, which also is equipped with an sc-smooth bundle action by $G_{x_z}$. 
Then $|TU(x_z)|$ and $|U(x_z)^1|$ are metrizable. The collection of sets ${(|U(x_z)^1|)}_{z\in |X^1|}$ is an open covering of $|X^1|$ and we can find a subordinate 
locally finite closed subcover denoted by ${(A_z)}_{z\in\Sigma}$, where $\Sigma\subset |X^1|$ consists of those points for which $A_z\neq \emptyset$. 
For $z\in\Sigma$ we find a  $G_{x_z}$-invariant closed subset of $U(x_z)^1$ and take the restriction of $TU(x_z)$ to this subset.
This defines, after dividing out by the $G_{x_z}$ a closed subset $\wt{A}_z$ of $|TX|$ which is closed and paracompact. Moreover
${(\wt{A}_z)}_{z\in\Sigma}$ is a locally finite covering of $|TX|$ by closed paracompact sets. Since as an ep-groupoid $|TX|$ is a regular Hausdorff space, it follows
that $|TX|$ is paracompact in view of the topological Proposition \ref{union_of_paracompact}. This completes the proof.
\qed \end{proof}

\section{Sc-Differential Forms on Ep-Groupoids}\label{SECT82}
Next we introduce the notion of a sc-differential form in the ep-groupoid context. It will be studied in later sections in more detail.
The notion has been already introduced in the context of M-polyfolds in Section \ref{subs_sc_differential}.
 \begin{definition}\index{D- Sc-Differential form}
 A {\bf sc-differential k-form on the ep-groupoid} \index{sc-differential k-form on the ep-groupoid}  $(X, \bm{X})$ is a  sc-smooth map $\omega:\oplus_k TX\rightarrow {\mathbb R}$ defined on the Whitney sum
 of the tangent of the object M-polyfold, which is
 linear in each argument and skew-symmetric,  and, moreover, possesses  the additional property that for every morphism $\phi: x\rightarrow y$ in $\bm{X}_1$, and its associated tangent $T\phi:T_xX\rightarrow T_yY$ it holds that 
 $$
 (T\phi)^\ast\omega_y=\omega_x.
 $$
\qed
 \end{definition}
   We denote by $\Omega_{ep}^k(X)$\index{$\Omega_{ep}^k(X)$} the vector space of sc-differential $k$-forms compatible with morphisms.
  If $\omega $ is a sc-differential form on $X^i$, then  $\omega$  will also be a differential form on $X^{i+1}$.  We can take the direct limit
 of the directed system
$$
\iota_{i,k}:\Omega^\ast_{ep}(X^i)\rightarrow\Omega^\ast_{ep}(X^{k})
 $$
 with $i\leq k$, induced by the pull-back via inclusions
 and obtain $\Omega^{\ast}_{ep,\infty}(X)$\index{$\Omega^{\ast}_{ep,\infty}(X)$}. With ${\mathbb N}=\{0,1,2,\ldots \}$,  the latter is defined by
 $$
 \Omega^{\ast}_{ep,\infty}(X)=\left(\bigoplus_{i\in {\mathbb N}} \Omega^\ast_{ep}(X^i)\right)/ D,
 $$
 where $D$ is the vector subspace generated by $\omega_i - \iota_{i, k}(\omega_i)$.
By definition,  the elements $[\omega]$ in  $\Omega^{\ast}_{ep,\infty}(X)$ are equivalence classes and can be represented 
by a differential form $\omega_i$ on $X^i$ for some $i$. Then $\omega_i$ and $\omega_k$ present the same element provided there exists an $l\geq k,i$ 
such  that the pullbacks of $\omega_i$ and $\omega_k$ to $X^l$ are the same.  
We define the exterior derivative of $[\omega]$ by
$$
d[\omega]:=[d\omega],
$$
where $d\omega$ is the exterior derivative of a sc-differential form on the object M-polyfold. 
We claim that 
$d[\omega]\in \Omega^{\ast}_{ep,\infty}(X)$, i.e. $d\omega$ is compatible with morphisms. This can be seen as follows.
  Given a morphism $\phi_0\colon x_0=s(\phi_0)\to y_0=t(\phi_0)\in \bm{X}_1$, we  choose  open neighborhoods $U(x_0)\subset X$ and $\bm{U}(\phi_0)\subset \bm{X}$  such  that $s:\bm{U}(\phi_0)\rightarrow U(x_0)$ and $t:U(\phi_0)\rightarrow U(y_0)$ are sc-diffeomorphisms. Defining the composition 
 $$
 f:U(x_0)\rightarrow U(y_0): x\rightarrow t((s\vert \bm{U}(\phi_0))^{-1}(x)), 
 $$
we note that with $\phi\in U(\phi_0)$, $x=s(\phi)\in U(x_0)\cap X_1$ and $y=t(\phi_0)$ the linear isomorphism
 $T\phi:T_xX\rightarrow T_yX$ is equal to $Tf(x)$.
 Therefore, choosing $x\in U(x_0)\cap X_1$,  we compute for  tangent vectors $h_i\in T_xX$,  $y=f(x)$, and $\phi\in \bm{U}(\phi_0)$ satisfying  $t(\phi)=y$,
 \begin{equation*}
 \begin{split}
(f^\ast\omega)_x(h_1,\ldots ,h_k)
&=\omega_y(Tf(x)h_1,\ldots ,Tf(x)h_k)\\
&=\omega_y(T\phi(h_1),\ldots ,T\phi(h_k))\\
&=\omega_x(h_1,\ldots h_k).
\end{split}
 \end{equation*}
Hence,  on $U(x_0)$, in view of $f^\ast \circ d=d\circ f^\ast$, 
 $$
 f^\ast d\omega = d f^\ast\omega = d\omega, 
 $$
so that the sc-differential form $d\omega$ is indeed compatible with morphisms.  From the properties of $d$ in the M-polyfold context we deduce that $d\circ d=0$.
  \begin{definition}\index{D- Sc de Rham complex}
 Given an ep-groupoid $X$,  we call $(\Omega^\ast_{ep,\infty}(X),d)$ the {\bf sc-de Rham complex} for the ep-groupoid $X$.
 We denote an element in $(\Omega^\ast_{ep,\infty}(X),d)$ again by $\omega$ though,  as an element in a direct limit,  it is rather an equivalence class.
 \qed
 \end{definition}

\section{Strong Bundles over Ep-Groupoids}\label{SST}

In this section we introduce the crucial notion of a strong bundle  over an ep-groupoid. Starting from a strong bundle $P: W\to X$ over the object M-polyfold 
$X$ of the ep-groupoid $(X, \bm{X} )$ we shall introduce the structure map 
$\mu: 
\bm{X}{{_s}\times_{P}}W\equiv \bm{W}\to W$.  Its properties allow to lift 
 the morphisms $\phi\in \bm{X}$ to linear isomorphisms 
 $\mu (\phi,\cdot ): W_{s(\phi)}\to W_{t(\phi)}$ between the corresponding fibers defining the morphisms in $\bm{W}$ between points of $W$ and ultimately leading to the sc-smooth functor $(P,\pi_1): (W,\bm{W})\to (X, \bm{X})$ between two ep-groupoids.
The  strong bundles over M-polyfolds have been introduced in Section \ref{section2.5_sb}.

We consider an ep-groupoid $X=(X, \bm{X})$ and a strong M-polyfold
bundle
$$
P : W\to X,
$$
where $P: W\to X$ is a surjective sc-smooth map between the M-polyfold $W$ onto the M-polyfold $X$ which is the object space of the ep-groupoid $(X, \bm{X})$. The M-polyfold $W$ possesses additional structures.  In particular, the fibers
$$
P^{-1}(x)=W_x
$$
over $x\in X$ carry the structure  of  Banach spaces. If $x\in X$ is a smooth point, the Banach space $W_x$ has a sc-structure.
 Since the source map $s : {\bf
X}\to X$ is by definition a local $\ssc$-diffeomorphism, the fibered product
$$
\bm{X}{{_s}\times_{P}}W=\{(\psi, w)\in \bm{X}\times W\vert \, s(\psi)=P(w)\}
$$
is an M-polyfold, and being viewed as  the pull-back of a strong M-polyfold bundle, is also a strong M-polyfold bundle.
Hence we have the strong bundle
$$
\pi_1:\bm{X}{{_s}\times_{P}}W\rightarrow \bm{X}.
$$
\begin{definition}[{\bf Strong bundle over an ep-groupoid}] \label{strong_bundle_ep}\index{D- Strong bundle over an ep-groupoid}
A strong bundle over the ep-groupoid $X$ is a pair $(P,\mu)$ \index{$(P,\mu)$} in which  $P: W\rightarrow X$ is a strong bundle over the object M-polyfold $X$
and $\mu: \bm{X}{{_{s}}\times_{P}}W\rightarrow W$ a strong bundle map covering the target map $t:\bm{X}\rightarrow X$, such  that the diagram  
\begin{equation*}
\begin{CD}
 \bm{X}{{_{s}}\times_{P}}W@>\mu>>W\\
@V\pi_1VV     @VVP V \\
\bm{X}@>t>> X,\\
\end{CD}
\end{equation*}
is commutative, so that 
$$
P\circ \mu (\psi, w)=t(\psi).
$$
Moreover, $\mu$ satisfies the  following additional 
properties.
\begin{Myitemize}
\item[(1)]\ $\mu$ is a surjective local $\ssc$-diffeomorphism
and linear on the fibers $W_x$.
\item[(2)]\ $\mu(1_x, w)=w$ for  all $x\in X$ and $w\in W_x$.
\item[(3)]\  $\mu (\psi\circ \gamma,  w)=\mu (\psi, \mu (\gamma, w))$
 for  all $\psi,\gamma \in \bm{X}$ and $w\in W$ satisfying  
 $$\text{$s(\gamma)=P(w)$\quad  and \quad $t(\gamma)=s(\psi)=P(\mu (\gamma,w))$.}$$
\end{Myitemize}
\qed
\end{definition}
Given a strong bundle $(P,\mu)$ over the ep-groupoid $X$ it follows, in particular, that
$$
\mu (\psi, \cdot ) : W_x\to W_y
$$
is a linear isomorphism if $\psi : x\to y$ is a morphism in $\bm{X}$.   In fact, the data of a strong bundle over an ep-groupoid  provide
via $\mu$ a prescription how to lift morphisms to linear isomorphisms in a compatible way with compositions.

We shall discuss the above definition now in detail. We shall  show that there is an  equivalent functorial picture in the background
which endows $W$ with morphisms resulting in an ep-groupoid  $(W,\bm{W})$ having  additional structures. In fact we shall obtain a sc-smooth functor 
$(P,\bm{P}): (W,\bm{W})\rightarrow (X,\bm{X})$, where $P$ and $\bm{P}$ are strong bundles.

We abbreviate $\bm{W}=\bm{X}{{_s}\times_{P}}W$. Now we can define the ep-groupoid
$$
W=(W, \bm{W})
$$
in which the object set $W$ is the M-polyfold $W$ we started with, and the morphism set $\bm{W}$ is the above fibered  product $\bm{X}{{_{s}}\times_{P}}W$.
Note, however, that we have more structure. The object set carries the structure of a strong bundle $W\rightarrow X$ and the morphism set
also carries the structure of a strong bundle,  namely  $\bm{W}\rightarrow \bm{X}$. For the following discussion we label the source and target maps
by the spaces they act in.  The source and target maps $s_W, t_W : \bm{W}\to W$ are defined as follows,
\begin{align*}
s_W(\psi, w)&=w\\
t_W(\psi, w)&=\mu (\psi, w)
\end{align*}
for $(\psi, w) \in \bm{X}{{_{s}}\times_{P}}W.$ These maps $s_W$ and $t_W$ are fiberwise linear   local $\ssc_\triangleleft$-diffeo\-morph\-isms covering the source and target maps $s_X,t_X\colon \bm{X}\to X$. Indeed,
\begin{align*}
P\circ s_W(\psi, w)&=P(w)=s_X(\psi)\\
P\circ t_W(\psi, w)&=P\circ \mu(\psi, w)=t_X(\psi).
\end{align*}
The action of a morphism $(\psi, w)\in  \bm{X}{{_{s}}\times_{P}}W$ is illustrated 
by the commutative diagram 
\begin{equation*}
\begin{CD}
s_{W}(\psi, w)@>(\psi, w)>>t_{W}(\psi, w)\\
@VPVV     @VVP V \\
s_{X}(\psi)@>\psi>> t_{X}(\psi).
\end{CD}
\end{equation*}
The orbit of a point $w\in W$ covers the orbit of $P(w)\in X$ and hence $P$ induces 
a surjective continuous map $\abs{P}: \abs{W}\to \abs{X}$ between the orbit spaces.

The identity morphism at $w\in W_x$ is the pair $(1_x, w)$ in $\bm{X}$ and if $\psi : x\to y$ is a morphism in $\bm{X}$, the inverse of the morphism  $(\psi, w)\in \bm{W}$ is the pair $i (\psi, w)=(\psi^{-1}, \mu (\psi, w))$. The multiplication map in $W$ is defined by
$$
(\psi, w')\circ (\gamma, w):=(\psi\circ \gamma, w)
$$
whenever  $w'=\mu (\gamma, w)\in W$.

The two sc-smooth projection maps $P : W\to X$ and $\bm{P}:=\pi_1 : \bm{W}\to \bm{X}$ together define an  sc-smooth functor denoted  by
$$
P=(P, \bm{P}) : W\to X
$$
between the two ep-groupoids $W=(W, \bm{W})$ and $X=(X, \bm{X})$.

As demonstrated above we can use the strong bundle $(P: W\to X, \mu)$ in order to construct additional data which allow to  view a strong bundle over the ep-groupoid $X$ alternatively as
the functor just constructed, which  in expanded notation is given by $(P,\bm{P}):(W,\bm{W})\rightarrow (X,\bm{X})$ with the projection  $\bm{P}:=\pi_1$.
We shall usually just call  $P:W\rightarrow X$ a strong bundle over the ep-groupoid $X$, whenever  the precise knowledge of the `structural'
strong bundle map $\mu$ is not needed.

For most purposes a strong bundle $(P,\mu)$ over the ep-groupoid $X$ should be viewed as a strong bundle over the object M-polyfold $X$, together
with a {\bf lift}\index{Lift of morphism} of the morphisms $\psi:x\rightarrow y$ to the  linear isomorphisms $\mu(\psi,\cdot ):W_x\rightarrow W_y$ satisfying some natural requirements.
For example,  $1_x$ is  lifted to the identity, and the lift associated to a composition is the composition of the lifts.
Moreover, the differential geometric side of things  requires that the data produce sc-smooth strong bundle maps.

Next we introduce the notion of a strong bundle map in the context of ep-groupoids. Again there are two equivalent viewpoints.
If we start with the original definition of strong bundles over ep-groupoids, say $(P,\mu)$ and $(P',\mu')$ then  a strong bundle map between
these two objects can be defined as follows.
\begin{definition}[{\bf Strong bundle map}]\label{def_strong_bundle_map}\index{D- Strong bundle functor}
Let $(P:W\rightarrow X,\mu)$ and $(P':W'\rightarrow  X',\mu')$ be two strong bundles over ep-groupoids. A {\bf strong bundle map}\index{Strong bundle map}
or {\bf strong bundle functor} between strong bundles over ep-groupoids consists of a strong bundle map $\Phi:W\rightarrow W'$ covering the sc-smooth functor $\varphi\colon X\rightarrow X'$ between ep-groupoids having the  property that
$$
\Phi(\mu(\gamma,w))=\mu'({ \varphi}(\gamma), \Phi(w)),
$$
for all $\gamma\in \bm{X}$ and $w\in W$ satisfying $s(\gamma)=P(w)$.
\qed
\end{definition}
While $\mu$ maps $w\in W_x$ to 
$\mu(\gamma,w)\in W_{t(\gamma)}$ for $\gamma\in \bm{X}$ satisfying  $x=s(\gamma)=P(w)$, 
the strong bundle map $\Phi$ maps $w\in W_x$ to $\Phi(w)\in W_{\varphi(x)}$ and  has to cover the sc-smooth map
$\varphi:X\rightarrow X'$ between the object M-polyfolds. Since $\varphi$ is a functor,  it follows
that the morphism $\gamma\in \bm{X}:s(\gamma)\rightarrow t(\gamma)$ is mapped to the morphism 
${\bm  \varphi}(\gamma): \varphi(s(\gamma))\rightarrow \varphi(t(\gamma))$ in $\bm{X}'$.
Hence the above required  relationship between  $\mu$ and $\mu'$ is well-defined.

There is a more functorial definition if we view $(P,\mu)$ as the functor $(P,\bm{P}):(W,\bm{W})\rightarrow (X,\bm{X})$ and similarly for $(P',\mu')$.
Namely,  given $\Phi$ as in the previous definition we can define the map $\bm{\Phi}:\bm{W}\rightarrow \bm{W}'$ by  
$$
\bm{\Phi}(\gamma,w)=({ \varphi}(\gamma),\Phi(w)) \ \ \text{for}\ \ (\gamma, w)\in \bm{W}=\bm{X}{{_s}\times_{P}}W.
$$
Then $\bm{\Phi}:\bm{W}\rightarrow \bm{W}'$ is a strong bundle map covering the map $\bm{\varphi}:\bm{X}\rightarrow \bm{X}'$ which is the  morphism part of  the functor 
$\varphi\colon X\to X'$ between the ep-groupoids and satisfying 
$$
\mu'\circ\bm{\Phi}(\gamma,w)=\mu'({\varphi}(\gamma),\Phi(w))=\Phi(\mu(\gamma,w))
$$
for $\gamma\in \bm{X}$.
Moreover, 
$$
s_{W'}(\bm{\Phi}(\gamma, w))=\Phi (w)=\Phi (s_{\bm{W}}(\gamma, w)).
$$
Recalling  that in the functorial picture $\mu$ and $\mu'$ are the target maps,
the two displayed equations show that the following diagrams commute 
 \begin{equation*}
  \begin{CD}
 \bm{W}@>\bm{\Phi}>>\bm{W}'\\
 @VsVV @VVsV\\
 W@>\Phi>> W'\\
 \end{CD}
 \qquad \qquad \qquad 
 \begin{CD}
 \bm{W}@>\bm{\Phi}>>\bm{W}'\\
 @VtVV @VVtV\\
 W@>\Phi>> W'.\\
 \end{CD}
 \end{equation*}
 Summarizing, we have ``lifted'' a strong bundle map $\Phi:W\rightarrow W'$ satisfying the requirements of  
Definition \ref{def_strong_bundle_map}  to a functor
 $(W,\bm{W})\rightarrow (W',\bm{W}')$.  The reverse is also true. Namely,  given such a functor covering the functor $\varphi\colon X\rightarrow X'$ between ep-groupoids, which is linear on the fibers, we obtain, using the strong bundle structures,  a strong bundle map $\Phi: W\rightarrow W'$ which has the properties listed in 
 Definition \ref{def_strong_bundle_map}. From the strong bundle map  we can recover the functor completely. Hence we obtain the following equivalent definition.

\begin{definition}\label{LinStBunM-def}\index{D- Strong bundle functor}
 A {\bf strong bundle functor} $\Phi : P\to P'$ between the two strong bundles $P : W\to X$ and $P' : W'\to X'$ over the ep-groupoids  $X$ and $X'$ consists of a functor $\Phi : W=(W, \bm{W})\to W'=(W', \bm{W}')$ between ep-groupoids which is linear on the fibers and which covers a sc-smooth functor 
 $(\varphi, {\bm  \varphi}) : 
 (X, \bm{X})\to (X', \bm{X}')$ between the bases. Moreover, the functor $\Phi$ induces strong bundle maps $\Phi : W\to W'$ and $\Phi : \bm{W}\to \bm{W}'$ between the object sets and morphism sets such that the diagram
\begin{equation*}
\begin{CD}
W@>\Phi>>W'\\
@VPVV     @VVP' V \\
X@>\varphi>> X' \\
\end{CD}
\end{equation*}
commutes.
\mbox{}\\[1ex]
\qed
\end{definition}

 In the following we shall either refer to a {\bf strong bundle map} between strong bundles over ep-groupoids if we appeal to the first definition,
or we shall refer to a {\bf strong bundle functor} if we use the second definition. In any case these are identical  notions.
For most purposes, however, the first view point is the most convenient one.
Namely,  if $W\rightarrow X$ and $W'\rightarrow X'$ are strong bundles over ep-groupoids
with respective structure maps $\mu$ and $\mu'$, then a strong bundle functor $\Phi$ is given by the following data. 
A sc-smooth functor $(\varphi,\bm{\varphi}):(X,\bm{X})\rightarrow (X',{\bf X'})$ between ep-groupoids 
and a strong bundle map $\Phi: W\rightarrow W'$
covering  $\varphi:X\rightarrow X'$, which is compatible with the structure maps in the sense that 
\begin{equation}\label{equation_section_functor}
\mu'(\varphi(\gamma),\Phi(w)) =\Phi(\mu(\gamma,w))
\end{equation}
for all $w\in W$ and $\gamma\in \bm{X}$ satisfying $s(\gamma)=P(w)$.

Having a strong bundle $(P:W\rightarrow X,\mu)$ over the ep-groupoid $X$ we can introduce  the notion of a section 
functor.

\begin{definition}\label{SECTION-FUNCTORS-X}
Let $(P: W\to X,\mu)$ be a strong bundle over the ep-groupoid $X$. 
\begin{itemize}
\item[(1)]\ A sc-smooth  {\bf section functor}\index{D- Section functor} $f$
is a  sc-smooth section $f: X\to W$ of the strong bundle $P:W\rightarrow X$ over the object M-polyfold satisfying 
 the  functorial condition
$$
f(y) =\mu(\phi,f(x))
$$
for  morphisms $\phi: x\rightarrow y$. 
We denote by $\Gamma(P,\mu)$\index{$\Gamma(P,\mu)$} the  vector space of sc-smooth section functors of $(P,\mu)$.
\item[(2)]\ A {\bf $\ssc^+$-section functor}\index{D- Sc$^+$-section functor} $f$ is an element in $\Gamma(P,\mu)$ satisfying  $f(x)\in W_{0,1}$ for all $x\in X$
and $f$,  viewed as a section of $W[1]\rightarrow X$,  is sc-smooth. By $\Gamma^+(P,\mu)\subset \Gamma (P, \mu)$\index{$\Gamma^+(P,\mu)$} we denote the vector subspace of $\ssc^+$-section functors.
\item[(3)]\  A {\bf sc-Fredholm section functor}\index{D- Sc-Fredholm section functor}  $f$ is an element in $\Gamma(P,\mu)$ such that, viewed as a section of the strong bundle
$W\rightarrow X$ over the object M-polyfold,  $f$ is a sc-Fredholm section. By $\text{Fred}(P,\mu)$\index{$\text{Fred}(P,\mu)$} we denote the collection  of sc-Fredholm section functors.
\end{itemize}
\qed
\end{definition}

These notions will be studied in more detail later on.

\section{Proper Covering Functors}\label{cov_g}
We need later on the notion of a proper  covering functor, which is inspired by the definition in the Lie groupoid case in  \cite{Mac}.
\begin{definition}\label{d-proper_covering}\index{D- Proper covering functor}
A sc-smooth functor $F:Y\rightarrow X$ between ep-groupoids is called a {\bf proper covering functor}
if the following holds true.
\begin{itemize}
\item[(1)]\ The sc-smooth map $F:Y\rightarrow X$ defined between the objects
is a surjective local sc-diffeomorphism.
\item[(2)]\ For every  object $x\in X$, the  preimage $F^{-1}(x)\subset Y$ is finite. There exists an  open neighborhood $U(x)\subset X$
and  there exist mutually disjoint open neighborhoods $U(y)\subset Y$ around  every point $y\in F^{-1}(x)$ such that the map
$F: U(y)\rightarrow U(x)$ is a sc-diffeomorphism and
$$
F^{-1}(U(x))=\bigcup_{y\in F^{-1}(x)} U(y).
$$
\item[(3)]\ The map 
$$
\bm{Y}\rightarrow \bm{X}{_{s}\times_F} Y,\quad \phi\mapsto (F(\phi),s(\phi))
$$
is a sc-diffeomorphism.
\end{itemize}
\qed
\end{definition}
The map
$$
\bm{Y}\rightarrow \bm{X}{_{s}\times_F} Y,\quad \phi\mapsto (F(\phi),s(\phi)), 
$$
being a  sc-diffeomorphism, preserves 
the degeneracy index, so that 
$$
d_{\bm{Y}}(\phi) = d_{  \bm{X}{_{s}\times_F} Y}(F(\phi),s(\phi))
$$
for the morphisms $\phi\in \bm{Y}$.
By construction, the fiber product   $ \bm{X}{_{s}\times_F} Y$ is a sub-M-polyfold of the product 
$ \bm{X}\times Y$. The latter has the degeneracy index 
$$
d_{ \bm{X}\times Y} (\phi,y)=d_{\bm{X}}(\phi) +d_Y(y).
$$
Since $F$ and $s$ are local sc-diffeomorphisms we conclude that
$$
d_{\bm{X}\times Y}(F(\phi),s(\phi))=d_{\bm{X}}(F(\phi))+d_Y(s(\phi))= d_{\bm{Y}}(\phi)+d_{\bm{Y}}(\phi)=2\cdot d_{\bm{Y}}(\phi).
$$
Combining this with $d_{\bm{Y}}(\phi)=d_{  \bm{X}{_{s}\times_F} Y}(F(\phi),s(\phi))$
we conclude that
$$
d_{  \bm{X}{_{s}\times_F} Y}(F(\phi),s(\phi))=\frac{1}{2}\cdot d_{\bm{X}\times Y}(F(\phi),s(\phi)),
$$
which shows that the fibered product lies in a particular position in the full product.

Next we study the composition of sc-smooth proper covering functors.

\begin{proposition}\index{P- Composition of coverings}
Assume that $X$, $Y$, and $Z$ are ep-groupoids. If $F:Y\rightarrow X$ and $H:Z\rightarrow Y$ are proper covering functors,
then the composition $F\circ H:Z\rightarrow X$ is a proper covering functor.
\end{proposition}
\begin{proof}
That $F\circ H$ satisfies the properties (1) and (2) follows easily. To verify (3) we have to show
that
$$
L_{F\circ H}:{\bf Z}\rightarrow \bm{X}{_{s}\times_{F\circ H}} Z, 
\quad \sigma\mapsto (F(H(\sigma)),s(\sigma))
$$
is a sc-diffeomorphism.
To this aim we introduce the sc-smooth map  $\wtilde{K}$ by
$$
\wtilde{K}:\bm{Y}{_{s}\times_H}Z\rightarrow  \bm{X}{_{s}\times_{F\circ H}}Z, 
\quad(\psi,z)\mapsto (F(\psi),z),
$$
and consider the sc-diffeomorphisms
\begin{align*}
L_H&:{\bf Z}\rightarrow \bm{Y}{_{s}\times_H} Z, 
\quad\sigma\mapsto  (H(\sigma),s(\sigma))&\\
L_F&:\bm{Y}\rightarrow \bm{X}{_{s}\times_F} Y, 
\quad\psi\mapsto  (F(\psi),s(\psi)).
\end{align*}
Then 
$$
L_{F\circ H} =\wtilde{K}\circ L_H.
$$
It suffices  to show that $\wtilde{K}$ is a sc-diffeomorphism, since we already know that the same is true fo $L_H$.
Considering  the  sc-smooth map
$$
\bm{X}{_{s}\times_{F\circ H}}Z\rightarrow \bm{Y} \times Z, 
\quad (\phi,z)\mapsto  (L_F^{-1}(\phi,H(z)),z),
$$
we note  that $L_F^{-1}(\phi,H(z))=\psi$ implies $F(\psi)=\phi$ and $s(\psi)=H(z)$. 
Hence 
$$
(L_F^{-1}(\phi,H(z)),z)\in \bm{Y}{_{s}\times_H} Z
$$
 and 
therefore we can define the sc-smooth map
$$
\wtilde{L}:\bm{X}{_{s}\times_{F\circ H}}Z\rightarrow \bm{Y}{_{s}\times_H} Z \quad \text{by}\quad \wtilde{L}(\phi,z)=(L_F^{-1}(\phi,H(z)),z).
$$
Then  $\wtilde{K}\circ \wtilde{L}(\phi,z)=\wtilde{K}(\psi,z)$, where $F(\psi)=\phi$ and $s(\psi)=H(z)$.
Hence $\wtilde{K}(\psi,z)=(\phi,z)$ implying $\wtilde{K}\circ\wtilde{L}=\id$. Moreover
$\wtilde{L}\circ \wtilde{K}(\psi,z)=\wtilde{L}(F(\psi),z)= (\psi,z)$ so that $\wtilde{L}\circ\wtilde{K}=\id$.
We have shown that $\wtilde{K}$ is a sc-diffeomorphism implying that the map $L_{F\circ H}$ is a sc-diffeomorphism.
\qed \end{proof}

The following theorem can be viewed  as a generalization
of the result about the natural action of the isotropy group.

\begin{theorem}[{\bf Proper Covering {I}}]\label{STRUC_1}\index{T- Proper covering {I}}
Let $F:Y\rightarrow X$ be a proper covering functor and let $x\in X$. Then for  a connected and sufficiently small open  neighborhood $U(x)$ in $X$,
admitting the natural $G_x$-action (given by the data $(\Phi,\Gamma)$), the following holds.
\begin{itemize}
\item[{\em (0)}]\  The subset  $F^{-1}(U(x))$ in $Y$  is a disjoint union of connected open neighborhoods $U(y)$ for all $y\in F^{-1}(x)$  such  that every map $F:U(y)\rightarrow U(x)$ is a sc-diffeomorphism. 
\item[{\em (1)}]\  The subset  $\Gamma(G_x,U(x))=\{\phi\in \bm{X}\, \vert \, s(\phi),t(\phi)\in U(x)\}=:\bm{U} (x)$ is the union of connected components $\bm{U}(g)$ with $g\in G_x$.
Moreover, the source and target maps  $s,t:\bm{U}(g)\rightarrow U(x)$ are sc-diffeomorphisms.
\item[{\em (2)}]\  For every pair $y,y'\in F^{-1}(x)$,  the set 
$$
\bm{U} (y, y')=\{\psi\in \bm{Y}\, \vert \, \text{$s(\psi)\in U(y)$ and $t(\psi)\in U(y')$} \}$$ 
is the disjoint union of connected open neighborhoods
$\bm{U}(\psi_0)\subset \bm{Y}$ centered at morphismsm  $\psi_0: y\to y'$ for which the following properties (3)-(5) hold. 
\item[{\em (3)}]\ The maps $s:\bm{U}(\psi_0)\rightarrow U(y)$ and $t:\bm{U}(\psi_0)\rightarrow U(y')$ are sc-diffeomorphisms.
\item[{\em (4)}]\ The map $\bm{F}: \bm{U}(\psi_0)\rightarrow \bm{U}(\bm{F}(\psi_0))$ is a sc-diffeomorphism. 
\item[{\em (5)}]\ The following diagram is a commutative diagram of sc-diffeomorphisms.
$$
\begin{CD}
\bm{U}(\psi_0)@>\bm{F}>> \bm{U}(\bm{F}(\psi_0))\\
@V s VV @V t VV\\
U(s(\psi_0)) @> F>> U(t(\bm{F}(\psi_0)).
\end{CD}
$$
\end{itemize}
\qed
\end{theorem}
\begin{remark}\index{R- On proper coverings}
(a) As $F$ is a functor we conclude for a morphism $\phi: y\to y'$ in $\bm{Y}$ connecting two objects $y, y'\in F^{-1}(x)$, that the morphism $\bm{F}(\phi): F(y)\to F(y')$ belongs to the isotropy group $G_x$.\par

\noindent (b) In Theorem \ref{STRUC_1}
 everything is determined by the  choice of a small connected open neighborhood $U(x)$
of an object $x\in X$, which  admits the natural $G_x$-action. As discussed in Section \ref{sec_1_1}, Theorem \ref{x-local-x},
the construction of the data $(\Phi,\Gamma)$, which determines the $G_x$-action,   is  natural. The main point in
proving the theorem is the interplay between the map $F$ restricted to suitable subsets,
and its compatibility with the source and target maps due to the functoriality.\qed
\end{remark}

\begin{proof}[Theorem \ref{STRUC_1}]
The set $A=F^{-1}(x)$ is finite,  and for a 
connected open neighborhood $U'(x)$,  the preimage $F^{-1}(U'(x))$ is the disjoint union
of connected open neighborhoods $U'(y)$, so that every map $F: U'(y)\rightarrow U'(x)$ is a sc-diffeomorphism. This is an immediate consequence of the definition of a sc-smooth proper
covering functor. We may assume that the set $U'(x)$ admits  the natural action of $G_x$, otherwise we replace $U'(x)$ by a smaller $G_x$-invariant and for the sets $U'(y)$ we take its preimages  $F^{-1}(U'(x))$.

We find around every point $y\in A$ 
a connected open neighborhood $V(y)$, contained in $U'(y)$ and admitting admitting the natural action of $G_{y}$, and,  for every morphism $\psi: y\to y'$,  a connected open neighborhood $\bm{V}(\psi)$ such that the source and target maps $s: \bm{V}(\psi)\to V(y)$ and $t: \bm{V}(\psi)\to V(y')$ are sc-diffeomorphims. We define the sc-diffeomorphism
$\wh{\psi}$ by 
\begin{eqnarray}\label{disp834}
\wh{\psi}=t\circ \bigl(s\vert \bm{V}(\psi)\bigr)^{-1}: V(y)\to V(y').
\end{eqnarray}
Denoting by $\bm{A}$ the set of all morphisms between objects in $A$,  the set 
$$
\bm{V}=\left\{\phi\in \bm{X}\, \vert \, s(\phi)\in \bigcup_{y\in F^{-1}(x)}V(y), 
t(\phi)\in \bigcup_{y'\in F^{-1}(x)}V(y')\right\}
$$
is a disjoint union of connected components $\bm{V}(\psi)$, 
$$
\bm{V}=\bigcup_{\psi \in \bm{A}}\bm{V} (\psi).
$$
We take a smaller connected open neighborhood $U(x)\subset U'(x)$ which is $G_x$-invariant and for which 
$(F\vert V(y))^{-1}(U(x))\subset V(y)$.  Since the sets $U'(y)$ are $G_y$-invariant, the same is true for the sets 
$U(y):=(F\vert V(y))^{-1}(U(x))$. We claim that the map 
$\wh{\psi}: V(y)\to V(y')$ defined in (\ref{disp834})
maps bijectively 
$U(y)$ onto $U(y')$, and hence it is a sc-diffeomorphism from $U(y)$ onto $U(y')$. To see this, we take $z\in U(y)$ and let $z'=\wh{\psi}(z)$. Then $z'\in V(y')$ and there is a morphism $\phi: z\to z'$. Hence $\bm{F}(\phi)$ is a morphism between $F(z)$ and $F(z')$.  Since $F(z)\in U(x)$ and $F(z')\in U'(x)$, and $U'(x)$ is $G_x$ invariant, there exists $g\in G_y$ for which $\bm{F}(\phi)=\Gamma (g, F(z))$. 
But $U(x)$ is a $G_y$-invariant subset of $U'(x)$ and hence 
$F(z')=s(\bm{F}(\phi))\in U(x)$. So, $z'\in U(y')$ and $\wh{\psi}(U(y))\subset U(y')$. Similar arguments shows that $\wh{\psi}^{-1}(U(y'))\subset U(y)$. Hence the map $\wh{\psi}$ restricted to $U(y)$ is a sc-diffeomorphim onto $U(y')$. 
Now it is clear that if $\bm{U} (\psi)=(s\vert \bm{V}(\psi))^{-1}(U(y))$, then the maps $s: \bm{U} (\psi)\to U(y)$ and 
$t: \bm{U} (\psi)\to U(y')$ are sc-diffeomorphisms. Moreover, 
if $\bm{U} =\{\phi\in \bm{X}\, \vert \, s(\phi)\in \bigcup_{y\in F^{-1}(x)}U(y), 
t(\phi)\in \bigcup_{y'\in F^{-1}(x)}U(y')\}$,
then 
$$\bm{U}=\bigcup_{\psi \in \bm{A}}\bm{U} (\psi).$$

At this point  the properties (0)-(3)  are satisfied and it remains to verify the property (4). We take a morphism $\psi_0: y\to y'$ between points $y, y'\in A$. Then 
$\bm{F}(\psi_0)=g\in G_x$.  Since $\bm{U} (\psi_0)$ is connected and the set $\bm{U} (g)$ is a connected component of the set $\bm{U} (x)$, it follows that $\bm{F}(\bm{U} (\psi_0))\subset \bm{U} (g).$ 
Using that the maps 
$F: U(y)\to U(x)$, $s: \bm{U} (g)\to U(x)$, and 
$s: \bm{U} (\psi_0)\to U(y)$ are  sc-diffeomorpism, we see that any morphism $\phi \in \bm{U} (g)$ is of the form $\phi=(s\vert \bm{U} (g))^{-1}\circ F\circ s(\psi)$ for a unique morphism $\psi\in \bm{U} (\psi_0)$. Both  morphisms $\bm{F}(\psi)$ and $\phi=$ belong to $\bm{U} (g)$ and have the same source,  so 
that $\bm{F}(\psi)=\phi=(s\vert \bm{U} (g))^{-1}\circ F\circ s(\psi).$
Hence $\bm{F}$ is a sc-diffeomorphism  as a  composition  of sc-diffeomorphisms,
$$
\bm{U}(\psi_0)\xrightarrow{s} U(s(\psi_0))\xrightarrow{F} U(x)\xrightarrow{(s\vert \bm{U} (g))^{-1}} \bm{U}(g).
$$
We have verified the  property (4) and the  property (5) is now obvious. The proof of Theorem \ref{STRUC_1} 
is complete.
\qed \end{proof}

Let $F:Y\rightarrow X$ be a proper covering functor.  Fixing  $x\in X$ and a sufficiently small connected open neighborhood $U(x)\subset X$, and open neighborhoods
$U(y)\subset Y$ for $y\in F^{-1}(x)$, such  that every map $F:U(y)\rightarrow U(x)$ is a sc-diffeomorphism,  we obtain,  on the object level,  the map
$$
F:\sqcup_{y\in F^{-1}(x)} U(y)\rightarrow U(x).
$$
Assuming  that $U(x)$ admits the natural $G_x$-action,  
we take the associated full subcategories and  obtain a local model  for the proper covering functor near $x$ and $F^{-1}(x)$.
The morphism structure incorporates the natural actions of the isotropy groups $G_y$ and their  relationships.
Next we shall study the structure of this local model.

Using the previous theorem, we shall  explain the local structure of a proper covering functor in more detail.
We    denote by $\text{Diff}_{\mathrm{sc}}$\index{$\text{Diff}_{\mathrm{sc}}$}
the  category whose objects are M-polyfolds and whose morphisms are sc-smooth diffeomorphisms between M-polyfolds.
\begin{definition}[{\bf Covering functor, geometric lift}]\label{DEFR845}
Let $D$ and $C$ be two categories having both finitely many objects and morphisms. In addition, $C$ has only one object, usually denoted by $x$, so that $C=\{x\}$, and all morphisms of $C$  and $D$ are isomorphisms. 
A {\bf covering functor}\index{D- Covering functor} 
$$f: D\rightarrow C$$
 is a functor such  that the map
$$
\bm{D}\rightarrow \bm{C}{_{s}\times_f} D,\quad \psi\rightarrow (f(\psi),s(\psi))
$$
is a bijection.  

A {\bf geometric lift}\index{D- Geometric lift} of a covering functor 
$f: D\to C$ consists of two functors $L_D:D\rightarrow \text{Diff}_{\mathrm{sc}}$ and $L_C:C\rightarrow \text{Diff}_{\mathrm{sc}}$ and a natural transformation
$$
\tau: L_D\rightarrow L_C\circ f.
$$
of the two functors $L_D, L_C\circ f: D\to \text{Diff}_{\mathrm{sc}}$.
\qed
\end{definition}
 
 We recall that the map $\tau: D\to \text{Diff}_{sc}$ associates with every object $y\in D$ a morphism 
 $$\tau (y): L_D(y)\to L_C\circ f(y)$$
 in $\text{Diff}_{sc}$ satisfying 
 $$\tau (y')\circ L_D(h)=(L_C\circ f)(h)\circ \tau (y)$$
 for every morphism $h: y\to y'$ in $\bm{D}$.
 
 If $f: D\to C$ is a covering functor and
 $\tau:L_D\rightarrow L_C\circ f$ a geometric lift of $f$, the unique object $x\in C$ is mapped to the M-polyfold $L_C(x)$, denoted by 
 $$
U(x):=L_C(x).
$$
For the finitely many objects $y\in D=f^{-1}(x)$ 
we obtain the M-polyfolds  $L_D(y)$ denoted by 
$$
U(y):=L_D(y),\quad y \in D=f^{-1}(x).
$$
Associated with  every objects $y \in D=f^{-1}(x)$ is the morphism 
$\tau (y): L_D(y)\to L_C(f(y))$  which, in view of $f(y)=x$, is the  sc-diffeomorphism  
$$
\tau(y):U(y)\rightarrow U(x)
$$
between the two M-polyfolds.

If $\psi:y\rightarrow y'$ is a morphism in $D$, then the  morphism $L_D(\psi): L_D(y)\to L_D(y')$ in $\text{Diff}_{sc}$ is the sc-diffeomorphism 
$$
\varphi_{yy'}(\psi):=L_D(\psi): U(y)\rightarrow U(y').
$$
The morphism $g\in \bm{C}=G_x$, $g: x\to x$, is mapped  to the morphism $L_C(g): L_C(x)\to L_C(x)$ which is the sc-diffeomorphism 
$$L_C(g): U(x)\to U(x)$$
of the M-polyfold $U(x)$. 
For two morphisms $\gamma: y\to y'$ and $\psi: y'\to y''$ in $\bm{D}$ so that 
$\psi\circ\gamma: y\to y''$, we obtain by functoriality the formula 
$$
\varphi_{yy''}(\psi\circ\gamma)=\varphi_{y'y''}(\psi)\circ\varphi_{yy'}(\gamma)
$$
which is a composition of sc-diffeomorphisms.

In the special case $y=y'$, the isotropy group $G_y\subset \bm{D}$ consists of the finitely many morphisms $y\to y$, and the $G_y$-action on the M-polyfold $U(y)$ is defined  by
$$
\text{$h\ast z:=L_D(h)(z)$, \quad $z\in U(y)$ and $h\in G_y$.}
$$
Similarly,  we denote the action of the isotropy group $G_x$ on $U(x)$ by $g\ast z=L_C(g)(z)$ for $z\in U(x)$ and $g\in G_x$.

Using that $\tau$ is a natural transformation we obtain for $h\in G_y$ the formula
\begin{equation*}
\begin{split}
\tau (y)(h\ast z)&=(\tau (y)\circ L_D(h))(z)\\
&=(L_C(f(h))\circ \tau (y))(z)\\
&=f(h)\ast \tau (y)(z)
\end{split}
\end{equation*}
for all $z\in U(y)$ where $\tau (y)(z)\in U(x)$. 
Therefore, the sc-diffeomorphisms $\tau (y): U(y)\to U(x)$ are equivariant with respect to the map $f: G_y\to G_x$.

\begin{proposition}[{\bf Geometric Lift}]\index{P- Geometric lift}\label{geo_lift}
Associated with a covering functor 
$f:D\rightarrow C$ and a geometric lift $\tau:L_D\rightarrow L_C\circ f$ of $f$  there exists a natural construction of a proper covering functor $F:Y\rightarrow X$ between the  two ep-groupoids $Y$ and $X$.
\end{proposition}
\begin{proof}
Let $x$ be the  unique object of $C$ and $G_x=\bm{C}$ its  isotropy group. 
Using  the previously introduced notations  we take  the translation groupoid 
$$X=G_x\ltimes U(x),$$
 whose  object M-polyfold is $U(x)$ and whose morphism
M-polyfold is $G_x\times U(x)$. The source and the target maps $s, t: G_x\times U(x) \to U(x)$ are defined as  $s(g,z)=z$, and $t(g,z)=g\ast z$.

The  object set  $Y$ of the  ep-groupoid  $(Y, \bm{Y})$  is defined as the disjoint 
union
$$
Y=\bigcup_{y\in D} U(y).
$$
It possesses a 
natural M-polyfold structure.
The  morphism set $\bm{Y}$  is the set 
$$
\bm{Y}=\bigcup_{(y,y')\in D\times D} [\text{mor}(y,y')\times U(y)],
$$ 
which has also  a natural M-polyfold structure.
The source and target maps $s, t: \bm{Y}\to Y$ are defined by 
$$s(\psi,z)=z\quad \text{and}\quad 
t(\psi,z) =\varphi_{yy'}(\psi)(z)
$$
for $\psi \in \text{mor}(y, y')$ and $z\in U(y).$
The maps  $s$ and $t$ are surjective local sc-diffeo\-mor\-phisms.
The composition of  the morphisms $(\psi,z')$ and $(\gamma,z)$ in $\bm{Y}$ is well defined provided  $z'=\varphi_{yy'}(\gamma)(z)$,  and is given by
$$
(\psi,z')\circ(\gamma, z)=(\psi\circ\gamma,z).
$$
The identity elements have the form $(1_y,z)$, where $z\in U(y)$ so that the unit map $u:Y\rightarrow \bm{Y}$ is the map 
$$
u(z)=(1_y,z),
$$
where $z\in U(y)$. The inversion map $\iota:\bm{Y}\rightarrow \bm{Y}$ is given by
$$
\iota(\psi,z)=(\psi^{-1},\varphi_{yy'}(\psi)(z)).
$$
For the M-polyfold structure on $\bm{Y}$  the inversion map $\iota$ is sc-smooth because  $\varphi_{yy'}(\psi)$ is a sc-diffeomorphism.
One verifies that $Y=(Y, \bm{Y})$ is an ep-groupoid, i.e the composition is  sc-smooth and the properness property holds.

The functor $F: (Y, \bm{Y})\rightarrow (X, \bm{X})$ is defined on the objects by 
$F|U(y):=\tau (y)$ for $y\in D$ which is a sc-diffeomorphism $U(y)\to U(x)$. On the morphism set ${\bf }$ we define $\bm{F}: \bm{Y}\to \bm{X}$ as follows. If $(\psi, z)\in \bm{Y}$, where $\psi : y\to y'$ is a morphism satisfying $s(\psi)=y$, and $z\in U(y)$, we set 
$$
\bm{F}(\psi,z)= (f(\psi),\tau (y)(z))\in \bm{X}.
$$
Since $F: U(y)\to U(x)$ is a sc-diffeomorphism
and $F^{-1}(U(x))=Y$, the map 
$F: X\to Y$ is a surjective local sc-diffeomorphism between the object sets. It follows from the construction that 
$$F: (Y, \bm{Y})\to (X, \bm{X})$$
is a sc-smooth functor between ep-groupoids.
Finally, we claim that the map 
\begin{equation}\label{map_star}
\bm{Y}\rightarrow \bm{X}{_{s}\times_F} Y,\quad (\psi,z)\rightarrow ((f(\psi),F(z)),z)
\end{equation}
 is a sc-diffeomorphism. Indeed, from 
$$
((f(\psi),F(z)),z)=((f(\psi'),F(z')),z').
$$
it follows that $z=z'$ 
 implying 
$$(f(\psi), \tau (y)(z))=(f'(\psi'), \tau (y') (z))$$
and hence $y=y'=s(\psi)=s(\psi')$. 
By the definition of a covering functor $f: D\to C$, the map 
$\bm{D}\to \bm{C}{_{s}\times_f}\bm{D}$, $\psi \mapsto (f(\psi), s(\psi))$ is a bijection. Therefore, 
we conclude that $\psi=\psi'$, which shows the injectivity of map \eqref{map_star}. In order to verify the surjectivity we  take an element 
$((g,q),z)\in \bm{X}{_{s}\times_F} Y$.
Then $F(z)=q$ and there exists a point $y\in Y$ for which $z\in U(y)$. Therefore, we find a unique  $\psi\in \bm{D}$ satisfying  $f(\psi)=g$ and $s(\psi)=y$ and consequently, 
$
((f(\psi),F(z)),z)=((g,q),z),
$
which  shows the surjectiviy of  the map \eqref{map_star}. Being a local sc-diffeomorphism, the map \eqref{map_star}  is therefore a sc-diffeomorphism as claimed.

The proof that the functor $F:Y\rightarrow X$ between the ep-groupoids is a proper covering functor according to Definition \ref{d-proper_covering}, is complete. \qed \end{proof}

Conversely,  combining Theorem \ref{STRUC_1} with  Proposition \ref{geo_lift} we shall construct, for  a proper covering functor $F:Y\rightarrow X$ between ep-groupoids, 
a local model near the object $x\in X$ and near the preimage $F^{-1}(x)\subset Y$, which is the geometric lift of the induced  covering functor
$$
f: D\rightarrow C,
$$
where $D$ is the full subcategory associated with  $F^{-1}(x)$ and $C$ the full subcategory associated with $x$. Moreover,  $f=F\vert D$.

\begin{theorem}[{\bf Proper Covering {II}}] \label{STRUC_2} \index{T- Proper covering {II}}
Let $F:Y\rightarrow X$ be a proper covering functor. Fixing an object $x\in X$, we let $C$ be the full subcategory of $X$ having the single object $x$ and the isotropy group $G_x$ as the morphism set. Moreover, we let $D$ be
the full subcategory having the finite object set $F^{-1}(x)$. Then the restriction $f=F|D:D\rightarrow C$ is a covering functor. Given a sufficiently small open neighborhood
$U(x)$ of $x$ in $X$ as described in Theorem \ref{STRUC_1}, there exists exists a geometric lift $A:K\rightarrow L$ of $f$ and there exist  fully faithful open 
embedding functors
$\Phi$ and $\Psi$, where the image of the functor $\Phi:L\rightarrow X$ is  the full subcategory associated with  $U(x)$, and the image of the functor $\Psi:K\rightarrow Y$ is  the full subcategory associated
with $F^{-1}(U(x))$, such  that the following diagram is commutative, 
$$
\begin{CD}
K @>\Psi >> Y\\
@V A VV @V F VV\\
L @> \Phi >> X.
\end{CD}
$$
\end{theorem}
\begin{proof}
Fixing the object $x\in X$ we obtain the category $C$  consisting of the single object $x$ and the morphism set $G_x$, and the full subcategory $D=(D, \bm{D})$ of $Y$ associated with the finite set of objects  $D=F^{-1}(x)$.
We denote by $f:D\rightarrow C$ the restriction of the proper covering functor $F$ to $D$, so that $f=F\vert D$.  From the covering functor $F$ we deduce that  $f$  is a covering functor between finite categories.
Next we employ Theorem \ref{STRUC_1} and choose the open neighborhood  $U(x)\subset X$ equipped  with its natural action of $G_x$ given by $(\Phi,\Gamma)$, so small that
the properties (0)-(5) of Theorem \ref{STRUC_1} hold. In particular,  we have open neighborhoods $U(x)\subset X$ and $U(y)\subset Y$ for $y\in D$ which are M-polyfolds.
Next we define the geometric lift of $f:D\rightarrow C$. The functor $L_C:C\to \text{Diff}_{sc}$ is defined by
$$
L_C(x)=U(x)\quad \text{and}\quad  L_C(g)=\Phi(g)\in \text{Diff}_{\textrm{sc}}(U(x)), \quad g\in G_x.
$$
On the object set $D$ the functor  $L_D:D\to \text{Diff}_{sc}$ is defined as 
$$
L_D(y)=U(y),\quad y\in D=F^{-1}(x).
$$
Recalling the statement (5) of Theorem  \ref{STRUC_1}, we define the functor $L_D:D\to \text{Diff}_{sc}$ on morphisms as follows. 
If  $\psi:y\to y'\in \bm{U}(\psi_0)$ for the distinguished morphism $\psi_0: y\to y'$, we put 
$$L_D(\psi )=t\circ (s\vert \bm{U}(\psi_0))^{-1}\in 
\text{Diff}_{sc}(U(y),U(y')).
$$
The natural transformation $\tau:L_D\to L_C\circ f$ associates with the object  $y\in F^{-1}(x)$  the morphism 
$\tau (t):L_D(y)\to L_C\circ f(y)$ which is the sc-diffeomorphism 
$$\tau (y):=F\vert U(y): U(y)\to U(x).
$$
One verifies that this defines a geometric lift of  the covering functor $f:D\rightarrow C$. The associated natural construction of a covering functor is denoted by
$A:K\rightarrow L$, where $L$ is the translation groupoid $L=G_x\ltimes U(x)$, and the ep-groupoid $(K,\bm{K})$  is defined by 
$$
K=\sqcup_{y\in F^{-1}(x)} U(y)\quad \text{and}\quad \bm{K}=\sqcup_{(y,y')\in D^2} [\text{mor}(y,y')\times U(y)].
$$
The proper covering functor $A:K\rightarrow L$ is defined on  the objects by $A|U(y):=F|U(y)$ and on the morphisms  $(\psi, z)\in \text{mor}(y,y')\times U(y)$ by  
$$
A(\psi,z)= (f(z), F(z)).
$$
The functor $\Phi:L\rightarrow X$ is the inclusion on objects and maps the morphism $(g,z)\in G_x\times U(x)$
to the morphism $\Gamma(g,z)\in \bm{X}$. The functor $\Psi:K\rightarrow Y$ is the inclusion on objects,  and maps  the morphism $(\psi,z) \in \text{mor}(y,y')\times U(y)$ where $\psi \in \bm{U}(\psi_0)$,  to the morphism $(s\vert \bm{U}(\psi_0)))^{-1}(z)\in \bm{Y}$. 
This way we obtain the commutative diagram
$$
\begin{CD}
K @>\Psi >> Y\\
@V A VV @V F VV\\
L @> \Phi >> X.
\end{CD}
$$
The image of $\Psi$ is the full subcategory associated with the object set  $F^{-1}(U(x))$, the image of $\Phi$ is the full subcategory associated with  $U(x)$, 
and the covering functor $A$ is  the restriction of $F$. The construction of the covering functor $A: K\rightarrow L$ defines, up to the fully faithful
open embeddings $\Psi$ and $\Phi$,  a local model for the proper covering functor $F$ near $x$ and $F^{-1}(x)$. The proof of Theorem \ref{STRUC_2} is complete.
\qed \end{proof}

For further considerations concerning the local structure of a proper covering functor we refer to  Appendix \ref{strustru}.

For applications we need a version of the previous discussion involving strong bundles. We recall  from Definition \ref{strong_bundle_ep} that the structure map $\mu$ of a strong bundle over an ep-groupoid $P:W\rightarrow Y$
is a strong bundle map $\mu:\bm{Y}{_{s}\times_P} W\rightarrow W$ covering the target map $t:\bm{Y}\rightarrow Y$. 
The pairs $(\psi,w)\in \bm{W}=\bm{Y}{_{s}\times_P} W$ are  viewed as the morphisms 
$$
(\psi,w): w\rightarrow \mu(\psi,w).
$$
of the object space $W$. Then $\bm{W}$ becomes  a strong bundle over $\bm{Y}$ via the projection map  
$$
\bm{P}: \bm{W}\rightarrow \bm{Y},\quad (\psi,w)\rightarrow \psi.
$$
The source and  the target maps $s, t: \bm{W}\rightarrow W$ are defined as $s(\psi,w)=w$ and $t(\psi,w)=\mu(\psi,w)$.  The reader  finds more details in Section \ref{SST} about strong bundles.
\begin{definition}\label{proper_sb_covering}
Let $(W,\mu)$ and $(V,\tau)$ be strong bundles over the ep-groupoids $Y$ and $X$, respectively, denoted  by  $P:W\rightarrow Y$ and $Q:V\rightarrow X$.
A {\bf proper strong bundle covering functor} is a strong bundle map $A$  (in particular fiberwise linear), for which the diagram
$$
\begin{CD}
W @> A>>  V\\
@V P VV @ V Q VV\\
Y @> F >> X
\end{CD}
$$
is commutative and which has 
the following additional properties.
\begin{itemize}
\item[(1)]\ $A$ is surjective and a local strong bundle isomorphism covering a sc-smooth proper covering functor $F:Y\rightarrow X$.
\item[(2)]\  The preimage $F^{-1}(x)$ of every object $x\in X$  is finite,  and there exist open neighborhoods $U(x)\subset X$ and $U(y)\subset Y$ for every $y\in F^{-1}(x)$, where the sets $U(y)$ are mutually disjoint, such  that the map
$$
A:W\vert U(y)\rightarrow V|U(x)
$$
 is a strong bundle isomorphism for every $y\in F^{-1}(x)$, and 
 $$
 F^{-1}(U(x))=\bigcup_{y\in F^{-1}(x)} U(y).
 $$
 \item[(3)]\ The map 
 $$
\bm{W}\rightarrow \bm{V}{_{s}\times_A} W:(\psi,w)\rightarrow ((F(\psi),A(w)),w)
 $$
 is a strong bundle isomorphism covering the sc-diffeomorphism $\bm{Y}\rightarrow \bm{X}{_{s}\times_F}Y:\psi\rightarrow (F(\psi),s(\psi))$. Here the strong bundle projection
 $$
 \bm{V}{_{s}\times_A} W\rightarrow \bm{X}{_{s}\times_F} Y
 $$
  is given by $((\phi,v),w)\rightarrow (\phi,P(w))$.
 \end{itemize}
 \end{definition}
The strong bundle data are  some kind of a lift of the data on the level of ep-groupoids. All structural statements about proper covering functors
 between ep-groupoids lift to the strong bundle situation.  We give some further discussion in the following appendix.
 \section{Appendix}

\subsection{Local Structure of Proper Coverings}\label{strustru}
We consider a proper covering functor 
$F:Y\rightarrow X$ between ep-groupoids 
as described in Definition \ref{d-proper_covering}. Fixing an object $x\in X$,  we denote by $X(x)$ the full subcategory whose object set consists of the single element $x$ and whose morphisms  are the elements of the isotropy group $G_x$. By $Y(x)$ we denote the full subcategory whose object set is the finite set $F^{-1}(x)\subset Y$ and denote the associated set of morphisms by 
$\text{mor}(F^{-1}(x))$. We obtain  the
functor
$$
F: Y(x)\to X(x).
$$
In view of the property (3) in Definition \ref{d-proper_covering} of a proper covering functor, the map 
$$
\text{mor}(F^{-1}(x))\rightarrow (G_x) {_{s}\times_F} F^{-1}(x)= G_x \times F^{-1}(x),\quad \phi\mapsto (F(\phi),s(\phi))
$$
is a bijection. This implies  
$$
\# \text{mor}(F^{-1}(x))= \# G_x\cdot \# F^{-1}(x).
$$

We decompose the finite set $F^{-1}(x)$ into isomorphism classes
$$
F^{-1}(x)=\Theta_1\cup\ldots \cup\Theta_q,
$$
where two objects are isomorphic if there exists an isomorphism in $\text{mor}(F^{-1}(x))$ between them. Clearly,  for any two elements 
$y,y'\in\Theta_j$, the isotropy groups $G_y$ and $G_{y'}$ are isomorphic.
\begin{lemma}\label{ABB}\index{L- Structure of ${\mathcal F}^{-1}_x$ (I)}
For every $j=1,\ldots,q$ we have the equality
$$
\# G_x\cdot \# \Theta_j = \#\{\psi\in \text{mor}(F^{-1}(x))\ |\ s(\psi),t(\psi)\in \Theta_j\}.
$$
Moreover, for every $j=1,\ldots ,q$ and $y_0\in \Theta_j$, 
$$
\#\{\psi\in \text{mor}(F^{-1}(x))\ |\ s(\psi),t(\psi)\in\Theta_j\} = \# G_{y_0}\cdot( \#\Theta_j)^2
$$
and 
$$
\# G_x =\# G_{y_0}\cdot \# \Theta_j.
$$
\end{lemma}
\begin{proof}
The map
\begin{equation}\label{x_burger}
\{\psi\in \text{mor}(F^{-1}(x))\ |\ s(\psi), t(\psi)\in \Theta_j\}\rightarrow G_x\times \Theta_j:\psi\rightarrow (F\psi),s(\psi))
\end{equation}
is a bijection by property (3) of a proper covering functor, implying the first claim.

We decompose the subset  $\{\psi\in \text{mor}(F^{-1}(x))\ |\ s(\psi),t(\psi)\in \Theta_j\}$ of morphisms  into the disjoint union
$$
\{\psi\in \text{mor}(F^{-1}(x))\ |\ s(\psi),t(\psi)\in \Theta_j\}=\sqcup_{(y,y')\in \Theta_j\times \Theta_j}\text{mor}(y,y').
$$
Since any two points in $\Theta_j$ are isomorphic, each of these sets $\text{mor}(y,y')$ is nonempty. If we choose  $y_0\in \Theta_j$ and  fix for every $y\in \Theta_j$ a morphism $\psi_y: y_0\rightarrow y$, 
then the map
$$
\text{mor}(y_0,y_0)\rightarrow \text{mor}(y,y'),\quad  \psi\mapsto  \psi_{y'}\circ\psi\circ\psi_{y}^{-1}
$$
is a bijection. Consequently, 
$$
\#\{\psi\in \text{mor}(F^{-1}(x))\ |\ s(\psi),t(\psi)\in\Theta_j\} = \# G_{y_0}\cdot( \#\Theta_j)^2.
$$
Combining this with the first statement we obtain the final claim of Lemma \ref{ABB}.
\qed \end{proof}

Let again $F:Y\rightarrow X$ be the above  proper covering functor between ep-groupoids and 
assume that $x,x'\in X$ are isomorphic.  We consider,  for given $y\in F^{-1}(x)$,  all objects  $y'\in F^{-1}(x')$ for which  there exists a morphism
$\phi:y\rightarrow y'$. If $y_1\in F^{-1}(x)$ is isomorphic to $y$ and $y_1'$  by the morphisms  $\psi:y_1\rightarrow y$, and $\psi':y_1\rightarrow y_1'$,  we see that
$y'$ and $y_1'$ are isomorphic. If  $y'$ is reachable by an isomorphism starting in $\Theta_j$ and $y_1'$ is isomorphic to
$y'$, then the latter is also reachable by an element starting in $\Theta_j$.  Hence we have proved the following lemma.
\begin{lemma}\label{LEM1}\index{L- Structure of ${\mathcal F}^{-1}_x$ (II)}
Given a proper covering functor $F:Y\rightarrow X$, we  assume that $x$ and $x'$ are isomorphic. Then there is a one-to-one correspondence 
between the isomorphism classes associated to $F^{-1}(x)$ and $F^{-1}(x')$, respectively.
If  $ F^{-1}(x)=\sqcup_{j=1}^q \Theta_j$  is  isomorphisms class decomposition, then there is the decomposition
$F^{-1}(x')=\sqcup_{j=1}^q\Theta_j'$ having  the property that if  $y\in F^{-1}(x)$ and $y'\in F^{-1}(x')$  are connected by an isomorphism, 
there exists a uniquely determined $j$ such  that $y\in\Theta_j$ and $y'\in\Theta_j'$.
\end{lemma}
Therefore, there is  a one-to-one correspondence between the  isomorphism classes for which  $\Theta_j$ corresponds to $\Theta_j'$ for $j=1,\ldots ,q$.
Two points $y_0\in \Theta_j$ and $y_0'\in\Theta_j'$ are isomorphic and therefore $\# G_{y_0}=\# G_{y_0'}$.
Moreover, since $x$ and $x'$ are isomorphic,  $\# G_x=\# G_{x'}$. Using the third claim of Lemma \ref{LEM1}, 
we infer that
$$
\#\Theta_j= \frac{\# G_x}{\# G_{y_0}}=\frac{\# G_{x'}}{\# G_{y_0'}}=\#\Theta_j'.
$$
We sum up the discussion in the 
following lemma.
\begin{lemma}\label{LEM2}\index{L- Structure of ${\mathcal F}^{-1}_x$ (III)}
Given a proper covering functor $F:Y\rightarrow X$, we  assume that $x$ and $x'$ are isomorphic. Then, with the decompositions
of $F^{-1}(x)$ and $F^{-1}(x')$ as described in Lemma \ref{LEM1}, 
we have the equality $\#\Theta_j=\#\Theta_j'$ for every  $j=1,\ldots ,q$.
In particular, 
$$
\# F^{-1}(x)=\sum_{i=1}^q \#\Theta_j=\sum_{i=1}^q\#\Theta_j'=\# F^{-1}(x').
$$
\end{lemma}
The structure of isomorphism classes of $F^{-1}(x)$ stays the same if we perturb $x'$. This follows immediately from the definition of a proper covering functor.
In particular, together with Lemma \ref{LEM2} we obtain the following result.
\begin{proposition}\index{P- Structure of ${\mathcal F}^{-1}_x$}
Let $F:Y\rightarrow X$ be a sc-smooth proper covering functor. Then for every connected component $c$ of the orbit space $|X|$ there exists a constant $k_c$ such  that the following holds.
If $z\in c$ and $\pi(x)=z$, then $\# F^{-1}(x) =k_c$. In other words,  viewing ${\mathbb N}=\{1,2,\ldots \}$ as an ep-groupoid with the identities as the morphisms,  
the map
$$
X\rightarrow {\mathbb N},\quad x\mapsto \# F^{-1}(x)
$$
is a sc-smooth functor.
\end{proposition}
\subsection{The Structure of Strong Bundle Coverings}
Recall Definition \ref{proper_sb_covering} and consider the diagram
$$
\begin{CD}
W @>A>> V\\
@V P VV @V QVV\\
Y @>F>> X,
\end{CD}
$$
where the properties of the maps have been specified in the before-mentioned definition.
A strong bundle $P:W\rightarrow Y$ over an ep-groupoid is by definition given as a strong bundle $W\rightarrow Y$ over the object M-polyfold
together with a map $\mu_Y$ which associates to a morphism $\phi:y\rightarrow y'$ a topological linear isomorphism $\mu_Y(\phi):P^{-1}(y)\rightarrow P^{-1}(y')$.
A morphism $w\rightarrow w'$ is given by a pair $(\phi,w)$ with $s(\phi)=P(w)$ so that
$$
(\phi,w):w\rightarrow \mu_Y(\phi)(w).
$$
The same holds for $Q:V\rightarrow X$.  The functor $A$ applied to morphisms has the form 
$$
A(\phi,w) = (F(\phi),A_{s(\phi)}(w)),
$$
where $A_y:P^{-1}(y)\rightarrow Q^{-1}(F(y))$ is the topological linear isomorphism induced by $A$ between the fibers of the object spaces.
By definition
$$
A(\phi,w):  A_{s(\phi)}(w)\rightarrow \mu_X(F(\phi))(A_{s(\phi)}(w))
$$
from which we deduce that 
\begin{eqnarray*}
\mu_X(F(\phi))\circ A_{s(\phi)}(w)&=&\mu_X(F(\phi))(A_{s(\phi)}(w))\\
&=&t(A(\phi,w))=A(t(\phi,w))\\
&=&A(\mu_Y(\phi)(w))\\
&=& A_{t(\phi)}\circ \mu_Y(\phi)(w)
\end{eqnarray*}
implying the relationship
$$
\mu_X(F(\phi))\circ A_{s(\phi)} = A_{t(\phi)}\circ \mu_Y(\phi)\ \ \text{for}\ \ \phi\in \bm{Y}.
$$
Denote by $\text{Ban}$ the category of Banach spaces with the morphisms being linear topological isomorphisms.
We can view $\mu_Y:Y\rightarrow \text{Ban}$ and $\mu_X:X\rightarrow \text{Ban}$ as functors. Then 
the assignment $Y\rightarrow \text{mor}(\text{Ban}):y\rightarrow A_y$ is nothing else, but a natural transformation
$$
\mu_Y\xrightarrow{A} \mu_X\circ F.
$$

\chapter{Branched \texorpdfstring{Ep$^+$}{plusplus}-Subgroupoids}\label{CHAPTER_9}
In this chapter we  introduce the notion of a branched ep$^+$-subgroupoid and study its properties. 
These objects will arise naturally
as solution spaces of sc-Fredholm sections provided there is enough transversality. If a branched ep$^+$-subgroupoid is oriented 
and has a decent boundary geometry, we shall show that we can integrate sc-differential forms over them. Moreover a version of Stokes theorem holds.
The chapter ends with an appendix raising some questions and describing some ideas pointing to possible useful generalizations.
\section{Basic Definitions}
Recall from Chapter \ref{chap4} that 
a (finite-dimensional) submanifold  $M$ of a M-polyfold $X$ is a subset which is locally an sc$^+$-retract. 
It inherits the structure of a M-polyfold and in fact this structure is equivalent to a M$^+$-polyfold structure.
Since, as we have seen, M$^+$-polyfolds essentially are manifolds modulo the fact that the boundary might be 
not particularly nice, we shall call M$^+$-polyfolds just (smooth) manifolds, see  the next 
 \ref{REM911} were we recall
some more details. 
 Using this language 
we can say that a submanifold $M$ of a M-polyfold  $X$
inherits a natural (smooth) manifold structure.
\begin{remark}\label{REM911}\index{R- On finite-dimensional submanifolds}
If $m\in M$
the dimension of $T_mM$ is locally constant and by definition the local dimension of $M$.  A local sc$^+$-retract inherits from $X$ 
the structure of a M$^+$-polyfold.  As was shown previously,
if $m_0\in M$ with $d_M(m_0)=0$, then there exists an open neighborhood of $m_0$ in $M$ which has a natural (classical) smooth manifold
structure. Near points $m_0$ with $d_M(m_0)=1$ we have a natural structural of a smooth manifold with boundary. Hence a M$^+$-polyfold is near points with $d_M(x)=0$ the same as a smooth manifold without boundary and near points with $d_M(x)=1$
a manifold with boundary. We also introduced previously the notion of a tame M$^+$-polyfold and showed that this object is the same as 
a smooth manifold with boundary with corners. In generally the behavior of a M$^+$-polyfold near boundary points with $d_M(x)\geq 2$ can be quite complicated.
\qed
\end{remark}

The basic definition is that of a branched ep$^+$-subgroup\-oid of an ep-groupoid and it will be given next. 
\begin{definition}\label{DEF912}\index{D- Branched ep$^+$-subgroupoid}
Let $X$ be an ep-groupoid. A {\bf branched ep$^+$-subgroup\-oid} is a functor  $\Theta:X\rightarrow {\mathbb Q}^+$  which has the following properties.
\begin{itemize}
\item[(1)]\ If $\Theta(x)>0$ then $x$ is smooth, i.e. the support $\supp(\Theta)=\{x\in X\ |\ \Theta(x)>0\}$ of $\Theta$ is contained in $X_\infty$.
\end{itemize}
In addition we require that for every  $x\in \supp(\Theta)$ the following holds.
\begin{itemize}
\item[(2)]\ There exists  an open neighborhood $U(x)$  and a finite family of submanifolds $M_i\subset U(x)$, $i\in I$, with $x\in M_i$, so that
the inclusion $M_i\rightarrow U(x)$ is proper.
\item[(3)]\  There exists for  every $i\in I$  a positive rational number $\sigma_i>0$, called {\bf weight}\index{Weight}.
\item[(4)]\ With the data in (ii) and (iii) we have the equality
$$
\Theta(y)= \sum_{\{i\in I\ |\ y\in M_i\}} \sigma_i \ \ \text{for}\ \ y\in U(x).
$$
We shall refer to ${(M_i)}_{i\in I}$ and ${(\sigma_i)}_{i\in I}$ as a {\bf local  branching structure}.\index{Local branching structure}
\end{itemize}
We shall call $\Theta$ {\bf  $\bm{n}$-dimensional}\index{$n$-dimensional} provided  for every point $x\in \supp(\Theta)$ the manifolds of a local 
branching  structure are all   $n$-dimensional. We say that $\Theta$ is {\bf locally  equidimensional}\index{Locally equidimensional} provided if  for every $x\in \supp(\Theta)$
the manifolds in a sufficiently small branching structure have all the same dimension.
\qed
\end{definition}
\begin{remark}\index{R- On finite-dimensional branched ep$^+$-subgroupoids}
a) We can view $\Theta$ as data which distinguishes a full subcategory of $X$, namely the one
associated to the objects in $\supp(\Theta)$, and in addition associates to the objects $x$ of $ \supp(\Theta)$
a rational positive weight $\Theta(x)$.\par

\noindent b) The local branching structures at a point $x$ are in general not unique, even close to $x$.
For example assume that $X={\mathbb R}^2$ with the morphisms being the identities.
Take two smooth functions $\beta_i:{\mathbb R}\rightarrow {\mathbb R}$, $i\in \{1,2\}$, having the following properties,
$\beta_1(s)=0$ for $s\leq 0$, $\beta_1(s)>0$ for $s>0$, $\beta_2(s)>0$ for $s<0$, and $\beta_2(s)=0$ for $s\geq 0$.
Associated to the  $\beta_i$ we have their graphs $G_1$ and $G_2$ and give each of them the weight $1$ and define
$$
\Theta(x,y)=|{\{i\in \{1,2\}\ |\ (x,y)\in G_i\}} |.
$$
We note that there is another branching structure with $\gamma_1(s)=0$, $s\in {\mathbb R}$,  and $\gamma_2(s)>0$ for $s\neq 0$ describing $\Theta$. It is clear from this example that local branching structures can be 
arbitrarily complicated, f.e. take oscillatory graphs with high numbers of tangencies.
\qed
\end{remark}
A basic  observation is given in the following lemma. Recall that a subset $A$ of a topological space $Z$ is {\bf locally closed}\index{Locally closed} provided every point $a\in A$ has an open neighborhood
$V(a)$ in $Z$ so that $V(a)\cap A$ is closed in $V(a)$. Also recall the definition of an open subset $U$ of $X$ having the properness
property, which will be frequently used, see Definition \ref{proper*}.
\begin{lemma}\index{L- Local closedness of $\Theta$}
Let $\Theta:X\rightarrow {\mathbb Q}^+$ be a branched ep$^+$-subgroupoid. 
Then $|\supp(\Theta)|$ is locally closed in $|X|$.
\end{lemma}
\begin{proof}
Let $z\in |\supp(\Theta)|$ and pick a representative $x\in X$. Take an open neighborhood $U(x)$
with the natural $G_x$-action, and the properness property, which in addition 
allows for a local branching structure ${(M_i)}_{i\in I}$, ${(\sigma_i)}_{i\in I}$, representing 
$\Theta$. Define $V(z):=|U(x)|$ as the open neighborhood of $z$ and pick $w\in V(z)$ with
$w\in \cl_{|X|}(|\supp(\Theta)|)$. Take a representative $y\in U(x)$ of $w$.
The submanifolds $M_i$ are properly mapped $M_i\rightarrow U(x)$. 
If there is an open neighborhood $U(y)\subset U(x)$ which does not intersect any of the $M_i$
it follows, since $\bigcup_{i\in I} M_i$ is $G_x$-invariant, that $|U(y)|$ does not intersect $\cl_{|X|}(|\supp(\Theta)|)$. Hence we find a sequence $x_k\in \bigcup_{i\in I} M_i$ converging to $y$.
After perhaps taking a subsequence we may assume  $x_k\in M_i$ for a suitable $i\in I$.
Since $M_i$ is properly mapped the preimage of the compact set $(x_k)\cup \{y\}$ is compact
and $x_k$ has a convergent subsequence in $M_i$, which implies $y\in M_i$. This completes the proof.
\qed \end{proof}
The following definitions will be  useful later on.
\begin{definition}\label{DEF915}\index{D- Closed or compact $\Theta$}\index{Closed $\Theta$}\index{Compact $\Theta$}
Let $\Theta:X\rightarrow {\mathbb Q}^+$ be a branched ep$^+$-subgroupoid.
\begin{itemize}
\item[(1)]\ We call $\Theta$ {\bf closed} provided $|\supp(\Theta)|$ is closed in $|X|$.
\item[(2)]\ We call $\Theta$ {\bf compact} provided $|\supp(\Theta)|$ is a compact subset of $|X|$.
\end{itemize}
\qed
\end{definition}
If $V(x')\subset U(x)$, $x$ smooth,  we can replace $M_i$ by $N_i=M_i\cap V(x')$ and obtain an induced representation using the same weights on $V(x')$
$$
\Theta(y)=\sum_{\{i\in I\ |\ y\in N_i\}} \sigma_i \ \ \text{for}\ \ y\in V(x').
$$
Note that some (or all) of the $N_i$ might be empty. 
We can always find symmetric representations.
Namely, take $V(x)\subset U(x)$ so that $V(x)$ admits the natural $G_x$-action. Without loss of generality
assume that for $V(x)$ the manifolds $(M_i)$ and weights $(\sigma_i)$ are given, where $i\in I$.
Then define a new index set $J:=I\times G_x$, and  $M_{i,g}:= g\ast M_i$.  With the weights $\sigma_{(i,g)}$ defined 
by
$$
\sigma_{(i,g)}=\frac{1}{|G_x|}\cdot \sigma_i,
$$
we obtain a new local branching structure ${(M_{(i,g)})}_{(i,g)\in J}$, ${(\sigma_{(i,g)})}_{(i,g)\in J}$. Indeed we compute, using that $\Theta$ is a functor.
\begin{eqnarray*}
\sum_{\{(i,g)\in J\ |\ y\in g\ast M_i\}} \sigma_{(i,g)}&=& \frac{1}{|G_x|}\cdot \sum _{\{(i,g)\in J\ |\ g^{-1}\ast y\in M_i\}}\sigma_i\\
&=& \frac{1}{|G_x|}\cdot \left(\sum_{g\in G_x} \left(\sum_{\{i\in I\ |\ g^{-1}\ast y\in M_i\}} \sigma_i\right)\right)\\
&=&\frac{1}{|G_x|}\cdot\sum_{g\in G_x} \Theta(g^{-1}\ast y)\\
&=&\Theta(y).
\end{eqnarray*}
The interesting point is that $G_x$ acts on $J$ via $(h,j)\rightarrow h(j)$,  and in addition  $h\ast M_{j}=M_{h( j)}$,  and  $\sigma_{j}=\sigma_{h(j)}$ for $h\in G_x$ and $j\in J$.
Hence we have established the following result.
\begin{lemma}
Let $\Theta:X\rightarrow {\mathbb Q}^+$ be a branched ep$^+$-subgroupoid of the ep-groupoid $X$. Then there exists for every 
object $x\in X$ with $\Theta(x)>0$ an open neighborhood $U(x)$ admitting the natural $G_x$-action,  submanifolds $(M_i)$, $i\in I$,
properly embedded in $U(x)$, a $G_x$-action on $I$,  and rational weights $\sigma_i>0$,  such that the following holds.
\begin{itemize}
\item[(1)]\ $g\ast M_i = M_{g(i)}$ for all $i\in I$.
\item[(2)]\ $\sigma_{i}=\sigma_{g(i)}$ for $g\in G_x$ and $i\in I$.
\item[(3)]\ $\Theta(y)=\sum_{\{i\in I\ |\ y\in M_i\}}\sigma_i$  for $y\in U(x)$. 
\end{itemize}
\qed
\end{lemma}
We can distinguish different types of branched ep$^+$-subgroupoids.
\begin{definition}[Types of $\Theta$]\label{DEF917}\index{D- Types of $\Theta$}
Let $X$ be an ep-groupoid and $\Theta:X\rightarrow {\mathbb Q}^+$ a branched ep$^+$-subgroupoid.
\begin{itemize}
\item[(1)]\ We say $\Theta$ is of {\bf manifold-type}\index{Manifold-type $\Theta$} provided $\Theta$ only takes the values $0$ and $1$,
and between objects  $x,y\in \supp(\Theta)$ there is at most one morphism.
\item[(2)]\ We say that $\Theta$ is of {\bf orbifold-type}\index{Orbifold-type $\Theta$} provided $\Theta$ only takes the values $0$ and $1$.
\end{itemize}
\qed
\end{definition}

The basic result about  branched ep-subgroupoids of the flavors described above is given in the following proposition. 
\begin{proposition}[Natural smooth structures]\label{PROPY918}\index{P- Smooth structures from types}
Let $X$ be an ep-groupoid and $\Theta:X\rightarrow {\mathbb Q}^+$ a branched ep$^+$-subgroupoid. Consider the orbit space $|\supp(\Theta)|\subset |X|$ with the induced
metrizable topology. 
\begin{itemize}
\item[{\em(1)}]\  If $\Theta$ is of manifold-type then the ep-groupoid $X$ induces on the orbit space 
$S=|\supp(\Theta)|$  a natural smooth manifold structure. 
\item[{\em(2)}]\  If $\Theta$ is of orbifold-type then $X$ induces on $S=|\supp(\Theta)|$ a natural smooth orbifold structure.
\end{itemize}
\end{proposition}
\begin{remark}\label{remark918}
Here the notion of manifold and orbifold  are the M$^+$-polyfold versions. In particular if $d_{|X|}(z)=0$ for a every $z\in |\supp(\Theta)|$ 
it holds that $S$ has a natural (classical) smooth manifold structure in case (1) and a natural (classical) smooth orbifold structure  in case (2).
\qed
\end{remark}
\begin{proof}
Denote by $Z$ the orbit space $|\supp(\Theta)|$ and equip it with the induced (metrizable) topology
from $|X|$. 
For every $x$ with  $\Theta(x)>0$ we find an open neighborhood
$U(x)$ with the natural  $G_{x}$-action so that $U(x)$ has the properness property, and properly embedded
manifolds $M_i$, $i\in I$, with $x\in M_i$, and rational weights $\sigma_i>0$ allowing us to write
$$
\Theta(y)=\sum_{\{i\in I\ |\ y\in M_i\}} \sigma_i\ \ \text{for}\ y\in U(x).
$$
Also $G_{x}$ acts on $I$ and $g\ast M_i=M_{g(i)}$.

First we draw a conclusion from the assumption that $\Theta$ only takes values in $\{0,1\}$.
Assuming that $M_i\neq M_j$ for some $i\neq j$
we find $y$ such that $y\in M_i$ and $y\not \in M_j$, which implies 
$$
0<\Theta(y)<\Theta(x)=1.
$$
 This is a contradiction
since $\Theta$ either takes the value $0$ or $1$.  Having established that all $M_i$ are identical we can replace our local branching structure by a single submanifold, denoted by $M_{x}$, properly embedded in $U(x)$ and the weight $\sigma=1$.
This implies the representation
$$
\Theta(y)= \delta_{M_{x}}(y)=\left[\begin{array}{cc}
1&\text{if}\ y\in M_x\\
0&\text{if}\ y\in U(x)\setminus M_{x}
\end{array}
\right.
$$
Define the maps $\Psi_x:M_x\rightarrow Z:w\rightarrow |w|$ for $x\in \supp(M)$.

In order to prove (1) we employ the assumption that between any two points in $\supp(\Theta)$ there is at most one 
isomorphism. The first consequence of this assumption is that every $\Psi_x$ is injective. 
Indeed, we must have that $G_{x}=\{1_{x}\}$ because otherwise there would be at least two isomorphisms 
from $x$ to itself.  Assuming $\Psi_x(w)=\Psi_x(w')$ for $w,w'\in M_x$ 
there exists $g\in G_{x}$ with $\Gamma(g,w)=w'$. Since  $g= 1_{x_z}$ we deduce that $w=w'$.

Assume that  $\Psi_x$ and $\Psi_{x'}$ have an intersecting image.   Pick $w\in M_x$ and $w'\in M_{x'}$ 
so that $\Psi_x(w)=\Psi_{x'}(w')$. Hence there exists a morphism $\phi: w\rightarrow w'$ in $\bm{X}$.
We find open neighborhoods $U(w)\subset U(x)$,  $U(w')\subset U(x')$, and $U(\phi)\subset \bm{X}$ so that
$s:U(\phi)\rightarrow U(w)$ and $t:U(\phi)\rightarrow U(w')$ are sc-diffeomorphisms, which we use to define 
the sc-diffeomorphism $f:U(w)\rightarrow U(w')$ by
$$
f(v) = t\circ (s|U(\phi))^{-1}(v).
$$
If $v\in U(w)\cap M_x$ it follows that $f(v)\in M_{x'}$ since both manifolds represent $\Theta$. If $v'\in U(w')\cap M_{x'}$ we similarly conclude that $f^{-1}(v')\in M_x$.
We note that $f(v) = \Psi_{x'}^{-1}\circ \Psi_x(v)$ which implies that $\Psi_x$ and $\Psi_{x'}$ are sc-smoothly compatible. 
Since the $\Psi_x$ are defined on smooth manifolds $M_x$  we can use them to construct ordinary M$^+$-polyfold charts for the metrizable space $Z$,
which then by the previous discussion are smoothly compatible. It is evident that the smooth manifold structure on $Z$ 
does not depend on the choices of the  specific $U(x)$ picked for $x$.  The easy argument follows the above line and is left to the reader.
This completes the proof of (1).

In order to establish (2) we consider the maps $\Psi_x : M_x\rightarrow Z$.  In this case $G_x$ acts on $M_x$
fixing $x$, and the maps $\Psi_x$ are the restrictions 
of the projections $\pi_x:U(x)\rightarrow |X|$ which induce homeomorphisms $_{G_x}\backslash U(x)\rightarrow |U(x)|$, i.e.
$$
\Psi_x =\pi_x|M_x.
$$
Hence $\Psi_x:M_x\rightarrow Z$ are continuous maps inducing homeomorphism when we pass to the $G_x$-quotient.
Given $x,x'\in \supp(\Theta)$ define ${\bm{M}}(x,x')$ by
$$
\bm{M}(x,x')=\{\phi\in\bm{X}\ |\ s(\phi)\in M_x,\ t(\phi)\in M_{x'}\}.
$$
This is submanifold of $\bm{X}$ and has therefore an induced smooth manifold structure for which $s:{\bm{M}}(x,x')\rightarrow M_x$ and $t:{\bm{M}}(x,x')\rightarrow M_{x'}$ are local (classical) diffeomorphisms.
Denote by $M$ the disjoint union of the $M_x$. We can view this manifold as the object space of a small category,
where the morphism set is the disjoint union of all $\bm{M}(x,x')$.  This category, denoted by $M$, turns out to be 
an \'etale proper Lie groupoid, see \cite{Mj}  (Note that our notion of smooth manifold is more general since we use M$^+$-polyfolds, see Remark
\ref{REM911}).  We have a natural functor 
$$
\beta:M\rightarrow Z: \left[\begin{array}{cc}
\beta(y) = \pi_x(y),& y\in M_x\\
\beta(\phi) =1_{\beta(s(\phi))}, & \phi\in \bm{M}(x,x')
\end{array}
\right.
$$
where $Z$ is viewed as a category  with only the identities as morphisms.
By construction, passing to the orbit space $|M|$ we obtain the homeomorphism $|\beta|:|M|\rightarrow |Z|=Z$.  Hence it equips $Z$ with an orbifold structure. Again different choices for $x\in \supp(\Theta)$
would lead to a different  \'etale proper Lie groupoid $M'$  (in our M$^+$-context) and $\beta':M'\rightarrow Z$. However, this data defines an equivalent orbifold structure.  The metrizable space $Z$ equipped with such an equivalence class
of orbifold structures is called an orbifold. The equivalence classes appear since in the construction choices have to be made.
Of course, the nature of such an equivalence class has to be investigated. We leave this to the reader. There are no surprises
and a reader knowledgeable with the notion of Morita equivalence of \'etale proper Lie groupoids will be able to carry this out.
For details in the Lie groupoid case,  see \cite{Mj} or \cite{AE}, and for the polyfold version for topological spaces see the constructions in Section \ref{SEC2} and Part IV. The relevant fibered product constructions are similar to those described in the next chapter \ref{SEC2}.
\qed \end{proof}

Next we show that a branched ep$^+$-subgroupoid $\Theta:X\rightarrow {\mathbb Q}^+$ has a natural dimensional decomposition.
\begin{proposition}\label{PROP919}\index{D- Dimensional decomposition of $\Theta$}
Let $X$ be an ep-groupoid and $\Theta:X\rightarrow {\mathbb Q}^+$ a branched ep$^+$-subgroupoid.
Then there is a uniquely determined sequence $\Theta_0,\Theta_1,..$ of branched ep$^+$-subgroupoids
$$
\Theta_k:X\rightarrow {\mathbb Q}^+,\ k\in \{0,1,...\},
$$
uniquely characterized by the following two properties.
\begin{itemize}
\item[{\em(1)}]\ $\Theta_k$ is $k$-dimensional.
\item[{\em(2)}]\ $\Theta=\Theta_0+\Theta_1+....$. 
\end{itemize}
We also note that the sum in (2) is locally finite.
\end{proposition}
\begin{proof}
  If $x$ does not belong to the support of $\Theta$ 
we define $\Theta_k(x)=0$ for every $k\in \{0,1,2,..\}$.
If  $x\in X$ is  a point in $\supp(\Theta)$ we take 
on a suitable open neighborhood $U(x)$, invariant under the $G_x$-action, and having a symmetric  local branching structure ${(M_i)}_{i\in I}$,
${(\sigma_i)}_{i\in I}$. We can decompose $I=I_0+....+I_n$, where $M_i$ for $i\in I_k$
is $k$-dimensional. Since $I$ is finite we see that $I_k=\emptyset $ for $k$ large enough.
Since the action of $G_x$ preserves dimension the sets $I_k$ are $G_x$-invariant.
Associated to ${(M_i)}_{i\in I_k}$ and ${(\sigma_i)}_{i\in I_k}$ we can define $\Theta_k$ on $U(x)$ in the obvious way by
$$
\Theta_k(y)=\sum_{\{i\in I_k\ |\ y\in M_i\}} \sigma_i.
$$
We can do this for every $x\in \supp(\Theta)$ and the local definitions are compatible
since morphisms preserve dimension.  By construction the sum is locally finite.
The uniqueness of the decomposition is obvious.
\qed \end{proof}

\section{The Tangent and Boundary  of \texorpdfstring{$\Theta$}{plu}}\label{SECX9.2}
In the preparation for defining  orientations for branched ep$^+$-subgroupoids we need to introduce the notion
of the tangent of $\supp(\Theta)$. As a by-product we shall also see that associated to $\Theta$
there is a branched ep-subgroupoid $T\Theta:TX\rightarrow {\mathbb Q}^+$, where, not surprisingly, the local branching structure 
on $TU(x)$ is given by ${(TM_i)}_{i\in I}$, ${(\sigma_i)}_{i\in I}$, if ${(M_i)}_{i\in I}$, ${(\sigma_i)}_{i\in I}$  is the one for $\Theta$ over $U(x)$.

Let $x$ be a smooth object in the ep-groupoid $X$. 
We consider finite formal sums 
$\sum_{i\in I} \sigma_i \cdot L_i$, where $I$ is a finite set, the $\sigma_i$ are positive rational numbers, and the $L_i$ are
smooth finite-dimensional linear subspaces of $T_xX$. 
With other words, we take first the free vector space over ${\mathbb Q}$ generated by all finite-dimensional smooth linear subspaces
and then take elements in the cone defined by sums with non-negative coefficients.
In particular  the trivial subspace denoted by
$O$ is one of the generators, so that for example  $1\cdot O + 3\cdot L$ is an admissible formal sum.
We identify
$$
\sigma_i\cdot L_i +\sigma_j \cdot L_j = (\sigma_i+\sigma_j)\cdot L
$$
provided $L_i=L_j=L$. This allows us to write any such sum as   $\sum \sigma_L\cdot L$,
where only a finite number of $\sigma_L$ are non-zero (and positive rational).  Sometimes we abbreviate
$$
\mathsf{L}:= \sum \sigma_L\cdot L.
$$
We denote the set of all such formal sums by  $\text{Gr}(x)$, and  define  $\text{Gr}(X)$ as the union
$$
\text{Gr}(X):=\bigcup_{x\in X_\infty} \text{Gr}(x).
$$
There is a natural projection map $\pi: \text{Gr}(X)\rightarrow X_\infty$.  We shall consider the points in the set $\text{Gr}(X)$ as the objects of a category.
The morphisms for this category are the pairs $(\mathsf{L},\phi)$ with $\phi\in \bm{X}_\infty$ and $\mathsf{L}\in \text{Gr}(s(\phi))$.
The source and target maps are given by 
$$
s(\mathsf{L},\phi)=\mathsf{L}
$$ 
 and 
$$
t\left(\sum\sigma_L\cdot L,\phi\right)=\sum\sigma_L\cdot T\phi(L),
$$
where we recall that for smooth $\phi$  we have the well-defined sc-isomorphism 
$$
T\phi:T_{s(\phi)}X\rightarrow T_{t(\phi)}X.
$$
The composition of morphisms $(\mathsf{L},\psi)$ and $\mathsf{K},\phi)$ with $s(\mathsf{L},\psi)=t(\mathsf{K},\phi)$ is defined by
$$
(\mathsf{L},\psi)\circ (\mathsf{K},\phi)=(\mathsf{K},\psi\circ\phi).
$$
The identity elements have the form $(\mathsf{L},1_{\pi({\mathsf{L}})})$, and
the projection $\pi:\text{Gr}(X)\rightarrow X_\infty$ extended to morphisms by $\pi(\mathsf{L},\phi)=\phi$ becomes a functor.

\begin{theorem}\label{P15.2*1}\index{T- Functor $\mathsf{T}_\Theta$}\index{$\mathsf{T}_\Theta$}
Let  $X$ be an ep-groupoid and $\Theta:X\rightarrow {\mathbb Q}^+$ a branched ep$^+$-sub\-group\-oid. Associated to $\Theta$
there is a natural section functor $\mathsf{T}_\Theta:X_\infty\rightarrow \text{Gr}(X)$, i.e. 
$$
\pi\circ T_\Theta = Id_{X_\infty},
$$
 which is characterized by the following properties.
\begin{itemize}
\item[{\em(1)}]\ If $x\in \supp(\Theta)$ it holds $\mathsf{T}_\Theta(x)=\sum_{i\in I} \sigma_i\cdot T_xM_i$, where $(M_i)$ and $(\sigma_i)$ is any local branching structure
for $\Theta$ at $x$.  If $x\not\in \supp(\Theta)$  it holds that $T_{\Theta} (x)$ is the zero-sum.
\item[{\em(2)}]\ For a morphism $\phi$ it holds that $\mathsf{T}_\Theta(\phi): \mathsf{T}_\Theta(s(\phi))\rightarrow \mathsf{T}_\Theta(t(\phi))$ is given by 
$\mathsf{T}_\Theta(\phi)=(\mathsf{T}_\Theta(s(\phi)),\phi)$.
\end{itemize}
\qed
\end{theorem}
The proof is given later after Lemma \ref{LEM1526}.
The difficult part is (1), where one needs to show that the definition does not depend on the local branching structure.  It suffices to prove the proposition in the equidimensional case and to define 
$$
\mathsf{T}_{\Theta} =\sum_{k=0}^\infty \mathsf{T}_{\Theta_k}.
$$
The proof needs a considerable amount of preparation and we follow the presentation in \cite{HWZ7}.
 
\begin{definition}\label{D15.2.2}\index{D- Good point}
Let $\Theta:X\rightarrow {\mathbb Q}^+$ be a branched equidimensional ep$^+$-subgroupoid and ${(M_i)}_{i\in I}$, ${(\sigma_i)}_{i\in I}$ be a local branching structure on $U=U(x)$, where $x\in \supp(\Theta)$.  Define $M_U = \bigcup_{i\in I} M_i$\index{$M_U$}. A point $y\in M_U$ is called a {\bf good point}\index{Good points}
(with the respect to the given branching structure) provided there exists an open neighborhood $V(y)\subset U(x)$ such that
$V(y)\cap M_i = V(y)\cap M_j$ for all $i,j\in I_y:=\{k\in I\ |\ y\in M_k\}$. 
\qed
\end{definition}
The notion of being a good point turns out to be independent of the local branching structure and an intrinsic notion
for $\Theta$ as shown in the next result.
\begin{lemma}[\cite{HWZ7}, Lemma 2.2]\label{L15.2.3}\index{L- Good points}
Assume that $\Theta:X\rightarrow {\mathbb Q}^+$ is a branched equidimensional ep$^+$-subgroupoid. For $x\in \supp(\Theta)$
consider on the open neighborhood  $U(x)$  two local branching structures resulting 
in $M_U$ and $N_U$ (which have to be equal). If $y\in M_U=N_U$ is a good point for the first 
local branching structure the same holds with respect to the  second local branching structure.
\end{lemma}
\begin{proof}
Abbreviate $A=M_U=N_U$ and assume that $y\in A$ is good for the first local branching structure.
We find an open neighborhood $V(y)\subset U(x)$ so that
\begin{eqnarray}\label{ART1}
&M_i\cap V(y)= M_{i'}\cap V(y) \ \ \text{for}\ \ i,i'\in I_y=\{i''\in I\ |\ y\in M_{i''}\}.&
\end{eqnarray}
We note that $i_0\not\in I_y$ implies $y\not\in\cl_X(M_{i_0})$. Indeed, that $M_{i_0}\subset U(x)$
is properly embedded  would imply that $y\in M_{i_0}$ and therefore $i_0\in I_y$.
Hence we can pick $V(y)$ so small that in addition to (\ref{ART1}) also the following holds
\begin{eqnarray}\label{ART2}
& M_i\cap V(y)=\emptyset\ \ \text{for}\ \ i\in I\setminus I_y.&
\end{eqnarray}
 Using the representation for $\Theta$ we see that 
for $z\in V(y)\cap A$
$$
\Theta(y) = \Theta(z).
$$
For the second branching structure ${(N_j)}_{j\in J}$, ${(\tau_j)}_{j\in J}$ we find an open neighborhood $W(y)\subset U(x)$
such that $J_z\subset J_y$ for $z\in W(y)\cap A$, which again is a consequence 
of the properness. This implies 
$$
\Theta(z) =\sum_{j\in J_z} \tau_j\leq \sum_{j\in J_y} \tau_j =\Theta(y)\ \text{for}\ z\in W(y)\cap A.
$$
On $V(y)\cap W(y)\cap A$ we must therefore have the equality $J_z=J_y$ because otherwise the left-hand side 
in the last expression would be strictly smaller, giving a contradiction.
Hence we have proved that with $Q(z)=V(y)\cap W(y)$ it holds that $N_j\cap Q(z)=N_{J'}\cap Q(z)$
for all $j,j'\in J_y$. This shows that $y$ is a good point for the second branching structure.
\qed \end{proof}
\begin{remark}\label{REMr3} \index{R- On good points}
(a)  In view of Lemma \ref{L15.2.3} being a good point for $\Theta$  is an intrinsic notion. 

\noindent (b)  At a good point $y$ we can take a suitable open neighborhood $V(y)$ and a submanifold
$M$ properly contained in $V(y)$ so that $\Theta(z)=\sigma\cdot \delta_M(z)$ for $z\in V(y)$,
where $\sigma>0$ is a rational weight and $\delta_M(y)=1$ if $y\in M$ and $0$ otherwise.

\noindent(c)  Since we can transport local branching structures around by morphisms we obtain the following result.
\qed
\end{remark}
\begin{lemma}\label{L15.2.4}\index{L- Good points}
Let $\Theta:X\rightarrow {\mathbb Q}^+$ be  a branched equidimensional ep$^+$-subgrou\-poid. 
Assume that $\phi\in \bm{X}$ so that $s(\phi)\in\supp(\Theta)$ (and then also $t(\phi)\in \supp(\Theta)$).
Then $s(\phi)$ is a good point if and only if $t(\phi)$ is a good point.
\end{lemma}

Following \cite{HWZ7}
we introduce the sets of good and bad points for $\Theta$.
\begin{definition}\index{$S^g_\Theta$}\index{$S^b_\Theta$}\index{D- Set of good points}\index{D- Set of bad points}
Let $\Theta:X\rightarrow {\mathbb Q}^+$ be a branched equidimensional ep$^+$-subgroupoid.  The {\bf set of good points} 
is denoted by $S^g_\Theta$ and consists of all points in $\supp(\Theta)$ which are good for $\Theta$
in the sense of Definition \ref{D15.2.2}. The {\bf set of bad points} is defined by $S^b_\Theta=\supp(\Theta)\setminus S^g_\Theta$.
\end{definition} 

In view of Remark \ref{REMr3} (c) following Lemma \ref{L15.2.3},   the notion of being a good point or bad point descends to orbit spaces.
With $|S^g_\Theta|$ and $|S^b_\Theta|$  being the induced subsets in $|X|$ it holds that 
$$
\pi^{-1}(|S^g_\Theta|)=S^g_\Theta\ \ \text{and}\ \ \pi^{-1}(|S^b_\Theta|)=S^b_\Theta,
$$
where $\pi:X\rightarrow |X|$ is the obvious map.

\begin{lemma}\label{LEM1526}
Let $\Theta:X\rightarrow {\mathbb Q}^+$ be a branched equidimensional ep$^+$-sub\-groupoid and $S^g_\Theta$ and $S^b_\Theta$
the set of good and bad points, respectively. 
\begin{itemize}
\item[{\em(1)}]\ The set $S^g_\Theta$ is open and dense, and $S^b_\Theta$ is closed and nowhere
dense in $\supp(\Theta)$. 
\item[{\em(2)}]\ Assume a bad point $x\in \supp(\Theta)$ and a local branching structure on $U(x)$,
say ${(M_i)}_{i\in I}$ and ${(\sigma_i)}_{i\in I}$ on $U(x)$ representing $\Theta$ on $U(x)$,
are given. Then for every $i\in I$ there exists a sequence of good points $(x_k)$ contained in
$M_i$ and converging to $x$.
\end{itemize}
\end{lemma}
\begin{proof}
\noindent(1) The set $\sg$ is open by  definition. Hence its complement $\sba$ is closed and we
show that it is nowhere dense. At this point the proof of (2) will imply the validity of (1).\par

\noindent(2) Fix a metric $d$ on object M-polyfold $X$ (Compatibility with morphisms is not required.). By taking a possibly smaller $U(x)$ we may assume that $M_i=M_{i'}$ in $U(x)$
provided it holds near $x$. By redefining the weights and the index set we may assume without loss of generality that
$M_i=M_{i'}$ if and only if $i=i'$.  Since $x$ is bad we must have at least two indices.
 
 Take   a monotonically decreasing sequence $(\varepsilon_k)$ of positive real numbers converging to 
$0$. Define $D_k(x)=\{y\in U(x)\ |\ d(x,y)\leq \varepsilon_k\}$. If  $\varepsilon_1$ is small enough
we may assume without loss of generality that every $D_k(x)$ is a closed subset of $X$.
For large $k$ the subsets $M_i\cap D_k(x)$ are compact. To see this recall that there exists
an sc$^+$-retraction $r:W(x)\rightarrow W(x)$ for a suitable open neighborhood $W(x)$ of $x\in X$, so that
$r(W(x))=W(x)\cap M_i$. If we take a sufficiently small closed $\varepsilon$-ball, which is contained in
$W(x)$, its image under $r$ will be compact using the compactness property of the inclusions
of the Banach spaces in the scale and the fact that $r$ is sc$^+$-smooth. In view of this we may assume
by taking perhaps $\varepsilon_1>0$ smaller that 
$M_i\cap D_k(x)$ is a compact metric space for all $k\geq 1$.
The same holds without loss of generality for
$$
A^k :=  \bigcup_{i\in I} (M_i\cap D_k(x)) \subset D_k(x).
$$
After fixing $i_0\in I$ and $k$,  consider for $i\in I\setminus \{i_0\}$ the closed subset
 of $M_{i_0}\cap D_k(x)$ defined
by
$$
\Sigma_k^i = M_i\cap M_{i_0}\cap D_k(x).
$$
 The
 complement 
$$
(M_{i_0}\cap D_k(x))\setminus \bigcup_{i\in I\setminus\{i_0\}} \Sigma^i_k
$$
it is an open subset of $M_{i_0}\cap D_k(x)$, which, of course, might be empty.
If it is nonempty  we we can pick a point $x_k$ in this set  having the additional property
that $d(x_k,x)<\varepsilon_k$, which then  necessarily is a good point. 
If the set is empty we have the identity
$$
M_{i_0}\cap D_k(x)=  \bigcup_{i\in I\setminus\{i_0\}} \Sigma^i_k.
$$
This means that at least one of the $\Sigma^i_k$ has non-empty interior, say
$\Sigma_{k}^{i_1}$. Pick a point $x_k'$ in the interior of the latter. If $x_k'$ lies for some $i\neq i_1$
in  $\Sigma_k^{i}$ but not in its interior,  we can move it slightly so that its does not belong to the set
$\Sigma_k^i$, but by staying
still in the first set $\Sigma_{k}^{i_1}$.  Continuing this way we find a $x_k\in M_{i_0}$
with $d(x,x_k)<\varepsilon_k$ so that for a nonempty $\Sigma_k^i$ we either do not belong to this set
or we belong to the interior. Clearly $x_k$ is a good point. By construction we have now a sequence 
of good points $(x_k)\subset M_{i_0}$ which converges to $x$.
This completes the proof of the lemma.
\qed \end{proof}
Now we are in the position to show that the functor $\mathsf{T}_\Theta$ is well-defined.
\begin{proof}[Theorem \ref{P15.2*1}]
Let $x\in \supp(\Theta)$. If $x$ is a good point the functor $\Theta$ can be  represented near $x$ by a single (near $x$ uniquely determined) submanifold
$M_x$ contained properly in an open neighborhood $U(x)$ and the rational weight $\sigma_x=\Theta(x)$.
Since $\Theta$ is constant on $M_x$ we have the identity 
$$
\Theta(y)=\Theta(x)\cdot\delta_{M_x}(y)\ \ \text{for}\ \ y\in U(x).
$$
Hence 
$$
\mathsf{T}_\Theta(y)=\Theta(y)\cdot T_yM_x\ \ \text{for}\ \ y\in U(x).
$$
This shows that $\mathsf{T}_\Theta$ is well-defined at good points.

Next we assume that $x\in\supp(\Theta)$ is a bad point. 
We take  a branching structure
on a suitable $U(x)$ representing $\Theta$, which we denote   by ${(M_i)}_{i\in I}$, ${(\sigma_i)}_{i\in I}$. 
Since $\Theta$ is equidimensional all occurring manifolds have the same dimension, say $n$.
 We define an  {\bf equivalence relation}\index{Equivalence relation $i\sim i'$} $\sim$  on $I$  by saying that $i\sim i'$ provided $T_xM_i=T_xM_{i'}$.
 An equivalence class is denoted by $[i]$.\index{$[i]$}
The next part of the construction is a finite induction which terminates  once we run out of indices. 
We formulate the needed assertion in the following proposition which is proved utilizing several lemmata.
\begin{proposition}\label{890}
Under the previously stated assumptions.
Starting with an $\sim$-equivalence class $\Xi$ of indices in $I$ there exists a finite sequence
of indices $i_1,..,i_p\in I$ such $[i_1]=...=[i_p]=\Xi$  with the following properties.
\begin{itemize}
\item[{\em(1)}]\ For every $e\in \{1,...,p\}$ there exists a  nonempty subset $I^e$ of $\Xi$ containing $i_e$. These subsets are mutually disjoint and $\Xi= I^1\sqcup..\sqcup I^p$. 
\item[{\em(2)}]\ For every $e\in \{1,...,p\}$ there exists a sequence of good points $(x_k^e)\subset U(x)$ converging to
$x$ such that $x_k^e\in M_i$ if and only if $i\in I^e$.
\end{itemize}
\end{proposition}
\begin{proof}
The consideration being local, 
we may assume that we are given a partial quadrant $C$ in the sc-Banach space $E$
and M$^+$-polyfolds ${(M_i)}_{i\in I}$ of dimension $n$ containing the point $x=0$.
The set $I$ is the disjoint union of equivalence classes and we fix one, say $\Xi$.
If $n=0$ the point $0$ is isolated and therefore a good point. Hence we may assume 
that $n\geq 1$.  Writing $\Xi=[i_1]$ define $L^1:=T_0M_{i_1}$ and consider
$L^1\cap C$. 
\begin{lemma}\label{1528}
Under the previously stated assumptions it holds $\text{int}_{L^1}(L^1\cap C)\neq \emptyset$.
\end{lemma}
\begin{proof} Take $M_{i}$ for $i\in \Xi$ and, with the considerations being local, we may assume
there exist a relatively open neighborhood $U\subset C$ of $0\in C$  and a sc$^+$-retraction
$r:U\rightarrow U$ with $r(U)=U\cap M_i= M_i$ and $r(0)=0$.

We  consider for every $h\in C$ with $h\in E_1$ the $C^1$-path 
$$
[0,\varepsilon)\rightarrow C:\tau\rightarrow r(\tau\cdot h).
$$
From this we conclude, differentiating at $\tau=0$,  the fact that $Dr(0)h\in C$.
Since $Dr(0)$ has  a continuous extension to $E_0$ we conclude that
$$
Dr(0)(C)\subset C.
$$
The operator $Dr(0):E_0\rightarrow L^1$ is surjective and therefore open.
Pick a nonzero vector $h_0$  in $C$ with $Dr(0)(h_0)=h_0$.
We find a close-by vector $h_1$ which lies in the interior of $C$ such that $Dr(0)h_1\neq 0$.
Hence we find a sufficiently small $\varepsilon >0$ such that $B_\varepsilon^{E_0}(h_1)\subset C$,
and   $Dr(0)(B_\varepsilon^{E_0}(h_1))\subset C\cap L^1$ consists of nonzero vectors,
and is an open neighborhood of $h_1$ in $C\cap L^1$. This completes the proof of the lemma.
\qed \end{proof}
Abbreviate $L:=L^1$ and pick  $h_0\in \text{int}_{L}(L\cap C)$ with $\norm{h_0}_{E_0}=1$.
Define for $\varepsilon>0$ 
$$
S_\varepsilon(h_0):=\{\tau\cdot h\ |\ \tau>0,\ h\in L,\  \norm{h}_{E_0}=1,\ \norm{h-h_0}_{E_0}<\varepsilon\},
$$
which lies for sufficiently small
 $\varepsilon>0$ in $\text{int}_{L}(L\cap C)$. We call $S_\varepsilon(h_0)$\index{$S_\varepsilon(h_0)$} the  {\bf open $\varepsilon$-cone}
 \index{$\varepsilon$-cone around $h_0$}
 around $h_0$ in $L$.  Moreover, in this case  the associated {\bf closed $\varepsilon$-cone}  $\bar{S}_\varepsilon(h_0):=\cl_L(S_\varepsilon(h_0))$ belongs to
 $L\cap C$ and is given by 
 $$
 \bar{S}_\varepsilon(h_0):=\cl_L(S_\varepsilon(h_0))=\{\tau\cdot h\ |\ \tau\geq 0,\ h\in L,\  \norm{h}_{E_0}=1,\ \norm{h-h_0}_{E_0}\leq\varepsilon\}.
 $$
\begin{lemma} \label{891}
Under the previously stated assumptions. Let $h_0\in \text{int}_L(L\cap C)$ with $\norm{h}_{E_0}=1$.
There exists $\delta >0$ and an sc-complement $Z$ of $L$ in $E$
such that suitable subsets of  $M_i$ with $i\in [i_1]$ can be written  as 
$$
a+B_i(a), \  \norm{a}_0<\delta,\ a\in \cl_L(S_\delta(h_0)),
$$
where $i\in [i_1]$, $B_i(0)=0$, $DB_i(0)=0$ and $B_i(a)\in Z$.
\end{lemma}
\begin{proof}
Pick an sc-complement $Z$ of $L$ so that $E=L\oplus Z$ and 
denote by $P:E\rightarrow E$ the associated sc-projection onto $L$.
We compute for $a\in L\cap C$ with $r=r_i$, the sc$^+$-retraction associated to $M_i$ near $0$
\begin{eqnarray}\label{892}
r(a)& =&\left(\int_{0}^{1} Dr(s\cdot a)ds\right)(a)\\
&=& Dr(0)(a) + \left(\int_0^1 (Dr(s\cdot a)-Dr(0))ds\right) (a)\nonumber\\
&=& a + \left(\int_0^1 (Dr(s\cdot a)-Dr(0))ds\right) (a)\nonumber\\
&=:& a + A(a)(a)\nonumber\\
&=& (a + PA(a)(a)) +(I-P)A(a)(a).\nonumber
\end{eqnarray}
We observe that by construction $A(0)=0$. 
The map  $a\rightarrow a+PA(a)a$ is defined on a suitable open neighborhood
$V(0)$ of $0$ in  $C\cap L$.   Hence we introduce the map
$$
\Phi: V(0)\rightarrow L:a\rightarrow a+PA(a)a.
$$
We need the following fact about $\Phi$.
\begin{lemma}\label{893}
Under the previously stated assumptions there exists a sufficiently small $\delta>0$ and $\sigma>0$ 
so that $ \cl_L (S_{\delta}(h_0))\subset C\cap L$, $\{a\in L\cap C\ |\ \norm{a}_{E_0}\leq   \sigma\}\subset V(0)$,  and in addition
the map
$$
\Phi:\{a\in \bar{S}_\delta(h_0)\ |\ \norm{a}_{E_0}\leq \sigma\}\rightarrow C\cap L
$$
is injective, has invertible derivatives,  and its image 
contains  the set 
$$
\{b\in \cl_L(S_{\delta/2}(h_0))\ |\ \norm{b}_{E_0}\leq \sigma/2\}.
$$
 The inverse of $\Phi$ restricted to the set $\{b\in \cl_L(S_{\delta/2}(h_0))\ |\ \norm{b}_{E_0}\leq \sigma/2\}$
defines a smooth map $\Psi$ with $\Psi(0)=0$
and $D\Psi(0)=Id_L$.
\end{lemma}
\begin{proof}
Given $h_0\in \text{int}_L(L\cap C)$ with $\norm{h_0}_{E_0}=1$ there exists $\varepsilon>0$
so that every $h\in L$ with $\norm{h}_{E_0}=1$ and $\norm{h-h_0}_{E_0}\leq \varepsilon$
satisfies $h\in C\cap L$. As a consequence, if $h\in L$ and $\norm{h- \norm{h}_{E_0}\cdot h_0}_{E_0}\leq \varepsilon\cdot \norm{h}_{E_0}$, then $h\in C\cap L$.  Now consider
for $a\in L\cap C$ (close enough to $0$) 
$$
\Phi(a)=a+PA(a)a =:a+\wt{A}(a).
$$
We have that $\wt{A}(0)=0$ and $D\wt{A}(0)=0$. Given $\tau\in (0,1)$ we therefore  find $\sigma>0$ so that for $a\in L\cap C$ with $\norm{a}_{E_0}\leq \sigma$ the following estimates hold
\begin{eqnarray}\label{1000}
&\norm{D\wt{A}(a)}_{L(E_0)}\leq \tau&\\
&\norm{\wt{A}(a)}_{E_0}\leq \tau\cdot \norm{a}_{E_0}&\nonumber \\
&\norm{\Phi(a)}_{E_0}\geq (1-\tau)\cdot \norm{a}_{E_0}.&\nonumber
\end{eqnarray}
Hence we estimate for $a\in L\cap C$ with $\norm{a}_{E_0}\leq \sigma$ and addition belonging to $\bar{S}_\delta(h_0)$
\begin{eqnarray*}
&&\norm{\Phi(a)-\norm{\Phi(a)}_{E_0}\cdot h_0}_{E_0}\\
&=& \norm{a +\wt{A}(a) - \norm{a +\wt{A}(a)}_{E_0}\cdot h_0}_{E_0}\\
&\leq & \norm{a-\norm{a}_{E_0}\cdot h_0}_{E_0} + 2\cdot \norm{\wt{A}(a)}_{E_0} \\
&\leq& \norm{a-\norm{a}_{E_0}\cdot h_0}_{E_0} + 2\tau \cdot \norm{a}_{E_0} \\
&\leq&(\delta +2\tau)\cdot\norm{a}_{E_0}\\
&\leq &(\delta +2\tau)\cdot (1-\tau)^{-1}\cdot  \norm{\Phi(a)}_{E_0}.
\end{eqnarray*}
Hence, if we take $\tau\in (0,1)$ and $\delta\in (0,\varepsilon)$ small enough,  we conclude for 
$a\in \bar{S}_\delta(h_0)$ with $\norm{a}_{E_0}\leq \sigma$  that $\Phi(a)\in \cl_L(S_\varepsilon(h_0))\subset L\cap C$.

Next we show that for sufficiently small $\sigma>0$, so that in particular the previous discussion  holds,
the image of $\Phi$ when restricted to points in $\bar{S}_\delta (h_0)$ with $\norm{a}_{E_0}\leq \sigma$
contains the points in $\bar{S}_{\delta/2}(h_0)$ with $\norm{a}_{E_0}\leq \sigma/2$.
 If $\sigma>0$ is small enough we may assume that (\ref{1000}) holds.

Pick $b\in \bar{S}_{\delta/2}(h_0)$ with $\norm{b}_{E_0}\leq \sigma/2$ and consider the fixed point problem
$$
 a=b-\wt{A}(a)=:T(a), 
 $$
where we look for a solution in $\bar{S}_\delta(h_0)$ with $\norm{a}_{E_0}\leq \sigma$.
The solution $a=a(b)$ is found by an iteration scheme starting with $b$ and applying successively 
$T$ to obtain the sequence $a_k:=T^k(b)$. We have to show that $b\in \bar{S}_{\delta/2}(h_0)$ with $\norm{b}_{E_0}\leq \sigma/2$
implies that $a=\lim T^k(b)$ exists and satisfies $a\in \bar{S}_{\delta}(h_0)$ with $\norm{a}_{E_0}\leq \sigma$.
We estimate
\begin{eqnarray}\label{999}
&&\norm{a_{k+\ell}-a_k}_{E_0} \\
&\leq &\sum_{e=1+k}^{k+\ell} \norm{a_e-a_{e-1}}_{E_0}\nonumber\\
&\leq& (\sum_{e=k}^{k+\ell-1} \tau^e)\norm{\wt{A}(b)}_{E_0}\nonumber\\
&\leq&
\tau^k\cdot (1-\tau)^{-1}\cdot \tau\cdot \norm{b}_{E_0}.\nonumber
\end{eqnarray}
We use this estimate in two ways. First of all we deduce that
\begin{eqnarray}\label{1001}
\norm{a_\ell-b}_{E_0}\leq (1-\tau)^{-1}\cdot \tau\cdot \norm{b}_{E_0}.
\end{eqnarray}
Hence, writing $a_\ell = b +o_\ell$
\begin{eqnarray}\label{1002}
&&\norm{a_\ell - \norm{a_\ell}_{E_0}\cdot h_0}_{E_0}\\
&=& \norm{b +o_\ell - \norm{b}_{E_0}h_0 +(\norm{b}_{E_0}-\norm{a_\ell}_{E_0})h_0}_{E_0}\nonumber\\
&\leq & \norm{b-\norm{b}_{E_0}h_0}_{E_0} +\norm{o_\ell}_{E_0} + |\norm{b}_{E_0}-\norm{a_\ell}_{E_0}|\nonumber\\
&\leq& \norm{b-\norm{b}_{E_0}h_0}_{E_0} + 2\norm{o_\ell}_{E_0}\nonumber\\
&\leq& (\delta/2) \norm{b}_{E_0} + 2 (1-\tau)^{-1}\tau \norm{b}_{E_0}\nonumber
\end{eqnarray}
If we take $\sigma$ small enough we can achieve that $\delta/2 + 2 (1-\tau)^{-1}\tau<\delta$.
From (\ref{1001}) we deduce using $\norm{b}_{E_0}\leq \sigma/2$,  if $\sigma$ is small enough (and consequently 
$\tau$ as small as we wish)
$$
\norm{a_\ell}_{E_0}\leq \norm{b}_{E_0}+\norm{a_\ell-b}_{E_0}\leq \sigma/2 + (1-\tau)^{-1}\tau \sigma/2\leq \sigma.
$$
From the previous discussion we conclude that for $\sigma$ small enough and $b\in \bar{S}_{\delta/2}(h_0)$ with $\norm{b}_{E_0}\leq \sigma/2$
the sequence $a_k:=T^k(b)$ is a Cauchy sequence which stays in $\bar{S}_\delta(h_0)$ and has norm at most $\sigma$.
Define $a(b)=\lim T^k(b)$ and observe $a(b)=b-\wt{A}(a_b)$.  It is clear that $b$ enters the discussion as a continuous parameter
so that the map $b\rightarrow a(b)$ is continuous.  We leave it to the reader to show that the map
$$
\Psi:\{b\in \bar{S}_{\delta/2}(h_0)\ |\ \norm{b}_{E_0}\leq \sigma/2\}\rightarrow \{a\in \bar{S}_{\delta}(h_0)\ |\ \norm{a}_{E_0}\leq \sigma\}:b\rightarrow a(b)
$$
is classically smooth.
\qed \end{proof}

Using (\ref{892}) and Lemma \ref{893} we can write on the subset of $\cl_L(S_{(\delta/2)}(h_0))$
consisting of points with $\norm{b}_{E_0}\leq \sigma/2$
$$
r\circ\Psi(b) = \Phi\circ \Psi(b) + (I-P)A(\Psi(b))(\Psi(b)) =: b+ B(b)
$$
as claimed in the Lemma \ref{891}.
\qed \end{proof}

 Next we can pick a sequence $x_k^1\rightarrow 0$ of good points lying on $M_{i_1}$
 in the $\delta/2$-sector 
 and we may assume that the associated set of indices $I_k^1$, consisting of those indices in $[i_1]$ to which $x_k^1$ belongs,
 is constant. We shall  denote it by
 $I^1$. Clearly $i_1\in I^1$. We can write $x_k^1 = \ell_k^1 + B_{i_1}(\ell_k^1)$.
 If $i\in [i_1]\setminus I^1$ we must have that $B_i(\ell_k^1)\neq  B_{i_1}(\ell_k^1)$.
 By perturbing $(\ell_k^1)$ slightly and perhaps after taking a subsequence 
 we my assume that $x_k^2:= \ell_k^2+B_i(\ell_k^2)$ are good points.
 They  do not lie in $M_{j}$ for $j\in I^1$ and consequently $I^2$ is disjoint. 
 If $i\in [i_1]\setminus (I^1\sqcup I^2)$ we perturb $\ell^2_k$ slightly 
 so that $x_k^3:= \ell_k^3+B_i(\ell_k^3)$ are good points. If the perturbation
 is small these do not lie in $M_i$  for $i\in I^1\sqcup I^2$ and we may assume 
 after taking a subsequence that $I^3_k$ is constant, i.e. $I^3$. Proceeding this way
 the induction will terminate and we are done proving the Proposition \ref{890}.
 \qed \end{proof}
 Consider the formal sum
$$
\sum_{e=1}^p \Theta(x^e_k)\cdot T_{x_k^e}M_{i_e}.
$$
This sum is  independent of the local branching structure since it only involves  
good points.  The coefficients are constant and the tangent  spaces converge as $k\rightarrow \infty$.
Hence we can talk about convergence of this formal sum.
Given a specific local branching structure, as  $k\rightarrow \infty$, it converges to 
$$
\sum_{e=1}^p \left(\sum_{i\in I^e}\sigma_i\right)\cdot T_{x}M_{i_e}
$$
which equals the
 sum
$$
\sum_{i\in [i_1]} \sigma_i \cdot T_xM_i.
$$
For any local branching structure describing $\Theta$ near $x$ the same result holds.
This holds for every $\sim$-equivalence class $\Xi\subset I$, 
which implies the main result, since the sum over these $\Xi$ equals 
$$
\sum_{i\in I} \sigma_i\cdot T_x M_i.
$$
This completes the proof of Theorem \ref{890}.
\qed \end{proof}

At this point we can define  the {\bf tangent}\index{Tangent of $\Theta$} of a $\Theta:X\rightarrow {\mathbb Q}^+$.
\begin{definition}\index{D- Tangent $T\Theta$}\label{DEFNX9212}
Let $X$ be an ep-groupoid and $\Theta:X\rightarrow {\mathbb Q}^+$ a branched 
ep$^+$-sub\-group\-oid. Then there exists a well-defined branched ep$^+$-subgroupoid
$T\Theta:TX\rightarrow {\mathbb Q}^+$ 
defined as follows. 
We put   $T\Theta(h)=0$  for $h\in T_xX$ with $\Theta(x)=0$. If $\Theta(x)>0$
we  take a local branching structure on $U(x)$,
say ${(M_i)}_{i\in I}$ and ${(\sigma_i)}_{i\in I}$, and define on $TU(x)$ the functor $T\Theta$ by
$$
T\Theta(h) =\sum_{\{i\in I\ |\ h\in TM_i\}} \sigma_i.
$$
We call $T\Theta$ the {\bf tangent} of $\Theta$.
\qed
\end{definition}
By the previous discussion this definition is independent of the local branching structure.
Moreover, we note that the local definitions fit together.  
Since we can transport local branching structures 
around using the local sc-diffeomorphims associated to the morphisms it suffices
to show that for two local branching structures on $U(x)$ with $x\in \supp(\Theta)$  the local definitions 
define the same object, i.e.
$$
\sum_{\{i\in I\ |\ h\in TM_i\}} \sigma_i = \sum_{\{j\in J\ |\ h\in TN_j\}} \tau_j\ \ \text{for}\ \ h\in TU(x).
$$
By the previous discussion it holds that $\mathsf{T}_\Theta (y) = \sum_{\{i\in I\ |\ y\in M_j\}} \sigma_i\cdot T_yM_i$
is well-defined and independent of the local branching structure taken. Hence for $y\in U(x)$ 
the number 
$$
\sum_{\{i\in I\ |\ h\in T_yM_i\}} \sigma_i,
$$
is intrinsic and independent of the local branching structure we picked. This shows that $T\Theta$ is well-defined.\par

Assume that $\Theta:X\rightarrow {\mathbb Q}^+$ is a branched ep$^+$-subgroupoid. There is a natural construction for a functor $\partial\Theta:X\rightarrow {\mathbb Q}^+$ called the {\bf boundary} \index{Boundary of $\Theta$}  of $\Theta$.
The support of $\partial\Theta$ is contained in $\partial X$. 
In general this functor has not a particularly good structure due to the fact that $\partial X$ might not have a good structure.
We shall later on describe situations in which $\partial \Theta (x)=\Theta(x)$ for all $x\in \partial X$ and $0$ otherwise.
One should note, however, that this is not(!) always the case. We start with a lemma discussing
the position of a submanifold intersecting $\partial X$. Recall for a M-polyfold as well as a M$^+$-polyfold
the degeneracy index $d_X$ which has been introduced in Definition \ref{M_polyfold_degeneracy _index}. We shall also utilize the notion
of the reduced tangent space $T^R_xX$ as well as that of the cone $C_xX$  (or depending on the situation the partial quadrant $C_xX$) at a smooth point  $x$.

\begin{lemma}[Characterization of $d_M(x)=0$]\label{15211}\label{15214}\index{L- Characterization of $d_M(x)=0$}
Let $M\subset X$ be a finite-dimensional submanifold of a M-polyfold $X$. The following two statements are equivalent.
\begin{itemize}
\item[{\em(1)}]\ $x\in M\cap \partial X$   with $d_M(x)=0$.
\item[{\em(2)}]\ $T_xM\subset T_x^RX$.
\end{itemize}
\end{lemma}
\begin{remark}\index{R- On points with $d_M(x)=0$}
Under the assumption of the lemma,  $M$ is a subset of $X$ which is locally an sc$^+$-retract. If $x\in M\cap \partial X$
it is possible that $x\not\in \partial M$. According to the lemma this precisely is the case when $T_xM\subset T^R_xX$.
Reformulating the statement,  $x\in \partial M$ if and only if $x\in\partial X$ and  $T_xM$ is not contained in the partial cone $C_xX$.
\qed
\end{remark}
\begin{proof}
The problem is entirely local. Hence we may assume $X=O$, where $(O,C,E)$ is a local model
and $x=0$. We may further assume that we are  given an open neighborhood $U$ of $0$ in $C$ containing $O$ and an sc-smooth retraction $r:U\rightarrow U$ with $r(U)=O$, and an sc$^+$-retraction $s:U\rightarrow U$ with $s(U)=M$.  \par

Assume first that (1) holds. Since $d_M(0)=0$ there exist an open neighborhood
$W\subset M$ of $0$ and an sc-diffeomorphism $\Phi:({\mathbb R}^n,0)\rightarrow (W,0)$, mapping $0$ to $0$.
Here $n=\dim(T_0M)$. For every $h\in {\mathbb R}^n$  it holds that $\Phi(h)=\Phi(h)-\Phi(0)\in C$
and consequently using the differentiability at $0$ we deduce that $T\Phi(0)(h)\subset C$
implying that $T_0M\subset T_0O \cap C=T^R_0O$. This shows that (2) holds.\par

Next assume (2). Without loss of generality $C=[0,\infty)^k\oplus W$, where $k=d_O(0)$,
and $E={\mathbb R}^k\oplus W$, in particular $W$ is an sc-Banach space. By assumption $T_0M\subset \{0\}\oplus W$. 
For $(0,h)\in T_0M$ small enough  we have that $(0,h)\in U$ and therefore $s(0,h)\in M \subset C$.
We compute similarly as before when we studied submanifolds that 
$$
s(0,h) = (0,h) + A(0,h) 
$$
with $A(0,0)=(0,0)$ and $DA(0,0)=0$, where $A$ is an sc-smooth map into $E$ defined on an open neighborhood $Q$
of $0$ in the smooth $T_0M$.  Take an sc-complement $Y$ of $T_0M$ so that $E=T_0M\oplus Y$
and denote by $P:E\rightarrow T_0M$ the associated sc-projection. We define
$$
\Psi: Q\rightarrow T_0M: (0,h)\rightarrow (0,h) +PA(0,h).
$$
Note that $D\Psi(0,0)=Id$ and the map is classically smooth so that by the implicit function theorem
$\Psi^{-1}$ is defined and classically smooth. Hence 
$$
s\circ \Psi^{-1}(0,h) = (0,h) + (Id-P)A(\Psi^{-1}(y))=: (0,h) + B(0,h),
$$
where $B(0,0)=0$, $DB(0,0)=0$, and $B(0,h)\in Y$.  By construction $s\circ\Psi^{-1}(0,h)\in C$. 
Define $\wt{s}(0,h)=s\circ\Psi^{-1}(0,h)$ for $(0,h)\in T_0M$ with $\norm{(0,h)}_{E_0}$ small.
Since $T\wt{s}(0,0)=Id: T_0M\rightarrow T_0M\subset T^R_0O\subset T_0O$
we can define for small $y\in Y$ the map
$$
((0,h),y)\rightarrow \wt{s}(0,h)+y,
$$
which is a local sc-diffeomorphism from an open neighborhood of $(0,0)$ in $T_0M\times Y$
to an open neighborhood of $0$ in $E$.  Denote the inverse by $F$ and note that the composition
$$
\wt{s}\circ P\circ F 
$$
is an sc$^+$-smooth map defined on an open neighborhood of $0$ in $E$ and has image 
$M$(near $0$).  Moreover, if defined on a suitable open neighborhood around $0$ it is a retraction.
Hence we see that $M$ near $0$ can be realized as an sc$^+$-retraction in $E$ with $C=E$.
This implies that $d_M(0)=0$.
\qed \end{proof}
Now we can state the crucial result needed to define $\partial\Theta$.
\begin{proposition}[Well-definedness of $\partial\Theta$]\index{P- Well-definedness of $\partial\Theta$}\label{PROPX9.2.15}
Let $\Theta:X\rightarrow {\mathbb Q}^+$ be a bran\-ched ep$^+$-subgroupoid and let $x\in \supp(\Theta)\cap \partial X$.
Suppose ${(M_i)}_{i\in I}$, ${(\sigma_i)}_{i\in I}$, and ${(N_j)}_{j\in J}$, ${(\tau_j)}_{j\in J}$ are two local branching structures
on a suitable open neighborhood $U(x)$ which represent $\Theta$ on $U(x)$. 
\begin{itemize}
\item[{\em(1)}]\ Then for every $y\in U(x)$ the equality
\begin{eqnarray}\label{polarizedx}
\sum_{\{i\in I\  |\ y\in \partial M_i\}}\sigma_i =\sum_{\{j\in J\ |\ y\in \partial N_j\}} \tau_j
\end{eqnarray}
holds.
\item[{\em(2)}]\ The sums in (\ref{polarizedx})  vanish if $y\in U(x)$ satisfies $d_X(y)=0$.
\item[{\em(3)}]\ There is a well-defined functor 
$\partial\Theta:X\rightarrow {\mathbb  Q}^+$
with
$\supp(\partial\Theta)$ $ \subset \partial X$, so that for $x\in \supp(\partial \Theta)$ it holds
$$
\partial \Theta(y) =\sum_{\{i\in I\ |\ y\in M_i,\ d_{M_i}(y)\geq 1\}}\sigma_i,\ \ y\in U(x).
$$
Here we take any local branching structure for $\Theta$ on some $U(x)$.
\end{itemize}
\end{proposition}
\begin{proof}
(1) Let ${(M_i)}_{i\in I}$, ${(\sigma_i)}_{i\in I}$ be a local branching structure on $U(x)$ and 
consider the sum $\sum_{\{i\in I\  |\ y\in \partial M_i\}}\sigma_i$. We also consider
$T_\Theta(y) = \sum_{\{i\in I\ |\ y\in M_i\}}\sigma_i\cdot T_yM_i$. The latter we know is independent of the local branching structure. We can rewrite this  as
\begin{eqnarray*}
T_\Theta(y) &=& \sum_{\{i\in I\ |\ y\in M_i\}}\sigma_i\cdot T_yM_i\\
&=&\sum_{\{i\in I\ |\ y\in M_i,\ d_{M_i}(y)\geq 1\}}\sigma_i\cdot T_yM_i  +
\sum_{\{i\in I\ |\ y\in M_i,\ d_{M_i}(y)=0\}}\sigma_i\cdot T_yM_i.
\end{eqnarray*}
In view of Lemma \ref{15211} and the remark thereafter,  this is a natural decomposition independent of the local branching structure
and from this the desired result (1) is easily obtained.\par

\noindent (2) Let  ${(M_i)}_{i\in I}$, ${(\sigma_i)}_{i\in I}$ be  a local branching structure on $U(x)$.
Assume that $d_X(x)=0$. Then, if $U(x)$ is small enough all the $d_{M_i}$ vanish identically
implying that $\sum_{\{i\in I\  |\ y\in \partial M_i\}}\sigma_i=0$.  
Hence for $y\in U(x)$ 
$$
\sum_{\{i\in I\ |\ y\in M_i,\ d_{M_i}(y)\geq 1\}}\sigma_i\cdot T_yM_i =0.
$$
(3) Follows immediately from (1) and (2). 

\qed \end{proof}

In view of Proposition \ref{15214} the boundary $\partial\Theta$ of $\Theta$ can be defined.
\begin{definition}\label{DEF9214}\index{D- Boundary of $\Theta$}\index{$\partial\Theta$}
Let $X$ be an ep-groupoid and $\Theta:X\rightarrow {\mathbb Q}^+$ a branched ep$^+$-groupoid.
The functor $\partial\Theta$, uniquely characterized by the following properties
\begin{itemize}
\item[(1)]\ $\supp(\partial\Theta)\subset \partial X$.
\item[(2)]\  For  $x\in \supp(\partial \Theta)$ it holds for any local branching structure on a suitable
$U(x)$
$$
\partial \Theta(y) =\sum_{\{i\in I\ |\ y\in M_i,\ d_{M_i}(y)\geq 1\}}\sigma_i,\ \ y\in U(x).
$$
\end{itemize}
 is called the {\bf boundary of $\Theta$}.
 \qed
\end{definition}
\begin{remark}\index{R- On the functor $\partial\Theta$}
The functor $\partial\Theta$ in general does not have particularly nice 
properties. In Section \ref{IandST},  Definition \ref{TAMERXX},
we shall introduce the notion of a tame branched ep$^+$-subgroupoid.
As we shall show their boundaries  $\partial\Theta$
are nice enough for an integration and orientation theory. 
\qed
\end{remark}

\section{Orientations}\label{SECXV93}
There is also  an oriented version of the $\text{Gr}(X)$-construction in the previous section, which we shall describe in the following.
\begin{remark}\index{R- Remark on \cite{HWZ7}} In  \cite{HWZ7} we used a particular definition of orientability 
of an ep$^+$-subgroupoid,  which is too special to be applicable to all occurring types of sc-Fredholm problems.
In \cite{HWZ7} we erroneously claimed that it does, see the  erratum \cite{HWZ7err} for more details.
However, the techniques from \cite{HWZ7} apply immediately to the general approach which 
is described below.
\qed
\end{remark}
 
An {\bf oriented}\index{Oriented} finite-dimensional linear vector space $L$ denoted by $\wh{L}$ is a pair $(L,o)$,
where $o$ is an orientation of $\Lambda^{\text{max}}L$, i.e. a choice of connected component in $\Lambda^{\text{max}}L\setminus\{0\}$.
In the case $L$ is the trivial space consisting only of the $0$-vector by definition $\Lambda^{\text{max}}L={\mathbb R}$.
If one of the components is denoted by $o$ we shall write for  the other one  $-o$.
For smooth $x$ denote by $\wh{\text{Gr}}(x)$ the set of formal sums
$$
\wh{\mathsf{L}}=\sum {\sigma_{\wh{L}}}\cdot\wh{L},
$$
where almost all $\sigma_{\wh{L}}$ vanish and the  $\wh{L}$  are  smooth, oriented, 
finite-dimensional linear subspaces of $T_xX$. We define the set
$$
\wh{\text{Gr}}(X):= \bigcup_{x\in X_\infty} \wh{\text{Gr}}(x),
$$
which again comes with a natural map into $X_\infty$, denoted by 
$$
\wh{\pi}:\wh{\text{Gr}}(X)\rightarrow X_\infty.
$$
 We can turn $\wh{\text{Gr}}(X)$ into a category
by viewing the elements $\wh{\mathsf{L}}$ as objects and defining the pairs $(\wh{\mathsf{L}},\phi)$,
with $s(\phi)=\wh{\pi}(\wh{\mathsf{L}})$, as morphisms. Again $s(\wh{\mathsf{L}},\phi)=\wh{\mathsf{L}}$  and 
$$
t(\wh{\mathsf{L}},\phi)=\sum \sigma_{(L,o_L)} \cdot (T\phi(L),T\phi_\ast(o_L)).
$$ 
There is a fiberwise 2-1 forgetful functor  fitting into the commutative diagram
$$
\begin{CD}
\wh{\text{Gr}}(X) @>\mathsf{f}>>\text{Gr}(X)\\
@VVV @VVV\\
X_\infty @=        X_\infty.
\end{CD}
$$
Note that $\sigma_{(L,o)} \cdot (L,o) +\sigma_{(L,-o)}\cdot (L,-o)$ is mapped to $(\sigma_{(L,o)}  +\sigma_{(L,-o)})\cdot L$.
If $\wh{\mathsf{L}}=\sum_{L\in\text{Gr}(x)} (\sigma_{(L,o)}\cdot (L,o) +\sigma_{(L,-o)}\cdot (L,-o))$ we denote $\mathsf{f}(\wh{\mathsf{L}})$ by $\mathsf{L}$. We observe that
$$
\mathsf{L}= \sum_{L\in \text{Gr}(x)} (\sigma_{(L,o)} +\sigma_{(L,-o)})\cdot L.
$$
\begin{definition}\index{D- Orientation for $\Theta$}\index{D- Orientability of $\Theta$}\label{DEFNC932}
An {\bf orientation}\index{Orientation} for the branched ep$^+$-subgroupoid $\Theta:X\rightarrow {\mathbb Q}^+$ is 
a choice of section functor  $\wh{\mathsf{T}}_\Theta:X_\infty\rightarrow \wh{\text{Gr}}(X)$ which lifts $\mathsf{T}_\Theta:X_\infty\rightarrow \text{Gr}(X)$, i.e.
$$
\mathsf{T}_\Theta=\mathsf{f}\circ \wh{\mathsf{T}}_\Theta
$$
 with the  property
that for every $x$ with $\Theta(x)>0$ there exists a local representation with oriented $(M_i,o_i)$ so that 
\begin{eqnarray}\label{15000}
\wh{\mathsf{T}}_\Theta(y) = \sum_{i\in I} \wh{\sigma}_i \cdot T_y(M_i,o_i).
\end{eqnarray}
We shall call $\wh{\mathsf{T}}_\Theta$ an {\bf oriented lift}\index{Oriented lift $\wh{\mathsf{T}}_\Theta$}
 of $T_\Theta$.
We say $\Theta$ is {\bf orientable}\index{Orientability} provided $T_\Theta$ admits an oriented lift.
\qed
\end{definition}

\begin{remark}
The local representation (\ref{15000}) associated to a local branching structure precisely means the local continuity
of the lift $x\rightarrow \wh{\mathsf{T}}_\Theta(x)$ on $\supp(\Theta)$. 
The reader should also observe that $\Theta$, if orientable, might allow more than two possible orientations
due to its `fractional nature'.
\qed
\end{remark}
We shall derive next some properties of an orientation $\wh{\mathsf{T}}_\Theta$ of $\Theta$.
Assume that we have fixed for every $x\in \supp(\Theta)$ an oriented local branching
structure ${(M_i^x,o^x_i)}_{i\in I_x}$ with weights ${(\sigma^x_i)}_{i\in I_x}$ so that
$$
\wh{\mathsf{T}}_\Theta(y) =\sum_{i\in I_x} {\sigma}^x_i \cdot T_y (M_i^x,o^x_i)\ \ \text{for}\ \ y\in U(x).
$$
In addition we have that $\mathsf{T}_\Theta(y)=\sum_{i\in I_x} {\sigma}_i^x\cdot T_{y}M_i$.
The interesting case is that some $L=T_xM_i$ occurs with the two different possible
orientations, i.e. there are two different indices $i,j\in I_x$ with $T_xM_i=T_xM_j$, but $o_i^x=-o^x_j$.
Hence  $\mathsf{T}_\Theta(x)$ has a term $\sigma\cdot L$ which comes
from  $\sigma^x_i \cdot (L,o^x_i) + \sigma^x_j\cdot (L,-o^x_j)$ with $\sigma^x_i+\sigma^x_j=\sigma$.

\begin{lemma}
Assume that $X$ is an ep-groupoid and $\Theta:X\rightarrow {\mathbb Q}^+$
a branched ep$^+$-subgroupoid. 
Then two orientations $\wh{\mathsf{T}}_\Theta^j$, $j=1,2$,  of $\Theta$, which agree at good points
 are the same. 
\end{lemma}
\begin{proof}
The branched ep$^+$-subgroupoid $\Theta$ defines $\mathsf{T}_\Theta$
and we assume that we are given oriented lifts $\wh{\mathsf{T}}_\Theta^j$, $j=1,2$,  satisfying (\ref{15000}), which
agree over the good points. Let $x\in \supp(\Theta)$ be a bad point 
and assume we have the oriented representations
$$
\wh{\mathsf{T}}^1_\Theta(y)= \sum_{i\in I_x}  \sigma_i\cdot T_y(M_i,o_i)\ \ \text{and}
\ \ \wh{\mathsf{T}}^2_\Theta(y)=\sum_{j\in J_x}\tau_j\cdot T_y(N_j,p_j)
$$
for $y\in U(x)$. 
We begin by studying the first orientation and associated branching structure.
Define an equivalence relation on $I_x$  by saying that $(M_i,o_i)$ and
$(M_{i'},o_{i'})$ are equivalent provided $T_xM_i=T_xM_{i'}$. Note that
at this point we ignore the orientations.
Denote the equivalence classes by $[i_1],...,[i_\ell]$. 

Employing Proposition \ref{890}, recall that 
for every such equivalence class $[i_h]$ there  exist sequences $(x_k^q)$ of good points converging to $x$, 
and a partition $[i_h]= I^1\sqcup...\sqcup I^e$ 
such that $x^q_k\in M_i$ if and only if $i\in I^q$.
Consider the formal sum
$$
\sum_{q=1}^e \wh{\mathsf{T}}_\Theta^1(x_k^q),
$$
where a term    $\wh{\mathsf{T}}_\Theta^1(x_k^q)$  has 
 the form
$$
\wh{\mathsf{T}}_\Theta^1(x_k^q) =\sum_{i\in I^q} \sigma_{i}\cdot T_{x_k^q}(M_i,o_i).
$$
We can define an equivalence relation on $I^q$ by saying that $i\approx i'$
provided 
$$
T_{x}(M_i,o_i)=T_{x}(M_{i'},o_{i'}).
$$
 Denote the two equivalence classes
by $I^q_\pm$ so that $I^1=I^q_+\sqcup I^q_{-}$. If $k$ is large enough 
$T_{x^q_k}(M_i,o_i)=T_{x^q_k}(M_{i'},o_{i'})$ if and only if $i,i'$ both lie in the same equivalence class.
Now we consider for all large $k$ the decomposition
$$
\wh{\mathsf{T}}_\Theta^1(x_k^q) =\sum_{i\in I^q_+} \sigma_{i}\cdot T_{x_k^q}(M_i,o_i)+
\sum_{i\in I^q_-} \sigma_{i}\cdot T_{x_k^q}(M_i,o_i).
$$
The  expression on the right would be the same for the second branching structure,
since the  $x^q_k$ are good points and the assumption $\wh{\mathsf{T}}_\Theta^1(x_k^q) =\wh{\mathsf{T}}_\Theta^2(x_k^q) $.
As $k\rightarrow \infty$ the right-hand side converges to
\begin{eqnarray}\label{16000}
\sum_{i\in I^q_+} \sigma_{i}\cdot T_{x}(M_i,o_i)+
\sum_{i\in I^q_-} \sigma_{i}\cdot T_{x}(M_i,o_i),
\end{eqnarray}
which is the part of $\wh{\mathsf{T}}_\Theta^1(x)$ associated to the tangent spaces $T_xM_i$, $i\in [i_h]$,
which are all equal, but in general have different orientations. 
If we apply the previous discussion to the second orientation we obtain the 
expression 
$$
\sum_{j\in J^q_+} \tau_{j}\cdot T_{x}(N_j,p_j)+
\sum_{j\in J^q_-} \tau_{j}\cdot T_{x}(N_j,p_j),
$$
which must equal (\ref{16000}). This sum is over indices $j$ so that $T_xN_j=T_xM_i$ with $i\in i_h]$ (or $j$ in the corresponding 
equivalence class $[j_{h'}]$.

The sum of these expressions define $\wh{\mathsf{T}}_\Theta^1(x)$ as well as $\wh{\mathsf{T}}_\Theta^2(x)$,  and since the ingredients are equal they coincide.
\qed \end{proof}

\begin{example}
We introduce four different smooth maps  $\beta^{\pm\pm}:{\mathbb R}\rightarrow {\mathbb R}$. 
The map $\beta^{++}$ satisfies  $ \beta^{++}(s)>0$ for $s\neq 0$, $\beta^{++}(0)=0$,  and has vanishing derivatives
of all order at $0$. Put $\beta^{--}:=-\beta^{++}$.
Define $\beta^{+-}$ by being equal to $\beta^{++}$ for $s\leq 0$ and equal to $\beta^{--}$
for $s>0$. Similarly we define $\beta^{-+}$. Associated to these maps 
we have their graphs denoted by $G^{++}, G^{--}, G^{+-}$, and $G^{-+}$.
Observe that
$$
|\{i \in \{++,--\}\ |\ (x,y)\in G^i\} | =|\{j\in \{+-,-+\}\ |\ (x,y)\in G^j\}|
$$
for all $(x,y)\in {\mathbb R}^2$. This defines $\Theta:{\mathbb R}^2\rightarrow {\mathbb Q}^+$.
Orient $G^{++}$ and $G^{+-}$  in such way that 
it is given at $(0,0)$ by the class of $(1,0)$ and orient $G^{--}$ and $G^{-+}$
so that the orientation at $(0,0)$ is given by the class of $(-1,0)$. 
Then $G^{++}$ and $G^{--}$ define a first orientation and $G^{+-}$ and $G^{-+}$ a second one.
They coincide at $(0,0)$, but are otherwise different.  In this case
we have in fact four different possibilities to orient. 
We also learn from this example, that given an oriented $\Theta$ and two local branching structures at $x\in \supp(\Theta)$ for which the first is oriented and the second is orientable,
the second cannot to be re-oriented to coincide with the orientation of the first one.
\qed
\end{example}

Next we would like to consider the boundary $\partial \Theta$ of an oriented $\wh{\Theta}$ and define an induced orientation.
Given $x\in \supp(\Theta)$
we find an open neighborhood $U(x)$ on which we have a local branching structure ${(M_i)}_{i\in I}$, ${(\sigma_i)}_{i\in I}$
such that
$$
\Theta(y)=\sum_{\{i\in I\ |\ y\in M_i\}} \sigma_i.
$$
Assuming that $\Theta$ is oriented we can pick an oriented branching structure such that
$$
\wh{\mathsf{T}}_\Theta(y) =\sum_{\{i\in I \ |\ y\in M_i\}} \sigma_i\cdot T_y(M_i,o_i).
$$
Assuming that $\Theta:X\rightarrow {\mathbb Q}^+$ is tame, the $(M_i,o_i)$ have the property that $M_i$ is tame, i.e.
in particular has the structure of a classical smooth manifold with boundary with corners.
In this case, for every $i\in I$, the set of points with $d_{M_i}(y)=1$ is open and dense in $\partial M_i$. 
We shall define an induced  boundary orientation for $\partial\Theta$ in the following.
First we state the usual orientation convention.
\begin{definition}\index{D- Orientation for $\partial M$}
Let $(M,o)$ be an oriented  smooth $n$-dimensional M$^+$-polyfold 
with $\partial M\neq \emptyset$. At a point $y\in \partial M$ with $d_M(y)=1$
an orientation $o_{\partial M,y}$ of $T_y\partial M$  is defined by a basis $(a_1,...,a_{n-1})$ of $T_y\partial M$
such that $(a_1,..,a_{n-1},e)$, where $e\in T_yM$ is inward pointing, defines $o_y$ on $T_yM$.
The map $o_{\partial M}$ which assigns to a point $y\in\partial M$ with $d_M(y)=1$ the orientation $o_{\partial M,y}$ is
called the {\bf orientation} of $\partial M$.
\qed
\end{definition}

Let $X$ be an ep-groupoid and $\wh{\Theta}:X\rightarrow {\mathbb Q}$ be a branched ep$^+$-subgroupoid.
\begin{definition}\index{D- Tame branched ep$^+$-subgroupoid}\label{TAMERXX}
We say that $\Theta:X\rightarrow {\mathbb Q}^+$ is {\bf tame}, provided for every $x\in \supp(\Theta)$ 
there exists an open neighborhood $U(x)$ admitting the natural $G_x$-action and having the properness
property and allowing for a local branching structure ${(M_i)}_{i\in I}$, ${(\sigma_i)}_{i\in I}$ representing $\Theta$ on $U(x)$,
 where 
every M$^+$-polyfold $M_i$ has an equivalent tame structure.  In particular, every $M_i$ has an equivalent structure of a smooth manifold with boundary with corners.
\qed
\end{definition}

Given the oriented, tame,  branched ep$^+$-subgroupoid $\wh{\Theta}:X\rightarrow {\mathbb Q}^+$ 
we find for every $x\in \supp(\Theta)$ an open neighborhood $U(x)$ and an associated oriented tame  branching structure
${(M_i,o_i)}_{i\in I}$, ${(\sigma_i)}_{i\in I}$ such that for $y\in U(x)$
\begin{eqnarray*}
&\Theta(y) = \sum_{\{i\in I\ |\ y\in M_i\}} \sigma_i&\\
& \wh{\mathsf{T}}_\Theta(y) = \sum_{\{i\in I\ |\ x\in M_i\}} \sigma_i\cdot T_y(M_i,o_i).&
\end{eqnarray*}
For $\partial\Theta:X\rightarrow {\mathbb Q}^+$ we define an orientation by 
\begin{eqnarray}\label{asdfgh}
&\wh{\mathsf{T}}_{\partial\Theta}(y) = 0\ \ \ \text{for}\ \ y\in X_\infty\setminus \supp(\Theta)&\\
& \wh{\mathsf{T}}_{\partial\Theta}(y) = \sum_{\{i\in I\ |\ y\in \partial M_i,\ d_{M_i}(y)=1\}} \sigma_i\cdot T_y(\partial M_i,o_{\partial M_i})\ \
\text{for}\ \ y\in U(x)\cap \partial X_\infty.&\nonumber
\end{eqnarray}
\begin{definition}[Boundary Orientation]\label{DEF936}\index{D- Orientation for $\partial\Theta$}
For the oriented, tame,  branched ep$^+$-subgroupoid $\wh{\Theta}:X\rightarrow {\mathbb Q}^+$ 
the map $\wh{\mathsf{T}}_{\partial\Theta}:X_\infty\rightarrow \wh{\text{Gr}}(X)$  defined in (\ref{asdfgh}) 
is called the orientation for $\partial\Theta$ induced from the orientation $\wh{\mathsf{T}}_\Theta$ of $\Theta$.
\qed
\end{definition}
In order to simplify notation we shall abbreviate 
$$
\wh{\Theta}:=(\Theta,\wh{\mathsf{T}}_\Theta)\ \ \text{and}\ \ 
\partial\wh{\Theta}:=(\partial\Theta,\wh{\mathsf{T}}_{\partial\Theta}).
$$
As we shall see later on, the following holds for a compact, oriented, tame, branched $n$-dimensional ep$^+$-subgroupoid  $\wh{\Theta}$ and
a degree $n-1$ sc-differential form $\omega$ for an appropriately defined {\bf branched integration} $\oint$:
$$
\oint_{\wh{\Theta}}d\omega = \oint_{\partial\wh{\Theta}}\omega.
$$
This is the  version of Stokes theorem for branched integration and details are carried out in Section \ref{IandST}.

\begin{remark}\index{R- Questions on non-tame branched ep$^+$-groupoids}
It would be interesting to carry out the orientation discussion for  branched ep$^+$-groupoids which are not necessarily
tame in order to see the generality under which Stokes theorem holds. 
The basic step would be to understand the boundary behavior of a M$^+$-polyfold. 
We refer the reader to Appendix \ref{SEC96} for a discussion of the issues involved.
\qed
\end{remark}

\section{Geometry of Local Branching Structures}\label{SEC15.4}
For the integration theory we need to understand the geometry of a local branching structure in more detail.
Of particular importance is the understanding of the branching sets. For this we need later on the following definitions.
\begin{definition}[\cite{HWZ7}]
Let $M$ be a smooth $n$-dimensional manifold with boundary and corners
and $K\subset M$. 
\begin{itemize}
\item[(1)]\ Associated to $z\in K$ there is a distinguished linear subspace 
$T^K_zM\subset T_zM$ defined by the following property.
Given a chart $\varphi:(V(z),z)\rightarrow ([0,\infty)^k\times {\mathbb R}^{n-k},0)$, where 
$k=d_M(z)$, the linear subspace $T\varphi(z)(T^K_zM)$\index{$T\varphi(z)(T^K_zM)$} of ${\mathbb R}^n$ is the linear hull
of all unit vectors in ${\mathbb R}^n$, which can be obtained by taking a sequence $(z_k)\subset K\setminus\{z\}$
with 
\begin{itemize}
\item[(a)]\ $\lim_{k\rightarrow \infty} z_k=z$.
\item[(b)]\ $\lim_{k\rightarrow\infty} \varphi(z_k)/|\varphi(z_k)| =e$.
\end{itemize}
The definition is independent of the choice of $\varphi$. 
\item[(2)]\
We say a point $z\in K$ is {\bf  essential}\index{D- Essential points} provided there exists a sequence $(z_k)\subset K$
such that
\begin{itemize}
\item[(a)]\ $\lim_{k\rightarrow \infty} z_k=z$.
\item[(b)]\ $T_{z_k}^K M= T_{z_k} M$.
\end{itemize}
If $K^e$ denotes the essential points in $K$ then the points in $K^{ne}:=K\setminus K^e$ are called 
{\bf non-essential}\index{D- Non-essential points} points.
\end{itemize}
Note that the statement that $z\in K$ is  essential is a statement about $z$ with respect to $K$ as a subset of $M$.
We need the knowledge of $M$ near $z$ to define the allowable $\varphi$.
\qed
\end{definition}

The following result is obvious from the definition.
\begin{lemma}
Let $M$ be a finite-dimensional manifold with boundary with corners and $K$ a subset.
Then $K^e$ is closed in $K$ and consequently $K\setminus K^e $ is open in $K$.
\qed
\end{lemma}
For the following we fix a local (oriented) branching structure ${(M_i,o_i)}_{i\in I}$, ${(\sigma_i)}_{i\in I}$, on an open neighborhood
$U(x)$. We impose the technical condition 
\begin{description}
\item[(T)] If for $i,i'\in I$ and some $y_0\in U(x)$  it holds $y_0\in M_i\cap M_{i'}$
and $T_{y_0}M_i=T_{y_0} M_{i'}$ and $o_{i,y_0}=\varepsilon \cdot o_{i',y_0}$ for some $\varepsilon \in \{-1,1\}$,
 then $T_y (M_i,o_i)=T_y(M_{i'},\varepsilon\cdot o_{i'})$ for all $y\in M_i\cap M_{i'}$ with $T_yM_i=T_yM_{i'}$.
\end{description}
Property {\em{\bf (T)}}  can  be achieved by restricting a given oriented branching structure on $U(x)$ to a smaller open neighborhood
$U'(x)$ as is shown in the following lemma.
\begin{lemma}
Let ${(M_i,o_i)}_{i\in I}$, ${(\sigma_i)}_{i\in I}$ be an oriented branching structure on an open neighborhood
$U(x)$ around the smooth point $x$. Then there exists $U'(x)\subset U(x)$ so that the restricted branching structure satisfies
the technical condition  {\em{\bf (T)}}.
\end{lemma}
\begin{proof}
The consideration is local. Without loss of generality we may assume that we are given  an open (rel $C$) neighborhood
$U$ of $0\in C\subset E$,  subsets $M_i$ which are sc$^+$-retracts, oriented and contain $0$.
We assume that the $M_i$ are properly embedded into $U$.  (There also exists an sc-smooth retract $O$ containing 
$0$ so that the $M_i$ lie in $O$, but this is not relevant for the argument).

Each of the $M_i$ is near $0$ the image of an sc$^+$-smooth retraction $r_i$. We consider
the maps 
$$
z\rightarrow Tr_i(r_i(z))E = T_{r_i(z)}M_i.
$$
for $z\in C$ close to $0$.  We can identify the smooth finite-dimensional spaces $T_{r_i(z)}M_i$ naturally with linear subspaces
of $E$.  Assume that our assertion is wrong.
Then we find  sequences $(y_k)$, $(z_k)$ of smooth points converging to zero,  indices $i_k^1\neq i_k^2$
and $\varepsilon_k\in \{-1,1\}$ such that
\begin{eqnarray}\label{PLKJ}
y_k\in M_{i_k^1}\cap M_{i^2_k},\ \ T_{y_k}M_{i_k^1}=T_{y_k}M_{i_k^2},\ \ o_{i_k^1,y_k}=\varepsilon_k\cdot o_{i_k^2,y_k}
\end{eqnarray}
and
\begin{eqnarray}\label{PLKJ2}
z_k\in M_{i_k^1}\cap M_{i^2_k},\ \ T_{z_k}M_{i_k^1}=T_{z_k}M_{i_k^2},\ \ o_{i_k^1,z_k}=-\varepsilon_k\cdot o_{i_k^2,z_k}.
\end{eqnarray}
After perhaps taking suitable subsequences we may assume that $i_k^1=:i_1$ and $i_k^2=:i_2$ are constant and $i_1\neq i_2$, and moreover
$\varepsilon_k$ is constant.  The points $y_k$ and $z_k$ are smooth and the tangent spaces $T_{y_k}M_{i_1}$ converge to $T_0M_{i_1}$
and similarly $T_{z_k}M_{i_2}\rightarrow T_0M_{i_2}$.   We also note that $T_0M_{i_1}=T_0M_{i_2}$. 
From (\ref{PLKJ}) we conclude  $o_{i_1,0}=\varepsilon\cdot o_{i_2,0}$ and from (\ref{PLKJ2}) that 
$o_{i_1,0}=-\varepsilon\cdot o_{i_2,0}$, which gives a contradiction.
\qed \end{proof}
We associate to   $i,i'\in I$ the closed subset of $M_i$ and $M_{i'}$  defined by
$$
K(i,i')=M_i\cap M_{i'}.
$$
We can view $K(i,i')$ as a subset of $M_i$ as well as a subset of $M_{i'}$.
Hence, a priori, we can talk about points $z\in K(i,i')$ which are essential with respect to
$K(i,i')$ as a subset of $M_i$ or as a subset of $M_{i'}$.
The detailed proof of the following statement, which shows that these two a priori notions of being essential
coincide, is left to the reader.
\begin{lemma}
If $z\in K(i,i')$ is essential for $K(i,i')\subset M_i$ then it is also essential for
$K(i,i')\subset M_{i'}$. Hence the subsets $K(i,i')^e$ and $K(i,i')^{ne}$ are well-defined
independent of the chosen reference $M_i$ or $M_{i'}$.
\end{lemma}
\begin{proof}
Given $z\in K\subset M\subset X$ the first observation is the following.
Take a M-polyfold chart $\psi:(U(z),z)\rightarrow (O,0)$ and note that $T_z^KM$ corresponds under $T\psi(z)$
to the linear subspace of $T_0O$ defined as follows.  Let $N=\psi(M\cap U(z))\subset O$ and $D=\psi(K\cap U(z))$.
Then define $T^D_0 N$ to be the span of all vectors in $T_0O$ such that there exists a sequence $(z_k)\subset D\setminus\{0\}$ with
$$
e=\text{lim}_{k\rightarrow\infty} \frac{z_k}{|{z_k}|_0}.
$$
As already said $T\psi(z)(T^K_zM)=T^D_0 N$.   From this description the proof of the lemma can be reduced
to the study of $0\in M_i,M_{i'}\subset O$ and $D=M_i\cap M_{i'}$ and the alternative description 
immediately implies that $T^K_zM_i =T^K_z M_{i'}$.
\qed \end{proof} 
In view of the lemma the subsets $K(i,i')^e$ and $K(i,i')^{ne}$ of $M_i$ and $M_{i'}$ are well-defined
independent of which reference we take.
The manifolds  $M_i$ and $M_{i'}$ come with orientations $o_i$  and $o_{i'}$. If $z\in K(i,i')^e$
and $T_{z}M_i$ and $T_{z}M_{i'}$ have the same orientation it holds in view of {\em{\bf (T)}} for all points
in $K(i,i')^{e}$, and similarly if the orientations differ.
In view of this observation we can associate to $(i,i')$ a map  $K(i,i')^{e}\rightarrow \{-1,1\}$.
If $K(i,i')^e=\emptyset$ it is the obvious map, and if the set is nonempty it is the constant map relating the orientations.

\begin{definition}\label{DEF944}\index{D- Measure zero}
We shall say that a set $K\subset M_i$ has {\bf measure zero} provided for every smooth chart 
$\varphi:M_i\supset U\rightarrow [0,\infty)^k\times {\mathbb R}^{n-k}$ the set $\varphi(U\cap K)$ has
Lebesgue-measure zero.
\qed
\end{definition}
  A basic fact is given by the following lemma,  which is proved similarly as  Lemma 2.8 in  \cite{HWZ7}.
\begin{lemma}
The sets $K(i,i')^{ne}$ have measure zero.
\end{lemma}
\begin{proof}
We
pick a point $x\in K (i,j)^{ne}$ and choose  a chart $\varphi:V(x)\subset
M_i\rightarrow O(0)\subset [0,\infty)^d\times {\mathbb R}^{n-d}$ satisfying $\varphi (x)=0$. We abbreviate the image set in $\R^n$ by $K=\varphi (K(i,j)^{ne}\cap V(x))$. It suffices to show that $0$ is a point of Lebesgue density $0$, i.e.,
$$
\lim_{\varepsilon\rightarrow 0}
\frac{1}{\varepsilon^n}\mu(B_\varepsilon\cap K)=0,
$$
where $\mu$ stands for the $n$-dimensional   Lebesgue measure  and  $B_{\varepsilon}$ is the  ball of radius $\varepsilon$ centered at the origin. We must have $T_x^{K(i, j)}M_i\neq T_xM_i$. Hence by
composing the chart with a rotation in the image,  we may assume
that
$$
\Sigma:=T\varphi(x)(T_x^{K(i,j)}M_i)\subset\R^{n-1}\times \{0\}.
$$
Now  take a number $\delta>0$, define the set
$\Gamma_\delta=\{(a, b)\in \R^{n-1}\times \R\vert \, \abs{b}\leq \delta \abs{a}\}$ and consider the subset
$(B_{\varepsilon}\setminus \Gamma_{\delta})\cap K$ of  $\R^n$. If, for a given  sequence $\varepsilon_k\to 0$,  there exists a sequence of points $x_k$ in
$(B_{\varepsilon_k}\setminus \Gamma_{\delta})\cap K$, then  we arrive at a contradiction with the definition of $\Sigma$. Consequently,   the set
$B_{\varepsilon}\cap K$  is contained in $\Gamma_{\delta}$ if $\varepsilon$ is sufficiently small. Hence
\begin{equation*}
\begin{split}
\limsup_{\varepsilon\to  0}\frac{1}{\varepsilon^n}\mu(B_\varepsilon\cap K)\leq
\limsup_{\varepsilon\rightarrow
0}\frac{1}{\varepsilon^n}\mu(B_\varepsilon\cap \Gamma_\delta)\leq
C(\delta)
\end{split}
\end{equation*}
for a constant $C(\delta)$ satisfying $C(\delta)\rightarrow 0$ as $\delta\rightarrow 0$. This shows that
$$
\lim_{\varepsilon\rightarrow 0}\frac{1}{\varepsilon^n}\mu(B_\varepsilon\cap K) =0.
$$
We have proved that the Lebesgue density at every point $x$ in the Borel set $K (i, j)^{ne}\subset M_i$ vanishes. This implies that this  set is of measure zero.
\qed \end{proof}
In view of the previous lemma the sets $K(i,i')^e$ have full measure in $K(i,i')$ and $K(i,i')^e\subset M_i\cap M_{i'}$
Denote by $M_U$ the union of the ${(M_i)}_{i\in I}$. Define for $y\in M_U$ the subset $I_y$ of $I$ by
$$
I_y=\{i\in I\ |\ y\in M_i\}.
$$
\begin{definition}\index{D- $y$-equivalence}
We say that $i,i'\in I_y$ are {\bf $y$-equivalent} ($i\sim_y i'$) provided there exists a chain of indices in $I_y$,
say $i_0=i,i_1,...,i_k=i'$ such that $y\in K(i_{p-1},i_p)^e   $ for $p=1,..,k$.  Note that if $i\sim_y i'$ then $T_y M_i=T_yM_{i'}$.
The equivalence relation $\sim_y$ defines a {\bf partition of $I_y$}\index{D- Partition $P_y$} given by the collection of equivalence classes
which we shall denote by $P_y:= I_y{/\sim_y}$. 
\qed
\end{definition}

\begin{lemma}
The set $\{y\in M_U\ |\ \sharp(I_y/\sim_y)\geq 2\}$ has measure zero.
\end{lemma}
\begin{proof}
Given $y\in M_U$ the equivalence relation $\sim_y$ defines a partition $P_y$ of the set $I_y$.
Define for a subset $J$ of $I$ and a finite partition $P$ of $J$ the set
$$
M_U^{(J,P)} =\{y\in M_U\ |\ I_y=J,\ P_y=P\}.
$$
The set $\{y\in M_U\ |\ \sharp(I_y/\sim_y)\geq 2\}$ can be written as a finite union
of such sets. Hence it suffices to show that $M_U^{(J,P)}$ has measure zero provided $\sharp P\geq 2$.

Assume that $\sharp P\geq 2$ and take two indices $i,i'$ which belong to different sets $\gamma_1$ and $\gamma_2$  of the partition $P$. If  $y\in M_U^{J,P}$ it follows that $i\not \sim_y i'$. That means that for every chain of indices
$i_0=i,...,i_k=i'$ in $I_y= J$ connecting $i_0$ with $i'$ it holds that $y$ belongs to some $K(i_{p-1},i_p)^{ne}$ with $p=p_y$.
Hence 
$$
y\in \bigcup_{j,j'\in J} K(j,j')^{ne}\subset \bigcup_{j,j'\in I} K(j,j')^{ne},
$$
where the right-hand side is a set of measure zero.
\qed \end{proof}
In view of the previous discussion we can take the subset $M_U^\ast$ of $M_U$ consisting of all
$y$ with $I_y{/\sim_y} \ = \{I_y\}$, i.e. all points $y$ so that all the elements in $I_y$ are $\sim_y$-equivalent.
By construction $M_U\setminus M_U^\ast$ has measure zero.
Denote for a subset $J\subset I$ by $M^{J}_U$ the subset of $M_U^\ast$ consisting of all
$y$ such that $I_y=J$ and $P_y=\{J\}$, i.e.
$$
M^J_U=\{y\in M_U^\ast\ |\ I_y =J,\ P_y=\{J\}\}.
$$
For $y\in M^{J}_U$ and every $i\in J$ we have the same
tangent space $T_yM_i$, but the orientations coming from the $o_i$ can differ.

This defines on $J$ an equivalence relation $\approx_J$, where by definition $i\approx_J i'$ provided
$T_y(M_i,o_i)=T_y(M_{i'},o_{i'})$  for all $y\in M^J_U$ with $i\sim_y i'$.
Recall that as a consequence of the technical condition {\em{\bf (T)}}  it holds that if   $i\sim_{y_0} i'$ then it is true for all $y$ with $i\sim_y i'$ in $M^J_U$. 

We obtain a partition of $J$ into possibly two nonempty subsets denoted by $A_J$ and $B_J$.
We define  for $y\in M^{J}_U$ 
$$
T^J_y:= T_yM_i,
$$
where we pick any element $i\in J$,
and 
$$
\sigma_J:=\sum_{i\in J} \sigma_i.
$$
With these definitions 
$$
\mathsf{T}_\Theta(y) = \sigma_J\cdot T_y^J\ \ \text{for}\ \ y\in M_U^{J}.
$$
 Associated to $A_J$, if nonempty,  we have the orientation 
$o_i$, where we pick some $i\in A_J$.  We  equip  $T_yM_i $ with $o_{A_J,y}=o_{i,y}$.    If  $B_J\neq \emptyset$ we define for $i'\in B_J$
the orientation $o_{B_J,y}=o_{i',y}$, which over the relevant set is just $-o_{A_J}$.
We define
$$
a_J:=\sum_{i\in A_J} \sigma_i\ \ \text{and}\ \ \ b_J:=\sum_{i'\in B_J} \sigma_{i'}
$$
and note that $\sigma_J=a_J+b_J$.
 Now consider the assignment
$$
M^J_U\ni y\rightarrow a_J\cdot (T^J_y,o_{a_J}) +b_J\cdot (T^J_y,o_{b_J})
$$
\begin{proposition}\label{PROP947}
Under the assumptions of the previous discussion the following holds.
The  collection
$$
M_U^\ast=\{ y \in M_U\ |\ \sharp P_y =1\}
$$
differs from $M_U$ by a set  of measure zero.  We can decompose $M_U^\ast$
by
$$
M_U^\ast = \sqcup_{J\subset I} M_U^{J}:=\sqcup_{J\subset I} \{y\in M_U^\ast\ |\ I_y=J\}.
$$
On $M_U^J$ we have the following identities
\begin{itemize}
\item[{\em(1)}]\ $\Theta(y)=\sum_{i\in J} \sigma_i=:\sigma_J$, $y\in M_U^J$.
\item[{\em(2)}]\ $\mathsf{T}_\Theta(y) =\sigma_J\cdot T_y^J$, $y\in M^J_U$.
\item[{\em(3)}]\  $\wh{\mathsf{T}}_\Theta(y) = a_J\cdot (T^J_y,o_{a_J}) +b_J\cdot (T^J_y,o_{b_J})$, $ y\in M_U^{J}$.
\item[{\em(4)}]\ $a_J+b_J = \sigma_J$ and $o_{a_J}=-o_{b_J}$.
\end{itemize}
\qed
\end{proposition}
The previous discussion allows us to compare two local branching structures on the open set
$U(x)$. Denote the second one by  ${(M_{i'},o_{i'})}_{i'\in I'}$,
${(\sigma_{i'})}_{i'\in I'}$. Then the union with the first one gives a new local branching structure indexed by $I''=I\sqcup I'$.
We abbreviate the three local branching structures by ${\mathcal B}$, ${\mathcal B}'$, and ${\mathcal B}''$. 
The local branching structure ${\mathcal B}''$ defines over  $U(x)$ the functor $2\cdot\Theta$ with orientation 
$$
\wh{\mathsf{T}}_{2\cdot\Theta}(y)=\sum_{\{i''\in I\sqcup I'\ |\ y\in M_{i''}\}} \sigma_{i''}\cdot T_y(M_{i''},o_{i''})
= 2\cdot \wh{\mathsf{T}}_\Theta(y).
$$
 The sets $M_U$, $M_U'$, and $M_U''$   associated to the three
local branching structures are all the same. However $M_U^\ast$, $M_U'^\ast$, and $M_U''^\ast$ generally are different, but 
their complements in $M_U$  have measure zero.
The set $M_U''^\ast$ is written as the disjoint union of sets 
$$
\{y\in M_U\ |\ I''_y = J'',\ \sharp P_y'' =1\}.
$$
where $J''\subset I''$. The intersection $M^\dagger_U= M_U^\ast\cap M_U'^\ast\cap M_U''^\ast$ has a complement of measure
zero in $M_U$. On $M^\dagger_U$ we can relate the  three local  branching structures in a way which is useful to define branched 
integration in the next subsection.  Define
\begin{eqnarray*}
&M^\dagger_U =\bigsqcup_{J\subset I} M^{\dagger J}_U :=\bigsqcup_{J\subset I} \{y\in M_U^\dagger\ |\ I_y = J\}&\\
&M^\dagger_U =\bigsqcup_{J'\subset I'} M^{\dagger J'}_U :=\bigsqcup_{J\subset I} \{y\in M_U^\dagger\ |\ I_y' = J'\}&\\
&M^\dagger_U =\bigsqcup_{J''\subset I''} M^{\dagger J''}_U :=\bigsqcup_{J''\subset I''} \{y\in M_U^\dagger\ |\ I_y'' = J''\}.&
\end{eqnarray*}
Pick $J\subset I$ and $J'\subset I'$ and consider the set $M^{\dagger J}_U \cap M^{\dagger J'}_U$. If $y$ belongs to this set
it holds that $I_y=J$ and $I_y'=J'$, with both being nonempty, since both branching structures define the same $\Theta$.
Since $y\in M^\dagger_U\subset M_U''^\ast$ it also satisfies $I''_y =J''$ and $\sharp P_y''=1$
for some $J''\subset I''$.   Clearly $J''=J\sqcup J'$ and we derive that 
\begin{eqnarray}\label{eq914}
M^{\dagger J}_U \cap M^{\dagger J'}_U \subset M^{\dagger J\sqcup J'}_U.
\end{eqnarray}
Since 
\begin{eqnarray}\label{eq915}
M^\dagger_U &=& \left( \bigsqcup_{J\subset I} M^{\dagger J}_U\right)\bigcap \left(\bigsqcup_{J'\subset I'} M^{\dagger J'}_U\right)\\
&=& \bigsqcup_{J\subset I,\ J'\subset J'} \left(M^{\dagger J}_U \cap M^{\dagger J'}_U\right),\nonumber
\end{eqnarray}
combining (\ref{eq914}) and (\ref{eq915}) we find the equality
\begin{eqnarray*}
M^{\dagger J}_U \cap M^{\dagger J'}_U = M^{\dagger J\sqcup J'}_U.
\end{eqnarray*}
From this discussion, since ${\mathcal B}$ and ${\mathcal B}'$ describe the same oriented local branching structure
and ${\mathcal B}''$ is their union  we infer that for $y\in  M^{\dagger J\sqcup J'}_U$ the following holds:
\begin{itemize}
\item[(1)]\ $\Theta(y) =\sum_{i\in J}\sigma_i =:\sigma_J = \sum_{i'\in J'} \sigma_{i'} =:\sigma_{J'}$
\item[(2)]\ $\mathsf{T}_{\Theta}(y) =\sigma_J\cdot T^J_y = \sigma_{J'}\cdot T^{J'}_y$
\item[(3)]\ $\wh{\mathsf{T}}_\Theta(y) = a_J\cdot (T^J_y,o_{a_J})+ b_J\cdot (T^J_y,o_{b_J}) = a_{J'}\cdot (T^{J'}_y,o_{a_{J'}})+ b_{J'}\cdot (T^{J'}_y,o_{b_{J'}})$.
\item[(4)]\ $a_J+b_J=\sigma_J=\sigma_{J'} =a_{J'}+b_{J'}$.
\end{itemize}
We can summarize this discussion in the following proposition.
\begin{proposition}\label{prop948}
Assume that $U(x)$ is an open neighborhood around the smooth point $x$ in a M-polyfold.  Suppose we are given two oriented local branching structures
on $U(x)$ denoted by ${(M_i,o_i)}_{i\in I}$, ${(\sigma_i)}_{i\in I}$, and ${(M_{i'},o_{i'})}_{i'\in I'}$, ${(\sigma_{i'})}_{i'\in I'}$ such that for $y\in U(x)$
the following holds
\begin{itemize}
\item[{\em(a)}]\ $\sum_{\{i\in I\ |\ y\in M_i\}} \sigma_i =  \sum_{\{i'\in I'\ |\ y\in M_{i'}\}} \sigma_{i'}$.
\item[{\em(b)}]\ $\sum_{\{i\in I\ |\ y\in M_i\}} \sigma_i\cdot T_y(M_i,o_i) =  \sum_{\{i'\in I'\ |\ y \in M_{i'}\}} \sigma_{i'}\cdot T_y(M_{i'},o_{i'})$.
\end{itemize}
Define $M_U:= \bigcup_{i\in I} M_i=\bigcup_{i'\in I'} M_{i'}$.  
Then there exists a subset $M^\dagger_U$ of $M_U$ and a finite partition of $M^\dagger_U=A_1\sqcup ...\sqcup A_\ell$ with the following properties.
\begin{itemize}
\item[{\em(1)}]\ $M_U\setminus M_U^\dagger$ has measure zero.
\item[{\em(2)}]\ For every $A_p$ the maps $y\rightarrow I_y=\{i\in I\ |\ y\in M_i\}$ and $y\rightarrow I_y'=\{i'\in I'\ |\ y\in M_{i'}\}$ are constant, and moreover
for $i,j\in I_y$ and $i',j'\in I_y'$
$$
T_yM_i =T_y M_j =T_y M_{i'}=T_y M_{j'}
$$
\item[{\em(3)}]\ Denote by $I_p$ the value of $y\rightarrow I_y$ on $A_p$ and similarly $I_p'$. Then 
$$
\sum_{i\in I_p} \sigma_i=\sum_{i'\in I_p'} \sigma_{i'}.
$$
\item[{\em(4)}]\ For $A_p$ with $I_p= J$ and $I'_p=J'$ there exist partitions $J=J_+\sqcup J_-$ and $J'=J'_+\sqcup J'_-$ such that
for $(i,j), (i',j')\in (J_+\times J_+')$ or $(i,j), (i',j')\in (J_-\times J_-')$ and all $y\in A_p$
$$
T_y (M_i,o_i)=T_y (M_j,o_j)= T_y(M_{i'},o_{i'})=T_y (M_{j'},o_{j'}).
$$
\item[{\em(5)}]\ For $y\in A_p$  and having fixed $i_\pm\in J_\pm$, $i'_\pm\in J'_\pm$ we have the identity
$$
[\sum_{j\in J_\pm}\sigma_j]\cdot T_y(M_{i_\pm},o_{i_\pm})= [\sum_{j'\in J_\pm'} \sigma_{j'} ]\cdot T_y(M_{i'_\pm},o_{i'_\pm})
$$
\end{itemize}
\qed
\end{proposition}
Proposition \ref{prop948} will be important for introducing the branched integration in the next 
section.
\section{Integration and Stokes}\label{IandST}
In the following we assume that $\Theta$ is an 
 oriented, tame, compact,  branched ep$^+$-subgroupoid of dimension $n$, see
 Definition \ref{TAMERXX}.
We denote by $S$ the orbit space $|\supp(\Theta)|$ which by assumption is a compact
subspace of the metrizable space $|X|$. 

Given a classically smooth $n$-manifold $M$ (perhaps with boundary with corners)   there is a natural class of measures equivalent to Lebesgue measure, i.e.
those which can be written in local coordinates as $f\cdot dx_1\wedge..\wedge dx_n$ for some function $f>0$.
This allows to introduce the notion of a subset $A\subset M$ to be of measure zero. We have introduced this notion already in
Definition \ref{DEF944}. Given a local branching structure ${(M_i)}_{i\in I}$, ${(\sigma_i)}_{i\in I}$ for $\Theta$ on
$U(x)$ we obtain 
$$
M_U=\bigcup_{i\in I}  M_i
$$
and $|M_U|\subset S$.  Given a subset $A$ of $S$  we can consider $A\cap |M_U|$ and the preimages
$A_i\subset M_i$ for $i\in I$ under the canonical map $M_U\rightarrow |M_U|$.
\begin{definition}\index{D- Measure zero}
We say a subset $A$ of $S=\supp(\Theta)$ has {\bf measure zero}, provided for every $x\in \supp(\Theta)$
and a choice of local branching structure on a suitable $U(x)$ the sets $A_i$ introduced above, have measure zero.
\qed
\end{definition}
According to the following lemma, see \cite{HWZ7}, Lemma 3.2, the definition does not depend on the choice of local branching structures.
\begin{lemma}\index{L- Measure zero}
The definition of a set of measure zero in $S$ does not depend on the choice of local branching structures.
\qed
\end{lemma}
After this preparation we can introduce the {\bf canonical $\sigma$-algebra ${\mathcal L}(S)$}, see also \cite{HWZ7}, Definition 3.3.
\begin{definition}[Canonical $\sigma$-Algebra ${\mathcal L}(S)$] \index{D- Canonical $\sigma$-Algebra ${\mathcal L}(S)$}
Let $X$ be an ep-groupoid and $\Theta: X\rightarrow {\mathbb Q}^+$ a tame, compact, branched ep$^+$-subgroupoid. Define
$S=|\text{supp}(\Theta)|$. The $\sigma$-algebra ${\mathcal L}(S)$\index{${\mathcal L}(S)$} is the smallest $\sigma$-algebra consisting of subsets of $S$,
which contains the Borel $\sigma$-algebra ${\mathcal B}(S)$\index{${\mathcal B}(S)$} and all subsets of $S$ of measure zero.
\qed
\end{definition}
It follows immediately from the definition that ${\mathcal L}(S)$ is obtained from ${\mathcal B}(S)$ by just adding the sets
of measure zero. Denote by ${\mathcal M}(S,{\mathcal L}(S))$\index{${\mathcal M}(S,{\mathcal L}(S))$} the vector space of finite signed measures on the Lebesgue $\sigma$-algebra  on $S$.

The basic first integration result concerns an sc-differential form and an oriented, tame, branched ep$^+$-subgroupoid.
\begin{theorem}\label{IandS}\index{T- Integration associated to $\wh{\Theta}$}
Let $X$ be an ep-groupoid and $\wh{\Theta}=(\Theta,\wh{\mathsf{T}}_\Theta):X\rightarrow {\mathbb Q}$
be an oriented, tame, compact,  branched ep$^+$-subgroupoid of dimension $n$.
There exists a linear map 
$$
\Phi_{\wh{\Theta}}:\Omega^n_\infty(X)\rightarrow {\mathcal M}(S,{\mathcal L}(S)):\omega\rightarrow \mu_\omega,
$$
characterized uniquely by the following property.
\begin{itemize}
\item  Given a point $x\in \text{supp}(\Theta)$ and an oriented local branching structure ${(M_i,o_i)}_{i\in I}$, ${(\sigma_i)}_{i\in I}$
on $U(x)$, which is assumed to admit the natural $G_x$-action, the following identity holds for every
measurable subset $K\subset S$ contained in $|U(x)|$ with $K_i\subset M_i$ the preimage of $K$ under
$\pi:U(x)\rightarrow |U(x)|:y\rightarrow |y|$
$$
\mu^\Theta_\omega(K) =\frac{1}{|G_x^{{\text{eff}}}|} \cdot \sum_{i\in I} \sigma_i\cdot \int_{K_i} \omega|(M_i,o_i).
$$
\end{itemize}
Here $G_x^{eff}$ is the effective isotropy group of $x$ and was previously defined in Definition \ref{DEF728}.
\qed
\end{theorem}
The proof of Theorem \ref{IandS} takes some preparation and we follow essentially \cite{HWZ7}.
In a first step we show that locally the definition does not depend on the choice of branching structure, see \cite{HWZ7}, Lemma 3.5.
\begin{lemma}[Independence]\label{LOP1}\index{L- Well definedness of $\mu^U_\omega$}
Assume that $x\in \supp(\Theta)$ and $U=U(x)$ is an open neighborhood invariant under the natural $G_x$-action
and having the properness property. Suppose further that $U$ harbors two oriented local branching structures,
say ${(M_i,o_i)}_{i\in I}$, ${(\sigma_i)}_{i\in I}$, and ${(M_{i'},o_{i'})}_{i'\in I'}$, ${(\sigma_{i'})}_{i'\in I'}$ describing the oriented $\wh{\Theta}$ over $U$.
Then for a measurable subset $K$ of $M_U$ with corresponding $K_i\subset M_i$ for $i\in I$ and $K_{i'}\subset M_{i'}$ for $i'\in I'$ it holds
$$
\frac{1}{|G_x^{{\text{eff}}}|} \cdot \sum_{i\in I} \sigma_i\cdot \int_{K_i} \omega|(M_i,o_i)
= \frac{1}{|G_x^{{\text{eff}}}|} \cdot \sum_{i'\in I'} \sigma_{i'}\cdot \int_{K_{i'}} \omega|(M_{i'},o_{i'}).
$$
Here the index sets are disjoint, i.e. $I\cap I'=\emptyset$. 
\end{lemma}
\begin{proof}
In order to prove the lemma we assume an oriented, compact, and tame  $\wh{\Theta}$ is given. Let $x\in \supp(\Theta)$ and $U(x)$ be
an open neighborhood equipped with the natural $G_x$-action, and the properness property.  On $U(x)$ we assume that
two local oriented branching structures are given, say
$$
{(M_i,o_i)}_{i\in I},\ {(\sigma_i)}_{i\in I}\ \ \text{and}\ \ {(M_{i'},o_{i'})}_{i'\in I'},\ \  {(\sigma_{i'})}_{i'\in I'},
$$
which represent $\Theta$ and $\wh{\mathsf{T}}_\Theta$, i.e. for $y\in U(x)$
\begin{eqnarray}
\sum_{\{i\in I\ |\ y\in M_i\}} \sigma_i &=&\sum_{\{i'\in I'\ |\ y\in M_{i'}\}}\sigma_{i'}\\
\sum_{\{i\in I\ |\ y\in M_i\}} \sigma_i \cdot T_y(M_i,o_i)&= &\sum_{\{i'\in I'\ |\ y\in M_{i'}\}} \sigma_{i'}\cdot T_y(M_{i'},o_{i'}).\nonumber
\end{eqnarray}
We assume that $\omega\in \Omega^n_\infty(Y)$ and denote by $\pi:U(x)\rightarrow |U(x)|$
the projection $y\rightarrow |y|$.   We show that for $K\subset |U(x)|$ with
$K_i= (\pi|M_i)^{-1}(K)$ and $K_{i'}:=(\pi|M_{i'})^{-1}(K)$  the equality
\begin{eqnarray}\label{XWEQ1}
\sum_{i\in I} \sigma_i\cdot \int_{K_i} \omega|(M_i ,o_i)=\sum_{i'\in I'} \sigma_{i'}\cdot \int_{K_{i'}}\omega|(M_{i'},o_{i'})
\end{eqnarray}
holds. 
If this is known to hold  the rest of the proof of Theorem \ref{IandS} follows easily. The difficulty in establishing the validity of 
(\ref{XWEQ1}) arises from the fact that in general there does not exist  a correspondence between the submanifolds
of the branching structures, so that one needs to analyze the branching sets in more detail.  
This is the place where the considerations of Section \ref{SEC15.4} enter, particularly Proposition \ref{prop948}.

Associated to the local branching structures we have with the notation in Section \ref{SEC15.4}
the sets $M_U=M_U'$ as well as $M_U^{\dagger}$, see Proposition \ref{prop948}. 
Denote by $A_1\sqcup..\sqcup A_p$  the decomposition of $M_U^\dagger$ according to Proposition \ref{prop948}.

If we subtract from the $K_i$ or $K_{i'}$ the set $M_U\setminus M_U^{\dagger}$, resulting in the sets $K_i^{\dagger}$ and $K_{i'}^{\dagger}$,
it holds, using that this is a measure zero modification
$$
\sum_{i\in I} \sigma_i\cdot\int_{K_i^{\dagger}}\omega|(M_i,o_i) =\sum_{i\in I} \sigma_i\cdot\int_{K_i}\omega|(M_i,o_i)
$$
and
$$
\sum_{i'\in I'} \sigma_{i'}\cdot\int_{K_{i'}^\dagger}\omega|(M_{i'},o_{i'}) =\sum_{i'\in I'} \sigma_{i'}\cdot\int_{K_{i'}}\omega|(M_{i'},o_{i'}).
$$
The set $K_i^\dagger $ decomposes
as 
$$
K^\dagger_i=\bigsqcup_{p=1}^\ell K^{\dagger p}_i,
$$
where $K^{\dagger p}_i:= K^\dagger_i\cap A_p$.
 We compute, using heavily Proposition \ref{prop948}:
\begin{eqnarray*}
&&\sum_{i\in I} \sigma_i\cdot\int_{K_i^{\dagger}}\omega|(M_i,o_i)= \sum_{i\in I} \sum_{p=1}^\ell\sigma_i\cdot\int_{K_i^{\dagger p}}\omega|(M_i,o_i)\\
&=&\sum_{p=1}^\ell \sum_{i\in I} \sigma_i\cdot\int_{K_i^{\dagger p}}\omega|(M_i,o_i)=\sum_{p=1}^\ell \sum_{i\in I_p}\sigma_i\cdot\int_{K_i^{\dagger p}}\omega|(M_i,o_i)\\
&=&\sum_{p=1}^\ell \left(\sum_{i\in I_p^+} \sigma_i-\sum_{i\in I_p^-}\sigma_i\right)\cdot \int_{K_{i_p^+}^{\dagger p}} \omega |(M_{i_p^+},o_{i_p^+})
\end{eqnarray*}
\begin{eqnarray*}
&=&\sum_{p=1}^\ell \left(\sum_{i'\in I'^+_p} \sigma_{i'}-\sum_{i'\in I'^-_p}\sigma_{i'}\right)\cdot \int_{K_{i'^+_p}^{\dagger p}} 
\omega |(M_{{i'}_p^+},o_{{i'}_p^+})\\
&=&\sum_{p=1}^\ell \sum_{i'\in I'_p}\sigma_{i'}\cdot\int_{K_{i'}^{\dagger p}}\omega|(M_{i'},o_{i'})=\sum_{p=1}^\ell \sum_{i'\in I'} \sigma_{i'}\cdot\int_{K_{i'}^{\dagger p}}\omega|(M_{i'},o_{i'})\\
&=& \sum_{i'\in I'}\sum_{p=1}^\ell \sigma_{i'}\cdot\int_{K_{i'}^{\dagger p}}\omega|(M_{i'},o_{i'})=\sum_{i'\in I'} \sigma_{i'}\cdot\int_{K_{i'}^{\dagger}}\omega|(M_{i'},o_{i'}).
\end{eqnarray*}
This shows that local integrals are well-defined and do not depend on the local branching structure.
\qed \end{proof}
The procedure  defines, associated
to an sc-differential form $\omega$ for a  local branching structure on $U(x)$,  which we assume to be invariant under the natural action of $G_x$
and having the properness property, a signed measure $\mu_\omega^{U(x)}$ for measurable subsets contained in $|U(x)|\subset S$, namely
$$
\mu^{U(x)}_\omega(K)=\frac{1}{|G^{\text{eff}}_x|} \cdot \sum_{i\in I}\sigma_i\cdot \int_{K_i} \omega|(M_i,o_i).
$$
From the construction it is clear that if $\tau=f\cdot \omega$ then 
\begin{eqnarray}\label{EQ918}
d\mu_\tau^{U(x)} = \bar{f}||U(x)|  \cdot d\mu_\omega^{U(x)},
\end{eqnarray}
with $\bar{f}\circ\pi =f$.
We shall use this fact later on.  
\begin{remark}\label{DFGX}
We point out that even if $U(x)$ and $U(x')$ are as before and $|U(x)|=|U(x')|$ with $|x|=|x'|$
it is still a possibility that $\mu^{U(x)}_\omega$ and $\mu^{U(x')}_\omega$ are not the same on $|U(x)|$.
This point is being addressed in Lemma \ref{LOP3}.
\qed
\end{remark}
  Since $S=|\supp(\Theta)|$ is compact the set $S$ can be 
covered by finitely many such constructions.  We want to show that this fact can be used to define a global signed measure, which boils down to showing that
the constructions coincide on the overlaps. We need Lemma 3.6 from \cite{HWZ7}. 

\begin{lemma}[Restrictions]\label{LOP2} \index{L- Restriction of $\mu^U_\omega$}
Let $X$ be an ep-groupoid and $\wh{\Theta}:X\rightarrow {\mathbb Q}^+$ an oriented, tame, branched ep$^+$-subgroupoid.   For $x\in \supp(\Theta)$
let $U=U(x)$ be an open subset invariant under the natural $G_x$-action, having the properness property,
and also supporting a local oriented  branching structure ${(M_i,o_i)}_{i\in I}$, ${(\sigma_i)}_{i\in I}$ for $\wh{\Theta}$.
Further let $y\in U$ with isotropy group $H_y$ and assume  $V=V(y)\subset U$ is an open neighborhood with the natural
$H_y$-action. For an sc-differential form $\omega$ denote the associated signed measures by $\mu^{U}_\omega$ and $\mu^{V}_\omega$.
Then, if a measurable subset $K\in {\mathcal L}(S)$ is contained in $|M_V|\subset |M_U|\subset S$ it holds 
$$
\mu^{U}_\omega(K)=\mu^{V}_\omega(K).
$$
\end{lemma}
\begin{proof}
The proof is given in \cite{HWZ7}.  Since the definition of the measures is independent of the actual 
local oriented branching structure we can take the one on $V$ which is induced from the one on $U$.
Then the result follows from a computation involving the orders of isotropy groups.
\qed \end{proof}

We address now the point raised in Remark \ref{DFGX},   see also \cite{HWZ7}, Lemma 3.7.
\begin{lemma}[Morphism Invariance]\label{LOP3}\index{L- Morphism invariance}
Let $X$ be an ep-groupoid and $\wh{\Theta}:X\rightarrow {\mathbb Q}^+$ a tame, compact, oriented branched ep$^+$-subgroupoid.
Assume that $x,x'\in \supp(\Theta)$ and $\phi:x\rightarrow x'$ is a morphism.  Let $U=U(x)$ and $U'=U(x')$
be open neighborhood allowing for the natural actions by $G_x$ and $G_{x'}$, respectively, and suppose the sets have the 
properness property. In addition we assume that the natural local sc-diffeomorpism associated to $\phi$, 
defines an sc-diffeomorphism $\wh{\phi}:U\rightarrow U'$. Given an sc-differential form $\omega$ on $X$ the equality
$$
\mu^U_\omega(K)=\mu^{U'}_\omega(K)
$$
holds for measurable subsets $K$ contained in $|M_U|=|M_{U'}|$. 
\end{lemma}
\begin{proof}
This is nothing else but change of variables.
\qed \end{proof}
Now we are in the position to prove Theorem \ref{IandS}.
\begin{proof}[Theorem \ref{IandS}]
Using the compactness of $S=|\supp(\Theta)|$ there exist finitely many points $x_1$,..,$x_\ell$
and open neighborhood $U_p:=U(x_p)$ admitting the natural $G_{x_p}$-action, having the properness property,
and having an oriented local branching structure ${(M^p_i,o^p_i)}_{i\in I^p}$, ${(\sigma_i^p)}_{i\in I^p}$ representing the oriented
$\wh{\Theta}$ over $U_p$, so that
$$
S=\bigcup_{p=1}^\ell |M_{U_p}|.
$$
 With the given sc-differential form $\omega$ we obtain the signed measures
$$
\mu^{U_p}_\omega
$$
on $|M_{U_p}|\subset S$.  Given a measurable subset $K$ of $S$ in ${\mathcal L}(S)$ we partition 
it in such a way that $K=\bigsqcup_{p=1}^\ell  K_p$ with $K_p\subset |M_{U_p}|$ and consider the sum
$$
\sum_{p=1}^\ell \mu^{U_p}_\omega(K_p).
$$
This definition does not depend on the choices involved, due to Lemmata \ref{LOP1}, \ref{LOP2}, and \ref{LOP3}.
It defines our desired signed measure $\mu_\omega^\Theta$. 
\qed \end{proof}
\begin{remark}\index{R- Questions about branched ep$^+$-subgroupoids}
(a) In \cite{HWZ7} the theorem was proved under the assumption that $X$ allows sc-smooth partitions of unity.
The modification of the proof presented above gets rid of this assumption.\par

(b) It should be possible to define $\mu_\omega$ also for compact, oriented, branched ep$^+$-subgroupoids, getting rid
of the tameness assumption. This would require to show that the boundaries of the occurring sub-M$^+$-polyfolds have measure
zero, which in all likelihood is the case. We leave such investigation to the reader.\par

(c) It might even be true that a compact sub-M$^+$-polyfold in a M-polyfold has a boundary which is a Lipschitz manifold,
which would have consequences for the the boundary structure of compact,  branched, ep$^+$-subgroupoids. 
It would in particular imply that the following considerations about Stokes theorem should be true for oriented 
$\Theta$ which are not tame. We also leave such investigations to the reader.
\qed
\end{remark}

Given an ep-groupoid $X$ and a compact, tame, oriented branched ep$^+$-subgroupoid $\wh{\Theta}:X\rightarrow {\mathbb Q}$
there is an oriented $\partial\wh{\Theta}:X\rightarrow {\mathbb Q}^+$ according to Definitions \ref{DEF9214} and \ref{DEF936}
and the associated discussions in previous sections.  Since $\Theta$ is tame the M$^+$-polyfolds occurring in the local branching 
structures have an equivalent structure as classically smooth manifolds $M$ with boundary with corners. Given such a M$^+$-polyfold
and a boundary point $x$ with degeneracy index $d=d_M(x)$, the boundary near $x$ consists of of the union of $d$-many local faces.
 
If $d\geq 2$ the intersection of two or more such faces,   as a subset of any such face, has $(n-1)$-dimensional Lebesgue measure $0$, where $n=\dim(M)$.
If $K\subset \partial M$ we have consequently a well-defined notion of being of {\bf $(n-1)$-dimensional measure zero}\index{Measure zero}. 
In view of this discussion we can define the Lebesgue $\sigma$-algebra ${\mathcal L}(\partial S)$, where $\partial S = |\supp(\partial\Theta)|$,
since there is a well-defined notion of a subset of $\partial S$ to be of $(n-1)$-dimensional measure zero.
\begin{definition}[Canonical $\sigma$-Algebra ${\mathcal L}(\partial S)$]\index{D- Canonical $\sigma$-Algebra ${\mathcal L}(\partial S)$}
Let $X$ be an ep-groupoid and $\Theta:X\rightarrow {\mathbb Q}^+$  a tame, compact, branched ep$^+$-subgroupoid.
The $\sigma$-algebra ${\mathcal L}(\partial S)$ \index{${\mathcal L}(\partial S)$} is the smallest $\sigma$-algebra containing the Borel $\sigma$-algebra ${\mathcal B}(\partial S)$\index{${\mathcal B}(\partial S)$}
and all subsets of $\partial S$ of measure zero.
\qed
\end{definition}
With the same techniques as used in the proof of Theorem \ref{IandS} we obtain the following theorem for boundary integration,
where ${\mathcal M}(\partial S,{\mathcal L}(\partial S))$ \index{${\mathcal M}(\partial S,{\mathcal L}(\partial S))$} is the vector spaces of signed measures.
\begin{theorem}[Canonical Boundary Measures]\label{THM9510}\index{T- Canonical boundary measures}
Assume $X$ is an ep-groupoid and $\wh{\Theta}:X\rightarrow {\mathbb Q}^+$ an oriented tame, compact, branched ep$^+$-subgroupoid of dimension $n$.
Then there exists a linear map 
$$
\Phi_{\partial\wh{\Theta}}:\Omega^{n-1}_\infty(X)\rightarrow {\mathcal M}(\partial S,{\mathcal L}(\partial S)):\omega\rightarrow \mu_\omega^{\partial\wh{\Theta}}
$$
which is uniquely characterized by the following property.
\begin{itemize}
\item Given $x\in \supp(\partial \Theta)$, an open neighborhood $U(x)$ in $X$ with the natural $G_x$-action and the properness property,
and an oriented branching structure on $U(x)$ representing $\wh{\Theta}$, say ${(M_i,o_i)}_{i\in I}$, ${(\sigma_i)}_{i\in I}$, the following identity
holds for a measurable subset $K\in {\mathcal L}(\partial S)$ with $K\subset |M_U^\partial| :=| \bigcup_{i\in I} \partial M_i|$
$$
\mu_{\omega}^{\partial\wh{\Theta}}(K) =\frac{1}{\sharp G_x^{\text{eff}}}\cdot \sum_{i\in I} \sigma_i\cdot \int_{K_i} \omega|(\partial M_i,o_{\partial M_i}),
$$
where $K_i\subset \partial M_i$ is the preimage of $\pi:\partial M_i\rightarrow |\supp(\partial\Theta)|\cap |U(x)|$.
\end{itemize}
\qed
\end{theorem}
\begin{remark}\index{R- Question concerning Stokes Theorem}
In our context we have a version of Stokes theorem, which can be formulated as follows.
\begin{theorem}[Stokes Theorem]\index{T- Stokes theorem}
Let $X$ be an ep-groupoid admitting sc-smooth partitions of unity, and $\wh{\Theta}:X\rightarrow {\mathbb Q}^+$ a tame, oriented, compact, branched ep$^+$-subgroupoid
of dimension $n$. Denote by $\partial\wh{\Theta}$ its oriented boundary and by $\omega$ an sc-differential form of degree $n-1$.
Then 
$$
\mu_\omega^{\partial\wh{\Theta}}(\partial S) =\mu_{d\omega}^{\wh{\Theta}}(S),
$$
or equivalently
$$
\oint_{\wh{\Theta}} d\omega   = \oint_{\partial\wh{\Theta}} \omega.
$$
\end{theorem}
The assumption of having an sc-smooth partition of unity is very likely not needed and enters since we use the standard ideas for proving 
Stokes theorem.
\end{remark}
\begin{proof}
We cover the compact $S$ by finitely many  open sets  $|U_k|$, where 
$U_k=U(x_k)$, $1\leq  k\leq \ell$ are open sets admitting the natural $G_{x_k}$-action, have the properness property
and are equipped with an oriented branching structure representing $\wh{\Theta}$.
We denote  the oriented  local  branching structures  by $(M_i,o_i)_{i\in I^k}$ (for $k=1,...,\ell$)  with the associated weights ${(\sigma_i)}_{i\in I^k}$. 
Then we have the total index set $I =I_1\sqcup ...\sqcup I_\ell$.

The $M_i$  are tame  finite dimensional submanifolds of $X$, all  of the same dimension $n$ and for $i\in I_k$
the $M_i$ are properly embedded in $U(x_k)$ and over $U(x_k)$ the functor $\Theta$ and associated orientation can be written as
\begin{eqnarray}
&\Theta(y)=\sum_{\{i\in I_k\ |\ y\in M_{i}\}} \sigma_i\ \ \text{for}\ \ y\in U(x_k)&\\
&\wh{\mathsf{T}}_{\Theta}(y) =\sum_{\{i\in M_i\ |\ y\in M_i\}} \sigma_i\cdot T_y(M_i,o_i)\ \ \text{for}\ \ y\in U(x_k).&\nonumber
\end{eqnarray}
Abbreviate $U_k:=U(x_k)$ for $k=1,...,\ell$ and 
 denote by $U_{k}^*$ the saturations of the sets $U_{k}$ in $X$ and add  another saturated open set $U^\ast_0$ so that
the sets $U^\ast_{k}$  cover $X$ and 
the sets $\abs{U^\ast_{k}\setminus \overline{U^\ast_0}}$  still cover $S$.
Given this finite open covering of $X$ by saturated set we can  take a subordinate sc-smooth partition of unity  $\beta_0,\ldots ,\beta_\ell$ on $X$ for the
ep-groupoid $X$.  Hence the $\beta_i$ are functors with values in $[0,1]$ and   $\supp \beta_k\subset U^{\ast}_k$.
By construction 
$$
\sum_{k=1}^\ell\beta_k =1\ \ \text{on}\ \  \bigcup_{k=1}^\ell M_{U_{k}}.
$$
We define   $\omega_k=\beta_k \omega$ so that its support is contained in $\supp \omega_k\subset U^*_k$ and
$\omega=\sum_{k=1}^n\omega_k$.

If  the set $|U_k|$ does not contain parts of $|\supp(\partial\Theta)|$,  we
conclude  by the standard Stokes theorem
$$
\int_{(M_i,o_i)} d\omega_{k} =0 \ \text{for}\ \ i\in I_k.
$$
In the other case, again by the standard Stokes theorem,
$$
\int_{(M_i,o_i)} d\omega_k=\int_{\partial (M_i,o_i)} \omega_{k}\ \ \text{for}\ \ i\in I_k.
$$
Summing  up  all these contributions we conclude the following.
\begin{eqnarray*}
\mu^{\wh{\Theta}}_{d\omega}(S)&=& \sum_{k=1}^\ell\mu^{\wh{\Theta}}_{d\omega_k}(S)\\
&=& \sum_{k=1}^\ell \frac{1}{\sharp G_{x_k}^{\text{eff}}}\cdot \sum_{i\in I_k}\sigma_i\cdot \int_{(M_i,o_i)}d\omega_k\\
&=& \sum_{k=1}^\ell \frac{1}{\sharp G_{x_k}^{\text{eff}}}\cdot \sum_{i\in I_k}\sigma_i\cdot \int_{\partial (M_i,o_i)}\omega_k\\
&=&\sum_{k=1}^\ell \mu_{\omega_k}^{\partial\wh{\Theta}}(\partial S)\\
&=& \mu^{\partial\wh{\Theta}}_\omega (\partial S).
\end{eqnarray*}
The proof of Stokes' theorem in the branched context is complete.
\qed \end{proof}

\section{Appendix}\label{SEC96}
In the following two subsections we shall discuss several questions, whose answers 
are related to possible generalizations of the results in this chapter.
Given an ep-groupoid $X$ and a compact oriented ep$^+$-subgroupoid $\wh{\Theta}:X\rightarrow {\mathbb Q}^+$
we have interpreted it as a full subcategory of $X$ with a weight function and we use some of the ambient
structure when we talk about local branching structures. Of course, it would be nice to define an intrinsic 
object, without reference to an ambient $X$, which would have similar properties, and would allow an integration theory, as well
as the notion of differential forms. Another question is concerned with the necessity of the tameness assumption.
It is possible that combining our discussion with results from a Lipschitz theory of manifolds would lead to an appropriate
generalization.

\subsection{Questions about M\texorpdfstring{$^+$}{pol}-Polyfolds}
A M$^+$-polyfold is a M-polyfold admitting an sc-smooth atlas, where the local models 
are sc$^+$-retracts $(O,C,E)$. If $d_M:M\rightarrow {\mathbb N}$  denotes  the associated degeneracy index, we know 
that near points $x\in M$ with $d_M(x)=0$ it has the structure of a smooth classical manifold.

In general the boundary behavior is complicated. However, as we have seen if $x\in \partial X$ and $X$ is tame near $x$,
then the fact that $T_xM$ is in good position to the partial quadrant $C_x$,  implies that $M$ near $x$ has the equivalent structure
of a classically smooth manifold with boundary with corners.  \par

\noindent {\bf Question 1:}  Since $O=r(U)$ for a relatively  open neighborhood of $0$ in the partial quadrant $C$ 
one can show that the sc$^+$-operator $Dr(0):E\rightarrow E$ has the property $Dr(0)(C)\subset C$. 
This distinguishes the smooth finite-dimensional linear subspaces $N$ of $E$ which can be written as the image 
of an sc$^+$-projection $P:E\rightarrow E$ with $P(C)\subset C$.  Let us refer to such $N$ as special $N$.
One can show that the interior of $N\cap C$ in $N$ is nonempty.
We can consider as local models for spaces open neighborhoods $V$ of $0\in N\cap C$, which we refer to as an open neighborhood
of $0$ in $(N,C,E)$, where $N$ is special. One can define in a standard way what a classically smooth map between two such
open sets is.  This allows us to define charts for metrizable spaces $M$ and smooth atlases.   This defines a class of
manifolds modeled on special $(N,C,E)$.  It is trivial to show that such manifolds are M$^+$-polyfolds using the property $P(C)\subset C$. 
The question is,  if every M$^+$-polyfold admits an equivalent atlas with special models?\par

\noindent{\bf Question 2:} Let $M$ be a compact M$^+$-polyfold of dimension $n$.  Does $\partial M$ have naturally the structure
of a Lipschitz-manifold of dimension $n-1$, see \cite{Heinonen,Rosen}. If that is the case one should be able 
to prove Stokes for every oriented, compact,  branched ep$^+$-subgroupoid $\wh{\Theta}$ using Stokes theorem in the Lipschitz context. 
In Remark \ref{rem116} we indicated even a possible generalization for defining M-polyfolds where $C$ is being replaced by more general
convex sets with nonempty interior. If it is indeed true that a Lipschitz version of Stokes can be applied, may be it even generalizes 
to such a more general setup.\par

The answers to these questions might allow to generalize some of the results about integration, but are also useful in studying the question
what is the intrinsic structure of an  ep$^+$-subgroupoid. We  raise some question about the latter in the next subsection.

\subsection{Questions about Branched Objects}
The definition of an ep$^+$-subgroupoid requires an ambient space. However, there seems to be an  underlying intrinsic object in the background.
The intrinsic object should be the metrizable space $|\supp(\Theta)|$ together with an additional structure. 
The following mentions some ideas to describe such an object, but the details have not been carried out.
Possibly these ideas need some modification.

Consider tuples $(M_U,w,U,C,E)$, where $E$ is an sc-Banach space,  $C\subset E$ a partial quadrant, $U\subset C$ a relatively open subset,
and $M_U$ a subset of $U$, and $w:M_U\rightarrow {\mathbb Q}^+$ is a map  with the following properties:
\begin{itemize}
\item[(1)]\ There exist finitely many  sc$^+$-retracts $(M_i,C,E)$ for $i\in I$ and positive rational weights $ {(\sigma_i)}_{i\in I}$
such that $M_i\subset U$ is properly embedded as submanifold.
\item[(2)]\ $M_U=\bigcup_{i\in I} M_i$.
\item[(3)]\ $w(y)= \sum_{\{i\in I\ |\ y\in M_i\}}\sigma_i$ for $y\in M_U$. 
\end{itemize}
We shall refer to $(M_U,w,U,C,E)$ as a branched manifold model. If $U'\subset C$ is relatively open and $U'\subset U$
we can define 
$$
(M_U,w,U,C,E)_{U'} =(M_{U'},w',U',C,E)
$$
 by restricting the data. We call this the restriction to $U'$.

We say that $(M_U^i,w^i,U^i,C^i,E^i)$ for $i=1,2$ are isomorphic provided there exists 
a homeomorphism $h:M_U^1\rightarrow M_U^2$ satisfying
\begin{itemize}
\item[(1)]\ $w^2\circ h =w^1$
\item[(2)]\ If $N\subset M_U$ is such that  $(N,C,E)$ is an sc$^+$-retract, then $(h(N),C,E)$ is an sc$^+$-retract.
The same for $h^{-1}$.
\item[(3)]\ $h:N\rightarrow h(N)$ is an sc-diffeomorphism.
\end{itemize}
The notion of equivalence accommodates the fact that $M_U$ can be written possibly  in different ways 
as union of M$^+$-polyfolds. Naturally the discussion of the geometry of local branching structures enters.
The local models have to be defined as subset of an ambient space to catch the variety of possible
branching structures. 

Let us refer to an object $\Xi=(M_U,w,U,C,E)$ as a branched M$^+$-polyfold model.  Following ideas in this book
one should be able to define the tangent $T\Xi$, and one should be able to define the notion of an orientation. 
There should be a notion of an sc-smooth map between such objects. Then one can define
 charts with the $\Xi$ as local models and sc-smooth atlases.  This would lead to the notion 
 of a branched M$^+$-polyfold.

 At this point we could employ the definition of an ep-groupoid, but take 
as local models the branched M$^+$polyfold models. Let us call such a space a branched ep$^+$-groupoid. 

The structure on $S=|\supp(\Theta)|$ would be given by an equivalence class $(Q,\beta)$, where $Q$ is a branched
ep$^+$-groupoid and $\beta:|Q|\rightarrow S$ is a homeomorphism. The notion of equivalence is parallel to the kind of Morita equivalence
used in giving a definition of orbifold using \'etale proper Lie groupoids, see \cite{AR,Mj,MM}. Such an idea is also applied
in Chapter \ref{chap11+}.
The boundaries of such spaces should be accessible to Lipschitz type methods mentioned in the previous subsection and there 
should be
Stokes theorem in this context.\par

\noindent{\bf Question 3:}  Is there a theory along the lines described above, i.e. build on branched local models, providing 
standard features expected from classical differential geometry or Lipschitz geometry (at least for the boundary portion)?

\chapter{Equivalences and Localization}\label{SEC2}
An equivalence between ep-groupoids is the sc-smooth version of an equivalence of categories.
In this chapter we shall study notions which are invariant under equivalences and introduce 
the notion of generalized maps. The ideas come from a well-known procedure in category theory, see \cite{GZ},
and also have been used in the classical Lie groupoid context, see \cite{Mj}.

\section{Equivalences}\label{SST_0}

The category ${\mathcal{EP}}$\index{${\mathcal{EP}}$} has  ep-groupoids as objects.  Its  class of morphisms $\textrm{mor}({\mathcal{EP}})$\index{$\textrm{mor}({\mathcal{EP}})$} consists of the sc-smooth functors between ep-groupoids and 
contains an interesting subclass  of morphisms called equivalences,  defined as follows.
\begin{definition}\label{equiv}\index{D- Equivalences}
An {\bf equivalence}  $F\colon X\rightarrow Y$ between the ep-groupoids  $X$ and $Y$ is a sc-smooth functor  having the following properties.
\begin{itemize}
\item[(1)]\ $F$ is faithful and full. 
\item[(2)]\  The induced map $|F|:|X|\rightarrow |Y|$ between the orbits spaces is a homeomorphism.
\item[(3)]\  $F$ is a local sc-diffeomorphism on objects.
\end{itemize}
We denote the {\bf class of equivalences}\index{Class of equivalences} by  ${\bf E}$\index{${\bf E}$}.
\qed
\end{definition}

We recall from the  category theory that a functor $F:X\to Y$ between two categories $X$ and $Y$ is called  {\bf faithful}\index{Faithful functor} if for every $x, x'\in X$, the induced map $\textrm{mor}_X(x, x')\to \textrm{mor}_Y(F(x), F(x'))$ is injective, 
and it is called {\bf full}\index{Full functor} if for every $x, x'\in X$, the induced map $\textrm{mor}_X(x, x')\to \textrm{mor}_Y(F(x), F(x'))$ is surjective.  If, in addition, the functor $F$ is {\bf essentially surjective}\index{Essentially surjective functor}, that is,  if 
for every object $y$ in $Y$ there exist an object $x$ in $X$ and an isomorphism $\varphi:F(x)\to y$, then the functor $F:X\to Y$ is called an {\bf equivalence}\index{Equivalence between categories}. If the functor $F:X\to Y$ between two categories $X$ and $Y$ is an equivalence in the sense of the category theory, then the induced map $\abs{F}:\abs{X}\to \abs{Y}$ is a bijection.

The following result characterizes equivalences between ep-groupoids.
\begin{lemma}\label{equivalence_in_the sense_of_category_theory}\index{L- Characterization of equivalences}
Assume that $F:X\rightarrow Y$ is a sc-smooth functor between ep-groupoids  which is an equivalence in the category sense 
and,  in addition,  is a local sc-diffeomorphism on the class of objects. Then  $F$ is an equivalence in the sense of Definition \ref{equiv}.
\end{lemma}
\begin{proof}
In view of the above remarks, it is enough to show that the induced map $\abs{F}:\abs{X}\to \abs{Y}$ is a homeomorphism.  Since the functor $F$ is sc-smooth and hence $\ssc^0$, the map 
$\abs{F}$ is  continuous between every level. Take a point $\abs{x}\in \abs{X}$ and its representative $x\in X$.  Then, by our assumption,  we find an open neighborhood $U$ of $x$ in $X$ so that 
$F:U\to F(U)$ is a sc-diffeomorphism. In particular, $F(U)$ is an open  subset of  $Y$. Since the  projection map $\pi:X\to \abs{X}$ is an open map and  $\abs{F}(\abs{U})=\abs{F(U)}$, it follows that $\abs{F(U)}$ is open in $\abs{Y}$. This holds for any open subset of $U$, and the proof is complete.
\qed \end{proof}

If $F\colon X\to Y$ and $G\colon Y\to Z$ are equivalences between ep-groupoids, then the composition $G\circ F\colon X\to Z$ is also an equivalence. In general, an equivalence is not invertible as a functor.

If $X$ is an ep-groupoid and ${\mathcal U}$ an open covering of the object M-polyfold, we define an ep-groupoid  $X_{{\mathcal U}}$ as follows. 
The object M-polyfold is the disjoint union of all the open sets $U\in {\mathcal U}$
$$
X_{{\mathcal U}}=\sqcup_{U\in {\mathcal U}} U.
$$
We shall write $(x,U)$ for its elements, where $x\in U$. The morphism M-polyfold ${\bm{X}}_{{\mathcal U}}$ consists of all triples  $(U,\phi,V)$ in which $U, V\in {\mathcal U}$ and $\phi\in {\bm{X}}$ is a morphism satisfying  $s(\phi)\in U$ and $t(\phi)\in V$.
The structural maps are the obvious ones. It is clear that this is a category whose object and morphism sets have M-polyfold structures induced from $X$ and ${\bm{X}}$, respectively.
We shall call $X_{{\mathcal U}}$ the {\bf refinement of the  ep-groupoid $X$}  associated with the  open covering ${\mathcal U}$.\index{Refinement  $X_{\mathcal U}$}\index{$X_{\mathcal U}$}
\begin{proposition}\label{fragmentation-x}\index{L- Fragmentation}
 $X_{{\mathcal U}}$ is an ep-groupoid and the map 
 $$
 F\colon X_{{\mathcal U}}\rightarrow X: (x,U)\mapsto  x
 $$
 is an equivalence of ep-groupoids.
 \end{proposition}
 \begin{proof}
 It is clear that object space and morphism space have natural M-polyfold structures coming from the ep-groupoid $X$.
The source and target maps,  defined by  $s(U,\phi,V)=(s(\phi),U)$ and $t(U,\phi,V)=(t(\phi),V)$,  are local sc-diffeomorphisms and $s,t\colon {\bm{X}}_{{\mathcal U}}\rightarrow X_{{\mathcal U}}$, 
 are by construction surjective. The smoothness properties of the structure  maps are obvious. Next we show that $X_{{\mathcal U}}$ is proper.  If  $(x,U)$ is  an object in $X_{{\mathcal U}}$,  we find an open neighborhood $V=V(x)$ of $x$ in $X$ such  that 
\begin{itemize}
\item[(1)]\ $\cl_X(V)\subset U$.
\item[(2)]\ $t:s^{-1}(\cl_X(V))\rightarrow X$ is proper. Here $s$ and $t$ are the source and target maps of the ep-groupoid  $X$.
\end{itemize}
The associated neighborhood of $(x,U)$ in $X_{{\mathcal U}}$ is then given by
$$
\wtilde{V}(x,U)=\{(y,U)\ |\ y\in V(x)\}.
$$
Note the closure of $\wtilde{V}\subset X_{{\mathcal U}}$ in $X_{{\mathcal U}}$ corresponds to the closure  of $V(x)$ in $X$.
Next, with $s$ and $t$ the source and target maps of $X_{{\mathcal U}}$ we show that 
$$
t:s^{-1}(\cl_{X_{{\mathcal U}}}(\wtilde{V}))\rightarrow X
$$ 
is proper.
We assume that $(\Phi_k)$ is a sequence in ${\bm{X}}_{|{\mathcal U}}$ with $s(\Phi_k)\in \cl(\wtilde{V})$ and $t(\Phi_k)$ belonging to a compact subset in $X_{|{\mathcal U}}$.
Writing $\Phi_k=(U_k,\phi_k,V_k)$ we know by assumption that $U_k=U$ and $s(\phi_k)\in \cl_X(V)$. Moreover $(t(\Phi_k))$ belongs to a compact subset 
$K$ of $X_{|{\mathcal U}}$. 
After perhaps taking a subsequence we may assume without loss of generality
that $V_k=W\in {\mathcal U}$ and $K\subset W$.
Hence we have a sequence $(\phi_k)\subset {\bm{X}}$ satisfying  $s(\phi_k)\in \cl_X(V)$ and $t(\phi_k)\in K$.
 Now using the properness of $X$ we may take another subsequence converging  to a morphism $\phi\in {\bm{X}}$ satisfying $s(\phi)\in \cl_X(V)$ and $t(\phi)\in K$. Then the convergence $\Phi_k\rightarrow (U,\phi,W)$ follows.
 Finally,  we note that $F$ is an equivalence of categories and a sc-smooth local sc-diffeomorphism, hence an equivalence of ep-groupoids in view of 
Lemma \ref{equivalence_in_the sense_of_category_theory}. 
  \qed \end{proof}

We recall from Section 
\ref{section_Tangent_Ep-Groupoid}, Definition \ref{natural_transf_def},  the notion of a natural equivalence between sc-smooth functors.
\begin{definition}[{\bf Natural equivalence}]\label{natural_equivalence}\index{D- Equivalence}
Two  sc-smooth functors $F:X\rightarrow Y$ and $G:X\to Y$ between the same ep-groupoids  are called {\bf naturally equivalent} if there exists 
 a sc-smooth map $\tau:X\rightarrow {\bm{Y}}$ having the following property.
 \begin{itemize}
 \item[(1)]\ For every object $x\in X$, $\tau (x)$ is a morphism,
$$
\tau (x)\colon F(x)\rightarrow G(x)\in {\bm{Y}}
$$
\item[(2)]\ The map $\tau$ is ``natural'' in the sense that 
$$
{\bf G}(\varphi)\circ \tau (x) = \tau (x')\circ {\bf F}(\varphi).
$$
for every morphism $\varphi:x\rightarrow x'$.  
\end{itemize}
The sc-smooth map $\tau\colon X\to {\bm{Y}}$ is called a {\bf natural transformation} of the two functors $F, G\colon X\to Y$ and symbolically referred to by 
$$\tau \colon F\to G.$$
\qed
\end{definition}

``Naturally equivalent'' is an equivalence relation for  sc-smooth functors $F, G\colon X\to Y$, which we denote by 
$$F\simeq G.$$
If $F\simeq G$ for two sc-smooth functors $F, G\colon X\to Y$ and if $\Phi\colon Y\to Z$ and $\Psi\colon Z\to X$ are sc-smooth functors, then 
$$F\circ \Psi\simeq G\circ \Psi\quad \text{and}\quad \Phi \circ F\simeq \Phi\circ G.$$
The following result is taken from \cite{HWZ3.5}, Proposition 2.7.
\begin{proposition}\label{equivalence_naturally_equivalent}\index{P- Equivalences {I}}
Assume that $F:X\to Y$ and $G:X\to Y$ are naturally equivalent  sc-smooth functors.  If one of them is an equivalence, so is the other.
\end{proposition}

Here is another useful property.
\begin{proposition}\label{another_property}\index{P- Equivalences {II}}
Assume that $F:X\rightarrow Y$ is an equivalence between ep-groupoids and $\Phi:Y\rightarrow Z$ a sc-smooth functor between ep-groupoids.
If $\Phi\circ F:X\rightarrow Z$ is an equivalence, so is $\Phi$.
\end{proposition}
\begin{proof}
Since on the object level $F$ and $\Phi\circ F$ are local sc-diffeomorphisms this holds for $\Phi$. 
Using that $\Phi\circ F$ is an equivalence it follows that it is essentially surjective, which then has to hold for $\Phi$.
Since $|F|$ and $|\Phi\circ F|$ are homeomorphisms we deduce via
$$
|\Phi| =|\Phi|\circ |F|\circ |F|^{-1} =|\Phi\circ F|\circ |F|^{-1}
$$
that $|\Phi|\colon \abs{X}\to \abs{Y}$ is a homeomorphism.

For two objects $y$ and $y'$ we have to show that the map $\Phi:{\bm{Y}}(y,y')\rightarrow{\bm{Z}}(\Phi(y),\Phi(y'))$ between the sets of morphisms is a bijection.
We assume that $\psi,\psi':y\rightarrow y'$ are morphisms satisfying $\Phi(\psi)=\Phi(\psi')$. Since $F$ is essentially surjective,  we find morphisms
$\phi:F(x)\rightarrow y$ and $\phi':F(x')\rightarrow y'$ and  consider the morphisms 
${(\phi')}^{-1}\circ\psi\circ \phi$ and  ${(\phi')}^{-1}\circ\psi'\circ \phi\colon F(x)\rightarrow F(x')$.  Because  $F$ is an equivalence,  there exist uniquely determined morphisms $\alpha,\alpha'\colon x\rightarrow x'$
solving 
$$
F(\alpha)={(\phi')}^{-1}\circ\psi\circ \phi\quad  \text{and}\quad F(\alpha')={(\phi')}^{-1}\circ\psi'\circ \phi.
$$
Applying the functor $\Phi$ to both expressions we compute
\begin{equation*}
\begin{split}
\Phi\circ F(\alpha)&=\Phi({(\phi')}^{-1}\circ\psi\circ \phi)\\
&=\Phi({(\phi')}^{-1})\Phi(\psi)\Phi(\phi)\\
&=\Phi({(\phi')}^{-1})\Phi(\psi')\Phi(\phi)\\
&=\Phi({(\phi')}^{-1}\circ\psi'\circ \phi)\\
&=\Phi\circ F(\alpha').
\end{split}
\end{equation*}
Since $\Phi\circ F$ is an equivalence,  we deduce $\alpha=\alpha'$,  implying
$$
{(\phi')}^{-1}\circ\psi\circ \phi={(\phi')}^{-1}\circ\psi'\circ \phi.
$$
Multiplying the inverses right and left,  we conclude that  $\psi=\psi'$. This proves the injectivity. 

A similar argument shows that  the map $\Phi:{\bm{Y}}(y,y')\rightarrow {\bm{Z}}(\Phi(y),\Phi(y'))$ is surjective.  Indeed, let  $\sigma:\Phi(y)\rightarrow\Phi(y')$ be a morphism in $ {\bm{Z}}(\Phi(y),\Phi(y'))$.  Since $F$ is essentially surjective, there exist 
$x,x'\in X$ and morphisms $\phi:F(x)\rightarrow y$ and $\phi':F(x')\rightarrow y'$. Then 
$$
\Phi(\phi')^{-1}\circ \sigma \circ\Phi(\phi) :\Phi\circ F(x)\rightarrow \Phi\circ F(x')
$$
belongs to $ {\bm{Z}}$,
and using that $\Phi\circ F$ is an equivalence we find a morphism $\alpha:x\rightarrow x'$ satisfying 
$$
\Phi\circ F(\alpha) = \Phi(\phi')^{-1}\circ \sigma \circ\Phi(\phi).
$$
We define the morphism  $\beta:y\rightarrow y'$ in ${\bm{Y}}$
by
$$
\beta=\phi'\circ F(\alpha) \circ  \phi^{-1}\in {\bm{Y}}(y, y')
$$
and compute,
\begin{equation*}
\begin{split}
\Phi(\beta)
&= \Phi(\phi')\circ \Phi\circ F(\alpha)\circ \Phi(\phi)^{-1}\\
&=\Phi(\phi')\circ \Phi(\phi')^{-1}\circ \sigma\circ \Phi(\phi)\circ \Phi(\phi)^{-1}\\
&=\sigma,
\end{split}
\end{equation*}
showing that $\Phi $ is surjective. 
This completes the proof of Proposition \ref{another_property}.
\qed \end{proof}
One of the basic results is the following theorem.
\begin{theorem}\label{gertrude}\index{T- Tangent of an equivalence}
If $F:X\rightarrow Y$ is an equivalence between ep-groupoids, then $TF:TX\rightarrow TY$ is an equivalence between ep-groupoids.
\end{theorem}

\begin{proof}
Since $F$ is a local sc-diffeomorphism, the same holds for the sc-smooth functor $TF:TX\to TY$. Hence, in view of Lemma \ref{equivalence_in_the sense_of_category_theory},  it suffices to show that $TF$ 
is an equivalence in the category theory sense. 
That $TF$ is a functor follows from the fact that $F$ is a functor.
Let us show next that $TF$ is essentially surjective.  We choose $a\in T_yY$. First of all we note that there exist 
$x\in X_1$ and $\phi\in {\bm{Y}}_1$ such  that 
$$
\phi:F(x)\rightarrow y.
$$
For  suitable open neighborhoods the source map  $s_Y:U(F(x))\rightarrow U(y)$ is a sc-diffeomorphism. We consider
for $x'$ near $x$ the map
$$
x'\mapsto  t_Y\circ (s_Y|U(\phi))^{-1}(F(x')).
$$
Differentiating at $x'=x$ in the direction $b\in T_xX$ gives the linear map
$$
b\mapsto Tt_Y(\phi)\circ (Ts_Y(\phi))^{-1}TF(x)(b)
$$
This is a linear isomorphism $T_xX\rightarrow T_yY$. Hence, we find a vector $b\in T_xX$ which is mapped to the given vector $a\in T_yY$.

Choosing $k=Ts_Y(\phi))^{-1}TF(x)(b)\in T_\phi{\bm{Y}}$, 
we obtain 
$$
Ts_Y(\phi)(k)=TF(x)(b)\ \ \text{and}\ \ Tt_Y(\phi)(k)=a.
$$
Therefore,  recalling Section \ref{section_Tangent_Ep-Groupoid}, the 
tangent vector $k\in T_\phi {\bm{Y}}$ is, by definition, the morphism 
$$
k\colon TF(b)\rightarrow a
$$
in $T{\bm{Y}}$.
This shows that $TF$ is essentially surjective. 
We leave the verification that $TF$ is full and faithful to the reader.
\qed \end{proof}

The next proposition shows that if two equivalences $F$ and $G$ satisfy $|F|=|G|$,   then  also $|TF|=|TG|$.
For the proof we shall need the following lemma.

\begin{lemma}\label{very_go}
Let $X$ be an  ep-groupoid, $x\in X$, and $U(x)$ a connected  open neighborhood of $x$ in $X$ equipped with the natural representation
$$
\Phi:G_x\rightarrow \text{Diff}_{sc}(U(x)).
$$
Suppose $\beta:U(x)\rightarrow X$ has an open image and  is an sc-diffeomorphism onto its image  satisfying $\beta(x)=x$, and 
$|\beta(y)|=|y|$ for all $y\in U(x)$. Then the following holds.
\begin{itemize}
\item[{\em(1)}]\ $\beta(U(x))=U(x)$ and $\beta:U(x)\rightarrow U(x)$ is an sc-diffeomorphism.
\item[{\em (2)}]\ There exists $g\in G_x$ such that $\beta(y)=g\ast y$ for all $y\in U(x)$.
\end{itemize}
\end{lemma}
\begin{proof}
Assume first that $X$ is effective.
Pick any point $y\in U(x)$ and take a continuous path $\gamma:[0,1]\rightarrow U(x)$ satisfying $\gamma(0)=x$ and $\gamma(1)=y$.
We consider the continuous map $\beta\circ \gamma:[0,1]\rightarrow X$. We show that the image of this map lies in $U(x)$.
Let $\varepsilon_0\in (0,1]$ be the largest number such that 
$$
\beta\circ\gamma([0,\varepsilon))\subset U(x).
$$
It is clear that $\varepsilon_0>0$.  Pick a sequence $(t_k)\subset [0,\varepsilon_0)$ converging monotonically 
to $\varepsilon_0$ and abbreviate $x_k:=\gamma(t_k)$.  Then $x_k\in U(x)$, $\beta(x_k)\in U(x)$ 
and $\beta(x_k) = g_k\ast x_k$ for a suitable $g_k\in G_x$. After perhaps taking a subsequence we may assume 
without loss of generality that $g=g_k$ for all $k$. Hence
$$
\beta(x_k)=g\ast x_k.
$$
By construction $x_k\rightarrow \gamma(\varepsilon_0)\in U(x)$ and therefore $\beta(x_k)\rightarrow g\ast \gamma(\varepsilon_0)\in U(x)$.
This implies $\beta(\gamma(\varepsilon_0))\in U(x)$. Therefore $\varepsilon_0=1$ and $\beta(U(x))\subset U(x)$.  

Next we show that $\beta(U(x))$ is open and closed in $U(x)$ which implies since $U(x)$ is connected that $\beta(U(x))=U(x)$.
By assumption $\beta(U(x))$ is open. We only need to show its closedness.  Pick $y\in U(x)$ so that there exists a sequence 
$(x_k)\subset U(x)$ with $\beta(x_k)\rightarrow y$. After perhaps taking a subsequence we find $g\in G_x$ such that 
$\beta(x_k)=g\ast x_k\rightarrow y$. This implies $x_k\rightarrow g^{-1}\ast y$. Hence $\beta(g^{-1}\ast y)=y$.
This completes the proof of (1). 

Pick any $y\in U_1(x)$, i.e. a point on level $1$. Since $X$ is metrizable and the neighborhood of $y$ is modeled on an sc-smooth retract
we can pick a metric $d$ which induces the topology (on level $0$) on $X$, but also 
for  a strictly decreasing sequence $\varepsilon_k$
converging to $0$ the spaces  
$$
\Sigma_k=\{z\in X\ |\ d(y,z)\leq \varepsilon_k\}\subset U(x)
$$
with the induced metric are complete.  We define the closed subsets 
$$
\Sigma_k^g=\{z\in \Sigma_k\ |\ \beta(z)=g\ast x\}.
$$
Then $\Sigma_k =\bigcup_{g\in G_x} \Sigma_k^g$ and by the Baire theorem there exists $g_k$ such that 
$\Sigma_k^{g_k}$ has a nonempty interior.  Hence we find $y_k$ with $d(y,y_k)<\varepsilon_k$ so
that $\beta(z)=g_k\ast z$ for $z$ near $y_k$. We may assume by slightly shifting $y_k$ that it is a smooth point.
We infer that $T\beta(y_k)=T\varphi_{g_k}(y_k)$. 
After perhaps taking a subsequence we may assume that $g_k=g$ and passing to the limit we find that
$$
T\beta(y)(h) = T\varphi_g(y)(h)\ \ \text{for all}\ \ h\in T_yX.
$$
We note that for large $k$ there cannot be a $g_k'\neq g$ with $\Sigma_k^{g'}$ having nonempty interior,
because it would imply by the same argument as above that also
$$
T\beta(y) =T\varphi_{g'}(y)
$$
for some $g'\neq g$. This, however is impossible by the current assumption that $X$ is effective.

Next take a large $k$ so that $\Sigma^g_k$ is the only set with nonempty interior. If $y'\in \Sigma_k$
with $d(y',y)<\varepsilon_k$ we can study by the same method a small neighborhood of $y'$
which is contained in $\Sigma^g_k$. By this argument we find an open set arbitrarily close to $y'$
on which $\beta(z)=g\ast z$.  Summing up we find an open dense subset of $\Sigma^g_k$ on which
$\beta(z)=g\ast z$. By continuity this must hold for all $z\in \Sigma_k^g$.  Hence we have proved
that every point $y\in U_1(x)$ has an open neighborhood $U(y)\subset U(x)$ on level $0$
so that $\beta(z)=g_y\ast z$ for $z\in U(y)$. Since $U_1(x)$ is also connected it follows that there is a single 
$g\in G_x$ such that $\beta(z)=g\ast z$ for all $z\in U_1(x)$. Using density of $U_1(x)$ in $U(x)$ and the continuity
of $\beta$ on level $0$ we finally conclude that 
$$
\beta(y)=g\ast y\ \ \text{for all}\ \ y\in U(x).
$$
Next we remove the effectivity assumption.  We can pass to $X_R$ which has the same object M-polyfold and  apply the argument
just established. This immediately implies the result. 
\qed \end{proof}
We use the previous lemma to prove the following useful result.
\begin{proposition}\label{prop4.14}
If  $F,G:X\rightarrow Y$ are two equivalences satisfying  
$$\abs{F}=\abs{G}\colon \abs{X}\rightarrow \abs{Y},$$
 then the following holds.
\begin{itemize}
\item[{\em (1)}]\ There exists for every $x\in X$   an sc-smooth map $U(x)\rightarrow \bm{Y}:z\rightarrow \alpha_z$
such that  
$$
\alpha_z:F(z)\rightarrow G(z).
$$
Here $U(x)$ is a connected open neighborhood of $x$ admitting the natural $G_x$-action.
\item[{\em (2)}]\ The tangent functors  $TF, TG:TX\rightarrow TY$ satisfy $\abs{TF}=\abs{TG}$.
\end{itemize}
\end{proposition}
\begin{remark}\index{R- Remark concerning Proposition \ref{prop4.14}}
In (1) we cannot prescribe $\alpha_x$, as the following example shows.
Let $X={\mathbb R}^2$ with the group ${\mathbb Z}_2$ acting with the nontrivial element
via $(x,y)\rightarrow (-x,-y)$. Then take the associated translation groupoid. We let $F$ be the identity functor and $G(x,y)=-(x,y)$.
Then we have the identity morphism $F(0,0)\rightarrow G(0,0)$, which however does not lie in an sc-smooth local family.
\end{remark}

Before we start  with the proof of the theorem we first recall from 
Section \ref{section_Tangent_Ep-Groupoid} that the tangent $T(X, { \bf X})=(TX, T{\bm{X}})$ of an ep-gropuoid $(X, {\bm{X}})$ is an ep-groupoid. The object space $TX$ is a M-polyfold and the morphism space $T{\bm{X}}$ is also a M-polyfold. The tangent vector $(\phi, k)\in T_{\phi}{\bm{X}}$ is viewed as a morphism
$$(\phi, h)\colon Ts(\phi )(h)\to Tt(\phi )(h)$$
between the two objects 
$Ts(\phi )(h)\in T_xX$ and $Tt(\phi )(h)\in T_yX$ in $TX$, if $\phi\in {\bm{X}}$ is a morphism $\phi\colon x\to y$.

\begin{proof}[Proposition \ref{prop4.14}]
Since $\abs{F}=\abs{G}$, given $x\in X$ there exists a morphism $\psi\colon F(x)\to G(x)$ belonging to ${\bm{Y}}$.  Since, by assumption $F$ and $G$ are equivalencies 
there exist open neighborhoods 
${\bm{U}}(\psi)$, $U(F(x))$, $U(G(x))$, $U(x)$, and $U'(x)$ such  that
$F\colon U(x)\rightarrow U(F(x))$, $G\colon U'(x)\rightarrow U(G(x))$, $s\colon {\bm{U}}(\psi)\rightarrow U(F(x))$ and $t\colon {\bm{U}}(\psi)\rightarrow U(G(x))$ are all
sc-diffeomor\-phisms. Hence we have the sequence of sc-diffeomorphisms
\begin{eqnarray}\label{unwrap}
U(x)\xrightarrow{F} U(F(x))\xrightarrow{s^{-1}}{\bm{U}}(\psi)\xrightarrow{t}U(G(x))\xrightarrow{G^{-1}}U'(x),
\end{eqnarray}
whose composition defines the sc-diffeomorphism 
$$
\beta:(U(x),x)\rightarrow (U'(x),x).
$$
Unwrapping the contents of the sequence  (\ref{unwrap}) there exist an sc-smooth map $U(x)\rightarrow \bm{Y}:z\rightarrow \psi_z$
with $\psi_x=\psi$, satisfying
$$
\psi_z :F(z)\rightarrow G(\beta(z))\ \ \text{for all}\ \ z\in U(x).
$$
Since $|F(z)|=|G(z)|$ we find $\phi_z:F(z)\rightarrow G(z)$, where at this point we do not know the precise dependence of $\phi_z$ on $z$.
Note that 
$$
\psi_z\circ \phi_z^{-1}:G(z)\rightarrow G(\beta(z))\ \ \text{for all}\ \ z\in U'(x).
$$
Since $F$ is an equivalence there exists $\sigma_z:z\rightarrow \beta(z)$ with $G(\sigma_z)=\psi_z\circ \phi_z^{-1}$. 
Hence $\beta:U(x)\rightarrow U'(x)\subset X$, $\beta(x)=x$, and $|z|=|\beta(z)|$. From  Lemma \ref{very_go}
we deduce that $U'(x)=U(x)$, $\beta:U(x)\rightarrow U'(x)$ is an sc-diffeomorphism, and there exists $g\in G_x$ such that
$$
\beta(z)=g\ast z\ \ \text{for all}\ z\in U(x).
$$
Define $\gamma_z:=\Gamma(g,z):z\rightarrow g\ast z$.  Then we define the sc-smooth map
$$
U(x)\rightarrow \bm{Y}:z\rightarrow \alpha_z:= G(\gamma_z^{-1})\circ \psi_z
$$
and note that 
$$
\alpha_z: F(z)\rightarrow G(z)\ \ \text{for all}\ \ z\in U(x).
$$

Next we consider (2).  Assume that $x\in X_1$.  Differentiating the sc-smooth expression
$$
\alpha_z : F(z)\rightarrow G(z)
$$
at $z=x$ and writing $\alpha:=\alpha_x$ we find that 
$$
T\alpha(x)\circ TF(x) =TG(x) :T_xX\rightarrow T_{G(x)}Y.
$$
This precisely means that $|TF(x)|=|TG(x)|$ for all $x\in X_1$ and therefore $|TF|=|TG|$.
\qed \end{proof}

We finally remark that equivalences between ep-groupoids descend to equivalences  between associated reduced ep-groupoids introduced at the end of Section \ref{subsec-Effective_Reduced}.
\begin{proposition}\label{prop2.16}
We assume that  $F:X\rightarrow Y$ is an equivalence of ep-groupoids.
We denote by $X_R$ and $Y_R$ the associated reduced ep-groupoids and define the map $F_R:X_R\to Y_R$ on the objects by $F_R=F$ and on the morphisms by $F_R([\phi])=[F(\phi)]$. 
Then $F_R$ is a well-defined functor and  $F_R:X_R\rightarrow Y_R$  is an equivalence. 
\end{proposition}
\begin{proof}
Since $F$ is an equivalence between ep-groupoids,  it follows immediately that $F(T_x)=T_{F(x)}$ and consequently $F_R$ is a well-defined functor.
Since $F=F_R$ on the objects,  the map $F_R$ is a local sc-diffeomorphism between the  objects. This implies that 
$F_R$ is a fully faithful functor. It is immediate that $F_R$ is essentially surjective, and hence an equivalence of categories.
All these facts together imply that $F_R$ is an equivalence of ep-groupoids.
\qed \end{proof}
\begin{example}
The assumption in Proposition \ref{prop2.16} that the functor $F\colon X\to Y$ is an equivalence cannot be omitted. Indeed, 
given ep-groupoids $X$ and $Y$ and a sc-smooth functor $\Phi:X\rightarrow Y$ it is in general not true that it passes to the reduced ep-groupoids as the following example illustrates.
Let $X={\mathbb Z}_4\ltimes {\mathbb R}$. Here the action of ${\mathbb Z}_4=\{0,1,2,3\}$ is given $0\ast x=x$, $1\ast x=-x$, $2\ast x=x$ and $3\ast x=-x$.
Let $Y={\mathbb Z}_4\ltimes {\mathbb R}^2$ with the action of ${\mathbb Z}_4$ by counter clock-wise rotation by $90$ degrees. We define
$\Phi$ on objects by $\Phi(x)=(0,0)$ and on morphisms by 
$$
\Phi(i,x)= (i,(0,0)),
$$
where $i=0,1,2,3$. Note that $T_x =\{(0,x),(2,x)\}$ for all $x\in {\mathbb R}$. For $(x,y)\in {\mathbb R}^2$ we have that $T_{(x,y)}=\{(0,(x,y))\}$ and note that $Y=Y_R$.
If $\phi_0=(i,x)$ and $\phi_1=(j,y)$ in ${\bm{X}}$ are equivalent it means that $x=y$ and $i=j$ mod $2$. Now we observe that
$\Phi((1,x))=(1,(0,0))$ and $\Phi(3,x)=(3,(0,0))$. These are however different points in $Y_R=Y$. Hence $\Phi$ does not pass to equivalence classes and consequently does not define
a functor between the reduced groupoids.
\qed
\end{example}

\section{The Weak Fibered Product}

 We continue with a useful  construction called the weak fibered product.

\begin{definition}\index{D- Weak fibered product}
Let  $F:X\rightarrow Y$ and $G:Z\rightarrow Y$ be 
sc-smooth functors between ep-groupoids.  The {\bf  weak fibered product}  $L=X\times_Y Z$\index{$X\times_Y Z$} is the category  whose objects are triples 
$(x,\varphi,z)$ in which $x\in X$, $z\in Z$, and $\varphi:F(x)\rightarrow G(z)$ is a morphism in ${\bm{Y}}$, i.e.
$$
L= X{_{F}\times_s}{\bm{Y}}{_{t}\times_G}Z.
$$
The morphism set ${\bm{L}}$ consists  of all triples  $(h,\varphi,k)\in {\bm{X}}\times{\bm{Y}}\times{\bm{Z}}$
satisfying  $s\circ F(h)=s(\varphi)$ and $s\circ G(k) = t(\varphi)$, so that 
$$
{\bm{L}}={\bm{X}}{_{s\circ F}\times_s}{\bm{Y}}{_{t}\times_{s\circ G}}{\bm{Z}}.
$$
The identity morphism $1_{(x, \varphi, z)}\in {\bm{L}}$ at the object $(x, \varphi, z)\in L$ is given by $1_{(x, \varphi, z)}:=(1_x, \varphi, 1_z)$. The inversion map $i:{\bm{L}}\to {\bm{L}}$ is defined by 
$i(h, \psi, k):=(h^{-1}, G(k)\circ \varphi\circ F(h)^{-1}, k^{-1})$. 
The source and target maps $s,t:{\bm{L}}\rightarrow L$ are defined by
$$
s(h,\varphi,k)=(s(h),\varphi,s(k))\ \ \text{and}\ \ t(h,\varphi,k)=(t(h), G(k)\circ\varphi\circ F(h)^{-1},t(k)).
$$
Finally, the multiplication $m\colon {\bm{L}}{_{s}\times_t}{\bm{L}}\to {\bm{L}}$ is defined by
$$m((h, \varphi, k), (h', \varphi', k'))=(h\circ h', \varphi', k\circ k').$$
\qed
\end{definition}
Since  $L=(L, {\bm{L}})$ is a small category, the weak fibered product $X\times_Y Z$ is a groupoid.  
Theorem 2.9 in \cite{HWZ3.5} shows that the weak fibered product is an ep-groupoid, provided at least one of the functors $F$ or $G$
is an equivalence.  Below we present a slightly modified proof of this fact.
\begin{theorem}\label{modified_weak_fibered_product}
\index{T- Weak fibered product}
Let $F:X\to Y$ be an equivalence of ep-groupoids, and let $G:Z\to Y$ be a sc-smooth functor.  Then the weak fibered product $X\times_Y Z$ is in a natural way an ep-groupoid and the projection
$\pi_Z: X\times_Y Z\rightarrow Z$ is an equivalence of  ep-groupoids.   Moreover, the sc-smooth functors 
$$F\circ \pi_X, G\circ  \pi_Z\colon X\times_Y Z\to Y$$
 are naturally equivalent and the degeneracy index $d_{X\times_Y Z}$ satisfies
$$
d_{X\times_Y Z}(x,\varphi,z)=d_Z(z),\quad \text{for every $(x,\varphi, z)\in X\times_Y Z$}.
$$
Moreover,  if $Z$ is a tame ep-groupoid so is the weak fibered product $ X\times_Y Z$. Further, if $Z$ has a paracompact orbit space, so has $X\times_Y Z$.
\qed
\end{theorem}

If $\pi_Z\colon X\times_YZ\to Z$ is an equivalence between ep-groupoids, then it is, by definition, a local sc-diffeomorphism between the object sets. Consequently, if one ep-groupoid is tame, so is the other, because tameness is invariant under sc-diffeomorphisms. Moreover, again by definition of an equivalence, the induced map $\abs{\pi_Z}\colon \abs{X\times_YZ}\to \abs{Z}$ is a homeomorphism between the orbit spaces. Consequently, if $\abs{Z}$ is paracompact, then also $\abs{X\times_YZ}$ is paracompact. Therefore, the last two statements of the theorem follow once we have shown that $\pi_Z$ is an equivalence.

If  we assume that $G:Z\to Y$ is an equivalence instead of $F:X\to Y$, then 
the fibered product $X\times_Y Z$ is again naturally an ep-groupoid and this time the projection 
$\pi_X:X\times _Y Z\rightarrow X$ is an equivalence of ep-groupoids. In this case,  the degeneracy index $d_{X\times_Y Z}$ satisfies 
$$d_{X\times_Y Z}(x,\varphi,z)=d_X(x),\quad \text{for every $(x,\varphi, z)\in X\times_Y Z$}.$$
If both $F:X\to Y$ and $G:Z\to Y$ are equivalences of ep-groupoids, then 
$$
d_X(x)=d_{X\times_Y Z}(x,\phi,z)=d_Z(z),\quad \text{for every $(x,\varphi, z)\in X\times_Y Z$}.
$$
The M-polyfold structures on $L$ and ${\bm{L}}$ in both cases are the natural ones induced
on $L\subset X\times{\bm{Y}}\times Z$ and ${\bm{L}}\subset {\bm{X}}\times{\bm{Y}}\times {\bm{Z}}$ as sub-M-polyfolds.

\begin{proof}[Proof of Theorem \ref{modified_weak_fibered_product}]
Clearly $ X\times{\bm{Y}}\times Z$ and ${\bm{X}}\times {\bm{Y}}\times {\bm{Z}}$  are  M-polyfolds.   We first show that $L$ is a  sub-M-polyfold.

Let  $(x,\varphi,z)\in L$.  Since $F:X\to Y$ is an equivalence and hence locally a sc-diffeomorphism, we find open neighborhoods $U=U(x)$ of $x$ in $X$ and $W=W(F(x))$ of $F(x)$ in $Y$, so that 
$$F\vert U:U\to W$$
is a sc-diffeomorphism. 
The source and the target maps $s, t:{\bm{Y}}\to Y$ are also local sc-diffeomorphisms.  Hence, replacing $U$ and $W$ by smaller neighborhoods if necessary,  we find open neighborhoods
 $W'=W'(G(z))$ of $G(z)$ in $Z$ and  ${\bm{U}}={\bm{U}}(\varphi)$ of $\varphi$  in ${\bm{Y}}$, such that the source and target maps 
$$
s\colon {\bm{U}}\to W\quad \text{and}\quad t\colon {\bm{U}}\to W'
$$
are sc-diffeomorphisms.  Moreover, since $G$ is sc-smooth and hence $\ssc^0$,  we find an open neighborhood $V=V(z)$ of $z$ in $Z$, so that 
$$G(V)\subset W'.$$
With these choices we define  the map 
$$
R:U\times {\bm{U}}\times V\to U\times {\bm{U}}\times V
$$ 
by 
\begin{equation}\label{retraction_weak_fibered_product_1}
R(x', \varphi', z')=\bigl( (F\vert U)^{-1}\circ  s\circ (t\vert {\bm{U}})^{-1}\circ G (z'), (t\vert {\bm{U}})^{-1}\circ G (z'), z').
\end{equation}
 Clearly $R$  is sc-smooth and $R\circ R=R$, so that $R$ is a sc-retraction.  The image of $R$ is equal to 
$$R(U\times {\bm{U}}\times V)=(U\times {\bm{U}}\times V)\cap L,$$
which implies that $L$ is a sub-M-polyfold of $X\times {\bm{Y}}\times Z$ and consequently inherits 
a M-polyfold structure. If $(V, \psi, (O, C, E))$ is a chart around $z$ on $Z$,  then 
the map $\Psi: (U\times {\bm{U}}\times V)\cap L\to O$,  defined by $\Psi ( x', \varphi', z')=\psi (z')$ defines a chart on $L$ around the point $(x, \varphi, z)$.
 We also observe that the projection map $\pi_Z:L\to Z$, defined by $\pi_Z(x, \varphi, z)=z$,  is in our charts the identity map. Consequently, $\pi_Z:L\to Z$ is a local sc-diffeomorphism. If $z\in Z$, then since $F$ is an equivalence, there exists an object $x\in X$ and a morphism $\varphi\in {\bm{Y}}$ between $F(x)$ and $G(z)$. Hence $(x, \varphi, z)\in L$ and $\pi_Z(x, \varphi, z)=z$, so that $\pi_Z$ is a surjection.

\noindent Next we shall show that ${\bm{L}}$ is also  M-polyfold.  We take a morphism  $(h, \phi, k)\in {\bm{L}}$. Since $F:X\to Y$ and $s:{\bm{X}}\to X$ are local sc-diffeomorphisms,  we find open neighborhoods ${\bm{U}}={\bm{U}}(h)$ of $h$ in ${\bm{X}}$ and $U=U(F(s(h)))$ of $F(s(h))$ in $Y$ such that 
$$F\circ s:{\bm{U}}\to U$$ 
is a sc-diffeomorphism.  Taking these neighborhoods smaller if necessary, we find, since $s, t:{\bm{Y}}\to Y$ are local sc-diffeomorphisms,  open neighborhoods ${\bm{V}}={\bm{V}}(\phi)$ of $\phi$ in ${\bm{Y}}$ and $V=V(G(s(k)))$ of $G(s(k))$ in $Y$, so that 
the maps 
$$s:{\bm{V}}\to U\quad \text{and}\quad t:{\bm{V}}\to  V$$
are sc-diffeomorphism. Finally, we choose an open neighborhood ${\bm{W}}={\bm{W}}(k)$ of $k$ in ${\bm{Z}}$ such that $(s\circ G)({\bm{W}})\subset {\bm{V}}$.  With these choices of open sets we define the 
map 
$${\bm{R}}\colon {\bm{U}}\times {\bm{V}}\times {\bm{W}}\to  {\bm{U}}\times {\bm{V}}\times {\bm{W}}$$ by 
$$
{\bm{R}}(h',\phi',k')=\bigl(   ((F\circ s)\vert {\bm{U}})^{-1}\circ s \circ (t\vert {\bm{V}})^{-1} \circ (s\circ G) (  k'),  (t\vert {\bm{V}})^{-1} \circ (s\circ G) (  k'), k'\bigr).
$$
Then ${\bm{R}}\circ {\bm{R}}={\bm{R}}$,  the map ${\bm{R}}$ is sc-smooth and has the image ${\bm{R}}( {\bm{U}}\times {\bm{V}}\times {\bm{W}})=( {\bm{U}}\times {\bm{V}}\times {\bm{W}})\cap {\bm{L}}$. Consequently, ${\bm{L}}$ is a sub-M-polyfold of the M-polyfold ${\bm{X}}\times {\bm{Y}}\times {\bm{Z}}$. Similarly as in the $L$-case, one defines charts on ${\bm{L}}$.  Moreover, the projection ${\bf \pi}_{\bm{Z}}:{\bm{L}}\to {\bm{Z}}$, defined by ${\bf \pi}_{\bm{Z}}(h, \varphi, k)=k$, is a local sc-diffeomorphism. 

Next we shall show that the structure maps are sc-smooth. We start with the source map $s:{\bm{L}}\to L$. Recall that 
$s(h, \varphi, k)=(s(h), \varphi, s(k))$ for every $(h, \varphi, k)\in {\bm{L}}$.  Locally, the source map $s$ is a composition of local sc-diffeomorphisms,  
$$s (h, \varphi, k)=\pi_Z^{-1}\circ s_Z\circ {\bf \pi}_{\bm{Z}}(h, \varphi, k),$$
and hence $s$ is a local sc-diffeomorphism. 
In our charts  on ${\bm{L}}$, the inversion map $i:{\bm{L}}\to {\bm{L}}$, given by $(h, \varphi, k)\mapsto (h, \varphi, k)^{-1}= (h^{-1},  G(k)\circ \varphi\circ F(h)^{-1} , k^{-1})$, takes the form $q\mapsto \psi'\circ i_Z\circ \psi^{-1}(q)$ where $\psi$ and $\psi'$ are suitable charts of $Z$. Since the inversion map $i_Z:{\bm{Z}}\to {\bm{Z}}$ is a local sc-diffeomorphism, the same holds for $i$. 

Now observe that the target map $t:{\bm{L}}\to {L}$ is equal to $t=i\circ s$ and so, $t$ is a local sc-diffeomorphism. 
The $1$-map $1:L\to {\bm{L}}$ is defined  by $1(x, \varphi, z)=(1_x, \varphi, 1_z)$ and,  in our charts on $L$ and ${\bm{L}}$,  is given by $q\mapsto  \psi' \circ 1_Z \circ \psi^{-1}(q)$ where $\psi$ and $\psi'$ are suitable charts on $Z$. Since the $1$-map $1_Z:Z\to {\bm{Z}}$ is sc-smooth, the same holds for $1_L:L\to {\bm{L}}$. Finally, we consider the multiplication map 
$m\colon {\bm{L}}{_{s}\times_t}{\bm{L}}\rightarrow {\bm{L}}$. We know that the fibered product 
$ {\bm{L}}{_{s}\times_t}{\bm{L}}$ is a  M-polyfold.  In local coordinates, the multiplication map $m$ takes the form 
$q\mapsto   \psi'\circ  m_Z(s^{-1}\circ t\cdot, \cdot )\circ  \psi^{-1}(q)$ for suitable charts $\psi$ and $\psi'$ on $Z$. This shows that the multiplication map $m$ is sc-smooth as claimed.

Next we will show that $L$ is proper. We take $(x, \varphi, z)\in L$ and choose neighborhoods $U=U(x)$ of $x$ in $X$ and $V=V(z)$ of $z$ in $Z$ 
such that the maps
$$t\colon s^{-1}\bigl(\ov{U}\bigr) \to X\quad \text{and}\quad t\colon s^{-1}\bigl(\ov{V}\bigr) \to Z$$
are proper.  We define the open neighborhood $W(x, \varphi, z)$ of $(x, \varphi, z)$ in $L$ by 
$$W(x, \varphi, z) =\{(x', \varphi', z')\in L\vert \, \text{$x'\in U,\, z'\in V, \,$  and $\varphi': F(x')\to G(z')$}\}$$
and claim that the restriction of the target map 
$$t:s^{-1}(\ov{W(x, \varphi, z)})\to L$$
is proper. To see this,  let $K$ be a compact subset of $\textrm{L}$ and let $((h_n, \psi_n, k_n))_{n\in \N}$ be a sequence of morphisms belonging to $ s^{-1}(\ov{W(x, \varphi, z')}$ so that $t(h_n, \psi_n, k_n)\in K$.  Since $K$ is compact, after taking a subsequence, the sequence  $t(h_n, \psi_n, k_n)$ converges in $\textrm{L}$ to the  morphism $(x', \psi', z')$. From 
$t(h_n, \psi_n, k_n)=(t(h_n), G(k_n)\circ \psi_n \circ F(h_n)^{-1}, t(k_n))$, it follows that $t(h_n)\to x'$, $t(k_n)\to z'$,  and $G(k_n)\circ \psi_n \circ F(h_n)^{-1}\to \psi'$.  Since 
$(h_n, \psi_n, k_n)\in s^{-1}(\ov{W(x, \varphi, z')})$, we conclude that $s(h_n)\in \ov{U}$ and $s(k_n)\in \ov{V}$. Then  the properness of the maps $t\colon s^{-1}\bigl(\ov{U}\bigr) \to X$ and $t\colon s^{-1}\bigl(\ov{V}\bigr) \to Z$ imply, after taking subsequences, the convergence  $h_n\to h$ in  $s^{-1}\bigl(\ov{U}\bigr)$. Hence 
the sequence of morphisms $(\psi_n)$ converges in ${\bm{Y}}$ to the morphism $G(k)^{-1}\circ \psi'\circ F(h)$.  Since $\psi_n$ is a morphism between $F(s(h_n))$ and $G(s(k_n))$, it follows that 
$\psi'$ is a morphism between $F(s(h))$ and $G(s(k))$. We have proved that the sequence $(h_n, \psi_n, k_n)$ has a subsequence converging to $(h, \psi, k)\in s^{-1}(\ov{W(x, \varphi, z')})$, and the proof of properness of $\textrm{L}$ is complete.

We already know that $\pi_Z$ is  surjective and a local sc-diffeomorphism, hence to show that $\pi_Z$ is an equivalence it suffices to show that for every two  objects
 $(x, \varphi, z)$, $ (x', \varphi', z')$ in $L$ the induced map 
$$\textrm{mor}_{L}((x, \varphi, z), (x', \varphi', z'))\to \textrm{mor}_{Z}(z, z')$$
is a bijection.  Take a morphism $k:z\to z'$. Then $\varphi'\circ G(k)\circ \varphi^{-1}$ is a morphism between $F(x)$ and $F(x')$.
Since $F:X\to Y$ is an equivalence and hence faithful and full, there exists a unique morphism $h$ in ${\bm{X}}$ between $x$ and $x'$ so that $F(h)=\varphi'\circ G(k)\circ \varphi^{-1}.$
Then $(h, \varphi, k)\in \textrm{mor}_{L}((x, \varphi, z), (x', \varphi', z'))$ and  $\pi_Z( (h, \varphi, k))=k$. Hence $\pi_Z$ is faithful and full, as claimed. The map $\pi_Z$ is also an essentially surjective. Indeed, let  $z\in Z$. Then since $F:X\to Y$ is an equivalence and hence, in particular, $\abs{F}:\abs{X}\to \abs{Y}$ is a bijection, there exists $x\in X$ and an isomorphism $\varphi\in {\bm{Y}}$ between  $F(x)$ and $G(z)$. Then $(x, \varphi, z)\in L$ and there exists a morphism between $\pi_Z(x, \varphi, z)=z$  and $z$, namely $1_z\in {\bm{Z}}$. Therefore $\pi_Z$ is an equivalence in the sense of category theory and hence $\abs{\pi_Z}:\abs{ X\times_Y Z}\to \abs{Z}$ is a bijection. Consequently, by Lemma \ref{equivalence_in_the sense_of_category_theory}, $\pi_Z$ is an equivalence in the sense of Definition \ref{equiv}. 

To see that the functors $F\circ \pi_X$ and $G\circ \pi_Z:X\times_YZ\to Y$ are naturally equivalent we define the map 
$\sigma\colon X\times_YZ\to {\bm{Y}}$ by 
$$\sigma (x, \varphi, z)=\varphi.$$
By definition of the triple $ (x, \varphi, z)\in X\times_YZ$, the morphism  $\varphi$ is a morphism between $F(x)$ and $G(z)$, hence between $F\circ \pi_X(x, \varphi, z)$ and $G\circ \pi_Z(x, \varphi, z).$  In order to  to show that $\sigma:F\circ \pi_X\to G\circ \pi_Z$ is a natural transformation we let $(x, \varphi, z)$ and $(x', \varphi', z')$ be two  triples in $X\times_YZ$ and $(h, \psi, k)$ a morphism between $ (x, \varphi, z)$ and $(x', \varphi', z')$. Then $h:x\to x'$, $k:z\to z'$, $\psi=\varphi$, and $\varphi'=G(k)\circ \varphi\circ F(h)^{-1}$. Consequently, 
\begin{equation*}
\begin{split}
 \sigma (x', \varphi', z') \circ (F\circ \pi_X)(h, \psi, k)&=\varphi'\circ F(h)\\
 &=G(k)\circ \varphi =(G\circ \pi_Z)(h, \psi, k)\circ  \sigma (x, \varphi, z).
 \end{split}
 \end{equation*}
It remains to verify that $\sigma: X\times_YZ\to {\bm{Y}}$ is sc-smooth. This follows from the fact that $\pi_Z:X\times_YZ\to Z$  is a local sc-diffeomorphism and hence we take it as a chart around $(x, \varphi, z)$. In this chart the composition $\sigma\circ \pi_{Z}^{-1}$ is equal to $(t\vert{\bm{U}})^{-1}\circ G$ where ${\bm{U}}$ is a suitable neighborhood of the  morphism $\varphi$. Since the maps $t$ and $G$ are sc-smooth, it follows that $\sigma$ is also sc-smooth. Hence we have proved the natural equivalence  $F\circ \pi_X\simeq G\circ \pi_Z$.
Finally, $\pi_Z\colon X\times_Y Z\to Z$ is a local sc-diffeomorphism and therefore,  
$$d_{X\times_Y Z}(x,\varphi,z)=d_Z(\pi_Z(x,\varphi, z))=d_Z(z)$$
for every $(x,\varphi, z)\in X\times_Y Z$. The proof of Theorem \ref{modified_weak_fibered_product} is complete.
\qed \end{proof}
\begin{remark}
Consider three ep-groupoids $A$, $B$, and $X$ and assume that $F:A\rightarrow X$ and $H:B\rightarrow X$ are  equivalences of ep-groupoids.
Then the weak fibered product $A\times_X B$  associated to the diagram $A\xrightarrow{F} X\xleftarrow{H} B$ is an ep-groupoid
and the projections $\pi_A: A\times_X B\rightarrow A$ and $\pi_B:A\times_X B\rightarrow B$ are equivalences of ep-groupoids in view of the previous discussion.  
By construction $F\circ \pi_A : A\times_X B\rightarrow X$ and $H\circ \pi_B: A\times_X B\rightarrow X$ are naturally equivalent
by $\tau:A\times_X B\rightarrow \bm{X}$ defined by
$$
\tau(a,\phi,b)=\phi.
$$
Given an equivalence $F:A\rightarrow X$ we may think of $A$ as a smaller model for $X$. The previous discussion can be understood
as saying that given smaller models for $X$, say $A$ and $B$, then there exists a model for $X$ which is smaller than $A$ and smaller than $B$.
\qed
\end{remark}
\section{Localization at the System of Equivalences}\label{section2.3_localization}
The category of  ep-groupoids ${\mathcal{EP}}$     \index{${\mathcal{EP}}$} has as  objects the ep-groupoids and as morphisms the  sc-smooth functors between ep-groupoids. They contain the distinguished class of morphisms ${\bf E}$\index{${\bf E}$},  called equivalences introduced in Definition \ref{equiv}.   Let us mention the following examples of equivalences. 
\begin{itemize}
\item[(E1)]\ \ \ Every isomorphism in the category ${\mathcal E}{\mathcal P}$ between two ep-groupoids is an equivalence.
\item[(E2)]\ \ \ If $F:X\rightarrow Y$ and $G:Y\rightarrow Z$ are equivalences, then $G\circ F$ is an equivalence. 
\item[(E3)]\ \ \ A functor which is naturally equivalent to an equivalence is an equivalence by Proposition 2.6.
\end{itemize}
We are going to construct a  new category
${\mathcal{EP}}({\bf E}^{-1})$\index{${\mathcal{EP}}({\bf E}^{-1})$} having  the same objects as ${\mathcal {EP}}$ for which there exists a functor
$$
i:{\mathcal{EP}}\rightarrow {\mathcal{EP}}({\bf E}^{-1})
$$
which is the identity on objects, but maps on morphisms every equivalence to an isomorphism. Moreover, this new category has some universal property.
The procedure is standard in category theory and called localization. In our case it has a quite explicit form.

From Theorem 
\ref{modified_weak_fibered_product} we recall that given an equivalence $F:X\rightarrow Y$ and a sc-smooth functor $\Phi:Z\rightarrow Y$ we can build the weak fibered product $X\times_Y Z$, which is an ep-groupoid. Moreover,  the projection $\pi_Z:X\times_Y Z\rightarrow Z$ is an equivalence  of ep-groupoids and the projection $\pi_X:X\times_Y Z\rightarrow X$ is a sc-smooth functor so that the sc-functors 
$\Phi\circ \pi_Z$ and $F\circ \pi_X$
are naturally equivalent via the natural equivalence $\tau\colon X\times_Y Z\rightarrow {\bm{Y}}$ defined by $\tau(x,\phi,z)=\phi$,
$$
F\circ\pi_X\xrightarrow{\tau} \Phi\circ\pi_Z.
$$
Next, given two ep-groupoids $X$ and $Y$, we  consider  {\bf diagrams} of the form
$$
 d\colon X\xleftarrow{F}Z\xrightarrow{\Phi} Y
 $$
\index{$d\colon X\xleftarrow{F}Z\xrightarrow{\Phi} Y$}
in which  $F\colon Z\to X$  is an equivalence and $\Phi\colon Z\to Y$ is a sc-smooth functor. 
\begin{definition}[{\bf Refinement}]
The diagram $d':X\xleftarrow{F'}Z'\xrightarrow{\Phi'} Y$ is called 
a {\bf refinement of the diagram}\index{D- Refinement of $d$} $d:X\xleftarrow{F}Z\xrightarrow{\Phi} Y$ if there exists an equivalence 
$$
H\colon Z'\rightarrow Z 
$$
of ep-groupoids such that  $F\circ H$ and $F'$ are naturally equivalent and $\Phi\circ H$ and $\Phi'$ are naturally equivalent.
\qed
\end{definition}

The situation, illustrated in the diagram
\begin{equation*}
\begin{CD}
d:X@<F<<Z@>\Phi>>Y\\
@. @AAHA @. \\
d':X@<F'<<Z'@>\Phi'>>Y,\\ 
\end{CD}
\end{equation*}
is sometimes abbreviated by
$$H:d'\rightarrow d.$$
For example, if 
$d: X\xleftarrow{F}Z\xrightarrow{\Phi} Y$ is a diagram and $H:Z'\rightarrow Z$ an equivalence, then 
$$d': X\xleftarrow{F\circ H}Z'\xrightarrow{\Phi\circ H} Y$$
is a refinement of  $d$. Another refinement of $d$ is the diagram 
$d':X\xleftarrow{G} Z\xrightarrow{\Psi} Y$
where $G\colon Z\rightarrow X$ is naturally equivalent to $F$ and $\Psi\colon Z\rightarrow Y$ is naturally equivalent to $\Phi$.  
The equivalence $H\colon Z\to Z$ is in this case $H=Id_Z$.
\begin{lemma}\index{L- Refinements of $d$}
Consider diagrams $d:X\xleftarrow{F}Z\xrightarrow{\Phi} Y$, $d':X\xleftarrow{F'}Z'\xrightarrow{\Phi'} Y$, and $d'':X\xleftarrow{F''}Z''\xrightarrow{\Phi''} Y$.  If $d'$ is a refinement of $d$ and $d''$ is a refinement of $d'$, then $d''$ is a refinement of $d$.
\end{lemma}
\begin{proof}
The verification  of the lemma is straightforward. Indeed,  we have the following diagram
\begin{equation*}
\begin{CD}
X@<F<<Z@>\Phi>>Y\\
@. @AAHA @. \\
X@<F'<<Z'@>\Phi'>>Y\\ 
@. @AAH'A @. \\
X@<F''<<Z@>\Phi''>>Y\\ \\
\end{CD}
\end{equation*}
in which  $F\circ H\simeq F'$, $F'\circ H'\simeq F''$ and $\Phi\circ H\simeq \Phi'$, $\Phi'\circ H'\simeq \Phi''$. 
Therefore,   $F\circ (H\circ H')\simeq F''$  and  $\Phi\circ (H\circ H')\simeq \Phi''$.
\qed \end{proof}
\begin{definition}[{\bf Common refinement}]\index{D- Common refinement}
Given two diagrams 
$$
d:X\xleftarrow{F} Z\xrightarrow{\Phi} Y\quad \text{and}\quad d':X\xleftarrow{F'} Z'\xrightarrow{\Phi'} Y,
$$
then a third diagram $d'':X\xleftarrow{F} Z''\xrightarrow{\Phi} Y$ is called  a {\bf common refinement} \index{Common refinement} of $(d, d')$, if $d''$ is a refinement of 
$d$ and a refinement of $d'$.
\qed
\end{definition}

Having a common refinement defines an equivalence relation on diagrams. This is a consequence of the following lemma.
\begin{lemma}\label{lemma_three}\index{L- Common refinements}
Consider the three diagrams 
$$d\colon X\xleftarrow{F}Z\xrightarrow{\Phi} Y,\quad d'\colon X\xleftarrow{F'}Z'\xrightarrow{\Phi'} Y\quad \text{and}\quad  d''\colon X\xleftarrow{F''}Z''\xrightarrow{\Phi''} Y.$$
 If $(d, d')$ has  a common refinement and  $(d', d'')$ has  a common refinement, then $(d, d'')$ has a common refinement.
\end{lemma}
\begin{proof}
Let $X\xleftarrow{G}W\xrightarrow{\Psi} Y $ be common refinement of $(d, d')$ so that we have the diagrams
\begin{equation*}
\begin{CD}
d:X@<F<<Z@>\Phi>>Y\\
@. @AAHA @. \\
\phantom{d:}X@<G<<W@>\Psi>>Y\\ \\
\end{CD}
\quad \text{and}\quad 
\begin{CD}
d':X@<F'<<Z'@>\Phi'>>Y\\
@. @AAH'A @. \\
\phantom{d':}X@<G<<W@>\Psi>>Y\\ \\
\end{CD}
\end{equation*}
in which $F\circ H\simeq G$,  $\Phi\circ H\simeq \Psi$ and $F'\circ H'\simeq G$,  $\Phi'\circ H'\simeq \Psi$. Similarly, let 
$X\xleftarrow{G'}W'\xrightarrow{\Psi} Y $ be  a common refinement of $(d', d'')$ illustrated by  the diagrams
\begin{equation*}
\begin{CD}
d':X@<F'<<Z'@>\Phi'>>Y\\
@. @AAKA @. \\
\phantom{d':}X@<G'<<W'@>\Psi'>>Y\\ \\
\end{CD}
\quad \text{and}\quad 
\begin{CD}
d'':X@<F''<<Z''@>\Phi''>>Y\\
@. @AAK'A @. \\
\phantom{d'':}X@<G'<<W'@>\Psi'>>Y\\ \\
\end{CD}
\end{equation*}
in which $F'\circ K\simeq G'$,  $\Phi'\circ K\simeq \Psi'$ and $F''\circ K'\simeq G'$,  $\Phi''\circ K'\simeq \Psi'$. 
Since $H':W\to Z'$ and $K:W'\to Z'$ are equivalences, we consider the weak  fibered  product $W''=W\times_{Z'}W'$. Then $G\circ \pi_W:W\times_{Z'}W'\to X$ is an equivalence and we claim that the digram
$$X\xleftarrow{G\circ \pi_W}W\times_{Z'}W'\xrightarrow{\Psi'\circ \pi_{W'}} Y $$
is a common refinement of the diagrams $d$ and $d''$.  Indeed,  $H\circ \pi_W:W\times_{Z'}W'\to Z$ is an equivalence as a composition of equivalences and 
\begin{equation*}
F\circ (H\circ \pi_W)=(F\circ H)\circ \pi_W\simeq G\circ \pi_W
\end{equation*}
and, using $\Phi'\circ H'\simeq \Phi$ and $H'\circ \pi_W\simeq K\circ \pi_{W'}$, 
\begin{equation*}
\begin{split}
\Phi\circ (H\circ \pi_W)&=(\Phi\circ H)\circ \pi_W\simeq\Psi\circ \pi_W\simeq (\Phi'\circ H')\circ \pi_W=\Phi'\circ (H'\circ \pi_W)\\
&\simeq\Phi'\circ (K\circ \pi_{W'})=(\Phi'\circ K)\circ \pi_{W'}\simeq \Psi'\circ \pi_{W'},
\end{split}
\end{equation*}
proving our claim.

Similarly, one shows that $X\xleftarrow{G\circ \pi_W}W\times_{Z'}W'\xrightarrow{\Psi'\circ \pi_{W'}} Y $ is refinement of $d''$ with an equivalence $K'\circ \pi_{W'}:W\times_{Z'}W'\to Z''$.
This completes the proof of Lemma \ref{lemma_three}.
\qed \end{proof}
\begin{definition}[{\bf Equivalence of diagrams}]\index{D- Common refinements as equivalence relation}
Two diagrams 
$$
d\colon X\xleftarrow{F}Z\xrightarrow{\Phi}Y\quad \text{and}\quad d'\colon X\xleftarrow{F'}Z'\xrightarrow{\Phi'}Y
$$
are called {\bf equivalent}\index{equivalence of diagrams}  if they have a {\bf common refinement}.
\qed
\end{definition}

This defines an equivalence relation on diagrams. We denote the {\bf equivalence class}  of $d$ by $[d]$ or by $[ X\xleftarrow{F}Z\xrightarrow{\Phi}Y]$ 
\index{$[ X\xleftarrow{F}Z\xrightarrow{\Phi}Y]$} and are going to  interpret $[d]$  as a morphism
$$
[d]\colon X\rightarrow Y
$$
between the ep-groupoids.\index{$[d]\colon X\rightarrow Y$}

We shall show shortly that such equivalence classes $[d]$  can be composed and the ep-groupoids with these type of morphisms define
a new category (denoted by ${\mathcal E}{\mathcal P}({\bf E}^{-1})$) \index{${\mathcal E}{\mathcal P}({\bf E}^{-1})$} 
possessing  useful properties.

Since the equivalence $F\colon Z\to X$ induces a homeomorphism $\abs{F}\colon \abs{Z}\to \abs{X}$ between the orbit spaces, we can associate with the diagram $d\colon X\xleftarrow{F}Z\xrightarrow{\Phi}Y$,  the map 
$\abs{d}\colon \abs{X}\rightarrow \abs{Y}$  between the orbit spaces, defined by 
$$\abs{d}:=\abs{\Phi}\circ \abs{F}^{-1}.$$
The  map $\abs{d}$ is continuous between all levels of the orbit spaces.  Let   $d'\colon X\xleftarrow{F'}Z'\xrightarrow{\Phi'}Y$  be   a refinement  of  the  digram  $d\colon X\xleftarrow{F}Z\xrightarrow{\Phi}Y$,  illustrated by the  diagram  
\begin{equation*}
\begin{CD}
X@<F<<Z@>\Phi>>Y\\
@. @AAHA @. \\
X@<F'<<Z'@>\Phi'>>Y. \\
\end{CD}
\end{equation*}
Then 
$$
\abs{F}\circ \abs{H} =\abs{F'}\quad \text{and}\quad \abs{\Phi}\circ\abs{H}=\abs{\Phi'},
$$
and therefore, since $\abs{F}$, $\abs{F'}$, and $\abs{H}$ are homeomorphisms, 
$$
\abs{d'}= \abs{\Phi'}\circ \abs{F'}^{-1}= \abs{\Phi}\circ \abs{H}\circ \abs{H}^{-1}\circ \abs{F}^{-1}= \abs{\Phi}\circ \abs{F}^{-1} =\abs{d}.
$$
Hence we have proved the following lemma.
\begin{lemma}\index{L- Induced map by a diagram}
If $H\colon d'\rightarrow d$ is a refinement of diagrams, where $d,d'\colon X\rightarrow Y$, then  the induced maps $|d|,|d'|\colon |X|\rightarrow |Y|$ between the orbit spaces are the same, i.e.,  $|d|=|d'|$.
\qed
\end{lemma}
This implies that we have a well-defined map
associated with  an equivalence class of diagrams $[d]:X\rightarrow Y$,  namely the map
$$
\abs{[d]}:\abs{X}\to \abs{Y}, 
$$
defined by 
$$
|[d]|:= |\Phi|\circ |F|^{-1},
$$
 where we have chosen  any representative $d:X\xleftarrow{F}Z\xrightarrow{\Phi} Y$ of the equivalence class $[d]$.
  We note that if in the equivalence class $[d]=[X\xleftarrow{F}Z\xrightarrow{\Phi} Y]$ both functors $F$ and $\Phi$ are equivalences, then the induced map $\abs{[d]}\colon \abs{X}\to \abs{Y}$ is a homeomorphism between the orbit spaces.
 
\begin{definition}\index{D- Induced map}
Given an equivalence class $[d]$ of diagrams $d\colon X\rightarrow Y$,  the {\bf induced map  between the orbit spaces}
is denoted by 
$$\abs{[d]}\colon \abs{X}\rightarrow \abs{Y}.$$
\qed
\end{definition}
Next we shall show  for two equivalence classes  $[d]:X\rightarrow Y$ and $[d']:Y\rightarrow Z$ that there is a well-defined associative 
composition $[d']\circ [d]$ defining an equivalence class $[d'']$ of diagrams $d''\colon X\rightarrow Z$. 

\begin{theorem}[{\bf Composition of equivalence classes}]\label{proppp}\index{T- Composition for $[d]$}
Let 
$$
d\colon X\xleftarrow{F} A\xrightarrow{\Phi}Y\quad \text{and}\quad d'\colon Y\xleftarrow{G} B\xrightarrow{\Psi} Z
$$
be representatives  of the equivalence classes $[d]\colon X\rightarrow Y$ and $[d']\colon Y\rightarrow Z$, respectively.  Then the equivalence class $[d'']$ of the diagram 
$$
d''\colon X\xleftarrow{F\circ\pi_A} A\times_Y B\xrightarrow{\Psi\circ\pi_B} Z
$$
is independent of the choices of representatives in the equivalence classes $[d]$ and $[d']$. The equivalence class $[d'']:X\to Z$ is called the {\bf composition}
of $[d]$ and $[d']$ and is denoted by 
$$[d'']=[d']\circ [d].$$
Moreover, the composition is associative; if $[d]:X\rightarrow Y$, $[d']:Y\rightarrow Z$,  and $[d'']:Z\rightarrow W$, then 
$$
([d'']\circ [d'])\circ [d]=[d'']\circ ([d']\circ [d]).
$$
In addition, 
$$
[d]\circ [1_X]=[d]\quad \text{and}\quad  [1_Y]\circ [d]=[d]
$$
for every equivalence class $[d]:X\rightarrow Y$.  Here $[1_X]$ is  the equivalence class 
$$[1_X]=[X\xleftarrow{1_X} X\xrightarrow{1_X}X].$$
\qed
\end{theorem}
The somewhat lengthy technical proof is postponed to  Appendix \ref{x-proppp}.

\begin{definition}[{\bf Generalized maps, generalized isomorphisms}]   \index{D- Generalized map}\index{D- Generalized isomorphism}
Let $X$ and $Y$ be two ep-groupoids.
An equivalence class $[a]=[X \xleftarrow{F}A\xrightarrow{G}Y]$ of diagrams in which $F$ is an equivalence and $G$ is a sc-smooth functor is called a {\bf generalized map}, and abbreviated by 
$$[a]\colon X\to Y$$
if the representation is irrelevant.
A generalized map $[a]\colon X\to Y$ is called a 
{\bf generalized isomorphism} or {\bf invertible}, if there exists a generalized map $[b]\colon Y\to X$ satisfying 
$$[b]\circ [a]=[1_X]\quad \text{and}\quad [a]\circ [b]=[1_Y].$$
In this case the equivalence class $[b]$ is called the {\bf inverse} of $[a]$ and denoted by 
$$[b]=[a]^{-1}.$$
\qed
\end{definition}
In the following theorem we characterize generalized isomorphisms.
\begin{theorem}\label{strong-iso}\index{T- Characterization of generalized inverse}
Let $X$, $Y$, and $A$  be ep-groupoids.
\begin{itemize}
\item[{\em(1)}]\
The equivalence class $[a]=[X\xleftarrow{F}A\xrightarrow{G}Y]$ in which $F$ is an equivalence and $G$ is a sc-smooth functor  is a generalized isomorphism if and only if the functor $G$ is an  equivalence. 
\item[{\em(2)}]\
The inverse of a generalized isomorphism $[a]=[X\xleftarrow{F}A\xrightarrow{G}Y]$ is the equivalence class 
$$[a]^{-1}=[Y\xleftarrow{G}A\xrightarrow{F}X].$$
\end{itemize}
\end{theorem}
The proof is postponed to  Appendix \ref{x-strong-iso}. \qed

Assume that $\mathfrak{f}:X\rightarrow Y$ and $\mathfrak{g}:Y\rightarrow Z$ are generalized maps between ep-groupouids.
Then their composition is defined giving $\mathfrak{g}\circ \mathfrak{f}:X\rightarrow Z$.  Assume that $\mathfrak{f}=[d]$ and $\mathfrak{g}=[d']$
with representatives $d$ and $d'$ given by
$$
X\xleftarrow{F} A \xrightarrow{\Phi} Y\ \ \text{and}\ \ Y\xleftarrow{G} B\xrightarrow {\Psi}  Z,
$$
respectively.   We know already that $|\mathfrak{f}|:|X|\rightarrow |Y|$ is given by $|\mathfrak{f}|= |\Phi|\circ |F|^{-1}$ 
and $|\mathfrak{g}|:|Y|\rightarrow |Z|$ by $|\mathfrak{g}|=|\Psi|\circ |G|^{-1}$. As a representative $d''$ for $[d']\circ [d]=\mathfrak{g}\circ \mathfrak{f}$ we can take 
the diagram 
$$
d''\colon X\xleftarrow{F\circ \pi_A} A\times_Y B\xrightarrow{\Psi\circ \pi_B} Z.
$$
Here $F\circ \pi_A$ is an equivalence.  We compute $|d''|=|[d'']|$
Given  $|x|\in |X|$  we find $(a,\varphi,b)\in A\times_Y B$ and a morphism $\sigma:F(a)\rightarrow x$. Then  by definition $|d''|(|x|) = |\Psi(b)|$.
Observe that $\varphi:\Phi(a)\rightarrow G(b)$ by definition of the weak fibered product.  Hence $|\Phi(a)|=|G(b)|$. This implies, using that $G$ is an equivalence
\begin{eqnarray*}
|[d'']|(|x|)&=& |\Psi(b)|\\
&=& |\Psi|\circ |G|^{-1}\circ |\Phi(a)| \\
&=& |\Psi|\circ |G|^{-1}\circ |\Phi|\circ |F|^{-1}(|x|)\\
&=&|[d']|\circ |[d]|(|x|).
\end{eqnarray*}
Hence we have proved the following result
\begin{proposition}\label{equal_prop}\index{P- Composition and induced maps}
If $\mathfrak{f}\colon X\to Y$ and $\mathfrak{g}\colon Y\to Z$ are generalized maps  between ep-groupoids,  then 
$$\abs{\mathfrak{g}\circ \mathfrak{f}}=\abs{\mathfrak{g}}\circ \abs{\mathfrak{f}}.$$
\qed
\end{proposition}

In view of the Theorems  \ref{proppp}  and \ref{strong-iso} we can now introduce the new category whose morphisms are the generalized maps.
\begin{definition}\index{D- Category ${\mathcal{EP}}({\bf E}^{-1})$}
The {\bf category ${\mathcal{EP}}({\bf E}^{-1})$}  has as  objects the objects of ${\mathcal{EP}}$, namely the ep-groupoids $X$,  and as  morphisms
the equivalence classes $[d]\colon X\to Y$ of diagrams (the generalized maps) between ep-groupoids. 
\qed
\end{definition}

The notation ${\mathcal{EP}}({\bf E}^{-1})$ indicates that we added inverses to the elements in ${\bf E}$.
  The basic result  relating  ${\mathcal E}{\mathcal P}$ with  ${\mathcal E}{\mathcal P}({\bf E}^{-1})$
is  the following theorem.
\begin{theorem}[{\bf Properties of the Localization at ${\bf E}$}]\label{localization}\index{T- Localization}
There is a  natural functor 
$$
i\colon {\mathcal{EP}}\rightarrow {\mathcal{EP}}({\bf E}^{-1}),
$$
which is the identity on objects and maps every  sc-smooth functor $\Phi\colon X\rightarrow Y$ between ep-groupoids to the equivalence class
$[d_\Phi]$ defined by 
$$
[d_\Phi]=[X\xleftarrow{1_X} X\xrightarrow{\Phi} Y].
$$
If the functor $F:X\rightarrow Y$ is an equivalence, then   the morphism $i(F)$ is invertible and
$$
i(F)^{-1} =[Y\xleftarrow{F} X\xrightarrow{1_X} X].
$$
\end{theorem}
\begin{proof}
We first verify  that 
$$
[d_\Psi]\circ [d_\Phi] = [d_{\Psi \circ\Phi}].
$$
for sc-smooth functors $\Phi:X\rightarrow Y$ and $\Psi:Y\rightarrow Z$.  Since the diagrams 
$$
d_\Phi \colon X \xleftarrow{1_X} X\xrightarrow{\Phi} Y\quad \text{and}\quad d_\Psi\colon Y\xleftarrow{1_Y} Y\xrightarrow{\Psi} Z, 
$$
are representative of the equivalence classes $[d_\Phi]$ and $[d_\Psi]$, the diagram 
\begin{equation}\label{diagram_1}
X\xleftarrow{\pi_X}  X\times_Y Y\xrightarrow{\Psi\circ \pi_Y} Z
\end{equation}
is a representative of the class $[d_\Psi]\circ [d_\Phi]$.
The class $[d_{\Psi\circ \Phi}]$ is represented by the diagram
$$d_{\Psi\circ \Phi}\colon  X\xleftarrow{1_X} X\xrightarrow{\Psi\circ \Phi} Z, $$
and we claim that the diagram $d_{\Psi\circ \Phi}$ is a refinement of the digram \eqref{diagram_1}.  For the equivalence between $X\times_Y Y$ and $X$ we take the  projection $\pi_X:X\times_Y Y\to X$. Then $1_X\circ \pi_X=\pi_X$ and, since $\Phi\circ \pi_X\simeq \pi_Y$ by  Theorem  \ref{modified_weak_fibered_product},   we find that
\begin{equation*}
(\Psi\circ \Phi)\circ \pi_X=\Psi\circ (\Phi\circ \pi_X)\simeq \Psi \circ \pi_Y.
\end{equation*}
Therefore, 
\begin{equation*}
i(\Psi\circ\Phi)=[d_{\Psi\circ \Phi}]=[d_{\Psi}]\circ [d_{\Phi}]=i(\Psi)\circ i(\Phi).
\end{equation*}
We already know that $i(1_X)=[X\xleftarrow{1_X}X\xrightarrow{1_X} X]=[1_X]$.

 Next we assume that $F:X\to Y$ is an equivalence and let $i(F)=[d_F]$ where $d_F$ is the diagram 
$$d_F\colon X\xleftarrow{1_X} X\xrightarrow{F}Y.$$ 
Considering the diagram $d\colon Y\xleftarrow{F} X\xrightarrow{1_X}X$, we claim that $[d]\circ [d_F]=[1_X]$ and $[d_F]\circ [d]=[1_Y]$.  The composition 
$[d]\circ [d_F]$ is represented by the diagram
\begin{equation}\label{invertible_1}
X\xleftarrow{\pi_X} X\times_Y X\xrightarrow{\pi_X} X.
\end{equation}
Since $F:X\to Y$ is an equivalence, the projection $\pi_X\colon X\times_Y X\to X$ is an equivalence by Theorem \ref{modified_weak_fibered_product}. Consequently, in view of
Theorem \ref{modified_weak_fibered_product},  the diagram in \eqref{invertible_1} is a representative of the equivalence class $[1_X]$, that is, $[d]\circ [d_F]=[1_X]$. The proof of 
$[d_F]\circ [d]=[1_Y]$ is similar. Therefore,  $i(F)$ is invertible and $i(F)^{-1}=[ Y\xleftarrow{F} X\xrightarrow{1_X}X]$. The proof of Theorem \ref{localization} is complete.
\qed \end{proof}
We summarize the fact that generalized maps induce continuous maps between orbit spaces, compatible with compositions,  as follows, where $\text{TOP}$ is the category of topological spaces.
\begin{theorem}\index{T- The functor $\vert .\vert$ on ${\mathcal{EP}}({\bf E}^{-1})$}
There exists a natural functor $|.|: {\mathcal{EP}}({\bf E}^{-1})\rightarrow \text{TOP}$ which associates to an ep-groupoid $X$ its orbit space $|X|$
and to a generalized map 
$$
\mathfrak{f}=[d]:X\rightarrow Y,
$$
 where the representative $d$ is given by $X\xleftarrow{F} A\xrightarrow{\Phi} Y$,
the continuous map $|[d]| =|\Phi|\circ |F|^{-1}:|X|\rightarrow |Y|$.
\qed
\end{theorem} 
\begin{remark}\index{R- On polyfolds}
Before we move to the strong bundle case we give an example why the previous considerations are important. These ideas will studied in more generality later on.
Assume that $Z$ is a metrizable space. Consider pairs $(X,\beta)$, where $\beta:|X|\rightarrow Z$ is a homeomorphism. 
In this case $|X|$ is metrizable as well. Given a second pair $(X',\beta')$ we say it is equivalent to $(X,\beta)$ provided there exists a generalized
isomorphism $\mathfrak{f}:X\rightarrow X'$ satisfying $\beta'\circ |\mathfrak{f}|=\beta$. The metrizable space $Z$ equipped with an equivalence class
of pairs $(X,\beta)$ will be called an sc-smooth polyfold.  It  is similar to an orbifold, with the difference that  as local models we allow $G\backslash O$, were $G$ is a finite group acting 
on the M-polyfold $O$.  Later we study which notions behave well with respect to equivalences. As a by-product one obtains the notions which make sense 
on polyfolds.
\qed
\end{remark}

\begin{theorem}\label{THMX10316}\index{T- Uniqueness theorem for generalized isomorphisms}
Assume that $X$ and $Y$ are ep-groupoids and $\mathfrak{f},\mathfrak{f}':X\rightarrow Y$ generalized isomorphisms 
satisfying $|\mathfrak{f}|=|\mathfrak{f}'|$. Then $\mathfrak{f}=\mathfrak{f}'$.
\end{theorem}
 \begin{proof}
 We take  representatives $d:X\xleftarrow{F} A\xrightarrow{H} Y$ for $\mathfrak{f}$  and
$d:X\xleftarrow{F'} A'\xrightarrow{H'} Y$ for $\mathfrak{f'}$.  Take the weak fibered product $A\times_X A'$
associated to the diagram 
$$
A\xrightarrow{F} X\xleftarrow{F'} A'.
$$
Then 
$$
X\xleftarrow{F\circ \pi_A} A\times_X A'\xrightarrow{H\circ\pi_A} Y
$$
 is a diagram of equivalences  equivalent to $d$.
The diagram 
$$
X\xleftarrow{F'\circ \pi_{A'}} A\times_X A' \xrightarrow{H'\circ \pi_{A'}} Y
$$
is equivalent to $d'$. We also not that $F\circ\pi_A$ and $F'\circ\pi_{A'}$ are naturally equivalent.
In view of this discussion
 we have shown that we can take the representatives for $\mathfrak{f}$ and $\mathfrak{f}'$ 
in such a way that  $A=A'$ and $F$ and $F'$ are naturally equivalent. This even allows us to replace $F'$ by $F$.
Hence 
\begin{eqnarray*}
& d\colon X\xleftarrow{F} A\xrightarrow{H} Y&\\
&d'\colon X\xleftarrow{F} A\xrightarrow{H'} Y.&
\end{eqnarray*}
The property that $|\mathfrak{f}|=|\mathfrak{f}'|$, or equivalently $|[d]|=|[d']|$,  then means that $|H|=|H'|$. 
To proceed with our argument we study the equivalences $H,H':A\rightarrow Y$ which have the property that the induced homemorphisms 
between the orbit spaces are the same. In view of Proposition \ref{prop4.14} we find for every object $a\in A$ and open neighborhood
$U(a)$ so that there exists a (sc-smooth) natural transformation $\alpha_a:U(a)\rightarrow \bm{X}$ 
such that 
$$
\alpha_a(b):H(b)\rightarrow H'(b)\ \ \text{for}\ \ b\in U(a).
$$
The collection ${\mathcal U}$ of all these $U(a)$, where $a$ varies over the objects in $A$ is an open covering and 
associated to it we have the ep-groupoid $A_{\mathcal U}$ together with natural  equivalences
$$
H_{\mathcal U}: A_{\mathcal U}\rightarrow Y\ \ \text{and}\ \ H_{\mathcal U}':A_{\mathcal U}\rightarrow Y.
$$
From the $(\alpha_a)$ we obtain a natural transformation $\alpha:H_{\mathcal U}\rightarrow H_{\mathcal U}'$.
We the natural equivalence $\pi: A_{\mathcal U}\rightarrow A$ the diagram
$$
X\xleftarrow{F\circ \pi} A_{\mathcal U}\xrightarrow{H_{\mathcal U}} Y
$$
is equivalent to $X\xleftarrow{F}A\xrightarrow{H} Y$ and 
$$
X\xleftarrow{F\circ \pi} A_{\mathcal U}\xrightarrow{H_{\mathcal U}'} Y
$$
is equivalent to  $X\xleftarrow{F}A\xrightarrow{H'} Y$. By construction the two displayed diagrams are equivalent 
which proves that $\mathfrak{f}=\mathfrak{f}'$.
  \qed \end{proof}

 \section{Strong Bundles and Equivalences}\label{STE__x}

Next we discuss the notion of equivalence between strong bundles over ep-groupoid. 
In Section \ref{SST} we have introduced the notion of a strong bundle $(P\colon W\to X, \mu)$ over the ep-groupoid $X$, where $P\colon W\to X$ is a strong bundle over the object M-polyfold $X$ and the structure map $\mu$ is a strong bundle isomorphism 
$$
\mu:{\bm{X}}{_{s}\times_P W}\rightarrow W
$$
covering the target map $t:{\bm{X}}\rightarrow X$ and possessing additional  properties.  In particular,  $\mu$ lifts a morphism $\phi\colon x\rightarrow y$
to a linear isomorphism $\mu(\phi,\cdot ):W_x\rightarrow W_y$ which is compatible with the composition of morphisms.  As we have seen Section \ref{SST}, $W$ can be viewed 
as the object space of a category $(W, {\bm{W}})$ whose morphisms space is ${\bm{W}}:={\bm{X}}{_{s}\times_P W}$. The projection 
${\bm{P}}:{\bm{W}}\rightarrow {\bm{X}}$, ${\bm{P}}(\phi,w)=\phi$,  defines a strong bundle. Moreover,  the pairs $(W[i],{\bm{W}}([i])$ for
$i=0,1$ are ep-groupoids and  the projections $(P,{\bm{P}})\colon (W[i],{\bm{W}}[i])\rightarrow (X,{\bm{X}})$ onto ep-groupoids are sc-smooth functors.

A strong bundle functor $\Phi \colon (P\colon W\to X, \mu)\rightarrow (P'\colon W'\to X', \mu')$  between strong bundles over ep-groupoids is a strong bundle map $\Phi:W\rightarrow W$
covering the  sc-smooth functor $\varphi\colon X\rightarrow X'$ between the ep-groupoids preserving the structural maps in the sense that
\begin{equation}\label{equation_section_functor_1}
\mu'(\bm{\varphi}(\gamma),\Phi(w))=\Phi(\mu(\gamma ,w))
\end{equation}
for all $(\gamma ,w)\in {\bm{W}}$. 
In the more functorial description,  a strong bundle functor is a  functor
$(\Phi,\bm{\Phi}):(W,{\bm{W}})\rightarrow (W,{\bm{W}}')$  in which  $\Phi$ and $\bm{\Phi}$ are strong bundle maps covering the sc-functor $\varphi\colon X\to X'$.

\begin{definition}\label{LinStBunEq-def}\index{D- Strong bundle equivalence}
A {\bf strong bundle equivalence} 
$$
\Phi \colon (P\colon W\to X, \mu)\rightarrow (P'\colon W'\to X', \mu')
$$
  between two  strong bundles over ep-groupoids is a linear strong bundle 
functor between strong bundles over ep-group\-oids satisfying the following properties:
\begin{itemize}
\item[(1)]\ The functors $\Phi \colon W[i]\to W'[i]$ for $i=0,1$ are  equivalences between  ep-group\-oids, covering the equivalence $\varphi \colon X\to X'$ between the underlying  ep-groupoids.
\item[(2)]\  The induced maps $\Phi \colon W\to W'$ and $\bm{\Phi} \colon {\bm{W}}\to {\bm{W}}'$ between the object sets and the morphism sets preserve the strong bundle structures and are locally strong bundle isomorphisms.
\end{itemize}
\qed
\end{definition}
If $(P\colon W\rightarrow X,\mu)$ is a strong bundle over the  ep-groupoid $X$,   we consider  $w\in W$ as an object in the ep-groupoid $(W, {\bm{W}})$. By definition,  
 $w$ and $w'\in W$ are isomorphic if there exists a morphism $\sigma \colon P(w)\to P(w')\in {\bm{X}}$ satisfying 
\begin{equation}\label{strong_bundle_equ_eq1}
\mu(\sigma,w)=w'.
\end{equation}
By $\abs{W}$ we denote, as usual, the  orbit space of $W$, whose elements are the  equivalence classes $|w|=\{w'\, \vert \,  w'\,  \text{is isomorphic to}\,   w\}$.
A  strong bundle map $\Phi \colon (P\colon W\to X, \mu)\rightarrow (P'\colon W'\to X', \mu')$ induces the continuous map  $|\Phi|\colon |W|\rightarrow |W'|$ between the orbit spaces.

Finally we recall from Section \ref{SST} that a  sc-smooth section functor $f\colon X\to W$ 
of the strong bundle $(P\colon W\to X, \mu)$
satisfies 
\begin{equation}\label{property_eq2}
f(t(\gamma))=\mu(\gamma,f(s(\gamma))).
\end{equation}
for every morphism $\gamma\in {\bm{X}}$. This relation 
encapsulates the functoriality of the section $f$. 
If  $\gamma\in {\bm{X}}$ belongs to the isotropy group $G_x$ at $x\in X$, then 
\begin{equation}\label{property_eq3}
f(x)=\mu(\gamma,f(x)).
\end{equation}
This relation has  an  interesting consequence for the pull-backs $\Phi^\ast(f')$ of section functors $f'$ of the bundle $W'\to X'$, defined at $x\in X$ by the formula
$$
\Phi^\ast(f')(x) =\Phi^{-1}\circ (f' (\varphi (x)))\in W_x,
$$
provided $\Phi$ is a strong bundle equivalence.
\begin{proposition}\label{Strong_bundle_equivalences_I}\index{P- Strong bundle equivalences {I}}
We assume that 
$\Phi \colon (P\colon W\to X, \mu)\rightarrow (P'\colon W'\to X', \mu')$
 is a strong bundle equivalence covering the equivalence $\varphi\colon X\to X'$ between ep-groupoids,  and $\Psi \colon (P\colon W\to X, \mu)\rightarrow (P'\colon W'\to X', \mu')$  is a 
strong bundle equivalence covering the equivalence $\psi\colon X\to X'$. 
If $|\Psi|=|\Phi|$, then 
$$
\Phi^\ast(f')=\Psi^\ast(f')
$$
for every sc-smooth section functor $f'\colon W'\to X'$.
\end{proposition}
\begin{proof}
The section functor 
$f'\colon X'\to W'$ satisfies, by definition, 
$$
f'(\psi (x)) = \mu'(\tau', f'(\varphi (x))) 
$$
for every morphism $\tau'\colon \varphi (x)\to \psi (x)$ in $\bx '$. 
Abbreviating
$$
\text{$g=\Psi^\ast (f')$\quad and \quad $f=\Phi^\ast (f')$},$$
we deduce 
\begin{equation}\label{pull_back_eq4}
\Psi (g(x))=f'(\psi (x))=\mu'(\tau', f'(\varphi (x)))=\mu'(\tau', \Phi (f(x)))
\end{equation}
for every morphism $\tau'\colon \varphi (x)\rightarrow \psi (x)$. In view of the assumption, 
that $\abs{\Phi}=\abs{\Psi}$  there exists  a morphism $\Sigma\colon \Phi (g(x))\to \Psi (g(x))$ in 
${\bf W'}$. The vector $\Phi (g(x))$ belongs to the fiber $W'_{\varphi (x)}$ and the vector 
$\Psi (g(x))$ belongs to the fiber $W_{\psi (x)}'$. Therefore, in view of \eqref{strong_bundle_equ_eq1}, there exists a morphism $\sigma\colon \varphi (x)\to \psi (x)$ in ${\bm{X}}'$ such that 
\begin{equation}\label{pull_back_eq5}
\Psi (g(x))=\mu'(\sigma, \Phi (g(x)).
\end{equation}
From \eqref{pull_back_eq4} and \eqref{pull_back_eq5} we deduce the relation 
$$
\mu'(\tau', \Phi (f(x)))=\mu'(\sigma, \Phi (g(x))),
$$
which implies, using the properties  of the structure map $\mu$ is Section \ref{SST},
\begin{equation}\label{pull_back_eq6}
\Phi (f(x))=\mu'((\tau')^{-1}\circ \sigma, \Phi (g(x))).
\end{equation}
We note that $(\tau')^{-1}\circ \sigma\colon \varphi (x)\to \varphi (x)$ in ${\bm{X}}'$. 
Since $\Phi$ is a strong bundle equivalence, which covers the equivalence  
$\varphi\colon X\to X'$ between ep-groupoids, there exists a morphism $\eta\colon x\to x$ in ${\bm{X}}$ satisfying 
$$
 \varphi(\eta)=(\tau')^{-1}\circ \sigma\in {\bm{X}}'.
 $$
The strong bundle equivalence $\Phi$ preserves the structure maps. Therefore, 
$$
\mu'(\varphi(\eta), \Phi (g(x)))=\Phi (\mu (\eta, g(x))).
$$
Combining this with \eqref{pull_back_eq6} leads to 
\begin{equation}\label{edi_zehnder}
\Phi(f(x))=\Phi (\mu (\eta, g(x))).
\end{equation}
Since $g$ is a section functor and $\eta\colon x\to x$ belongs to $G_x$, we obtain, in view of \eqref{property_eq3}, 
$g(x)=\mu (\eta, g(x))$, and hence conclude from the previous relation that 
$$\Phi(f(x))=\Phi(g(x)).$$
Now, $\Phi$ is fiberwise an isomorphism implying that $g(x)=f(x)$,  as we wanted to prove. 
\qed \end{proof}

We shall prove the analogous result for push-forwards.  For  a strong bundle equivalence
$$
\Phi\colon (P\colon W\to X, \mu)\rightarrow (P'\colon W'\to X',\mu')
$$
as in Proposition \ref{Strong_bundle_equivalences_I},   we define the push-forward 
$$\Phi_\ast:\Gamma(P,\mu)\rightarrow \Gamma(P',\mu')$$
of sc-smooth section functors as follows. By definition, the equivalence $\Phi$ is essentially surjective. Hence fixing 
$x'\in X'$ we find  $x\in X$ and a morphism $\sigma:\varphi(x)\rightarrow x'$ in ${\bm{X}}'$ and define the push-forward $\Phi_\ast (f)\colon X'\to W'$ of the section functor $f\colon X\to W$  of the bundle $P\colon W\to X$ at the point $x'\in X'$ by 
$$
(\Phi_\ast f)(x') = \mu'(\sigma,\Phi(f(x))).
$$
Let us first make sure that this is well-defined. If $y\in X$ and $\sigma':\varphi (y)\rightarrow x'$, then 
$\sigma^{-1}\circ\sigma':\varphi(y)\rightarrow \varphi(x)$ and since $\varphi$ is an equivalence between ep-groupoids 
we find a morphism $\tau:y\rightarrow x$ satisfying  ${\bls \varphi}(\tau)=\sigma^{-1}\circ\sigma'$, so that $\sigma\circ {\bls \varphi}(\tau)=\sigma'$.
We compute using the functoriality of $\Phi$ and the properties of the structure map $\mu'$, 
\begin{equation}\label{equation_independent}
\begin{split}
\mu'(\sigma',\Phi(f(y)))
&=\mu'(\sigma\circ{\bls \varphi}(\tau),\Phi(f(y)))\\
&=\mu'(\sigma,\mu'({\bls \varphi}(\tau),\Phi(f(y)))\\
&=\mu'(\sigma,\Phi(\mu(\tau,f(y))))\\
&=\mu'(\sigma,\Phi(f(x))).
\end{split}
\end{equation}
This proves that our definition is independent of the choices. The sc-smoothness
of $\Phi_\ast f$ is easily verified. 

\begin{proposition}\label{Strong_bundle_equivalences_II}\index{P- Strong bundle equivalences {II}}
If  $\Phi$ and $\Psi$ are the strong bundle equivalences of Proposition \ref{Strong_bundle_equivalences_I} satisfying 
$|\Phi|=|\Psi|$, 
then 
$$\Phi_\ast f=\Psi_\ast f$$
for every $f\in \Gamma (P, \mu)$.
\end{proposition}
\begin{proof}
Let $f\in \Gamma (P, \mu)$ and $x'\in X'$. We take $x\in X$ and a morphism $\sigma:\varphi(x)\rightarrow x'\in {\bm{X}}'$. From $|\Psi|=|\Phi|\colon \abs{W}\to \abs{W'}$  we conclude that $|\varphi|=|\psi|\colon \abs{X}\to \abs{X'}$ implying the existence of a morphism $\sigma':\psi(x)\rightarrow x'$ in ${\bm{X}}'$. Then, 
by definition, 
$$
\text{$(\Phi_\ast f)(x') = \mu'(\sigma,\Phi(f(x)))$\quad  and\quad 
$(\Psi_\ast f)(x') = \mu'(\sigma',\Psi(f(x)))$}.
$$
Since $\abs{\Phi}=\abs{\Psi}$, we conclude the 
existence of  a morphism $\Phi(f(x))\rightarrow \Psi(f(x))$ in ${\bm{W}}'$.  Hence, arguing as in Proposition \ref{Strong_bundle_equivalences_I}, there exists a morphism 
$\tau':\varphi(x)\rightarrow \psi(x)$ in ${\bm{X}}'$ such that  $\Psi(f(x))=\mu'(\tau',\Phi(f(x)))$. Consequently, 
$$
(\Psi_\ast f)(x')=\mu'(\sigma',\mu'(\tau',\Phi(f(x))))=\mu'(\sigma'\circ\tau',\Phi(f(x))).
$$
Then $\sigma ^{-1}\circ \sigma'\circ \tau'\colon \varphi(x)\rightarrow \varphi(x)$ and,  since $\varphi$ is an equivalence between ep-groupoids, 
we find a morphism $\eta\colon x\rightarrow x$ in ${\bm{X}}$ satisfying  ${\bls \varphi}(\eta)=\sigma ^{-1}\circ \sigma'\circ \tau'$ implying 
$$
\sigma\circ {\bls \varphi}(\eta)=\sigma'\circ \tau'.
$$
Hence,  
\begin{equation*}
\begin{split}
(\Psi_\ast f)(x')
&=\mu'(\sigma\circ {\bls \varphi} (\eta),\Phi(f(x)))\\
&=\mu'(\sigma,\mu'({\bls \varphi} (\eta),\Phi(f(x))))\\
&=\mu'(\sigma,\Phi(\mu(\eta,f(x))))\\
&=\mu'(\sigma,\Phi(f(x)))\\
&=(\Phi_\ast f)(x').
\end{split}
\end{equation*}
The proof of Proposition \ref{Strong_bundle_equivalences_II} is complete.
\qed \end{proof}
\begin{proposition}\label{composition_push_pull}
If $\Phi\colon (P,\mu)\to (P', \mu')$ and $\Psi\colon (P',\mu')\to (P'',\mu'')$ are strong bundle equivalences, then 
$$(\Psi\circ \Phi)^\ast=\Phi^\ast \circ \Psi^\ast\quad \text{and}\quad (\Psi\circ \Phi)_\ast=\Psi_\ast \circ \Phi_\ast.$$
\end{proposition}
\begin{proof}
We assume that the equivalence $\Phi\colon ((P\colon W\to X),  \mu)\to (P'\colon W'\to X'), \mu')$ covers the equivalence $\varphi\colon X\to X'$ between ep-groupoids and the equivalence $\Psi \colon ((P'\colon W'\to X'),  \mu)\to (P''\colon W''\to X''), \mu')$ covers the equivalence $\psi\colon X'\to X''$ between ep-groupoids. If $f''\colon X''\to W''$ is a section functor and $x\in X$, then, by definition,
\begin{equation*}
\begin{split}
((\Psi\circ \Phi)^\ast f'')(x)&=(\Psi\circ \Phi)^{-1}(f''(\psi \circ \varphi (x))\\
&=\Phi^{-1}\circ \Psi^{-1}(f''(\psi (\varphi (x)))\\
&=\Phi^{-1}(\Psi^\ast f'' (\varphi (x))=(\Phi^\ast\circ \Psi^\ast)f''(x).
\end{split}
\end{equation*}
Next we consider push-forwards.  We fix  $x''\in X''$. Using the fact that  the composition of equivalences is again an equivalence, we choose $x\in X$ and a morphism $\sigma''\colon \psi \circ \varphi (x)\to x''$ in ${\bf X'}$. Then, by definition,  
$$(\Psi\circ \Phi)_\ast (f)(x'')=\mu''(\sigma'' , \Psi\circ \Phi (f(x))).$$
Moreover, again using that  $\phi$  and $\psi$ are equivalences, we find a point  $x'\in X'$ and a morphism $\sigma'\colon \psi (x')\to x''$ in ${\bm{X}}''$ and a $y\in X$ and a morphism $\sigma\colon \varphi (y)\to x'.$
Then 
$$(\Phi_\ast f)(x')=\mu'(\sigma, \Phi (f(y)))$$
and, from \eqref{equation_section_functor_1}  
and the properties of the structure map $\mu$, we deduce
\begin{equation*}
\begin{split}
\bigl(\Psi_\ast (\Phi_\ast f)\bigr)(x'')&=\mu''\bigl(\sigma', \Psi ((\Phi_\ast f)(x')\bigr)\\
&=\mu''(\sigma', \Psi (\mu'(\sigma , \Phi (f(y)))\\
&=\mu''(\sigma',\mu''( {\bls \psi}(\sigma) ,\Psi( \Phi (f(y))))\\
&=\mu''(\sigma'\circ {\bls \psi}(\sigma),\Psi \circ \Phi (f(y)).
\end{split}
\end{equation*}
Since $\sigma\colon \varphi (y)\to x'$, we have ${\bls \psi}(\sigma)\colon \psi \circ \varphi (y)\to \psi (x')$,  and $\sigma'\colon \psi (x')\to x''$, the composition $\sigma'\circ {\bls \psi}(\sigma)$ is a morphism between $\psi \circ \varphi (y)$ and $x''$. 
We conclude, in view of the independence of the choices involved in the definition of push-forward, \eqref{equation_independent},
that 
$(\Psi\circ \Phi)_\ast (f)(x'')=\bigl(\Psi_\ast (\Phi_\ast f)\bigr)(x'')$, as claimed.
\qed \end{proof}

It follows from the definitions that  $\Phi_\ast$ and $\Phi^\ast$ are mutual inverses. The transformations
$\Phi_\ast:\Gamma(P,\mu)\rightarrow \Gamma(P',\mu')$ and $\Phi^\ast:\Gamma(P',\mu')\rightarrow \Gamma(P,\mu)$
map the subspaces of $\ssc^+$-section functors in one bundle to the corresponding ones in the other bundle.
If the  sc-smooth section functor $f$ of $(P,\mu)$, just viewed as a sc-smooth section
of the strong bundle $W\rightarrow X $ over the object M-polyfold,  is sc-Fredholm, then the same is true for $\Phi_\ast f$ and $\Phi^\ast f$. 
This  follows immediately from the definition of a sc-Fredholm section  in Part I.

We  summarize the above results as follows.
\begin{theorem}\label{useful}\index{T- Strong bundle equivalences and sections}
If $(P:W\rightarrow X,\mu)$ and $(P':W'\rightarrow X',\mu')$ are strong bundles over ep-groupoids and
$\Phi:W\rightarrow W'$ is a strong bundle equivalence, then the pull-back $\Phi^\ast$ and the push-forward $\Phi_\ast$
define mutually inverse bijections
$$
\Phi_\ast :\Gamma(P,\mu)\rightarrow \Gamma(P',\mu')\quad  \text{and}\quad \Phi^\ast:\Gamma(P',\mu')\rightarrow \Gamma(P,\mu).
$$
If  $\Phi$ and $\Psi\colon W\rightarrow W'$ are two equivalences of strong bundles satisfying $|\Phi|=|\Psi|$, then 
$\Phi_\ast=\Psi_\ast$ and $\Phi^\ast=\Psi^\ast$. Moreover,  $\Phi_\ast$ and $\Phi^\ast$ preserve $\ssc^+$-sections and 
sc-Fredholm sections, so that 
$$
\Phi_\ast:\Gamma^+(P,\mu)\rightarrow \Gamma^+(P',\mu')\ \text{and}\ \Phi^\ast:\Gamma(P',\mu')\rightarrow \Gamma(P,\mu)
$$
and
$$
\Phi_\ast:\textrm{Fred}(P,\mu)\rightarrow \textrm{Fred}(P',\mu')\ \text{and}\ \ \Phi^\ast:\text{Fred}(P',\mu')\rightarrow \text{Fred}(P,\mu).
$$
\qed
\end{theorem}

In order to generalize the notion of natural equivalence to strong bundle maps we 
 consider the two strong bundles 
 $$
 \text{$(P \colon W\to X,\mu)$\quad  and \quad  $(P' \colon W'\to X',\mu')$}
 $$
 over ep-groupoids and the two strong  bundle functors
$$
\Phi, \Psi \colon (P,\mu)\to (P',\mu').
$$
\begin{definition}\label{NatEquiv-def}\index{D- Natural equivalence}
The 
strong  bundle functors $\Phi$ and $\Psi$  are called {\bf naturally equivalent}, $\Phi\simeq \Psi$, if there exists a natural transformation
$$
T \colon W\to {\bm{W}}'
$$
associating with every object $w\in W$ a morphism $\tau (w)\colon \Phi (w)\to \Psi (w)$ in 
${\bm{W}}'$ satisfying 
$$
T (w')\circ {\bf \Phi}(h)={\bf \Psi}(h)\circ T (w)
$$
for every morphism $h\colon w\to w'$ in ${\bm{W}}$. In addition, there exists a natural transformation 
$$
\tau\colon X\to {\bm{X}}'
$$
of the two underlying functors $\varphi$ and $\psi\colon X\to X'$ between ep-groupoids such that 
$$
\tau\circ P(w)=P'\circ T(w),\quad w\in W.
$$
where $P\colon W\to X$ and $P'\colon W'\to X'$ are the projection functors. 
\qed
\end{definition}
\begin{remark}\index{R- Sc-smoothness of natural transformations}
The natural transformation $\tau$ between the sc-smooth functors $\varphi$ and $\psi$ is sc-smooth by definition,
see Definition \ref{natural_equivalence}. We seemingly did not spell out any smoothness assertion about $T$.
However, as we shall see, $T$ is completely detemined by $\tau$ and automatically has the right sc-smoothness properties.
\qed
\end{remark}
Recall that a strong bundle over an ep-groupoid $(P,\mu)$ consists of a strong bundle $P:W\rightarrow X$ over the object 
space and a lift of the morphisms to strong linear bundle maps. In particular  $\bm{W}$ is a natural construction 
associated to $(P,\mu)$. It consists of the pairs $(\phi,w)\in \bm{X}\times W$ with $s(\phi)=P(w)$.
It follows from Definition \ref{NatEquiv-def} that the natural transformation 
$$
T \colon W\rightarrow {\bm{W}}' := {\bm{X}}'{{_s}\times_{P'}}W',
$$
covering $\tau \colon X\rightarrow {\bm{X}}'$ has the form
$$
T(w)=(\tau(P(w)),A(w)).
$$
In view of the definition of the morphisms in ${\bm{W}}'$ we obtain $\Phi(w)=s(T(w))=A(w)$ and $\Psi(w)=t(T(w))=\mu(\tau(P(w)),A(w))$, so that 
\begin{align}\label{Kalign}
T(w)&=(\tau(P(w)),\Phi(w))\\
\Psi(w)&=\mu(\tau(P(w)),\Phi(w)).\nonumber
\end{align}
We see that the strong bundle functor $\Psi$ and the natural transformation $T$ are uniquely determined by the strong bundle functor $\Phi$ and the natural transformation $\tau$. We see from (\ref{Kalign}) and the sc-smoothness properties of $P$, $\Phi$ and $\tau$, that $T$ is necessarily an sc-smoothness strong bundle map.
This allows us to define a natural transformation  between the two strong bundle functors equivalently as follows. 
\begin{definition}\label{NatTranBM}\index{D- Natural transformation}
A {\bf natural transformation} $T$ between the strong bundle functors $\Phi \colon P\to P'$ and $\Psi \colon P\to P'$
(covering the underlying functors $\varphi \colon X\to X'$  and $\psi \colon X\to X'$)
is a strong bundle map
$$
T \colon W\rightarrow {\bm{X}}'{{_s}\times_{P'}}W'
$$
covering a natural transformation $\tau \colon X\rightarrow {\bm{X}}'$ between  the functors $\varphi$ and $\psi$ of the form
$$
T(w)=(\tau(P(w)),\Phi(w)), 
$$
and satisfying
$$
\Psi(w)=\mu(\tau(P(w)),\Phi(w)).
$$
\qed
\end{definition}
From the natural equivalence $\Phi\simeq \Psi$ we deduce that their induced maps 
$$
\abs{\Phi},\abs{\Psi}\colon \abs{W}\to \abs{W'}
$$
 between the orbit spaces agree,
$$\abs{\Phi}=\abs{\Psi}.$$
At this point we can introduce a category whose objects are strong bundles over ep-groupoids.
\begin{definition}[{\bf The category ${\mathcal{SEP}}$}] \index{D- Category of strong bundles}
The category  ${\mathcal{SEP}}$\index{${\mathcal{SEP}}$} has as objects the 
strong bundles over ep-groupoids. The morphisms are the strong bundle  functors.
The morphism set has a distinguished subset $ {\bf F}$\index{${\bf F}$} of strong bundle equivalences.
\qed
\end{definition}

\section{Localization in the  Strong Bundle Case} 
Next we carry out the previous localization procedure but now in the context of strong bundle equivalences.
As in the case of  the category $\mathcal{E}\mathcal{P}$
we have a distinguished class of morphisms which are the strong bundle equivalences.
The preparatory material for this subsection can be found in the Sections \ref{SST} and \ref{SST_0}.
The only new thing is that we need to discuss the behavior of sections under generalized strong bundle isomorphisms.
There will be  the distinguished class ${\bf F}$\index{${\bf F}$} \index{Class of strong bundle equivalences} of strong bundle equivalences and,  following the procedure carried out
for ep-groupoids,  we shall define the  new category ${\mathcal S}{\mathcal E}{\mathcal P}({\bf F}^{-1})$.\index{${\mathcal S}{\mathcal E}{\mathcal P}({\bf F}^{-1})$}

  We consider two strong bundles
$(P \colon W\to X,\mu)$ and $(P' \colon W'\to X',\mu')$ over ep-groupoids and study the
diagrams
$$
W\xleftarrow{\Phi} W''\xrightarrow{\Psi}W',
$$
in which  $P'' \colon W''\to X''$ is a  third strong bundle over an ep-groupoid, and where $\Phi$ is a strong bundle equivalence and $\Psi$ a sc-smooth strong bundle functor.  Following our  earlier construction of a generalized map between ep-groupoids, one introduces the notion of a refinement of a diagram and then calls two diagrams equivalent, if they possess a common refinement. The equivalence class of the  diagram $D\colon W\xleftarrow{\Phi} W''\xrightarrow{\Psi}W',$ is then, by definition, a {\bf generalized strong bundle map},\index{Generalized strong bundle map}\index{$D\colon W\xleftarrow{\Phi} W''\xrightarrow{\Psi}W'$}  and denoted by
$$
[D] \colon W\rightarrow W'.\index{$[D] \colon W\rightarrow W'$}
$$
A diagram $D:W\xleftarrow{\Phi}W''\xrightarrow{\Psi}W'$ always covers a diagram $d\colon X\xleftarrow{\phi} X''\xrightarrow{\psi} X'$ of the underlying ep-groupoids. Hence the generalized strong bundle map  $[D]$ covers the generalized map 
$$
[d] \colon X\rightarrow X'
$$
between ep-groupoids. 
{The generalized strong bundle map $[D]=[W\xleftarrow{\Phi}W''\xrightarrow{\Psi}W']$ induces the map 
$$\abs{[D]}\colon \abs{W}\to \abs{W'}$$
between the orbit spaces, defined by 
$\abs{[D]}=\abs{\Psi}\circ \abs{\Phi}^{-1}$. The definition does not depend on the choice of the representative diagram for 
$[D]$.  
\begin{remark}\index{R- On the interpretation of certain diagrams}
Given the previous discussion it is tempting to introduce the diagram $\mathscr{D}$
\begin{eqnarray}\label{EQN1051}
{\mathscr{D}}\colon\ \ \ \ \ \ &
\begin{CD}
W @>[D]>>  W'\\
@V P VV   @V P' VV\\
X@>[d]>> X
\end{CD}
\end{eqnarray}
which literally, however,  does not make too much sense, since compositions are not defined.  Nevertheless, we might view (\ref{EQN1051}) 
as a {\bf symbolic diagram} and this is the view point we shall take.
Clearly passing to orbit spaces it becomes a commutative diagram $|\mathscr{D}|$ for continuous maps
\begin{eqnarray*}
|{\mathscr{D}}|\colon\ \ \ \ \ \ &
\begin{CD}
\abs{W} @>\abs{[D]}>> \abs{ W'}\\
@V\abs{P} VV   @V \abs{P'} VV\\
\abs{X}@>\abs{[d]}>> \abs{X}
\end{CD}
\end{eqnarray*}
There are  additional local properties of such diagrams, which we shall discuss below.
\qed
\end{remark} 
There is another property of the symbolic diagram. Namely pick $x\in X$ and $x'\in X'$ such that $|[d]|(|x|)=|x'|$. 
Taking a representative $d  \colon X\xleftarrow{\phi} X''  \xrightarrow{\psi} X'$  of $[d]$, we find $x''$ and morphisms 
$$
x\xleftarrow{a} \phi(x'') \ \text{and}\ \  \psi(x'')\xrightarrow{b} x'.
$$
These diagrams define a group homomorphism 
$$
\gamma\colon G_{x}\rightarrow G_{x'}'\colon g\rightarrow \gamma (g),
$$
obtained as a composition of several homomorphisms.
\begin{itemize}
\item $g\rightarrow h =  a\circ g \circ a^{-1}$\ \ (group isomorphism).
\item  $h\rightarrow k=\phi^{-1}(h)$\ \  (group isomorphism).
\item  $k\rightarrow \ell= \psi(k)$ \ \ (group homomorphism).
\item $\ell\rightarrow b\circ \ell \circ b^{-1}$ \ \ (group isomorphism).
\end{itemize}
We also find  suitable open neighborhoods $U(x)$, $U(x'')$, $U(x')$, $U(\phi(x''))$, and $U(\psi(x'')$, and a germ of equivariant sc-diffeomorphisms $q:{\mathcal O}(X,x)\rightarrow {\mathcal O}(X',x')$
such that
$$
q(g\ast z) = \gamma(g)\ast q(z).
$$
The germ $q$ is obtained 
from the following diagram where all arrows are sc-diffeomorphisms with the exception of $\psi$
are sc-diffeomorphisms
$$
U(x)\xleftarrow{t} U(a)\xrightarrow{s} U(\phi(x''))\xleftarrow{\phi} U(x'')\xrightarrow{\psi} U(\psi(x''))\xleftarrow{s} U(b)\xrightarrow{t} U(x')
$$
Associated to the germ $q$ there is a natural lift $Q:W|{\mathcal O}(X,x)\rightarrow W|{\mathcal O}(X',x')$ which is a strong bundle map fitting into the commutative diagram
of equivariant maps
$$
\begin{CD}
W|{\mathcal O}(X,x)@> Q>> W|{\mathcal O}(X',x')\\
@V P VV @V P' VV\\
{\mathcal O}(X,x) @> q>> {\mathcal O}(X',x')
\end{CD}
$$
 These are the local expressions 
induced by $\mathscr{D}$ and these are compatible with the projections $P$ and $P'$ and define the  strong bundle structures.

\begin{definition}\index{D- Category  ${\mathcal S}{\mathcal B}({\bf F}^{-1})$}
The {\bf category ${\mathcal S}{\mathcal E}{\mathcal P}({\bf F}^{-1})$}  has the strong bundles  over ep-group\-oids as objects and the generalized 
strong bundles maps as morphisms.
\qed
\end{definition}
 As in the ep-groupoid case a generalized strong bundle isomorphism has a particular form.
\begin{definition}\index{D- Generalized strong bundle isomorphism}
A {\bf generalized strong bundle isomorphism}
 $[D]\colon W\rightarrow W'$ is an equivalence class 
 $$
 [D]=[W\xleftarrow{\Phi} W''\xrightarrow{\Psi} W']
 $$
 of diagrams  $W\xleftarrow{\Phi} W''\xrightarrow{\Psi} W'$ in which  both $\Phi$ and $\Psi$ are strong bundle equivalences. 
 \qed
\end{definition}
So the objects are morphisms can be viewed as the following diagrams.
$$
\begin{array}{cc}
\begin{CD} 
W\\
@V P VV\\
X
\end{CD}
&\ \ \ \ \ \ \ \ \ 
\begin{CD}
W @>>[D]>  W'\\
@V P VV   @V P' VV\\
X@>[d]>> X
\end{CD}
\end{array}
$$
or in short form $[D]:W\rightarrow W'$.

Theorem \ref{useful} allows to push-forward or pull-back sc-smooth section functors by  strong bundle
equivalences.   We assume that we have a generalized strong bundle isomorphism
$$
[D]:(P,\mu)\rightarrow (P',\mu')
$$
and take a representative $W\xleftarrow{\Phi} W''\xrightarrow{\Psi} W'$ in which  both maps are strong bundle equivalences.
For a sc-smooth section functor $f$ of the strong bundle $(P,\mu)$,  we define  the push-forward by 
$$
[D]_\ast f = \Psi_\ast\Phi^\ast f.
$$
This definition does not depend on the choice of the diagram representing  $[D]$. To see this we  take a second diagram
$W\xleftarrow{\Phi'} W'''\xrightarrow{\Psi'} W'$ refining the previous one and let $H:W'''\rightarrow W''$ be the refining equivalence.
Then $\Phi\circ H\simeq \Phi'$ and $\Psi\circ H\simeq \Psi'$. Using the 
Propositions \ref{Strong_bundle_equivalences_I}-\ref{composition_push_pull} we obtain,  
\begin{equation*}
\begin{split}
\Psi_\ast\Phi^\ast f
&=\Psi_\ast H_\ast H^\ast \Phi^\ast f\\
&=(\Psi\circ H)_\ast (\Phi\circ H)^\ast f\\
&=\Psi_\ast'  (\Phi')^\ast f.
\end{split}
\end{equation*}
Since our equivalence relation is  the common refinement we see that $[D]_\ast$ is well-defined. Similarly we define
$$
[D]^\ast g = \Phi_\ast \Psi^\ast g
$$
and again this is well-defined.  Again by Theorem \ref{useful},  we see that $[D]_\ast$ and $[D]^\ast$ preserve $\ssc^+$-section functors
and sc-Fredholm functors, see Definition \ref{SECTION-FUNCTORS-X} for the different classes of sc-smooth section functors. Hence we obtain the following theorem.
\begin{theorem}\label{Push-Forw-prop}\index{T- Generalized bundle isomorphisms}
Let $(P \colon W\to X,\mu)$ and $(P' \colon W'\to X',\mu')$ be two strong bundles over ep-groupoids and $[D]:W\rightarrow W'$ a generalized strong bundle isomorphism.
Denote by $\Gamma(P,\mu)$ and $\Gamma(P',\mu')$ the vector spaces of sc-smooth section functors. Then $[D]$ induces a well-defined isomorphism
$$
[D]_\ast\colon\Gamma(P,\mu)\rightarrow \Gamma(P',\mu').
$$
Its inverse is given by the inverse diagram and equals the pull-back $[D]^\ast$.  The same assertion hold for $sc^+$-section functors 
$$
[D]_\ast\colon\Gamma^+(P,\mu)\rightarrow \Gamma^+(P',\mu')\colon
$$
and sc-Fredholm sections functors
$$
[D]_\ast\colon \text{Fred}(P,\mu)\rightarrow \text{Fred}(P',\mu').
$$
Moreover, if $[D], [D']: W\rightarrow W'$ are two generalized strong bundle isomorphisms inducing the same maps between orbit spaces, i.e. $|[D]|=|[D']|$, then
$[D]_\ast =[D']_\ast$ and $[D]^\ast =[D']^\ast$.
\qed
\end{theorem}

We end this section by stating a result which generalizes Theorem \ref{THMX10316} to the strong bundle situation.
\begin{theorem}\index{T- Uniqueness result for generalized strong bundle isomorphisms}
Assume that $(P \colon W\to X,\mu)$ and $(P' \colon W'\to X',\mu')$ are  two strong bundles over ep-groupoids and $[D], [D'] :W\rightarrow W'$ are generalized strong bundle isomorphisms satisfying 
$|[D]|=|[D']|$. Then $[D]=[D']$.\qed
\end{theorem}
\begin{proof}
The proof follows essentially along the lines of the proof of Theorem \ref{THMX10316}.
\qed \end{proof}
\section{Appendix}

\subsection{Proof of Theorem \ref{proppp}}\label{x-proppp}
We recall the statement of the theorem for the convenience of the reader.\par

{\bf Theorem \ref{proppp}\textcolor{red}{.}}
Let 
$$
d\colon X\xleftarrow{F} A\xrightarrow{\Phi}Y\quad \text{and}\quad d'\colon Y\xleftarrow{G} B\xrightarrow{\Psi} Z
$$
be representatives  of the equivalence classes $[d]\colon X\rightarrow Y$ and $[d']\colon Y\rightarrow Z$, respectively.  Then the equivalence class $[d'']$ of the diagram 
$$
d''\colon X\xleftarrow{F\circ\pi_A} A\times_Y B\xrightarrow{\Psi\circ\pi_B} Z
$$
is independent of the choices of representatives in the equivalence classes $[d]$ and $[d']$. The equivalence class $[d'']:X\to Z$ is called the {\bf composition}
of $[d]$ and $[d']$ and is denoted by 
$$[d'']=[d']\circ [d].$$
Moreover, the composition is associative; if $[d]:X\rightarrow Y$, $[d']:Y\rightarrow Z$,  and $[d'']:Z\rightarrow W$, then 
$$
([d'']\circ [d'])\circ [d]=[d'']\circ ([d']\circ [d]).
$$
In addition, 
$$
[d]\circ [1_X]=[d]\quad \text{and}\quad  [1_Y]\circ [d]=[d]
$$
for every equivalence class $[d]:X\rightarrow Y$.  Here $[1_X]$ is  the equivalence class 
$$[1_X]=[X\xleftarrow{1_X} X\xrightarrow{1_X}X].$$
\par

\begin{proof}
Assume that $X\xleftarrow{F'} A'\xrightarrow{\Phi'}Y$ is a refinement of $X\xleftarrow{F} A\xrightarrow{\Phi}Y$ and $Y\xleftarrow{G'} B'\xrightarrow{\Psi'} Z$ is a refinement of 
$Y\xleftarrow{G} B\xrightarrow{\Psi} Z, 
$ so that we have the following two diagrams 
\begin{equation*}
\begin{CD}
X@<F<<A@>\Phi>>Y\\
@. @AAHA @. \\
X@<F'<<A'@>\Phi'>>Y\\ \\
\end{CD}
\qquad \text{and}\qquad 
\begin{CD}
X@<G<<B@>\Psi>>Y\\
@. @AAKA @. \\
X@<G'<<B'@>\Psi'>>Y\\ \\
\end{CD}
\end{equation*}
in which $H$ and $K$ are equivalences and $F\circ H\simeq F'$, $\Phi\circ H\simeq \Phi'$ and $G\circ K\simeq G'$, $\Psi\circ K\simeq \Psi'$. 
We claim that 
the diagram 
$$
X\xleftarrow{F'\circ\pi_{A'}} A'\times_Y B'\xrightarrow{\Psi'\circ\pi_{B'}} Z
$$ 
is a refinement of $d''$, so that 
\begin{equation*}
\begin{CD}
d'':X@<F\circ \pi_A<<A\times_YB@>\Psi\circ \pi_B>>Y\\
@. @AALA @. \\
\phantom{d'':}X@<F'\circ \pi_{A'}<<A'\times_{Y}B'@>\Psi'\circ \pi_{B'}>>Y\\ \\
\end{CD}
\end{equation*}
for an equivalence $L$ satisfying $F'\circ \pi_{A'}\simeq(F\circ \pi_A)\circ L$ and 
$\Psi'\circ \pi_{B'}\simeq(\Psi\circ \pi_B)\circ L$.

Indeed, take  $(a, \varphi, b)\in  A'\times_Y B'$. Then $\varphi$ is a morphism in ${\bm{Y}}$ between the points $\Phi' (a)$ and $G'(b)$. Since  $\Phi\circ H\simeq \Phi'$, there exists a sc-smooth map $\tau_{A'}:A'\to {\bm{Y}}$ such  that $\tau_{A'}(a)$ is a morphism between the points $\Phi'(a)$ and $(\Phi \circ H)(a)=\Phi (H(a))$ in ${\bm{Y}}$. Similarly, there  exists a sc-smooth map $\tau_{B'}:B'\to {\bm{Y}}$ so that $\tau_{B'}(b)$ is a morphism between $G'(b)$ and $(G \circ K)(b)=G(K(b))$. Consequently, $\tau_B(b)\circ \varphi \circ (\tau_{A'}(a))^{-1}$ is a morphism in ${\bm{Y}}$ between $H(a)$ and $K(b)$. So we define the map 
$$L: A'\times_Y B'\to  A\times_Y B$$
by 
$$L(a, \varphi, b)=(H(a), \tau_{B'}(b)\circ \varphi \circ (\tau_{A'}(a))^{-1}, K(b)).$$
The map $L$ is sc-smooth and 
\begin{equation*}
\begin{split}
(F\circ \pi_A)\circ L(a, \varphi, b)&=F\circ (\pi_A\circ H)(a, \varphi, b)=F\circ H(a)\\
&=F'(a)=(F'\circ \pi_{A'})(a, \varphi, b)
\end{split}
\end{equation*}
so that $(F\circ \pi_A)\circ L=F'\circ \pi_{A'}$. Similarly, $(\Phi \circ \pi_B)\circ L =\Psi'\circ \pi_{B'}$.
Since $F\circ \pi_A$ and $F'\circ \pi_{A'}$ are local sc-diffeomorphism, also $L$ is a local sc-diffeomorphism. Moreover, the induced maps between orbit spaces satisfy
 $\abs{(F\circ \pi_A)\circ L}=\abs{F\circ \pi_A}\circ \abs{L}=\abs{F'\circ \pi_{A'}}$ and since $\abs{F\circ \pi_A}$ and $\abs{F'\circ \pi_{A'}}$ are homeomorphisms, the same is true for the map $\abs{L}: \abs{A'\times_Y B'}\to  \abs{A\times_Y B}$.  Next we show that $L$ is faithful and full.  Take two objects  $(a, \varphi, b)$ and $(a', \varphi', b')$ in $A'\times_Y B'$. We have to verify that the induced map 
 $$\textrm{mor}_{A'\times_Y B'}( (a, \varphi, b), (a', \varphi', b'))\to \textrm{mor}_{A\times_YB}( L(a, \varphi, b), L(a', \varphi', b'))$$
 is a bijection. We take a morphism $(h', \psi, k')\in \textrm{mor}_{A\times_YB}( L(a, \varphi, b), L(a', \varphi', b'))$. In view of the definition of $L$, the triple $(h', \psi, k')$ is a morphism between $(H(a), \tau_{B'}(b)\circ \varphi \circ (\tau_{A'}(a))^{-1}, K(b))$ and $(H(a'),  \tau_{B'}(b')\circ \varphi' \circ (\tau_{A'}(a'))^{-1}, K(b'))$ so that  we have a following diagram 
 \begin{equation*}
\begin{CD}
H(a)@>\psi= \tau_{B'}(b)\circ \varphi \circ (\tau_{A'}(a))^{-1}>>K(b)\\
 @ Vh'VV   @VVk'V \\
H(a')@> \psi'= \tau_{B'}(b')\circ \varphi' \circ (\tau_{A'}(a'))^{-1}  >>K(b').\\
\end{CD}
\end{equation*}
Since $H$ and $K$ are equivalences, and hence faithful and full, we find unique morphisms $h:a\to a'$ in ${\bf A}'$ and $k:b\to b'$ in ${\bf B}'$ such that $H(h)=h'$ and $K(k)=k'$. Then 
$$L(h, \psi, k)=(h', \psi, k')$$
showing that $L$ is faithful and full. 

We have proved that $L$ is an equivalence and that $(F\circ \pi_A)\circ L\simeq F'\circ \pi_{A'}$ and $(\Phi \circ \pi_B)\circ L\simeq \Psi'\circ \pi_{B'}$. Consequently, the diagram 
$
X\xleftarrow{F'\circ\pi_{A'}} A'\times_Y B'\xrightarrow{\Psi'\circ\pi_{B'}} Z
$
is a refinement of $d''$ as claimed. Hence the equivalence class $[d'']$ is independent of the choice of the  representatives in the  equivalence classes $[d]$ and $[d']$.  

To see that the composition is associative, we consider three diagrams
$$d\colon X\xleftarrow{F}A\xrightarrow{\Phi}Y,\quad d'\colon Y\xleftarrow{F'}A'\xrightarrow{\Phi'}Z,\quad \text{and}\quad d''\colon Z\xleftarrow{F''}A''\xrightarrow{\Phi''}W.$$
The composition $[d'']\circ [d']$ is represented  by the diagram 
\begin{equation*}
Y\xleftarrow{F'\circ\pi_{A'} } A' \times_Z A'' \xrightarrow{\Phi''\circ\pi_{A''} } W, 
\end{equation*}
so that the composition $([d'']\circ [d'])\circ [d]$ is represented by the diagram 
\begin{equation}\label{associativity_1}
X\xleftarrow{F\circ\wt{\pi}_{A}}A\times_Y (A'\times_Z A'')\xrightarrow{(\Phi''\circ\pi_{A''})\circ \wt{\pi}_{A'\times_ZA''} } W
\end{equation}
where we have denoted by $\wt{\pi}$ the projection  onto the corresponding factor of the weak fibered product $A\times_Y (A'\times_Z A'')$. 
Similarly,  the composition $[d']\circ [d]$ is represented  by the diagram 
$$
X\xleftarrow{F\circ\pi_{A}} A\times_Y A'\xrightarrow{\Phi''\circ\pi_{A'}} Z
$$
and  the composition $[d'']\circ ( [d'])\circ [d])$ is represented by the diagram 
\begin{equation}\label{associativity_2}
X\xleftarrow{(F\circ {\pi}_{A})\circ \wt{\pi}_{A\times_YA'} } (A\times_Y A')\times_Z A''\xrightarrow{(\Phi''\circ \wt{\pi}_{A''} } W.
\end{equation}
where this time $\wt{\pi}$ denotes  projection onto the corresponding  factor of the weak fibered product $(A\times_Y A')\times_Z A''$. 
Note that $(A\times_Y A')\times_Z A''$ is the sub  M-polyfold of $A\times{\bm{Y}}\times A'\times {\bm{Z}}\times A''$
consisting of points $(a,\phi,a',\psi,a'')$ in which  $\phi:\Phi(a)\to F'(a')$ and $\psi:\Phi' (a')\to F'' (a'')$ are morphisms.
Moreover, $A\times_Y (A'\times_Z A'')$ is the same space.   Then 
\begin{equation*}
\begin{aligned}
&(F\circ\wt{\pi}_{A} )(a,\phi,a',\psi,a'')=F(a')\\
&(F\circ {\pi}_{A})\circ \wt{\pi}_{A\times YA'} (a,\phi,a',\psi,a'')= (F\circ {\pi}_{A})(a, \phi, a')=F(a)
\end{aligned}
\end{equation*}
and
\begin{equation*}
\begin{aligned}
&(\Phi''\circ\pi_{A''})\circ \wt{\pi}_{A'\times_ZA''}   (a, \phi, a', \psi, a'')\Phi''\circ\pi_{A''}(a', \psi, a'')=\Phi''(a'')\\
&(\Phi''\circ \wt{\pi}_{A''}) (a,\phi,a',\psi,a'')=\Phi''(a'').
\end{aligned}
\end{equation*}
Hence,  taking the identity map $i_{(A\times_Y A')\times_Z A''}:(A\times_Y A')\times_Z A''\to A\times_Y (A'\times_Z A'')$ as the equivalence, the diagrams \eqref{associativity_1} and \eqref{associativity_2} are refinements of each other which implies  that $([d'']\circ [d'])\circ [d]=[d'']\circ ( [d']\circ [d])$ as claimed.

Next we shall show that $[1_X]\circ [d]=[d]$ where $d\colon X\xleftarrow{F}A\xrightarrow{\Phi} Y$. To see this it suffices to show that the diagram $X\xleftarrow{\pi_X}X\times_XA\xrightarrow{\Phi\circ \pi_A}$ is a refinement of the diagram $X\xleftarrow{F}A\xrightarrow{\Phi} Y$. The map $\pi_A\colon X\times_XA\to A$ is an equivalence and we only  have to show that $F\circ \pi_A\simeq \pi_X$.  We define $\tau:X\times_XA\to {\bm{X}}$  by $\tau (x, \varphi, a)=\varphi$. Clearly the map $\tau$ is sc-smooth and for every $ (x, \varphi, a)\in X\times_XA$, $\tau (a, \varphi, a)=\varphi$ is a morphism between $x$ and $(F\circ \pi_A) (x, \varphi, a)=F(a).$ If $(h, \psi, k)$ is a morphism in $X\times_XA$ between 
$(x, \varphi, a)$ and $(x', \varphi', a')$, then $h:x\to x'$, $k:F(a)\to F(a')$, and $\psi=\phi$ and $\varphi'=F(k)\circ \psi\circ h^{-1}=F(k)\circ \varphi\circ h^{-1}$.  Consequently, 
$$\tau (x', \varphi', a')\circ \pi_X(h, \psi, k)=(F\circ \pi_A)(h, \psi, k)\circ \tau (x, \varphi, a)$$
showing that  the map $\tau$  is natural  and that $F\circ \pi_A$ and $\pi_X$ are naturally equivalent. Hence $[1_X]\circ [d]=[d]$. Similarly, one shows that $[d]\circ [1_Y]=[d]$. The proof of Theorem \ref{proppp}  is complete.
\qed \end{proof}

\subsection{Proof of Theorem \ref{strong-iso}}\label{x-strong-iso}
We recall the theorem for the convenience of the reader.\par

{\bf Theorem \ref{strong-iso}\textcolor{red}{.}}
Let $X$, $Y$, and $A$ be ep-groupoids.
The equivalence class $[a]=[X\xleftarrow{F}A\xrightarrow{G}Y]$ in which $F$ is an equivalence and $G$ is a sc-smooth functor  is a generalized isomorphism if and only if the functor $G$ is an  equivalence. 

The inverse of a generalized isomorphism $[a]=[X\xleftarrow{F}A\xrightarrow{G}Y]$ is the equivalence class 
$$[a]^{-1}=[Y\xleftarrow{G}A\xrightarrow{F}X].
$$
\par

\begin{proof}
(1) We first verify that the generalized map $[a]\colon X\to Y$ given by the equivalence class $[a]=[X\xleftarrow{F}A\xrightarrow{G}Y]$ in which $F$ and $G$ are equivalences is a generalized isomorphism whose inverse is equal to $[a]^{-1}=[Y\xleftarrow{G}A\xrightarrow{F}X]$.

Indeed, the composition 
$$[a]^{-1}\circ [a]= [Y\xleftarrow{G}A\xrightarrow{F}X]\circ [X\xleftarrow{F}A\xrightarrow{G}Y]$$
is represented by the diagram 
$X\xleftarrow{F\circ \pi_1}A\times_YA\xrightarrow{F\circ \pi_2}X$ which is refined by the diagram $X\xleftarrow{F}A\xrightarrow{F}X$ by means of the equivalence $A\to A\times_YA$ defined by $a\mapsto (a, a)$. The diagram $X\xleftarrow{F}A\xrightarrow{F}X$ refines also the identity diagram $X\xleftarrow{1_X}X\xrightarrow{1_X}X$ via the equivalence $F\colon A\to X$.
Therefore, $[a]^{-1}\circ [a]=[1_X]$. Similarly one sees that $[a]\circ [a]^{-1}=[1_Y]$ so that  $[a]\colon X\to Y$ is indeed a generalized isomorphism.

(2)\, Conversely we assume that $[a]\colon X\to Y$ is a generalized isomorphism and $[b]\colon Y\to X$ is the inverse so that $[b]\circ [a]=[1_X]$ and  $[a]\circ [b]=[1_Y]$.
If $[a]=[X\xleftarrow{F}A\xrightarrow{\Phi}Y]$ and 
$[b]=[Y\xleftarrow{G}A\xrightarrow{\Psi}X]$, we shall verify that 
$\Phi\colon A\to Y$ is an equivalence.

In view of Theorem \ref{proppp},  the diagram
$$X\xleftarrow{F\circ \pi_A} A\times_Y B\xrightarrow{\Psi\circ\pi_B} X$$
is a representative of the class $[b]\circ [a]=[1_X]$ 
and the diagram
$$Y\xleftarrow{G\circ  \pi'_B} B\times_X A\xrightarrow{\Phi\circ\pi'_A} Y$$
is a  a representative of the class $[a]\circ [b]=[1_Y]$. By Theorem \ref{proppp},  the sc-functor $\Psi\circ\pi_B: A\times_Y B\to  X$ is an equivalence and the same holds for the sc-functor $\Phi\circ\pi'_A:B\times_X A\to  Y$.  

\noindent {\bf (C1)} We claim that the functors $\Phi:A\to Y$ and $\Psi:B\to X$ are equivalences in the sense of the category theory. \\[0.5ex]
We first show that  $\Phi$ is full, that is, given $a, a'\in A$, the induced map $\textrm{mor}_{A}(a, a')\to \textrm{mor}_{Y}(\Phi (a), \Phi (a'))$ is surjective.  With  $a, a'\in A$, the objects $\Phi (a), \Phi(a')$ belong to $Y$ and since $G:B\to Y$ is essentially surjective, there are $b, b'\in B$ and morphisms $\varphi\colon G(b)\to  \Phi (a)$ and $\varphi'\colon G(b')\to  \Phi (a')$. In view of the definition of the weak fibered product $B\times_XA$, the triples $(b, \varphi, a)$ and $(b', \varphi, a')$ belong to $B\times_XA$. We know that the map 
\begin{equation}\label{faithful_morphism_1}
\textrm{mor}_{B\times_XA}((b,\varphi, a), ( b', \varphi', a'))\to \textrm{mor}_{Y}((\Phi \circ \pi'_A)(a, \varphi, b), (\Phi \circ \pi'_A)(a', \varphi', b')),
\end{equation}
induced by the equivalence  $\Phi\circ \pi'_A$, is surjective. We note that  the sets $\textrm{mor}_{Y}((\Phi \circ \pi'_A)(a, \varphi, b), (\Phi \circ \pi'_A)(a', \varphi', b'))$ and $\textrm{mor}_{Y}(\Phi(a), \Phi(a'))$ are equal.  Hence given the morphism 
morphism $h':\Phi(a)\to \Phi(a')$, we find  a morphism  
$$
(k, \psi, h)\in \textrm{mor}_{B\times_XA}((b,\varphi, a), ( b', \varphi',a'))
$$
 such that $\Phi\circ \pi'_A (k, \psi, h)=h'$. Since $\Phi\circ \pi'_A (k, \psi, h)=\Phi (h)$ and,  by definition of the morphism $(k, \psi, h)$, $h\in \textrm{mor}_{A}(a, a')$, we conclude that the map  $\textrm{mor}_{A}(a, a')\to \textrm{mor}_{Y}(\Phi (a), \Phi (a'))$ is surjective, as claimed. Next we show that this map is also injective.  Take two distinct elements $h, h'$ in $\textrm{mor}_{A}(a, a')$. Then the morphisms $\Phi (h)$ and $\Phi (h')$ belong to 
$\textrm{mor}_{Y}(\Phi(a), \Phi(a'))$ and we claim that they are distinct. Arguing as above we find two objects $b$ and $b'\in B$ and two morphisms $\varphi\colon G(b)\to  \Phi (a)$ and $\varphi'\colon G(b')\to  \Phi (a')$ such  that $(b, \varphi, a), (b', \varphi', a')\in B\times_XA$. Setting $\wt{k}:=(\varphi')^{-1}\circ \Phi (h)\circ \varphi$ and $\wt{k}':=(\varphi')^{-1}\circ \Phi (h')\circ \varphi$, we obtain two morphisms in ${\bm{Y}}$ between $G(b)$ and $G(b')$. Since the equivalence $G:B\to Y$ is full, there are morphisms $k$ and $k'$ in ${\bf B}$ between $b$ and $b'$ such that $G(k)=\wt{k}$ and $G(k')=\wt{k}'$. 
Then $(k, \varphi, h)$ and $(k', \varphi', h')$ are distinct morphisms in $\textrm{mor}_{B\times_XA}((b,\varphi, a), ( b', \varphi', a'))$ and since the map \eqref{faithful_morphism_1} induced by $\Phi\circ \pi'_A$ is injective, the morphisms $(\Phi\circ \pi'_A)(k, \varphi, h)$ and $(\Phi\circ \pi'_A)(k', \varphi', h')$ are distinct. From  $(\Phi\circ \pi'_A)(k, \varphi, h)=\Phi (h)$ and $(\Phi\circ \pi'_A)(k', \varphi', h')=\Phi (h')$  it follows that $\Phi (h)$ and $\Phi (h')$ are distinct elements of $\textrm{mor}_{Y}(\Phi(a), \Phi(a'))$, proving that the functor   $\Phi:A\to Y$ is faithful.  Finally,  since the functor 
$\Phi\circ\pi'_A:B\times_XA\to Y$ is essentially surjective,  we conclude for a given object $y\in Y$ that  there exists an object $(b, \varphi, a)\in B\times_XA$ and a morphism $\psi$ in ${\bm{Y}}$ between $(\Phi\circ\pi'_A)(b, \varphi, a)=\Phi (a)$ and $x$.  But this implies that the functor $\Phi:A\to Y$ is essentially surjective. Hence $\Phi$ is an equivalence in the sense of the category theory,  and similar arguments show the same for the functor $\Psi$. This finishes the proof of  the claim (C1).

It remains, in view of Lemma \ref{equivalence_in_the sense_of_category_theory}, to prove the following claim.\\[0.5ex]
\noindent {\bf (C2)}  The sc-functors $\Phi:A\to Y$ and $\Psi:B\to X$ are local sc-diffeomorphisms. \\[0.5ex]
 The arguments for the functors $\Phi$ and $\Psi$ are similar and hence we only consider the $\Phi$-case. 
 Fix a point $a\in A$ and choose an object $(b, \varphi, a)\in B\times_XA$.  Here $a\in A$, $b\in B$, and $\varphi:\Psi (b)\to F(a)$ is a morphism in ${\bm{X}}$.
Since the source and the target maps $s, t\colon {\bm{X}}\rightarrow X$ and also the equivalence $F\colon A\rightarrow X$ are  local sc-diffeomorphisms,  we find open neighborhoods
$U=U(\Psi  (b))$ of $\Psi (b)$ in $X$, $V=V(F(a))$ of $F(a)$ in $X$, ${\bm{U}}={\bm{U}}(\varphi)$ of $\varphi$ in ${\bm{X}}$, and $W=W(a)$ of $a$ in $A$ such that  the maps 
$$
s\vert {\bm{U}}\colon {\bm{U}}\to U,\quad t\vert {\bm{U}}\colon {\bm{U}}\to V,\quad \text{and}\quad  F\vert W\colon W\to V
$$ 
are sc-diffeomorphism.  Let $U':=\Psi^{-1}(U)\subset B$ and define for $b'\in U'$
$$\Gamma (b')=(b', (s\vert {\bm{U}})^{-1}(b'), (F\vert W)^{-1}\circ t\circ  (s\vert {\bm{U}})^{-1}(b')).$$
Then $\Gamma (b')\in B\times_XA$ so that $\Gamma:U'\to B\times_XA$. The map $\Gamma$ is a local sc-diffeo\-mor\-phism and its inverse is equal to a suitably restricted projection 
$\pi'_B:B\times_XA\to B$. 
Since $\Phi\circ \pi'_A$ is an equivalence, the composition $(\Phi\circ \pi'_A)\circ \Gamma:U'\to Y$, 
$$(\Phi\circ \pi'_A)\circ \Gamma (b')= \Phi\circ  (F\vert W)^{-1}\circ t\circ  (s\vert {\bm{U}})^{-1}(b')),$$
is a local sc-diffeomorphism. In addition the map $\Gamma':U'\to A$, defined by 
$$\Gamma'(b')= (F\vert W)^{-1}\circ t\circ  (s\vert {\bm{U}})^{-1}(b'),$$
is a local sc-diffeomorphism. Consequently, the map $\Phi:A\to Y$ is a local sc-diffeomorphism and the claim (C2) is proved.

Having  proved  the claims (C1) and (C2), it follows now from Lemma \ref{equivalence_in_the sense_of_category_theory}
that $\Phi$ is an equivalence. The proof of Theorem \ref{strong-iso} is complete.
\qed \end{proof}

\chapter{Geometry  up to Equivalences}\label{CHAPTER_11}

We now study the category ${\mathcal{EP}}({\bf E}^{-1})$ in more detail.  In particular we investigate how certain concepts
behave with respect to equivalences and generalized maps.

\section{Ep-Groupoids and Equivalences}
We recall that the objects are the ep-groupoids and the morphisms are the generalized maps $[d]\colon X\to Y$ abbreviating the equivalence  classes  
$$
[d]=[X\xleftarrow{F}A\xrightarrow{\Phi} Y]
$$
of diagrams, in which $F$ is an equivalence and $\Phi$ a sc-smooth functor between ep-groupoids. 
The invertible generalized maps are the equivalence classes
$$
[d]=[X\xleftarrow{F}A\xrightarrow{\Phi} Y]
$$
in which both functors $F$ and $\Phi$ are equivalences, in view of Theorem \ref{strong-iso}. Often we denote a generalized map
by $\mathfrak{f}:X\rightarrow Y$, where $\mathfrak{f}=[d]$.

\begin{theorem}\label{TANGENTXXX}\index{T- Tangent functor in ${\mathcal{EP}}({\bf E}^{-1})$}
The tangent functor induces in a natural way a functor
$$
T:{\mathcal{EP}}({\bf E}^{-1})\rightarrow {\mathcal{EP}}({\bf E}^{-1}).
$$
It maps an ep-groupoid $X$ to the ep-groupoid $TX$ and a morphism  $[d]:X\rightarrow Y$ to the morphism  $T[d]:TX\rightarrow TY$
defined by 
$$
T[X\xleftarrow{F}A\xrightarrow{\Phi}Y]=[TX\xleftarrow{TF} TA\xrightarrow{T\Phi}TY].
$$
\end{theorem}
\begin{proof}
On objects we just take the tangent $TX$ of the ep-groupoid $X$.  We  already know from Theorem \ref{main_KK} that $TX$  is an ep-groupoid. 
Let  $X\xleftarrow{F}A\xrightarrow{\Phi}Y$ 
be  a representative of the equivalence class  $[d]:X\rightarrow Y$. Since $F$ is an equivalence,   the sc-smooth functor $TF$ is an equivalence by Theorem \ref{gertrude} and we claim that if a diagram $X\xleftarrow{G}B\xrightarrow{\Psi}Y$ is a refinement of the diagram $X\xleftarrow{F}A\xrightarrow{\Phi}Y$, then 
the diagram $TX\xleftarrow{TG} TB\xrightarrow{T\Psi}TY$ is a refinement of the diagram $TX\xleftarrow{TF} TA\xrightarrow{T\Phi}TY$. 
So,   we consider the diagram 
\begin{equation*}
\begin{CD}
X@<F<<A@>\Phi>>Y\\
@. @AAHA @. \\
X@<G<<B@>\Psi>>Y\\ \\
\end{CD}
\end{equation*}
in which $H:B\to A$ is an equivalence and $F\circ H\simeq G$ and  $\Phi\circ H\simeq \Psi$. We claim that $T(F\circ H)=TF\circ TH\simeq TG$ and $T(\Phi\circ H)\simeq T\Psi$.
In view of Theorem \ref{gertrude} the tangent functors   $TF$, $TG$,  and $TH$ are  equivalences. 
Since $F\circ H\simeq G$, there exists a sc-smooth map $\sigma:B\to {\bm{X}}$ such that, for every $x\in B$, $\sigma (x)$ is a morphism $G(x)\to F\circ H(x)$ in ${\bm{X}}$. Hence  
$T\sigma \colon TB\to T{\bm{X}}$ is sc-smooth and by Proposition \ref{169} a natural transformation $TG\rightarrow TF\circ TH$ and consequently $TF\circ TH\simeq TG$.
The remaining argument is similar.
\qed \end{proof}

If $X$ is an ep-groupoid we denote by  $d_X$  the degeneracy index map for the underlying object M-polyfold.
If $\phi:x\rightarrow y$ is a morphism,  then  $d_X(x)=d_X(y)$ so that we obtain the  {\bf induced degeneracy map}\index{Induced degeneracy map}
$$
d_{|X|}:|X|\rightarrow {\mathbb N}
$$
on the orbit space $|X|$.
\begin{proposition}\label{deg-x}\index{P- Invariance of $d_{|X|}$}
We assume that $X$ and $Y$ are ep-groupoids and $d_X$ and $d_Y$ the respective degeneracy index maps, and
$d_{|X|}$  and $d_{|Y|}$ the induced degeneracy maps.
Let  $[d]\colon X\rightarrow Y$ be  an isomorphism in ${\mathcal{EP}}({\bf E}^{-1})$ and 
$|[d]|\colon |X|\rightarrow|Y|$
 the induced homeomorphism.
If $|[d]|(\abs{x})=\abs{y}$, then 
$$
d_X(x')=d_Y(y')
$$
for every $x'\in \abs{x}$ and $y'\in \abs{y}$ and $d_{|X|}(\abs{x})=d_{|Y|}(\abs{y})$.
\end{proposition}
\begin{proof}
If $x'\in \abs{x}$, then $d_{X}(x')=d_{X}(x)$ since $s, t:{\bm{X}}\to X$ are local sc-diffeo\-mor\-phisms. Hence it suffices to show that $d_{X}(x)=d_{Y}(y)$.
Since $[d]:X\to Y$ is a generalized  isomorphism, the equivalence class  $[d]$ is,  in view of Proposition \ref{strong-iso}, represented by a diagram 
$$
X\xleftarrow{F} A\xrightarrow{G} Y
$$
 in which $F$ and $G$ are equivalences. 
By assumption, 
$$
\abs{[d]}(\abs{x})=\abs{G}\circ \abs{F}^{-1}(\abs{x})=\abs{y}.
$$
  Let $\abs{F}^{-1}(\abs{x})=\abs{a}$. Then $\abs{F}(\abs{a})=\abs{F(a)}=\abs{x}$ and $\abs{G}(\abs{a}]=\abs{G(a)}=\abs{y}$ and, since  equivalences are essentially surjective, 
there are morphisms  $\varphi\in {\bm{X}}$ and $\psi\in {\bm{Y}}$ so that
$$
\varphi\colon F(a)\rightarrow x \quad \text{and}\quad \psi\colon G(a)\rightarrow y.
$$
Using again the fact that $s$ and $t$ are local sc-diffeomorphisms ,  we find  that 
$$d_{X}(x)=d_{X}(F(a))\quad \text{and}\quad d_{Y}(y)=d_{Y}(G(a)).$$
Since $F$ and $G$ are local sc-diffeomorphisms, 
$$d_{X}(x)=d_{X}(F(a))=d_{A}(a)=d_{Y}(G(a))=d_{Y}(y).$$
Consequently, $d_{X}(x)=d_{Y}(y)$ as claimed.
\qed \end{proof}

If $X$ is a  tame M-polyfold we have introduced the notion of a face. It is by definition the closure of a connected component of 
the set $\{x\in X\ |\ d_X(x)=1\}$, see Definition \ref{DEF248}. The tame M-polyfold was called face-structured provided a point $x\in X$ lies in precisely
$d_X(x)$-many faces. Recall that $x$ always lies in at most $d_X(x)$-many faces.  When we study an ep-groupoid we can talk about the faces
of the underlying M-polyfold, which is useful for certain discussions. If we study ep-groupoids up to equivalences we need to take the 
morphisms into consideration and the definition of a face in the ep-groupoid sense should accommodate this fact. 
If $X$ is an ep-groupoid and $\phi\in \bm{X}$ 
 we have the equality 
 $$
 d_{X}(s(\phi))=d_X(t(\phi)),
 $$
 which we already used to define $d_{|X|}:|X|\rightarrow {\mathbb N}$. This allows us to give a more appropriate definition  of
 a face in the case of an ep-groupoid. We assume that the ep-groupoid is tame so that $\partial X$ has reasonable geometrical properties.
 If $X$ is an ep-groupoid and $\theta\subset |X|$,  there exists an associated full subcategory $X_\theta$\index{$X_\theta$} of $X$ associated to all
 objects $x\in X$ with $|x|\in\theta$.

\begin{definition}\label{DEF1113}\index{D- Face of a tame ep-groupoid}\index{D- Face structured ep-groupoid}
Let $X$ be a tame ep-groupoid. A {\bf face} of $X$ (in the ep-groupoid sense) is the full-subcategory $X_\theta$
associated to the closure $\theta$ in $|X|$ of a connected component $\theta^\circ$ of $\{z\in |X|\ |\ d_{|X|}(z)=1\}$. We shall 
call $\theta$ a face of $|X|$.
  A tame ep-groupoid is called {\bf face structured} (in the ep-groupoid sense)
provided a point $z\in |X|$ lies in precisely $d_{|X|}(z)$-many faces of $|X|$.  We say the tame ep-groupoid $X$ is {\bf weakly face structured}
provided the underlying object M-polyfold is face-structured.
\qed
\end{definition}
\begin{example}
To illustrate the definition consider the manifold with boundary with corners $X=C=[0,\infty)^2\subset {\mathbb R}^2$
and take as morphisms only the identities. Then $X$ is a face structed ep-groupoid.    In the second example take 
the same object space and we have ${\mathbb Z}_2=\{-1,1\}$ act as follows. The element $1$ acts through the identity
and $-1$ interchanges coordinates.  Then ${\mathbb Z}_2\ltimes C$ is an ep-groupoid. It is not face structured, but weakly face structured.
\qed
\end{example}
\begin{remark}\index{R- On being face-structured}
In certain contexts it is even true that face structured ep-groupoids have an additional structure, 
namely they are {\bf ordered face structured}\index{Ordered face structured}. The additional structure 
is that on the set of faces there is a partial ordering, which induces for an object $x\in X$ a total ordering on 
the set of faces containing $x$.  This additional structure occurs in applications to Floer-theoretic problems.
\qed
\end{remark}

There is a comprehensive discussion of faces of ep-groupoids in \cite{HWZ3.5}. Although in \cite{HWZ3.5} we still discussed only spaces built on splicing cores, the discussion applies verbatim to the tame ep-groupoids considered here. Recall that a M-polyfold build on  splicing cores is automatically tame.
According to Proposition 3.15 in \cite{HWZ3.5}, an equivalence $\Phi\colon X\to Y$ between ep-groupoids defines the natural bijections 
\begin{align*}
\Phi_\ast&\colon \text{Faces}_X\to \text{Faces}_{Y}\colon\quad \Phi_\ast(F)=\text{saturation}(\Phi (F))\\
\Phi^\ast&\colon \text{Faces}_Y\to \text{Faces}_{X}\colon\quad \Phi^\ast(G)=\Phi^{-1} (G))
\end{align*}
and $\Phi^\ast =(\Phi_\ast)^{-1}$. Moreover, if $X$ is face structured, then also $Y$ is face structured and vice versa.
Consequently, considering a generalized isomorphism
$$[d]=[X\xleftarrow{\Phi}A\xrightarrow{\Psi}Y]$$
we can define the bijections 
\begin{equation*}
[d]_\ast\colon \text{Faces}_X\to \text{Faces}_{Y}\quad \text{and}\quad 
[d]^\ast\colon \text{Faces}_Y\to \text{Faces}_{X}
\end{equation*}
by
\begin{equation*}
[d]_\ast=\Psi_\ast\circ \Phi^\ast\quad \text{and}\quad [d]^\ast=\Phi_\ast\circ \Psi^\ast.
\end{equation*}
They satisfy $[d]^\ast=([d]_\ast)^{-1}.$

\begin{proposition}\index{P- Boundary of face-structured $X$}
If the tame ep-groupoid $X$ is face structured, then each face $F$ has an object set $\text{obj}_F$ which is a sub-M-polyfold of $X$ 
on which the induced structure is tame.  
The same holds for the morphism set $\text{mor}_X$. Equipped with the induced structures the face $F$ is in a natural way
an ep-groupoid. If the intersection of two faces $F$ and $G$ is  nonempty, then it is a face  in $F$ as well as in $G$.
 If $X$ and $Y$ are tame ep-groupoids and $[d]:X\rightarrow Y$ is a generalized  isomorphism, then
$[d]$ induces the  bijection
$$
[d]_\ast\colon \text{Faces}_X\rightarrow{Faces}_Y.
$$
Moreover,  for every face $F$ in $X$ and its corresponding face $G=[d]_\ast (F)$ in $Y$ the generalized isomorphism $[d]$ induces a natural equivalence class
$$
[d]_F:F\rightarrow G
$$
which is a generalized  isomorphism.
\end{proposition}

\begin{proof}
Pick a connected component $\theta^\circ $ in $\{z\in |X|\ |\ d_{|X|}(z)=1\}$ and denote by $\theta$ its closure in $|X|$.
A point $x\in X_\theta$ lies in exactly $d=d_X(x)$-many local faces of the M-polyfold $X$. 
If $x$ would lie in less than $d$-many global faces of the M-polyfold $X$ it would immediately follow
that there are less than $d$-many connected components $\theta_i^\circ$  of $\{z\in |X|\ |\ d_{|X|}(z)=1\}$
with their closures $\theta_i$ in $|X|$ containing $|x|$. 

Let $\theta$ be any face of $|X|$. In view of the previous discussion we find that an object $x\in X_\theta$
lies in precisely one global face $F$ of the M-polyfold $X_\theta$ so that $F\subset X_\theta$. 
The global faces $F$ of the M-polyfold $X$  are closed in the M-polyfold $X$,  and $X_\theta$ consists of a disjoint 
union of these. It has been proved in Theorem \ref{FACE_XX} that $F$ is a sub-M-polyfold with an induced tame structure.
Since $X_\theta$ is a disjoint union of such $F$ we see that the object space $X_\theta$ is a sub-M-polyfold
and consequently the associated full subcategory is an ep-groupoid.

If $[d]:X\rightarrow Y$ is a generalized  isomorphism  between the tame ep-groupoids $X$ and $Y$, then it is represented by a diagram 
$$
X\xleftarrow{\Phi}A\xrightarrow{\Psi} Y
$$
in which $\Phi$ and $\Psi$ are equivalences.  If $\theta^\circ$ is a connected component of $\{z\in |X|\ |\ d_{|X|}(z)=1\}$ it follows for the homeomorphism
$|[d]|:|X|\rightarrow |Y|$ that the image 
$|[d]|(\theta^\circ)$ is a connected component of $\{q\in |Y|\ |\ d_{|Y|}(q)=1\}$. Since $|[d]|$ is a homeomorphism 
$$
\cl_{|Y|}([d](\theta^\circ)) =|[d]|(\theta)
$$
We therefore define
$$
[d]_\ast X_\theta =Y_{|[d]|(\theta)}\ \ \text{and}\ \  [d]^\ast Y_{\theta'}= X_{|[d]|^{-1}(\theta')}.
$$
Clear $[d]_\ast$ and $[d]^\ast$ are mutually inverse to each other.

Alternatively, if $F$ is a face of $X$ one verifies easily that  $B=\Phi^{-1}(F)$ is a face of $A$ and in turn the saturation $C$ of $\Psi(B)$  is a face of $Y$.
In  other words,  there are  natural bijections between the faces of $A$ and those of $X$ and $Y$. Also it is clear by the above construction, that the restrictions of $\Phi$ and $\Psi$ define 
the diagram
\begin{equation}\label{walter_}
F\xleftarrow{\Phi\vert B} B\xrightarrow{\Psi\vert B} C.
\end{equation}
The equivalence class  $[e]:=[d]\vert F$ of this diagram is a generalized  isomorphism $F\rightarrow C$ by Theorem \ref{strong-iso}. Indeed,  the restricted functors in \eqref{walter_} are essentially surjective,
full and faithful and local sc-diffeomorphism, and hence equivalences by 
Lemma \ref{equivalence_in_the sense_of_category_theory}.
\qed \end{proof}

\begin{remark}\label{REM1116}\index{R- Face-structured and generalized isomorphisms}
The reader should note the following. Assume that $X$ and $Y$ are ep-groupoids and $\mathfrak{f}:X\rightarrow Y$
is a generalized isomorphism. If one these ep-groupoids is tame and face-structured the same holds for the other one.
\qed
\end{remark}

The following somewhat stronger result holds true.
\begin{proposition}\index{P- Correspondence of faces}
Let $X$ and $Y$ be tame face-structured ep-groupoids and $\mathfrak{f},\mathfrak{g}:X\rightarrow Y$ two generalized isomorphisms
satisfying $|\mathfrak{f}|=|\mathfrak{g}|$.
Then the induced bijections between the collection of faces are the same, i.e. $\mathfrak{f}_\ast=\mathfrak{g}_\ast$.
Further, the induced generalized isomorphisms $F\rightarrow \mathfrak{f}_\ast F$  and $G\rightarrow \mathfrak{f}_\ast G$ induce the same map $\abs{F}\to \abs{\mathfrak{f}_\ast F}$ on the orbit spaces.
\end{proposition}
\begin{proof}
This follows from the previous discussion.
\qed \end{proof}

\section{Sc-Differential Forms and Equivalences}
We next investigate the behavior of sc-differential forms on ep-groupoids (introduced in Section 
\ref{SECT82}) under equivalences. The main result in this section is the following 
theorem. 

 \begin{theorem}\label{OOmaine}\index{T- Generalized sc-smooth maps and forms}
 We assume that $\mathfrak{f}\colon X\rightarrow Y$ is a generalized sc-smooth map between ep-groupoids, and that $\mathfrak{h},\mathfrak{k}\colon X\rightarrow Y$
 are two sc-smooth generalized isomorphisms. Then the following holds.
 \begin{itemize}
 \item[{\em (1)}]\ $d\colon \Omega^{\ast}_{ep,\infty}(X)\rightarrow \Omega^{\ast}_{ep,\infty}(X)$ and $d\circ d=0$.
\item[{\em (2)}]\ There is a well-defined operation of pull-back $\mathfrak{f}^\ast\colon \Omega^{\ast}_{ep,\infty}(Y)\rightarrow \Omega^{\ast}_{ep,\infty}(X)$
 such that $d_X\circ \mathfrak{f}^\ast = \mathfrak{f}^\ast \circ d_Y$.
 \item[{\em (3)}]\ For the generalized isomorphism $\mathfrak{h}\colon X\to Y$ there is a well-defined push-forward $\mathfrak{h}_\ast\colon \Omega^{\ast}_{ep,\infty}(X)\rightarrow \Omega^{\ast}_{ep,\infty}(Y)$ satisfying $\mathfrak{h}_\ast\circ d_X=d_Y\circ \mathfrak{h}_\ast$. Moreover, with the pull-back $\mathfrak{h}^\ast$  we have that $\mathfrak{h}_\ast\circ \mathfrak{h}^\ast =Id$ and $\mathfrak{h}^\ast\circ \mathfrak{h}_\ast=Id$.
 \item[{\em (4)}]\ If the generalized isomorphisms $\mathfrak{h}$ and $\mathfrak{k}$ induce the same  maps between orbit spaces,  so that  $|\mathfrak{h}|=|\mathfrak{k}|$, then $\mathfrak{h}_\ast=\mathfrak{k}_\ast$ and $\mathfrak{h}^\ast=\mathfrak{g}^\ast$.
 \end{itemize}
 \qed
 \end{theorem}
For statement  (1) we refer to Section 
\ref{SECT82}.  The other statements  of the  theorem will follow from several results of independent interest, where we make precise
how $\mathfrak{f}^\ast$ and $\mathfrak{h}_\ast$ are defined. We recall that the generalized map $\mathfrak{f}:X\rightarrow Y$ is an equivalence class
of diagrams 
 $$
 X\xleftarrow{F} A\xrightarrow{\Phi}Y
 $$
in which  $F$ is an equivalence of ep-groupoids and $\Phi$ a sc-smooth functor.
The generalized isomorphism is given by the same  diagram in which, however,  $\Phi$ is also an equivalence.
 
 \begin{lemma}\label{differential_form_1}\index{L- Differential forms {I}}
 We consider an equivalence 
$F:X\rightarrow Y$   between ep-groupoids and a  sc-differential form
$\omega$   on $X$.  Given $y\in Y_1$ and $h_1,\ldots,h_k\in T_yY$, we let  $\psi$ and $\psi'\in {\bm{Y}}_1$
be two morphisms  satisfying $\psi:F(x)\rightarrow y$ and $\psi'\colon F(x')\rightarrow y$ for $x$ and $x'\in X_1$.
If $l_1,\ldots ,l_k\in T_xX$ and $l_1',\ldots ,l_k'\in T_{x'}X$
satisfy
\begin{equation}\label{diff_sc_forms}
T\psi(TF(x)l_j )=h_j=T\psi'(TF(x')l_j'),
\end{equation}
then 
$$
 \omega_x(l_1,\ldots,l_k)=\omega_{x'}(l_1',\ldots,l_k').
 $$
 \end{lemma}
 \begin{proof}
 We start with the morphism
 $$
\tau:= \psi^{-1}\circ\psi'\colon F(x')\rightarrow F(x)
 $$
in ${\bm{Y}}_1$. Its linearization at $x'$ satisfies, in view of the assumption \eqref{diff_sc_forms}, 
 \begin{equation}\label{oip}
 T \tau(TF(x')l_j') = TF(x)l_j
 \end{equation}
 for all $j=1,\ldots ,k$.
  Since $F$ is an equivalence there exists a uniquely determined morphism $\sigma\colon x'\rightarrow x$ satisfying
 $$
 F(\sigma)= \tau.
 $$
This implies  $\sigma\in {\bm{X}}_1$ since $\tau\in {\bm{Y}}_1$. In order to relate $T\sigma$ and $T\tau$ we    take,  for suitable open neighborhoods, 
 the sc-diffeomorphisms 
 $s_X\colon {\bm{U}}(\sigma)\rightarrow U(x')$ and $t_X\colon {\bm{U}}(\sigma)\rightarrow U(x)$ and recall the  definition of $T\sigma$,
$$
T\sigma = T(t_X\circ (s_X|U(\sigma))^{-1})(x').
$$
By the functoriality of $F$, 
\begin{equation*}
\begin{split}
F(t_X\circ (s_X|U(\sigma))^{-1}(q))
 =t_Y\circ F((s_X|U(\sigma))^{-1}(q))
 =t_Y\circ (s_Y|U(\tau))^{-1}(F(q)).
 \end{split}
 \end{equation*}
Taking the derivative in $q$ at $q=x'$ leads to 
$$
TF(x)\circ T\sigma = T\tau\circ TF(x').
$$
Consequently, in view of \eqref{oip},
$$
TF(x)l_j = T \tau(TF(x')l_j') = TF(x)(T\sigma l_j').
$$
This implies, since $TF(x)\colon T_xX\to T_{F(x)}Y$ is an isomorphism,  that
$$
l_j=T\sigma l_j'\quad \text{for all $ j$.}
$$
Therefore, using the compatibility of sc-differential forms with morphisms, 
$$
\omega_x(l_1,\ldots,l_k)=\omega_x(T\sigma l_1',\ldots,T\sigma l_k')=\omega_{x'}(l_1',\ldots,l_k'),
$$
as claimed in Lemma 
\ref{differential_form_1}.
 \qed \end{proof}
 Lemma \ref{differential_form_1} allows to define the push-forward $F_\ast\omega\in \Omega_{ep,\infty}(X)$ under an equivalence $F\colon X\to Y$.
 \begin{definition}\label{d_push_forward}\index{D- Push-forward of forms}
The {\bf push-forward}\index{D- Push-forward of sc-differential form} of the sc-differential form $\omega\in \Omega^k_{ep,\infty}(X)$ by an equivalence  $F\colon X\rightarrow Y$ between ep-groupoids, denoted by $$F_\ast\omega\in \Omega^k_{ep,\infty}(Y),$$ is defined as follows.
If $y\in Y_1$ and $h_1,\ldots, h_k\in T_yY$, we  take $x\in X_1$ (using that $F$ is essentially surjective) such  that there is a morphism $\psi\colon F(x)\rightarrow y$. Then we take  
 $l_j\in T_xX$ satisfying $h_j=T\psi(TF(x)l_j)$,  and define
 $$
 (F_\ast\omega)_y(h_1,\ldots ,h_k)=\omega_x(l_1,\ldots ,l_k).
 $$
\qed
 \end{definition}
 In view of Lemma \ref{differential_form_1} this is well-defined and independent of the choice of the morphism $\psi$.
 In order to see that the push-forward  form $F_\ast \omega$ on $Y$ as defined above is sc-smooth, we recall that the equivalence $F\colon X\to Y$ is, by definition, locally a sc-diffeomorphism. Moreover, the morphism $\psi\colon F(x)\to y$ has a unique sc-smooth extension $\wh{\psi}$ to a sc-diffeomorphism between properly chosen open neighborhoods 
$U(F(x))=F(U(x))$ and $U(y)$ in $Y$
$$
\wh{\psi}:U(F(x))\rightarrow U(y).
$$
 The composition  $\wh{\psi}\circ F\colon U(x)\to U(y)$ is a sc-diffeomorphism satisfying $(\wh{\psi}\circ F)(x)=y$. Hence the push-forward of the sc-smooth form 
$\omega$ on $X$ by the equivalence $F\colon X\to Y$ is the pull-pack of the form $\omega$ by the sc-diffeomorphism 
$(\wh{\psi}\circ F)^{-1}\colon U(y)\to U(x)$,
\begin{equation}\label{push-forward_eq1}
F_\ast \omega=\bigl((\wh{\psi}\circ F)^{-1}\bigl)^\ast \omega.
\end{equation}
This shows sc-smoothness.
We also note that  standard arguments now imply that 
\begin{equation}\label{push-forward_eq2}
d_Y(F_\ast \omega)=F_\ast (d_X\omega).
\end{equation}

 \begin{lemma}\label{equivalence_forms}\index{C- Equivalences and forms}
Assume that $F:X\rightarrow Y$ and $G:Y\rightarrow Z$ are equivalences between ep-groupoids.
Then 
\begin{itemize}
\item[{\em (1)}]\ $(G\circ F)^\ast \omega =F^\ast \circ G^\ast\omega$.
\item[{\em (2)}]\ $(G\circ F)_\ast\omega = G_\ast \circ F_\ast\omega$
\item[{\em (3)}]\ $F_\ast \circ F^\ast \omega=\omega\quad  \text{and}\quad  F^\ast \circ F_\ast\tau=\tau.$
\end{itemize}
\end{lemma}
\begin{proof} The well-known formula (1) is an immediate consequence of the definition of the pull-back of a form $\omega$ on $Y$, which is the form $F^\ast\omega$ on $X$ given by
$$(F^\ast \omega)_x(k_1,\ldots, k_n)=\omega_{F(x)}(TF(x)k_1,\ldots, TF(x)k_n).$$
The statements (2)  and (3) will  follow from 
Definition \ref{d_push_forward}, the properties of equivalences and the invariance of the forms.

We consider the sc-form $\omega$ on $X$ and fix $z\in Z_1$ and $h_1,\ldots ,h_k\in T_zZ$. The composition $G\circ F\colon X\to Z$ is an equivalence and hence essentially surjective. Therefore,  there exist $x\in X_1$ and a morphism $\psi\colon G\circ F(x)\to z$ in ${\bm{Z}}_1$. Taking the unique solutions $l_j$ of 
$T\psi\circ T(G\circ F)(x)l_j=h_j$, we obtain, by definition, 
$$
((G\circ F)_\ast\omega)_z(h_1,\ldots ,h_k)=\omega_x(l_1,\ldots ,l_k).$$
On the other hand, if $y=F(x)$, we have the morphism $\psi\colon G(y)\to z$ and taking the solutions $T\psi\circ TG(y)l_j'=h_j$ which are equal to $l_j'=TF(x)l_j$ we obtain, by definition,
$$
(G_\ast (F_\ast \omega))_z(h_1,\ldots ,h_k)=(F_\ast \omega)_y(TF(x)l_1,\ldots ,(TF(x)l_k).$$
Since $y=F(x)$, we have the $1$- morphism $1_y\colon y\to y$ in ${\bm{Y}}$ whose linearization is the identity map and get for the solution of 
$TF(x)l_j'=TF(x)l_j$, which is given by $l_j'=l_j$, again by definition,
$$
(F_\ast \omega)_y(TF(x)l_1,\ldots , TF(x) l_k)=\omega_x(l_1,\ldots ,l_k).
$$
We have verified that $G_\ast (F_\ast\omega)=(G\circ F)_\ast \omega$.

In order to verify $F_\ast \circ F^\ast \omega=\omega$ for forms  $\omega$ on $Y$, we fix $y\in Y$ and $h_1,\ldots, h_k\in T_yY$, find $x \in X$ and a morphism $\phi\colon F(x)\to y$ in ${\bm{Y}}$. With the solution $l_j$ of $T\phi\circ  TF(x)l_j=h_j$ we get 
\begin{equation*}
\begin{split}
(F_\ast(F^\ast \omega))_y(h_1, \ldots ,h_k)&=(F^\ast \omega)_x(l_1,\ldots ,l_k)\\
&=(F^\ast \omega)_{F(x)}(TF(x)l_1,\ldots ,TF(x)l_k)\\
&=\omega_y (h_1,\ldots ,k_k),
\end{split}
\end{equation*}
where we have used that 
$l_j=(TF(x))^{-1}\circ (T\phi)^{-1}h_j$ and also have used the invariance of $F^\ast \omega$ under morphisms. The second formula in (3) is verified the same way and the proof of Lemma \ref{equivalence_forms} is complete. 
\qed \end{proof}

Next we study the compatibility with natural equivalences.
 \begin{lemma}\label{Goddard}\index{L- Differential forms {II}}
Let $F,G\colon X\rightarrow Y$ be equivalences between ep-groupoids which are naturally equivalent by the natural transformation
$$
\tau\colon F\rightarrow G.
$$
If $\omega$ is a sc-differential form on $X$ and $\omega'$ a sc-differential form on $Y$, then
$$
F_\ast\omega =G_\ast\omega\quad  \text{and}\quad  \ F^\ast\omega'= G^\ast\omega'.
$$
 If $\Phi,\Psi:X\rightarrow Y$ are sc-smooth functors which are naturally equivalent,
then $\Psi^\ast\omega=\Phi^\ast\omega$.
\end{lemma}
 \begin{proof}
 By assumption,  $\tau:X\rightarrow {\bm{Y}}$ is a sc-smooth map associating with every $x\in X$ the morphism
 $$
\tau(x)\colon  F(x)\rightarrow G(x)
 $$
such that  for every morphism $\phi\colon x\rightarrow x'$ in ${\bm{X}}$, the diagram 
 $$
 \begin{CD}
 F(x)@>\tau(x)>> G(x)\\
 @V F(\phi)VV @V G(\phi)VV\\
 F(x')@>\tau(x') >> G(x')
 \end{CD}
 $$
 is commutative.
 We choose $h_1,\ldots, h_k\in T_yY$. We take  $\psi\colon F(x)\rightarrow y$ and $\psi'\colon G(x')\rightarrow y$ which have to be on level $1$ 
 and take  $l_1,\ldots,l_k\in T_xX$ solving  $T\psi\circ TF(x)l_j=h_j$ and similarly take $l_j'$ in $T_{x'}X$ satisfying  $T\psi'\circ TG(x')l_j'=h_j$.
 By definition, 
 $$
 F_\ast\omega(h_1,\ldots ,h_k)=\omega(l_1,\ldots,l_k)\quad \text{and}\quad G_\ast\omega(h_1,\ldots,h_k)=\omega(l_1',\ldots,l_k').
 $$
Recalling the morphism $\tau(x):F(x)\rightarrow G(x)$, we consider the composition
 $$
 G(x')\xrightarrow{\psi'}y\xrightarrow{\psi^{-1}} F(x)\xrightarrow{\tau(x)} G(x).
 $$
 Since $y$ is on level $1$,  this is true for the data above as well. The functor  $G$ is an equivalence, hence  we find a morphism
 $
 \sigma\colon x'\rightarrow x
 $
 satisfying  $G(\sigma)= \tau(x)\circ\psi^{-1}\circ \psi'$. A computation, similar to that in the previous lemma, leads to 
 $$
 T\sigma(l_j')=l_j
 $$
 which implies that 
 $$
 \omega_x(l_1,\ldots,l_k)=\omega_{x'}(l_1',\ldots,l_k').
 $$
Consequently,  $(F_\ast\omega)_y(h_1,\ldots,h_k)=(G_\ast\omega)_y(h_1,\ldots,h_k)$  proving the statement for the push-forward. 
 
 For the pull-back we start again with  the morphism
 $$
\tau(x)\colon  F(x)\rightarrow G(x)\quad \text{in ${\bm{Y}}$}, 
 $$
so that  $F(x)=s(\tau(x))$ and $G(x)=t(\tau(x))$ for all $x\in X$. We choose suitable open neighborhoods such  $s:{\bm{U}}(\tau(x))\rightarrow U(F(x))$ is a sc-diffeomorphism. Then 
 $$
 t\circ (s\vert {\bm{U}}(\tau(x)))^{-1}\circ F(q) = G(q).
 $$
We differentiate at $q=x$ which we assume to be on level $1$ and, abbreviating $\phi=\tau(x)$, we obtain
 $$
 TG(x) = T\phi\circ TF(x).
 $$
 By the compatibility of the sc-differential form $\omega$ on $Y$ with morphisms we conclude,
  \begin{equation*}
  \begin{split}
G^\ast \omega(h_1,\ldots, h_k)&= \omega(TG(x)h_1,\ldots,TG(x)h_k)\\&=\omega(T\phi\circ TF(x)h_1,\circ ,T\phi\circ TF(x)h_k)\\
 &=\omega(TF(x)h_1,\ldots,TF(x)h_k)=F^\ast \omega(h_1,\ldots, h_k).
 \end{split}
\end{equation*}
 This proves the assertion about the pull-backs.  Note that we did not use 
 the fact that $F$ and $G$ are equivalences.
 Hence we can apply the same idea to the natural transformation $\tau:\Phi\rightarrow\Psi$ in order to complete the proof of Lemma \ref{Goddard}.
\qed \end{proof}

\begin{lemma}\label{compatibilityyy}\index{C- Pull-back of forms}
We consider two diagrams $d\colon X\xleftarrow{F} A\xrightarrow{G}Y$  and $d'\colon X\xleftarrow{F'} A'\xrightarrow{G'}Y$, in which $F$ and $F'$ are equivalences and $G, G'$ are sc-smooth functors between ep-groupoids, and assume  
that  $H:A'\rightarrow A$ is an equivalence, as illustrated in the diagram

\begin{equation*}
\begin{CD}
d:X@<F<<A@>G>>Y\\
@. @AAHA @. \\
d':X@<F'<<A'@>G'>>Y.\\ 
\end{CD}
\end{equation*}
\mbox{}\\[1ex]
If $F\circ H\simeq F'$ (naturally equivalent) and $G\circ H\simeq G'$, then 
$$
F_\ast \circ G^\ast\omega' = (F')_{\ast} \circ (G')^{\ast}\omega'.
$$
for every sc-differential form $\omega'$ on $Y$.
If, in addition, $G$ and $G'$ are equivalences, then 
$$G_\ast \circ F^\ast\omega = (G')_{\ast} \circ (F')^{\ast}\omega$$
for every sc-differential form $\omega$ on $X$.
\end{lemma}
\begin{proof} Recalling that the composition of equivalences is an equivalence, and using Lemma \ref{differential_form_1} and Lemma \ref{d_push_forward} we compute,
\begin{equation*}
\begin{split}
{F'}_\ast \circ (G')^\ast\omega &=(F\circ H)_\ast \circ (G')^\ast\omega\\
&=F_\ast \circ H_\ast  \circ (G\circ H)^\ast\omega\\
&= F_\ast \circ H_\ast \circ H^\ast \circ G^\ast\omega\\
&= F_\ast \circ  G^\ast\omega.
\end{split}
\end{equation*}
The second part is proved the same way.
 \qed \end{proof}
We recall that two diagrams in the equivalence class $[d]=[X\xleftarrow{F}A\xrightarrow{G}Y]$, in which $F$ is an equivalence and $G$ a sc-smooth functor between ep-groupoids, are equivalent if they possess a common refinement. Therefore, we are able, in view of the previous lemma, to introduce the pull-back $[d]^\ast \omega$ 
for generalized maps $[d]\colon X\to Y$ and the push-forward $[d]_\ast\omega$ for generalized isomorphisms $[d]\colon X\to Y$ between ep-groupoids.
\begin{definition}[{\bf Pull-back, push-forward by generalized maps}]
If $[d]=[X\xleftarrow{F}A\xrightarrow{G}Y]\colon X\to Y$ is a generalized map between ep-groupoids then the {\bf pull-back} $[d]^\ast\omega'$ of a  
sc-differential form $\omega'$ on $Y$ is the sc-differential form on $X$,  defined by 
$$[d]^\ast\omega'=F_\ast\circ G^\ast \omega'\in \Omega_{\textrm{ep},\infty}(X).$$ 
If $[d]=[X\xleftarrow{F}A\xrightarrow{G}Y]\colon X\to Y$ is a generalized 
isomorphism (i.e., both $F$ and $G$ are equivalences), the {\bf push-forward} $[d]_\ast\omega$ of a sc-differential form
 $\omega$ on $X$ is the sc-differential form on $Y$,  defined by 
$$[d]_\ast\omega=G_\ast\circ F^\ast \omega\in \Omega_{\textrm{ep},\infty}(Y).$$ 
\qed
\end{definition}
In view of Lemma \ref{compatibilityyy} the push-forward and the pull-back by generalized isomorphism are well-defined, i.e., independent of the choice of the diagram representing the equivalence class $[d]$.

With Definition \ref{differential_form_1}, the statements (2) and (3) of Theorem \ref{OOmaine} now follow from Lemma \ref{differential_form_1} and Lemma \ref{equivalence_forms}.

The proof of the  statement (4) in Theorem \ref{OOmaine} requires a preparation.
 \begin{proposition}\label{dependence_on_F}\index{P- Dependence on $F$}
We assume that $F,G:X\rightarrow Y$ are equivalences between ep-groupoids satisfying $|F|=|G|\colon \abs{X}\to \abs{Y}$.
 Then for  the  sc-differential forms $\omega$ on $X$ and $\omega'$ on $Y$ we have  the identities  
 $$
 F_\ast\omega=G_\ast\omega\ \text{and}\ \ F^\ast\omega'=G^\ast\omega'.
 $$
 \end{proposition}
 \begin{proof}
 Since $|F|=|G|$,  Proposition \ref{prop4.14} implies for every $x\in X_1$ that 
 there exists a morphism $\psi\colon F(x)\rightarrow G(x)$ satisfying
 $$
 T\psi(TF(x)h)=TG(x)h\quad  \text{for all}\quad  h\in T_xX.
 $$
 If $\omega'$ is a sc-differential form on $Y$,  we compute, using 
 $\omega'_{G(x)}\circ (T\psi\oplus\ldots \oplus T\psi)=\omega_{F(x)}'$, 
 \begin{equation*}
 \begin{split}
(G^\ast\omega')_x(h_1,\ldots,h_k)
&=\omega'_{G(x)}(TG(x)h_1,\ldots,TG(x)h_k)\\
&=\omega'_{G(x)}(T\psi(TF(x)h_1),\ldots,T\psi(TF(x)h_k))\\
&=\omega'_{F(x)}(TF(x)h_1,\ldots,TF(x)h_k)\\
&=(F^\ast\omega')_x(h_1,\ldots,h_k).
\end{split}
\end{equation*}
Therefore,
$$
G^\ast\omega'=F^\ast\omega'.
$$
Next we investigate the 
push-forward.
We consider $F_\ast\omega$ and $G_\ast\omega$ and fix  $y\in Y_1$.
For $h_1,\ldots,h_k\in T_yY$ we compute $F_\ast\omega$ by taking  $x\in X_1$ and a morphism $\phi\in {\bm{Y}}_1$ such that
$\phi:F(x)\rightarrow y$. Then we take the solutions  $l_1,\ldots, l_k\in T_xX$ of  $T\phi(TF(x)l_j)=h_j$, and  obtain, by definition, 
$$
\omega_x(l_1,\ldots ,l_k)=:(F_\ast\omega)_y(h_1,\ldots,h_k).
$$
The definition for $G$ is similar. 
Recalling the morphism 
$\psi\colon F(x)\to G(x)$, we introduce  the morphism  
$$\phi':=\phi\circ \psi^{-1}\colon G(x)\to y$$
 and take 
$l_j'\in T_{x'}X$ satisfying  $T\phi'(TG(x)l_j')=h_j$. Then,   by definition, 
$$
\omega_{x}(l_1',\ldots,l_k')=:(G_\ast\omega)_y(h_1,\ldots,h_k).
$$
By construction,
$
T\phi(TF(x)l_j)=T\phi'(TG(x)l_j'
$
or equivalently,
$$
TF(x)l_j = T(\phi^{-1}\circ\phi')(TG(x)l_j').
$$
By Proposition \ref{prop4.14} there exists a morphism
$\sigma\colon F(x)\rightarrow G(x)$ such that 
$$
T\sigma\circ TF(x)h=TG(x)h\quad \text{for all $h\in T_xX$}.
$$
Hence
$$
TG(x)l_j=T\sigma\circ TF(x)l_j= T\sigma \circ T(\phi^{-1}\circ\phi')(TG(x)l_j')
$$
for all $j$. Since  $\sigma\circ \phi^{-1}\circ\phi'\colon G(x)\rightarrow G(x)$ 
in ${\bm{Y}}$ and since $G$ is an equivalence
there exists a unique morphism 
$\tau\colon x\rightarrow x$ satisfying 
$$
G(\tau)=\sigma\circ \phi^{-1}\circ\phi'.
$$
We consider the germs $t_Y\circ (s_Y\vert {\bm{U}}(\sigma\circ \phi^{-1}\circ\phi'))^{-1}$ and $t_X\circ (s_X\vert {\bm{U}}(\tau))^{-1}$. 
Because  $G$ is an equivalence we must have  the relationship 
$$
t_Y\circ (s_Y\vert {\bm{U}}(\sigma\circ \phi^{-1}\circ\phi'))^{-1}\circ G(b) = G(t_X\circ (s_X\vert {\bm{U}}(\tau))^{-1}(b))
$$
for $b$ near $x$. Differentiating in $b$ at $b=x$  leads to 
$$
T(\sigma\circ \phi^{-1}\circ\phi')\circ TG(x)h= TG(x)T\tau(h)
$$
for all $h\in T_xX$. In view of \eqref{oip} we arrive at
$
TG(x)l_j = TG(x)\circ T\tau(l_j')
$ for all $j$. Since $TG(x)$ is an isomorphism,  
$$
l_j= T\tau(l_j')\quad  \text{for all $j=1,\ldots ,k$.}
$$
From the compatibility of $\omega$ with morphisms we  conclude
$$
\omega(l_1,\ldots ,l_k)=\omega(T\tau l_1',\ldots ,T\tau l_k')=\omega(l_1',\ldots ,l_k'),
$$
which implies the desired equality
$
F_\ast\omega=G_\ast\omega.
$
The proof of Proposition \ref{dependence_on_F} is complete.
 \qed \end{proof}

 \begin{proposition}\label{dependence_f_ast}\index{C- Dependence of $\mathfrak{f}_{\ast}$}\index{C- Dependence of $\mathfrak{f}^{\ast}$}
 Assume that $\mathfrak{f},\mathfrak{g}:X\rightarrow Y$ are sc-smooth generalized isomorphism satisfying $|\mathfrak{f}|=|\mathfrak{g}|$.
 Then 
 $$\text{$\mathfrak{f}_\ast=\mathfrak{g}_\ast$\quad   and \quad $\mathfrak{f}^\ast=\mathfrak{g}^\ast$}.$$
 \end{proposition}
 \begin{proof}
 Let $\mathfrak{f}=[X\xleftarrow{F}A\xrightarrow{G}Y]$ and 
$\mathfrak{g}=[X\xleftarrow{F'}A'\xrightarrow{G'}Y]$. We consider the two equivalences 
$F\colon A\to X$ and $F'\colon A'\to X$.
According to Theorem \ref{modified_weak_fibered_product}, the associated weak fibered  
 product $A\times_X A'$ is an ep-groupoid and the two functors
 $$F\circ \pi_A,F'\circ \pi_{A'}\colon 
A\times_X A'\to X$$
are naturally equivalent. Moreover,
the functors 
$\pi_A\colon A\times_XA'\to A$ and 
$\pi_{A'}\colon A\times_XA'\to A'$ are equivalences. Hence the compositions $F\circ \pi_A$ and $F'\circ \pi_{A'}$ are equivalences.
Since two equivalences between ep-groupoids which are naturally equivalent induce the same maps on the orbit spaces, we have 
\begin{equation}\label{equation2.33}
\abs{F\circ \pi_A}=
\abs{F'\circ \pi_{A'}}\colon \abs{A\times_XA'}\to \abs{X}.
\end{equation}
In addition, in view of Lemma \ref{equivalence_forms},
$$(F\circ \pi_A)^\ast=(F'\circ \pi_{A'})^\ast.$$
Using  that $(\pi_A)_\ast\circ (\pi_A)^\ast=Id$, by Lemma \ref{equivalence_forms} (3), we compute for a given sc-differential form $\omega$ on $X$,
\begin{equation}\label{equation2.34}
\begin{split}
G_\ast \circ F^\ast \omega&=
G_\ast \circ (\pi_A)_\ast\circ (\pi_A)^\ast\circ F^\ast \omega\\
&=(G \circ \pi_A)_\ast\circ (F\circ \pi_A)^\ast\omega\\
&=(G \circ \pi_A)_\ast\circ (F'\circ \pi_{A'})^\ast\omega.
\end{split}
\end{equation}
Considering now $G'\circ \pi_{A'}$, we use the assumption that 
$|G|\circ |F|^{-1}= |G'|\circ |F'|^{-1}$, use \eqref{equation2.33}
and compute,
\begin{equation*}
\begin{split}
\abs{G'\circ \pi_{A'}}&=\abs{G'}\circ \abs{\pi_{A'}}\\
&=\abs{G}\circ \abs{F}^{-1}\circ \abs{F'}\circ \abs{\pi_{A'}}\\
&=\abs{G}\circ \abs{F}^{-1}\circ \abs{F'\circ\pi_{A'}}\\
&=\abs{G}\circ \abs{F}^{-1}\circ \abs{F\circ \pi_{A}}\\
&=\abs{G}\circ \abs{F}^{-1}\circ \abs{F}\circ \abs{\pi_{A}}\\
&=\abs{G\circ \pi_A}.
\end{split}
\end{equation*}
Consequently, we can apply Proposition \ref{dependence_on_F} and obtain using \eqref{equation2.34}, Lemma \ref{equivalence_forms} and Lemma \ref{differential_form_1},
\begin{equation*}
\begin{split}
G_\ast F^\ast\omega&= (G\circ \pi_A)_\ast \circ (F'\circ\pi_{A'})^\ast\omega\\
&=(G'\circ\pi_{A'})_\ast \circ (F'\circ\pi_{A'})^\ast\omega\\
&=(G')_\ast \circ (\pi_{A'})_\ast \circ (\pi_{A'})^\ast \circ (F')^\ast\omega\\
& = (G')_\ast \circ (F')^\ast\omega.
\end{split}
\end{equation*}
We have verified the first assertion of Proposition \ref{dependence_f_ast}. The  second assertion is proved the same way and the proof of Proposition \ref{dependence_f_ast} is complete.
  \qed \end{proof}
 With Proposition \ref{dependence_f_ast} also the last statement (4) in Theorem \ref{OOmaine} is verified, so that the proof of Theorem \ref{OOmaine} is  finally complete.\hfill \qed

   \section{Branched Ep\texorpdfstring{$^+$}{Ep}-Subgroupoids and Equivalences}
 
   We shall study the behavior of branched ep$^+$subgroupoid under equivalences and generalized isomorphisms.
  We begin with a categorical definition.
   \begin{definition}\index{D- Pushforward and Pullback of $\Theta$}
   Let $F:X\rightarrow Y$ be an equivalence between ep-groupoids and $\Theta:X\rightarrow{\mathbb Q}^+$ and $\Theta':Y\rightarrow {\mathbb Q}^+$  just functors
   \begin{itemize}
   \item[(1)]\ The {\bf pullback} $F^\ast\Theta':X\rightarrow {\mathbb Q}^+$ is defined by $F^\ast\Theta' :=\Theta'\circ F$.
   \item[(2)]\ The {\bf push-forward} $F_\ast\Theta:X'\rightarrow {\mathbb Q}^+$ is defined as follows. For $y\in Y$ pick
   $x\in X$ and $\psi\in \bm{Y}$ with $s(\psi)=F(x)$ and $t(\psi)=y$. Then set $(F_\ast\Theta)(y):=\Theta(x)$.
   \end{itemize}
   \end{definition}
   Clearly the pull-back is well-defined as a functor, but the push-forward might depend on the choices involved.
   \begin{lemma}
  {\em(1)} Let $F:X\rightarrow Y$ be an equivalence and $\Theta:X\rightarrow {\mathbb Q}^+$ a functor 
   Then $F_\ast\Theta$ is well-defined as a functor.  Moreover $F^\ast F_\ast \Theta=\Theta$.
   
   \noindent {\em(2)} Let $\Theta':Y\rightarrow {\mathbb Q}^+$ be a functor.  Then $F_\ast F^\ast\Theta'=\Theta'$.
   
   \noindent{\em(3)} Let $\Theta$ and $\Theta'$ be as in (b) and assume that $F,G:X\rightarrow Y$ are equivalences which are naturally equivalent.
   Then $F_\ast\Theta=G_\ast\Theta$ and $F^\ast\Theta'=G^\ast\Theta'$.
   \end{lemma}
   \begin{proof}
   In order to show that $F_\ast\Theta$ is well-defined, assume that $x,x'\in X$ and $\psi,\psi'\in\bm{Y}$ satisfy
   $$
   \psi:F(x)\rightarrow y\ \ \text{and}\ \ \ \psi':F(x')\rightarrow y.
   $$
   Then $\psi^{-1}\circ\psi':F(x')\rightarrow F(x)$. Since $F$ is an equivalence we find $\sigma:x'\rightarrow x$ satisfying $F(\sigma)=\psi^{-1}\circ\psi'$.
   From this we conclude that
   $$
   \Theta(x)=\Theta(x'),
   $$
   which show that $F_\ast\Theta$ is well-defined.  
    The properties $F_\ast F^\ast\Theta'=\Theta'$ and $F^\ast F_\ast\Theta=\Theta$ are immediate (b) as well 
    as the properties stated in (3).
   \qed \end{proof}
   If $\mathfrak{f}:X\rightarrow Y$ is a generalized isomorphism, it is represented by a diagram $X\xleftarrow{F}A\xrightarrow{G}Y$,
   where $F$ and $G$ are equivalences. We define $\mathfrak{f}_\ast = G_\ast F^\ast$ and $\mathfrak{f}^\ast = F_\ast G^\ast$ and as we shall see   this is well-defined independent of the diagram we picked to represent $\mathfrak{f}$.
   Indeed, if we have to equivalent diagrams representing $\mathfrak{f}$ there is a common refinement resulting in 
   $$
   \begin{CD}
   X @<F <<  A @>G>> Y\\
   @.   @A H AA  @.\\
   X@<F''<<  A''  @> G''>> Y\\
   @.  @V H'VV  @.\\
   X @<F' <<  A' @> G' >> Y
   \end{CD}
   $$
   Here $F\circ H \simeq F''$, $G\circ H\simeq G''$, $F'\circ H'\simeq F''$, and $G'\circ H'\simeq G''$.
   This implies that 
   \begin{eqnarray*}
   && G_\ast F^\ast\Theta\\
   &=& G_\ast H_\ast H^\ast F^\ast\Theta\\
   &=& G''_\ast F''^\ast \Theta\\
   &=& (G'\circ H')_\ast (F'\circ H')^\ast\Theta\\
   &=& G'_\ast F'^\ast\Theta,
   \end{eqnarray*}
   so that $\mathfrak{f}_\ast\Theta$ is well-defined. Similarly $\mathfrak{f}^\ast\Theta'$ is well-defined.
   It is easily verified that for $\mathfrak{f}:X\rightarrow Y$ and $\mathfrak{g}:Y\rightarrow Z$ we have the properties
   \begin{eqnarray*}
   &(\mathfrak{g}\circ\mathfrak{f})_\ast=\mathfrak{g}_\ast\circ \mathfrak{f}_\ast&\\
   &(\mathfrak{g}\circ\mathfrak{f})^\ast=\mathfrak{f}^\ast\circ \mathfrak{g}^\ast&\\
   &\mathfrak{f}_\ast^{-1}=\mathfrak{f}^\ast&.
   \end{eqnarray*}
   The considerations so far were purely algebraic for arbitrary functors $\Theta$ and $\Theta'$.
   Next we consider $\Theta$ and $\Theta'$ which are branched ep$^+$-subgroupoids.
   As we shall see this property is preserved by the generalized isomorphisms.
   The first basic result in this section is the following theorem.
    \begin{theorem}\label{THM1133}\index{T- Branched ep$^+$-subgroupoids under generalized isomorphisms}
    Let $\mathfrak{f}:X\rightarrow Y$ be a generalized isomorphism between two ep-groupoids. 
    \begin{itemize}
    \item[{\em(1)}]\ If $\Theta:X\rightarrow {\mathbb Q}^+$ is a branched  ep$^+$-subgroupoid then $\mathfrak{f}_\ast\Theta:Y\rightarrow {\mathbb Q}^+$ is a branched ep$^+$-subgroupoid.
    \item[{\em(2)}]\ If If $\Theta':Y\rightarrow {\mathbb Q}^+$ is a branched ep$^+$-subgroupoid then $\mathfrak{f}_\ast\Theta:X\rightarrow {\mathbb Q}^+$ is  a branched ep$^+$-subgroupoid.
    \end{itemize}
    Moreover, the following holds for the maps $\mathfrak{f}_\ast$ and $\mathfrak{f}^\ast$.
    \begin{itemize}
   \item[{\em (3)}]\  If $\Theta$ and $\Theta'$ are closed, see Definition \ref{DEF915}(i), so are $\mathfrak{f}_\ast\Theta$ and $\mathfrak{f}^\ast\Theta'$.
  \item[{\em (4)}]\  If $\Theta$ and $\Theta'$ are compact, see Definition \ref{DEF915}(ii), so are $\mathfrak{f}_\ast\Theta$ and $\mathfrak{f}^\ast\Theta'$.
  \item[{\em(5)}]\ The maps $\mathfrak{f}_\ast$ and $\mathfrak{f}^\ast$ also preserve the property of being of manifold-type or orbifold-type, respectively.
  For the definitions see Definition \ref{DEF917}.
  \item[{\em (6)}]\ The maps $\mathfrak{f}_\ast$ and $\mathfrak{f}^\ast$ also preserve the dimensional decomposition given in Proposition \ref{PROP919}.
   \end{itemize}
    \end{theorem}
   \begin{proof}
   The basic considerations which lead to the desired results are concerned with the behavior of branched 
   ep$^+$-subgroupoids with respect to equivalences $F:X\rightarrow Y$.   
   
   First consider $F^\ast\Theta'= \Theta'\circ F$.  Let $x\in \supp(F^\ast\Theta')$. We find open neighborhoods
   $U(x)$ and $U(F(x))$ so that $F:U(x)\rightarrow U(F(x))$ is an sc-diffeomorphism and $\Theta'$ is represented on $U(F(x))$
   by a local branching structure ${(M_i')}_{i\in I}$ , ${(\sigma_i)}_{i\in I}$, i.e.
   $$
   \Theta'(y)=\sum_{\{i\in I\ |\ y\in M_i'\}}\sigma_i\ \ \text{for}\ \ y\in U(F(x)).
   $$
   Using $F|U(x)$ we can pull back  the $M_i'$, defining $M_i:= (F|U(x))^{-1}(M_i')$. We keep the weights.
   Then 
   $$
   (F^\ast\Theta')(z) = \Theta'\circ F(z)=\sum_{\{i\in I\ |\ z\in M_i\}}\sigma_i\ \ \text{for}\ \ z\in U(x).
   $$
   This shows that $F^\ast\Theta'$ is a branched ep$^+$-subgroupoid.
   
   Next we consider $F_\ast\Theta$. Pick $y\in \supp(F_\ast\Theta)$. We find $x\in\supp(\Theta)$ and a morphism
   $\psi:F(x)\rightarrow y$.  The local sc-diffeomorphism $\wh{\psi}$ between neighborhoods of $F(x)$ and $y$ 
   will map a local branching structure around $F(x)$ to a local branching structure around $y$.
   Since $F$ is a local sc-diffeomorphisms a local branching structure around $x$ is mapped to one around $F(x)$.
   Using $F$ and $\wh{\psi}$ we can map a local branching structure around $x$ for $\Theta$ to one 
   around $y$ for $F_\ast\Theta$.
   
 If $\mathfrak{f}$ is represented  by
 $$
d\colon  X\xleftarrow{F} A\xrightarrow{G} Y
 $$
 then $\mathfrak{f}_\ast \Theta =G_\ast F^\ast\Theta$ which implies that $\mathfrak{f}_\ast\Theta$ is a branched ep$^+$-subgroupoid.
 Similarly for $\mathfrak{f}^\ast\Theta'= F_\ast G^\ast\Theta'$ in view of the previous discussion. This completes the proofs of (1) and (2).
 
 A branched ep$^+$-subgroupoid $\Theta:X\rightarrow {\mathbb Q}^+$ defines a map
 $|\Theta|:|X|\rightarrow {\mathbb Q}^+$ defined by $|\Theta|(|x|):=\Theta(x)$. We note the identity
 $$
 |\mathfrak{f}_\ast\Theta|= |\Theta|\circ |\mathfrak{f}|^{-1}.
 $$
If  $\Sigma:=\{z\in |X|\ |\ |\Theta|(z)>0\}$ is closed or compact, the same holds for $|\mathfrak{f}|(\Sigma)=\{q\in |Y|\ |\ |\mathfrak{f}_\ast\Theta|(q)>0\}$. A similar argument holds for $\mathfrak{f}^\ast\Theta'$.  This proves (3) and (4).

If $\Theta$ or $\Theta'$  only takes the values $0$ and $1$ the same holds for $\mathfrak{f}_\ast\Theta$ and $\mathfrak{f}^\ast\Theta'$.
This proves that $\mathfrak{f}_\ast$ and $\mathfrak{f}^\ast$ preserve the property of having orbifold-type.
Assume that between two objects $x,x'\in \supp(\Theta)$ there are is at most one morphism.  Let $y,y'\in \supp(\mathfrak{f}_\ast\Theta)$.
If there exist two different morphisms $y\rightarrow y'$ it follows that there exists a nontrivial morphism $\phi':y\rightarrow y$.
If $\psi:F(x)\rightarrow y$. It follows that $\psi\circ\phi'\circ\psi^{-1}:F(x)\rightarrow F(x)$ is a nontrivial isomorphism.
Then $\sigma:x\rightarrow X$ with $F(\sigma)= \psi\circ\phi'\circ\psi^{-1}$ is nontrivial and $x\in \supp(\Theta)$, which gives 
a contradiction. The same argument works for $F^\ast\Theta'$. Hence $\mathfrak{f}_\ast$ and $\mathfrak{f}^\ast$ preserves
manifold-type.   The point (4) is trivial. This completes the proof of the theorem.
   \qed \end{proof}

   Next we consider the behavior of orientations. We begin with a quite algebraic consideration.
   Consider the ep-groupoids $X$ and $Y$ with associated $\pi_X:\wh{\text{Gr}}(X)\rightarrow X_\infty$ and 
   $\pi_Y:\wh{\text{Gr}}(y)\rightarrow Y_\infty$. Denote by $\Gamma(\pi_X)$ and $\Gamma(\pi_X)$ the space of section functors  (no continuity requirement).
   Assume that $F:X\rightarrow Y$ is an equivalence of ep-groupoids. We define maps 
   $$
   F_\ast:\Gamma(\pi_X)\rightarrow \Gamma(\pi_Y)\ \ \text{and}\ \   F^\ast:\Gamma(\pi_Y)\rightarrow \Gamma(\pi_X)
   $$
   as follows.  Given $h\in \Gamma(\pi_X)$ we set  $h'=F_\ast h\in\Gamma(\pi_Y)$ by
\begin{eqnarray}\label{EQN118}
   h'(y)= T\psi \circ TF(x) h(x),
  \end{eqnarray}
   where $\psi:F(x)\rightarrow y$. Since $F$ is an equivalence we find for given $y\in Y_\infty$ such an $x$ and $\psi$. 
   The map $T\psi\circ TF(x):T_xX\rightarrow T_yY$ is a linear sc-isomorphism, which maps the formal finite sum 
   in the obvious way. Since $h$ is functorial and $F$ an equivalence the definition does not depend on the choices involved.
   If $F,G:X\rightarrow Y$ are naturally equivalent functors one immediately deduces that $G_\ast h=F_\ast h$. 
   
   Similarly we can define for $h'\in\Gamma(\pi_Y)$ the element $h=F^\ast h'$ by
\begin{eqnarray}\label{EQN119}
   TF(x) h(x) =h'(F(x)),
\end{eqnarray}
   and again for naturally equivalent $F$ and $G$ it holds that $F^\ast h'=G^\ast h'$.
   
   Next consider a generalized isomorphism $\mathfrak{f}:X\rightarrow Y$. It is represented by an equivalence of diagrams, where the equivalence
   is defined by `common refinement'. From our discussion above it follows immediately that if $\mathfrak{f}$ is represented by 
   a diagram 
   $$
   d\colon X\xleftarrow{F}A\xrightarrow{G} Y.
   $$
   then $\mathfrak{f}_\ast h:= G_\ast F^\ast h$ and $\mathfrak{f}^\ast h':= F_\ast G^\ast h'$ are well-defined and do not depend
   on the specific choice of representative diagram. Hence we have proved.
   \begin{lemma}
   Given a generalized isomorphism $\mathfrak{f}:X\rightarrow Y$ between ep-groupoids
   there are well-defined maps 
   $$
   \mathfrak{f}_\ast:\Gamma(\pi_X)\rightarrow \Gamma(\pi_Y)\ \ \text{and}\ \ \ \mathfrak{f}^\ast:\Gamma(\pi_Y)\rightarrow \Gamma(\pi_X)
   $$
   such hat $\mathfrak{f}_\ast^{-1}=\mathfrak{f}^\ast$.  If $\mathfrak{g}:Y\rightarrow Z$ is another generalized isomorphism
   then $(\mathfrak{g}\circ\mathfrak{f})_\ast=\mathfrak{g}_\ast\circ \mathfrak{f}_\ast$ and $(\mathfrak{g}\circ\mathfrak{f})^\ast=\mathfrak{f}^\ast\circ\mathfrak{g}^\ast$.  If $\mathfrak{f}$ is represented by the diagram
   $$
   d\colon X\xleftarrow{F} A\xrightarrow{G} Y
   $$
   then $\mathfrak{f}_\ast=G_\ast F^\ast$ and $\mathfrak{f}^\ast = F_\ast G^\ast$, where for an equivalence
   the push-forward and pullback are defined by (\ref{EQN118}) and (\ref{EQN119}), respectively.
   \end{lemma}

   Given an oriented $\wh{\Theta}:X\rightarrow {\mathbb Q}^+$ the associated $\wh{\mathsf{T}}_\Theta$ is a section functor
   in $\Gamma(\pi_X)$. If $\mathfrak{f}:X\rightarrow Y$ is a generalized isomorphism one would like to equip
   the branched ep$^+$-subgroupoid $\mathfrak{f}_\ast\Theta$ with $\mathfrak{f}_\ast \wh{\mathsf{T}}_\Theta\in\Gamma(\pi_Y)$.
   The only thing one needs to check is that it has the local representation property, so that it is indeed an orientation.
  Similarly there is a construction involving $\mathfrak{f}^\ast$.
   Our second main result in this section is the following theorem.
   \begin{theorem}\label{THM1135}\index{T- Transformations of orientation}
   Let $\mathfrak{f}:X\rightarrow Y$ be a generalized isomorphism and $(\Theta,\wh{\mathsf{T}}_\Theta)$ an oriented branched ep$^+$subgroupoid on $X$ 
   and $(\Theta',\wh{\mathsf{T}}_{\Theta'})$ an oriented branched ep$^+$-subgroupoid on $Y$.
   Then $\mathfrak{f}_\ast\wh{\mathsf{T}}_\Theta$ is an orientation for $\mathfrak{f}_\ast\Theta$ and $\mathfrak{f}^\ast\wh{\mathsf{T}}_{\Theta'}$
   an orientation for $\mathfrak{f}^\ast\Theta'$.
     \end{theorem}
   \begin{proof}
   The following arguments show that $\mathfrak{f}_\ast \wh{\mathsf{T}}_\Theta$ and $\mathfrak{f}^\ast\wh{\mathsf{T}}_{\Theta'}$ are locally representable
   by the tangents of oriented local branching structures.
   
    Let $\wh{\mathsf{T}}_\Theta:X_\infty\rightarrow \wh{\text{Gr}}(X)$ be an orientation for $\Theta:X\rightarrow {\mathbb Q}^+$. 
   We already know that $F_\ast\Theta:Y\rightarrow {\mathbb Q}^+$ is a branched ep$^+$-subgroupoid, when $F:X\rightarrow Y$ is an equivalence.
   We define $F_\ast\wh{\mathsf{T}}_\Theta:Y_\infty\rightarrow \wh{\text{Gr}}(Y)$ as follows. For $y\in Y_\infty$ we take $x\in X_\infty$ 
   so that $|F(x)|=|y|$. Pick an isomorphism $\psi:F(x)\rightarrow y$ and pick open neighborhoods $U(x)$, $V(F(x))$, and $V(y)$
   so that
   $$
   F:U(x)\rightarrow V(F(x))\ \ \text{and}\ \ \wh{\psi}:V(F(x))\rightarrow V(y)
   $$
   are sc-diffeomorphisms and $U(x)$ supports an oriented branching structure 
   $$
   {(M_i,o_i)}_{i\in I},\ \ {(\sigma_i)}_{i\in I}
   $$
    representing $\Theta$ and
   $\wh{\mathsf{T}}_\Theta$. We map this data by the sc-diffeo\-mor\-phism   $\wh{\Psi}\circ (F|U(x))$ to obtain
   $(M_i',o_i') := \wh{\psi}\circ F(M_i,o_i)$, and we keep the weights. Then for $z\in V(y)$ we define 
   $$
   (F_\ast\wh{\mathsf{T}}_\Theta)(z)=\sum_{\{i\in I\ |\ z\in M_i\}} \sigma_i\cdot T_z(M_i',o_i').
   $$
 Again this definition involves choices, and we have to show that the definition does not depend on them.  Hence assume $\psi':F(x')\rightarrow y$
 and pick a local branching structure on $U(x')$ representing near $x'$ the oriented $\wh{\Theta}$. We denote it by ${(N_j,o_j)}_{j\in J}$, ${(\tau_j)}_{j\in J}$.
Then $x$ and $x'$ are isomorphic by $\gamma:x'\rightarrow x$ defined by $\psi\circ\gamma=\psi'$.
 We  can map the latter local branching structure  via $\wh{\gamma}$ to a neighborhood 
 of $x$, where it will represent the oriented $\wh{\Theta}$. After perhaps restricting $U(x)$ we  map it via $\wh{\psi}$ and obtain another local branching 
 structure near $y$.  It is clear that near $y$ it defines the same expression. Of course, the result is the same as using $\wh{\psi}'\circ F|U(x')$
 and this shows that $F_\ast\wh{\mathsf{T}}_\Theta$ is well-defined. 
 
  Assume next that $F,G:X\rightarrow Y$ are equivalences which are naturally isomorphic $\tau:F\rightarrow G$. Clearly
 $F_\ast\Theta=G_\ast\Theta$ and by an argument as above we also see that $F_\ast\wh{\mathsf{T}}_\Theta=G_\ast\wh{\mathsf{T}}_\Theta$.
 
 Assume that $\Theta':Y\rightarrow {\mathbb Q}^+$ is a branched ep$^+$-subgroupoid and $F:X\rightarrow Y$ an equivalence. 
 The orientation of $F^\ast\Theta'$ is defined by $F^\ast\wh{\mathsf{T}}_\Theta:X_\infty\rightarrow \wh{\text{Gr}}(X)$ as follows. 
 Given $x\in X_\infty$ we pick open neighborhoods $U(x)$ and $V(F(x))$, so that $F:U(x)\rightarrow V(F(x))$ is an sc-diffeomorphism.
 If $V(y)$ is small enough we take an oriented local branching structure and its preimage. This will define 
 $F^\ast\wh{\mathsf{T}}_\Theta$ over $U(x)_\infty$. Assume that $F,G:X\rightarrow Y$ are equivalences which are naturally equivalent by a natural transformation.
 Clearly $G^\ast\Theta'=F^\ast\Theta'$ and $F^\ast\wh{\mathsf{T}}_\Theta=G^\ast\wh{\mathsf{T}}_\Theta$.
 
 If $\mathfrak{f}:X\rightarrow Y$ is a generalized isomorphism  it follows from the previous discussion
 that taking any representative diagram
 $$
 d\colon X\xleftarrow{F} A\xrightarrow{G} Y
 $$
 the expressions $G_\ast F^\ast\wh{\mathsf{T}}_\Theta$ and $G^\ast F_\ast\wh{\mathsf{T}}_{\Theta'}$ are well-defined and do not depend
 on the particular  choice of the representatives.  Hence $\mathfrak{f}_\ast\wh{\mathsf{T}}_\Theta$ and $\mathfrak{f}^\ast\wh{\mathsf{T}}_{\Theta'}$
 are well-defined.   Our discussion shows that orientations can be transformed by generalized isomorphisms. 
 \qed \end{proof}
 In view of Theorem \ref{THM1135} we can make the following definition.
 \begin{definition}\index{D- Transformation of $\wh{\Theta}$}
 Let $\wh{\Theta}=(\Theta,\wh{\mathsf{T}}_\Theta):X\rightarrow {\mathbb Q}^+$ and $\wh{\Theta}':Y\rightarrow {\mathbb Q}^+$
 be oriented, branched ep$^+$-subgroupoids. Assume that $\mathfrak{f}:X\rightarrow Y$ is a generalized isomorphism.
 Then $\mathfrak{f}_\ast\wh{\Theta}$ and $\mathfrak{f}^\ast\wh{\Theta}'$ are defined by
 \begin{eqnarray*}
 &\mathfrak{f}_\ast\wh{\Theta}=(\mathfrak{f}_\ast\Theta,\mathfrak{f}_\ast\wh{\mathsf{T}}_\Theta)&\\
 &\mathfrak{f}^\ast\wh{\Theta}' =(\mathfrak{f}^\ast\Theta',\mathfrak{f}^\ast\wh{\mathsf{T}}_{\Theta'}).&
 \end{eqnarray*}
 \qed
 \end{definition}
 If $\mathfrak{f}:X\rightarrow Y$ and $\mathfrak{g}:Y\rightarrow Z$ are generalized isomorphisms its is clear from the definition that
 the operations on the oriented, branched ep$^+$-subgroupoids have the following functorial property
\begin{eqnarray}
& (\mathfrak{g}\circ\mathfrak{f})_\ast \wh{\Theta}=\mathfrak{g}_\ast\circ \mathfrak{f}_\ast\wh{\Theta}\ \ \text{and}\ \ \  (\mathfrak{g}\circ\mathfrak{f})^\ast\wh{\Theta}'=\mathfrak{f}^\ast\circ \mathfrak{g}^\ast\wh{\Theta}'&\\
&\mathfrak{f}^\ast\mathfrak{f}_\ast\wh{\Theta}=\wh{\Theta}\ \ \text{and}\ \ \ \mathfrak{f}_\ast\mathfrak{f}^\ast\wh{\Theta}'=\wh{\Theta}'.&\nonumber
\end{eqnarray}
Next we consider the tangent of $\Theta:X\rightarrow {\mathbb Q}^+$ and its behavior under equivalences.
\begin{theorem}\index{T- Transformations of $T\Theta$}
Let $\mathfrak{f}:X\rightarrow Y$ be a generalized isomorphism and $\Theta:X\rightarrow {\mathbb Q}^+$ and $\Theta':Y\rightarrow {\mathbb Q}^+$
branched ep$^+$-subgroupoids. Then 
$$
(T\mathfrak{f})_\ast T\Theta = T(\mathfrak{f}_\ast \Theta)\ \ \text{and}\ \ \ (T\mathfrak{f})^\ast(T\Theta')= T(\mathfrak{f}^\ast\Theta').
$$
\end{theorem}
\begin{proof}
If $\mathfrak{f}:X\rightarrow Y$ is a generalized isomorphism which is representable by the diagram
$$
d\colon X\xleftarrow{F} A\xrightarrow{G} Y,
$$
then we obtain the  generalized isomorphism $T\mathfrak{f}:TX\rightarrow TY$ by
$$
T\mathfrak{f}:= [Td]:= [TX\xleftarrow{TF} TA\xrightarrow {TG} TY].
$$
Recalling that $T\Theta$ is locally represented by the tangents of the manifolds occurring in the local branching structure 
the result is easily established. We leave the further details to the reader.
\qed \end{proof}

Next we study the relationship between the boundary construction $\partial\Theta$ and
the maps $\mathfrak{f}_\ast$ and $\mathfrak{f}^\ast$. We shall prove the following result
\begin{theorem}\label{THMX1138}\index{T- Boundary operation and equivalences}
Let $\mathfrak{f}:X\rightarrow Y$ be  a generalized isomorphism between ep-groupoids.
Assume that $\Theta:X\rightarrow {\mathbb Q}^+$ and $\Theta':Y\rightarrow {\mathbb Q}^+$ are branched ep$^+$-subgroupoids
with associated  boundaries $\partial\Theta$ and $\partial\Theta'$.
Then  we have the equalities
\begin{eqnarray}
&\mathfrak{f}_\ast(\partial\Theta) =\partial (\mathfrak{f}_\ast\Theta)&\\
&\mathfrak{f}^\ast(\partial\Theta') = \partial(\mathfrak{f}^\ast\Theta').&\nonumber
\end{eqnarray}
\end{theorem}
\begin{proof}
In order to do so,  we have to understand as before the behavior with respect to equivalences $F:X\rightarrow Y$,
as well as $G:X\rightarrow Y$ provided $F$ and $G$ are naturally equivalent. 
The basic idea is as before that the local branching structures can be moved around by equivalences,
and the local sc-diffeomorphism $\wh{\psi}$ associated to a morphism $\psi$.
If $M\subset X$ is a submanifold and $x\in M$ we can move small neighborhoods around in the same way
and not that $d_M(x)$ is being preserved. 
The proof of the theorem follows from these considerations and is left to the reader.
\qed \end{proof}

Let $\wh{\Theta}:X\rightarrow {\mathbb Q}^+$ be an oriented, tame, branched ep$^+$-subgroupoid. 
Then we have seen that $\partial\Theta$ has a natural  orientation $\wh{\mathsf{T}}_{\partial\Theta}$, so that we can define
$\partial\wh{\Theta}$, see Definition \ref{DEF936}.  The section $\wh{\mathsf{T}}_{\partial\Theta}$ can be pushed forward 
or pulled back by a generalized isomorphism. According to the following theorem  these procedures define the corresponding
orientations for the transformed $\partial\Theta$.

\begin{theorem}\index{T- Transformations and boundaries}\label{THMX11.3.9}
Let $\mathfrak{f}:X\rightarrow Y$ be a generalized isomorphism between ep-groupoids and $\wh{\Theta}:X\rightarrow {\mathbb Q}^+$
and $\wh{\Theta}':Y\rightarrow {\mathbb Q}^+$ oriented, tame, branched ep$^+$-subgroupoids.  Then the following formulae hold.
\begin{eqnarray*}
&\partial(\mathfrak{f}_\ast\wh{\Theta}) =\mathfrak{f}_\ast(\partial\wh{\Theta})&\\
&\partial(\mathfrak{f}^\ast\wh{\Theta}')=\mathfrak{f}^\ast(\partial\wh{\Theta}').&
\end{eqnarray*}\
\qed
\end{theorem}
We leave the proof, which follows the lines of previous proofs in this section,  to the reader.

 \section{Equivalences and Integration}\label{SECRTY114}
The ingredients for integration, sc-differential forms and oriented, branched ep$^+$-subgroupoids behave well under equivalences.
Recall the statement of Theorem \ref{IandS}, which asserts that given an ep-groupoid $X$, we can associate
to an oriented, tame, compact, branched ep$^+$-subgroupoid $\wh{\Theta}:X\rightarrow {\mathbb Q}^+$ of dimension $n$
a natural linear map
$$
\Phi_{\wh{\Theta}}:\Omega^n_\infty(X)\rightarrow {\mathcal M}(S,{\mathcal L}(S)),
$$
where $S=|\supp(\Theta)|$, associating to an sc-differential $n$-form $\omega$ a signed measure $\mu_\omega^\Theta$.
\begin{theorem}[Transformation of canonical measures]\label{THM1141} \index{T- Transformation of canonical measures}
Assume that $\mathfrak{f}:X\rightarrow Y$ is a generalized isomorphism between two ep-groupoids and
$\wh{\Theta}:X\rightarrow {\mathbb Q}^+$ and $\wh{\Theta}':Y\rightarrow {\mathbb Q}^+$ are oriented, tame,  compact, branched
ep$^+$-subgroupoids of dimension $n$. Suppose further that $\omega$ and $\omega'$ are sc-differential forms 
on $X$ and $Y$, respectively, of degree $n$. Then, for a measurable subset $K$ of $|\supp(\Theta)|$
$$
\mu^{\wh{\Theta}}_{\omega} (K)  = \mu^{\mathfrak{f}_\ast\wh{\Theta}}_{\mathfrak{f}_\ast\omega}( |\mathfrak{f}|(K))\ \ \text{and}\ \ \
\mu^{\mathfrak{f}^\ast\wh{\Theta}'}_{\mathfrak{f}^\ast\omega'}(K)   = \mu^{\wh{\Theta}'}_{\omega'}( |\mathfrak{f}|(K)).
$$
With the obvious notation we can rewrite this more suggestively as
$$
|\mathfrak{f}|_\ast\mu^{\wh{\Theta}}_{\omega}   = \mu^{\mathfrak{f}_\ast\wh{\Theta}}_{\mathfrak{f}_\ast\omega}\ \ \text{and}\ \ \
\mu^{\mathfrak{f}^\ast\wh{\Theta}'}_{\mathfrak{f}^\ast\omega'}   = |\mathfrak{f}|^\ast \mu^{\wh{\Theta}'}_{\omega'}.
$$
\qed
\end{theorem}
The theorem will be proved later.
\begin{remark}\label{REM1142}\index{R- Remarks on polyfold structures}
An obvious consequence of this result is the following fact, which will be elaborated on further
in Chapter \ref{chap11+}.  Let $Z$ be a metrizable space and consider a pair  $(X,\beta)$, where $X$ is an ep-groupoid
and $\beta:|X|\rightarrow Z$ is a homeomorphism. We shall refer to such a pair $(X,\beta)$ as a {\bf polyfold structure}\index{Polyfold structure on $Z$}
on the metrizable space $Z$. Given another such pair $(X',\beta')$ we say that it defines an {\bf equivalent polyfold structure}
\index{Equivalent polyfold structure}
on $Z$ provided there exists a generalized isomorphism $\mathfrak{f}:X\rightarrow X'$ satisfying
$$
\beta'\circ |\mathfrak{f}| =\beta.
$$
With other words $|\mathfrak{f}|={(\beta')}^{-1}\circ \beta$. Hence, $(X,\beta)$ and $(X',\beta')$ are equivalent if the {\it transition map}
$ {(\beta')}^{-1}\circ\beta :|X|\rightarrow |X'|$ is induced by a generalized isomorphism $\mathfrak{f}:X\rightarrow X'$. Recall that 
$|\mathfrak{f}|=|\mathfrak{g}|$ implies that $\mathfrak{f}=\mathfrak{g}$, see Theorem \ref{THMX10316}, which implies that $\mathfrak{f}$ is unique.

We explain this concept in a special case. Let $X$ be an \'etale proper Lie groupoid, so that between two objects 
there are is at most one morphism. It is an easy exercise that $|X|$ has a (classical)
natural smooth manifold structure. Assume that $Z$ is a metrizable topological space and $\beta:|X|\rightarrow Z$ a homeomorphism.
Then we can equip $Z$ with the unique smooth manifold structure making $\beta$ a diffeomorphism. 
If $(X',\beta')$ is similar and the transition map is induced from a generalized isomorphism $\mathfrak{f}$ it means that the smooth manifold 
structure
induced on $Z$ is the same as that coming from $(X,\beta)$. Hence the formalism in the more general case can be viewed as 
a generalization of the definition of a smooth manifold. In fact it is a generalization of the concept of an orbifold with boundary and corners
(a manifold is a special case).

A metrizable space equipped with an equivalence class of polyfold structures is called a {\bf polyfold}\index{Polyfold}.
A polyfold $Z$ is the generalization of an orbifold to the sc-smooth world.

There are several concepts which can be carried over to polyfolds.   For example
a differential form on $Z$ is an equivalence class $[(X,\beta),\omega]$. Here $\omega$ is an sc-differential form
on $X$ and $((X,\beta),\omega)\sim ((X',\beta',\omega')$ provided for a generalized isomorphism $\mathfrak{f}:X\rightarrow X'$
with $\beta'\circ|\mathfrak{f}|=\beta$ it holds that $\mathfrak{f}^\ast\omega'=\omega$. 
We can also define branched sub$^+$-polyfold of $Z$ as an equivalence class $[(X,\beta),\Theta]$ with the obvious notion of equivalence.
Similarly the notions of being compact, closed, tame, and oriented carry over. 

If $[(X,\beta),\wh{\Theta}]$ is an oriented, compact, tame branched $n$-dimensional sub$^+$-polyfold and $[(X,\beta),\omega]$
an sc-differential $n$-form on $Z$, we can consider the support $S$ of $|\Theta|\circ \beta^{-1}:Z\rightarrow {\mathbb Q}^+$. 
which call the support of $[(X,\beta),{\Theta}]$, and which is defined independently of the representative which we picked.
Then we can defined the Lebesgue $\sigma$-algebra ${\mathcal L}(S)$ and the vector spaces of signed measures
${\mathcal M}(S,{\mathcal L}(S))$. We can associate 
to $[(X,\beta),\wh{\Theta}]$ and $[(X,\beta),\omega]$ the signed measure $\mu_{[(X,\beta),\omega]}^{[(X,\beta),\wh{\Theta}]}$.
The upshot is that all notions which transform well under equivalences define intrinsic objects on the polyfold $Z$.
\qed
\end{remark}

There is a version of Theorem \ref{THM1141} for the boundary case, which relies on Theorem \ref{THM9510}.
\begin{theorem}[Transformation of canonical boundary measures]\label{THM1143} \index{T- Transformation of canonical boundary measures}
Assume that $\mathfrak{f}:X\rightarrow Y$ is a generalized isomorphism between two ep-groupoids and
$\wh{\Theta}:X\rightarrow {\mathbb Q}^+$ and $\wh{\Theta}':Y\rightarrow {\mathbb Q}^+$ are oriented, tame,  compact, branched
ep$^+$-subgroupoids of dimension $n$. Suppose further that $\omega$ and $\omega'$ are sc-differential forms 
on $X$ and $Y$, respectively, of degree $n-1$. Then on $|\supp(\partial \Theta)|$
$$
\mu^{\partial\wh{\Theta}}_{\omega}   = \mu^{\mathfrak{f}_\ast\partial \wh{\Theta}}_{\mathfrak{f}_\ast\omega}\circ |\mathfrak{f}|\ \ \text{and}\ \ \
\mu^{\mathfrak{f}^\ast\partial\wh{\Theta}'}_{\mathfrak{f}^\ast\omega'}   = \mu^{\partial\wh{\Theta}'}_{\omega'}\circ |\mathfrak{f}|
$$
With the obvious notation we can rewrite this more suggestively as
$$
|\mathfrak{f}|_\ast\mu^{\partial\wh{\Theta}}_{\omega}   = \mu^{\mathfrak{f}_\ast\partial\wh{\Theta}}_{\mathfrak{f}_\ast\omega}\ \ \text{and}\ \ \
\mu^{\mathfrak{f}^\ast\partial\wh{\Theta}'}_{\mathfrak{f}^\ast\omega'}   = |\mathfrak{f}|^\ast \mu^{\partial\wh{\Theta}'}_{\omega'}.
$$
\qed
\end{theorem}
\begin{proof}[Theorem \ref{THM1141}]
Consider the case where $F:X\rightarrow Y$ is an equivalence between two ep-groupoids.
Let $x\in \supp(\Theta)$ and consider open neighborhoods $U(x)$ and $U(F(x))$ such that $F:U(x)\rightarrow U(F(x))$
is an sc-diffeomorphism, and both neighborhoods have the properness property, allow the natural actions.
We assume that $U(x)$ is equipped with an oriented local branching structures ${(M_i,o_i)}_{i\in I}$, ${(\sigma_i)}_{i\in I}$ so that
$$
\Theta(y)=\sum_{\{i\in I\ |\ y\in M_i\}}\sigma_i\ \ \text{for}\ \ y\in U(x).
$$
We equip $U(F(x))$ with the branching structure obtained via $F|U(x)$ and denote the associated submanifolds by
${(M_i',o_i')}_{i\in I}$.
We compute
\begin{eqnarray}\label{KLM1141}
\mu_\omega^{\wh{\Theta}}(K)&=& \frac{1}{\sharp G^{\text{eff}}_x}\cdot \sum_{i\in I} \sigma_i\cdot \int_{K_i} \omega\\
&=& \frac{1}{\sharp G^{\text{eff}}_{F(x)}}\cdot \sum_{i\in I} \sigma_i\cdot \int_{F(K_i)} F_\ast \omega\\
&=&  \mu_{F_\ast\omega}^{F_\ast\wh{\Theta}}(|F|(K)).\nonumber
\end{eqnarray}
This shows that the homeomorphism $|F|:|X|\rightarrow |Y|$ relates locally the measures 
$\mu^{\wh{\Theta}}_{\omega}$ and $\mu^{F_\ast\wh{\Theta}}_{F_\ast\omega}$. Since measures
are $\sigma$-additive it follows that
$$
\mu^{\wh{\Theta}}_{\omega} =|F|^\ast \mu^{F_\ast\wh{\Theta}}_{F_\ast\omega}.
$$
There is a similar relationship involving pull-backs and the proof follows the same lines
$$
\mu^{F^\ast\wh{\Theta}'}_{F^\ast\omega'} = |F|^\ast \mu^{\wh{\Theta}'}_{\omega'}.
$$
Given the generalized isomorphism $\mathfrak{f}:X\rightarrow Y$ take a representative diagram
$$
X\xleftarrow{F} A\xrightarrow{G} Y
$$
We compute
\begin{eqnarray*}
|\mathfrak{f}|_\ast \mu^{\wh{\Theta}}_\omega = G_\ast F^\ast\mu^{\wh{\Theta}}_\omega
= G_\ast  \mu^{F^\ast\wh{\Theta}}_{F^\ast\omega}
=\mu^{G_\ast F^\ast\wh{\Theta}}_{G_\ast F^\ast\omega}
= \mu^{\mathfrak{f}_\ast\wh{\Theta}}_{\mathfrak{f}_\ast\omega}
\end{eqnarray*}
The pull-back formula is established in the same way.
\qed \end{proof}

\begin{proof}[Theorem \ref{THM1143}]
The considerations are very similar to those in the roof of Theorem \ref{THM1141}.
Using a local branching structure at a point $x$ in $\supp(\partial\Theta)$ we establish
for equivalences that the measures agree locally. This step corresponds to (\ref{KLM1141}).
The rest is formally as in the proof of Theorem \ref{THM1141}.
The details are left to the reader.
\qed \end{proof}

\section{Strong Bundles up to Equivalence}
For the following developments we have to study the behavior of various concepts, introduced above, under generalized strong bundle maps.
In  other words,  we have to carry out some `strong bundle geometry up to equivalence'. This is similar to the case of ep-groupoids;  added are the sc-smooth section functors.  We begin by recalling Theorem \ref{Push-Forw-prop} from Chapter \ref{SEC2}.

{\bf Theorem \ref{Push-Forw-prop}\textcolor{red}{.}}
Let $(P \colon W\to X,\mu)$ and $(P' \colon W'\to X',\mu')$ be two strong bundles over ep-groupoids and $[D]:W\rightarrow W'$ a generalized strong bundle isomorphism.
Denote by $\Gamma(P,\mu)$ and $\Gamma(P',\mu')$ the vector spaces of sc-smooth section functors. Then $[D]$ induces a well-defined isomorphism
$$
[D]_\ast\colon\Gamma(P,\mu)\rightarrow \Gamma(P',\mu').
$$
Its inverse is given by the inverse diagram and equals the pull-back $[D]^\ast$.  The same assertion hold for $sc^+$-section functors 
$$
[D]_\ast\colon\Gamma^+(P,\mu)\rightarrow \Gamma^+(P',\mu')
$$
and sc-Fredholm sections functors
$$
[D]_\ast\colon \text{Fred}(P,\mu)\rightarrow \text{Fred}(P',\mu').
$$
Moreover, if $[D], [D']: W\rightarrow W'$ are two generalized strong bundle isomorphisms inducing the same maps between orbit spaces, i.e. $|[D]|=|[D']|$, then
$[D]_\ast =[D']_\ast$ and $[D]^\ast =[D']^\ast$.
\qed

Next we show that generalized isomorphisms map auxiliary norms to auxiliary norms.
\begin{proposition}\label{prop_auxiliary_norm}\index{P- Auxiliary norms}
Let $(P:W\rightarrow X,\mu)$ and $(P':W'\rightarrow X',\mu')$ be strong bundles over ep-groupoids and $[D]:W\rightarrow W'$ a generalized strong bundle isomorphism.
Given an auxiliary norm $N'$ for $(P':W'\rightarrow X',\mu')$ there is a well-defined pull-back $N=[D]^\ast N'$ which is an auxiliary norm for $(P,\mu)$.
Similarly there is a push forward $[D]_\ast N$ defining an auxiliary norm on $(P',\mu')$. Moreover, 
$$
[D]^\ast[D]_\ast N=N\ \text{and}\ [D]_\ast [D]^\ast N' =N'.
$$
In addition,  if $[D]$ and  $[D']:W\rightarrow W'$ satisfy $|[D]|=|[D']|$,  then $[D]_\ast N=[D']_\ast N$ and $[D]^\ast N'=[D']^\ast N'$.
\end{proposition}
\begin{proof}
In order to  define the auxiliary norm $N:=[D]^\ast N'$we choose  a representative $D\colon W\xleftarrow{\Phi}W''\xrightarrow{\Psi} W'$ and take $e\in W_{0,1}$. Since $\Phi$ is a strong bundle equivalence, we find $e''\in W''_{0,1}$ and a morphism $\Phi(e'')\rightarrow e$. Then we define
$$
N(e):= N'(\Psi(e'')).
$$
As one easily verifies,  the definition does not depend on the choices involved since $N'$ is a functor $W'\rightarrow {\mathbb R}^+$. Moreover it follows from the definition that
$[D]^\ast N=[D']^\ast N'$ provided $|[D]|=|[D']|$.
We can also go in the reverse direction to define the push forward and derive the desired results.
\qed \end{proof}

The $\ssc^+$-multisections behave nicely too under  generalized strong bundle isomorphisms, as we shall see next.

In order to define the pull-back and push forward of $\ssc^+$-multisections we consider the equivalence $\Phi\colon W\to W'$ between the strong bundles 
$(P:W\rightarrow X,\mu)$ and $(P':W'\rightarrow X',\mu')$ covering the equivalence $\phi\colon X\to X'$ between the underlying ep-groupoids. The {\bf pull back of the $\ssc^+$-multisection} \index{pull back of the $\ssc^+$-multisection}
$\Lambda'\colon W'\to \Q^+$ is the multisection $\Lambda\colon W\to \Q^+$ defined by 
$$\Lambda (w)=\Phi^\ast \Lambda' (w)=\Lambda' (\Phi (w)),\quad w\in W.$$
The push-forward $\Phi_\ast \Lambda$ of the $\ssc^+$-multisection $\Lambda^\ast\colon W\to \Q^+$ at the point $w'\in W'$ is defined as follows. Using that the equivalence is essentially surjective we find a point $w\in W$ and a morphism $(\sigma', \Phi (w))\in {\bm{W}}'$,
$$(\sigma', \Phi (w))\colon \Phi (w)\to w'=\mu'(\sigma', \Phi (w))\in W'.$$
Here $\sigma'\in {\bm{X}}'$ is the underlying morphism $\sigma'\colon \phi (x)\to x'$, where $x=Pw\in X$ and $x'=P'w'\in X'$. The {\bf push-forward of the $\ssc^+$-multisection}\index{push-forward of a $\ssc^+$-multisection} $\Lambda'\colon W'\to \Q^+$ is defined as the $\ssc^+$-multisection $\Lambda'\colon W'\to \Q^+$,
$$\Lambda'(w')=\Phi_\ast \Lambda (w')=\Lambda (w).$$
It is independent of the choices of $w$ and morphism $\Phi (w)\to w'$. Indeed, if $w_0\in W$ and $\Phi (w_0)\to w'$ is a second morphism in ${\bm{W}}'$, we obtain the  morphism $\Phi (w)\to \Phi (w_0)$ in $W'$ and hence  a morphism $w\to w_0$ in ${\bm{W}}$ and conclude by the functoriality of 
$\Lambda$ that $\Lambda (w)=\Lambda (w_0)$.
\begin{lemma}
$\Lambda'=\Phi_\ast \Lambda$ is a functor.
\end{lemma}
\begin{proof}
Given a morphism $w_1'\to w_2'$ in ${\bf W'}$,  we have to show that $\Lambda' (w_1')=\Lambda' (w_2')$. There are $w_1, w_2\in W$ and morphisms 
$\Phi (w_1)\to w_1'$ and $\Phi (w_2)\to w_2'$ in ${\bm{W}}'$. Since $w_1'\to w_2'$ is a morphism, there exists also a morphism $\Phi (w_1)\to \Phi (w_2)$. Since  $\Phi$ is an equivalence,  there exists a morphism $w_1\to w_2$ in ${\bm{W}}$ and therefore $\Lambda (w_1)=\Lambda (w_2)$, in view of the functoriality of $\Lambda$. Hence, by definition of $\Phi_\ast \Lambda$, 
$$\Lambda' (w_1')=\Lambda (w_1)=\Lambda (w_2)=\Lambda' (w_2').$$
\qed \end{proof}

To obtain a local section structure $(s'_i, \sigma_i')$ around $x'\in X'$ we take a point $x\in X$ and the morphism $\sigma'\colon \phi (x)\to x'$ in ${\bm{X}}'$ from above, take the local section structure $(s_i, \sigma_i)$ of $\Lambda$ around $x$, define $\sigma_i'=\sigma_i$ and $s_i'(x')=\mu (\sigma', \Phi (s_i (x))$, and use that $\sigma'$ extends to a unique local sc-diffeomorphism by $t\circ s^{-1}$.

By construction, 
$$\Phi^\ast\circ \Phi_\ast \Lambda =\Lambda$$
for every $\ssc^+$-mulitisection $\Lambda\colon W\to \Q^+$.

If $[D]=[W\xleftarrow{\Phi}W''\xrightarrow{\Psi}W']$ is a generalized strong bundle isomorphism covering the generalized isomorphism $[d]=[X\xleftarrow{\phi}X''\xrightarrow{\psi}X']$, we choose a representative diagram and define the push forward $[D]_\ast \Lambda$ as the $\ssc^+$-multisection 
$$[D]_\ast \Lambda=\Psi_\ast\circ \Phi^\ast\colon W'\to \Q^+.$$
The definition does not depend on the choice of the diagram.

\begin{proposition}\index{P- Push forward and pull back of $\ssc^+$-multi sections}
Let $(P:W\rightarrow X,\mu)$ and $(P':W'\rightarrow X',\mu')$ be strong bundles over ep-groupoids and $[D]:W\rightarrow W'$ a generalized strong bundle isomorphism.
Given a $\ssc^+$-multisection $\Lambda$ for $(P,\mu)$ there exists a well-defined push forward $[D]_\ast\Lambda$ and for a $\ssc^+$-multisection $\Lambda'$ for $(P'\,\mu')$
a well-defined pull-back $[D]^\ast\Lambda'$. These operations are mutually inverse. Moreover, if $[D],[D']:W\rightarrow W'$ are generalized strong bundle maps
satisfying $|[D]|=|[D']|$, then the push-forward and pullback operations for $\ssc^+$-multisections are the same. 
\qed
\end{proposition}

The proof follows from previously established properties.

\section{Coverings and Equivalences}\label{SQWERTY116}
We study proper coverings up to equivalence. We recall that a  proper covering functor $F\colon Y\rightarrow X$ between ep-groupoids is a sc-smooth functor such that 
 \begin{itemize}
 \item[(1)]\  The sc-smooth map $F$ between the object M-polyfolds is a surjective local sc-diffeomorphism.
 \item[(2)]\ The preimage $F^{-1}(x)$  of every object $x\in X$ is finite and there exist open neighborhoods
 $U(x)$ of $x$ and $U(y)$ for every $y\in F^{-1}(x)$ such  that $F\colon U(y)\rightarrow U(x)$ is a sc-diffeomorphism 
 and
 $$
 F^{-1}(U(x))=\bigcup_{y\in F^{-1}(x)} U(y).
 $$
 \item[(3)]\ The map
 $$
\alpha\colon  {\bm{Y}}\rightarrow \bx{_{s}\times_F}Y, \quad \phi\mapsto ({\bf F}(\phi),s(\phi))
 $$
 is a sc-diffeomorphism.
 \end{itemize}

Denote for an object $x\in X$ by $Y(x)$ the full subcategory of $Y$ generated by the finitely many 
objects in $F^{-1}(x)$.
 We can define a group action of $G_x$ on $Y(x)$ as follows. For $g\in G_x$ and $y\in Y(x)$
we define $g\ast y := t(\psi)$, where $\alpha(\psi)=(g,y)$. 
Assuming that $\alpha(\psi_1)=(g_1,y)$ and $\alpha(\psi_2)=(g_2,t(\psi_1))$ we see that 
$F(\psi_1)=g_1$, $s(\psi_1)=y$,  $F(\psi_2)=g_2$,  and $s(\psi_2)=t(\psi_1)$. 
There exists a unique $\psi_3$ such that $F(\psi_3)=g_2\circ g_1$ and $s(\psi_3)=y$.
We note that $s(\psi_2\circ \psi_1)=s(\psi_1)=y$ which implies that $\psi_3=\psi_2\circ \psi_1$
by uniqueness.  Consequently defining $g\ast y:= t(\alpha^{-1}(g,y))$ we compute
$$
(g_2\circ g_1)\ast y = t(\psi_3) =t(\psi_2\circ\psi_1)= t(\psi_2) = g_2\ast t(\psi_1) =g_2\ast(g_1\ast y).
$$
Moreover $1\ast y =t(1_y)=y$. Hence we can build the translation groupoid $G_x\ltimes Y_x$.
The morphisms are the pairs $(g,y)$ with $s(g,y)=y$ and $t(g,y)=g\ast y$. The composition
of $(g_2,y_2)$ and $(g_1,y_1)$ is defined if $y_2=g\ast y_1$ and in this case
$(g_2,y_2)\circ (g_1,y_1)=(g_2\circ g_1,y_1).$

The category $Y(x)$ is isomorphic to the category $G_x\ltimes Y_x$ by the functor 
which on objects is the identity $Y(x)\rightarrow Y(x)$ and on morphisms $\alpha$.  Indeed, by construction
$\alpha(1_y)= (1_x,y)$ and 
$$
\alpha(\psi_2\circ \psi_1)= (F(\psi_2\circ\psi_1),s(\psi_1))=(F(\psi_2)F(\psi_1),s(\psi_1)).
$$
We note that  $F(\psi_1)\ast s(\psi_1) = t(\psi_1)=s(\psi_2)$. Hence 
$$
\alpha(\psi_2\circ \psi_1) =(F(\psi_2),s(\psi_2))\circ (F(\psi_2),s(\psi_1))= \alpha(\psi_2)\circ \alpha(\psi_1).
$$

Hence we have proved.
\begin{lemma}
Let $F:Y\rightarrow X$ be a proper covering functor between ep-groupoids. Then for every object $x\in X$
the full subcategory $Y(x)$ generated by the objects in $F^{-1}(x)$ can be naturally identified 
with the translation groupoid $G_x\ltimes Y_x$ for a naturally defined action
$$
G_x\times Y_x\rightarrow Y_x : (g,y)\rightarrow g\ast y:= t(\alpha^{-1}(g,y)).
$$
The identification of $(Y(x),\bm{Y}(x))\rightarrow G_x\ltimes Y(x)$ is given by the identity on objects
and by the map $\alpha$ on morphisms.
\qed
\end{lemma}

Every morphism $\phi:x\rightarrow x'$ defines an equivalence of categories as follows.
\begin{lemma}
Let $\phi:x\rightarrow x'$ be a morphism in $\bm{X}$. Then $\phi$ induces an equivariant map
$$
B_\phi:Y(x)\rightarrow Y(x')
$$
and consequently an equivalence of translation groupoids, which is bijective on objects and morphisms
$$
\wt{\phi}: G_x\ltimes Y(x)\rightarrow G_{x'}\ltimes Y(x'): (g,y)\rightarrow (\phi\circ g\circ \phi^{-1},B_\phi(y)).
$$
The map $B_\phi:Y(x)\rightarrow Y(x')$ is defined by
$$
B_\phi(y)=t(\psi),
$$
where $F(\psi)=\phi$ and $s(\psi)=y$. The map
$$
\alpha^{-1}\circ\wt{\phi}\circ\alpha: (Y(x),\bm{Y}(x))\rightarrow (Y(x),\bm{Y}(x'))
$$
is a natural  equivalence of categories associated to $\phi:x\rightarrow x'$.
\end{lemma}
\begin{proof}
Given $\phi:x\rightarrow x'$ there is an induced equivariant bijection $B_\phi:Y_x\rightarrow Y_{x'}$ 
defined by 
$$
B_\phi(y) = t(\alpha^{-1}(\phi,y)).
$$
The equivariance is obtained as follows. We consider $(g,y)\in G_x\times Y(x)$ and
$\phi:x\rightarrow x'$ and compute $B_{\phi}(g\ast y)$ as follows. We take $\psi_1\in \bm{Y}$
with $F(\psi_1)=\phi$ and $s(\psi_1)=g\ast y =t(\psi_2)$ where $F(\psi_2)=g$ and $s(\psi_2)=y$.
Then $\psi_1\circ \psi_2$ is well-defined and $s(\psi_1\circ\psi_2)=s(\psi_2)=y$.
Moreover $F(\psi_1\circ\psi_2)=\phi\circ g$.
Hence 
$$
B_{\phi}(g\ast y)= t(\psi_1)=t(\psi_1\circ\psi_2) = B_{\phi\circ g}(y).
$$
The morphism $\phi:x\rightarrow x'$ determines the  group isomorphism $G_x\rightarrow G_{x'}$
given by
$$
\gamma_\phi(g) = \phi\circ g\circ \phi^{-1}.
$$
Hence $\gamma_\phi(g)\circ \phi= \phi\circ g$. Consequently
$$
B_{\phi}(g\ast y) = B_{\gamma_\phi(g)\circ \phi}(y).
$$
By definition $B_{\gamma_\phi(g)\circ \phi}(y)=t (\psi_3)$ with $F(\psi_3)=\gamma_\phi(g)\circ\phi$
and $s(\psi_3)=y$. Pick $\psi_4$ with $F(\psi_4)=\gamma_\phi(g)$ and $s(\psi_4)=B_{\phi}(y)$.
We also take $\psi_5$ with $F(\psi_5)=\phi$ and $s(\psi_5)=y$. 
Then $\psi_3=\psi_4\circ \psi_5$ which implies
$$
B_{\gamma_\phi(g)\circ \phi}(y)=t (\psi_3) = t(\psi_4)=\gamma_\phi(g)\ast B_\phi(y).
$$
Bringing everything together we obtain
$$
\gamma_\phi(g)\ast B_\phi(y) = B_\phi( g\ast y).
$$
\qed \end{proof}

Given two proper covering functors $F:Y\rightarrow X$ and $F':Y'\rightarrow X'$ assume we are 
 given a pair of functors $\Psi$ and $\Phi$ making the following diagram commutative
 $$
 \begin{CD}
 Y @>\Psi>> Y'\\
 @V FVV @V F'VV\\
 X@>\Phi>> X'.
 \end{CD}
 $$
 For objects $x\in X$ and $x'\in X'$
 denote by $Y(x)$ and $Y'(x')$ the full subcategories generated by the finitely many objects 
 in $F^{-1}(x)$ and $F'^{-1}(x')$, respectively.  The maps $\alpha$ and $\alpha'$ define bijective equivalences
 of categories
 $$
\bm{Y}(x)\rightarrow G_x\ltimes Y(x):\psi\rightarrow (F(\psi),s(\psi))
$$
and
$$
\bm{Y}'(x')\rightarrow G_{x'}\ltimes Y'(x'):\psi'\rightarrow (F'(\psi'),s(\psi')).
$$
From the commutative diagram we infer that $\Psi$ induces a functor $Y(x)\rightarrow Y'(x')$.
\begin{definition}
Let $F:Y\rightarrow X$ and $F':Y'\rightarrow X'$ be two proper covering functors.
\begin{itemize}
\item[(1)]\ For a pair $(\Psi,\Phi)$ of sc-smooth functors fitting into the commutative diagram
 $$
 \begin{CD}
 Y @>\Psi>> Y'\\
 @V FVV @V F'VV\\
 X@>\Phi>> X', 
 \end{CD}
 $$
we say that $(\Psi,\Phi)$ is an {\bf sc-smooth functor between proper coverings}\index{D- Functor between proper coverings} $F$ and $F'$, written
as 
$$
(\Psi,\Phi):F\rightarrow F'.
$$
\item[(2)]\ An {\bf equivalence between proper covering functors}\index{D- Equivalence between proper coverings} is an sc-smooth functor $(\Psi,\Phi):F\rightarrow F'$
so that $\Phi$ and $\Psi$  are  equivalences of ep-groupoids and for every object $x\in X$ with $x'=\Phi(x)$ the functor $\Psi:Y(x)\rightarrow Y'(x')$
is a bijection on objects.
\end{itemize}
We note that there is a subclass of functors in (1) for which $\Psi:Y(x)\rightarrow Y(\Phi(x))$
is a bijection on objects for every $x\in X$.
\qed
\end{definition}

Assume that $\Gamma=(\Psi,\Phi):F\rightarrow F'$ is an equivalence between proper coverings 
and $(\zeta,\xi):F''\rightarrow F'$ an sc-smooth functor. Take the weak fibered products 
of ep-groupoids associated to the diagrams
\begin{eqnarray*}
& Y\xrightarrow{\Psi} Y'\xleftarrow{\zeta} Y''&\\
&X\xrightarrow{\Phi} X'\xleftarrow{\xi} X''.&
\end{eqnarray*}
We obtain $P:=Y \times_{Y'} Y''$ and $Q:=X\times_{X'} X''$.  For  $(y,\psi',y'')\in Y\times_{Y'} Y''$
we consider $(F(y),F'(\psi'),F''(y''))$ and note that 
$$
\Phi(F(y))= F'(\Psi(y))= F'(s(\psi'))=s(F'(\psi'))
$$
and
$$
\xi(F''(y''))= F'(\zeta(y''))=F'(t(\psi'))=t(F'(\psi')),
$$
so that $(F(y),F'(\psi'),F''(y''))\in X\times_{X'}X''$. We abbreviate 
$$
H:P\rightarrow Q,
$$
where $H(y,\psi',y'')=(F(y),F'(\psi'),F''(y''))$. 
Before we give the main result we need to introduce the notion of 
a natural transformation in our context.
\begin{definition}
Let $F$ and $F'$ be proper covering functors and $(\zeta,\xi),(\zeta',\xi'):F\rightarrow F'$ a
functor between proper coverings. A natural transformation $\tau$ is a pairs $(\tau_Y,\tau_X)$
of natural transformations 
$$
\tau_Y:Y\rightarrow \bm{Y}'\ \ \text{and}\ \ \tau_X:X\rightarrow \bm{X}'
$$
of ep-groupoids $\tau_Y: \zeta\rightarrow \zeta'$ and $\tau_X:\xi\rightarrow \xi'$ so that
$$
\tau_X\circ F = \bm{F}'\circ \tau_Y.
$$
We shall write $\tau:(\zeta,\xi)\rightarrow (\zeta',\xi')$. 
\qed
\end{definition}
If $y\in Y$ we have the morphisms $\tau_Y(y):\zeta(y)\rightarrow \zeta'(y)$ so that 
for any morphism $\psi:y_1\rightarrow y_2$ the diagram
$$
\begin{CD}
\zeta(y_1)@>\zeta(\psi)>> \zeta(y_2)\\
@V \tau_Y(y_1) VV   @V\tau_Y(y_2) VV\\
\zeta'(y_1) @>\zeta'(\psi)>> \zeta'(y_2)
\end{CD}
$$
Applying $ \bm{F}'$ to the diagram we obtain the following diagram where $x_i= F(y_i)$
$$
\begin{CD}
\xi(x_1)@>\xi(\bm{F}(\psi))>> \xi'(x_2)\\
@V \tau_X(x_1)VV @V \tau_X(x_2)VV\\
\xi'(x_1) @>\xi'(\bm{F}(\psi))>> \xi'(x_2).
\end{CD}
$$

The following is the main result in this section.
\begin{theorem}\label{THM1165}
Assume that $\Gamma=(\Psi,\Phi):F\rightarrow F'$ is an equivalence between proper coverings and 
$(\zeta,\xi):F''\rightarrow F'$ an sc-smooth functor between proper coverings. With the ep-groupoids 
$P$ and $Q$ defined as the weak fibered products 
$P=Y\times_{Y'}$ and $Q=X\times_{X'} X''$ 
the functor $H:P\rightarrow Q$  defined by
$$
H(y,\psi',y'')=(F(y),F'(\psi'),F(y''))
$$
the following holds.
\begin{itemize}
\item[{\em(1)}]\ The functor $H$ is well-defined and a proper covering functor.
\item[{\em(2)}]\ The functor  pair $(\pi_{Y''},\pi_{X''}):H\rightarrow F''$ is an equivalence of proper covering functors.
\item[{\em(3)}]\ The functor pairs $(\Psi\circ \pi_{Y},\Phi\circ \pi_{X}):H\rightarrow F'$ and 
$(\zeta\circ \pi_{Y''},\xi\circ\pi_{X''}):H\rightarrow F'$ between proper coverings are naturally equivalent.
\end{itemize}
Moreover, if $X''$ is a tame ep-groupoid (so that trivially $Y''$ is tame) it holds that $P$ and $Q$ are tame ep-groupoids.
If $(\zeta,\xi)$ is also an equivalence between proper coverings, then the pair of projections 
$(\pi_{Y},\pi_X)$ is an equivalence as well.
\end{theorem}

\begin{proof}
The lengthy proof is divided into several steps.\par

\noindent {\bf (A)}   $H$ is well-defined on morphisms. Hence $H$ is a well-defined functor between
ep-groupoids $H:P\rightarrow Q$. Moreover $H$ is sc-smooth.\par

The morphism space $\bm{P}$ associated to $P=Y\times_{Y'} Y''$ consists of all tuples
$(\psi,\psi',\psi'')\in \bm{Y}\times\bm{Y}'\times\bm{Y}''$ satisfying
$$
s(\Psi(\psi))\xrightarrow{\psi'} s(\zeta(\psi'')).
$$
Applying $H$ which by definition is $(F(\psi),F'(\psi'),F''(\psi''))$ we compute that
\begin{eqnarray*}
&s(\Phi(F(\psi)))=s(F'(\Psi(\psi)))=s(F'(\psi'))&\\
&s(\xi(F''(\psi'')))=s(F'(\zeta(\psi'')))=F'(s(\zeta(\psi'')))=F'(t(\psi'))=t(F'(\psi'')).&\nonumber
\end{eqnarray*}
Hence
$$
s(\Phi(F(\psi)))\xrightarrow{F'(\psi')}s(\xi(F''(\psi''))).
$$
This shows that $H:P\rightarrow Q$ is a well-defined functor. 
Since $F,F'$ and $F''$ are sc-smooth the same holds for $H$.\par

\noindent {\bf (B)} $H$ is surjective on objects.\par

Pick $(x,\phi',x'')\in Q$. Then $\phi:\Phi(x)\rightarrow \xi(x'')$. Since $F''$ is a proper covering functor 
we can pick $y''\in Y''$ such that $F''(y'')=x''$. Since $F'\circ \zeta(y'')= \xi\circ F''(y'')=\xi(x'')=t(\phi')=s(\phi'^{-1})$
there exists a unique $\psi'^{-1}$ with
$$
F'(\psi'^{-1})=\phi'^{-1}\ \ \text{and}\ \ s(\psi'^{-1})=\zeta(y'').
$$
Clearly
$$
F'(\psi')=\phi'\ \ \text{and}\ \ t(\psi')=\zeta(y'').
$$
Consider $s(\psi')$ and note that $F'(s(\psi'))=s(F'(\psi'))=s(\phi')=\Phi(x)$. 
Since $(\Psi,\Phi):F\rightarrow F'$ is an equivalence between proper covering functors 
we know that $\Psi:Y(x)\rightarrow Y'(\Phi(x))$ is a bijection. Since $s(\psi')\in Y'(\Phi(x))$
we find a unique $y\in Y(x)$ with $\Psi(y)=s(\psi')$. Consequently
$$
(y,\psi',y'')\in P
$$
and $H(y,\psi',y'')=(x,\phi,x'')$. \par

Next we consider $\wt{\alpha}: \bm{P}\rightarrow \bm{Q}{_{s}\times_H} P$ defined by
$$
\wt{\alpha}(\psi,\psi',\psi'')= (H(\psi,\psi',\psi''),(s(\psi),\psi',s(\psi''))).
$$

\noindent{\bf (C)} The map $\wt{\alpha}$ is an sc-diffeomorphism.\par

It is clear that $\wt{\alpha}$ is sc-smooth. 
We can write $\wt{\alpha}$ as 
$$
\wt{\alpha}(\psi,\psi',\psi'')=((F(\psi),F'(\psi'),F''(\psi'')),(s(\psi),\psi',s(\psi''))).
$$
Note that $s(F(\psi))=F(s(\psi))$, $s(F'(\psi'))=F'(s(\psi'))$, and $s(F''(\psi''))=F''(s(\psi''))$.
Let us first show that $\wt{\alpha}$ is a bijection.  So assume 
we are given $((\phi,\phi',\phi''),(y,\psi',y''))$ in $\bm{Q}{_{s}\times_H} P$. 
Hence 
$$
\phi':s(\Phi(\phi))\rightarrow s(\xi(\phi''))\ \ \text{and}\ \ \psi':\Psi(y)\rightarrow \zeta(y''),
$$
and moreover
$$
(F(y),F'(\psi'),F(y''))=(s(\phi),\phi',s(\phi'')).
$$
Since $F(y'')=s(\phi'')$ there exists a unique $\psi''$ with $F''(\psi'')$ and $s(\psi'')=y''$, and similarly 
there exists a unique $\psi$ with $F(\psi)=\phi$ and $s(\psi)=y$ since $F(y)=s(\phi)$. 
We compute 
$$
s(\Psi(\psi))=\Psi(s(\psi))=\Psi(y) =t(\psi')\ \ \text{and}\ \ s(\zeta(\psi''))=\zeta(s(\psi''))=\zeta(y'')=t(\psi').
$$
Consequently $(\psi,\psi',\psi'')\in \bm{P}$ and 
$$
\wt{\alpha}(\psi,\psi',\psi'')=((\phi,\phi',\phi''),(y,\psi',y'')).
$$
The uniqueness  of choices of $\psi$ and $\psi''$ in the previous construction shows that
$\wt{\alpha}$ is also injective.

Next we need to show that $\wt{\alpha}$ is a local sc-diffeomorphism.  For this it suffices 
to analyze the sc-smoothness of the previous construction in the surjectivity proof.
Assume that  $((\phi_0,\phi'_0,\phi''_0),(y_0,\psi'_0,y''_0))$ in $\bm{Q}{_{s}\times_H} P$
is given. Varying $y''$ sc-smoothly near $y''_0$ the data $\zeta(y'')$ varies sc-smoothly.
Using that $t$ is an sc-diffeomorphism there is a unique $\psi'(y'')$ near $\psi'_0$
with 
$$
t(\psi'(y''))=\zeta(y'').
$$
Then s$(\psi'(y''))$ varies sc-smoothly in $y''$ and since $\Psi$ is an equivalence we find
$y=y(y'')$ near $y_0$ varying sc-smoothly in $y''$ with
$$
\Psi(y(y''))=s(\psi'(y'')).
$$
Hence we obtain the sc-smooth map defined near $y''_0$ and given by
$$
y''\rightarrow (y(y''),\psi'(y''),\zeta(y'')).
$$
This map is also a local sc-diffeomorphism onto an open neighborhood in $P$ of the point $(y_0,\psi'_0,y''_0)$.
Moreover, the map $y''\rightarrow H(y(y''),\psi'(y''),\zeta(y''))$ is sc-smooth into $P$
and maps $y''_0$ to $(y_0,\psi'_0,y''_0)$. The source map $s:\bm{Q}\rightarrow Q$
is a local sc-diffeomorphism and we obtain uniquely defined sc-smooth maps defined for $y''$ near $y''_0$
$$
y''\rightarrow (\phi(y''),\phi'(y''),\phi''(y''))\in Q
$$
with $(\phi(y''_0),\phi'(y''_0),\phi''(y''_0))=(\phi_0,\phi'_0,\phi''_0)$, so that
$$
s( (\phi(y''),\phi'(y''),\phi''(y'')))=H(y(y''),\psi'(y''),\zeta(y'')).
$$
Hence we have the sc-smooth map defined for $y''$ near $y''_0$ with image in $\bm{Q}{_{s}\times_H}P$
defined by
$$
y''\rightarrow ((\phi(y''),\phi'(y''),\phi''(y'')),(y(y''),\psi'(y''),y'')).
$$
Since $F(y(y''))=s(\phi(y''))$ we define $\psi(y'')=\alpha^{-1}(\phi(y''),y(y''))$ and observe that 
it is depending sc-smoothly on $y''$.  Similarly using that $F''(y'')= s(\phi''(y''))$ we obtain the sc-smooth
$y''\rightarrow \psi''(y'')=\alpha''^{-1}(\phi''(y''),y'')$. The same gives also $\psi'(y'')$ satisfying
$F'(\psi'(y))=\phi'(y)$ and $s(\psi'(y))= \Psi(y(y''))$. Hence we obtain the sc-smooth map
$$
y''\rightarrow ((\psi(y''),\psi'(y''),\psi''(y''))
$$
in $\bm{Y}\times\bm{Y}'\times\bm{Y}''$. By construction the image
lies
in the sub-M-polyfold $\bm{P}$. Define 
$$
\wt{\beta}(((\phi(y''),\phi'(y''),\phi''(y'')),(y(y''),\psi'(y''),y'')))= ((\psi(y''),\psi'(y''),\psi''(y'')).
$$
We note that
\begin{eqnarray*}
&&\wt{\alpha}\circ \wt{\beta}((\phi(y''),\phi'(y''),\phi''(y'')),(y(y''), \psi'(y''),y''))\\
&=&((\phi(y''),\phi'(y''),\phi''(y'')),(y(y''), \psi'(y''),y'')).
\end{eqnarray*}
where the right-hand side is an sc-smooth coordinate  depending on $y''$ and $\wt{\beta}$ is sc-smooth.
This shows that $\wt{\alpha}$ is an sc-diffeomorphism.\par

\noindent {\bf (D)}  $H:P\rightarrow Q$ is a local sc-diffeomorphism on objects.\par

Let $(x_0,\phi'_0,x''_0)\in Q$. For $x''$ near $x''_0$ we have the sc-smooth coordinates
$x''\rightarrow (x(x''),\phi'(x''),x'')$, where $\phi'(x''_0)=\phi'_0$ and $t(\phi'(x''))=\xi(x'')$.
With $x''\rightarrow s(\phi'(x''))$ being sc-smooth and $\Phi$ being an equivalence we find $x(x'')$ such that
$\Phi(x(x''))=s(\phi'(x''x))$. Pick any point $(y_0,\psi'_0,y''_0)$ with $H(y_0,\psi'_0,y''_0)=(x_0,\phi'_0,x''_0)$.
Using $F''$ we find sc-smooth $x''\rightarrow y''(x'')$ with $F''(y''(x''))=x''$ and similarly 
using $F$ we find $y(x'')$ with $F(y(x''))=x(x'')$. We note that
$$
s(\phi'(x'')) =\Phi(x(x''))=\Phi\circ F(y(x''))=F'\circ \Psi(y(x'')).
$$ 
We find $\psi'(x'')$ such that
$$
F'(\psi'(x''))=\phi'(x'')\ \ \text{and}\ \ s(\psi'(x''))=\Psi(y(x'')).
$$
Then $(y(x''),\psi'(x''),y''(x''))$ depends sc-smoothly on $x''$ and is a local inverse
for $H$.\par

\noindent {\bf (E)}   $(\pi_{Y''},\pi_{X''}): F\times_{F'} F''\rightarrow F''$ is an equivalence of proper covering functors.\par

We know from the constructions in the ep-groupoid framework that $\pi_{Y''}: Y\times_{Y'} Y''\rightarrow Y''$
and $\pi_{X''}: X\times_{X'}X''\rightarrow X''$ are equivalences of ep-groupoids.  We compute on objects
$$
\pi_{X''}\circ H(y,\psi',y'')= \pi_{X'}\circ (F(y),F'(\psi'),F''(y''))= F''(y'') = F''\circ \pi_{Y'}(y,\psi',y'')
$$
and similarly on morphisms, so that we have the commutative diagram of  functors
$$
\begin{CD}
P=Y\times_{Y'} Y'' @> \pi_{Y''} >>Y''\\
@V H VV   @V F'' VV\\
Q=X\times_{X'} X'' @>\pi_{X''}>> X''.
\end{CD}
$$
Consider $q\in Q$ and the associated $P(q)$ and pick any $y''\in Y''(\pi_{X''}(q))$.
With $q=(x,\phi',x'')$ and $y''\in Y''(x'')$.  Then $F''(y'')=x''$ and 
 we find since $t(\phi')=\xi(x'')=\xi\circ F''(y'')=F'\circ \zeta(y'')$ a unique $\psi'\in \bm{Y}'$
 with $F'(\psi')=\phi'$ and $t(\psi')= \zeta(y'')$.  Since $\Phi(x)=s(\phi')=s(F'(\psi'))=F'(s(\psi'))$
 we can consider $\Psi:Y(x)\rightarrow Y'(\Phi(x))$ an by the properties of $\Psi$ there 
 exists a unique $y\in Y(x)$ with $\Psi(y)= s(\psi')$. Then $(y,\psi',y'')\in P(q)$ and
 $$
 \pi_{Y''}(y,\psi',y'')=y''.
 $$
 This shows that $\pi_{Y''}:P(q)\rightarrow Y''(\pi_{X''}(q))$ is surjective, and the uniqueness assertions
 during the proof show that the map is a bijection. \par
 
  \noindent{\bf (F)}  The  functor pairs $(\Psi\circ\pi_Y,\Phi\circ\pi_X)$ and $(\zeta\circ \pi_{Y''},\xi\circ \pi_{X''})$
 are naturally equivalent as functors $H\rightarrow F'$.\par
 
 We define $\tau_{P}:P\rightarrow \bm{Y}'$ and $\tau_{Q}:Q\rightarrow\bm{X}'$ as
 follows
 \begin{eqnarray*}
 &\tau_{P}(y,\psi',y'') = \psi'&\\
 &\tau_{Q}(x,\phi',x'')=\phi'.&
 \end{eqnarray*}
 We compute 
 $$
 \bm{F}'\circ \tau_P(y,\psi',y'')=\bm{F}'(\psi')
 =\tau_{Q}(F(y),\bm{F}'(\psi'),F(y''))=\tau_{Q}\circ H(y,\psi',y'').
 $$
 Hence $(\tau_P,\tau_Q)$ is compatible with $H$ and $F$ in the sense that 
 $\bm{F}'\circ \tau_P =\tau_Q\circ H$.  We compute further that
\begin{eqnarray*}
\Psi\circ \pi_Y(y,\psi',y'')=\Psi(y) \xrightarrow{\psi'} \zeta(y'')=\zeta\circ \pi_{Y''}(y,\psi',y''),
\end{eqnarray*}
which means since $\tau_Y(y,\psi',y'')=\psi'$, that $\tau_P: \Psi\circ\pi_Y\rightarrow \zeta\circ \pi_{Y''}$.
 Similarly we show that 
 $$
 \tau_Q: \Phi\circ \pi_X\rightarrow \xi\circ \pi_{X''}.
 $$
 Hence $(\tau_P,\tau_Q):(\Psi\circ\pi_Y,\Phi\circ \pi_X) \rightarrow (\zeta\circ \pi_{Y''},\xi\circ \pi_{X''})$ as claimed.\par

 The assertions about tameness are trivial.  If $X''$ is tame, i.e. the objects M-polyfold is tame, then 
 with $s$ and $t$ being local sc-diffeomorphsims the morphism space $\bm{X}''$ is tame. 
 Since the object space $X\times_{X'}X''$ is locally sc-diffeomorphic to $X''$ it is tame and then the same holds
 for the associated morphism space.  Since $F'':Y''\rightarrow X''$ is a proper covering functor it is a surjective local sc-diffeomorphism implying the tameness for the object space $Y''$ and then by the standard argument for the morphism space. The same argument which we used for $X\times_{X'} X''$ implies that $Y\times_{Y'} Y''$ is tame.
\qed \end{proof}

Consider the category ${\mathcal PC}$\index{${\mathcal{PC}}$} of proper coverings between ep-groupoids and the functors defined above as morphisms. Then we have the special collection  ${\bf E}={\bf E}_{{\mathcal PC}}$\index{${\bf E}_{{\mathcal PC}}$} of equivalences,
and as in the ep-groupoid case we can pass to the category ${\mathcal PC}({\bf E}^{-1})$, i.e. carry out the localization
at ${\bf E}_{\mathcal PC}$. The procedure is similar as in the ep-groupoid case and we just sketch the 
basic points and leave details to the reader.

As before the new category ${\mathcal PC}({\bf E}^{-1})$ has the same objects as the category ${\mathcal PC}$.
In order to define the morphisms in the new category we
consider  diagrams 
$$
d\colon F\xleftarrow{\Gamma} A\xrightarrow{\varepsilon} F'
$$
where $F$, $A$, and $F'$ are proper coverings of ep-groupoids, $\Gamma$ an equivalence of proper coverings, and $\varepsilon$ a functor between proper coverings.
We shall call $d$  a {\bf diagram from $F$ to $F'$} and we shall write $d\colon F\rightarrow F'$.

\begin{definition}\index{D- Common refinement}
Assume we are given two diagrams $d,d'\colon F\rightarrow F'$ between proper coverings.
A {\bf common refinement} for the diagrams $d,d':F\rightarrow F'$ is given 
by equivalences $H:A''\rightarrow A$ and $H':A''\rightarrow A'$ giving the diagram
$$
\begin{CD}
F @< \Gamma<< A @>\varepsilon >> F'\\
  @.         @A H AA   @.\\
@.     A''     @. \\
@.  @V H' VV @.\\
F @<\Gamma' << A' @>\varepsilon' >> F'
\end{CD}
$$
so that $\Gamma \circ H\simeq \Gamma'\circ H'$ and $\varepsilon\circ H\simeq \varepsilon'\circ H'$.
\qed
\end{definition}
Having a common refinement defines an equivalence relation as the following proposition shows.
\begin{proposition}
Assume that $d,d',d'':F\rightarrow F'$ are diagrams between proper coverings.  If $d$ and $d'$ have a common refinement and $d', d''$ have a common refinement,
then $d, d''$ have a common refinement.
\end{proposition}
\begin{proof}
We consider the diagrams 
$$
\begin{array}{cc}
\begin{CD}
F @< \Gamma<< A @>\varepsilon >> F'\\
  @.         @A H AA   @.\\
@.     A''     @. \\
@.  @V H' VV @.\\
F @<\Gamma' << A' @>\varepsilon' >> F'
\end{CD}
&
\begin{CD}
F @< \Gamma'<< A' @>\varepsilon' >> F'\\
  @.         @A K AA   @.\\
@.     B''     @. \\
@.  @V K' VV @.\\
F @<\Gamma''<< B'' @>\varepsilon'' >> F'.
\end{CD}
\end{array}
$$
By assumption 
$$
\Gamma\circ H\simeq \Gamma'\circ H',\ \varepsilon\circ H\simeq \varepsilon'\circ H',\
\Gamma'\circ K\simeq \Gamma''\circ K',\ \varepsilon'\circ K\simeq \varepsilon''\circ K'.
$$
We  take the weak fibered product associated to 
$$
A''\xrightarrow{H'} A' \xleftarrow{K} B''
$$
 resulting in
$A''\times_{A'} B''$. In view of Theorem \ref{THM1165}  the  projections $\pi_{A''}$ and $\pi_{B''}$ are equivalences
between proper coverings  since $H'$ and $K$ are equivalences. Moreover, again by Theorem \ref{THM1165}
$$
H'\circ \pi_{A''}\simeq K \circ \pi_{B"}.
$$
Consider
$$
\begin{CD}
F@<\Gamma<< A @>\varepsilon>> F'\\
@.   @A  H\circ \pi_{A''}AA@.\\
@.    A''\times_{A'} B'' @.\\
@.    @V K'\circ \pi_{B''}VV   @.\\
F @<\Gamma''<< B'' @ >\varepsilon''>> F'.
\end{CD}
$$
We note that 
$$
\Gamma\circ (H\circ \pi_{A''})\simeq \Gamma'\circ H'\circ \pi_{A''}\simeq \Gamma'\circ K\circ \pi_{B''}\simeq\Gamma''\circ (K'\circ \pi_{B''}),
$$
and
$$
\varepsilon\circ (H\circ\pi_{A''})\simeq \varepsilon'\circ H'\circ \pi_{A''}\simeq \varepsilon'\circ (K\circ \pi_{B''}).
$$
This concludes the proof that having a common refinement defines an equivalence relation.
\qed \end{proof}
The morphisms in our new category ${\mathcal PC}({\bf E}^{-1})$ will be the equivalence classes
$[d]$ of diagrams. We need to define a composition. As in the ep-groupoid case 
we define for $d:F\rightarrow F'$ and $d':F'\rightarrow F''$ the composition 
$[d']\circ [d]:F\rightarrow F''$ as follows. Writing more precisely
$$
F\xleftarrow{\Gamma} A\xrightarrow{\varepsilon} F'\xleftarrow{\Gamma'} A'\xrightarrow{\varepsilon'} F''
$$
we take the weak fibered product $A\times_{F'} A'$ associated to the part $A\xrightarrow{\varepsilon} F'\xleftarrow{\Gamma'} A'$ 
and define a new diagram $d''$ by
$$
d''\colon F\xleftarrow{ \Gamma\circ \pi_{A}}  A\times_{F'} A'\xrightarrow{\varepsilon'\circ \pi_{A'}} F''.
$$
As in the ep-groupoid case the equivalence class of the diagram $d''$ is independent of the representatives
picked in $[d]$ and $[d']$. Hence we define 
$$
[d']\circ [d]:= [d''].
$$
The identities in the new category are given by $[F\xleftarrow{1_F} F\xrightarrow{1_F} F]$ and the inverses
to $[F\xleftarrow{\Gamma} A\xrightarrow{\Gamma'} F']$ by $[F'\xleftarrow{\Gamma'} A\xrightarrow{\Gamma}F]$.
The verifications are purely formal and follow the scheme used in the ep-groupoid case.
There exists a natural functor
$$
I: {\mathcal PC}\rightarrow {\mathcal PC}({\bf E}^{-1})
$$
which is the identity on objects and on morphisms $\varepsilon:F\rightarrow F'$  is defined by 
$$
I(\varepsilon) = [F\xleftarrow{1_F} F\xrightarrow{\varepsilon} F'].
$$
\begin{remark}\index{R- On strong polyfold bundles}
We can go one step further and incorporate strong bundles.  The following discussion
uses the results from Section \ref{cov_g} and Section \ref{STE__x}.

Recall from Section \ref{cov_g}
the definition of proper strong bundle covering functors, see Definition \ref{proper_sb_covering}.

Let $(W,\mu)$ and $(V,\tau)$ be strong bundles over ep-groupoids $Y$ and $X$, respectively.
A {\bf proper strong bundle covering functor} is a strong bundle map $A$  (in particular fiberwise linear), for which the diagram
\begin{eqnarray}\label{EQN1111}
\begin{CD}
W @> A>>  V\\
@V P VV @ V Q VV\\
Y @> F >> X
\end{CD}
\end{eqnarray}
is commutative and which has 
the following additional properties.
\begin{itemize}
\item[(1)]\ $A$ is surjective and a local strong bundle isomorphism covering a sc-smooth proper covering functor $F:Y\rightarrow X$.
\item[(2)]\  The preimage $F^{-1}(x)$ of every object $x\in X$  is finite,  and there exist open neighborhoods $U(x)\subset X$ and $U(y)\subset Y$ for every $y\in F^{-1}(x)$, where the sets $U(y)$ are mutually disjoint, such  that the map
$$
A:W\vert U(y)\rightarrow V|U(x)
$$
 is a strong bundle isomorphism for every $y\in F^{-1}(x)$, and
 $$
 F^{-1}(U(x))=\bigcup_{y\in F^{-1}(x)} U(y).
 $$
 \item[(3)]\ The map 
 $$
\bm{W}\rightarrow \bm{V}{_{s}\times_A} W:(\psi,w)\rightarrow ((F(\psi),A(w)),w)
 $$
 is a strong bundle isomorphism covering the sc-diffeomorphism $\bm{Y}\rightarrow \bm{X}{_{s}\times_F}Y:\psi\rightarrow (F(\psi),s(\psi))$. Here the strong bundle projection
 $$
 \bm{V}{_{s}\times_A} W\rightarrow \bm{X}{_{s}\times_F} Y
 $$
  is given by $((\phi,v),w)\rightarrow (\phi,P(w))$.
 \end{itemize}
 For the economy of notation we shall introduce the following notation.  A {\bf proper strong bundle covering} $A$
 refers to $A:W\rightarrow V$ and the rest of the data in 
 the diagram (\ref{EQN1111}) is implicitly assumed.   Sometimes we shall write instead of $A$ more comprehensively
 $A\rightarrow F$ to emphasize the underlying proper covering between ep-groupoids.
 If we need to reveal for an argument more data we use the full diagram (\ref{EQN1111}).
 
 There is a category ${\mathcal SBPC}$ where the objects are proper coverings between strong bundles $A$
 and the morphisms are (pairs of ) functors $a: A\rightarrow A'$. There is a again a natural notion of equivalence
 and one can carry out the associated localization procedure resulting in 
 ${\mathcal SBPC}[{\bf E}_{{\mathcal SBPC}}^{-1}]$. The details are left to the reader.
 \qed
  \end{remark}

%
%
%

\begin{partbacktext}
\part{Fredholm Theory in Ep-Groupoids}
\noindent 
In this part we bring the sc-Fredholm theory into the context of ep-groupoids.
We introduce  sc$^+$-multisections and prove associated extension results.
We  introduce the integration results in the ep-groupoid context and discuss
orientation and transversality questions.  The core of  Part {III} consists of the Chapters \ref{CHAPH12} to \ref{CHAPX16}
and extends the work in \cite{HWZ3,HWZ3.5,HWZ7,HWZ7err}.

In Chapter \ref{chap11+} we explore the previous study of the  notions,  which are invariant under equivalences.
This leads to the generalization of the usual finite-dimensional theory of orbifolds and their bundles to the M-polyfold framework.
The results are immediate consequences of the discussion in the core chapters and we allow ourselves to be brief.
In fact, we  restrict ourselves to ourselves to the basic notions and  some sample results.
The reader is encouraged to bring other results into this framework.
The origin of \ref{chap11+} is in the paper \cite{HWZ3.5}.

\end{partbacktext}

\chapter{Sc-Fredholm Sections}\label{CHAPH12}
We shall develop the sc-Fredholm theory for sc-Fredholm section functors $f$. 
  This means we shall discuss their compactness properties
 and develop a transversality and perturbation theory. The functorial property of $f$ means that $f$ locally respects certain symmetries.
 It is a known fact that symmetry and transversality are antagonistic concepts, i.e. in general it is impossible to achieve transversality by using perturbations 
 respecting symmetries. For this reason we have to develop a multi-valued perturbation theory.

 \section{Introduction and Basic Definition}

The definition of an sc-Fredholm section follows immediately from the definitions in the M-polyfold case,
see Definition \ref{oi}.
 We shall consider a strong bundle over an  ep-groupoid $(P:W\rightarrow X,\mu)$ and an sc-smooth section functor 
 $f:X\rightarrow W$. 
 \begin{definition}\index{D- Sc-Fredholm functor}
 The sc-smooth section functor $f$ of $P$ is said to be a {\bf sc-Fredholm section functor}  provided $f$ as an sc-smooth section 
 of the strong bundle $W\rightarrow X$ over the  M-polyfold $X$, the object space, is sc-Fredholm.
 \qed
 \end{definition}
 This means in particular that $f$ is {\bf regularizing} in the following sense.  If $x$ is an object in $X_m$ and $f(x)\in W_{m,m+1}$, then
 $x\in X_{m+1}$. The regularizing property for morphisms holds as well.
To see this, assume that $\phi:x\rightarrow y$ is a morphism in $\bm{X}_{m}$ and $f(\phi)\in \bm{W}_{m,m+1}$.
Using  $f(s(\phi))= s(f(\phi))$ we see that $s(\phi)\in X_{m+1}$ implying that $\phi\in \bm{X}_{m+1}$.
 Moreover, for every smooth object $x\in X$ the germ $(f,x)$ has a filled version $(g,0)$ so that $(g-s,0)$
 for a suitable local sc$^+$-section is conjugated to a basic germ. 
 
In the following discussion we shall introduce the functor versions of the concepts already appearing in the 
M-polyfold version of the sc-Fredholm theory.

 \section{Auxiliary Norms}\label{anorms-ssect}
In  Part I  we have introduced  the notion of an  auxiliary norm for  an M-polyfold. In this section we generalize this concept  to ep-groupoids by incorporating  morphisms.
The additional requirement is the compatibility with morphisms.

\begin{definition}[{\bf Auxiliary norm}] \label{auxuilary_norm_def}\index{D- Auxiliary norm for $(P, \mu)$}
Let $(P:W\rightarrow X,\mu)$ be a strong bundle over an ep-groupoid according to Definition \ref{strong_bundle_ep}. We denote by $W_{0,1}$ the topological space defined by 
$$
W_{0,1}=\{h\in W\ |\ h\ \text{has bi-regularity}\ \ (0,1)\}, 
$$
and by $P:W_{0,1}\rightarrow X$ the projection onto the object space $X$. The fibers $P^{-1}(x)$  are,  as usual,  Banach spaces.
An {\bf auxiliary norm $N$\index{$N$} for $(P,\mu)$}  is a continuous map $N:W_{0,1}\rightarrow [0,\infty)$ having the following properties.
\begin{itemize}
\item[(1)]\ The restriction of $N$ to each fiber $P^{-1}(x)$ is a complete norm.
\item[(2)]\ If $(h_k)$ is a sequence in $W_{0,1}$ such  that $(P(h_k))$ converges in $X$ to some $x$, and $N(h_k)\rightarrow 0$, then $h_k\rightarrow 0_x$ in $W_{0,1}$.
\item[(3)]\ If  $h=\mu(\phi,k)$ for $\phi\in {\bf X}$ and $k\in W_{0,1}$ satisfying  $s(\phi)=P(k)$,  then 
$$N(h)=N(k).$$
\end{itemize}
\qed
\end{definition}
If we ignore the compatibility with morphisms for the moment,  it is easy to construct an auxiliary norm locally, see Part I.
However, in order to construct an auxiliary norm globally we need partitions of unity and we have to incorporate the compatibility with morphisms. 
These partitions of unity do not have to be sc-smooth, since we do not require $N$ to have smoothness properties. Nevertheless
we need their existence and therefore  have to assume paracompactness of the orbit space $|X|$.
The main result in this section guarantees
an auxiliary norm for strong bundles over ep-groupoids whose orbit spaces  are paracompact.

\begin{theorem}\label{ANorm1-prop}\index{T- Existence of auxiliary norms}
Every strong  bundle $P \colon W\rightarrow X$  over an ep-groupoid $X$ with paracompact orbit space $|X|$ admits
an auxiliary norm.
\qed
\end{theorem}
The proof requires some preparation and we start with continuous partitions of unity. Recall that a subset $A$ of an ep-groupoid $X$  is saturated
provided $A=\pi^{-1}(\pi(A))$, where $\pi:X\rightarrow |X|$ is the quotient map.
\begin{definition}\index{D- Partitions of unity on ep-groupoids}
Let $X$ be an ep-groupoid with paracompact orbit space $|X|$. A {\bf continuous partition of unity of the ep-groupoid}  $X$ 
consists of a family of continuous functions $\beta_\lambda:X\rightarrow [0,1]$, $\lambda\in\Lambda$, defined on the object space $X$ and  possessing the  following properties.
\begin{itemize}
\item[(1)]\ If there exists a morphism $\phi:x\rightarrow y$, then $\beta_\lambda(x)=\beta_\lambda(y)$ for all $\lambda$.
\item[(2)]\ For every  point $x\in X$ there exists a saturated open neighborhood $U(x)$ such  that there are only finitely many indices $\lambda$ for which 
$\supp(\beta_\lambda)\cap U(x)\neq \emptyset$.
\item[(3)]\ $\sum_{\lambda\in\Lambda}\beta_\lambda(x)=1$  for every $x\in X$. The  sum is locally finite in view of  property  (2).
\end{itemize}
\qed
\end{definition}
The existence of a continuous partition of unity for the ep-groupoid $X$ is equivalent to  the existence of such a partition on  the orbit space $|X|$.
Indeed,  if $(\what{\beta}_\lambda)$ is a continuous partition of unity for  $|X|$, then $(\beta_\lambda=\what{\beta}_{\lambda}\circ \pi)$ is one for  $X$. Conversely, a continuous partition of unity $(\beta_\lambda)$
for  $X$ descends, in view of (1),  to one on the quotient space $|X|$.
\begin{definition}
Given an open cover ${(U_\lambda)}_{\lambda\in\Lambda}$ of saturated sets for the ep-groupoid $X$, a {\bf subordinate continuous partition of unity}
is a continuous partition of unity for the ep-groupoid $X$ indexed by the same set , say ${(\beta_\lambda)}_{\lambda\in\Lambda}$ (some of the functions might vanish),
such that $\supp(\beta_\lambda)\subset U_\lambda$ for all $\lambda\in\Lambda$.
\qed
\end{definition}

\begin{proposition}\index{P- Existence of continuous  partitions of unity on ep-groupoids}
Consider an ep-groupoid $X$ with paracompact orbit space $|X|$ . Then for every open cover by saturated sets there exists an subordinate continuous partition of unity for the object M-polyfold  $X$.
\end{proposition}
\begin{proof}
If $(U_\lambda)$ is an open cover of saturated sets, then 
$(\pi(U_\lambda))$ is an open cover of $|X|$ and we can take a subordinate continuous partition of unity $(\what{\beta}_\lambda)$. Then $\beta_\lambda:=\what{\beta}_{\lambda}\circ \pi$ is the desired partition of unity on $X$ where some of the functions $\beta_\lambda$ may vanish identically.
\qed \end{proof}

In order to construct the auxiliary norm, the neighborhoods have to be chosen carefully.
This is accomplished by the following theorem, which will be used quite frequently.
We shall refer to the theorem as the good neighborhood theorem.
\begin{theorem}[{\bf Good Neighborhood Theorem}]\label{river}\index{T- Good neighborhood}
Let $P:W\rightarrow X$ be a strong bundle over an ep-groupoid $X$.
 Then every object $x\in X$ possesses open neighborhoods $U'(x)$ and $U''(x)$ having  the following properties.
\begin{itemize}
\item[{\em (1)}]\ $\cl_X(U'(x))\subset U''(x)$.
\item[{\em (2)}]\ $U''(x)$ admits the  natural $G_x$-action (defined by $(\Phi, \Gamma)$),  which leaves $U'(x)$ invariant,
 i.e., $\Phi (g)(U'(x))=U'(x)$ for all $g\in G_x$. 
\item[{\em (3)}]\ $W\vert U''(x)$ is strong bundle isomorphic to a local model $K\rightarrow O$, where we use the local models  $(O,C,E)$ and $(K, U\triangleleft  F, E\triangleleft F)$.
\item[{\em (4)}]\ A  sequence of morphisms $(\phi_k)\subset {\bf X}$ satisfying $s(\phi_k)\in \cl_X(U''(x))$ and $t(\phi_k)\rightarrow z\in X$ possesses a subsequence converging to 
$\phi\in {\bf X}$ satisfying  $s(\phi)\in\cl_X(U''(x))$ and $t(\phi)=z$.
 \item[{\em (5)}]\ If $U(x)=\pi^{-1}(\pi (U'(x)))$ is the saturation of $U'(x)$, then  $$
(U(x)\setminus U'(x) )\cap \cl_X(U''(x)) =\emptyset.
$$
\end{itemize}
\end{theorem}
\begin{proof}
If $x\in X$ we choose by Theorem \ref{x-local-x} an open neighborhood $V(x)$ on which we have the natural $G_x$-action defined by $(\Phi, \Gamma)$  such that the map  
 $t\colon s^{-1}(\cl_X(V(x)))\rightarrow X$ is proper, and there exists a strong bundle isomorphism between $W|V(x)$ and a local strong bundle model.
 Then we  take open neighborhoods $U''(x)$ and $U'(x)$ satisfying 
$$
U'(x)\subset \cl_X(U'(x))\subset U''(x)\subset \cl_X(U''(x))\subset V(x)
$$
and such that $U'(x)$ and $U''(x)$ are invariant under the diffeomorphisms $\Phi (g)$, $g\in G_x$.  At this point (1)-(4) hold and it remains to prove  (5).

Arguing by contradiction we assume that $(U(x)\setminus U'(x))\cap \cl_X(U''(x))\neq \emptyset$. Hence we  find  a morphism $\phi:y\rightarrow z$ from 
$y\in \cl_X(U''(x))\setminus U'(x)$ to  $z\in U'(x)$.  Since $y,z\in V(x)$ there exists a uniquely determined $g\in G_x$ satisfying
$\phi=\Gamma(g,y)$. Consequently,  $z=t(\Gamma(g,y))=\Phi (g)(y)$. Since $z\in U'(x)$ and $U'(x)$ is invariant,  it follows that $y\in U'(x)$ contradicting $y\not \in U'(x)$.
\qed \end{proof}
Now we can prove the main result.
\begin{proof}[Proof of Theorem \ref{ANorm1-prop}]
Let  $U(x)=\pi^{-1}(\pi(U'(x)))$ be the saturation of $U'(x)$ where $U'(x)$ has the properties listed in Theorem \ref{river}. The sets $|U(x)|=|U'(x)|=\pi(U'(x))$, where $x$ varies over $X$,
define an open cover of $|X|$. Since $\abs{X}$  is assumed to be paracompact, we find a subordinate partition of unity  $(\what{\beta}_x)$ for $|X|$, and define for the open cover $(U(x))$ of $X$ the partition of unity 
$\beta_x=\what{\beta}_x\circ\pi$.
For every $x\in X$ we can take $U'(x)$ and fix a strong bundle isomorphism $\Psi:W|U'(x)\rightarrow K$ covering $\psi:U'(x)\rightarrow O$.
Here we use models $(O,C,E)$ and $(K,C\triangleleft F,E\triangleleft F)$.

On $W|U'(x)\rightarrow U'(x)$ we define the action of $G_x$ by the 
\begin{equation}\label{eq_star}
g\ast w= \mu(\Gamma(g, P(w)),w)\quad \text{for $g\in G_x$, $w\in W\vert U'(x)$}
\end{equation}
recalling that $s( \Gamma(g, P(w)))=P(w)$. We note  that $P(g\ast w)= P(\mu(\Gamma(g,y),w))=t(\Gamma (g, P(w))=\Phi(g)(P(w))$ and hence,
$$P(g\ast w)=\Phi(g)(P(w)).$$
In order to verify that \eqref{eq_star} defines a group action we have to recall from the natural representation in Theorem \ref{x-local-x}, the homeomorphism $\Phi$ and the map $\Gamma$, and from Definition \ref{strong_bundle_ep}, the properties of the strong bundle map $\mu\colon {\bf X}{_{s}\times_P} W\to W$. If $g=1_x$, then $\Gamma (1_x, P(w))=1_x$ and $\mu (1_x, w)=w$ for all $w$, and hence $1_x\ast w=w$. 

Moreover, 
\begin{equation*}
\begin{split}
\gamma\ast (g\ast w)&=\mu (\Gamma (\gamma , P(g\ast w)), g\ast w)=
\mu (\Gamma (\gamma , P(g\ast w), \mu (\Gamma (g, P(w)), w)\\
&=
\mu (\Gamma (\gamma , P(g\ast w)\circ \Gamma (g, P(w)), w).
\end{split}
\end{equation*}
Since $t(\Gamma (g, P(w))=\Phi (g)(P(w))$ and $s(\Gamma (\gamma, P(g\ast w)))=\Phi (\gamma )\Phi (g)(P(w))=\Phi (\gamma \circ g)(P(w))$, the composition above is a morphism
$$P(w)\to \Phi (\gamma \circ g)(P(w)),$$
and hence is equal to 
$\Gamma (\gamma \circ g, P(w))$ in view of the properties of $\Gamma$. Consequently, 
$$\gamma \ast (g\ast w)=\mu (\Gamma (\gamma \circ g, P(w)), w)=(\gamma \circ g)\ast w$$
and \eqref{eq_star} defines indeed a group action of $G_x$.  It is the {\bf natural lift} of the natural $G_x$-action on $U'(x)$ to the bundle $W\vert U'(x)$.

We define on $W_{0,1}|U'(x)$ the map $N'_x\colon W_{0,1}|U'(x)\rightarrow [0,\infty)$ by 
$$
N_x'(w)=\frac{1}{\abs{G_x}}\cdot \sum_{g\in G_x}\norm{\pr_2\circ\Psi(g\ast w)}_{F_1}.
$$

The map $N_x'$ is invariant under the group action of $G_x$, i.e.,  
$$
N_x'(h\ast w)=N_x'(w)
$$
for all $h\in G_x$ and 
$w\in W_{0,1}\vert U'(x)$.  

Next we extend the definition to $W_{0,1}|U(x)$.
If  $w\in W_{0,1}\vert U(x)$ with base point $y=P(w)\in U(x)$, there is a  morphism $\phi$ satisfying  $s(\phi)=y$ and $t(\phi)\in U'(x)$,  and we 
define 
$$N''_x(w)=N'_x(\mu(\phi,w)).$$
Here, as before,  $\mu$ is the structural map of the strong bundle structure. 
This is well-defined and the extension is compatible with the initial definition of $N_x'$.
 Next we define the function $N_x:W_{0,1}\rightarrow [0,\infty)$ as follows.
$$
N_x(w)=\begin{cases}
0 &\quad  \text{if $P(w)\not\in U(x)$}\\
\beta_x(P(w))\cdot N''_x(w)&\quad  \text{if $P(w)\in U(x)$}.
\end{cases}
$$
In order to show  that $N_x$ is a continuous function we assume that $w_k\rightarrow w$ in $W_{0,1}$.  If $P(w)=y\in U(x)$, then   $y_k=P(w_k) \in U(x)$   for large 
$k$.  Since $U(x)$ is the saturation of $U'(x)$,  we find a morphism $\phi\colon y\rightarrow z$ with $z\in U'(x)$.
We define $w'=\mu (\phi, w)$ so that 
$P(w')=P\mu (\phi, w)=t(\phi)=z\in U'(x)$. For $k$ large we similarly define for the sequence $w_k$ satisfying $P(w_k)=y_k\in U(x)$ the sequence $w_k'=\mu (t\circ s^{-1}(y_k), w_k)$. 
Since $w_k\to w$ in $W_{0,1}$ and $\mu$ and $t\circ s^{-1}$ are local sc-diffeomorphisms, we conclude $w_k'\to w'$.  Therefore, $N_x'(w_k')\to N_x'(w')$ and so, $N_x(w_k)\to N_x(w).$

If $P(w)\not \in U(x)$ we know that $N_x(w)=0$. If $P(w)$ does not belong to the closure of $U(x)$,  the continuity assertion is trivial.
Hence we assume that $P(w)\in \cl_X(U(x))\setminus U(x)$. By construction,  $\supp(\beta_x)$ is a closed subset in $X$ contained in $U(x)$.
We find an open neighborhood $V(P(w))$ in $X$ not intersecting the support of $\beta_x$. If now $w_k\rightarrow w$, then $P(w_k)\in V(P(w))$ for large 
$k$,  and the continuity assertion follows.

Clearly, if $\beta_x(y)\neq 0$,  then $N_x''$,  restricted to the fiber of $W_{0,1}$  over $y$,  is a complete norm.
Assume that $(w_k)$ is a sequence so that $y_k=P(w_k)$ has the properties
\begin{itemize}
\item[(1)]\  $y_k\rightarrow y$ with $\beta_x(y)>0$.
\item[(2)]\ $N_x(w_k)\rightarrow 0$.
\end{itemize}
Then we  find as above an isomorphic sequence $(w_k')$ for which  $P(w'_k)=y_k'\rightarrow y'\in U_x'$. It follows from the construction of $N_x'$
that $w_k'\rightarrow 0_{y'}$ and the same holds evidently for $(w_k)$. 

Finally,  we define the function $N\colon W_{0,1}\to [0,\infty)$ by 
$$
N=\sum_{x\in X} N_x.
$$
Given any $y\in X$ there exists,  by construction,  a saturated open neighborhood $V=V(y)$ in $X$ so that there are only finitely many
$x\in X$ for which $N_x\neq 0$ for the fibers above $y$. This shows that $N$ is continuous. It clearly induces fiber-wise a complete norm on $W_{0,1}$.
Moreover, if $N(w_k)\rightarrow 0$ and $P(w_k)\rightarrow y$ there is at least one $x$ for which  $\beta_x(y)>0$. We conclude
that $w_k\rightarrow 0_y$ in $W_{0,1}$ by the previous discussion.

We have verified that the function $N$ is an auxiliary norm for $(P,\mu)$ and the proof of Theorem \ref{ANorm1-prop} is complete.
\qed \end{proof}

In certain applications we shall need to extend given auxiliary norms
over $\partial X$ to auxiliary norms over
 $X$. The basic result is the following theorem, where the reader should note that we require 
 $X$ to be a tame ep-groupoid.
\begin{theorem}[{\bf Extension of Auxiliary Norms}]\label{EXTT}\index{T- Extension of auxiliary norms}
Let $(P\colon W\rightarrow X,\mu)$ be a strong bundle over a  tame   ep-groupoid $X$ with paracompact orbit space $|X|$.
Suppose $N$ is an auxiliary norm given for $W|\partial X$, i.e., 
$$
N:W_{0,1}|\partial X\rightarrow [0,\infty)
$$
has the properties required in Definition \ref{auxuilary_norm_def}. Then there exists an auxiliary norm $\wtilde{N}:W_{0,1}\rightarrow [0,\infty)$
which coincides over $\partial X$ with $N$.
\end{theorem}
\begin{proof}
It is clear from the proof of Theorem \ref{ANorm1-prop} that we only need to show, for every $x\in \partial X$, that $N$ has
a local extension to a saturated small open neighborhood. The proof of Theorem  \ref{ANorm1-prop} shows how to define 
auxiliary norms on saturated open neighborhoods $U(x)$ for which we assume that $\cl_X(U(x))\cap \partial X=\emptyset$.
Then a partition of unity argument,  as previously used, gives an auxiliary norm for $W_{0,1}\rightarrow X$.

Since the construction is local we can work in local coordinates.  Hence we consider  a strong local bundle $P:K\rightarrow O$,
where $(O,C,E)$ is a tame 
local M-polyfold model and $K\subset C\triangleleft F$ is the image of a strong bundle retraction
$R:U\triangleleft F\rightarrow U\triangleleft F$ covering the sc-smooth retraction $r:U\rightarrow U$ satisfying  $r(U)=O$,  and $U\subset C$ is open in the partial quadrant $C=[0,\infty)^n\oplus H$ in the sc-Banach space $E$. Since $r$ is tame it satisfies $d_C(r(x))=d_C(x)$ for $x\in U$, which will be important. 
We may assume that we work near a point $x_0=(0,h_0)\in O$. (We cannot assume that $h_0\in E$ is zero, since $x_0$ may be on level $0$).
We shall first carry out a local construction, in order to construct $\wtilde{N}$ on $W_{0,1}|U'(x_0)$. Here $U'(x_0)\subset O$ is an open neighborhood
invariant under $G_x$ and $\tilde{N} =N$ on $W_{0,1}|(\partial O\cap U'(x_0))$. 

We  now turn to the details. We abbreviate  $e_0=(1,1,1,\ldots ,1,0)\in C=[0,\infty)^n\oplus H$. Then there exists a $\ssc^0$-function  $\tau\colon C\rightarrow [0,\infty)$ so that for given $(t,h)\in [0,\infty)^n\oplus H$ the vector $(t,h)-\tau(t,h)e_0$ has at least 
one of its first $n$-coordinates vanishing.  The map $\tau$ only depends on $t$, and is given by 
$$
\tau(t,h) =\min\{t_1,\ldots ,t_n\}.
$$
Then the map 
$$
\text{$q\colon C\rightarrow C$, \quad  defined by \quad $q(t,h)=(t,h)-\tau(t)e_0$}, 
$$ 
is a $\ssc^0$-retraction onto $\partial C$.  If  $(t,h)\in  [0,\infty)^n\oplus H$  is near $(0,h_0)$, then  $(t,h)\in U$
and $r(q(t,h))\in U$. Using the tameness property of $r$ we see that  $r(q(t,h))\in O\cap  \partial C =\partial O$. Also for $(t,h)\in \partial O$
it holds that $r(q(t,h))=r(t,h)=(t,h)$
Consequently, 
$$
(t,h)\mapsto r(q(t,h))
$$
is a $\ssc^0$-retraction defined near $(0,h_0)$, which fixes the points in $\partial O$ near $(0,h_0)$.
 For  $((t,h),w)\in K$ in which  $(t,h)$ is close to $(0,h_0)$ and $w$ is on level $1$,   we define $N'$ by 
$$
N'((t,h),w)= N(r(q(t,h)),A(r(q(t,h)))w) + \norm{w-A(r(q(t,h)))w}_{F_1}, 
$$
where $R(u,w)=(r(u),A(u)w)$ is the strong bundle retraction. If $((t,h),w)\in K$ with $(t,h)\in \partial O$ near $(0,h_0)$, then $R((t, h), w)=(r(t, h), A(t, h)w)=((t, h), w)$ and $r(q(t, h))=r(t, h)=(t, h)$. Consequently, 
$$
N'((t,h),w)=N((t,h),w).
$$
Therefore $N'$ extends $N$. At this point we have constructed for $W_{0,1}\vert U'(x_0)$ an extension of $N$, where $U'(x_0)$
is a sufficiently small open neighborhood of $x_0$ in $O$, which we may assume to be invariant under the $G_{x_0}$-action.
Now we average $N'$ over the group $G_{x_0}$ to obtain as a result $N''$. Since $N$ was already compatible with morphisms,  the averaged $N''$ is still an extension 
of $N$. 

At this point we leave the local coordinates. Pulling back the local data by a strong bundle chart we have proved the following local extension result.
Under the assumptions of the theorem there exists for every $x\in \partial X$ an open neighborhood $U'(x)\subset X$ invariant under $G_x$
and an auxiliary norm $N_x'' :W_{0,1}|U'(x)\rightarrow [0,\infty)$ so that $N=N''_x$ on $W_{0,1}|(\partial X\cap U'(x))$.  Moreover,
if $\phi\in {\bf X}$ satisfies  $s(\phi),t(\phi)\in U'(x)$ and $w'=\mu(\phi,w)$, then $N''_x(w)=N_x''(w')$.  By taking $U'(x)$ smaller but still invariant under
$G_x$, we may assume that  every sequence of morphisms $(\phi_k)$ with $s(\phi_k)\in \cl_X(U'(x))$ and $(t(\phi_k))$ in a compact subset, has a convergent
subsequence. We can extend, as already demonstrated  previously in the proof of Theorem \ref{ANorm1-prop}, $N''_x$ to the saturation $U(x)$ of $U''(x)$. 
So far  we have constructed for every $x\in \partial X$ an extension of $N$ to a saturated open neighborhood $U(x)\subset X$.
We already know how to define auxiliary norms for saturated open neighborhoods of points $x\in X\setminus\partial X $.
In this case we may assume that the saturated open neighborhoods are small enough, so that their closures do not intersect $\partial X$. 
Finally,  using a partition of unity argument,  we obtain an auxiliary norm possessing the desired properties. The proof of Theorem \ref{EXTT} 
is complete.
\qed \end{proof}

The strong bundle  $P\colon W\rightarrow X$ over an ep-groupoid is said  to have  {\bf reflexive
$1$-fibers}\index{Reflexive $1$-fibers} if there exists a strong bundle atlas with local strong bundle models $K\subset C\triangleleft F$ in which 
 $F_1$ is a reflexive Banach space, and $C\subset E$ a partial quadrant in a sc-Banach space. Following \cite{HWZ3}
 we introduce the notion of mixed convergence.  
  \begin{definition}[{\bf Mixed convergence}]\label{DEFP1228}\index{D- Mixed convergence}
 A sequence $(w_k)\subset W$ of bi-regularity $(0,1)$, i.e. $w_k\in W_{0,1}$, is said to be {\bf mixed convergent} to
 an element $w\in W_{0,1}$ provided $P(w_k)\rightarrow P(w)=:x$ in $X_0$, and there exists a strong bundle chart
 $\Psi:W|U(x)\rightarrow K\subset E\triangleleft F$ ($F_1$ being reflexive),  
 such  that for  $\Psi(w_k)=(a_k,h_k)$ and $\Psi(w)=(a,h)$  the sequence $h_k$ converges weakly to $h$ in $F_1$,  denoted by $h_k\rightharpoonup h$ in $F_1$.
 We shall write
 $$
 w_k\stackrel{m}{\longrightarrow} w
 $$
if  $w_k$ is mixed convergent to $w$.
\qed
\end{definition}
\begin{remark}\index{R- On mixed convergence}
 The definition does not depend on the choice of the chart $\Psi$. Indeed, if $\Phi:K\rightarrow K'$
is a strong bundle isomorphism (arising as a transition map) and $\Phi(a,h)=(\phi(a),A(a,h))$, we consider a sequence $w_k=(a_k,h_k)\in K$ satisfying $a_k\rightarrow a$ in $E_0$ and 
$h_k\rightharpoonup h$ in $F_1$. The latter implies that $h_k\rightarrow h$ in $F_0$.
Then $\Phi(a_k,h_k)=:(b_k,l_k)$ satisfies $b_k\rightarrow b$ in $E_0'$ and $l_k\rightarrow  l$ in $F_0'$.
For large $k$,  the operator norms of continuous linear operators $A(a_k,\cdot )\colon F_1\rightarrow F_1'$ are uniformly bounded.  Therefore,  the sequence 
$(l_k)=(A(a_k,h_k))$ is bounded in $F_1'$.
From  $h_k\rightarrow h$ in $E_0$ we conclude the convergence  $(b_k,l_k)=\Phi(a_k,h_k)\rightarrow \Phi(a,h)=:(b,l)$ in $E_0'\oplus F_0'$. The boundedness of $(l_k)$ in $F_1'$ and the convergence of $(b_k,l_k)$ on level $0$ 
implies the weak convergence $l_k\rightharpoonup l$ in $F_1'$. 
\qed
\end{remark}
In the context of strong bundles with reflexive fibers we consider special auxiliary norms.
\begin{definition}\index{D- Reflexive auxiliary norm}
Let $P:W\rightarrow X$ be a a strong bundle over an ep-groupoid with reflexive $1$-fibers . A {\bf reflexive auxiliary norm}
$N:W_{0,1}\rightarrow [0,\infty)$ is a continuous map possessing  the following properties.
\begin{itemize}
\item[(0)]\ $N$ is a functor, i.e. if $h'=\mu(\phi,h)$ for $h,h'\in W_{0,1}$ and $\phi:P(h)\rightarrow P(h')$ a morphism in ${\bf X}$, then
$N(h)=N(h')$.
\item[(1)]\ For every $x\in X$ the function $N$,  restricted to the fiber of $W_{0,1}$ over $x$,  is a complete norm.
\item[(2)]\ If $(w_k)$ is a sequence in $W_{0,1}$  satisfying  $P(w_k)\rightarrow x$ in $X$ and $N(w_k)\rightarrow 0$,
then $w_k\rightarrow 0_x$ in $W_{0,1}$.
\item[(3)]\ If $(w_k)\subset W_{0,1}$ is mixed convergent to $w\in W_{0,1}$, then 
$$
N(w)\leq \text{liminf}_{k\rightarrow 0} N(w_k).
$$
\end{itemize}
\qed
\end{definition}
The next  result guarantees the existence of reflexive auxiliary norms for strong bundles $P:W\rightarrow X$
with reflexive $1$-fibers and paracompact orbit spaces $|X|$. 

\begin{theorem}[{\bf Reflexive auxiliary norms}]\label{EXTTT}\index{T- Reflexive auxiliary norms}
We assume 
$ P \colon  W\rightarrow X$ is a strong bundle  with reflexive $1$-fibers, over an ep-groupoid $X$ with paracompact orbit space $|X|$. Then there exists a  reflexive auxiliary norm $N$.  Moreover,  for every
functorial continuous map $f:X\rightarrow (0,\infty)$ i.e., $f$ is invariant under the morphisms, the map $h\rightarrow f(P(h))N(h)$ is a reflexive auxiliary norm.
For every auxiliary norm $N'$, there exist two reflexive auxiliary norms $N_1$ and $N_2$ 
satisfying
$$
N_1\leq N'\leq N_2.
$$
\qed
\end{theorem}
\begin{remark}\index{R- On the constructions of auxiliary norms}
Most of the constructions are similar to the  constructions in Theorem \ref{EXTT}. Before we begin with the proof 
let us note the following facts which are easily verified.
\begin{itemize}
\item[(1)]\ If $N$ is a reflexive auxiliary norm, so is $(f\circ P) \cdot N$, if  $f\colon X\rightarrow (0,\infty)$ is a functorial continuous map.
\item[(2)]\ If  $N_1$ and $N_2$ are two auxiliary norms, there is a function  $f$ as above such  that $N_1\leq (f\circ P) \cdot N_2$.
\item[(3)]\ If $U$ is a saturated open subset of $X$, and $N$ a reflexive auxiliary norm for $W_{0,1}|U$, then,  given a functorial
continuous map $f:X\rightarrow [0,\infty)$ with support $\supp(f)$ contained in $U$, the map $\wtilde{N}:W_{0,1}\rightarrow [0,\infty)$,  
defined by $\wtilde{N}(w)=0$ if $P(w)\not\in U$,  and by $\wtilde{N}(w)=f(P(w))N(w)$ if $P(w)\in U$, has the property that for a mixed convergent sequence 
$w_k\stackrel{m}{\longrightarrow}w$ the estimate $N(w)\leq \text{liminf}_{k\rightarrow \infty} N(w_k)$  holds. 
\end{itemize}
\qed
\end{remark}
\begin{proof}[Theorem \ref{EXTTT}]
It suffices to construct for every $x\in X$ an open saturated neighborhood $U$ and an reflexive auxiliary norm for $W_{0,1}|U$. 
Then the above facts and the construction in Theorem \ref{EXTT} will lead to the desired result.
Around  an object  $x\in X$ we  take a strong bundle chart $\Psi:W|U'(x)\rightarrow K$, where $K\subset C\triangleleft F$,
and $C$ is a partial quadrant in a sc-Banach space. Here $F_1$ is a reflexive Banach space.
We may assume that $U'(x)$ admits the natural $G_x$-action, which has a natural lift 
to $W|U'(x)$, and   
define 
$$
N_{U'(x)}(h) =\frac{1}{|G_x|} \sum_{g\in G_x} \norm{\text{pr}_2\circ \Psi(g\ast h)}_{F_1}.
$$
 We extend $N_{U'(x)}$ to the saturation ${U}(x)=\pi^{-1}(\pi (U'(x))$ of $U'(x)$ as follows.
 If $w\in W_{0,1}|{U}(x)$ there exists a morphism $\phi\colon P(w)\rightarrow y\in U'(x)$
 and we define $N_{{U}(x)}(w)=N_{U'(x)}(w')$, where $w'=\mu(\phi, w)$.
 The definition does not depend on the choice of the morphism $\phi$ as long as it has the stated properties.
 Then $N_{U(x)}$ has the desired properties and the proof is completed by a previous partition of unity argument.
 \qed \end{proof}

We also need for later inductive constructions an extension theorem for reflexive auxiliary norms.

\begin{theorem}[{\bf Extension of reflexive auxiliary norms}]\label{THMOP12213}\index{T- Extension of auxiliary norms}
Let $W$ be a strong bundle with reflexive 1-fibers over the tame  ep-groupoid $X$. We assume that a reflexive auxiliary norm
$N$ is given for $W_{0,1}\vert \partial X$. Then there exists a reflexive auxiliary $\wtilde{N}:W_{0,1}\rightarrow [0,\infty)$ which extends $N$, i.e.
$$
\wtilde{N}|(W_{0,1}|\partial X)= N.
$$
\end{theorem}
\begin{proof}
The local extension defined in Theorem \ref{EXTT} is well-behaved under mixed limits which follows immediately by inspection.
We can define local auxiliary norms over saturated neighborhoods of points $x\in X$ not belonging to the boundary, which are also 
well-behaved with respect to mixed convergence,  and the  argument is completed by partitions of unity,  as in  Theorem \ref{EXTTT}.
\qed \end{proof}

\section{\texorpdfstring{$\bssc^+$}{sssc}-Section  Functors}\label{SECR12.3}
The main point of this subsection is to show that  there is a significant supply of $\ssc^+$-sections under suitable conditions. In addition,  we show that $\ssc^+$-sections defined over the boundary can be extended to the whole ep-groupoid.  For the notion of a sc-smooth bump function in the M-polyfold framework we refer to   Part I.
\begin{definition}\label{d-ep-groupoid_bump}
An ep-groupoid $X$ {\bf admits sc-smooth bump functions}\index{D- Ep-groupoid admitting sc-smooth bump functions} provided for every
object $x$ and every saturated open neighborhood $U(x)$ there exists a sc-smooth map $f\colon X\rightarrow {\mathbb R}$ not identically $0$,  having its support in $U(x)$ and 
satisfying $f(y)=f(z)$ if  there exists a morphism $\phi:y\rightarrow z$.
\qed
\end{definition}
We point out  that in our definition we do not require $f$ to take values in $[0,1]$. 
As we shall see next there is an easy criterion to decide whether  the ep-groupoid admits sc-smooth bump function.

\begin{theorem}\label{t-bump_functions}\mbox{}\index{T- Bump functions}
Let $X$ be an ep-groupoid.
\begin{itemize}
\item[{\em (1)}]\ The ep-groupoid $X$ admits sc-smooth bump functions if and only if the underlying object M-polyfold admits sc-smooth bump functions. 
\item[{\em (2)}]\ If the ep-groupoid $X$ admits sc-smooth bump functions, then there exists for every object $x\in X$ and every saturated open neighborhood
$U(x)$ a sc-smooth bump function $f:X\rightarrow [0,1]$ satisfying $f(x)=1$ and $\supp(f)\subset U(x)$.
\end{itemize}
\end{theorem}
The subtle point  of the statement in (1) is the following. For  ep-groupoids the bump functions are compatible with the morphisms, whereas for the underlying object space $X$ viewed 
as a M-polyfold there is no compatibility condition.  Concerning (2) we note that for a given $x$ all the guaranteed sc-smooth bump functions might not take values in $[0,1]$, and even worse
it could be that $f(x)=0$. However,  the theorem claims that one can construct  sc-smooth bump functions possessing  the additional properties once there is the guaranteed supply of sc-smooth bump functions.
\begin{proof} [Theorem \ref{t-bump_functions}]
(1) We 
first  assume that the ep-groupoid $X$ admits sc-smooth bump functions. We consider the object M-polyfold $X$ and choose a point $x\in X$ and an open neighborhood $U(x)\subset X$ of $x$ in the  M-polyfold $X$.  Then we take a smaller open neighborhood 
$V(x)$ satisfying  $\cl_X(V(x))\subset U(x)$ such that $V(x)$ is equipped with a natural $G_x$-action $\Phi$ as defined in Theorem \ref{x-local-x} and for which, moreover, the map   $t:s^{-1}(\cl_X(V(x)))\rightarrow X$ is proper.
Next we take another open neighborhood $W(x)$ having the following properties.
\begin{itemize}
\item $W(x)\subset \cl_X(W(x))\subset V(x)$.
\item $W(x)$ is invariant under the $G_x$-action $\Phi$.
\end{itemize}
On   the saturation $\wt{W}(x)=\pi^{-1}(\pi(W(x)))$ of $W(x)$ we find,  by assumption,  a sc-smooth bump function $f:X\rightarrow {\mathbb R}$ satisfying $\supp(f)\subset \wtilde{W}(x)$ and $f(y')=f(y'')$ if there is a morphism $\phi\colon y'\to y''$. 
There is a point $z'\in \wtilde{W}(x)$ at which  $f(z')\neq 0$ and there exists   a morphism $\phi\colon z'\rightarrow y$ with $y\in W(x)$, and hence $f(y)=f(z')\neq 0$. Therefore $f$ does not vanish 
on $W(x)$ and we 
define $\beta\colon X\rightarrow {\mathbb R}$ by
\begin{equation}\label{con_xx}
\beta(z)=\begin{cases}
f(z)&\quad  \text{for $z\in W(x)$}\\
0&\quad  \text{for $z\in X\setminus W(x)$}.
\end{cases}
\end{equation}
Then $\beta\not\equiv 0$ and 
$$
\supp(\beta)\subset \cl_X(W(x))\subset V(x)\subset \cl_X(V(x))\subset U(x).
$$
We shall show that $\beta$ is a sc-smooth function, so that that $\beta$ is a sc-smooth bump function having its  support in $U(x)$. 

By construction,  $\beta$ is sc-smooth at all points which do not belong to 
$\cl_X(W(x))\setminus W(x)$.  Hence we only have to consider points $z \in 
\cl_X(W(x))\setminus W(x)$. In the case 
$z\not\in \wtilde{W}(x)$, using that $\supp(f)$ is closed in $X$ and a subset of $\wtilde{W}(x)$, we find an open neighborhood $Q(z)$ not intersecting $\supp(f)$  and hence $\beta=0$ on $Q(z)$.    In the remaining case where $z\in \wtilde{W}(x)$, there exists a morphism $\phi\colon z_0\to z$ for some $z_0\in W(x)$. Since $z\in\cl_X(W(x))\setminus W(x)$, we conclude that 
$z\in V(x)$. But also $z_0\in V(x)$ and therefore by Theorem \ref{x-local-x}, for a suitable $g\in G_x$ we have  
$
\Phi (g)(z_0)=z.
$
From the invariance  of $W(x)$ under the action of $G_x$ we deduce that $z\in W(x)$ contradicting $z\not \in W(x)$ and showing that this case does not occur. We have proved that $\beta$ is a  sc-smooth function with $\supp (\beta)\subset U(x)$.  

Conversely,  we assume that the object M-polyfold $X$ admits sc-smooth bump functions.  For a given $x\in X$ and a saturated open neighborhood
$U(x)$ we choose  an open neighborhood $V(x)$ satisfying  $\cl_X(V(x))\subset U(x)$ such  that on $V(x)$ we have the natural $G_x$-action
and $t:s^{-1}(\cl_X(V(x)))\rightarrow X$ is proper.  We take another open neighborhood $W(x)$ which is invariant under the $G_x$-action
and satisfies $\cl_X(W(x))\subset V(x)$. Then we take a sc-smooth bump function $f\colon X\rightarrow {\mathbb R}$ having its support in $W(x)$ and 
consider $f\colon W(x)\rightarrow {\mathbb R}$. We can average with respect to the group $G_x$ and may assume without loss of generality that
$f$,  on $W(x)$,  is $G_x$-invariant.  Now we define the function $\beta\colon X\rightarrow {\mathbb R}$ by
$$
\beta(y)=
\begin{cases}
f(z)&\quad  \text{if there exists a morphism $\phi\colon z\to y$ with $z\in W(x)$}\\
0&\quad \text{if $y\not\in \pi^{-1}(\pi(W(x))$}.
\end{cases}
$$
The function $\beta$ is well-defined, i.e., it coincides 
with $f$ on $W(x)$. 
Indeed, if $y\in W(x)$ and $\phi\colon z\to y$ for $z\in W(x)$, then $z=\Phi (g)(y)$ for some $g\in G_x$ and hence, in view of the invariance of $W(x)$, 
$f(z)=f(y)$ so that $\beta=f$ on $W(x)$.
Moreover, by construction $\beta\neq 0$ and $\supp(\beta)\subset U(x)$. Further $\beta$ is compatible with morphisms.

It remains to show that $\beta$ is sc-smooth. This is clear at point a  $z\in X\setminus \wt{W}(x)$ and a point $z\in \wt{W}(x)=\pi^{-1}(\pi (W(x)))$.
Hence we only have to consider points $z\in \cl_X(\wtilde{W}(x))\setminus \wtilde{W}(x)$. We shall show that such a point $z$ has an open neighborhood $Q(z)$ not intersecting the support of $\beta$, which will imply the result.  Otherwise,  arguing indirectly,  we find a convergent sequence $z_k\rightarrow z$ with $\beta(z_k)\neq 0$. Then, by the construction
of $\beta$, there exist morphisms $\phi_k$ satisfying  $t(\phi_k)=z_k$ and $s(\phi_k)\in W(x)$. Using the properness, after perhaps taking a subsequence,
we may assume that $\phi_k\rightarrow\phi$. Hence $t(\phi)=z$ and $s(\phi)\in\cl_X(W(x))$. Since $z\not\in \wtilde{W}(x)$
we must have $s(\phi)\not\in W(x)$. We find arbitrarily small open neighborhoods of $s(\phi)$ on which $\beta$ vanishes.
Taking the image under the local sc-diffeomorphism $t\circ s^{-1}$ we see that $z$ has such a neighborhood as well. This leads  to  a contradiction and
the proof of (1)  is complete.\par

(2) If the ep-groupoid $X$ admits sc-smooth bump functions we know from  part (1) that the object M-polyfold admits 
sc-smooth bump functions.  In view of a result in Part I, (Proposition \ref{prop-x5.36}) we find on the object M-polyfold sc-smooth bump functions 
which take values in $[0,1]$, satisfy $f(x)=1$, and have support in the prescribed open neighborhood $U(x)$.  Now carrying out the construction
in (1) for this function  $f$ we first average and then extend to obtain a sc-smooth bump function  on the ep-groupoid 
which,  in addition,  takes values in $[0,1]$ and at the point $x$ the value $1$. This completes the proof of (2). The proof of Theorem \ref{t-bump_functions} is complete.
\qed \end{proof}
Next we use the previous discussion to prove  the existence of suitable $\ssc^+$-sections for strong bundles $(P, \mu)$ over ep-groupoids $X$ admitting sc-smooth bump functions.  Sections $\sigma\colon X\to W$ of strong bundles over ep-groupoids have, of course, to be compatible with morphisms, i.e., they  satisfy
$$
\sigma (t(\phi))=\mu(\phi,\sigma (s(\phi)))\quad  \text{for all $\phi\in {\bf X}$.}
$$
\begin{theorem}\label{LLLOP}\index{T- Existence of $\ssc^+$-sections}
We assume that $(P:W\rightarrow X,\mu)$ is a strong bundle over the ep-groupoid $X$ which admits sc-smooth bump functions, and which has, in addition, a paracompact orbit space $\abs{X}$.
Let $N$ be an auxiliary norm for $(P,\mu)$. 
Given a smooth point $x\in X$ with isotropy $G_x$, a smooth point $e\in W_x$  satisfying $e=\mu (g, e)$ for all $g\in G_x$,  an open neighborhood
$\wh{O}$ of $\pi(x)$ in $|X|$, and $\varepsilon>0$, then there exists a $\ssc^+$-section functor $\sigma\colon X\to W$ possessing  the following properties.
\begin{itemize}
\item[{\em (1)}]\ The support of $\sigma$ is contained  in $\pi^{-1}(\wh{O})$.
\item[{\em (2)}]\ $\sigma(t(\phi))=\mu (\phi, \sigma(s(\phi))$  for all $\phi \in {\bf X}$. 
\item[{\em (3)}]\ $\sigma(x)=e$.
\item[{\em (4)}]\ $N(\sigma(y))\leq N(e)+\varepsilon$ for all $y\in X$.
\end{itemize}
\end{theorem}

The hypothesis that $|X|$ is paracompact is used to insure the existence of auxiliary norms $N$.
The existence of sc-smooth bump functions is used for the actual construction of the functor $\sigma$.

\begin{proof}
We choose  an open neighborhood $U(x)$ of $x$ contained in $\pi^{-1}(\wh{O})$ and having  the following properties.
\begin{itemize}
\item[$\bullet$] The natural $G_x$-action is defined on $U(x)$. 
\item[$\bullet$] $t:s^{-1}(\cl_X(U(x)))\rightarrow X$ is proper.
\item[$\bullet$] There exists a strong bundle isomorphism $\Psi\colon W\vert U(x)\rightarrow K$ which covers  the sc-diffeomorphism $\psi\colon U(x)\rightarrow  O$ satisfying $\psi (x)=0$.  Here we use the local model $(O,C,E)$ for the object M-polyfold $X$ and the model $(K, C\triangleleft F, E\triangleleft F)$ for the local strong bundle 
retract $K=R(U\triangleleft F)$ defined by the retraction $R(u, f)=(r(u), \rho (u)f)$ for $(u, f)\in U\oplus F$. The sc-smooth map $r\colon U\to U$ is a retraction onto $O=r(U)$. The diagram 
$$
\begin{CD}
P^{-1}(U(x))@>\Psi>> K\\
@VV PV  @VV p V\\
U(x)@>\psi>>O
\end{CD}
$$
commutes.
\end{itemize}
We choose a smaller  neighborhood $V(x)$, still invariant under the action of $G_x$,  such  that $\cl_X(V(x))\subset U(x)$ and denote by $\wtilde{V}(x)$ the saturation of $V(x)$.
By construction, $\wtilde{V}(x)\subset \pi^{-1}(\wh{O})$. Then  we take on  the object M-polyfold $X$ a sc-smooth bump
function $f\colon U(x)\rightarrow [0,1]$ satisfying  $f(x)=1$ and having its support in $V(x)$. We also assume that $f$ is invariant under the $G_x$-action.
As used already in the  previous proof  there is a canonical way to extend it to a sc-smooth bump function for the ep-groupoid $X$.
We denote this extension by $\beta$. 

The construction of the desired section $\sigma$ of 
$P\colon W\vert U(x)\to U(x)$ starts in the bundle chart 
$\Psi$. Recalling that $\psi (x)=0$, we denote the image of the distinguished smooth point $e\in W_x$ under the strong bundle isomorphism $\Psi\colon
 W\vert U(x)\to K$ by $\Psi (e)=(0, e')$. We define the section $s'$ of the local bundle $p\colon K\to O$ by 
 $$s'(m)=R(m, e')=(m, \rho (m)e')$$
 for all $m\in O$, so that $s'(0)=R(0, e')=(0, e')$ because $R$ retracts onto $K$.
 
 The strong bundle isomorphism $\Psi\colon  W\vert U(x)\to K$ induces the $G_x$-action on $K$, defined by 
 $$g\ast \Psi (w)=\Psi (g\ast w)\quad \text{and}\quad 
 g\ast \psi (y)=\psi (g\ast y),$$
 for all $g\in G_x$ and $w\in  W\vert U(x)$ and $y\in O$, where $g\ast w=\mu \bigl(\Gamma (g, P(w)), w\bigr)$. Recalling the assumption $e=\mu (g, e)$ on $e\in W_x$ and the property $\Gamma (g, x)=g$ for $g\in G_x$, we deduce that $g\ast e=e$ and 
 $g\ast \Psi (e)=\Psi (g\ast e)=\Psi (e)$ and hence 
 $g\ast (0, e')=(0, e')$ for all $g\in G_x$.  
 
 If $h(y)=\Psi^{-1}(s'(\psi (y)))$ is the pull-back section $h\colon U(x)\to W\vert U(x)$, we introduce for every $g\in G_x$ the section $h_g\colon U(x)\to W\vert U(x)$ by 
 $$h_g(y)=\mu \bigl(\Gamma (g, g^{-1}\ast y), h(g^{-1}\ast y)\bigr)
 $$
 for all $y\in U(x)$. It satisfies $h_g(e)=e$. Averaging 
 these sections over $G_x$ we obtain the section 
 $\sigma''\colon U(x)\to W\vert U(x)$, given by 
 \begin{equation*}
 \begin{split}
 \sigma''(y)&=\dfrac{1}{\abs{G_x}}\sum_{g\in G_x}h_g(y)\\
&=\dfrac{1}{\abs{G_x}}\sum_{g\in G_x}\mu \bigl(\Gamma (g, g^{-1}\ast y), h(g^{-1}\ast y)\bigr).
\end{split}
\end{equation*}
Then $\sigma'' (e)=e$ and we claim that $\sigma''$ is $G_x$-equivariant in the sense that 
\begin{equation}\label{sdouble_equiv}
\sigma''(\gamma\ast y)=\gamma\ast \sigma'' (y)=
\mu \bigl(\Gamma (\gamma, y),\sigma'' (y)\bigr)
\end{equation}
for all $y\in U(x)$ and $\gamma\in G_x$. Indeed, recalling that $g^{-1}\ast (\gamma \ast y)=(g^{-1}\circ \gamma)\ast y$, the left-hand side becomes
\begin{equation}\label{sdouble_equiv_1}
\sigma''(\gamma\ast y)=\dfrac{1}{\abs{ G_x}}\sum_{g\in G_x}
\mu \bigl(\Gamma (g, (g^{-1}\circ \gamma)\ast y), h((g^{-1}\circ \gamma)\ast y)\bigr).
\end{equation}
By the properties of $\Gamma$,
\begin{equation*}
\Gamma (g, (g^{-1}\circ \gamma)\ast y)=
\Gamma (\gamma \circ \gamma^{-1}\circ g, (g^{-1}\circ \gamma)\ast y)=
\Gamma (\gamma, y)\circ 
\Gamma (\gamma^{-1}\circ g, (g^{-1}\circ \gamma)\ast y)
\end{equation*}
Hence, by the properties of the structure map $\mu$,
\begin{equation}\label{sdouble_equiv_2}
\begin{split}
\mu &\bigl(
\Gamma (g, (g^{-1}\circ \gamma)\ast y), h((g^{-1}\circ \gamma)\ast y)\bigr)\\
&=\mu \bigl(
\Gamma (\gamma, y)\circ 
\Gamma (\gamma^{-1}\circ g, (g^{-1}\circ \gamma)\ast y), h((g^{-1}\circ \gamma)\ast y)\bigr)\\
&=\mu \bigl(
\Gamma (\gamma, y), \mu\bigl(
\Gamma (\gamma^{-1}\circ g, (g^{-1}\circ \gamma)\ast y), h((g^{-1}\circ \gamma)\ast y)\bigr)\bigr).
\end{split}
\end{equation}
Therefore, inserting \eqref{sdouble_equiv_2} into \eqref{sdouble_equiv_1}, changing to the summation index $g^{-1}\circ \gamma$ and using the linearity of the map 
$\mu \bigl(\Gamma (\gamma, y), \cdot  \bigr)\colon W_y\to W_{\gamma\ast y}$ we arrive at the formula 
$\sigma'' (\gamma \ast y)=\mu \bigl(\Gamma (\gamma, y), \sigma'' (y)\bigr)$ as claimed. 

Since every morphism $\phi$ connecting two points in $U(x)$ has the form $\phi=\Gamma (g, y)$ for $s(\phi )=y$ and $t(\phi )=g\ast y$ and a unique $g\in G_x$, we deduce from \eqref{sdouble_equiv_2} that 
\begin{equation}\label{sdouble_equiv_4}
\sigma'' (t(\phi))=\mu\bigl(\phi, \sigma'' (s(\phi))\bigr)
\end{equation}
for every morphism $\phi$ between points in $U(x)$. 

The section $\sigma'\colon U(x)\to W\vert U(x)$ defined by $\sigma' (y)=\beta (y)\cdot \sigma'' (y)$ is again a $G_x$-equivariant $\ssc^+$-section. In addition, 
$\sigma'$ has its support in $V(x)$ and satisfies $\sigma' (e)=e$. Finally, we proceed as in the proof of 
Theorem \ref{t-bump_functions} and use the morphisms in $\bx$ to extend the section $\sigma'$ from $U(x)$ to the desired $\ssc^+$-section functor $\sigma\colon X\to W$. It satisfies the properties (1)-(3). The property (4) is achieved by choosing $V(x)$ small enough,
see also Remark \ref{extension-Q} below. The proof of Theorem \ref{LLLOP} is complete.
\qed \end{proof}

\begin{remark}\index{R- Transversality versus functoriality}
As we shall see later on in the Fredholm theory on ep-groupoids or polyfolds, the property of a section of a strong bundle of being functorial
and the property of being transversal are competing and  do quite often  not hold simultaneously. This forces us
to use multivalued sections. Note that in the previous theorem there would be no such section functors $\sigma$ if we would not have required
that $e$ has $G_x$ as an automorphism group. Later on we shall see that such (restricted) section functors $\sigma$ used as perturbations will not be
sufficient to achieve transversality. 
\qed
\end{remark}

Another remark is concerned with the above proof, whose ingredients allow to define $\ssc^+$-sections having  properties required in 
situations in which we try to achieve  transversality.
\begin{remark}\label{extension-Q}\index{R- Extensions of sc$^+$-sections}
Given a smooth point $x\in X$ and a $G_x$-invariant open neighborhood $U(x)$ of $x$, one constructs a $\ssc^+$-section $s$ with desired properties
near $x$ and compatible with the $G_x$-action, i.e., compatible with morphisms whose targets and sources belong to  $U(x)$. Then we can take 
an invariant sc-smooth function  $f:X\rightarrow [0,1]$ with support in $U(x)$ and satisfying $f=1$ near $x$ so that the product $f\cdot s$ is a  $\ssc^+$-section on $U(x)$.
The proof of Theorem \ref{LLLOP} shows how to extend $f\cdot s$ to a $\ssc^+$-section defined on all of $X$. The construction, i.e.,  the choice of $f$
allows to control the size of $f\cdot s$ with respect to the auxiliary norm $N$ and therefore controls the extension $\wtilde{s}$, which is entirely determined by $f\cdot s$.
As this remark shows,  local constructions can always be extended in a controlled way. In applications we shall see local constructions,
where the section attains particular values, or where the linearization at $x$ has certain properties. Basically any finite-dimensional idea
about a perturbation can be implemented in the polyfold set-up.\qed
\end{remark}

We come to the sc-smooth extension of a $\ssc^+$-section defined over the boundary $\partial X$ of an ep-groupoid.
We recall the discussion of local faces in Part I. Namely,  if $(O,C,E)$ is a tame local model,
then there exists an open subset $U$ of $C$ and a sc-smooth retraction $r:U\rightarrow U$ preserving the degeneracy index $d_C$, having $O$ as its image, and for every smooth point $x\in O$ there exists a sc-subspace $A$ of $E_x$ satisfying $T_xO\oplus_{sc} A=E_x$.
\begin{definition}\index{D- $\bssc^+$-section over $\partial X$}
We consider a strong bundle $(P:W\rightarrow X,\mu)$ over the tame ep-groupoid $X$ and assume that  the map $s:\partial X\rightarrow W$ satisfies $P\circ s=Id$. 
Then  $s$ is called  {\bf $\ssc^+$-section defined over the boundary}  provided for every point $x\in \partial X$ and every local face $F$ containing $x$, the restriction 
$s\vert F$ is a $\ssc^+$-section of the bundle $W|F\to F$. 
\qed
\end{definition}
The boundary  $\partial X$ is not a M-polyfold. However,  locally it is the union of a finite number of faces due to the tameness assumption.
These faces are M-polyfolds, which allows us to define some kind of sc-smoothness via the properties of all local restrictions. This control is enough to guarantee
local $\ssc^+$-smooth extensions. 

Recall that for a tame retract $(O,C,E)$ the equality $d_O(x)=d_C(x)$ holds
for $x\in O$, see Proposition \ref{tame_equality}. In particular $\partial O=O\cap \partial C$.
\begin{proposition}\label{hucky}\index{P- Extension of $\ssc^+$-sections}
Let $K\rightarrow O$ be a tame  local strong bundle model, where $K=R(U\triangleleft  F)$ is a strong bundle retract and the strong bundle retraction $R$ covers  the tame retraction
$r:U\rightarrow U$ of the  local model $(O,C,E)$. We assume that $s$ is a $\ssc^+$-section defined on an open neighborhood in $\partial O$
of some $x_0\in \partial O$. Then there exists an open neighborhood $V(x_0)\subset O$ and a $\ssc^+$-section 
$\wtilde{s}$ of $W|V(x_0)$ such that $\wtilde{s}=s$ on $\partial O$ near $x_0$.
\end{proposition}
\begin{proof}
There is no  loss of generality to assume  that $E={\mathbb R}^n\oplus H$, $C=[0,\infty)^n\oplus H$, and $x_0=(0,h_0)$. Possibly $h_0\neq 0$ 
if $x_0$ is not smooth. For a subset $I$ of $\{1,\ldots ,n\}$ and a point $(t,h)$ near $(0,h_0)$, we define 
$$
(t,h)_I = r(t_I,h),
$$
where  $t_I$ is obtained from $t$ by putting  $t_i=0$ if $i\in I$. We define $\wtilde{s}$ on an open neighborhood $V(x_0)$ in $O$ by the following formula
$$
\wtilde{s}(t,h) = \sum_{j=1}^n (-1)^{j+1} \sum_{|I|=j} R((t,h),s((t,h)_I)).
$$
The maps $(t,h)\rightarrow s((t,h)_I)$ are $\ssc^+$-sections of the bundle $O\triangleleft F\to O$ near $x_0$,   implying that $\wtilde{s}$ is a local $\ssc^+$-section.
The following straightforward computation shows that $\wtilde{s}=s$ on $\partial O$ near $x_0$. Namely, considering $(t,h)\in O$ with $t_{i_0}=0$ near $x_0$ and using 
$(t,h)_{\{i_0\}}=r(t_{\{i_0\}},h)=r(t,h)=(t,h)$, we compute 
\begin{equation*}
\begin{split}
&\wt{s}(t, h)=\quad \sum_{j=1}^n (-1)^{j+1} \sum_{|I|=j} s((t,h)_I)\\
&\quad=\sum_{j=1}^n (-1)^{j+1} \sum_{|I|=j, i_0\in I} s((t,h)_I)+ \sum_{j=1}^n (-1)^{j+1} \sum_{|I|=j, i_0\not\in I} s((t,h)_I)\\
&\quad=\sum_{j=1}^n (-1)^{j+1} \sum_{|I|=j-1, i_0\not\in I} s((t,h)_{I\cup\{i_0\}})+ \sum_{j=1}^n (-1)^{j+1} \sum_{|I|=j, i_0\not\in I} s((t,h)_I)\\
&\quad=\sum_{j=0}^{n-1} (-1)^{j} \sum_{|I|=j, i_0\not\in I} s((t,h)_{I})+ \sum_{j=1}^n (-1)^{j+1} \sum_{|I|=j, i_0\not\in I} s((t,h)_I)\\
&\quad=\sum_{j=0}^{n-1} (-1)^{j} \sum_{|I|=j, i_0\not\in I} s((t,h)_{I})+ \sum_{j=1}^{n-1} (-1)^{j+1} \sum_{|I|=j, i_0\not\in I} s((t,h)_I)\\
&\quad= s((t,h)_{\{i_0\}})=s(r(t,h))=s(t,h).
\end{split}
\end{equation*}
\qed \end{proof}
There is a  corollary in the presence of symmetries. 
\begin{corollary}\label{ext_t_t}\index{C- Equivariant extension}
Let  $(O,C,E)$ be  as described in Proposition \ref{hucky},  and $x_0\in \partial O=\partial C\cap O$  a point such  that $G_{x_0}$ acts by sc-diffeomorphisms 
on $O$ with a lift to a strong bundle model $P:K\rightarrow O$. We assume that  $s$ is a $\ssc^+$-section defined on an $G_{x_0}$-invariant open neighborhood of $x_0$ in $\partial O$ and satisfying
$s(g\ast x) =g\ast s(x)$ for $x\in \partial O$ near to $x_0$. Then on a suitable $G_{x_0}$-invariant neighborhood of $x_0$ in $O$ there exists a $\ssc^+$-section $\wh{s}$, which is invariant under the $G_{x_0}$-action and which restricts to $s$ on  $\partial O$ near $x_0$.
\end{corollary}
\begin{proof}
Averaging the extension $\wt{s}$ guaranteed by Proposition \ref{hucky},  
we obtain the local section $\wh{s}$ near $x_0$ in $O$, defined by 
$$
\wh{s}(x)=\frac{1}{|G_{x_0}|} \sum_{g\in G_{x_0}} g^{-1}\ast \left( \wt{s}(g\ast x)\right).
$$
Since $s$  is   by assumption  equivariant,  it follows that $\wh{s}=s$ on $\partial O$ near $x_0$.
\qed \end{proof}

Using this local result we  prove the following extension result.
\begin{theorem}[{\bf Extension of $\ssc^+$-section functors}]\label{extension-ssc-section functors}\index{T- Extensions of $\ssc^+$-functors}
Let  $(P:W\rightarrow X,\mu)$ be  a strong bundle over the  tame ep-groupoid $X$ whose  
orbit space $|X|$ is paracompact.
Suppose that the ep-groupoid admits sc-smooth partitions of unity.  Given a $\ssc^+$-section functor $s$ of $W|\partial X$
and a saturated open neighborhood $U$ of  $\partial X$ in $X$, there exists exists a $\ssc^+$-section functor $\wtilde{s}$ of $P$
supported  in $U$,  which on $\partial X$ restricts to $s$. Moreover,  if an auxiliary norm $N$  and $\delta>0$  are given such that
$N(s(x))<\delta $ for all $x\in \partial X$, then the extension $\wtilde{s}$ can be constructed in such a way that $N(\wtilde{s}(x))\leq\delta $ for all $x\in X$.
\end{theorem}
\begin{proof}
For every $x\in \partial X$ we take the open neighborhood $U(x)\subset X$ guaranteed by Corollary \ref{ext_t_t} and contained in the open neighborhood $U\subset X$ of $\partial X$, such that the section $\wt{s}$ on $U(x)$ is an extension of $s\vert (U(x)\cap \partial X)$.
We extend $\wt{s}$ in an equivariant way to the saturation  $\wt{U}(x)=\pi^{-1}(\pi (U(x))$ of $U(x)$ and denote this extension by  $\wtilde{s}_x$. The complement  $U_\infty=X\setminus\partial X$ of the boundary is saturated and so we obtain a covering of $X$  by the saturated open sets $U_\infty$ and $\wt{U}(x)$, $x\in \partial X$.
Using Theorem \ref{t-bump_functions} (2), we choose a subordinated sc-smooth partition of unity consisting of the sc-smooth bump functions  $\beta_\infty$ and $\beta_x$  on $\wt{U}(x)$ for every $x\in \partial X$ according to Definition \ref{d-ep-groupoid_bump},  where $\beta_x\colon X\to [0,1]$ satisfies  $\beta_x(x)=1$ and $\supp (\beta_x)\subset \wt{U}(x)$.
Using  only the functions $\beta_x$ (and not $\beta_\infty$),  we  define the section 
$\wt{s}$ by 
$$
\wt{s}=\sum_{x\in \partial X}\beta_x \wt{s}_x,
$$
which is a sc-smooth section functor extending $s$ and supported in $U$.

The map  $z\mapsto  N(\wtilde{s}(z))$ is continuous and 
we find a saturated open neighborhood $V$ of $\partial X$, whose closure is contained in $U$ such  that $N(\wtilde{s}(z))<\delta$ on $V$. Then we 
take another saturated open neighborhood $V'$ of $\partial X$ whose  closure is contained in $V$. The two sets $V$ and $X\setminus \cl_X(V')$ constitute 
an open saturated cover  of  $X$. 
Arguing as before we find 
an associated sc-smooth partition of unity in which the sc-smooth bump function $\beta\colon X\to [0,1]$ has its support in $V'$ and satisfies $\beta=1$ on $\partial X$. Then the product $\beta\cdot \wt{s}$ is an extension of $s$ which satisfies $N((\beta\cdot \wt{s})(x))\leq \delta$ for all $x\in X$. The proof of Theorem 
\ref{extension-ssc-section functors} is complete.
\qed \end{proof}

\section{Compactness Properties}\label{SEC114}
The standing assumption for this section is that $(P:W\rightarrow X,\mu)$ is a strong bundle over an ep-groupoid. 
We assume that we have fixed an auxiliary norm $N:W_{0,1}\rightarrow {\mathbb R}^+$ (a functor). One can extend $N$ to $W$ by setting
$N(h)=+\infty$ if $h$ is not of bi-regularity $(0,1)$. Hence we may view $N$ as a functor $W\rightarrow {\mathbb R}^+\cup\{+\infty\}$.
Most of the results we shall describe here are straight forward generalizations of the results in Chapter \ref{sec_package}.
We shall denote for a section $f$ of $P$ by $S_f$ the full subcategory of $X$ associated to all zero-objects for $f$, i.e. all $x\in X$ with $f(x)=0$.
\begin{definition}\index{D- Compactness notions}\label{DEFN12.4.1}
Let $f$ be an sc-smooth section functor of $P$, where $(P:W\rightarrow X,\mu)$ is a strong bundle over 
a tame ep-groupoid and $N$ an auxiliary norm. 
\begin{itemize}
\item[(1)]\ We say $f$ has {\bf compact solution set} (in the ep-groupoid sense) \index{Compact solution set, ep-groupoid case} if the orbit space $|S_f|$ associated to $S_f$ is compact.
\item[(2)]\ The section $f$ is said to be {\bf  proper}\index{Proper, ep-groupoid case} if $|S_f|$ is compact and for every auxiliary norm $N$
there exists an open saturated neighborhood $U$ of $S_f$ with the property that the orbit space of the set
$$
\{x\in U\ |\ N(f(x))\leq 1\}
$$
has compact closure.
\end{itemize}
\qed
\end{definition}
These two compactness notions are equivalent for sc-Fredholm sections functors.
The ep-groupoid version of Theorem \ref{x-cc} is given by the following result.
\begin{theorem}\label{THMB1242}
Let $f$ be an sc-Fredholm section of the strong bundle $(P:W\rightarrow X,\mu)$ over an ep-groupoid with paracompact orbit space
$|X|$, which guarantees the existence of auxiliary norms.  The following two statements are equivalent.
\begin{itemize}
\item[{\em(1)}]\ $f$ is proper.
\item[{\em(2)}]\ $f$ has a compact solution set.
\end{itemize}
\end{theorem}
\begin{proof}
(1) $\Rightarrow$ (2).  It is clear that properness implies a compact solution set in the ep-groupoid sense.\par

\noindent (2) $\Rightarrow$ (1). Assume next that $f$ has a compact solution set and an auxiliary norm $N$ is given.  Then we take for every $z\in |S_f|$ 
a representative $x_z\in X$ and we find, employing Theorem \ref{x-cc},  an open neighborhood $U(x_z)$ admitting the natural action
and having the properness property,   so that  
$\Sigma_z:=\{y\in U(x_z)\ |\ N(f(y))\leq 1\}$ satisfies
\begin{itemize}
\item $\cl_{X}(\Sigma_z)$ is compact.
\end{itemize}
This property persists if we replace $U(x_z)$ by a smaller open neighborhood, for which we may assume that it admits the 
natural $G_{x_z}$-action and that it has the property that $\cl_{|X|}(|U|)=|\cl_{X}(U)|$, see Proposition \ref{Prop7.1.9}.
Then the  orbit set 
$$
V_z := |U(x_z)|
$$
is an open neighborhood of $z$. We note that 
$\cl_{|X|}(|\Sigma_z|)$ is compact and contained in $|\cl_{X}(U(x_z))|$.

Since $|S_f|$ is compact we find finitely many $z_1,...,z_k$ 
so that the corresponding $V_{z_i}$ cover $|S_f|$. 
Denote by $U$ the  saturated open neighborhood of $S_f$ defined by
$$
U=\bigcup_{i=1}^k \pi^{-1}(V_{z_i}).
$$
Then
$$
\cl_{|X|}(|\{x\in U\ |\ N(f(x))\leq 1\}|)\subset 
\bigcup_{i=1}^k \cl_{|X|}( |\Sigma_{z_i}|),
$$
 which proves the claim.
\qed \end{proof}

In order to formulate an ep-groupoid version 
of Theorem \ref{THM528} we need the notion of the reflexive local compactness property in the ep-groupoid context.
\begin{definition}\label{DEFR1243}\index{D- Reflexive local compactness property, ep-groupoid context}
Assume that $P:Y\rightarrow X$ is a strong bundle over an ep-groupoid with reflexive $(0,1)$-fiber and let $f$ be an sc-smooth section functor of $P$, which in addition is regularizing.
We say that $f$ has the {\bf ep-groupoid reflexive local compactness property}  provided for every reflexive auxiliary norm $N$ (a functor)
and every point $x\in X$ there exists an open neighborhood $U(x)$ in $X$, i.e. on level $0$, so that 
$\cl_X\left(  \{y\in U(x)\ |\ N(f(y))\leq 1\}\right) $ is a compact subset of $X$. 
\qed
\end{definition}
Let us discuss this definition in more detail. We need the following fact.
\begin{lemma}\label{LEMM1244}
Let $X$ be an ep-groupoid. The following holds.
\begin{itemize}
\item[{\em(1)}]\ For a saturated subset $\Sigma$ of $X$ the equality $\cl_{|X|}(|\Sigma|) =|\cl_X(\Sigma)|$ holds.
\item[{\em(2)}]\ For a subset $\Sigma$ of $X$  with compact closure the equality $\cl_{|X|}(|\Sigma|) =|\cl_X(\Sigma)|$ holds.
\end{itemize}
\end{lemma}
\begin{proof}
For any subset $\Sigma$ the continuity of $\pi:X\rightarrow |X|:x\rightarrow |x|$ implies the inclusion
 $|\cl_X(\Sigma)|\subset \cl_{|X|}(|\Sigma|) $.\par
 
 \noindent (2)  Assume that $z\in \cl_{|X|}(|\Sigma|)$. We pick a representative $x$ of $z$ and an open neighborhood
 $U(x)$ with the natural $G_x$-action and the properness property. Let $(z_k)\subset |\Sigma|$ with $z_k\rightarrow z$.
 There is no loss of generality assuming that $(z_k)\subset |U(x)|$. We can pick $(x_k)\subset U(x)$ with
 $z_k=|x_k|$ and must have $x_k\rightarrow x$. Since $\Sigma$ is saturated $x_k\in\Sigma$ and we conclude
 that $x\in \cl_X(\Sigma)$, i.e. $z\in |\cl_X(\Sigma)|$.\par
 
 \noindent (2) Pick $z\in \cl_{|X|}(|\Sigma|)$ and  $(z_k)\subset |\Sigma|$
 converging to $z$. We can pick $(x_k)\subset \Sigma$ with $|x_k|=z_k$.
 After perhaps taking a subsequence and using the compactness of $\cl_X(\Sigma)$ we may assume
 that $x_k\rightarrow x\in\cl(\Sigma)$. Then $z=|x|\in |\cl_X(\Sigma)|$.
\qed \end{proof}

Assume that  $f$ has the reflexive local compactness property.
 Given a point $x\in X$ we find an open neighborhood $U(x)$, which is invariant under the $G_x$-action,  having the properness property, i.e.
 $t:s^{-1}(\cl_X(U(x)))\rightarrow X$ is proper, so that 
 $$
 \cl_X(\{y\in U(x)\ |\ N(f(y))\leq 1\})\ \ \text{is compact}.
 $$
Denoting by $\wt{U}(x)$ the saturation of $U(x)$ and using that $f$ and $N$ are functors we observe that 
\begin{eqnarray}\label{EQNNH}
&& \{\wt{y}\in \wt{U}(x)\ |\ N(f(\wt{y}))\leq 1\}\\
& =& \{y\in X\ |\ \exists \ \phi\in\bm{X},\ s(\phi)=y,\ t(\phi)\in U(x)\ \text{and}\ \ N(f(t(\phi)))\leq 1\}\nonumber
\end{eqnarray}
are saturated sets and that 
$$
\cl_X(\{y\in U(x)\ |\ N(f(y))\leq 1\}) \ \ \text{is compact}.
$$
Hence 
$$
\cl_{|X|}(|\{y\in \cl_X(U(x))\ |\ N(f(y))\leq 1\}|) = |\cl_X(\{y\in U(x)\ |\ N(f(y))\leq 1\})|.
$$
Consequently, using (\ref{EQNNH}) 
\begin{eqnarray*}
&& |\cl_X(\{\wt{y}\in \wt{U}(x)\ |\ N(f(\wt{y}))\leq 1\})|\\
 &=&\cl_{|X|}(|\{\wt{y}\in \wt{U}(x)\ |\ N(f(\wt{y}))\leq 1\}|)\\
 &=& |\cl_X(\{y\in U(x)\ |\ N(f(y))\leq 1\})|.
\end{eqnarray*}
Now  the following proposition is easily obtained.
\begin{proposition}
Assume that $P:Y\rightarrow X$ is a strong bundle over an ep-groupoid with reflexive $(0,1)$-fiber and let $f$ be an sc-smooth section functor of $P$, which in addition is regularizing. Then the following three statements are equivalent.
\begin{itemize}
\item[{\em(1)}]\ $f$ has the ep-groupoid reflexive local compactness property.
\item[{\em(2)}]\ For every reflexive auxiliary norm $N$  for $(P:W\rightarrow X,\mu)$ and $z\in|X|$ there exists a saturated open subset $\wt{U}$
with $z\in |\wt{U}|$ so that $|\cl_X(\{y\in \wt{U}\ |\ N(f(y))\leq 1\})|$ is compact.
\item[{\em(3)}]\  For every reflexive auxiliary norm $N$  for $(P:W\rightarrow X,\mu)$ and $z\in|X|$ there exists a saturated open subset $\wt{U}$
with $z\in |\wt{U}|$ so that $\cl_{|X|}(|\{y\in \wt{U}\ |\ N(f(y))\leq 1\}|)$ is compact.
\end{itemize}
\qed
\end{proposition}

A basic result, which generalizes Theorem \ref{THM528} is the following.

\begin{theorem}\index{T- Extension of controlling $(U_\partial, N)$, ep-groupoid context}\label{THMXXX1246}
Assume that $(P:W\rightarrow X,\mu)$ is a strong bundle over the tame ep-groupoid $X$ with paracompact orbit space $|X|$, where $W$ has reflexive $(0,1)$-fibers. Denote by $N$ a reflexive auxiliary norm and let $f$ be an sc-Fredholm functor with compact solution set
having the reflexive local compactness property. Defining  $S_{\partial f}=\{x\in \partial X\ |\ f(x)=0\}$  let $U_\partial\subset \partial X$ be an open saturated neighborhood in $\partial X$ of $S_{\partial f}$ so that $\cl_{|X|}(|\{y\in U_\partial\ |\ N(f(y))\leq 1\}|)$ is compact in $|X|$. Then there exists a saturated open neighborhood  $U$  of $S_f=\{x\in X\ |\ f(x)=0\}$ with the following properties.
\begin{itemize}
\item[{\em(1)}]\ $U\cap\partial X= U_\partial$.
\item[{\em(2)}]\ The closure of $|\{x\in U\ |\ N(f(x))\leq 1\}|$ in $|X|$ is compact.
\end{itemize}
\end{theorem}
\begin{proof}
For $z\in |U_{\partial} |$ we take a representative $x_z$ of $z$. Following the argument of Theorem \ref{THM528}
we find an open neighborhood $U(x_z)$  in $X$ such that one of the following properties holds.
\begin{itemize}
\item[(a)]\ \ $N(f(x))>1$ for $x\in U(x_z)$ and $U(x_z)\cap \partial X\subset U_\partial$.
\item[(a')]\ \  $\cl_X(\{x\in U(x_z)\ |\ N(f(x))\leq 1\})$ is compact and $U(x_z)\cap \partial X\subset U_\partial$.
\end{itemize}
In addition to one of these properties we can always assume the following
\begin{itemize}
\item[(b)] \ $U(x_z)$ admits the natural $G_{x_z}$-action.
\item[(c)]\  $\cl_{|X|}(|U(x_z)|) = |\cl_X(U(x_z))|$.
\end{itemize}
By construction, we obtain for every $z\in |U_{\partial}|$ an open neighborhood $U(x_z)$
with either the properties (a), (b), (b) or the properties (a'), (b), (c).

The set  $\cl_{|X|}(|\{x\in U_\partial\ |\ N(f(x))\leq 1\}|)$ is compact.  We find $\bar{z}_1,...,\bar{z}_\ell\in |U_{\partial} |$ so that 
$$
|\cl_{X}(\{x\in U_\partial\ |\ N(f(x))\leq 1\})|\subset |U(\bar{x}_1)|\cup..\cup |U(\bar{x}_\ell)|,
$$
where $\bar{x}_i=x_{\bar{z}_i}$. Denoting by $\wt{U}(\bar{x}_i)$ the saturation of $U(\bar{x}_i)$ we obtain, since $U_{\partial}$ is saturated that
$$
(\wt{U}(\bar{x}_1)\cup...\cup \wt{U}(\bar{x}_\ell))\cap \partial X\subset U_\partial.
$$
 If $z\in |U_\partial| \setminus  |(\wt{U}(\bar{x}_1)\cup..\cup \wt{U}(\bar{x}_\ell))|$
we must have $N(f(x_z))>1$. Consider the union 
$$
\wt{U}:= \wt{U}(\bar{x}_1)\cup..\cup \wt{U}(\bar{x}_\ell) \cup \bigcup_{z\in |U_\partial| \setminus |(\wt{U}(\bar{x}_1)\cup..\cup \wt{U}(\bar{x}_\ell))|} \wt{U}(x_z).
$$
Here $\wt{U}(x_z)$ is the saturation of $U(x_z)$. Clearly $\wt{U}(x_z)$ does not depend on the choice of $x$ and we shall define
$$
\wt{U}(z):=\wt{U}(x_z).
$$
Then $|\wt{U}(z)|=|U(x_z)|$ and we write with our modified notation
$$
\wt{U} = \wt{U}(z_1)\cup...\cup \wt{U}(z_\ell) \cup \bigcup_{z\in  |U_\partial| \setminus |(\wt{U}(z_1)\cup..\cup \wt{U}(z_\ell))|} \wt{U}(z).
$$
This is an open subset of $X$ satisfying $\wt{U}\cap \partial X=U_\partial$.   Moreover, 
the closure of the orbit space of $\{x\in \wt{U}\ |\  N(f(x))\leq 1\}$ is compact. Indeed, take a sequence $(x_k)\subset \wt{U}$ satisfying $N(f(x_k))\leq 1$.
We may assume up to isomorphism
that   $(x_k)\subset U(\bar{x}_1)\cup..\cup U(\bar{x}_\ell)$ and since the closure of each of the finitely many sets $\{x\in U(\bar{x}_i)\ |\ N(f(x))\leq 1\}$
is compact, we see that $(x_k)$ has a convergent subsequence. The orbit space of the set $\{x\in X\setminus \wt{U}\ |\ f(x)=0\}$ is compact 
and we find finitely many $U(x_i)\subset X\setminus \partial X$, say $i=\ell+1,...,e$, whose orbit spaces cover this compact set so that in addition the closure of each set $\{x\in U(x_i)\ |\ N(f(x))\leq 1\}$ is compact. Moreover, we assume as usual that these sets have the natural action on them.
Finally define with $z_i=|x_i|$ and $\wt{U}(z_i)$ the saturation of $U(x_i)$
$$
U = \wt{U}\cup \wt{U}(z_{\ell+1})\cup...\cup \wt{U}(z_e).
$$
Then $U\cap \partial X=U_\partial$ and $|\cl_X(\{x\in U\ |\ N(f(x))\leq 1\})|$ is compact. The latter
satisfies
$$
|\cl_X(\{x\in U\ |\ N(f(x))\leq 1\})|=\cl_{|X|}(|\{x\in U\ |\ N(f(x))\leq 1\}|)
$$
 in view of Lemma \ref{LEMM1244} since $U$ is saturated.
The proof is complete.
\qed \end{proof}

\section{Orientation Bundle}\label{SEC115}

For many applications of the sc-Fredholm theory orientation questions play an important role.
We already discussed these issues in the context of the M-polyfold theory and we shall generalize the discussion to cover the case
where we deal with an sc-Fredholm section functor of a strong bundle over an ep-groupoid.
We shall denote by $(P:W\rightarrow X,\mu)$ a strong bundle over an ep-groupoid and by $f$ an sc-Fredholm section functor.
Given a smooth $x\in X$ we consider the affine  space of linearizations $\text{Lin}(f,x)$ introduced in Definition \ref{LINUXX}.
The elements of $\text{Lin}(f,x)$ are all linear sc-operators $S:T_xX\rightarrow W_x$ obtained by taking a local sc$^+$-section $s$
satisfying $s(x)=f(x)$ and putting $S:= (f-s)'(x)$. We do not require $s$ to be compatible with the action of the isotropy group $G_x$ of $x$.
Given a morphism $\phi:x\rightarrow x'$ between smooth points we obtain the sc-isomorphism $T\phi:T_xX\rightarrow T_{x'}X$
and define
$$
\phi_\ast :\text{Lin}(f,x)\rightarrow \text{Lin}(f,x'): \phi_\ast S :=  \mu(\phi, S\circ (T\phi)^{-1}).
$$
Of course, we need to verify that this is well-defined.
\begin{lemma}
The following holds for smooth points $x$ and $x'$.
\begin{itemize}
\item[{\em(1)}]\ Given a morphism $\phi:x\rightarrow x'$ and $S\in \text{Lin}(f,x)$ the sc-operator $\phi_\ast S$ belongs to
$\text{Lin}(f,x')$. 
\item[{\em(2)}]\ The map $\phi_\ast :\text{Lin}(f,x)\rightarrow \text{Lin}(f,x')$ is a bijection.
\item[{\em(3)}]\  $({1_x})_\ast = Id_{\text{Lin}(f,x)}$. For morphisms $\phi$ and $\psi$ in $\bm{X}$ with $s(\psi)=t(\phi)$ 
it holds that $(\psi\circ \phi)_\ast =\psi_\ast \circ \phi_\ast$.
\end{itemize}
\end{lemma}
\begin{proof}
Assume that $\phi:x\rightarrow x'$ is a morphism between smooth points. Then $\phi$ extends locally to an sc-diffeomorphism
$\wh{\phi}:(U(x),x)\rightarrow (U(x'),x')$, which can be written as
$$
\wh{\phi}(y) = t\circ (s|U(\phi))^{-1}(y),
$$
where $U(x),U(y),U(\phi)$ are suitable open neighborhoods so that $s:U(\phi)\rightarrow U(x)$ and $t:U(\phi)\rightarrow U(x')$ are sc-diffeomorphisms. The tangent $T\phi: T_xX\rightarrow T_{x'}X$ is an sc-isomorphism and by definition the tangent map associated 
to $\wh{\phi}$ taken  at $x$.
If we take $U(x)$ small enough and in such a way that we have the natural action of $G_x$ 
defined on $U(x)$ the sc-diffeomorphism $\wh{\phi}:U(x)\rightarrow U(x')$ is equivariant with respect to the group homomorphism
$$
\gamma\colon G_x\rightarrow G_{x'}\colon g\rightarrow \gamma(g):= \phi\circ g\circ \phi^{-1}.
$$
We can lift the equivariant $\wh{\phi}$ to a strong bundle isomorphism
$$
\wh{\Phi}: W|U(x)\rightarrow W|U(x'): w\rightarrow \mu((s|U(\phi))^{-1}(P(w)),w)
$$
which at $x$ maps $W_x$ via $\mu(\phi,.)$ to $W_{x'}$. 

If $s_0$ is a local sc$^+$-section defined near $x$ with $s_0(x)=f(x)$. Then $s_1= \wh{\Phi}_\ast s_0$ is a local sc$^+$-section at $x'$
satisfying $s_1(x')=f(x')$. We note that 
$$
\mu(\phi, (f-s_0)'(x)\circ T\phi^{-1}) = (f-s_1)'(x').
$$
From this we see that $\phi_\ast S:= \mu(\phi, S\circ (T\phi)^{-1})$ is well-defined as a map
$$
\text{Lin}(f,x)\rightarrow \text{Lin}(f,x').
$$
This proves (1). The functorial properties in (3) are obvious from the previous discussion.
From (1) and (3) the bijectivity statement (2) is obvious.
\qed \end{proof}
Recall from Section \ref{ORIENTXX}  that for a smooth object $x$ in the object M-polyfold $X$ 
the set $\text{Lin}(f,x)$ is a convex subset of ${\mathcal L}(T_xX,W_x)$, 
which we equip with the induced topology. As such it is a contractible space of Fredholm operators and we can consider 
the determinant bundle over it, which is denoted by $\text{DET}(f,x)$
$$
\text{DET}(f,x)\rightarrow \text{Lin}(f,x).
$$
This topological  line bundle has 
two different orientations which we refer to as the two possible orientations of $(f,x)$. 
There is in general not a preferred orientation, although this might happen in certain cases, for example in Gromov-Witten theory
the complex orientation.
If $\mathfrak{o}$ is one of the orientations we denote by $-\mathfrak{o}$ the other one.
An orientation 
$\mathfrak{o}$ of $(f,x)$ is by definition an orientation for the topological line bundle $\text{DET}(f,x)$.

Assume that $S\in\text{Lin}(f,x)$ and an orientation $\mathfrak{o}$ has been picked for $\text{DET}(f,x)$.
Then $\text{det}(S)$ obtains an orientation $\mathfrak{o}_S$. 
Given a morphism $\phi:x\rightarrow x'$ between smooth points the  orientation $\mathfrak{o}_S$ for $\text{det}(S)$ 
defines naturally an orientation $\phi_\ast\mathfrak{o}_S$ for $\text{det}(\phi_\ast S)$.  This in turn defines an orientation
for $\text{DET}(f,x')$. It does not depend on the specific choice of $S$.   Consequently given 
$\phi:x\rightarrow x'$ an orientation $\mathfrak{o} $ for $\text{DET}(f,x)\rightarrow \text{Lin}(f,x)$
is mapped to an orientation $\phi_\ast\mathfrak{o}$ for $\text{DET}(f,x)\rightarrow \text{Lin}(f,x')$.

We recall from Chapter \ref{ORIENTXX} that,  under the assumption 
that $X$ is tame, given a smooth point $x$ and an orientation for $(f,x)$ there exists a natural orientation  for all $(f,y)$ where $y$ is smooth and belongs to a sufficiently small open neighborhood of $x$. This means given $\mathfrak{o}_x$ there exists a preferred local germ
denoted by $\mathfrak{o}_{(x)}$, which associates to a smooth point $y$ near $x$ an orientation $\mathfrak{o}_{(x)}(y)$ which is obtained
by local propagation from $\mathfrak{o}_x$. There are precisely two such germs, namely $\pm \mathfrak{o}_{(x)}$. 

\begin{definition}\label{DEFNG1252}\index{D- Orientation bundle $\mathscr{O}_f$}
Let $(P:W\rightarrow X,\mu)$ be a strong bundle over a tame ep-groupoid and $f$ an sc-Fredholm functor.
 The set $\mathscr{O}_f$  consists of all tuples $(x, (\text{DET}(f,x),\mathfrak{o}))$ where $x\in X_\infty$, 
 $\text{DET}(f,x)$ stands for the topological line bundle $\text{DET}(f,x)\rightarrow \text{Lin}(f,x)$ and 
 $\mathfrak{o}$ is an orientation for the latter. 
By $\sigma$ we denote the  natural projection
$$
\sigma\colon \mathscr{O}_f\rightarrow X_\infty,
$$
which is a surjective $(2:1)$-map.  $\mathscr{O}_f$ is turned into a topological space 
by equipping it with the topology ${\mathcal T}$  making every germ $\mathfrak{o}_{(x)}$ continuous. 
We call $(\mathscr{O}_f,{\mathcal T})$ the {\bf orientation bundle} associated to $f$.
\qed
\end{definition}
From the construction it is evident that given $\phi:x\rightarrow x'$ between smooth points the associated local family
$\wh{\Phi}$ covering $\wh{\phi}$ can be used to push forward a germ $\mathfrak{o}_{(x)}$ which is precisely the germ $\mathfrak{o}_{(x')}$,
which at $x'$ is the push forward of the orientation $\mathfrak{o}_x$. 
Let us summarize the previous discussion which follows essentially from the facts established in Chapter \ref{ORIENTXX}.

\begin{theorem}\index{T- Orientation bundle}
Let $(P:W\rightarrow X,\mu)$ be a strong bundle over a tame ep-groupoid and $f$ an sc-Fredholm section functor. 
The orientation bundle $\mathscr{O}_f$ associated to the sc-Fredholm section $f$ of $P:W\rightarrow X$ is a topological space
for which 
$$
\sigma: \mathscr{O}_f\rightarrow X_\infty
$$
is a local homeomorphism and the local sections $\mathfrak{o}_{(x)}$ are continuous.
Moreover given $\phi:x\rightarrow x'$   between smooth points the associated local $\wh{\Phi}$ pushes a  local continuous section 
 forward to a continuous local section. 
 \qed
\end{theorem}
Now we are in the position to define orientability of an sc-Fredholm section $f$ of a strong bundle over a tame ep-groupoid
$(P:W\rightarrow X,\mu)$.
\begin{definition}\label{DEFN1254}\index{D- Orientation for sc-Fredholm section functor}
An {\bf orientation} for the  sc-Fredholm section functor $f$ of $P:W\rightarrow X$, where $X$ is tame,  consists of a global continuous section $\mathfrak{o}$ of $\sigma:\mathscr{O}_f\rightarrow X_\infty$,
say
$$
x\rightarrow \mathfrak{o}_x
$$
having the property that for every smooth morphism $\phi\in \bm{X}$ it holds that $\phi_\ast\mathfrak{o}_{s(\phi)}=\mathfrak{o}_{t(\phi)}$.
\qed
\end{definition}
\begin{remark}\index{R- On orientations}
Having orientations for an sc-Fredholm section we obtain natural orientations for the solution spaces assuming
a certain amount of transversality.  Before we study more general cases, 
we consider as an example  a strong bundle $(P:W\rightarrow X,\mu)$ over an ep-groupoid with $\partial X=\emptyset$.
Let $f$ be a compact sc-Fredholm section and assume that for every $x\in X$ with $f(x)=0$ the linearization
is onto.  Further assume that $f$ is oriented by $\mathfrak{o}$. In this case, as  consequences of the compactness assumption
and the regularity assumption $L:=\{x\in X\ |\ f(x)=0\}$ is an \'etale proper  Lie groupoid, see \cite{AR},  in the classical sense having a compact orbit space.
The orientation $\mathfrak{o}$ induces a natural orientation for $T_xL =\text{ker}(f'(x))$ so that in fact $L$ has a natural orientation
and the morphisms respect these orientations. The orbit space $|L|$ is in a natural way a compact oriented orbifold.  
We shall consider generalizations of this result when starting with a $f$ not assuming the transversality condition.
In the case where functorial perturbations are not sufficient to achieve transversality we shall show that after a suitable perturbation
we obtain a branched weighted orbifold, see Chapter \ref{CHAPX16}.
\qed
\end{remark}
The construction of the orientation bundle $\mathscr{O}_f\rightarrow X_\infty$ as well as orientations for them (if they exist) behave well under generalized strong bundle isomorphisms, provided a necessary condition for orientability  holds.
We shall discuss this in more detail in the following.   Assume we are given a strong bundle over a tame ep-groupoid, say $(P:W\rightarrow X,\mu)$,
and an sc-Fredholm section functor $f$ of $P$. Associated to this data we obtain the orientation bundle $\sigma:\mathscr{O}_f\rightarrow X_\infty$ 
which is a topological bundle having as fibers two point sets. Moreover a smooth morphism $\phi:x\rightarrow x'$
defines a bijection between the corresponding fibers. A necessary condition for  orientability  of $f$, i.e. the existence of a continuous 
section functor $\mathscr{o}$ of $\sigma$ is, of course, that elements in $G_x$ for smooth objects $x$ act trivially on $\sigma^{-1}(x)$.
We leave it to the reader to prove the following sufficient and necessary condition of orientability for an sc-Fredholm section functor $f$.
The proof is straight forward.
\begin{proposition}\index{P- Criteria for  orientability}\label{PROPT1256}
An sc-Fredholm section $f$ of the strong bundle over a tame ep-groupoid  is orientable if and only if the following two conditions hold.
\begin{itemize}
\item[{\em(1)}]\ Passing from $\sigma:\mathscr{O}_f\rightarrow X_\infty$ to orbit spaces the continuous map $|\sigma|:|\mathscr{O}_f|\rightarrow |X_\infty|$ has above every class $|x|$ two preimages, i.e. it is $(2:1)$.
\item[{\em(2)}]\ $|\sigma|:|\mathscr{O}_f|\rightarrow |X_\infty|$ admits a global continuous section.
\end{itemize}
\qed
\end{proposition}
The construction of the orientation bundle has some naturality properties.
 \begin{theorem}[Naturality properties of $\mathscr{O}_f$]\index{T- Naturality properties of $\mathscr{O}_f$}\label{THMS1257}
Assume $(P:W\rightarrow X,\mu)$  and $(P':W'\rightarrow X',\mu')$  are strong bundles over tame ep-groupoids and $
\bar{\mathfrak{f}}:W\rightarrow W'$ is  a generalized strong bundle isomorphism covering the generalized isomorphism
$\mathfrak{f}:X\rightarrow X'$.  We assume that $f$ is an sc-Fredholm section functor of $P$ and denote by $f'$ the push-forward sc-Fredholm section.
\begin{itemize}
\item[{\em(1)}]\  With $\sigma$ and $\sigma'$ being the orientation bundles, the generalized 
strong bundle isomorphism induces a fiber-preserving homeomorphism $|\bar{\mathfrak{f}}|_\ast$ fitting into the following commutative diagram
$$
\begin{CD}
|{\mathscr{O}}_f| @> {|\bar{\mathfrak{f}}|}_\ast>> |{\mathscr{O}}_{f'}|\\
@V |\sigma| VV  @V |\sigma'| VV\\
|X_{\infty}| @>|\mathfrak{f}|>> |X'_{\infty}|.
\end{CD}
$$
The construction of $|\bar{\mathfrak{f}}|_\ast$ is functorial.  
\item[{\em(2)}]\ Assume that 
 $\mathsf{o}$ is a continuous section functor of $\sigma$, so that in particular the isotropy groups $G_x$, $x\in X_\infty$, act trivially on $\sigma^{-1}(x)$.
Then $\bar{\mathfrak{f}}$ defines naturally a push-forward $\mathsf{o}'=\bar{\mathfrak{f}}_\ast\mathsf{o}$,
which is a continuous section functor of $\mathscr{O}_{f'}\rightarrow X_\infty'$. 
\end{itemize}
\end{theorem}
\begin{proof}
A representative for $\bar{\mathfrak{f}}$ is a diagram of strong bundle equivalences of the form
$$
\begin{CD}
W @< \Phi<<  V @>\Phi'>> W'\\
@V P VV     @V QVV @V P' VV\\
X @<F<< A @> F'>> X'.
\end{CD}
$$
Then the push-forward by $\Phi'$ of the pull-back of $f$ by $\Phi$ is precisely $f'$.  We 
denote by $f'':A\rightarrow V$ the pull-back of $f$ by $(\Phi,F)$.

Pick any smooth point $x'\in X_\infty'$.  We find a point $a\in A$ and a morphism $\phi'$ resulting in the diagram
$$
F(a)\xleftarrow{F} a\xrightarrow{F'} F'(a)\xrightarrow{\phi'} x'.
$$
Using previously introduced constructions this diagram defines a composition of bijections
$$
\sigma^{-1}(F(x))\rightarrow \sigma''^{-1}(a)\rightarrow \sigma'^{-1}(F'(a))\rightarrow\sigma'^{-1}(x').
$$
Different choices would give
$$
F(b)\xleftarrow{F} b\xrightarrow{F'} F'(b)\xrightarrow{\psi'} x'.
$$
Since $F$ and $F'$ are equivalences the two diagrams are related by morphisms where
the vertical arrows can be canonically filled in using that $F$ and $F'$ are as equivalences fully faithful.
$$
\begin{CD}
F(a)@< F<< a @>F'>> F'(a)@> \phi'>> x'\\
          @V VV                           @V VV       @V\psi'^{-1}\circ\phi'VV    @|\\
F(b)@<F<< b @>F'>> F'(b)@>\psi'>> x'
\end{CD}
$$
These diagrams have canonical lifts to maps between the bundle fibers. 
If the actions of the isotropy groups $G_x$ on the $\sigma^{-1}(x)$ are trivial and in addition we have
a continuous section functor $\mathsf{o}$ of $\sigma$ we can define $\mathsf{o}'$ as follows.
We consider $\mathsf{o}_{F(a)} = (F(a),\text{DET}(f,F(a)),\mathfrak{o}_{F(a)})$. Using the above diagram
and its lift to bundles we obtain an orientation for $\text{DET}(f',x')$ which we denote by $\mathsf{o}'_{x'}$.
Since $\mathsf{o}$ is compatible with morphisms the resulting $\mathsf{o}'_{x'}$ is the same 
if we start with $\mathsf{o}_{F(b)}$. The data in the diagram moves sc-smoothly if we start moving $x'$
which implies that our construction is not only well-defined, but also is continuous. 
We leave it to the reader that representing $\bar{\mathfrak{f}}$ and $\mathfrak{f}$ by a different 
diagram results in the same definition. This is, of course, a consequence of the definition when two diagrams
represent the same generalized isomorphism or bundle map.

If do not impose an assumption on the isotropy groups we see that the above procedure is still
well-defined on the level of orbit spaces and we obtain (1).
\qed \end{proof}


\chapter{\texorpdfstring{Sc$^+$}{sccc}-Multisections}\label{CHAPS13}
In this section we shall introduce the notion of structurable as well as structured sc$^+$-multisections, which will
be needed for a sophisticated perturbation theory.  We first recall some of the relevant properties of ep-groupoids.

\section{Structure Result}
The following fact about ep-groupoids will be used frequently and follows from a previously established result, see Proposition \ref{Prop7.1.9}. Let $X$ be an ep-groupoid, and denote 
as usual the associated morphism M-polyfold by $\bm{X}$.
For every  $x\in X$,  there exists an open neighborhood  $U(x)$ admitting  the natural $G_x$-action and  the target map $t :s^{-1}(\cl_X(U(x)))\rightarrow X$ is proper. We also know that $|U(x)|$ is an open subset of $|X|$ and 
the map $U(x)\rightarrow |U(x)| :y\rightarrow \abs{y}$,  induces a homeomorphism $_{G_{x}}\backslash U(x)\rightarrow |U(x)|$.

\begin{lemma}\label{simple_lemma}
Suppose that $U(x)$ is as described above and let $V$ be a $G_x$-invariant subset of $U(x)$. Then the following statement  are equivalent:
\begin{itemize}
\item[{\em(1)}]\  $\cl_X(V)\subset U(x)$.
\item[{\em(2)}]\  $\cl_{|X|}(|V|)\subset |U(x)|$.
\end{itemize}
\end{lemma}
\begin{proof}
Assuming  (1) we take $z\in \cl_{|X|}(|V|)$ and a sequence $(z_k)\subset  |V|$  such that $z_k\rightarrow z$ in $|X|$. 
Then choosing  points $y_k\in V$ and $y\in X$ satisfying $|y_k|=z_k$ and $|y|=z$, we  find a sequence of morphisms
 $(\phi_k)$ such that  $s(\phi_k)=y_k$ and $t(\phi_k)\rightarrow y$.
Using the properness property of the target map $t\colon s^{-1}(\cl_{X}(V))\to X$,  we may assume,  after perhaps taking a subsequence,  that $\phi_k\rightarrow \phi$ in $\bm{X}$. This implies the convergence $y_k\rightarrow y'$ where $y'$,  in view of the assumption (i), belongs to  $\cl_{X}(V)\subset U(x)$. 
Since $\phi$ is a morphism between $y'$ and $y$, we have $|y'|=\abs{y}=z\in |U(x)|$, proving the statement  (1).

In the other direction,  let $y\in \cl_X(V)$ and $(y_k)$ be a sequence in $V$ satisfying $y_k\rightarrow y$. Then $|y_k|\rightarrow |y|$ so that, by assumption, $|y|\in |U(x)|$. Hence there is a morphism $\phi\colon z\to y$ where $z\in U(x)$. Using that $s$ and $t$ are local sc-diffeomorphisms, we find open neighborhoods $W(z)$, $W(y)$, and $W(\phi)$ such that 
the source and target maps $s\colon W(\phi)\rightarrow W(z)$ and $t\colon  W(\phi)\rightarrow W(y)$ are sc-diffeomorphisms. 
For large $k$, we find morphisms $\phi_k\in W(\phi)$ satisfying $s(\phi_k)=z_k\in W(y)$ and $t(\phi_k)=y_k$. From $y_k\to y$, we conclude the convergence  $\phi_k\to \phi$ and $z_k\to z$.  Since $y_k, z_k\in U(x)$, there are $g_k\in G_x$ such that $\Gamma (g_k, z_k)=\phi_k$.
After perhaps taking a subsequence we may assume that $\Gamma (g, z_k)=\phi_k$. This implies that $y_k=g\ast z_k\to y=g\ast z$ and  shows 
that $y\in U(x)$.

\qed \end{proof}

A crucial concept is that of a good system of open neighborhoods defined next.
\begin{definition}\label{DEFF1}\index{D- Good system}\index{$\mathfrak{U}$}
A collection  $\mathfrak{U}={(U(x))}_{x\in X}$ of open neighborhoods of the points in the ep-groupoid $X$ is called
a {\bf good system} of open neighborhoods if the following holds.
\begin{itemize}
\item[(1)]\ For every $x\in X$, the target map $t\colon s^{-1}(\cl_X(U(x)))\rightarrow X$ is proper.
\item[(2)]\ Every open neighborhood  $U(x)$ is equipped with the natural $G_x$-action.
\item[(3)]\ For every two points $x,x'\in X$,  the open subset $\bm{U}(x, x')$\index{$\bm{U}(x,x')$} of $\bm{X}$, defined by 
$$
\bm{U}(x,x'):=\{\phi\in \bm{X}, \vert \,  \text{$s(\phi)\in U(x)$ and $t(\phi)\in U(x')$}\}, 
$$
has the following property. For every connected component $\Sigma\subset \bm{U}(x,x')$, 
the source map $s\colon \Sigma\rightarrow U(x)$ and the target map $t\colon \Sigma\rightarrow U(x')$ are sc-diffeo\-mor\-phisms onto open subsets of $U(x)$ and $U(x')$, respectively.
\end{itemize}
\qed
\end{definition}
The following obvious result allows us to introduce the notion of refinement of a good system of open neighborhoods.
\begin{proposition}\label{prop_obvious}
Let ${(U(x))}_{x\in X}$ be  a good system of open neighborhoods and $(\wt{U}(x))_{x\in X}$ a family of $G_x$-invariant open neighborhoods satisfying $\wt{U}(x)\subset U(x)$ for all $x\in X$. Then ${(\wt{U}(x))}_{x\in X}$ is a good system of open neighborhoods.
\qed
\end{proposition}
Given ${(U(x))}_{x\in X}$ and open neighborhoods $V(x)$ for every point $x\in X$ we can pick for every $x\in X$ 
a $G_x$-invariant open neighborhood $W(x)$ which is contained in $U(x)\cap V(x)$. As a consequence of the proposition
${(W(x))}_{x\in X}$ is a good system of open neighborhoods.
This prompts the following definition.
\begin{definition} \index{D- Refinement of a good system}\index{D- Subordinate good system}
Let $X$ be an ep-groupoid and $\mathfrak{U}={((U(x))}_{x\in X}$ a good system of open neighborhoods. We say that another good system
$\mathfrak{W}={(W(x))}_{x\in X}$ is a {\bf refinement} provided $W(x)\subset U(x)$ for all $x\in X$.
Given an arbitrary open covering of $X$ be open sets we say that a good system $\mathfrak{U}$  is {\bf subordinate}
provided every $U(x)$ is contained in a set of the open covering.
\qed
\end{definition}
In the following we shall explore the properties of a good system of open neighborhoods $\mathfrak{U}$.
We begin with the study of open invariant neighborhoods $U(x)$ and $U(x')$ equipped with the natural $G_x$ and $G_{x'}$-actions, respectively.  We denote by $\bm{U}(x,x')$ the collection of morphisms in $\bm{X}$ starting in $U(x)$ and ending in $U(x')$.
We define the 
sc-smooth map $G_{x'}\times \bm{U}(x,x')\times G_x\rightarrow \bm{U}(x,x')$  by 
$$
g'\ast\phi \ast g:= \Gamma'(g', t(\phi))\circ\phi\circ \Gamma(g,g^{-1}\ast s(\phi)).
$$
By definition,  $t(g'\ast\phi\ast g) = g'\ast t(\phi)$ and $s(g'\ast\phi\ast g) = g^{-1}\ast s(\phi)$,  so that
$$
g'\ast\phi\ast g\colon  g^{-1}\ast s(\phi)\rightarrow g'\ast t(\phi).
$$
Moreover, if  $h,g\in G_x$ and $g',h'\in G_{x'}$, then 
\begin{equation*}
\begin{split}
&\phantom{=.}(g'h')\ast\phi\ast (gh)\\
&= \Gamma'(g'h',t(\phi))\circ \phi\circ \Gamma(gh,(h^{-1}g^{-1})\ast s(\phi))\\
&= \Gamma'(g',h'\ast t(\phi))\circ\Gamma'(h',t(\phi))\circ\phi\circ  \Gamma(g,g^{-1}\ast s(\phi))\circ \Gamma(h,(h^{-1}g^{-1})\ast s(\phi))\\
&= \Gamma'(g',h'\ast t(\phi))\circ (h'\ast\phi \ast g)\circ \Gamma(h,h^{-1}\ast (g^{-1}\ast s(\phi)))\\
&=\Gamma'(g',t(h'\ast\phi \ast g))\circ (h'\ast\phi \ast g)\circ \Gamma(h, h^{-1}\ast s(h'\ast\phi \ast g))\\
&= g'\ast (h'\ast\phi\ast  g)\ast h.
\end{split}
\end{equation*}
This together with   $1_{x'}\ast \phi\ast 1_x=\phi$ implies that 
$$
G_{x'}\times \bm{U}(x,x')\times G_x\rightarrow \bm{U}(x,x')
$$
defines a right-action of $G_x$ and  a left-action of $G_{x'}$ on $\bm{U}(x,x')$, which are commuting. 
The source map $s:\bm{U}(x,x')\rightarrow U(x)$ is $G_x$ equivariant since
$$
s(\phi\ast g) = s(\phi\circ \Gamma(g,g^{-1}\ast(\phi))) = g^{-1}\ast s(\phi)
$$
when we let $G_x$ act on the right of $U(x)$. Similarly $t:\bm{U}(x,x')\rightarrow U(x')$ satisfies
$$
t(g'\ast \phi) = \Gamma'(g',t(\phi))= g'\ast t(\phi),
$$
where $G_{x'}$ acts from the left, i.e. the natural action of $G_{x'}$.
\begin{definition}\index{D- Natural bi-action on $\bm{U}(x,x')$}
Given $x,x'\in X$ and open neighborhoods $U(x)$ and $U(x')$ with the natural actions by $G_x$ and $G_{x'}$, respectively,
we shall call 
$$
G_{x'}\times \bm{U}(x,x')\times G_x\rightarrow \bm{U}(x,x')
$$
the {\bf natural bi-action}.
\qed
\end{definition}

Before we give a basic existence result for systems of good neighborhoods we shall derive some of their properties with respect to our actions. 
\begin{lemma}
Let  ${(U(x))}_{x\in X}$ be  a good system of open neighborhoods and $\Sigma$ a connected component of $\bm{U}(x,x')$. Then the following assertions hold true.
\begin{itemize}
\item[{\em(1)}]\ If $g'\in G_{x'}$ and $g'\ast\Sigma=\Sigma$,  then $g'=1_{x'}$.
\item[{\em(2)}]\ If $g\in G_x$ and $\Sigma\ast g=\Sigma$,  then $g=1_x$.
\item[{\em(3)}]\ If \  $\Sigma_1$ and $\Sigma_2$ are connected components of $\bm{U}(x, x')$,  then either 
$s(\Sigma_1)=s(\Sigma_2)$ or $s(\Sigma_1)\cap s(\Sigma_2)=\emptyset$.
 The same assertion holds for the target map $t$.
\end{itemize}
\end{lemma}
\begin{proof}
(1) The target map $t\colon \Sigma\rightarrow U(x')$ is a sc-diffeomorphism onto an open subset $O'$.  From $g'\ast\Sigma=\Sigma$ it follows that 
$g'\ast O'=O'$. For  $\phi \in \Sigma$ consider $g'\ast \phi\in \Sigma$ and note that  the morphisms $\phi$ and $g'\ast \phi$ have the same source. Since 
the source map  $s\colon \Sigma\rightarrow U(x)$ is a sc-diffeomorphism onto an open subset of $U(x)$,   we deduce
 that $g'\ast\phi=\phi$. This means  that $\Gamma'(g',t(\phi))\circ \phi=\phi$ and therefore 
$\Gamma'(g',t(\phi))=1_{t(\phi)}$. This implies  that $g'=1_{x'}$, as claimed.  Part (2) is proved similarly.

(3) Assuming that  $s(\Sigma_1)\cap s(\Sigma_2)\neq \emptyset$,  we find morphisms $\phi_1\in \Sigma_1$ and $\phi_2\in\Sigma_2$ 
satisfying $s(\phi_1)=s(\phi_2)$, and a unique $g'\in G_{x'}$ for which 
$$
\phi_2 = \Gamma'(g',t(\phi_1))\circ \phi_1.
$$
Then, considering the sc-smooth map 
$$
\Sigma_1\rightarrow \bm{U}(x,x'), \quad  \phi\mapsto  \Gamma'(g',t(\phi))\circ \phi, 
$$
we deduce that $g'\ast \Sigma_1=\Sigma_2$ since both $\Sigma_1$ and $\Sigma_2$ are connected components of $\bm{U}(x, x')$.
Since $s(g'\ast \Sigma_1)=s(\Sigma_1)$, we see that $s(\Sigma_1)=s(\Sigma_2)$, which proves our claim.  
A similar argument applies to the target map $t$.
\qed \end{proof}

\begin{theorem}[{\bf Extension of good partial systems}]\index{T- Extension of a good partial system}\label{THMX1317}
Let $Y$ be an ep-groupoid and $Q$ an open subset of the object space so that $|Q|=|Y|$.
Equip the   full subcategory associated to $Q$ with the induced structure of an ep-groupoid. 
Assume we are given a good system ${(U(q))}_{q\in Q}$ of open neighborhoods for the ep-groupoid $Q$.
Then given for  every $y\in Y\setminus Q$  an open neighborhood $W(y)\subset Y$, we can find 
open neighborhoods $U(y)$ for $y\in Y\setminus Q$ so that $U(y)\subset W(y)$ and $\mathfrak{U}$
$$
\mathfrak{U}={(U(y))}_{y\in Y}
$$
is a good system of open neighborhoods for $Y$.  
\end{theorem}
\begin{proof}
For $y\in Y\setminus Q$ pick $\psi_y\in \bm{Y}$ with $q_y:= s(\psi_y)\in Q$ and $t(\psi_y)=y$. We find open neighborhoods
$O(\psi_y)$, $O(q_y)$, and $U(y)$ such that the following holds.
\begin{itemize}
\item $s:O(\psi_y)\rightarrow O(q_y)$ and $t:O(\psi_y)\rightarrow U(y)$ are sc-diffeomorphisms.
\item $O(q_y)\subset U(q_y)$ and $O(q_y)$ is invariant under the action on $U(q_y)$.
\item $U(y)$ admits the natural $G_y$-action and has the properness property and is contained in $W(y)$.
\end{itemize}
We note that we can pick these three neighborhoods arbitrarily small. Using the properness property 
of $U(q_y)$, which contains $O(q_y)$,  we can pick $U(y)$ so small that the following holds in addition.
\begin{itemize}
\item For every connected component $\Sigma$ of $\bm{U}(q_y,y)=\{\phi\in \bm{Y}\ |\ s(\phi)\in U(q_y),\ t(\phi)\in U(y)\}$ the maps
$s:\Sigma\rightarrow U(q_y)$ and $t:\Sigma\rightarrow U(y)$ are sc-diffeomorphisms onto their image.
\end{itemize}

At this point we have defined a collection $\mathfrak{U}={(U(y))}_{y\in Y}$ of open neighborhoods, where for $q\in Q$
the open set $U(q)$ is the original given one.  Then $\mathfrak{U}$
clearly satisfies Properties (1) and (2) of Definition \ref{DEFF1}. Next we consider the  Property (3).
If $q_1,q_2\in Q$ we already know that $\bm{U}(q_1,q_2)$ satisfies (3). 
We need to investigate the following three cases for Property (3) which involve new sets.
\begin{itemize}
\item[(1)]\  $\bm{U}(q,y)$ for $q\in Q$ and $y\in Y\setminus Q$.
\item[(2)]\ $\bm{U}(y,q)$ for  $y\in Y\setminus Q$ and $q\in Q$.
\item[(3)]\ $\bm{U}(y,z)$ for $y,z\in Y\setminus Q$.
\end{itemize}
\noindent{\bf Case (1)}  We define an sc-smooth  map $A: \bm{U}(q,y)\rightarrow \bm{U}(q,q_y)$
by 
$$
A(\phi) =[(t|O(\psi_y))^{-1}(t(\phi))]^{-1}\circ\phi.
$$
We note that $A$ is a local sc-diffeomorphism. 
 We  verify
\begin{eqnarray}\label{EQN13X}
s(A(\phi))=   s(\phi)
\end{eqnarray}
and
\begin{eqnarray}\label{EQN13XX}
t(A(\phi)) = t ([(t|O(\psi_y))^{-1}(t(\phi))]^{-1}\circ\phi) =  \wh{\psi}_y^{-1}(t(\phi)),
\end{eqnarray}
where $\wh{\psi}_y:O(q_y)\rightarrow O(y)$ is the sc-diffeomorphism associated to $\psi_y$.
Now we show that $A$  is injective, and hence an sc-diffeomorphism onto its (open) image.
To see the injectivity assume that $A(\phi_1)=A(\phi_2)$. Then we derive from  (\ref{EQN13XX})
that  $t(\phi_1)=t(\phi_2)$ which implies, using the definition of $A$, that $\phi_1=\phi_2$. 
We rewrite (\ref{EQN13XX}) as 
\begin{eqnarray}\label{EQN13Y}
t(\phi)= \wh{\psi}_y\circ t(A(\phi)).
\end{eqnarray}
If $\Sigma\subset \bm{U}(q,y)$ is a connected component, then $A(\Sigma)$ lies in a connected component
$\Sigma'$ of $\bm{U}(q,q_y)$. Since $s|\Sigma'$ and $t|\Sigma'$ are sc-diffeomorphisms onto their images
it follows that $s|\Sigma$ as well as $t|\Sigma$ have the same property.\par

\noindent{\bf Case (2)}  The map $B: \bm{U}(y,q)\rightarrow \bm{U}(q,y):\phi\rightarrow \phi^{-1}$ is an sc-diffeomorphism. Using that $t\circ B =s$ and $s\circ B=t$ the assertion follows from {\bf Case (1)}.\par

\noindent{\bf Case (3)}  The map 
$C:  \bm{U}(y,z)\rightarrow \{\phi\in\bm{Y}\ |\ s(\phi)\in O(q_{y}),\ t(\phi)\in O(q_{z})\}
$
defined by 
$$
C(\phi) =\left ((t|O(\psi_z))^{-1}(t(\phi))\right)^{-1}\circ \phi\circ  ((t|O(\psi_y))^{-1}(s(\phi)))
$$
is an sc-diffeomorphism and 
$$
t(C(\phi))= \wh{\psi}_z^{-1}(t(\phi))\ \ \text{and}\ \ s(C(\phi)) = \wh{\psi}_y^{-1}(s(\phi)).
$$
Moreover $\{\phi\in\bm{Y}\ |\ s(\phi)\in O(q_{y}),\ t(\phi)\in O(q_{z})\}\subset \bm{U}(q_y,q_z)$.
If $\Sigma$ is a connected component in  $\bm{U}(y,z)$ ist image under $C$ lies in a connected component
$\Sigma'$ of $\bm{U}(q_y,q_z)$ on which $t$ and $s$ have the desired properties.
This immediately with the previous discussion implies the desired result.
\qed \end{proof}

The basic result is the following theorem which is concerned with an ep-groupoid  $X$ having a paracompact orbit space.
The paracompactness of $|X|$ is used  in a crucial way. The theorem  shows 
that there exists a system of good open neighborhoods $\mathfrak{U}$ with the $U(x)$ being arbitrarily small.
\begin{theorem}[{\bf Existence of good systems}]\label{xxxx-structure}\index{T- Existence of a good system}
Let $X$ be an ep-groupoid with paracompact  orbit space  and let ${(W(x))}_{x\in X}$ be  collection of open neighborhoods around the points in $X$.
Then there exists a good system of open neighborhoods ${(U(x))}_{x\in X}$ satisfying $U(x)\subset W(x)$ for all $\in X$.
\end{theorem}

\begin{proof} In view of the extension theorem Theorem \ref{THMX1317} it suffices to verify the following assertion.\par

\noindent {\bf Assertion:} There exists an open subset $V$ of $X$ having the property $|V|=|X|$, so that the ep-groupoid 
associated to $V$ admits a good system of open neighborhoods ${(U(x))}_{x\in V}$ with $U(x)\subset V$.
The paracompactness of $|X|$ is important for this part of the construction.\par

By taking suitable smaller subsets for $W(x)$ we may assume without loss of generality that the following holds true.
\begin{description}
\item[(1)]  For every $x\in X$ the map $t:s^{-1}(cl_X(W(x)))\rightarrow X$ is proper.
\item[(2)]  $W(x)$ is equipped with the natural $G_x$-action. 
\end{description}
These will be our standing assumptions for the following constructions.\par

 We start with the open covering ${(|W(x)|)}_{x\in X}$ of $|X|$.
Using the paracompactness of $|X|$ we know that $|X|$ is metrizable by Theorem \ref{ATHOME}.
We can find open subsets $\wt{V}_x$ of $|W(x)|$ so that ${(\wt{V}_x)}_{x\in X}$ is a locally finite open covering of
$|X|$ and $\cl_{|X|}(\wt{V}_x)\subset |W(x)|$. Of course, many of these sets might   be empty. Let $A\subset X$ be the collection of all $x\in X$ with $\wt{V}_x\neq \emptyset$.
 Define for $a\in A$
\begin{eqnarray}
V_a =\pi^{-1}(\wt{V}_a) \cap W(a).
\end{eqnarray}
Note that we write $V_a$ instead of $V(a)$, to indicate that $V_a$ need not to contain the point $a$.
The following holds by construction.
\begin{description}
\item[(3)] $V_a\subset W(a)$ and $V_a$ is invariant under $G_a$ for $a\in A$.
\item[(4)] $\cl_{|X|}(|V_a|)\subset |W(a)|$ for $a\in A$.
\item[(5)] ${(|V_a|)}_{a\in A}$ is a locally finite covering of $|X|$.
\end{description}
It follows that $\cl_X(V_a)\subset W(a)$ for $a\in A$. Consider the open subset $V$ of $X$ defined by
$$
V=\bigcup_{a\in A} V_a.
$$
Clearly $|V|=|X|$.  If $x\in V$ there exist a $G_x$-invariant open neighborhood $U_0(x)$ so that $U_0(x)$ intersects only a finite number
of $\cl_X(V_a)$ nontrivially. The set of such $a$  can be written as
$$
\{a^x_1,...,a^x_{\ell_x}\}.
$$
We may also assume the following after perhaps a renumbering of the $a^x_i$.
\begin{description}
\item[(6)] $\cl_X(U_0(x))\subset V_{a^x_1}$ and $U_0(x)\subset W(x)$.
\item[(7)]  $\cl_X(U_0(x))\subset W(a^x_i)$ for $i\in \{1,...,\ell_x\}$.
\end{description}
For $x\in V$ we consider the set of all morphisms starting in $\cl_X(V)$ and ending in $x$. Since a small open neighborhood
of $|x|$ only intersects a finite number of the $|V_a|$, and using the properness of the associated $t:s^{-1}(\cl_X(W(a)))\rightarrow X$
we deduce that this  set is finite and denote it by $\Theta_x$
$$
\Theta_x =\{\phi\in \bm{X}\ |\ s(\phi)\in \cl_X(V),\ t(\phi)=x\}.  
$$
Given an open neighborhood  ${\mathcal O}$ of $\Theta_x$ we find an open $G_x$-invariant  neighborhood $U_1(x)\subset U_0(x)$
such that $\{\phi\in \bm{X}\ |\ s(\phi)\in \cl_X(V),\ t(\phi)\in U_1(x)\}$ is contained in ${\mathcal O}$.  
This implies, taking ${\mathcal O}$ sufficiently small, and $U_1(x)$ appropriately, that for a connected component $\Sigma$
of $\{\phi\in \bm{X}\ |\ s(\phi)\in V,\ t(\phi)\in U_1(x)\}$ the maps $s|\Sigma$ and $t|\Sigma$ are 
sc-diffeomorphisms onto open subsets of $V$ and $U_1(x)$, respectively.  Define the collection ${(U(x))}_{x\in V}$ by
$U(x)=U_1(x)$. 

At this point we have constructed open neighborhoods $U(x)$ around the points $x\in V$ so that the following holds.
\begin{itemize}
\item[(i)]\ \ \ Each $U(x)$ admits the natural $G_x$-action.
\item[(ii)]\ \ \  $t:s^{-1}(\cl_X(U(x)))\rightarrow X$ is proper.
\item[(iii)]\ \ \ For $x,x'\in V$ the M-polyfold $\bm{U}(x,x')$ has the property that for a connected component $\Sigma$ the maps $s|\Sigma$ and $t|\Sigma$
are sc-diffeomorphisms onto open subsets of $U(x)$ and $U(x')$, respectively.
\end{itemize}
 Hence, we satisfy the definition of a good
system of open neighborhoods but only on $V\subset X$. However, it is important to recall that $|V|=|X|$.
Now employing Theorem \ref{THMX1317} this family of sets can be extended to $X$ and the proof is complete.
\qed \end{proof}

We discuss some of the consequences of having a good system of neighborhoods ${(U(x))}_{x\in X}$ for the ep-groupoid $X$ with paracompact orbit space $|X|$.
Assume that $\Sigma$ is a connected component in ${\bf U}(x,x')$.
Define $O=s(\Sigma)$ and $O'=t(\Sigma)$ and denote by $G_{x,\Sigma}$ the subgroup of $G_x$ consisting of the elements $g\in G_x$ satisfying $ g\ast O =O$.
Similarly we define $G_{x',\Sigma}$ as the subgroup of $G_{x'}$ consisting of all $g'$ with $g'\ast O'=O'$. 
Define the sc-diffeomorphism $D_\Sigma:O\rightarrow O'$ by
$$
D_\Sigma = t\circ (s|\Sigma)^{-1}.
$$
For $g\in G_{x,\Sigma}$ and $z\in O$ we consider $D_\Sigma(g\ast z)$ and $D_\Sigma(z)$.  We pick $\psi, \phi \in \Sigma$
with $s(\psi)=z$ and  $s(\phi)=g\ast z$, i.e. $\psi=(s|\Sigma)^{-1}(z)$ and $\phi=(s|\Sigma)^{-1}(g\ast z)$. Then 
$$
t(\psi)= t\circ (s|\Sigma)^{-1}(z)=D_\Sigma(z)\ \text{and}\ t(\phi)= t\circ (s|\Sigma)^{-1}(g\ast z)= D_\Sigma(g\ast z).
$$
For every $z\in O$ we define the morphism
$$
\Theta_g(z):=[(s|\Sigma)^{-1}(g\ast z)] \circ \Gamma(g,z)\circ [(s|\Sigma)^{-1}(z)]^{-1}: D_\Sigma(z)\rightarrow D_\Sigma(g\ast z),
$$
 which has source and target in $O'$.
For $z\in O$ there exists  a unique $g'_z\in G_{x'}$ that that
$$
\Gamma'(g_z',D_\Sigma(z))= \Theta_g(z).
$$
Since the right-hand side changes smoothly in $z$ it follows that $g_z'$ is independent of $z\in O$. It will be denoted by $g'(g)$ and
it is also clear that $g'(g)\in G_{x',\Sigma}$. Hence we obtain for every $z\in O$ and $g\in G_{x,\Sigma}$  the identity
$$
\Gamma'(g'(g),D_\Sigma(z))=\Theta_g(z).
$$
Applying $t$ we arrive at
$$
g'(g)\ast D_\Sigma(z)= D_\Sigma(g\ast z).
$$
Using the uniqueness of the map $G_{x,\Sigma}\rightarrow G_{x',\Sigma}:g\rightarrow g'(g)$ it follows from 
\begin{eqnarray*}
&& (g'(g)g'(h))\ast D_\Sigma(z)\\
&=&g'(g)\ast(g'(h)\ast D_\Sigma(z))\\
&=& g'(g)\ast (D_\Sigma(h\ast z))\\
&=&D_\Sigma(g\ast (h\ast z))\\
&=& D_\Sigma((gh)\ast z)\\
&=& g'(gh)\ast D_\Sigma(z)
\end{eqnarray*}
that $g'(gh)=g'(g)g'(h)$.
For later reference we note the following proposition.
\begin{proposition}
Let ${(U(x))}_{x\in X}$ a good system of open neighborhoods and  $\Sigma\subset {\bf U}(x,x')$ be a connected component.
Then there exists a uniquely determined group isomorphism $G_{x,\Sigma}\rightarrow G_{x',\Sigma}:g\rightarrow g'(g)$ such that for all $z\in s(\Sigma)$
$$
(s|\Sigma)^{-1}(g\ast z)= g'(g)\ast (s|\Sigma)^{-1}(z)\ast g^{-1}.
$$
\end{proposition}
\begin{proof}
By definition, with $\psi_z = (s|\Sigma)^{-1}(z)$ for $z\in O$ it holds that $D_\Sigma(z)=t(\psi_z)$.  Recall that
$$
\psi_{g\ast z} \circ \Gamma(g,z)\circ \psi_z^{-1}=\Gamma'(g'(g),t(\psi_z)).
$$
Hence
$$
\psi_{g\ast z}\circ \Gamma(g,z) = \Gamma'(g'(g),t(\psi_z))\circ \psi_z = g'(g)\ast \psi_z.
$$
Moreover, we compute 
$$
\psi_{g\ast z}\ast g = \psi_{g\ast z} \circ \Gamma(g,g^{-1}\ast s(\psi_{g\ast z})) = \psi_{g \ast z}\circ \Gamma(g,z).
$$
Putting these facts together we arrive at
$$
\psi_{g\ast z}= g'(g)\ast \psi_z \ast g^{-1}.
$$
\qed \end{proof}

\section{General \texorpdfstring{$\bssc^+$}{bss}-Multisections}\label{multioo}

We shall introduce the notion of a $\ssc^+$-multisection.
Multisections are a particular case of  set-valued sections. Set-valued operators have been used in nonlinear functional analysis for a long time and we refer the reader to \cite{AE}. Our definition is along the lines of the definition introduced  in \cite{CRS} and  related to \cite{FO}. The following definition from \cite{HWZ3.5} brings a definition in  \cite{CRS} into the groupoid framework.

We view the non-negative rational numbers ${\mathbb Q}^+=\Q\cap [0,\infty)$ as a category having only the identities as morphisms.

\begin{definition}[{\bf $\bssc^+$-multisection}] \label{sc+-section-functor}\index{D- $\bssc^+$-multisection}
Let $(P:W\rightarrow X,\mu)$ be a strong bundle over the ep-groupoid $X$. Then an {\bf $\ssc^+$-multisection} is a 
functor 
$$
\Lambda:W\rightarrow {\mathbb Q}^+
$$
 such that the following local representation 
(called {\bf local section structure})\index{Local section structure} holds true.  Around  every object $x_0\in X$,  there exists an open neighborhood $U(x_0)$ on which the isotropy group $G_{x_0}$ acts by its natural representation, and  finitely many $\ssc^+$-sections
$s_1,\ldots ,s_k :U(x_0) \to W$ (called {\bf local sections}) with associated  positive rational numbers $\sigma_1,\ldots ,\sigma_k$ (called {\bf weights})\index{Weights} satisfying the following properties:
\begin{itemize}
\item[(1)]\ $\sum_{i=1}^k \sigma_i =1.$
\item[(2)]\ $\Lambda(w)=\sum_{i\in \{1,\ldots k\  |\ w=s_i(P(w))\}} \sigma_i$ 
\end{itemize}
for all $w\in W|U(x_0)$, where the empty sum has by definition the value $0$.
We shall refer to the right hand side of the identity (2) 
as a {\bf local representation}\index{Local representation of $\Lambda$} near $x_0$.
\qed
\end{definition}

The functoriality of $\Lambda$ implies that 
$\Lambda (w')=\Lambda (w)$, if there exists a morphism $w\to w'$ in $\bm{W}$. Explicitly,
$$
\Lambda (\mu (\phi, w))=\Lambda (w)
$$ 
for all $(\phi, w)\in \bm{X}_{_{s}\times_P}W$. Hence $\Lambda$ induces a map $\abs{\Lambda}: \abs{W}\to \Q^+$ on the orbit space of $W$.
\begin{definition}\label{DEFNX1322}
\begin{itemize}
\item[(1)]\ The  {\bf domain support}\index{D- Domain support of $\Lambda$} of $\Lambda$ \index{$\text{dom-supp}(\Lambda)$} is the subset of $X$, defined by 
$$
\text{dom-supp}(\Lambda)=
\text{cl}_X\left(\{x\in X\, \vert \, \ \exists
\text{ $w\in W_x\setminus\{0\}$ for which  $\Lambda(w)>0$}\}\right).
$$
\item[(2)]\ The {\bf support of $\bm{\Lambda}$} \index{D- Support of $\Lambda$} is the subset $\supp(\Lambda)$\index{$\supp(\Lambda)$} of $W$ defined by
$$
\supp(\Lambda)=\{w\in W\, \vert \, \Lambda(w)>0\}.
$$
\qed
\end{itemize}
\end{definition}
For every $x\in X$, the set 
$$
\supp(\Lambda)(x)=\{w\in W\,\vert \, \text{$P(w)=x$ and $\Lambda (w)>0$}\}\index{$\supp(\Lambda)(x)$}
 $$
  is finite and 
$$\sum_{w\in \supp(\Lambda)(x)}\Lambda (w)=1.$$
Moreover, if $x\in U(x_0)$, then 
$$
\supp(\Lambda)(x)=\{s_1(x),\ldots ,s_k(x)\}.
$$
\begin{definition}\label{DEF1223}\index{D- Pointwise norm $N(\Lambda)$}
If $N$ is an auxiliary norm for the strong bundle $(P,\mu)$ we define for $x\in X$, 
$$
N(\Lambda)(x)=\max\{N(w)\, \vert \, \text{$w\in\supp(\Lambda)$ satisfying  $P(w)=x$}\}.\index{$N(\Lambda)(x)$}\index{$N(\Lambda)(x)$}.
$$
We shall call $N(\Lambda)$ the {\bf pointwise norm} of the sc$^+$-multisection with respect to the auxiliary norm $N$.
\qed
\end{definition} 
Denoting by  $0_W$  the zero-section of $W$,  we obtain the relations
$$
\text{dom-supp}(\Lambda)=\cl_X(P(\supp(\Lambda)\cap (W\setminus 0_W))) =\cl_X(\{x\in X\ |\ N(\Lambda)(x)>0\}).
$$
\begin{definition}
The $\ssc^+$-multisection is called {\bf trivial  
on the set $V\subset X$} \index{D- $\Lambda$ trivial over $V$} if $\Lambda$ is identically equal to $1$ on the zero-section  over $V$, i.e., $\Lambda (0_x)=1$ for all $x\in V$ (and hence $\Lambda (w_x)=0$ for all $w_x\neq 0_x$).
\qed
\end{definition}
We note that the domain-support of $\Lambda$ is the smallest closed set in $X$ outside of which $\Lambda$ is trivial.
  In general the local representation is not unique. Among the local representations 
there are certain preferable ones. This is discussed next.

Let $\Lambda:W\rightarrow {\mathbb Q}^+$ be an sc$^+$-multisection,  assume $U(x)$ is an open $G_x$-invariant neighborhood
(natural action) so that the properness property holds, 
 and ${(s_j)}_{i\in J}$ is  the associated local section structure with weights $(\sigma_j)$.

We take  a suitable rational number $\sigma=q/p\in (0,1]$ and, for every $j\in J$, a positive integer $k_j$ such that 
$k_j (q/p)=\sigma_j$.  We introduce  a new index set
$\wh{J}$ consisting of all pairs $(j,i)\in J\times \N$ where $1\leq i\leq q k_j$, and new sections $(s_{(j,i)})_{(j, i)\in \wh{J}}$ on $U(x)$ by
$$
s_{(j,i)}=s_j.
$$
The cardinality of the index set $\wh{J}$ is equal to  $p$, 
$$
|\wh{J}|= \sum_{j\in J} q\cdot k_j =\sum_{j\in J} p\cdot \sigma_j = p.
$$
We have replaced the original section $s_j$ having  weight $\sigma_j$ by $(k_j\cdot q)$-many sections $s_{(j, i)}$, 
each having weight $\frac{\sigma_j}{k_j\cdot q} =\frac{1}{p}=\frac{1}{\wh{J}}$. Hence we have constructed a local $\ssc^+$-section structure over $U(x)$
of $p$ many sections having equal weights $1/p$.

Consequently, we may assume that we are given an $\ssc^+$-section structure $(s_i)_{i\in I}$ over $U(x)$ each having weight $1/\abs{I}$. Then
$$
\Lambda(w)=\frac{1}{\abs{I} }\abs{   \{i\in I \, \vert \, s_i(P(w))=w\}    }
$$
for $ w\in W$ satisfying $P(w)\in U(x)$.

 Abbreviating $\wh{I}=I\times G_x$, 
we define, for $(i,g)\in \wh{I}$,  a new $\ssc^+$-section $s_{(i,g)}$  on $U(x)$ by
$$
s_{(i,g)}(g\ast y) =\mu(\Gamma(g,y),s_i(y)).
$$
We have defined a section structure $(s_{(i,g)})_{(i,g)\in\wh{I}}$ over $U(x)$, where each section has the weight $1/\abs{\wh{I}}$, and such that 
$$
\Lambda(w) =\frac{1}{\abs{\wh{I}}} \abs{\{(i, g)\in \wh{I}\, \vert  \,   s_{(i, g)}(P(w))=w\}}
$$
for $w\in W$ satisfying $P(w)\in U(x)$.  Defining the action of $G_x$ on $\wh{I}$ by $g(i, h):=g\ast (i, h)=(i, gh)$ for $(i, h)\in \wh{I}$,  we compute

\begin{equation*}
\begin{split}
s_{g(i,h)}( g\ast y) &=s_{(i,gh)}(g\ast y)\\
&=\mu(\Gamma(gh,h^{-1}\ast y),s_i(h^{-1}\ast y))\\
&=\mu(\Gamma(g, y),\mu(\Gamma(h,h^{-1}\ast y),s_i(h^{-1}\ast y)))\\
&=\mu(\Gamma(g, y),s_{(i,h)}( y))
\end{split}
\end{equation*}

In  other words we can always find a $\ssc^+$-section structure ${(s_i)}_{i\in I}$, where every section carries the same weight $1/|I|$,
the isotropy group
$G_x$ acts on the index set $G_x\times I\rightarrow I$, $(g,i)\rightarrow g(i)$,  and the family of sections has an invariance property with respect to the action
$$
s_{g(i)}(g\ast y)=\mu(\Gamma(g,y),s_i(y)).
$$
This prompts the following definition. 

\begin{definition}\label{XDEF1}
A {\bf symmetric $\ssc^+$-section structure}\index{D- Symmetric sc$^+$-section structure} for 
$$
\Lambda: W\to \Q^+,
$$
a $\ssc^+$-multisection,   on $U(x)$ is a pair $(U(x), (s_i))_{i\in I_x}$  consisting of the 
$G_x$-invariant open neighborhood $U(x)$ such that the target map $t: s^{-1}(\cl_X(U(x)))\rightarrow X$ is proper, a finite family $(s_i)_{i\in I_x}$ of $\ssc^+$-sections  defined on $U(x)$,  and an action of $G_x$ on the index set $I_x$ such that the following holds:
\begin{itemize}
\item[(1)]\ $s_{g(i)}(g\ast y) = \mu(\Gamma(g,y),s_i(y))$ for all $y\in U(x)$, $i\in I_x$, and $g\in G_x$.
\item[(2)]\ $ \Lambda(w) =\frac{1}{|I_x|}|\{i\in I_x\ |\ s_i(P(w))=w\}|$ for all $w\in W$ satisfying $P(w)\in U(x)$.
\end{itemize}
\qed
\end{definition}
It is important to notice that we allow for different $i,i'\in I_x$ that $s_i=s_{i'}$.
From the previous discussion we obtain the following result.
\begin{lemma}
 Every sc$^+$-multisection functor locally admits a symmetric sc$^+$-sec\-tion structure.
 \qed
\end{lemma}
As a first new object we consider tuples $(\Lambda,\mathfrak{U},\mathfrak{S})$, where $\mathfrak{U}$ is a a good system of
open neighborhoods, and $\mathfrak{S}$ a family of sc$^+$-section structures.
\begin{definition}\index{D- System of symmetric sc$^+$-section structures}
A sc$^+$-multisection with a {\bf system of symmetric  section structures} is a tuple $(\Lambda,\mathfrak{U},\mathfrak{S})$,
where $\Lambda:W\rightarrow {\mathbb Q}^+$ is an sc$^+$-multisection, $\mathfrak{U}={(U(x))}_{x\in X}$ is 
a good system of open neighborhoods, and $\mathfrak{S}$ is a family $\left( \mathfrak{s}_x\right)_{x\in X}$,
where for every $x\in X$  there is a given a finite set $I_x$ and a symmetric sc$^+$-section structure $\mathfrak{s}_x=(s_i^x)$ parametrized by $I_x$ and defined on $U(x)$
so that 
$$
\Lambda(w)=\frac{1}{|I_x|}\cdot \sharp \left( \{i\in I_x\ |\ s_i^x(P(w))=w\}\right) .
$$
\qed
\end{definition}
In general the precise parametrization of the $(s_i^x)$ does not matter, and moreover only the germ near $x$ matters.
This leads us to define a notion of equivalence  between sc$^+$-multisections with a system of symmetric section structures.
\begin{definition}\index{D- Decorated sc$^+$-multisection}
Let $(\Lambda,\mathfrak{U},\mathfrak{S})$
and $(\Lambda',\mathfrak{U}',\mathfrak{S}')$ be sc$^+$-multisections with a system of symmetric section structures.
We say that both are equivalent provided $\Lambda=\Lambda'$, and there exists a good system of open neighborhoods
$\mathfrak{U}''$ and for every $x\in X$ an equivariant bijection $b_x:I_x\rightarrow I_{x}'$ such that the following holds.
\begin{itemize}
\item[(1)]\ $U''(x)\subset U(x)\cap U'(x)$ for all $x\in X$.
\item[(2)]\ $s'_{b_x(i)}(y) =  s_i(y)$ for $i\in I_x$ and $y\in U''(x)$.
\end{itemize}
An equivalence class will be denoted by $[\Lambda,\mathfrak{U},\mathfrak{S}]$ and will be called a {\bf decorated
sc$^+$-multisection}. By abuse of notation we shall refer to $\Lambda$ as a decorated sc$^+$-multisection equating
$\Lambda\equiv [\Lambda,\mathfrak{U},\mathfrak{S}]$ if the context is clear.
\qed
\end{definition}

\begin{definition}\index{D- Sum of sc$^+$-multisections}
Let $(P: W\rightarrow X,\mu)$ be a strong bundle over the ep-grouoid $X$, and $\Lambda'$, $\Lambda''$ two decorated $\ssc^+$-multisections.
Then the  sum $\Lambda''=\Lambda'\oplus \Lambda''$\index{$\Lambda'\oplus \Lambda''$} defined as follows
 is a naturally  decorated
sc$^+$-multisection $[\Lambda'',\mathfrak{U}'',\mathfrak{S}'']$, where 
\begin{itemize}
\item[(1)]\  $(\Lambda'\oplus\Lambda'')(w) =\sum_{\{(w',w'')\ |\ w',w''\in W|P(w),\ w'+w''=w\}} \Lambda'(w')\cdot\Lambda''(w'').$
\item[(2)]\  $\mathfrak{U}''$ by $U''(x)=U(x)\cap U'(x)$ for $x\in X$.
\item[(3)]\ $I_x'' =I_x\times I_{x}'$ and $s^{''x}_{(i,i')} := s_i^x + s_{i'}^{'x}$ on $U''(x)$.
\end{itemize}
\qed
\end{definition}
The sum operation is associative and commutative. However, in general, there are no inverses. 
\begin{definition}\index{D- Product $a\odot\Lambda$}\label{D- Scalar product}
Given $(P:W\rightarrow X,\mu)$ and a decorated $\ssc^+$-multisection $\Lambda$ we define for a real number $a$ the {\bf scalar product}
$a\odot\Lambda$\index{$a\odot\Lambda$}, which is another  decorated $\ssc^+$-multisection,  as follows. If $a=0$, then 
$(0\odot \Lambda)(w)=0$ if $w\neq 0$ and $(0\odot \Lambda)(0)=1$. Moreover,  for $a\neq 0$, 
$$
(a\odot \Lambda)(w):=\Lambda\left(\frac{1}{a}w\right).
$$
\qed
\end{definition}
We note that $a\odot\Lambda$ is decorated as well. Namely taking $\mathfrak{U}$ and $\mathfrak{S}={(\mathfrak{s}_x)}_{x\in X}$
from $[\Lambda,\mathfrak{U},\mathfrak{S}]$ the new $\mathfrak{S}'$ is given by the families
$$
\mathfrak{s}_x' ={(a\cdot s_i^x)}_{i\in I_x}.
$$
\begin{lemma}\label{lemma_scalar_prod}\index{L- Scalar product $\beta\odot\Lambda$}
Let $(P:W\rightarrow X,\mu)$ be a strong bundle over an ep-groupoid and $\Lambda$ a decorated  $\ssc^+$-multisection.
We assume that $\beta\colon X\rightarrow {\mathbb R}$ is a sc-smooth functor, i.e., $\beta (x)=\beta (x')$ if there is a morphism $\phi: x\to x'$.  Then,  $\beta\odot\Lambda$\index{$\beta\odot\Lambda$},  defined by
$$
(\beta\odot\Lambda)(w)=(\beta(P(w))\odot \Lambda)(w),
$$
is a decorated $\ssc^+$-multisection.
\end{lemma}
\begin{proof}
The proof is trivial.
\qed \end{proof}
Next we consider locally finite sums of decorated sc$^+$-multisections. Recall that a subset $O$ of the ep-groupoid
$X$ is called {\bf saturated}\index{Saturated} provided $O=\pi^{-1}(\pi(O))$.
\begin{definition}\label{locally-finite-family}\index{D- Locally finite family}
We assume that $(P:W\rightarrow X,\mu)$ is a strong bundle over an ep-groupoid. The family  $(\Lambda_j)_{j\in J}$  of decorated $\ssc^+$-multisections is called 
{\bf locally finite}, provided for every $x\in X$ there exists an open saturated neighborhood $O(x)$ having  the property
that there are only finitely many indices $j\in J$ for which there exist $w\in W\vert O(x)$ satisfying  $w\neq 0$ and $\Lambda_j(w)>0$.
\qed
\end{definition}
We can take the saturated $O(x)$ always in such a way that for every other saturated $O'(x)\subset O(x)$ the 
following equality holds.
\begin{eqnarray*}
\{j\in J\ |\ \exists w\in W|O(x),\ w\neq 0,\ \Lambda_j(w)>0\}\\
=\{j\in J\ |\ \exists w\in W|O'(x),\ w\neq 0,\ \Lambda_j(w)>0\}.
\end{eqnarray*}
In the following we always assume that such a $O(x)$ is taken. 
Assume we are given a family $(\Lambda_j)$ as described in the definition. 
Pick $x\in X$ and take the associated saturated $O(x)$. Denote by $J_x$ the finite set of indices
$j$ for which there exists $w\neq 0$, $w\in W|O(x)$ with $\Lambda_j(w)>0$. For every $j\in J_x$ 
we consider $\Lambda_j$ and $\mathfrak{U}_j$ and $\mathfrak{S}_j$.  We pick $U_j(x)$ and ${(s_{j,i}^x)}_{i\in I_x^j}$.
Then define $U(x)$ by
$$
U(x)=O(x) \bigcap \left(\bigcap_{j\in J_x} U_j(x)\right).
$$
The  index set $I_x$ is defined by $I_x = \prod_{j\in J_x} I_x^j$. The elements of $I_x$ are written as
${(i_j)}_{j\in J_x}$ or $(i_j)$ for short.
We define  $\mathfrak{s}_x ={(s_{i_j})}_{(i_j)\in I_x}$ by
$$
s^x_{(i_j)} = \sum_{j\in J_x} s^{jx}_{i_j}\ \ \text{on}\ \ U(x).
$$
With $\mathfrak{S}={(\mathfrak{s}_x)}_{x\in X}$ and  $\mathfrak{U}={(U(x)}_{x\in X}$ 
the following result is trivial.
\begin{lemma}\index{L- Well-definedness $\oplus_{i\in I}\Lambda_i$}
Let $(P:W\rightarrow X,\mu)$ be a strong bundle over an ep-groupoid. If $(\Lambda_i)_{i\in I}$ is a locally finite set of 
decorated $\ssc^+$-multisections, then  the 
sum $\Lambda =\oplus_{i\in I} \Lambda_i$ is a well-defined  decorated $\ssc^+$-multisection $[\oplus_{j\in J}\Lambda_j,\mathfrak{U},\mathfrak{S}]$.
\qed
\end{lemma}

\section{Structurable Sc\texorpdfstring{$^+$}{pp}-Multisections}\label{SECT133}

Assume that $P:W\rightarrow X$ is a strong bundle over an ep-groupoid and $\Lambda:W\rightarrow {\mathbb Q}^+$
a decorated sc$^+$-multisection.
For more complicated constructions, for example extension theorems, it is important that the different symmetric
sc$^+$-section structures are coherent, which means they are related in a suitable sense. This will be discussed in the following.

Assume that $[\Lambda,\mathfrak{U},\mathfrak{S}]$ is a decorated sc$^+$-multisection.
Taking a representative $(\Lambda,\mathfrak{U},\mathfrak{S})$ we have the following data.
\begin{itemize}
\item A good system of open neighborhoods ${(U(x)}_{x\in X}$, see Definition \ref{DEFF1}.
\item For $x\in X$ a symmetric section structure $\mathfrak{s}_x={(s_i^x)}_{i\in I_x}$ for $\Lambda$, see Definition \ref{XDEF1}, where the $s^x_i$ are sc$^+$-sections defined for $W|U(x)$.
\end{itemize}
For the moment we shall only consider the index sets ${(I_x)}_{x\in X}$.
In order to formalize what it means that the local choices of symmetric $\ssc^+$-section structures are coherent  we  
consider diagrams of maps,  written as
\begin{equation}\label{just}
d:I_x \stackrel{a}{\twoheadleftarrow} I \stackrel{b}{\twoheadrightarrow} I_{x'}, 
\end{equation}
where the maps $a$ and $b$ are surjective and the number of preimages in $a^{-1}(\{i\})$ is independent of $i\in I_x$
and the number of preimages  in $b^{-1}(\{j\})$ is independent of $j\in I_{x'}$.  Reading the  diagram from left to right we refer to \eqref{just} as a {\bf diagram from $I_x$ to $I_{x'}$} \index{$d:I_x\rightarrow I_{x'}$} and with the obvious abbreviation we shall write $d: I_x\rightarrow I_{x'}$.

\begin{definition}\index{D- Equivalence of diagrams}\index{D- Correspondence from $I_x$ to $I_{x'}$}
Two diagrams $d: I_x \stackrel{a}{\twoheadleftarrow} I \stackrel{b}{\twoheadrightarrow} I_{x'}$ and $d': I_x \stackrel{a'}{\twoheadleftarrow} I \stackrel{b'}{\twoheadrightarrow} I_{x'}$ are called {\bf equivalent} if there exists a bijection $c: I\to I'$ such that
the following diagram is commutative
$$
\begin{CD}
I_x @<a<< I @>b>> I_{x'}\\
@|  @V c VV @|\\
I_x @<a'<< I' @>b'>> I_{x'}.
\end{CD}
$$
We shall refer to the equivalence class $[d]$ of the diagram 
$d: I_x \stackrel{a}{\twoheadleftarrow} I \stackrel{b}{\twoheadrightarrow} I_{x'}$ as a {\bf correspondence from $I_x$ to $I_{x'}$} and often write $[d]: I_x\to I_{x'}$. 
\qed
\end{definition}
The set of  correspondences $[d]: I_x\to I_{x'}$
will be denoted by 
$$
\mathfrak{J}(x,x'):=\{[d]:I_x\rightarrow I_{x'}\}.\index{$\mathfrak{J}(x,x')$}
$$
\begin{definition}\index{D- Pseudo-inverse $\sharp$}
The {\bf pseudo-inverse} is a  bijection $\sharp: \mathfrak{J}(x,x')\rightarrow \mathfrak{J}(x',x)$ defined by 
$$
[I_x \stackrel{a}{\twoheadleftarrow} I \stackrel{b}{\twoheadrightarrow} I_{x'}]^\sharp := [I_{x'} \stackrel{b}{\twoheadleftarrow} I \stackrel{a}{\twoheadrightarrow} I_{x}].
$$
\qed
\end{definition}

In view of the  $G_x$ and $G_{x'}$-actions  on the index sets  $I_x$ and $I_{x'}$, respectively, we obtain a bi-action on $\mathfrak{J}(x,x')$
\index{D- Bi-action on $\mathfrak{J}(x,x')$}
$$
G'_{x'}\times \mathfrak{J}(x,x')\times G_{x} \rightarrow \mathfrak{J}(x,x'). 
$$
defined by
$$
g'\ast [I_x \stackrel{a}{\twoheadleftarrow} I \stackrel{b}{\twoheadrightarrow} I_{x'}]\ast g
:= [I_x \stackrel{g^{-1}\circ a}{\twoheadleftarrow} I \stackrel{g'\circ b}{\twoheadrightarrow} I_{x'}]
$$
Here we use the action of $G_x$ on $I_x$ which means that $g\in G_x$ induces a permutation of $I_x$ written as $g :I_x\rightarrow I_x$.
Similarly for $g'\in G_{x'}$ we have $g':I_{x'}\rightarrow I_{x'}$.

\begin{example}\label{example-simple}
Here is a simple  example. Assume that the index set $I_x$ consists of one point, $I_{x}=\{\ast\}$,  and let 
$$
[d]=[\{\ast\}\twoheadleftarrow I \stackrel{b}{\twoheadrightarrow} I_{x'}] :  \{\ast\}\rightarrow I_{x'}.
$$
Then we have 
$$
g'\ast[d]\ast g= [\{\ast\}\twoheadleftarrow I \stackrel{g'\circ b}{\twoheadrightarrow} I_{x'}]
$$
Using that $g'$ acts on $I_{x'}$ as a permutation, and the map $b: I\to I_{x'}$ is a surjection having the property that the preimages $b^{-1}(\{i\})$ 
have the same cardinality for every $i\in I_{x'}$, one can construct a bijection $c: I\to I$ satisfying $b\circ c = g' \circ b$. Consequently, 
the diagram $\{\ast\}\twoheadleftarrow I \stackrel{g'\circ b}{\twoheadrightarrow} I_{x'}$ is equivalent 
to the diagram   $\{\ast\}\twoheadleftarrow I \stackrel{b}{\twoheadrightarrow} I_{x'}$  and 
$$
g'\ast [\{\ast\}\twoheadleftarrow I \stackrel{b}{\twoheadrightarrow} I_{x'}]\ast g=[\{\ast\}\twoheadleftarrow I \stackrel{b}{\twoheadrightarrow} I_{x'}]
$$
for all $g\in G_x$ and $g'\in G_{x'}$.
\qed
\end{example}

\begin{definition}\index{D- Correspondence}
A {\bf correspondence} between the  two local symmetric sc$^+$-section structures $(U(x),{(s_i^x)}_{i\in I_x})$ at $x$ and $(U(x'),{(s_i^{x'})}_{i\in I_{x'}})$ at $x'$ is a map 
$$
\tau_{x,x'}:\bm{ U}(x,x')\rightarrow \mathfrak{J}(x,x')
$$
having the following properties:
\begin{itemize}
\item[(1)]\ $\tau_{x,x'}$ is locally constant and takes only finitely many values. 
\item[(2)]\ $\tau_{x,x'}(g'\ast\phi\ast g) =g'\ast\tau_{x,x'}(\phi)\ast g$ for all $\phi\in \bm{U}(x,x')$, $g\in G_x$, and $g'\in G_{x'}$.
\item[(3)]\ If $\phi\in \bm{U}(x, x')$ and $I_x \stackrel{a}{\twoheadleftarrow} I \stackrel{b}{\twoheadrightarrow} I_{x'}$ is a representative of 
 $\tau_{x, x'}(\phi)$, then 
$$
s^{x'}_{b(k)}(t(\psi))=\mu(\psi,s^{x}_{a(k)}(s(\psi)).
$$
for all $k\in I$ and all morphisms $\psi$ belonging to the connected component containing $\phi.$
\end{itemize}
\qed
\end{definition}
\begin{remark}\index{R- On corresponences}
A correspondence is an important piece of information how two different local sc$^+$-section structures are related
over open subsets $O\subset U(x)$ and $O'\subset U(x')$ which are given by the images (via $s$ and $t$)  of a connected component 
$\Sigma$. 
\qed
\end{remark}
The formulation in (3) depends on the choice of representative  $I_x \stackrel{a}{\twoheadleftarrow} I \stackrel{b}{\twoheadrightarrow} I_{x'}$
of $\tau_{x,x'}(\phi)$, but the required property does not.  To see this take a second representative for $\tau_{x,x'}(\phi)$
which together with a bijection $d:I'\rightarrow I$ fits into the following commutative diagram
$$
\begin{CD}
I_x @<a'<< I'@>b'>>I_{x'}\\
@|    @V d VV  @|\\
I_x  @<a<< I@>b>> I_{x'}.
\end{CD}
$$
It holds that $a\circ d=a'$ and $b\circ d=b'$.  This implies in view of the displayed property in (3)
that 
$$
s^{x'}_{b'(k)}(t(\psi))=s^{x'}_{b(d(k))}(t(\psi))=\mu(\psi,s^{x}_{a(d(k))}(s(\psi))=\mu(\psi,s^{x}_{a'(k)}(s(\psi)).
$$

Next we introduce  two  crucial concepts. The first is  that of a structurable $\ssc^+$-multisection functor and the second that of a structured
sc$^+$-multisection functor.
\begin{definition}\index{D- Structurable sc$^+$-multisection}\label{DEFN13.3.6}
Let $P:W\rightarrow X$ be a strong bundle over an ep-groupoid. 
A decorated $\ssc^+$-multisection functor $\Lambda:  W\rightarrow {\mathbb Q}^+$ is said to be {\bf structurable} if the following holds,
where $\Lambda=[\Lambda,\mathfrak{U},\mathfrak{S}]$. There exists a representative $(\Lambda,\mathfrak{U},\mathfrak{S})$
and for every $x,x'\in X$  there is a correspondence $\tau_{x,x'}: \bm{U}(x,x')\rightarrow\mathfrak{J}(x,x')$ 
between local section structures $(U(x), {\mathfrak s}_x)$ and $(U(x'), {\mathfrak s}_{x'})$ which, in addition, satisfies
\begin{itemize}
\item[(1)]\ If $x=x'$ it is required that
$\tau_{x,x}(\Gamma(g,y))=[I_x\stackrel{\id}{\twoheadleftarrow}I_x\stackrel{g}{\twoheadrightarrow} I_x]$ for all $g\in G_x$ and $y\in U(x)$.
\item[(2)]\ $\tau_{x',x}(\phi) = {(\tau_{x,x'}(\phi^{-1}))}^\sharp$ for $\phi\in \bm{ U}(x', x)$.
\end{itemize}
\qed
\end{definition}
\begin{remark}\index{R- On being structurable}
If $\Lambda$ is structurable the definition guarantees the existence of an additional finer structure  and we can consider
$(\Lambda,\mathfrak{U},\mathfrak{S},\tau)$, 
where $\tau$ stands for the collection of all $\tau_{x,x'}$, and $(x,x')$ varies over $X\times X$, and $\mathfrak{I}={(I_x)}_{x\in X}$. 
The notion of being structurable is the key notion. When carrying out constructions one usually takes an associated 
structured version and usually it does not matter which one.
\qed
\end{remark}

The tuple $(\Lambda,\mathfrak{U},\mathfrak{S},\tau)$ is a useful object in its own right.  Given another such tuple
 $(\Lambda',\mathfrak{U}',\mathfrak{S}',\tau')$ we shall introduce a notion 
 of equivalence.
 \begin{definition}\index{D- Equivalence of $(\Lambda,\mathfrak{U},\mathfrak{S},\tau)$}\label{DEFNX13.3.8}
 We call $(\Lambda,\mathfrak{U},\mathfrak{S},\tau)$ and  $(\Lambda',\mathfrak{U}',\mathfrak{S}',\tau')$ {\bf equivalent} provided the following holds.
 \begin{itemize}
 \item[(1)]\ $\Lambda=\Lambda'$.
 \item[(2)]\ There exists a good system of open neighborhoods $\mathfrak{U}''$ such that $U''(x)\subset U(x)\cap U'(x)$ for all $x\in X$.
 \item[(3)]\ For every $x\in X$ it holds $I_x=I_x'$.
 \item[(4)]\ For $x,x'\in X$ it holds on $\bm{U}''(x,x')=\{\phi\in \bm{X}\ |\ s(\phi)\in U''(x),\ t(\phi)\in U''(x')\}$
 that $\tau_{x,x'}=\tau'_{x,x'}$.
  \end{itemize}
An {\bf equivalence class} is denoted by $[\Lambda,\mathfrak{U},\mathfrak{S},\tau]$.\index{$[\Lambda,\mathfrak{U},\mathfrak{S},\tau]$}
\qed
\end{definition}
\begin{remark}
It would have been possible to define a different equivalence relation 
accommodating different choices of sets $I_x$. Namely instead of requiring 
that $I_x=I_x'$ we could have stipulated the existence of $G_x$-equivariant bijections
$b_x:I_x\rightarrow I_x'$ and in addition that $\tau_{x,x'}$ and $\tau'_{x,x'}$
are related by 
$$
b_x\circ a = a'\circ c\ \ \text{and}\ \ b_{x'}\circ b=b'\circ c,
$$
where 
 $$
 \tau_{x,x'}(\phi)=[I_x\stackrel{a}{\twoheadleftarrow} I\stackrel{b}{\twoheadrightarrow}I_{x'}]\ \text{and}\ \tau_{x,x'}'(\phi)=[I_x'\stackrel{a'}{\twoheadleftarrow} I'\stackrel{b'}{\twoheadrightarrow}I_{x'}']
 $$
 and  $c:I\rightarrow I'$ is a bijection.  However, in applications it does not seem 
 too much of a difference in general, and we picked the more restrictive formulation.
 \qed
\end{remark}

\begin{definition}\label{DEF13.3.8}\index{D- Structured ac$^+$-multisection}
An equivalence class $[\Lambda,\mathfrak{U},\mathfrak{S},\tau]$ is called a  {\bf structured sc$^+$-multisection}.
We shall often abbreviate $\Lambda\equiv [\Lambda,\mathfrak{U},\mathfrak{S},\tau]$ and call it a structured sc$^+$-multisection
and $(\Lambda,\mathfrak{U},\mathfrak{S},\tau)$ a {\bf representative} of the structured sc$^+$-multisection.
\qed
\end{definition}

With these definitions we can consider for a strong bundle $P:W\rightarrow X$ over the ep-groupoid $X$ various classes of sc$^+$-multisection functors.
\begin{definition}\index{D- Classes $\Gamma_m^+(P),\Gamma^+_s(P),\Gamma^+_{st}(P)$}\index{$\Gamma_m^+(P),\Gamma^+_s(P),\Gamma^+_{st}(P)$}
\begin{itemize}
\item[(1)]\ $\Gamma_m^+(P)$ denotes the set of sc$^+$-multisection functors.
\item[(2)]\  The set of structurable sc$^+$-multisection functors is denoted by $\Gamma^+_{s}(P)$.
\item[(3)]\  The set of structured sc$^+$-multisection functors is denoted by 
$\Gamma^+_{st}(P)$.
\end{itemize}
\qed
\end{definition}
Clearly,  we have a forgetful map 
$$
\mathsf{forget}:\Gamma^+_{st}(P)\rightarrow \Gamma^+_s(P),\quad \wt{\Lambda}\mapsto \Lambda.
$$
It is clear from the definition that every $\Lambda\in\Gamma^+_s(P)$ can lifted to an element in $\Gamma^+_{st}(P)$ (usually not uniquely).
In fact the definition of a structurable sc$^+$-multisection $\Lambda$ precisely says that there exists $\wt{\Lambda}\in \Gamma^+_{st}(P)$
with $\mathsf{forget}(\wt{\Lambda})=\Lambda$.

The next proposition shows that being structurable is preserved under the sum operation previously introduced 
for two $\ssc^+$-multisection functors. Moreover if $\Lambda_1$ and $\Lambda_2$ are structured so is $\Lambda_1\oplus\Lambda_2$.

\begin{proposition}
Let $P:W\rightarrow X$ be  a strong bundle over the ep-groupoid $X$. 
If $\Lambda_1,\Lambda_2: E\rightarrow {\mathbb Q}^+$ are structurable $\ssc^+$-multisections functors, then $\Lambda_1\oplus \Lambda_2$ is also a structurable $\ssc^+$-multisection  functor. If $\Lambda_1$ and $\Lambda_2$ are structured, the sc$^+$-multisection $\Lambda_1\oplus\Lambda_2$ is naturally structured.
\end{proposition}
\begin{proof}
For $i=1,2$, we lift  the structurable $\ssc^+$-multisection $\Lambda_i$ to 
a structured $\ssc^+$-mutisection and take representatives 
 $(\Lambda_i,\mathfrak{U}_i,\mathfrak{S}_i,\tau_i)$ . 
 In view Proposition \ref{prop_obvious} we may assume that $\mathfrak{U}:=\mathfrak{U}_1=\mathfrak{U}_2$.  We define 
  $I_x=I_x^1\times I^2_{x}$  for $x\in X$. Introduce the family ${\mathfrak I}$ of index sets 
by 
$$
\mathfrak{I} =\{I_x\, \vert \, \ x\in X\},
$$
and the family $\tau$ of maps $\tau_{x, x'}: {\bf U}(x, x')\to {\mathfrak J}(x, x')$ by 
$$
\tau_{x,x'}=\tau^1_{x,x'}\times \tau^2_{x,x'}, $$
where 
the right-hand side is defined as
$$
[I_x^1\stackrel{a^1}{\twoheadleftarrow}I^1\stackrel{b^1}{\twoheadrightarrow} I_{x'}^1]\times [I_x^2\stackrel{a^2}{\twoheadleftarrow}I^2\stackrel{b^2}{\twoheadrightarrow }I_{x'}^1]:=[I_x^1\times I^2_x\stackrel{a^1\times a^2}{\twoheadleftarrow}I^1\times I^2\stackrel{b^1\times b^2}{\twoheadrightarrow} I_{x'}^1\times I^2_{x'}].
$$
Next we  define the  family $\mathfrak{S}={(\mathfrak{s}_x)}_{x\in X}$ of $\ssc^+$-section structures  ${\mathfrak s}^x=( s^x_{(i,j)})_{(i, j)\in I_x}$  on $U(x)$ by  
$$
s_{(i,j)}^x= s_i^{1,x}+s_j^{2,x}.
$$
With the action of $G_x$ on the index set $I_x=I^1_x\times I^2_x$ defined by $g(i, j)=(g(i), g(j))$, the $\ssc^+$-sections structures $s^x_{i, j}$ are symmetric over $U(x)$. 
Moreover, if $x,x'\in X$ and the map $\tau_{x,x'}(\phi): I_{x}\to I_{x'}$ is represented by the diagram 
$$
I_x=I^1_{x}\times I^2_{x}\stackrel{a^1\times a^2}{\twoheadleftarrow}I^1\times I^2\stackrel{b^1\times b^2}{\twoheadrightarrow} I_{x'}^1\times I^2_{x'}=I_{x'},
$$
then  
\begin{equation*}
\begin{split}
s^{x'}_{(b^1(i),b^2(j))}(t(\phi))&= s^{1,x'}_{b^1(i)}(t(\phi)) + s^{2,x'}_{b^2(j)}(t(\phi))\\
&=\mu(\phi,s^{1,x}_{a^1(i)}(s(\phi))) +\mu(\phi,s^{2,x}_{a^2(j)}(s(\phi)))\\
&=\mu(\phi, s^x_{(a^1(i),a^2(j))}(s(\phi)))
\end{split}
\end{equation*}
for all $(i, j)\in I^1\times I^2$, 
showing  the compatibility of the local section structures. Clearly, the maps $\tau_{x, x'}$ are locally constant  and satisfy
$$
\tau_{x,x}(\Gamma (g, y))=[I_x\stackrel{\id }{\twoheadleftarrow}I_x\stackrel{g}{\twoheadrightarrow}I_{x}]
$$
 for all $g\in G_x$ and $y\in U(x)$.  This completes the proof of the proposition.
\qed \end{proof}
\begin{remark}
Observe that in the above construction the structured $\Lambda_1\oplus\Lambda_2$ is in general different from 
$\Lambda_2\oplus \Lambda_1$, since for example  the index sets are different, i.e.  $I_x^1\times I_x^2$ versus $I_x^2\times I^1_x$.
Of course, one could introduce an appropriate notion of equivalence to equate both sums, since the underlying functors
are the same. However, since being structured is an auxiliary concept we shall not study this further. For us the basic notion is being structurable and as such $\Lambda_1\oplus \Lambda_2=\Lambda_2\oplus \Lambda_1$.\qed
\end{remark}

We recall that  given $r\in \R$  and a $\ssc^+$-multisection functor $\Lambda: E\to \Q^+$, the 
$\ssc^+$-multisection functor $r\odot\Lambda$ is defined as follows.  If $r=0$,  it is determined by the requirement
 $(0\odot\Lambda)(0_x)=1$ for $x\in X$, and if $r\neq0$, then $(r\odot \Lambda)(e):=\Lambda(\frac{e}{r})$.
 The following result is obvious.
 \begin{lemma}
Let $P:W\rightarrow X$ be a strong bundle over an ep-groupoid. 
If  $\Lambda:W\rightarrow {\mathbb Q}^+$ is a structurable $\ssc^+$-multisection, then $r\odot\Lambda$ is a structurable $\ssc^+$-multisection for every $r\in {\mathbb R}$. If $\beta:X\rightarrow {\mathbb R}$ is an sc-smooth functor and  $\Lambda$ is a structurable    sc$^+$-multisection
 so is $\beta\odot\Lambda$.
\qed
\end{lemma}
Let $P:W\rightarrow X$ be a strong bundle over the ep-groupoid $X$. 
We call a $\ssc^+$-multisection functor $\Lambda: W\rightarrow {\mathbb Q}^+$ {\bf trivial near a point $x\in X$}\index{Trivial near a point} if  there exists an open neighborhood $O(x)$ of $x$ such that 
$$
\Lambda(0_y)=1\quad \text{for all $y\in O(x)$}.
$$
  This means that $\Lambda$ represents the zero-section 
over $O(x)$.  

Assume that $[\Lambda,\mathfrak{U},\mathfrak{S}]$ is structurable and pick a structured version represented by
$(\Lambda,\mathfrak{U},\mathfrak{S},\tau)$.  Fix for every $x\in X$, where $\Lambda$ is trivial near $x$, a new $G_x$-invariant open neighborhood $U'(x)$ contained in $U(x)$ so that $\Lambda(0_y)=1$ for $y\in U'(x)$. Otherwise keep $U(x)$. This defines a new
$\mathfrak{U}'$. 
 Denote by $T_\Lambda$ the collection of all $x$ so that
$\Lambda$ is trivial near $x$.  By  assumption on $\mathfrak{U}'$ it holds that  $\Lambda(0_y)=1$ for $y\in U'(x)$ for $x\in T_\Lambda$.

For $x\in T_\Lambda$ we redefine $I_x'=\{*\}$, i.e. a one-point set and take as section structure $\mathfrak{s}'_x$, which consists
of the single zero-section over $U(x)$ with weight $1$.
If $x\in X\setminus T_\Lambda$ we define $I_x':=I_x$.
Consider for $x,x'\in X$ 
$$
\tau_{x,x'}:\bm{U}(x,x')\rightarrow \mathfrak{F}(x,x').
$$
Assume that $\tau_{x,x'}$ on a connected component is represented by
$$
I_x\stackrel{a}{\twoheadleftarrow} I\stackrel{b}{\twoheadrightarrow} I_{x'}.
$$
If $x\in T_\Lambda$  we replace $a$ by the unique map $a'$ onto $I_x'$, and if $x'\in T_\Lambda$  we replace  $b$ by the unique map
$b'$ onto $I_{x'}'$.  Through these modifications we obtain $\tau'$.
We replace the representative $(\Lambda,\mathfrak{U},\mathfrak{S},\tau)$ by 
$(\Lambda,\mathfrak{U},'\mathfrak{S}',\tau')$.  We note that $\tau'$ is a correspondence between the local section structures.
Hence we obtain a new structured $[\Lambda,\mathfrak{U}',\mathfrak{S}',\tau']$,  which has the following  properties
\begin{description}
\item[(1)]  $U'(x)\subset U(x)$ for all $x\in T_\Lambda$ and $U'(x)=U(x)$ for $x\in X\setminus T_\Lambda$.
\item[(2)]  For every point  $x\in T_\Lambda$   we have that $\Lambda(0_y)=1$ for $y\in U'(x)$ and $I_x=\{*\}$.
\end{description}
For the following considerations we need to modify this even further, in fact somewhat drastically.
Pick for $x\in T_\Lambda$ any $U''(x)$ with the natural $G_{x}$-action such that $U''(x)\subset T_\Lambda$.
 For $x\in X\setminus T_\Lambda$ we define $U''(x):=U(x)$.  If $x,x'\in X\setminus T$ define $\tau_{x,x'}':=\tau_{x,x'}$ on $\bm{U}''(x,x')=\bm{U}(x,x')$. If $x\in T_\Lambda$ and $x'\in X\setminus T_\Lambda$ define $\tau_{x,x'}'$ on $\bm{U}''(x,x')$ by
 $$
 \tau_{x,x'}''(\phi) = [I_x=\{*\} \twoheadleftarrow I_{x'} \stackrel{Id}\twoheadrightarrow I_{x'}].
 $$
If $x\in X\setminus T_\Lambda$ and $x'\in T_\Lambda$ we define 
$$
 \tau_{x,x'}''(\phi) = [I_x\stackrel{Id}\twoheadleftarrow I_{x}\twoheadrightarrow I_{x'}=\{*\}].
 $$
In the case of $x,x'\in T_\Lambda$ we define 
$$
 \tau_{x,x'}''(\phi)=[\{\ast\}  \stackrel{Id}\twoheadleftarrow \{*\} \stackrel{Id}\twoheadrightarrow \{*\}].
 $$
With this data we obtain $(\Lambda,\mathfrak{U}'',\mathfrak{S}'',\tau'')$.   The previous construction proves the following result.
\begin{proposition}\label{prop2.34}
Let $(P:W\rightarrow X,\mu)$ be a strong bundle over the ep-groupoid $X$. Let $[\Lambda,\mathfrak{U},\mathfrak{S},\tau]$
be a  structured sc$^+$-multisection and $(\Lambda,\mathfrak{U},\mathfrak{S},\tau)$ a representative.
  Then there exists a structured  sc$^+$-multisection $[\Lambda,\mathfrak{U}'',\mathfrak{S}'',\tau'']$ where the representative 
   $(\Lambda,\mathfrak{U}'',\mathfrak{S}'',\tau'')$ has the following properties.
   \begin{itemize}
   \item[{\em(1)}]\ For $x\in T_\Lambda$ we have  $I_x''=\{*\}$ and we can take any open neighborhood $U''(x)$ with the property
   that it admits the $G_x$-action and $U''(x)\subset T_\Lambda$. The section structure consists of the single zero-section.
   \item[{\em(2)}]\ For $x\in X\setminus T_\Lambda$ it holds  $U''(x)=U(x)$ and $I_x''=I_x$.
   \item[{\em(3)}]\ For $x,x'\in X\setminus T_\Lambda$ it holds $\tau_{x,x'}''=\tau_{x,x'}$.
   \item[{\em(4)}]\ For $x\in T_\Lambda$ and $x'\in X\setminus T_\Lambda$ we have $\tau_{x,x'}''(\phi)=[\ast\twoheadleftarrow I_{x'}''\stackrel{Id}\twoheadrightarrow I_{x'}'']$ and for $x\in X\setminus T_\Lambda$ and $x'\in T_\Lambda$ we take $\tau_{x,x'}''(\phi)=[I_x\stackrel{Id}\twoheadleftarrow I_x\twoheadrightarrow \{*\}]$.
   \item[{\em(5)}]\ For $x,x'\in T_\Lambda$ we have $\tau_{x,x'}''(\phi)=[\{\ast\}\twoheadleftarrow\{*\}\twoheadrightarrow \{*\}]$.
   \end{itemize}
   \qed
   \end{proposition}
\begin{remark}
The fact that $U''(x)$ for $x\in T_\Lambda$ can be taken arbitrarily subject to the requirements $U''(x)\subset T_\Lambda$
and $U''(x)$ admits the natural $G_x$-action will be important.
\qed
\end{remark}
In view of Proposition \ref{prop2.34} we can make the following definition.
\begin{definition}\index{D- Tightly structured}
Let $(P:W\rightarrow X,\mu)$ be a strong bundle over an ep-groupoid  with paracompact $|X|$.  
A structured $[\Lambda,\mathfrak{U},\mathfrak{S},\tau]$ is called {\bf tightly structured} provided there exists a representative
$(\Lambda,\mathfrak{U},\mathfrak{S},\tau)$ having the following properties.
\begin{itemize}
\item[(1)]\ For $x\in T_\Lambda$ we have  $I_x=\{*\}$ and $\mathfrak{s}_x$ consists of the  single zero-section.
 \item[(2)]\ For $x\in T_\Lambda$ and $x'\in X\setminus T_\Lambda$ we have $\tau_{x,x'}(\phi)=[\ast\twoheadleftarrow I_{x'}\stackrel{Id}\twoheadrightarrow I_{x'}]$ and for $x\in X\setminus T_\Lambda$ and $x'\in T_\Lambda$ we have $\tau_{x,x'}(\phi)=[I_x\stackrel{Id}\twoheadleftarrow I_x\twoheadrightarrow \{*\}]$.
   \item[(3)]\ For $x,x'\in T_\Lambda$ we have $\tau_{x,x'}''(\phi)=[\{\ast\}\twoheadleftarrow\{*\}\twoheadrightarrow \{*\}]$.
   \end{itemize}
\qed
\end{definition}
With this definition and employing Propsoition \ref{prop2.34} we obtain the following corollary.
\begin{corollary}
Let $(P:W\rightarrow X,\mu)$ a strong bundle over the ep-groupoid $X$.
For every structurable $\Lambda$ there exists a tightly structured version.
\qed
\end{corollary}

A basic result is the following.  
\begin{proposition}\label{sum_of_structurable_is_structurable}
Let $P: W\rightarrow X$ be a strong bundle over an ep-groupoid $X$ and let 
 ${(\Lambda_l)}_{l\in L}$ be  a family of structurable $\ssc^+$-multisection functors having locally finite domain supports. 
Then  the sum $\bigoplus_{l\in L}\Lambda_l$ is a structurable $\ssc^+$-multisection functor.
\end{proposition}
\begin{proof}
We consider representatives $(\Lambda_l,\mathfrak{U}^l, \mathfrak{S}^l,\tau^l)$, $l\in L$,
 of   structured versions. For every $l\in L$ set $T_l:=T_{\Lambda_l}$. 
We define a map $X\rightarrow 2^L:x\rightarrow  J_x$ by associating to $x$ the finite set of all $l\in L$ such that
for every open neighborhood $O(x)$ it holds $O(x)\cap (X\setminus T_l)\neq \emptyset$.
In particular we find for every $x\in X$ an open neighborhood $V(x)$ satisfying
\begin{description}
\item[(1)] $V(x)$ admits the natural $G_x$-action.
\item[(2)] $V(x) \subset  T_l $  for all $l\in L\setminus J_x$.
\end{description}
 For $x\in X$ define $U(x)$ by
 $$
 U(x):= V(x)\bigcap\left (\bigcap_{l\in J_x} U^l(x)\right).
 $$
The following holds for $\mathfrak{U}={(U(x))}_{x\in X}$:
\begin{description}
\item[(3)] For $x\in X$ and $l\not\in J_x$ we have  that $U(x)\subset T_l$. 
\item[(4)] If $l\in J_x$ then $U(x)\subset U^l(x)$.
\end{description}
With the above properties guaranteed we can pass, following Proposition \ref{prop2.34}, for every $l\in L$ to a tightly structured version
using the same $\mathfrak{U}$ as a good system of open neighborhoods.  Without loss of generality we may assume therefore
that the $(\Lambda_l,\mathfrak{U}^l, \mathfrak{S}^l,\tau^l)$ already have these properties.  In particular
$\mathfrak{U}=\mathfrak{U}^l$ for $l\in L$.   Define
$$
\Lambda=\bigoplus_{l\in L}\Lambda_l.
$$
We use $\mathfrak{U}$  as a  good system. Define $K^l_x = I^l_x$ if $l\in J_x$ and otherwise $K^l_x=\{*\}$.
The product $K_x$ defined as 
$$
K_x:=\prod_{l\in L} K^l_x
$$
is a finite set. We can view an element $\mathsf{k}\in K_x$ as a map which assigns to $l\in L$ an element $\mathsf{k}(l)\in K^i_x$.
For $\mathsf{k}\in K_x$ we define  $\mathfrak{s}_x ={(s^x_{\mathsf{k}})}_{\mathsf{k}\in K_x}$  on $U(x)$ by
$$
s^x_{\mathsf{k}} = \sum_{l\in L} s^{l,x}_{\mathsf{k}(l)}.
$$
Note that this is a finite sum since up to finitely many $l$ the other terms are the zero-sections.
This gives $\mathfrak{S}={(\mathfrak{s}_x)}_{x\in X}$.  
Next we define $\tau_{x,x'}(\phi)$ on $\bm{U}(x,x')$. The obvious idea is to take the product 
 product of the ${ (\tau^l_{x,x'}(\phi))}_{l\in L}$, but one has to make sure that this is well-defined.
Every $\tau^l_{x,x'}(\phi)$ for $l\in L$ has a representative of the form
$$
I^l_x \twoheadleftarrow \bar{I}^l_x\twoheadrightarrow I^l_{x'}.
$$
Moreover, this diagram takes the following forms given $x$ and $x'$.
\begin{itemize}
\item[(1)]\ With the exception of finitely many $\ell\in L$ : $\{\ast\}\twoheadleftarrow \{\ast\}\twoheadrightarrow \{\ast\}$ for $x,x'\in T_l$ 
\item[(2)]\ For finitely many values of $l\in L$ : $ \{\ast\}\twoheadleftarrow I^l_{x'} \twoheadrightarrow I^l_{x'}$ for $x\in T_l$.
\item[(3)]\ For finitely many values $l\in L$ : $ I^l_x\twoheadleftarrow I^l_{x} \twoheadrightarrow \{\ast\}$ for $x'\in T_l$.
\item[(4)]\  For finitely many values $l\in L$ : $I^l_x \twoheadleftarrow \bar{I}^l_x\twoheadrightarrow I^l_{x'}$ 
\end{itemize}
From this it follows that the product of the $\tau_{x,x'}^l$ is defined and a correspondence.
 \qed \end{proof}

\section{Equivalences, Coverings and Structurability}\label{SECT134}
Next we shall study the properties of structured sc$^+$-multisections with respect to equivalences of strong bundles
and generalized strong bundle isomorphisms, but also with respect to proper covering functors.

\begin{theorem}\label{pullbackX}\index{T- Pull-back of structurable $\Lambda$}
Let $(P:W\rightarrow X,\mu)$ and $(P':W'\rightarrow X',\mu')$ be strong bundles over ep-groupoids.
Suppose $\Phi:W\rightarrow W'$ is a strong bundle equivalence covering the equivalence $F:X\rightarrow X'$ and $\Lambda':W'\rightarrow {\mathbb Q}^+$ is a
structurable sc$^+$-multisection functor.
Then the sc$^+$-multisection functors $\Phi^\ast\Lambda':=\Lambda'\circ\Phi$ is  structurable. 
If $\Lambda'$ is structured there is a canonical way to structure $\Phi^\ast\Lambda'$.
\end{theorem}
\begin{proof}
We start by showing that $\Phi^\ast\Lambda' = \Lambda'\circ \Phi$ is structurable.
Fix a structured $[\Lambda',\mathfrak{U}',\mathfrak{S}',\tau']$ with representative $(\Lambda',\mathfrak{U}',\mathfrak{S}',\tau')$.
Recall that we can take the sets in $\mathfrak{U}'$ as small as we wish, so that the other associated new data are just the restrictions of the old data.

For every $x\in X$ we find $U(x)$ with the natural $G_x$-action, so that
\begin{description}
\item[(1)]  $t:s^{-1}(\cl_X(U(x)))\rightarrow X$ is proper, and 
\item[(2)] $F:U(x)\rightarrow F(U(x))$ is an sc-diffeomorphism. 
\item[(3)] $F(U(x))\subset  U'(F(x))$, and without loss of generality $F(U(x))=U'(F(x))$.
\end{description}

We define $I_{x}:=I'_{F(x)}$ and $\mathfrak{U}={(U(x))}_{x\in X}$. Given $x,y\in X$ the equivalence 
$F$ defines an sc-diffeomorphism ${\bf U}(x,y)\rightarrow {\bf U}'(F(x),F(y))$ and we infer that $s:{\bf U}(x,y)\rightarrow U(x)$
and $t:{\bf U}(x,y)\rightarrow U(y)$, when restricted to a connected component,  are sc-diffeomorphisms onto open subsets. 
With other words $\mathfrak{U}$ is a good system of open neighborhoods.

We define $\tau$ by
$$
\tau_{x,y}(\phi):=\tau'_{F(x),F(y)}(F(\phi)).
$$
Finally the symmetric sc$^+$-section structures $\mathfrak{s}_x={(s^x_i)}_{i\in I_x}$ are defined via
$$
\Phi(s^x_i(z)) = {s'}^{F(x)}_i(F(z))\ \text{for}\ z\in U(x),\ i\in I_x,
$$
and determine $\mathfrak{S}$. 
Hence $\Phi^\ast\Lambda'$ is structurable and $[\Phi^\ast\Lambda',\mathfrak{U},\mathfrak{S},\tau]$ is the canonically 
structured version when starting with $[\Lambda',\mathfrak{U}',\mathfrak{S}',\tau']$.
\qed \end{proof}
There is a similar result for the push-forward $F_\ast\Lambda$. However, there is not a canonical structuring 
for the latter. Here we shall also need an extension result which we formulate now.
\begin{proposition}\index{P- Extension of $\Lambda$}\label{PROPEXT}
Let $(P:W\rightarrow X,\mu)$ be a strong bundle over an ep-groupoid and assume that $V$ is an open subset of $X$ satisfying
$|X|=|V|$. We view $V$ as an ep-groupoid and assume $\Lambda:W|V\rightarrow {\mathbb Q}^+$ is a given structured sc$^+$-multisection functor.
Then the following holds.
\begin{itemize}
\item[(1)]\ $\Lambda$ has a canonical extension to $W$ denoted by $\bar{\Lambda}$.
\item[(2)]\  Given open neighborhoods $O(x)$ for $x\in X\setminus V$ the sc$^+$-multisection $\bar{\Lambda}$ has a structured version
so that the induced data to $V$ is the original structured $\Lambda$ and the sets $U(x)$ for $x\in X\setminus V$ are contained in $O(x)$.
\end{itemize}
\end{proposition}
\begin{proof}
(1) For $w\in W$ we can pick find $w'\in W|V$ and $\phi:P(w')\rightarrow P(w)$ such that $\mu(\phi,w')=w$.  We define
$$
\bar{\Lambda}(w):=\Lambda(w').
$$
One easily verifies that this definition does not depend on the choices involved. The extension $\bar{\Lambda}$ is an sc$^+$-multisection functor.
Namely we can take open neighborhoods $U(\phi)$, $U(P(w))$ and $U(P(w'))$ so that
$$
t:U(\phi)\rightarrow U(P(w))\ \ \text{and}\ \ s:U(\phi)\rightarrow U(P(w'))
$$
are sc-diffeomorphisms. We may assume that $U(P(w'))$ is small enough so that it supports a local sc$^+$-section structure.
With the data at hand we can push it forward to $U(P(w))$, where it represents $\bar{\Lambda}$. \par

\noindent (2) For every $q\in V$ we have an open $U(q)\subset V$ invariant under the $G_q$-action, having the properness property
and supporting a symmetric sc$^+$-section structure $\mathfrak{s}_q$. In addition, the collection $\mathfrak{U}={(U(q))}_{q\in V}$,
is a good system for $V$.  In addition we have correspondences $\tau_{q,q'}$ relating the local sc$^+$-section structures.

We employ the extension theorem for $\mathfrak{U}$ and obtain ${(U(x))}_{x\in X}$ which extends $\mathfrak{U}$
and where $U(x)\subset O(x)$ for $x\in X\setminus V$. We also need to extend the correspondences and define specific 
local sc$^+$-section structures. For this we first have a closer look at the proof of Theorem \ref{THMX1317} to see how
data was being moved around. Without loss of generality we assume that $U(x)=O(x)$, where $U(x)$ will be fixed arbitrarily small.

For every $x\in X\setminus V$ we fix $\psi_x\in \bm{X}$ with $q_x:= s(\psi_x)$ and $t(\psi_x)=x$. Moreover, we find open neighborhoods
$O(\psi_x)$, $O(q_x)$, and $U(x)$ having the following properties.
\begin{itemize}
\item $s:O(\psi_x)\rightarrow O(q_x)$ and $t:O(\psi_y)\rightarrow U(x)$ are sc-diffeomorphisms.
\item $O(q_x)\subset U(q_x)$ and $O(q_x)$ is invariant under the $G_{q_x}$-action on $U(q_x)$.
\item  $U(x)$ admits the natural $G_x$-action and has the properness property.
\item  For every connected component $\Sigma$ of $\bm{U}(q_x,x)$ the maps $s:\Sigma\rightarrow U(q_x)$ and $t:\Sigma\rightarrow U(x)$ are sc-diffeomorphisms onto their image.
\end{itemize}
We have seen in Theorem \ref{THMX1317} that ${(U(x))}_{x\in X}$ is a good system. We need to define local section structures on the new sets 
and correspondences for the following types of sets 
\begin{itemize}
\item[(i)]\ \ \ $\bm{U}(q,x)$ for $q\in V$ and $x\in X$.
\item[(ii)]\ \ \ $\bm{U}(x,q)$ for $x\in X$ and $q\in V$.
\item[(iii)]\ \ \ $\bm{U}(x,y)$ for $x, y\in X$.
\end{itemize}
We start with the local sc$^+$-section structures.  Pick $x\in X\setminus V$ and define the  equivariant sc-diffeomorphism
$$
D_x\colon O(q_x) \rightarrow U(x)\colon D(q) =t\circ (s|U(\psi_x))^{-1}(q),
$$
where $\gamma_x:G_{q_x}\rightarrow G_x$ is defined by $\gamma_x(\phi)  \psi_x\circ \phi\circ \psi_x^{-1}$.
With $O(q_x)\subset U(q_x)$ we push forward the restriction of $\mathfrak{s}_{q_x}$ to $O(q_x)$. More precisely with
$\mathfrak{s}_{q_x}={(s_i^{q_x})}_{i\in I_{q_x}}$ we define $I_x:= I_{q_x}$ and $s_i^x$ by
$$
s_i^x(D_x(q)) =\mu((s|O(\psi_x))^{-1}(q),s^{q_x}_i(q))\ \ \text{for}\ \ q\in O(q_x).
$$
We let $G_x$ act on $I_x=I_{q_x}$ by $g(i):= (\gamma_x^{-1}(g))(i)$.  This prescription defines for every $x\in X\setminus V$
a symmetric local sc$^+$-section structure on $U(x)$ representing $\bar{\Lambda}$.

In a next step we have to define the correspondences. In case (i) we take the sc-smooth embedding
$$
A: \bm{U}(q,x)\rightarrow \bm{U}(q,q_x)\colon A(\phi) = [(t|O(\psi_x))^{-1}(t(\phi))]^{-1}\circ \phi,
$$
whose image is $\{\psi\in\bm{X}\ |\ s(\psi)\in U(q),\ t(\psi)\in O(q_x)\}$. We define 
$\tau_{q,x}:\bm{U}(q,x)\rightarrow \mathfrak{I}(q,x)$ by
$$
\tau_{q,x}(\phi):= \tau_{q,q_x}\circ A(\phi).
$$
This map is obviously locally constant since this did hold for $\tau_{q,q_x}$.  Moreover, it defines a correspondence between
the local sc$^+$-section structures as we shall show now. Take $\phi\in \bm{U}(q,x)$ and define 
$$
\psi:= A(\phi) = [(t|O(\psi_x))^{-1}(t(\phi))]^{-1}\circ \phi \in \{\sigma\in\bm{X}\ |\ s(\sigma)\in U(q),\ t(\sigma)\in O(q_x)\}.
$$
Assume $\tau_{q,q_x}(\psi)$ is represented by $I_x\twoheadleftarrow I\twoheadrightarrow I_{q_x}$ with maps $a$ and $b$ so that 
$$
s_{b(i)}^{q_x}(t(\psi))=\mu(\psi, s_{a(i)}^q(s(\psi))).
$$
By construction $\tau_{q,x}(\phi)$ is represented by the same diagram where $I_x=I_{q_x}$ by definition. We compute
\begin{eqnarray*}
s_{b(i)}^x(t(\phi))&=& \mu((s|O(\psi_x))^{-1}(D_x^{-1}(t(\phi))), s_{b(i)}^{q_x}(D_x^{-1}(t(\phi))))\\
&=& \mu((s|O(\psi_x))^{-1}(t(\psi)), s_{b(i)}^{q_x}(t(\psi)))\\
&=&\mu((s|O(\psi_x))^{-1}(t(\psi)), \mu(\psi, s_{a(i)}^q(s(\psi))))\\
&=&\mu((s|O(\psi_x))^{-1}(t(\psi)), \mu(\psi, s_{a(i)}^q(s(\phi))))\\
&=&\mu(\phi, s_{a(i)}^q(s(\phi))).
\end{eqnarray*}
This shows the compatibility.  The case (ii) is similar, where we use the pseudo-inverse of $\tau_{q,x}$. 
In case (iii) we define $\tau_{x,y}$ as follows.  Recall that we have the  sc-diffeomorphism
$$
C:\bm{U}(x,y)\rightarrow \{\phi\in\bm{X}\ |\ s(\phi)\in O(q_x),\ t(\phi)\in O(q_y)\}
$$
defined by
$$
C(\phi) =\left ((t|O(\psi_y))^{-1}(t(\phi))\right)^{-1}\circ \phi\circ \left((t|O(\psi_x))^{-1}(s(\phi))\right).
$$
In this case we define $\tau_{x,y}$ by 
$$
\tau_{x,y}=\tau_{q_x,q_y}\circ C.
$$
Through a straight forward computation as in case (i) one verifies the assertion.
\qed \end{proof}
\begin{theorem}\label{pushforwardX}\index{T- Push-forward of structurable $\Lambda$}
Let $P:W\rightarrow X$ and $P':W'\rightarrow X'$ be strong bundles over ep-groupoids. Assume that $\Phi:W\rightarrow W'$ is a strong bundle equivalence covering the equivalence $F:X\rightarrow X'$ and $\Lambda:W\rightarrow {\mathbb Q}^+$ 
a structurable $\ssc^+$-multisection functor. Then the push-forward $\Phi_\ast\Lambda$  is a  structurable $\ssc^+$-multisection functor. 
\end{theorem}
\begin{proof}
Define  $\Lambda':W'\rightarrow {\mathbb Q}^+$
as follows.  If $w'\in W'$, using the fact that $\Phi$ is a strong bundle equivalence covering the equivalence of ep-groupoids $F$,
we find $w\in W$, $\psi\in \bm{X}'$ such that
$$
w' =\mu'(\psi,\Phi(w))\ \text{and}\ \ \psi:F(P(w))\rightarrow P'(w').
$$
Define $\Lambda'(w'):=\Lambda(w)$ and  note that the definition does not depend on the choices involved.
Hence we have defined, as a functor
$$
\Lambda' \colon W'\rightarrow {\mathbb Q}^+.
$$
We pick a structured version   $[\Lambda, \mathfrak{U},\mathfrak{S},\tau]$  represented by $(\Lambda, \mathfrak{U},\mathfrak{S},\tau)$ of the structurable 
$\ssc^+$-multisection functor $\Lambda: W\to \Q^+$.  Since $F$ is an equivalence, and hence a local sc-diffeomorphism, we may assume, in view of Theorem  \ref{xxxx-structure}, that for every open neighborhood $U(x)\in {\mathfrak U}$ the map $F: U(x)\to F(U(x))$ is a sc-diffeomorphism
and $F(U(x))$ is contained in a set having the properness property.

For every $y\in F(X)$, we choose a point $x_y\in X$ such that $F(x_y)=y$ and  set 
$$
U'(y):= F(U(x_y))\quad \text{and}\quad  I'_y:= I_{x_y}.
$$
We define a $\ssc^+$-section structure $\mathfrak{s}_y'=(s'^{y}_i)_{i\in I'_y}$ on $U'(y)$ as the  
push-forward of $(s^{x_y}_i)_{i\in I_{x_y}}$  by $\Phi$ and obtain
$$
s'^{y}_{i}:= \Phi\circ s^{x_y}_i\circ (F\vert U(x_y))^{-1}\quad \text{for all $i\in I'_{y}$}.
$$
Using the fact that $F$ is faithful and full, the  map $F: \morp (x_{y}, x_{y'})\to \morp (y, y')$ is a bijection. Next we define the action of $G_y'$ on the index set $I_{y}'$ and the maps 
$$
\tau_{y, y'}':  \bm{U}'(y, y')\to {\mathfrak I}'(y, y').
$$
We put $I_{y}':=I_{x_y}$ and define  the action of the isotropy group $G_y'$ on the index set $I_{y}'$ by 
$$
g'(i)=((F|\text{mor}(x_y,x_y))^{-1}(g'))(i), \quad  \text{$g\in G_y'$,\, $i\in I'_{y}$}.
$$
The equivalence $F$ induces an sc-diffeomorphism
$$
\bm{U}(x_y,x_{y'})\rightarrow \bm{U}'(y,y'):\phi \rightarrow F(\phi).
$$
In order to see this, recall that the map is clearly a local sc-diffeomorphism.  
The map is surjective. Indeed, let $\phi'\in \bm{U}'(y,y')$ pick the uniquely determined
points $x\in U(x_y)$ and $x'\in U(x_{y'})$ such that $F(x)=s(\phi')$ and $F(x')=t(\phi')$.
Since $F$ is an equivalence there exists $\phi$ with $s(\phi)=x$ and $t(\phi)=x'$ such that
$F(\phi)=\phi'$. The map is also injective. If $F(\phi_1)=F(\phi_2)$ we see that 
\begin{eqnarray*}
&F(s(\phi_1))= s(F(\phi_1))=s(F(\phi_2))= F(s(F(\phi_2)))\ \ \text{and}&\\
&F(t(\phi_1))= t(F(\phi_1))=t(F(\phi_2))= F(t(F(\phi_2))),&
\end{eqnarray*}
from which we conclude that $\phi_1$ and $\phi_2$ have the same source and target. Since $F$ is faithful
this implies $\phi_1=\phi_2$.

Now we define the maps  $\tau'_{y,y'}: \bm{U}'(y, y')\to {\mathfrak I}'(y, y')$   by   
$$
\tau_{y,y'}'(\phi) = \tau_{x_y,x_{y'}}((F|\bm{U}(x_y,x_{y'}))^{-1}(\phi)), \quad \text{$\phi \in \bm{U}'(y, y')$ }.
$$
One can easily check that with these definitions  the $\ssc^+$-section structure $(s^y_i)_{i\in I'_y}$ is symmetric and  the maps $\tau_{y, y'}'$ are correspondences for all $y, y'\in F(X)$.   

We note that $(\Lambda'|(W'|F(X)),\mathfrak{U}',\mathfrak{S}',\tau')$ 
is a representative for a structured sc$^+$-multisection for the strong bundle $W'|F(X)\rightarrow F(X)$. Note that $F(X)$ is an ep-goupoid
for the induced structure as an open subset of $X'$, and moreover
$|F(X)|=|X'|$. Now we apply Propsition \ref{PROPEXT} and obtain that the extension $\Lambda'$ is structurable.

\qed \end{proof}

The structurable $\ssc^+$-multisections behave nicely with respect to pull-backs by proper covering functors.
\begin{theorem}\label{OTHM1344}\index{T- Pullback by proper covering functors}
Let $(P: E\rightarrow Y,\nu)$, $(Q:W\rightarrow X,\mu)$  be  strong bundles over  ep-groupoids  and $\Phi:E\rightarrow W$ a proper strong bundle covering map covering the proper covering functor $F:Y\rightarrow X$. If 
$\Lambda:W\rightarrow {\mathbb Q}^+$ is a structurable $\ssc^+$-multisection functor, then the pull-back
$\Phi^\ast\Lambda$ is a structurable $\ssc^+$-multisection functor.
\end{theorem}
\begin{proof}
Fix a structured version 
$[\Lambda,\mathfrak{U},\mathfrak{S},\tau]$  of $\Lambda$.
In view of the definition of a proper covering functor, without loss of generality we may assume that the open neighborhoods $U(x)\in {\mathfrak U}$ are sufficiently small 
so that the following holds:
\begin{itemize}
\item[$\bullet$]\  There are finitely many mutually disjoint open neighborhoods $V(y)\subset Y$ around every point 
$y\in F^{-1}(x)$ such that the restrictions $F: V(y)\to U(x)$  are sc-diffeomorphisms and 
 $$
 F^{-1}(U(x))=\bigcup_{z\in F^{-1}(x)}V(z),
 $$
 \item[$\bullet$]\ the open neighborhoods $\mathfrak{V}={(V(y))}_{y\in Y}$ form a good system of open neighborhoods on $Y$. 
 \end{itemize} 
In addition, since $F$ is a proper covering functor, the map $\bm{Y} \to \bm{X}{_{s}\times_F}Y$, defined by 
\begin{equation}\label{prop_cover_sc_diff}
\phi\mapsto (F(\phi),s(\phi)), 
\end{equation}
is a sc-diffeomorphism.

On every open neighborhood $U(x)$,  we are given the symmetric $\ssc^+$-section structure $\mathfrak{s}_x = {(s^x_i)}_{i\in I_x}$ and for every two points $x,x'\in X$ we have correspondences 
$$
\tau_{x,x'}:\bm{U}(x,x')\rightarrow \mathfrak{I}(x,x')
$$
 relating the local section structures ${\mathfrak s}_x$ and ${\mathfrak s}_{x'}$.

For every $y\in F^{-1}(x)$, we take as index set $I_y$ the index set $I_x$ and define the $\ssc^+$-section structure ${\mathfrak s}_y$ on $V(y)$ as the pull-back of
of ${\mathfrak s}_x$, 
\begin{equation}\label{compatibility_prop_cover_pull_back}
\Phi (s^y_i(z))= s^x_i(F(y)),\quad z\in V(y), i\in I_y.
\end{equation}
The action of  the isotropy group $G_y$ on the index set $I_y$ is defined by $g(i)=F(g)(i)$ for $g\in G_y$ and $i\in I_y$. 
Since the map defined in \eqref{prop_cover_sc_diff} is a sc-diffeomorphim, we deduce (for the natural actions) that 
$F(\Gamma^Y (g, z))=\Gamma^X(F(g), F(z))$ for all $g\in G_y$ and $z\in V(y)$.  Using this fact  we conclude that  
the sections defined by \eqref{compatibility_prop_cover_pull_back} are symmetric. Indeed,  we compute for $g\in G_y$, $i\in I_y$, and $z\in V(y)$, 
\begin{equation*}
\begin{split}
\Phi(s^y_{g(i)}(g\ast z))&=s^x_{g(i)}(F(g\ast z)=s^{x}_{F(g)(i)}(F(g)\ast F(z))\\
&=\mu(\Gamma^X(F(g), F(z)), s_i^x(F(z)))\\
&=\mu(F(\Gamma^Y(g, z)), \Phi (s_i^y(z)))\\
&=\Phi (\mu (\Gamma^Y (g, z), s_i^y(z))).
\end{split}
\end{equation*}
Since $\Phi$ is a fiber-wise isomorphism, we conclude that 
$$
s^y_{g(i)}(g\ast z)=\mu (\Gamma^Y (g, z), s_i^y(z))
$$
 which proves our claim.

Having defined the good system of  open neighborhoods $\mathfrak{V}={(V(y))}_{y\in Y}$, the family $\mathfrak{I}=(I_y)_{y\in Y}$ of index sets as well as the family  $\mathfrak{S}={(\mathfrak{s}_y)}_{y\in Y}$ of symmetric $\ssc^+$-section structures we define the family ${\mathfrak \tau}=(\tau_{y, y'})_{y, y'\in X}$ of maps $\tau_{y, y'}: V(y, y')\to \mathfrak{I}(y, y')$ as follows. 
If $y\in F^{-1}(x)$ and $y'\in F^{-1}(x')$,  we put
$$
\tau_{y,y'}(\phi):= \tau_{x,x'}(F(\phi))
$$
for all $\phi \in V(y, y')$. Clearly, the maps $\tau_{y,y'}$ are locally finite and satisfy  $\tau_{y',y}(\phi)=\tau_{y, y'}(\phi^{-1})^\sharp$ for $\phi\in V(y',y)$ and 
$\tau_{y, y}(\Gamma^Y (g, z))=[I_{y}\stackrel{\id}\twoheadleftarrow I_{y}\stackrel{g}\twoheadrightarrow I_{y}]$ for every $g\in G_{y}$ and $z\in V(y)$. Finally, we shall verify the compatibility condition.  
If  $y\in F^{-1}(x)$, $y'\in F^{-1}(x')$, $\phi\in {\bf V}(y,y')$,  and the diagram $I_{y}\stackrel{a}\twoheadleftarrow I\stackrel{b}\twoheadrightarrow I_{y'}$  is a representative of $\tau_{y,y'}(\phi)$, we compute 
\begin{equation*}
\begin{split}
\Phi\circ \nu(\phi,s^y_{a(i)}(s(\phi)))&=\mu(F(\phi),\Phi(s^y_{a(i)}(s(\phi)))\\
&=\mu(F(\phi),s^x_{a(i)}(s(F(\phi))))\\
&=s^{x'}_{b(i)}(t(F(\phi)))\\
&=\Phi\circ s^{y'}_{b(i)}(t(\phi)), 
\end{split}
\end{equation*}
which implies that 
$$
\nu(\phi,s^y_{a(i)}(s(\phi)))=s^{y'}_{b(i)}(t(\phi))
$$\
for every $i\in I$. The proof of the theorem is complete.
\qed \end{proof}
In a next step we introduce the notion of two structured $\ssc^+$-multisection functors $\Lambda_1,\Lambda_2$ being commensurable.
More precisely let $[\Lambda_i,\mathfrak{U}_i,\mathfrak{S}_i,{\tau}_i]$, $i=1,2$,  be structured sc$^+$-mulitsection functors for the strong bundle $P:W\rightarrow X$ over the ep-groupoid $X$ with paracompact orbit space.
\begin{definition}\index{D- Commensurable sum}
Let $P:W\rightarrow X$ be the strong bundle  over the ep-groupoid $X$ having paracompact orbit space and let 
$\Lambda_i\equiv[\Lambda_j,\mathfrak{U}_j,\mathfrak{S}_j,{\tau}_j]$, where  $j=1,2$,  be structured $\ssc^+$-mulitsection functors. 
\begin{itemize}
\item[(1)]\ The functors ${\Lambda}_1$ and ${\Lambda}_2$ are called {\bf  commensurable} if we can take 
representatives for which   $\mathfrak{U}_1=\mathfrak{U}_2$, $\tau_1 =\tau_2$, and 
$I^1_x=I^2_x$ for all $x\in X$. Denote such a choice of compatible data by $\mathfrak{c}$.
\item[(2)]\ The {\bf commensurable sum} of the commensurable $\ssc^+$-multisection functors ${\Lambda}_1$ and ${\Lambda}_2$, for
the data $\mathfrak{c}$, 
denoted by 
$$
{\Lambda}:={\Lambda}_1\boxplus_{\mathfrak{c}} {\Lambda}_2 :W\to 
 {\mathbb Q}^+, 
$$
is defined by $[\Lambda, \mathfrak{U},\mathfrak{S}, \tau]$, where 
the  underlying $\ssc^+$-multisections functor $\Lambda$ is given as follows.  
The symmetric $\ssc^+$-section structure ${\mathfrak{s}}_x=(s^x_i)_{i\in I_x}$ on $U(x)$ is defined by 
$$s^x_i = s_i^{1,x}+s_i^{2,x},$$
and if $x\in X$ and $w\in W$ satisfies $P(w)\in U(x)$, then
$$
\Lambda(w)=\frac{1}{\abs{I_x}}\abs{\{ i\in I_x\, \vert \, s_i^x(P(w))=w\}}.
$$
\end{itemize}
\qed
\end{definition}
The commensurable sum can also be taken when we are given an infinite family 
${([{\Lambda}_j,\mathfrak{U},\mathfrak{S}_j,\tau])}_{j\in J}$, where every two are commensurable 
and in addition the underlying family $(\Lambda_j)$ is locally finite, i.e. there exists for every $x\in X$ an open neighborhood $O=O(x)$
so that with the exception of finitely many indices we have for $w\in W_y$, $y\in O$, that $\Lambda_j(w)=0$ provided $w\neq 0_y$. In this case
$$
\Lambda := \boxplus_{j\in J}^\mathfrak{c} \Lambda_j\ \ \ \text{(Here $\mathfrak{c}$ denotes as before a choice.)}
$$
is well-defined and a structured sc$^+$-multisection functor with data $(\mathfrak{U},\mathfrak{S},\tau)$.
 This is evident and does not require a proof.  The commensurable sum is not intrinsic and depends on choices.
 Nevertheless it is very important in constructions utilizing partition of unity arguments.

\section{Constructions of  Sc\texorpdfstring{$^+$}{yuo}-Multisections}
In this subsection we demonstrate how sc$^+$-multisection with certain properties  can be constructed.

\begin{definition}\label{KDEF}\index{D- Atomic sc$^+$-multisection}
A $\ssc^+$-multisection functor $\Lambda:W\rightarrow {\mathbb Q}^+$ is called  {\bf atomic} if the following conditions are satisfied: 
\begin{itemize}
\item[(1)]\  There exist a point $x\in X$ and $G_x$-invariant open neighborhoods  $U(x)$ and $V(x)$ such  that 
$\cl_X(V(x))\subset U(x)$ and $t: s^{-1}(\cl_X(U(x)))\rightarrow X$ is proper. 
\item[(2)]\ There exists a symmetric $\ssc^+$-section structure $(s_i)_{i\in I_x}$ on $U(x)$ vanishing 
outside of $V(x)$.  
\item[(3)]\ The functor $\Lambda$ is related to the local section structure as follows:
\begin{itemize}
\item[$\bullet$]  If $w\in W$ and $w'=\mu(\phi,w)$, where $\phi : P(w)\to t(\phi)$ is a morphism  such that  $t(\phi)\in U(x)$, then 
$$
 \Lambda(w) = \frac{1}{\abs{I_x}}\abs{ \{i\in I_x\ |\ s_i(P(w'))=w'\}}.
 $$
\item[$\bullet$] If $w\in W$ and there is no morphism $\phi$ satisfying  $s(\phi)=P(w)$ and $t(\phi)\in U(x)$, then 
 $$
 \Lambda(w)=\begin{cases}
1, &\quad  \text{if $w=0$,}\\
0, &\quad  \text{otherwise.}
 \end{cases}
 $$
 \end{itemize}
 \end{itemize}
 \end{definition}
\begin{remark}\index{R- Atomic perturbations}
  In the perturbation theory one often introduces the atomic $\ssc^+$-multi\-section $\Lambda$ by first constructing one section which takes 
 at a prescribed smooth point $x$ a particular value, or has a tangent map with a prescribed property. Then one moves this section around 
 with the action of $G_x$ and finally  extends by a standard procedure. If one perturbs over a compact solution set, then  a finite number of such perturbations added together are sufficient  to achieve transversality. In order to carry out such a procedure one needs the existence of sc-smooth   bump functions but not necessarily sc-smooth partitions of unity.  
 \qed
  \end{remark}

 \begin{theorem}\label{THME1353}\index{T- Structurability}
Let $P:W\rightarrow X$ be a strong bundle over an ep-groupoid 
 If  $\Lambda:W\rightarrow {\mathbb Q}^+$ is an  atomic $\ssc^+$-multisection functor, then 
it  is structurable. 
\end{theorem}
\begin{proof}
We consider  an  atomic $\ssc^+$-multisection functor $\Lambda:W\rightarrow {\mathbb Q}^+$. 
We have  to construct a system of good neighborhoods ${\mathfrak U}=(U(x))_{x\in X}$, the family ${\mathfrak I}=(I_x)_{x\in X}$ of index sets,  the family of correspondences ${\mathfrak \tau}=(\tau_{x, x'})_{x, x'\in X}$, and the family ${\mathfrak S}=({\mathfrak s}_x)_{x\in X}$ of symmetric $\ssc^+$-section structures on $U(x)$ satisfying the compatibility conditions.

Since $\Lambda$ is atomic   there are $G_{x_0}$-invariant open neighborhoods $U(x_0)$ and $V(x_0)$ such that $\cl_{X}V(x_0)\subset U(x_0)$ and the target map $t: \cl_{X}(s^{-1}(U(x_0))\to X$ is proper. Moreover, on $U(x_0)$, there exists a symmetric $\ssc^+$-section structure ${\mathfrak s}_{x_0}=(s^{x_0}_i)_{i\in I_{x_0}}$ such that the sections $s^{x_0}_i$ vanish on 
$U(x_0)\setminus V(x_0)$ and $\Lambda$ can be written in terms of this data
$$
\Lambda(w)=\frac{1}{|I_{x_0}|}\cdot |\{i\in I_{x_0}\ |\ s_i(P(w))=w\}|\ \ \text{if}\ \ P(w)\in U(x_0).
$$

We take $U(x_0)$ as our first set with index set $I_{x_0}$, and sc$^+$-section structure $\mathfrak{s}_{x_0}={(s^{x_0}_i)}_{i\in I_{x_0}}$.
We define $\tau_{x_0,x_0}:\bm{U}(x_0,x_0)\rightarrow \mathfrak{I}(x_0,x_0)$ by 
$$
\tau_{x_0,x_0}(\Gamma(g,y))=[I_{x_0}\stackrel{Id}\twoheadleftarrow I_{x_0}\stackrel{g}\twoheadrightarrow I_{x_0}]\ \ \text{for}\ \ y\in U(x_0),\ g\in G_{x_0}.
$$
Next consider the points $x\in U(x_0)$. For each of them we pick an open $G_x$-invariant neighborhood
$U(x)$ satisfying $U(x)\subset U(x_0)$ and take $\mathfrak{s}_x:={(s^{x_0}_i|U(x))}_{i\in I_x}$ with $I_x:=I_{x_0}$.
We define a group homomorphism 
$$
\gamma_x:G_x\rightarrow G_{x_0}: g\rightarrow \gamma_x(g) 
$$
where $\gamma_x(g)$ is uniquely determined by the equality  $\Gamma^x(g,x)=\Gamma^{x_0}(\gamma_x(g),x)$.
Here $\Gamma^x$ and $\Gamma^{x_0}$ are associated to the natural representations on $U(x)$ and $U(x_0)$, respectively.
Then $G_x$ acts on $I_x$ by
$$
g(i) :=( \gamma_x(g))(i),\ i\in I_{x}=I_{x_0}.
$$
Given $x,x'\in U(x_0)$ the set $\bm{U}(x,x')$ is an open subset of $\bm{U}(x_0,x_0)$ and
we define $\tau_{x,x'}$ as the restriction of $\tau_{x_0,x_0}$.  

Among the  points $x\in X\setminus U(x_0)$ consider first the points $x$ which have an open neighborhood $O(x)$ 
with the property that $\Lambda(0_y)=1$ for $y\in O(x)$.  In this case fix an open $G_x$-invariant neighborhood
$U(x)$ with the properness property and take as sc$^+$-section structure $\mathfrak{s}_x$ the single zero-section
$s^x_\ast$ parameterized by $I_x=\{*\}$.  

If $x\in X\setminus U(x_0)$ but such a neighborhood $O(x)$ does not exist, we can pick a morphism
$\phi_x$ with $s(\phi_x)=x$ and $t(\phi_x)\in \cl_X(V(x_0))$. We find open neighborhoods $U(\phi_x)$, $U(x)$, and $U'(t(\phi_x))$
such that 
\begin{itemize}
\item $\sigma_x:=t\circ (s|U(\phi_x))^{-1}:U(x)\rightarrow U'(t(\phi_x))$ is an sc-diffeomorphism.
\item  $U'(t(\phi_x))\subset U(t(\phi_x))$.
\end{itemize}
We define $I_x:=I_{t(\phi_x)}$ and take as section structure $\mathfrak{s}_x$ the pull-back of the restriction 
of $\mathfrak{s}_{t(\phi_x)}$ to $U'(t(\phi_x))$.  We define a group isomorphism $\gamma_x:G_x\rightarrow G_{t(\phi_x)}$
by
$$
\gamma_x(g)  =  \phi_x \circ g\circ \phi^{-1}_x
$$
which we use to define the $G_x$-action on $I_x$ by 
$$
g(i):= (\gamma_x(g))(i),\ i\in I_x=I_{t(\phi_x)}=I_{x_0}.
$$
At this point we have defined a family $\mathfrak{U}={(U(x))}_{x\in X}$, an index set $I_x$ for every $x\in X$, together 
with group actions of $G_x$ on $I_x$.  In addition we have defined for every $x\in X$ a sc$^+$-section structure
$\mathfrak{s}_x$.   We also have defined the correspondences $\tau_{x,x'}$ for $x,x'\in U(x_0)$.
We need to define $\tau_{x,x'}$ for the remaining cases. Assume first that $x,x'$ are given outside of $U(x_0)$
such that the previously introduced $\phi_x$ and $\phi_{x'}$ with targets in $\cl_X(V(x_0))$ exist.
Using the sc-embeddings $\sigma_x$ and $\sigma_{x'}$ we can embed 
$\bm{U}(x,x')$ into $\bm{U}(t(\phi_x),t(\phi_{x'}))$ via the map
$$
\phi \rightarrow ( (s|U(\phi_{x'})^{-1}(t(\phi) ) \circ \phi  \circ {(( s|U(\phi_x))^{-1}(s(\phi)))}^{-1},
$$
and define 
$$
\tau_{x,x'}(\phi) := \tau_{t(\phi_x),t(\phi_{x'})} \left(( (s|U(\phi_{x'})^{-1}(t(\phi) ) \circ \phi  \circ {(( s|U(\phi_x))^{-1}(s(\phi)))}^{-1}\right).
$$
Assume next that $|x|\not\in |U(x_0)|$ and $x'$ is arbitrary.
We define
$$
\tau_{x,x'}(\phi)= [\{\ast\}\twoheadleftarrow I_{x'} \stackrel{Id}\twoheadrightarrow I_{x'}].
$$
If $x $ is arbitrary and $x'\in X$ satisfies $|x'|\not\in |U(x_0)|$ we define 
$$
\tau_{x,x'}(\phi) =[I_x\stackrel{Id}\twoheadleftarrow I_x \twoheadrightarrow \{*\}].
$$
The latter two definitions are compatible if $|x|,|x'|\not\in |U(x_0)|$.  On verifies easily that the 
sc$^+$-section structures are compatible via the correspondences.
With this data $[\Lambda,\mathfrak{U},\mathfrak{S},\tau]$ is a structured version of $\Lambda$.
\qed \end{proof}

 The following two Theorems  together show that one can construct atomic $\ssc^+$-multisections 
possessing  special properties.
\begin{theorem}[Construction I]\index{T- Construction of $\Lambda$}\label{THMXX1354}
Let $(P:W\rightarrow X,\mu)$  be  a strong bundle over an ep-groupoid admitting smooth bump functions.  We assume that $x$ is a smooth object  and $U(x)$ is a $G_x$-invariant   open neighborhood having the properness property and admitting the natural
$G_x$-action.  Let $V(x)$ be an another open neighborhood such that $\cl_X(V(x))\subset U(x)$ and $V(x)$
is invariant under $G_x$. Then the following statements hold. 
\begin{itemize}
\item[{\em(1)}]\ If $w\in W$ is a smooth vector satisfying  $P(w)=x$,  
then there exists a finite family of $\ssc^+$-section $s_i$, $i\in I$, defined on $U(x)$ vanishing on $U(x)\setminus V(x)$
 and an action of $G_x$ on $I$ satisfying 
$s_{g(i)}(g\ast z) =\mu(\Gamma (g, z),s_i(z))$ for $z\in U(x)$,  and $s_{i_0}(x)=w$ for some $i_0\in I$.
\item[{\em(2)}]\ If $s$ is a $\ssc^+$-section near $x\in U(x)$  satisfying  $s(x)=0$,  then there exists a finite family of $\ssc^+$-sections $s_i$ for $i\in I$, supported  in $U(x)$, vanishing on $U(x)\setminus V(x)$, and an action of  $G_x$ on $I$ such that 
$s_{g(i)}(g\ast z) =\mu(\Gamma (g,z), s_i(z))$ for $z\in U(x)$ and  $i\in I$. 
Moreover,  
$s_{i_0}(y)=s(y)$ for $y$ near  $x$ for some $i_0\in I$.
\end{itemize}
\end{theorem}
\begin{proof}
The proof  is straightforward. In case (1) we fix an open neighborhood $V(x)$ with $\cl_X(V(x))\subset U(x)$
so that $V(x)$ is invariant under the $G_x$-action. We construct in local coordinates a $\ssc^+$-section $s$ satisfying  $s(x)=w$. We multiply 
$s$ with a  sc-smooth bump function $\beta$ supported in $V(x)$ and satisfying $\beta (x)=1$, so that the product $\beta\cdot s$ 
is a $\ssc^+$-section supported in $U(x)$. Recall such $\beta$ exists with arbitrarily small support.
Now,  applying  the $G_x$ action,  we define for every $g\in G_x$ the local $\ssc^+$-section $s_g: U(x)\to W$ by $s_{g}(g\ast z)=\mu(\Gamma (g, z), s_0(z))$ for $z\in U'(x)$, where $s_0=\beta\cdot s$.  Using the properties of the sc-diffeomorphism $\Gamma: G_x\times U(x)\to \{\phi \in {\bf X}\, \vert \, \text{$s(\phi)$ and $t(\phi)\in U(x)$}\}$ and of the structure map $\mu$  one sees that if $g=\text{id}\in G_x$ then $s_{\text{id}}(z)=s_0(z)$. Moreover, for every $g,h\in G_x$ the formula 
$s_{g\circ h}(g\ast z)=\mu (\Gamma (g, z), s_h(z))$ holds for every $z\in U(x)$. Consequently, introducing the finite index set $I=G_x$ and the group isomorphism $\tau: G_x\to \text{Per}(I)$ into the permutations of $I$, defined by $\tau (g)(h)=g\circ h$, we obtain the equation
$$
s_{\tau (g)(h)}(g\ast z)=\mu (\Gamma (g, z), s_h(z))
$$
for every $z\in U(x)$ and $g, h\in G_x$.
This way we obtain a called symmetric section structure near $x$ defined on $U(x)$, where the $s_g$ have support 
in $V(x)$ with $\cl_X(V(x))\subset U(x)$. \par

In case (2) we use a sufficiently concentrated  bump function and define $s_0=\beta\cdot s$.
\qed \end{proof}
The previous proposition tells us how to  obtain locally a $G_x$-invariant family of $\ssc^+$-section supported in a suitable
neighborhood $U(x)$ for which there exists $V(x)$ with $\cl_X(V(x))\subset U(x)$, so that the support is contained in $V(x)$.  Such families can be used to construct atomic $\ssc^+$-multisections as the next result illustrates. 
\begin{theorem}[Construction II]\index{T- Construction of $\Lambda$}\label{THMXX1355}
Let $(P:W\rightarrow X,\mu)$ be  a strong bundle over an ep-groupoid admitting smooth bump functions. 
We assume that $x$ is a smooth point and $U(x)$ an open neighborhood  which is invariant under the natural $G_x$-action having the properness property.  Further, we  let $s_i$, $i\in I$, be a finite number of $\ssc^+$-sections of the strong bundle $W\rightarrow X$ supported in $U(x)$, where we view $X$ just as a M-polyfold.
We assume in  addition that   $G_x$ acts on $I$ such  that $s_{g(i)}(g\ast z)=\mu(\Gamma (g, z),s_i(z))$ for all $z\in U(x)$ and $i \in I$.
Then there exists a uniquely determined atomic $\ssc^+$-multisection $\Lambda$ of  $P$ which is supported in the saturation of $U(x)$ and which over $U(x)$ can be written 
$$
\Lambda(w)=\frac{1}{|I|}\cdot |\{i\in I\ |\ s_i(P(w))=w\}|.
$$
Moreover $\Lambda$ is structurable.
\end{theorem}
\begin{proof}
Let $U(x)$ be an open neighborhood admitting the natural $G_x$-action, so that $t\colon s^{-1}(\cl_X(U(x)))\rightarrow X$ is proper.
We assume we are given a finite set $I$ with an action by $G_x$ and a family of sc$^+$-sections ${(s_i)}_{i\in I}$ of $W|U(x)$
having the following properties.
\begin{itemize}
\item[(i)]\ \ $s_{g(i)}(g\ast y)= \mu(\Gamma(g,y),s_i(y))$ for $y\in U(x)$.
\item[(ii)]\  \ There exists $V(x)$ with $\cl_X(V(x))\subset U(x)$ and invariant under $G_x$, such that
$s_i(y)=0$ for $y\in U(x)\setminus V(x)$. 
\end{itemize}
First we define $\Lambda: W\rightarrow {\mathbb Q}^+$. If $w\in W$ and there exists a morphism 
$\phi$ with $t(\phi)\in U(x)$ and $s(\phi)=P(w)$ we define
$$
\Lambda(w)=\frac{1}{|I|}\cdot |\{i\in I\ |\ s_i(P(w))=w\}|.
$$
If no such $\phi$ exists we define 
$$
\Lambda(w)=\left[\begin{array}{cc}
1&\ \ \text{if}\ w=0_{P(w)}\\
0&\ \ \text{otherwise.}
\end{array}\right.
$$
We need to show that $\Lambda$ can be represented locally by sc$^+$-smooth section structures.\\

In the case that there exist $\phi:P(w)\rightarrow y$ with $y\in U(x)$ we can use the associated local sc-diffeomorphism
$\wh{\phi}$ to pull-back the local section structure and it will represent $\Lambda$ on $W$ restricted to a suitable 
small neighborhood of $P(w)$.\par

In the case that there does not exist a morphism starting at $y=P(w)$ and ending in $U(x)$
we find an open $G_y$-invariant neighborhood $U(y)$ so that $|U(y)|\cap |V(x)|=\emptyset$.  In this case
we take as section structure the single zero section on $U(y)$ which we may assume to possess the properness property.\par

After this discussion we see that $\Lambda$ is an sc$^+$-multisection functor. We have already established
that they are structurable.

\qed \end{proof}
The basic existence result deduced  from the two previous results  is the following theorem which provides in addition  a control  of the size of the sections. 
\begin{theorem}\label{existence_lambda}\index{T- Controlled atomic sc$^+$-multisection}
Let $(P:W\rightarrow X,\mu)$ be a strong bundle over the ep-groupoid $X$ with paracompact orbit space $|X|$,  which admits sc-smooth bump functions.
Let  $N$ is a given auxiliary norm for $(P,\mu)$. {\em (}It  exists due to  the paracompactness of $|X|${\em)}. 
For every  smooth vector  $w\in W\setminus 0_W$ satisfying $P(w)=x$ and every  open neighborhood $U(x)$ around $x=P(w)$ equipped with the natural $G_x$-action, there exists an atomic $\ssc^+$-multisection  satisfying  the following properties, 
\begin{itemize}
\item[{\em(1)}]\ $\text{dom-supp}(\Lambda)\subset \pi^{-1}(U(x))$.
\item[{\em(2)}]\ $\Lambda(g\ast w)>0$ for all $g\in G_x$.
\end{itemize}
Moreover,  if $s$ is a $\ssc^+$-section defined in the  open neighborhood $U(x)$ and if $\varepsilon>0$, then there exists an atomic $\ssc^+$-multisection $\Lambda$ supported   in $\pi^{-1}(U(x))$ and  satisfying $N(\Lambda)(y)\leq N(s(x))+\varepsilon$ for all $y\in X$,
whose  local section structure near $x$ contains  the 
$\ssc^+$-section $s$ having a positive weight.
\end{theorem}

The importance of atomic $\ssc^+$-section functors comes from the following result.
It states that every structurable $\ssc^+$-multisection functor $\Lambda:W\rightarrow {\mathbb Q}^+$ associated to $P:W\rightarrow X$,
where $|X|$ is paracompact and $X$  admits sc-smooth partitions of unity, can be obtained as   a locally finite commensurable
sum of atomic $\ssc^+$-multisection functors.
\begin{theorem}\index{T- Atomic decomposition}\label{THMX1357}
Let $(P:W\rightarrow X,\mu)$ be a strong bundle over the ep-groupoid $X$ which has a paracompact orbit space $|X|$ and admits sc-smooth partitions of unity.
We assume that $\Lambda:W\rightarrow {\mathbb Q}^+$ is a  structurable $\ssc^+$-multisection functors. Then there exists a locally finite family of commensurable structured $\ssc^+$-multisection functors 
${([\Lambda_j,\mathfrak{U},\mathfrak{S}_j,\tau])}_{j\in J}$, each having an underlying atomic $\ssc^+$-multisection functor $\Lambda_j$, such that  
$$
[\Lambda,\mathfrak{U},\mathfrak{S},\tau] =\boxplus_{j\in J}[ {\Lambda}_j,\mathfrak{U},\mathfrak{S}_j,\tau].
$$
\end{theorem}
\begin{proof}
We take a structured lift  of $\Lambda$ represented by $(\Lambda,\mathfrak{U},\mathfrak{S},\tau)$.
Setting $\wt{U}(x):=\pi^{-1}(|U(x)|)$ for every $x\in X$, the family ${(\wt{U}(x))}_{x\in X}$ is an open covering of $X$ by saturated sets. By assumption the orbit space $\abs{X}$ is paracompact and $X$ admits sc-smooth partition of unity functors as an ep-groupoid. Consequently, there exists a locally finite refinement $(V_j)_{j\in J}$ of ${(\wt{U}(x))}_{x\in X}$ and,  for every $j\in J$,  a point $x_j\in X$ such that $\cl_{\abs{X}}(V_j)\subset \abs{U(x_j)}$, and, moreover, a sc-smooth partition  of unity $(\beta_j)_{j\in J}$ of functors $\beta_j: X\to [0,1]$ satisfying 
$\supp (\abs{\beta_j})\subset V_j$.

Fixing $j\in J$, we define the $\ssc^+$-section structure $\mathfrak{s}^{j,x}=(s^{j,x}_{i})_{i\in I_x}$ on every $U(x)$ by 
$$s^{j,x}_i(y)=\beta_j(y)\cdot s^x_i(y), $$
which in turn defines the $\ssc^+$-multisection functor  $\Lambda_j$ and the structured 
$$
[\Lambda_j,\mathfrak{U},\mathfrak{S_j},\tau].
$$ 
We claim that,  for every $j\in J$, the functor $\Lambda_j$ is atomic. Since 
$\supp \abs{\beta_j}\subset  V_j \subset |U(x_j)|$, it follows that $\wt{\beta}_j:=\beta_j\vert U(x_j): U(x_j)\rightarrow [0,1]$
vanishes outside of the set $\pi^{-1}(V_j) \cap U(x_j)$ whose closure in $X$ is contained in $U(x_j)$ in view of Lemma \ref{simple_lemma}.
Hence $\Lambda_j$ is atomic and by construction the family $(\Lambda_j)$ is locally finite. By construction
$$
[\Lambda,\mathfrak{U},\mathfrak{S},\tau] =\boxplus_{j\in J} [\Lambda_j,\mathfrak{U},\mathfrak{S}_j,\tau].
$$
The proof is complete.
\qed \end{proof}

\chapter{Extension of Sc\texorpdfstring{$^+$}{pl}-Multisections}\label{CHAPH14}

The main  result in this chapter  concerns an extension  of a structurable  $\ssc^+$-multi\-sec\-tion functor  defined over $\partial X$, i.e.
$\Lambda:W\vert \partial X\rightarrow {\mathbb Q}^+$, 
to  a structurable $\ssc^+$-multisection  functor $\Lambda':W\rightarrow {\mathbb Q}^+$ over all of $W$,   satisfying some additional properties:
$$
[\Lambda:W\vert \partial X\rightarrow {\mathbb Q}^+]\ \   \rightsquigarrow\ \ [\Lambda':W\rightarrow {\mathbb Q}^+].
$$
The additional properties are concerned with the size of the extension and its support.
First we shall provide the appropriate definitions and state the main result, and then we outline 
 the proof so that the reader will be able to follow the sizable construction.
\section{Definitions and Main Result}\label{SECR14.1}

In order to formulate the extension result we need some preparation. We assume that $X$ is a tame ep-groupoid with paracompact orbit space $|X|$ admitting sc-smooth partitions of unity. We note that in general the boundary $\partial X$ is not an ep-groupoid, however, it has enough structure to talk about smooth data on it. 

Recall that a tame M-polyfold  $M$ has for every point $x$ an arbitrarily small open neighborhood $V=V(x)$, so that there are
precisely $d(x)$-many connected components in $\{y\in V\ |\ d(y)=1\}$ having $x$ in its closure, and each such closure
is a sub-M-polyfold, see Corollary \ref{cor_2.42} and Section \ref{subsec_tame_m_polyfolds} for more detail. If we take such a neighborhood small enough we may assume that the properness property holds.
Of course, in addition, we may always assume that we have the natural $G_x$-action on $V(x)$.
Hence we make the following definition. 
\begin{definition}\label{DEF_good_open}\index{D- Good neighborhood of $x\in \partial X$}
Let $X$ be a tame  ep-groupoid and $x\in \partial X$. A {\bf good} open neighborhood $U(x)$ in $\partial X$
is a subset of $\partial X$ which is given as the intersection $U(x)=V(x)\cap \partial X$, where $V(x)$ has the following properties.
\begin{itemize}
\item[(1)]\  $V(x)$ is an open neighborhood of $x$ in $X$ and admits the natural $G_x$-action.
\item[(2)]\  $t:s^{-1}(\cl_X(V(x)))\rightarrow X$ is proper.
\item[(3)]\  There are $d_X(x)$-many  connected components $\theta_i$  in $\{y\in V(x)\ |\ d_X(y)=1\}$ so that their closures $\bar{\theta}_i$
in $V(x)$ contain $x$, and for $y\in V(x)$ it holds that
$$
d_X(y)=\sharp\{i\in \{1,..,d_X(x)\}\ |\ y\in \bar{\theta}_i\}.
$$
\end{itemize}
We shall refer to $V(x)$ as a {\bf good} open neighborhood of $x\in\partial X$ in $X$. Note that the $\bar{\theta}_j$ are tame M-polyfolds. \qed
\end{definition}
\begin{remark}
For the existence of such neighborhoods $V(x)$ see Proposition \ref{FACE_XXXX} and the discussion preceding it.
\end{remark}
Next we shall introduce a class of sc$^+$-multisections over $\partial X$.

\begin{definition}\label{symm}\index{D- Sc$^+$-multisection over $\partial X$}
Let $(P:W\rightarrow X,\mu)$ be a strong bundle over the tame ep-groupoid $X$ with paracompact $|X|$. A functor 
$\Lambda:W|\partial X\rightarrow {\mathbb Q}^+$ is an {\bf sc$^+$-multisection (functor) over $\partial X$} if for every $x\in \partial X$
 there exists 
a good neighborhood $U(x)$ in $\partial X$, an index set $I_x$ with $G_x$-action and a family of section ${(s^x_i)}_{i\in I_x}$
defined for $W|U(x)$ so that 
\begin{itemize}
\item[(1)]\ The restriction of $s_i^x$ to every $\bar{\theta}_j$ is an sc$^+$-section of $W|\bar{\theta_j}$.
\item[(2)]\  $s^x_{g(i)}(g\ast y)=\mu(\Gamma(g,y),s^x_i(y))$ for all $y\in U(x)$.
\item[(3)]\  $\Lambda(w)=\frac{1}{|I_x|}\cdot |\{i\in I_x\ |\ s_i^x(P(w))=w\}|$\ \ for $w\in W$ with $P(w)\in U(x)$.
\end{itemize}
We call $\mathfrak{s}_x={(s^x_i)}_{i\in I_x}$ a {\bf symmetric sc$^+$-section structure}\index{D- Symmetric sc$^+$-section structure} for $\Lambda$ near $x$ (in $\partial X$).
\qed
\end{definition}
\begin{remark}[Extensions]\index{R- Remark concerning extensions}
Let for $x\in \partial X$ be the data as in Definition  \ref{symm} and assume that $X$ admits sc-smooth partitions of unity.
Then it is an easy exercise to show that there exists an open neighborhood $V(x)$ of $x$ in $X$
being invariant under the natural $G_x$-action and sc$^+$-sections ${(\overline{s}_i^x)}_{ i\in I}$ of $W|V(x)$
having the following properties:
\begin{itemize}
\item[(a)]\ $V(x)\cap \partial X=U(x)$.
\item[(b)]\ $\overline{s}_i|U(x)= s_i$ for all $i\in I$.
\item[(c)]\  $\overline{s}^x_{g(i)}(g\ast y) =\mu(\Gamma(g,y),\overline{s}^x_i(y))$ for $y\in V(x)$ and $g\in G$.
\end{itemize}
To see this, we first construct in a straight forward many $V(x)$. 
Pick $i\in I$ and  employ Proposition \ref{hucky}  and a partition of unity argument
to construct $\overline{s}_i^x$.  Then we use the formula (c) to define $s^x_{g(i)}$ for every $g\in G$.
In a next step we pick $i\in I\setminus G(I)$ and repeat the previous step. We are done after a finite number of iterations.
\qed
\end{remark}

Next we shall introduce the notion of good systems of neighborhoods for $\partial X$. 

\begin{definition}\index{D- Good system of neighborhoods ($\partial X$-case}
Let ${(U(x))}_{x\in\partial X}$ be a collection of open neighborhoods $U(x)\subset \partial X$ of points $x\in \partial X$
in $\partial X$.  Then ${(U(x))}_{x\in\partial X}$ is called a {\bf good system} of neighborhoods for $\partial X$ if the following conditions are satisfied:
\begin{itemize}
\item[(1)]\  The $U(x)$ are good neighborhoods in the sense of Definition \ref{DEF_good_open}
\item[(2)]\  For two points $x,x'\in \partial X$,  the  subset $\bm{U}(x,x')=\{\psi\in \bm{X}\ |\ s(\psi)\in U(x),\ t(\psi)\in U(x')\}$ of $\bm{X}$ has  the following property.  There exist good open neighborhoods $V(x)$ and $V(x')$ of $x$ and $x'$ in $X$, see
Definition \ref{DEF_good_open}, restricting to  $\partial X$
as $U(x)$ and $U(x')$ so that
for every connected component $\Sigma\subset \bm{V}(x,x')$ the source map $s:\Sigma\rightarrow V(x)$ and the target map
$t:\Sigma\rightarrow V(x')$ are sc-diffeomorphisms onto open subsets. 
\end{itemize}
\qed
\end{definition}
Given an open covering ${(O(x))}_{x\in\partial X}$ of $\partial X$ ($O(x)$ open in $\partial X$), we can find a good system of open neighborhoods
${(U(x))}_{x\in\partial X}$, $U(x)\subset O(x)$,  by employing Theorem \ref{xxxx-structure} as outlined in the following proof.
\begin{theorem}\label{THMX14.1.15}
Let $X$ be an ep-groupoid with paracompact orbit space $|\partial X|$ and ${(O(x))}_{x\in\partial X}$ an  open cover
of $\partial X$ by open neighborhoods of its points (Here $O(x)$ is open in $\partial X$).  Then there exists a good system of neighborhoods  $ {(U(x))}_{x\in\partial X}$
for $\partial X$
with $U(x)\subset O(x)$ for all $x\in \partial X$.
\end{theorem}
\begin{proof}
Extend the elements of the collection ${(O(x))}_{x\in\partial X}$ to open neighborhoods $\wt{O}(x)$ of $x\in \partial X$
with $\wt{O}(x)\cap \partial X=O(x)$.  By adding open neighborhoods $\wt{O}(x)$ for $x\in X\setminus \partial X$
we obtain the collection  ${(\wt{O}(x))}_{x\in X}$. In view of Theorem \ref{xxxx-structure}
we find a good system of neighborhoods ${(V(x))}_{x\in X}$ with $V(x)\subset \wt{O}(x)$ for $x\in X$.
Then define for $x\in \partial X$ the set $U(x)$ by $U(x):=\partial X\cap V(x)$.
\qed \end{proof}

Now we are in the position to define structured $\ssc^+$-multisection functors for $W\vert \partial X$.
\begin{definition}\label{corr-xxxx}\index{D- Structured $\ssc^+$-multisection ($\partial X$-case)}
Let $(P:W\rightarrow X,\mu)$ be a strong bundle over a tame ep-groupoid with paracompact orbit space $|X|$. A {\bf structured} $\ssc^+$-multisection  over $\partial X$ is an equivalence class of tuples $[\Lambda,\mathfrak{U},\mathfrak{S},\tau]$ in which  $\Lambda:W\vert \partial X\rightarrow {\mathbb Q}^+$ is a $\ssc^+$-multi\-sec\-tion functor, and
\begin{itemize}
\item[(1)]\  $\mathfrak{U}={(U(x))}_{x\in \partial X}$ is a good system of open neighborhoods for $\partial X$,
\item[(2)]\   $\tau={(\tau_{x,x'})}_{(x,x')\in \partial X\times \partial X}$ is a collection of correspondences
$$
\tau_{x,x'}:\bm{U}(x,x')\rightarrow \mathfrak{J}(x,x'),
$$
Every $\tau_{x,x'}$ only takes finitely many values. 
We assume in addition that 
\begin{itemize}
\item[(a)]\ $\tau_{x',x}(\phi^{-1})=(\tau_{x,x'}(\phi))^\sharp$.
\item[(b)]\  $\tau_{x,x}(\Gamma(g,z))=[I_x\stackrel{Id}{\twoheadleftarrow} I_x\stackrel{g}{\twoheadrightarrow} I_x]$.
\item[(c)]\ $g'\ast\tau_{x,x'}(\phi)\ast g=\tau_{x,x'}(g'\ast\phi\ast g)$ for $g'\in G_{x'}$, $g\in G_x$,
and $\phi\in \bm{U}(x,x')$.
\end{itemize}
\item[(3)]\  $\mathfrak{S}= {(\mathfrak{s}_x)}_{x\in \partial X} $, where $\mathfrak{s}_x$  is a symmetric $\ssc^+$-section structure on $U(x)$ indexed by $I_x$ and representing $\Lambda$:
$$
\Lambda(w)=\frac{1}{|I_x|}\cdot |\{i\in I_x\ |\ s^x_i(P(w))=w\}|\ \ \ \text{if}\ \ P(w)\in U(x).
$$
\item[(4)]\ For $x,x'\in \partial X$ the compatibility condition $s^{x'}_{b(k)}(t(\phi))=\mu(\phi,s^x_{a(k)}(s(\phi)))$, $k\in I$,
for all $\phi\in \bm{U}(x,x')$ holds. Here $a,b$ is the usual data associated to $\tau_{x,x'}(\phi)$, i.e. 
$$
\tau_{x,x'}(\phi) = [I_x\stackrel{a}\twoheadleftarrow I\stackrel{b}\twoheadrightarrow I_{x'}].
$$
\end{itemize}
A $\ssc^+$-multisection functor $\Lambda:W|\partial X\rightarrow {\mathbb Q}^+$ is called  {\bf structurable}\index{D- Structurable $\ssc^+$-multisection ($\partial X$-case)} if there exists 
a structured $\ssc^+$-multisection  of the form $[\Lambda,\mathfrak{U},\mathfrak{S},\tau]$, where the notion
of equivalence is similar to the $X$-case, see Definition \ref{DEF13.3.8} and before.
\qed
\end{definition}
For a latter discussion in Section \ref{SECTX16.5}    we note the following result.
\begin{proposition}\label{prop14.1}
Let $(P:W\rightarrow X,\mu)$ and $(P':W'\rightarrow X',\mu')$ be strong bundles over tame ep-groupoids with paracompact orbit spaces
and assume that $\overline{\mathfrak{f}}:W\rightarrow W'$ is a generalized strong bundle isomorphism covering $\mathfrak{f}:X\rightarrow X'$.
Further let $\Lambda:W|\partial X\rightarrow {\mathbb Q}^+$ and $\Lambda':W'|\partial X'\rightarrow {\mathbb Q}^+$ be structurable 
sc$^+$-multisection functors in the sense of Definition \ref{corr-xxxx}.  Then $\bar{\mathfrak{f}}^\ast\Lambda'$ defines a structurable 
sc$^+$-multisection functor defined on $W|\partial X$ and $\bar{\mathfrak{f}}_\ast \Lambda$ defines a structurable sc$^+$-multisection functor
on $W'|\partial X'$.
\qed
\end{proposition}
\begin{proof}
Given $\overline{\mathfrak{f}}:W\rightarrow W'$ we take a representative which gives the commutative diagram
$$
\begin{CD}
W @< \Phi << V @> \Psi>> W'\\
@V P VV  @V Q VV @V P' VV\\
X@<E<<  A@> F>>X, 
\end{CD}
$$
where $Q: V\rightarrow A$ is a strong bundle over a tame ep-groupoid, and the horizontal arrows
give strong bundle equivalences. 

Let $\Lambda:W\vert \partial X\rightarrow {\mathbb Q}^+$ be a structurable sc$^+$-multisection functor. 
By an argument  similar to that used in the proof of Theorem \ref{pullbackX} it is a straight forward  exercise to show that $E^\ast\Lambda$ is a  structurable sc$^+$-multisection functor.   Then employing an argument as
in the proof of Theorem \ref{pushforwardX} we see that $\overline{\mathfrak{f}}_\ast \Lambda=F_\ast E^\ast \Lambda 
$ is structurable as well.  The argument for the pullback via $\overline{\mathfrak{f}}$ is similar.
\qed
\end{proof}

One of the  main results in this chapter is the following extension theorem.
\begin{theorem}[$\partial$-Extension Theorem]\label{p-main-p}\index{T- $\partial$-Extension theorem}
Assume  $(P:W\rightarrow X,\mu)$ is   a strong bundle over a tame  ep-groupoid $X$ with paracompact orbit space $|X|$ and $X$, as the M-polyfold,   admits  sc-smooth partitions of unity. 
Let 
$$
\Lambda:W|\partial X\rightarrow {\mathbb Q}^+
$$
be  a structurable $\ssc^+$-multisection functor defined over the boundary and let 
$N:W\rightarrow [0,\infty]$ be an auxiliary norm. 
If $\wt{U}$ is a saturated open neighborhood 
of $\text{dom-supp}(\Lambda)$ in $X$ and $f:X\rightarrow [0,\infty)$ is a continuous functor supported in $\wt{U}$ and satisfying
$N(\Lambda)(x)<f(x)$ for all $x\in \text{dom-supp}(\Lambda)$, then there exists  a structurable $\ssc^+$-multisection functor $\Lambda':W\rightarrow {\mathbb Q}^+$
having  the following properties:
\begin{itemize}
\item[\em{(1)}]\ $N(\Lambda')(x)\leq f(x)$ for all $x\in X$.
\item[\em{(2)}]\ $\text{dom-supp}(\Lambda')\subset \wt{U}$.
\item[\em{(3)}]\ $\Lambda'|(W|\partial X)=\Lambda$.
\end{itemize}
\qed
\end{theorem}
\begin{remark}
The proof of this theorem requires some preparations and is rather involved.
Extension theorems for sc$^+$-section functors are, given sc-smooth partitions of unity,
not very difficult. However, the complexity in the case of multisections is much higher.
\qed
\end{remark}
\section{Good Structured Version of \texorpdfstring{$\Lambda$}{Lam}}
We assume that $\Lambda$  is given as described in Theorem \ref{p-main-p}.
Since   $\Lambda: W|\partial X\rightarrow {\mathbb Q}^+$ is structurable we can take a structured lift
represented by
\begin{eqnarray}\label{display7}
 (\Lambda,\mathfrak{U}',\mathfrak{S}',\tau')
\end{eqnarray}
 In detail this gives 
us the following list of data with certain properties, which will be used in the following.
\begin{description}
\item[{\bf (P1)}]\ \ For every $x\in \partial X$ a good neighborhood $U'(x)$ of $x$ in $\partial X$ which by definition has the following properties
\begin{itemize}
\item[(i)]\ \ $U'(x)=V(x)\cap \partial X$, where $V(x)$ is an open neighborhood of $x$ in $X$ admitting the natural $G_x$-action.
\item[(ii)]\ \ $t:s^{-1}(\cl_X(V(x)))\rightarrow X$ is proper.
\item[(iii)]\ \ There are $d_X(x)$-many connected components $\theta_i$ in 
$$
\{y\in V(x)\ |\ d_X(y)=1\}\subset \partial X,
$$
 so that their closures $\bar{\theta}_i$ in $V(x)$ are tame M-polyfolds containing $x$ and 
$$
d_X(y)=\sharp\{i\in \{1,...,d_X(x)\}\ |\ y\in\bar{\theta}_i\}.
$$
\end{itemize}
\item[{\bf (P2)}]\ \  For every $x\in \partial X$ a finite index set $I_x$ and a symmetric sc$^+$-section structure
$\mathfrak{s}_x'$ given by $s_i'^x$, $i\in I_x$, on $U'(x)$, which represents $\Lambda$ over $U'(x)$,
see Definition \ref{symm}. We abbreviate $\mathfrak{J}={(I_x)}_{ x\in \partial X}$.
\item[{\bf (P3)}]\ \  For $x,x'\in \partial X$ correspondences $\tau_{x,x'}':\bm{U}(x,x')\rightarrow \mathfrak{J}(x,x')$ with the properties listed in Definition \ref{corr-xxxx}. We denote by $v'(x,x')$
the image of $\tau_{x,x'}'$.
\end{description}
We shall modify the given data 
later on. The index sets in $\mathfrak{J}$  will stay unchanged.
In general $\mathfrak{U}'$ will be modified to achieve certain goals.
We denote by $\bm{U}'(x,x')$ the subset of $\partial \bm{X}$ defined by
$$
\bm{U}'(x,x')=\{\phi\in \partial\bm{X} \ |\ s(\phi)\in U'(x),\ t(\phi)\in U'(x')\}.
$$
The first important fact is that 
$$
\tau_{x,x'}':\bm{U}'(x,x')\rightarrow \mathfrak{J}(x,x')
$$
only takes finitely many values.  This defines a subset $v'(x,x')$ of $\mathfrak{J}(x,x')$.  If we replace the $U'(x)$ by smaller $U(x)$
the associated $\tau_{x,x'}$ obtained by restriction takes  values defining a set $v(x,x')\subset v'(x,x')$.

In a first step replace every $U'(x)$ by a smaller $U(x)$ with the same properties, but in addition satisfying
$$
\cl_{|\partial X|}(|U(x)|)\subset |U'(x)|.
$$
This implies $\cl_{X}(U(x))\subset U'(x)$ and the following holds true.
\begin{lemma}
Defining $\bm{U}(x,x'):=\{\phi\in \partial\bm{X}\ |\ s(\phi)\in U(x),\ t(\phi)\in U(x)\}$ the closures of the preimages in $\bm{X}$ under $\tau_{x,x'}|\bm{U}(x,x')$  of two different $[d]\neq [d']$  are disjoint.
\end{lemma}
\begin{proof}
Let  $(\phi_k),(\psi_k)$ be sequences in $\bm{U}(x,x')$ satisfying 
$$
\tau_{x,x'}'(\phi_k)=[d]\ \ \text{and}\ \  \tau_{x,x'}'(\psi_k)=[d'].
$$
Suppose further that $\phi_k\rightarrow \phi_0$ and $\psi_k\rightarrow \phi_0$ in $\bm{X}$ (which implies $\phi_0\in\partial\bm{X}$).
Then $\phi_0\in {\bf U}'(x,x')$ and define  $[d_0]:=\tau_{x,x'}'(\phi_0)$. Since $\tau_{x,x'}':\bm{U}'(x,x')\rightarrow \mathfrak{J}(x,x')$ is locally constant
it follows for large enough $k$
$$
[d]=\tau_{x,x'}'(\phi_k)=[d_0]=\tau_{x,x'}'(\psi_k)=[d'].
$$
This completes the proof.
\qed \end{proof}

In view of the lemma, replacing $\mathfrak{U}'$ by a good system of neighborhoods with smaller sets, we may assume without loss of generality  that $\mathfrak{U}'$ has this property already. It will persist when we replace
these sets by even smaller ones. Hence we have proved.
\begin{proposition}[Good Structured Version]\label{basic-data}\index{D- Good structured version}
Under the assumption of Theorem \ref{p-main-p} we can pick a structured lift $[\Lambda,\mathfrak{U}',\mathfrak{S}',\tau']$
and representative $(\Lambda,\mathfrak{U}',\mathfrak{S}',\tau')$ 
such that the following  property, besides {\bf (P1)-(P3)},  holds in addition.
\begin{itemize}
\item[{\bf (P4)}]\ \ \ For different $[d]\neq [d']$ in $v'(x,x')$ the closures of $\tau_{x,x'}'^{-1}([d])$ and $\tau_{x,x'}'^{-1}([d'])$ in $\bm{X}$ are disjoint.
\end{itemize}
\qed
\end{proposition}

From now on we assume that $(\Lambda,\mathfrak{U}',\mathfrak{S}',\tau')$ 
has, besides {\bf (P1), (P2), (P3)}  the additional property {\bf (P4)}.

\section{Extension of Correspondences}
Under the assumptions of Theorem \ref{p-main-p} we have picked a representative 
$$
(\Lambda,\mathfrak{U}',\mathfrak{S}',\tau')
$$
 of a suitable structured version
  satisfying 
{\bf (P1), (P2), (P3), (P4)}. 
In particular we have for every pair $x,x'\in\partial X$ the map
\begin{eqnarray}\label{EQNXX14}
\tau_{x,x'}':\bm{U}'(x,x')\rightarrow \mathfrak{J}(x,x'),
\end{eqnarray}
which defines a correspondence between the section structures. In this section we 
shall show first that there is a canonical extension of every $\tau_{x,x'}'$ 
to an open neighborhood $\wh{O}_{x,x'}$ of $\bm{U}'(x,x')\subset \partial\bm{X}$
in $\bm{X}$.

There exists for every $x\in\partial X$ an open neighborhood $\wt{U}'(x)\subset X$
which has the following property
\begin{itemize}
\item[{\bf (P5)}]\ \ \ \ (i)  $\wt{U}'(x)$ admits the natural action by $G_x$ and has the properness property.
\item[]\ \ \ (ii)  $\wt{U}'(x)\cap \partial X=U'(x)$.
\end{itemize}

It is clear that ${(\wt{U}'(x))}_{x\in\partial X}$ is a covering of $\partial X$ by
sets open in $X$. We define
$$
\wt{\bm{ U}}'(x,x'):=\{\phi\in \bm{X}\ |\ s(\phi)\in \wt{U}'(x),\ t(\phi)\in \wt{U}'(x')\}.
$$
The following proposition is crucial. It is the first step aiming at  extending the maps $\tau_{x,x'}':\bm{U}'(x,x')\rightarrow\mathfrak{J}(x,x')$ to open subsets in $\bm{X}$. The sets $\wt{O}_{x,x'}$ we are going to construct next are an intermediate step
in constructing the sets $\wh{O}_{x,x'}$ previously mentioned (\ref{EQNXX14}).
\begin{lemma}\label{pp-prop-pp}
Under the assumptions of Theorem \ref{p-main-p} consider a representative $(\Lambda,\mathfrak{U}',\mathfrak{S}',\tau')$  of a 
structured lift of $\Lambda$ having the properties {\bf (P1)--(P4)}.  Also let 
${(\wt{U}'(x))}_{x\in\partial X}$ be a covering of $\partial X$ by
sets open in $X$ satisfying {\bf (P5)}.
Then there exists for every $[d]\in v'(x,x')$ an open neighborhood $\wt{O}_{x,x',[d]}$ of $\tau_{x,x'}'^{-1}([d])$ in $\bm{X}$ satisfying the following.
\begin{itemize}
\item[{\em(1)}]\  $\wt{O}_{x,x',[d]}\subset \wt{\bm{U}}'(x,x')$ for all $[d]\in v'(x,x')$.
\item[{\em(2)}]\  $\cl_{\bm{X}}(\wt{O}_{x,x',[d]})\cap \cl_{\bm{X}}(\wt{O}_{x,x',[d']})=\emptyset$ if $[d],[d']\in v'(x,x')$ satisfy $[d]\neq [d']$.
\item[{\em(3)}]\  $h\ast \wt{O}_{x,x',[d]}\ast g=\wt{O}_{x,x',h\ast [d]\ast g}$ for $g\in G_x$ and $h\in G_{x'}$.
\item[{\em(4)}]\  Defining  $\wt{O}_{x,x'}=\bigcup_{[d]\in v'(x,x')} \wt{O}_{x,x',[d]}$ we have  $\wt{O}_{x,x'}\cap\partial\bm{ X}=\bm{U}'(x,x')$.
\end{itemize}
\end{lemma}
\begin{proof}
The sets $\tau_{x,x'}'^{-1}([d])$ for $[d]\in v'(x,x')$ have  mutually disjoint closures in $\bm{X}$, which
lie in $\partial\bm{X}$. Recalling that $\wt{U}'(x)\cap\partial X=U'(x)$ for $x\in \partial X$,  and that
$\bm{U}'(x,x')$ decomposes into the disjoint union 
$$
\bm{U}'(x,x') =\bigcup_{[d]\in v'(x,x')} \tau_{x,x'}'^{-1}([d]),
$$
we note that each $\tau_{x,x'}'^{-1}([d])$ is open and closed in $\bm{U}'(x,x')$.
Of course, $\cl_{\bm{X}}(\tau_{x,x'}'^{-1}([d]))$ is usually larger than $\tau_{x,x'}'^{-1}([d])$.

We can find open neighborhoods 
$\wt{O}'_{x,x',[d]}$ of $\tau_{x,x'}'^{-1}([d])$  in $\bm{X}$ which are contained in $\wt{\bf U}'(x,x')$ and have mutually disjoint closures in $\bm{X}$. 
These sets therefore  satisfy (1) and (2). We also have by construction
$$
\tau'^{-1}_{x,x'}([d])\subset \wt{O}'_{x,x',[d]}\cap \partial \bm{X}.
$$
On the other hand, for  $\phi\in \wt{O}'_{x,x',[d]}\cap \partial\bm{X}$ it follows that $\phi\in \bm{U}'(x,x')$ 
and $\tau_{x,x'}'(\phi)=[d]$ and therefore
$$
\tau_{x,x'}'^{-1}([d]) = \wt{O}'_{x,x',[d]}\cap\partial\bm{X}.
$$
From this and the properties of $\tau'_{x,x'}$ we deduce that for $(h,[d'],g)\in G_{x'}\times v'(x,x')\times G_x$ with $h\ast [d']\ast g=[d]$
it holds that $h\ast \wt{O}'_{x,x',[d']}\ast g$ is an open neighborhood of $\tau_{x,x'}^{-1}([d])$.
Define 
$$
\wt{O}_{x,x',[d]}:=\bigcap_{(h,[d'],g)\in G_{x'}\times v'(x,x')\times G_x,\ h\ast[d']\ast g=[d]} h\ast \wt{O}'_{x,x',[d']}\ast g
$$
Since $\wt{\bm{U}}'(x,x')$ is invariant under the bi-action this defines an open neighborhood of $\tau'^{-1}_{x,x'}([d])$ in ${\bf X}$ contained in  $\wt{\bm{U}}'(x,x')$. Again it holds
$$
\wt{O}_{x,x',[d]} \cap\partial{\bf X} =\tau'^{-1}_{x,x'}([d])\ \text{for all}\ [d]\in v'(x,x').
$$
Then $\wt{O}_{x,x',[d]}$ satisfies (1) and (2) and from the construction it follows that $h\ast\wt{O}_{x,x',[d]}\ast g= \wt{O}_{x,x',h\ast[d]\ast g}$,
so that (3) holds. Finally $\wt{O}_{x,x'}$ has obviously property (4).
\qed \end{proof}
Now we can define the sets $\wh{O}_{x,x'}$ we are interested in.
We denote by $\wh{O}_{x,x'}$\index{$\wh{O}_{x,x'}$} the union of the connected components in $\wt{O}_{x,x'}$ intersecting $\partial\bm{X}$.
This set is bi-invariant and 
\begin{eqnarray}\label{OEQNXXX}
\bm{U}'(x,x')\subset \wh{O}_{x,x'}\subset \wt{O}_{x,x'}\subset  \wt{\bm{U}}'(x,x').
\end{eqnarray}
Each  connected component in $\wh{O}_{x,x'}$ intersects $\partial\bm{X}$, which allows
us to take for $\phi\in \wh{O}_{x,x'}$ an element $\psi\in \wh{O}_{x,x'}\cap \partial\bm{X}\subset \bm{U}'(x,x')$
in the same connected component.
 We  assign to the element
$\phi$ the value $[d]_\phi=\tau_{x,x'}'(\psi)\in v'(x,x')$.  
\begin{lemma}
The map $\phi\rightarrow [d]_\phi$ is well-defined. 
\end{lemma}
\begin{proof}
Assume that $\Sigma$ is a connected component in $\wh{O}_{x,x'}$ and $\psi,\psi'\in\Sigma\cap \partial{\bm{X}}$.
Then there exist $[d],[d']$ such that $\phi\in \wt{O}_{x,x',[d]}$ and $\phi'\in \wt{O}_{x,x',[d']}$. If $[d]\neq [d']$
we know that the closures of these sets are disjoint implying that $\psi$ and $\psi'$ cannot be on the same connected component.
\qed \end{proof}
In view of this lemma we can define 
$$
\wh{O}_{x,x',[d]}=\{\phi\in \wh{O}_{x,x'},\ [d]_\phi=[d]\}, \index{$\wh{O}_{x,x',[d]}$}
$$
and the  set $\wh{O}_{x,x'}$ decomposes as disjoint union 
\begin{eqnarray}\label{req1}
\wh{O}_{x,x'}=\bigcup_{[d]\in v'(x,x')} \wh{O}_{x,x',[d]}.
\end{eqnarray}
We can extend $\tau_{x,x'}':{\bf U}'(x,x')\rightarrow \mathfrak{J}(x,x')$ in a unique way to a locally constant map
$$
\wh{\tau}_{x,x'}:\wh{O}_{x,x'}\rightarrow \mathfrak{J}(x,x'),
$$
which restricted to ${\bf U}'(x,x')$ is $\tau_{x,x'}'$. Namely we define 
\begin{eqnarray}\label{req2}
 \wh{\tau}_{x,x'}(\phi) =[d]\ \text{if}\  \phi\in \wh{O}_{x,x',[d]}.
\end{eqnarray}
We summarize the construction so far as follows.
\begin{proposition}\label{CORR13.6.12}
Under the assumptions of Theorem \ref{p-main-p} consider a representative $(\Lambda,\mathfrak{U}',\mathfrak{S}',\tau')$  of a 
structured lift of $\Lambda$ having the properties {\bf (P1)--(P4)}.  Also let 
${(\wt{U}'(x))}_{x\in\partial X}$ be a covering of $\partial X$ by
sets open in $X$ satisfying {\bf (P5)}. Then there exist for $x,x'\in \partial X$ open subsets $\wh{O}_{x,x'}$ of $\bm{X}$
satisfying
\begin{itemize}
\item[{\em(1)}]\ $\bm{U}'(x,x')\subset \wh{O}_{x,x'}\subset \wt{\bm{U}}'(x,x')$.
\item[{\em(2)}]\  The set $\wh{O}_{x,x'}$ is invariant under the bi-action of $(G_{x'},G_x)$.
\item{\em[(3)}]\ $\tau'_{x,x'}:\bm{U}'(x,x')\rightarrow \mathfrak{J}(x,x')$ has a unique locally constant extension 
$$
\wh{\tau}_{x,x'}:\wh{O}_{x,x'}\rightarrow \mathfrak{J}(x,x')
$$
taking the same values as $\tau_{x,x'}'$.
This extension has the additional properties
$$
\wh{\tau}_{x,x'}(g'\ast\phi\ast g) = g'\ast\wh{\tau}_{x,x'}(\phi)\ast g\ \text{for}\ g\in G_x,\ g'\in G_{x'},\ \phi\in \wh{O}_{x,x'},
$$
and
$$
\wh{\tau}_{x,x}(\Gamma(g,y)) =\tau'_{x,x}(\Gamma(g,z))
$$
for $z\in U'(x)$ and $y\in \wt{U}'(x)$.
\end{itemize}
Importantly $\wh{O}_{x,x'}$ can be written as a disjoint union 
\begin{eqnarray*}
\wh{O}_{x,x'} =\bigcup_{[d]\in v'(x,x')} \wh{O}_{x,x',[d]},
\end{eqnarray*}
with
\begin{eqnarray*}
&\wh{O}_{x,x',[d]}\cap\partial \bm{X} =\tau_{x,x'}'^{-1}([d])\subset \bm{U}'(x,x')&\\
&\tau_{x,x'}'^{-1}([d'])\cap \wh{O}_{x,x',[d]}=\emptyset\ \ \text{for}\  [d']\neq [d].&\nonumber
\end{eqnarray*}
\end{proposition}
\begin{proof}
Follows from Lemma \ref{pp-prop-pp} and the discussion thereafter.
\qed \end{proof}
\begin{remark}\label{remark14.3.4}
At this point we have constructed an extension $\wh{\tau}_{x,x'}:\wh{O}_{x,x'}\rightarrow \mathfrak{J}(x,x')$ for every choice of $x,x'\in\partial X$. Note that $\wh{O}_{x,x'}$ is not(!)  of the form 
$\wt{\bm{U}}(x,x')=\{\phi\in \bm{X}\ |\ s(\phi)\in \wt{U}(x),\ t(\phi)\in \wt{U}(x')\}$
for suitable open sets $\wt{U}(x)$ and $\wt{U}(x')$. This means that at this point $\wh{\tau}_{x,x'}$ has not the form required from a correspondence.  When we finally succeed in constructing a structurable extension $\Lambda'$ of $\Lambda$ we shall obtain a correspondence of the usual form between the local section structures. However, in general, the restriction of the data defining the structure will have little in common
with the original data coming from a structured choice over $\partial X$.  This being said they will however be related on some level. \qed
\end{remark}

Next we pick  for every $x\in \partial X$ an open neighborhood $\wt{U}(x)$ of $x$ in $X$ with the following property.
\begin{itemize}
\item[{\bf (P6)}] \ \ \ (i)  \ $\wt{U}(x)$ admits the natural $G_x$-action
and has the properness property.
\item[]\ \ \ (ii) \ $\cl_{|X|}(|\wt{U}(x)|)\subset |\wt{U}'(x)|\ \text{for} \ x\in \partial X.$
\end{itemize}

In particular it holds that
$$
\cl_X(\wt{U}(x))\subset \wt{U}'(x)\ \text{for}\ x\in \partial X.
$$
We shall write 
 $\wt{\bm{U}}(x,x')$ for $\{\phi\in\bm{X}\ |\ s(\phi)\in \wt{U}(x),\ t(\phi)\in \wt{U}(x')\}$.
 We define 
\begin{eqnarray}
&U(x):=\wt{U}(x)\cap \partial X&\\
&\bm{U}(x,x')=\{\phi\in \bm{X}\ |\ s(\phi)\in U(x),\ t(\phi)\in U(x')\}.\nonumber
\end{eqnarray}
For the following we shall rely on the data already constructed, which satisfies {\bf (P1)---(P6)}.
Together with (\ref{OEQNXXX}) the following holds for the various sets we currently consider
\begin{eqnarray}
&\bm{U}(x,x')\subset cl_{\bm{X}}(\bm{U}(x,x'))\subset \bm{U}'(x,x')\subset \wh{O}_{x,x'}\subset \wt{O}_{x,x'}\subset  \wt{\bm{U}}'(x,x')&\\
&\bm{\wt{U}}(x,x')\subset \cl_{\bm{X}}(\bm{\wt{U}}(x,x'))\subset \bm{\wt{U}}'(x,x').&\nonumber
\end{eqnarray}

The next lemma is the place where we use the paracompactness of $|X|$.
\begin{lemma}\label{Lemma2-62}
Under the assumption that $X$ is an ep-groupoid with paracompact $|X|$ 
there exists a subset $A$ of $\partial X$ and for every $a\in A$  a nonempty open subset $\wt{V}_a$ of $\wt{U}(a)$
with the following properties.
\begin{itemize}
\item[\em{(1)}]\ \  $\cl_{|X|}(|\wt{V}_a|)\subset |\wt{U}(a)|$ for all $a\in A$.
\item[\em{(2)}]\ \  $\wt{V}_a$ is invariant under the natural action of $G_a$ on $\wt{U}(a)$.
\item[\em{(3)}]\ \  ${(|\wt{V}_a|)}_{a\in A}$ is a locally finite covering of $|\partial X|$ by sets open in $|X|$.
\end{itemize}
Note that in general $a\not\in \wt{V}_a$, which is indicated by writing $\wt{V}_a$ instead of $\wt{V}(a)$.
\end{lemma}
\begin{proof}
The collection of open subsets ${(|\wt{U}(x)|)}_{x\in \partial X}$ is an open covering of $|\partial X|$ by sets open in $|X|$.
If we add $|X\setminus \partial X|$ we obtain an open covering of $|X|$.
   Using that $|X|$ is paracompact, and therefore even metrizable, we find a subordinate open, locally finite covering
${(W_a)}_{a\in X}$, where,  of course, most sets are empty. We also may assume that $\cl_{|X|}(W_a)\subset |\wt{U}(a)|$.
We denote by $A\subset \partial X$ the set indexing 
nonempty sets.   Then $\emptyset\neq \cl_{|X|}(W_a)\subset |\wt{U}(a)|$  for $a\in A$.  We define
$$
\wt{V}_a: = \wt{U}(a) \cap \pi^{-1}(W_a).
$$
  Then $\cl_{|X|}(|\wt{V}_a|)\subset | \wt{U}(a)|$ for all $a\in A$ and ${(|\wt{V}_a|)}_{a\in A}$ is a locally finite  covering of $|\partial X|$ by sets open in $|X|$.
\qed \end{proof}

We note that for $a\in A$ the  various sets are related by
\begin{eqnarray}
\emptyset\neq \wt{V}_a\subset  \cl_{X}(\wt{V}_a) \subset \wt{U}(a)\subset\cl_X(\wt{U}(a))\subset  \wt{U}'(a).
\end{eqnarray}
In view of Proposition \ref{CORR13.6.12} we also have for $x,x'\in\partial X$ open subsets $\wh{O}_{x,x'}$ of $\bm{X}$ which satisfy 
\begin{eqnarray}
\wh{O}_{x,x'}\cap\partial \bm{ X} =\bm{U}'(x,x')\subset \wh{O}_{x,x'}\subset \wt{\bm{U}}'(x,x').
\end{eqnarray}
For the next result we use the basic data  $(\Lambda,\mathfrak{U}',\mathfrak{S}',\tau')$ from Proposition \ref{basic-data},  the open sets ${(\wt{U}'(x))}_{x\in \partial X}$ defined 
thereafter satisfying {\bf (P5)}, and the open sets ${(\wt{U}(x))}_{x\in \partial X}$ satisfying 
{\bf (P6)}.  In addition recall  ${(|\wt{V}_a|)}_{a\in A}$ obtained in Lemma \ref{Lemma2-62} and 
\begin{eqnarray}
U(x):=\wt{U}(x)\cap \partial X\ \ \text{and}\ \  U'(x)=\wt{U}'(x)\cap \partial X.
\end{eqnarray}
  The goal in the remaining part of this section 
is to construct for $a\in A$ open subsets $\wh{V}_a\subset\wt{U}'(a)$ of $X$ with suitable properties 
such that in particular
$$
\wh{\bm{V}}_{a,a'}:=\{\phi\in\bm{X}\ |\ s(\phi)\in\wh{V}_a,\ t(\phi)\in \wh{V}_{a'}\}
$$
satisfies $\wh{\bm{V}}_{a,a'}\subset \wh{O}_{a,a'}$ for $a,a'\in A$. This implies as a consequence of Proposition \ref{CORR13.6.12} that,
in particular, the extension $\wt{\tau}_{a,a'}$ is defined on $\wh{\bm{V}}_{a,a'}$. 
We observe that with respect to Remark \ref{remark14.3.4} we are moving in the right direction.
However, we note that $\wh{V}_a$ is in general not a neighborhood of $a$, and the points $a$ 
only belong to the subset $A$ of $\partial X$, rather than allowing all the points in $\partial X$.

\begin{lemma}\label{hoofers}
Fixing $a\in A$, for every  $y\in U'(a)$ there exists a finite subset $o(y)\subset A$ and a $G_y$-invariant open neighborhood
$\wt{O}(y)\subset \wt{U}'(a)$ with the following properties.
\begin{itemize}
\item[\em{(1)}]\ \ The orbit set $|\wt{O}(y)|$ intersects $|\wt{V}_b|$ only if $b\in o(y)$.
\item[\em{(2)}]\ \ $g\ast \wt{O}(y) =\wt{O}(g\ast y)$ and $o(y)=o(g\ast y)$ for $g\in G_a$.
\item[\em{(3)}]\ \  For every $b\in o(y)$ it holds that 
$$
\{\phi\in {\bf X}\ |\ s(\phi)\in\wt{O}(y)\cap\wt{U}(a),\ t(\phi)\in \wt{U}(b)\}\subset \wh{O}_{a,b}.
$$
\end{itemize}
\end{lemma}
\begin{proof}
(1) Fix $a\in A$ and pick $y\in U'(a)$. Since ${(|\wt{V}_a|)}_{a\in A}$ is a locally finite  covering of $|\partial X|$ by open subsets in $|X|$
we find a finite subset $o(y)$ of $A$ so that every neighborhood of $|y|$ intersects $|\wt{V}_b|$ for $b\in o(y)$ 
and there exists a sufficiently small neighborhood intersecting these sets only. Hence we find an open neighborhood $\wt{O}(y)$
which is $G_y$-invariant and contained in $ \wt{U}'(a)$ so that (1) holds for $|\wt{O}(y)|$.  \par

\noindent (2) We define for $g\in G_x$ the set 
$\wt{O}(g\ast y):= g\ast \wt{O}(y)$.
Clearly $|\wt{O}(g\ast y)|=|\wt{O}(y)|$, and  $o(g\ast y) =o(y)$ which proves (2). 
We can take $\wt{O}(y)$ as small as we wish.\par

\noindent(3)  We would like to show that (3) can be achieved if $\wt{O}(y)$ is small enough. Here $y\in U'(a)$. Arguing indirectly we find 
  a sequence of morphisms $(\phi_k)\subset\bm{X}$ with the following properties.
  \begin{itemize}
  \item $s(\phi_k)\in \wt{U}(a)$  and  $(s(\phi_k))$ converging to $y$.
\item $t(\phi_k)\in \wt{U}(b_k)$, where $b_k\in o(y)$.
\item $\phi_k\not\in \wh{O}_{a,b_k}$.
\end{itemize}
From the  properness properties of the finitely many $\wt{U}(b_k)$ we deduce that each subsequence of $(\phi_k)$ has a convergent subsequence, where the limit
then satisfies $s(\phi_0)=y\in \cl_X(U(a))$ and $t(\phi_0)\in \cl_X(U(b))$ for some $b\in o(y)$. 
Hence $\phi_0\in \bm{U}'(a,b)$ and therefore $\phi_0\in \wh{O}_{a,b}$. This contradiction implies that for $k$ large we must have that 
$\phi_k\in \bigcup_{b\in o(y)} \wh{O}_{a,b}$. Hence for $\wt{O}(y)$ sufficiently small (3) holds. Taking $\wt{O}(y)$ even smaller
the same must hold for the finitely many $\wt{O}(g\ast y)$.
\qed \end{proof}
We define for $a\in A$ the open subset $\wh{V}_a$ of $X$ by
\begin{eqnarray}
\wh{V}_a =\left(\bigcup_{y\in U'(a)} \wt{O}(y)\right)\cap  \wt{V}_a.
\end{eqnarray}
The collection of open sets ${(\wh{V}_a)}_{a\in A}$  has the following properties.
\begin{proposition}\label{important-xx}
\begin{itemize}
\item[\em{(1)}]\ \ For $a\in A$ the open subset $\wh{V}_a$ of $X$ is $G_a$-invariant and contained in $\wt{V}_a\subset \wt{U}(a)$
\item[\em{(2)}]\ \ $\wh{V}_a\cap \partial X= \wt{V}_a\cap \partial X$ for $a\in A$.
\item[\em{(3)}]\ \ The collection ${(|\wh{V}_a|)}_{a\in A}$ is a locally finite covering of $|\partial X|$ by subsets open in $|X|$.
\item[\em{(4)}]\ \ Defining for $a,a'\in A$ the sets 
$$
\wh{\bm{V}}_{a,a'}:=\{\phi\in {\bm{X}}\ |\ s(\phi)\in \wh{V}_a,\ t(\phi)\in \wh{V}_{a'}\}
$$
 it holds that
$$
\wh{\bm{V}}_{a,a'}\subset \wh{O}_{a,a'}\ \text{for}\ a,a'\in A. 
$$
In particular we can restrict the locally constant map
$\wh{\tau}_{a,a'}:\wh{O}_{a,a'}\rightarrow \mathfrak{J}(a,a')$ to $\wh{\bm{ V}}_{a,a'}$.
\end{itemize}
\end{proposition}
\begin{proof}
By construction $\wh{V}_a$ is open, $G_a$-invariant, and contained in $\wt{V}_a$. This implies (1).\par

Point (2) follows since by construction $\bigcup_{y\in U'(x)} \wt{O}(y)$ contains $\wt{V}_a\cap \partial X$. Hence 
$$
\wh{V}_a\cap \partial X=\partial X\cap \left(\left(\bigcup_{y\in U'(a)} \wt{O}(y)\right)\cap  \wt{V}_a\right)= \wt{V}_a\cap\partial X.
$$
The collection ${(|\wh{V}_a|)}_{a\in A}$ is an open covering of $\partial X$. Since $\wh{V}_a\subset \wt{V}_a$ for $a\in A$
it is in addition locally finite. This gives (3). \par

Point (4)  holds by construction of the sets $\wh{V}_a$ for $a\in A$, see specifically Lemma \ref{hoofers} (3).   With $a,a'\in A$ 
assume that $\phi\in \wh{\bm{V}}_{a,a'}$.  By construction $\wh{V}_a\subset \wt{U}(a)$ and there exists $y\in U'(x)$
such that $s(\phi)\in \wt{O}(y)$. By Lemma \ref{hoofers} (1)  and (3) we must have that $t(\phi) \in \wt{V}_b$ for some 
$b\in o(y)$. Since $\wh{V}_b\subset \wt{V}_b$ we conclude that $a'\in o(y)$ and $\phi\in \wh{O}_{a,a'}$.
\qed \end{proof}
At this point we have constructed for $a\in A$ open subsets $\wh{V}_a$ of $X$ so that
\begin{eqnarray}
\wh{V}_a\subset \wt{V}_a\subset\cl_X(\wt{V}_a)\subset \wt{U}(a)\ \text{for}\ a\in A,
\end{eqnarray}
and
\begin{eqnarray}
|\partial X|  =\bigcup_{a\in A} |\wh{V}_a|.
\end{eqnarray}
Moreover ${(|\wh{V}_a|)}_{a\in A}$ is a locally finite covering of $|\partial X|$ by sets open in $|X|$.
In addition, for $a,a'\in A$ it holds
\begin{eqnarray}
\wh{\bm{V}}_{a,a'}=\{\phi\in \bm{X}\ |\ s(\phi)\in \wh{V}_a,\ t(\phi)\in \wh{V}_{a'}\}\subset \wh{O}_{a,a'},
\end{eqnarray}
so that $\wh{\tau}_{a,a'}$ defines a map 
\begin{eqnarray}
\wh{\tau}_{a,a'}|\wh{\bm{ V}}_{a,a'}: \wh{\bm{ V}}_{a,a'}\rightarrow \mathfrak{J}(a,a').
\end{eqnarray}
Define for $a\in A$ the open subset $V_a$ of $U(a)$ by 
\begin{eqnarray}\label{LLLX1}
V_a =U(a)\cap  \wh{V}_a =U(a)\cap \wt{V}_a
\end{eqnarray}
and note that
\begin{eqnarray}\label{LLLX2}
\cl_{|X|}(|V_a|)\subset |U(a)|.
\end{eqnarray}
Recall also that
$$
U(a)=\partial X\cap \wt{U}(a).
$$
The data we have constructed will be used in the next section.

\section{Implicit Structures and Local Extension}
The set-up is again as described in Theorem \ref{p-main-p} and we assume we have
taken a representative $(\Lambda,\mathfrak{U}',\mathfrak{S}',\tau')$ satisfying 
{\bf (P1)--(P4)}.  Also let 
${(\wt{U}'(x))}_{x\in\partial X}$ be a covering of $\partial X$ by
sets open in $X$ satisfying {\bf (P5)}, and let ${(\wt{U}(x))}_{x\in\partial X}$ 
satisfy {\bf (P6)}.
We start this section with an important remark.
\begin{remark}\label{I_REMARK}\index{R- Important remark about correspondences}
For the following part of the construction it is {\bf important (!)} to be aware 
of a subtlety, 
which might not be apparent at first sight. For $x,x'\in\partial X$ the maps $\tau'_{x,x'}:{\bf U}'(x,x')\rightarrow \mathfrak{J}(x,x')$, 
 define correspondences between different sc$^+$-section structures.  It is crucial to note 
 that in general the symmetric sc$^+$-section structure 
${(s'^x_i)}_{i\in I_x}$ might implicitly have additional properties over certain subsets of $U'(x)$, which have to be taken into account
to define a structurable extension. For example assume that 
we have morphisms $\psi,\phi,\sigma\in \partial \bm{X}$ so that $s(\psi)=s(\phi)$ and $s(\sigma)=t(\psi)$ and $t(\sigma) = t(\phi)$
so that in addition 
$$
\sigma\circ \psi=\phi,
$$
where $s(\psi), s(\phi)\in U'(x)$, $t(\psi)\in U'(x')$, and $t(\phi)\in U'(x'')$. 
Associated to these three morphisms there are correspondences $\tau'_{x,x'}(\psi)$, $\tau'_{x,x''}(\phi)$, and $\tau'_{x',x''}(\sigma)$. Assume they are represented by the diagrams
$$
I_{x}\stackrel{\alpha}{\twoheadleftarrow} I\stackrel{\beta}{\twoheadrightarrow} I_{x'},\ 
I_{x}\stackrel{\alpha'}{\twoheadleftarrow} I'\stackrel{\beta'}{\twoheadrightarrow} I_{x''},\ \text{and}\ 
I_{x'}\stackrel{\alpha''}{\twoheadleftarrow} I''\stackrel{\beta''}{\twoheadrightarrow} I_{x''}.
$$
It may happen that for $i_1\not\in G_x(i_2)$, where $i_1,i_2\in I_{x}$,  there exist $k\in I, k'\in I', k''\in I''$ such that
$$
i_1 =\alpha(k),\ \beta(k)=\alpha''(k''),\ \beta''(k'') =\beta'(k'),\ \alpha'(k') =i_2.
$$
This forces the equality $s_{i_1}^x(y)=s_{i_2}^x(y)$ for $y\in U'(x)$ near $s(\psi)$. So the global structure described by $\tau'$
forces  additional local conditions. The constructions needed to prove an extension theorem for a structurable sc$^+$-multisection functor has to accommodate this fact. As preparation we introduce some additional structures to keep track of the implicit conditions.
\qed
\end{remark}
In view of this important remark we shall capture the relevant local information from  the global structure 
provided by $\tau'$. As it turns out we only need to study  for most parts 
\begin{eqnarray*}
\left\{\wh{\tau}_{x,x'}:\wh{O}_{x,x'}\rightarrow \mathfrak{J}(x,x')\right\}_{(x,x')\in A\times A}\ \text{and}\ \left\{ \wh{V}_a\right\}_{a\in A}.
\end{eqnarray*}
Recall that 
\begin{eqnarray*}
\cl_X(\wh{V}_a)\subset \wt{U}(a)\ \text{for}\ x\in A, \text{and}\ \wh{V}_a\cap\partial X=V_a,\ a\in A,
\end{eqnarray*}
and according to Proposition \ref{important-xx} with $\wh{\bm{V}}_{a,a'}:=\{\phi\in\bm{X}\ |\ s(\phi)\in \wh{V}_a,\ t(\phi)\in \wh{V}_{a'}\}$ it holds that
$$
\wh{\bm{V}}_{a,a'}\subset \wh{O}_{a,a'}\ \ \text{for}\ \  a,a'\in A.
$$
The following definitions are needed to address the issues raised in Remark \ref{I_REMARK}.
\begin{definition}\index{D- $z$-loop}
 Fix $a\in A$ and 
and let $z\in \wh{V}_a$. A $z$-{\bf loop (for $a$)} of morphisms is given by the following data.
\begin{itemize}
\item[(1)]\   A finite sequence $a_0,...,a_n$ in $A$ satisfying $a_0=a=a_n$.
\item[(2)]\  Morphisms $\phi_0,...,\phi_{n-1}$ with the properties
\begin{itemize}
\item[(a)]\  $t(\phi_j)=s(\phi_{j+1})$ for $j=0,..,n-2$, $s(\phi_0)=z=t(\phi_{n-1})$.
\item[(b)]\  $s(\phi_j)\in \wh{V}_{a_j}$ and $t(\phi_j)\in \wh{V}_{a_{j+1}}$ for $j=0,..,n-1$.
\item[(c)]\  $1_z=\phi_{n-1}\circ..\circ \phi_0$.
\end{itemize}
\end{itemize}
\qed
\end{definition}
Given a $z$-loop $\phi_0,..,\phi_{n-1}$ we note that $\phi_j\in \wh{\bm{V}}_{a_j,a_{j+1}}$.
Consequently we  can consider $\wh{\tau}_{a_j,a_{j+1}}(\phi_j)$ for $j=0,..,n-1$, where the latter are represented by diagrams
$$
I_{a_j}\stackrel{\alpha_j}{\twoheadleftarrow} I_j\stackrel{\beta_j}{\twoheadrightarrow} I_{a_{j+1}}.
$$

\begin{definition}\index{D- $z$-arc}
Fix $a\in A$ and pick $z\in \wh{V}_a$.
A {\bf $z$-arc (for $a$)} associated to the $z$-loop consists of  a sequence of indices $(i_0,i_1,...,i_n)$ satisfying
\begin{itemize}
\item[(1)]\  $i_j\in I_{a_j}$ for $j\in \{0,..,n\}$.
\item[(2)]\   $(i_j,i_{j+1})=(\alpha_j(k_j),\beta_j(k_j))$ for $j\in \{0,...,k-1\}$ and suitable $k_j\in I_j$.
\end{itemize}
\qed
\end{definition}
Note that $i_0,i_n\in I_a$. In general they do not need to be the same and one might view $i_n$ as the `index holonomy' 
of $i_0$ along the $z$-loop. However,  this `holonomy' is not unique, since even using the same $z$-loop and starting with  the same $i_0$,
possible different choices of the $ k_j$ might lead to different endpoints.
\begin{definition}\index{D- $z$-equivalent}
Fix $a\in A$ and 
let $z\in \wh{V}_a$ and $i,i'\in I_a$. We say that $i$ and $i'$ are {\bf $z$-equivalent (for $a$)} if there exists a $z$-arc with
$i_0=i$ and $i_n=i'$.  In this case we  shall write $i\sim_z i'$.
\qed
\end{definition}
We shall denote the $\sim_z$-equivalence class of some index $i\in I_a$ by $[i]^{\sim_z}$,
where
$$
[i]^{\sim_z}:=\{i'\in I_a\ |\ i'\sim_z i\}.
$$
The equivalence relations $\sim_z$ have the following  properties summarized in two results.
The first statement is essentially obvious.
\begin{lemma}
Let $a\in A$ and assume  that $z\in \wh{V}_a$,   $i,i'\in I_a$, and $g\in G_a$. Then 
$i\sim_z i'\ \text{if and only if}\   g(i)\sim_{g\ast z} g(i')$. In particular
$$
[g(i)]^{\sim_{g\ast z}} = g([i]^{\sim_z})\ \ \text{for}\ \ i\in I_a,\ g\in G_a,\ z\in \wh{V}_a.
$$
\end{lemma}
\begin{proof}
If $i\sim_z i'$ there exists a $z$-loop relating $i$ with $i'$. Denote the morphisms in the loop by
$$
\phi_0,...,\phi_{n-1}.
$$
Then $\Gamma(g^{-1},g\ast z),\phi_0,...,\phi_{n-1},\Gamma(g,z)$ is a $g\ast z$-loop and it relates the indices $g(i)$ and $g(i')$. Hence 
$g(i)\sim_{g\ast z} g(i')$. That $g(i)\sim_{g\ast z} g(i')$ implies $i\sim_z i'$ follows similarly.
\qed \end{proof}

\begin{proposition} \label{approx}
The equivalence relations $\sim_z$, $z\in \wh{V}_a$, on $I_a$  have the following properties.
\begin{itemize}
\item[\em{(1)}]\  Let  $z_0\in \wh{V}_{a}$ and $i,i'\in I_a$.   If $i\sim_{z_0} i'$ the same holds for all
$z\in \wh{V}_a$  close to $z_0$.
\item[\em{(2)}]\   If $i,i'\in I_a$ and $(z_k)\subset  V_a$  is a sequence for which $i\sim_{z_k} i'$ and
$z_k\rightarrow z_0\in \cl_X(V_a)$, then there exists an open neighborhood $Q(z_0)\subset U(a)$
for which $s'^a_i(y)=s'^a_{i'}(y)$ if $y\in Q(z_0)$.
\end{itemize}
\end{proposition}
We note that (2) is a particularly important observation.   The proof needs some preparation and the following lemma will be important.
  \begin{lemma}\label{L14.3.7}
 Denote by $X$ an ep-groupoid and  let $x_1,...,x_\ell$ be a finite collection of isomorphic objects. Given open neighborhoods $Q(x_i)$ in $X$ there exist open neighborhoods
  $O(x_i)\subset Q(x_i)$ with the following properties.
  \begin{itemize}
  \item[\em{(1)}]\  Each $O(x_i)$ is equipped with the natural $G_{x_i}$-action.
  \item[\em{(2)}] \ For every morphism $\psi:x_i\rightarrow x_j$ there exists an open neighborhood $W(\psi)$ 
  such that $s:W(\psi)\rightarrow O(x_i)$ and $t:W(\psi)\rightarrow O(x_j)$ are sc-diffeomorphisms. 
  \item[\em{(3)}]\ For $\psi:x_i\rightarrow x_j$ and $\phi:x_j\rightarrow x_k$ the map $m:W(\phi){_{s}\times_t}W(\psi)\rightarrow \bm{X}:(\phi',\psi')\rightarrow \phi'\circ \psi'$
  defines an sc-diffeomorphism onto $W(\phi\circ\psi)$.
  \end{itemize}
   \end{lemma}
\begin{proof}
We prove this result by induction with the respect to the cardinality of the set of isomorphic objects.
The statement for one object is equivalent to the existence of arbitrarily small neighborhoods
with the natural action. Assume we are given isomorphic objects $x_1,...,x_{\ell+1}$. By induction hypothesis 
we find $O'(x_i)$ for $i=1,...,\ell$ with the natural $G_{x_i}$-action, and for every morphism $\psi:x_i\rightarrow x_j$ between these objects an open neighborhood $W(\psi)$ such that $s:W(\psi)\rightarrow O'(x_i)$ and 
$t:W(\psi)\rightarrow O'(x_j)$ are sc-diffeomorphisms, so that (3) holds.  Add the point $x_{\ell+1}$
to our considerations. We pick for every $i=1,..,\ell$ a morphism starting at $x_{\ell+1}$ and ending at $x_i$
$$
\psi_i:x_{\ell+1}\rightarrow x_i.
$$
Associated to $\psi_i$ we have a local sc-diffeomorphism $\wh{\psi}_i$. It is obtained in the usual way  as
$\wh{\psi}_i= t\circ (s|W(\psi_i))^{-1}$, where $W(\psi_i)$ is a suitable open neighborhood so that $s:W(\psi_i)\rightarrow O(x_{\ell+1})$ and $t:W(\psi_i)\rightarrow O(x_i)$ are sc-diffeomorphisms.
We pick $O(x_{\ell+1})$
admitting the $G_{x_{\ell+1}}$-action and so small that $\wh{\psi}_i(O(x_{\ell+1}))\subset O'(x_i)$
and $\wh{\psi}_i:O(x_{\ell+1})\rightarrow \wh{\psi}_{i}(O(x_{\ell+1}))$ is an sc-diffeomorphism.
Define $O(x_i)= \wh{\psi}_{i}(O(x_{\ell+1}))$. 
For $g\in G_{x_{\ell+1}}$ we define $W(g)= \Gamma(g,O(x_{\ell+1}))$, and for  $\psi:x_{\ell+1}\rightarrow x_i$
we define $W(\psi)$  as follows, where we use the unique decomposition $\psi=\psi_i\circ g$ for some $g\in G_{x_{\ell+1}}$.
$$
W(\psi)=\{ (s|W(\psi_i))^{-1}(g\ast z)\circ\Gamma(g,z)\ |\ z\in O(x_{\ell+1})\}.
$$
Let us show that the new set of data is the desired construction associated to the isomorphic objects 
$x_1,..,x_{\ell+1}$.

First we consider $s:W(\psi)\rightarrow O(x_{\ell+1})$ for a $\psi:x_{\ell+1}\rightarrow x_i$.
It follows immediately from the construction that $s$ is an sc-diffeomorphism. Similarly
$t:W(\psi)\rightarrow O(x_i)$ takes the form 
$$
t(\sigma) =\wh{\psi}_i(g\ast s(\sigma))
$$
which implies that $t$ is an sc-diffeomorphism.  This completes the proof of (2).

In order to verify (3) assume that $\psi:x_{\ell+1}\rightarrow x_i$ and $\phi:x_i\rightarrow x_j$.
Consider  
$$
W(\phi){_{s}\times_t} W(\psi)\rightarrow W(\phi\circ\psi).
$$
We first construct an sc-smooth map
$\Delta: W(\psi)\rightarrow W(\phi){_{s}\times_t} W(\psi)$ by
$$
\sigma \rightarrow ((s|W(\phi))^{-1}(t(\sigma)), \sigma)
$$
This map is an sc-diffeomorphism with inverse being the projection onto $W(\psi)$. 
We also note that 
$$
\delta: W(\psi)\rightarrow W(\phi\circ\psi) :\sigma\rightarrow ((s|W(\phi))^{-1}(t(\sigma)))\circ  \sigma)
$$
 is an sc-diffeomorphism
onto $W(\phi\circ \psi)$. Since $m\circ \Delta =\delta$ the result follows.
\qed \end{proof}
\begin{proof}[Proposition \ref{approx}]
(1) For $z_0\in \wh{V}_a\subset \wt{U}(a)$ there exists an open neighborhood of $|z_0|$ only intersecting finitely many $|\cl_X(\wh{V}_y)|$, $y\in A$.
Denote the corresponding $y$ by $a_1,...,a_\ell$. We may assume that $a_1=a$.
A morphism $\psi$ starting 
in $\wh{V}_a$ near $z_0$ and ending in one of the $\wh{V}_y$, $y\in A$,  has to end in one of the $\wh{V}_y$ with $y\in \{a_1,..,a_\ell\}$.
In each of the $\wh{V}_{a_i}$ there can be only a finite number of points which can be reached from $z_0$. 
Denote the collection of all these points, where we vary over the $a_1,...,a_\ell$   by $x_1,...,x_b$, and note that they are all isomorphic. 
Take a $z_0$-loop of morphisms starting and ending at $z_0$ so that the composition is $1_{z_0}$ and the corresponding
$z_0$-arc starts at $i$ and connects it with $i'$. As a consequence of Lemma \ref{L14.3.7}, when slightly moving $z_0$ to $z$ 
the morphisms deform to a $z$-loop and since the $\tau_{a,a'}'$ are locally constant we infer that $i\sim_z i'$.\par

Before we prove (2) we need the following consideration which already was discussed in Lemma \ref{simple_lemma}. Recall from (\ref{LLLX1}) and (\ref{LLLX2}) that $V_a\subset U(a)$ and $\cl_{|X|}(|V_a|)\subset |U(a)|$. In addition $U(a)$
is a good open neighborhood of $a$ in the sense of Definition \ref{DEF_good_open} and it has 
properness property,  which means that $t:s^{-1}(\cl_X(U(a)))\rightarrow X$ is a proper map.
The lemma then asserts 
that $\cl_X(V_a)\subset U(a)$.\par

(2)  In order to prove (2) we pick $z_0\in \cl_X(V_a)$ and with the help
of Lemma \ref{L14.3.7}, brought into the right context,
 we can argue in a similar fashion as before.
 For the point  $z_0\in \cl_X(V_a)$ consider all objects which belong to some $\cl_X(V_b)$, $b\in A$,
and which are isomorphic to $z_0$. We note that there are only finitely many such objects
due to the fact that $(|V_b|)$ is a locally finite covering of $|\partial X|$ and $V_a\subset U(a)$,
and $U(a)$ has the properness property.
Denote these objects (including $z_0)$ by $x_1=z_0,...,x_\ell$ and apply Lemma \ref{L14.3.7}.

Assume that $(z_k)\subset V_a$, $z_k\rightarrow z_0$, and there exists a $z_k$-loop 
with an associated $z_k$-arc relating $i$ and $i'$ such that $i\sim_{z_k} i'$.  If $z_k$ is close enough
to $z_0$ it follows that the occurring morphisms for $z_k$ lie in the $W(\phi)$ for suitable 
$\phi$ as given in Lemma \ref{L14.3.7}. We may assume, starting with sufficiently small $O(x_i)$
that $\tau'_{x_i,x_j}$ is constant on the associated $W(\phi)$. Just taking $k$ large enough  we can take
a loop of morphisms 
starting and ending at $z_0$ and composing to $1_{z_0}$ (having the same combinatorial structure as the $z_k$-loop) so that the index $i$ is related to the index $i'$. This implies  for a suitable open neighborhood $Q(z_0)$,
 that $s'^a_i(y)=s'^a_{i'}(y)$ for $y\in Q(z_0)$. The proof is complete.
\qed \end{proof}

In next step we define a new equivalence relation $\approx_z$ for points  $z\in \cl_X(\wh{V}_a)$.
\begin{definition}\index{$\Theta_z$}
For $a\in A$ and $z\in \cl_X(\wh{V}_a)$
we define the set $\Theta_z$   to consist of all
$(i,i')\in I_a\times I_a$ such that there exists 
$(z_k)\subset \wh{V}_a$ with $i\sim_{z_k} i'$ and $z_k\rightarrow z$. 
We define $i\approx_z i'$\index{$\approx_z$}  provided there exist $i_0,...,i_k\in I_a$  with $(i_p,i_{p+1})\in \Theta_z$ for
$p=0,..,k-1$, and $i_0=i$, $i_p=i'$.
\qed
\end{definition}

It is clear that for $z\in \wh{V}_a$ the equivalence $i\sim_z i'$ implies $i\approx_z i'$.  
Moreover, if $z\in \cl_X(V_a)=\cl_{\partial X}(V_a)$, and $i\approx_z i'$ it follows that $s'^a_i(y)=s'^a_{i'}(y)$ for $y\in Q(z)$, where $Q(z)$ is an open neighborhood
in $\partial X$ contained in $U(z)$, see Proposition \ref{approx} (2).

For the issue discussed in Remark \ref{I_REMARK} the equivalence relation $\sim_z$ encapsulates the requirements.
However, the equivalence relation $\approx_z$ is more convenient for the necessary constructions.
We shall define for given $i\in I_a$ and $z\in \cl_X(\wh{V}_a)$
$$
 [i]^{\approx_z}= \{i'\in I_a\ |\ i'\approx_z i\}.
$$
\begin{lemma}\label{true-y}\index{L- Properties of $\sim_z$ and $\approx_z$}
The following holds true.
\begin{itemize}
\item[\em{(1)}]\  If $i\sim_z i'$ for $z\in\wh{V}_a$, then $i\approx_z i'$.
\item[\em{(2)}]\  If $z\in \cl_X(\wh{V}_a)$ and $g\in G_a$ then $i\approx_z i'$ implies $g(i)\approx_{g\ast z} g(i')$. In particular
$g([i]^{\approx_z})= [g(i)]^{\approx_{g(z)}}$.
\item[\em{(3)}]\   If $z\in \cl_X({V}_a)$ and $i\approx_z i'$ then $s'^{a}_{i}(y)=s'^{a}_{i'}(y)$ for $y\in {U}(a)$ close to $z$.
\item[\em{(4)}]\  With  $i\in I_a$ and  $z_0\in \cl_X(\wh{V}_a)$ there exists an open neighborhood
$\wt{Q}(z_0)$ in $\wt{U}(a)$ so that $[i]^{\sim_z}\subset [i]^{\approx_{z_0}}$ for all $z\in \wh{V}_a\cap \wt{Q}(z_0)$.
\end{itemize}
\end{lemma}
\begin{proof}
(1) and (2) are trivial, and  (3) follows from Proposition  \ref{approx} (2). 
For the  proof of (4) we argue indirectly.  Let  $z_0\in \cl_X(\wh{V}_a)$ and $i\in I_a$.
 Assume there exists
a sequence $(z_k)\subset \wh{V}_a$ satisfying $z_k\rightarrow z_0$ and $i_0\in I_a$ satisfying
$i_0\not\in [i]^{\approx_{z_0}}$ and $i\sim_{z_k} i_0$. Then by definition $i\approx_{z_0} i_0$ implying that
$i_0\in [i]^{\approx_{z_0}}$ giving a contradiction.
\qed \end{proof}

The previous discussion allows us to define open neighborhoods around the points $z\in \cl_X(\wh{V}_a)$ 
with suitable properties, which is accomplished in the next proposition. Recall that for $a\in A$ it holds
that $\cl_X(\wh{V}_a)\subset \wt{U}(a)$.
\begin{proposition}\label{neighbor}
Fix $a\in A$. Then there exists for every $z\in \cl_X(\wh{V}_a)$ an open neighborhood $\wt{\mathsf{O}}(z)$ with the following properties.
\begin{itemize}
\item[\em{(1)}]\ \ $\wt{\mathsf{O}}(z)\subset \wt{U}(a)$.
\item[\em{(2)}]\ \ $\wt{\mathsf{O}}(g\ast z)=g\ast \wt{\mathsf{O}}(z)$ for $g\in G_a$ and $\wt{\mathsf{O}}(z)\cap \wt{\mathsf{O}}(g\ast z)=\emptyset$ if $z\neq g\ast z$.
\item[\em{(3)}] \ \ If $z\in \cl_X(\wh{V}_a)$ and $i\in I_a$  then $[i]^{\sim_y}\subset [i]^{\approx_z}$ for $y\in \wt{\mathsf{O}}(z)\cap \wh{V}_a$.
\item[\em{(4)}]\ \ If $z\in \cl_X(V_a)$ it holds $s_i'^a(y)=s_{i'}'^a(y)$ for $i\approx_z i'$ and $y\in \wt{\mathsf{O}}(z)\cap U(a)$.
\end{itemize}
\end{proposition}
\begin{proof}
This is an immediate consequence of the previous discussion, particularly Lemma \ref{true-y}.
\qed \end{proof} 

The proposition prompts the following definition.
\begin{definition}\label{GOODNEI}
A {\bf good neighborhood system} for $\cl_X(\wh{V}_a)\subset \wt{U}(a)$ is given by a family 
$$
{ \left(\wt{\mathsf{O}}(z)\right)}_{z\in \cl_X(\wh{V}_a)}
$$
of open (in $X$) neighborhoods  with the following properties. 
\begin{itemize}
\item[(1)]\ \  $\wt{\mathsf{O}}(z)\subset \wt{U}(a)$.
\item[(2)]\ \ $\wt{\mathsf{O}}(g\ast z)=g\ast \wt{\mathsf{O}}(z)$ for $g\in G_a$ and $\wt{\mathsf{O}}(z)\cap \wt{\mathsf{O}}(g\ast z)=\emptyset$ if $z\neq g\ast z$.
\item[(3)]\ \ If $z\in \wh{V}_a$ and $i\in I_a$  then $[i]^{\sim_y}\subset [i]^{\approx_z}$ for $y\in \wt{\mathsf{O}}(z)\cap \wh{V}_a$.
\item[(4)]\ \ If $z\in \cl_X(V_a)$ it holds  $s_i'^a(y)=s_{i'}'^a(y)$ for $i\approx_z i'$ and $y\in \wt{\mathsf{O}}(z)\cap U(a)$.
\end{itemize}
\qed
\end{definition}
\begin{remark}
The good neighborhood system captures part of the  implicit information of having a system of correspondences.
Note that the initial choice of ${(V_a)}_{a\in A}$ and the extensions $\wt{V}_a,a\in A$ are important. The existence of such a collection
$(V_a)$ is a consequence of the assumption that $|X|$ is  paracompact.
\qed
\end{remark}
We shall use the previous discussion to define suitable local extensions of some of the local section structures, which are compatible 
with the implicit relations. They will be used later on for the proof of the global extension theorem.

The following constructions are concerned with the $\mathfrak{s}'_a$, $a\in A$,
$$
\mathfrak{s}_a' ={(s^a_i)}_{i\in I_a},
$$
where $s^a_i$ is an sc$^+$-section of $W|U'(a)$.  Recall that 
\begin{eqnarray*}
&\wh{V}_a\subset\cl_{X}(\wh{V}_a)\subset  \wt{U}(a)\subset \wt{U}'(a).
&
\end{eqnarray*}
 Note that $\Lambda$ 
is as an sc$^+$-multisection (without the structuring) completely determined by the collection
of all $\mathfrak{s}_a'$, where $a$ varies in $A$. This is, of course, the case since the collection
$(|V_a|)_{a\in A}$ covers $|\partial X|$, where $V_a=\wh{V}_a\cap\partial X$. 

For the further constructions we shall need extensions of the $s^a_i$ to $\wh{V}_a$ having suitable properties.
The following discussion prepares the construction of such extensions.  Let $z\in \cl_X(\wh{V}_a)\cap \partial X$
and denote by $[z]$ its orbit under the $G_a$-action.
We observe that $I_a$ has different decompositions 
into equivalence classes associated to the different elements $z'\in [z]$ defined by $\approx_{z'}$. 
Fix a $z'\in [z]$ and denote by ${[i_1]}^{\approx_{z'}},...,{[i_k]}^{\approx_{z'}}$ the associated different equivalence 
classes of $I_a$. If we take $z''\in [z]$ and $g\in G_a$ with $g\ast z'=z''$ then 
the collection 
$$
g({[i_j]}^{\approx{z'}}):=\{ g(i)\ |\ i\in {[i_1]}^{\approx{z'}}\}
$$
where $j=1,...,k$ is the set of $\approx_{z''}$ equivalence classes. With other words
$g({[i_j]}^{\approx_{z'}})= {[g(i_j)]}^{\approx_{z''}}$ and 
${[g(i_1)]}^{\approx_{z''}},...,{[g(i_k)]}^{\approx_{z''}}$ are the $\approx_{z''}$-equivalence classes.
Taking a different $g'$ with $g'\ast z'=z''$ only gives a different presentation of these equivalence classes.

\begin{lemma}\label{LEMMX14413}
For fixed $a\in A$ and $z\in \cl_X(\wh{V}_a)\cap \partial X$ define $\wt{\mathsf{O}}([z])$ as the union of all
$\wt{\mathsf{O}}(z')$ with $z'\in [z]$, where $[z]$ is the $G_a$-orbit of $z$.  Fix $z'\in [z]$ 
and let $G_a(z')=\{g\in G_a\ |\ g\ast z'=z'\}$. For $I_a$ let ${[i_1]}^{\approx_{z'}},.., {[i_k]}^{\approx_{z'}}$
be the decomposition of $I_a$ into $\approx_{z'}$-equivalence classes
$$
I_a = {[i_1]}^{\approx_{z'}}\sqcup ......\sqcup {[i_k]}^{\approx_{z'}}.
$$
Then the following holds.
\begin{itemize}
\item[\em{(1)}]\ The action of $G_a(z')$ on $I_a$ preserves the decomposition into the equivalence classes assocoiated
to $\approx_{z'}$.
\item[\em{(2)}]\  For an element $ {[i_j]}^{\approx_{z'}}$ and an element $g$ in the subgroup of $G_a(z')$ consisting of elements 
fixing ${[i_j]}^{\approx_{z'}}$ it holds
$$
s_{i_j}'^a(g\ast y)=\mu(\Gamma(g,y),s_{i_j}'^a(y))\ \ \text{for}\ \  y\in \wt{\mathsf{O}}(z')\cap\partial X,
$$
where we recall that   $s_{i'}'^a(y)=s_{i_j}'^a(y)$ for $y\in \wt{\mathsf{O}}(z')\cap\partial X$ and $i'\in {[i_j]}^{\approx_{z'}}$.
\item[\em{(3)}]\  If for some $g\in G_a(z')$ it holds that $g( {[i_j]}^{\approx_{z'}} )   \neq  [i_j]^{\approx_{z'}}$, which means that these
sets are in fact disjoint, then $g( {[i_j]}^{\approx_{z'}} )= {[g(i_j)]}^{\approx_{z'}}$ and
$$
s_{g(i_j)}'^a(g\ast y) = \mu(\Gamma(g,y), s_{i_j}'^a(y))\ \ \text{for}\ \ y\in \wt{\mathsf{O}}(z')\cap\partial X,
$$
where we recall that   $s_{i'}'^a(y)=s_{i_j}'^a(y)$ for $y\in \wt{\mathsf{O}}(z')\cap\partial X$ and $i'\in {[i_j]}^{\approx_{z'}}$.
\item[\em{(4)}]\  For $z''\in [z]$ different from $z'$ the sections over $\wt{\mathsf{O}}(z'')\cap \partial X$ are determined completely
by the sections over $\wt{\mathsf{O}}(z')\cap \partial X$.
\end{itemize}
\end{lemma}
\begin{proof} (1)  In view of Lemma \ref{true-y} the equivalence $i\approx_z i'$ implies the equivalence
$g(i)\approx_{g(z)} g(i')$ for every $g\in G_a$.  Assuming next that $g\in G_a(z')$ and 
 $i\in {[i_j]}^{\approx{z'}}$ the fact that $g(i) \in {[i_j]}^{\approx{z'}}$ implies $g({[i_j]}^{\approx{z'}})= {[i_j]}^{\approx{z'}}$
 If $g(i)\not\in {[i_j]}^{\approx{z'}}$ the image of ${[i_j}|^{\approx{z'}}$ must be disjoint.\par
 
 \noindent(2)  Assuming that $g\in G_a(z')$ and $g({[i_j]}^{\approx{z'}})= {[i_j]}^{\approx{z'}}$ we compute
 since for $i,i'\in [i_j]$ the associated sections coincide of $\wt{\mathsf{O}}(z)\cap \partial X$
 \begin{eqnarray*}
 s_{i_j}'^a(g\ast y) &=& s_{g(i_j)}'^a(g\ast y)\\
 &=&\mu(\Gamma(g,y),s_{i_j}'^a(y))
 \end{eqnarray*}
 for $y\in \wt{\mathsf{O}}(z)\cap \partial X$.\par
 
 \noindent(3) Assume that $g( {[i_j]}^{\approx_{z'}} )   \neq  [i_j]^{\approx_{z'}}$.  Then 
 \begin{eqnarray*}
 s_{g(i_j)}'^a(g\ast y)&=& \mu(\Gamma(g,y),s_{i_j}'^a(y))
 \end{eqnarray*}
 for $y\in \wt{\mathsf{O}}(z')\cap\partial X$ by the fact that the section structure is symmetric.\par

\noindent(4) The previous points imply that the section structure is determined by knowing it over
$\wt{\mathsf{O}}(z')\cap \partial X$.  Indeed, for every class ${[i_j]}^{\approx_{z'}}$ we have the section
$s_{i_j}'^a |\wt{\mathsf{O}}(z)\cap \partial X$ which we shall denote by $s_j^{z'}$. This section is $g$-equivariant
if $g\ast z'=z'$ and $g({[i_j]}^{\approx_{z'}})={[i_j]}^{\approx_{z'}}$. If $g ({[i_j]}^{\approx_{z'}})={[i_{j'}]}^{\approx_{z'}}$
for some $j'\neq j$ and $g\in G_a(z')$, then the section $s_{j'}^{z'}$ is determined by $s_j^{z'}$ which is statement (3).

Now fix $z''\in [z]$ different from $z'$ and pick $g\in G_a$ satisfying $g\ast z'=z''$.
The decomposition of $I_a$ according to $\approx_{z''}$ is given by 
$$
I_a = {[g(i_1)]}^{\approx_{z''}}\sqcup...\sqcup {[g(i_k)]}^{\approx_{z''}}.
$$
Pick any of these classes, say $ {[g(i_j)]}^{\approx_{z''}}$. Then 
\begin{eqnarray*}
s_{g(i_j)}'^{a}(g\ast y) &=& \mu(\Gamma(g,y),s_{i_j}'^a(y))\ \ y\in \wt{\mathsf{O}}(z')\cap \partial X,
\end{eqnarray*}
which determines $s_{g(i_j)}'^{a}$ on $\wt{\mathsf{O}}(z'')\cap \partial X$ as well as those sections which have index
in $ {[g(i_j)]}^{\approx_{z''}}$.
\qed \end{proof}

What follows is a very important local extension result. 
\begin{proposition}\label{PROP14.4.3}
For fixed $a\in A$ there exists a local sc$^+$-smooth symmetric section structure ${(\wh{s}^a_{i})}_{i\in I_a}$ defined on the open set $\wh{V}_a$
and having the following properties.
\begin{itemize}
\item[\em{(1)}]\  $\wh{s}^a_{i}(y)=s'^a_i(y)$ for all $y\in V_a$ and $i\in I_a$.
\item[\em{(2)}] \  If $z\in \wh{V}_a$ and $i,i'\in I_a$ such that $i\sim_z i'$ then $\wh{s}^a_{i}(z)=\wh{s}^a_{i'}(z)$.
\end{itemize}
\end{proposition}
\begin{proof}
For the fixed $a\in A$ take  a good neighborhood system 
$$
{ \left(\wt{\mathsf{O}}(z)\right)}_{z\in \cl_X(\wh{V}_a)}.
$$
We denote by $[z]$ the $G_a$-orbit generated by $z\in \cl_X(\wh{V}_a)$. We define a $G_a$-invariant subset 
of $\wt{U}(a)$ by
$$
\wt{\mathsf{O}}([z]) =\bigcup_{z'\in [z]} \wt{\mathsf{O}}(z').
$$
Recall that by definition of a good neighborhood system, see Definition \ref{GOODNEI}, the sets  $\wt{\mathsf{O}}(z')$ and $\wt{\mathsf{O}}(z'')$ for $z',z''\in [z]$
are  disjoint unless $z'=z''$.
For every $z'\in [z]$ the index set  $I_a$ is written as disjoint union of equivalence classes $[i]^{\approx_{z'}}$.
We recall the property $g([i]^{\approx_z}) = [g(i)]^{\approx_{g(z)}}$.
For  $y\in \wt{\mathsf{O}}([z])\cap\partial X$ we find $g\in G_a$  such that $y\in \wt{\mathsf{O}}(g\ast z)$. Then, with $z'=g\ast z$
we have that 
\begin{eqnarray}\label{PQEQNX}
s'^{a}_{i}(y)=s'^{a}_{i'}(y) \ \text{if}\ y\in \wt{\mathsf{O}}(z')\cap \partial X\ \text{and}\ i\approx_{z'} i'.
\end{eqnarray}
In the next step we pick for every $[z]$ with $z\in \cl_X(V_a)$ and $i\in I_a$  sc$^+$-smooth $\wh{s}^{[z]}_i:\wh{\mathsf{O}}([z])\rightarrow W$
so that  the following holds.

\begin{itemize}
\item[(A)]\ \ If $[z]\subset  \cl_X(\wh{V}_a)\cap\partial X$ we take sc$^+$-sections $\wh{s}^{[z]}_i$, $i\in I_a$, defined on $\wt{\mathsf{O}}([z])$ such that the following holds.
\begin{itemize}
\item[(i)] \ \ $\wh{s}^{[z]}_i(y)=s'^{a}_i(y)$ for $y\in \partial X\cap \wt{\mathsf{O}}([z])$.
\item[(ii)]\ \  $\wh{s}^{[z]}_{g(i)}(g\ast y)= \mu(\Gamma(g,y),\wh{s}^{[z]}_i(y))$ for all $y\in \wt{\mathsf{O}}([z])$, $g\in G_a$, and $i\in I_a$.
\item[(iii)]\ \ $\wh{s}^{[z]}_i(y)=\wh{s}^{[z]}_{i'}(y)$ for $y\in \wt{\mathsf{O}}(z')$ and $i\approx_{z'} i'$, where $z'\in [z]$.
\end{itemize}
\item[(B)]\ \
For $[z]\subset \cl_X(\wh{V}_a)\setminus \partial X$  we take sc$^+$-sections $\wh{s}^{[z]}_i$, $i\in I_a$, defined on $\wt{\mathsf{O}}([z])$ so that the following holds.
\begin{itemize}
\item[(i)]\ \ $\wh{s}^{[z]}_{g(i)}(g\ast y) =\Gamma(g,\wh{s}^{[z]}_i(y))\ \text{for all}\ y\in \wt{\mathsf{O}}([z]),\ i\in I_a,\ g\in G_a.$
\item[(ii)]\ \  $\wh{s}^{[z]}_{i}(y)=\wh{s}^{[z]}_{i'}(y)$ for $y\in \wt{\mathsf{O}}(z')$ and $i\approx_{z'} i'$, where $z'\in [z]$.
\end{itemize}
\end{itemize}
Note that sections with the above properties can indeed be picked.  To see this recall Lemma \ref{LEMMX14413},
and  assume first
that $[z]\subset \cl_X(\wh{V}_a)\cap\partial X$. Pick a representative $z'\in [z]$ and consider $\wt{\mathsf{O}}(z')$.
For $y\in \wt{\mathsf{O}}(z')\cap \partial X$ (\ref{PQEQNX}) holds. Then $z'$ defines a partition of $I_a$ into equivalence 
classes, say $I_a = {[i_1]}^{\approx_{z'}}\sqcup... \sqcup {[i_k]}^{\approx_{z'}}$ and 
$$
g({[i_j]}^{\approx_{z'}})={[g(i_j)]}^{\approx_{g(z')}}= {[g(i_j)]}^{\approx_{z'}}\ \ \text{if}  \ \ g\ast z'= z'.
$$
  We pick the class ${[i_1]}^{\approx_{z'}}$ and we take an extension of $s_{i_1}'^a|(\wt{\mathsf{O}}(z')\cap \partial X)$
  to $\wt{\mathsf{O}}(z')$. We average this extension with respect to all $g\in G_a(z')$
  which fix ${[i_1]}^{\approx_{z'}}$. The averaged section is still an extension since we extend a section with the same symmetry.
At this point we have constructed 
\begin{eqnarray}\label{EQNC14}
\wh{s}_i^{z'}\ \ \text{of}\ \ W|\wt{\mathsf{O}}(z')\ \ \text{for}\ \ i \in {[i_1]}^{\approx_{z'}}.
\end{eqnarray}
Next we consider ${[i_2]}^{\approx_{z'}}$. If there exists $g\in G_a(z')$ with $g({[i_1]}^{\approx_{z'}})={[i_2]}^{\approx_{z'}}$
we take such a $g$ and define $\wh{s}_{g(i)}^{z'}$ by
$$
\wh{s}_{g(i)}^{z'}(g\ast y) =\mu(\Gamma(g,y),\wh{s}_i^{z'}(y)).
$$
First of all we note that all sections are the same independently of the choice of $i\in {[i_1]}^{\approx_{z'}}$.
The definition also does not depend on the choice of $g$. In fact if there are two such elements $g_1,g_2\in G_a(z')$
we note that $ g_2^{-1}(g_1({[i_1]}^{\approx_{z'}}))={[i_1]}^{\approx_{z'}}$ and since the sections defined in 
(\ref{EQNC14}) have the invariance property with repespect to such elements $g_2^{-1}\circ g_1$ the claim follows.
Hence we obtain
\begin{eqnarray}
\wh{s}_i^{z'}\ \ \text{of}\ \ W|\wt{\mathsf{O}}(z')\ \ \text{for}\ \ i \in {[i_2]}^{\approx_{z'}}.
\end{eqnarray}
If there is no $g\in G_a(z')$ mapping ${[i_1]}^{\approx_{z'}}$ to ${[i_2]}^{\approx_{z'}}$ we deal with the situation
in the same way we obtained (\ref{EQNC14}). Then we take ${[i_3]}^{\approx_{z'}}$ and distinguish the case
if it is related by a group element to $ {[i_1]}^{\approx_{z'}}$ or ${[i_2]}^{\approx_{z'}}$. If that is the case
push forward the data.  If no such relationship exists we deal with the situation as in the case of 
${[i_1]}^{\approx_{z'}}$.  After a finite number of steps we have constructed 
\begin{eqnarray}
\wh{s}_i^{z'}\ \ \text{of}\ \ W|\wt{\mathsf{O}}(z')\ \ \text{for}\ \ i \in I_a.
\end{eqnarray}
This data has the following properties.
\begin{itemize}
\item[(i)]\ \ \ If $g\in G_a(z')$ then $\wh{s}_{g(i)}^{z'}(g\ast y)=\mu(\Gamma(g,y),\wh{s}_{i}^{z'}(y))$ for $y\in \wt{\mathsf{O}}(z')$.
\item[(ii)]\ \ \  If $i\approx_{z'} i'$ then $ \wh{s}_{i}^{z'}=\wh{s}^{z'}_{i'}$.
\item[(iii)]\ \ \  $\wh{s}^{z'}_i(y)=s_i'^a(y)$ for $y\in \wh{\mathsf{O}}(z')\cap X$ and $i\in I_a$.
\end{itemize}

If we take $z''\neq z'$ we fix $g\in G_a$ with $g\ast z'=z''$ and push the data forward. Again the choice of $g$ does not matter
since for two possible choices $g_1$ and $g_2$ we have that $g_2^{-1}\circ g_a\in G_a(z')$ and the data associated
with $z'$ was constructed in such a way that it has an invariance property with respect to the action of $G_a(z')$
and the equivalence relation $\approx_{z'}$.  At this point we have defined $\wh{s}_{i}^{z'}$ for all $z'\in [z]$ which defines
$\wh{s}^{[z]}_i$ of $W|\wt{\mathsf{O}}$.  This gives us a system of sections satisfying $(\mathsf{A})$.
The construction in the case $(\mathsf{B})$ is similar, but easier, since we do not have to extend  sections given on the boundary.

At this point we have constructed for $a\in A$ and a given good neighborhood system for $\cl_X(\wh{V}_a)\subset \wt{U}(a)$,
for every $\wt{\mathsf{O}}([z])$ a system of sections ${(\wh{s}_i^{[z]})}_{i\in I_a}$ with the properties listed above.
The sets 
$$
{(\wt{\mathsf{O}}([z])}_{[z]\in |\cl_X(\wh{V}_a)|}
$$
 define an open covering of $\cl_X(\wh{V}_a)$ by $G_a$-invariant subsets contained in
$\wt{U}(a)$. We take a sub-ordinate functorial sc-smooth partition of unity ${(\sigma_{[z]})}_{\{[z]\}}\cup\{\sigma_\infty\}$ 
indexed by the set $\{[z]\} \cup \{\infty\}$. Here $\infty$ is the index of a map supported in $\wt{U}(a)\setminus \wh{V}_a$.
We disregard this map and only use the  ${(\sigma_{[z]})}_{[z]}$, which satisfy
$$
\text{supp}(\sigma_{[z]})\subset \wt{\mathsf{O}}([z]).
$$
Of course most of these maps vanish. Then we introduce  the maps $\wh{s}^a_i$, $i\in I_a$, which are defined on $\wh{V}_a$ by
$$
\wh{s}^a_i(y)=\sum_{[z]} \sigma_{[z]}(y)\cdot \wh{s}^{[z]}_i(y),\ y\in \wh{V}_a.
$$
By construction it is obvious that the following two properties hold true.
\begin{itemize}
\item[(i)] \ \ \ $\wh{s}^a_{g(i)}(g\ast y) =\mu(\Gamma(g,y),\wh{s}^a_i(y)),\ y\in \wh{V}_a,\ i\in I_a$.
\item[(ii)] \ \ \ $\wh{s}^a_i(y)=s'^{a}_i(y)$ for $y\in \partial X\cap \wh{V}_a$ and $i\in I_a$.
\end{itemize}
Assume next that $i\sim_y i'$ for some $y\in \wh{V}_a$. There are only  finitely many $\sigma_{[z_1]},..,\sigma_{[z_\ell]}$
which are nonzero on
 $y$ and by  construction
$$
\text{supp}(\sigma_{[z_j]})\subset \wt{\mathsf{O}}([z_j]).
$$
For every $[z_j]$ pick $z_j'\in [z_j]$ such that $y\in \wt{\mathsf{O}}(z_j')$.
Since $[i]^{\sim_y}\subset [i]^{\approx_{z_j'}}$ for every $z_j'$, see Definition \ref{GOODNEI} (3),  and $\wh{s}^{[z_j]}_{i}(y)=\wh{s}^{[z_j]}_{i'}(y)$ for
$i\approx_{z_j'} i'$ it follows that
$$
\wh{s}^a_i(y) = \wh{s}^a_{i'}(y) \ \text{provided}\ i\sim_y i'.
$$
\qed \end{proof}
The next section uses the this local result to construct global extensions.

\section{Extension of  the Sc\texorpdfstring{$^+$}{ppl}-Multisection}
Let us summarize the data we are dealing with. We started with an sc$^+$-multisection functor
$$
\Lambda:W|\partial X\rightarrow {\mathbb Q}^+,
$$
and we have fixed a representative  
\begin{eqnarray}\label{ERX1450}
(\Lambda,\mathfrak{U}',\mathfrak{S}',\tau')
\end{eqnarray}
 of a structured version.
We fixed $({\wt{U}(x))}_{x\in\partial X}$ admitting the natural $G_x$-action and satsifying $\wt{U}(x)\subset \wt{U}'(x)$.
We have constructed ${(\wh{V}_a)}_{a\in A}$, where $\wh{V}_a$ is open in $X$ and $\cl_X(\wh{V}_a)\subset \wt{U}(a)$.
Moreover, the whole collection ${(|\wh{V}_a|)}_{a\in A}$ is a locally finite open covering of $|\partial X|$. On the associated $\wh{\bm{V}}_{a,a'}$
we have defined a canonical $\wh{\tau}_{a,a'}:\wh{\bm{V}}_{a,a'}\rightarrow \mathfrak{J}(a,a')$ obtained from
the data provided by $\tau'$.
In addition for every $a\in A$ we are given a symmetric sc$^+$-section structure 
\begin{eqnarray}\label{ERX14500}
{(\wh{s}^a_{i})}_{i\in I_a}\ \text{on}\ \ \wh{V}_a
\end{eqnarray}
so that $\wh{s}^a_{i}(y)=s'^a_i(y)$ for $y\in V_a$ and $i\in I_a$ and $\wh{s}^a_{i}(y)=\wh{s}^a_{i'}(y)$ for $y\in \wh{V}_a$
and $i\sim_y i'$.  We shall use this data to define in a suitable way a structurable  extension $\Lambda':W\rightarrow {\mathbb Q}^+$
of $\Lambda$. In this section we first shall construct $\Lambda'$ in a quite specific way and the structurability is addressed in the next section. In  a first step  we  shall decompose $\Lambda:W|\partial X\rightarrow {\mathbb Q}^+$ for which we need the following lemma.

\begin{lemma}\label{ABBXZ0}
We work under the assumptions of Theorem \ref{p-main-p}. In particular $X$ is tame and admits sc-smooth partitions of unity as an ep-groupoid. Then there is a collection of  sc-smooth functors $\beta_a:X\rightarrow [0,1]$, $a\in A$,  having the following properties.
\begin{itemize}
\item[\em{(1)}]\  For $a\in A$ the sc-smooth map $\beta_a$ has support in $\pi^{-1}(\pi(\wh{V}_a))$.
\item[\em{(2)}]\   For every $y\in X$ there exists a saturated open neighborhood $Q(y)$ in $X$ such that there are only finitely many $a\in A$ for which
$\beta_a|Q(y)\neq 0$.
\item[\em{(3)}]\  For every $y\in \partial X$ it holds $\sum_{a\in A} \beta_a(y)=1$.
\end{itemize}
\end{lemma}
\begin{proof} 
The collection of sets ${(\pi^{-1}(|\wh{V}_a|))}_{a\in A}$
 together with the open set $\wh{V}_\infty:= X\setminus\partial X$
defines a locally finite open covering of $X$ and we take a subordinate family of partition of unity functors. From these 
we only use the ${(\beta_a)}_{a\in A}$.
\qed \end{proof}
We note that for $a\in A$ it holds
$$
\text{supp}( |\beta_a|)\subset  |\wh{V}_a|,
$$
and consequently
\begin{eqnarray}
\cl_X(\{ y\in \wt{U}(a)\ |\ \beta_a(y)\neq 0\}) \subset \wh{V}_a\ \text{for}\ a\in A.
\end{eqnarray}
The $(\beta_a)$ restricted to $\partial X$ define an sc-smooth partition of unity.  We define structured ${\Lambda}_a:W|\partial X\rightarrow {\mathbb Q}^+$
for $a\in A$  by 
$$
\Lambda_a =(\beta_a|\partial X)\odot \Lambda.
$$
The  structuring data we are going to use comes from (\ref{ERX1450}).
 We take as good system $\mathfrak{U}'$, as index set $\mathfrak{I}$, and as correspondence $\tau'$. The section structure
$\mathfrak{S}_a' ={(\mathfrak{s}_{a,x})}_{x\in \partial X}$ is given by
$$
\mathfrak{s}_{a,x}:\ \ s'^{x}_{a,i}(y)=\beta_a(y)\cdot s'^{x}_i(y)\ \text{for}\ x\in \partial X,\ i\in I_{x},\ y\in U'(x).
$$
We note the following obvious fact.
\begin{lemma}
The structured  $\wt{\Lambda}:=(\Lambda,\mathfrak{U}',\mathfrak{S}',\tau')$ is the commensurable locally finite sum
of the $\wt{\Lambda}_a:=({\Lambda}_a,\mathfrak{U}',\mathfrak{S}'_a,\tau')$
\begin{eqnarray}\label{LEMMX1452}
\wt{\Lambda}=\boxplus_{a\in A} \wt{\Lambda}_a.
\end{eqnarray}
\qed
\end{lemma}
\begin{remark} It is also important to notice that $\Lambda_a$ is completely determined by the symmetric sc$^+$-section structure 
$$
{((\beta_a|U'(a))\cdot s'^a_i)}_{i\in I_a}.
$$
In fact, it is even determined by ${((\beta_a|V_a)\cdot s'^a_i)}_{i\in I_a}$. Of course, as a structured sc$^+$-section functor 
the global knowledge of $\tau'$ is needed. For example we have used it for taking the commensurable sum in (\ref{LEMMX1452}).
\qed
\end{remark} 

We begin now with the important constructions needed  to extend $\Lambda$ to a saturated open neighborhood
of $\partial X$. We fix $a\in A$ and consider $\Lambda_a:W|\partial X\rightarrow {\mathbb Q}^+$.
The local section structure for $\Lambda_a$ is given at $x\in \partial X$ on the set $U'(x)$ by
$$
{(\beta_a\cdot s'^x_i)}_{i\in I_x}.
$$
First we shall only  work on the open sets $\wh{V}_b$, where $b$ varies in $A$. Most of the construction
takes place on $\wh{V}_a$, and we can transport data to other $\wh{V}_b$ with the help of the $\wh{\tau}_{a,b}$. Of course, transporting data around relies on using a choice of morphisms, as well as the correspondences
$\wh{\tau}$. The careful analysis of the implicit structures is needed to avoid ambiguities.

\begin{proposition}
Recall the family ${(\wh{s}_i^a)}_{i\in I_a}$, $a\in A$, see (\ref{ERX14500}), obtained through the local extension procedure in the previous section.
For each given  $a\in A$ there exists for every $b\in A$  a symmetric sc$^+$-section structure ${(\wt{s}_{a,i}^b)}_{i\in I_b}$
where $\wt{s}_{a,i}^b$ is defined on $\wh{V}_b$ and the following holds.
\begin{itemize}
\item[\em{(1)}]\  $\wt{s}_{a,i}^a(y) =\beta_a(y)\cdot \wh{s}^a_{i}(y)$ for $i\in I_a$ and $y\in \wh{V}_a$.
\item[\em{(2)}]\   $\text{supp}(\wt{s}_{a,i}^a)$ is closed in $X$ and contained in $\wh{V}_a$ for $i\in I_a$.
\item[\em{(3)}] \  $\wt{s}^b_{a,i}(y)=\beta_a(y)\cdot s'^{b}_{i}(y)$ for $y\in V_b$,  $b\in A$.
\item[\em{(4)}]\   For  $b,b'\in A$ and $\phi\in \wh{\bf V}_{b,b'}$ with $\wh{\tau}_{b,b'}(\phi)$ given by $[I_{b}\stackrel{\alpha(k)}{\twoheadleftarrow} I\stackrel{\beta(k)}{\twoheadrightarrow} I_{b'}]$ it holds for $k\in I$
$$
\wt{s}^{b'}_{a,\alpha(k)}(t(\phi))=\mu(\phi,\wt{s}^{b}_{a,\beta(k)}(s(\phi))).
$$
\end{itemize}
\end{proposition}
The proof of this result is subtle mainly due to the requirement (4). In fact (4) required the lengthy preparations introducing the equivalence relations  $\sim_z$ and $\approx_z$ in the previous section.
\begin{proof}
In Proposition \ref{PROP14.4.3} we have already constructed a symmetric sc$^+$-section structure $\wh{s}^a_{i}$, $i\in I_a$, on $\wh{V}_a$ which extends 
$s'^a_i$, $i\in I_a$.   We define
$\wt{s}^a_{a,i}$ for $i\in I_a$ by
$$
\wt{s}^a_{a,i}(y)= \beta_a(y)\cdot \wh{s}^a_{i}(y),\ y\in \wh{V}_a.
$$
This defines $\mathfrak{s}^a_a$ and 
 (1) and (2) hold, and in addition (3) for the special case $b=a$.

In a next step we define, keeping the initial $a\in A$ fixed related section structures for all $\wh{V}_b$, where $b\in A$. 
If $b\in A$ such that there is no morphism
starting in $\wh{V}_a$ and ending in $\wh{V}_b$ we associate to every $i\in I_b$ the zero section on $\wh{V}_b$ which defines the symmetric sc$^+$-section structure
$\mathfrak{s}^b_{a}$.  
Suppose that $b\in A$ is such that there exists a morphism starting in $\wh{V}_a$ and ending in $\wh{V}_b$.
In a first step we define $s^b_{a,i'}$ for $i'\in I_b$ on the subset of $\wh{V}_b$ which can be reached by morphisms
starting  in  $\wh{V}_a$ and ending in $\wh{V}_b$. Subsequently we shall see that the extension by zero over the rest 
of the set $\wh{V}_b$ is sc$^+$-smooth.
Pick $\phi\in \wh{\bm{V}}_{a,b}$ and take 
$$
I_a\stackrel{\alpha}{\twoheadleftarrow} I\stackrel{\beta}{\twoheadrightarrow} I_b
$$
 representing
$\wh{\tau}_{a,b}(\phi)$. For given $i'\in I_b$ pick $k\in I$ satisfying $\beta(k)=i'$ and with $i=\alpha(k)$ define
$$
\wt{s}^b_{a,i'}(t(\phi))=\mu(\phi,\wt{s}^a_{a,i}(s(\phi))).
$$
We have to verify that the definition is independent of the choices. Suppose $\psi\in \wh{\bm{V}}_{a,b}$ with $t(\psi)=t(\phi)$ and $\wh{\tau}_{a,b}(\psi)$ is represented
by 
$$
 I_a\stackrel{\gamma}{\twoheadleftarrow} J\stackrel{\delta}{\twoheadrightarrow} I_b.
 $$
  We assume that we have picked $j\in J$ satisfying
$\delta(j)=i'$. We define $e=\gamma(j)$. Using this data the alternative definition for $\wt{s}^b_{a,i'}(t(\phi))$ is given by 
$\mu(\psi,\wt{s}^a_e(s(\psi)))$. We need to show  the equality 
\begin{eqnarray}\label{meretx}
\mu(\phi,\wt{s}^a_{a,i}(s(\phi))) = \mu(\psi,\wt{s}^a_{a,e}(s(\psi))).
\end{eqnarray}
We can write $\psi^{-1}\circ \phi =\Gamma(g,s(\phi))$ for a uniquely determined $g\in G_a$. Hence
$$
\phi=\psi\circ\Gamma(g,s(\phi)).
$$
Since $\phi=\psi\circ\Gamma(g,s(\phi))=\psi\circ\Gamma(g,g^{-1}\ast s(\psi))$ we have that $\phi=\psi\ast g$ implying
$$
\wh{\tau}_{a,b}(\phi)=\wh{\tau}_{a,b}(\psi\ast g) =\wh{\tau}_{a,b}(\psi)\ast g.
$$
Therefore
$$
[I_a\stackrel{\alpha}{\twoheadleftarrow} I\stackrel{\beta}{\twoheadrightarrow} I_b] =[ I_a\stackrel{g^{-1}\circ \gamma}{\twoheadleftarrow} J\stackrel{\delta}{\twoheadrightarrow} I_b],
$$
so that without loss of generality we may assume that $\alpha=g^{-1}\circ \gamma$ and $\beta=\delta$.  With $z=s(\psi)$ we note that
$\psi , \phi^{-1}, \Gamma(g,s(\phi))$ is a $z$-loop and with 
$$
\wh{\tau}_{a,a}(\Gamma(g,s(\phi)))=[I_a\stackrel{Id}{\twoheadleftarrow} I_a\stackrel{g}{\twoheadrightarrow} I_a]
$$
we see that we obtain the $z$-arc $e,i',i,g(i)$ via
$$
e =\gamma(j),\ \delta(j) = i' = \beta(k),\  \alpha(k) = i,\ g(i).
$$
 so that 
\begin{eqnarray}\label{meretxXX}
 e\sim_z g(i).
\end{eqnarray}
Using that 
$$
i' = \beta(k),\ i =\alpha(k) = g^{-1}(\gamma(k)), \ \text{and}\ e=\gamma(j).
$$
we compute
\begin{eqnarray*}
&& \mu(\phi,\wt{s}^a_{a,i}(s(\phi)))\\
&=& \mu(\psi\circ\Gamma(g,s(\phi)),\wt{s}^a_{a,i}(s(\phi)))\\
&=& \mu(\psi,\mu(\Gamma(g,s(\phi)),\wt{s}^a_{a,i}(s(\phi))))\\
&=& \mu(\psi,\wt{s}^a_{a,g(i)}(g\ast s(\phi)))\\
&=&\mu(\psi,\wt{s}^a_{a,g(i)}(s(\psi))).
\end{eqnarray*}
We need to show that $ \mu(\psi,\wt{s}^a_{a,g(i)}(s(\psi)))=\mu(\psi,\wt{s}^a_{a,e}(s(\psi)))$. However, since $e\sim_{s(\psi)} g(i)$, and $z=s(\psi)$, 
see (\ref{meretxXX}),  it follows
that 
\begin{eqnarray*}
\wt{s}^a_{a,g(i)}(s(\psi))&=&\beta_a(s(\psi))\cdot \wh{s}^a_{g(i)}(s(\psi))\\
&=&\beta_a(s(\psi))\cdot \wh{s}^a_{e}(s(\psi))\\
&=&\wt{s}^a_{a,e}(s(\psi))
\end{eqnarray*} 
which proves this  identity, and therefore that $\wt{s}^b_{a,i'}$ is well-defined on all points of $\wh{V}_b$  reachable from
$\wh{V}_a$.  

Next let $z'\in \wh{V}_b$  be in the closure of reachable points from $a$. Hence we find a sequence of morphisms 
$\psi_k$ satisfying $s(\phi_k)\in \wh{V}_a$, $t(\phi_k)\in \wh{V}_b$ and $z_k':= t(\phi_k)\rightarrow z'$. Using properness
we may assume with out loss generality that $\phi_k\rightarrow \phi$ with $s(\phi)\in \cl_X(\wh{V}_a)\subset \wt{U}(a)$ and
$t(\phi)=z'$. If $s(\phi)\in \wh{V}_a$ the point $z'$ is reachable and $\wt{s}^b_{a,i'}$ is already defined. So assume $s(\phi)\not\in \wh{V}_a$.
Then for large $k$ we  already have that $\wt{s}^a_{a,i}$ vanishes in an open neighborhood of $s(\phi_k)$.
Hence it follows that  $\wt{s}^b_{a,i'}$ where defined vanishes on the reachable points near $z'$. Hence we can extend
$\wt{s}^b_{a,i'}$ over the rest  of $\wh{V}_b$ by $0$.  This defines ${(\wt{s}^b_{a,i'})}_{i'\in I_b}$. 
At this point we have constructed for a given $a\in A$
symmetric sc$^+$-section structures ${(\wt{s}_{a,i}^b)}_{i\in I_b}$ for every $b\in A$  such that
$$
\wt{s}_{a,i'}^b(t(\psi))=\mu(\psi,\wt{s}_{a,i}^a(s(\psi)))
$$
provided $i'$ and $i$ correspond under $\wh{\tau}_{a,b}(\psi)$.   The last statement proves (4) for a special case.

In order to prove (4) in generality, still for fixed $a$,  consider $b, b'\in A$ and let
$\psi\in \wh{\bm{V}}_{b,b'}$.   We have to consider several cases.
Assume first that $s(\psi)$ or $t(\psi)$ is reachable from $\wh{V}_a$. In this case both are reachable 
and we find $\phi_1,\phi_2$ with
$$
s(\phi_1)=s(\phi_2)\in\wh{V}_a\ \text{and}\ \psi\circ \phi_1=\phi_2.
$$
In this case, if $i\in I_{b}$ and $i'\in I_{b'}$ correspond under $\wh{\tau}_{b,b'}(\psi)$, we can pick $j\in I_a$
corresponding under $\wh{\tau}_{a,b}(\phi_1) $ to $i$, and  pick $j'\in I_a$ corresponding under $\wh{\tau}_{a,b'}(\phi_2)$ to $i'$.
We note that $j\sim_{s(\phi_1)} j'$. This implies using previous definitions 
\begin{eqnarray*}
&& \wt{s}^{b'}_{a,i'}(t(\psi))\\
&=& \mu(\phi_2,\wt{s}^a_{a,j'}(s(\phi_2)))\\
&=&\mu(\psi\circ\phi_1,\wt{s}^a_{a,j'}(s(\phi_2)))\\
&=&\mu(\psi,\mu(\phi_1,\wt{s}^a_{a,j'}(s(\phi_2))))\\
&=& \mu(\psi,\mu(\phi_1,\wt{s}^a_{a,j'}(s(\phi_1))))\\
&=&\mu(\psi,\mu(\phi_1,\wt{s}^a_{a,j}(s(\phi_1))))\\
&=&\mu(\psi,\wt{s}^b_{a,i}(s(\psi))).
\end{eqnarray*}
At this point we have constructed for  given $a$  a symmetric sc$^+$-section structures  ${(\wt{s}^b_{a,i})}_{i\in I_b}$ on $\wh{V}_b$
for every $b\in A$. Moreover, these sc$^+$-section structures are compatible under the correspondences
$\wh{\tau}_{b,b'}:\wh{\bm{V}}_{b,b'}\rightarrow \mathfrak{J}(b,b')$ and it holds that
$$
\wt{s}^{b}_{a,i}(y) =\beta_a(y)\cdot s'^{b}_i(y)\ \text{for}\ y\in V_b,\ i\in I_b.
$$
\qed \end{proof}
 Define for every $b\in A$
the symmetric sc$^+$-section structure on $\wh{V}_b$
$$
{(\wt{s}^b_i)}_{i\in I_b}
$$
by setting
$$
\wt{s}^b_i=\sum_{a\in A} \wt{s}_{a,i}^b.
$$
This is a locally finite sum and therefore well-defined.  It also extends the section structure of $\Lambda$ on $V_b$.
By construction it is compatible with $\wh{\tau}$. 

At this point we can define for every $b\in A$ 
an sc$^+$-multisection functor ${\Lambda}'_b$ on $W|\wh{V}_b$ by
$$
{\Lambda}'_b(w)=\frac{1}{|I_b|} \cdot |\{i\in I_b\ |\ \wt{s}_i^b(P(w))=w\}|.
$$
If $\phi$ is a morphism starting in $\wh{V}_b$ and ending in $ \wh{V}_{b'}$ it follows with $w'=\mu(\phi,w) $ that
${\Lambda}'_{b'}(w')={\Lambda}'_b(w)$ so that the sc$^+$-multisection functors are compatible.

The symmetric sc$^+$-section structures ${(\wt{s}^b_i)}_{i\in I_b}$ are related by the $\wh{\tau}_{b,b'}:\bm{\wh{V}}_{b,b'}\rightarrow\mathfrak{J}(b,b')$. We shall construct from this data a structurable  sc$^+$-section functor $\Lambda':W\rightarrow {\mathbb Q}^+$
which extends $\Lambda:W|\partial X\rightarrow {\mathbb Q}^+$.  We note that the structure data for $\Lambda'$
will in general not restrict to the one originally fixed for $\Lambda$.

Recall that ${(|\wh{V}_b|)}_{b\in A}$ is an open covering of $|\partial X|$. For $|x|\in |X|$ not contained in 
$\bigcup _{b\in A} |\wh{V}_b|$ we can take a representative $x$ and an open neighborhood $\wh{V}_x$ invariant under the natural action,
with the properness property, so that $|\wh{V}_x|$ only contains isomorphism classes of points on which the sections vanish.
We equip $\wh{V}_x$ with the zero section $(s_1^x)$ as section structure and define $I_x=\{1\}$. The correspondence
for $b\in A$ and such a $x$ we define by 
$$
\wh{\tau}_{b,x}(\psi)=[I_b\stackrel{Id}{\twoheadleftarrow} I_b\twoheadrightarrow \{1\}].
$$
Similarly for the reversed pair $(b,x)$. For $x,x'$ as just described we define $\wh{\tau}_{x,x'}(\psi)=[\{1\}\twoheadleftarrow \{1\}\twoheadrightarrow \{1\}]$.

At this point we have defined for a suitable $A'\subset X$, containing $A$, open subsets $\wh{V}_x$
with the natural $G_x$-action and symmetric sc$^+$-section structures compatible for some correspondence.
Moreover
$$
\bigcup_{x\in A'} |\wh{V}_x| =|X|.
$$
Hence the data defines a global sc$^+$-multisection $\Lambda':E\rightarrow {\mathbb Q}^+$ extending
$\Lambda$. By Theorem \ref{criterionn}, which is proved below,  $ \Lambda'$ is structurable. The domain support of $\Lambda'$
has an orbit space  contained in the open neighborhood $\bigcup_{a\in A} |\wh{V}_a|$ of $|\partial X|$.
Define
$$
\wt{U}= \pi^{-1}\left(\bigcup_{a\in A} |\wh{V}_a|\right).
$$
Define the map $h:|X|\rightarrow {\mathbb R}^+$ by
$$
h(|x|) = \text{max} \left(\{N(s_i^x(x))\ |\ i\in I_x\}\right).
$$
Then $h$ is a continuous map and satisfies for $|x|\in |\partial X|$
$$
h(|x|) =N(\Lambda)(x) < f(x).
$$
We can find an sc-smooth functor $\sigma:X\rightarrow {\mathbb R}^+$ such that $\sigma(x)=1$ for $x\in \partial X$
having support in $\wt{U}$ such that 
$$
h(|x|)\cdot \sigma(x)\leq f(x)\ \text{for all}\ x\in X.
$$
Now consider $\sigma\odot \Lambda'$ which is the required extension. Since $\Lambda'$ is structurable according to Theorem \ref{criterionn} the same is true
for $\sigma\odot\Lambda'$. This completes modulo the structurability assertion the proof of Theorem \ref{p-main-p}.

The following result is a very useful criterion for structurability and applies to the situation in the previous section.
\begin{theorem}\label{criterionn}
Let  $P:W\rightarrow X$ be a strong bundle over the ep-groupoid $X$ having a paracompact orbit space $|X|$ and let
$\Lambda:W\rightarrow {\mathbb Q}^+$  be  A $\ssc^+$-multisection functor. We assume that there exists a subset $A$ of $X$ and for every $a\in A$ an open neighborhood $U(a)$ and an open subset $V_a\subset U(a)$ with the following properties (a)--(c): 
\begin{itemize}
\item[{\em (a)}] \ \ $U(a)$ is $G_a$-invariant and the target map $t: s^{-1}(\cl_X (U(a)) \to X$ is proper. 
\item[{\em (b)}] \ \ $V_a$ is $G_a$-invariant and $\cl_{|X|}(|V_a|)\subset |U(a)|$.
\item[{\em (c)}]\ \ $\bigcup_{a\in A} |V_a| =|X|$.
\end{itemize}
In addition we assume that the following two properties hold true.
\begin{itemize} 
\item[{\em (d)}]\ \ For every $a\in A$, there exists a symmetric $\ssc^+$-section structure $\mathfrak{s}^a={(s_i^a)}_{i\in I_a}$ on $V_a$ indexed by a finite set $I_a$ such that if $w\in W$ satisfies $P(w)\in V_a$, then 
$$
\Lambda(w)=\frac{1}{|I_a|} |\{i\in I_a\, \vert \,  s^a_i(P(w))=w\}|.
$$
\item[{\em (e)}] \ \ For every $a,a'\in A$, with $\bm{V}_{a,a'}:=\{\phi\in\bm{X}\ |\ s(\phi)\in V_a,\ t(\phi)\in V_{a'}\}$, 
there exists a 
 correspondence 
 $$
 \tau_{a,a'}:\bm{V}_{a,a'} \rightarrow \mathfrak{J}(a,a')
 $$
  relating the section structures $\mathfrak{s}^a$ and $\mathfrak{s}^{a'}$ over $V_a$ and $V_{a'}$ and only taking finitely many values.
\end{itemize}
Then  $\Lambda$ is a $\ssc^+$-structurable multisection functor.
\end{theorem}
\begin{proof}
Define  the open subset $V=\bigcup_{a\in A} V_a$ of $X$.  It has the  property that $|V|=|X|$  and inherits the structure of an ep-groupoid
for which the inclusion functors $F: V\rightarrow X$ is an equivalence of ep-groupoids. In the same way the inclusion functor
$\Phi: W|V\rightarrow W$ is a strong bundle equivalence. Clearly 
$$
\Phi_\ast (\Lambda|(W|V)) =\Lambda
$$
and in view of this fact, employing Theorem \ref{pushforwardX},   it suffices to show that $\Lambda|(W|V)$ is structurable.

In order to show the structurability of $\Lambda|(W|V)$ there is no loss of generality viewing $\Lambda$
as an sc$^+$-multisection functor on $W\rightarrow V$. Then the data given in the hypothesis of the theorem provides
us with an open covering ${(V_a)}_{a\in A}$. First we shall pick a map $V\rightarrow A:x\rightarrow a_x$ so that 
$x\in V_{a_x}$ and we pick an open neighborhood $Q(x)$ so that for the natural $G_{a_x}$-action
on $U(a_x)$ it holds that $g\ast x=g$ implies $g\ast Q(x)=Q(x)$. We define a local symmetric sc$^+$-section structure 
on $Q_x$ by restricting  the data from $V_{a_x}$. More precisely we define $\bar{I}_x:= I_{a_x}$
$$
\bar{\mathfrak{s}}_x\colon\ \ \  \bar{s}_i^x := s^{a_x}_i\ \text{for}\ \ i\in I_{a_x}.
$$
The group $G_x$ acts on $\bar{I}_x$  as follows. For $g\in G_x$ there exists a unique $h_g\in G_{a_x}$ with $\Gamma(h_g,x)=g$
and we define $g(i):= h_g(i)$, so that we act through a subgroup of $G_{a_x}$.   

Given $x,x'\in V$ we consider $\bm{Q}(x,x')$ which is an open subset of $\bm{V}(a_x,a_{x'})$ which allows us to restrict
$\tau_{a_x,a_{x'}}$. Hence we define
$$
\bar{\tau}_{x,x'} :\bm{Q}(x,x')\rightarrow {\mathfrak{J}}(x,x')
$$
as follows.   We note that by construction ${\mathfrak{J}}(x,x'):=\mathfrak{J}(a_x,a_{x'})$ and 
$$
\bar{\tau}_{x,x'}(\phi) :=\tau_{a_x,a_{x'}}(\phi).
$$
We gather the data as follows. We define $\bar{\mathfrak{U}}={(Q(x))}_{x\in V}$,  $\bar{\mathfrak{S}}={(\bar{\mathfrak{s}}_x)}_{x\in V}$, and
$\bar{\tau} ={(\bar{\tau}_{x,x'})}_{(x,x')\in V\times V'}$.

We show that $(\Lambda,\bar{\mathfrak{U}}, \bar{\mathfrak{S}},\bar{\tau})$ is a structured version of $\Lambda$, which will
complete the proof of the theorem. First of all we have for every $x\in V$ an open neighborhood $Q(x)$ with the natural $G_x$-action
and a symmetric sc$^+$-section structure $\bar{\mathfrak{s}}_x$ representing $\Lambda$. Moreover, given 
$x,x'\in V$ we have a correspondence 
$$
\bar{\tau}_{x,x'}:\bm{Q}(x,x')\rightarrow \mathfrak{J}(x,x')
$$
which only takes finitely many values,  and is locally constant.  Moreover, by construction if $\phi\in \bm{V}(x,x')$
then $\bar{\tau}_{x,x'}(\phi)= \tau_{a_x,a_{x'}}(\phi)$. Assume the latter is represented by the diagram
$$
I_{a_x}\stackrel{\alpha}{\twoheadleftarrow} I\stackrel{\beta}{\twoheadrightarrow} I_{a_{x'}}.
$$
Then by hypothesis $s^{a_{x'}}_{\beta(k)}(t(\phi))=\mu(\phi,s^{a_x}_{\alpha(k)}(s(\phi)))$ which precisely implies
$$
\bar{s}_{\beta(k)}^{x'}(t(\phi))=\mu(\phi,\bar{s}^{x}_{\alpha(k)}(s(\phi))).
$$
Moreover, using previous definitions
\begin{eqnarray*}
\bar{\tau}_{x,x}(\Gamma_x(h,y))=\tau_{x,x}(\Gamma_{a_x}(g,y))
= \left[I_{a_x}\stackrel{Id}{\twoheadleftarrow} I_{a_x}\stackrel{g}{\twoheadrightarrow} I_{a_x}\right]
= \left[\bar{I}_x\stackrel{Id}{\twoheadleftarrow} \bar{I}_{x}\stackrel{h}{\twoheadrightarrow} \bar{I}_{x}\right]
\end{eqnarray*}
\qed \end{proof}

\section{Remarks on Inductive Constructions}
The main result Theorem \ref{p-main-p} proved in the previous sections 
stated that a structurable sc$^+$-multisection defined on the boundary of a tame ep-groupoid  can be extended
as a structurable sc$^+$-multisection. In applications, very often, the perturbations 
are constructed inductively and additional problems arise. For example, if the underlying ep-groupoid
has a boundary with corners and it is face-structured the inductive procedure might define
perturbations over some of the faces so that the (global) sc$^+$-multisection on the boundary has to be assembled
from the pieces defined on the faces. Of course, that can only be accomplished provided the sc$^+$-multisections
 on the intersection of faces 
coincide. In order to extend the sc$^+$-multisection it is important to have a structuring of the assembled
sc$^+$-multisection.  The structuring data for  the sc$^+$-multisections on the faces
 is in general a consequence of the inductive
scheme, and a matter of the concrete situation and will not be discussed here. 
Here we shall restrict ourselves to highlight the compatibility issue of the data coming from
the different faces if we have a boundary with corners. We shall also restrict ourselves
to the case where $X$ is a face-structured ep-groupoid. More general
situations are possible, but left to future discussions.

Let $X$ be a face-structured tame ep-groupoid and $P:W\rightarrow X$ a strong bundle.  
In this case a face $F$ is saturated and has the structure of a sub-M-polyfold.
Equipping $F$ with the induced M-polyfold structure, it becomes with the associated morphism set
an ep-groupoid.  If $x\in \partial X$ and $U(x)$ an open neighborhood in $X$ 
with the natural $G_x$-action it follows that $g\ast (F\cap U(x))=F\cap U(x)$ for all faces $F$ 
which contain $x$. Recall that there are exactly $d_X(x)$-many of such faces, where $d_X$ is the degeneracy index.
Denote by ${\mathcal F}_X$\index{${\mathcal F}_X$} the set of faces, where by assumption every $F$ is a sub-ep-groupoid
with the induced structure being that of a tame face-structured ep-groupoid. 
We denote by ${\mathcal F}_{X,x}$\index{${\mathcal F}_{X,x}$} the collection of all faces containing $x$. As already pointed out it holds that 
$d_X(x)=\sharp{\mathcal F}_{X,x}$.

We assume that for every face $F$ we are given a structured sc$^+$-multisection functor 
$$
\wt{\Lambda}_F=[\Lambda_F,\mathfrak{U}_F,\mathfrak{S}_F,\tau_F]
$$
for $W|F$. A representative of $\wt{\Lambda}_F$ is given by $(\Lambda_F,\mathfrak{U}_F,\mathfrak{S}_F,\tau_F$,
where $\Lambda_F:W|F\rightarrow {\mathbb Q}^+$ is a functor and for every $x\in F$ 
a  good neighborhood $U_F(x)\subset F$ has been given. This means $U_F(x)$ admits the natural $G_x$-action 
and $t:s^{-1}(\cl_F(U_F(x)))\rightarrow F$ is proper. Moreover we are given for every $x\in F$ an index set $I_x^F$ 
with an action by $G_x$ and a symmetric section structure $\mathfrak{s}_x^F$, given by
$$
{(s^x_{F,i})}_{i\in I_x^F}.
$$
Moreover, for $x,x'\in F$ a correspondence $\tau_{F,x,x'}:\bm{U}_F(x,x')\rightarrow \mathfrak{J}_F(x,x')$
is given relating in the usual way the section structures associated to $x$ and $x'$, respectively.
We note that we can always shrink the sets $U_F(x)$ restrict $\tau_{F,x,x'}$.

We consider the assignment ${\mathcal F}_X\ni F\rightarrow \wt{\Lambda}_F$ and shall formulate 
a compatibility condition.
\begin{definition}\index{D- Restriction $\wt{\Lambda}_{F|F'}$}
Given a face $F'\in {\mathcal F}_X$ the restriction of $\wt{\Lambda}_F$ to $W|(F\cap F')$ denoted by
$\wt{\Lambda}_{F|F'}$ is defined as follows, where 
$$
\wt{\Lambda}_{F|F'} =[\Lambda_{F|F'}, \mathfrak{U}_{F|F'},\mathfrak{S}_{F|F'},\tau_{F|F'}]
$$
and the following data is obtained by taking a representative $(\Lambda_F,\mathfrak{U}_F,\mathfrak{S}_F,\tau_F)$ 
of $\wt{\Lambda}_F$
\begin{itemize}
\item[(1)] \  $\Lambda_{F|F'}=\Lambda_F|(W|(F\cap F'))$.
\item[(2)] \  $\mathfrak{U}_{F|F'} = {(U_F(x)\cap F')}_{x\in F\cap F'}$, where $\mathfrak{U}_F ={(U_F(x))}_{x\in F}$ is a good system of open neighborhoods.
\item[(3)] \  $\mathfrak{S}_{F|F'}= {(\mathfrak{s}^{F|F'}_x)}_{x\in F\cap F'}$, where 
$\mathfrak{s}^{F|F'}_x$ is the symmetric sc$^+$-section structure obtained from $\mathfrak{s}^F_x$
by restricting the $s^x_i, i\in I_x^F$ to $U_F(x)\cap F'$.
\item[(4)] \  The correspondences $\tau_{F|F',x,x'}$ are the restrictions of $\tau_{F,x,x'}$ for $x,x'\in F\cap F'$
to $\{\phi\in \bm{X}\ |\ s(\phi)\in U_F(x)\cap F',\ t(\phi)\in U_F(x)\cap F'\}$.
\end{itemize}
\qed
\end{definition}
The compatibility condition is now given as follows.
\begin{definition}\index{D- Compatibility for $(\wt{\Lambda}_F)$}
Let $P:W\rightarrow X$ be a strong bundle over a tame face-structured ep-groupoid
and denote by ${\mathcal F}_X$ the set of faces. Assume that $F\rightarrow \wt{\Lambda}_F$
is an assignment which associates to a face $F\in {\mathcal F}_X$ a structure sc$^+$-multisection
on $W|F$.  Then the assignment is said to be compatible provided for every choice 
of faces $F,F'$ the equality
$$
\wt{\Lambda}_{F|F'}=\wt{\Lambda}_{F'|F}
$$
holds.
For a compatible assignment we define $\wt{\Lambda}_{F\cap F'}=\wt{\Lambda}_{F|F'}$.
\qed
\end{definition}
With this definition we see that $\wt{\Lambda}_{F\cap F'}$ is a structured sc$^+$-multisection 
of $W|(F\cap F')$.  A face of $F\cap F'$ is obtained by intersecting $F\cap F'$ with a suitable
$F''\in {\mathcal F}_X$. Then the compatibility allows us to consider 
$$
\wt{\Lambda}_{F\cap F'},\ \wt{\Lambda}_{F\cap F''},\ \text{and}\ \ \wt{\Lambda}_{F'\cap F''}.
$$
We note that $F\cap F'\cap F''$ is a common face of these three ep-groupoids and the restriction 
of the structures sc$^+$-multisections defines a common $\wt{\Lambda}_{F\cap F'\cap F''}$.
This can be done for an arbitrary number of faces.  

We note that as a consequence of the compatibility condition it holds that 
for 
$$
I^F_x=I^{F'}_x \ \ \text{for}\  F, F'\in {\mathcal F}_{X,x}
$$
We shall abbreviate $I_x:= I_x^F$ and also write $\mathfrak{J}(x,x')$ instead of $\mathfrak{J}_F(x,x')$.
 Moreover on $U_F(x)\cap U_{F'}(x)$ we have the equality
$s^x_{F,i}(y)=s^x_{F',i}(y)$ for $i\in I_x$.  We shall take for every $x\in \partial X$ a good open neighborhood
$V(x)$ of $x$ in $X$ such that  $U(x)=V(x)\cap \partial X$ satisfies
$$
U(x)\cap F\subset U_F(x)\ \ \text{for all}\ F\in {\mathcal F}_{F,x}.
$$

\begin{theorem}\index{T- Extensions for compatible $(\wt{\Lambda}_F)$}
Let $P:W\rightarrow X$ be a strong bundle over a tame face-structured ep-groupoid
and ${\mathcal F}_X\ni F\rightarrow \wt{\Lambda}_F$ be a an assignment of compatible
structured sc$^+$-multisection functors defined on the various $W|F$. Then 
this data defines a naturally structured $\wt{\Lambda}: W|\partial X\rightarrow {\mathbb Q}^+$.
In particular the underlying $\Lambda$ has an extension to $W|X$. The extension 
can be chosen to satisfy certain size restrictions, see Theorem \ref{p-main-p}.
\end{theorem}
\begin{proof}
The compatibility condition implies immediately that $\Lambda:W|F\rightarrow {\mathbb Q}^+$ 
given by
$$
\Lambda(w) =\Lambda_F(w)\ \ \text{for}\ P(w)\in F,\ F\in {\mathcal F}_{X,P(w)}
$$
is a well-defined functor. As already explained we have fixed 
for every $F\in {\mathcal F}_X$ a representative 
$$
(\Lambda_F,\mathfrak{U}_F,\mathfrak{S}_F,\tau_F)
$$
so that for every $x\in \partial X$ there exists an open neighborhood $U(x)\subset \partial X$
coming from a good open neighborhood $V(x)$ in $X$ such that $U(x)=V(x)\cap\partial X$ and 
$U(x)\cap F\subset U_F(x)$ for all $F\in {\mathcal F}_{X,x}$.  Moreover
it holds that 
$$
s^x_{F',i}(y)=s^x_{F,i}(y)\ \text{for} \ i\in I_x,\  y\in U(x)\cap F\cap F'.
$$ 
This allows us to define a symmetric sc$^+$-section structure $\mathfrak{s}_x={(s^x_i)}_{i\in I_x}$ on $U(x)$ by
$$
s^x_i(y)=s^x_{F,i}(y)\ \text{for}\ y\in F\cap U(x).
$$
We do this for every $x\in \partial X$. Finally we define the correspondences as follows.  
For $x\in \partial X$ and $x'\in \partial X$ consider $\bm{U}(x,x')$ which consists of all morphisms 
$\phi$ with $s(\phi)\in U(x)$ and $t(\phi)\in U(x')$. If $s(\phi)\in F\cap U(x)$ then necessarily $t(\phi)\in F\cap U(x')$.
Hence $\phi\in \bm{ U}_F(x,x')$. We can define 
$$
\tau_{x,x'}(\phi) =\tau_{F,x,x'}(\phi).
$$
This is well-defined since $\tau_{F,x,x'}(\phi)=\tau_{F',x,x'}(\phi)$ if $s(\phi)\in F\cap F'$.
Moreover, it follows from the definition that the map $\tau_{x,x'}:\bm{U}(x,x')\rightarrow \mathfrak{J}(x,x')$
is locally constant. It is also easy to verify that $\tau_{x,x'}$ relates the two local section structures.
Our desired structured sc$^+$-multisection for $W|\partial X$ is given by
$[\Lambda,\mathfrak{U},\mathfrak{S},\tau]$.  At this point we can apply Theorem \ref{p-main-p}
to construct the extension.
\qed \end{proof}

\chapter{Transversality and Invariants}\label{CHAPX16}
In this chapter we shall describe some elements of a global perturbation theory.

\section{Natural Constructions}\label{SEC141}
We assume that $(P:W\rightarrow X,\mu)$ is a strong bundle over an ep-groupoid and  consider for a smooth $x\in X$ the set of all sc-Fredolm operators $L\colon T_xX\rightarrow W_x$.  We define 
the set $\text{Gr}_F(x)$ whose elements are finite formal sums with non-negative rational coefficients of  sc-Fredholm operators.
An element of $\text{Gr}_F(x)$ has the form  
$$
\mathsf{L}=\sum\sigma_L\cdot L,
$$
where all $\sigma_L$ are non-negative and only a finite number are nonzero, and $L:T_xX\rightarrow W_x$.
Finally we put
$$
\text{Gr}_F(X)=\bigcup_{x\in X_\infty} \text{Gr}_F(x). \index{${\text{Gr}}_F(x)$}
$$
and write  $\pi:\text{Gr}_F(X)\rightarrow X_\infty$ for the  natural projection.
Note that for a smooth $x$, a sc-Fredholm operator $L:T_xX\rightarrow W_x$, and a morphism $\phi$ with $s(\phi)=x$ 
we can define $\phi_\ast L:T_{t(\phi)}X\rightarrow W_{t(\phi)}$ by 
$$
\phi_\ast L \circ T\phi (h) =\mu(\phi, L(h)),\ h\in T_xX.
$$
Hence $\phi:x\rightarrow y$ defines a bijection $\text{Gr}_F(x)\rightarrow \text{Gr}_F(y)$ via
$$
\phi_\ast(\sum\sigma_L\cdot L) =\sum\sigma_L\cdot \phi_\ast L.
$$
Elements in 
$\text{Gr}_F(X)$ can be viewed as objects of a category, where the morphism set consists of all pairs 
$(\phi,\mathsf{L})$ with $s(\phi)=\pi(\mathsf{L})$. The source map and target map are given by
$$
s(\phi,\mathsf{L})=\mathsf{L}\ \ \text{and}\ \  t(\phi,\mathsf{L})=\phi_\ast \mathsf{L}.
$$
The first result is the following.
\begin{theorem}[Well-defined $\mathsf{T}_{(f,\Lambda)}$]\label{THM1511} \index{T- Well-defined  $\mathsf{T}_{(f,\Lambda)}$}
Assume $(P:W\rightarrow X,\mu)$ is a strong bundle over an ep-groupoid,  $f$ an sc-Fredholm section, and 
 $\Lambda:W\rightarrow {\mathbb Q}^+$ an sc$^+$-multisection. Then there exists an associated section functor
 $\mathsf{T}_{(f,\Lambda)}$ of $\text{Gr}_F(X)\rightarrow X_\infty$ uniquely characterized by the following properties.
 \begin{itemize}
 \item[{\em(1)}] \  $\mathsf{T}_{(f,\Lambda)}(x)=0$ for $x\in X_\infty\setminus \supp(\Lambda\circ f)$.
 \item[{\em(2)}]\  For $x\in \supp(\Lambda\circ f)$ and any local sc$^+$-section structure ${(s_i)}_{i\in I}$, ${(\sigma_i)}_{i\in I}$ describing 
 $\Lambda$ around $x$ it holds 
 $$
 \mathsf{T}_{(f,\Lambda)}(y) = \sum_{\{i\in I\ |\ f(y)=s_i(y)\}} \sigma_i\cdot (f-s_i)'(y)\ \ \text{for smooth}\ \ y \ \text{near}\ x.
 $$
 \end{itemize}
\end{theorem}
\begin{remark}\index{R- On linearizations}
If $f(y)-s_i(y)\neq 0$ the expression $(f-s_i)'(y)$ is not defined, but it is also not counted towards the sum.
Setting the empty sum equal to $0$ property (1) follows from (2).
\qed
\end{remark}
\begin{proof}
Let $(P:W\rightarrow X,\mu)$ be a strong bundle over an ep-groupoid, $f$ an sc-Fredholm section and $\Lambda$ an sc$^+$-multisection.
Assume that $x\in \supp(\Lambda\circ f)$  and pick an open neighborhood $U(x)$ with the natural $G_x$-action and 
a local sc$^+$-section structure ${(s_i)}_{i\in I}$, ${(\sigma))}_{i\in I}$.
We need to analyze the formal sum
$$
\sum_{\{i\in I\ |\ f(y)=s_i(y)\}} \sigma_i\cdot (f-s_i)'(y)\ \ \text{for}\ \  y\in U(x)\cap X_\infty.
$$
Without loss of generality we may assume that $X=O$, where $(O,C,E)$ is a local model and  $x=0\in O$.
Also $W\rightarrow O$ is a strong bundle model and the local sc$^+$-section structure 
${(s_i)}_{i\in I}$, ${(\sigma_i)}_{i\in I}$, is defined on $O$.
\begin{proposition}\label{PROP1561}
There exists a sequence $(x_k)\subset  O$
with $x_k\neq 0$ and $x_k\rightarrow 0$ with the following property.
\begin{itemize}
\item For $i,i'\in I_e$ it holds that $s_i=s_{i'}$ near $x_k$.
\end{itemize}
\end{proposition}
The proof of the proposition relies on several lemmata.
\begin{lemma}
 Equip $O$ with the metric $d$ induced from $E_0$ and pick $y\in O$. 
Define $B_\varepsilon(y)=\{z\in O\ |\ d(z,y) < \varepsilon\}$ for some $\varepsilon>0$, and 
assume that $B_\varepsilon(y)$ can be written as a finite union 
$$
B_\varepsilon(y)=\bigcup_{i\in I} \Sigma_i
$$
of closed subsets $\Sigma_i$ of $B_\varepsilon(y)$. Then at least one of the sets $\Sigma_i$ has a nonempty interior.
\end{lemma}
\begin{proof}
Pick any $z\in B_\varepsilon(y)$ and, since every $\Sigma_i$ is closed, we can pick a $\delta$ such that
\begin{itemize}
\item $B_\delta(z)\subset B_\varepsilon(y)$.
\item For every $i\in I$ it holds that  $B_\delta(z)\cap\Sigma_i=\emptyset$ if and only if $z\not\in\Sigma_i$.
\end{itemize}
Denote by $K\subset I$ the subset consisting of the indices $i$ with $z\not\in \Sigma_i$. 
We note that moving $z$ slightly still $z\not\in\Sigma_i$ for $i\in K$.
We use this as follows. Take $i_1\in I\setminus K$ so that $z\in\Sigma_{i_1}$. If $z$ does not lie in the interior 
we can move it slightly to obtain $z_1$ so that $z_1\not\in \Sigma_i$ for $i\in K':=K\cup\{i_1\}$.
Define $K_1=\{i\in I\ |\ z_1\not\in\Sigma_i\}$. 
 We continue this process which at each step increases the set $K$ by at least one element and note that before we run out of indices
we find one for which $z$ lies in the interior. 
\qed \end{proof}
\begin{lemma}
Let $y\in O$ and $\varepsilon>0$. Then there exists a point $z\in O$ with $d(z,y)<\varepsilon$,
a  $\delta>0$, and a partition $I_1,...,I_p$ of $I$ such that 
\begin{itemize}
\item[{\em (1)}]\  $B_\delta(z)\subset B_\varepsilon(y)$.
\item[{\em(2)}]\ $s_i(w)=s_{i'}(w)$ for $w\in B_\delta(z)$, $i,i'\in I_e$ and $e=1,..,p$.
\item[{\em(3)}] \  $s_i(w)\neq s_j(w)$ for $w\in B_\delta(z)$ and $i\in I_e$, $j\in I_{e'}$ for $e\neq e'$.
\end{itemize}
\end{lemma}
\begin{proof}
Take a point $y\in O$ and fix an index $i_1\in I$. Consider on ${B}_\varepsilon(y)$ 
the sets $\Sigma_{j}:=\{z\in {B}_\varepsilon(y)\ |\ s_j(z)=s_{i_1}(z)\}$ and define $\bigcup_{j\neq i_1} \Sigma_j$. The latter set is closed 
in $B_\varepsilon(y)$. 
If its complement in ${B}_\varepsilon(y)$ is nonempty we can pick a point $z_1$ which has a neighborhood $B_{\delta_1}(z_1)$
on which every $s_j$ is point-wise different from $s_{i_1}$. We define $I_1=\{i_1\}$.
If $ \bigcup_{j\neq i_1} \Sigma_j={B}_\varepsilon(y)$
then there exists $j\neq i_1$ for which $\Sigma_{j}$  has a nonempty interior. We can take a point $z$ in this interior
and by slightly moving it we may assume that either $z$ belongs to the interior of some $\Sigma_k$ or it does not belong to
$\Sigma_k$. We pick such a $z_1$ and define $I_1=\{j\in I\ |\ s_j(z_1)=s_{i_1}(z_1)\}$. We also find $\delta >0$ such that
\begin{itemize}
\item $s_i(w)=s_{i'}(w)$ for $w\in B_{\delta_1}(z_1)$ and $i,i'\in I_1$.
\item  $B_{\delta_1}(z_1)\subset B_\varepsilon(y)$.
\item  $s_k(w)\neq s_{i}(w)$ for $w\in B_{\delta_1}(z_1)$, $i\in I_1$ and $k\not\in I_1$.
\end{itemize}
We note that moving $z_1$ less than $\delta_1$  the properties listed above are still true for a suitably
chosen smaller $\delta_2$.

In a next step we take an index $i_2\not\in I_1$ (if $I_1\neq I$)  and find $z_2$ with $d(z_2,z_1)<\delta_1/4$
and argue as before to obtain a set $I_2$ disjoint from $I_1$ which has the above properties 
for some $\delta_2< \delta_1/4$ and $I_2$. Note that the properties above persist for $I_1$ and
$\delta_2$ as well.  Proceeding inductively we obtain $\delta>0$ and $z$ having the properties
stated in the lemma.
\qed \end{proof}
\begin{proof}[Proof of Proposition \ref{PROP1561}]
We pick a sequence of points $(y_k)$  in $O$ with $y_k\neq 0$ and $y_k\rightarrow 0$.
We also take a sequence $\varepsilon_k\rightarrow 0$. We find sequences $(z_k)$ and $\delta_k>0$
such that $d(y_k,z_k)<\varepsilon_k$ and  the following holds:
\begin{itemize}
\item For every $k$ there is a partition of $I$, say $I=I^k_1\sqcup....\sqcup I^k_{p_k}$
so that $s_i(w)=s_{i'}(w)$ for $d(w,z_k)<\delta_k$ and $i,i'\in I^k_e$, where $e\in \{1,..,p_k\}$.
\item $s_i(w)\neq s_j(w)$ for $d(z_k,w)<\delta_k$, $i\in I^k_e$, $j\in I^k_{e'}$ and $e\neq e'$.
\end{itemize} 
After perhaps taking a first subsequence of $(z_k)$ we may assume that $I_1^k$ is constant.
Then we can take  of this subsequence another  subsequence so that $I^k_2$ is constant.
Proceeding this way we find after a finite sumber of steps a sequence $(x_k)$ converging to $0$
so that the assertion of our proposition holds. 
\qed \end{proof}

Next take a second local sc$^+$-section  structure ${(t_j)}_{j\in J}$, ${(\tau_j)}_{j\in J}$, where 
we again work on the local model $W\rightarrow O$ with $x=0$.  Define a partition of $J$
denoted by $J_1^k,...,J_p^k$ by setting
$$
J_e^k=\{j\in J\ |\ t_j(x_k)=s_{i_e}(x_k)\},
$$
where $i_e\in I_e$.  After taking a subsequence of $(x_k)$  we may assume that this partition
is independent of $k$.  Hence, without loss of generality  we  conclude for  $(x_k)$
that there are partitions 
$$
I=I_1\sqcup..\sqcup I_e\ \ \text{and}\ \ J=J_1\sqcup...\sqcup J_e
$$
where $I_e $ and $J_e$ correspond via the requirement that $s_{i_e}(x_k)=t_{j_e}(x_k)$. Moreover
$s_i=s_{i'}$ near $x_k$ if and only if they belong to the same $I_e$ and the same for $(t_j)$ with respect to
the partition of $J$. Also note that $t_j(z)=t_{j'}(z)$ for $z$ near $x_k$ provided $j,j'$ belong to the same $J_e$.
We deduce for $z$ near $x_k$
$$
\sum_{i\in I_e} \sigma_i =\Lambda(s_{i_e}(z))=\sum_{j\in J_e} \tau_j.
$$
We can say for $z$ near $x_k$  that the sc$^+$-section  $z\rightarrow s_{i_e}(z)$ occurs with multiplicity $\sigma^e:=\sum_{i\in I_e}\sigma_i$
and similarly $t_{j_e}$ which equals $s_{i_e}$ for $z$ near $z_k$ occurs with the multiplicity $\tau^e:=\sum_{j\in J_e}\tau_j$
which equals $\sigma^e$.  Of course, $s_{i_e}$ and $t_{j_e}$ are equal near $x_k$.
We note that the formal sums
$$
\sum_{e=1}^p \sigma^e\cdot T(f-s_{i_e})(x_k)
$$
and
$$
\sum_{e=1}^p \tau^e\cdot T(f-t_{j_e})(x_k)
$$
are equal and converge to $ \sum_{e=1}^p \sigma^e\cdot T(f-s_{i_e})(0)$ and 
$\sum_{e=1}^p \tau^e\cdot T(f-t_{j_e})(0)$ respectively, which have to be equal. Moreover, it follows
that
$$
\sum_{i\in I} \sigma_i\cdot T(f-s_i)(0)=\sum_{j\in J} \tau_j\cdot T(f-t_j)(0).
$$
Since $f(0)=s_i(0)$ and $f(0)=t_j(0)$ we conclude that 
$$
\sum_{i\in I} \sigma_i\cdot (f-s_i)'(0)=\sum_{j\in J} \tau_j\cdot (f-t_j)'(0).
$$
This completes the proof of Theorem \ref{THM1511}.
\qed \end{proof}
In view of the fact that 
$$
\mathsf{T}_{(f,\Lambda)}:X_\infty\rightarrow \text{Gr}_F(X)
$$
is well-defined we can make the following definition.
\begin{definition}\index{$\mathsf{T}_{(f,\Lambda)}$}\index{D- Linearization of $\Lambda\circ f$}
Let $(P:W\rightarrow X,\mu)$  be a strong bundle over an ep-groupoid,  $f$ an sc-Fredholm section, and 
 $\Lambda:W\rightarrow {\mathbb Q}^+$ an sc$^+$-multisection. We shall call 
 $\mathsf{T}_{(f,\Lambda)}$ the {\bf linearization of $\Lambda\circ f$} at the points in $\supp(\Lambda\circ f)$.
 Also we shall call for $x\in\supp(\Lambda\circ f)$ the expression  $\mathsf{T}_{(f,\Lambda)}(x)$ the linearization 
 of $\Lambda\circ f$ at $x$.
\end{definition}
As we shall see transversality questions are related to the behavior of $\mathsf{T}_{(f,\Lambda)}$.
We note that there is also an oriented version of the previous discussion which, however, is postponed to Section 
\ref{SEC153}.

\section{Transversality and Local Solution Sets}

In applications involving complicated moduli spaces it often occurs that perturbations  have to adhere to additional structures, limiting the transversality which can be achieved.
  The properties of  $\mathsf{T}_{(f,\Lambda)}$  control the geometry of 
$\supp(\Lambda\circ f)$ to first order.  Before we go deeper into transversality questions it is worthwhile to point out that as a rule of thumb
every transversality question within the polyfold theory can be resolved, provided  related finite-dimensional questions in an orbibundle context can be solved. Hence, the reader familiar with standard finite-dimensonal transversality questions, will be able to apply her/his insights in the polyfold context.

\begin{definition}\index{D- Transversality}
Let $(P:W\rightarrow X,\mu)$ be a strong bundle over an ep-groupoid, 
$f$ an sc-Fredholm section functor, and $\Lambda$ an sc$^+$-multisection.
\begin{itemize}
\item[(1)]\  $(f,\Lambda)$ is said to be {\bf transversal}\index{D- Transversal pair $(f,\Lambda)$} provided for every $x\in \supp(\Lambda\circ f)$
the linearization $\mathsf{T}_{(f,\Lambda)}$ at the point $x$ is a nonzero combination of surjective sc-Fredholm operators, i.e. if ${(s_i)}_{i\in I}$
is a local sc$^+$-section structure for $\Lambda$, and $f(x)=s_i(x)$ for $i\in I$, then all the $(f-s_i)'(x):T_xX\rightarrow W_x$ are surjective.
\item[(2)]\   If $x\in X$ with $\Lambda\circ f(x)>0$ we say that $(f,\Lambda)$ is {\bf transversal at $x$}\index{D- Transversal at $x$} provided
$\mathsf{T}_{(f,\Lambda)}(x)$ is a nonzero combination of surjective sc-Fredholm operators.
\end{itemize}
We shall say that $\mathsf{T}_{(f,\Lambda)}(x)$ is {\bf onto} or {\bf surjective} \index{Surjectivity of  $\mathsf{T}_{(f,\Lambda)}(x)$} at a point $x$  in $\supp(\Lambda\circ f)$ provided every
occurring sc-Fredholm operator with positive weight is onto.
\qed
\end{definition}

This notion of transversality has strong implications provided $d_X(x)=0$. 
\begin{theorem}[Interior local perturbation]\label{THM1524}\index{T- Interior local perturbation}
Assume  $(P:W\rightarrow X,\mu)$ is a strong bundle over an ep-groupoid, 
$f$ an sc-Fredholm section functor, and $\Lambda$ an sc$^+$-multisection functor.
Suppose   $x\in \supp(\Lambda\circ f)$ with $d_X(x)=0$  and $(f,\Lambda)$ is transversal  at $x$.  
Then there exists an open neighborhood $U(x)$ equipped with the natural $G_x$-action so that the following holds true
with $M_{U(x)}=\{y\in U(x)\ |\ \Lambda\circ f(y)>0\}$.
\begin{itemize}
\item[{\em (1)}]\  $d_X(y)=0$ for all $y\in U(x)$.
\item[{\em(2)}] \ $\mathsf{T}_{(f,\Lambda)}(y)$ is onto for every $y\in M_{U(x)}$.
\item[{\em (3)}] \  There exist submanifolds  ${(M_i)}_{i\in I}$ and positive rational weights ${(\sigma_i)}_{i\in I}$ 
with the following properties:
\begin{itemize}
\item[{\em (a)}] \ Every $M_i$ has $d_{M_i}\equiv 0$ and therefore has a natural (classical) smooth manifold structure.
\item[{\em (b)}]\ Every $M_i$ is properly embedded in $U(x)$.
\item[{\em (c)}]\ $M_{U(x)} =\bigcup_{i\in I} M_i$. 
\item[{\em (d)}]\ For every $y\in U(x)$ the identity 
$\Lambda\circ f(y)=\sum_{\{i\in I\ |\ y\in M_i\}} \sigma_i$
holds.
\item[{\em(e)}] \  $T_{(f,\Lambda)}(y) =\sum_{\{i\in I\ |\ f_i(y)=s_i(y)\}} \sigma_i\cdot (f-s_i)'(y)$ for $y\in U(x)$ and $T_{(f,\Lambda)}(y)$ is surjective
provided $y\in \supp(\Lambda\circ f)$.
\end{itemize}
\end{itemize}
\end{theorem}
\begin{proof}
This follows  from Theorem \ref{IMPLICIT0}. Pick an open neighborhood $U'(x)$
which allows the natural $G_x$-action and a local sc$^+$-section structure ${(s_i)}_{i\in I}$, ${(\sigma_i)}_{i\in I}$ on $U'(x)$
representing $\Lambda$ over $U'(x)$.

By assumption $f(x)=s_i(x)$ for $i\in I$.  Since $\mathsf{T}_{(f,\Lambda)}(x)$ is onto it follows that $(f-s_i)'(x)$ is onto
and the implicit function theorem gives for every $i\in I$ a solution set $M_i$ containing $x$ so that there exists an open neighborhood
$U_i(x)$ satisfying
\begin{itemize}
\item $M_i'=\{y\in U_i(x)\ |\ f(y)=s_i(y)\}$, $d_{M_i'}\equiv  0$, and $M_i'$ is a submanifold and consequently has 
an equivalent smooth manifold structure without boundary.
\item $(f-s_i)'(y):T_yX\rightarrow W_y$ is surjective for every $y\in M_i'$.
\end{itemize}
We find an open neighborhood $U(x)\subset U'(x)\bigcap\left(\bigcap_{i\in I} U_i(x)\right)$ which is invariant under $G_x$.
Then we define $M_i:=M_i'\cap U(x)$, and  by construction for $y\in M_{U(x)}:=\{y\in U(x)\ |\ \Lambda\circ f(y)>0\}$ the assertion 
(3) holds.
\qed \end{proof}
\begin{remark}\index{R- On transversality theory using sc$^+$-multisections}
We observe that on the ep-groupoid $U(x)$ the functor $\Theta:U(x)\rightarrow {\mathbb Q}^+$
defined by
$$
U(x)\rightarrow {\mathbb Q}^+:z\rightarrow \Lambda\circ f(z)
$$
is a branched, tame  ep$^+$-subgroupoid of $U(x)$ with support $M_{U(x)}$. We may think of transversality theory
as finding small $\Lambda:W\rightarrow {\mathbb Q}^+$ so that $\Lambda\circ f$ becomes a branched ep$^+$-subgroupoid, see Definition \ref{DEF912},
or at least a functor with reasonable properties.
\qed
\end{remark}

The situation becomes more complicated at points $x$ with $\Lambda\circ f(x)>0$ and $d_X(x)\geq 1$.
In this case we need additional information. For example if  $\partial X$ has some structure
near $x$ and $(\Lambda,f)$  lies in a reasonable position to $\partial X$ the solution space associated to $\Lambda\circ f(y)>0$ will have a good structure near $x$.
A version of Definition \ref{new_good_position_def} reads as follows.
\begin{definition}\index{D- Good position of $(f,\Lambda)$}
Let $(P:W\rightarrow X,\mu)$ be a strong bundle over an ep-groupoid, $f$ an sc-Fredholm section and $\Lambda$
an sc$^+$-multisection functor. 
Assume that $x\in \partial X$ satisfies $\Lambda\circ f(x)>0$.
We say that $(f,\Lambda)$ is in {\bf good position} at $x$ provided
\begin{itemize}
\item[(1)] \  $X$ is tame near $x\in X$.
\item[(2)] \ $(f,\Lambda)$ is transversal  at $x$, i.e. $\mathsf{T}_{(f,\Lambda)}(x)$ is surjective.
\item[(3)]\   The surjective sc-Fredholm operators $L_i$  of $\mathsf{T}_{(f,\Lambda)}(x)=\sum_{i\in I}\sigma_i\cdot L_i$
have kernels $\ker(L_i)$ which are in good position to the partial quadrant $C_xX$, see Definition \ref{mission1}.
\end{itemize}
\qed
\end{definition}
Next we describe the solution set near a point $x$ with $d_X(x)\geq 1$ and $\Lambda\circ f(x)>0$.
\begin{theorem}[Boundary local perturbation]\label{THMXXC1525} \index{T- Boundary local perturbation}
 Let $(P:W\rightarrow X,\mu)$ be a strong bundle over an ep-groupoid, $f$ an sc-Fredholm section and $\Lambda$
 a sc$^+$-multisection functor.
 Suppose $x\in \partial X$ with $\Lambda\circ f(x)>0$ and $(f,\Lambda)$ is in good position at $x$.  
 There exist an open neighborhood $U(x)$ with the natural $G_x$-action, tame submanifolds  ${(M_i)}_{i\in I}$ and positive rational weights ${(\sigma_i)}_{i\in I}$ 
with the following properties, where $M_{U(x)}=U(x)\bigcap\left(\bigcup_{i\in I} M_i\right)$.
\begin{itemize}
\item[{\em (1)}]\  Every $M_i$ is properly embedded in $U(x)$.
\item[{\em (2)}]\ $\mathsf{T}_{(f,\Lambda)}(y)$ is surjective for all $y\in M_{U(x)}$ and the kernels of the 
surjective sc-Fredholm operators with positive weights are in good position to $C_yX$.
\item[{\em (3)}] \ For every $y\in U(x)$ the identity 
$\Lambda\circ f(y)=\sum_{\{i\in I\ |\ y\in M_i\}} \sigma_i$
holds.
\end{itemize}
Note that the tame $M_i$ have equivalent smooth structures as manifolds with boundary with corners.
We also note that by definition $\Lambda\circ f: U(x)\rightarrow {\mathbb Q}^+$ is a tame (see Definition \ref{TAMERXX}) branched ep$^+$-subgroupoid,
since it is represented by tame $M_i$. The same conclusion holds if (2) is replaced by the general position requirement 
that the restrictions of $T_{(f,\Lambda)}(y)$ to $T_y^RX$ are surjective, see also the upcoming Definition \ref{DEFRT1535}.
 \end{theorem}
 \begin{proof}
 The proof follows along the lines of the proof of Theorem \ref{THM1524} using as ingredient
 the implicit function theorem for sc-Fredholm sections in the M-polyfold context, when
 the linearization is onto and the kernel is in good position to a tame $\partial X$, see Theorem \ref{IMPLICIT0}.
 We leave the details to the  reader.
 \qed \end{proof}
Before we start with the global perturbation theory we make a remark concerning orientations.
The details will be carried out in a larger context in Section \ref{SEC153}.
\begin{remark}[Formalism for orientations]\index{R- On the formalism for orientations}
Using results from Section \ref{SEC141} we can extend the previous discussion and incorporate 
orientations. 
In  this case we start with a strong bundle over an ep-groupoid $(P:W\rightarrow X, \mu)$, and oriented 
sc-Fredholm section functor $(f,\mathfrak{o})$ and an sc$^+$-multisection $\Lambda:W\rightarrow {\mathbb Q}^+$.
Assume that $(f,\Lambda)$ is transversal  at $x$,  and if $d_X(x)\geq 1$ we assume in addition that it is in good position.
Then we  can write with  the help of a local branching structure on $U(x)$, involving tame $M_i$ and suitable weights,
$$
\Lambda\circ f(y) = \sum_{\{i\in I\ |\ y\in M_i\}} \sigma_i,\ \ y\in U(x).
$$
This time, however, each $M_i$ inherits an orientation $o_i$ from $\mathfrak{o}$, which we have to keep track of.
For this we consider the formal sum
$$
\mathsf{T}_{((f,\mathfrak{o}),\Lambda)}(y) =\sum_{\{i\in I\ |\ y\in M_i\}} \sigma_i\cdot ((f-s_i)'(y),o_{(f-s_i)'(y)})
$$
for $y\in M_{U(x)}$. Since the $(f-s_i)'(y)$ are surjective and oriented we obtain an orientation $o_{i,y}$ for $T_yM_i$ 
so that we can pass to the formal sum
$$
\wh{\mathsf{T}}_{((f,\mathfrak{o}),\Lambda)}(y)=\sum_{\{i\in I\ |\ y\in M_i\}} \sigma_i\cdot (T_yM_i,o_{i,y}).
$$
 A priori all this might depend on the choice of the local sc$^+$-section structure, but as we shall see later it does not.
 We note that  $\wh{\mathsf{T}}_{((f,\mathfrak{o}),\Lambda)}$ can be viewed as an oriented lift of 
 the tangent $\mathsf{T}_{\Lambda\circ f}$ of the ep$^+$-subgroupoid $\Lambda\circ f$.
  The upshot will be that $\wh{\mathsf{T}}_{((f,\mathfrak{o}),\Lambda)}$
 is the natural orientation of $\Lambda\circ f$ coming from $(f,\mathfrak{o})$.  
 \qed
 \end{remark}

\section{Perturbations}
In this section we shall consider a strong bundle $(P:W\rightarrow X,\mu)$ over a tame ep-groupoid.   Assume that $f$ is an sc-Fredholm section functor
 of $P$. Then the associated solution category is $S_f=\{x\in X\ |\ f(x)=0\}$ and we assume as a standing assumption that
 the orbit space $|S_f|$ is compact, see Section \ref{SEC114}. Our aim is to study perturbations of $f$. As already previously pointed out,  as a rule of thumb, whatever would be possible in finite dimensions by multisection perturbations, will be possible in the ep-groupoid set-up as well (provided we have
 sc-smooth bump functions).
 Moreover, in most cases, what would be possible in a Sard-Smale perturbation theory without symmetries, see \cite{Smale}, 
 is possible to achieve by sc$^+$-multisection perturbations in cases of symmetry, i.e. the ep-groupoid setup.
 The price to be paid, since the latter is a perturbation theory over ${\mathbb Q}$, is that on a homological level, when studying solution spaces, 
 only a theory over ${\mathbb Q}$ can be used, unless, of course, one can use in a given instance less general perturbations, f. e. single-valued perturbations.
 
We  equip $W$  with an auxiliary 
 norm $N:W_{0,1}\rightarrow {\mathbb R}^+$, see Definition \ref{auxuilary_norm_def}. We can always extend 
 an auxiliary norm to all of $W$ by defining $N(h)=+\infty$ if $h$ is not of bi-regularity $(0,1)$. Hence without loss of generality 
 we always assume that $N:W\rightarrow {\mathbb R}^+\cup\{+\infty\}$.
As a consequence of the discussion in Section \ref{SEC114} there exists an open saturated  neighborhood $U$ of the {\bf solution category}\index{Solution category $S_f$}\index{$S_f$}
$S_f=\{x\in X\ |\ f(x)=0\}$ so that 
the orbit space of the set $\{x\in U\ |\ N(f(x))\leq 1\}$ has compact closure in $|X|$. 
\begin{definition}\index{D- Compactness control $(N,U)$}
Assume that $(P:W\rightarrow X,\mu)$ is a strong bundle over the ep-groupoid $X$ having a paracompact orbit space $|X|$.
Let $N:W\rightarrow [0,\infty]$ be an auxiliary norm and $f$ an sc-Fredholm section functor. Denote by $S_f$  the solution category associated to $f$.
Suppose that $U$ is a saturated open neighborhood of $S_f$.  
 We shall say that $(N,U)$ {\bf controls
the compactness} of $f$ provided the set $\{x\in U\ |\ N(f(x))\leq 1\}$ has an orbit space in $|X|$ which has a compact closure.
\qed
\end{definition}

In the following the relevant data for us is $(P:W\rightarrow X,\mu)$, $f$, and $(N,U)$.
Here $(P,\mu)$ is a strong bundle, $f$ is an sc-Fredholm section, and $(N,U)$  controls compactness, where $U$ is saturated open neighborhood of $S_f$.
We also assume that $X$ admits sc-smooth bump functions so that we can construct sc$^+$-multisections.
One of the aims  of the following considerations is to show the existence of many sc$^+$-multisection functors $\Lambda:W\rightarrow {\mathbb Q}^+$,
having domain support in $U$, having small point-wise norm, say  $N(\Lambda)(x)<1$ for $x\in X$, such that 
$$
\Lambda\circ f:X\rightarrow {\mathbb Q}^+
$$
is a tame branched ep$^+$-subgroupoid, and possibly having other additional properties.  The point-wise norm $N(\Lambda)$ is introduced in
Definition \ref{DEF1223}.
\begin{remark}\index{R- Bump functions or geometric perturbations}
In general we need the existence of sc-smooth bump functions to construct the sc$^+$-multisections for our perturbations.
In concrete situations,  one might be able to use more geometric maps to construct these perturbations. In these cases
compactness after perturbation might not follow from general principles and would have to be checked. 
\qed
\end{remark}

\begin{definition}\index{D- Admissible perturbation}
Let $(P:W\rightarrow X,\mu)$ be a strong bundle over an ep-groupoid, $f$ an sc-Fredholm section functor 
with $|S_f|$ being compact. Assume  $N$ is  an auxiliary norm and $U$ a saturated open neighborhood of $S_f$, so that $(N,U)$ controls compactness.
 An {\bf $(N,U)$-admissible perturbation} for $f$  is an sc$^+$-multisection $\Lambda:W\rightarrow {\mathbb Q}^+$
so that $N(\Lambda)(x)\leq 1$ for all $x\in X$ and the set
$$
\{x\in X\ |\ \exists\ t\in[-1,1]\ \text{with}\ (t\odot \Lambda)(f(x))>0\}
$$
is a subset of $U$.
\qed
\end{definition}
For example if the domain support of $\Lambda$ belongs to $U$ and $N(\Lambda)(x)\leq 1$ for all $x\in X$, then 
$\Lambda$ is admissible.   Indeed, assume that  $x\not\in U$ and solves $\Lambda\circ f(x)>0$.   Then  $\Lambda(h)=0$ for $h\in W_x\setminus \{0_x\}$, which implies
that $f(x)=0$ so that $x\in S_f\subset U$ giving a contradiction.
Admissible perturbations behave well under products. To see this assume we have two sets of data and take the product
so that $f\times f'$ is an sc-Fredholm section of $W\times W'\rightarrow X\times X'$. The auxiliary norm is defined
by $(N\times N')(h,h') =\text{max}\{N(h),N'(h')\}$. Given for the  two problems sc$^+$-section functors $\Lambda$ and $\Lambda'$, 
we define $\Lambda\times\Lambda':W\times W'\rightarrow {\mathbb Q}^+$ by
$$
(\Lambda\times\Lambda')(h,h')=\Lambda(h)\cdot \Lambda'(h'),
$$
which defines an sc$^+$-multisection functor.
Assume that $\Lambda$ is admissible with respect to $(N,U)$ and $\Lambda'$ with respect to $(N',U')$.
We observe that 
$$
(N\times N')(\Lambda\times\Lambda')(x,x')=\text{max}\{N(\Lambda)(x),N'(\Lambda')(x')\}\leq 1.
$$
Moreover assume that $(t\odot (\Lambda\times\Lambda'))\circ (f\times f')(x,x')>0$. This implies  that
$((t\odot \Lambda)\circ f(x))\cdot ((t\odot\Lambda')\circ f'(x'))>0$ from which we deduce that $(x,x')\in U\times U'$.
Note that we used the identity $t\odot(\Lambda\times\Lambda') = (t\odot\Lambda)\times (t\odot\Lambda')$.
\begin{remark}\index{R- On fibered products of sc$^+$-multisection functors}
The reader should note that if $\Lambda$ has domain support in $U$, and $\Lambda'$ domain support in $U'$,
the domain support of $\Lambda\times \Lambda'$ lies in $(U\times X' )\cup (X\times U')$ and not necessarily in 
$U\times U'$.  
\qed
\end{remark}
Next we define what it means to be in general position.
\begin{definition}\index{D- General position}\label{DEFRT1535}
Let $f$ be an sc-Fredholm section of $P$ where $(P:W\rightarrow X,\mu)$ is a strong bundle over a tame ep-groupoid $X$.
Assume further that $\Lambda:W\rightarrow {\mathbb Q}^+$ is an sc$^+$-multisection. 
We say that $(f,\Lambda)$ is in {\bf general position} provided for every $x\in \supp(\Lambda\circ f)$
the linearization $\mathsf{T}_{(f,\Lambda)}(x)$ restricted to the reduced tangent space at $x$ is surjective, i.e.
if $\mathsf{T}_{(f,\Lambda)}(x)=\sum_{i\in I} \sigma_i\cdot L_i$, where $\sigma_i>0$, then 
$$
L_i|T_x^RX:T^R_xX\rightarrow W_x
$$
is surjective for every $i\in I$. This is a natural extension of the Definition \ref{DEF_539}.
We shall call the restriction of $\mathsf{T}_{(f,\Lambda)}(x)$ to $T^R_xX$ the {\bf reduced linearization}\index{D- Reduced linearization} and denote it by
$\mathsf{T}_{(f,\Lambda)}^R(x)$. Hence $\mathsf{T}_{(f,\Lambda)}^R(x)=\sum_{i\in I} \sigma_i\cdot L_i|T^R_xX$.
\qed
\end{definition}
\begin{remark}\index{R- On good position and general position}
If $(f,\Lambda)$ is in general position then for every $x\in \supp(\Lambda\circ f)$ the linearization  $\mathsf{T}_{(f,\Lambda)}(x)$
is onto. If $x$ also happens to belong to the boundary $\partial X$ we have that $\mathsf{T}_{(f,\Lambda)}(x)$ is in good position, in fact even in general position.
In particular it follows that $\Lambda\circ f:X\rightarrow {\mathbb Q}^+$ is a tame branched ep$^+$-groupoid.\qed
\end{remark}
The next result shows that we always can force by a small perturbation a general position. Of course, if the perturbations are restricted 
by additional requirements that might not be possible. The reader should recall Theorem \ref{p:=}, Theorem \ref{thm_pert_and_trans} and Theorem \ref{SARD}
which deal with the case of strong bundles over a tame M-polyfold. The ideas used there are straight forward generalizations of the finite-dimensional case.
The following result is an adaption to the corresponding ep-groupoid case. The proof of this theorem follows along the lines of a similar result in \cite{HWZ3.5}. 

\begin{theorem}[General Position]\label{THM1536}\index{T- General position}
Assume $(P:W\rightarrow X,\mu)$  is a strong bundle over a tame ep-groupoid admitting sc-smooth bump functions, and 
 $f$ a  sc-Fredholm section functor with compact solution set $|S_f|$. Let $N$ be an 
 auxiliary norm and $U$ be a saturated subset of $X$ containing $S_f$
 so that $(N,U)$ controls compactness. Given  any continuous functor $h:X\rightarrow (0,1]$, 
 there exists an sc$^+$-multisection functor $\Lambda:W\rightarrow {\mathbb Q}^+$ with the following properties.
 \begin{itemize}
 \item[{\em (1)}] \  $N(\Lambda)(x)<h(x)$ for all $x\in X$.
 \item[{\em (2)}]\   The domain support of $\Lambda$ is contained in $U$.
 \item[{\em (3)}]\ $\mathsf{T}_{(f,\Lambda)}$ is surjective on $\supp(\Lambda\circ f)$.
 \item[{\em (4)}] \  For $x\in \supp(\Lambda\circ f)$ the kernels of the operators from $\mathsf{T}_{(f,\Lambda)}(x)$ are in general position
 to $\partial X$.
 \end{itemize}
 In particular $\Lambda\circ f:X\rightarrow {\mathbb Q}^+$ is a compact, tame, branched ep$^+$-subgroupoid.
 Also the perturbation $\Lambda$ can be taken  to be structurable (!).
 \end{theorem}
\begin{proof}
 We shall construct a structurable perturbation 
with the claimed properties. Point (4) is equivalent to the requirement that for $x\in \supp(\Lambda\circ f)$ the 
reduced linearization $T^R_{(f,\Lambda)}(x)$ is onto. The basic idea of the proof is as follows. One constructs for a large parameter space ${\mathbb R}^h$, i.e. 
$h>>0$, a family $\Lambda_t:W\rightarrow {\mathbb R}^{2n}$, $t\in {\mathbb R}^h$,  of structured sc$^+$-multisections, so that $\Lambda_0$ is the zero-section with weight $1$
and $\wt{\Lambda}:{\mathbb R}^h\times W\rightarrow {\mathbb Q}^+$ is an sc$^+$-multisection functor for the pull-back of $W$ by
${\mathbb R}^h\times X\rightarrow X$. The $\Lambda_t$ are constructed from local pieces with suitable properties so that
the  set $\{0\}\times S_f$ belongs to $\supp(\wt{\Lambda}\circ \wt{f})$, where $\wt{f}$ is defined by $\wt{f}(t,x)=(t,f(x))$, and moreover
$\mathsf{T}_{(\wt{f},\wt{\Lambda})}(0,x)$ is surjective for $x\in S_f$ when restricted to ${\mathbb R}^h\times T^R_xX$, or equivalently
the reduced $\mathsf{T}_{(\wt{f},\wt{\Lambda})}^R(0,x)$ are surjective.
This then implies that the solution set consisting of all $(t,y)$ near $\{0\}\times S_f$ is a finite-dimensional tame 
branched ep-subgroupoid. Applying an obvious  branched ep$^+$-groupoid version of Theorem \ref{SARD} the preimage of a suitable(!) small regular value $t_0$
for the projection
to ${\mathbb R}^h$ is the desired solution space and $\Lambda_{t_0}$ the structurable perturbation. Suitable here means that it is a regular value for restrictions to intersection of faces. The details are as follows.  \par

\noindent {\bf Local constructions at $x\in S_f$ with $d_X(x)=0$:} We pick smooth $w^1,...,w^\ell$ in $W_x$
so that $R(f'(x))$ together with $w^1,...,w^\ell$ span $W_x$. 
Pick an open neighborhood $U(x)$ which admits the $G_x$-action,
 has the properness property, so that in addition every $y\in U(x)$ satisfies $d_X(y)=0$. Let $V(x)$ be an open neighborhood
which is invariant under $G_x$ satisfying $\cl_X(V(x))\subset U(x)$. 
For every $w^j$ we can find an sc$^+$-section $s^j_1$ having the properties
\begin{itemize}
\item $s_1^j$ is defined on   $U(x)$  and $s_1^j(y)=0$ for $y\in U(x)\setminus V(x)$.
\item $s^j_1(x)=w^j$.
\item $N(s_1^j(y))\leq \text{min}\{1,h(y)\}$ for $y\in U(x)$.
\end{itemize}
Define  $s^{x,1}_{t}$ for $t\in {\mathbb R}^\ell$ by
$$
s^{x,1}_t(y)=\sum_{j=1}^\ell t_j\cdot s^j_1(y)\ \ \text{for}\ \ y\in U(x),\ t\in {\mathbb R}^\ell.
$$
We consider the local sc-Fredholm section $(t,y)\rightarrow f(y) -s^{x,1}_t(y)$ and note that the linearization at $(0,x)$ is surjective. 
Using the group elements  $g\in G_x$ we define for $y\in U(x)$ and $t\in {\mathbb R}^\ell$ 
$$
s^{x,g}_t(y) =\mu(\Gamma(g,g^{-1}\ast y),s^{x,1}_t(g^{-1}\ast y)).
$$
Of course, with $e\in G_x$ being the unit element we have that $s^{x,e}_t=s^{x,1}_t$.
Since $f(g\ast y) =\mu (\Gamma(g,y),f(y))$ it follows immediately  that the linearization of
$$
(t,y)\rightarrow f(y)-s^{x,g}_t(y)
$$
at $(0,x)$ is also surjective for every $g\in G$. The collection ${(s^{x,g}_t)}_{g\in G_x}$ is for fixed $t$ a symmetric local sc$^+$-section structure
and defines (for fixed $t$)  an atomic sc$^+$-multisection functor $\Lambda_t$. Hence for fixed $t$ it is structurable. 
Moreover, we can view ${\mathbb R}^\ell\times W$ as a strong bundle over ${\mathbb R}^\ell\times X$
and define $\wt{\Lambda}:{\mathbb R}^\ell\times W\rightarrow {\mathbb Q}^+$ by $(t,y)\rightarrow \Lambda_t(y)$.
We  define $\wt{f}(t,y)=(t,f(y))$, which is also an sc-Fredholm section. Clearly, since $x\in S_f$, it follows
that $(0,x)\in \supp(\wt{\Lambda}\circ \wt{f})$.  The linearization $\mathsf{T}_{(\wt{f},\wt{\Lambda})}(0,x)$  is by the previous 
discussion surjective  and from the local implicit function theorem we deduce the following.\par

\begin{itemize}
\item[(Interior)] \ \ \ \ \ \ \ \ \ \ There exists $\varepsilon_x>0$, and an open neighborhood $Q(x)$, invariant under $G_x$, and   contained in $U(x)$ 
such that $\{t\ |\ |t|<\varepsilon_x\}\times Q(x):(t,y)\rightarrow \Lambda_t^x\circ f(x)$ is a tame ep$^+$-subgroupoid. In addition, 
for $(t,y)$ in its support the linearization $\mathsf{T}_{(\wt{f}, \wt{\Lambda})}(t,y)$ is surjective. In particular, for every $x\in S_f\cap Q(x)$ the linearization 
$T_{(\wt{f},\wt{\Lambda})}(0,x)$ is surjective. The same assertion also holds for points $x$ in the saturation of $Q(x)$.
\end{itemize}

\noindent {\bf Local constructions at $x\in S_f$ with $d_X(x)\geq 1$:}  This part is similar to  the previous one, but we need
to be more careful in the construction, since ultimately we need to move the perturbed solution set into a general position.
We can find $U(x)$ and $V(x)$ as before, so that in addition $U(x)\cap \partial X$ has $d_X(x)$-many local faces, say $\Theta_1,...,\Theta_d$,
where $d=d_X(x)$. Then any finite intersection $\Theta_I=\bigcap_{i\in I} \Theta_i$  of such faces is a tame M-polyfold $\Theta_I$ containing $x$. Moreover $T_x^RX \subset T_x\Theta_I$.  We pick $w^1,...,w^\ell\in W_x$ such that the image of $f'(x)|T_X^RX$ together with the $w_i$
spans $W_x$. Of course, we can pick the $w_i$ as small as we wish. Now we proceed as in the previous case and construct
$\Lambda_t^x$, which again for fixed $t$ is structurable since it is atomic. By construction the associated $(\wt{f},\wt{\Lambda})$ 
is in general position at $(0,x)$. By the corresponding implicit function theorem we deduce the following.\par

\begin{itemize}
\item[(Boundary)] \ \ \ \ \ \ \ \ \ \ There exists $\varepsilon_x>0$, and an open neighborhood $Q(x)$ invariant under $G_x$  contained in $U(x)$ 
such that $\{t\ |\ |t|<\varepsilon_x\}\times Q(x):(t,y)\rightarrow \Lambda_t^x\circ f(x)$ is a tame ep$^+$-subgroupoid so that
for $(t,y)$ in its support the linearization $\mathsf{T}_{(\wt{f}, \wt{\Lambda})}(t,y)$ when restricted to the reduced tangent space
at $(t,y)$ is surjective. In particular, for every $x\in S_f\cap Q(x)$ the reduced linearization $\mathsf{T}^R_{(\wt{f},\wt{\Lambda})}(0,x)$ is surjective.
The same holds for $x\in S_f$ belonging in addition to the saturation of $Q(x)$.
\end{itemize}

\par

\noindent{\bf From local to global:} The collection ${|Q(x)|}_{x\in S_f}$ is an open covering of the compact set $|S_f|$.
We find finitely many points $x_1,...,x_k$ with the property
$$
|S_f|\subset \bigcup_{i=1}^k |Q(x_i)|.
$$
We abbreviate by $Q$ the union of the $Q(x_i)$ and by $\wt{Q}$ the saturation of $Q$. Then 
$$
S_f\subset \wt{Q}\subset U.
$$
We take the data associated to the $x_i$ previously constructed. Define $\bar{\ell}=\ell_{x_1}+...+\ell_{x_k}$
and ${\mathbb R}^{\bar{\ell}}={\mathbb R}^{\ell_{x_1}}\oplus...\oplus {\mathbb R}^{\ell_{x_k}}$ and for $t=(t_1,...,t_k)\in {\mathbb R}^{\bar{\ell}}$
$$
\Lambda_t (w) =\left(\bigoplus_{i=1}^k \Lambda_{t_i} \right)(w).
$$
Then every $\Lambda_t:W\rightarrow {\mathbb Q}^+$ is a structurable sc$^+$-multisection functor. In addition the 
associated $\wt{\Lambda}:{\mathbb R}^{\bar{\ell}}\times W\rightarrow {\mathbb Q}^+$  is an sc$^+$-multisection functor.
If $x\in Q$ we find $i$ so that an isomorphic $x'$ belongs to $Q(x_i)$. The properties of $(\wt{f},\wt{\Lambda})$ at $(0,x')$ are the same 
as those of $(0,x)$. From the  construction it follows that 
\begin{eqnarray}
\mathsf{T}_{(\wt{f},\wt{\Lambda})}^R(0,x)\ \ \text{is onto for all}\ \ x\in S_f.
\end{eqnarray}
We also note that there exists $\varepsilon_0>0$ so that $N(\Lambda_t)(x)<1$ for all $x\in X$, provided $|t|\leq \varepsilon_0$.
Since $(N,U)$ controls compactness and the domain support of every $\Lambda_t$ belongs to $U$ it follows that
\begin{eqnarray}
\{(t,x)\in {\mathbb R}^{\bar{\ell}}\times X\ |\ |t|\leq \varepsilon_0,\ \Lambda_t(f(x))>0\}\ \ \text{is compact.}
\end{eqnarray}
As a consequence of the implicit function theorem, f.e. use Theorem \ref{THMXXC1525},  we find $\varepsilon\in (0,\varepsilon_0]$ such that 
$$
\wt{\Theta}: \{t\in {\mathbb R}^{\bar{\ell}}\ | |t|<\varepsilon\}\times X\rightarrow {\mathbb Q}^+:(t,y)\rightarrow \Lambda_t\circ f(y)
$$
is a tame branched ep$^+$-subgroupoid, having the property that for $(t,y)\in \supp (\wt{\Theta})$ the reduced linearization 
$\mathsf{T}^R_{(\wt{f},\wt{\Lambda})}(t,y)$ is surjective. As a consequence $\wt{\Theta}$ can be locally written 
as with weights and a local branching structure consisting of tame M$^+$-polyfolds, or equivalently classically smooth manifolds
with boundaries with corners.  In order to find a specific value $t_0$ for which the perturbation $\Lambda_{t_0}$ has
the desired property we use a modification of the argument given in the proof of Theorem \ref{thm_pert_and_trans}. 
This is being discussed next.\par

\noindent{\bf Achieving global transversality:}  
At this point we can consider the map
\begin{eqnarray}\label{EQNXX1537}
\left\{(t,y)\in {\mathbb R}^{\bar{\ell}}\times X\ |\ |t|< {\varepsilon},\ \bar{\Lambda}_t(f(y))>0\right\}\rightarrow {\mathbb R}^{\bar{\ell}}:(t,y)\rightarrow t.
\end{eqnarray}
For $(t_0,x_0)$ in this set we can consider $d_X(x_0)$ and find an open neighborhood of the form $\wt{U}(t_0,y_0)=\{t\ |\ |t-t_0|<\delta\}\times U(x_0)$,
having the following properties.
\begin{itemize}
\item $U(x_0)$ admits the natural $G_{x_0}$-action and the properness property.
\item $U(x_0)$ contains $d_X(x_0)$-many local faces containing $x_0$.
\item On $\wt{U}(t_0,x_0)$ the branched ep$^+$-subgroupoid can be written with a local branching structure as
$$
\Lambda_t(f(x))=\frac{1}{|I|}\cdot |\{i\in I\ |\ (t,y)\in M_i\}|.
$$
\item Each $M_i$ is a smooth manifold with boundary with corners properly embedded in $\wt{U}(t_0,x_0)$
and the intersection of $M_i$ with an arbitrary number of local faces is a smooth manifold with boundary with corners.
\item For $(t,y)\in M_i$ it holds that $\mathsf{T}^R_{(\wt{f},\wt{\Lambda})}(t,x)$ is surjective. 
\end{itemize}
Denote by ${\mathcal F}$ the intersection of a finite number of local faces in $\wt{U}(t_0,x_0)$. The empty intersection 
we define to be $\wt{U}(t_0,x_0)$. Associated to ${\mathcal F}$ we have the smooth manifolds with boundary with corners
${\mathcal F}\cap M_i$, where $i\in I$.  Now we consider for ever ${\mathcal F}$ (there are a finite number of such intersections)
and $i\in I$ the map
$$
{\mathcal F}\cap M_i\rightarrow {\mathbb R}^{\bar{\ell}}:(t,x)\rightarrow t.
$$
For each such map there is a set of measure zero $Z_{i,{\mathcal F}}$ so that ${\mathbb R}^{\bar{\ell}}\setminus Z_{i,{\mathcal F}}$ 
consists of regular values. Denote by $Z$ the finite union of the $Z_{i,{\mathcal F}}$ which again has measure zero.
For $t$ in the complement the collection of all $(t,y)\in M_i$ is a smooth manifold with boundary with corners 
in general position to the boundary of ${\mathbb R}^{\bar{\ell}}\times X$. See the proof of Theorem 5.21 in \cite{HWZ3.5} for 
the argument proving the latter well-known fact of a parametrized transversality theory. A countable number of such sets $\wt{U}(t_0,x_0)$ have saturations 
which 
cover the set in (\ref{EQNXX1537}). Hence we find a set of measure zero in ${\mathbb R}^{\bar{\ell}}$ so that for $t$ in its complement 
(which contains points arbitrarily close to $0$) the perturbation $\Lambda_t$ has the desired properties. By slightly modifying the proof,
using the compactness of the closure of the orbit space of the set in (\ref{EQNXX1537}) even a discussion of a finite number of sets as
above suffices.
\qed \end{proof}

The standard example of a transversal extension theorem is given by the following result.

\begin{theorem}[Transversal Extension Theorem] \label{THM1537}\index{T- Transversal extension}
Assume that $(P:W\rightarrow X,\mu)$  is a strong bundle over a tame ep-groupoid and 
 $f$ a  sc-Fredholm section functor with compact solution set, i.e. $|S_f|$ is compact. We assume that $|X|$ is paracompact and $X$ admits sc-smooth partitions of unity.
  Let $N$ be an auxiliary norm and $U$ a saturated subset of $X$
 so that $(N,U)$ controls compactness.
 Suppose further that  $\Lambda:W|\partial X\rightarrow {\mathbb Q}^+$ is an $(U,N)$-admissible,  structurable 
 sc$^+$-multisection so that there exists a continuous functor $h:X\rightarrow (0,1]$ 
 with 
 $$
 N(\Lambda)(x)<h(x)\ \ \text{for all}\ \ x\in\partial X,
 $$
  and $\mathsf{T}_{(f|\partial X,\Lambda)}(x)$
 is surjective with respect to intersections of local faces for every smooth $x$ in $\partial X$. Then there exists a structurable  sc$^+$-multisection
 $\wt{\Lambda}:W\rightarrow {\mathbb Q}^+$ with the following properties.
 \begin{itemize}
 \item[{\em (1)}] \  $\wt{\Lambda}|(W|\partial X)=\Lambda$.
 \item[{\em (2)}] \  The domain support of $\wt{\Lambda}$ is contained in $U$ and $N(\wt{\Lambda})(x)<h(x)$ for all $x\in X$, i.e. $\wt{\Lambda}$ is $(N,U)$ admissible.
 \item[{\em (3)}] \  $\mathsf{T}_{(f,\Lambda)}$ is onto on $\supp(\Lambda\circ f)$.
 \end{itemize}
 Moreover, $\wt{\Lambda}\circ f:X\rightarrow {\mathbb Q}^+$ is a compact,  tame, branched ep$^+$-subgroupoid.
 \end{theorem}
\begin{proof}
We are given the saturated open subset $U$ of $X$ so that $(N,U)$ controls compactness.
Over the boundary we have the sc$^+$-multisection functor $\Lambda: W|\partial X\rightarrow {\mathbb Q}^+$ which is structurable and
satisfies for $x\in\partial X$ the inequality $N(\Lambda)(x)<h(x)$. We first take a continuous functor $h':X\rightarrow [0,1]$ satisfying 
$N(\Lambda)(x)<h'(x)<h(x)$ for $x\in \partial X$ and which is supported in $U$.
We apply Theorem \ref{p-main-p} and obtain a structurable extension 
$\Lambda':W\rightarrow {\mathbb Q}^+$ of $\Lambda$ which satisfies $N(\Lambda')(x)\leq h(x)$ for all $x\in X$ and has the domain support in $U$. 
The sc$^+$-multisection $\Lambda'$ extends $\Lambda$. By the assumptions on $\Lambda$ we find for every $x\in \partial X$ with $\Lambda\circ f(x)>0$
an open neighborhood $Q(x)\subset U$ admitting the natural $G_x$-action and having the properness property so that 
$$
\mathsf{T}_{(f,\Lambda')}^R(y)\ \ \text{is onto for}\ \ y\in U(x).
$$
The saturation of the union of all these $Q(x)$ defines a saturated open neighborhood of $\{y\in\partial X\ |\ \Lambda\circ f(y)>0\}$ and is contained
in $U$. At this point we have transversality for our solution set near $\partial$. For $x\in U\setminus Q$ with $\Lambda'\circ f(x)>0$ we note that $d_X(x)=0$
and we can argue as in Theorem \ref{THM1536} to construct parameter-depending perturbations we achieve local transversality. Using the compactness
of the orbit space of the support of $\Lambda'\circ f$ we can construct a structurable $\Lambda_t$ with domain support away from $\partial X$, but contained in $U$,
so that the structurable $\Lambda_t\oplus\Lambda'$ where $t\in {\mathbb R}^\ell$ has the property that
$$
(t,y)\rightarrow (\Lambda_t\oplus\Lambda')\circ f(x)
$$
has near $(0,x)$ with $x\in S_f$ the property that $\mathsf{T}_{(\wt{f},\wt{\Lambda})}^R(0,x)$ is onto. Then for arbitrarily small
$t$ we obtain the perturbation $\wt{\Lambda}:=\Lambda_t\oplus \Lambda'$ with the desired properties.
\qed \end{proof} 

There is an obvious corollary dealing with homotopy invariance. 
\begin{corollary}\label{CORRX1539}
Assume that $(P:W\rightarrow X,\mu)$ is a strong bundle over an ep-groupoid with $d_X\equiv 0$ and paracompact orbit space $|X|$.  Let $f$ be an sc-Fredholm functor
for which $|S_f|$ is compact.  
\begin{itemize}
\item[{\em (1)}]\   Fix an auxiliary norm $N$ and a saturated open neighborhood $U$ of $S_f$  so that $(N,U)$ controls compactness.
Then there exists a structurable  sc$^+$-multisection functor $\Lambda$ with domain support in $U$ and satisfying $N(\Lambda)(x)<1$ for all $x\in X$
having the property that $\mathsf{T}_{(f,\Lambda)}(x)$ is onto for all $x\in \supp(\Lambda\circ f)$. In particular $\Theta=\Lambda\circ f$ is a tame weighted 
ep$^+$-subgroupoid for which the orbit space of the support is compact. Let us call $(N,U,\Lambda)$ a structurable regular perturbation of $f$.
\item[{\em (2)}]\   Assume that $(N^i,U^i,\Lambda^i)$ for $i=0,1$ are regular perturbations of $f$.  Define $\wt{X}=[0,1]\times X$ and $\wt{W}=[0,1]\times W$, so that 
$\wt{W}\rightarrow\wt{X}$ is a strong bundle. Also define the sc-Fredholm section $\wt{f}$ by $\wt{f}(t,x)=(t,f(x))$. Then there exists an auxiliary
norm $N$ for $\wt{W}$ which over $i=0,1$ restricts to $N^i$ and a saturated open neighborhood $U$ of $S_{\wt{f}}$ restricting over $i=0,1$ to $U^i$
so that $(N,U)$ controls compactness.
\item[{\em (3)}] \  There exists a structurable regular perturbations $(N,U,\Lambda)$, where $\Lambda$ over $i=0,1$ is identical to $\Lambda^i$.
\end{itemize}
\end{corollary}
\begin{remark}\index{R- Topological content and orientations}
We note that the topological content is possibly empty since we do not have an orientation and we do a perturbation theory over ${\mathbb Q}$.
However, immediately, if we are in the case that $f$ is oriented it follows that $\Lambda\circ \wt{f}$ has a natural orientation and we obtain invariants.
For this see the next section.
\qed
\end{remark}
\begin{proof}[Corollary \ref{CORRX1539}]
By applying a special case of Theorem \ref{THM1536} we find a structurable  sc$^+$-multisection  $\Lambda_0$ with domain support in $U$
satisfying $N(\Lambda_0)(x)<1$ for all $x\in X$, so that $\mathsf{T}_{(f,\Lambda_0)}$ is onto on $\supp(\Lambda_0\circ f)$.
Then $\Theta_0=\Lambda_0\circ f$ is a compact branched ep$^+$-subgroupoid.   Assume we have also picked a second structurable
$\Lambda_1$ with similar properties.  Consider the tame ep-groupoid $[0,1]\times X$ equipped with the obvious projection
$[0,1]\times X\rightarrow X$. We can pull-back the strong bundle and obtain the strong bundle $\wt{W}\rightarrow [0,1]\times X$.
The fiber over $(t,x)$ is $W_x$ and $N$ defines an auxiliary norm on $\wt{W}$ which denote by $\wt{N}$. Define $\wt{U}=[0,1]\times X$.
Then $(\wt{N},\wt{U})$ controls compactness for the sc-Fredholm section $\wt{f}$ defined by $\wt{f}(t,x)=f(x)$. 
Applying Theorem \ref{THM1537} we find a structurable $\wt{\Lambda}$ for $\wt{W}$ with $\wt{N}(\wt{\Lambda})(t,x)<1$ for $(t,x)\in [0,1]\times X$
so that restricted to $\wt{W}|(\{i\}\times X)=W$ for $i=0,1$ we obtains $\Lambda_i$.  Hence $\wt{\Theta}:=\wt{\Lambda}\circ \wt{f}$ defines a tame, compact, branched ep$^+$-subgroupoid which restricts to $\{i\}\times X=X$ as $\Theta_i$.  This means we obtain a compact cobordism $\wt{\Theta}$  between $\Theta_0$ and $\Theta_1$.

 Assume that $(N,U)$ and $(N',U')$ control compactness and we have  associated sc$^+$-multisections
$\Lambda$ and $\Lambda'$.  Then we define  $U''=U\cap U'$ and $N''=\text{max}\{N,N'\}$ and note that $(N'',U'')$ also controls compactness.
We can take an associated $\Lambda''$. Then $\Lambda''$ is controlled  by $(N,U)$ as well as $(N',U')$. Hence we can connect in a controlled way
$\Lambda$ with $\Lambda''$ and $\Lambda''$ with $\Lambda'$, i.e. $\Lambda$ with $\Lambda'$.   
\qed \end{proof}

\section{Orientations and Invariants}\label{SEC153}
Next  we consider an oriented version of the previous  discussion.
We note that trans\-ver\-sality and orien\-tation are separate issues. Hence we may assume that transversality is already achieved
before we consider the orientation question.
The set-up is given by the following data.
\begin{eqnarray}\label{ASSUMP}
\text{\bf Assumptions}
\end{eqnarray}
\begin{itemize}
\item  A strong bundle over a tame ep-groupoid $(P:W\rightarrow X,\mu)$.
\item  An oriented sc-Fredholm section  functor $(f,\mathfrak{o})$, where we assume that $f$ has a compact solution set, i.e. 
$|S_f|$ is compact.  
\item A sc$^+$-multisection $\Lambda:W\rightarrow {\mathbb Q}^+$ such that $\mathsf{T}_{(f,\Lambda)}$ is onto
on $\supp(\Lambda\circ f)$,  and at points  $x\in \partial X\cap \supp(\Lambda\circ f)$ we assume that the kernels
of $\mathsf{T}_{(f,\Lambda)}$ are in good position to $\partial X$.
\item $\Lambda\circ f:X\rightarrow {\mathbb Q}^+$ is a branched ep$^+$-subgroupoid. (Note that this condition
says that at the boundary the behavior is resonable.)
\end{itemize}
Such data can be obtained by starting with an oriented sc-Fredholm section $(f,\mathfrak{o})$ 
and by applying perturbation theory to obtain $\Lambda$ with the above properties. 
Then $\Theta:X\rightarrow {\mathbb Q}^+$ defined by 
$$
\Theta:=\Lambda\circ f
$$
is a tame,  branched ep$^+$-subgroupoid. Hence we have the tangent functor
$$
\mathsf{T}_\Theta:X_\infty\rightarrow \text{Gr}(X).
$$
In order to orient $\Theta$ we need a suitable lift $\wh{\mathsf{T}}_\Theta:X_\infty\rightarrow \wh{\text{Gr}}(X)$.
The interesting fact, proved here, is that the orientation $\mathfrak{o}$ of $f$ defines a natural lift
$$
\wh{\mathsf{T}}_{(\Theta,\mathfrak{o})}:X_\infty\rightarrow \wh{\text{Gr}}(X).
$$
 In order to prove this 
we need some preparation.
Starting point is the data from Assumption \ref{ASSUMP}.
Then for  every smooth $x$ the convex set of sc-Fredholm operators
$(f-s)'(x):T_xX\rightarrow W_x$, where $s$ is a local sc$^+$-section with $s(x)=f(x)$, have induced orientations with suitable properties, see Section \ref{SEC115}.

For a smooth point $x$ denote by $\wh{\text{Gr}}_F(x)$ \index{$\wh{\text{Gr}}_F(x)$} the set of finite formal
sums $\sum_{\wh{L}} \sigma_{\wh{L}}\cdot \wh{L}$, where $\wh{L}=(L,o)$ is an oriented 
sc-Fredholm operator, $L:T_xX\rightarrow W_x$.  That means $o$ is an orientation of $\text{det}(L):= \Lambda^{max}\ker(L)\otimes (\Lambda^{max}\text{coker}(L))^\ast$. In the formal sum  the weights are nonnegative rational numbers and only a finite number of them are nonzero.
We define 
$$
\wh{\text{Gr}}_F(X):=\bigcup_{x\in X_\infty} \wh{\text{Gr}}_F(x).
$$
There is a natural projection $\wh{\pi}:\wh{\text{Gr}}_F(X)\rightarrow X_\infty$.
As before we can view $ \wh{\text{Gr}}_F(X)$ as the object set of a category.
The morphism set has the elements $(\phi,\wh{\mathsf{L}})$ where $\wh{\pi}(\wh{\mathsf{L}})=s(\phi)$.
The source and target maps are defined by $s(\phi,\wh{\mathsf{L}})=\wh{\mathsf{L}}$ and 
$t(\phi,\wh{\mathsf{L}})=\phi_\ast (\wh{\mathsf{L}})$. For an oriented $(L,o)$ we define
$\phi_\ast(L,o)=(\phi_\ast L,\phi_\ast\mathfrak{o})$.  

Forgetting orientations defines a 2-1 forgetful functor $\mathsf{f}_F$ fitting into the following diagram
$$
\begin{CD}
\wh{\text{Gr}}_F(X) @>\mathsf{f}_F>> \text{Gr}_F(X)\\
@V\wh{\pi} VV @V\pi VV\\
X_\infty @= X_\infty
\end{CD}
$$
We have seen that associated to $(f,\Lambda)$ there is a well-defined section functor 
$\mathsf{T}_{(f,\Theta)}$ of $\pi$.  As it turns out if $f$ is oriented there is a canonical lift
$\wh{\mathsf{T}}_{(f,\Lambda)}$ to a section functor of $\wh{\pi}$ such that
$$
\mathsf{f}_F\circ \wh{\mathsf{T}}_{(f,\Lambda)}={\mathsf{T}}_{(f,\Lambda)}.
$$
Let $x\in \supp(\Lambda\circ f)$, take an open neighborhood $U(x)$ with the natural $G_x$-action
and a local sc$^+$-section structure ${(s_i)}_{i\in I}$, ${(\sigma_i)}_{i\in I}$. Then 
$$
\mathsf{T}_{(f,\Lambda)}(y) =\sum_{\{i\in I\ |\ f(y)=s_i(y)\}} \sigma_i\cdot (f-s_i)'(y)\ \ \text{for}\ \ y\in U(x).
$$
This representation does not depend on the choice of the local sc$^+$-section structure. 
Since $f$ is oriented by $\mathfrak{o}$ every $(f-s_i)'(y)$ inherits an orientation $o_{i,y}$.
We define
$$
\wh{\mathsf{T}}_{((f,\mathfrak{o}),\Lambda)}(y)=\sum_{\{i\in I\ |\ f(y)=s_i(y)\}} \sigma_i\cdot ((f-s_i)'(y),o_{i,y})\ \ \text{for}\ \ y\in U(x)
$$
This definition does not depend on the choice of the local sc$^+$-section structure since the occurring
operators in the formal sum are independent of such choices.

Since for $y\in \supp(\Lambda\circ f)$  the operators occurring in the linearization $\mathsf{T}_{(f,\Lambda)}(y)$
are surjective we obtain an orientation $o_{i,y}$ for $\ker((f-s_i)'(y))$. If we vary $y$ in 
 $M_i = \{y\in U(x)\ |\ f(y)=s_i(y)\}$ the orientations $o_{i,y}$ vary continuously defining an orientation $o_i$ for $M_i$.
 Now we define 
 $$
 \wh{\mathsf{T} }_{(\Theta,\mathfrak{o})}(y) =\sum_{\{i\in I\ |\ y\in M_i\}} \sigma_i\cdot T_y(M_i,o_i).
 $$
 Note that the right-hand sum equals
 $$
 \sum_{\{i\in I\ |\ f(y)=s_i(y)\}} \sigma_i\cdot (\ker((f-s_i)'(y)),o_{i,y}),
 $$
 where $o_{i,y}$ is the orientation of $\Lambda^{max}\ker((f-s_i)'(y))$ which can be canonically identified with the
 orientation of $\det((f-s_i)'(0))$ via  identifying $\Lambda^{max}\ker((f-s_i)'(y))$
  with $(\Lambda^{max}\ker((f-s_i)'(y)))\otimes {\mathbb R}^\ast=\det((f-s_i)'(0))$.
\begin{theorem}[Canonical Orientation $\wh{\mathsf{T}}_{(\Lambda\circ f,\mathfrak{o})}$]\label{THM1541}
\index{T- Canonical orientation $\wh{\mathsf{T}}_{(\Lambda\circ f,\mathfrak{o})}$}
Let $(P:W\rightarrow X,\mu)$ be a strong bundle over an ep-groupoid, $(f,\mathfrak{o})$ an oriented sc$^+$-Fredholm section functor, and $\Lambda:W\rightarrow {\mathbb Q}^+$ an sc$^+$-multisection so that Assumption \ref{ASSUMP} holds. Then 
with $\Theta:=\Lambda\circ f:X\rightarrow {\mathbb Q}^+$ the functor 
 $\mathsf{T}_{\Theta}:X_\infty\rightarrow \text{Gr}(X)$ has a canonical lift 
$$
\wh{\mathsf{T}}_{(\Theta,\mathfrak{o})}:X_\infty\rightarrow \wh{\text{Gr}}(X),
$$
characterized by by
$$
 \wh{\mathsf{T} }_{(\Theta,\mathfrak{o})}(y)=\sum_{\{i\in I\ |\ f(y)=s_i(y)\}} \sigma_i\cdot (\ker((f-s_i)'(y)),o_{i,y}),
 $$
 where we use the identifications explained above.
 \qed
 \end{theorem}

\begin{definition}\index{D- Canonical orientation for $\Lambda\circ f$}\label{DEF1542}
Under Assumption \ref{ASSUMP} we
 call $\wh{\mathsf{T}}_{(\Lambda\circ f,\mathfrak{o})}$ defined in Theorem \ref{THM1541} the {\bf canonical orientation}
for $\Theta=\Lambda\circ f$.
\qed
\end{definition}

Next we incorporate differential forms into the picture.  We obtain the following result, which is a consequence of the preceding discussion.
\begin{theorem}
Let $(P:W\rightarrow X,\mu)$ be a strong bundle over the tame ep-groupoid $X$ and $(f,\mathfrak{o})$ an oriented
sc-Fredholm section of Fredholm index $k$ so that $|S_f|$ is compact. Assume that $N$ is an auxiliary norm and 
$U$ a saturated open neighborhood of $S_f$ so that $(N,U)$ controls compactness. We assume that $\Lambda:W\rightarrow {\mathbb Q}^+$
is an sc$^+$-multisection functor with domain support in $U$ satisfying
\begin{itemize}
\item[{\em (1)}] \ $N(\Lambda)(x)\leq 1$ for $x\in X$.
\item[{\em (2)}]\   $\mathsf{T}_{(f,\Lambda)}(x)$ is onto for all $x\in \supp(\Lambda\circ f)$. 
\item[{\em (3)}]\  $\mathsf{T}_{(f,\Lambda)}(x)$ for $x\in \partial X\cap \supp(\Lambda\circ f)$ is in good position to $\partial X$
\end{itemize}
Then $\Theta:=\Lambda\circ f$ is a compact, tame,  branched ep$^+$-subgroupoid, which is naturally oriented by $\mathfrak{o}$ via 
Theorem \ref{THM1541} and Definition \ref{DEF1542}, and $\partial\Theta$ has an induced orientation, Definition \ref{DEF936}. Then for every sc-differential form $\omega$ on $X$ of degree $(k-1)$
it holds that
$$
\oint_{(\Theta,\mathfrak{o})}d\omega =\oint_{\partial(\Theta,\mathfrak{o})} \omega.
$$
\qed
\end{theorem}
In view of Corollary \ref{CORRX1539} and the previous theorem we obtain the following result,
whose easy proof is left to the reader.
\begin{corollary}
Let $(P:W\rightarrow X,\mu)$ be a strong bundle over an ep-groupoid $X$ with paracompact orbit space and $d_X\equiv 0$.
Assume that $f$ is an sc-Fredholm section with $|S_f|$ being compact and $\mathfrak{o}$ is an orientation for $f$.
Then there exists a well-defined map
$$
\Phi_{(f,\mathfrak{o})}: H^\ast_{dR}(X)\rightarrow {\mathbb R}
$$
uniquely characterized for a homogenous element $[\omega]$ by
$$
\Phi_{(f,\mathfrak{o})} ([\omega]) = \oint_{(\Lambda\circ f,\mathfrak{o})}\omega,
$$
where $(N,U,\Lambda)$ is a regular perturbation and $\Lambda\circ f$ is oriented via the canonical orientation 
coming from $(f,\mathfrak{o})$, also denoted by $\mathfrak{o}$.
\qed
\end{corollary}

\begin{remark}\index{R- Invariants by triangulations}
There is another way to associate invariants to sc-Fredholm functors with compact moduli space. Namely for carefully 
chosen perturbations it is possible to triangulate the weighted solution set and to produce a chain with rational coefficients.
This is not carried out here and is left to the reader.
\qed
\end{remark}

\chapter{Polyfolds}\label{chap11+}
  We start by defining the notion of a polyfold, which is the generalization of 
   an orbifold with boundary and corners to the sc-smooth framework.
   All the results follow from the ep-groupoid case and the study of the concepts which behave well
   under generalized isomorphisms. As a consequence we allow our arguments to be quite brief.

\section{Polyfold Structures}\label{polyy}

This section capitalizes on the nice transformation properties of  previously introduced objects under equivalences,
and on the  compatibility  with the localization process. This allows to collect several results using a general  formalism.
Here is the basic definition.
\begin{definition}\label{SECDEF121}
Let $Z$ be a topological space. A {\bf polyfold structure}\index{D- Polyfold structure} for $Z$ is  a pair $(X,\beta)$ consisting of a ep-groupoid  $X$ and a homeomorphism
$\beta\colon |X|\rightarrow Z$.  
\qed
\end{definition}
The  polyfold structure $(X,\beta)$ relates the  ep-groupoid $X$ with the  topological space $Z$ via a homeomorphism
$\beta:|X|\rightarrow Z$.  One should view it as a additional structure on a topological space, which turns $Z$ into some kind 
of smooth space. In some sense $(X,\beta)$ is a generalization of a smooth atlas. 

The orbit space  $\abs{X}$ of an ep-groupoid  $X$ is a locally metrizable, regular, Hausdorff space (
Section \ref{section1.3_top_prop}) and therefore $Z$ must have the same  properties in order 
to admit a polyfold structure. If we know that $Z$ is a paracompact topological space  which  admits a polyfold structure, then  it follows via the Nagata-Smirnov Theorem that
$Z$ is metrizable. We summarize this in the following proposition.
\begin{proposition}
If a topological space $Z$ admits a polyfold structure its topology is locally metrizable and regular Hausdorff.
If a paracompact topological space admits a polyfold structure, then its topology is metrizable.
\qed
\end{proposition}
We distinguish different types of polyfold structures depending on properties of the ep-groupoid $X$.
The polyfold structure $(X,\beta)$  is called {\bf tame},\index{Tame polyfold structure} if the ep-groupoid $X$ is tame.
The polyfold structure $(X,\beta)$ is called {\bf face-structured}\index{Face-structured polyfold structure} if  $X$ is face-structured,
and it is called {\bf weakly face-structured}\index{Weakly face-structured polyfold structure} provided $X$ is weakly face-structured,
see Definition \ref{DEF1113}.

In order to introduce the notion of equivalent polyfold structures we recall from 
Section  \ref{section2.3_localization} that a generalized isomorphism 
$$
\mathfrak{f}=[X\xleftarrow{F}A\xrightarrow{G}Y]\colon X\to Y
$$
 between ep-groupoids induces the canonical homeomorphism 
$$\abs{\mathfrak{f}}=\abs{G}\circ \abs{F}^{-1}\colon \abs{X}\to \abs{Y}$$
between the orbit spaces. It does not depend on the choice of the diagram from the equivalence class. The inverse generalized isomorphism, 
$$
\mathfrak{f}^{-1}=[Y\xleftarrow{G}A\xrightarrow{F}X]\colon Y\to X,
$$
is  also a generalized isomorphism  and $\abs{\mathfrak{f}^{-1}}=\abs{\mathfrak{f}}^{-1}\colon \abs{Y}\to \abs{X}.$
The composition $\mathfrak{g}\circ \mathfrak{f}\colon X\to Z$ of two generalized isomorphisms 
$\mathfrak{f}\colon X\to Y$ and $\mathfrak{g}\colon Y\to Z$ is again a generalized isomorphism, and the induced homeomorphism $\abs{\mathfrak{f}\circ \mathfrak{g}}$ between orbit spaces satisfies 
$$\abs{\mathfrak{g}\circ \mathfrak{f}}=\abs{\mathfrak{g}}\circ \abs{\mathfrak{f}}\colon \abs{X}\to \abs{Z}$$
in view of Proposition \ref{equal_prop}. We also recall Theorem \ref{THMX10316} which asserts that for two generalized
isomorphisms $\mathfrak{f}:\mathfrak{f}':X\rightarrow Y$ the equality $|\mathfrak{f}|=\mathfrak{f}'|$ implies that $\mathfrak{f}=\mathfrak{f}'$.
We  define the  equivalence  of two polyfold  structures as follows.

\begin{definition}[{\bf Equivalent polyfold structures}] 
Two polyfold structures $(X,\alpha)$ and $(Y,\beta)$ for the topological space $Z$ are {\bf equivalent}\index{D- Equivalent  polyfold structures}
 provided there exists a generalized isomorphism
$\mathfrak{f}:X\rightarrow Y$ satisfying $\alpha=\beta \circ |\mathfrak{f}|$.  (As pointed out, $\mathfrak{f}$ is necessarily unique.)
\qed
\end{definition}
Now we can introduce the notion of a polyfold.
\begin{definition}[{\bf Polyfold}]
A {\bf polyfold} \index{D- Polyfold}  is a pair $(Z,c)$ consisting of a topological space $Z$ equipped with an equivalence class $c$
of polyfold structures.
\qed
\end{definition}
Most of the time we shall suppress $c$ in the notation and talk simply  about the polyfold $Z$. Of course, a topological space $Z$ might be a polyfold in different ways.
Next we shall study  properties which behave well with respect to equivalent polyfold structures.

\begin{proposition}\index{P- Existence of bump functions under equivalence}
For two equivalent polyfold  structures $(X,\alpha)$ and $(Y,\beta)$  of the topological space $Z$ the following holds.
\begin{itemize}
\item[{\em (1)}] \ $X$ admits sc-smooth bump functions if and only the same holds for $Y$.
\item[{\em (2)}]\   $X$ admits sc-smooth partitions of unity (as ep-groupoid)  if and only if the same holds for $Y$.
\end{itemize}
\end{proposition}
\begin{proof}
By assumption there exists a generalized isomorphism $[d]:X\rightarrow Y$. Taking a representative we obtain the diagram
$$
d\colon X\xleftarrow{F}A\xrightarrow{G}Y
$$
 involving  only  equivalences. 
If $X$ admits sc-smooth bump functions, there exists, by definition, a sc-smooth bump function in every open set. Since $F$ is a local sc-diffeomorphism, the same is true for the ep-groupoid $A$, and using that also $G$ is a local sc-diffeomorphism this also  holds for the image of $G$. If $y\not \in G(A)$, we find using that $G$ is essential surjective  a point $a\in A$ and a morphism $\psi\colon G(a)\to y$ in $\bm{Y}$. It has a unique extension to a local sc-diffeomorphism denoted by $\wh{\psi}$. The composition 
$\wh{\psi}\circ G$ is a local sc-diffeomorphism satisfying $\wh{\psi}\circ G(a)=y$ which allows to construct bump functions in a small neighborhood of $y$ outside of the image of $G$.

In order to see (2) assume that $X$ admits sc-smooth partitions of unity.  Take any open covering of the object space $Y$ by saturated open
sets, say $(U_\lambda)$. Denote by $\mathfrak{f}:X\rightarrow Y$ the generalized isomorphism satisfying $\beta\circ |\mathfrak{f}|=\alpha$ and by 
 $(V_\lambda)$ the covering of $X$ by saturated open subsets defined by
$$
|\mathfrak{f}|(|V_\lambda|)= |U_\lambda|.
$$
By assumption we can pick a subordinate sc-smooth partition of unity $(\sigma_\lambda)$ for $(V_\lambda)$. 
Define $\tau_\lambda$ by by $\tau_\lambda(y) := \sigma_\lambda(x)$, where $|\mathfrak{f}|(|x|)=|y|$.
It is an easy exercise that each $\tau_\lambda$ is sc-smooth and that $(\tau_\lambda)$ is an sc-smooth partition of unity subordinate 
to $(U_\lambda)$.
\qed \end{proof}
Sc-smooth bump functions or sc-smooth partitions of unity are important for constructions. In view of the proposition
the following definitions can be made.

\begin{definition}
We consider the  polyfold $(Z,c)$.
\begin{itemize}
\item[(1)]\  The polyfold $Z$ is called {\bf paracompact}\index{D- Paracompact polyfold}  if  the underlying topological space $Z$ is paracompact.
\item[(2)]\   The polyfold $Z$ is said to {\bf admit sc-smooth bump functions}\index{D- Polyfold admitting bump functions} provided there exists a  polyfold structure $(X,\alpha)\in c$ whose ep-groupoid $X$ admits sc-smooth bump functions. Note that also the  equivalent polyfold structures $(Y,\beta)\in c$ will have this property.
\item[(3)]\ The polyfold $Z$ is said to {\bf admit sc-smooth partitions of unity}\index{D- Polyfold admitting partions of unity} provided this  is true for a model $(X,\alpha)\in c$. Then this is true for the other models in $c$ as well.
\end{itemize}
\qed
\end{definition}
Our next aim is to construct a category ${\mathcal P}$ whose objects are polyfolds $(Z,c)$ and whose morphisms are what we shall call sc-smooth maps between them.  The definition of the morphisms needs some preparation.
For this let us  first consider another category,  called the {\bf category of polyfold structures}\index{Category of polyfold structures}, and denoted by  ${\mathcal P}{\mathcal S}$.\index{${\mathcal P}{\mathcal S}$}
The objects  in ${\mathcal P}{\mathcal S}$ are the  polyfold structures $(X,\alpha)$.
This means that $X$ is an ep-groupoid and 
$$
\alpha:|X|\rightarrow Z_\alpha
$$
 is a homeomorphism to the topological space $Z_\alpha$.
The {\bf morphisms}
$$
(\mathfrak{f},f)\colon (X,\alpha)\rightarrow (Y,\beta)
$$
 consist of pairs $(\mathfrak{f},f)$ in which  $f:Z_\alpha\rightarrow Z_\beta$ is a continuous map
and $\mathfrak{f}:X\rightarrow Y$ is a generalized (sc-smooth) map between ep-groupoid, i.e. $\mathfrak{f}$ is an equivalence class of diagrams
$$
d\colon X\xleftarrow{F}A\xrightarrow{\Phi} Y,
$$
where  $F$ is an equivalence of ep-groupoids and $\Phi$ a sc-smooth functor. Moreover,  
 $$
 f\circ\alpha = \beta\circ |\mathfrak{f}|.
 $$
We can visualize the structure by the following diagram, where the left-hand side stands for the diagram on right-hand side. 
 Of course, the top part does not make sense as a commutative diagram of maps, since $\mathfrak{f}$ is only a generalized isomorphism.
($\mathfrak{f}$  establishes locally a correspondence up to isomorphism between open subsets in the object spaces $X$ and $Y$, so that one has to interpret the upper part of the diagram accordingly.)
 \begin{equation*}
 (\mathfrak{f},f)\colon (X,\alpha)\rightarrow (Y,\beta) \quad \quad  \quad 
 \begin{CD}
 X @>\mathfrak{f}>> Y\\
 @V\pi_X VV @V \pi_Y VV\\
 |X|       @>|\mathfrak{f}|>> |Y|\\
 @V\alpha VV  @V\beta VV\\
 Z_{\alpha} @>f >>  Z_{\beta}.
 \end{CD}
 \end{equation*}
 The upper diagram is the part which defines the sc-smooth structure.
The composition of two morphisms is defined  by 
$$
(\mathfrak{f},f)\circ (\mathfrak{g},g)=(\mathfrak{f}\circ\mathfrak{g},f\circ g),
$$
and the identity elements have the form $1_{(X,\alpha)}=(1_{[1_X]},Id_{Z_\alpha})$.
The {\bf  isomorphisms} in ${\mathcal P}{\mathcal S}$ are  the pairs $(\mathfrak{f},f):(X,\alpha)\rightarrow (Y,\beta)$ in which $f\colon Z_\alpha\rightarrow Z_\beta$  is a homeomorphism 
and $\mathfrak{f}:X\rightarrow Y$ a generalized isomorphism between the ep-groupoids satisfying $|\mathfrak{f}|= \beta^{-1}\circ f\circ\alpha$.
As a consequence of Theorem \ref{THMX10316} the generalized isomorphism $\mathfrak{f}$ is uniquely determined by $f$. Of course, not every homeomorphism $f$ is induced by such
a $\mathfrak{f}$. 

We now consider a polyfold $(Z,c)$ and choose a polyfold structure  $(X,\alpha)\in c$.  Then $c$ consists of all objects $(Y,\beta)$ in ${\mathcal P}{\mathcal S}$
for which  there exists an isomorphism $( \mathfrak{f},f)\colon (X,\alpha)\rightarrow (Y,\beta)$ 
between the polyfold structures of the special form $(\mathfrak{f},Id_Z)$. Roughly speaking we consider
for a given $(X,\alpha)$ all the diagrams of the form below, where $\mathfrak{f}$ is a generalized isomorphism.
$$
\begin{CD}
X @>\mathfrak{f}>> Y\\
@V \alpha\circ \pi_X VV @ V\beta\circ \pi_Y VV\\
Z @=  Z
\end{CD}
$$
The diagram has to be understood in the sense that there exists a generalized isomorphism 
for which $|\mathfrak{f}|$ fits into the obvious commutative diagram. Since $|\mathfrak{f}|$
is completely determined we obtain, as already mentioned before, the uniqueness of $\mathfrak{f}$.

We might identify $c$ with this particular subcategory (which is not full). 
This allows us, given polyfolds $(Z,c)$ and $(W,d)$,  to define the notion of equivalence of morphisms $(\mathfrak{f},f)\colon (X,\alpha)\rightarrow (Y,\beta)$ between the polyfold structures $(X,\alpha)\in c$ and $(Y,\beta)\in d$.

\begin{definition}[{\bf Equivalence of morphisms}] \index{D- Equivalence of morphisms in ${\mathcal P}{\mathcal S}$}
Let $(Z,c)$ and $(W,d)$ be polyfolds and assume that $(X,\alpha),(X',\alpha')$ are  two polyfold structures in $c$ and $(Y,\beta),(Y',\beta')$ two polyfold structures in $d$.   The two morphisms
$(\mathfrak{f},f):(X,\alpha)\rightarrow (Y,\beta)$ and $(\mathfrak{f}',f'):(X',\alpha')\rightarrow (Y',\beta')$ in  ${\mathcal P}{\mathcal S}$ {\bf are 
equivalent} if  there exist  two isomorphisms $(\mathfrak{h},Id_Z):(X,\alpha)\rightarrow (X',\alpha')$ and $(\mathfrak{k},Id_W):(Y,\beta)\rightarrow (Y',\beta')$ in ${\mathcal P}{\mathcal S}$ satisfying $(\mathfrak{f}',f')\circ (\mathfrak{h},Id_Z)=(\mathfrak{k},Id_W)\circ (\mathfrak{f},f),$ as illustrated in the diagram
\begin{equation*}
\begin{CD}
(X,\alpha)@>(\mathfrak{f},f)>>(Y,\beta)\\
@V(\mathfrak{h},Id_Z)VV  @V(\mathfrak{k}.Id_W)VV\\
(X',\alpha')@>(\mathfrak{f}',f')>>(Y,\beta').
\end{CD}
\end{equation*}
\qed
\end{definition}
An equivalent pair $(\mathfrak{f},f)\sim (\mathfrak{f}',f')$ necessarily satisfies  $f=f'\colon Z\rightarrow W$.
If $(Z,c)$ and $(W,d)$ are two polyfolds and $(X,\alpha)\in c$ and $(Y,\beta)\in d$ two polyfold structures, then  a morphism
$(\mathfrak{f},f)\colon (X,\alpha)\rightarrow (Y,\beta)$ determines for two other  polyfold structures $(X',\alpha')\in c$ and $(Y',\beta')\in d$ the equivalent morphism 
$(\mathfrak{k}^{-1}\circ \mathfrak{f}\circ \mathfrak{h},f)\colon (X',\alpha')\rightarrow (Y',\beta')$, where $(\mathfrak{h},Id_Z)\colon (X',\alpha')\rightarrow (X,\alpha)$ 
and $(\mathfrak{k},Id_W)\colon (Y',\beta')\rightarrow (Y,\beta)$ are morphisms between the polyfold structures. 
\begin{remark}
If $f:Z\rightarrow W$ is a homeomorphism between polyfolds and $(X,\alpha)$ and $(Y,\beta)$
are polyfold structures, then there is at most one generalized isomorphism $\mathfrak{f}:X\rightarrow Y$
satisfying $f\circ\alpha =\beta\circ |\mathfrak{f}|$. However if $f$ is just a continuous map
there might be none, or perhaps many generalized sc-smooth map $\mathfrak{f}:X\rightarrow Y$
such that $f\circ\alpha =\beta\circ |\mathfrak{f}|$. So the uniqueness is special for homeomorphisms.
This is also a known fact in the orbifold world.
\qed
\end{remark}

\begin{definition}[{\bf Sc-smooth maps}] \label{sc-smooth_ map_pol}\index{D- Sc-smooth map between polyfolds}
A  sc-smooth  map $\what{f}:(Z,c)\rightarrow (W,d)$ between two polyfolds
is  an equivalence class $\what{f}=[(\mathfrak{f},f)]$ of morphisms between representatives
of the structures $c$ and $d$.  
\qed
\end{definition}

\begin{remark}
From a practical perspective one can  think of a sc-smooth map $\what{f}:(Z,c)\rightarrow (W,d)$
as a morphism $(\mathfrak{f},f):(X,\alpha)\rightarrow (Y,\beta)$ consisting of  a continuous map  $f:Z\rightarrow W$ 
and a generalized sc-smooth map $\mathfrak{f}:X\rightarrow Y$  satisfying  $f=\beta\circ |\mathfrak{f}|\circ\alpha^{-1}$.
Since $f:Z\rightarrow W$ is independent of the elements $(\mathfrak{f},f)$ in $\what{f}$
we shall often identify $\what{f}\equiv f$. However, the actual overhead (over $f$) is important for constructions
because  only on the level of the ep-groupoids   the sc-smoothness is captured. 
\qed
\end{remark}

\begin{definition}[{\bf Composition}] \index{D- Composition of sc-smooth polyfold maps}
The {\bf composition} $\wh{g}\circ \wh{f}\colon (Z, c)\to  (V,e)$ of the sc-smooth maps
$$
\text{$\what{f}:(Z,c)\rightarrow (W,d)$\quad  \text{and}\quad  $\what{g}:(W,d)\rightarrow (V,e)$}
$$
between polyfolds is defined via representatives as the  equivalence class
$$
\what{g}\circ\what{f}=[(\mathfrak{g}\circ\mathfrak{f},g\circ f)]
$$
in which
$(\mathfrak{f},f)\colon (X,\alpha)\rightarrow (Y,\beta)$ and $(\mathfrak{g},g)\colon (Y,\beta)\rightarrow (U,\gamma)$ are representatives of  $\what{f}$ and $\what{g}$, respectively.
The definition does not depend on the choices involved.
\qed
 \end{definition}

We are ready to introduce the polyfold category  ${\mathcal P}$.

\begin{definition}[{\bf Polyfold category ${\mathcal P}$}] \label{polyfold_category_def} The {\bf polyfold category}\index{D- Polyfold category} ${\mathcal P}$\index{${\mathcal P}$} has as objects the polyfolds $(Z,c)$ and as morphisms the sc-smooth maps $\what{f}:(Z,c)\rightarrow (W,d)$.
Further,  $1_{(Z,c)}=[([1_X],Id_Z)]$, where $(X,\alpha)\in c$. 
\qed
\end{definition}
The {\bf degeneracy index} \index{$d_{(X,\beta)}$}
$
d_{(X,\beta)}\colon (Z, c)\to  {\mathbb N}$ of the polyfold structure $(X, \beta)\in c$ is defined by  
$$
d_{(X,\beta)}(z):=d_X(x),\quad z\in Z,
$$
where $x\in X$ satisfies  $\beta(|x|)=z$. 
Since $s$ and $t$ are local sc-diffeomorphisms,  the definition does not depend on the
choice of $x$ as long as $\beta(|x|)=z$.  Even more is true.
\begin{lemma}
If $(X,\beta)$ and $(X',\beta')$ are two equivalent polyfold structures on $Z$, then 
$$
d_{(X,\beta)}=d_{(X',\beta')}.
$$
\end{lemma}
\begin{proof}
By the definition of equivalence there exists a generalized isomorphism $\mathfrak{f}:X\rightarrow X'$ satisfying  $\beta'\circ|\mathfrak{f}|=\beta$.
Let $z\in Z$ and choose points  $x\in X$ and $x'\in X'$ satisfying $\beta(|x|)=z=\beta(|x'|)$. We can represent $\mathfrak{f}$ by a diagram of equivalences
$$
X\xleftarrow{F}A\xrightarrow{F'} X'.
$$ 
Using that $F$ and $F'$ are equivalences,  we find a point $a\in A$ and morphisms $\phi:F(a)\rightarrow x$ and $\phi':F'(a)\rightarrow x'$.
This implies the equalities 
$$
d_{(X,\beta)}(z)=d_X(x)=d_X(F(a))=d_A(a)=d_{X'}(F'(a))=d_{X'}(x')=d_{(X',\beta')}(z)
$$
and the proof is complete. 
\qed \end{proof}
As a consequence we can define the degeneracy index map  for a polyfold.
\begin{definition}[{\bf Degeneracy index for polyfolds}] \index{D- Degeneracy index}
For  a polyfold $(Z,c)$ there exists a well-defined map $d_Z:Z\rightarrow {\mathbb N}$, called {\bf degeneracy index},\index{$d_Z$}
which has the property that for every polyfold structure $(X,\beta)\in c$  the identity
$$
d_Z(z)=d_X(x),
$$
holds, provided $\beta(|x|)=z$.
\qed
\end{definition}

Similarly  as we defined sub-M-polyfolds we define sub-polyfolds. 
\begin{definition}[{\bf Sub-polyfold}] \index{Sub-polyfold} A subset $B\subset Z$ of a polyfold $(Z, c)$ is called a {\bf sub-polyfold}\index{D- Subpolyfold} of $(Z,c)$ if there exists 
a polyfold structure $(X,\alpha)\in c$ such  that the saturated subset $A$ of $X$,  defined by 
$$
A= (\alpha\circ \pi)^{-1}(B), 
$$
is an ep-subgroupoid. Here $\pi\colon X\to \abs{X}$ is the usual projection onto the orbit space.
\qed
\end{definition}
\begin{remark}
If  $(Y,\beta)$ is another representative of $c$, then there exists  a generalized isomorphism $\mathfrak{f}:X\rightarrow Y$
satisfying $\beta\circ |\mathfrak{f}|=\alpha$. 
Let $\pi:X\rightarrow |X|$ and $\pi:Y\rightarrow |Y|$ be  the quotient maps. If the saturated subset $A = (\alpha\circ \pi)^{-1}(B)$ of $X$  is an ep-subgroupoid of $X$, then 
$A$ is a sub-M-polyfold of the object space $X$, and,  as we have shown,  the associated full subcategory is in a natural way an ep-groupoid.
If the diagram  $X\xleftarrow{F} C\xrightarrow{G} Y$ represents  $\mathfrak{f}$, the preimage of $A$ under $F$ is an ep-subgroupoid of $C$ 
and the saturation of its image in $Y$ is an ep-subgroupoid as well. Let us denote the latter by $A'$. Then obviously, 
$$
A' = (\beta\circ \pi')^{-1}(B).
$$
Hence, in the above definition,  the choice of a representative in $c$ does not matter.  We note that the generalized isomorphism
$\mathfrak{f}:X\rightarrow Y$ induces a generalized isomorphism $\mathfrak{f}_{A}:A\rightarrow A'$ between the ep-groupoids satisfying  
$$
(\beta| (\pi'(A')))\circ |\mathfrak{f}_A|=\alpha|(\pi(A)).
$$
Therefore with $|A|=\pi(A)$ and $|A'|=\pi'(A')$ it follows
that  
$$
(A,\alpha|(|A|))\ \ \text{and}\ \ (A',\beta\vert(|A'|))
$$
 are equivalent polyfold structures on $B$.
This shows that a subpolyfold inherits in a natural way a polyfold structure.
\qed
\end{remark}

Having the degeneracy index $d_Z$ we shall define the boundary of $Z$.
\begin{definition} \label{boundary_polyfold}\index{D- Boundary of a polyfold}
The {\bf boundary of  a polyfold } $Z$,  denoted by $\partial Z$,   is the collection of all points $z$ satisfying $d_Z(z)\geq 1$. 
\end{definition}
In general the boundary of a polyfold is just a set without any obvious meaningful geometric structure.
This changes for tame polyfolds, a concept we will introduce below.
\begin{example}
If $X={\mathbb R}^2$  we have the group action of $G={\mathbb Z}_2$, where the nontrivial element acts by $(x,y)\rightarrow (x,-y)$. 
We obtain the translation groupoid $G\ltimes {\mathbb R}^2$, which is an ep-groupoid (even tame). The degeneration index vanishes identically. We can identify the orbit space homeomorphically
with $Z={\mathbb R}\times [0,\infty)$ via $\beta:|(x,y)|\rightarrow (x,|y|)$. Hence $(X,\beta)$ defines a polyfold structure on $Z$. The boundary of $Z$ as a polyfold is empty, though as a manifold
it would be non-empty. Hence one has to be somewhat cautious talking about boundaries, since the geometric intuition might be misleading.
\qed
\end{example}

In order to introduce tame polyfolds we first prove the following simple lemma.
\begin{lemma}
We assume that the   polyfold structures $(X,\alpha)$ and $(Y,\beta)$ for $Z$ are equivalent. If the ep-groupoid $X$ is tame so is $Y$.
\end{lemma}
\begin{proof}
By definition of the equivalence of polyfold structures, there is a generalized isomorphism $\mathfrak{f}:X\rightarrow Y$ between ep-groupoids satisfying
$$
\beta\circ|\mathfrak{f}|=\alpha.
$$
The generalized isomorphism $\mathfrak{f}$ can be represented by a diagram 
$$
X\xleftarrow{F} A\xrightarrow{G} Y.
$$
 On the object level $F$ and $G$ are local sc-diffeomorphisms.
If any of the three ep-groupoids is tame, then all the others are as well.
\qed \end{proof}
The lemma allows to define a tame polyfold.

\begin{definition}[{\bf Tame polyfold}] \index{D- Tame polyfold}
A polyfold $(Z,c)$ is said to be {\bf tame} if  at least one (and then all) of the defining polyfold structures $(X,\beta)\in c$
has $X$ as  a tame ep-groupoid.
\qed
\end{definition}

If  $F:X\rightarrow Y$ is an  equivalence  between tame ep-groupoids or more generally a generalized isomorphism, we know from Remark \ref{REM1116}  that if one of the ep-groupoids is face-structured, so is the other.
This has an immediate consequence for generalized isomorphisms for which the same result holds.
\begin{definition}
A polyfold $(Z,c)$ is said to be {\bf  face-structured}\index{D- Face-structured polyfold} provided it is tame 
and for an admissible polyfold structure $(X,\alpha)\in c$ the ep-groupoid   $X$  is face-structured.
A {\bf face} $A$ of $Z$ is  by definition  the image of a face $F$ of $X$ under the map $\alpha\circ \pi:X\rightarrow Z$.
For a point $z\in Z$ the degeneracy index $d_Z(z)$ equals precisely the number of faces $z$ belongs to.
We call $(Z,c)$  {\bf weakly face-structured}\index{D- Weakly face-structured polyfold}  if   it is tame
and if for $(X,\alpha)\in c$ the ep-groupoid $X$ has  the property that the underlying object M-polyfold is face-structured.
\qed
\end{definition}

The following result is a consequence of a previous result saying that the property of being face-structured behaves well under generalized
isomorphisms.

\begin{proposition}
Every face $A$ of a face structured polyfold $(Z,c)$ has a natural polyfold structure $(A,c\vert A)$.
\qed
\end{proposition}

\section{Tangent of a Polyfold}
Let us show next that a polyfold $(Z,c)$ has a tangent space which again is a polyfold, denoted by $T(Z,c)=(TZ,Tc)$ and equipped 
with a sc-smooth map $p:TZ\rightarrow Z$.  The tangent space at a point will not be a vector space, but a quotient
of a vector space by a linear group action. 

We need some preparation and note that if  $(Z,x)$ is a polyfold, we can employ the degeneracy index
$d_Z:Z\rightarrow {\mathbb N}$ to define the subsets $Z_i$ by
$$
Z_i = \{z\in Z\ |\ d_Z(z)\geq i\}
$$
for $i\in {\mathbb N}$.
Clearly $Z_0=Z$ and $...\subset Z_{i+1}\subset Z_i..\subset Z_0=Z$.
The subsets $Z_i$ have natural polyfold structures $c^i$ defined as follows. If $(X,\alpha)\in c$
then the homeomorphism $\alpha:|X|\rightarrow Z$ defines a bijection
$\alpha^i: |X^i|\rightarrow Z_i$ and we can equip $Z_i$ with the unique topology making $\alpha^i$
a homeomorphism. Since $X^i$ is an ep-groupoid we see that $(X^i,\alpha^i)$ defines a polyfold structure
on $Z_i$. If $(X,\alpha),\ (Y,\beta)\in c$ it follows that $(X^i,\alpha^i)$ and $(Y^i,\beta^i)$ are equivalent
polyfold structures for $Z_i$. We shall write $Z^i$ for the set $Z_i$ equipped with the polyfold structure
$c^i$ and refer to the polyfold $(Z^i,c^i)$ or $Z^i$ for short. We shall call it the polyfold obtained
from $(Z,c)$ by raising the index by $i$
$$
(Z,c)^i:= (Z^i,c^i).
$$

Now we are in the position to carry out the tangent construction. Applying the tangent functor $T$ to the ep-groupoid $X$ we obtain the tangent ep-groupoid $TX$. The tangent 
$$
T\mathfrak{f}\colon TX\rightarrow TY
$$
of a generalized isomorphism $\mathfrak{f}=[X\xleftarrow{F} A\xrightarrow{H} Y]$ between ep-groupoids is defined as the generalized isomorphism
$$
T\mathfrak{f}=T[X\xleftarrow{F} A\xrightarrow{H} Y]:=[TX\xleftarrow{TF} TA\xrightarrow{TH} TY]
$$
between the tangent ep-groupoids, $T\mathfrak{f}\colon TX\to TY$. Since $TF$ and $TH$ are equivalences between ep-groupoids in view of Theorem \ref{gertrude}, the generalized map $T\mathfrak{f}\colon TX\to TY$ is indeed a generalized isomorphism. In view of Theorem \ref{TANGENTXXX} it is independent of the choice of the representative in the equivalence class $\mathfrak{f}$ of diagrams.
The following Lemma is a trivial consequence of Theorem \ref{THMX10316}.

\begin{lemma}\label{dafe}
If the two morphisms  $( \mathfrak{f},f)$ and $(\mathfrak{g},g)\colon (X,\alpha)\rightarrow (Y,\beta)$ between polyfold structures of $Z$ satisfy $f=g=\text{id}_Z$,  then 
$$
\mathfrak{f}=\mathfrak{g}\ \ \text{and}\ \ T\mathfrak{f}=T\mathfrak{g}.
$$
\end{lemma}
\begin{proof}
We assume that the diagram $d\colon X\xleftarrow{F}A\xrightarrow{G} Y$ represents $\mathfrak{f}$ and $d'\colon X\xleftarrow{F'}B\xrightarrow{G'} Y$ represents
$\mathfrak{g}$. From  $\beta\circ \abs{\mathfrak{f}}=\alpha =\beta\circ \abs{\mathfrak{g}}$ we conclude 
\begin{equation}\label{iiii}
|\mathfrak{f}|=|\mathfrak{g}|,
\end{equation}
which implies in view of Theorem \ref{THMX10316} that $\mathfrak{f}=\mathfrak{g}$ and therefore
$T\mathfrak{f}=T\mathfrak{g}$.
\qed \end{proof} 

We fix an  admissible polyfold structure 
$(X,\alpha)$ of the polyfold $(Z, c)$ and  take the tangent $TX$ of the ep-groupoid $X$. If $(Y,\beta)$ is an equivalent  polyfold structure, meaning there exists the isomorphism $(\mathfrak{f},Id_Z):(X,\alpha)\rightarrow (Y,\beta)$, Lemma \ref{dafe} tells us that $\mathfrak{f}$ and therefore 
$$
T\mathfrak{f}\colon TX\rightarrow TY
$$
between the tangent ep-groupoids are well-defined and unique. 

Now we shall define the tangent space $T(Z,c)$ of the polyfold $(Z,c)$.
For a point $z\in Z$ on level $1$, i.e. $z\in Z_1$,  and a polyfold structure $(X,\alpha)\in c$, we introduce  equivalence classes of tuples $[(z,(X,\alpha),h)]$ in which 
 $h\in T_xX$ and $\alpha(|x|)=z$. The equivalence is   defined as follows. 

If $(z,(X,\alpha),k)$ with $k\in T_yX$ is another such tuple
we call  $(z,(X,\alpha),h)$ equivalent to $(z,(X,\alpha),k)$ if  there exists 
a morphism $\phi\colon x\rightarrow y$ satisfying $T\phi (h)=k$. 
(The linear map $T\phi$ is defined since $\phi\in {\bf X}_1$.)
We denote the collection of these  equivalence classes by
$$
T^{(X,\alpha)}(Z,c) =\bigcup_{z\in Z_1} \{z\}\times \{[(z,(X,\alpha),h)]\}.
$$
Associated with  the generalized  isomorphism $\mathfrak{f}\colon X\rightarrow Y$ satisfying $\beta\circ|\mathfrak{f}|=\alpha$, we obtain the  uniquely determined canonical bijection
$$
T^{(X,\alpha)}(Z,c)\rightarrow T^{(Y,\beta)}(Z,c),\quad [(z,(X,\alpha),h)]\rightarrow [(z,(Y,\beta),|T\mathfrak{f}|(|h|))].
$$
Identifying the points via the canonical bijections we obtain a  set which we 
denote by 
$$
T(Z,c)\index{$T(Z,c)$}.
$$ 
This set is,  by definition,   the {\bf tangent space}\index{Tangent space of a polyfold} of $(Z,c)$,   as a set. Its elements are equivalence classes of equivalence classes
$$
[[(z,(X,\alpha),h)]].
$$
For every $(X,\alpha)\in c$ we have a natural map $T\alpha: |TX|\rightarrow T(Z,c)$ defined by
$$
T\alpha(|h|) = [[(z,(X,\alpha),h)]],
$$
where $h\in T_xX$ and $x\in X_1$ with $\alpha(|x|)=z$. This map is a bijection and fits into the commutative diagram
$$
\begin{CD}
|TX| @>T\alpha >>  T(Z,c)\\
@V |\tau_X| VV   @V\tau_ZVV\\
|X^1| @>\alpha >>      Z^1.
\end{CD}
$$
Given $(X,\alpha)$ and $(Y,\beta)$ in $c$ there exists a unique $\mathfrak{f}:X\rightarrow Y$
satisfying $\beta\circ |\mathfrak{f}|=\alpha$, defining $T\mathfrak{f}:TX\rightarrow TY$.
We note that 
$$
T\beta\circ |T\mathfrak{f}|=T\alpha.
$$
From these facts and the paracompactness discussion of ep-groupoids the following result follows immediately.
\begin{lemma}\label{l_para}
The set $T(Z,c)$ has a unique natural topology ${\mathcal T}$ for which for every $(X,\alpha)$
the map $T\alpha:|TX|\rightarrow T(Z,c)$ is a homeomorphism. The topology is regular and Hausdorff.
In particular if $(Z,c)$ is paracompact the same holds
for ${\mathcal T}$.
\qed
\end{lemma}

Next we define the polyfold structure on the topological space $T(Z,c)$.  
With the polyfold structure $(X,\alpha)$ for $(Z,c)$ we associate a polyfold structure $(TX, T\alpha)$ for the topological space $T(Z, c)$ by taking the tangent $TX$ of the ep-groupoid $X$ and the  homeomorphism
$$
T\alpha\colon |TX|\rightarrow T(Z,c).
$$
The pair $(TX,T\alpha)$ defines a polyfold structure on the topological space  $T(Z, c)$. 
Indeed, if $(Y, \beta)\in c$ and if  $\mathfrak{f}:X\rightarrow Y$ is the generalized isomorphism satisfying $\beta\circ|\mathfrak{f}|=\alpha$, then
$T\mathfrak{f}:TX\rightarrow TY$ satisifies
$$
T\beta\circ |T\mathfrak{f}|=T\alpha.
$$
Therefore,  the polyfold structures $(TX, T\alpha)$ for  $T(Z,c)$ associated with the polyfold structures $(X, \alpha)$ for $(Z, c)$ are all equivalent.

\begin{definition}[{\bf Tangent of $(Z, c)$}]
By $Tc$ we denote the equivalence class of  all polyfold structures for the topological space $T(Z,c)$
equivalent to those of the form  $(TX,T\alpha)$ for $T(Z, c)$ associated with the polyfold structures $(X,\alpha)$ of the polyfold $(Z,c)$. Then the pair $(T(Z, c), Tc)$ is a polyfold,  called the {\bf tangent of the polyfold $(Z,c)$},  and 
abbreviated by 
$T(Z,c)\index{D- Tangent of a polyfold}$.
\qed
\end{definition}
The continuous canonical projection map 
$$
p:T(Z, c)\rightarrow (Z^1, c^1)
$$
is defined, for $(X, \alpha)\in c$, by 
$$
p([[(z,(X,\alpha),h]])=z.
$$  
The  map $p$  comes from the  sc-smooth map
$$
[(P^{(X,\alpha)},p)]:T(Z,c)\rightarrow (Z^1,c^1)
$$
between the two polyfolds. Indeed, if $(TX, T\alpha)$ is the polyfold structure of $T(Z, c)$ 
associated with the polyfold structure $(X,\alpha)$ of $(Z, c)$, the generalized sc-smooth 
map  $P^{(X,\alpha)}\colon TX\rightarrow X^1$ between the ep-groupoids is the equivalence class 
$$P^{(X,\alpha)}=[TX\xleftarrow{1_{TX}}TX\xrightarrow{P}X^1]$$
where $P\colon TX\to X^1$ is the  canonical sc-smooth functor,   defined by 
$P(h)=x,$
if $h\in T_xX$ and $(X,\alpha)\in c$. 
It satisfies 
$$
p\circ T\alpha=\alpha\circ |P^{(X,\alpha)}|. 
$$
Hence, the sc-smooth map 
$[(P^{(X,\alpha)},p)]\colon T(Z,c)\rightarrow (Z^1,c^1)$ between the polyfolds induces,  on the topological level, the continuous projection  $p$. It  is the canonical sc-smooth bundle projection
to which we refer to sometimes as $p:TZ\rightarrow Z^1$,  if $c$ is understood.

\begin{theorem}\index{T- Tangent functor for polyfolds}
In the polyfold  category ${\mathcal P}$ introduced in Definition \ref{polyfold_category_def} there  exists a natural functor $T$, called the  {\bf tangent functor}\index{Tangent functor for polyfolds}\, which associates with the polyfold $(Z, c)$ its tangent $T(Z, c)$ and with a sc-smooth map $\wh{f}\colon (Z, c)\to (W, d)$ between polyfolds according to Definition \ref{sc-smooth_ map_pol}, its tangent $T\wh{f}\colon T(Z, c)\to T(W, d)$.  
\end{theorem}
\begin{proof}
We have already constructed $T$ on the objects of ${\mathcal P}$ and now consider a sc-smooth map  $\wh{f}\colon (Z,c)\rightarrow (W,d)$ between polyfolds. By definition, $\wh{f}=[( \mathfrak{f},f)]$ where $f\colon Z\to W$ is a continuous map and $\mathfrak{f}\colon X\to Y$ is a generalized map between ep-groupoids satisfying $f\circ \alpha =\beta \circ \abs{\mathfrak{f}}$  for the polyfold structures $(X, \alpha)\in c$ and $(Y, \beta)\in d$. We define its tangent  map 
$$T\wh{f}=[( T\mathfrak{f},Tf)]\colon T(Z, c)\to T(W, d)$$
by 
$$
Tf([[(z,(X,\alpha),h)]])=[[(f(z),(Y,\beta),k)]]  
$$
where $|k|=|T\mathfrak{f}|(|h|)$. 
One easily verifies that $T$ is a functor.
\qed \end{proof}

Next we study sc-differential forms on the polyfold $(Z, c)$ and choose 
a polyfold structure $(X, \alpha)$ in $ c$, so that $\alpha:|X|\rightarrow Z$ is a homeomorphism.
If $(X',\alpha')\in c$ is another polyfold structure for $(Z, c)$, then there exists a unique  sc-smooth generalized isomorphism $\mathfrak{h}:X\rightarrow X'$ between ep-groupoids satisfying
\begin{equation}\label{fish}
\alpha'\circ |\mathfrak{h}|=\alpha.
\end{equation}
 If $\omega\in \Omega^\ast_{ep,\infty}(X)$ is a sc-form on $X$, then the push-forward sc-form $\mathfrak{h}_\ast\omega\in 
\Omega^\ast_{ep,\infty}(Y)$ is well-defined,  and by Theorem \ref{OOmaine} does not depend  on $\mathfrak{h}$ provided \eqref{fish} holds.
The same is true  for the pull-back of a form $\omega'\in \Omega^\ast_{ep,\infty}(Y)$.

As a consequence of the previous discussion we shall  define the notion of a sc-differential form on a polyfold $Z$.
\begin{definition}\label{def_2.70}
 A {\bf sc-differential form on the polyfold $(Z,c)$}  is an equivalence class $[(\omega,(X,\alpha)]$  where $(X,\alpha)$ is a polyfold structure in $c$
and $\omega\in \Omega_{ep,\infty}(X)$ a sc-differential form on $X$. The equivalence is defined as follows,  
$$(\omega,(X,\alpha))\sim(\tau,(Y,\beta))$$
for $(X, \alpha)$ and $(Y, \beta)$ in $c$, if  there exists a generalized isomorphism 
 $\mathfrak{h}:X\rightarrow Y$ between ep-groupoids satisfying 
 $\beta \circ \abs{\mathfrak{h}}=\alpha$ and 
$\tau=\mathfrak{h}_\ast\omega$.
\end{definition}
We denote by $\Omega_{sc}^\ast(Z,c)$ the collection of sc-differential forms on $(Z,c)$.  This is a vector space with the operations defined by 
$$
[(\omega,(X,\alpha))]+\lambda[(\tau,(X,\alpha))]=[(\omega+\lambda\tau,(X,\alpha))].
$$
The exterior differential $d=d_Z$ is defined as
$$
d[(\omega,(X,\alpha))]=[(d\omega,(X,\alpha))].
$$
If 
$$
\what{f}:(Z,c)\rightarrow (Z',c')
$$
is a sc-smooth map between the two polyfolds,  and 
if $(X,\alpha)\in c$ and $(X',\alpha')\in c'$ are polyfold structures, then there exists a representative $(\mathfrak{f},f)$, where $f:Z\rightarrow Z'$ is a continuous map
and $\mathfrak{f}:X\rightarrow X'$ a generalized map between ep-groupoids, and we define the {\bf pull-back} of  $[(\omega',(X',\alpha')]\in \Omega_{sc}^\ast(Z',c')$ by 
$$
(\what{f})^\ast [(\omega',(X',\alpha'))] = [(\mathfrak{f}^\ast\omega',(X,\alpha))].
$$
By construction,  
$$
d_Z\circ \what{f}^\ast =\what{f}^\ast\circ d_{Z'}.
$$
We can summarize the discussion, which reduces it to the already discussed
ep-groupoid case,  in  the following theorem.
\begin{theorem}
Given a polyfold $(Z,c)$,  there is an associated deRham complex $(\Omega^\ast_{sc}(Z,c),d)$. A  sc-smooth map
$\what{f}:(Z,c)\rightarrow  (Z',c')$ between polyfolds induces  the co-chain map
$$
\text{$\what{f}^\ast\colon \Omega^\ast_{sc}(Z',c')\rightarrow \Omega^\ast_{sc}(Z,c)$ \quad satisfying $\ d_Z\circ \what{f}^\ast=\what{f}^\ast\circ d_{Z'}$}.
$$
Moreover,  the usual functorial properties hold.
\end{theorem}

\section{Strong Polyfold Bundles}
We shall introduce the notion of a strong polyfold bundle. The basic building blocks are again 
equivalence classes of strong bundles over ep-groupoids, quite in the spirit of the previous discussions.
Since there are no new ideas needed we allow ourselves to be somewhat sketchy.

The polyfold $(Z, c)$ is a topological space $Z$ equipped with an equivalence class $c$ of polyfold structures $(X, \alpha)$. Similarly we are going to define the strong polyfold bundle $(p\colon Y\to Z, \ov{c})$ over the polyfold $(Z, c)$ as a continuous and surjective map $p$ between the topological spaces equipped with an equivalence class $\ov{c}$ of strong polyfold bundle structures for $p$.
\begin{definition}[{\bf Strong polyfold bundle structures}]
Let $Y$  and $Z$ be two topological spaces and let 
$$p\colon Y\to Z$$
be 
a continuous and surjective map.
 A {\bf strong polyfold bundle structure}\index{D- Strong polyfold bundle structure} for $p\colon Y\rightarrow Z$ is a tuple 
 $$
 ((P\colon W\rightarrow X,\mu),\Gamma,\gamma)
 $$ 
 consisting of a 
strong bundle  $(P\colon W\rightarrow X,\mu)$ over the  ep-groupoid $X$, a the homeomorphism $\Gamma\colon |W|\rightarrow Y$ covering the homeomorphism $\gamma\colon |X|\rightarrow Z$ so that $p\circ \Gamma = \gamma\circ |P|$, as illustrated in the diagram
\begin{equation*}
\begin{CD}
W@>>> \abs{W}@>\Gamma>>Y \\
@V PVV @V\abs{P}VV   @VVp V\\
X@>>> \abs{X} @>\gamma>>Z,\\
\end{CD}
\end{equation*}
where the horizontal arrows on the left are the maps passing to  orbits.
\qed
\end{definition}
The smoothness properties of $p:Y\rightarrow Z$ are described by $(P:W\rightarrow X,\mu)$.  In general the latter 
is only taken up to (some notion of) equivalence, which is explained now.
\begin{definition}[{\bf Equivalence of strong polyfold bundle structures}]\index{D- Equivalence of strong polyfold bundle structures}
 Two strong polyfold bundle structures 
 $$
 ((P\colon W\to X, \mu), \Gamma, \gamma)\ \ \text{and}\ \ 
((P'\colon W'\to X', \mu'), \Gamma', \gamma')
$$ 
for the map 
 $$p\colon Y\to Z$$
 are {\bf equivalent} if there exists a generalized strong bundle isomorphism
 $$
\mathfrak{F}\colon W\to W'
 $$
  covering the generalized isomorphism $\mathfrak{f}\colon X\to X'$ between ep-groupoids and satisfying  
 \begin{align*}
 \gamma'\circ \abs{\mathfrak{f}}=\gamma\quad \text{and}\quad \Gamma'\circ \abs{\mathfrak{F}}&=\Gamma.
  \end{align*}
 The situation is illustrated in the diagrams 
 \begin{equation*}
 \begindc{\commdiag}[500]
\obj(1,1)[aa]{$W$}
\obj(3,1)[ap]{$W'$}
\obj(1,0)[bb]{$X$}
\obj(3,0)[bp]{$X'$}
\mor{aa}{ap}{$\mathfrak{F}$}
\mor{aa}{bb}{$P$}[\atright,\solidarrow]
\mor{ap}{bp}{$P'$}
\mor{bb}{bp}{$\mathfrak{f}$}
\enddc
\qquad \qquad  \qquad 
\begindc{\commdiag}[500]
\obj(2,2)[yy]{$Y$}
\obj(1,1)[ww]{$\abs{W}$}
\obj(3,1)[wp]{$\abs{W'}$}
\obj(1,0)[xx]{$\abs{X}$}
\obj(3,0)[xp]{$\abs{X'}.$}
\obj(2,-1)[zz]{$Z$}
\mor{ww}{yy}{$\Gamma$}
\mor{wp}{yy}{$\Gamma'$}[\atright,\solidarrow]
\mor{ww}{wp}{$\abs{\mathfrak{F}}$}
\mor{ww}{xx}{$\abs{P}$}[\atright,\solidarrow]
\mor{wp}{xp}{$\abs{P'}$}
\mor{xx}{xp}{$\abs{\mathfrak{f}}$}
\mor{xx}{zz}{$\gamma$}[\atright,\solidarrow]
\mor{xp}{zz}{$\gamma'$}
\enddc
\end{equation*}
\qed
 \end{definition}
 An equivalence class $\ov{c}$ of strong bundle structures for $p\colon Y\to Z$ induces an equivalence class $c$  of polyfold structures for the polyfold $(Z, c)$.

 \begin{definition}
 A {\bf strong polyfold bundle}\index{D- Strong polyfold bundle}
  $(p:Y\rightarrow Z,\bar{c})$ consists of a surjective continuous map $p:Y\rightarrow Z$  between topological spaces and  an equivalence class $\bar{c}$ of strong polyfold bundle structures for $p$.
\qed
 \end{definition}

\begin{definition}
The strong polyfold bundle  $(p\colon Y\rightarrow Z,\bar{c})$ is called  
\begin{itemize}
\item[(1)]  \ {\bf tame}\index{D- Tame strong polyfold bundle} if the  polyfold $(Z,c)$ is tame.
\item[(2)]\  {\bf weakly face-structured}\index{D- Weakly face-structured} if $(Z,c)$ is weakly face-structured.
\item[(3)]\  {\bf face-structured} \index{D- Face structured} if$(Z,c)$ is face-structured.
\end{itemize}
\qed
\end{definition}
Next we introduce the category ${\mathcal S}{\mathcal P}{\mathcal B}{\mathcal S}$ \index{Category  ${\mathcal S}{\mathcal P}{\mathcal B}{\mathcal S}$} of {\bf strong polyfold bundle structures}. \index{Category of strong bundle structures} 
The objects of the category are the tuples
$$
\mathsf{S}=((P:W\rightarrow X,\mu), \Gamma,\gamma), 
$$
consisting of a strong bundle $(P:W\rightarrow X,\mu)$ over an ep-groupoid and homeomorphisms 
$\Gamma:|W|\rightarrow Y$ and $\gamma:|X|\rightarrow Z$, where $Y$ and $Z$ are topological spaces (depending on $\mathsf{S}$).
We note that $\gamma \circ |P|\circ \Gamma^{-1}:Y\rightarrow Z$ is a surjective continuous map and we shall denote it by $p$ so that
$$
p:Y\rightarrow Z.
$$
View this as a kind of bundle situation. We note that $Y=Y_{\mathsf{S}}$, $Z=Z_{\mathsf{S}}$, and $p=p_{\mathsf{S}}$.
The following discussion essentially is parallel to the discussion in the polyfold case. 
We shall call $\mathsf{S}$ a {\bf strong polyfold bundle structure} on $p_{\mathsf{S}}$, i.e. $p:Y\rightarrow Z$.
It is sometimes convenient to view $\mathsf{S}=((P:W\rightarrow X,\mu), \Gamma,\gamma)$ as 
$\mathsf{S}=((P:W\rightarrow X,\mu), \Gamma,\gamma,p)$, where $p$ is, of course, redundant information.

\begin{definition}\index{D- The category ${\mathcal S}{\mathcal P}{\mathcal B}{\mathcal S}$}
A {\bf morphism}  in the category ${\mathcal S}{\mathcal P}{\mathcal B}{\mathcal S}$ between two objects $\mathsf{S}$ and $\mathsf{S}'$
with  $\mathsf{S}= ((P\colon W\rightarrow X,\mu), \Gamma,\gamma)$ and $\mathsf{S}'=
((P'\colon W' \rightarrow X',\mu'), \Gamma',\gamma')$
is a $4$-tuple $( \mathfrak{A}, \mathfrak{a},A,a)$ in which $A\colon Y\to Y'$ and $a\colon Z\to Z'$ are continuous maps between topological spaces satisfying 
$$
p_{\mathsf{S}'}\circ A=a\circ p_{\mathsf{S}},
$$
i.e., $A\colon Y\to Y'$ is a bundle map covering the map $a\colon Z\to Z'$. Here $p_{\mathsf{S}}:Y\rightarrow Z$ and $p_{\mathsf{S}'}:Y'\rightarrow Z'$ is the data associated to $\mathsf{S}$ and $\mathsf{S}'$, respectively.
Moreover, $\mathfrak{A}\colon W\to W'$ is a generalized strong bundle map covering the generalized map $\mathfrak{a}\colon X\to X'$    between the underlying ep-groupoids and satisfying the compatibility conditions
$$
\text{$\Gamma' \circ |\mathfrak{A}| = A\circ \Gamma$\quad \text{and}\quad
$\gamma'\circ |\mathfrak{a}| =a\circ \gamma.$}
$$
The morphism $( \mathfrak{A}, \mathfrak{a},A,a)$  is {\bf invertible} (i.e. an isomorphism) if the continuous maps $A\colon Y\to Y'$ and $a\colon Z\to Z'$ are homeomorphism, and $\mathfrak{A}\colon W\to W'$ is a generalized strong bundle isomorphism covering the generalized isomorphism $\mathfrak{a}\colon X\to X'$ between ep-groupoids.
\qed
\end{definition}
The definition is illustrated by the following diagrams.
$$
\begin{array}{ccc}
\begin{CD}
W @>\mathfrak{A}>> W'\\
@V PVV @V P' VV\\
X @>\mathfrak{a}>> X'
\end{CD}
&
\begin{CD}
|W|@ >|\mathfrak{A}| >> |W'|\\
@V \Gamma VV     @V \Gamma' VV\\
Y @> A>> Y'
\end{CD}
&
\begin{CD}
|X| @>|\mathfrak{a}|>> |X'|\\
@V\gamma VV  @V \gamma' VV\\
Z @> a>> Z'
\end{CD}
\end{array}
$$
The composition of two morphisms is defined by
$$
(\mathfrak{A}',\mathfrak{a}',A',a')\circ (\mathfrak{A},\mathfrak{a},A,a)=(\mathfrak{A}'\circ\mathfrak{A},\mathfrak{a}'\circ\mathfrak{a},A'\circ A, a'\circ a).
$$
To simplify the notation we shall abbreviate the morphism 
$$
(\mathfrak{A},A):=(\mathfrak{A}, \mathfrak{a},A,a),
$$
keeping in mind the underlying maps $(\mathfrak{a},a).$
Generalizing the polyfold map we next introduce the notion of a sc-smooth {\bf strong bundle map} between two strong polyfold bundles 
$$
(p\colon Y\rightarrow Z, \ov{c})\quad \text{and}\quad (p'\colon Y'\rightarrow Z', \ov{c}').
$$
Starting with  the strong polyfold bundle $(p\colon Y\rightarrow Z, \ov{c})$ the equivalence class $\ov{c}$ of polyfold structures consists of all the objects 
$$
\mathsf{S}= ((P\colon W\to X, \mu), \Gamma, \gamma)
$$
 of the category ${\mathcal S}{\mathcal P}{\mathcal B}{\mathcal S}$, for which there exists an isomorphism $(\mathfrak{A},A)$ between them of the special form $(\mathfrak{A},Id_Y)$, where $p= \gamma\circ |P|\circ \Gamma^{-1}$.
We now proceed as in the polyfold case.  Take two morphisms
$(\mathfrak{A}_0,A_0)$  and $( \mathfrak{A}_1,A_1)$  in the category 
${\mathcal S}{\mathcal P}{\mathcal B}{\mathcal S}$, where 
\begin{equation*}
\begin{split}
( \mathfrak{A}_0,A_0)\colon& ((P_0\colon W_0\to X_0, \mu_0), \Gamma_0, \gamma_0, p\colon Y\to Z)\\
&\quad \to ((P'_0\colon W'_0\to X'_0, \mu'_0), \Gamma'_0, \gamma'_0, p'\colon Y'\to Z')
\end{split}
\end{equation*}
and 
\begin{equation*}
\begin{split}
( \mathfrak{A}_1,A_1)\colon &((P_1\colon W_1\to X_1), \mu_1), \Gamma_1), \gamma_1, p\colon Y\to Z)\\
&\quad  \to 
((P'_1\colon W_1'\to X_1', \mu_1'), \Gamma_1', \gamma_1', p'\colon Y'\to Z')
\end{split}
\end{equation*}
for two polyfold polyfold bundle structures for $p\colon Y\to Z$, and two polyfold bundle structures for $p'\colon Y'\to Z'$.

\begin{definition}[{\bf Equivalence of morphisms in ${\mathcal S}{\mathcal P}{\mathcal B}{\mathcal S}$}]
The two morphisms $( \mathfrak{A}_0,A_0)$  and $( \mathfrak{A}_1,A_1)$ are {\bf equivalent} if there exists an isomorphism $(\mathfrak{H},Id_Y)$ between the polyfold bundle structures of $p$ and an isomorphism $(\mathfrak{H}',Id_{Y'})$ between the polyfold bundle structures of $p'$ satisfying 
$$
(\mathfrak{A}_1,A_1)\circ ( \mathfrak{H},Id_{Y})=( \mathfrak{H}',Id_{Y'})\circ (\mathfrak{A}_0,A_0).
$$
\qed
\end{definition}
From the definition it follows that $A_1=A_0:Y\rightarrow Y'$.  We also have that $a_1=a_0:Z\rightarrow Z'$. 
In particular the following  diagram holds.
$$
\begin{CD}
|W_0|@>\Gamma_0>> Y @< \Gamma_1 << |W_1|\\
@V|\mathfrak{A}_0| VV   @ V AVV @ V |\mathfrak{A}_1| VV\\
|W_0'|     @>\Gamma_0' >>   Y' @< \Gamma_1' << |W_1'|
\end{CD}
$$
In addition the generalized isomorphisms $\mathfrak{H}:W_0\rightarrow W_1$ and $\mathfrak{H}':W_0'\rightarrow W_1'$ satisfy
$$
\mathfrak{A}_1\circ \mathfrak{H} = \mathfrak{H}'\circ \mathfrak{A}_0.
$$
These diagrams induce similar diagrams for the base spaces.
\begin{definition}[{\bf Strong polyfold bundle maps}] 
A sc-smooth strong polyfold bundle  map
$$\wh{F}\colon (p\colon Y\to Z, \ov{c})\to  (p'\colon Y'\to Z', \ov{c}')$$
is an equivalence class $( \mathfrak{A},A)$ of morphism in 
${\mathcal S}{\mathcal P}{\mathcal B}{\mathcal S}$.
\qed
\end{definition}

\begin{definition}[{\bf Auxiliary norm on strong polyfold bundles}]
\index{D- Auxliary norms for srong polyfold bundles}
An {\bf auxiliary norm} $n$  of a strong polyfold bundle  $(p\colon Y\to  Z, \ov{c})$ over the polyfold $Z$  is a continuous map
$$n\colon Y_{(0,1)}\rightarrow [0,\infty)
$$
with the property,  that for every strong bundle structure $m=((P\colon W\rightarrow X, \mu),\Gamma,\gamma)$  there exists  an auxiliary norm
$N\colon W_{0,1}\rightarrow [0,\infty)$ satisfying  
$$
n(\Gamma (\abs{e}  ))= N (e)
$$
 for every $e\in E_{0,1}$. Here the subset $Y_{0,1}\subset Y$ is given by 
 $Y_{0,1}=\Gamma (\abs{W_{0,1}})$. It is independent of the choice of the strong bundle structure $m\in \ov{c}$ used.
 \qed
\end{definition}

We assume that the polyfold structure $((P'\colon W'\to X',\mu'), \Gamma', \gamma')\in \ov{c}$ possesses the auxiliary norm $N'$. Then every equivalent polyfold structure 
 $((P\colon W\to X,\mu), \Gamma, \gamma)\in \ov{c}$ possesses, in view of Theorem \ref{prop_auxiliary_norm}, 
the auxiliary norm $N=[D]^\ast N'$ which is defined as follows. We choose a representative diagram $W\xleftarrow{\Phi}W''\xrightarrow{\Psi}W'$ for $[D]$. Fixing $e\in W$, we find, using that $\Phi$ is an equivalence, a point $e''\in W''$ and  morphism $\Phi (e'')\to e$ and define 
$$
N(e)=([D]^\ast N') (e):=N'(e'),
$$
 where $e'=\Psi (e'')$. Since $N'\colon W'\to \R^+$ is a functor, the definition does not depend of the choice of the representative diagram. From 
$\abs{[D]}(\abs{e})=\abs{\Psi}\circ \abs{\Phi}^{-1}(\abs{e})=\abs{\Psi}(\abs{e''})=\abs{e'}$ and $\Gamma'\circ \abs{[D]}=\Gamma$, we conclude that 
$$n(\Gamma (\abs{e}  )=N(e)=N'(e')=n(\Gamma' (\abs{e'}  ).$$

Therefore, the auxiliary norm $n$ on the strong polyfold bundle $p$ is well-defined if one of the strong polyfold structures in $\ov{c}$ possesses an auxiliary norm. According to Theorem \ref{ANorm1-prop} this happens if $p\colon Y\to Z$ is a paracompact polyfold $Z$. 

\begin{proposition}[{\bf Existence of an auxiliary norm}]\label{existence_auxiliary_norm} 
A strong polyfold bundle $(p\colon Y\rightarrow Z, \ov{c})$ over a paracompact base polyfold $Z$ possesses an auxiliary norm $n\colon Y_{0,1}\to \R^+$.
\qed
\end{proposition}

\begin{remark}\index{R- Remark on auxiliary norms}
We should remark that an auxiliary norm $n$ is, in general, not a norm on a fiber which is not necessary a vector space.
\qed
\end{remark}

\begin{definition}[{\bf Reflexive $1$-fibers}]\label{DEF1228}\index{D- Strong polyfold bundle with reflexive $1$-fiber}
A strong polyfold bundle $p\colon Y\rightarrow Z$ is said to have {\bf reflexive $1$-fibers}  if there exists a strong polyfold  structure 
$\mathfrak{m}=((P\colon W\rightarrow X, \mu),\Gamma,\gamma)$ in which the fibers  $W_{0,1}$ are reflexive Banach spaces. (The equivalent strong polyfold bundle structures then have also reflexive $1$-fibers).
\qed
\end{definition}

\begin{definition}
[{\bf Reflexive auxiliary norm}]\label{DEF1229}\index{Reflexive auxiliary norm}
A {\bf reflexive auxiliary norm}  of strong polyfold bundle $(p\colon Y\to Z, \ov{c})$ possessing reflexive fibers is an auxiliary norm $n\colon Y_{0,1}\to \R^+$ for which there exists a {\bf reflexive} auxiliary norm $N\colon W_{o,1}\to \R^+$ in a model $m=((P\colon W\to X,\mu), \Gamma, \gamma)\in \ov{c}$ satisfying $n(\Gamma (\abs{w}))=N(w)$ for all $w\in W_{0,1}$.
\qed
\end{definition}
\begin{definition}[{\bf Mixed convergence for $p$}]
Let $(p\colon Y\to Z, \ov{c})$ be a strong polyfold bundle having reflexive $1$-fibers. A sequence $(y_k)\subset Y_{0,1}$ is called {\bf mixed convergent} to $y\in Y_{0,1}$, if there exists a reflexive auxiliary norm $N\colon W_{0,1}\to \R^+$ in the model $m\in \ov{c}$ such that the sequence $\Gamma^{-1}(y_k)\subset \abs{W_{0,1}}$ and the point $\Gamma^{-1}(y)$ have representative $w_k$ resp. $w$ in $W_{0,1}$, for which $w_k\xrightarrow{m}w$ in $W_{0,1}$.
\qed
\end{definition}

From Theorem \ref{EXTTT} we deduce the following result.
\begin{theorem}\index{T- Reflexive auxiliary norms}
Let $p\colon Y\rightarrow Z$ be a strong polyfold bundle over the paracompact polyfold $Z$. If  $p$ has reflexive $1$-fibers, then 
the following holds.
\begin{itemize}
\item[{\em (1)}] \ \ The bundle $p$ admits a reflexive auxiliary norm $n$.
\item[{\em (2)}] \ \ For every  auxiliary norm $n$  for $p$ there exist reflexive auxiliary norms $n_1$ and $n_2$ for $p$ 
satisfying
$$
n_1\leq n\leq n_2.
$$
\end{itemize}
\qed
\end{theorem}

The theorem guarantees a large supply of reflexive auxiliary norms for $p$. Let us collect  some of the properties of such an auxiliary norm.
\begin{proposition}\index{P- Mixed convergence in polyfold bundles}
 If $(p\colon Y\rightarrow Z, \ov{c})$ is a strong polyfold bundle over the paracompact polyfold $Z$ possessing reflexive $1$-fibers, then the following holds for a reflexive auxiliary norm $n$ for $p$.
\begin{itemize}
\item[{\em (1)}] \  If $(y_k)\subset Y_{0,1}$ satisfies $p(y_k)\rightarrow z$ and $n(y_k)\rightarrow 0$,
then $y_k\rightarrow 0_{z}$ in $Y_{0,1}$.
\item[{\em (2)}]\  If $(y_k)$ is mixed convergent to $y\in Y_{0,1}$, then 
$$
n(y)\leq \liminf_{k\rightarrow\infty} n(y_k).
$$
\end{itemize}
\qed
\end{proposition}
In the case that the polyfold $Z$ is paracompact and hence, as we have pointed out, metrizable, an extension results for auxiliary norms defined on the boundary $\partial Z$ to the whole of $Z$ are available.  These results are based on Theorem \ref{EXTT} and 
Theorem \ref{THMOP12213} in the case of reflexive auxiliary norms, and use the fact that 
auxiliary norms behave well under generalized strong bundle isomorphisms. Pullbacks and push-forwards under the  latter
also preserve the relexiveness property.
\begin{theorem}\index{T- Extension of auxiliary norms}\label{THMX2213X}
Let $(p:Y\rightarrow Z,\bar{c})$ be a strong polyfold bundle over the tame paracompact
polyfold $Z$ and assume that $n: Y_{0,1}|\partial Z\rightarrow {\mathbb R}^+$ is an auxiliary norm. Then the following holds true.
\begin{itemize}
\item[(1)]\ Then there exists an auxiliary norm $\overline{n}:Y_{0,1}\rightarrow {\mathbb R}^+$
which extends $n$.
\item[(2)]\ If $n$ is reflexive the extension $\overline{n}$ be be taken to be reflexive as well.
\end{itemize}
\qed
\end{theorem}

A strong polyfold bundle $(p\colon Y\to Z, \ov{c})$ possesses a canonical zero section $z\mapsto 0_z$ induced from the zero sections of the strong bundles $(P\colon W\to X,\mu)$ in the overhead $\ov{c}$.
In view the discussion of the behavior of sc-smooth section of strong bundles over ep-groupoids under generalized strong bundle maps in Section \ref{STE__x}
we can define sc-smooth sections of strong polyfold bundles. The constructions  turn out to be compatible with sc-Fredholm sections and $\ssc^+$-sections.

We start with the  strong polyfold bundle $p\colon Y\rightarrow Z$ and choose a strong polyfold bundle structure  $((P\colon W\to X, \mu), \Gamma, \gamma),$ where 
$P\colon W\rightarrow X$ is a strong bundle over the ep-groupoid $X$ and  $\Gamma\colon \abs{W}\rightarrow Y$ is  a homeomorphism  covering the homeomorphism $\gamma\colon \abs{X}\rightarrow Z$ satisfying $p\circ \Gamma=\gamma \circ \abs{P}.$

In this bundle structure  a sc-smooth section is, by definition, a pair $(f,F)$ in which $f\colon Z\to Y$ is  a continuous section of the strong polyfold bundle $p$,  and $F\colon X\to W$ is a sc-smooth functor of the strong bundle $P$, such that 
$$
\Gamma\circ |F|= f\circ \gamma\quad  \text{on}\ \ |X|.
$$
We represent these  data as a tuple $(f,\Gamma, W,F)$. If we take  a different strong polyfold bundle chart 
$((P'\colon W'\to X', \mu'), \Gamma', \gamma')$ for $p$ we obtain the tuple 
$(f',\Gamma',W',F')$ and introduce the equivalence relation 
$$(f,\Gamma, W,F) \sim(f',\Gamma',W',F'),$$ defined by the requirement $f=f'$ and there exists a generalized strong bundle isomorphism $\mathfrak{A}\colon W\rightarrow W'$ satisfying
$$
\Gamma' \circ |\mathfrak{A}| =\Gamma\quad  \text{and}\quad   \mathfrak{A}_\ast F= F'.
$$
We recall that if $\mathfrak{A}=[W\xleftarrow{\Phi}W''\xrightarrow{\Psi}W']$, then $\mathfrak{A}_\ast=\Psi_\ast\circ \Phi^\ast$.
This defines an equivalence relation $\sim$.

\begin{definition}\index{D- Strong bundle section}
A {\bf sc-smooth section of the strong polyfold bundle} $p\colon Y\rightarrow Z$ consists of an equivalence class
$[(f,\Gamma,W,F)]$, where $(W,\Gamma)$ is an admissible  strong bundle structure for $p$,  $F$ is a sc-smooth section functor
of the strong bundle $P\colon Y\to X$ over the ep-groupoid $X$, and $f$ is a continuous section of $p$ satisfying 
$$f\circ \gamma= \Gamma\circ |F|.$$
\qed
\end{definition}

Any admissible strong bundle structure $(P:W\rightarrow X,\Gamma,\gamma)$ of a strong polyfold bundle $p\colon Y\rightarrow Z$
defines a bi-level structure $Y_{m,k}$ with $0\leq k\leq m+1$ on $Y$.  In particular,  $Y_{0,1}$ is a topological space with a surjective continuous map  $p=p_{0,1}\colon Y_{0,1}\rightarrow Z$. 

\begin{definition}\index{D- $\ssc^+$-polyfold section}
A {\bf $\ssc^+$-smooth section of the strong polyfold bundle} $p\colon Y\rightarrow Z$ consists of an equivalence class
$[(f,\Gamma,W,F)]$, where $(W,\Gamma)$ is an admissible  strong bundle structure for $p$,  $F$ is a $\ssc^+$-smooth section functor
of $P\colon W\rightarrow X$, and $f$ is a continuous section of $p\colon Y_{0,1}\rightarrow Z$ satisfying 
 $$f\circ\gamma= \Gamma\circ \abs{F}.$$
 \qed
\end{definition}
\section{Branched Finite-Dimensional Orbifolds}
The definition of a polyfold utilizes  the notion of a generalized isomorphism.
Our aim 
is to introduce the concept of a branched (weighted) suborbifold of $Z$ as well as  related notions.
That this is possible is rooted in the fact that  a branched ep$^+$-subgroupoid
is well-behaved under generalized isomorphisms. This allows to combine the discussions
in Chapters \ref{CHAPTER_9} and \ref{CHAPTER_11}.  The following results
are easily reduced to already established facts. 

Assume that $(Z,c)$ is a polyfold and $(X,\alpha)\in c$ a polyfold structure. 
As shown in Theorem \ref{THM1133} a generalized isomorphism $\mathfrak{f}:X\rightarrow X'$
can be used to push-forward or pull-back branched ep$^+$-subgroupoids $\Theta$.  Moreover the 
operations $\mathfrak{f}_\ast$ and $\mathfrak{f}^\ast$ preserve many of the properties
such a branched ep$^+$-subgroupoid may have. Namely being closed or compact, being of manifold-type or being of orbifold-type, as well as the dimensional decomposition,
see Definition \ref{DEF915} (closed or compact), Definition \ref{DEF917} (manifold-type or orbifold-type), and Proposition \ref{PROP919} (dimensional decomposition).
Moreover, if $\Theta$ is tame this property is preserved under pull-backs and push-forwards.

\begin{definition}\index{D- Branched suborbifold}
A {\bf branched (weighted) suborbifold} of $(Z,c)$ is an equivalence class
$\theta\equiv [\Theta,(X,\alpha)]$, where $(X,\alpha)\in c$ and $\Theta:X\rightarrow {\mathbb Q}^+$
is a branched ep$^+$-subgroupoid. Here $(\Theta,(X,\alpha))$ and $(\Theta',(X',\alpha'))$
are equivalent provided there exists a generalized isomorphism $\mathfrak{f}:X\rightarrow X'$ 
satisfying
\begin{itemize}
\item[(1)]\ $\alpha'\circ |\mathfrak{f}|=\alpha$.
\item[(2)]  \ $\mathfrak{f}^\ast\Theta'=\Theta$.
\end{itemize}
\qed
\end{definition}
Given $\theta$ and two representatives $(\Theta,(X,\alpha))$ and $(\Theta',(X',\alpha'))$,
 we note that $|\Theta'|\circ |\mathfrak{f}| =|\Theta|$, which implies
that 
$$
|\Theta|\circ \alpha^{-1} = |\Theta'|\circ |\mathfrak{f}|\circ \alpha^{-1} =|\Theta'|\circ \alpha'^{-1}.
$$
This means that we obtain a well-defined map $Z\rightarrow {\mathbb Q}^+$, again denoted by $\theta$, by defining, using a representative $(\Theta,(X,\alpha))$
$$
\theta= |\Theta|\circ \alpha^{-1}.
$$
Of course, there is a lot of overhead associated to this map.

\begin{definition}\index{D- Types of branched suborbifolds}
Let $\theta$ be a branched suborbifold of $(Z,c)$, and let $(\Theta,(X,\alpha))$ 
be a representative of $\theta$. 
\begin{itemize}
\item[(1)] \ $\theta$ is called {\bf closed} provided $\Theta$ is closed.
\item[(2)]  \ $\theta$ is called {\bf compact} provided $\Theta$ is compact.
\item[(3)]  \ $\theta$ is said to be of {\bf manifold-type} provided $\Theta$ is of manifold-type.
\item[(4)] \ $\theta$ is said to be of {\bf orbifold-type} provided $\Theta$ is of orbifold-type.
\item[(5)]   \ $\theta$ is said to be {\bf tame} provided $\Theta$ is tame.
\end{itemize}
\qed
\end{definition}
As a consequence of previously established results
we obtain the following proposition.
\begin{proposition}
If $\theta$ is a branched suborbifold of the polyfold $(Z,c)$
let $A=\supp(\theta)$ be the subset of all $z\in Z$ with $\theta(z)>0$
which is a topological subspace of $Z$. If $\theta$ is of manifold-type,
then $A$ has a natural smooth manifold structure, and if $\theta$ is of orbifold-type,
then $A$ has a natural smooth orbifold structure.  Both these structures
are  induced from the overhead
of $(Z,c)$.
\end{proposition}
\begin{proof}
This follows from Proposition \ref{PROPY918} and Theorem \ref{THM1133}.
\qed \end{proof}
We can also define orientations for branched suborbifolds using Theorem \ref{THM1135}.
\begin{definition}\label{DEFNX16.31}
Let $(Z,c)$ be a polyfold. An {\bf oriented branched suborbifold} \index{D- Oriented branched suborbifold} of $(Z,c)$
is given by an equivalence class $[(\Theta,\wh{\mathsf{T}}_\Theta),(X,\alpha)]$,
where $(X,\alpha)\in c$ and $(\Theta, \wh{\mathsf{T}}_\Theta)$ is an oriented
branched ep$^+$-subgroupoid. The equivalence 
$$
((\Theta,\wh{\mathsf{T}}_\Theta),(X,\alpha))\sim ((\Theta',\wh{\mathsf{T}}_{\Theta'}),(X',\alpha'))
$$
is defined by requiring the existence of a generalized isomorphism $\mathfrak{f}:X\rightarrow X'$
satisfying
$$
\mathfrak{f}^\ast(\Theta',\wh{\mathsf{T}}_{\Theta'})=(\Theta,\wh{\mathsf{T}}_{\Theta})\ \text{and}\ \alpha'\circ|\mathfrak{f}|=\alpha.
$$
\qed
\end{definition}
We shall simplify notation by setting $\wh{\theta}=[(\Theta,\wh{\mathsf{T}}_\Theta),(X,\alpha)]$.

Next we introduce the notion of differential form on a polyfold $(Z,c)$.
\begin{definition}
Given $(X,\alpha)\in c$ a {\bf sc-differential form} is given by an equivalence class
$[\omega,(X,\alpha)]$, where $\omega\in \Omega^\ast_{ep,\infty}(X)$.
Here
$$
(\omega,(X,\alpha))\sim (\omega',(X',\alpha'))
$$
provided there exists a generalized isomorphism $\mathfrak{f}:X\rightarrow X'$
satisfying $\alpha'\circ |\mathfrak{f}|=\alpha$ and $\mathfrak{f}\omega'=\omega$.
\qed
\end{definition}
Given an sc-smooth map $\wh{g}:(Z,c)\rightarrow (Z',c')$ 
the pull-back $\wh{g}^\ast [(\omega',(X',\alpha'))]$ is defined by 
$$
\wh{g}^\ast [(\omega',(X',\alpha'))] = [(\mathfrak{g}^\ast\omega',(X,\alpha))].
$$
One easily verifies that this is well-defined.
Next we introduce the exterior derivative.
\begin{definition}
The exterior derivative $d$ is defined by
\begin{eqnarray}
d([(\omega,(X,\alpha))]) =[(d\omega,(X,\alpha))]
\end{eqnarray}
\qed
\end{definition}
In view of Theorem \ref{OOmaine} this is well-defined and
\begin{eqnarray}
d_{(Z,c)} \left(\wh{g}^\ast  [(\omega',(X',\alpha'))]\right)= \wh{g}^\ast \left(d_{(Z',c')}  [(\omega',(X',\alpha'))]\right).
\end{eqnarray}
As we have seen in Section \ref{SECX9.2}, Proposition \ref{PROPX9.2.15}
a branched ep$^+$-subgroupoid $\Theta:X\rightarrow {\mathbb Q}^+$
 has a well-defined boundary $\partial\Theta:X\rightarrow {\mathbb Q}^+$.
 In general $\partial\Theta$ does not have pleasant properties unless
 $\Theta$ lies in a sufficiently nice position to the boundary, for example
 if $\Theta$ is tame. In the case that $\Theta$ is tame 
 we not only can define $\partial\Theta$, but it also makes sense to talk
 about its orientation.  Namely of given the tame, oriented branched 
 ep$^+$-groupoid $(\Theta,\wh{\mathsf{T}}_\Theta)$  the boundary $\partial\Theta$
 has an induced orientation $\wh{\mathsf{T}}_{\partial\Theta}$ and we define
 $$
 \partial(\Theta,\wh{\mathsf{T}}_\Theta):=(\partial\Theta,\wh{\mathsf{T}}_{\partial\Theta}).
 $$
In view of Theorem \ref{THMX11.3.9} this data behaves well with respect to generalized
isomorphisms and this will allow us ultimately to state a version of Stokes' Theorem.
Indeed, we can define for an oriented, tamed,
branched ep-groupoid $\wh{\Theta}$
$$
\partial([(\Theta,\wh{\mathsf{T}}_\Theta),(X,\alpha)])
=[(\partial\Theta,\wh{\mathsf{T}}_{\partial\Theta}),(X,\alpha)],
$$
and this definition is compatible with the notion of our equivalence relation,
which one easily derives
employing Theorem \ref{THMX11.3.9}.  With other words 
the following definition makes sense.
\begin{definition}
Let $(Z,c)$ be a polyfold and $\wh{\theta}\equiv [(\Theta,\wh{\mathsf{T}}_\Theta),(X,\alpha)]$
an oriented, tame, branched sub-orbifold, i.e. $\wh{\Theta}$ is an oriented, tame, branched
ep$^+$-subgroupoid. Then the {\bf boundary} \index{D- Boundary of a tame, oriented, branched orbifold} $\partial\wh{\theta}$ is defined
by $[\partial(\Theta,\wh{\mathsf{T}}_\Theta),(X,\alpha)]$,
where $\partial(\Theta,\wh{\mathsf{T}}_\Theta)$ is defined in Definition \ref{DEF936}.
\qed
\end{definition}

Finally we use the constructions and results in Section \ref{SECRTY114}
to formulate the Stokes Theorem.  Recall  Remark \ref{REM1142} and 
 assume that $Z$ is a polyfold,
$[(\Theta,\wh{\mathsf{T}}_{\Theta}),(X,\alpha)]$ is an oriented, tame, compact, branched sub$^+$-polyfold of dimension $n$, and $[\omega, (X,\alpha)]$
an sc-differential form of degree $n-1$. In order to simplify notation we put 
$$
\wh{\theta}=[(\Theta,\wh{\mathsf{T}}_{\Theta}),(X,\alpha)]\ \ \text{and}\ \ \tau=[\omega,(X,\alpha)].
$$
We can define $d\tau:=[d\omega,(X,\alpha)]$ and $\partial\wh{\theta}:=
[(\partial\Theta,\wh{\mathsf{T}}_{\partial\Theta}),(X,\alpha)]
$.
These are all objects associated to $Z$, but one needs the overhead to describe the structures.
We can define branched integration as follows.  Let $S\subset Z$ be the support of $\theta$ 
defined as the subset consisting of points $z\in Z$ with $|\Theta|\circ \alpha^{-1}(z)>0$. This definition 
is independent of the representative for $\wh{\theta}$. We denote similarly by $\partial S$ the support
of $\partial\theta$. Pushing forward the Lebesgue $\sigma$-algebras for $|\supp(\Theta)$ and
$|\supp(\partial\Theta)|$ we obtain the corresponding $\sigma$-algebras for $S$ and $\partial S$, respectively. The definition of the latter does not depend on the representative which was taken.  

Given a $n$-form $\tau$ on $Z$ and an oriented, tame, compact, branched sub-orbifold $\wh{\theta}$ of dimension $n$
we define
$$
\oint_{\wh{\theta}}\tau:= \oint_{\wh{\Theta}}\omega = \mu_{\omega}^{\wh{\Theta}}(|\supp(\Theta)|).
$$
The definition does not depend on the choice of the representatives to define it.  There is a similar
definition for a $(n-1)$-form $\tau'$ for the boundary case.
$$
\oint_{\partial\wh{\theta}} \tau' :=\oint_{\partial\wh{\Theta}} \omega' =\mu_{\omega'}^{\partial\wh{\Theta}}(|\supp(\partial\Theta|)).
$$
On the ep-groupoid level Stokes' theorem holds and consequently it holds
in the polyfold context as well.
\begin{theorem}[Stokes]\index{T- Polyfold Stokes}
Given a polyfold $(Z,c)$, a differential $(n-1)$-form $\tau=[\omega,(X,\alpha)]$
with $(X,\alpha)\in c$, and an oriented, tame, compact, branched sub-orbifold
$\wh{\theta}=[(\Theta,\wh{\mathsf{T}}_{\Theta}),(X,\alpha)]$ of dimension $n$
the Stokes type formula
$$
\oint_{\wh{\theta}} d\tau = \oint_{\partial\wh{\theta}}\tau
$$
holds.
\qed
\end{theorem}

\section{Sc\texorpdfstring{$^+$}{pl}-Multisections}\label{SECTX16.5}
In this section we shall employ previously established results concerning sc$^+$-multi\-sec\-tion functors, see Section \ref{SECT133}.
From Theorem \ref{pullbackX} and Theorem \ref{pushforwardX} we know that structurable sc$^+$-multisection functors 
are well-behaved under pullbacks and push-forwards.  The  previously established
fact that structurability behaves well under generalized
strong bundle isomorphisms follows immediately.  Consequently, we obtain a notion of sc$^+$-multisections for strong polyfold bundles.
 As we shall see, also  the pullback of a structurable sc$^+$-multisection  by a proper strong polyfold bundle covering  is structurable.

\begin{definition}\index{D- Sc$^+$-multisection}
Let $(p:Y\rightarrow Z,\bar{c})$ be a strong polyfold bundle.  An {\bf sc$^+$-mul\-tisec\-tion} 
$$
\lambda:Z\rightarrow {\mathbb Q}^+
$$
is  given by an equivalence class
$\lambda\equiv[(\lambda,\Gamma,W,\Lambda)]$, where $(W,\Gamma)$ is an admissible strong bundle structure
$$
\begin{CD}
|W| @>\Gamma>> Y\\
@V |P| VV     @ Vp VV\\
|X| @>\gamma >>  Z
\end{CD}
$$
and $\Lambda: W\rightarrow {\mathbb Q}^+$ a sc$^+$-multisection functor having the property
$$
\lambda \circ \Gamma| (|w|) = \Lambda(w)\ \ \text{for all}\ \ w\in W.
$$
Here $(\lambda,\Gamma,W,\Lambda)$ is {\bf equivalent} to $(\lambda',\Gamma',W',\Lambda')$ provided
$\lambda=\lambda'$, and there exists a generalize strong bundle isomorphism $\bar{\mathfrak{f}}:W\rightarrow W'$,
covering the generalized isomorphism $\mathfrak{f}:X\rightarrow X'$, fitting into the commutative diagram
$$
\begin{CD}
|W| @>|\bar{\mathfrak{f}}|  >> |W'|\\
@V |P| VV  @V |P'| VV\\
|X| @>|\mathfrak{f}| >> |X'|
\end{CD}
$$
and satisfying $\Lambda =\bar{\mathfrak{f}}^\ast\Lambda'$.
\qed
\end{definition}
In Definition \ref{DEFNX1322} we introduced the notion of the domain support for a sc$^+$-multi\-sec\-tion functor,
and we can introduce a related notion here.
\begin{definition}
For the sc$^+$-multisection  $\lambda:Y\rightarrow {\mathbb Q}^+$ the  {\bf domain support}\index{D- Domain support of $\lambda$} is defined by
$$
\text{dom-supp}(\lambda)=\cl_Z(\{z\in Z\ |\ \exists \ y\in p^{-1}(z),\ \text{with}\ \lambda(z)>0\}).
$$ 
\qed
\end{definition}
\begin{remark}
Incidentally, if we consider a representative $(\lambda,\Gamma,W,\Lambda)$, then $\Lambda$
induces a homeomorphism 
$\Gamma : \vert\text{dom-supp}(\Lambda)\vert\rightarrow\text{dom-supp}(\lambda)$.
\qed
\end{remark}

Given an auxiliary norm $n:Y\rightarrow [0,\infty]$ we can measure the point-wise  size of a sc$^+$-multisection
$\lambda:Y\rightarrow {\mathbb Q}^+$ by defining
\begin{eqnarray}
n(\lambda)(z) =\text{max}\{ n(y)\ |\ y\in p^{-1}(z),\ \lambda(y)\neq 0\}.
\end{eqnarray}
Then $n(\lambda):Z\rightarrow {\mathbb R}$ is a continuous map.

Since structurability is well-behaved under generalized strong bundle isomorphisms the following definition makes sense.
\begin{definition}\index{D- Structurable sc$^+$-multisection}
Let $(p:Y\rightarrow Z,\bar{c})$ be a strong polyfold bundle and $\lambda:Z\rightarrow {\mathbb Q}^+$ a 
sc$^+$-multisection.  We say that $\lambda$ is {\bf structurable} if one, and therefore all of its representatives $\Lambda:W\rightarrow {\mathbb Q}^+$ are structurable.  
\qed
\end{definition}

Assume that we are given a a strong polyfold bundle $p:Y\rightarrow X$ 
with a tame paracompact base $X$, an auxiliary norm
$n:Y_{0,1}\rightarrow {\mathbb R}^+$, and a functor 
$\lambda:Y|\partial Z\rightarrow {\mathbb Q}^+$.  
In a first step we define what is means that $\lambda$ is a sc$^+$-multisection.
\begin{definition}\index{D- Sc$^+$-multisection over $\partial Z$}
A sc$^+$-multisection of $\lambda$ of $Y|\partial Z$ is an equivalence class
$$
\lambda\equiv [(\lambda,\Gamma,W,\Lambda)],
$$
 where $P:W\rightarrow X$
is an sc-smooth local model for $p:Y\rightarrow Z$, and $\Lambda:W|\partial X\rightarrow
{\mathbb Q}^+$ an sc$^+$-multisection functor satisfying
$$
\lambda\circ \Gamma(\vert w\vert) =\Lambda(w)\ \ \text{for all}\ w\in W,\ P(w)\in \partial X.
$$
Two tuples $(\lambda,\Gamma,W,\Lambda)$ and $(\lambda',\Gamma',W',\Lambda')$ are said to be equivalent provided $\lambda=\lambda'$ and the there exists a generalized strong bundle isomorphism 
$$
\begin{CD}
W @>\overline{\mathfrak{f}} >> W'\\
@ V P VV @V P' VV\\
X @>\mathfrak{f} >>  X'
\end{CD}
$$
(the diagram interpreted with the usual care) such that 
 $$
 \Lambda =\overline{\mathfrak{f}}^\ast\Lambda'\ \ \text{and}\ \ \ \Gamma'\circ \vert\overline{\mathfrak{f}}\vert = \Gamma.
 $$
\qed
\end{definition}
Using Definition \ref{corr-xxxx} and Proposition \ref{prop14.1}
we can define structurable sc$^+$-multi\-sec\-tions over $\partial Z$.
\begin{definition}\index{D- Structurable sc$^+$-multisection over $\partial Z$}
Let $p:Y\rightarrow Z$ be a strong bundle over the tame paracompact polyfold $Z$.
Assume that $\lambda: Y|\partial Z\rightarrow {\mathbb Q}^+$ is an sc$^+$-multisection 
over the boundary $\partial Z$. We say that $\lambda$ is {\bf structurable} provided
one of its representative $(\lambda,\Gamma,W,\Lambda)$ has a structurable $\Lambda:W|\partial X\rightarrow {\mathbb Q}^+$.
\end{definition}

Using Theorem \ref{p-main-p} we can state an extension theorem for structurable sc$^+$-multi\-sect\-ions.

\begin{theorem}
Let $(p:Y\rightarrow Z,\bar{c})$ be a strong polyfold bundle over the tame paracompact polyfold $Z$.
We assume that $Z$ admits sc-smooth partitions of unity. Denote by $n:Y\rightarrow [0,\infty]$ an auxiliary norm
and $\lambda :Y|\partial Z\rightarrow {\mathbb Q}^+$ a structurable sc$^+$-multisection. Assume that
$\wt{U} \subset Z$ is an open neighborhood of $\text{dom-supp}(\lambda)$ and $f:Z\rightarrow[0,\infty)$ a continuous map
with support in $\wt{U}$ satisfying $n(\lambda)(z)<f(z)$ for $z\in \text{dom-supp}(\lambda)$. 
Then there exists a structurable sc$^+$-multisection $\lambda':Y\rightarrow {\mathbb Q}^+$ with the following properties.
\begin{itemize}
\item[{\em(1)}]\ $n(\lambda')(z)\leq f(z)$ for all $z\in Z$.
\item[{\em(2)}]\ $\text{dom-supp}(\lambda')\subset \wt{U}$.
\item[{\em(3)}]\ $\lambda'|(Y|\partial Z)=\lambda$.
\end{itemize}
\qed
\end{theorem}
As a corollary of 
Theorem \ref{OTHM1344} and the study of proper coverings up to equivalences in Section \ref{SQWERTY116}
one is  able to define the pull-back of a sc$^+$-multisections by a proper covering.
We leave the details to the reader.  

In practice, when one uses the theory to define invariants via sc-Fredholm theory, 
 it is usually easier to work in  an ep-groupoid model.  One producea the invariants 
 in this context, and then  concludes from the transformation properties
 of generalized isomorphisms, that the invariants are independent of the choices.
 Particularly, invariants involving proper coverings, need the fine structure provided
 by the ep-groupoid description.

\section{Fredholm Theory}
In this section we develop the sc-Fredholm theory for sc-smooth sections of strong polyfold bundles.
Essentially this is an exercise about the sc-Fredholm theory in ep-groupoids and its behavior with respect to equivalences.

With   $p:Y\rightarrow Z$ being a strong polyfold bundle recall that an sc-smooth section 
is given by an equivalence class $[(f,\Gamma,W,F)]$, where $(W,\Gamma)$ is a strong bundle structure for
$p$, i.e. 
$$
\begin{CD}
|W|  @>\Gamma>>  Y\\
@V|P| VV   @V p VV\\
|X| @>\gamma>>   Z,
\end{CD}
$$
$F$ is an sc-smooth section of $P:W\rightarrow X$ and $f$ is a continuous section of $p$  such that 
$$
\Gamma \circ |F|  = f\circ \gamma.
$$
Here $\Gamma$ and $\gamma$ are homeomorphisms.  
  Assume that $(f,\Gamma,W,F)$ is a representative of the equivalence class and that $F$ is an sc-Fredholm section.
If $(f,\Gamma',W',F')$ is another representative we can conclude that $F'$ is an sc-Fredholm section 
of $P':W'\rightarrow X'$.  In order to see this note that $F$ and $F'$ are related by a strong bundle isomorphism
$$
\begin{CD}
W@> \mathfrak{F}>>  W'\\
@V P VV @V P'VV\\
X @>\mathfrak{f} >> X",
\end{CD}
$$
where one has to exercise the usual caution when interpreting this diagram.
We have already seen that sc-Fredholm sections transform well under push-forwards and 
pull-backs of strong bundle maps.  Since the action of $\mathfrak{F}$  on $F$ can be represented 
as a pull-back followed by a push-forward and since $F'=\mathfrak{F}_\ast F$ it follows that 
$F'$ is sc-Fredholm. See Theorem \ref{useful} and Theorem \ref{Push-Forw-prop} for the necessary background material. As a consequence of this discussion we can give the definition of an sc-Fredholm section 
of a strong polyfold bundle.

\begin{definition}\index{D- Polyfold Fredholm section}
A {\bf sc-smooth Fredholm section} of the strong polyfold bundle  $p:Y\rightarrow Z$ consists of an equivalence class
$[(f,\Gamma,W,F)]$, where $(W,\Gamma)$ is an admissible  strong bundle structure for $p$,  $F$ is a sc-smooth Fredholm section functor
of the strong bundle $P\colon W\rightarrow X$, and $f$ is a continuous section of $p$ satisfying  
$$f\circ \gamma= \Gamma\circ \abs{F}.
$$
\qed
\end{definition}

The following compactness definition is basic.
\begin{definition}\index{D- Compact polyfold Fredholm section}\index{Compact solution set}
Let $f\equiv [(f,\Gamma,W,F)]$ be an sc-Fredholm section of the strong polyfold bundle $p:Y\rightarrow Z$.
Then $f$ is said to have a {\bf compact solution} set provided the collection of all $z\in Z$ such that $f(z)=0$ is a compact subset of $Z$.
\end{definition}

Due to the fact that for an sc-Fredholm section in the ep-groupoid setting different compactness notions are equivalent we obtain the following 
stability result as a consequence of Theorem \ref{THMB1242}.
\begin{theorem}
Let $f\equiv [(f,\Gamma,W,F)]$ be an sc-Fredholm section of the strong polyfold bundle $p:Y\rightarrow Z$, where $Z$ is paracompact.
Suppose that $n:Y_{0,1}\rightarrow {\mathbb R}^+$ is an auxiliary norm and $f$ has a compact solution set.
Then there exists an open neighborhood $U$ of $S_f=\{z\in Z\ |\ f(z)=0\}$ such that the closure of the set of all
$z\in U$ with $n(f(z))\leq 1$ is compact.
\end{theorem}
\begin{proof}
We take a model $P:W\rightarrow X$, which by assumption will be a strong bundle over the ep-groupoid $X$.
By assumption $X$ has a paracompact orbit space $|X|$.  Moreover, we are given an sc-Fredholm section functor $F$ of $P$
which fits into the commutative diagram
$$
\begin{CD}
|W| @> \Gamma >> Y\\
@A |F| AA     @A f AA\\
|X| @>\gamma>> Z,
\end{CD}
$$
where $\Gamma$ and $\gamma$ define the strong polyfold bundle structure. The auxiliary norm $n$ is represented in 
the model by $N:W_{0,1}\rightarrow {\mathbb R}^+$.  Employing Theorem \ref{THMB1242} we conclude that not only
 $S_F=\{x\in X\ |\ F(x)=0\}$ has a compact orbit space $|S_F|$, but that also there exists a saturated open neighborhood
 $U'$ of $S_F$ so that the closure of the orbit space of $\{x\in U'\ |\ N(F(x))\leq 1\}$ is compact.  
 Then $U=|U'|$ has the desired properties.
\qed \end{proof}
We can consider oriented sc-Freholm sections of strong polyfold bundles, since orientations of sc-Fredholm section also behave well under generalized strong bundle isomorphisms.  
We begin by reviewing the orientation discussion in the ep-groupoid setting.
Assume that $P:W\rightarrow X$ is a strong bundle over the tame ep-groupoid $X$ and $F$ an sc-Fredholm functor.
As constructed in Section \ref{SEC115} there is a well-defined orientation bundle, see Definition \ref{DEFNG1252}, associated to $f$
and denoted by
$$
\sigma:\mathscr{O}_F\rightarrow X_\infty.
$$
Here $\mathscr{O}_F$ is a topological space, $\sigma$ is a local homeomorphism and 
the map $\sigma$ is a $2:1$ covering.  As discussed in Section \ref{SEC115},
given a  point $x\in X_\infty$, there exist precisely two local continuous
section germs denoted by $\pm\mathfrak{o}_{(x)}$. Of course, there is no preferred local choice.
 The value at $x$ is denoted by $\pm\mathfrak{o}_x$.
Given a smooth morphism $\phi\in \bm{X}$, i.e. $\phi\in \bm{X}_\infty$, we obtain a local
germ of sc-diffeomorphism 
$$
\wh{\phi}:U(s(\phi))\rightarrow U(t(\phi))
$$
 between small open neighborhoods in $X$
and this germ can be lifted to a local strong bundle isomorphism $\wh{\Phi}: W|U(s(\phi))\rightarrow W|U(t(\phi))$,
which can used to push $\mathfrak{o}_{(x)}$ forward.
In particular $\phi$ defines a push forward $\phi_\ast\mathfrak{o}_x$.
Recall that $F$ is orientable provided $\sigma:\mathscr{O}_F\rightarrow X_\infty$ admits a global continuous
section $\mathsf{o}$ with the property that
$$
\phi_\ast \mathsf{o}_{s(\phi)} = \mathsf{o}_{t(\phi)}\ \ \text{for all  morphisms in}\ \ \bm{X}_\infty.
$$
With other words $\mathsf{o}$ is a continuous section functor of $\sigma$.
The choice of a specific section functor  $\mathsf{o}$ is called an orientation for $F$ and denoted by $\mathsf{o}_F$.
If $|X_\infty|$ is connected, we see that there is at most one orientation, and if $F$  is orientable there are precisely two possible orientations.  

Relevant for our orientation discussion in the polyfold framework is Theorem \ref{THMS1257}
which we recall for convenience. We have adapted the notation to the current situation\par

\noindent{\bf Theorem \ref{THMS1257}\textcolor{red}{.}}
Assume $(P:W\rightarrow X,\mu)$  and $(P':W'\rightarrow X',\mu')$  are strong bundles over tame ep-groupoids and $
\bar{\mathfrak{f}}:W\rightarrow W'$ is  a generalized strong bundle isomorphism covering the generalized isomorphism
$\mathfrak{f}:X\rightarrow X'$.  We assume that $F$ is an sc-Fredholm section functor of $P$ and denote by $F'$ the push-forward sc-Fredholm section.
\begin{itemize}
\item[(1)]\  With $\sigma$ and $\sigma'$ being the orientation bundles, the generalized 
strong bundle isomorphism induces a fiber-preserving homeomorphism $|\bar{\mathfrak{f}}|_\ast$ fitting into the following commutative diagram
$$
\begin{CD}
|{\mathscr{O}}_F| @> {|\bar{\mathfrak{f}}|}_\ast>> |{\mathscr{O}}_{F'}|\\
@V |\sigma| VV  @V |\sigma'| VV\\
|X_{\infty}| @>|\mathfrak{f}|>> |X'_{\infty}|.
\end{CD}
$$
The construction of $|\bar{\mathfrak{f}}|_\ast$ is functorial.  
\item[(2)]\ Assume that 
 $\mathsf{o}$ is a continuous section functor of $\sigma$, so that in particular the isotropy groups $G_x$, $x\in X_\infty$, act trivially on $\sigma^{-1}(x)$.
Then $\bar{\mathfrak{f}}$ defines naturally a push-forward $\mathsf{o}'=\bar{\mathfrak{f}}_\ast\mathsf{o}$,
which is a continuous section functor of $\mathscr{O}_{F'}\rightarrow X_\infty'$. 
\end{itemize}
\qed\par
Now we are ready for the discussion of orientations for sc-Fredholm sections of strong bundles over tame polyfolds.
From Proposition \ref{PROPT1256} we can conclude the following. Let $P:W\rightarrow X$ be a strong bundle over a tame ep-groupoid with a paracompact orbit space and $F$ is a sc-Fredholm section. We denote by $\sigma:\mathscr{O}_F\rightarrow X_\infty$ the orientation bundle. 
Then $F$ is orientable provided the following two conditions hold.
\begin{itemize}
\item[(1)]\ The continuous map $\vert\sigma\vert:\vert\mathscr{O}_F\vert\rightarrow \vert X_\infty\vert$ has above every class $\vert x\vert$ precisely two preimages.
\item[(2)]\ $\vert\sigma\vert:\vert\mathscr{O}_F\vert\rightarrow \vert X_\infty\vert$ admits a global continuous section.
\end{itemize}
A choice of section of $\vert\sigma\vert$ then by definition is an orientation of $F$. An orientation for $ [(f,\Gamma,W,F)]$ can be defined as a choice of orientation for a representative $(W,F)$. If we take another representative 
$(W',F')$ we obtain via Theorem \ref{THMS1257} an induced orientation. Hence we see that there is a well-defined notion of orientation.
\begin{definition}\index{D- Orientation for $f$}
An {\bf orientation} (if it exists) for a sc-Fredholm section $f$ of the tame strong polyfold bundle $p:Y\rightarrow Z$ consists
of the choice of an orientation for a representative $(W,F)$.
\qed
\end{definition}
As a simple application we note the following result where we use generic perturbations of $f$ by sc$^+$-multi\-sections $\lambda$.
For such a $\lambda$ the maps $\theta =\lambda\circ f:Z\rightarrow {\mathbb Q}^+$ defines branched suborbifold
of $Z$ and if $f$ is oriented $\theta$ has  a natural orientation. This follows from the discussion in Section \ref{SEC153}.
See in particular Theorem \ref{THM1541} and the subsequent Definition \ref{DEF1542}.
\begin{theorem}
Let $p:Y\rightarrow Z$ be a strong polyfold bundle over the paracompact polyfold $Z$  with $\partial Z=\emptyset$.
Assume that $f$
is an oriented sc-Fredholm section of $p$ with compact solution set. Suppose the Fredholm index is $k$.
Let $n:Z_{0,1}\rightarrow {\mathbb R}^+$ be an auxiliary norm and $U$ an open neighborhood of $S_f=\{z\in Z\ |\ f(z)=0\}$
so that $(n,U)$ control compactness. Then, for given sc$^+$-multisection $\lambda$ controlled by $(n,U)$,
so that $f,\lambda)$ is in general position, the integral $\oint_{\wh{\theta}} \tau$ 
for  given  closed sc-differential $k$-form $\tau$ on $Z$ is well-defined and independent of the perturbation $\lambda$.
Moreover, the integral only depends on the cohomology class $[\tau]$. Here $\wh{\theta}$ is $\theta=\lambda\circ f$ equipped
with the inherited orientation from $f$.
\qed
\end{theorem}

We leave it to the reader to bring more of the results from Part III into the current context. We just note that 
that as along as the definitions and results are invariant under generalized strong bundle isomorphisms 
a result in the polyfold context can be obtained.  In application it is usually advantageous to work in 
the ep-groupoid context, however, sometimes the polyfold formulation allows for more compact statements, in particular
if the overhead is suppressed in the statements.

%
%
%

\begin{partbacktext}
\part{Fredholm Theory in Groupoidal Categories}
\noindent  In this part we shall bring the ep-groupoid theory into the more general framework
of groupoidal categories. First we shall develop the theory of sc-smooth structures on a certain class of groupoidal categories,
This is followed by a theory of sc-Fredholm section functors. Finally we develop in Chapter \ref{CHAPTER19X}  useful ideas to construct polyfold structures in applications. For example the papers \cite{FH2,FH3} are build on this framework.  In \cite{H2} an application of the ideas to Gromov-Witten theory
is outlined.
\end{partbacktext}

\chapter{Polyfold Theory for  Categories}\label{CHAPX17}
In this part we shall develop a theory of Fredholm functors for certain categories.
The hard work has already been done in the previous parts.  The resulting theory
is very convenient in applications since it provides a transparent language 
with a large body of results. The construction of symplectic field theory is an excellent example for the use of this `categorical  polyfold theory'. 
We refer the reader to \cite{FH2}  for the illustration of the methods
in this highly nontrivial example.

The  motivated reader is also invited to turn the Gromov-Witten (GW) example into the categorical framework.
The reference \cite{H2} gives an outline of such a theory. The paper \cite{HWZ5} gives a complete and detailed construction of GW, however,  the motivation of this paper was to 
discuss the issues and to show  how they can be phrased and dealt with in the polyfold framework.  The aim was not to give an abstract stream-lined
proof in which important issues might not be as visible, since they are taken care of by the `machinery'.   Of course, once one has enough experience and knows the inner workings 
of the machinery, it perfectly makes sense to use the technology in a slick and stream-lined way, which is the approach to SFT in \cite{FH2}.
Using ideas from  \cite{FH2}, for example the  construction functors for polyfolds, and some of the technical work in \cite{HWZ5},   the 
theory  in the current  book allows a rather fast and complete construction of the Gromov-Witten invariants.

In the following we assume the reader familiar with the material of Part I to Part III and we allow ourselves at time 
to be sketchy if the constructions and arguments are straight forward.

\section{Polyfold Structures and Categories}

The starting point is a category $\mathscr{C}$ with additional  structures. Usually we shall denote the object class
associated to $\mathscr{C}$ by $\text{obj}(\mathscr{C})$ and the morphism class by $\text{mor}({\mathscr{C}})$.
However, we sometimes find it convenient to abbreviate $C=\text{obj}(\mathscr{C})$ and $\bm{C}=\text{mor}({\mathscr{C}})$
if the context is clear.  We shall  consider categories where every morphism is an isomorphism.
In this case we can consider for every object $\alpha$ the class $|\alpha|$  of objects isomorphic to $\alpha$.
The {\bf orbit space}\index{Orbit space} $|\mathscr{C}|$ of $\mathscr{C}$ is the class consisting of all isomorphism classes of objects. 

\begin{definition}\label{DEF_GCT}\index{D- GCT}
A {\bf GCT}, where GC stands for groupoidal category and T for topology, is given by a pair
$(\mathscr{C},{\mathcal T})$, where $\mathscr{C}$ is a category having the following three properties.
\begin{itemize}
\item[(1)]\  The orbit space ${\abs{\mathscr{C}}}$ of the category $\mathscr{C}$ is a set.
\item[(2)]\ The set of morphisms between any two objects is finite.
\item[(3)]\  All morphisms are isomorphisms.
\end{itemize}
In addition  ${\mathcal T}$ is a topology on the orbit set $|{\mathscr{C}}|$ which is metrizable.
Instead of saying $(\mathscr{C},{\mathcal T})$ is  a GCT we sometimes say it is groupoidal category with (metrizable) topology, where it is understood
that the topology is defined on the orbit space.
 We also sometimes refer to the GCT $\mathscr{C}$ and suppress the topology in the notation.
 \qed
\end{definition}

Our main goal is to define the notion of an sc-smooth structure on $(\mathscr{C},{\mathcal T})$. It will turn out that 
this is the starting point for a very rich theory, which is a mixture of a generalized differential geometry,
blended with nonlinear functional analysis and category theory. 

\begin{definition}\index{D- Translation groupoid}
 Given a M-polyfold $O$,  a finite group $G$, and a group homomorphism 
 $\gamma:G\rightarrow \text{Diff}_{sc}(O)$, denote the associated group action by sc-diffeomorphisms by
$$
G\times O\rightarrow O: (g,x)\rightarrow g\ast x,\ \ \text{where}\ g\ast x=\gamma(g)(x).
$$
The  {\bf translation groupoid}\index{Translation groupoid} $G\ltimes O$ associated to $(O,G,\gamma)$ is the category
with $\text{obj}(G\ltimes O)=O$ and $\text{mor}(G\ltimes O)=G\times O$. The source and target maps are defined by
$$
s(g,x)=x\ \text{and}\ \ t(g,x) =g\ast x,
$$
so that $(g,x)$ is viewed as a morphism $(g,x):x\rightarrow g\ast x$.
\qed
\end{definition}
One easily verifies that $G\ltimes O$ is an ep-groupoid. For the following we note that ${_G}\backslash O=|G\ltimes O|$ is a metrizable space. Namely take any metric $d'$ on $O$ and average to obtain the metric $d$ defined by
$$
d(q,q') =\frac{1}{|G|}\cdot  \sum_{g\in G} d'(g\ast q,g\ast q').
$$
Then define $\bar{d}$ on the quotient by 
$$
\bar{d}([q],[q']):= \text{min}_{g,g'\in G} d(g\ast q,g'\ast q').
$$
The building blocks for our upcoming theory are given in the following definition.
\begin{definition}\index{D- Local uniformizer}\index{D- Tame local unifomizer}
Let $\mathscr{C}$ be a GCT and $\alpha$ an object in ${\mathscr{C}}$. A {\bf local uniformizer}\index{Local uniformizer} for $\mathscr{C}$ at $\alpha$ is a covariant functor $\Psi:G\ltimes O\rightarrow \mathscr{C}$ such that the following holds.
\begin{itemize}
\item[(1)] \  $G\ltimes O$ is a translation groupoid associated to the M-polyfold $O$ equipped with an sc-smooth action by a finite group $G$, i.e. the data $(O,G,\gamma)$. 
\item[(2) ]\   $O$ is  paracompact (and hence metrizable).
\item[(3)] \  The functor $\Psi$ is injective on objects and there exists $q_0\in O$ with $\Psi(q_0)=\alpha$.
Moreover $\Psi$ is  full and faithful, i.e. fully faithful.
\item[(4)] \  The induced map $|\Psi|:{_{G\backslash}O}\rightarrow |\mathscr{C}|$ is a homeomorphism onto an open neighborhood $U$ of $|\alpha|$ in $|\mathscr{C}|$.
\end{itemize}
Associated to a local uniformizer $\Psi$ we have its {\bf footprint}\index{Footprint} defined by
$$
\text{footprint}(\Psi)=|\Psi(O)|,
$$
which is an open subset of $|\mathscr{C}|$.   If $O$ is tame we shall refer to $\Psi$ as a {\bf tame local uniormizer}.
\qed
\end{definition}
Denote local uniformizers  at objects $\alpha$ and $\alpha'$ by $\Psi$ and $\Psi'$. 
\begin{definition}\index{D- Transition set}
The {\bf transition  set}  ${\bm{M}}(\Psi,\Psi')$  is defined as the weak fibered product associated to the diagram 
$$
O\xrightarrow{\Psi} \mathscr{C}\xleftarrow{\Psi'} O'.
$$
It consists of all tuples $(q,\phi,q')\in O\times \bm{C} \times O'$ so that $\phi:\Psi(q)\rightarrow \Psi'(q')$.
\qed
\end{definition}
We first note that given three objects we have the following {\bf structure maps}\index{Structure maps} for associated local uniformizers
$\Psi,\Psi'$ and $\Psi''$.
\begin{itemize}
\item[(1)] \ \ \ The {\bf source and target maps}\index{Source map}\index{Target map} $s:{\bm{M}}(\Psi,\Psi')\rightarrow O:(q,\phi,q')\rightarrow q$ and $t:{\bm{M}}(\Psi,\Psi')\rightarrow O':(q,\phi,q')\rightarrow q'$.
\item[(2)]\ \ \  The $1$-map or {\bf unit map} \index{Unit map} $O\rightarrow {\bm{M}}(\Psi,\Psi):q\rightarrow (q,1_{\Psi(q)},q)$
\item[(3)] \ \ \ The {\bf inversion map}\index{Inversion map} ${\bm{M}}(\Psi,\Psi')\rightarrow {\bm{M}}(\Psi',\Psi):(q,\phi,q')\rightarrow (q',\phi^{-1},q)$.
\item[(4)] \ \ \ The {\bf multiplication map}\index{Multiplication map} ${\bm{M}}(\Psi',\Psi''){_{s}\times_t}{\bm{M}}(\Psi,\Psi')\rightarrow{\bm{M}}(\Psi,\Psi'')$ defined by
$$
m((q',\phi',q''),(q,\phi,q'))= (q,\phi'\circ\phi,q'').
$$
\end{itemize}

Denote the category of sets by $\text{SET}$. We shall write $\mathscr{C}^-$ for the category obtained from $\mathscr{C}$ by keeping all objects,
but only allowing the identities as morphisms.
\begin{definition}\index{D- Uniformizer construction}
A {\bf uniformizer construction} is given by  a functor  $F:\mathscr{C}^-\rightarrow \text{SET}$  
which associates
to an object $\alpha$ a set of local uniformizers $F(\alpha)$ at $\alpha$.
\qed
\end{definition}
\begin{remark}\label{REM1716} \index{R- On uniformizer constructions}
Here are some remarks about functors $F$ and the way they occur in applications.\par

\noindent  a) A unformizer construction appears in application usually as a well-defined recipe for a construction involving
for a given object $\alpha$ set-many choices resulting in the set $F(\alpha)$. \par

\noindent b) In many circumstances we have more structure, namely
$F$ is in fact a functor $F:\mathscr{C}\rightarrow \text{SET}$ reflecting the often occurring feature 
that given an isomorphism $\phi:\alpha\rightarrow \alpha'$ one can 
 establish a  bijective correspondence between the choices available for $\alpha$ and $\alpha'$, 
resulting in a bijection $F(\phi):F(\alpha)\rightarrow F(\alpha')$. \par

\noindent c)   In many cases,  $F(\phi)$ even has a geometric interpretation and the following describes the additional structures which one might see.
A frequent situation is as follows.  Given $\Psi\in F(\alpha)$ and $\phi:\alpha\rightarrow \alpha'$ define $\Psi'=F(\phi)(\Psi)$.  By construction $\Psi:G\ltimes O\rightarrow \mathscr{C}$ and $\Psi':G'\ltimes O'\rightarrow \mathscr{C}$.  Assume  there exists a well-defined group isomorphism $\gamma_{(\phi,\Psi)}:G\rightarrow G'$ and an equivariant
sc-diffeomorphism $\sigma_{(\phi,\Psi)}:O\rightarrow O'$ defining an sc-diffeomorphic functor
$$
\bar{\sigma}_{(\phi,\Psi)}:=\gamma_{(\phi,\Psi)}\ltimes \sigma_{(\phi,\Psi)}:G\ltimes O\rightarrow G'\ltimes O'.
$$
Having  the functors $\Psi$ and $\Psi'\circ \bar{\sigma}_{(\phi,\Psi)}$ defined on $G\ltimes O$ we assume that 
there is a natural construction of a natural transformation 
$$
\tau_{(\phi,\Psi)}: \Psi\rightarrow \Psi'\circ \bar{\sigma}_{(\phi,\Psi)}.
$$
In particular, for $q\in O$ the object $\Psi(q)$ and the object $\Psi'(\sigma_{(\phi,\Psi)}(q))$ are isomorphic via
$\tau_{(\phi,\Psi)}(q)$. For $g\in G$ the morphism $(g,q)$ which induces $\Psi(g,q):\Psi(q)\rightarrow \Psi(g\ast q)$
corresponds to $\Psi'(\bar{\sigma}_{(\phi,\Psi)}(g,q)):\Psi'(\sigma_{(\phi,\Psi)}(q))\rightarrow \Psi'(\sigma_{(\phi,\Psi)}(g\ast q))$.
The latter comes from the commutative diagram
$$
\begin{CD}
\Psi(q) @> \Psi(g,q)>> \Psi(g\ast q)\\
@V \tau_{(\phi,\Psi)}(q)VV @V \tau_{(\phi,\Psi)}(g\ast q)VV\\
\Psi'({\sigma}_{(\phi,\Psi)}(q) )@>\Psi'(\bar{\sigma}_{(\phi,\Psi)}(g,q))>>\Psi'(\sigma_{(\phi,\Psi)}(g\ast q))
\end{CD}
$$
In many cases the constructions  $(\phi,\Psi)\rightarrow \sigma_{(\phi,\Psi)}$, $(\phi,\Psi)\rightarrow \gamma_{(\phi,\Psi)}$,
and $(\phi,\Psi)\rightarrow \tau_{(\phi,\Psi)}$ satisfy the following compatibilities
\begin{itemize}
\item[(1)] \  $\sigma_{(\phi',F(\phi)(\Psi))}\circ\sigma_{(\phi,\Psi)}=\sigma_{(\phi'\circ\phi,\Psi)}$.
\item[(2)] \  $\gamma_{(\phi',F(\phi)(\Psi))}\circ \gamma_{(\phi,\Psi)}=\gamma_{(\phi'\circ \phi,\Psi)}$.
\item[(3)] \   $\tau_{(\phi',F(\phi)(\Psi))}(\sigma_{(\psi,\Psi)}(q))\circ\tau_{(\phi,\Psi)}(q)=\tau_{(\phi'\circ \phi,\Psi)}(q)$ for $q\in O$.
\end{itemize}
\qed
\end{remark}

\begin{definition}\index{D- Polyfold structure for $(\mathscr{C},{\mathcal T})$}
A {\bf polyfold structure}   for the GCT  $\mathscr{C}$ consists of two constructions:
\begin{itemize}
\item[(1)]\ \ \ A uniformizer construction  $F:\mathscr{C}^-\rightarrow \text{SET}$.
\item[(2)]\ \ \  A construction of  a M-polyfold structure for every transition set $\bm{M}(\Psi,\Psi')$, where $\Psi\in F(\alpha)$ and $\Psi'\in F(\alpha')$.
\end{itemize}
These two constructions satisfy the following compatibility conditions:
\begin{itemize}
\item[(a)] \ \ \ The source and target maps 
$$
O\xleftarrow{s} \bm{M}(\Psi,\Psi')\xrightarrow{t} O'
$$
 are local sc-diffeomorphisms.
\item[(b)] \ \ \ The unit map $u: O\rightarrow \bm{M}(\Psi,\Psi)$ and the inversion map $\iota:\bm{M}(\Psi,\Psi')\rightarrow \bm{M}(\Psi',\Psi)$ are sc-smooth.
\item[(c)] \ \ \ The multiplication map ${\bm{M}}(\Psi',\Psi''){_{s}\times_t}{\bm{M}}(\Psi,\Psi')\rightarrow{\bm{M}}(\Psi,\Psi'')$ is sc-smooth.
\end{itemize}
 We shall refer to $(F,{\bm{M}})$ as a {\bf polyfold structure}\index{Polyfold structure}. If all local uniformizers associated to $F$ are tame, we shall  call it a {\bf tame polyfold construction}\index{Tame polyfold construction}.  A {\bf polyfold}\index{D- Polyfold} $\mathscr{C}$ consists of a GCT $\mathscr{C}$
 together with a polyfold structure $(F,\bm{M})$ 
\qed
\end{definition}
We can draw the following consequences from this definition. 
\begin{lemma}\label{LEMMX17.1.8}
Given $\Psi:G\ltimes O\rightarrow {\mathcal S}$ in $F(\alpha)$ with $\Psi(q)=\alpha$
the map 
\begin{eqnarray}\label{EQNX17.1}
\hat{\Psi}: G\ltimes O\rightarrow \bm{M}(\Psi,\Psi):(g,p)\rightarrow (p,\Psi(g,p),g\ast p)
\end{eqnarray}
 is 
a sc-diffeomorphism. Moreover, for given $(q,g,q)\in \bm{M}(\Psi,\Psi)$, $g\in G=G_\alpha$,
there exists an open neighborhood $U(q,g,q)$ characterized by the property that
$$
s:U(q,g,q)\rightarrow O
$$
is an sc-diffeomorphism. Further it holds that $U(q,g,q)\cap U(q,g',q)=\emptyset $ for $g\neq g'$
and 
$$
\bm{M}(\Psi,\Psi)=\bigcup_{g\in G} U(q,g,q).
$$
In addition it holds that $t:U(q,g,q)\rightarrow O$ is given by $t(p,\phi,p')= g\ast p$.
\end{lemma}
\begin{proof}
By assumption $\Psi: G\ltimes O\rightarrow {\mathcal S}$ is injective on objects and fully faithful.
If $g\in G$ and $p\in O$ then $\Psi(g,p):\Psi(p)\rightarrow \Psi(g\ast p)$ which implies
that $(p,\Psi(g,p),g\ast p)\in \bm{M}(\Psi,\Psi)$. Moreover the map  in (\ref{EQNX17.1})
is injective since $\Psi$ is injective on objects and $\Psi$ is faithful. 
Given $p, p'\in O$ and $(p,\phi,p')\in \bm{M}(\Psi,\Psi)$ the fact that $\Psi$
is full implies that there exists $(g,p)$ with $\Psi(g,p) = \phi$ which also implies that $g\ast p=p'$.
This shows that the map in (\ref{EQNX17.1}) is a bijection. 
From $s\circ \hat{\Psi}(g,p)= p$ we conclude that $\hat{\Psi}$ is a local sc-diffeomorphism, since $s$ has this property.
Of course, with $\hat{\Psi}$ being a bijection we conclude that it is a sc-diffeomorphism.
Since $\{g\}\times O$ is open in $G\times O$ we see that the map 
$$
O\rightarrow \bm{M}(\Psi,\Psi):p\rightarrow (p,\Psi(g,p),g\ast p)
$$
is an sc-diffeomorphism onto some open subset. which we shall denote by $U(q,g,q)$.
We note that $(q,g,q)\in U(q,g,q)$.  It is clear that by construction $U(q,g,q)$ and $U(q,g',q)$ are disjoint if $g\neq g'$ and further that
$$
\bm{M}(\Psi,\Psi)=\bigcup_{g\in G} U(q,g,q).
$$
Clearly $s:U(q,g,q)\rightarrow O$ has the form $s(p,\phi,p')=p$. From
$$
\hat{\Psi}(g,p)=(p,\Psi(g,p),g\ast p)
$$
 it follows that $t\circ \hat{\Psi}(g,p)= g\ast p$ implying that if 
$(p,\phi,p')\in U(q,g,q)$ we must have $p'=g\ast p$ and consequently $t(p,\phi,p')=g\ast p$.
\qed \end{proof}

Assume we are given a polyfold structure $(F,{\bm{M}})$ for $\mathscr{C}$. 
We can connect this structure to the theory of ep-groupoids as follows.
Using that $|\mathscr{C}|$ is a set and ${\mathcal T}$ a metrizable topology we can pick a set  of uniformizers $\bm{\Psi}:={(\Psi_\lambda)}_{\lambda\in\Lambda}$ covering the orbit set of $\mathscr{C}$, i.e. 
$$
|\mathscr{C}| =\bigcup_{\lambda\in\Lambda} \text{footprint}(\Psi_\lambda).
$$
Recall that the footprints are open subsets.
We shall refer to $\bm{\Psi}$ as a {\bf covering set of uniformizers}\index{Covering set of uniformizers}, or just simply 
as a {\bf covering set}.
By definition $\Psi_\lambda : G_\lambda\ltimes O_\lambda\rightarrow \mathscr{C}$, and 
we associate to this family the M-polyfolds
$$
X=\bigsqcup_{\lambda\in\Lambda} O_\lambda\ \ \text{and}\ \ {\bm{X}}=\bigsqcup_{(\lambda,\lambda')\in \Lambda\times\Lambda} {\bm{M}}(\Psi_\lambda,\Psi_{\lambda'}).
$$
We  define source and target maps $s,t:{\bm{X}}\rightarrow X$ by associating 
to $(q_\lambda,\phi,q_{\lambda'})$ the objects
$$
s(q_\lambda,\phi,q_{\lambda'})=q_\lambda\ \ \text{and}\ \ t(q_\lambda,\phi,q_{\lambda'})=q_{\lambda'}.
$$
These maps are local sc-diffeomorphisms, since $(F,{\bm{M}})$ defines a polyfold structure.
It follows immediately that the global inversion map
$\iota:{\bm{X}}\rightarrow {\bm{X}}:(q_\lambda,\phi,q_{\lambda'})\rightarrow (q_{\lambda'},\phi^{-1},q_{\lambda})$ 
is an sc-diffeomorphism and the unit map $u:X\rightarrow {\bm{X}}$ is sc-smooth. Since $s$ and $t$ are
local sc-diffeomorphisms,  ${\bm{X}}{_{s}\times_t}{\bm{X}}$ has a natural M-polyfold structure 
and it follows that the multiplication map is sc-smooth. In summary the family $\bm{\Psi}$
defines a small category $X=X_{\bm{\Psi}}$ where the object and morphism set carry M-polyfold structures
for which the standard structure maps are sc-smooth and $s$ and $t$ are local sc-diffeomorphisms.
However, more is true.

\begin{theorem}\index{T- Category associated to $\bm{\Psi}$}
Let $(F,\bm{M})$ be a polyfold construction for the GCT $\mathscr{C}$. Then the  small category $X_{\bm{\Psi}}$ associated to the covering set of uniformizers $\bm{\Psi}$
has the structure of an ep-groupoid. If $(F,\bm{M})$ is tame then $X_{\bm{\Psi}}$ is a tame ep-groupoid.
\end{theorem}
\begin{proof}
It suffices to verify the properness property.
We find an open neighborhood $V=V(q_\lambda)\subset O_\lambda$ so that
$V$ admits the natural action of the stabilizer group $G_{q_\lambda}$ of $g_\lambda$
and $|\Psi_\lambda| : {_{G_{q_\lambda}}\backslash} V\rightarrow |\Psi_\lambda(V)|$  is a homeomorphism.
Since $|\mathscr{C}|$ is metrizable it is in particular normal and we find an open neighborhood $W$ of $|\Psi_\lambda(q_\lambda)|$ with $\cl_{|\mathscr{C}|}(W)\subset |\Psi_\lambda(V)|$. Define $U\subset O_\lambda$ to consist
of all $p_\lambda$ satisfying
$$
|\Psi_\lambda(p_\lambda)|\in W.
$$
Then $U$ is an open neighborhood of $q_\lambda$ in $O_\lambda$ and we shall show that
it has the desired properties.  Assume that $(q_k,\phi_k,p_k)$ is a sequence of elements in ${\bm{X}}$
such that $q_k\in \cl_X(U)$ and $(p_k)$ belongs to a compact subset $K$ in $X$.
After perhaps taking a subsequence we may assume that $p_k\in O_{\lambda'}$ and $p_k\rightarrow p_0$.
Since $|\Psi_{\lambda}(q_k)|=|\Psi_{\lambda'}(p_k)|\rightarrow |\Psi_{\lambda'}(p_0)|$ it follows
that $|\Psi_{\lambda'}(p_0)|\in \cl_{|\mathscr{C}|}(W)\subset V$. Using that
$$
|\Psi_\lambda| : {_{G_{q_\lambda}}\backslash} V\rightarrow |\Psi_\lambda(V)|
$$
is a homeomorphism we conclude that a subsequence of $(q_k)$ is convergent.
\qed \end{proof} 
 In view of the theorem we can associate to a covering set of uniformizers $\bm{\Psi}:={(\Psi_\lambda)}_{\lambda\in\Lambda}$
 an ep-groupoid $X_{\bm{\Psi}}$. In addition, we can define a functor
 $$
 \Gamma_{\bm{\Psi}}:X_{\bm{\Psi}}\rightarrow \mathscr{C}
 $$
 by associating to $q_\lambda\in O_\lambda$ the object $\Psi_\lambda(q_\lambda)$ and to a morphism
 $(q_\lambda,\phi,q_{\lambda'})$ the morphism $\phi$.

 \begin{lemma}
 The functor $\Gamma_{\bm{\Psi}}:X_{\bm{\Psi}}\rightarrow \mathscr{C}$ is an equivalence of categories
 and $|\Gamma_{\bm{\Psi}}|:|X_{\bm{\Psi}}|\rightarrow |\mathscr{C}|$ is a homeomorphism.
 \end{lemma}
 \begin{proof}
 Since the footprints cover $|\mathscr{C}|$ it follows that $\Gamma_{\bm{\Psi}}$ is {\bf essentially surjective}\index{Essentially surjective}, i.e.
 for every object $\alpha$ in $\mathscr{C}$ there exists $x\in X_{\bm{\Psi}}$ and $\phi\in \bm{C}$ satisfying
 $$
 \Gamma_{\bm{\Psi}}(x)\xrightarrow{\phi} \alpha.
 $$
 Given two objects $q_\lambda, q_{\lambda'}$ the map
 $$
 \Gamma_{\bm{\Psi}}:{\bm{X}}(q_\lambda,q_{\lambda'})\rightarrow {\bf C}(\Psi_\lambda(q_\lambda),\Psi_{\lambda'}(q_{\lambda'}))
 $$
 is trivially a bijection.  This shows that  $\Gamma_{\bm{\Psi}}$ is an equivalence of categories.
 From this it follows immediately that $| \Gamma_{\bm{\Psi}}|$ is a bijection.
 Since the $\Psi_\lambda$ are uniformizers we infer that $|\Gamma_{\bm{\Psi}}|$  is a homeomorphism.
\qed \end{proof}

Given a polyfold structure $(F,\bm{M})$ for the GCT $\mathscr{C}$ the above discussion shows that we can construct
an sc-smooth small category  $X_{\bm{\Psi}}$ with some additional features, i.e. an ep-groupoid, and an equivalence of categories 
to $\mathscr{C}$. The construction of $X_{\bm{\Psi}}$ involved choices. However, as we shall see next, the results 
of two different choices are (sc-smoothly) Morita-equivalent in a canonical way, i.e. there exists a canonical generalized isomorphism
 $$
 \mathfrak{f}:X_{\bm{\Psi}}\rightarrow X_{\bm{\Psi}'}
 $$
 compatible, in some sense, with the equivalences $\Gamma_{\bm{\Psi}}$ and $\Gamma_{\bm{\Psi}'}$.
In order to see this assume that $\bm{\Psi}$ and $\bm{\Psi}'$ are two covering sets of uniformizers associated
to $(F,{\bm{M}})$. Then the union $\bm{\Psi}''$
also has footprints covering the orbit space. We obtain natural inclusion functors
$$
X_{\bm{\Psi}}\xrightarrow{A} X_{\bm{\Psi}''}\ \ \text{and}\ \ X_{\bm{\Psi}'}\xrightarrow{A'}X_{\bm{\Psi}''}.
$$
One readily verifies that $A$ and $A'$ are sc-smooth equivalences between ep-groupoids. We also note that
$$
\Gamma_{\bm{\Psi}}= \Gamma_{\bm{\Psi}''}\circ A\ \ \text{and}\ \ 
\Gamma_{\bm{\Psi}'}= \Gamma_{\bm{\Psi}''}\circ A'.
$$
We take  the weak fibered product $X_{\bm{\Psi}}\times_{X_{\bm{\Psi}''}}X_{\bm{\Psi}'}$ associated to the diagram 
$$
X_{\bm{\Psi}}\xrightarrow{A} X_{\bm{\Psi}''}\xleftarrow{A'} X_{\bm{\Psi}'}.
$$
This is an ep-groupoid 
and the projections onto the factors give equivalences of ep-groupoids resulting in the diagram
\begin{eqnarray}\label{polq1}
d\colon X_{\bm{\Psi}}\xleftarrow{\pi_1}   X_{\bm{\Psi}}\times_{X_{\bm{\Psi}''}}X_{\bm{\Psi}'}\xrightarrow{\pi_2}                      X_{\bm{\Psi}'}.
\end{eqnarray}
This diagram $d$ defines a generalized isomorphism $[d]:X_{\bm{\Psi}}\rightarrow X_{\bm{\Psi}'}$.
One also verifies easily that $|\Gamma_{\bm{\Psi}'}|\circ |[d]| =|\Gamma_{\bm{\Psi}}|$ and we can summarize the discussion as follows,
where we appeal to the definition of a polyfold structure on a topological  space as given in Definition \ref{SECDEF121}.
\begin{theorem}\label{THMX17110}\index{T- Morita equivalence}
The construction associating to a covering set of uniformizers $\bm{\Psi}$ 
the ep-groupoid $X_{\bm{\Psi}}$ produces for different families $\bm{\Psi}$ and $\bm{\Psi}'$ 
ep-groupoids and a natural generalized isomorphism 
$$
\mathfrak{f}\colon X_{\bm{\Psi}}\rightarrow X_{\bm{\Psi}'},
$$
which is the identity if $\bm{\Psi}=\bm{\Psi}'$.
If $\mathfrak{f}$ is associated to $(\bm{\Psi},\bm{\Psi}')$ and $\mathfrak{f}'$ to $(\bm{\Psi}',\bm{\Psi}'')$, then 
$\mathfrak{f}'\circ\mathfrak{f}$ is associated to $(\bm{\Psi},\bm{\Psi}'')$.
The generalized isomorphism $\mathfrak{f}:X_{\bm{\Psi}}\rightarrow X_{\bm{\Psi}'}$ associated to (\ref{polq1})
satisfies $|\Gamma_{\bm{\Psi}'}|\circ |\mathfrak{f}| =|\Gamma_{\bm{\Psi}}|$.
Consequently the pairs $(X_{\bm{\Psi}},|\Gamma_{\bm{\Psi}}|)$ define
equivalent polyfold structures on the metrizable space $|\mathscr{C}|$ and consequently $|\mathscr{C}|$ is naturally a paracompact polyfold.
\qed
\end{theorem}
In view of this theorem we can associate to a polyfold  structure for the GCT $\mathscr{C}$
a polyfold structure for the underlying metrizable space $|\mathscr{C}|$. With other words we have a forgetful functor
$$
(\text{GCT-polyfold}\ \  \mathscr{C})\ \ \ \rightsquigarrow \ \ (\text{Paracompact topological polyfold}\  |\mathscr{C}|).
$$
This procedure
loses the algebraic knowledge of the underlying category structure, but keeps
some information about it in the form of the Morita equivalence class of ep-groupoids.
Theorem \ref{THMX17110} is important for constructions in the categorical context since it allows to transfer
any notion which behaves well with respect to generalized isomorphisms. Since the arguments are usually straight forward
we allow ourselves to be sketchy at times.
\begin{remark}\index{R- Unions of covering families}
Assume that $X$ and $X'$ are ep-groupoids constructed from covering families $\bm{\Psi}$ and $\bm{\Psi}'$.
We denote by $X''$ the ep-groupoid associated to the union of the two covering families.
The weak fibered product $X\times_{X''} X'$ has as objects the tuples $(q_\lambda,(q_\lambda,\phi,q_{\lambda'}'),q_{\lambda'}')$,
where $q_\lambda$ is an object in $O_\lambda$, $q_{\lambda'}'\in O_{\lambda'}'$, and $(q_\lambda,\phi,q_{\lambda'}')$ is a morphism
in $X''$  between $q_\lambda$ and $q_{\lambda'}'$.
We define  
$$
\tau: X\times_{X''} X'\rightarrow \bm{C}: (q_{\lambda},(q_{\lambda},\phi,q_{\lambda'}'),q_{\lambda'}')\rightarrow \phi
$$
which associates to an object in $X\times_{X''} X'$ a morphism in $\mathscr{C}$. It turns out that $\tau$ 
is a natural transformation between $\Gamma\circ \pi_1$ and $\Gamma'\circ\pi_2$ which define equivalences 
$$
X\times_{X''} X'\rightarrow \mathscr{C}.
$$
\qed
\end{remark}
In the following  we shall carry over parts of the discussion of ep-groupoids and the concepts
which are compatible with generalized isomorphisms.
With $\mathscr{C}$ being equipped with a polyfold structure we obtain a filtration 
$\mathscr{C}_i$ for $i\in {\mathbb N}$. Namely we can pick for an object $\alpha$
a uniformizer $\Psi\in F(\alpha)$,  and writing $\alpha=\Psi(q)$ with $q\in O$ we can say
that $\alpha$ has {\bf regularity}\index{Regularity of an object} $i\in {\mathbb N}\cup \{\infty\}$ provided $q\in O_i$.
Of course, if $i<k$ and $q\in O_k$ it also has regularity $i$. Usually, for constructions it only matters that an object has a sufficient amount of regularity, or is smooth, i.e. belongs to all $O_i$. We can define $\text{reg}(\alpha)\in {\mathbb N}\cup\{+\infty\}$ as the maximal $i$ (including $\infty$)
such that $\alpha$ has this regularity. This all is well-defined, independent of the choice of $\Psi\in F(\alpha)$.  If $\alpha$ and $\alpha'$
are isomorphic they have the same regularity. This follows from the fact that $s$ and $t$ are local sc-diffeomorphisms.
As a consequence we obtain a {\bf regularity filtration}\index{Regularity filtration of $\mathscr{C}$}
$$
\mathscr{C}_\infty\rightarrow...\rightarrow \mathscr{C}_i\rightarrow \mathscr{C}_{i-1}...\rightarrow \mathscr{C}_0=\mathscr{C},
$$
which also descends to isomorphism classes.
As in the case of M-polyfolds or ep-groupoids the full subcategory $\mathscr{C}_i$ for a  $i\in {\mathbb N}$
has a natural polyfold structure as well. To see this we start with the equivalence 
$$
\Gamma_{\bm{\Psi}}:X_{\bm{\Psi}}\rightarrow \mathscr{C},
$$
which induces the homeomorphism $|\Gamma_{\bm{\Psi}}|\colon |X_{\bm{\Psi}}|\rightarrow |\mathscr{C}|$.
As we already have seen $|X|$ is metrizable and in particular paracompact. This implies that $|X^i|$ is paracompact and therefore metrizable
and consequently $|\mathscr{C}_i|$ has a natural metrizable topology, see Section \ref{section1.3_top_prop} for the relevant discussion.
In particular $\mathscr{C}_1$ is a GCT.
Define 
$$
F^1:\mathscr{C}_1\rightarrow \text{SET}
$$
by associating to $\alpha$ the collection $F^1(\alpha)$ of all $\Psi^1:G\ltimes O^1\rightarrow \mathscr{C}_1$.
Here $\Psi^1$ is obtained from $\Psi$ by lifting the index by $1$. We define for $\Psi^1$ and $\Psi^{'1}$
the M-polyfold ${\bm{M}}_{F^1}(\Psi^1,\Psi^{'1})$ as follows
$$
{\bm{M}}_{F^1}(\Psi^1,\Psi^{'1}):=({\bm{M}}(\Psi,\Psi'))^1,
$$
with the natural identification,  observing that the right-hand side has a natural M-polyfold structure.
The following result is obvious.
\begin{theorem}\label{THMX17111} \index{T- Natural polyfold structure for $\mathscr{C}_1$} 
Given the polyfold structure $(F,{\bm{M}})$ for $\mathscr{C}$ the full subcategory $\mathscr{C}_1$ is in a natural way a GCT
and, moreover, has the natural polyfold structure $(F^1,{\bm{M}}_{F^1})$. 
\qed
\end{theorem}
If $\mathscr{C}$ is equipped with a polyfold structure we shall write $\mathscr{C}^1$ for the category
$\mathscr{C}_1$ equipped with the induced natural polyfold structure. We also define
$$
\mathscr{C}^{i+1}:={(\mathscr{C}^i)}^1.
$$ 
\begin{definition}
Assume $\mathscr{C}$ is equipped with a polyfold structure. The degeneracy functor
$$
d:\mathscr{C}\rightarrow {\mathbb N}
$$
is defined by $d(\alpha):=d_O(q)$, where $\Psi\in F(\alpha)$ and $\Psi(q)=\alpha$. 
Here $d_O$ is the degeneracy index on $O$, where $\Psi:G\ltimes O\rightarrow \mathscr{C}$.
The {\bf boundary} $\partial\mathscr{C}$  of $\mathscr{C}$ is the full subcategory associated to objects of degeneracy at least $1$.
\qed
\end{definition}

\section{Tangent Construction}\label{TANGENTX172}
In this section we shall  introduce the tangent construction, which associates to a GCT $\mathscr{C}$ a  $GTC$ $T\mathscr{C}$ together with a projection
functor
$$
P\colon T\mathscr{C}\rightarrow \mathscr{C}^1
$$
and equips $T\mathscr{C}\rightarrow \mathscr{C}^1$ with the structure of a so-called  polyfold bundle, which induces the already given structure on $\mathscr{C}^1$. A polyfold bundle is a weaker notion than that of a strong bundle in the categorical context,  introduced later in Section \ref{CHAP176-}.

Given an object $\alpha$ in $\mathscr{C}_1$ we pick $\Psi, \Psi'\in F(\alpha)$, written as 
$$
\Psi\colon G\ltimes O\rightarrow \mathscr{C}\ \ \text{and}\ \ \Psi'\colon  G\ltimes O'\rightarrow \mathscr{C},
$$
and consider tuples $(\alpha,\Psi,h)$, with $h\in T_qO$ and   $\Psi(q)=\alpha$, and similarly 
$(\alpha,\Psi',h')$, with  $h'\in T_{q'}O'$, where by  $q'\in O'$ we denote the point satisfying $\Psi'(q')=\alpha$.
We shall define a notion of equivalence for two such tuples $(\alpha,\Psi,h)$ and $(\alpha,\Psi',h')$.

In order to do so we fix the tuple $(q,1_\alpha,q')\in \bm{M}(\Psi,\Psi')$,  and using that the source and target maps $s$ and $t$ are local sc-diffeomorphisms 
we find open neighborhoods $U(q)$, $U(q')$ and $U(q,1_\alpha,q')$ such that 
$$
s\colon U(q,1_\alpha,q')\rightarrow U(q)\ \ \text{and}\ \ t\colon U(q,1_\alpha,q')\rightarrow U(q')
$$
are sc-diffeomorphisms defining an sc-diffeomorphism $\sigma:U(q)\rightarrow U(q')$ by $\sigma(p) = t\circ (s|U(q,1_\alpha,q'))^{-1}(p)$.
We declare $(\alpha,\Psi,h)$ and $(\alpha,\Psi',h')$ to be {\bf equivalent}, i.e. 
$$
(\alpha,\Psi,h)\sim (\alpha,\Psi',h'),
$$
 provided 
$$
T\sigma(q)h=h'.
$$
\begin{lemma}\index{L- Equivalence relation $\sim$}
 ``$\sim$" defines an equivalence relation on the tuples $(\alpha,\Psi,h)$.
 \end{lemma}
\begin{proof}
Clearly $(\alpha,\Psi,h)\sim (\alpha,\Psi,h)$ since for $(q,1_\alpha,q)$ the local sc-diffeomor\-phism $\sigma$ is the identity.
The symmetry of the relation, i.e. $(\alpha,\Psi,h)\sim (\alpha,\Psi',h')$ implies $(\alpha,\Psi',h')\sim(\alpha,\Psi,h)$,
follows from the fact that the local sc-diffeomorphisms associated to $(q,1_\alpha,q')$ 
and $(q',1_\alpha,q)$ are inverse to each other.  For the transitivity property consider $(\alpha,\Psi,h)\sim (\alpha,\Psi',h')$
and $(\alpha,\Psi',h')\sim (\alpha,\Psi'',h'')$. Then with $\sigma$ corresponding to $(q,1_\alpha,q')$ and $\sigma'$ to $(q',1_\alpha,q'')$ 
we have 
$$
h'= T\sigma(q)(h)\ \ \text{and}\ \ h''=T\sigma(q')(h').
$$
This implies $h'' = T\sigma(q')\circ T\sigma(q)(h)$. We note that for $p$ near $q$ it holds that 
\begin{eqnarray*}
\sigma'\circ \sigma (p) &=& t\circ (s|U(q',1_\alpha,q''))^{-1} \circ t\circ  (s|U(q,1_\alpha,q'))^{-1}(p)\\
&=& t\circ (s|U(q,1_\alpha,q''))^{-1}(p)\\
&=&\sigma''(p)
\end{eqnarray*}
which implies that $h'' = T\sigma''(q)(h)$ and proves our assertion.  
\qed \end{proof}
We shall write $[\alpha,\Psi,h]$\index{$[\alpha,\Psi,h]$} for the equivalence class containing $(\alpha,\Psi,h)$.
\begin{definition}\index{D- Tangent vector at $\alpha$}
Let $\mathscr{C}$ be a GCT equipped with a polyfold structure. Given an object $\alpha$ in $\mathscr{C}_1$
we call an equivalence class $[\alpha,\Psi,h]$, where $\Psi\in F(\alpha)$, say $\Psi:G\ltimes O\rightarrow \mathscr{C}$,
and $h\in T_qO$ with $\Psi(q)=\alpha$ a 
 {\bf tangent vector}
at $\alpha$.
\qed
\end{definition}
 Note that $\alpha$ needs to be of regularity one in order to be able to talk about tangent vectors.
\begin{definition}\index{D- Tangent space $T_\alpha\mathscr{C}$}
Given a GCT $\mathscr{C}$ equipped with a polyfold structure,  we define for an object $\alpha$ of regularity $1$ the {\bf tangent space} $T_\alpha\mathscr{C}$
as the set of equivalence classes $[\alpha,\Psi,h]$.
\qed
\end{definition}
 The set $T_\alpha\mathscr{C}$ carries the structure  of a real vector space by setting
$$
[\alpha,\Psi,h]+\lambda[\alpha,\Psi,k]=[\alpha,\Psi,h+\lambda k].
$$
We can equip $T_\alpha\mathscr{C}$  with the structure of a Banach space by requiring 
that
the linear bijection
$$
S_\Psi\colon T_\alpha\mathscr{C}\rightarrow T_qO\colon [\alpha,\Psi,h]\rightarrow h
$$
is a topological linear isomorphism. The Banach space structure on $T_\alpha\mathscr{C}$ does not depend on the choice 
of $\Psi\in F(\alpha)$, since  $S_{\Psi'}\circ S_{\Psi}^{-1}(h)=T\sigma(q)(h)$ is a topological linear isomorphism.
We also note that for smooth objects $T_\alpha\mathscr{C}$   is naturally a sc-Banach space. 
Hence we obtain the following result.
\begin{proposition}\index{P- Banach space $T_\alpha\mathscr{C}$}
Given a GCT $\mathscr{C}$ equipped with a polyfold structure $(F,\bm{M})$, there exists for every $\alpha$ in $\mathscr{C}_1$ a natural associated Banach space $T_\alpha\mathscr{C}$
called the tangent space at $\alpha$. In case that $\alpha$ is smooth $T_\alpha\mathscr{C}$ has a natural sc-structure.
\qed
\end{proposition}
At this point we have a construction which, starting with a GTC $\mathscr{C}$ equipped with a polyfold structure,
associates to every objects $\alpha\in \mathscr{C}_1$  a Banach space $T_\alpha\mathscr{C}$. We shall show that 
this extends to a functor $\mathscr{C}_1\rightarrow \text{Ban}$, where $\text{Ban}$ \index{$\text{Ban}$} is the {\bf category of Banach spaces}\index{Category of Banach spaces: Ban} with the morphisms being topological linear isomorphisms.

For this assume  that $\phi:\alpha\rightarrow \alpha'$ is an isomorphism belonging to level $1$ and let $[\alpha,\Psi,h]$
be an element in $T_\alpha\mathscr{C}$.
We consider the transition M-polyfold $\bm{M}(\Psi,\Psi')$ and pick $q_0\in O$ satisfying $\Psi(q_0)=\alpha$ and $q_0'$
with $\Psi'(q_0')=\alpha'$. Take open neighborhoods $U(q_0)$, $U(q_0,\phi,q_0')$, and $U(q_0')$ such that $s:U(q_0,\phi,q_0')\rightarrow U(q_0)$ and $t:U(q_0,\phi,q_0')\rightarrow U(q_0')$ are sc-diffeomorphisms. Then $\sigma:U(q_0)\rightarrow U(q_0')$
defined by 
$$
\sigma(q) =t\circ (s|U(q_0,\phi,q_0'))^{-1}(q)
$$
is an sc-diffeomorphism satisfying $\sigma(q_0)=q_0'$. We define $T\phi:T_\alpha\mathscr{C}\rightarrow T_{\alpha'}\mathscr{C}$ 
\index{$T\phi:T_\alpha\mathscr{C}\rightarrow T_{\alpha'}\mathscr{C}$} by
$$
T\phi([\alpha,\Psi,h])=[\alpha',\Psi',T\sigma(q_0)(h)].
$$
Assume next that $\psi:\alpha'\rightarrow \alpha''$ and pick $\Psi''\in F(\alpha'')$. Then, with the obvious notation
\begin{eqnarray*}
&&T\psi(T\phi([\alpha,\Psi,h]))=T\psi([\alpha',\Psi',T\sigma(q_0)(h)])\\
& =&[\alpha'',\Psi'',T\sigma'(q_0')(T\sigma(q_0)(h))]=[\alpha'',\Psi'',T\sigma''(q_0)(h)].
\end{eqnarray*}
Summarizing, we obtain the following result.
\begin{proposition}\index{P- The functor $\mu_{\mathscr{C}}$}  \label{PROPU1725}
Given a GCT $\mathscr{C}$ equipped with a polyfold structure $(F,\bm{M})$ there exists a well-defined functor 
$$
\mu_{\mathscr{C}}: \mathscr{C}_1\rightarrow \text{Ban},
$$
 which associates to an object $\alpha$ on level $1$ the Banach space $T_\alpha\mathscr{C}$ and to a morphism
 $\phi:\alpha\rightarrow \alpha'$ the linear isomorphism $T\phi$.
 \qed
 \end{proposition}
 As we shall see next the tangent spaces $T_\alpha\mathscr{C}$ fit in some sense sc-smoothly together.
 The category $T\mathscr{C}$\index{$T\mathscr{C}$}, which we are going to define, will have as objects all the equivalence classes $[\alpha,\Psi,h]$,
 where $\alpha$ varies over $\mathscr{C}_1$. We shall use the functor $\mu_{\mathscr{C}}$ to carry out the construction.
 As obvious data we have the   projection functor $P:T\mathscr{C}\rightarrow \mathscr{C}_1$,  which for the moment
 is only defined on the object level by
\begin{eqnarray}\label{EQN172-}
P([\alpha,\Psi,h])=\alpha.
\end{eqnarray}
As morphisms we take the tuples $([\alpha,\Psi,h],\phi,T\phi([\alpha,\Psi,h]))$, where 
$\alpha$ is an object in $\mathscr{C}_1$, $[\alpha,\Psi,h]\in T_{\alpha}\mathscr{C}$ and $\phi$ is a morphism in
$\mathscr{C}$ with $s(\phi)=\alpha$ and $t(\phi)=\alpha'$. We view the morphism as
$$
([\alpha,\Psi,h],\phi,T\phi([\alpha,\Psi,h])): [\alpha,\Psi,h]\rightarrow T\phi([\alpha,\Psi,h]),
$$
so that 
\begin{eqnarray*}
& s([\alpha,\Psi,h],\phi,T\phi([\alpha,\Psi,h]))=[\alpha,\Psi,h]& \\
& t([\alpha,\Psi,h],\phi,T\phi([\alpha,\Psi,h]))=T\phi([\alpha,\Psi,h]).&
\end{eqnarray*}
For convenience of notation we shall often write a tangent vector $[\alpha,\Psi,h]$ as $\bar{h}$ and a morphism
$([\alpha,\Psi,h],\phi,[\alpha',\Psi',h'])$ as $\Phi$ or $\Phi:\bar{h}\rightarrow \bar{h}'$, or $(\bar{h},\phi,\bar{h}')$, or $(\bar{h},\phi,T\phi(\bar{h}))$,
which depends on the situation. 
We define $P$ on morphisms as
$$
P(\bar{h},\phi,\bar{h}')=\phi,
$$
which together with (\ref{EQN172-}) defines a functor 
\begin{eqnarray}
P:T\mathscr{C}\rightarrow \mathscr{C}_1.
\end{eqnarray}
Consider any object $\alpha$ belonging to $\mathscr{C}_1$ and pick $\Psi\in F(\alpha)$ which we write as 
$$
\Psi:G\ltimes O\rightarrow \mathscr{C}.
$$
We define with $\tau:TO\rightarrow O^1$
being the sc-smooth bundle projection the tangent $T\Psi$, say 
\begin{eqnarray}
T\Psi: G\ltimes TO \rightarrow T\mathscr{C}
\end{eqnarray}
by 
$$
T\Psi(h) = [\Psi(\tau(h)), \Psi', T\sigma(h)].
$$
 Here $\Psi'\in F(\tau(h))$ and $\sigma$ is the local sc-diffeomorphism 
associated to the element $(\tau(h),1_{\Psi(\tau(h))},q')\in \bm{M}(\Psi,\Psi')$, where $\Psi'(q')=\Psi(\tau(h))$. The definition involves the choice 
of a $\Psi'$ for every $p= \tau(h)$, but as we shall see  it is independent of such a choice since it is compensated for by $\sigma$, which changes 
with $\Psi'$.  Given a morphism $(g,h):h\rightarrow g\ast h$ in $G\ltimes TO$, say $h\in T_qO$,  we define $T\Psi(h)\rightarrow T(\Psi(g,q))(T\Psi(h))$
by
$$
T\Psi(g,h) = (T\Psi(h),\Psi(g,q),T(\Psi(g,q))(T\Psi(h))).
$$
Note that the morphism $\Psi(g,q):\Psi(q)\rightarrow \Psi(g\ast q)$ for $q\in O_1$ has the associated tangent operator
$$
T(\Psi(g,q)):T_{\Psi(q)}\mathscr{C}\rightarrow T_{\Psi(g\ast q)}\mathscr{C},
$$
introduced in Proposition \ref{PROPU1725}.
\begin{lemma}
$T\Psi: G\ltimes TO \rightarrow T\mathscr{C}$ is well-defined.
\end{lemma}
\begin{proof}
Assume $p\in O^1$ is fixed and pick $\Psi',\Psi''\in F(\Psi(p))$. With $G'$ being the isotropy group of $\Psi(p)$
and writing 
 $$
 \Psi':G'\ltimes O'\rightarrow \mathscr{C}\ \ \text{and} \ \ \Psi'':G'\ltimes O''\rightarrow \mathscr{C}
 $$
consider the points  $q'\in O'$ and $q''\in O''$
satisfying  $\Psi'(q')=\Psi''(q'')=\Psi(p)$. The local sc-diffeomorphism $\sigma$ with $\sigma(p)=q'$ is associated by the usual construction
to $(p,1_{\Psi(p)},q')\in \bm{M}(\Psi,\Psi')$ and $\sigma'$ with $\sigma'(p)=q''$ to $(p,1_{\Psi(p)},q'')\in \bm{M}(\Psi,\Psi'')$.
Consider with $\beta=\Psi(p)$ the tuples 
$$
(\beta,\Psi',T\sigma(p)(h))\ \ \text{and}\ \ (\beta,\Psi'',T\sigma'(p)(h)).
$$
 These tuples are equivalent 
since with $\sigma''$ associated to $(q',1_\beta,q'')\in \bm{M}(\Psi',\Psi'')$ we have that 
$$
\sigma''\circ \sigma =\sigma' \ \ \text{near}\ \ p.
$$
Hence $T\sigma''(q') T\sigma(p) = T\sigma'(p)$, which implies that 
$$
[\beta,\Psi',T\sigma(p)(h)] =[\beta,\Psi'', T\sigma'(p)(h)].
$$
Therefore $T\Psi$ is well-defined.
\qed \end{proof}
Assume that $\Psi_1\in F(\alpha)$ and $\Psi_2\in F(\alpha')$, which can be written
as
$$
\Psi_1:G\ltimes O\rightarrow \mathscr{C}\ \ \text{and}\ \ \Psi_2:G'\ltimes O'\rightarrow \mathscr{C}
$$
are given.
For $(o,\phi,o')\in \bm{M}(\Psi_1,\Psi_2)$ we can pick open neighborhoods
$$
U(o),\ U(o,\phi,o')\ \ \text{and}\ \ U(o')
$$
 such that  
$$
s:U(o,\phi,o')\rightarrow U(o)\ \ \text{and}\ \ t:U(o,\phi,o')\rightarrow U(o')
$$
are sc-diffeomorphisms. There exists a well-defined map $U(o)\ni q\rightarrow \phi_q$ such that
$$
(s|U(o,\phi,o'))^{-1}(q) = (q,\phi_q,\sigma(q)),
$$
where $\sigma(q) =t\circ (s|U(o,\phi,o'))^{-1}(q)$.  The existence of this well-defined map is, of course,
a consequence of having a polyfold structure on $\mathscr{C}$. Hence we 
have the assignment 
$$
q\rightsquigarrow \left[\Psi_1(q)\xrightarrow{\phi_q} \Psi_2(\sigma(q))\right]\ \ \ \text{for}\ \ q\in U(o).
$$
Using this assignment we lift it to the tangent level by defining 
\begin{eqnarray}\label{EQNXX67}
k\rightsquigarrow\left[T\Psi_1(k)\xrightarrow{T\phi_q} T\Psi_2(T\sigma(k))\right]\ \ \ \text{for}\ \  k\in T_qU(o),\ q\in U(o)^1.
\end{eqnarray}
We note that all the ingredients are well-defined and the diagrams, which depend on $q$, make sense.  
We record (\ref{EQNXX67}) as having the family of morphisms in $T\mathscr{C}$ given as
$$
k\rightarrow\Phi_k:=(T\Psi_1(k),\phi_{\tau(k)},T\Psi_2(T\sigma(k)))\ \ \text{for}\ \ k\in TU(o).
$$
The  tuples
$(k,\Phi_k,T\sigma(k))$ belong to $\bm{M}(T\Psi_1,T\Psi_2)$ and  we obtain the map
\begin{eqnarray}\label{LEMMX1727}
TU(o)\rightarrow \bm{M}(T\Psi_1,T\Psi_2): k\rightarrow (k,\Phi_k,T\sigma(k)).
\end{eqnarray}

\begin{lemma}
The following holds.
\begin{itemize}
\item[{\em (1)}] \  The  construction of a map as in (\ref{LEMMX1727}) can be done for all $(o,\phi,o')\in \bm{M}(\Psi_1,\Psi_2)$ on level $1$, and 
consequently every  element in $\bm{M}(T\Psi_1,T\Psi_2)$ lies in the image of such a map.
\item[{\em (2)}]\   Consider two such maps $TU(o)\rightarrow \bm{M}(T\Psi_1,T\Psi_2): k\rightarrow (k,\Phi_k,T\sigma(k))$
and $TU(o_1)\rightarrow \bm{M}(T\Psi_1,T\Psi_2): k\rightarrow (k,\Phi_k^1,T\sigma'(k))$.
If the images of these maps intersect nontrivially,  the transition map has the form $k\rightarrow k$.
\end{itemize}
\end{lemma}
\begin{proof}
(1) is obvious. In order to prove (2) assume the images intersect. In this case we find $k_0\in TU(o)\cap TU(o')$ such that 
 $$
 (k_0,\Phi_{{k_0}},T\sigma(k_0))=(k_0,\Phi^1_{k_0},T\sigma'(k_0)).
 $$
Underlying is the  identity $(q_0,\phi_{q_0},\sigma(q_0))=(q_0,\phi_{q_0}',\sigma'(q_0))$, where $q_0=\tau(k_0)$.
This implies that the data near $q_0$ coincides as well (using the properties of a polyfold structure on $\mathscr{C}$) and the definitions. Consequently  the transition map is the identity.
\qed \end{proof} 
As a consequence we can view the 
maps $k\rightarrow (k,\Phi_k, T\sigma(k))$ as inverses of charts which are sc-smoothly compatible.
There is a uniquely defined topology ${\mathcal T}$  on $\bm{M}(T\Psi_1,T\Psi_2)$ making these maps homeomorphisms
onto open sets. 
\begin{definition}\index{D- Natural topology on $\bm{M}(T\Psi_1,T\Psi_2)$}
Let $\Psi_i\in F(\alpha_i)$ for $i=1,2$. Then the topology ${\mathcal T}$ defined on $\bm{M}(\Psi_1,\Psi_2)$ is called the {\bf natural topology}.
\end{definition}
The natural topology has very good properties.
\begin{proposition}
The natural  topology ${\mathcal T}$ on $\bm{M}(T\Psi_1,T\Psi_2)$  is metrizable.
The source and target maps are local homeomorphisms. The inversion maps are  homeomorphisms
and the unit maps as well as the multiplication maps are continuous.
\end{proposition}
\begin{proof}
We first shall proof that ${\mathcal T}$ is metrizable.
Using the Nagata-Smirnov Theorem it suffices to show that ${\mathcal T}$ is Hausdorff, regular, locally metrizable,
and paracompact.\par

\noindent{\bf Hausdorff:}  For $j=1,2$ let $(k_j,\Phi_j,k_j')$ be two different points in
$\bm{M}(T\Psi_1,T\Psi_2)$. We find  two natural maps 
$$
TU(o_i)\rightarrow \bm{M}(T\Psi_1,T\Psi_2): k\rightarrow (k,\Phi_k^i,T\sigma_i(k))
$$
for $i=1,2$. If the images are disjoint we obtain immediately disjoint open neighborhoods.
If the images intersect the transition map is of the form $k\rightarrow k$ and we may assume
that the two natural maps coincide. Since $k_1\neq k_2$ we find disjoint open neighborhoods in 
$TU(o)$ and the images of these neighborhoods provided us with disjoint open neighborhoods.\par

\noindent{\bf Local Metrizibility:} Since $\bm{M}(T\Psi_1,T\Psi_2)$ is locally homeomorphic to a M-polyfold via the source map $s$  it is locally metrizable.\par

\noindent{\bf Regularity:}  
Take a point $(k,\Phi,k')$ and a closed subset $A$ of $\bm{M}(T\Psi_1,T\Psi_2)$ not containing this point.
We find an open neighborhood $U(k,\Phi,k')$ and an open neighborhood $U(k)$ so that 
$$
s:U(k,\Phi,k')\rightarrow U(k) \ \ \text{is an sc-diffeomorphism.}
$$
 Since $TO$ is metrizable we can find an open neighborhood $W(k)\subset TO$
with $ \cl_{TO}(W(k))\subset U(k)$ and we can take a continuous function $\beta:U(k)\rightarrow [0,1]$
which takes the value $0$ at $k$ and the value $1$ on $U\setminus \cl_{TO}(W)$. 
If $(h,\Phi',h')\in U(k,\Phi,k')$ we define $\bar{\beta}(h,\Phi',h')=\beta(h)$ and otherwise it takes the value $1$.
This map is continuous. \par

\noindent{\bf Paracompactness:}  Consider the local homeomorphism $s:\bm{M}(\Psi_1,\Psi_2)\rightarrow O$.
For every $q\in O$ there exist only a finite number of points lying above it. Hence we find finitely many open sets
$U_1,...,U_{k(q)}$ in $\bm{M}(\Psi_1,\Psi_2)$ and an open neighborhood $U(q)$ so that
$$
s:U_i\rightarrow U(q)
$$
is an sc-diffeomorphism. Since $O$ is metrizable and therefore paracompact we find a locally finite covering ${(K_\lambda)}_{\lambda\in\Lambda}$
of $O$ consisting of closed sets so that $K_\lambda\subset U(q_\lambda)$. Then we take the preimages
and obtain $\wt{K}_{\lambda,i}$, where $\lambda$ varies in $\Lambda$ and $i\in \{1,...,k(q_\lambda)\}$.
Then the $\wt{K}_{\lambda,i}$ are paracompact closed subsets of $\bm{M}(T\Psi_1,T\Psi_2)$.
The whole collection $(\wt{K}_{\lambda,i})$ is a covering of $\bm{M}(T\Psi_1,T\Psi_2)$ by closed sets.
The covering is also locally finite. Hence the Hausdorff topological space $\bm{M}(T\Psi_1,T\Psi_2)$
has been written as a locally finite covering of closed paracompact sets. This implies paracompactness of
$\bm{M}(T\Psi_1,T\Psi_2)$ since it is Hausdorff and regular.
\par

\noindent  The verification of the properties of the structure maps is straight forward. 
\qed \end{proof}
Having the topology in place we can equip $\bm{M}(T\Psi_1,T\Psi_2)$ with a M-polyfold structure.

\begin{proposition}
The metrizable spaces $\bm{M}(T\Psi,T\Psi')$  have natural M-polyfold structures for which 
the maps $s:\bm{M}(T\Psi,T\Psi')\rightarrow TO$ and $t:\bm{M}(T\Psi,T\Psi')\rightarrow TO'$
are local sc-diffeomorphisms, the inversion map 
$$
\bm{M}(T\Psi,T\Psi')\rightarrow \bm{M}(T\Psi',T\Psi)
$$
is an sc-diffeomorphism and the unit map and multiplication maps are sc-smooth.
\end{proposition}
\begin{proof}
We use the maps $TU(o)\rightarrow \bm{M}(T\Psi_1,T\Psi_2): k\rightarrow (k,\Phi_{\tau(k)}, T\sigma(k))$
as charts and since the transition maps are identities we obtain a M-polyfold structure. 
With this definition all source maps are sc-smooth and local sc-diffeomorphisms.  A target map $t$ in local coordinates
has the form
$$
k\rightarrow T\sigma(k)
$$
and this expression is sc-smooth and a local sc-diffeomorphism. The other maps have also local expressions 
which are sc-smooth. We leave the details to the reader.
\qed \end{proof}

Starting with a a GCT $\mathscr{C}$ equipped with a polyfold structure $(F,\bm{M})$ we described several constructions.
The first construction defines a groupoidal category $T\mathscr{C}$ together with a functor $P:T\mathscr{C}\rightarrow \mathscr{C}_1$
so that the fibers $P^{-1}(\alpha)$ are Banach spaces. We shall show soon that $|T\mathscr{C}|$ carries a natural metrizable topology.
A second construction associates to $\Psi\in F(\alpha)$ with $\alpha\in \mathscr{C}_1$, say $\Psi:G\ltimes O\rightarrow \mathscr{C}$ 
a functor 
$$
T\Psi: G\ltimes TO\rightarrow  T\mathscr{C}.
$$
A third construction associates to $\Psi$ and $\Psi'$ as just described a natural M-polyfold structure on $\bm{M}(T\Psi,T\Psi')$.
The canonical map 
$$
\bm{M}(T\Psi,T\Psi')\rightarrow \bm{M}(\Psi^1,\Psi'^1):(k,(T\Psi(k),\phi,T\Psi'(k')),k')\rightarrow(\tau(k),\phi,\tau'(k))
$$
is sc-smooth since in suitable local coordinates it takes the form 
$$
k\rightarrow \tau(k).
$$
This is a bundle situation for which we have a natural system of charts. 
Pick $(o,\phi,o')$ on level $1$. Then we find open neighborhoods $U(o,\phi,o')$ and $U(o)$ in $\bm{M}(\Psi,\Psi')$ and $O$
such that $s:U(o,\phi,o')\rightarrow U(o)$ is an sc-diffeomorphism. We lift the level by one and obtain the sc-diffeomorphism
$s:U(o,\phi,o')^1\rightarrow U(o)^1$ which we may take as (M-polyfold valued) chart.
Then we obtain the commutative diagram
$$
\begin{CD}
\bm{M}(T\Psi,T\Psi')|U(o,\phi,o')^1 @> s>>  TU(o)\\
@VVV @V \tau VV\\
U(0,\phi,o')^{1} @>s>>  U(o)^1,
\end{CD}
$$
which defines a bundle chart.   We also have the commutative diagram
$$
\begin{CD}
G\ltimes TO  @> T\Psi>> T\mathscr{C}\\
@VVV  @VVV\\
G\ltimes O^1    @> \Psi^1 >>  \mathscr{C}^1,
\end{CD}
$$
where on the bottom we have a unformizer for the polyfold $\mathscr{C}^1$. In order to see this diagram as a bundle uniformizer
we need to equip $|T\mathscr{C}|$ with a metrizable topology.
The collection of isomorphism classes in $T\mathscr{C}$ is a set,  and as we shall see next 
the orbit space $|T\mathscr{C}|$ carries a natural metrizable topology. 
\begin{proposition}
The orbit space $|T\mathscr{C}|$ carries a natural metrizable topology. In particular $T\mathscr{C}$ is in a natural way a GCT.
Moreover $|P|:|T\mathscr{C}|\rightarrow |\mathscr{C}_1|$ is continuous. (Recall that $\mathscr{C}_1$ is a GCT.)
\end{proposition}
\begin{proof}
Pick a covering family $\bm{\Psi}$ and construct the associated ep-groupoid $X=X_{\bm{\Psi}}$ together with the natural equivalence 
$$
\Gamma:=\Gamma_{\bm{\Psi}}\colon X\rightarrow \mathscr{C},
$$
which induces a homeomorphism between the orbit spaces. Then $X$ is a metrizable ep-groupoid
implying that $TX$ is a metrizable ep-groupoid, see Section \ref{section_Tangent_Ep-Groupoid}.  We define 
\begin{eqnarray}\label{EQNW17211}
T\Gamma:TX\rightarrow T\mathscr{C}
\end{eqnarray}
as follows. Recall that the object M-polyfold is the disjoint union of $O_\lambda$,  and $q_\lambda\in O_\lambda$ is mapped
by $\Gamma$ to the object $\Psi_\lambda(q_\lambda)$. With $\alpha=\Psi_\lambda(q_\lambda)$ pick $\Psi\in F(\alpha)$ and assume that $\alpha$ is on level one.
Assuming that  $\Psi(q)=\alpha$ we obtain the element $(q_{\lambda},1_{\alpha},q)\in \bm{M}(\Psi_\lambda,\Psi)$ and  as before
the sc-diffeomorphism 
$$
\sigma:(U(q_\lambda),q_\lambda)\rightarrow (U(q),q).
$$
Then we define for $h_\lambda \in T_{q_\lambda}O_\lambda$ 
$$
T\Gamma( h_\lambda) = [\Psi_{\lambda}(q_\lambda), \Psi, T\sigma(q_\lambda)(h_\lambda)].
$$
Passing to orbit spaces we obtain the bijection
$$
|T\Gamma| \colon |TX|\rightarrow |T\mathscr{C}|.
$$
The orbit space $|TX|$ has a natural metrizable topology since $X$ is metrizable, see Theorem \ref{paraXCV}, and consequently defines a metrizable topology
on $|T\mathscr{C}|$. It is easily verified that the definition of the topology on $|T\mathscr{C}|$ does not depend on the choice of
the covering family $\bm{\Psi}$. 
\qed \end{proof}

Hence, starting with $(F,\bm{M})$ for the GCT $\mathscr{C}$ we have constructed a GCT $T\mathscr{C}$ together with
$P:T\mathscr{C}\rightarrow \mathscr{C}_1$ which induces a continuous map between orbit spaces.
We also constructed a functor $TF:\mathscr{C}_1^-\rightarrow \text{SET}$
which associates to $\alpha$ the set $TF(\alpha)$ consisting of all $T\Psi$ with $\Psi\in F(\alpha)$. These we view as the commutative diagrams
$$
\begin{CD} 
G\ltimes TO @> T\Psi>> T\mathscr{C}\\
@VVV @V P VV\\
G\ltimes O^1 @>\Psi^1 >>  \mathscr{C}_1.
\end{CD}
$$
The top horizontal map is fiber-wise a topological linear isomorphism. Passing to orbit spaces 
$|T\Psi|$ and $|\Psi^1|$ become homeomorphisms onto open sets. In addition we have a construction $\bm{M}$
which associates to $T\Psi_1$ and $T\Psi_2$ the structure of a bundle over an M-polyfold 
$$
\bm{M}(T\Psi_1,T\Psi_2)\rightarrow \bm{M}({(\Psi_1)}^1,{(\Psi_2)}^1).
$$
We  associate to $\alpha$
the collection $TF(\alpha)$ of all $T\Psi$, where $\Psi$ varies in $F(\alpha)$. In this way we keep track of the bundle structure
and we obtain the construction $(TF,\bm{M})$ for the diagram $T\mathscr{C}\rightarrow \mathscr{C}_1$ between GCT's
and equip this with a polyfold bundle structure. Then $TF:\mathscr{C}_1^-\rightarrow \text{SET}$ is a functor.

A slightly weaker construction produces the functor $TFP:=TF\circ P: T\mathscr{C}\rightarrow \text{SET}$. Here we forget 
about the bundle structure and just consider the GCT $T\mathscr{C}$.
 The polyfold structure $TFP$ associates to $k\in T_\alpha\mathscr{C}$ the set $TF(P(k))$ consisting of all $T\Psi$, where
$\Psi$ varies over $\Psi\in F(\alpha)$. Hence, if $k,k'\in T_\alpha\mathscr{C}$ it holds that $TF(P(k))=TF(P(k'))$ and consequently
$TFP:T\mathscr{C}\rightarrow \text{SET}$.

Summarizing the construction we obtain the following result.
\begin{theorem}
Let $\mathscr{C}$ be a GCT equipped with a polyfold structure $(F,\bm{M})$. Then there exists a natural construction 
of a category $T\mathscr{C}$ called the {\bf tangent category}\index{Tangent category} and a simultaneous polyfold structure
for $T\mathscr{C}\rightarrow \mathscr{C}_1$ denoted by $(TF,\bm{M})$. 
\end{theorem}
 We can visualize the two different constructions as follows
$$
\mathscr{C} \ \ \ \  \rightsquigarrow  \ \ \ \   T\mathscr{C}\ \ \ \text{This associates to a polyfold another one}.
$$
The bundle view-point is represented by 
 the following diagram
$$
\mathscr{C} \ \ \ \  \rightsquigarrow \ \ \ \    \begin{CD} 
T\mathscr{C}\\
@V P VV\\
\mathscr{C}^1
\end{CD}
\ \ \ \text{This associates to a polyfold a polyfold bundle.}
$$
\begin{remark}\label{REMARG238}\index{R-  Remark on how to transport notions and results over to the categorical setting}
There are many constructions and notions we can carry over from the M-polyfold or ep-groupoid situation.
Let $\mathscr{C}$ be a polyfold and $\alpha$ a smooth object. Then the tangent $T_\alpha\mathscr{C}$ is
an sc-Banach space.  We can define the {\bf  partial cone}\index{Partial cone $C_\alpha\mathscr{C}\subset T_\alpha\mathscr{C}$} $C_\alpha\mathscr{C}\subset T_\alpha\mathscr{C}$ and the {\bf reduced tangent space}\index{Reduced tangent $T^R_\alpha\mathscr{C}$} $T^R_\alpha\mathscr{C}$ following
Definition \ref{reduced_cone_tangent}. These can be defined in the obvious way for an sc-smooth model
$X$ and pushed forward via $T\Gamma$, where $\Gamma:X\rightarrow \mathscr{C}$ is the associated 
equivalence.  In the case that $\mathscr{C}$ is tame the partial cone will be a partial quadrant
and we can define what it means that a finite-dimensional linear subspace is in {\bf good position}\index{Good position to $C_\alpha\mathscr{C}$} to $C_\alpha\mathscr{C}$.
\qed
\end{remark}
\section{Subpolyfolds}
Let the GCT $\mathscr{C}$ be equip\-ped with a polyfold structure $(F,\bm{M})$.
\begin{definition}\index{D- Subpolyfold}
A {\bf subpolyfold} $\mathscr{A}$ of $\mathscr{C}$ is a saturated full subcategory such that for 
an object $\alpha$ in $\mathscr{A}$ and $\Psi\in F(\alpha)$, say $\Psi:G\ltimes O\rightarrow \mathscr{C}$,
the set of all objects $A_\Psi\subset O$ with $\Psi(A)\subset \mathscr{A}$ is a sub-M-polyfold of $O$.
\qed
\end{definition}
Although the definition seems to require to check the properties of $\mathscr{A}$ with respect to all $\Psi$, a collection which is not even a set in general, this is not true. In fact, it only has to be 
checked with respect any covering family $\bm{\Psi}$. In this spirit we give an alternative definition.
\begin{definition}[Alternative Definition of subpolyfold]  Let the GCT $\mathscr{C}$ be equip\-ped with a polyfold structure
$(F,\bm{M})$. A saturated full subcategory $\mathscr{A}$ is subpolyfold provided for a suitable sc-smooth local model
$X_{\bm{\Psi}}$ with natural equivalence $\Gamma_{\bm{\Psi}}:X_{\bm{\Psi}}\rightarrow \mathscr{C}$
the preimage $A$ under $\Gamma_{\bm{\Psi}}$ is an ep-subgroupoid in the sense of Definition \ref{Xep-subgroupoid}.
\qed
\end{definition}
\begin{remark}\index{R- On the equivalence of subpolyfold notions}
Note that the choice of $X_{\bm{\Psi}}$ in the definition does not matter.  We leave it to the reader to verify that 
the two given definitions of a subpolyfold are equivalent.
\qed
\end{remark}
Given a subpolyfold $\mathscr{A}$ of $\mathscr{C}$ we shall show next that it inherits 
from the structure on $\mathscr{C}$ a natural polyfold structure.
For every object $\alpha$ and $\Psi\in F(\alpha)$
$$
\Psi:G\ltimes O\rightarrow \mathscr{C}
$$
 define $\Psi_{\mathscr{A}}:G\ltimes A_\Psi\rightarrow \mathscr{A}$
as the restriction of $\Psi$. Cleary $\Psi_{\mathscr{A}}$ is injective on objects and fully faithful. Moreover,
we have the commutative diagram
$$
\begin{CD}
|G\ltimes A_\Psi| @>|\Psi_\mathscr{A}|>>| \mathscr{A}|\\
@VVV  @VVV\\
|G\ltimes O| @>|\Psi| >>  |\mathscr{C}|
\end{CD}
$$
where the horizontal maps are homeomorphic maps onto the open images, and the vertical maps
are topological embeddings, which in fact is even true level-wise. Hence every $\Psi_{\mathscr{C}}$ is a uniformizer.
We obtain a functor $F_{\mathscr{A}}:\mathscr{A}^-\rightarrow \text{SET}$ which associates to an object $\alpha$ the 
set of unifomizers $F_{\mathscr{A}}(\alpha)$ consisting of all $\Psi_\mathscr{A}$, where $\Psi$ varies in $F(\alpha)$.
Given two uniformizers $\Psi_\mathscr{A}$ and $\Psi_{\mathscr{A}}'$ associated to objects $\alpha$ and $\alpha'$ in $\mathscr{A}$
we consider the associated transition set $\bm{M}(\Psi_\mathscr{A},\Psi_\mathscr{A}')$ which is a subset of 
$\bm{M}(\Psi,\Psi')$.
\begin{lemma}
$\bm{M}(\Psi_\mathscr{A},\Psi_\mathscr{A}')$ is a sub-M-polyfold of $\bm{M}(\Psi,\Psi')$.
\end{lemma}
\begin{proof}
For $(o,\phi,o')\in \bm{M}(\Psi_{\mathscr{A}},\Psi_{\mathscr{A}}')$ we find open neighborhoods $U(o,\phi,o')$ in 
$\bm{M}(\Psi,\Psi')$ and $U(o)$ in $O$ such that $s:U(o,\phi,o')\rightarrow U(o)$ is an sc-diffeomorphism. Then $A_\Psi\cap U(o)$ is a sub-M-polyfold of $U(o)$ and also of $O$ and the preimage of $U(o)\cap A_\Psi$ under $s$ is precisely 
$U(o,\phi,o')\cap \bm{M}(\Psi_{\mathscr{A}},\Psi_{\mathscr{A}}')$. This implies that $\bm{M}(\Psi_{\mathscr{A}},\Psi_{\mathscr{A}}')$ is a sub-M-polyfold of $\bm{M}(\Psi,\Psi')$.
\qed \end{proof}
Since the $\bm{M}(\Psi_\mathscr{A},\Psi_\mathscr{A}')$ are sub-M-polyfolds, it follows that the associated 
structural maps are sc-smooth, since they are restrictions of the ambient structural maps for the $\bm{M}(\Psi,\Psi')$.
Summarizing we obtain a polyfold construction $(F_\mathscr{A},\bm{M})$ for the subpolyfold $\mathscr{A}$.
Hence we have proved.
\begin{proposition}
Given a GCT $\mathscr{C}$ with polyfold structure $(F,\bm{M})$ and a subpolyfold $\mathscr{A}$, the 
full subcategory $\mathscr{A}$ as stand-alone GCT inherits a natural polyfold structure $(F_\mathscr{A},\bm{M})$.
\qed
\end{proposition}

A typical example of subpolyfolds are  local faces or faces if the category is face-structured.
First we introduce the notion of faces into the current context.
\begin{definition}\index{D- Face of $\mathscr{C}$}
Assume the GCT $\mathscr{C}$ is equipped with a tame polyfold structure. A {\bf face} of $\mathscr{C}$ 
is the full subcategory $\mathscr{C}_\theta$ associate to the closure $\theta$ in $|\mathscr{C}|$ 
of a connected component $\theta^\circ$ in $\{z\in |\mathscr{C}|\ |\ d_{|\mathscr{C}|}(z)=1\}$.
We say that $\mathscr{C}$ is {\bf face-structured} provided every $z\in |\mathscr{C}|$ belongs to precisely $d(z)$-many $\theta$.
\qed
\end{definition}
Assume that $X_{\bm{\Psi}}$ is an sc-smooth model of $\mathscr{C}$ and $\Gamma_{\bm{\Psi}}: X_{\bm{\Psi}}\rightarrow \mathscr{C}$ the natural equivalence.
Passing to orbit spaces we obtain the homeomorphism
$$
|\Gamma_{\bm{\Psi}}|: |X_{\bm{\Psi}}|\rightarrow |\mathscr{C}|,
$$
which preserves the induced degeneracy index. Hence the closures of connected components
in the subspaces defined by $d_{|X_{\bm{\Psi}}|}$ and $d_{|\mathscr{C}|}$ correspond.
In view of Definition \ref{DEF1113} the faces of $X_{\bm{\Psi}}$ and 
$\mathscr{C}$ correspond. 
In particular,  the tame $\mathscr{C}$ is face-structured if and only if any of its sc-smooth models $X_{\bm{\Psi}}$ is face-structured. This follows from the discussion after Definition \ref{DEF1113}.  
\begin{proposition}
Assume the GCT $\mathscr{C}$ is equipped with a tame polyfold structure $(F,\bm{M})$.
   If $\mathscr{C}$ is face structured, every face is a subpolyfold.\qed
\end{proposition}

\section{Branched Ep\texorpdfstring{$^+$}{plll}-Subcategories}
The branched ep$^+$-subcategories are important since they appear as the solution spaces
for sc-Fredholm section functors.
Assume that $\mathscr{C}$ is a GCT equipped with a polyfold structure.
\begin{definition}
A functor $\Theta:\mathscr{C}\rightarrow {\mathbb Q}^+$ is called a branched ep$^+$-subcategory
provided for an object $\alpha$ and $\Psi\in F(\alpha)$ the functor $\Theta\circ\Psi:G\ltimes O\rightarrow {\mathbb Q}^+$
is a branched ep$^+$-subgroupoid in the sense of Definition \ref{DEF912}.
\qed
\end{definition}
Again, the property of being a branched ep$^+$-subcategory only has to be checked 
for a covering set of uniformizers. 
More precisely, if $\bm{\Psi}={(\Psi_\lambda)}_{\lambda\in\Lambda}$  is a covering set of uniformizers and for every $\Psi_\lambda$
the functor $\Theta\circ\Psi_\lambda$ is a a branched ep$^+$-subgroupoid, the same property holds for an arbitrary $\Psi$ coming from 
$F$.  Consider $\Gamma:X\rightarrow \mathscr{C}$, where $\Gamma$ and $X$ are associated to the covering family $\bm{\Psi}$.
Given $\Theta:\mathscr{C}\rightarrow {\mathbb Q}^+$ we define the {\bf pull-back}\index{Pull-back $\Gamma_{\bm{\Psi}}^\ast\Theta$} by 
$$
\Gamma^\ast\Theta :=\Theta\circ\Gamma,
$$
which defines a functor $X\rightarrow {\mathbb Q}^+$.  Given $\Theta:X\rightarrow {\mathbb Q}^+$ we can also define a {\bf push forward}
$$
\Gamma_\ast\Theta:\mathscr{C}\rightarrow {\mathbb Q}^+
$$
as follows. Given an object $\alpha$ of $\mathscr{C}$ pick an object $x\in X$ so that there exists a morphism $\phi:\Gamma(x)\rightarrow \alpha$
and define $(\Gamma_\ast\Theta)(\alpha):=\Theta(x)$. This is possible since $\Gamma$ as an equivalence is essentially surjective.
 One easily verifies that this  is well-defined.
The following facts are obvious.
\begin{proposition}\label{XPROP1742}
Assume $X_{\bm{\Psi}}$ is a sc-smooth model for $\mathscr{C}$ (equipped with $(F,\bm{M})$) and $\Gamma_{\bm{\Psi}}:X_{\bm{\Psi}}\rightarrow \mathscr{C}$ is the canonical equivalence. Then the following holds.
\begin{itemize}
\item[{\em (1)}] \ If $\Theta:X_{\bm{\Psi}}\rightarrow {\mathbb Q}^+$
is a branched ep$^+$-subgroupoid, then $(\Gamma_{\bm{\Psi}})_\ast\Theta$ is a branched ep$^+$-subcategory.
\item[{\em (2)}]\   If $\Theta:\mathscr{C}\rightarrow {\mathbb Q}^+$ is a branched ep$^+$-subcategory, then $\Gamma_{\bm{\Psi}}^\ast\Theta$
is a branched ep$^+$-subgroupoid.
\end{itemize}
\qed
\end{proposition}
Given two covering families $\bm{\Psi}$ and $\bm{\Psi}'$ with associated equivalences $\Gamma:X\rightarrow \mathscr{C}$ and $\Gamma':X'\rightarrow \mathscr{C}$ denote by $\mathfrak{f}:X\rightarrow X'$ the canonical generalized isomorphism. Recall that we can use $\mathfrak{f}$ to push forward or pull back branched ep$^+$-subgroupoids. The reader will easily verify the identities
$$
\mathfrak{f}^\ast (\Gamma')^\ast \Theta = \Gamma^\ast\Theta \ \text{and}\ \ \Gamma_\ast\Theta' = \Gamma'_\ast\mathfrak{f}_\ast\Theta',
$$
for $\Theta:\mathscr{C}\rightarrow {\mathbb Q}^+$ and $\Theta':X\rightarrow {\mathbb Q}^+$.

\begin{definition}\index{D- Compact $\Theta$}
Assume that $\mathscr{C}$ is a polyfold and $\Theta:\mathscr{C}\rightarrow {\mathbb Q}^+$ a branched ep$^+$-subcategory.
\begin{itemize}
\item[(1)] \  We say $\Theta$ is {\bf compact} provided the orbit space $|\supp(\Theta)|$ is a compact subset
of $|\mathscr{C}|$.
\item[(2)] \  We say $\Theta$ is {\bf tame}\index{D- Tame $\Theta$} if for a suitable covering family $\bm{\Psi}$ with associated
$\Gamma:X\rightarrow \mathscr{C}$ the pull-back $\Gamma^\ast\Theta$ which is a branched ep$^+$-subgroupoid is tame. We note that this is well-defined since for two different choices of covering family
the generalized isomorphism $\mathfrak{f}:X\rightarrow X'$ via pull-back or push-forward of
branched ep$^+$-subgroupoids preserves tameness.
\end{itemize}
\qed
\end{definition}

In view of Proposition \ref{XPROP1742} we can define following Definition \ref{DEF917} 
branched ep$^+$-subcategories which are of manifold-type or orbifold type. 
\begin{definition}\label{GHDEF}
Let $\mathscr{C}$ be a polyfold and $\Theta:\mathscr{C}\rightarrow {\mathbb Q}^+$ a branched
ep$^+$-subcategory. 
\begin{itemize}
\item[(1)] \  We say that $\Theta$ is of {\bf manifold-type} provided between any two objects 
$\alpha$ and $\alpha'$ in $\supp(\Theta)$ there is at most one morphism, and in addition
$\Theta$ only takes the values $0$ and $1$.
\item[(2)]\   $\Theta$ is said to be of {\bf orbifold type} provided $\Theta$ only takes the values $0$ and $1$.
\item[(3)] \ We say $\Theta$ is tame provided for $\alpha$ being in the support of $\Theta$ and $\Psi\in F(\alpha)$
the functor $\Theta\circ \Psi:G\ltimes O\rightarrow {\mathbb Q}^+$ is a tame branched ep$^+$-subgroupoid. Equivalently
for an sc-smooth model $X$ the functor $\Theta\circ \Gamma:X\rightarrow {\mathbb Q}^+$ is tame branched ep$^+$-subgroupoid.
\end{itemize}
\qed
\end{definition}
As a consequence of Proposition \ref{PROPY918} we obtain the following result.
\begin{proposition}
Let the GCT $\mathscr{C}$ be equipped with a polyfold structure $(F,\bm{M})$ and
assume that $\Theta:\mathscr{C}\rightarrow {\mathbb Q}^+$ is a branched ep$^+$-subcategory.
\begin{itemize}
\item[{\em (1)}] \  If $\Theta$ is of manifold-type then $S=|\supp(\Theta)|$ has in a natural way the structure
of a smooth finite-dimensional manifold (the M$^+$-polyfold version).
\item[{\em (2)}] \ If $\Theta$ is of orbifold-type then $S=|\supp(\Theta)|$ has a natural smooth orbifold structure
(the M$^+$-polyfold version).
\end{itemize}
\end{proposition}
See Remark \ref{remark918} concerning our notion of smooth manifold or orbifold.
\begin{proof}
We reduce the consideration to Proposition \ref{PROPY918}.  Let $\bm{\Psi}$ be a covering set of uniformizers and $X=X_{\bm{\Psi}}$ be the associated ep-groupoid and 
$\Gamma:X\rightarrow \mathscr{C}$ the corresponding equivalence. If follows from the 
definition that $\Theta\circ\Gamma$ is of manifold or orbifold-type in the ep-groupoid 
sense, see Definition \ref{DEF917}, if and only if $\Theta$ is of the corresponding type
defined in Definition \ref{GHDEF}.  Assume that $\bm{\Psi}'$ is a second 
covering set defining $X'$ and $\Gamma':X'\rightarrow \mathscr{C}$. There exists a 
natural generalized isomorphism $\mathfrak{f}:X\rightarrow X'$ satisfying 
\begin{eqnarray}\label{OPEQN007}
|\Gamma|=|\Gamma'|\circ |\mathfrak{f}|,
\end{eqnarray}
and
$$
\Theta\circ\Gamma' = \mathfrak{f}_\ast (\Theta\circ\Gamma).
$$
In the case that $\Theta$ is of manifold-type the generalized isomorphism $\mathfrak{f}$ induces
a classically smooth diffeomorphism 
$$
|\supp(\Theta\circ \Gamma)|\xrightarrow{|\mathfrak{f}|} |\supp(\Theta\circ\Gamma')|
$$
There exists a unique manifold structure on $|\supp(\Theta)|$
for which $|\Gamma|:|\supp(\Theta\circ\Gamma)|\rightarrow |\supp(\Theta)|$ is a diffeomorphism
and this structure is independent of the choice of the covering family in view of
(\ref{OPEQN007}).  

In the case of $\Theta$ being of orbifold type the argument is similar and $|\mathfrak{f}|$
is a smooth orbifold maps. Details are left to the reader.
\qed \end{proof}

As we shall see a branched ep$^+$-subcategory $\Theta:\mathscr{C}\rightarrow {\mathbb Q}^+$ has, as in  the ep$^+$-subgroupoid case,
a tangent functor $T\Theta:T\mathscr{C}\rightarrow {\mathbb Q}^+$ as well as the associated section functor 
$\mathsf{T}_\Theta$, which is being used to define orientations via suitable lifts, see Definition \ref{DEFNX9212} and Theorem \ref{P15.2*1}.

The tangent $T\Theta:T\mathscr{C}\rightarrow {\mathbb Q}^+$ is defined as follows.  We pick a covering family $\bm{\Psi}$
and obtain the associated equivalence $\Gamma:X\rightarrow \mathscr{C}$. We can take the tangent $T(\Gamma^\ast\Theta)$
which is a branched ep$^+$-subgroupoid $T(\Gamma^\ast\Theta):TX\rightarrow {\mathbb Q}^+$. Then we use $T\Gamma:TX\rightarrow T\mathscr{C}$ to push it forward to obtain
$$
(T\Gamma)_\ast (T(\Gamma^\ast\Theta)):T\mathscr{C}\rightarrow {\mathbb Q}^+.
$$
This construction might depend a prior on the choice of $\bm{\Psi}$, but this turns out to be not the case. 
Namely with $\mathfrak{f}:X\rightarrow X'$ being the canonical generalized isomorphism it follows that $T\mathfrak{f}:TX\rightarrow TX'$
is the canonical isomorphism for the uniformizer construction for $T\mathscr{C}$. Using that
$$
\mathfrak{f}_\ast (\Gamma^\ast\Theta)=(\Gamma'^\ast\Theta)
$$
it follows that
$$
T\mathfrak{f}_\ast (T(\Gamma^\ast\Theta))= T(\Gamma'^\ast\Theta).
$$
Pushing forward by $T\Gamma'$ we obtain
\begin{eqnarray*}
   (T\Gamma')_\ast (T({(\Gamma')}^\ast\Theta))= (T\Gamma')_\ast T\mathfrak{f}_\ast (T(\Gamma^\ast\Theta))= (T\Gamma)_\ast (T(\Gamma^\ast\Theta)).
\end{eqnarray*}
\begin{definition}\index{D- Tangent of $\Theta$}
Let $\mathscr{C}$ be a polyfold and $\Theta:\mathscr{C}\rightarrow {\mathbb Q}^+$ a branched ep$^+$-subcategory.
The well-defined $T\Theta: T\mathscr{C}\rightarrow {\mathbb Q}^+$ given via a covering set $\bm{\Psi}$ and $\Gamma:X\rightarrow \mathscr{C}$ by
$$
T\Theta:=(T\Gamma)_\ast (T(\Gamma^\ast\Theta)),
$$
is called the {\bf tangent} of $\Theta$.
\qed
\end{definition}
Next we introduce $\mathsf{T}_\Theta$.
Starting with the polyfold $\mathscr{C}$ we define a functor on the full subcategory associated to the smooth objects denoted by
$$
\nu\colon \mathscr{C}_\infty\rightarrow \text{SET}.
$$
We associate to $\alpha$ the set $\nu(\alpha)$ consisting of finite formal sums $\sum \sigma_L\cdot L$, where $\sigma_L\geq 0$ is a rational
number and for almost all $L$ the number is zero. Moreover $L\subset T_\alpha\mathscr{C}$ is a smooth finite-dimensional linear subspace. If $\phi:\alpha\rightarrow \alpha'$ we define 
$$
\nu(\phi)\colon\nu(\alpha)\rightarrow \nu(\alpha') : \sum \sigma_L\cdot L\rightarrow \sum \sigma_L\cdot \phi_\ast L,
$$
where $\phi_\ast L:= T\phi (L)$ and $T\phi:T_\alpha\mathscr{C}\rightarrow T_{\alpha'}\mathscr{C}$ is the sc-operator 
associated to $\phi$. Associated to $\nu$ we have the category ${\mathcal L}_\nu$ whose objects are the pairs
$(\alpha,\mathsf{L})$, with $\mathsf{L}\in\nu(\alpha)$. This category fibers over $\mathscr{C}_\infty$
$$
{\mathcal L}_\nu\rightarrow \mathscr{C}_\infty:(\alpha,\mathsf{L})\rightarrow \alpha.
$$
In order to stay consistent with the notation introduced in the ep-groupoid case
we set $\text{Gr}(\mathscr{C}):={\mathcal L}_\nu$,
so that
$$
\text{Gr}(\mathscr{C})\rightarrow \mathscr{C}_\infty.
$$
To a given  branched ep$^+$-subcategory $\Theta:\mathscr{C}\rightarrow {\mathbb Q}^+$ we can associate 
a section functor $\mathsf{T}_{\Theta}$ of $\text{Gr}(\mathscr{C})\rightarrow \mathscr{C}_\infty$ as follows.
Take an sc-smooth model $X$ associated to a covering family $\bm{\Psi}$ and the canonical equivalence 
$$
\Gamma:X\rightarrow \mathscr{C} 
$$
with associated tangent functor $T\Gamma:TX\rightarrow T\mathscr{C}$, see (\ref{EQNW17211}). 
Then $\Gamma^\ast\Theta$ is a branched ep$^+$-subgroupoid which has an associated section functor
$\mathsf{T}_{\Gamma^\ast\Theta}$ of $\text{Gr}(X)\rightarrow X_\infty$.  Since $\Gamma^\ast\Theta$ is a functor and has near a point $x$ in the support
a particular local representation it follows that $\mathsf{T}_{\Gamma^\ast\Theta}(x)$ is invariant under all
$T\phi$, where $\phi\in G_x$, i.e. with $\sum\sigma_L\cdot L= \mathsf{T}_{\Gamma^\ast\Theta}(x)$ it holds
\begin{eqnarray}
T\phi(\sum\sigma_L\cdot L) =\sum \sigma_L\cdot L\ \ \text{for all}\ \phi\in G_x.
\end{eqnarray}
Observe that for a smooth $x\in X$ with $\Gamma(x)=\alpha$ a smooth finite-dimensional subspace 
$L$ of $T_xX$ is mapped by $T\Gamma$ to a finite-dimensional linear subspace of $T_\alpha\mathscr{C}$.
This also defines a map for the finite sums $\sum \sigma_L\cdot L$.  Hence we can transport for smooth $x\in X$
the finite combination via $T\Gamma$ to a finite combination of smooth finite-dimensional subspaces of $T_\alpha\mathscr{C}$.
We also have to make a definition for an arbitrary smooth $\alpha$ which might not lie in the image of $\Gamma$.
Here the fact that $\Gamma$ is essentially surjective will be important.
Assume that $x\in X$ is smooth and $\alpha$ is an object in $\mathscr{C}$ and $x\in X$  such that there exists a morphism
$\phi:\Gamma(x)\rightarrow \alpha$. Then the image of the finite sum $\mathsf{T}_{\Gamma^\ast\Theta}(x)$
under first applying $T\Gamma$ and then $T\phi$ is not depending on the choice of $\phi$
since the formal sum is invariant under the action of the isotropy group.
Therefore there is a well-defined push-forward operation $T\Gamma_\ast(\mathsf{T}_{\Gamma^\ast\Theta})$
which defines a section of $\text{Gr}(\mathscr{C})\rightarrow \mathscr{C}_\infty$.
\begin{proposition}
Let the GCT $\mathscr{C}$  be equipped with the polyfold structure $(F,\bm{M})$ and
let $\Theta:\mathscr{C}\rightarrow {\mathbb Q}^+$ be a branched ep$^+$-subcategory.
Then there exists a well-defined  section $\mathsf{T}_\Theta$ of $\text{Gr}(\mathscr{C})\rightarrow \mathscr{C}_\infty$
uniquely defined as follows.
Given a covering family $\bm{\Psi}$ with associated ep-groupoid $X$ and equivalence
$\Gamma:X\rightarrow \mathscr{C}$ the identity
$\mathsf{T}_\Theta = T\Gamma_\ast\mathsf{T}_{\Theta\circ\Gamma}$
holds.
\qed
\end{proposition}
The straight forward proof is left to the reader.
Next we introduce the important notion of an orientation. Again the basic constructions needed are those carried out in the ep-groupoid context, see  Section \ref{SECXV93}. 
We assume that the GCT $\mathscr{C}$ has been equipped with  a polyfold structure $(F,\bm{M})$. 
For a smooth object $\alpha$ in $\mathscr{C}$ we consider the tangent space $T_\alpha\mathscr{C}$ 
and consider finite sums 
$$
\wh{\mathsf{L}}=\sum\sigma_{\wh{L}} \cdot \wh{L},
$$
where all $\sigma_{\wh{L}}$ are non-negative rational numbers and almost all of these numbers are zero, and
$\wh{L}$ are smooth oriented finite-dimensional subspaces of the sc-Banach space $T_\alpha\mathscr{C}$.
The collection of such finite formal sums associated to $\alpha$ is denoted by $\wh{\text{Gr}}(\alpha)$ and the union
of all such objects is denoted by $\wh{\text{Gr}}(\mathscr{C})$, where the objects vary over all smooth ones.
We obtain
$$
\wh{\text{Gr}}(\mathscr{C})\rightarrow \mathscr{C}_\infty.
$$
There is a forgetful functor $\mathsf{f}:\wh{\text{Gr}}(\mathscr{C})\rightarrow \text{Gr}(\mathscr{C})$
which just forgets the orientations.  Given a covering set $\bm{\Psi}$ with associated equivalence 
$\Gamma:X\rightarrow \mathscr{C}$ and section functor $\wh{\mathsf{T}}$ of $\wh{\text{Gr}}(\mathscr{C})\rightarrow \mathscr{C}_\infty$
we can define a pull-back denoted by $T\Gamma^\ast\mathsf{T}$, as follows.  If $x\in X$ is smooth first evaluate $\mathsf{T}$ on 
the image $\Gamma(x)$ giving the formal sum $\sum \sigma_{\wh{L}} \wh{L}$, and then taking the preimages under the linear sc-operator
$T\Gamma(x):T_xX\rightarrow T_{\Gamma(x)}\mathscr{C}$. We define accordingly
$$
(T\Gamma^\ast\wh{\mathsf{T}})(x):= T\Gamma(x)^{-1}\left(\wh{\mathsf{T}}(\Gamma(x))\right).
$$
\begin{definition}\label{DEFX1746}
Let $\mathscr{C}$ be a polyfold and $\Theta:\mathscr{C}\rightarrow {\mathbb Q}^+$ 
be a branched ep$^+$-subcategory with the section  functor $\mathsf{T}_\Theta$ of $\text{Gr}(\mathscr{C})\rightarrow \mathscr{C}_\infty$. An orientation 
for $\Theta$ is a section functor $\wh{\mathsf{T}}_\Theta$ of $\wh{\text{Gr}}(\mathscr{C})\rightarrow \mathscr{C}_\infty$ such that
$$
\mathsf{f}\circ \wh{\mathsf{T}}_\Theta = \mathsf{T}_\Theta,
$$
having the additional property that for a covering set $\bm{\Psi}$ with associated
$\Gamma:X\rightarrow \mathscr{C}$ the pull-back $T\Gamma^\ast \wh{\mathsf{T}}_\Theta$
is an orientation for $\Gamma^\ast \Theta$ in the ep-groupoid sense according to Definition \ref{DEFNC932}.
\qed
\end{definition}
\begin{remark}
The part of Definition \ref{DEFX1746} requiring that he pull-back $T\Gamma^\ast \wh{\mathsf{T}}_\Theta$
is an orientation for $\Theta\circ\Gamma$ can be replaced by requiring that for an object $\alpha$ in $\supp(\Theta)$ and $\Psi\in F(\alpha)$
the pull back $T\Psi^\ast \wh{\mathsf{T}}_\Theta$ is an orientation for $\Theta\circ\Psi.$
\qed
\end{remark}

\section{Sc-Differential Forms and Stokes}
Let $\mathscr{C}$ be a polyfold. Our next goal is the definition of sc-differential forms on $\mathscr{C}$ and the associated de Rham complex.
Denote by $P:T\mathscr{C}\rightarrow \mathscr{C}^1$ the tangent bundle. On the categorical level we can build the $k$-fold
Whitney sum
$$
\oplus^k T\mathscr{C}\rightarrow \mathscr{C}_1.
$$
\begin{proposition}\index{P- Whitney sum polyfold bundle structure}
If $\mathscr{C}$ is a polyfold the $k$-fold Whitney sum has a natural structure as polyfold bundle
$$
\oplus^k T\mathscr{C}\rightarrow \mathscr{C}^1.
$$
\end{proposition}
\begin{proof}
Pick an object $\alpha$ in $\mathscr{C}_1$ and $\Psi\in F(\alpha)$, say $\Psi:G\ltimes O\rightarrow \mathscr{C}$.   Then $T\Psi$ fits into the commutative diagram
$$
\begin{CD}
TO @> T\Psi >> T\mathscr{C}\\
@V p VV   @V P VV\\
O^1 @> \Psi^1>> \mathscr{C}^1.
\end{CD}
$$
The collection of such $T\Psi$ defines the polyfold bundle structure for $T\mathscr{C}\rightarrow \mathscr{C}^1$ inducing 
the existing structure on $\mathscr{C}^1$.  

Given a morphism $(g,h)$ in $G\ltimes TO$, where $g\in G$ and $h\in T_qO$,  we have that 
$$
T\Psi(g,h)= (T\Psi(h),\Psi(g,q),T(\Psi(g,q))(T\Psi(h))).
$$
Here $T(\Psi(g,q))$ is the tangent of $\Psi(g,q):q\rightarrow g\ast q$.
 The $k$-fold Whitney sum $\oplus^k TO$ has the structure of a M-polyfold bundle 
over $O^1$, see Section \ref{subs_sc_differential}.
 We shall define $\oplus^kT\Psi$ covering $\Psi^1$ so that this data defines the structure for $
\oplus^k T\mathscr{C}\rightarrow \mathscr{C}_1$.  These functors
$$
\oplus^k T\Psi: G\ltimes \oplus ^k TO\rightarrow \oplus^k T\mathscr{C}
$$
are given as follows.  On objects 
$(h_1,...,h_k)$, where $h_i\in T_qO$,
$$
(\oplus^k T\Psi)(h_1,...,h_k)=(T\Psi(h_1),...,T\Psi(h_k)).
$$
A morphism $(g,(h_1,...,h_k)):(h_1,..,h_k)\rightarrow (g\ast h_1,...,g\ast h_k)$ is mapped as
\begin{eqnarray*}
&&(\oplus^k T\Psi)(g,(h_1,...,h_k))\\
& =&((T\Psi(h_1),\Psi(g,x),T(\Psi(g,x))(T\Psi(h_1))),......,\\
&&\ \ \ \ \ \ \ \ \ \ \ (T\Psi(h_k),\Psi(g,x),T(\Psi(g,x))(T\Psi(h_k)))).
\end{eqnarray*}
Given an object $\alpha$ of regularity $1$ we can associate to it the set consisting of all
$\oplus^k T\Psi$ where $\Psi$ varies over $F(\alpha)$. These cover $\Psi^1$ and fit into the commutative diagram
$$
\begin{CD}
\oplus^k TO @> \oplus^k T\Psi >> \oplus^k T\mathscr{C}\\
@V \oplus^k p VV   @V \oplus^k P VV\\
O^1 @> \Psi^1 >> \mathscr{C}^1.
\end{CD}
$$
The transition sets 
$\bm{M}(\oplus^k T\Psi,\oplus^k T\Psi')$ can be identified with a $k$-fold Whitney sum 
of $\bm{M}(T\Psi,T\Psi')\rightarrow\bm{M}(\Psi^1,\Psi'^1)$ which allows to define the bundle structure. This equips $\bm{M}(\oplus^k T\Psi,\oplus^k T\Psi')\rightarrow \bm{M}(\Psi^1,\Psi'^1)$ with the structure of a M-polyfold bundle and metrics for the orbit spaces.
It follows immediately that the source and target maps are sc-smooth local bundle isomorphisms.
 Clearly
$\oplus^k T\mathscr{C}$ is a GCT and it is easily verified  that the  $\oplus^kT\Psi$ together with the $\Psi^1$
have the polyfold bundle uniformizer properties.  
 \qed \end{proof}
Next we consider functors $\omega:\oplus^k T\mathscr{C}\rightarrow {\mathbb R}$ which are skew-symmetric, 
so that for $T\Psi$ the functor $\omega_{\,\Psi}:= \omega\circ (T\Psi\oplus..\oplus T\Psi)$ is sc-smooth 
$$
\omega_{\,\Psi}: TO\oplus..\oplus TO\rightarrow {\mathbb R}.
$$
Notationally we shall use the abbreviation
$$
\Psi^\ast \omega := \omega\circ \oplus^k T\Psi.
$$
If $\omega$ is such that $\omega_{\,\Psi}$ is sc-smooth for a covering set of uniformizers, then the smoothness holds
for any uniformizer coming from $F$. 
\begin{definition}
Assume a functor $\omega:\oplus^k T\mathscr{C}\rightarrow {\mathbb R}$, which is skew-symmetric and linear in each argument, 
is given. We say that $\omega$ is {\bf sc-smooth} provided there exists a covering set of uniformizers $\bm{\Psi}$
so that $\omega_{\,\Psi}$ is sc-smooth for every $\Psi\in \bm{\Psi}$.  We denote by $\Omega^k(\mathscr{C})$ the real vector space
consisting of the functors $\omega$ just defined.
\qed
\end{definition}
If $\mathscr{C}$ is equipped with a polyfold structure the same is true for $\mathscr{C}^i$, where $i=0,1,....$,
and the inclusion functors 
$$
\mathscr{C}^{i+1}\rightarrow \mathscr{C}^i
$$
are sc-smooth. Indeed, for  a uniformizer $\Psi$ we obtain the commutative diagram
$$
\begin{CD}
O^{i+1}  @>\Psi^{i+1}>> \mathscr{C}^{i+1}\\
@V\text{incl}VV   @V\text{incl} VV\\
O^{i} @>\Psi^{i}>>  \mathscr{C}^{i}
\end{CD}
$$
from which we see that the local representative for the inclusion $\text{incl}$ is an sc-smooth inclusion $O^{i+1}\rightarrow O^i$.
We can pass to tangents giving us
$$
\begin{CD}
T(O^{i+1}) @> T(\Psi^{i+1}) >> T(\mathscr{C}^{i+1})\\
@V VV @VVV\\
T(O^i) @> T(\Psi^i) >>  T(\mathscr{C}^i)
\end{CD}
$$
Here the vertical maps are inclusions. Given $\omega$ on $\oplus^k T(\mathscr{C}^i)$ we can pull it back to $\oplus^k T(\mathscr{C}^{i+1})$.
Hence we obtain 
$$
\Omega^k(\mathscr{C})\rightarrow\Omega^k(\mathscr{C}^1)\rightarrow....
$$
and can pass to the direct limit denoted by $\Omega^k_\infty(\mathscr{C})$. This can be done for every integer $k\geq 0$ and we define
$\Omega^\ast_\infty(\mathscr{C})$ to be the direct sum of all the $\Omega^k_\infty(\mathscr{C})$.
An element is denoted by $[\omega]$. We shall call $[\omega]$ {\bf homogeneous}\index{Homogeneous element} provided it belongs to some $\Omega^k_\infty(\mathscr{C})$. 
An element $[\omega]\in \Omega^k_\infty(\mathscr{C})$ has representatives $\omega$ 
defined on $\Omega^k(\mathscr{C}^i)$ which are of class sc$^i$. 

Recall from Section \ref{subs_sc_differential} that the exterior derivative
$d:\Omega^k(O^{i+1})\rightarrow \Omega^{k+1}(O^i)$ is well-defined and as a consequence 
we obtain
$$
d:\Omega^k_\infty(O)\rightarrow \Omega^{k+1}_\infty(O)
$$
by having $d$ act on $[\omega]$ via $d([\omega])=[d\omega]$. Using this fact we can define
$d$ on $\Omega^\ast_\infty(\mathscr{C})$ via uniformizers as $[\omega]\rightarrow [d\omega]$, where,
assuming the homogeneous case, 
for a representative $\omega$ in $\Omega^k(\mathscr{C}^{i+1})$ the representative $d\omega\in \Omega^{k+1}(\mathscr{C}^i)$ is given by
\begin{eqnarray}
 (d\omega)\circ \oplus^{k+1} T\Psi =  d(\omega\circ \oplus^kT\Psi)
\end{eqnarray}
for every uniformizer $\Psi$ or equivalently  (viewed on the right $i$-level) 
$$
\Psi^\ast d\omega = d(\Psi^\ast\omega).
$$
Again one only needs a covering set of unifomizers.  As a consequence, alternatively we can defined $d$ by taking a covering set $\bm{\Psi}$ of uniformizers and assume that $[\omega]$ is a homogeneous 
element in $\Omega^k_\infty(\mathscr{C})$. Associated to the covering set we have the sc-smooth model $X$ and equivalence
$$
\Gamma:X\rightarrow \mathscr{C}.
$$
We can represent a homogeneous $[\omega]$ for a $i\geq 2$ by some $\omega\in \Omega^k(\mathscr{C}^i)$. Then $\Gamma^\ast\omega =\omega\circ \oplus^kT\Gamma$ is at least sc$^1$ and defined on $T(X^i)\oplus..\oplus T(X^i)\rightarrow X^{i+1}$.
On the ep-groupoid we can apply the exterior derivative to obtain $d(\Gamma^\ast\omega)$ and we can push the result forward.
Then the definition $d[\omega] =[\Gamma_\ast d(\Gamma^\ast\omega)]$ is independent of the choices involved and it follows  that $d\circ d=0$ since it holds for ep-groupoids.
\begin{definition}\index{D- de Rham complex}
Let $\mathscr{C}$ be a GCT equipped with a polyfold structure $(F,\bm{M})$. Then 
$(\Omega^\ast_\infty(\mathscr{C}),d)$ is called the {\bf de Rham complex} associated to the polyfold $\mathscr{C}$.
\end{definition}

Having introduced sc-differential forms on polyfolds and oriented branched ep$^+$-subcategories we can introduce $\oint_{\wh{\Theta}}[\omega]$,
which is the integration of an sc-differential form over an oriented tame branched ep$^+$-subcategory provided the latter has a support with a compact orbit space. The definition of such an integration is being done via the use of an sc-smooth model, so that we can employ the results from Section \ref{IandST}. 
\begin{theorem}\label{THJK1754}
Assume that $[\omega]\in \Omega^k_\infty(\mathscr{C})$ and $\wh{\Theta}:\mathscr{C}\rightarrow {\mathbb Q}^+$ is a tame, oriented, branched ep$^+$-subcategory so that $\supp(\Theta)$ has a compact orbit space. Given a covering set $\bm{\Psi}$ with associated 
sc-smooth model $X$ and equivalence $\Gamma:X\rightarrow \mathscr{C}$ introduce the oriented $\Theta\circ\Gamma$, $\Gamma^\ast\wh{\mathsf{T}}_\Theta$, abbreviated by $\wh{\Theta}_{\,\bm{\Psi}}$
and the sc-smooth differential form $[\omega_{\,\bm{\Psi}}]=[\Gamma^\ast\omega]$. Then the branched integral in the ep-groupoid sense 
$$
\oint_{{\wh{\Theta}}_{\,\bm{\Psi}}}[\omega_{\,\bm{\Psi}}]
$$
has a value which is independent of the choice of the covering family $\bm{\Psi}$.  In particular the definition 
$$
\oint_{\wh{\Theta}} [\omega]:= \oint_{{\wh{\Theta}}_{\,\bm{\Psi}}}[\omega_{\,\bm{\Psi}}]
$$
is independent of the choice of $\bm{\Psi}$.
\end{theorem}
\begin{proof}
The basic input comes from Section \ref{SECRTY114}, where the compatibility of the relevant notions 
with respect to equivalences was discussed in depth.  Assume that we have taken two different covering sets of uniformizers $\bm{\Psi}$ and $\bm{\Psi}'$ resulting in two sc-smooth models and 
equivalences 
$$
\Gamma:X\rightarrow \mathscr{C}\ \ \text{and}\ \ \Gamma':X'\rightarrow \mathscr{C}.
$$
Denote by $\mathfrak{f}:X\rightarrow X'$ the associated natural generalized isomorphism. 
Considering $\wh{\Theta}_\Gamma:=(\Gamma^\ast\Theta$, $\Gamma^\ast\wh{\mathsf{T}}_\Theta)$ and
$\wh{\Theta}_{\Gamma'}:=(\Gamma'^\ast\Theta$, $\Gamma'^\ast\wh{\mathsf{T}}_\Theta)$ 
we note that 
$$
\mathfrak{f}^\ast\wh{\Theta}_{\Gamma'} =\wh{\Theta}_\Gamma.
$$
Similarly with $[\omega_{\,\bm{\Psi}}]$ and $[\omega_{\,\bm{\Psi}'}]$ being the pull-back forms
it holds hat
$$
\mathfrak{f}^\ast[\omega_{\,\bm{\Psi}'}]=[\mathfrak{f}^\ast\omega_{\,\bm{\Psi}'}]=[\omega_{\,\bm{\Psi}}].
$$
Applying Theorem \ref{THM1141} we compute with $K=|\supp(\Theta)|$
\begin{eqnarray*}
 \oint_{\wh{\Theta}_{\,\bm{\Psi}}}[\omega_{\,\bm{\Psi}}]&:=& \oint_{\wh{\Theta}_{\,\bm{\Psi}}}\omega_{\,\bm{\Psi}}
=\mu_{\omega_{\,\bm{\Psi}}}^{\wh{\Theta}_{\,\bm{\Psi}}}(|\Gamma|^{-1}(K))\\
&=& \mu_{\mathfrak{f}_\ast \omega_{\,\bm{\Psi}}}^{\mathfrak{f}_\ast\wh{\Theta}_{\,\bm{\Psi}}}(|\mathfrak{f}|\circ |\Gamma|^{-1}(K))=\mu_{\omega_{\,\bm{\Psi}'}}^{\wh{\Theta}_{\,\bm{\Psi}'}}(|\Gamma'|^{-1}(K))\\
&=&\oint_{\wh{\Theta}_{\,\bm{\Psi}'}}\omega_{\,\bm{\Psi}'}
=\oint_{\wh{\Theta}_{\,\bm{\Psi}'}}[\omega_{\,\bm{\Psi}'}].
\end{eqnarray*}
In view of this argument we can define $\oint_{\wh{\Theta}}[\omega]$ by picking a covering set of uniformizers and setting
\begin{eqnarray}
\oint_{\wh{\Theta}}[\omega] := \oint_{\wh{\Theta}_{\,\bm{\Psi}}}[\omega_{\,\bm{\Psi}}].
\end{eqnarray}
\qed \end{proof}
The  definition of the integral $\oint_{\wh{\Theta}}[\omega]$ is reduced to the computation of an integral
in the ep-groupoid setting, where we know that Stokes' Theorem holds. Hence it is not unreasonable to
expect that in our more general set-up a Stokes' Theorem holds as well, which turns out to be true.
In order to formulate the result we need to define a boundary integral based on a similar construction 
for ep-groupoids. 

First we need to introduce the boundary of $\Theta:\mathscr{C}\rightarrow {\mathbb Q}^+$.
We start Given a  branched ep$^+$-subcategory $\wh{\Theta}:\mathscr{C}\rightarrow {\mathbb Q}^+$
defined on the polyfold $\mathscr{C}$. We take a covering set $\bm{\Psi}$ of uniformizers with associated sc-smooth model $X$ and equivalence $\Gamma:X\rightarrow{\mathbb Q}^+$.   Then $\Gamma^\ast\Theta:X\rightarrow{\mathbb Q}^+$
is a branched ep$^+$-subgroupoid and following Definition \ref{DEF9214} there is an associated
functor $\partial(\Gamma^\ast\Theta):X\rightarrow {\mathbb Q}^+$ with support in $\partial X$.
We define 
$$
\partial\Theta:\mathscr{C}\rightarrow {\mathbb Q}^+: \alpha\rightarrow \Gamma_\ast(\partial(\Gamma^\ast\Theta).
$$
The definition does not depend on the choice of $\bm{\Psi}$. If $\bm{\Psi}'$ is a second covering family defining $\Gamma':X'\rightarrow \mathscr{C}$ we obtain the natural generalized isomorphism
$\mathfrak{f}:X\rightarrow X'$ which satisfies, as we recall $|\Gamma'|\circ |\mathfrak{f}|=|\Gamma$.
Since $\mathfrak{f}_\ast (\partial(\Gamma^\ast\Theta))=\partial(\Gamma'^\ast\Theta)$ in view of Theorem 
\ref{THMX1138} it follows that 
$$
\Gamma_\ast(\partial(\Gamma^\ast\Theta)=\Gamma'_\ast(\partial(\Gamma'^\ast\Theta)).
$$
\begin{definition}\index{D- Boundary $\partial\Theta$}
Assume that $\Theta:\mathscr{C}\rightarrow {\mathbb Q}^+$ is a branched ep$^+$-subcategory,
where $\mathscr{C}$ is a polyfold.  Then $\partial\Theta:\mathscr{C}\rightarrow {\mathbb Q}^+$ defined 
via a covering family $\bm{\Psi}$ and associated $\Gamma:X\rightarrow \mathscr{C}$
by
$$
\partial\Theta : = \Gamma_\ast(\partial(\Gamma^\ast\Theta))
$$
is well-defined and independent of the choices and called the {\bf boundary} of $\Theta$.
\qed
\end{definition}
As already explained in the context of ep-groupoids, although $\partial\Theta$ is well-defined
it usually does not have very good differential geometric properties.  One needs to impose additional
requirements to make sure that $\partial\Theta$ is sufficiently nice.  
We recall that for a branched ep$^+$-subcategory $\Theta:\mathscr{C}\rightarrow {\mathbb Q}^+$
the notions of being compact, tame, and oriented are well-defined.
As a consequence of  Theorem \ref{THM9510} and using that our relevant notions are well-behaved
under generalized isomorphisms 
we obtain the following theorem by arguments similar to those in the proof of Theorem 
\ref{THJK1754}.
\begin{theorem}\label{SOS1756}\index{T- Boundary integration}
Consider a polyfold  $\mathscr{C}$ and $\wh{\Theta}:\mathscr{C}\rightarrow {\mathbb Q}^+$ 
which is an oriented, tame, compact, branched ep$^+$-subcategory of dimension $n$. Let the orientation of $\Theta$ be given by $\wh{\mathsf{T}}_{\Theta}:\mathscr{C}_\infty\rightarrow \wh{\text{Gr}}(\mathscr{C})$.
Denote by $\partial\wh{\Theta}$ the boundary of $\Theta$ equipped with the induced
orientation $\wh{\mathsf{T}}_{\partial\Theta}$.  Given $[\omega]\in \Omega^{n-1}_\infty(\mathscr{C})$ there is a well-defined integral $\oint_{\partial\wh{\Theta}}[\omega]$. This integral is defined by taking 
a covering set $\bm{\Psi}$ with associated $\Gamma:X\rightarrow \mathscr{C}$ via
$$
\oint_{\partial\wh{\Theta}}[\omega] := \oint_{\partial(\Gamma^\ast\wh{\Theta})}[\omega_{\,\bm{\Psi}}],
$$
and the definition does not depend on the choice of the covering family.
\qed
\end{theorem}
Having Theorem \ref{THJK1754} and Theorem \ref{SOS1756} in place we can state the
Stokes Theorem. The proof is immediate since the statement is equivalent to the the ep-groupoid version.
\begin{theorem}[Stokes Theorem]\index{T- Stokes Theorem in polyfolds}
Let $\mathscr{C}$ be a polyfold and $\wh{\Theta}:\mathscr{C}\rightarrow {\mathbb Q}^+$ 
a tame, compact, oriented, branched ep$^+$-subcategory of dimension $n$. For  $[\omega]\in\Omega^{n-1}_\infty(\mathscr{C})$ the following identity holds.
\begin{eqnarray}
\oint_{\wh{\Theta}}d[\omega] =\oint_{\partial\wh{\Theta}} [\omega].
\end{eqnarray}
\end{theorem}
The Stokes Theorem plays an important role when defining invariants of moduli spaces obtained
as from solution categories of sc-Fredholm section functors after a slight perturbation by an sc$^+$-multisection functor.

 \section{Strong Bundle Structures}\label{CHAP176-}

The strong bundles over polyfolds, which we are going to consider arise in the following  way.
We denote by $\mathscr{C}$ a GCT equipped with a polyfold structure $(F,\bm{M})$ and we assume 
we are given a functor $\mu:\mathscr{C}\rightarrow \text{Ban}$. We shall abbreviate $\wh{\phi}:=\mu(\phi)$ 
and consider a category $\mathscr{E}=\mathscr{E}_\mu$, whose objects are pairs $(\alpha,e)$, where $e\in\mu(\alpha)$.
A morphism in $\mathscr{E}$ is a pair $(\phi,e)$, where $e\in \mu(s(\phi))$ and $(\phi,e)$ is  seen as
$$
(\phi,e):(s(\phi),e)\rightarrow (t(\phi),\wh{\phi}(e)),
$$
covering the underlying $\phi$
$$
\begin{CD}
(s(\phi),e) @> (\phi,e)>> (t(\phi),\wh{\phi}(e))\\
@VVV @VVV\\
s(\phi) @>\phi >> t(\phi).
\end{CD}
$$
We shall denote by $P:\mathscr{E}\rightarrow \mathscr{C}$ the natural projection functor which on objects takes the form
$P(\alpha,e)=\alpha$ and on morphisms $P(\phi,e)=\phi$.  The preimage of an object $\alpha$ is a Banach space, where the vector space operations are defined
$$
(\alpha,e)+\lambda\cdot (\alpha,e')=(\alpha,e+\lambda\cdot e').
$$
We also note that for a  given morphism the preimage $P^{-1}(\phi)$ is a Banach space, where the operations are defined similarly by
$$
(\phi,e)+\lambda\cdot (\phi,e')=(\phi,e+\lambda\cdot e').
$$
\begin{definition}[Vector Bundle GCT]\index{D- Vector bundle GCT}
Let $\mathscr{C}$ be a GCT. A {\bf vector bundle GCT} over $\mathscr{C}$ consists of a functor 
$\mu:\mathscr{C}\rightarrow \text{Ban}$ and a metrizable topology ${\mathcal T}$ on $|\mathscr{E}_\mu|$ such that the following holds. 
\begin{itemize}
\item[(1)]\ \ For every object $\alpha$ the natural map $P^{-1}(\alpha)\rightarrow |\mathscr{E}|$ is continuous.
\item[(2)] \ \ The map $|P|:|\mathscr{E}|\rightarrow |\mathscr{C}|$ is surjective, continuous and open.
\end{itemize}
\qed
\end{definition}
We would like to equip $P:\mathscr{E}\rightarrow \mathscr{C}$ with what we shall call the structure of a strong bundle over a polyfold for which the sc-smooth models are
strong bundles over ep-groupoids, which have been introduced previously.
To explain this we consider as the building blocks the  strong local bundles
$K\rightarrow O$ equipped with an action by strong bundle isomorphisms of a finite group $G$.  Hence we obtain translation groupoids fibering over each other
$$
p\colon G\ltimes K\rightarrow G\ltimes O.
$$
Then $p$ is a functor which on objects maps $h\in K_x$ to $x$ and maps the morphism $(g,h)$ to the morphism $(g,x)$.
Given an object $\alpha$ in $\mathscr{C}$, we are interested in functors  
$$
\bar{\Psi}:G\ltimes K\rightarrow \mathscr{E},
$$
where $p:K\rightarrow O$ is a strong bundle over the M-polyfold $O$ and the $G$-action on $O$ lifts to an action by strong bundle 
maps $K\rightarrow K$. Moreover, $\bar{\Psi}$ is injective on objects, fully faithful and fiberwise a topological linear isomorphism.
In addition we require that for a suitable $q\in O$ we have that $\bar{\Psi}(0_q)=(\alpha,0)$.
By mapping $q\in O$ first to $0_q\in K_q$ and then to $P(\bar{\Psi}(0_q))$ we obtain a functor
$$
\Psi: G\ltimes O\rightarrow \mathscr{C}
$$
which is injective on objects, and fully faithful. 
We obtain the commutative diagram
$$
\begin{CD}
G\ltimes K @>\bar{\Psi} >>   \mathscr{E}\\
@V p VV   @V P VV\\
G\ltimes O  @>\Psi >>   \mathscr{C}.
\end{CD}
$$  
Recall that by assumption we have a polyfold construction $(F,\bm{M})$ for $\mathscr{C}$. We can therefore say
that $\bar{\Psi}$ is compatible with the polyfold structure for $\mathscr{C}$ provided the underlying $\Psi$ belongs to $F(\alpha)$.

\begin{definition}\index{D- Strong polyfold bundle}
Let $P:\mathscr{E}\rightarrow \mathscr{C}$ be a vector bundle GCT, which is associated to a polyfold $\mathscr{C}$, a functor $\mu:\mathscr{C}\rightarrow \text{Ban}$ and a topology ${\mathcal T}$ on $|\mathscr{E}|$.
 A {\bf strong bundle uniformizer} for $P$ at the object $\alpha$ consists of a covariant functor
$\bar{\Psi}\colon G\ltimes K\rightarrow \mathscr{E}$ covering a $\Psi\in F(\alpha)$ that the following holds.
\begin{itemize}
\item[(1)]\   $G\ltimes K$ is the translation groupoid associated to the strong bundle $p:K\rightarrow O$ equipped with the action of a finite group by strong bundle isomorphisms.
\item[(2)] \  The functor $\bar{\Psi}$ is fully faithful and fiberwise a topological linear Banach space isomorphism, i.e.
for every $q\in O$ the map
$$
K_q \rightarrow P^{-1}(\Psi(\alpha)): k\rightarrow \bar{\Psi}(k)
$$
is a Banach space isomorphism.
\item[(3)] \ The induced map $|\bar{\Psi}|:|K|\rightarrow |\mathscr{E}|$ is a homeomorphism onto an open subset of the form $|P|^{-1}(U)$, covering the homeomorphism
$|\Psi|:|O|\rightarrow U$.
\end{itemize}
\qed
\end{definition}
The {\bf footprint}\index{Footprint} of $\bar{\Psi}$ is $|P|^{-1}(U)$ and it covers the {\bf base footprint}\index{Base footprint} $U$ of the underlying $\Psi$.  

\begin{definition}
Let $\mathscr{C}$ by a GCT equipped with a polyfold construction $(F,\bm{M})$ and assume $P:\mathscr{E}\rightarrow \mathscr{C}$
is a vector bundle GCT associated to $\mu:\mathscr{C}\rightarrow \text{Ban}$ and the topology ${\mathcal T}$. A {\bf compatible strong bundle uniformizer construction}
is given by a functor $\bar{F}:\mathscr{C}^-\rightarrow \text{SET}$ which associates to an object $\alpha$ a set $\bar{F}(\alpha)$ of strong bundle uniformizers
$$
\bar{\Psi}: G\ltimes K\rightarrow \mathscr{C},
$$
so that $(\alpha,0)$ is in the image and the induced $\Psi:G\ltimes O\rightarrow \mathscr{C}$ belongs to $F(\alpha)$.
\qed
\end{definition}
It is clear that given two strong bundle uniformizers $\bar{\Psi}$ and $\bar{\Psi}'$ we can form the transition set
$\bm{M}(\bar{\Psi},\bar{\Psi}')$ defined as a the weak fibered product associated to the diagram
$$
K\xrightarrow{\bar{\Psi}}\mathscr{E}\xleftarrow{\bar{\Psi}'} K'.
$$
Hence an element takes the form $(k,(\phi,\bar{\Psi}(k)),k')$ with $\phi:\Psi(p(k))\rightarrow \Psi'(p'(k'))$ being a morphism in $\mathscr{C}$ and $\bar{\Psi}'(k') =\mu(\phi)(\bar{\Psi}(k))$.
With $\Psi$ and $\Psi'$ being the underlying uniformizers we obtain a natural map
$$
\bm{M}(\bar{\Psi},\bar{\Psi}')\rightarrow \bm{M}(\Psi,\Psi'): (k,(\phi,\bar{\Psi}(k)),k')\rightarrow (p(k),\phi,p'(k')),
$$
where the fibers are Banach spaces in a natural way.  
We leave it to the reader to write down all the usual structural maps 
which are assumed to be fiberwise linear. For example the source maps fit into the commutative diagram
$$
\begin{CD}
\bm{M}(\bar{\Psi},\bar{\Psi}') @> s>>  K\\
@VVV @V p VV\\
\bm{M}(\Psi,\Psi') @> s>> O.
\end{CD}
$$
Finally we can define what it means to have a strong bundle construction. 
\begin{definition}
Let $P:\mathscr{E}\rightarrow \mathscr{C}$ be a vector bundle GCT, which is associated to a polyfold $\mathscr{C}$
and a functor $\mu:\mathscr{C}\rightarrow \text{Ban}$ and a topology ${\mathcal T}$ on $|\mathscr{E}|$.  A {\bf strong polyfold bundle construction}
for $P:\mathscr{E}\rightarrow \mathscr{C}$ compatible with the polyfold structure $(F,\bm{M})$ for the GCT $\mathscr{C}$ consists of a functor
$\bar{F}:\mathscr{C}^-\rightarrow \text{SET}$ associating to an object $\alpha$ a set of strong bundle uniformizers 
$$
\bar{\Psi}:G\ltimes K\rightarrow \mathscr{E}
$$
at $\alpha$ so that the underlying functors $\Psi:G\ltimes O\rightarrow \mathscr{C}$ belong to $F(\alpha)$. In addition we are given a strong bundle construction 
for every 
$$
\bm{M}(\bar{\Psi},\bar{\Psi}')\rightarrow \bm{M}(\Psi,\Psi'),
$$
 where the structure on the base is the already existing one 
coming from the polyfold structure for $\mathscr{C}$. The data $(\bar{F},\bm{M})$ has the property that for the strong bundle structures
all structural maps are strong bundle maps and have the usual properties, listed below for the reader, otherwise. We shall call $(\bar{F},\bm{M})$ a {\bf strong polyfold bundle structure}\index{D- Strong polyfold bundle structure}
for $P:\mathscr{E}\rightarrow \mathscr{C}$ compatible with the polyfold structure $(F,\bm{M})$ for $\mathscr{C}$.
\qed
\end{definition}
Assume we are given the polyfold $\mathscr{C}$, the functor $\mu:\mathscr{C}\rightarrow\text{Ban}$
and the topology ${\mathcal T}$ for $|\mathscr{E}|$.  If $(\bar{F},\bm{M})$ is the strong polyfold bundle construction (covering $(F,\bm{M})$) we are in the following situation.  Every transition set 
$\bm{M}(\bar{\Psi},\bar{\Psi}')$ comes equipped with a strong bundle structure for 
$\bm{M}(\bar{\Psi},\bar{\Psi}')\rightarrow \bm{M}(\Psi,\Psi')$, with the induced structure on the base being
the already existing one.  The source and target maps form sc-smooth strong bundle maps, which are local strong bundle isomorphisms
$$
\begin{CD}
K @< s<<  \bm{M}(\bar{\Psi},\bar{\Psi}')@>t>>  K'\\
@V p VV    @VVV @V p'VV\\
O @< s<<   \bm{M}(\Psi,\Psi') @> t>> O'.
\end{CD}
$$
The inversion map $\iota$ defined by $\iota(k,(\phi,\bar{\Psi}(k)),k')=(k',(\phi^{-1},\bar{\Psi}'(k')),k)$
covers the inversion map on the base and is an sc-smooth strong bundle isomorphism
$$
\begin{CD}
\bm{M}(\bar{\Psi},\bar{\Psi}')@>\iota>> \bm{M}(\bar{\Psi}',\bar{\Psi})\\
@VVV @VVV\\
\bm{M}(\Psi,\Psi') @>\iota>> \bm{M}(\Psi',\Psi).
\end{CD}
$$
The unit map $u$ is an sc-smooth strong bundle map fitting into the diagram
$$
\begin{CD}
K @> u>>  \bm{M}(\bar{\Psi},\bar{\Psi})\\
@V p VV   @VVV \\
O@>u>>   \bm{M}(\Psi,\Psi).
\end{CD}
$$
Finally the multiplication map is an sc-smooth strong bundle map fitting into the commutative diagram
$$
\begin{CD}
\bm{M}(\bar{\Psi}',\bar{\Psi}''){_{s}\times_t}\bm{M}(\bar{\Psi},\bar{\Psi}')@>m>> \bm{M}(\bar{\Psi},\bar{\Psi}'')\\
@VVV  @VVV\\
\bm{M}({\Psi}',{\Psi}''){_{s}\times_t}\bm{M}({\Psi},{\Psi}')@>m>> \bm{M}({\Psi},{\Psi}'').
\end{CD}
$$

\section{Proper Covering Constructions}\label{SECV17.7}
In applications the covering construction is important and is closely related to the construction for ep-groupoids given in 
Section \ref{SQWERTY116}.
\begin{definition}\index{D- Proper covering functor}\label{DEFNX17.7.1}
Let $\mathscr{A}$ and $\mathscr{B}$ be two GCT's. We say 
 that  the functor $\mathsf{P}:\mathscr{A}\rightarrow \mathscr{B}$  is a {\bf proper covering  functor} of GCT's provided the following holds
 \begin{itemize}
 \item[(1)]\  $|\mathsf{P}|:|\mathscr{A}|\rightarrow |\mathscr{B}|$ is continuous and surjective.
 \item[(2)]\      $\mathsf{P}$   on objects is finite-to-one and surjective.
 \item[(3)]\   For every $z\in |\mathscr{B}|$ and the points $y_1,...,y_k\in |\mathsf{P}|^{-1}(z)$
 there exists an open neighborhood $U(z)$ and mutually disjoint open neighborhoods $U(y_i)$ 
 such that
 $$
 |\mathsf{P}|^{-1}(U(z))= \bigcup_{i=1}^k U(y_i).
 $$
 \item[(4)] \  The map $\bm{A}\rightarrow \bm{B}{_{s}\times_{\mathsf{P}}}A:\phi\rightarrow (\mathsf{P}(\phi),s(\phi))$ is a bijection.
 \end{itemize}
Here $\bm{A}$ and $\bm{B}$ stand for the morphism classes and $A$ for the object class.
\qed
\end{definition}
Our aim is to equip $\mathsf{P}:\mathscr{A}\rightarrow \mathscr{B}$ with an sc-smooth structure
which reflects the covering property, and  in addition defines polyfold structures for $\mathscr{A}$ and $\mathscr{B}$.  

\begin{remark}\index{R- On uniformizer constructions for GCT's}
Recall that a uniformizer construction for a GCT $\mathscr{C}$ can be viewed as follows.  Given an object $\alpha$
we obtain the associated category having this single object $\alpha$ and as morphisms the isotropy group $G=G_\alpha$.  We write it as $G\ltimes\{\alpha\}$.
Then the uniformizer construction consists of a method of associating to the object $\alpha$ a M-polyfold
$O$ and to the morphisms a representation in $\text{Diff}_{sc}(O)$, together with a functor from 
the associated translation groupoid $G\ltimes O$ into $\mathscr{C}$. The construction $\bm{M}$
then `certifies' the compatibility of the various uniformizers. The construction we need for $\mathsf{P}$ generalizes this picture,  
and the basic building blocks have already been  described in Definition \ref{DEFR845} with the basic idea
of a geometric lift.
\qed
\end{remark}

Starting with $\mathsf{P}:\mathscr{A}\rightarrow \mathscr{B}$ we pick an object $\beta$ in $\mathscr{B}$ and note
that $\mathsf{P}^{-1}(\beta)$ consists of finitely many objects $\alpha_1,..,\alpha_k$.  Consider the associated
full subcategory denoted by $\mathscr{A}_\beta$ with object set $A_\beta=\{\alpha_1,..,\alpha_k\}$.
The functor $\mathsf{P}$ induces the covering functor
\begin{eqnarray}\label{EQN1713}
\mathsf{P}_\beta:\mathscr{A}_\beta\rightarrow G_\beta\ltimes \{\beta\}.
\end{eqnarray}
In order to define the uniformizers we have to describe the type of  domains we would like to have. The domains should be geometric lifts of the covering functor $\mathsf{P}_{\beta}$ in (\ref{EQN1713}) in the sense of Definition \ref{DEFR845}. 
The only difference is that we shall use pointed spaces as objects.
\begin{definition}
We denote by $\text{Diff}_{sc}^\ast$
the category whose objects are M-polyfolds together with a distinguished point and the morphisms are sc-diffeomorphisms preserving the distinguished points.
More precisely, an object in $\text{Diff}_{sc}^\ast$\index{$\text{Diff}_{sc}^\ast$} is a pair $(O,o)$ with $O$ being a M-polyfold and $o\in O$. 
A morphism $f:(O,o)\rightarrow (O',o')$ is an sc-diffeomorphism $f:O\rightarrow O'$ with $f(o)=o'$.
\qed
\end{definition}\index{D- Category $\text{Diff}_{sc}^\ast$}

The relevant domains for our upcoming functor constructions are associated 
to the following functorial data.
\begin{itemize}
\item[(i)] \ \ A choice of  functors $\mathsf{A}: \mathscr{A}_\beta\rightarrow \text{Diff}^\ast_{sc}$ and $\mathsf{B}:G_\beta\ltimes \{\beta\}\rightarrow \text{Diff}_{sc}^\ast$.
\item[(ii)]\ \   A natural transformation $\tau:\mathsf{A}\rightarrow \mathsf{B}\circ \mathsf{P}_\beta$.
\end{itemize}
The functor $\mathsf{B}$  associates to the (single) object $\beta$  a pointed M-polyfold $(O_\beta,o_\beta)$ with distinguished point $o_\beta\in O_\beta$.
To  a morphism $g\in G_\beta$ it associates an sc-diffeomorphism
$$
g \ast : (O_\beta,o_\beta)\rightarrow (O_\beta,o_\beta): q\rightarrow g\ast q
$$
preserving $o_\beta$. The functoriality guarantees that $g\ast(h\ast q)=(g\circ h)\ast q$.
This data will define the ep-groupoid $G_\beta\ltimes O_\beta$, which contains the distinguished object $o_\beta$.

 Associated to an object $\alpha$ in $\mathscr{A}_\beta$ we have a pointed M-polyfold
$({O}_\alpha,{o}_a)$ and associated to a morphism ${\phi}:\alpha\rightarrow \alpha'$ an sc-diffeomorphism
$$
\phi\ast:(O_\alpha,o_\alpha)\rightarrow (O_{\alpha'},o_{\alpha'}):o\rightarrow \phi\ast o.
$$
Again functoriality implies if $s(\phi')=t(\phi)$ that $\phi'\ast (\phi\ast q)=(\phi'\circ\phi)\ast q$
for $q\in O_{s(\phi)}$.
The natural transformation $\tau:\mathsf{A}\rightarrow \mathsf{B}\circ \mathsf{P}_\beta$ defines an
sc-diffeomorphisms for every object $\alpha$ in $\mathscr{A}_\beta$
$$
\tau_\alpha:(O_\alpha,o_\alpha)\rightarrow (O_\beta,o_\beta),
$$
so  that we obtain  for every morphism $\phi:\alpha\rightarrow \alpha'$ in $\mathscr{A}_\beta$  a commutative diagram of sc-diffeomorphisms
$$
\begin{CD}
(O_\alpha,o_\alpha) @>\phi\ast >> (O_{\alpha'},o_{\alpha'})\\
@V\tau_\alpha VV   @V \tau_{\alpha'} VV\\
(O_\beta,o_\beta) @> (\mathsf{P}_\beta(\phi))\ast >> (O_\beta,o_\beta).
\end{CD}
$$
We can take the disjoint union of the $O_\alpha$ denoting it by $E$, i.e.
\begin{eqnarray}
E=\bigsqcup_{\alpha\in \mathsf{P}^{-1}(\beta)} O_\alpha,
\end{eqnarray}
 and view it as
the object M-polyfold of an ep-groupoid, where the morphisms are pairs
$(\phi,o)$, $\phi\in \bm{A}_\beta$ with $o\in O_{s(\phi)}$, $s(\phi,o)=o$ and $t(\phi,o)=\phi\ast o$, i.e.
$$
(\phi,o):o\rightarrow \phi\ast o.
$$
We denote the morphism set by $\bm{E}$. It is the disjoint union
of all $\{\phi\}\times O_{s(\phi)}$, where $\phi$ varies over the morphisms 
in $\mathscr{A}_\beta$
\begin{eqnarray}
\bm{E} =\bigsqcup_{\phi\in \bm{A}_\beta} \left(\{\phi\}\times O_{s(\phi)}\right).
\end{eqnarray}
Each of the sets $\{\phi\}\times O_{s(\phi)}$ carries
a natural M-polyfold structure coming from the identification with $O_{s(\phi)}$. 
Hence $\bm{E}$ has a natural M-polyfold structure.
It is a trivial exercise to equip $E\equiv (E,\bm{E})$
with the structure of an ep-groupoid. The source map is given by $s(\phi,o)=o$
and the target map by $t(\phi,o)=\phi\ast o$.  We have the natural local sc-diffeomorphism on the object level
\begin{eqnarray}\label{EQNX11.75}
\wh{\tau}:E\rightarrow O_\beta
\end{eqnarray}
mapping $o\in O_\alpha$ to $\tau_\alpha(o)$. The preimages of the distinguished point $o_\beta$ in $O_\beta$ are
the distinguished points $o_\alpha\in O_\alpha$ for $\alpha\in A_{\beta}=\mathsf{P}^{-1}(\beta)$. 
 Then we can view $\wh{\tau}$ as a map
$$
\wh{\tau}:(E,A_\beta)\rightarrow (O_\beta,\beta).
$$
The morphisms 
$(\phi,o)$ are mapped by 
\begin{eqnarray}
\bm{\wh{\tau}}:\bm{E}\rightarrow G_\beta\times O_\beta: \bm{\wh{\tau}}(\phi,o)=(\mathsf{P}_\beta(\phi),\wh{\tau}(o)),
\end{eqnarray}
 which obviously is a local sc-diffeomorphism.  This together with (\ref{EQNX11.75}) shows that
 $\wh{\tau}\equiv (\wh{\tau},\bm{\wh{\tau}})$ defines a functor
\begin{eqnarray}
 \wh{\tau}:E\equiv (E,\bm{E})\rightarrow G_\beta\ltimes O_\beta.
\end{eqnarray}
 Finally we note that the map
$$
\bm{E}\rightarrow \text{mor}(G_\beta\ltimes O_\beta){_{s}\times_{\wh{\tau}}} E: (\phi,o)\rightarrow ((\mathsf{P}_\beta(\phi),\wh{\tau}(o)),o)
$$
is an sc-diffeomorphism. With other words $\wh{\tau}:E\rightarrow G_\beta\ltimes O_\beta$ is a proper covering map
between ep-groupoids of a particular form and it can be considered as a geometric lift of $\mathsf{P}_\beta:\mathscr{A}_\beta\rightarrow G_\beta\ltimes \{\beta\}$.  
The previous discussion motivates the next definition.
\begin{definition}\label{DEFQ1773}\index{D- Proper covering model}
A {\bf proper ep-groupoid covering model} consists of a translation groupoid $G\ltimes O$ associated to a finite group $G$ acting by sc-diffeomorphisms on a M-polyfold $O$, a distinguished point $\bar{o}\in O$ which is fixed by 
$G$, an ep-groupoid $E$, and a sc-smooth  functor $\wh{\tau}:E\rightarrow G\ltimes O$
having the following properties.
\begin{itemize}
\item[(1)]\  $\wh{\tau}:E\rightarrow O$ is a surjective local sc-diffeomorphism on the object level.
\item[(2)]\  For every $y\in \wh{\tau}^{-1}(\bar{o})$ there exists an open neighborhood $U(y)$
such that the $U(y)$ are mutually disjoint and
$$
\wh{\tau}^{-1}(O) = \bigcup_{y\in \wh{\tau}^{-1}(\bar{o})} U(y),
$$
and  for every such $y$ the map $\wh{\tau}:U(y)\rightarrow O$ is an sc-diffeomorphism.
\item[(3)] \  The map $\Gamma:\bm{E}\rightarrow (G\ltimes O){_{s}\times_{\wh{\tau}}}E$ given by
$$
\psi\rightarrow (\wh{\tau}(\psi),s(\psi))
$$
is an sc-diffeomorphism.
\end{itemize}
\qed
\end{definition}
\begin{remark}\label{REM17.7.5}
Compare the previous definition with Definition \ref{d-proper_covering} and the structural results
Theorem \ref{STRUC_1} and Theorem \ref{STRUC_2}.
\qed
\end{remark}
Let us show that the definition describes the previously discussed model in the case that the base M-polyfold is  connected.
Let $\wh{\tau}:E\rightarrow G\ltimes O$ be a proper ep-groupoid covering model in the sense 
of Definition \ref{DEFQ1773} with distinguished point $\bar{o}$. Suppose further that $O$ is connected. 
We consider the category $G\ltimes \{\bar{o}\}$ and denote by $\mathscr{A}_{\bar{o}}$
the full subcategory of $E$ associated to the object set $A_{\bar{o}}=\wh{\tau}^{-1}(\bar{o})$.
Denote by $\mathsf{P}:\mathscr{A}_{\bar{o}}\rightarrow G\ltimes\{\bar{o}\}$ the restriction of $\wh{\tau}$.
Using that $O$ is connected we infer that every $U(y)$ is connected. We associate to a morphism $\psi\in \bm{E}$ the pair $(y_{s(\psi)},y_{t(\psi)})\in A_{\bar{o}}^2$, so that
$$
s(\psi)\in U(y_{s(\psi)})\ \ \text{and}\ \ t(\psi)\in U(y_{t(\psi)}).
$$
\begin{lemma}
The map $\bm{E}\rightarrow A_{\bar{o}}^2:\psi\rightarrow (y_{s(\psi)},y_{t(\psi)})$ 
is locally constant, i.e. continuous.
\end{lemma}
\begin{proof}
The maps $\psi\rightarrow s(\psi)$ and $\psi\rightarrow t(\psi)$ are continuous and the sets $U(y)$
are open.  Hence the set of all $\psi$ with $(s(\psi),t(\psi))\in U(y)\times U(y')$
is open. 
\qed \end{proof}
Recall that the map
$$
\Gamma:\bm{E}\rightarrow (G\times O){_{s}\times_{\wh{\tau}}} E:\psi\rightarrow (\wh{\tau}(\psi),s(\psi))
$$
is a sc-diffeomorphism. For the following denote by $\tau_y:U(y)\rightarrow O$ the 
restriction of $\wh{\tau}$.  We obtain a map
$$
A_{\bar{o}}\rightarrow \text{mor}(\text{Diff}^\ast_{sc}):y\rightarrow \tau_y.
$$
For $\phi\in \bm{A}_{\bar{o}}$ denote by $\bm{E}_\phi$ the connected component
of $\bm{E}$ containing $\phi$. 
\begin{lemma}
If $\phi,\phi'\in \bm{A}_{\bar{o}}$ are two different elements then $\bm{E}_\phi\cap \bm{E}_{\phi'}=\emptyset$. Further we have the identity
$$
\bm{E} =\bigcup_{\phi\in \bm{A}_{\bar{o}}} \bm{E}_\phi.
$$
\end{lemma}
\begin{proof}
Take $\psi\in \bm{E}$ with $s(\psi)\in U(y)$ and $t(\psi)\in U(y')$. 
The map $\Gamma:{\bf E}\rightarrow (G\ltimes O){_{s}\times_{\wh{\tau}}}E:\psi\rightarrow (\wh{\tau}(\psi),s(\psi)) $ is a sc-diffeomorphism and denote the image of $\psi$ 
 by $((g_\psi,\wh{\tau}(s(\psi))),s(\psi))$.
Take a continuous path $\gamma:[0,1]\rightarrow O$ connecting
$\wh{\tau}(s(\psi))$ at time $a=0$ with $\bar{o}$ at time $1$, and consider the continuous path 
$$
[0,1]\rightarrow (G\times O)\times U(y): a\rightarrow ((g_\psi,\gamma(a)),\tau_y^{-1}(\gamma(a))).
$$
Then we obtain the continuous path
$$
[0,1]\rightarrow \bm{E}:a\rightarrow \Gamma^{-1}((g_\psi,\gamma(a)),\tau_y^{-1}(\gamma(a)))
$$
starting at $\psi$ and ending at the element $\phi_\psi:=\Gamma^{-1}((g_\psi,\bar{o}),y)$.
This shows that every $\psi$ belongs to some $\bm{E}_\phi$.

Assume next that $\phi,\phi'\in A_{\bar{o}}$ belong to the same path component. Take 
a continuous path $\wh{\gamma}$ connecting them. Then 
$$
\Gamma\circ \wh{\gamma}(a)= ((g_a,o_a),q_a).
$$
Clearly $g_a$ has to be independent of $a$ and we write $g=g_a$.
Then $\Gamma(\phi)=((g,\bar{o}),y)$ and $\Gamma(\phi')=((g,\bar{o}),y)$ which implies $\phi=\phi'$.
\qed \end{proof}
The maps $\psi\rightarrow y_{s(\psi)}$, $\psi\rightarrow y_{t(\psi)}$, and $\psi\rightarrow g_\psi$
are constant on $\bm{E}_\phi$.  Hence, if $\phi:y\rightarrow y'$, then 
$y_{s(\psi)}=y$, $y_{t(\psi)}=y'$, and $g_\psi= g_\phi$, where $\Gamma(\phi)=((g_\phi,\bar{o}),y)$.
\begin{lemma}
The following identity holds on $\bm{E}_\phi$, where $\phi:y\rightarrow y'$
$$
t(\psi)   = \tau_{y'}^{-1}(g_\phi \ast\tau_y(s(\psi))).
$$
\end{lemma}
\begin{proof}
If $\psi\in \bm{E}_\phi$, where $\phi:y\rightarrow y'$,  then
$\wh{\tau}(\psi) = (g_\phi,\wh{\tau}(s(\psi)))$ and consequently
$$
\tau_{y'}(t(\psi))=\wh{\tau}(t(\psi))=t(\wh{\tau}(\psi))= g_\phi\ast \wh{\tau}(s(\psi))=g_\phi\ast \tau_y(s(\psi)).
$$
This implies
$$
t(\psi)=\tau_{y'}^{-1}(g_\phi\ast \tau_y(s(\psi))).
$$
\qed \end{proof}
Given $\phi\in \bm{A}_{\bar{o}}$ we define the sc-diffeomorphism
$$
\phi\ast :U(s(\phi))\rightarrow U(t(\phi)): \phi\ast q =\tau_{t(\phi)}^{-1}(g_\phi\ast \tau_{s(\phi)}(q)).
$$
If $\phi,\phi'\in \bm{A}_{\bar{o}}$ with $t(\phi)=s(\phi')$ we compute
$$
(\phi'\circ\phi)\ast q = \phi'\ast(\phi\ast q)
$$
using that $(g_{\phi'\circ\phi},\bar{o})=\wh{\tau}(\phi'\circ\phi)=\wh{\tau}(\phi')\circ\wh{\tau}(\phi)=(g_{\phi'}g_\phi,\bar{o})$. We define the functor $\mathsf{A}:\mathscr{A}\rightarrow \text{Diff}^\ast_{\text{sc}}$ by
$y\rightarrow (U(y),y)$ and $\phi\rightarrow \phi\ast$. The functor $\mathsf{B}$ associates to $\bar{o}$ the pointed M-polyfold $(O,\bar{o})$.  The natural transformation $\tau:\mathsf{A}\rightarrow \mathsf{B}\circ\mathsf{P}$ associates to $y\in A$ the sc-diffeomorphism ${\tau}_y$. 
It is an easy exercise that the proper covering functor $\wh{\tau}: F\rightarrow G\ltimes O$ associated to the data $\mathsf{A},\mathsf{B},\tau$
is isomorphic to the $\wh{\tau}:E\rightarrow G\ltimes O$ we started with.

Next we are in the position to define the notion of uniformizer in the current context.
\begin{definition}\label{DEFR1774}\index{D- Proper covering uniformizer}
Let $\mathsf{P}:\mathscr{A}\rightarrow \mathscr{B}$ be a proper covering functor between two GCT's. 
A {\bf proper covering uniformizer} at the object $\beta$  in $\mathscr{B}$ with isotropy $G=G_\beta$, consists of a pair of functors
$\wh{\Psi}:E\rightarrow \mathscr{A}$ and $\Psi:G\ltimes O_\beta\rightarrow \mathscr{B}$, where 
$\wh{\tau}:E\rightarrow G\ltimes O_\beta$ is a proper ep-groupoid covering model, see Definition \ref{DEFQ1773}, such that the following holds.
\begin{itemize}
\item[(1)]\  The distinguished point $o_\beta\in O_\beta$ is mapped by $\Psi$ to $\beta$, and the distinguished points $o_\alpha\in O_\alpha$ with $\alpha\in \mathsf{P}^{-1}(\beta)$ 
are mapped bijectively by $\wh{\Psi}$ to the points in $\mathsf{P}^{-1}(\beta)$.
\item[(2)]\   $\Psi$ and $\wh{\Psi}$ are fully faithful, injective on objects,  and fit into the commutative functor diagram
$$
\begin{CD}
E @>\wh{\Psi} >> \mathscr{A}\\
@V\wh{\tau} VV  @V \mathsf{P}VV\\
G\ltimes O_\beta@>\Psi>>  \mathscr{B}.
\end{CD}
$$
\item[(3)]\   On the object level the preimage of $\Psi(O_\beta)$ under $\mathsf{P}$ is the image of $\wh{\Psi}$.
\item[(4)]\  There exists an open neighborhood $U$ of $|\beta|$ in $|\mathscr{B}|$
such that $|\Psi|:|G\ltimes O_\beta|\rightarrow U$ and $|\wh{\Psi}|:|E|\rightarrow |\mathsf{P}|^{-1}(U)$
are homeomorphisms.
\end{itemize}
\qed
\end{definition}
We note that passing to orbit spaces we obtain the commutative diagram
$$
\begin{CD}
|E| @>|\wh{\Psi}| >> |\mathsf{P}|^{-1}(U)\\
@V|\wh{\tau}| VV  @V |\mathsf{P}|VV\\
|G\ltimes O_\beta|@>|\Psi|>> U.
\end{CD}
$$
Here the horizontal arrows are homeomorphisms. The local singularity structure of $|\mathsf{P}|$ is
that given by $|\wh{\tau}|$.
\begin{remark}
Let us also note that $\Psi:G\ltimes O_\beta\rightarrow \mathscr{B}$ is a uniformizer in the usual sense,
but that $\wh{\Psi}$ is not.  In fact, the latter is defined on a special ep-groupoid.
Note, however that the restrictions of $\wh{\Psi}$ to $G_{\alpha}\ltimes O_\alpha$ for $\alpha\in \mathsf{P}^{-1}(\beta)$
are uniformizers.  Nevertheless we shall refer to $\wh{\Psi}$ also as a uniformizer.
\qed
\end{remark}

\begin{definition}\label{DEFF1775}\index{D- Proper covering uniformizer construction}
Let $\mathsf{P}:\mathscr{A}\rightarrow \mathscr{B}$ be a proper covering functor between two GCT's.
A {\bf proper covering uniformizer construction} for $\mathsf{P}$ consists of a functor $F:\mathscr{B}^-\rightarrow \text{SET}$ which associates to an object $\beta$ a set $F(\beta)$ of proper covering uniformizers, see Definition \ref{DEFR1774},
$$
\begin{CD}
E @>\wh{\Psi} >> \mathscr{A}\\
@V\wh{\tau} VV  @V \mathsf{P}VV\\
G\ltimes O_\beta@>\Psi>>  \mathscr{B}
\end{CD}
$$
with $\Psi(o_\beta)=\beta$. 
\qed
\end{definition}
Given a proper covering uniformizer construction for $\mathsf{P}:\mathscr{A}\rightarrow \mathscr{B}$ 
we  consider as usual the associated transition sets. In our case, since a uniformizer is a pair,
we obtain the diagram
\begin{eqnarray}\label{EQWE1714}
\bm{M}(\wh{\Psi},\wh{\Psi}')\xrightarrow{\mathsf{p}} \bm{M}(\Psi,\Psi').
\end{eqnarray}
Here the  set $\bm{M}(\wh{\Psi},\wh{\Psi}')$ consists of all tuples $(e,\phi,e')$, where $\phi$ is a morphism
in $\mathscr{A}$ between $\wh{\Psi}(e)$ and $\wh{\Psi}'(e')$
$$
\bm{M}(\wh{\Psi},\wh{\Psi}')=\left\{(e,\phi,e')\ |\ e\in E,\ e'\in E',\ \phi\in \text{mor}(\wh{\Psi}(e),\wh{\Psi}'(e'))\right\},
$$
the set $\bm{M}(\Psi,\Psi')$
is defined as usual, and the map $\mathsf{p}$ is defined as follows
\begin{eqnarray}
\mathsf{p}(e,\phi,e')=(\wh{\tau}(e),\mathsf{P}(\phi),\wh{\tau}'(e')).
\end{eqnarray}
\begin{definition}\label{DEFG1776}\index{D- Proper polyfold covering structure}
Let $\mathsf{P}:\mathscr{A}\rightarrow \mathscr{B}$ be a proper covering functor between two GCT's.
A {\bf proper  polyfold covering structure} for $\mathsf{P}$ consists of a pair $(F,\bm{M})$,
where $F$ is a proper covering uniformizer construction and $\bm{M}$ is a construction
of M-polyfold structures for the transition sets $\bm{M}(\Psi,\Psi')$ and $\bm{M}(\wh{\Psi},\wh{\Psi}')$
so that $\mathscr{A}$ and $\mathscr{B}$ for the induced structures become polyfolds and the map
$\mathsf{p}$ defined in (\ref{EQWE1714}) becomes a local sc-diffeomorphism.
A {\bf proper polyfold covering functor} consists of a  proper covering functor between two GCT's,
say $\mathsf{P}:\mathscr{A}\rightarrow \mathscr{B}$, equipped with a proper polyfold covering  structure.
\qed
\end{definition}
Next we shall study properties of $\mathsf{p}:\bm{M}(\wh{\Psi},\wh{\Psi}')\rightarrow\bm{M}(\Psi,\Psi')$ in more detail.
\begin{lemma}\label{LEMM1776}
Assume we are given a proper covering uniformizer construction $(F,\bm{M})$  for $\mathsf{P}:\mathscr{A}\rightarrow \mathscr{B}$. 
For uniformizer pairs $(\wh{\Psi},\Psi)$ and $(\wh{\Psi}',\Psi')$ and the associated transition sets, we  have a well-defined bijective map
\begin{eqnarray}\label{GEQ1714}
\bm{M}(\wh{\Psi},\wh{\Psi}')\rightarrow \bm{M}(\Psi,\Psi'){_{s}\times_{\wh{\tau}}}E
\end{eqnarray}
mapping $(e,\phi,e')$ to $(\mathsf{p}(e,\phi,e'),s(e,\phi,e'))= ((\wh{\tau}(e),\mathsf{P}(\phi),\wh{\tau}'(e')),e)$. 
Here $\wh{\tau}:E\rightarrow G_\beta\ltimes O_\beta$, with $\wh{\Psi}:E\rightarrow \mathscr{A}$ and $\Psi: G_\beta\ltimes O_\beta\rightarrow \mathscr{B}$, and similarly for $\wh{\tau}'$.
\end{lemma}
\begin{proof}
To see that the map defined in (\ref{GEQ1714}) is bijective
 let us first show that it is injective.
Assume that  $(e,\phi,e')$ and $(e_1,\phi_1,e'_1)$ are mapped to the same point, implying
$$
((\wh{\tau}(e),\mathsf{P}(\phi),\wh{\tau}'(e')),e)=((\wh{\tau}(e_1),\mathsf{P}(\phi_1),\wh{\tau}'(e_1')),e_1).
$$
This implies that $e=e_1$ and moreover
$$
s(\phi)=\wh{\Psi}(e)=\wh{\Psi}(e_1)=s(\phi_1)\ \ \text{and}\ \ 
 \mathsf{P}(\phi)=\mathsf{P}(\phi_1).
$$
Since $\mathsf{P}$ is a proper covering map between categories it follows from these latter identities
that $\phi=\phi_1$.  Since $\wh{\Psi}'$ is injective on objects it follows that  $e'=e'_1$ and
 injectivity is established.

In order to establish surjectivity assume that 
$$
((\wt{e},\wt{\phi},\wt{e}'),e)\in \bm{M}(\Psi,\Psi'){_{s}\times_{\wh{\tau}}}E
$$
 is given. Then  $\wh{\tau}(e)=\wt{e}$  and $\wt{\phi}:\Psi(\wt{e})\rightarrow \Psi'(\wt{e}')$ hold.
From Definition \ref{DEFF1775} we deduce that  
$$
\mathsf{P}(\wh{\Psi}(e)) =\Psi (\wh{\tau}(e))=\Psi(\wt{e}),
$$
implying $s(\wt{\phi}) = \Psi(\wt{e}) =\mathsf{P}(\wh{\Psi}(e))$.
 Since $\mathsf{P}$ is a proper covering functor, there exists a unique
morphism $\phi$ with $s(\phi)=\wh{\Psi}(e)$ so that $\mathsf{P}(\phi)=\wt{\phi}$. The object $t(\phi)$ satisfies
$$
\mathsf{P}(t(\phi))=t(\mathsf{P}(\phi))=t(\wt{\phi})= \Psi'(\wt{e}').
$$
In view of property (3) in Definition \ref{DEFR1774} 
$t(\phi)$ belongs to the image of the objects  $E'$ under $\wh{\Psi}'$ and consequently for 
a uniquely determined $e'\in E'$ with $\wh{\Psi}'(e') =t(\phi)$ we  obtain the object $(e,\phi,e')$
in $\bm{M}(\wh{\Psi},\wh{\Psi}')$ which satisfies 
$$
(e,\phi,e')\rightarrow (\mathsf{p}(e,\phi,e'), e)=((\wt{e},\wt{\phi},\wt{e}'),e).
$$
The surjectivity has been established.
\qed \end{proof}
In view of Lemma \ref{LEMM1776} the map
$$
\bm{M}(\wh{\Psi},\wh{\Psi}')\rightarrow \bm{M}(\Psi,\Psi'){_{s}\times_{\wh{\tau}}}E
$$
is a bijection. However, more is true.
\begin{proposition}
Assume that $\mathsf{P}:\mathscr{A}\rightarrow \mathscr{B}$ is a proper polyfold covering functor.
Then for two uniformizers the map
$$
\bm{M}(\wh{\Psi},\wh{\Psi}')\rightarrow \bm{M}(\Psi,\Psi'){_{s}\times_{\wh{\tau}}}E
$$
is an sc-diffeomorphism.
\end{proposition}\index{P- Pullback $\mathsf{P}^\ast\Theta$}
\begin{proof}
 From the sc-smoothness of $\mathsf{p}$ it follows that the above map is sc-smooth. Since
$\mathsf{p}$ is a local sc-diffeomorphism and the same holds for $\wh{\tau}$ it follows that our map has to be a local sc-diffeomorphism, which of course is together with the other properties only possible if the map is a global sc-diffeomorphism. 
\qed \end{proof}
Given a proper polyfold covering functor $\mathsf{P}:\mathscr{A}\rightarrow \mathscr{B}$ 
a {\bf covering family of uniformizers}\index{Covering family of uniformizers} ${(\wh{\Psi}_\lambda,\Psi_\lambda)}_{\lambda\in\Lambda}$ for $\mathsf{P}$
is a family such that ${(\Psi_\lambda)}_{\lambda\in\Lambda}$ is a covering family for $\mathscr{B}$.
Associated to a such a family we have the covering of ep-groupoids
$\wh{\tau}_\lambda:E_\lambda\rightarrow G_\lambda\ltimes O_\lambda$ fitting into the diagram
$$
\begin{CD}
E_\lambda@>\wh{\Psi}_\lambda>> \mathscr{A}\\
@V\wh{\tau}_\lambda VV @V \mathsf{P} VV\\
G_\lambda\ltimes O_\lambda@>\Psi_\lambda >>  \mathscr{B}.
\end{CD}
$$
Denote by $E$ the disjoint union of all the object M-polyfolds $E$  and similarly by $X$ the union of the $O_\lambda$, i.e.
$$
E=\bigsqcup_{\lambda\in\Lambda} E_\lambda\ \ \text{and}\ \ X=\bigsqcup_{\lambda\in\Lambda} O_\lambda.
$$
Then we obtain a surjective  local sc-diffeomorphism $\wh{\tau}:E\rightarrow X$ by setting
$$
\wh{\tau}(e_\lambda) = \wh{\tau}_\lambda(e_\lambda)\ \ \text{for}\ \ e_\lambda\in E_\lambda.
$$
We can turn $E$ and $X$ into ep-groupoids
and $\wh{\tau}$ into a proper covering functor of ep-groupoids by defining $\bm{E}$ and $\bm{X}$ as follows
$$
\bm{E}=\bigsqcup_{(\lambda,\lambda')\in\Lambda\times\Lambda} \bm{M}(\wh{\Psi}_\lambda,\wh{\Psi}_{\lambda'})\ \ \text{and}\ \ \bm{X}=\bigsqcup_{(\lambda,\lambda')\in\Lambda\times\Lambda} \bm{M}({\Psi}_\lambda,{\Psi}_{\lambda'}).
$$
Then, as previously seen $E\equiv (E,{\bf E})$ and $X\equiv(X,{\bf X})$ are ep-groupoids. One easily verifies 
that $\wh{\tau}$ extends naturally to a proper covering functor, denoted by  $\wh{\tau}$,  between these ep-groupoids.
Of course, a different covering family produces a different proper covering functor and the different choices
 are related by some kind Morita equivalence. We leave it to the reader to work out the straight-forward,
but somewhat lengthy details. We call $\wh{\tau}:E\rightarrow X$ an sc-smooth model for the proper polyfold covering functor
$\mathsf{P}:\mathscr{A}\rightarrow \mathscr{B}$. Of course, it comes with a natural pair of equivalences
\begin{eqnarray}\label{ERQN1779}
\begin{CD}
E @>\wh{\Gamma}>> \mathscr{A}\\
@V \wh{\tau} VV @V \mathsf{P} VV\\
X @>\Gamma>>  \mathscr{B}.
\end{CD}
\end{eqnarray}
Next we consider the behavior of branched ep$^+$-subcategories with respect to proper polyfold covering functors.
\begin{proposition}
Assume that $\mathsf{P}:\mathscr{A}\rightarrow \mathsf{B}$ is a proper polyfold covering functor
and $\Theta:\mathscr{B}\rightarrow {\mathbb Q}^+$ a branched ep$^+$-subcategory. 
View $\mathscr{A}$ and $\mathscr{B}$ with the induced structures as polyfolds. The pull-back functor
$$
\mathsf{P}^\ast\Theta :
\mathscr{A}\rightarrow {\mathbb Q}^+
$$
defined by $\mathsf{P}^\ast\Theta:=\Theta\circ \mathsf{P}$ is a branched ep$^+$-subcategory.
Moreover the following holds.
\begin{itemize}
\item[{\em(1)}]\  If $\Theta$ is of manifold type also $\mathsf{P}^\ast\Theta$ is of manifold type.
\item[{\em(2)}]\   If $\Theta$ is of orbifold type also $\mathsf{P}^\ast\Theta$ is of orbifold type.
\item[{\em (3)}]\   If $\Theta$ is tame also  $\mathsf{P}^\ast\Theta$ is tame.
\item[{\em (4)}]\   If $\Theta$ is compact also $\mathsf{P}^\ast\Theta$ is compact.
\end{itemize}
\end{proposition}
\begin{proof}
We take a sc-smooth  model $\wh{\tau}:E\rightarrow X$ as in (\ref{ERQN1779}). Then, by definition,
$\Theta$ is a branched ep$^+$-subcategory
if and only if $\Gamma^\ast\Theta$ is a branched ep$^+$-subgroupoid.  It is an easy exercise that the pull-back
of a branched ep$^+$-subgroupoid by a proper covering functor $\wh{\tau}:E\rightarrow X$ between ep-groupoids
is a branched ep$^+$-subgroupoid of $E$.  Consequently $\wh{\tau}^\ast(\Gamma^\ast\Theta)$ is a branched ep$^+$-subgroupoid of $E$, which implies that $\wh{\Gamma}_\ast(\wh{\tau}^\ast(\Gamma^\ast\Theta))$ is a branched ep$^+$-subcategory of $\mathscr{A}$.  From the identity
\begin{eqnarray*}
\mathsf{P}^\ast\Theta&=& \wh{\Gamma}_\ast (\wh{\tau}^\ast(\Gamma^\ast\Theta))
\end{eqnarray*}
we obtain the desired result that $\mathsf{P}^\ast\Theta$ is a branched ep$^+$-subcategory.
All the procedures in the proof of the previous result preserve the tameness property, which yields (3).
Assume that $\Theta$ only takes the values $\{0,1\}$. Then the same is true for $\mathsf{P}^\ast\Theta$.
This proves (2). 

 In order to show (1) assume that $\alpha$ and $\alpha'$ are two objects 
with $\Theta\circ\mathsf{P}(\alpha)=1$ and $\Theta\circ\mathsf{P}(\alpha')=1$. Consider 
two  isomorphisms $\phi,\phi':\alpha\rightarrow \alpha'$. Then  $\mathsf{P}(\phi), \mathsf{P}(\phi'):\mathsf{P}(\alpha)\rightarrow \mathsf{P}(\alpha')$ and it follows from the fact that $\Theta$ is of manifold-type that
$\mathsf{P}(\phi)=\mathsf{P}(\phi')$. Define the element  $\sigma:=\phi'\circ\phi^{-1}$ in $G_{\alpha'}$
and note that  $\mathsf{P}(\sigma)=1_{\mathsf{P}(\alpha')}$.
Since $\mathsf{P}$ is a proper covering functor the map
$$
\bm{A}\rightarrow \bm{B}{_{s}\times_{\mathsf{P}}}A:\psi\rightarrow (\mathsf{P}(\psi),s(\psi))
$$
is a bijection. As we have shown it holds that $\sigma$ is mapped as
$$
\sigma\rightarrow (\mathsf{P}(\sigma),s(\sigma))= (1_{\mathsf{P}(\alpha')},\alpha').
$$
Since $1_{\alpha'}$ is also mapped to $(1_{\mathsf{P}(\alpha')},\alpha')$ we conclude that
$\sigma=1_{\alpha'}$ or equivalently $\phi=\phi'$.
This proves (1).

By the assumption in (4)  $|\supp(\Theta)|$ is compact. Consider $|\supp(\Theta\circ\mathsf{P})|$
and take a sequence $(\wt{z}_k)$ in the latter and define $z_k=|\mathsf{P}|(\wt{z}_k)$. There is no loss of generality assuming the $z_k\rightarrow z\in |\supp(\Theta)|$. We take a proper covering uniformizer around $\beta$ with $|\beta|=z$,
 and we can represent the sequences, after perhaps passing to subsequences by $(o_k)\subset O_b$ with
$\Psi(o_\beta)=\beta$ and $o_k\rightarrow o_\beta$. Moreover, we find  $(e_k)\subset E$ with $\wh{\tau}(e_k)=o_k$. 
We may assume that $(e_k)$ lies in a connected component sc-diffeomorphic to $O_\beta$ from which
it follows that $(e_k)$ converges.
\qed \end{proof}

Finally, we shall briefly discuss the situation of strong bundles in the covering context.
Assume that  $\mathsf{P}:\mathscr{A}\rightarrow \mathscr{B}$
is a proper polyfold covering functor. We have to consider strong bundles over $\mathscr{A}$ and $\mathscr{B}$
and in order that they fit into our covering scheme we need to impose some additional structure.
As before we denote by $\text{Ban}$ the category whose objects are Banach spaces and the morphisms 
are topological linear isomorphisms.  We assume we are given functors
$$
\mu:\mathscr{B}\rightarrow \text{Ban}\ \ \text{and}\ \ \wt{\mu}:\mathscr{A}\rightarrow \text{Ban}
$$
together with a natural transformation $\gamma: \wt{\mu}\rightarrow \mu\circ\mathsf{P}$.  This implies 
that for every object $\alpha$ in $\mathscr{A}$ we are given Banach spaces $\wt{\mu}(\alpha)$ and
$\mu(\mathsf{P}(\alpha))$ together with a linear topological isomorphism
$$
\gamma_\alpha:\wt{\mu}(\alpha)\rightarrow \mu(\mathsf{P}(\alpha)),
$$
so that a morphism $\phi:\alpha\rightarrow \alpha'$ produces a topological linear isomorphism
$\wt{\mu}(\phi):\wt{\mu}(\alpha)\rightarrow \wt{\mu}(\alpha')$ and $\mu(\mathsf{P}(\phi)):\mu(\mathsf{P}(\alpha))\rightarrow
\mu(\mathsf{P}(\alpha'))$ fitting into the commutative diagram 
$$
\begin{CD}
\wt{\mu}(\alpha) @>\wt{\mu}(\phi)>> \wt{\mu}(\alpha')\\
@V\gamma_\alpha VV @V \gamma_{\alpha'} VV\\
\mu(\mathsf{P}(\alpha))@>\mu(\mathsf{P}(\phi))>> \mu(\mathsf{P}(\alpha')).
\end{CD}
$$
Associated to $\wt{\mu}$ we can define as before $\wt{\mathscr{W}}$ having as objects 
$(\alpha,\wt{w})$ defining the functor $\wt{P}:\wt{\mathscr{W}}\rightarrow \mathscr{A}$. The base is already equipped 
with a polyfold structure. Associated to $\mu $ we obtain $P:\mathscr{W}\rightarrow \mathscr{B}$, where again the base
has a polyfold structure.  We define $\wt{\mathsf{P}}:\wt{\mathscr{W}}\rightarrow \mathscr{W}$ on objects  by
$$
 (\alpha,\wt{e})\rightarrow (\mathsf{P}(\alpha),\gamma_\alpha(\wt{e})),
$$
and on morphisms by
$$
(\phi,\wt{w})\rightarrow (\mathsf{P}(\phi),\gamma_{s(\phi)}(\wt{w})).
$$
With these definitions we obtain the commutative functor diagram
\begin{eqnarray}\label{XXXX4567}
\begin{CD}
\wt{\mathscr{W}}@ >\wt{P}>>  \mathscr{A}\\
@V \wt{\mathsf{P}} VV      @V\mathsf{P}VV\\
\mathscr{W} @> P>> \mathscr{B},
\end{CD}
\end{eqnarray}
where the vertical arrows are proper covering functors, and $\mathsf{P}$ already is a proper polyfold covering functor.
One needs to introduce structures which turn the first vertical arrow into a proper strong bundle covering functor.
This necessary structure is indicated by the following diagram, where  in the bottom plane 
we have the diagram (\ref{XXXX4567})
$$
\begin{tikzcd}[row sep=scriptsize, column sep=scriptsize]
& \wt{W} \arrow[dl] \arrow[rr] \arrow[dd] & & E \arrow[dl] \arrow[dd] \\
W \arrow[rr, crossing over] \arrow[dd] & & X \\
& \wt{\mathscr{W}} \arrow[dl] \arrow[rr] & & \mathscr{A} \arrow[dl] \\
\mathscr{W} \arrow[rr] & &\mathscr{B} \arrow[from=uu, crossing over]\\
\end{tikzcd}
$$
and the vertical arrows stand for uniformizers. The  right-hand face depicts the proper polyfold covering structure 
where the uniformizers are defined on $E\rightarrow X$. In  the top plane we need to pick strong bundles
$\wt{W}\rightarrow X$ and $W\rightarrow X$ over-ep-groupoids and a strong bundle map which is also is  covering map,
but isomorphic between the fibers.  Since the data in the vertical face on the right is induced by data on the left
we only have to explain the requirements for the diagram
\begin{eqnarray}
\begin{CD}
\wt{W} @>>>\wt{\mathscr{W}}\\
@V F VV @V\wt{\mathsf{P}}VV\\
W@>>>\mathscr{W}.
\end{CD}
\end{eqnarray}
With the objects introduced previously we are given the GCT $\wt{\mathscr{W}}$ whose objects
are pairs $(\alpha,\wt{w})$ with $\alpha$ an object in $\mathscr{A}$ and $\wt{w}\in \wt{\mu}(\alpha)$.
The morphisms have the form $(\phi,\wt{w})$ with $\wt{w}\in \wt{\mu}(s(\phi) )$.  Similarly
$\mathscr{W}$ has as objects $(\beta,w)$ with $w\in \mu(\beta)$ and morphisms $(\psi,w)$
satisfying $w\in \mu(s(\psi))$.  The proper covering  functor
$$
\wt{\mathsf{P}}\colon \wt{\mathscr{W}}\rightarrow \mathscr{W}
$$
is defined on objects by $(\alpha,\wt{w})\rightarrow (\mathsf{P}(\alpha),\gamma_{\alpha}(\wt{w}))$ and on morphisms 
by $ (\phi,\wt{w})\rightarrow (\mathsf{P}(\phi),\gamma_{s(\phi)}(\wt{w}))$.  We note that for fixed $\alpha$ the map
$(\alpha,\wt{w})\rightarrow (\mathsf{P}(\alpha),\gamma_\alpha(\wt{w}))$ is a linear topological isomorphism 
between Banach spaces. We would like to model the situation with strong bundles over ep-groupoids
denoted by $\wt{W}$ and $W$. Hence $\wt{W}$ and $W$ are strong bundles over ep-groupoids and 
$A:\wt{W}\rightarrow W$ is a proper strong bundle covering functor in the sense of Definition \ref{proper_sb_covering},
where we suppress here the underlying base data. A uniformizer construction for our current situation 
gives for an object $\beta$ in the base of $\mathscr{W}$ a functor $\Psi:W\rightarrow \mathscr{W}$
which maps a uniquely determined object $(o_\beta,0)$ to $(\beta,0_{\mu(\beta)})$ and is a linear isomorphism
between fibers.  Similarly $\wt{\Psi}:\wt{W}\rightarrow \wt{\mathscr{W}}$. The other requirements are similar to those 
explained in the case of proper covering functors. In our context all occurring maps in the transition sets
have to be in addition strong bundle maps.   We leave the remaining  details to the reader.

\chapter{Fredholm Theory in Polyfolds}\label{CHAPX18}
This chapter is concerned with the sc-Fredholm theory, which is the main topic of this book. We have discussed sc-Fredholm section functors
in great detail in the context of strong bundles over ep-groupoids and we shall carry the ideas over to the categorical context.

We assume we are given a GCT $\mathscr{C}$ equipped with a polyfold structure $(F,\bm{M})$
and a functor $\mu:\mathscr{C}\rightarrow \text{Ban}$ defining the bundle $P:\mathscr{E}\rightarrow \mathscr{C}$.
Suppose  that $\mathscr{E}$ is equipped with a strong bundle structure $(\bar{F},\bm{M})$ inducing the 
previously defined polyfold structure on $\mathscr{C}$.  Recall that $|\mathscr{C}|$ and $|\mathscr{E}|$
have as consequence of the definition of  a GCT metrizable topologies and $|P|:|\mathscr{E}|\rightarrow |\mathscr{C}|$ is continuous.

We shall study sc-Fredholm section functors of $P$ and 
we shall describe how some of the material developed in the ep-groupoid context can be generalized to the categorical context.
Many of the concepts can be carried over just using an sc-smooth model $X_{\bm{\Psi}}$ and the natural equivalence 
of categories $\Gamma_{\bm{\Psi}}:X_{\bm{\Psi}}\rightarrow \mathscr{C}$ together with the bundle version
$$
\begin{CD}
W_{\bm{\Psi}} @>\bar{\Gamma}_{\bm{\Psi}}>>  \mathscr{E}\\
@V p VV @V P VV\\
X_{\bm{\Psi}}  @>\Gamma_{\bm{\Psi}}>> \mathscr{C}.
\end{CD}
$$
Alternatively we can study, more locally, the diagrams
$$
\begin{CD}
G\ltimes K @>\bar{\Psi}>>   \mathscr{E}\\
@V p VV    @V P VV\\
G\ltimes O @>\Psi >>   \mathscr{C},
\end{CD}
$$
where $\bar{\Psi}\in \bar{F}(\alpha)$ covering $\Psi\in F(\alpha)$.  If $f:\mathscr{C}\rightarrow \mathscr{E}$ is a section functor,
we obtain for every choice $\bar{\Psi}$ a local representative $f_{\bar{\Psi}}: O\rightarrow K$, which is $G$-equivariant.
The properties relevant for our studies are those which can be defined by imposing requirements on ${(f_{\bar{\Psi}})}_{\bar{\Psi}\in\bm{\bar{\Psi}}}$,
where $\bm{\bar{\Psi}}$ is a covering family, and where the definition does not depend on the covering family which has been chosen. 
The sc-Fredholm property is of this kind.

\section{Basic Concepts}
Assume that $P:\mathscr{E}\rightarrow \mathscr{C}$ is a strong bundle over a polyfold with the structure defined by $(\bar{F},\bm{M})$.
\begin{definition}\label{DEFG1811}
A section functor $f$ of $P:\mathscr{E}\rightarrow \mathscr{C}$ is called an {\bf sc-Fredholm section functor} provided for any object $\alpha$ and $\bar{\Psi}\in \bar{F}(\alpha)$
a strong bundle uniformizer
$\bar{\Psi}:G\ltimes K\rightarrow \mathscr{C}$ covering $\Psi$ the representative $f_{\bar{\Psi}}:O\rightarrow K$ is sc-Fredholm in the sense
of Definition \ref{DEFX3116}.
\qed
\end{definition}
\begin{remark}\index{R- On sc-Fredholm  section functors}
 (a) Assume we are given a covering family $\bm{\bar{\Psi}}$ for $P$, i.e. the underlying $\bm{\Psi}$ is a covering family for $\mathscr{C}$.
If $f$ is a section functor such that all $f_{\bar{\Psi}}$, $\bar{\Psi}\in \bm{\bar{\Psi}}$,  are sc-Fredholm sections it follows that for an arbitrary object $\alpha$ and $\bar{\Psi}\in \bar{F}(\alpha)$
the representative $f_{\bar{\Psi}}$ is sc-Fredholm.  This is an easy exercise. So the sc-Fredholm property is checkable with a covering family.
As a consequence the pull-back of $f$ by $\bar{\Gamma}_{\bm{\bar{\Psi}}}:W_{\bm{\bar{\Psi}}}\rightarrow \mathscr{E}$ denoted by
$f_{\bm{\bar{\Psi}}}$ is an sc-Fredholm section of a strong bundle over an ep-groupoid if and only if $f$ is sc-Fredholm.\par

\noindent(b)  If $f$ is an sc-Fredholm section functor it follows that $f$ is regularizing. Indeed, if $\alpha$ is an object
in $\mathscr{C}_m$ and $f(\alpha)\in \mathscr{E}_{m,m+1}$, then $\alpha\in \mathscr{C}_{m+1}$.\par

\noindent(c) A definition,  which is equivalent to Definition \ref{DEFG1811}, requires $f$ to be regularizing
as described in (b) and being sc-smooth, which can be checked by a covering family. Moreover, it stipulates that for every smooth $z\in |\mathscr{C}|$ there exists an object $\alpha$
with $|\alpha|=z$ and $\bar{\Psi}\in\bar{F}(\alpha)$, with underlying $\Psi$ satisfying $\Psi(\bar{q})=\alpha$,
so that $(f_{\bar{\Psi}},\bar{q})$ is an sc-Fredholm germ in the sense of Definition \ref{oi}.
\qed
\end{remark}
Assume that $f$ is an sc-Fredholm section functor and $\alpha$ an object satisfying $f(\alpha)=0$. The tangent $T_\alpha\mathscr{C}$ is an sc-Banach space
and there is a well-defined linearization $f'(\alpha):T_\alpha\mathscr{C}\rightarrow \mathscr{E}_\alpha$ which is an sc-Fredholm functor.
It is defined as follows.  Take $\bar{\Psi}\in\bar{F}(\alpha)$ covering $\Psi$ and consider a tangent vector $[\alpha,\Psi,h]$. The linearization is defined by
\begin{eqnarray}\label{linXXXc}
f'(\alpha)([\alpha,\Psi,h])= (\alpha,\bar{\Psi}(q_0)(f_{\bar{\Psi}}'(h))),
\end{eqnarray}
where $f_{\bar{\Psi}}'(q_0):T_{q_0}O\rightarrow K_{q_0}$ is the linearization of $f_{\bar{\Psi}}$ at $q_0$, which satisfies $\Psi(q_0)=\alpha$.
The definition does not depend on the choice of $\bar{\Psi}$.
\begin{definition}\label{QWDEF1812}\index{D- Linearization $f'(\alpha)$}
The sc-Fredholm operator $f'(\alpha):T_\alpha\mathscr{C}\rightarrow \mathscr{E}_\alpha$ which via (\ref{linXXXc}) 
is well-defined, is called the {\bf linearization} of the sc-Fredholm section functor $f$ at the solution object $\alpha$.\qed
\end{definition}
Note that we  can also define a set of linearizations at smooth objects, which are not necessarily solution objects.
For this we take a uniformizer $\bar{\Psi}$ at $\alpha$ and obtain $f_{\bar{\Psi}}:O\rightarrow K$.
Denote by   $q_0$ the point which satisfies $\Psi(q_0)=\alpha$. For every local sc$^+$-section $s$ defined 
near $\bar{q}$, satisfying $s(q_0)=f_{\bar{\Psi}}(q_0)$ we can take the linearization
$(f_{\bar{\Psi}}-s)'(q_0):T_{q_0}O\rightarrow K_{q_0}$. The collection of all
$$
f_s'(\alpha)([\alpha,\Psi,h])=(\alpha,\bar{\Psi}(q_0)((f_{\bar{\Psi}}-s)'(q_0)(h)))
$$
consists of linear sc-Fredholm operators,  which all differ by sc$^+$-operators and therefore 
have the same index. The collection is independent of the choice of $\bar{\Psi}$ at $\alpha$.
\begin{definition}\index{D- Set of lineraizations}
The {\bf set of linearizations} of the sc-Fredholm section functor $f$ of $P:\mathscr{E}\rightarrow \mathscr{C}$
at the smooth object $\alpha$ is denoted by $\text{Lin}(f,\alpha)$. The index $\text{ind}(f,\alpha)$
is the Fredholm index of any of its elements in $\text{Lin}(f,\alpha)$.
\qed
\end{definition}
We note that if $\alpha$ is a zero object the set of linearizations consists of one element, namely the one defined
in Definition \ref{QWDEF1812}.

If $\phi:\alpha\rightarrow \alpha'$ is a morphism between two smooth objects and $f(\alpha)=0$, then $f(\alpha')=0$ as well.  Moreover,
with $T\phi: T_\alpha\mathscr{C}\rightarrow T_{\alpha'}\mathscr{C}$ and $\wh{\phi}:=\mu(\phi):\mu(\alpha)\rightarrow \mu(\alpha')$ being the lift by $\mu:\mathscr{C}\rightarrow \text{Ban}$, it holds that 
$$
\wh{\phi}\circ f'(\alpha) = f'(\alpha')\circ T\phi.
$$
Given  $\phi:\alpha\rightarrow \alpha'$ between smooth objects we obtain an induced map
$$
\phi_\ast :\text{Lin}(f,\alpha)\rightarrow\text{Lin}(f,\alpha'): \phi_\ast(L)= \wh{\phi}\circ L\circ T\phi^{-1}.
$$
Clearly $\phi'_\ast\circ\phi_\ast=(\phi'\circ\phi)_\ast$ and $({1_\alpha})_\ast=Id$.  As we have seen in the discussion
of sc-Fredholm section functors, in the context of ep-groupoids,  the set of linearizations is important when
discussion orientation questions.

Let us state a version of the implicit function theorem in the categorical context, which is based 
on results on ep-groupoids. We denote for a polyfold $\mathscr{C}$ and an open subset $U$ of $|\mathscr{C}|$
by $\mathscr{C}_U$ the
full subcategory associated to objects with $|\alpha|\in U$. 
\begin{theorem}
Let $P:\mathscr{E}\rightarrow \mathscr{C}$ be a strong bundle over a tame polyfold and $f$ an sc-Fredholm section functor. Assume that the smooth object $\alpha$ satisfies $f(\alpha)=0$ and the linearization
$f'(\alpha): T_\alpha\mathscr{C}\rightarrow \mathscr{E}_\alpha$ is surjective and the kernel $\ker(f'(\alpha))$
is in good position to the boundary (see Remark \ref{REMARG238}). Then there exists an open neighborhood $U$ of $|\alpha|$ in $|\mathscr{C}|$ so that the following holds.
\begin{itemize}
\item[{\em(1)}]\   The functor 
$\Theta\colon \mathscr{C}_U\rightarrow \{0,1\}$ defined by 
$$
\Theta(\beta)=\left[\begin{array}{cc}
1&\ \text{if}\ \  f(\beta)=0\\
0&\ \text{if}\ \  f(\beta)\neq 0
\end{array}
\right.
$$ 
is a tame orbifold-type ep$^+$-subcategory. 
\item[{\em(2)}]\    For  every object $\beta$ in $\mathscr{C}_U$ with
$\Theta(\beta)=1$ the linearization $f'(\beta)$ is surjective and the kernel lies in good position to the boundary.
\end{itemize}
Note that $|\supp(\Theta)|$ is in a natural way a smooth orbifold with boundary with corners.
That means that the local model is smooth manifold with boundary with corners divided out by a finite 
group action.
\qed
\end{theorem}
\begin{proof}
The result is a consequence of Theorem  \ref{THM1524}, Theorem \ref{THMXXC1525}, and Proposition 
\ref{PROPY918}.
\qed \end{proof}
\begin{remark}\index{R- Philosophy of perturbations}
Define $\Lambda:\mathscr{E}\rightarrow \{0,1\}$ by associating to the zero vectors the value $1$
and to other vectors the value $0$. Then, as we shall see   later on, $\Lambda$  is an sc$^+$-multisection  functor
and $(f,\Lambda)$ will be in good position at $\alpha$. 
With this definition it holds that $\Theta =\Lambda\circ f$ and the good position of $(f,\Lambda)$ at $\alpha$
implies that $\Theta$ is an ep$^+$-subcategory on a suitable $\mathscr{C}_U$.
 As we already have seen in the ep-groupoid case
it holds more generally that for a pair $(f,\Lambda)$ in sufficiently generic position (to be made precise)
the associated $\Lambda\circ f$ is a branched  ep$^+$-subcategory.
\qed
\end{remark}
Next we import the notion of an auxiliary norm  from the ep-groupoid context into the categorical framework. We assume that we have a strong polyfold bundle $P:\mathscr{E}\rightarrow \mathscr{C}$.
Recall that $\mathscr{E}$ has double filtration $\mathscr{E}_{m,k}$ for $0\leq k\leq m+1$.
\begin{definition}\label{FGH1890}\index{D- Auxiliary norm}
Let $P:\mathscr{E}\rightarrow \mathscr{C}$ be a strong bundle over a polyfold. 
An {\bf auxiliary norm} is a functor $N:\mathscr{E}\rightarrow [0,\infty]$ having the following properties.
\begin{itemize}
\item[(1)] \  If $e\in \mathscr{E}\setminus\mathscr{E}_{0,1}$ then  $N(e)=+\infty$.
\item[(2)]  \   $N$ defines on each fiber of $\mathscr{E}_{0,1}$  a complete norm.
\item[(3)] \  Given a sequence $(|h_k|)\subset |\mathscr{C}_{0,1}|$ satisfying $|P(h_k)|\rightarrow  |\alpha|$ and $|N|(|h_k|)\rightarrow 0$ it holds  $|h_k|\rightarrow |0_\alpha|$ in $|\mathscr{E}_{0,1}|$. 
\end{itemize}
\qed
\end{definition}
An easy lemma whose proof is left to the reader gives an alternative description of an auxiliary norm.
\begin{lemma}
Assume that $P:\mathscr{E}\rightarrow \mathscr{C}$ is a strong bundle over a polyfold and $N:\mathscr{C}\rightarrow [0,\infty]$ a functor.
Then $N$ is an auxiliary norm if and only if for a covering subset of strong bundle uniformizers $\bm{\bar{\Psi}}$ the functor 
$N_{\bm{\bar{\Psi}}}=N\circ\bar{\Gamma}_{\bm{\bar{\Psi}}}: W_{\bm{\bar{\Psi}}}\rightarrow [0,\infty]$ is an auxiliary 
norm for the strong bundle over the ep-groupoid $W_{\bm{\bar{\Psi}}}\rightarrow X_{\bm{\Psi}}$.
\qed
\end{lemma}
In order to introduce the notion of a reflexive auxiliary norm we distinguish a particular class of strong polyfold bundles $P:\mathscr{E}\rightarrow \mathscr{C}$.
\begin{definition}
We say that $P$ has reflexive $(0,1)$-fibers provided for every object $\alpha$
the fiber $(P^{-1}(\alpha))_{(0,1)}$ is a reflexive Banach space.
\qed
\end{definition}
We note that this can be checked on the orbit space in the following sense.  
If $z\in |\mathscr{C}|$ take an $\alpha$ representing $z$. If the $(0,1)$-fiber over
$\alpha$ is reflexive the same  will hold for any other $\alpha'$ isomorphic to $\alpha$.
If $p: W\rightarrow X$ is a local model for $P$ 
$$
\begin{CD}
W @>\bar{\Gamma}>>  \mathscr{E}\\
@V p VV   @V P VV\\
X @>\Gamma>> \mathscr{C}
\end{CD}
$$
and $|\Gamma(x)|=|\alpha|$, it follows that ${(W_x)}_{(0,1)}$ is a reflexive Banach space.
With other words any sc-smooth model has reflexive $(0,1)$-fibers. Following Definition \ref{DEFP1228}
we introduce the notion of mixed convergence.
\begin{definition}\index{D- Mixed convergence}
Let $P:\mathscr{E}\rightarrow \mathscr{C}$ be a strong polyfold bundle with reflexive $(0,1)$-fibers.
We say a sequence $(y_k)\subset |\mathscr{E}_{(0,1)}|$ is {\bf mixed convergent} to $y\in |\mathscr{E}_{0,1}|$ provided
$|P|(y_k)|\rightarrow |P|(y)$ (on level $0$) and for a choice of uniformizer $\bar{\Psi}$ covering $\Psi$
at a representative $\alpha$ for $|P|(y)$, say
$$
\begin{CD}
G\ltimes K@>\bar{\Psi}>> \mathscr{E}\\
@Vp VV  @V P VV\\
G\ltimes O @>\Psi >>  \mathscr{C}
\end{CD}
$$
with $\Psi(x)=\alpha$, we can pick for large $k$ representatives $e_k\in K$ with $|\bar{\Psi}(e_k)|=y_k$
and $e$ with $p(e)=x$ and $|\bar{\Psi}(e)|=y$ such that $e_k\xrightarrow{m} e$ (mixed convergent, see Definition \ref{DEFP1228}).
\qed
\end{definition}
Having defined mixed convergence we introduce the notion of a reflexive auxiliary norm.
\begin{definition}\index{D- Reflexive auxiliary norm}
Let $P:\mathscr{E}\rightarrow \mathscr{C}$ be a strong polyfold bundle with reflexive $(0,1)$-fibers and $N$ an auxiliary norm for $P$.  We say that $N$ is a {\bf reflexive auxiliary norm}, if besides the usual properties
of an auxiliary norm given in Definition \ref{FGH1890}, the following property holds.
If $(y_k)\subset |\mathscr{E}_{(0,1)}|$ is mixed convergent to $|y|\in |\mathscr{E}_{0,1}|$, then for suitable representatives 
$(e_k)$ and $e$ for $(y_k)$ and $y$ it holds
$$
N(e)\leq \text{liminf}_{k\rightarrow\infty} N(e_k)
$$
Alternatively we could require $|N|(y)\leq \text{liminf}_{k\rightarrow\infty} |N|(y_k)$.
\qed
\end{definition}
We can transport all results about auxiliary norms in the ep-groupoid context over to the categorical 
case, where $P:\mathscr{E}\rightarrow \mathscr{C}$ is a strong polyfold bundle. For this 
we just take an sc-smooth local model
$$
\begin{CD} 
W @>\bar{\Gamma} >>  \mathscr{E}\\
@V p VV   @V PVV\\
X@>\Gamma>> \mathscr{C}
\end{CD}
$$
and note that given an auxiliary norm $n$ for $p$ the push forward $N:= \bar{\Gamma}_\ast n$
is an auxiliary norm for $P$.  If $n$ is reflexive so is $N$.  Since $|\mathscr{C}|$ is by definition paracompact and therefore metrizable it follows that $|X|$ is metrizable and therefore paracompact. The choice of local
sc-smooth model does not matter since the notions around auxiliary norms are well-behaved 
under generalized strong bundle isomorphisms in the ep-groupoid context.
Consequently we can apply
Theorem \ref{EXTTT} and obtain the existence of an reflexive auxiliary norm on $P$, provided it has
reflexive $(0,1)$-fibers. Applications of  Theorem \ref{EXTT} and Theorem \ref{THMOP12213} provide 
extension theorems. 
\begin{theorem}\label{THG18112}\index{T- Existence of auxiliary norms}
Let $P:\mathscr{E}\rightarrow \mathscr{C}$ be a strong polyfold bundle. Then the following holds.
\begin{itemize}
\item[{\em (1)}]\   $P$ admits an auxiliary norm.
\item[{\em(2)}] \ If $P$ has reflexive $(0,1)$-fibers, there exists for every auxiliary norm $N'$ for $P$ two reflexive 
auxiliary norms $N_1$ and $N_2$ satisfying $N_1\leq N'\leq N_2$.
\item[{\em (3)}]\  If $P$ is defined over a tame polyfold an auxiliary norm defined on $\mathscr{E}|\partial\mathscr{C}$
can be extended to $P$.  In the case that $P$ has reflexive $(0,1)$-fibers a reflexive auxiliary norm
on $\mathscr{E}|\partial\mathscr{C}$ can be extended as a reflexive auxiliary norm.
\end{itemize}
\qed
\end{theorem}
We leave it to the reader to transport other interesting results around auxiliary norms 
into the categorical context.

\section{Compactness Properties}
The starting point is a strong polyfold bundle $P:\mathscr{E}\rightarrow \mathscr{C}$
and an sc-Fredholm section functor $f$ of $P$. Given an sc-smooth model associated to a covering family $\bm{\bar{\Psi}}$ we obtain $p:W\rightarrow X$ and the sc-Fredholm section $f_{\bm{\bar{\Psi}}}$ which fit into the commutative diagram
$$
\begin{CD}
W @>\bar{\Gamma}>> \mathscr{E}\\
@A f_{\bm{\bar{\Psi}}} AA    @ A f AA\\
X @>\Gamma>>  \mathscr{C}.
\end{CD}
$$
Many compactness properties are statements which can be formulated in terms of orbit spaces. We note that $|\Gamma|$ and $|\bar{\Gamma}|$ are homeomorphisms
fitting into the commutative diagram
$$
\begin{CD}
|W |@>|\bar{\Gamma}|>> |\mathscr{E}|\\
@A |f_{\bm{\bar{\Psi}}}| AA    @ A |f| AA\\
|X| @<|\Gamma|^{-1}<<  |\mathscr{C}|
\end{CD}
$$
resulting in the relationship
\begin{eqnarray}
|f| =|\bar{\Gamma}|\circ |f_{\bm{\bar{\Psi}}}|\circ |\Gamma|^{-1}
\end{eqnarray}
which allows us to carry over compactness concepts discussed in the ep-groupoid context in Section 
\ref{SEC114}.  If $f$ is an sc-Fredholm section functor of $P$ we define the associated {\bf solution category}\index{Solution category $S_f$}
$S_f$ to be the full subcategory of $\mathscr{C}$ associated to objects $\alpha$ satisfying $f(\alpha)=0$. Since $f$ is regularizing a solution object is smooth.
\begin{definition} \index{D- Reflexive local compactness property, category context}
 We say that the sc-Fredholm section functor $f$ of $P$ defines a {\bf compact moduli space}\index{D- Compact moduli space}
provided the orbit space $|S_f|$ of the solution category is compact. We refer to $|S_f|$ as the {\bf coarse moduli space} associated to $f$.\index{D- Coarse moduli space}
\qed
\end{definition}
Having a compact moduli space for an sc-Fredholm section functor implies the seemingly stronger 
properness property, which follows from Theorem \ref{THMB1242}.
\begin{theorem}[Stability of compactness]\index{T- Stability of compactness}\label{THMX18.2.2}
Let $f$ be an sc-Fredholm section functor of the strong polyfold bundle $P:\mathscr{E}\rightarrow \mathscr{C}$ having a compact (coarse) moduli space $|S_f|$. Then given an auxiliary norm $N:\mathscr{E}\rightarrow [0,\infty]$
there exists an open neighborhood $U$ of $|S_f|$ in $|\mathscr{C}|$ such that the closure of the set
$\{z\in U\ | |N|\circ |f|(z)\leq 1\}$ (on level $0$) is compact.
\qed
\end{theorem}
This result is important, since as a consequence compactness is not being destroyed by small perturbations, and prompts the
following definition.
\begin{definition}\index{D- Controls compactness}
Let $P:\mathscr{E}\rightarrow \mathscr{S}$ be a strong polyfold bundle, $N:\mathscr{E}\rightarrow [0,+\infty]$ be an auxiliary norm,
and $f$ an sc-Fredholm section functor for which $|S_f|$ is compact. Given an open neighborhood $U$ f $|S_f|$ we say that
$(N,U)$ {\bf controls compactness} for $f$ provided  $\cl_{|\mathscr{C}|}(\{z\in U\ |\ |N|\circ |f|(z)\leq 1\})$ is compact.
\end{definition}

In many applications a stronger form of compactness holds with important additional properties, which were discussed in the ep-groupoid context in Section \ref{SEC114}. We first discuss a version of the reflexive local compactness property as given in Definition \ref{DEFR1243}. We note that Definition \ref{DEFR1243} could have formulated equivalently  in terms of orbit spaces, and the following definition is a straight forward generalization 
of such a reformulation. 
\begin{definition}
Let $P:\mathscr{E}\rightarrow \mathscr{C}$ be a strong polyfold bundle with reflexive $(0,1)$-fibers and 
$f$ an sc-Fredholm section functor. We say that $f$ has the (categorical) {\bf reflexive local compactness property}
provided for every reflexive auxiliary norm $N:\mathscr{C}\rightarrow [0,\infty]$ and every point $z\in |\mathscr{C}|$
there exists an open neighborhood $U(z)$ in $|\mathscr{C}|$ such that the closure of the 
set $\{y\in U(z)\ |\ |N|\circ |f|(y)\leq 1\}$ is compact in $|\mathscr{C}|$. 
\qed
\end{definition}
\begin{remark}\index{R- On the categorical reflexive local compactness property}
Note that $z\in |\mathscr{C}|$ is not(!) required to satisfy $|f|(z)=0$. 
For example if $z\in |\mathscr{C}|_0\setminus |\mathscr{C}|_1$ then $|f|(z)\not\in |\mathscr{E}_{0,1}|$ since otherwise, by the regularizing property, we conclude that $z\in |\mathscr{C}_1|$.
Consequently, we infer that $|N|\circ |f|(z)=+\infty$. Using the reflexive local compactness property we see
that  for  a suitable open neighborhood $U(z)$ it holds that $|N|\circ |f|(y)>1$ for $y\in U(z)$. 
Indeed, otherwise there exists a sequence $z_k\rightarrow z$ with $|N|\circ |f|(z_k)\leq 1$.
From this we deduce that $|N|\circ |f|(z)\leq \text{liminf}_{k\rightarrow\infty} |N|\circ |f|(z_k)\leq 1$ giving a contradiction.
Consequently  $\cl_{|\mathscr{C}|}(\{y\in U(z)\ |\ N\circ |f|(z)\leq 1\})=\emptyset$, which is compact.
\qed
\end{remark}

From Theorem \ref{THMXXX1246} we obtain the following extension theorem which frequently appears 
in inductive constructions for moduli spaces with some Floer-theoretic structure.
\begin{theorem}\label{THEMG1824}\index{T- Controlled extension of neighborhoods}
Assume that $P:\mathscr{E}\rightarrow \mathscr{C}$ is a strong polyfold bundle  over the tame
$\mathscr{C}$ with reflexive $(0,1)$-fibers
and $f$ an sc-Fredholm section functor.  Assume that $N:\mathscr{C}\rightarrow [0,\infty]$
is a reflexive auxiliary norm and $f$ has a compact moduli space and also the reflexive local compactness property.
Let $S_{\partial f}$ be the solution category associated to $\partial f:=f|\partial\mathscr{C}$ and $U_\partial\subset |\partial\mathscr{C}|$
an open neighborhood of the compact $|S_{\partial f}|$  in $|\partial\mathscr{C}|$ such that the closure 
of $\{z\in U_\partial\ |\ |N|\circ |f|(z)\leq 1\}$ is compact in $|\mathscr{C}|$. Then there exists an open neighborhood
$U$ of the coarse moduli space $|S_f|$ with the following properties.
\begin{itemize}
\item[{\em (1)}]\ \ $U\cap |\partial\mathscr{C}|=U_\partial$.
\item[{\em(2)}]\ \ The closure of $\{z\in U\ |\ |N|\circ |f|(z)\leq 1\}$ in $|\mathscr{C}|$ is compact, i.e. $(N,U)$ controls compactness for $f$.
\end{itemize}
\qed
\end{theorem}
\begin{remark}
 In applications a situation where $N$ at the beginning is only defined for $\mathscr{E}|\partial\mathscr{C}$ arises frequently. 
In this case we can make statements about compactness properties of $\partial f$ involving $N$.
In a next step one would apply Theorem \ref{THG18112} (3) to extend $N$ over $\mathscr{E}$.
Theorem \ref{THEMG1824} then allows to extend $U$ by keeping track of compactness properties.
\qed
\end{remark}

\section{Sc\texorpdfstring{$^+$}{qq}-Multisection Functors}
At this point we consider strong polyfold bundles $P:\mathscr{E}\rightarrow \mathscr{C}$.
We shall use auxiliary norms and if $P$  has reflexive $(0,1)$-fibers  we may assume the auxiliary norm
to be reflexive.
In this section we shall derive important results about sc$^+$-multisections. Since the relevant notions 
behave well under generalized isomorphisms all results which are available in the ep-groupoid
context are also available in the categorical context. We shall exploit this for 
the main notions and leave it to the reader to deal with the remaining concepts.

\begin{definition}\index{D- Sc$^+$-mulisection functor}
A {\bf sc$^+$-multisection functor} for $P:\mathscr{E}\rightarrow \mathscr{C}$ is a functor 
$\Lambda:\mathscr{E}\rightarrow {\mathbb Q}^+$ such that for a sc-smooth model
$p:W\rightarrow X$ and associated equivalences $\bar{\Gamma}$ and $\Gamma$ the functor 
$$
\bar{\Gamma}^\ast\Lambda := \Lambda\circ\bar{\Gamma}:W\rightarrow {\mathbb Q}^+
$$
is an sc$^+$-multisection functor in the ep-groupoid sense, see Definition \ref{sc+-section-functor}.
\qed
\end{definition}
\begin{remark}
The definition does not depend on the choice of sc-smooth model which follows 
from the behavior of sc$^+$-multisection functors in the ep-groupoid context under 
generalized strong bundle isomorphisms.
\qed
\end{remark}

Given an auxiliary norm $N:\mathscr{E}\rightarrow [0,\infty]$ we can describe the size
of an sc$^+$-multisection functor by the functor $N(\Lambda):\mathscr{C}\rightarrow [0,\infty)$ which we call the 
{\bf point-wise norm}\index{Point-wise norm of $\Lambda$} of $\Lambda$. It is defined for an object $\alpha$ by
$$
N(\Lambda)(\alpha):= \text{max}\{N(\alpha,e)\ |\ (\alpha,e)\in P^{-1}(\alpha),\ \Lambda(\alpha,e)>0\}.
$$
Of course, $|N(\Lambda)|:|\mathscr{C}|\rightarrow [0,\infty)$ provides the same information and 
we shall also refer to it as the point-wise norm. 

For applications the structurable sc$^+$-multsection functors are important and they are defined as follows.
\begin{definition}\label{KLOP000}\index{D- Structurable sc$^+$-multisection functor}
Let $P:\mathscr{E}\rightarrow \mathscr{C}$ be a strong polyfold bundle.  A  sc$^+$-multisection functor 
$\Lambda:\mathscr{E}\rightarrow {\mathbb Q}^+$  is said to be {\bf structurable} provided
for a suitable sc-smooth model $p:W\rightarrow X$ the pull-back $\bar{\Gamma}^\ast\Lambda$ is structurable.
\qed
\end{definition}
As we have seen in Section \ref{SECT134}, in the ep-groupoid context,  the notion of being structurable is being preserved by pull-backs and push-forwards by strong bundle maps and therefore also by
generalized strong bundle isomorphisms. Therefore it does not matter in Definition \ref{KLOP000}
which sc-smooth model is being taken.  If $\Lambda_i:\mathscr{E}\rightarrow {\mathbb Q}^+$, $i=1,2$,
we define the functor $\Lambda_1\oplus\Lambda_2$ by
$$
(\Lambda_1\oplus\Lambda_2)(\alpha,e)=\sum_{e'+e''=e} \Lambda_1(\alpha,e')\cdot\Lambda_2(\alpha,e''),
$$
where $e,e',e''\in \mu(\alpha)$. For a functor  $\beta:\mathscr{C}\rightarrow {\mathbb R}$ we
define $f\odot\Lambda:\mathscr{E}\rightarrow {\mathbb Q}^+$ by 
$$
(\beta\odot\Lambda)(\alpha, e):=\left[ \begin{array}{cc}
\Lambda(\alpha, \beta(\alpha)^{-1}\cdot e) & \text{if}\ \beta(\alpha)\neq 0\\
1 &\text{if}\ \beta(\alpha)=0,\ e=0\\
0&\text{if}\ \beta(\alpha)=0,\ e\neq 0
\end{array}
\right.
$$
From the results in Section \ref{SECT133} we obtain
the following theorem.
\begin{theorem}
Let $P:\mathscr{E}\rightarrow \mathscr{C}$ be a strong polyfold bundle. If
$\Lambda_1,\Lambda_2:\mathscr{E}\rightarrow {\mathbb Q}^+$ are structurable sc$^+$-multisection vectors
and $\beta:\mathscr{C}\rightarrow {\mathbb R}$ is an sc-smooth functor, then 
$$
\Lambda_1\oplus \Lambda_2\ \ \text{and}\ \ \beta\odot\Lambda_1
$$
are structurable sc$^+$-multisection functors.
\qed
\end{theorem}
Structurability is  preserved by pull-backs via proper covering functors. 
This gives the following result which is a consequence of Theorem \ref{OTHM1344}.
\begin{theorem}
Let  $\wt{\mathsf{P}}:\wt{\mathscr{W}}\rightarrow \mathscr{W}$ be an sc-smooth strong bundle proper covering functor between the strong polyfold bundles $\wt{P}:\wt{\mathscr{W}}\rightarrow \wt{\mathscr{C}}$
and $P:\mathscr{W}\rightarrow \mathscr{C}$, or more precisely
$$
\begin{CD}
\wt{\mathscr{W}}@>\wt{\mathsf{P}}>> \mathscr{W}\\
@V\wt{P}VV @V P VV\\
\wt{\mathscr{C}}@> \mathsf{P}>> \mathscr{C}.
\end{CD}
$$
If $\Lambda:\mathscr{W}\rightarrow {\mathbb Q}^+$ is a structurable sc$^+$-multisection functor,
then also $\wt{\mathsf{P}}^\ast\Lambda$ is a structurable sc$^+$-multisection functor.
\qed
\end{theorem}
At this point we have discussed the basic results about sc$^+$-multisections and specifically structurable ones.
In the next section we discuss the construction of sc$^+$-multisection functors as well as extension theorems.

\section{Constructions and Extensions}\label{SERC184}
In order to construct sc$^+$-mulisection functors we need  sc-smooth partitions of unity or at least sc-smooth bump functions.
\begin{definition}\index{D- Admitting sc-smooth partitions of unity}\index{D- Admitting sc-smooth bump functions}
We say that the polyfold $\mathscr{C}$ admits {\bf sc-smooth partitions of unity} or {\bf sc-smooth bump functions}
provided it holds for an sc-smooth model.  
\qed
\end{definition}
\begin{remark}\index{R- On sc-smooth partitions of unity}
The existence of sc-smooth partitions of unity or sc-smooth bump functions is being preserved
under generalized isomorphisms between ep-groupoids and therefore the specific sc-smooth model
occurring in the previous definition does not matter.  If the sc-smooth model $X$, which has a metrizable 
orbit space $|X|$, admits as a M-polyfold sc-smooth partitions of unity then it admits an sc-smooth 
partition of unity in the ep-groupoid sense, i.e. the maps are functors, see the discussion in Subsection \ref{SCPART}.
\qed
\end{remark}
For the next consideration we assume that we  are given a strong polyfold bundle $P:\mathscr{E}\rightarrow \mathscr{C}$ admitting sc-smooth bump functions.  Consider a smooth object  $\alpha$ with isotropy group $G$ and   let $\bar{\Psi}$ be a strong bundle
uniformizer in $\bar{F}(\alpha)$. This gives us the diagram
$$
\begin{CD}
G \ltimes K @>\bar{\Psi}>> \mathscr{E}\\
@V p VV   @V P VV\\
G\ltimes O @>\Psi>> \mathscr{C}.
\end{CD}
$$
Denote by $\bar{q}\in O$ the unique point satisfying $\Psi(\bar{q})=\alpha$.  Since $|\mathscr{C}|$ is metrizable we can pick an open neighborhood $V$ of $\bar{q}$ in $O$ which is invariant under the $G$-action such that $\cl_{|\mathscr{C}|}(|\Psi(V)|) \subset |\Psi(O)|$. Assume that $(\alpha,e)$ is a smooth vector in $\mathscr{E}$ and let $\bar{h}$ be the preimage
in $K$ under $\bar{\Psi}^{-1}$.
Using sc-smooth bump functions  we can construct an sc$^+$-section $s_1$ on $O$ which vanishes outside of $V$
and satisfies $s_1(\bar{q})=\bar{h}$. Applying the $g$-action we obtain $s_g:O\rightarrow K$ defined by
$$
s_g(g\ast q)=g\ast s_1(q)\ \ \text{for}\ \ q\in O.
$$
Then $s_{1_G}=s_1$ and ${(s_g)}_{g\in G}$ is a symmetric sc$^+$-section structure which vanishes
outside of $V$.
Define $\Lambda:\mathscr{C}\rightarrow {\mathbb Q}^+$ as follows. For $\bar{\Psi}(h)$ with $h\in K$
we define 
$$
\Lambda(\bar{\Psi}(h))=\frac{1}{|G|}\cdot \sharp{\{g\in G\ |\ s_g(p(h))=h\}}.
$$
If $(\beta,e)$ is an object in $\mathscr{E}$ and there exists $(\phi,e):(\beta,e)\rightarrow (t(\phi),\wh{\phi}(e))$
with $t(\phi)\in \Psi(O)$ we define with $\bar{\Psi}(h')=(t(\phi),\wh{\phi}(e))$
$$
\Lambda(\beta,e):=\Lambda(\bar{\Psi}(h')),
$$
where we note that the right-hand side has been previously defined.  
The definition also does not depend on the choice of $\phi$
as long as $t(\phi)\in\Psi(O)$. In the remaining cases
there does not exists $\phi:\beta\rightarrow t(\phi)$ with $t(\phi)\in \Psi(O)$. In this case we define
$$
\Lambda(\beta,e)=\left[\begin{array}{cc}
1 & \ \ \text{if}\ \ e=0\\
0&\ \ \text{otherwise}.
\end{array}
\right.
$$
\begin{proposition}
The previous construction defines an sc$^+$-multisection functor $\Lambda:\mathscr{E}\rightarrow {\mathbb Q}^+$ with $\Lambda(\alpha,e)>0$.  In a suitable sc-smooth model $p:W\rightarrow X$ the pull-back $\bar{\Gamma}^\ast
\Lambda$ is atomic and consequently structurable. Hence $\Lambda$ is a structurable sc$^+$-multisection functor 
for $P$.
\end{proposition}
\begin{proof}
The easy proof is left to the reader and follows immediately from the theory in Chapter \ref{CHAPS13}
discussing sc$^+$-multisection functors in the ep-groupoid context, specifically Theorem \ref{THME1353}.
\qed \end{proof}
The previous result shows how we can construct global structurable sc$^+$-multisec\-tions 
using one uniformizer.  Of course, we could apply the procedure at a finite number of objects 
producing $\Lambda_1,...,\Lambda_k$ which could be summed up to produce $\Lambda=\Lambda_1\oplus..\oplus\Lambda_k$.
More generally,  using an sc-smooth model 
$$
\begin{CD}
W @>\bar{\Gamma} >> \mathscr{E}\\
@V p VV @V P VV\\
X@>\Gamma>> \mathscr{C}
\end{CD}
$$
we can construct for $P:\mathscr{E}\rightarrow \mathscr{C}$ sc$^+$-multisection functors
with a variety of properties by just making the appropriate constructions for $P:W\rightarrow X$ in the ep-groupoid context and 
pushing the data forward by $\bar{\Gamma}$. Chapter \ref{CHAPS13} demonstrated in detail
the wealth of possibilities. We leave the details to the reader and only mention the following extension result
which is the categorical version of Theorem \ref{p-main-p} and can be proved by using an sc-smooth model.
\begin{theorem}
Let $P:\mathscr{E}\rightarrow \mathscr{C}$ be a strong polyfold bundle over a tame polyfold
and assume that $\mathscr{C}$ admits sc-smooth partitions of unity and $N:\mathscr{E}\rightarrow [0,\infty]$ is a given auxiliary norm. Assume that $\Lambda: \mathscr{E}|\partial \mathscr{C}\rightarrow {\mathbb Q}^+$
is a structurable sc$^+$-multisection functor with domain support $\text{dom-supp}(\Lambda)$.
Let $V$ be an open neighborhood of $|\text{dom-supp}(\Lambda)|\subset |\partial\mathscr{C}|$ 
in $|\mathscr{C}|$ and $\beta :\mathscr{C}\rightarrow [0,\infty)$ a continuous functor with 
$|\supp(\beta)|\subset V$ and satisfying
$$
N(\Lambda)(\alpha)< \beta(\alpha)\ \text{for}\ \alpha\ \text{belonging to}\ \text{dom-supp}(\Lambda).
$$
Then there exists a structurable sc$^+$-multisection functor $\Lambda':\mathscr{E}\rightarrow {\mathbb Q}^+$
such that 
\begin{itemize}
\item[{\em (1)}]\ $N(\Lambda')(\alpha)\leq \beta(\alpha)$ on $\mathscr{C}$.
\item[{\em(2)}]\   $|\text{dom-supp}(\Lambda')|\subset V$.
\item[{\em (3)}]\   $\Lambda'|(\mathscr{E}|\partial\mathscr{C})=\Lambda$.
\end{itemize}
\qed
\end{theorem}

\section{Orientations}
In order to use the moduli space $|S_f|$ associated to an sc-Fredholm section functor
to produce invariants we need to perturb it in general by a small generic sc$^+$-multisection functor 
$\Lambda$ to obtain a branched ep$^+$-subcategory $\Theta=\Lambda\circ f$ from which we are  going 
to extract  invariants. Unless $\Lambda$ is single-valued one cannot expect to obtain meaningful
invariants,  unless we also have  an orientation for $\Theta$. As discussed in the ep-groupoid case 
such orientations for $\Lambda\circ f$ are  naturally induced by an orientation for $f$.  In this section 
we import the ideas from the ep-groupoid setting to the categorical set-up, see Section \ref{SEC153}
for the ep-groupoid situation.

Define for a smooth object $\alpha$ the set $\wh{\text{Gr}}_F(\alpha)$ to consist of finite, non-negative rational,
formal
sums $\wh{\mathsf{L}}=\sum_{\wh{L}}\sigma_{\wh{L}}\cdot \wh{L}$, where $\wh{L}=(L,o)$ is an oriented sc-Fredholm operator, and 
$$
L\colon T_\alpha\mathscr{C}\rightarrow \mathscr{E}_\alpha,
$$
with $\mathscr{E}_\alpha=P^{-1}(\alpha)$. That means that only a finite number of the $\sigma_{\wh{L}}$ are non-zero,
and in this case positive rational numbers.
 The collection of all $(\alpha,\wh{\mathsf{L}})$
defines the objects of a category denoted by $\wh{\text{Gr}}_F(\mathscr{C})$ which fibers over
the full subcategory $\mathscr{C}_\infty$ associated to smooth objects in $\mathscr{C}$
$$
\wh{\pi}\colon \wh{\text{Gr}}_F(\mathscr{C})\rightarrow\mathscr{C}_\infty: (\alpha,\wh{\mathsf{L}})\rightarrow \alpha.
$$
The morphisms in $\wh{\text{Gr}}_F(\mathscr{C})$ are pairs $(\phi,\wh{\mathsf{L}})$,
where $\phi$ is a smooth morphism in $\mathscr{C}$ and $s(\phi)=\wh{\pi}(\alpha,\wh{\mathsf{L}})$.
Here $s(\phi,\wh{\mathsf{L}})=(s(\phi),\wh{\mathsf{L}})$ and $t(\phi,\wh{\mathsf{L}})= (t(\phi),\phi_\ast\wh{\mathsf{L}})$ where
$$
\phi_\ast\wh{\mathsf{L}}=\phi_\ast\left(\sum_{\wh{L}}\sigma_{\wh{L}}\wh{L}\right)=\sum_{\wh{L}}\sigma_{\wh{L}}\cdot
\phi_\ast\wh{L},
$$
with  $\phi_\ast(L,o)=(\phi_\ast L,\phi_\ast o)$ and  $\phi_\ast L =\mu(\phi)\circ L\circ T\phi^{-1}$,
so that
$$
(\phi,\wh{\mathsf{L}}):(s(\phi),\wh{\mathsf{L}})\rightarrow (t(\phi),\phi_\ast(\wh{\mathsf{L}})).
$$
In a next step we define the notion of an orientation for an sc-Fredholm section functor of a strong polyfold bundle
$P:\mathscr{E}\rightarrow \mathscr{C}$. For a smooth object $\alpha$ we can consider the convex set of linearizations
$\text{Lin}(f,\alpha)$ consisting of previously defined sc-Fredholm operators $L:T_\alpha\mathscr{C}\rightarrow \mathscr{E}_\alpha$. Just viewing the operators as classical linear Fredholm operators between the Banach spaces on level $0$
we obtain a contractible convex space of Fredholm operators so that the associated determinant bundle 
has two possible orientations.  Following  the ep-groupoid case, see Definition \ref{DEFNG1252}, and properly translated into our set-up we can build a category whose objects are $(\alpha,(\text{DET}(f,\alpha),\mathfrak{o}))$
where $\text{DET}(f,\alpha)$ stands for the topological line bundle $\text{DET}(f,\alpha)\rightarrow \text{Lin}(f,\alpha)$
and $\mathfrak{o}$ for an orientation of the latter.  If we have picked an orientation
$\mathfrak{o}$ the only other possible orientation is then $-\mathfrak{o}$.
The class of all such objects is denoted by ${\mathscr{O}}_f$ and it fibers over $\mathscr{C}_\infty$ via
$$
\sigma\colon {\mathscr{O}}_f\rightarrow \mathscr{C}_\infty:(\alpha,(\text{DET}(f,\alpha),\mathfrak{o}))\rightarrow \alpha.
$$
The morphisms have the form $(\phi,(\text{DET}(f,\alpha),\mathfrak{o}))$ with $s(\phi)=\alpha$,
where we view such a morphisms as 
$$
(\phi,(\text{DET}(f,\alpha),\mathfrak{o})):(\alpha, (\text{DET}(f,\alpha),\mathfrak{o}))\rightarrow
(t(\phi),  (\text{DET}(f,t(\phi)),\phi_\ast(\mathfrak{o}))).
$$
\begin{definition}
Let $P:\mathscr{E}\rightarrow \mathscr{C}$ be a strong polyfold bundle and $f$ an sc-Fredholm section functor.
The category $\sigma: \mathscr{O}_f\rightarrow \mathscr{C}_\infty$ fibering over the smooth object category 
is called the {\bf orientation category} associated to $f$.
\qed
\end{definition}
We note that a smooth morphism $\phi:\alpha\rightarrow \alpha'$ lifts to a bijection
$\phi_\ast$ between two-point sets
$$
 \{(\text{DET}(f,\alpha),\mathfrak{o}),(\text{DET}(f,\alpha),-\mathfrak{o})\}\rightarrow
\{(\text{DET}(f,\alpha'),\mathfrak{o}'),(\text{DET}(f,\alpha'),-\mathfrak{o}' )\}.
$$
If $\alpha$ has nontrivial isotropy  it can in principle happen that for a suitable
$g\in G_\alpha$ it holds that $g_\ast (\text{DET}(f,\alpha),\pm\mathfrak{o})=(\text{DET}(f,\alpha),\mp\mathfrak{o})$,
which obstructs orientability of $f$.

Next we introduce the notion of an orientation, which defined by several properties. We shall begin 
with the more algebraic requirements. An orientation, if it exists, associates in particular to a smooth object $\alpha$ one of the two possible orientations 
for $\text{DET}(f,\alpha)$, i.e. it is  a section functor $\mathsf{o}$ of $\sigma$ and can be written as 
$$
\alpha\rightarrow \mathsf{o}_\alpha = (\alpha,(\text{DET}(f,\alpha),\mathfrak{o}_\alpha)).
$$
Moreover, $\mathsf{o}$ is assumed to be a functor, so that the choice of orientations is compatible with
morphisms
$$
\mathsf{o}_{t(\phi)}=\phi_\ast\mathsf{o}_{s(\phi)}\ \ \text{for smooth morphisms}\ \phi.
$$
Besides these purely functorial properties one needs a version of local propagation, i.e. local continuity
as in the ep-groupoid case, see Section \ref{SECTX65} and specifically Definition \ref{SDEF6511}.
We need some preparation.  Given $\bar{\Psi}\in \bar{F}(\alpha)$ for a smooth object $\alpha$ we obtain the commutative 
functor diagram
$$
\begin{CD}
G_\alpha\ltimes K@>\bar{\Psi} >> \mathscr{E}\\
@Af_{\bar{\psi}}AA @A fAA\\
G_\alpha\ltimes O @ >\Psi>>  \mathscr{C}.
\end{CD}
$$
For a smooth element $q\in O$   the map
$$
\bar{\Psi}_\ast\colon {\mathcal L}(T_qO,K_q)\rightarrow {\mathcal L}(T_{\Psi(q)}\mathscr{C},\mathscr{E}_{\Psi(q)}):L\rightarrow \bar{\Psi}(q)\circ L\circ (T\Psi(q))^{-1}
$$
defines a bijection between sc-operators and preserves several types of maps.
For example sc-Fredholm operators are mapped to sc-Fredholm operators
and sc$^+$-operators are mapped to sc$^+$-operators.  Given an sc-Fredholm operator $L$ above $q$,
an orientation $\mathfrak{o}$ for $\det(L)$ is mapped canonically to an orientation $\bar{\Psi}(q)_\ast \mathfrak{o}$
 of $\det(\bar{\Psi}_\ast L)$ using the linear isomorphisms
$T\Psi(q)$ and $\bar{\Psi}(q)$ which establish  correspondences between the kernels and co-kernels.
Hence we have a push-forward operation for orientations which we denote by $\bar{\Psi}_\ast$ and its inverse is
the pull-back $\bar{\Psi}^\ast$.
As a consequence, given a section functor $\mathsf{o}$ for $\sigma:\mathscr{O}_f\rightarrow \mathscr{C}_\infty$ 
we can consider the associated $\mathsf{o}_{\bar{\Psi}} =\bar{\Psi}^\ast\mathsf{o}$, which associates to
$\text{DET}(f_{\bar{\Psi}},q)$ an orientation ${(\mathsf{o}_{\bar{\Psi}})}_q$.

Given a smooth morphism $\phi:\alpha\rightarrow \alpha'$, $\bar{\Psi}\in \bar{F}(\alpha)$ and $\bar{\Psi}'\in\bar{F}(\alpha')$
we obtain associated to $\Phi=(\bar{q},\phi,\bar{q}')$, where $\Psi(q)=\alpha$, $\Psi'(\bar{q}')=\alpha'$,  a local sc-diffeomorphism 
$\wh{\Phi}: (U(\bar{q}),\bar{q})\rightarrow (U(\bar{q}'),\bar{q}')$.  This local sc-diffeomorphism together 
with its lift $\bar{\Phi}$ to the bundles allows to pull-back the germ of $\mathsf{o}_{\bar{\Psi}'}$ near $\bar{q}'$ 
and we obtain near $\bar{q}$ the following identity of germs
$$
\bar{\Phi}^\ast \mathsf{o}_{\bar{\Psi}'} =\mathsf{o}_{\bar{\Psi}}\ \ \text{near}\ \ \bar{q}.
$$
The pull-back operation $\bar{\Phi}$ has clearly a continuity property. 
\begin{definition}\index{D- Continuity property of $\mathsf{o}$}
Let $P:\mathscr{E}\rightarrow \mathscr{C}$ be a strong polyfold bundle and $f$ an sc-Fredholm section functor.
Denote by $\sigma:\mathscr{O}_f\rightarrow \mathscr{C}_\infty$ the associated orientation category.
Let $\mathsf{o}$ be a section functor of $\sigma$. We say that $\mathsf{o}$ has the {\bf continuity property},
provided for every class $z\in |\mathscr{C}_\infty|$, representative $\alpha_z$ of $z$ and $[\bar{\Psi}:G\ltimes K\rightarrow \mathscr{C}]\in\bar{F}(\alpha_z)$, where $p:K\rightarrow O$, $\bar{o}\in O$ with $\Psi(\bar{o})=\alpha$, the following holds.
If $f_{\bar{\Psi}}:O\rightarrow K$ is the local representative of $f$ equipped with the orientations $\bar{\Gamma}^\ast\mathsf{o}$, then the orientations near $\bar{o}$ are related by continuation.
\qed
\end{definition}
\begin{remark}\index{R- On the continuity property of $\mathsf{o}$}
The above definition requires an orientation check for a set of objects $\{\alpha_z\ |\ z\in |\mathscr{C}_\infty|\}$,
where we pick for every $z$ a $\bar{\Psi}_z\in \bar{F}(\alpha_z)$. If we take now an arbitrary $\bar{\Psi}$
in some $\bar{F}(\alpha)$ the sc-smoothness of the associated strong transition bundle guarantees
that also with respect to $\bar{\Psi}$ we have local continuity. This follows from the previous discussion.
\qed
\end{remark}
\begin{definition}[Orientation for $f$]\index{D- Orientation for $f$}
Let $P:\mathscr{E}\rightarrow \mathscr{C}$ be a strong polyfold bundle and $f$ an sc-Fredholm section functor.
An {\bf orientation} for $f$ is given by a  section functor $\mathsf{o}$ of $\sigma:{\mathcal O}_f\rightarrow \mathscr{C}_\infty$  which has the continuity property. If such a section $\mathsf{o}$ exists we call $f$ {\bf orientable}.\index{D- Orientable sc-Fredholm section functor}
\qed
\end{definition}
\begin{remark}
We leave the verification of the following to the reader. If $f$ is an orientable sc-Fredholm section functor of $P$
then $\sigma:{\mathcal O}_f\rightarrow \mathscr{C}_\infty$, when passing to orbit spaces, defines 
a ${\mathbb Z}_2=\{-1,1\}$-bundle over $|\mathscr{C}_\infty|$ which we denote by 
$|\sigma|:|{\mathcal O}_f|\rightarrow |\mathscr{C}_\infty|$.  An orientation  $\mathsf{o}$ defines a continuous section $|\mathsf{o}|$ of $|\sigma|$.
Moreover, $f$ is orientable provided $|\sigma|:|{\mathcal O}_f|\rightarrow |\mathscr{C}_\infty|$ is a $(2:1)$-map and admits a continuous section.
See also Proposition \ref{PROPT1256} for the ep-groupoid case.
\qed
\end{remark}
Assume that $f$ is an sc-Fredholm section of the strong polyfold bundle $P:\mathscr{E}\rightarrow \mathscr{C}$.
Given two different covering families $\bm{\Psi}$ and $\bm{\Psi}'$ we obtain
$$
\begin{array}{cc}
\begin{CD}
W @>\bar{\Gamma}>> \mathscr{E}\\
@A f_{\bm{\bar{\Psi}}} AA @A f AA\\
X @>\Gamma>> \mathscr{C}
\end{CD}\ \ \ 
&\ \ \ 
\begin{CD}
W' @>\bar{\Gamma}'>> \mathscr{E}\\
@A \bar{f}_{\bm{\bar{\Psi}}'} AA @A f AA\\
X' @>\Gamma'>> \mathscr{C}.
\end{CD}
\end{array}
$$
The following holds true.
\begin{proposition} 
Assume that $f$ is an sc-Fredholm section functor of the strong polyfold bundle $P:\mathscr{E}\rightarrow \mathscr{C}$.
\begin{itemize}
\item[{\em (1)}] \  The orientability of $f$ is equivalent to the orientability of $f_{\bm{\bar{\Psi}}}$.
\item[{\em (2)}]\   An orientation for $f$ pulls-back to an orientation for $f_{\bm{\bar{\Psi}}}$ and an orientation
for $f_{\bm{\bar{\Psi}}'}$.
\item[{\em (3)}]\  Given orientations $\mathsf{o}$ for $f_{\bm{\bar{\Psi}}}$ and $\mathsf{o}'$
for $f_{\bm{\bar{\Psi}}'}$ they correspond to the same orientation for $f$ provided the natural generalized
strong bundle isomorphism $\bar{\mathfrak{f}}:W_{\bm{\bar{\Psi}}}\rightarrow W_{\bm{\bar{\Psi}}'}$
covering $\mathfrak{f}:X_{\bm{\Psi}}\rightarrow X_{\bm{\Psi}'}$  pulls back $\mathsf{o}'$ to $\mathsf{o}$.
\end{itemize}
\end{proposition}
\begin{proof}
The proof goes along the lines of the proof of Theorem \ref{THMS1257} and is left to the reader.
\qed \end{proof}

\section{Perturbation Theory}
We assume that $P:\mathscr{E}\rightarrow \mathscr{C}$ is a strong polyfold bundle over a tame $\mathscr{C}$, and $N:\mathscr{E}\rightarrow [0,\infty]$  an auxiliary norm. We study sc-Fredholm section functors $f$ of $P$ which have compact moduli spaces
$|S_f|$. We can fix an open neighborhood $U$ of $|S_f|$ so that 
$$
\cl_{|\mathscr{C}|}(\{z\in U\ |\ |N|\circ |f|(z)\leq 1 \})\ \ \text{is compact}.
$$
We shall use the latter fact to show that there are many structurable sc$^+$-multisection functors $\Lambda:\mathscr{E}\rightarrow {\mathbb Q}^+$, which are small with respect to $(N,U)$, which controls compactness,
so that $\Theta:\mathscr{C}\rightarrow {\mathbb Q}^+$ defined by
$$
\Theta(\alpha) =\Lambda\circ f(\alpha)
$$
is a  ep$^+$-subcategory with the standard properties:
\begin{itemize}
\item[(1)]\  $|\supp(\Theta)|\subset U$.
\item[(2)]\  $\Theta$ is compact.
\item[(3)]\ $f'(\alpha):T_\alpha\mathscr{C}\rightarrow P^{-1}(\alpha)$ is surjective for every $\alpha$ belonging to the support of $\Theta$.
\item[(4)]\ In addition we can force good boundary behavior, for example $\Theta$ is in general position
or in in good position, both implying that $\Theta$ is tame.
\end{itemize}
In applications, depending on the situation, we might be able to take a $\Lambda$ which only takes the values
$\{0,1\}$ to achieve the above.  In this case $|\supp(\Theta)|$ has naturally the structure of a compact orbifold with boundary and corners.  If we know in addition that between any two objects in the support of $\Theta$ there is at most one 
morphism it follows that $|\supp(\Theta)|$ has naturally the structure of a smooth compact manifold with boundary with corners. If in addition $f$ is oriented via $\mathsf{T}_f$ it follows that $|\supp(\Theta)|$ is naturally oriented either as an orbifold or manifold.
Of course, as in the ep-groupoid case, a general perturbation for which $\Theta$ is a branched ep$^+$-subcategory is always possible.
The following results are immediate corollaries of results proved in the ep-groupoid context.

\begin{theorem}
Assume that $P:\mathscr{E}\rightarrow \mathscr{C}$ is a strong polyfold bundle over a tame $\mathscr{C}$ admitting sc-smooth bump functions and $f$ is an sc-smooth section functor of $P$ with compact $|S_f|$. Assume that $N:\mathscr{E}\rightarrow [0,\infty]$
is an auxiliary norm and $U$ an open neighborhood of $|S_f|$ so that $(N,U)$ controls compactness. Let $h:\mathscr{C}\rightarrow (0,1]$
be a continuous functor. Then there exists an sc$^+$-multisection functor $\Lambda:\mathscr{C}\rightarrow {\mathbb Q}^+$ with the following properties.
\begin{itemize}
\item[{\em (1)}]\   $N(\Lambda)(\alpha) < h(\alpha)$ for all objects in $\mathscr{C}$.
\item[{\em (2)}]\   The domain support of $\Lambda$ is contained in $\mathscr{C}_U$.
\item[{\em(3)}]\   $\mathsf{T}_{(f,\Lambda)}$ is surjective for objects in $\supp(\Lambda\circ f)$.
\item[{\em (4)}]\   For an object in $\supp(\Lambda\circ f)$ the kernels of the operators $\mathsf{T}_{(f,\Lambda)}(\alpha)$
are in general position to $\partial\mathscr{C}$.
\end{itemize}
In particular $\Lambda\circ f:\mathscr{C}\rightarrow {\mathbb Q}^+$ is a compact, tame, branched ep$^+$-subcategory.
Further the perturbation $\Lambda$ can be taken structurable.
\end{theorem}
\begin{proof}
Follows directly from Theorem \ref{THM1536} by studying the question in a smooth model.
\qed \end{proof}
We can also easily bring Theorem \ref{THM1537} into the polyfold context. As an exercise the reader  might also prove a version of  Theorem \ref{THMX5.3.12} for the  ep-groupoid case, and its polyfold generalization. Another good and useful exercise is the generalization of Theorem \ref{MORSE-type} to an ep-groupoid as well as a categorical context.

\chapter{General Constructions}\label{CHAPTER19X}
In this chapter we shall describe several basic constructions which are very useful in applications, i.e. 
the construction of moduli spaces in concrete situations.  In  \cite{FH2} we shall describe new tools in the context
of symplectic geometry, which focus on developing a modular approach to the construction of moduli spaces in symplectic geometry. The procedure 
leads to polyfold structures on certain groupoidal categories $\mathscr{S}$, strong bundle structures for $P:{\mathcal E}\rightarrow \mathscr{S}$,
and a Fredholm theory for certain section functors of $P$. The theory in Chapters \ref{CHAPX17} and \ref{CHAPX18} then allows to construct moduli spaces.

\section{The Basic Constructions}\label{SECX19.1}
We start with a familiar definition.
\begin{definition}\index{D- Groupoidal category}\label{DEFNX19.1.1}
A category  $\mathscr{S}$ is said to be a {\bf groupoidal  category} provided the following holds.
\begin{itemize}
\item[(1)]\   Every morphism is an isomorphism.
\item[(2)]\   Between any two objects 
there are at most only a finite number of morphisms.
\item[(3)]\  The class of isomorphism classes of objects is a set.
\end{itemize}
\qed
\end{definition}
We shall obtain the polyfold structure on $\mathscr{S}$, and later on similarly the strong bundle structure 
on $P:{\mathcal E}\rightarrow \mathscr{S}$,  as a consequence of a particular kind of construction, which we shall introduce now.
This type of  construction will produce in the case of a groupoidal category $\mathscr{S}$ the following at the same time.
\begin{itemize}
\item[$\bullet$] A natural topology ${\mathcal T}$ on the orbit space $|\mathscr{S}|$. In case ${\mathcal T}$ is metrizable 
the pair $(\mathscr{S},{\mathcal T})$ is a GCT, see Definition \ref{DEF_GCT}.
 \item[$\bullet$] A  structure for the GCT $(\mathscr{S},{\mathcal T})$, which is a polyfold structure provided ${\mathcal T}$ is metrizable.
\end{itemize}

That ${\mathcal T}$ is metrizable has to be studied in any given context by an adhoc method. It will turn out that ${\mathcal T}$ is always locally metrizable, i.e.
every point in $|\mathscr{S}|$ will have an open neighborhood on which the topology is metrizable. This, of course, does not imply that the topology is Hausdorff.
There are criteria, which we shall
describe later,  which are easy to apply in practice, and which  guarantee that ${\mathcal T}$ is Hausdorff and completely regular. 
If one can show that in addition the topology is second countable, Urysohn's metrizability theorem implies that ${\mathcal T}$ is metrizable ({\bf Urysohn's metrization theorem}\index{Urysohn's metrization theorem}: A second countable, regular Hausdorff space is metrizable.).
In many applications this is not too difficult to show.  Alternatively,  if ${\mathcal T}$ is paracompact (rather than second countable) in addition to being
regular and Hausdorff,  the Nagata-Smirnov theorem implies metrizability as well. ({\bf Nagata-Smirnov metrization theorem}: A topological space is metrizable if and only if it is regular, Hausdorff and has a countably locally finite basis, i.e the topology has a basis which is a countable union of sets of open sets which are locally finite.  In particular, a paracompact regular Hausdorff space, which is locally metrizable is metrizable, see \cite{RS})\index{Nagata-Smirnov metrization theorem}.

We need the following definition.
\begin{definition}
Let $\mathscr{S}$ be a groupoidal category and $\alpha$ an object with isotropy group $G$.
We shall call a functor  
$$
\Psi:G\ltimes O\rightarrow \mathscr{S}
$$
a {\bf uniformizer}\index{D- Uniformizer  at $\alpha$} at $\alpha$ provided it has  the following properties.
\begin{itemize}
\item[(1)] \  $G\ltimes O$ is the translation groupoid consisting of the M-polyfold $O$ with  $G$ acting by sc-diffeomorphism.
\item[(2)] \  $\Psi$ is injective on objects and fully faithful.
\item[(3)] \ There exists $\bar{o}\in O$ with $\Psi(\bar{o})=\alpha$.
\end{itemize}
\qed
\end{definition}

We denote by $\mathscr{S}^-$ the category which has the same objects as $\mathscr{S}$,
but only has the identities as morphisms.
 The first part of the data is
given by a functor 
$$
F:\mathscr{S}^-\rightarrow \text{SET}
$$
 associating to every object $\alpha$ (with automorphism group $G$)
a set $F(\alpha)$ consisting of uniformizers at $\alpha$.
\begin{definition}\index{D- Uniformizer construction}
Given a groupoidal category $\mathscr{S}$ a functor $F:\mathscr{S}^-\rightarrow \text{SET}$,
which associates to an object $\alpha$ a set of uniformizers $F(\alpha)$ at $\alpha$ is called
a {\bf uniformizer construction}.
\qed
\end{definition}
\begin{remark} \index{R- Functoriality of uniformizer constructions}
In most applications $F$ is in fact a functor defined on $\mathscr{S}$. That means  that 
given a morphism $\Phi:\alpha\rightarrow \alpha'$ we obtain a bijection $F(\Phi):F(\alpha)\rightarrow F(\alpha')$.
As  part of the constructions in many cases  even more is true and  $F(\Phi)$ will have  a geometric meaning, see 
for example Remark \ref{REM1716}.
\qed
\end{remark}

Once $F$ is given, we can consider  the transition sets ${\bm{M}}(\Psi,\Psi')$ associated
to the functors $\Psi$ and $\Psi'$. 
For $\Psi\in F(\alpha)$ and $\Psi'\in F(\alpha')$ the  {\bf transition set}
${\bm{M}}(\Psi,\Psi')$\index{Transition set}  is given as the weak fibered product associated to the diagram
$$
O\xrightarrow{\Psi} \mathscr{S} \xleftarrow{\Psi'} O'.
$$
More precisely 
$$
{\bm{M}}(\Psi,\Psi')=\{(o,\Phi,o')\ |\ o\in O,\ o'\in O',\ \Phi:\Psi(o)\rightarrow\Psi'(o')\}.
$$
This construction comes with several maps called {\bf structure maps}\index{Structure maps}. Namely 
\begin{itemize}
\item $s:{\bm{M}}(\Psi,\Psi')\rightarrow O:(o,\Phi,o')\rightarrow o$, the {\bf source map}.
\item  $t:{\bm{M}}(\Psi,\Psi')\rightarrow O':(o,\Phi,o')\rightarrow o'$, the {\bf target map}.
\item  $u:O\rightarrow {\bm{M}}(\Psi,\Psi):o\rightarrow (o,1_{\Psi(o)},o)$, the {\bf unit map}.
\item $\iota:{\bm{M}}(\Psi,\Psi')\rightarrow {\bm{M}}(\Psi',\Psi):(o,\Phi,o')\rightarrow (o',\Phi^{-1},o')$, the {\bf inversion map}.
\item $m:{\bm{M}}(\Psi',\Psi''){_{s}\times_t}{\bm{M}}(\Psi,\Psi')\rightarrow {\bm{M}}(\Psi,\Psi'') :m((o',\Phi',o''),(o,\Phi,o'))=(o,\Phi'\circ\Psi,o'')$, the {\bf multiplication map}.
\end{itemize}
Given  $F$, we have for every object $\alpha$ a set of fully faithful functors $\Psi:G\ltimes O\rightarrow \mathscr{S}$ 
 having $\alpha$ in the image. Of course,
there is at this point no compatibility requirement between different uniformizers. Also note that we cannot talk about continuity properties
of $|\Psi|$ since $|\mathscr{S}|$ is not even equipped with a topology yet.

The second set of data, denoted by ${\mathcal F}$, consists of  a construction involving the transition sets.
It  associates to $h=(o,\Phi,o')\in {\bm{M}}(\Psi,\Psi')$ a germ of map 
$$
F_h:{\mathcal O}(O,o)\rightarrow ({\bm{M}}(\Psi,\Psi'),h)
$$
having the following properties, were we abbreviate $f_h:= t\circ F_h$, which is a germ of map $f_h:{\mathcal O}(O,o)\rightarrow {\mathcal O}(O',o')$.
We can think of $F_h$ as a map defined on an open neighborhood of $o$.

\noindent{\bf (A) Diffeomorphism Property:}  The germs $f_h:{\mathcal O}(O,o)\rightarrow {\mathcal O}(O',o')$ are local sc-diffeomorphisms, and 
$s(F_h(q))=q$ for $q$ near $o$. If $\Psi=\Psi'$ and $h=(o,\Psi(g,o),g\ast o)$ then $F_h(q) = (q,\Psi(g,q),g\ast q)$ for $q$ near $o$, so that $f_h(q)=g\ast q$.

\noindent{\bf (B) Stability Property:}  $F_{F_{h}(q)}(p)= F_h(p)$ for $q$ near $o=s(h)$ and $p$ near $q$.

\noindent{\bf (C) Identity Property:}  $F_{u(o)}(q) =u(q)$ for $q$ near $o$.

\noindent{\bf (D) Inversion Property:} $F_{\iota(h)}(f_h(q)) = \iota(F_h(q))$ for $q$ near $o=s(h)$.

\noindent{\bf (E) Multiplication  Property:}   If $s(h')=t(h)$ then $f_{h'}\circ f_h(q)=f_{m(h',h)}(q)$ for $q$ near $o=s(h)$ and 
 $m(F_{h'}(f_h(q)),F_h(q))=F_{m(h',h)}(q)$ for $q$ near $o=s(h)$.

\noindent{\bf (F) $\bm{M}$-Hausdorff Property:}  For different $h_1,h_2\in {\bm{M}}(\Psi,\Psi')$ with $o=s(h_1)=s(h_2)$  the images under $F_{h_1}$ and $F_{h_2}$ of small neighborhoods $o$ are disjoint.

\begin{definition}\label{DEFNX19.1.5}
The construction ${\mathcal F}$ is called a {\bf transition construction}\index{Transition construction}  associated to the uniformizer construction $F$.
\qed
\end{definition}

Usually, in applications,  the constructions $F$ and  ${\mathcal F}$ are natural and the listed properties are consequences of this naturality.
\begin{definition}\label{DEF1916}\index{D- Basic construction}
Given a groupoidal category $\mathscr{S}$ we shall refer to $(F,{\mathcal F})$, where the  functor $F:\mathscr{S}^-\rightarrow\text{SET}$ is a uniformizer
construction and ${\mathcal F}$ is an associated transition construction,  \index{Transition construction} 
 as a {\bf basic construction}\index{Basic construction}.
 \qed
\end{definition}
\begin{remark}\index{R- Remarks on transition constructions}
The Properties (A) and  (D)  say that there exist for $h\in \bm{M}(\Psi,\Psi')$ with $h=(o,\Phi,o')$ open neighborhoods $V(o)$ and $V(o')$
with the property that $f_h:(V(o),o)\rightarrow (V(o'),o')$ is an sc-diffeomorphism, and $F_{\iota(h)}(f_h(q))=\iota(F_h(q))$ for $q\in V(o)$.
Of course if $\Psi=\Psi'$ the transition set $\bm{M}(\Psi,\Psi)$ has a special form and consists of all $(q,\Psi(g,q),g\ast q)$, with $q\in O$ and $g\in G$.
In this special case (A) requires that on a suitable open neighborhood $V(o)$ it holds  $f_{(o,\Psi(g,o),g\ast o)}(q)=g\ast q$ and $F_{(o,\Psi(g,o),g\ast o)}(q)=(q,\Psi(g,q),g\ast q)$.

The Property (C) stipulates the existence of $V(o)$ with $F_{u(o)}(q)=u(q)$ for $q\in V(o)$. From the multiplication property with $h'=(o',\Phi',o'')$ 
we have the existence of three open neighborhoods $V(o)$, $V(o')$, and $V(o'')$ such that $f_h:(V(o),o)\rightarrow (V(o'),o')$ and $f_{h'}:(V(o'),o')\rightarrow (V(o''),o'')$
are sc-diffeomorphisms  and the following holds $f_{h'}\circ f_h(q)=f_{m(h',h)}(q)$ for $q\in V(o)$ and $m(F_{h'}(f_h(q)),F_h(q))=F_{m(h',h)}(q)$ for $q\in V(o)$.
These are, given the underlying structures, natural properties.  All these properties persist if we pass to smaller neighborhoods.

The properties we have discussed  allow us to consider $F_h(q)$ as an element close to $h$ in $\bm{M}(\Psi,\Psi')$ provided $q$ is close to $s(h)$.
The Stability Property (B) then says that the map  $F_{h^\ast}$ for $h^\ast$ near $h$ is the same as $F_h$, i.e. a stability property of the constuction. More precisely
(B) guarantees given $h=(o,\Phi,o')\in \bm{M}(\Psi,\Psi')$ an open neighborhood, again denoted by $V(o)$,
such that $F_{F_h(q)}(p)=F_h(p)$ for $q\in V(o)$ provided $p\in V(o)$ and $p$ is close to $q$.

The Hausdorff property just says that two different elements $h_1$ and $h_2$ with the same source should be considered as far apart
given our notion of `nearness' coming from $F_h$.
\qed
\end{remark}

In the following we shall derive the important  properties of a basic construction. The two important facts are that a basic construction defines
a natural topology ${\mathcal T}$ on the orbit space $|\mathscr{S}|$. This topology in general does not need to be well-behaved, and its properties 
in concrete situation have to be investigated. Moreover a basic construction defines also a metrizable topology on every $\bm{M}(\Psi,\Psi')$.

\section{The Natural Topology \texorpdfstring{${\mathcal T}$}{T} for \texorpdfstring{$\vert\mathscr{S}\vert$}{S}}   \label{SECNX19.2}
The starting point is a groupoidal category $\mathscr{S}$ and a basic construction  $(F,{\mathcal F})$.
Associated to the basic construction  there is a natural topology ${\mathcal T}$ on $|\mathscr{S}|$, which we shall construct now.
\begin{definition}\label{DEF1921}
Assume $\mathscr{S}$ is a groupoidal category and $(F,{\mathcal F})$ a basic construction.
The collection ${\mathcal T}$ of subsets of $|\mathscr{S}|$ consists of all subsets $U$ of $|\mathscr{S}|$ so that for every point $z\in U$
there exists an object $\alpha_z$ with $|\alpha_z|=z$,  $\Psi\in F(\alpha_z)$, and an open neighborhood
$V(\bar{o})$, where $\Psi(\bar{o})=\alpha_z$, such that $|\Psi(V(\bar{o}))|\subset U$.
\qed
\end{definition}
The basic topological result is given by the following theorem.
\begin{theorem}\label{TOPOLOGY}
Assume that we are given for the groupoidal category $\mathscr{S}$ a basic construction $(F,{\mathcal F})$.
Then the following holds.
\begin{itemize}
\item[{\em (1)}]\  The set ${\mathcal T}$ defined in Definition \ref{DEF1921} defines a topology on $|\mathscr{S}|$.
\item[{\em (2)}]\   With $\Psi\in F(\alpha)$, say $\Psi:G\ltimes O\rightarrow \mathscr{S}$, and an open subset $V$ of $O$
the set $|\Psi(V)|$ is open in $|\mathscr{S}|$.
\item[{\em (3)}]\ With $\Psi$ as above $|\Psi|: {_G\backslash O} \rightarrow |\mathscr{S}|$ is a homeomorphism onto an open neighborhood of $|\alpha|$.
\item[{\em (4)}]\    For the topology ${\mathcal T}$ every point has a countable neighborhood basis.
\item[{\em (5)}]\   Every point in $|\mathscr{S}|$ has (for ${\mathcal T}$) an open neighborhood which is metrizable.
({\em (5)} implies {\em(4)}).

\end{itemize}
\end{theorem}
\begin{remark}\index{R- On Hausdorffness}
There is in general no reason that ${\mathcal T}$ is Hausdorff.  We shall address this question and the related 
metrizability   question for ${\mathcal T}$ later on.
\end{remark}
\begin{proof}
\noindent (1) By the definition  of ${\mathcal T}$ it is clear that $\emptyset$ and $|\mathscr{S}|$ belong to $\mathscr{S}$.
If ${(U_\lambda)}_{\lambda\in\Lambda}$ is a family of sets in ${\mathcal T}$ it follows immediately from the construction
that $U=\bigcup_{\lambda\in\Lambda} U_\lambda$ belongs to ${\mathcal T}$. 

Let $U_1,...,U_k$ be a finite
family of elements in ${\mathcal T}$ and $z\in U:=\bigcap_{i=1}^k U_i$. 
For $z\in U_i$ we find 
$\Psi_i:G_i\ltimes O_i\rightarrow \mathscr{S}$ with $\Psi_i(\bar{o}_i)=\alpha_i$ and $|\alpha_i|=z$ 
so that for a suitable open neighborhood $V_i'$ of $\bar{o}_i$ we have  
$$
|\Psi_i(V_i')|\subset U_i.
$$
 We find  $h_i=(\bar{o}_1,\Phi_i,\bar{o}_i)\in {\bm{M}}(\Psi_1,\Psi_i)$ and take the associated germs $f_{h_i}:{\mathcal O}(O_1,\bar{o}_1)\rightarrow {\mathcal O}(O_i,\bar{o}_i)$.
By replacing $V_1' $ by a possibly  smaller  open neighborhood $V_1$ of $\bar{o}_1$ we may assume that all $V_i:=f_{h_i}(V_1)$ are open and contained in $V_i'$. In addition we may assume that $F_{h_i}$ is defined on $V_1$. Since $f_{h_i}=t\circ F_{h_i}$
it follows that $|\Psi_i(f_{h_i}(v))|=|\Psi_1(v)|$ for $v\in V_1'$ and consequently
$|\Psi_1(V_1)|\subset U$, which completes the proof that ${\mathcal T}$ is a topology.\par

\noindent (2) Assume that $\Psi(\bar{o})= \alpha$.  Let  $V\subset O$ be open and consider
the set $|\Psi(V)|$.  Take $z'\in |\Psi(V)|$ and  an object $\alpha'$ with $|\alpha'|=z'$.
Pick $\Psi'\in F(\alpha')$ with  $\Psi'(\bar{o}')=\alpha'$. For a suitable $o\in V$ we find an element
$h=(\bar{o}',\Phi,o)\in{\bm{M}}(\Psi',\Psi)$. Take the associated germ $f_h$.  We find an open neighborhood
$V'$ of $\bar{o}'$ which is sc-diffeomorphically mapped onto an open neighborhood of $o$ which we may assume
to be contained in $V$. This shows that $z'\in |\Psi'(V')|\subset |\Psi(V)|$, using that $F_h(q')=(q',\Phi_h(q'), f_h(q'))$,
implying $|\Psi(f_h(q'))|=|\Psi'(q')|$. \par

\noindent (3)  In view of (2) the set $U=|\Psi|({_G\backslash} O)=|\Psi(O)|$ is an open neighborhood of $|\alpha|$.
Consider $|\Psi|:{_G\backslash} O\rightarrow U$. This map is a bijection since $\Psi$ is fully faithful.
We first show the continuity. Pick a point $G\ast o$ in ${_G\backslash} O$ and define $z'=|\Psi(G\ast o)|$.
Pick an open neighborhood $V$ around $z'$. We find $\alpha'$ representing $z'$ and $\Psi'\in F(\alpha')$ 
so that $\Psi'(\bar{o}')=\alpha'$ and $|\Psi'(W')|\subset V$ for a suitable open neighborhood $W'$ of $\bar{o}'$.
We find an element $h=(o,\Phi,\bar{o}')$ and take the corresponding $f_h$. For a sufficiently small
open neighborhood $W$ of $o$  the open image $f_h(W)$ is contained in $W'$. Moreover for $w\in W$
$$
s(F_h(w))=w\ \text{and}\ \ t(F_h(w))=f_h(w),
$$
where $F_h$ has the form  $F_h(w)=(w,\Phi_w,f_h(w))$. Therefore $\Phi_w:\Psi(w)\rightarrow \Psi'(f_h(w))$ and
$$
|\Psi(w)| = |\Psi'(f_h(w))|
$$
implying $|\Psi(W)|\subset |\Psi'(W')|$. Then $|\Psi(G\ast W)|\subset V$. This proves continuity.
Since for an open subset $A$ in $O$ the image $|\Psi(A)|$ is open we conclude that $|\Psi|$ is open
and this completes the proof of (3). \par

\noindent (4) For given $z\in |\mathscr{S}|$ pick an object $\alpha$ satisfying $z=|{\alpha}|$ and $\Psi\in F(\alpha)$.
Let $\Psi(\bar{o})=\alpha$ and note that $\bar{o}$ has a countable neighborhood basis in $O$, say
$$
 V_1\supset V_2\supset V_3\supset...
 $$
In view of (2)  every $|\Psi(V_i)|$ is an open neighborhood 
of $z$. Pick any open neighborhood $U$ of $z$. We have to show that for sufficiently large $i$ 
it holds $|\Psi(V_i)|\subset U$.  We find $\Psi'\in F(\alpha')$ with $z=|\alpha'|$ and an open neighborhood
$V'$ of $\bar{o}'$, where $\Psi'(\bar{o}')=\alpha'$, with $|\Psi'(V')|\subset U$.  Since we have a morphism $\Phi:\Psi(\bar{o})\rightarrow \Psi'(\bar{o}')$
can take $f_h$, where $h=(\bar{o},\Phi,\bar{o}')$. We find $i$ such that $f_h(V_i)\subset V'$. 
This proves that $z\in |\Psi(V_i)|\subset U$.\par

\noindent (5) For $z\in |\mathscr{S}|$ pick an object $\alpha$ with $|\alpha|=z$ and $\Psi\in F(\alpha)$.  Then 
$$
\Psi:G\ltimes O\rightarrow \mathscr{S}
$$
with $\Psi(\bar{o})=\alpha$.  Then $|\Psi|(|G\ltimes O|)|=|\Psi(O)|$ is an open neighborhood of $z$ in view of (3).
Since the M-polyfold $O$ is metrizable we can pick a metric $\bar{d}$, and average this metric with respect to $G$
to obtain an invariant metric $d$ for $O$. This induces a metric $\rho$ for ${_G\backslash} O$ from which we conclude that
the open neighborhood $|\Psi(O)|$ of $z$ is metrizable.
\qed \end{proof}
\begin{remark}\index{R- On topological questions}
Starting with a groupoidal category $\mathscr{S}$ and a basic construction $(F,{\mathcal F})$ we obtain a  natural 
locally metrizable topology ${\mathcal T}$ for the orbit space $|\mathscr{S}|$. 
If by an adhoc method one can verify that ${\mathcal T}$ is Hausdorff, regular, and second countable it follows
that $(\mathscr{S},{\mathcal T})$ is a GCT and the theory developed in Chapters   \ref{CHAPX17} and \ref{CHAPX18} is applicable.
The metrizability question will be investigated in more depth in Section \ref{SECX194}.
\qed
\end{remark}

\section{The Natural Topology for \texorpdfstring{${\mathbf M}(\Psi,\Psi')$}{MPP}}\label{SECXN19.3}
As before we assume that we are given a basic construction $(F,{\mathcal F})$ for the groupoidal category $\mathscr{S}$.
We have shown that $|\mathscr{S}|$ carries a natural topology ${\mathcal T}$, which is locally metrizable.
We shall  show next that the transition sets ${\bm{M}(\Psi,\Psi')}$ carry natural topologies as well.

\begin{definition}\index{D- Topology ${\mathcal T}_{\Psi,\Psi'}$}
Let ${\mathcal T}_{\Psi,\Psi'}$ be the set consisting of subsets $W$ of ${\bm{M}}(\Psi,\Psi')$
having the following property. Given $h=(o,\Phi,o') \in W$ there exists an open neighborhood
$V$ of $o$ in $O$ such that $F_h(V)\subset W$.
\qed
\end{definition}

The basic result is the following.
\begin{theorem}\label{bold_M}
Let $(F,{\mathcal F})$ be a basic construction for the groupoidal category $\mathscr{S}$. Then the set
 ${\mathcal T}_{\Psi,\Psi'}$ is a topology on ${\bm{M}}(\Psi,\Psi')$. This topology is metrizable. Furthermore
the following properties hold:
\begin{itemize}
\item[{\em (1)}] \  The source map $s:{\bm{M}}(\Psi,\Psi')\rightarrow O$ is  a local homeomorphism.
\item[{\em (2)}]\   The target map $t:{\bm{M}}(\Psi,\Psi')\rightarrow O'$ is  a local homeomorphism.
\item[{\em (3)}]\   The unit map $u:O\rightarrow {\bm{M}}(\Psi,\Psi)$ is continuous.
\item[{\em (4)}]\   The inversion map $\iota:{\bm{M}}(\Psi,\Psi')\rightarrow {\bm{M}}(\Psi',\Psi)$ is a homeomorphism.
\item[{\em (5)}]\ The multiplication map
$m\colon \bm{M}(\Psi',\Psi''){_{s}\times_t}\bm{M}(\Psi,\Psi')\rightarrow \bm{M}(\Psi,\Psi'')$
 is continuous.
\end{itemize}
\qed
\end{theorem}
The proof of Theorem \ref{bold_M} requires some preparation and will be a consequence of the following discussion.
\begin{proposition}
The collection of subsets ${\mathcal T}_{\Psi,\Psi'}$ of ${\bm{M}}(\Psi,\Psi')$ is a topology and the following holds.
\begin{itemize}
\item[{\em (1)}]\   Given $h=(o,\Phi,o')$ there exists an open neighborhood $V$ of $o$ in $O$ so that
for every open subset $W$ of $V$ the image $F_h(W)$ belongs to ${\mathcal T}_{\Psi,\Psi'}$.
\item[{\em (2)}]\   If $W$ is as in (1) the map $F_h:W\rightarrow {\bm{M}}(\Psi,\Psi')$ is continuous.
\item[{\em (3)}]\   The topology  ${\mathcal T}_{\Psi,\Psi'}$  is Hausdorff.
\end{itemize}
\end{proposition}
\begin{proof}
By construction $\emptyset$, ${\bm{M}}(\Psi,\Psi')$,  and any union of elements in ${\mathcal T}_{\Psi,\Psi'}$ belong to
${\mathcal T}_{\Psi,\Psi'}$. Let $W_1,...,W_k$ be a finite collection of elements in ${\mathcal T}_{\Psi,\Psi'}$.
Assume that $h=(o,\Phi,o')\in W:=\bigcap_{i=1}^k W_i$.  For every $i$ we find an open neighborhood 
$V_i$ of $o$ such that $F_h(V_i)\subset W_i$. Consequently $F_h(\bigcap_{i=1}^k V_i)\subset W$.
This proves the  assertion that ${\mathcal T}_{\Psi,\Psi'}$ is a topology on $\bm{M}(\Psi,\Psi')$.\par

\noindent (1)   We note that there exist open neighborhoods $V$ of $o$ and $V'$ of $o'$ 
so that $f_h:V\rightarrow V'$ is an sc-diffeomorphism mapping $o$ to $o'$ and in addition
$$
F_{F_h(q)}(p)=F_h(p)\ \ \text{for}\ q\in V,\ p\in V.
$$
We also have that 
$$
s(F_h(q))=q\ \ \text{and}\ \ t(F_h(q))=f_h(q)\ \ \text{for all}\ q\in V.
$$
Let $W\subset V$ be an open subset and consider $F_h(W)$.
Take any $h'\in F_h(W)$ and note that there is a unique way to write it as
$F_h(q)=h'$.
Recall that 
$$
F_{F_h(q)}(p)=F_h(p)\ \ \text{for all}\ p\ \text{near} \ q.
$$
Hence the image of a small open neighborhood of $q$ is contained in $F_h(W)$. This proves that $F_h(W)$ is open. \par

\noindent (2) We start with $F_h:V(o)\rightarrow \bm{M}(\Psi,\Psi')$, where $o=s(h)$, 
and consider for  $q\in W\subset V$ the element  $k=F_h(q)$. Assume an open neighborhood $P$ of $k$ is given. We can take 
a small open neighborhood  of $k$ of the form $F_{k}(U)$, where $U$ is a sufficiently small open neighborhood of $o'=s(k)$ such that 
$F_k(U)\subset P$ and 
 $F_k(p)=F_h(p)$ for $p$ close to $o'$. Hence $F_h(U)\subset P$ and this proves continuity.\par

\noindent (3) Assume next that $h_i=(o_i,\Phi_i,o_i')$, for $i=1,2$, are two different elements in ${\bm{M}}(\Psi,\Psi')$.
Using the $\bm{M}$-Hausdorff property  and the previous discussion we conclude that these points have disjoint open neighborhoods.
\qed \end{proof}

Given a basic construction for the groupoidal category $\mathscr{S}$ we have just shown that  every transition set 
$\bm{M}(\Psi,\Psi')$ carries a natural Hausdorff topology.
We shall show next that the structural maps have the required continuity/homeomorphism properties, which will complete 
the proof of Theorem \ref{bold_M}, with the exception that  the metrizability of the topology ${\mathcal T}_{\Psi,\Psi'}$ still has to be shown.
\begin{proposition}
Let $\alpha,\alpha'$ and $\alpha''$ be three objects and $\Psi\in F(\alpha)$, $\Psi'\in F(\alpha')$, and $\Psi''\in F(\alpha'')$.
The following holds.
\begin{itemize}
\item[{\em (1)}]\   The source map $s:{\bm{M}}(\Psi,\Psi')\rightarrow O$ is  a local homeomorphism.
\item[{\em(2)}]\   The target map $t:{\bm{M}}(\Psi,\Psi')\rightarrow O'$ is  a local homeomorphism.
\item[{\em (3)}]\   The unit map $u:O\rightarrow {\bm{M}}(\Psi,\Psi)$ is continuous.
\item[{\em(4)}]\ The inversion map $\iota:{\bm{M}}(\Psi,\Psi')\rightarrow {\bm{M}}(\Psi',\Psi)$ is a homeomorphism.
\item[{\em(5)}]\  The multiplication map $m$ is continuous.
\end{itemize}
\end{proposition}
\begin{remark}\index{R- On transitions}
\noindent (a) Assuming that the proposition holds we note that the unit map $u$ is a homeomorphism onto its image, which is an open subset 
of $\bm{M}(\Psi,\Psi)$. Indeed, the map $u:O\rightarrow \bm{M}(\Psi,\Psi)$ is an injection.
Pick $o\in O$ and recall that $F_{u(o)}$ maps a small open neighborhood of $o$ to an open neighborhood of $u(o)$.
Since $F_{u(o)}(q)=u(q)$ for $q$ near $o$ we infer that the image of $u$ is open. 
Then $s|u(O)\rightarrow O$ as  a bijective local homeomorphism is a homeomorphism and the inverse to $u:O\rightarrow u(O)$.\par

\noindent (b) The multiplication is a local homeomorphism. Indeed, assume that we are given $h'\in \bm{M}(\Psi',\Psi'')$ and $h\in \bm{M}(\Psi,\Psi')$,
with $s(h')=t(h)$ and define $h''=m(h',h)$.  We find an open neighborhoods $U(s(h))$ of $s(h)$, $U(t(h))$ of $t(h)$,   $U(h)$ of $h$, $U(h')$ of $h'$ so that 
the following maps are homeomorphisms
\begin{itemize}
\item[] $s:U(h)\rightarrow U(s(h))$, $t:U(h)\rightarrow U(t(h))$, $f_h:U(s(h))\rightarrow U(t(h))$.
\item[] $s:U(h')\rightarrow U(t(h))$.
\end{itemize}
Then the fibered product $ {\bm{M}}(\Psi',\Psi''){_{s}\times_t}{\bm{M}}(\Psi,\Psi')$ near $(h',h)$ is homeomorphically
parametrized by the map 
$$
q\rightarrow ((s|U(t(h)))^{-1}(f_h(q)),(s|U(s(h)))^{-1}(q)),
$$
 where $q\in U(s(h))$. Then
$$
s(m(((s|U(t(h)))^{-1}(f_h(q)),(s|U(s(h)))^{-1}(q))) =q
$$
which implies the assertion.
\qed
\end{remark}
\begin{proof}
(1) Let $h=(o,\Phi,o')\in {\bm{M}}(\Psi,\Psi')$. Then $s(h)=o$ and 
$$
s(F_h(q))=q\ \text{for}\ q\ \ \text{near}\ o.
$$ 
Given a sufficiently small open neighborhood $U$ of $o$ we know that $V=F_h(U)$ is an open neighborhood 
of $h$ and trivially $s(V)\subset U$. This proves continuity of $s$.  Next consider $s:F_h(U)\rightarrow U$. 
We also know that $F_h:U\rightarrow F_h(U)$ is injective in view of $s\circ F_h(q)=q$ for $q\in U$ and by construction
the map is surjective, and hence a bijection. Clearly $F_h:U\rightarrow F_h(U)$ is a continuous inverse for $s$. \par

\noindent (2) Since $t=s\circ \iota$ the assertion (2) follows from (4).\par

\noindent (3) Let $o\in O$ and $h=u(o)$. Take a sufficiently small open neighborhood $V$ around $o$.
Then $F_{u(o)}(V)=\{(q,1_{\Psi(q)},q)\ |\ q\in V\}$ is open and $u(V)=F_{u(o)}(V)$ implying continuity.\par

\noindent (4) Clearly $\iota$ is a bijection. If $h\in {\bm{M}}(\Psi,\Psi')$  let $h'=\iota(h)$. 
Pick a small neighborhood around $h'$. We may assume that it has the form $F_{h'}(V')$ with
$V'$ an open neighborhood of $o'=s(h')=t(h)$. If $V'$ is small enough we find $V=V(o)$ so that
$f_h:V\rightarrow V'$ is an sc-diffeomorphism.  Using that 
$$
\iota(F_h(q)) = F_{h'}(f_h(q))\ \ \text{for}\ \ q\ \text{near}\ o
$$
it follows that $\iota(F_h(V))= F_{h'}(V')$ showing continuity. Interchanging the role of $\Psi$ and $\Psi'$
we see that the inverse of $\iota$ is also continuous.\par

\noindent (5) Consider $m\colon \bm{M}(\Psi',\Psi''){_{s}\times_t}\bm{M}(\Psi,\Psi')\rightarrow \bm{M}(\Psi,\Psi'')$
and fix $h=(o,\Phi,o')$ and $h'=(o',\Phi',o'')$ so that $m(h',h)=(o,\Phi'\circ\Phi,o'')=:h''$. Assume we are given an open neighborhood of $h''$.
We may assume that it has the form $W'':=F_{h''}(U)$ for a suitable open neighborhood $U$ of $o$.  If we take $U$ small enough we may 
assume that $F_{h''}:U\rightarrow W''$ is a homeomorphism. We also may assume that 
$f_h:(U,o)\rightarrow (U',o')$ is a homeomorphism for a suitable open neighborhood $U'$ of $o'$ so that in addition
defining  $W= F_{h}(U)$ and $W'=F_{h'}(U)$ 
$$
F_h:U\rightarrow W \ \ \ \text{and}\ \ \ F_{h'}:U'\rightarrow W'
$$
are homeomorphisms. Given an  element 
$$
(k',k)\in W'{_{s}\times_t}W 
$$
we find a unique $u\in U$ such that $F_h(u)=k$ and $t(k)=t\circ F_h(u))=f_h(u)=s(k') = F_{h'}(f_h(u))$. This implies
$$
 k'\circ k = F_{h'}(f_h(u))\circ F_h(u)= F_{h''}(u)\in W.
$$ 
\qed \end{proof}

The desired metrizability property is proved in the next proposition.
\begin{proposition}
For every pair $\Psi\in F(\alpha)$ and $\Psi'\in F(\alpha')$ the topology ${\mathcal T}_{\Psi,\Psi'}$
on the transition set ${\bm{M}}(\Psi,\Psi')$ is metrizable.
\end{proposition}
\begin{proof}
First of all we note that ${\bm{M}}(\Psi,\Psi')$ is locally metrizable since it is locally homeomorphic to a metrizable space.
We have already proved that it is Hausdorff.  It suffices to show that the space is paracompact.
Then the Nagata-Smirnov metrization theorem implies that the space is metrizable.

Let ${(V_\lambda)}_{\lambda\in\Lambda}$ be an open covering of ${\bm{M}}(\Psi,\Psi')$. We have to show that it admits 
a subordinate locally finite open cover. Given an element $o\in O$,  there are at most $\sharp G'$-many 
$h\in {\bm{M}}(\Psi,\Psi')$ with $s(h)=o$.  Denote them by $h^o_1,..,h^o_{k_o}$. 
We find a sufficiently small open neighborhood $V(o)$ of $o$ and $U(h^o_i)$ so that
\begin{itemize}
\item[(i)] \ \ \ $s:U(h^o_i)\rightarrow V(o)$ is a homeomorphism.
\item[(ii)]\  \ \ The $U(h^o_i)$ are mutually disjoint.
\item[(iii)] \ \ \ For every $i$ there exists a $\lambda_i$ such that $U(h^o_i)\subset V_{\lambda_i}$.
\end{itemize}
Since $O$ is metrizable we find a locally finite subcover subordinate to $(V(o))_{o\in O}$.
We find a subset $K\subset O$ and nonempty open sets $W_o$ for $o\in K$ so that
$({W_o})_{o\in K}$ covers $O$, $W_o\subset V(o)$, and the cover is locally finite.
We do not require $W_o$ to be an open neighborhood of $o$.
We define 
$$
V(o,i):= (s|U(h^o_i))^{-1}(W_o)\ \text{for}\ o\in K,\ i\in \{1,..,k_o\}.
$$
Then $(V(o,i))$, which is indexed by $o\in K$ and $i\in \{1,...,k_o\}$,  consists of open sets and is subordinate to $(V_\lambda)$. If $h\in {\bm{M}}(\Psi,\Psi')$
consider $q=s(h)$. We find $o\in K$ such that $q\in W_o$ and $h$ belongs to one of the $V(o,i)$.
This shows that $(V(o,i))$ is an open covering of ${\bm{M}}(\Psi,\Psi')$ subordinate to the original $(V_\lambda)$.

Let us show that the covering is locally finite.  Pick $h\in {\bm{M}}(\Psi,\Psi')$ and define $q=s(h)$.
We find an open neighborhood $A$ of $q$ which only intersects finitely many of the $W_o$. 
Taking $A$ small enough we find an open neighborhood $B$ of $h$ so that $s:B\rightarrow A$ is a homeomorphism.
Assume $V(o,i)$ intersects $B$. Then $s(V(o,i))$ and $s(B)$ intersect, implying that $W_o$ and $A$ intersect.
This implies that there are only finitely many $V(o,i)$ intersecting $B$. This shows that ${\bm{M}}(\Psi,\Psi')$
is paracompact.
\qed \end{proof}
\begin{remark}
Given a groupoidal category $\mathscr{S}$ and a basic construction $(F,{\mathcal F})$ for $\mathscr{S}$
we have shown that the orbit space $|\mathscr{S}|$ carries a natural topology ${\mathcal T}$, which is locally metrizable,
and every transition set $\bm{M}(\Psi,\Psi')$ carries a natural metrizable topology, so that all structure maps 
are continuous, $s$ and  $t$ are local homeomorphisms, and $\iota$ is a homeomomorphism.
We investigate the metrizability properties for ${\mathcal T}$ in the next section.
\qed
\end{remark}

\section{Metrizability Criterion for \texorpdfstring{${\mathcal T}$}{TT}}\label{SECX194}
Next we consider  the metrizability question  for  ${\mathcal T}$ which we have mentioned before.
For this we introduce a useful criterion, called the  properness property.
We recall that the $\bm{M}(\Psi,\Psi')$ are metrizable topological spaces. 
\begin{definition}\label{DEF1941}\index{D- Properness property}
Assume that the groupoidal category $\mathscr{S}$ admits the basic construction  $(F,{\mathcal F})$.
\begin{itemize}
\item[(1)]\   We say that {\bf properness property (v1)} holds, 
provided for   $\Psi\in F(\alpha)$  and given $o\in O$ there exists an open neighborhood $U=U(o)$ in $O$
so that for any  $\Psi'\in F(\alpha')$, $\Psi':G'\ltimes O'\rightarrow \mathscr{S}$,  the metrizable space 
${\bm{M}}(\Psi,\Psi')$ has the property  that every sequence $(h_k)$ in ${\bm{M}}(\Psi,\Psi')$ with $s(h_k)\in U$ and $t(h_k)\in O'$ convergent in $O'$
has a convergent subsequence in $\bm{M}(\Psi,\Psi')$.
\item[(2)]\  We say that {\bf properness property (v2)} holds, provided for every uniformizer $\Psi$ the following holds with respect to ${\mathcal T}$.
Given $o\in O$ there exists an open neighborhood  $U=U(o)$ in $O$
so that $\cl_{|\mathscr{S}|}(|\Psi(U(o))|)\subset |\Psi(O)|$. 
\item[(3)]\  We say that {\bf properness property (v3)} holds provided for every $\Psi:G\ltimes O\rightarrow \mathscr{S}$
the following property holds true.  Given $o\in O$ there exists an open neighborhood $U=U(o)\subset O$ such that 
for a sequence $(q_k)\subset U(o)$ the convergence of $|\Psi(q_k)|$ in $(|\mathscr{S}|,{\mathcal T})$ to some 
$z$ implies that a suitable subsequence of $(q_k)$ converges to an element in $\cl_O(U(o))\subset O$.
\end{itemize}
The following proposition shows that the three formulations of a properness property are equivalent.
In the future, if $(F,{\mathcal F})$, satisfies one of the properties (and therefore all three) we shall say that the {\bf properness property} holds.
\qed
\end{definition}

\begin{proposition}\index{P- Equivalence of properness criteria}
Assume that  $(F,{\mathcal F})$ is a basic construction for the groupoidal category $\mathscr{S}$ and let ${\mathcal T}$ be the associated
topology for $|\mathscr{S}|$.  Then the three properness properties formulated in Definition \ref{DEF1941} are equivalent.
\end{proposition}
\begin{proof}
\noindent (1) implies (2):  By assumption, given $o\in O$ there exists an open neighborhood
$U(o)\subset O$ which satisfies the properness property (v1).  We shall show that $U(o)$ satisfies (v2).
Assuming  that (2) does not hold,  arguing indirectly  we may assume that $\cl_{|\mathscr{S}|}(|\Psi(U(o))|)\not\subset |\Psi(O)|$.
This implies that there exists a sequence $(q_k)\subset U(o)$ such that $|\Psi(q_k)|\rightarrow z\in |\mathscr{S}|$,
where $z\not\in |\Psi(O)|$. We pick $\alpha'$ with $|\alpha'|=z$ and $\Psi'\in F(\alpha')$,  and find for large 
$k$ an element $h_k\in \bm{M}(\Psi,\Psi')$ with $s(h_k)=q_k$ and $t(h_k)=p_k\rightarrow
\bar{o}'$ with $\Psi'(\bar{o}')=\alpha'$. In view of (1) we may assume after perhaps taking a subsequence that $h_k\rightarrow h$ in $\bm{M}(\Psi,\Psi')$. The element $h$ satisfies
$s(h)\in cl_O(U(o))$ and $t(h)=\bar{o}'$. Hence $z=|\alpha'|\in |\Psi(O)|$, which is a contradiction.\par

\noindent  (2) implies (3): We take the open neighborhood $U=U(o)$ guaranteed by (2) and show
that $U$ has the properties stated in (3). Pick any sequence $(q_k)\subset U$ and assume that $|\Psi(q_k)|\rightarrow z$  with respect to ${\mathcal T}$.
Then, in view of (2), we deduce that $z\in \cl_{|\mathscr{S}|}(\Psi(U)|)\subset |\Psi(O)|$. Hence we can find $\bar{q}\in O$ satisfying $|\Psi(\bar{q})|=z$.
Since $|\Psi|:|G\ltimes O|\rightarrow |\Psi(O)|$ is a homeomorphism we find a sequence $(g_k)\subset G$ satisfying
$$
g_k\ast q_k\rightarrow \bar{q}.
$$
After perhaps taking a subsequence we may assume that $(g_k)$ is a constant sequence associated to $g\in G$,  and therefore, without loss of generality
$$
q_k\rightarrow g^{-1}\ast \bar{q}.
$$
Hence  $(q_k)$ converges to an element in $\cl_O(U)\subset O$.\par

\noindent (3) implies (1):   Given $\Psi:G\ltimes O\rightarrow \mathscr{S}$ let $o\in O$ be a given point. 
By assumption we  find an open neighborhood $U=U(o)$ in $O$ such that (iii) holds. 
Assume $\Psi':G'\ltimes O'\rightarrow \mathscr{S}$ and $(h_k)\subset \bm{M}(\Psi,\Psi')$
is a sequence satisfying $q_k:=s(h_k)\in U$ and $(p_k:=t(h_k))$ is a convergent sequence in $O'$, say converging to $p$.
Then $|\Psi'(t(h_k))|$ converges to some element $z\in |\mathscr{S}|$. In view (3) we may assume without loss of generality 
that $q_k\rightarrow q $ in $O$, where $q\in \cl_O(U)$.  Hence $(h_k)\subset \bm{M}(\Psi,\Psi')$
has the property that $s(h_k)\rightarrow q\in O$ and $t(h_k)\rightarrow p\in O'$. 

Consider the set of morphisms 
in $\mathscr{S}$ having source $\Psi(q)$ and target in $\Psi'(O')$. This collection is a finite set and contains at most 
$\sharp G'$-many elements. Denote them by $\Phi_1,...,\Phi_\ell$. We find with $\bar{h}_i=(q,\Phi_i,{(\Psi')}^{-1}(t(\Phi_i)))$ open neighborhoods 
$U(q)$ and $U(\bar{h}_i)$ for $i=1,..,\ell$ such that
$$
s:U(\bar{h}_i)\rightarrow U(q)
$$
is a homeomorphism.
We know that for large $k$ the point $q_k=s(h_k)$ belongs to $U(q)$. Fix a sufficiently large  $k$ so that $s(h_k)=q_k\in U(q)$.
We fix $i$ and  we find an element  $\wt{h}_k\in  U(\bar{h}_i)$ satisfying 
$$
s(h_k) =s(\wt{h}_k).
$$
Write $\wt{h}_k= (q_k,\wt{\Phi}_k,r_k)$, where $\Psi'(r_k)=t(\wt{\Phi}_k)$. We note that $t(h_k)$ and $t(\wt{h})$ only differ by 
the action of a suitable $g'_k$, i.e. 
$$
h_k  =m((r_k,\Psi'(g'_k,r_k),g'_k\ast r_k), \wt{h}_k).
$$
After perhaps taking a subsequence we can replace the $g_k'$ by a constant sequence $g'$ and obtain
$$
h_k=m((r_k,\Psi'(g',r_k),g'\ast r_k), \wt{h}_k).
$$
Since $\wt{h}_k = (s|U(\bar{h}_i))^{-1}(q_k)\rightarrow (s|U(\bar{h}_i))^{-1}(q)$ and therefore $r_k\rightarrow r$ 
we find that $(r_k,\Psi'(g',r_k),g'\ast r_k)\rightarrow (r,\Psi'(g',r),g'\ast r)$. From the continuity of $m$ we conclude the convergence 
of $(h_k)$.
\qed \end{proof}
Next we consider the ramifications of the properness condition.
\begin{proposition}
Assume that $(F,{\mathcal F})$ is a basic construction for the groupoidal category $\mathscr{S}$,
and ${\mathcal T}$ is the natural topology for $|\mathscr{S}|$, which is locally metrizable.
 If the construction satisfies the properness property, then ${\mathcal T}$ is Hausdorff and completely regular.
 \end{proposition}
 \begin{proof}
 For  $|\alpha|\neq |\alpha'|$ we pick $\Psi\in F(\alpha)$ and $\Psi'\in F(\alpha')$
 so that $\Psi(\bar{o})=\alpha$ and $\Psi'(\bar{o}')=\alpha'$. Around the points 
 $\bar{o}$ and $\bar{o}'$ we pick monotonic neighborhood bases $U_k=U_k(\bar{o})$ and
 $U'_k = U'_k(\bar{o}')$. Then $(|\Psi(U_k)|)$ and $(|\Psi(U_k)|)$ are neighborhood bases
 for $|\alpha|$ and $|\alpha'|$, respectively.  We show that they do not intersect for large $k$.
 Arguing indirectly we find sequences
 $o_k\rightarrow \bar{o}$ with $o_k\in U_k$ and $o'_k\in U'_k$ with $o'_k\rightarrow \bar{o}'$
 such that $|\Psi(o_k)|=|\Psi'(o_k')|$. Hence we find $h_k\in \bm{M}(\Psi,\Psi')$
 with $s(h_k)=o_k$ and $t(h_k)=o'_k$.   Using the properness property we may
 assume without loss of generality that $h_k\rightarrow h$ which implies that $|\alpha|=|\alpha'|$
 which contradicts our initial assumption.  This proves that ${\mathcal T}$ is Hausdorff.
 
 Next we show that ${\mathcal T}$ is completely regular. For this take a point 
 $z$ and a closed subset $A$ of $|\mathscr{S}|$ not containing $z$.  We pick $\alpha$
 with $|\alpha|=z$ and $\Psi\in F(\alpha)$. Then $\Psi:G\ltimes O\rightarrow \mathscr{S}$
 and $\Psi(\bar{o})=\alpha$ for a suitable $\bar{o}$.  We can  pick an open $G$-invariant neighborhood $U$ of $\bar{o}$ such that $|\Psi(U)|\cap A=\emptyset$ and, moreover
 we can pick the $G$-invariant $V=V(\bar{o})$ with $\cl_O(V)\subset U$ so that
 $V$ satisfies $\cl_{|\mathscr{S}|}(|\Psi(V)|) \subset |\Psi(U)|$. 
 Then $\cl_O(V)\subset U\subset O$ and we can take a continuous map
 $\beta:O\rightarrow [0,1]$ with $\beta(\bar{o})=0$ and $\beta|(O\setminus V)\equiv 1$.
 Also, perhaps using an averaging procedure we map assume that $\beta$ is invariant 
 with respect to the $G$-action. Then $\beta$ descends to a continuous function
 $\bar{\beta}:|\Psi(O)|\rightarrow [0,1]$, which we extend by the value $1$. This map
 $\bar{\beta}:|\mathscr{S}|\rightarrow [0,1]$ is continuous and has the desired properties.
 Indeed, we can argue using sequences, since every point has a countable neighborhood basis. 
 Let $z\in |\mathscr{S}|$ and assume that $z_k\rightarrow z$. If $z\in |\Psi(O)|$ we see that for large $k$
 it holds that $z_k\in |\Psi(O)|$. We can pick $q_k,q\in O$ with $q_k\rightarrow q$, $|\Psi(q_k)|=z_k$, and $|\Psi(q)|=z$.
 Then $\bar{\beta}(z_k)-\bar{\beta}(z)=\beta(q_k)-\beta(q)\rightarrow 0$. If $z\not\in |\Psi(O)|$, we find an open neighborhood
 $W=W(z)$ such that $|\Psi(V)|\cap W=\emptyset$. This immediately implies the continuity assertion since for large $k$
 we must have $z_k\in W$ so that $\bar{\beta}(z_k)=1=\bar{\beta}(z)$.  If we cannot find such a $W$ there exists a sequence 
 $(q_k)\subset V$ with $|\Psi(q_k)|\rightarrow z$. Applying the properness property we may assume without loss of generality that
 $q_k\rightarrow q\in \cl_O(V)\subset U\subset O$ and therefore $z=|\Psi(q)|$ implying that $z\in |\Psi(O)|$ giving a contradiction.
 This shows that the topology ${\mathcal T}$ is completely regular.
  \qed \end{proof}

The main result in this section is given by the following theorem which summarizes the essence of the previous discussion.
Recall that the natural topology is always locally metrizable and the properness
assumption implies in addition that it is Hausdorff and completely regular.
\begin{theorem}\label{MET}
Let $\mathscr{S}$ be a groupoidal category and $(F,{\mathcal F})$ a basic construction
for which the properness assumption holds.
\begin{itemize}
\item[{\em (1)}]\  If  ${\mathcal T}$ is second countable, then it is metrizable.
\item[{\em (2)}]\  If ${\mathcal T}$ is paracompact, then it is metrizable.
\end{itemize}
\end{theorem}
\begin{proof}
\noindent (1)  By Urysohn's metrization theorem a topological space which is Hausdorff,
completely regular and second countable is metrizable.\par

\noindent (2) A  locally metrizable space, which in addition is Hausdorff, completely regular
and paracompact is metrizable by the Nagata-Smirnov metrization theorem.
\qed \end{proof}

Finally we note the following result.
\begin{proposition}\index{P- Metrizability characterization}
Let $\mathscr{S}$ be a groupoidal category and $(F,{\mathcal F})$ a basic construction.
Denote by ${\mathcal T}$ the natural topology on $|\mathscr{S}|$. Then the following statements are equivalent.
\begin{itemize}
\item[{\em(1)}] \  ${\mathcal T}$ is metrizable.
\item[{\em(2)}]\  $(F,{\mathcal F})$ satisfies the properness property and ${\mathcal T}$ is paracompact.
\end{itemize}
\end{proposition}
\begin{proof}
\noindent (2) implies (1):  This follows from the previous discussion.\par

\noindent (1) implies (2):  A metrizable space is paracompact. We have to show that the metrizability of ${\mathcal T}$ requires the basic construction 
to satisfy one of the properness conditions.  Let $\Psi:G\ltimes O\rightarrow \mathscr{S}$ be a uniformizer at $\alpha$. 
Denote by $\bar{o}\in O$ the element mapped to $\alpha$. We know that $|\Psi|: {_G\backslash}O\rightarrow |\Psi(O)|$ is a homeomorphism and $|\alpha|\in |\Psi(O)|$.
Since ${\mathcal T}$ is metrizable we can pick an open neighborhood $V$ of $|\alpha|$ such that $\cl_{|\mathscr{S}|}(V)\subset |\Psi(O)|$. 
Denote by $U\subset O$ the preimage of $V$ under $\Psi$, which is a $G$-invariant open subset.

We verify  that (v3) holds. Assume that $(q_k)\subset U$ and $|\Psi(q_k)|\rightarrow z$ in $|\mathscr{S}|$. By construction 
$|\Psi(q_k)|\in V$ so that $z\in \cl_{|\mathscr{S}|}(V)\subset |\Psi(O)|$. Hence we find $\bar{q}\in O$ with $|\Psi(\bar{q})|=z$.
Since $|\Psi|:{_G\backslash}O\rightarrow |\Psi(O)|$ is a homeomorphism we find a sequence $(g_k)\subset G$ such that
$$
g_k\ast q_k\rightarrow \bar{q}.
$$
After perhaps taking a subsequence we may assume that $(g_k)$ is a constant sequence and without loss of generality
$q_k\rightarrow g^{-1}\ast \bar{q}=:q$. This completes the proof.
\qed \end{proof}

\section{The Polyfold Structure for \texorpdfstring{$(\mathscr{S},{\mathcal T})$}{ST}}\label{SECNX19.5}
We assume that $\mathscr{S}$ is a groupoidal category and $(F,{\mathcal F})$ a basic construction.
The associated natural topology for the orbit space $|\mathscr{S}|$ is denoted by ${\mathcal T}$.
We have seen that every ${\bm{M}}(\Psi,\Psi')$ has a natural metrizable topology for every pair
$(\Psi,\Psi')$. For this topology the source and target maps are local homeomorphisms and the other structure maps
are continuous.  We also equipped $|\mathscr{S}|$ with a topology ${\mathcal T}$ so that every point has
a countable neighborhood basis of open sets.
Recall that  for every  object $\alpha$ and $\Psi\in F(\alpha)$ the map on orbit space
$$
|\Psi|:{_G\backslash} O\rightarrow |\mathscr{S}|
$$
is a homeomorphism onto an open neighborhood of $|\alpha|$ in $|\mathscr{S}|$.
Using this fact and the previous discussion we can prove the following theorem.
\begin{theorem}\label{THMX19.5.1}
The metrizable  space ${\bm{M}}(\Psi,\Psi')$ has a natural  M-polyfold structure for which $s$ and $t$ are local sc-diffeomorphisms. 
\end{theorem}
\begin{proof}
As we shall see there is a sc-smooth atlas for ${\bm{M}}(\Psi,\Psi')$ which has transition maps of a simple form.
That this is the case is immediate apparent in the case $\bm{M}(\Psi,\Psi)$. In this case we have a bijection
$\tau:G\times O\rightarrow \bm{M}(\Psi,\Psi):(g,o)\rightarrow (o,\Psi(g,o),g\ast o)$. The sc-smooth atlas for
$\bm{M}(\Psi,\Psi)$ will define a M-polyfold structure for which the bijection $\tau$ is an sc-diffeomorphism
and $G\times O$ is equipped with the natural M-polyfold structure coming from $O$ (and the trivial structure on the finite set $G$).

Recall that with $h=(o,\Psi(g,o),g\ast o)$ the germ $F_h$ as the form 
$$
F_h(q)=(q,\Psi(g,q),g\ast q)
$$
 and for a small open neighborhood
$U$ of $o$ the map $F_h:U\rightarrow F_h(U)$ is a homeomorphism and $F_h(U)\in {\mathcal T}_{\Psi,\Psi}$.
We see that $s_h:F_h(U)\rightarrow U: s_h(o,\Psi(g,o),g\ast o)=o$ is precisely the inverse of $F_h$. We also 
note that if $F_h(U)\cap F_{h'}(U')\neq \emptyset$ the transition map $s_{h'}\circ s_h^{-1}$ has the form
$$
s_{h'}\circ s_h^{-1}(q)=q.
$$
This shows that the system $(s_h)$ is an sc-smooth atlas.  Moreover
$$
s_h\circ \tau(g,q)  = q
$$
proving that $\tau$ is an sc-diffeomorphism. 

We next consider the general case and as before we shall show that the source maps define an sc-smooth atlas. Again it will turn out that the 
transition maps take an easy form.
Pick $h\in {\bm{M}}(\Psi,\Psi')$ and pick open neighborhoods $U(h)$, $U(o)$, where $o=s(h)$, with the following properties
\begin{itemize}
\item[(i)] \ \ \  $s:U(h)\rightarrow U(o)$ is a homeomorphism, denoted by $s_h$.
\item[(ii)] \ \ \ $F_h:U(o)\rightarrow U(h)$ is the inverse of $s_h$. 
\item[(iii)] \ \ \ $t(F_h(q))=f_h(q)$ for $q\in U(o)$.
\end{itemize}
We shall show that the collection $(s_h)$, where $h$ varies over ${\bm{M}}(\Psi,\Psi')$,
has sc-smooth overlaps, i.e. $s_{h'}\circ s_h^{-1}:s_h(U(h)\cap U(h'))\rightarrow s_{h'}(U(h)\cap U(h'))$
is sc-smoth. This, of course, will imply that ${\bm{M}}(\Psi,\Psi')$ has a natural  M-polyfold structure. 

Assuming  $q\in s_h(U(h)\cap U(h'))$ we compute 
$$
s_{h'}\circ s_h^{-1}(q)= s_{h'}\circ F_h (q)=s(F_h(q))=q.
$$
This shows that ${\bm{M}}(\Psi,\Psi')$ has a tame M-polyfold structure and by construction
$s$ is a local sc-diffeomorphism.  The same holds for $t$. To see this fix $h\in {\bm{M}}(\Psi,\Psi')$
and consider $s_h:U(h)\rightarrow U(s(h))$. Then 
$$
t\circ s_h^{-1} (q) =t\circ F_h(q)=f_h(q)
$$
which is sc-smooth and a local sc-diffeomorphism. 
\qed \end{proof}

\begin{corollary}
Given $\Psi,\Psi',\Psi''$ associated to the objects $\alpha,\alpha',\alpha''$, consider
$$
{\bm{M}}(\Psi,\Psi'),\ {\bm{M}}(\Psi',\Psi''), \  \text{and}\ \ {\bm{M}}(\Psi,\Psi'')
$$
with their natural M-polyfold structures.
Then all structure maps are sc-smooth.
\end{corollary}
\begin{proof}
We begin with the unit map $u:O\rightarrow {\bm{M}}(\Psi,\Psi)$. Then $s\circ u(o)=o$ and since 
$s$ is a local sc-diffeomorphism we conclude that $u$ is sc-smooth.

Next we consider the inversion map $\iota:{\bm{M}}(\Psi,\Psi')\rightarrow {\bm{M}}(\Psi',\Psi)$. Fix $h=(o,\Phi,o')$
and let $h'=(o',\Phi^{-1},o)=\iota(h)$. Then $F_h:U(o)\rightarrow U(h)$ for suitable open neighborhoods 
is an sc-diffeomorphism and  for $q$ near $o$ we compute
$$
t(\iota(F_h(q))= t(F_{h'}(f_h(q)).
$$
Since $t$,$f_h$, and  $F_{h'}$ are local sc-diffeomorphisms we see that the left-hand side is sc-smooth.
Since $t$ and $F_h$ are local sc-diffeomorphisms we concluded that $\iota $ is sc-smooth.

Finally we consider the multiplication map $m$.  Assume that $h\in \bm{M}(\Psi,\Psi')$ and $h'\in \bm{M}(\Psi',\Psi'')$
satisfying $t(h)=s(h')$. 
Using that for $q$ near $o=s(h)$ we have the relationship
$$
m(F_{h'}(f_h(q)),F_h(q))=F_{m(h',h)}(q)
$$
we see that the right-hand side is a local sc-diffeomorphism, which implies the same for the expression
on the left. We note that
$$
q\rightarrow (F_{h'}(f_h(q)),F_h(q)) = m((f_{h}(q),\Phi_{f_h(q)}',f_{m(h',h)}(q)),(q,\Phi_q,f_h(q)))
$$
is the inverse of a  local chart in the fibered product ${\bm{M}}(\Psi',\Psi''){_{s}\times_t}{\bm{M}}(\Psi,\Psi')$.
Hence $m$ is sc-smooth.
\qed \end{proof}

In view of the discussion in this section we can conclude the following.
\begin{theorem}
Given a basic construction $(F,{\mathcal F})$ for $\mathscr{S}$, which satisfies the properness property and for which  the natural topology ${\mathcal T}$ is either 
second countable or paracompact there is an underlying natural  polyfold structure $(F,{\bm{M}})$ for $\mathscr{S}$. 
If the uniformizers $\Psi:G\ltimes O\rightarrow \mathscr{S}$ associated to $F$ have tame domains $O$ the polyfold structure is tame.
\qed
\end{theorem}

\begin{remark}\index{R- Comment on polyfold constructions}
In applications the category $\mathscr{S}$ very often is constructed having a geometric purpose in mind
and more often than not there is a natural  uniformizer construction $F$,  
and a transition construction ${\mathcal F}$.  In a first step
one obtains the natural topology ${\mathcal T}$ for $|\mathscr{S}|$
as well as the metrizable topologies for the transition sets $\bm{M}(\Psi,\Psi')$.
In fact,  the latter even carry natural M-polyfold structures.  We recall that ${\mathcal T}$ is locally metrizable.
After checking that
$(F,{\mathcal F})$ has, indeed, the properness property one concludes that ${\mathcal T}$ 
is locally metrizable, Hausdorff, and completely regular.  These facts imply that for a given 
compact subset $K$ of  $|\mathscr{S}|$ there exists an open neighborhood $U=U(K)$ which is paracompact and hence metrizable,
given the already established properties of ${\mathcal T}$.
This sometimes is all what one needs in applications.  In case that ${\mathcal T}$ is also second countable 
we note that ${\mathcal T}$ is metrizable and for such basic constructions $(F,{\mathcal F})$ we have just 
shown that there is an associate polyfold construction $(F,\bm{M})$. The same  holds if ${\mathcal T}$ is paracompact.  We also note that if all uniformizers $\Psi$ associated to $F$ are tame, i.e. the domains are tame, the resulting polyfold construction is tame.
\qed
\end{remark}

\section{A Strong Bundle Version}\label{SECX19.6}
In order to equip $P:{\mathcal E}\rightarrow \mathscr{S}$ with a strong bundle structure, there is a similar scheme as in the 
$\mathscr{S}$-case,  and we allow ourselves to be somewhat sketchy. Here $\mathscr{S}$ and ${\mathcal E}$ are  groupoidal categories, so that in addition the fibers $P^{-1}(\alpha)$ are equipped with Banach space structures.
The morphisms in ${\mathcal E}$ are assumed to be  linear lifts associated to the underlying morphisms in $\mathscr{S}$. More precisely, for two objects $\alpha$ and $\alpha'$ in $\mathscr{S}$
and an isomorphism $\Phi:\alpha\rightarrow \alpha'$ there exists a well-defined topological linear isomorphism
$$
\wh{\Phi}:P^{-1}(\alpha)\rightarrow P^{-1}(\alpha').
$$ 
We shall call $\wh{\Phi}$ the {\bf lift}  of $\Phi$ \index{$\wh{\Phi}$}\index{Lift of $\Phi$} and we assume that the procedure 
$\Phi\rightarrow \wh{\Phi}$ has the obvious properties
\begin{itemize}
\item[(i)] \ \ $\wh{1_\alpha} = Id_{P^{-1}(\alpha)}$.
\item[(ii)] \ \ If $s(\Phi)=t(\Phi')$ then $\wh{\Phi'\circ \Phi}= \wh{\Phi}'\circ \wh{\Phi}$. 
\end{itemize}
A morphism in ${\mathcal E}$ is assumed to be of the form $(\Phi,e)$ with $s(\Phi)=P(e)$
so that 
$$
s(\Phi,e)=e\ \ \text{and}\ \ t(\Phi,e) =\wh{\Phi}(e).
$$
Of course, there is no difference if the morphisms are taken as tuples $(e,\wh{\Phi},e')$, where $e'=\wh{\Phi}(e)$, since
the two descriptions are canonically isomorphic. We take the shorter one, i.e. the first one.
Hence 
$$
e\xrightarrow{(\Phi,e)} \wh{\Phi}(e)\ \ \text{or}\ \ e\xrightarrow{(\Phi,e)} e'\ \ (\text{provided}\ \ e'=\wh{\Phi}(e)\ ).
$$
Having these assumptions on ${\mathcal E}$ there is an alternative description of our data, which is the most useful.
We have seen this already in the ep-groupoid context.
Denote by $\text{Ban}$ the category whose objects are Banach spaces and where the morphisms 
are linear topological isomorphisms.  Then $P:{\mathcal E}\rightarrow \mathscr{S}$, with the structures
described above, gives a covariant
functor 
$$
\mu:\mathscr{S}\rightarrow \text{Ban}.
$$
Indeed, on objects $\mu(\alpha):=P^{-1}(\alpha)$, and on morphisms 
$$
\mu(\Phi) =\wh{\Phi}.
$$
The original $P:{\mathcal E}\rightarrow \mathscr{S}$ can be,  up to natural isomorphism, reconstructed from $(\mathscr{S},\mu)$.

Given a groupoidal category $\mathscr{S}$ and a functor
$\mu:\mathscr{S}\rightarrow \text{Ban}$ we can define a category ${\mathcal E}_\mu$
as follows. The objects are the pairs $(\alpha,e)$ with $e\in \mu(\alpha)$, and the morphisms 
are the pairs $(\Phi,e)$  with $e\in \mu(s(\Phi))$, viewed as morphisms 
$$
(\Phi,e):e\rightarrow \mu(\Phi)(e).
$$
In view of this discussion we study $P_\mu:{\mathcal E}_\mu\rightarrow \mathscr{S}$ associated to a functor $\mu:\mathscr{S}\rightarrow \text{Ban}$ and abbreviated by $P:{\mathcal E}\rightarrow \mathscr{S}.$ 

\begin{definition}\index{D- Strong bundle unifomizer}\label{DERG19.9}
Consider $(\mathscr{S},\mu)$, where $\mathscr{S}$ is a groupoidal category and $\mu :\mathscr{S}\rightarrow \text{Ban}$ a functor. Denote the associated functor diagram by $P:{\mathcal E}\rightarrow \mathscr{S}$. A {\bf strong bundle uniformizer} at the object $\alpha$ is a functor 
$$
\bar{\Psi}: G\ltimes K\rightarrow {\mathcal E}
$$
with the following properties.
\begin{itemize}
\item[(1)] \  $p:K\rightarrow O$ is a strong bundle over  a M-polyfold,  and $G$ acts by strong bundle diffeomorphisms. Moreover $G\ltimes K$ is the associated translation groupoid.  The functor 
$\bar{\Psi}$ is a  topological linear isomorphism between fibers. In addition there exists
a uniformizer at $\alpha$ denoted by $\Psi:G\ltimes O\rightarrow \mathscr{S}$ fitting into the commutative diagram
$$
\begin{CD}
G\ltimes K @>\bar{\Psi}>> {\mathcal E}\\
@V Id\ltimes p VV     @V P VV\\
G\ltimes O   @>\Psi >>   \mathscr{S}.
\end{CD}
$$
Note that $\Psi$ is completely determined by $\bar{\Psi}$.
\item[(2)] \   $\bar{\Psi}$  is injective on objects and fully faithful.
\item[(3)] \   There exists $0_{\bar{o}}\in K_{\bar{o}}$ with $\Psi(\bar{o})=\alpha$ and
$\bar{\Psi}(0_{\bar{o}}) = (\alpha,0_{\mu(\alpha)})$.
\end{itemize}
\qed
\end{definition}
\begin{remark}\label{remark19.8}
We can also require, more generally, that $p:W\rightarrow E$ is a strong bundle over an ep-groupoid
and that $\bar{\Psi}$ and $\Psi$ are fully faithful functors, which both are injective on objects and fit into the commutative diagram
$$
\begin{CD}
W @>\bar{\Psi}>> {\mathcal E}\\
@V p VV     @V P VV\\
E   @>\Psi >>   \mathscr{S},
\end{CD}
$$
and that a distinguished point $e\in E$ is mapped to the given object $\alpha$ and consequently $\bar{\Psi}(0_e)=(\alpha,0_{\mu(e)})$.
Of course, we also assume that $\bar{\Psi}$ is fiber-wise a topological  linear isomorphism.
Hence, rather than taking $G\ltimes K\rightarrow G\ltimes O$ as local models one could take $W\rightarrow E$.
\qed
\end{remark}
From the definition it follows that a strong bundle uniformizer $\bar{\Psi}:G\ltimes K\rightarrow {\mathcal E}$
at $\alpha$ induces a strong bundle uniformizer $\Psi:G\ltimes O\rightarrow \mathscr{S}$.
\begin{definition}\index{D- Strong bundle uniformizer construction} 
A {\bf strong bundle unformizer construction} is a functor $\bar{F}:\mathscr{S}^-\rightarrow \text{SET}$ associating to every object $\alpha$ a set of strong bundle uniformizers $ \bar{F}(\alpha)$.
\qed
\end{definition}
Associated to $\bar{F}$ we have the underlying uniformizer construction $F$ for $\mathscr{S}$.
Given $\bar{F}$ we can build for $\bar{\Psi}\in \bar{F}(\alpha)$ and $\bar{\Psi}'\in \bar{F}(\alpha')$
the transition set $\bm{M}(\bar{\Psi},\bar{\Psi}')$ associated to the functor diagram
$$
K\xrightarrow{\bar{\Psi}}  {\mathcal E} \xleftarrow{\bar{\Psi}'} K',
$$
which covers 
$$
O\xrightarrow{\Psi} \mathscr{S} \xleftarrow{\Psi} O'.
$$
The elements in the transition set have the form $(a,(\Phi,e),a')$, where  $e\in \mu(s(\Phi))$, $\bar{\Psi}(a)=(s(\Phi),e)$
and $\bar{\Psi}'(a') = (t(\Phi),\wh{\Phi}(e))$. Since $e$ is determined by $a$ via $\bar{\Psi}$, and $\Phi$ acts as a linear map,
i.e. $\wh{\Phi}=\mu(\Phi)$,  we just write
$(a,\Phi,a')$ instead of $(a,(\Phi,e),a')$.
There exists a natural map 
$$
\bm{m}\colon \bm{M}(\bar{\Psi},\bar{\Psi}')\rightarrow \bm{M}(\Psi,\Psi')
$$
defined by
$$
\bm{m}(a,\Phi,a')= (p(a),\Phi,p'(a')).
$$
The fibers of $\bm{m}$ are Banach spaces. Indeed, with $(o,\Phi,o')$ being the base point,
the elements above it consist of all $(a,\Phi,{(\bar{\Psi}')}^{-1}\circ \wh{\Phi}\circ \bar{\Psi}(a))$, where $a\in p^{-1}(o)$.  The structure maps associated to strong bundle unifomizers are given as follows and are linear on the fibers
and cover those for the underlying unifomizer construction. 
\begin{description}
\item[$\bullet$]  The {\bf source map}\index{Bundle source map}
$$
\begin{array}{cc}
\begin{CD}
\bm{M}(\bar{\Psi},\bar{\Psi}') @> s >>  K \\
@V\bm{m} VV                                   @ V p VV    \\     
\bm{M}(\Psi,\Psi') @> s >>                    O
 \end{CD}
&\ \ \ \ 
\begin{CD}
 (a, \Phi,a') @>>>   a\\
 @V VV      @V VV \\
     (p(a),\Phi, p'(a')) @>>>   p(a).
\end{CD}
\end{array}
$$
\item[$\bullet$] The {\bf target map}\index{Bundle target map}
$$
\begin{array}{cc}
\begin{CD}
\bm{M}(\bar{\Psi},\bar{\Psi}') @> t >>  K' \\
@V\bm{m} VV                                   @ V p' VV    \\     
\bm{M}(\Psi,\Psi') @> t >>                    O'
 \end{CD}
&\ \ \ \ 
\begin{CD}
 (a, \Phi,a') @>>>   a'\\
 @V VV      @V VV \\
     (p(a),\Phi, p'(a')) @>>>   p'(a').
\end{CD}
\end{array}
$$
\item[$\bullet$] The {\bf unit map}\index{Bundle unit map}
$$
\begin{array}{cc}
\begin{CD}
K @> u >> \bm{M}(\bar{\Psi},\bar{\Psi})\\
@V p VV        @V\bm{m} VV\\
O @> u>>   \bm{M}(\Psi,\Psi')
\end{CD}
&\ \ \ \
\begin{CD}
a @>>>   (a,{1}_{\Psi(p(a))},a)\\
@VVV @VVV\\
p(a) @>>>   (p(a),1_{\Psi(p(a))},p(a)).
\end{CD}
\end{array}
$$
\item[$\bullet$] The {\bf inversion map}\index{Bundle inversion map}
$$
\begin{CD}
\bm{M}(\bar{\Psi},\bar{\Psi}') @>\iota>>  \bm{M}(\bar{\Psi}',\bar{\Psi})\\
@V\bm{m} VV    @V\bm{m} VV \\
\bm{M}(\Psi,\Psi') @>\iota >> \bm{M}(\Psi',\Psi)
\end{CD}
$$
$$
\begin{CD}
(a,{\Phi},a')  @>>> (a',{\Phi}^{-1}, a)\\
@ VVV @VVV \\
(p(a),\Phi,p'(a'))  @>>> (p'(a'),\Phi^{-1},p(a)).
\end{CD}
$$
\item[$\bullet$] The {\bf multiplication map}\index{Bundle multiplication map}, where $o=p(a), o'=p'(a'), o''=p''(a'')$
$$
\begin{CD}
\bm{M}(\bar{\Psi}',\bar{\Psi}''){_{s}\times_t} \bm{M}(\bar{\Psi},\bar{\Psi}') @> m>> \bm{M}(\bar{\Psi},\bar{\Psi}'')\\
@V \bm{m}{_{s}\times_t}\bm{m}  VV @ V \bm{m} VV\\
\bm{M}(\Psi',\Psi''){_{s}\times_t} \bm{M}({\Psi},{\Psi}') @> m>> \bm{M}({\Psi},{\Psi}'')
\end{CD}
$$
$$
\begin{CD}
((a',{\Phi}',a''),(a,{\Phi},a')) @>>>   (a,{\Phi}'\circ {\Phi},a'')\\
@VVV @VVV\\
((o',\Phi',o''),(o,\Phi,o')) @>>>   (o,\Phi'\circ\Phi,o'').
\end{CD}
$$
\end{description}
\mbox{}

Next we need a construction $\bar{\mathcal F}$ which extends ${\mathcal F}$ in the $\mathscr{S}$-case.  
Recall that ${\mathcal F}$  associates to $h =(o,\Phi,o')\in \bm{M}(\Psi,\Psi')$ a germ
of map
$$
F_h :{\mathcal O}(O,o)\rightarrow ( \bm{M}(\Psi,\Psi'),h).
$$
We can write $F_h$ as
$$
F_h(u) =(u,\Phi_u, f_h(u)),
$$
where $\Phi_u:\Psi(u)\rightarrow \Psi'(f_h(u))$. This morphism lifts to the linear isomorphism
$$
\wh{\Phi}_u: P^{-1}(\Psi(u))\rightarrow P^{-1}(\Psi'(f_h(u)))
$$
and suggests to define  $\bar{F}_h$  as follows
$$
\bar{F}_h(a) = (a, {\Phi}_{p(a)},  a'(a)),
$$
where $a'(a)$ is uniquely determined by $\bar{\Psi}'(a'(a))= \wh{\Phi}_{p(a)} (a)$.
The germ  $\bar{F}_h$ is defined on $K|{\mathcal O}(O,o)$ and covers $F_h$.
$$
\begin{CD}
K|{\mathcal O}(O,o) @> \bar{F}_h >> \bm{M}(\bar{\Psi},\bar{\Psi}')\\
@Vp VV @V \bm{m} VV\\
{\mathcal O}(O,o) @ > F_h >> \bm{M}(\Psi,\Psi')
\end{CD}
$$
We observe that 
$$
\bm{m}(\bar{F}_h(a)) = F_h(p(a)), 
$$
so that the definition of $\bar{F}_h$ is a natural way to lift $F_h$.  We can define 
$\bar{f}_h: K|{\mathcal O}(O,o)\rightarrow K'|{\mathcal O}(O',o')$ covering $f_h$
by
$$
\bar{f}_h (a)= t\circ \bar{F}_h(a).
$$
We require the following properties:
\begin{definition}\index{D- Transition requirements for strong bundles}
The following properties imposed on $\bar{\mathcal F}$ and ${\mathcal F}$
are referred to as {\bf transition requirements} for strong bundles.
\begin{itemize}
\item[(A)]\ {\bf Diffeomorphism Property}:  The germs $\bar{f}_h=t\circ\bar{F}_h:K|{\mathcal O}(O,o)\rightarrow K'|{\mathcal O}(O',o')$ are local strong bundle isomorphisms covering $f_h$
$$
\begin{CD}
K|{\mathcal O}(O,o)@>\bar{f}_h >> K'|{\mathcal O}(O',o')\\
@V p VV   @V p'VV\\
{\mathcal O}(O,o)@> f_h >>{\mathcal O}(O',o')
\end{CD}
$$
 and $s(\bar{F}_h(a))=a$ for $a\in K$ with $p(a)$ near $o$. If $\bar{\Psi}=\bar{\Psi}'$
and $h=(o,\Psi(g,o),g\ast o)$ it holds $\bar{F}_h(a)=(a,{\Psi}(g,p(a)),g\ast a)$ if $p(a)$ is near $o$. In particular
$\bar{f}_h(a)=g\ast a$ (covering $f_h(p(a))= g\ast (p(a))=p(g\ast a)$.).

\item[(B)]\ {\bf Stability Property}:  $\bar{F}_{F_{h}(q)}(a)= \bar{F}_h(a)$ for $q$ near $o=s(h)$ and $a$ such that $p(a)$
 is near $q$.

\item[(C)]\  {\bf Identity Property}:  $\bar{F}_{u(o)}(a) = (a,{1}_{\Psi(p(a))},a) =u(a)$ for $p(a)$ near $o$.

\item[(D)]\  {\bf Inversion Property}: $\bar{F}_{\iota(h)}(\bar{f}_h(a)) = \iota(\bar{F}_h(a))$ for $a$ such that $p(a)$ is near $o=s(h)$.

\item[(E)]\ {\bf Multiplication  Property}:   If $s(h')=t(h)$ then $\bar{f}_{h'}\circ \bar{f}_h(a)=\bar{f}_{m(h',h)}(a)$ for $p(a)$ near $o=s(h)$ and 
 $m(\bar{F}_{h'}(\bar{f}_h(a)),\bar{F}_h(a))=\bar{F}_{m(h',h)}(a)$ for $p(a)$ near $o=s(h)$.

\item[(F)]\ $\bm{M}$-{\bf Hausdorff Property}: The $F_h$ underlying $\bar{F}_h$ has the property 
that for different $h_1,h_2\in {\bm{M}}(\Psi,\Psi')$ with $o=s(h_1)=s(h_2)$  the images under $F_{h_1}$ and $F_{h_2}$ of small neighborhoods $o$ are disjoint.
\end{itemize}
\qed
\end{definition}
\begin{definition}\index{D- Basic strong bundle construction}
Given a pair $(\mathscr{S},\mu)$, where $\mathscr{S}$ is a groupoidal category, and $\mu:\mathscr{S}\rightarrow \text{Ban}$ is a covariant functor with associated  $P:{\mathcal E}\rightarrow \mathscr{S}$,
we shall refer to the construction $(\bar{F},\bar{\mathcal F})$, where $\bar{F}$ is a strong bundle uniformizer construction, and $\bar{\mathcal F}$ is an associated strong bundle transition construction as
as {\bf basic strong bundle construction}  for $(\mathscr{S},\mu)$.
\qed
\end{definition}

At this point we can run through the construction in the $\mathscr{S}$-case and equip 
$P:{\mathcal E}\rightarrow \mathscr{S}$ with a strong bundle structure, provided the underlying 
basic construction for $\mathscr{S}$ induces a polyfold construction for $\mathscr{S}$.

\subsection*{Natural Topology $\bar{\mathcal T}$}

For the following discussion we  assume that starting with $(\mathscr{S},\mu)$ we have a basic bundle construction $(\bar{F},\bar{\mathcal F})$ for $P:{\mathcal E}\rightarrow \mathscr{S}$ covering $(F,{\mathcal F})$.
In a first step we define a topology $\bar{\mathcal T}$ for $|{\mathcal E}|$. 
\begin{definition}\label{DEF1964}
A subset $\wh{U}$ of $|{\mathcal E}|$ belongs to $\bar{\mathcal T}$ provided for every $\bar{z}\in \wh{U}$
there exists an object $\alpha$ in $\mathscr{S}$ with $|\alpha|=z:=|P|(\bar{z})$, $\bar{\Psi}\in\bar{F}(\alpha)$
$$
\bar{\Psi}:G\ltimes K\rightarrow {\mathcal E},
$$
a point $a\in K$ with $|\bar{\Psi}(a)|=\bar{z}$, and an open neighborhood 
$\wh{V}$ of $a$ in $K$ such that $|\bar{\Psi}(\wh{V})|\subset \wh{U}$.
\qed
\end{definition}
The underlying basic construction $(F,{\mathcal F})$ defines the topology ${\mathcal T}$
for $|\mathscr{S}|$.  The following is easily verified using the continuity of the strong bundle maps
$p:K\rightarrow O$. The following result is an analogue to Theorem \ref{TOPOLOGY}.

\begin{theorem}\index{T- Natural topology $\bar{\mathcal T}$}
Assume we are given a basic bundle construction $(\bar{F},\bar{\mathcal F})$ for $(\mathscr{S},\mu)$
with underlying basic construction $(F,{\mathcal F})$. Then the following holds.
\begin{itemize}
\item[{\em (1)}]\  The set $\bar{\mathcal T}$ is a topology on $|{\mathcal E}|$. 
\item[{\em (2)}]\  For $\bar{\Psi}\in \bar{F}(\alpha)$, $\bar{\Psi}:G\ltimes K\rightarrow {\mathcal E}$,  and $\wh{W}\subset K$ open the set $|\bar{\Psi}(W)|$
is open in $|{\mathcal E}|$.
\item[{\em (3)}]\   With $\bar{\Psi}$ as in (ii) the induced $|\bar{\Psi}|:|G\ltimes K|\rightarrow |{\mathcal E}|$ is a homeomorphism onto an open set.
\item[{\em (4)}] \   For every point $z\in |\mathscr{S}|$ there exists an open neighborhood $V$ of $z$ in $|\mathscr{S}|$
so that $|P|^{-1}(V)$ with the induced topology from $\bar{\mathcal T}$ is metrizable.
\item[{\em (5)}]\  The map $|P|:|{\mathcal E}|\rightarrow |\mathscr{S}|$ is continuous.
\end{itemize}
\qed
\end{theorem}

The metrizability of $\bar{\mathcal T}$ only depends on the metrizability of the underlying ${\mathcal T}$ 
on $|\mathscr{S}|$.
\begin{proposition}\index{P- Metrizability for $\bar{\mathcal T}$}
Assume we are given a basic bundle construction $(\bar{F},\bar{\mathcal F})$ for $(\mathscr{S},\mu)$
with underlying basic construction $(F,{\mathcal F})$.  If the natural topology ${\mathcal T}$ for
$|\mathscr{S}|$ is metrizable the topology $\bar{\mathcal T}$ for $|{\mathcal E}|$ is metrizable as well.
\qed
\end{proposition}

\subsection*{Natural Strong Bundle  Structures for $\bm{M}(\bar{\Psi},\bar{\Psi}')$}
Under the same assumptions as before we can equip $\bm{M}(\bar{\Psi},\bar{\Psi}')$ with a metrizable topology.
\begin{definition}\index{D- Topology $\bar{\mathcal T}_{\bar{\Psi},\bar{\Psi}'}$}
For $\bar{\Psi}\in \bar{F}(\alpha)$ and $\bar{\Psi}'\in \bar{F}(\alpha')$ denote by $\bar{\mathcal T}_{\bar{\Psi},\bar{\Psi}'}$ the collection of subsets $\wh{W}$ of $\bm{M}(\bar{\Psi},\bar{\Psi}')$ having the following property.
Given $\bar{h}=(a,\wh{\Phi},a')$ there exists an open neighborhood $\wh{V}$ of $\bar{h}$ in $K$
such that $\bar{F}_h(\wh{V})\subset \wh{W}$.
\qed
\end{definition}\index{T- Topological properties of $\bm{M}(\bar{\Psi},\bar{\Psi}')$}
The basic topological  results are  stated in the following theorem.
\begin{theorem}
Assume we are given a basic bundle construction $(\bar{F},\bar{\mathcal F})$ for $(\mathscr{S},\mu)$
with underlying basic construction $(F,{\mathcal F})$. The set $\bar{\mathcal T}_{\bar{\Psi},\bar{\Psi}'}$
defines a topology on $\bm{M}(\bar{\Psi},\bar{\Psi}')$.  The topological spaces
$(\bm{M}(\bar{\Psi},\bar{\Psi}'), \bar{\mathcal T}_{\bar{\Psi},\bar{\Psi}'})$ are metrizable and have the following properties.
\begin{itemize}
\item[{\em (1)}]\   The source map $s:\bm{M}(\bar{\Psi},\bar{\Psi}')\rightarrow K$ is a local homeomorphism.
\item[{\em (2)}]\  The target  map $t:\bm{M}(\bar{\Psi},\bar{\Psi}')\rightarrow K'$ is a local homeomorphism.
\item[{\em (3)}]\  The unit map $u:K\rightarrow \bm{M}(\bar{\Psi},\bar{\Psi})$ is continuous.
\item[{\em (4)}]\   The inversion map 
$\iota:\bm{M}(\bar{\Psi},\bar{\Psi}')\rightarrow \bm{M}(\bar{\Psi}',\bar{\Psi})$ is a homeomorphism.
\item[{\em (5)}]\  The multiplication map
$$
m:\bm{M}(\bar{\Psi}',\bar{\Psi}''){_{s}\times_t} \bm{M}(\bar{\Psi},\bar{\Psi}') \rightarrow \bm{M}(\bar{\Psi},\bar{\Psi}'')
$$ 
is continuous.
\item[{\em (6)}]\   $\bm{m}: \bm{M}(\bar{\Psi},\bar{\Psi}')\rightarrow \bm{M}(\Psi,\Psi')$ 
for the natural topologies is continuous.
\end{itemize}
In addition the projection map $\bm{m}:\bm{M}(\bar{\Psi},\bar{\Psi}')\rightarrow \bm{M}(\Psi,\Psi')$ is continuous.
\qed
\end{theorem}

\subsection*{Strong Bundle Structure for $P:{\mathcal E}\rightarrow \mathscr{S}$}
Given $(\bar{F},\bar{\mathcal F})$ it follows that $\bm{m}:\bm{M}(\bar{\Psi},\bar{\Psi}')\rightarrow \bm{M}(\Psi,\Psi')$ has in a natural way the structure of a strong bundle over a M-polyfold.

\begin{theorem}\index{T- Natural strong bundle structure for $\bm{m}$}
Assume we are given a basic bundle construction $(\bar{F},\bar{\mathcal F})$ for $(\mathscr{S},\mu)$
with underlying basic construction $(F,{\mathcal F})$. Then, with the spaces 
equipped with their natural topologies $\bm{m}:\bm{M}(\bar{\Psi},\bar{\Psi}')\rightarrow \bm{M}(\Psi,\Psi')$ admits a natural structure of a M-polyfold bundle, where the structure on the base
coincides with the previously defined natural M-polyfold structure. For this structure the 
structural maps $s$ and $t$ are locally diffeomorphic strong bundle maps. The unit map 
and multiplication are strong bundle maps, and the inversion map is a strong bundle isomorphism.
\qed
\end{theorem}
From this we deduce the following basic result.
\begin{theorem}
Assume we are given a basic bundle construction $(\bar{F},\bar{\mathcal F})$ for $(\mathscr{S},\mu)$
with underlying basic construction $(F,{\mathcal F})$.  Assume that the natural topology ${\mathcal T}$ for 
$|\mathscr{S}|$ is metrizable. Then the natural topology $\bar{\mathcal T}$ for $|{\mathcal E}|$ is metrizable.
Further for objects $\alpha,\alpha'$ in $\mathscr{S}$ and $\bar{\Psi}\in\bar{F}(\alpha)$, $\bar{\Psi}'\in \bar{F}(\alpha')$
$$
\bm{M}(\bar{\Psi},\bar{\Psi}')\rightarrow \bm{M}(\Psi,\Psi')
$$
has in a natural way the structure of a strong bundle over a M-polyfold. Hence there is an associated
strong bundle construction $(\bar{F},\bm{M})$ for $P:{\mathcal E}\rightarrow \mathscr{S}$.
\qed
\end{theorem}
\begin{remark}
The results in Chapter \ref{CHAPTER19X} allow to construct polyfold structures on certain groupoidal categories $\mathscr{S}$ as well as certain bundles
over them $P:{\mathcal E}\rightarrow \mathscr{S}$. Once $P$ is equipped with the structure as a strong bundle over a polyfold we can study
section functors $S$ of $P$. Of particular interest are sc-Fredholm functors, which can be studied by the methods introduced in this book, see in particular Chapters \ref{CHAPX17} and \ref{CHAPX18}.
\qed
\end{remark}

\section{Covering Constructions}\label{SECX19.7}
The main point in this section is the description of a basic type of construction, which provides a given 
finite-to-one covering functor $\mathsf{P}:\mathscr{A}\rightarrow \mathscr{B}$, defined below, with an adapted uniformizer construction, where the uniformizers respect the structure functor $\mathsf{P}$, and in addition are sc-smoothly compatible.
 As such this is a generalization of the ideas in Section \ref{SECX19.1} and quite close to the corresponding strong bundle version 
given in Section \ref{SECX19.6}.
Recall Definition \ref{DEFNX19.1.1} for the notion of a groupoidal category.  We shall denote by $A$ the class
of objects associated to $\mathscr{A}$ and by $\bm{A}$ the class of morphisms. Similarly for $\mathscr{B}$ they are denoted by $B$ and
$\bm{B}$.

\begin{definition}\index{D- Finite to one covering functor}\label{DEFNX19.7.1}
A {\bf finite-to-one  covering functor} $\mathsf{P}:\mathscr{A}\rightarrow \mathscr{B}$ 
 is a functor between  two groupoidal categories,
 which has the following properties:
 \begin{itemize}
 \item[(1)]\   $\mathsf{P}$ is surjective on objects and  the preimage of every object is finite.
 \item[(2)] \ The assignment $\bm{A}\rightarrow \bm{B}{_s\times_{\mathsf{P}}} A$ which maps $\phi$ to
 $(\mathsf{P}(\phi),s(\phi))$ is a bijection.
 \end{itemize}
 \qed
\end{definition}
\begin{remark}
Note that there is no topology in play and compare with Definition \ref{DEFNX17.7.1}.
Essentially the above definition reflects the properties (ii) and (iii) from \ref{DEFNX17.7.1}.
We shall discuss constructions which shall also define topologies on $|\mathscr{A}|$ and $|\mathscr{B}|$.
These additional structures will turn under suitable assumptions the groupoidal categories $\mathscr{A}$  and $\mathscr{B}$  into GCT's 
and $\mathsf{P}$ into a proper covering functor.
\qed
\end{remark}
 Next we shall introduce the notion 
of uniformizer but, for the lack of  having a topology on the orbit spaces, the definitions 
will not involve a topology for the moment.
\begin{definition}\label{DEFNX19.7.2}
Given a finite-to-one covering functor $\mathsf{P}:\mathscr{A}\rightarrow \mathscr{B}$, a {\bf proper covering uniformizer}
(or uniformizer for short)  for $\mathsf{P}$ at the object $\beta$ in $\mathscr{B}$ with isotropy $G$ consists of a pair of functors
$\wh{\Psi}:E\rightarrow \mathscr{A}$ and  $\Psi: G\ltimes O_\beta\rightarrow \mathscr{B}$, where $\wh{\tau}:E\rightarrow G\ltimes O_\beta$ is a proper ep-groupoid covering model, see Definition \ref{DEFQ1773},  so that the following holds
\begin{itemize}
\item[(1)] \ The distinguished point $o_\beta$ is mapped by $\Psi$ to $\beta$. 
\item[(2)] \ The functors $\Psi$ and $\wh{\Psi}$ are fully faithful and injective on objects and fit into the commutative  diagram
\begin{eqnarray}\label{EQNY19.1}
\begin{CD}
E @>\wh{\Psi} >>  \mathscr{A}\\
@V\wh{\tau} VV    @V \mathsf{P} VV\\
G\ltimes O_\beta @>\Psi>> \mathscr{B}
\end{CD}
\end{eqnarray}
\item[(3)]  \ The diagram (\ref{EQNY19.1}) is cartesian (see the next definition).
\end{itemize}
\qed
\end{definition}
\begin{definition}\index{D- Cartesian diagram}
The functor diagram (\ref{EQNY19.1}) is {\bf cartesian} provided on the object level the  image of the functor pair
$\langle\wh{\tau},\wh{\Psi}\rangle :E\rightarrow O_\beta\times A$ precisely consists of all pairs 
$(o,a)$ such that $\mathsf{P}(a) =\Psi(o)$.  
On the morphism level  the following is required.
 If $\wh{\Phi}$ belongs to $\bm{A}$
and $\mathsf{P}(\wh{\Phi})= \Psi(g,o)$ then there exists $\phi\in \bm{E}$
with $\wh{\tau}(\phi)=(g,o)$ and $\wh{\Psi}(e)=\wh{\Phi}$.
\qed
\end{definition}
Let us make the following observations concerning Definition \ref{DEFNX19.7.2}. Recalling 
Definition \ref{DEFQ1773}  we see that on the object level
$E$ is the disjoint union of finitely many open subsets $\wh{O}_a$, where $a$ varies in $\wh{\tau}^{-1}(o_\beta)$, so that
$\wh{\tau}:(O_a,a)\rightarrow (O_\beta,o_\beta)$ is an sc-diffeomorphism. This implies that on the object level $\Psi:O_b\rightarrow \mathscr{B}$ is determined by $\wh{\Psi}$ via
$$
\Psi(o) =\mathsf{P}\circ\wh{\Psi}\circ (\wh{\tau}|\wh{O}_a)^{-1}.
$$
Moreover, $\Psi$ is determined on morphisms as follows.   Since $\wh{\tau}$ is a proper ep-groupoid covering model the map
$$
\bm{E}\rightarrow \text{mor}(G\ltimes O_\beta){_{s}\times_{\wh{\tau}}}E: \phi\rightarrow (\wh{\tau}(\phi),s(\phi))
$$
is an sc-diffeomorphism.  Given a morphism $(g,o)\in G\times O_\beta$  we can fix a point $\wh{o}\in \wh{O}_a$ satisfying $\wh{\tau}(\wh{o})=o$.
We find a unique $\phi\in \bm{E}$ such that
$$
\wh{\tau}(\phi) =(g,o)\ \ \text{and}\ \ \ s(\phi)=\wh{o}.
$$
Then by commutativity  of the functor diagram
$$
\Psi(g,o) = \Psi\circ \wh{\tau}(\phi) = \mathsf{P}\circ \wh{\Psi}(\phi).
$$
Hence we have proved the following lemma.
\begin{lemma}\label{LEMMX19.7.3}
Let  $\wh{\tau}:E\rightarrow G
\ltimes O_\beta$ be a proper ep-groupoid covering model,  $\mathsf{P}:\mathscr{A}\rightarrow \mathscr{B}$
 a finite-to-one covering functor, and  $\wh{\Psi}:E\rightarrow \mathscr{A}$  having the property 
that  there exists a uniformizer $\Psi:G\ltimes O_\beta\rightarrow \mathscr{B}$
so that $(\wh{\Psi},\Psi)$ is a proper covering uniformizer. Then $\Psi$ is uniquely determined by $\wh{\Psi}$.
\qed
\end{lemma}
\begin{remark}
As a consequence, rather than  referring to Definition \ref{DEFNX19.7.2} we shall just refer to a uniformizer $\wh{\Psi}$
for $\mathsf{P}$. This then implicitly implies  that there exists a $\Psi$  fitting into the obvious commutative diagram. As shown such a $\Psi$ is unique, if it exists.
\qed
\end{remark}

As before we can introduce the notion of a uniformizer construction which fits the current context.
\begin{definition}
Given a finite-to-one  covering functor $\mathsf{P}:\mathscr{A}\rightarrow \mathscr{B}$, a {\bf proper covering uniformizer construction}
or simply just a uniformizer construction is given by a functor $\bar{F}:\mathscr{B}^-\rightarrow \text{SET}$
which associates to an object $\beta$ a set $\bar{F}(\beta)$ of unifomizers for $\mathsf{P}$ at $\beta$
in the sense of Definition \ref{DEFNX19.7.2} or the  remark after  Lemma \ref{LEMMX19.7.3}.
\qed
\end{definition}
\begin{remark}
The construction associates to $\beta$ a set of uniformizers $\wh{\Psi}$ for $\mathsf{P}$ and 
from this data we can construct the associated $\Psi$.
As already pointed out in previous discussions $\bar{F}$ usually has more structure and very often is a functor 
on $\mathscr{B}$.
\qed
\end{remark}
In a next step we note that given $\bar{F}$, there are transition sets associated $\wh{\Psi}$ and $\wh{\Psi}'$.
Since this data  also determines $\Psi$ and $\Psi'$ we obtain transition data for the latter. The sets 
$\bm{M}(\wh{\Psi},\wh{\Psi}')$ and $\bm{M}(\Psi,\Psi')$ are naturally related by the map $\mathsf{p}$
\begin{eqnarray}
\mathsf{p}:\bm{M}(\wh{\Psi},\wh{\Psi}')\rightarrow \bm{M}(\Psi,\Psi'):\mathsf{p}(e,\phi,e')=(\wh{\tau}(e),\mathsf{P}(\phi),\wh{\tau}'(e')).
\end{eqnarray}
We note that the map $\mathsf{p}$ is surjective and the preimage of each point is a finite set. 
Using the same arguments as in Lemma \ref{LEMM1776} we obtain the following result.
\begin{lemma}
The map $\bm{M}(\wh{\Psi},\wh{\Psi}')\rightarrow \bm{M}(\Psi,\Psi'){_{s}\times_{\wh{\tau}}}E$ 
defined by 
$$
(e,\phi,e')\rightarrow (\mathsf{p}(e,\phi,e'),s(e,\phi,e'))= ((\wh{\tau}(e),\mathsf{P}(\phi),\wh{\tau}'(e')),e)
$$
is a bijection.
\qed
\end{lemma}
We shall list the various transition sets which we can build. The source and target maps fit into the following diagram
$$
\begin{CD}
E@< s<<  \bm{M}(\wh{\Psi},\wh{\Psi}') @>t >>  E'\\
@V \wh{\tau} VV @V \mathsf{p} VV  @V \wh{\tau}' VV \\
O_\beta @< s <<   \bm{M}(\Psi,\Psi') @> t>> O_{\beta'}'.
\end{CD}
$$  
The inversion maps and unit maps give the following two diagrams.
$$
\begin{array}{cc}
\begin{CD}
\bm{M}(\wh{\Psi},\wh{\Psi}') @>\iota >>  \bm{M}(\wh{\Psi}',\wh{\Psi})\\
@V \mathsf{p} VV  @V \mathsf{p}VV\\
\bm{M}(\Psi,\Psi') @>\iota >> \bm{M}(\Psi',\Psi)
\end{CD}
\ \ \ &\ \ \
\begin{CD}
E @> u >> \bm{M}(\wh{\Psi},\wh{\Psi})\\
@V \wh{\tau} VV   @V \mathsf{p} VV\\
O_\beta @> u >>   \bm{M}(\Psi,\Psi).
\end{CD}
\end{array}
$$
Finally the multiplication is compatible with the structural $\mathsf{p}$-maps
$$
\begin{CD}
\bm{M}(\wh{\Psi}',\wh{\Psi}''){_{s}\times_t}\bm{M}(\wh{\Psi},\wh{\Psi}')@> m>> \bm{M}(\wh{\Psi},\wh{\Psi}'')\\
@V\mathsf{p}{_{s}\times_t}\mathsf{p} VV @V \mathsf{p}VV\\
\bm{M}({\Psi}',{\Psi}''){_{s}\times_t}\bm{M}({\Psi},{\Psi}')@> m>> \bm{M}({\Psi},{\Psi}'').
\end{CD}
$$
As in Section \ref{SECX19.1} we can consider transition constructions ${\mathcal F}$, see Definition
\ref{DEFNX19.1.5},  for the transition sets $\bm{M}(\wh{\Psi},\wh{\Psi}')$
and $\bm{M}(\Psi,\Psi')$. Of course, since we have the structural maps $\mathsf{p}$ these constructions
should be in some sense $\mathsf{p}$-compatible, and we shall introduce the needed compatibility requirements.
\begin{definition}\index{D- $\mathsf{P}$-Transtion construction}
Assume $\mathsf{P}:\mathscr{A}\rightarrow \mathscr{B}$ is a finite-to-one covering functor
and   $\bar{F}$ is a uniformizer construction for $\mathsf{P}$.
Given  transition  constructions  $\bar{\mathcal F}$  for sets of the form  $\bm{M}(\wh{\Psi},\wh{\Psi}')$
and ${\mathcal F}$ for sets of the form $\bm{M}(\Psi,\Psi')$ we say that they are {\bf $\mathsf{P}$-compatible} if the following holds.
Given $\wh{h}\in \bm{M}(\wh{\Psi},\wh{\Psi}')$
we have a germ 
$$
\wh{F}_{\wh{h}}:{\mathcal O}(E,s(\wh{h}))\rightarrow (\bm{M}(\wh{\Psi},\wh{\Psi}'),\wh{h}):\wh{u}\rightarrow \wh{F}_{\wh{h}}(\wh{u}),
$$
and for $h\in \bm{M}(\Psi,\Psi')$
$$
F_h:{\mathcal O}(O_\beta,s(h))\rightarrow (\bm{M}(\Psi,\Psi'),h):u\rightarrow F_h(u).
$$
We require that given $\wh{h}\in \bm{M}(\wh{\Psi},\wh{\Psi}')$
the following holds
\begin{eqnarray}\label{EQNX19.2}
\mathsf{p}\circ \wh{F}_{\wh{h}}(\wh{u}) = F_{\mathsf{p}(\wh{h})}\circ \wh{\tau}(\wh{u})\ \text{for}\ \wh{u}\ \text{near to}\ s(\wh{h}).
\end{eqnarray}
\qed
\end{definition}
Each of these transition constructions $\bar{\mathcal F}$ and ${\mathcal F}$ satisfy a list of properties which are stated before Definition \ref{DEFNX19.1.5}.
The fact that $\mathscr{A}$ and $\mathscr{B}$ are related by the covering functor $\mathsf{P}$ requires us to demand
a compatibility between $\bar{\mathcal F}$ and ${\mathcal F}$, which we stated above.   
Recall that $\wh{\tau}:E\rightarrow G\ltimes O_\beta$ on the object level is a local sc-diffeomorphism and it maps sc-diffeomorphically the
neighborhood germs
$$
{\mathcal O}(E,s(\wh{h}))\rightarrow {\mathcal O}(O_\beta,s(\mathsf{p}(\wh{h}))).
$$
Hence it relates with $h=\mathsf{p}(\wh{h})$ the germ $\wh{F}_{\wh{h}}$ with the germ $F_h$. 
Recall that $\wh{F}_{\wh{h}}$ and $F_h$ can be written as follows.
\begin{eqnarray}
F_h(u)=(u,\Phi_{h,u},f_h(u))\ \ \text{and}\ \  \wh{F}_{\wh{h}}(\wh{u})=(\wh{u},\wh{\Phi}_{\wh{h},\wh{u}},\wh{f}_{\wh{h}}(\wh{u}).
\end{eqnarray}
 We explore the equality given in (\ref{EQNX19.2})
and infer that for $\wh{u}$ near $\wh{o}=s(\wh{h})$ with $f_h = t\circ F_h$ and $\wh{f}_{\wh{h}}=t\circ \wh{F}_{\wh{h}}$
\begin{eqnarray}
\mathsf{P}(\wh{\Phi}_{\wh{u}}) = \Phi_{\wh{\tau}(\wh{u})}\ \ \text{and}\ \  \wh{\tau}\circ  \wh{f}_{\wh{h}}(\wh{u})= f_h\circ \wh{\tau}(\wh{u}).
\end{eqnarray}
This also implies  that the (germs of) sc-diffeomorphisms  $f_h$ are uniquely determined by the $\wh{f}_{\wh{h}}$. 
Indeed, for $a\in\wh{\tau}^{-1}(o_\beta)$ the restriction $\tau_a:(O_a,a)\rightarrow (O_\beta,o_\beta)$
of $\wh{\tau}$ is a sc-diffeomorphism. Hence, picking an $a\in \wh{\tau}^{-1}(o_\beta)$ we can write
$$
f_h(u) = \wh{\tau}\circ \wh{f}_{\wh{h}}\circ \tau_a^{-1}(u)
$$
for $u$ near $s(h)$, where $h=\mathsf{p}(\wh{h})$.
Moreover,  we must have
\begin{eqnarray}
\Phi_{h,u} = \mathsf{P}(\wh{\Phi}_{\wh{h},\tau_a^{-1}(u)})\ \ \text{for}\ \ u\ \text{near}\ s(h)
\end{eqnarray}
and therefore the germ $F_h$ is given by
\begin{eqnarray}
F_h(u) =(u,\mathsf{P}(\wh{\Phi}_{\wh{h},\tau_a^{-1}(u)}), \wh{\tau}\circ \wh{f}_{\wh{h}}\circ \tau_{a}^{-1}(u)).
\end{eqnarray}
Having the transition construction for $\mathsf{P}$ we obtain a natural topology $\bar{\mathcal T}$ for $|\mathscr{A}|$
and a natural topology ${\mathcal T}$ for $|\mathscr{B}|$, see Section \ref{SECNX19.2}.
Further $\bm{M}(\Psi,\Psi')$ and 
$\bm{M}(\wh{\Psi},\wh{\Psi}')$ obtain natural metrizable topologies as discussed in Section \ref{SECXN19.3}, and natural M-polyfold structures as shown in Theorem \ref{THMX19.5.1} of Section  \ref{SECNX19.5}.
\begin{proposition}
Considering  $\bm{M}(\wh{\Psi},\wh{\Psi}')$ and $\bm{M}(\Psi,\Psi')$ with  their natural M-poly\-fold structures 
the map 
$$
\mathsf{p}:\bm{M}(\wh{\Psi},\wh{\Psi}')\rightarrow \bm{M}(\Psi,\Psi')
$$
 is sc-smooth and
a surjective  local sc-diffeomorphism. Moreover the map
\begin{eqnarray}\label{EQNX19.3}
\Gamma:\bm{M}(\wh{\Psi},\wh{\Psi}')\rightarrow \bm{M}(\Psi,\Psi'){_{s}\times_{\wh{\tau}}} E:
(\wh{o},\wh{\Phi},\wh{o}')\rightarrow (\mathsf{p}(\wh{o},\wh{\Phi},\wh{o}'),\wh{o})
\end{eqnarray}
is a sc-diffeomorphism.
\end{proposition}
\begin{proof}
We already know that $\mathsf{p}$ is surjective. Pick any $\wh{h}\in \bm{M}(\wh{\Psi},\wh{\Psi}')$ and define
$h=\mathsf{p}(\wh{h})$. Defining  $\wh{o}=s(\wh{h})$ and $o=s(h)$ we find an open neighborhoods
$U(\wh{o})$ and $U(o)$ so that the following holds.
\begin{itemize}
\item[(1)]\ $\wh{\tau}:U(\wh{o})\rightarrow U(o)$ is a sc-diffeomorphism.
\item[(2)]\  $\bar{F}_{\wh{h}}:U(\wh{o})\rightarrow \bm{M}(\wh{\Psi},\wh{\Psi}')$ is a sc-diffeomorphism
onto an open neighborhood $U(\wh{h})$.
\item[(3)]\ $F_h:U(o)\rightarrow \bm{M}({\Psi},{\Psi}')$ is a sc-diffeomorphism onto an open neighborhood
$U(h)$.
\end{itemize}
From the property $\mathsf{p}\circ \bar{F}_{\wh{h}}(\wh{u})=F_h\circ \wh{\tau}(\wh{u})$ we see
that 
$$
\mathsf{p} (\wh{k}) = F_h\circ \wh{\tau}\circ \bar{F}^{-1}_{\wh{h}}(\wh{k})\ \ \wh{k}\in U(\wh{h}).
$$
This implies that the image of $\mathsf{p}$ equals $U(h)$ and that as a composition of sc-diffeomorphisms 
the map $\mathsf{p}:U(\wh{h})\rightarrow U(h)$ is a sc-diffeomorphism.   
Since $s:\bm{M}(\Psi,\Psi')\rightarrow O_\beta$ is a local sc-diffeomorphism  it follows that 
$\bm{M}(\Psi,\Psi'){_{s}\times_{\wh{\tau}}} E$ has a natural M-polyfold structure which is obtained
from the fact that $\bm{M}(\Psi,\Psi'){_{s}\times_{\wh{\tau}}} E$ is a sub-M-polyfold
of $\bm{M}(\Psi,\Psi')\times E$. The map in (\ref{EQNX19.3}) is a bijection and 
since $s: \bm{M}(\wh{\Psi},\wh{\Psi}')\rightarrow E$ is sc-smooth and $\mathsf{p}$ is sc-smooth 
we infer that the map in (\ref{EQNX19.3}) is sc-smooth.  To complete the proof it suffices to show that this map is a local sc-diffeomorphism.  Pick $\wh{h}\in \bm{M}(\wh{\Psi},\wh{\Psi}')$ and define $h=\mathsf{p}(\wh{h})$.
Introducing $o=s(h)$ and $\wh{o}=s(\wh{h})$ we  find sufficiently small open neighborhoods $U(\wh{o})$, $U(\wh{h})$, $U(o)$, and $U(h)$ such that $\wh{\tau}:U(\wh{o})\rightarrow U(o)$, $F_h:U(o)\rightarrow U(h)$, and
$\wh{F}_{\wh{h}}:U(\wh{o})\rightarrow U(\wh{h})$  are  sc-diffeomorphisms. We note that the following map
is a local sc-diffeomorphism
$$
\Theta: U(o)\rightarrow \bm{M}(\Psi,\Psi'){_{s}\times_{\wh{\tau}}} E: u\rightarrow \Theta(u):=(F_h(u),\wh{\tau}^{-1}(u)).
$$
Then 
$$
\Gamma^{-1}\circ \Theta(u) = \wh{F}_{\wh{h}}\circ \wh{\tau}^{-1}(u),
$$
 which shows that 
$\Gamma^{-1}$ near $(h,\wh{o})$ is sc-smooth. This completes the proof.
\qed \end{proof}

In order to obtain polyfold structures for $\mathscr{A}$ and $\mathscr{B}$, it is important that the natural topologies
satisfy the properness property,  and that the topologies are paracompact. Since the data for $\mathscr{A}$ and 
$\mathscr{B}$ is intertwined it follows that the properness property for ${\mathcal T}$ and $\bar{\mathcal T}$
are equivalent requirements, and the same holds for the paracompactness property.

\begin{lemma}\index{L- Continuity and openness of $\vert\mathsf{P}\vert$}
Let $|\mathscr{A}|$ and $|\mathscr{B}|$ be equipped with the natural topologies. 
Then the map $|\mathsf{P}|$ is continuous and open.
\end{lemma}
\begin{proof}
Let $|\beta|\in |\mathscr{B}|$ be given.  Take a representative $\beta$
and a proper covering uniformizer giving the diagram
\begin{eqnarray}\label{ERTP}
\begin{CD}
E @>\wh{\Psi} >> \mathscr{A}\\
@V\wh{\tau}VV   @V\mathsf{P} VV\\
G\ltimes O_\beta @>\Psi>>   \mathscr{B}
\end{CD}
\end{eqnarray}
Passing to orbit spaces we obtain the following  commutative diagram, where the horizontal
maps are topological embeddings
$$
\begin{CD}
|E|  @> |\wh{\Psi}| >> |\mathscr{A}|\\
@V |\wh{\tau}|  VV   @V |\mathsf{P}| VV\\
_G\backslash O_\beta @>|\Psi| >>  |\mathscr{B}|
\end{CD}
$$
Take an open neighborhood $U(|\beta|)$ contained in the image of $|\Psi|$
and let $U(|o_\beta|)$ be the preimage under $|\Psi|$. Since $|\wh{\tau}|$ is continuous 
we see that 
$$
|\mathsf{P}|^{-1}(U(|\beta|))= |\wh{\Psi}|(|\wh{\tau}|^{-1}(U(|o_\beta|)
$$
which is an open set.  This proves continuity.

In order to obtain the openness we just need to show that $|\wh{\tau}|$ is open.
Indeed, if $|\alpha|\in |\mathscr{A}|$ we can take a representative $\alpha$ and define $\beta=\mathsf{P}(\alpha)$.
We need to show that the image of an arbitrarily small open neighborhood of $|\alpha|$ is open.
We take a proper covering uniformizer, which as usual is $\wh{\Psi}$ with underlying $\Psi$.
We take $a\in E$ with $\wh{\Psi}(a)=\alpha$ and take a small open neighborhood 
$U(a)\subset E$ so that $\wh{\tau}:U(a)\rightarrow U(b)$ is an sc-diffeomorphism 
for a suitable open neighborhood $U(o_\beta)$. By a standard fact $|\wh{\Psi}(U(a))|$ and $|\Psi(U(o_\beta))|$ are open.
Now we notice that
$$
|\mathsf{P}|(|\wh{\Psi}(U(a))|)= |\Psi(U(o_\beta))|,
$$
which follows from the commutative diagram (\ref{ERTP}).
\qed \end{proof}
We shall need some more information about $|\mathsf{P}|$.
\begin{lemma}\label{LEMMX19.7.10}
For every  $z\in |\mathscr{B}|$ and given open neighborhoods $U(z)$ and $\wh{U}$ of $|\mathsf{P}|^{-1}(z)$
there exists  a proper covering model
$\wh{\tau}:E\rightarrow G\ltimes O_\beta$  with $|\Psi(o_\beta)|=z$ so that the following holds.
\begin{itemize}
\item[{\em (1)}]\  There exists an open $G$-invariant neighborhood $V(o_\beta)$ in $O_\beta$ having the properness property and satisfying  $|\Psi(V(o_\beta))|\subset U(z)$.
\item[{\em (2)}] \    There exist points $a_1,...a_k\in \wh{\tau}^{-1}(o_\beta)$ and open neighborhoods $V(a_i)$ invariant 
under the $G_{a_i}$-action and having the properness property so that the $|\wh{\Psi}(V(a_i))|$ are mutually 
disjoint and 
\begin{eqnarray}\label{KLOPs9}
|\mathsf{P}|^{-1}(V(z))=  \bigsqcup_{i=1}^k |\wh{\Psi}(V(a_i))| \subset \wh{U}.
\end{eqnarray}
\end{itemize}
\end{lemma}
\begin{proof}
Given $z$ we pick an object $\beta$ representing it and a proper covering model $\wh{\tau}$.
Pick a maximal collection of points $a_1,...,a_k\in E$ with $\wh{\tau}(a_i)=b$
so that any two of them are not connected by an isomorphism. Then denote by $a_{k+1},...a_{k+\ell}$
the other preimages. 
We can find for $i\in \{1,...,k\}$ arbitrarily small open neighborhoods
$V(a_i)$ invariant under the $G_{a_i}$-action so that the restrictions of $\wh{\tau}$ are sc-diffeomorphisms
onto suitable open neighborhoods of $o_\beta$. We can arrange, by taking possibly smaller sets, that their images
are all the same, and using the properness property for $E$ that there are no morphisms 
starting in $V(a_i)$ and ending in $V(a_j)$ for $i\neq j$, where $i,j\in \{1,...,k\}$. Denote the common image
of these sets by $U(o_\beta)$.  Then we have the sc-diffeomorphisms $\wh{\tau}:V(a_i)\rightarrow U(o_\beta)$
and the sets $|\wh{\Psi}(V(a_i))|$ are mutually disjoint.  Perhaps, passing to smaller sets we may assume
that
$$
|\wh{\Psi}(V(a_i))|\subset \wh{U}\ \text{for}\ i\in \{1,...,k\},\ \text{and}\ \ |\Psi(U(o_\beta))|\subset U(z).
$$
It holds for every $i\in \{1,...,k\}$ that
\begin{eqnarray}\label{EQNX19.5}
|\mathsf{P}| ( |\wh{\Psi}(V(a_i))|) = |\Psi|(|\wh{\tau}(V(a_i))|)=|\Psi|(U(o_\beta))=: V(z).
\end{eqnarray}
  With the definitions of $V(a_i)$  and $V(z)$ in (\ref{EQNX19.5})
it follows that
$$
|\mathsf{P}|^{-1}(V(z)) = \bigcup_{i=1}^k |\wh{\Psi}(V(a_i))|\subset \wh{U}.
$$
\qed \end{proof}
Next we are ready to give a result relating the paracompactness for the natural topologies 
on $|\mathscr{A}|$ and $|\mathscr{B}|$ as well as the properness properties.

\begin{proposition}
Assume that $\mathsf{P}:\mathscr{A}\rightarrow \mathscr{B}$ is a proper covering functor, and we are given a uniformizer
construction $\bar{F}$ as well as a transition construction $\bar{\mathcal F}$ with underlying $F$ and ${\mathcal F}$. Then for the natural topologies it holds that ${\mathcal T}$ has the properness property if and only $\bar{\mathcal T}$ has the properness property. 
\end{proposition}
\begin{proof}
  Assume first that $\bar{\mathcal T}$ has the properness property.
Take any uniformizer $\wh{\Psi}:E\rightarrow \mathscr{A}$ for $\mathsf{P}$. Recall that the object space $E$ is a disjoint union
of finitely many $O_{a_1},...,O_{a_k}$ each of which admits the natural $G_{a_i}$-action. Part of the structure
are morphisms which relate the different $O_{a_i}$.  In any case we can consider a $\wh{\Psi}|O_{a_i}$.
By assumption, given $\wh{o}\in O_{a_i}$ there exists an open neighborhood $U(\wh{o})$ such that  
a sequence $(\wh{q}_\ell)\subset U(\wh{o})$ with $|\wh{\Psi}(\wh{q}_\ell)|\rightarrow \wh{z}$ in $(|\mathscr{A}|,\bar{\mathcal T})$  must have a converging subsequence in $O_{a_i}$, say without loss of generality
$\wh{q}_\ell\rightarrow \wh{q}$ with $|\wh{\Psi}(\wh{q})|=\wh{z}$.

By taking $U(\wh{o})$ small enough and assuming that $\wh{o}\in O_{a_i}$ we find an open neighborhood
$U(o)$, $o=\wh{\tau}(\wh{o})$ such that $\wh{\tau}:U(\wh{o})\rightarrow U(o)$ is a sc-diffeomorphism.
It suffices to show that 
$$
\cl_{|\mathscr{B}|}(|\Psi(U(o))|)\subset |\Psi(O_\beta)|.
$$
 Take a sequence 
$(z_j)\subset |\Psi(U(o))|$ converging to some $z\in |\mathscr{B}|$. We can pick $b_j\in U(o)$ with $|\Psi(b_j)|=z_j$
and find $\wh{b}_j\in U(\wh{o}_i)$ with $\wh{\tau}(\wh{b}_j)=b_j$. It holds 
$$
|\mathsf{P}|(|\wh{\Psi}(\wh{b}_j)|) = |\Psi(b_j)|\rightarrow z.
$$
If we can show that $(|\wh{\Psi}(\wh{b}_j)|)$ has a convergent subsequence, say without loss 
converges to some $\wh{z}\in|\mathscr{A}|$, it follows from the assumption  that $\wh{z}\in |\wh{\Psi}(O_{a_i})|$.
This would imply that $z\in |\Psi(O_\beta)|$ and the proof would be completed.
We take a proper covering model around a $\beta'$ representing $z$ and for large $j$
we can express the given data with respect to this model. Employing 
Lemma \ref{LEMMX19.7.10} we infer that a subsequence of $(|\wh{\Psi}(b_j)|)$ converges and the proof
of the first part of the assertion is complete.

Assume next that ${\mathcal T}$ has the properness property.  Given a proper covering model
and uniformizers we find for given $o\in O_\beta$ an open neighborhood $U(o)$
such that $\cl_{|\mathscr{B}|}(|\Psi(U(o))| )\subset |\Psi(O_\beta)|$.  We can take $U(o)$ so small
that we have for every $a\in \wh{\tau}^{-1}(o)$ an open neighborhood $U(a)$ so that
$\wh{\tau}:U(a)\rightarrow U(o)$ is an sc-diffeomorphism. We show that given such a $U(a)$
it holds that $\cl_{|\mathscr{A}|}(|\wh{\Psi}(U(a))|)\subset |\wh{\Psi}(E)|$. Assuming that
$a_j\in U(a)$ and $|\wh{\Psi}(a_j)|\rightarrow \wh{z}$ in $|\mathscr{A}|$ we see that 
$$
|\Psi (\wh{\tau}(a_j))|=|\mathsf{P}|(|\wh{\Psi}(a_j)|)\rightarrow |\mathsf{P}|(\wh{z})=:z.
$$
By assumption $z\in |\Psi(O_\beta)|$ and it follows from Lemma \ref{LEMMX19.7.10} that 
for large $j$ the points $a_j$ belong to an arbitrary neighborhood
of $\mathsf{P}^{-1}(z)$. This implies our assertion since then a subsequence 
of the $a_j$ has to converge to an element $a'\in E$ with $|\wh{\Psi}(a')|=\wh{z}$.
This completes the proof of the properness assertion.
\qed \end{proof}

\begin{proposition}
Assume that $\mathsf{P}:\mathscr{A}\rightarrow \mathscr{B}$ is a proper covering functor, and we are given a uniformizer
construction $\bar{F}$ as well as a transition construction $\bar{\mathcal F}$ with underlying $F$ and ${\mathcal F}$. Then ${\mathcal T}$ is paracompact  if and only $\bar{\mathcal T}$ is paracompact.
\end{proposition}
\begin{proof}
The basic input we shall need is Lemma \ref{LEMMX19.7.10}. 
Assume first that ${\mathcal T}$ is paracompact.  We start with an open covering 
${(U_{\lambda})}_{\lambda\in\Lambda}$ of $|\mathscr{A}|$. In a first step we take a suitable refinement
which has the following property. For every $z\in |\mathscr{B}|$ the preimage $|\mathsf{P}|^{-1}(z)$ is given by finitely many points $\wh{z}_1,...,\wh{z}_k\in |\mathscr{A}|$. Applying the lemma
we find $V(z)$ and $V_z(\wh{z}_i):= |\wh{\Psi}(V(a_i))|$ such that
$$
|\mathsf{P}|^{-1}(V(z))=\sqcup _{i=1}^k V_z(\wh{z}_i)
$$
and every $V_z(\wh{z}_i)$ is contained in one of the sets $U_\lambda$. Then collection 
of all $V_z(\wh{z}_i)$ with $z$ varying in $|\mathscr{B}|$ and $i\in \{1,...,k_z\}$ is an open covering of $|\mathscr{B}|$ and  a refinement of $(U_\lambda)$. The collection ${(V(z))}_{z\in |\mathscr{B}|}$ is an open covering 
of $|\mathscr{B}|$ and since its topology ${\mathcal T}$ is paracompact we find a locally finite 
open covering with the same index set (many sets may be empty) denoted by ${(V_z)}$ such that
$V_z\subset V(z)$.  Taking the preimage by $|\mathsf{P}|$ we obtain $V_{z,\wh{z}_i}$ with $i\in \{1,...,k_z\}$.
This collection is clearly a refinement of $(U_\lambda)$. Let us first show that $(V_{z,\wh{z}_i})$
is an open covering. The openness is clear, and  moreover
\begin{eqnarray*}
|\mathscr{A}|&=&|\mathsf{P}|^{-1}(|\mathscr{B}|) \\
&=& \bigcup_{z} |\mathsf{P}|^{-1}(V_z)\\
&=& \bigcup_{z,\wh{z}_i} V_{z,\wh{z}_i}
\end{eqnarray*}
Next let us show that the covering is locally finite.  For this let $\wh{o}\in |\mathscr{A}|$
and define $o:=|\mathsf{P}|(\wh{o})$. We find an open and sufficiently small  neighborhood $W(o)$
such that there are only finitely many $z$ with $V_z\cap W(o)\neq \emptyset$.
Take  $W(\wh{o}):=|\mathsf{P}|^{-1}(W(o))$ as the open neighborhood of $\wh{o}$.
If $V_{z,\wh{z}_i}\cap W(\wh{o})\neq \emptyset$ then $W(o)\cap V_z\neq \emptyset$.
This implies that there are finitely many such sets. This shows that $\bar{\mathcal T}$ is paracompact.

Next assume that $\bar{\mathcal T}$ is paracompact.  Let $(U_\lambda)$ be an open covering of $|\mathscr{B}|$.
Again we apply Lemma \ref{LEMMX19.7.10} and find $(V(z))_{z\in |\mathscr{B}|}$ which is a refinement
of $(U_\lambda)$ and has the property 
$$
|\mathsf{P}|^{-1}(V(z))= \sqcup_{i=1}^{k_z}  V_z(\wh{z}_i).
$$
Then the collection of all $V_z(\wh{z}_i)$ is an open covering of $|\mathscr{A}|$. 
Using the same index set we find open $V_{z,\wh{z}_i}\subset V_z(\wh{z}_i)$ such that
the new collection $(V_{z,\wh{z}_i})$ is locally finite. Then for every $z$ and $i\in \{1,...,k_z\}$
the sets $|\mathsf{P}|(V_{z,\wh{z}_i})\subset V(z)$ are open. Define 
$$
W_{z,i}:= |\mathsf{P}|(V_{z,\wh{z}_i})\ \ z\in |\mathscr{B}|,\ i\in \{1,...,k_z\}.
$$
Then the whole collection $(W_{z,i})$ is an open covering and refines $(U(z))$.
We finish the proof by showing that the collection is locally finite. For any  point $o\in |\mathscr{B}|$
we find an open neighborhood $\wh{W}:=\wh{W}(|\mathsf{P}|^{-1}(o))$ which is only intersected
by finitely many $V_{z,\wh{z}_i}$ and denote by $W$ the image of $\wh{W}$ under $|\mathsf{P}|$.
  If for some $(z,i)$ it holds that $W_{z,i}\cap W\neq \emptyset$ it holds that $V_{z,\wh{z}_i}\cap \wh{W}\neq \emptyset$.
  However, there are only finitely many indices $(z,i)$ for which the latter holds. This completes the proof
  of the second part.
\qed \end{proof}

Finally we show among other things that $\mathsf{P}$ is a sc-smooth functor
\begin{theorem}
Assume that $\mathsf{P}:\mathscr{A}\rightarrow \mathscr{B}$ is a proper covering functor 
and $\bar{F}$ and $\bar{\mathcal F}$ are a uniformizer construction  as well as a transition construction for $\mathsf{P}$.
Suppose further that the natural topology for $|\mathscr{A}|$ or for $|\mathscr{B}|$ has the properness property and is paracompact.  Then $\mathscr{A}$ and $\mathscr{B}$ have natural polyfold structures. 
The local representatives of $\mathsf{P}$ with respect to a proper covering model $\wh{\tau}:E\rightarrow G\ltimes O_\beta$
is $\wh{\tau}$. Hence $\mathsf{P}$ is an sc-smooth functor.  Moreover, for $\wh{\Psi}$ and $\wh{\Psi}'$
with $\wh{\Gamma}:\bm{M}(\wh{\Psi},\wh{\Psi}')\rightarrow \bm{A}:(e,\wh{\Phi},e')\rightarrow \wh{\Phi}$
and $\Gamma:\bm{M}(\Psi,\Psi')\rightarrow \bm{B}:(o,\Phi,o')\rightarrow \Phi$ the natural maps
we have the commutative diagram
$$
\begin{CD}
\bm{A} @> \mathsf{P}\times s  >>  \bm{B}{_{s}\times_\mathsf{P}} A\\
@A\wh{\Gamma} AA      @A\Gamma\times \Psi AA\\
\bm{M}(\wh{\Psi},\wh{\Psi}')@>\mathsf{p}\times s>>   \bm{M}(\Psi,\Psi') {_{s}\times_{\wh{\tau}}} E,
\end{CD}
$$
where the image of $\mathsf{P}\times s$ restricted to the image of $\wh{\Gamma}$ is the image of
$\Gamma\times\wh{\tau}$. Hence $\mathsf{P}\times s:\bm{A}\rightarrow \bm{B}{_{s}\times_{\mathsf{P}}}A$
is an sc-smooth functor.
\end{theorem}
\begin{proof}
We have shown that the paracompactness for ${\mathcal T}$ is equivalent to that of $\bar{\mathcal T}$. The same holds for the properness property.  Hence, the assumption implies that our data equips $\mathscr{A}$ as well as $\mathscr{B}$ with a polyfold structure. The proper covering functor $\mathsf{P}:\mathscr{A}\rightarrow \mathscr{B}$ 
is with respect to any $\wh{\Psi}$ and the underlying $\Psi$ in the local proper covering model $\wh{\tau}:E\rightarrow G\ltimes O_\beta$ represented by $\wh{\tau}$, which follows immediately from the commutative diagram
$$
\begin{CD}
E @> \wh{\Psi} >> \mathscr{A}\\
@V \wh{\tau} VV    @V\mathsf{P} VV\\
G\ltimes O_\beta @> \Psi >>   \mathscr{B}.
\end{CD}
$$
Hence $\mathsf{P}$ is an sc-smooth functor locally (with respect $\wh{\Psi}$) representable by a proper covering functor between ep-groupoids. Given $\wh{\Psi},\wh{\Psi}'$
consider
$$
\begin{CD}
\bm{A} @> \mathsf{P}\times s  >> \bm{B}{_{s}\times_\mathsf{P}} A\\
@A\wh{\Gamma} AA      @A\Gamma\times\Psi AA\\
\bm{M}(\wh{\Psi},\wh{\Psi}')@>\mathsf{p}\times s>>   \bm{M}(\Psi,\Psi') {_{s}\times_{\wh{\tau}}} E,
\end{CD}
$$
where the first vertical arrow maps $(e,\wh{\Phi},e')$ to $\wh{\Phi}$, and the second vertical arrow maps
$((o,\Phi,o'),e)$ to $(\Phi,\wh{\Psi}(e))$. Here the bottom horizontal arrow is a sc-diffeomorphism which implies that the top arrow   is sc-smooth represented on the 
image of $\Gamma$. 
\qed
 \end{proof}

\section{Covering Constructions for Strong Bundles}
We shall finish this chapter by discussing strong bundle constructions associated to $\mathsf{P}:\mathscr{A}\rightarrow \mathscr{B}$.
The different ingredients have already been discussed in previous sections and the main point of our consideration is to point out
the compatibility requirement between these different constructions, when they have to fit together.

The  starting point is a finite-to-one covering functor $\mathsf{P}:\mathscr{A}\rightarrow \mathscr{B}$ 
and $\text{Ban}$-valued functors 
$$
\mu_{\mathscr{A}}:\mathscr{A}\rightarrow \text{Ban}\ \ \text{and}\ \ 
\mu_{\mathscr{B}}:\mathscr{B}\rightarrow \text{Ban}.
$$
 Recall that $\text{Ban}$ is the category of Banach space
with the morphisms being topological isomorphisms.  
We assume these functors to be related via a natural transformation
\begin{eqnarray}
\Theta: \mu_{\mathscr{A}}\rightarrow \mu_{\mathscr{B}}\circ\mathsf{P}.
\end{eqnarray}
The natural transformation provides us for every object $\alpha$ in $\mathscr{A}$ with a topological linear isomorphism
\begin{eqnarray}
\Theta_\alpha :\mu_{\mathscr{A}}(\alpha)\rightarrow \mu_{\mathscr{B}}(\mathsf{P}(\alpha)).
\end{eqnarray}
In addition, for every morphism $\wh{\Phi}$ in $\mathscr{A}$, say $\wh{\Phi}:\alpha\rightarrow \alpha'$,
we obtain the commutative diagram
$$
\begin{CD}
\mu_{\mathscr{A}}(\alpha)@>\Theta_{\alpha}>>  \mu_{\mathscr{B}}(\mathsf{P}(\alpha))\\
@V \mu_{\mathscr{A}}(\wh{\Phi}) VV    @V \mu_{\mathscr{B}}(\mathsf{P}(\wh{\Phi})) VV \\
\mu_{\mathscr{A}}(\alpha') @>\Theta_{\alpha'} >>  \mu_{\mathscr{A}}(\mathsf{P}(\alpha')).
\end{CD}
$$
Denote by $\mathscr{E}_{\mathscr{A}}$ the obvious bundle category with objects $(\alpha,\wh{v})$,
where $\wh{v}\in \mu_{\mathscr{A}}(\alpha)$ and $\alpha$ is an object in $\mathscr{A}$.  The morphisms in $\mathscr{E}_{\mathscr{A}}$
are the pairs $(\wh{\Phi},\wh{v})$ with $\wh{v}\in \mu_{\mathscr{A}}(s(\wh{\Phi}))$ and $\wh{\Phi}$ is a morphism in $\mathscr{A}$ 
$$
 (s(\wh{\Phi}),\wh{v})\xrightarrow{(\wh{\Phi},\wh{v})} (t(\wh{\Phi}),\mu_{\mathscr{A}}(\wh{\Phi})(\wh{v})).
$$
Similarly we define   $\mathscr{E}_{\mathscr{B}}$ to consist of the objects $(\beta,v)$ with morphisms $(\Phi,v)$
$$
(s({\Phi}),{v})\xrightarrow{(\Phi,v)} (t({\Phi}),\mu_{\mathscr{B}}(\Phi)(v)).
$$
We define a functor $\wh{\mathsf{P}}:\mathscr{E}_{\mathscr{A}}\rightarrow \mathscr{E}_{\mathscr{B}}$ covering
$\mathsf{P}:\mathscr{A}\rightarrow \mathscr{B}$ by
\begin{eqnarray}
&\wh{\mathsf{P}}(\alpha,\wh{v}) =(\mathsf{P}(\alpha), \Theta_\alpha(\wh{v}))&\\
&\wh{\mathsf{P}}(\wh{\Phi},\wh{v})=({\mathsf{P}(\wh{\Phi})},\Theta_{s(\wh{\Phi})}(\wh{v}))&\nonumber\\
&P_{\mathscr{A}}(\alpha,\wh{v})=\alpha&\nonumber\\
&P_{\mathscr{B}}(\beta,v)=\beta&\nonumber
\end{eqnarray}
This fits into the following commutative diagram, abbreviated by ${\bf E}$
\begin{eqnarray}\label{EQNX19.12}
\begin{array}{ccc}
{\bf E}:\ \ \ \ \ \ &\ \ \ \ \ 
\begin{CD}
\mathscr{E}_{\mathscr{A}}@>\wh{\mathsf{P}}>> \mathscr{E}_{\mathscr{B}}\\
@V P_{\mathscr{A}} VV    @V P_{\mathscr{B}}VV\\
\mathscr{A} @>\mathsf{P} >>  \mathscr{B}
\end{CD}
\ \ \ \ \ & \ \ \ \ \
\begin{CD}
(\alpha,\wh{v}) @>\wh{\mathsf{P}}>>  (\mathsf{P}(\alpha),\Theta_{\alpha}(\wh{v}))\\
@V P_{\mathscr{A}}VV @V P_{\mathscr{B}}VV\\
\alpha @>\mathsf{P}>>   \mathsf{P}(\alpha),
\end{CD}
\end{array}
\end{eqnarray}
where we note that for fixed $\alpha$ the map $\mu_{\mathscr{A}}(\alpha)\ni\wh{v}\rightarrow \Theta_\alpha(\wh{v})\in \mu_{\mathscr{B}}(\mathsf{P}(\alpha))$ is a topological linear isomorphism. The morphisms are mapped as follows
\begin{eqnarray*}
\begin{CD}
(\wh{\Phi},\wh{v}) @>\wh{\mathsf{P}}>>  (\mathsf{P}(\wh{\Phi}),\Theta_{s(\wh{\Phi})}(\wh{v}))\\
@V P_{\mathscr{A}}VV @V P_{\mathscr{B}}VV\\
\wh{\Phi} @>\mathsf{P}>>   \mathsf{P}(\wh{\Phi}).
\end{CD}
\end{eqnarray*}
\begin{definition}\index{D- Bundle covering square}
A {\bf bundle covering square} is given by the diagram ${\bf E}$ in (\ref{EQNX19.12}) associated to the 
finite-to-one covering functor $\mathsf{P}:\mathscr{A}\rightarrow \mathscr{B}$ between groupoidal categories,
and the natural transformation $\Theta:\mu_{\mathscr{A}}\rightarrow \mu_{\mathscr{B}}\circ\mathsf{P}$,
where $\mu_{\mathscr{A}}:\mathscr{A}\rightarrow \text{Ban}$ and $\mu_{\mathscr{B}}:\mathscr{B}\rightarrow \text{Ban}$ are functors.
\qed
\end{definition}

In Section \ref{SECX19.6} we discussed how a general construction will provide 
$P_{\mathscr{A}}:\mathscr{E}_{\mathscr{A}}\rightarrow \mathscr{A}$ and $P_{\mathscr{B}}:\mathscr{E}_{\mathscr{B}}\rightarrow \mathscr{B}$ with the structure of strong bundles over polyfolds.
Previously we also discussed a general construction to equip $\mathsf{P}:\mathscr{A}\rightarrow \mathscr{B}$ with the structure 
of a proper covering functor. 
If these three constructions show a certain compatibility we
obtain a general construction of proper strong bundle covering functors.  

Our aim is to describe a procedure which would equip ${\bf E}$ with a uniformizer
as well as a transition construction. Since we have to consider a certain number of symbols
at the same time we use a notation which differs somewhat from the previous sections. It should not cause any confusion.

The local sc-smooth models for our situation, abbreviated by ${\bf K}$,  take the form 
\begin{eqnarray}\label{EQNX19.13}
\begin{array}{cc}
{\bf K}: \ \ \ \ \ \ &\ \ \ \ \
\begin{CD}
\wh{W} @>{\Xi}>>  W\\
@V \wh{Q}  VV   @V Q VV\\
\wh{E}@>{\tau}>>  E
\end{CD}
\end{array}
\end{eqnarray}
Here $\wh{Q}:\wh{W}\rightarrow \wh{E}$ and $Q:W\rightarrow E$  are  strong bundles over an ep-groupoids.
We assume that $E=G\ltimes O_\beta$, where $O_\beta$ contains the distinguished point $o_\beta$, which is mapped to the object $\beta$.
Moreover, ${\tau}:\wh{E}\rightarrow  E$ is a proper covering functor 
of ep-groupoids induced by the strong covering functor ${\Xi}:\wh{W}\rightarrow W$
between strong bundles. 
The functor  $\Xi:\wh{W}\rightarrow W$ in (\ref{EQNX19.13})  defines a topological linear isomorphism
$$
\wh{Q}^{-1}(\wh{e}) \rightarrow Q^{-1}(\tau(\wh{e})): \wh{w}\rightarrow\Xi_{\wh{e}}(\wh{w}).
$$
A morphism $(\wh{\phi},\wh{w}):\wh{w}\rightarrow \mu_{\wh{E}}(\wh{\phi})(\wh{w})$ 
is mapped by $\Xi$ to 
$$
(\tau(\wh{\phi}),\Xi(\wh{w})):\Xi(\wh{w})\rightarrow \Xi(\mu_{\wh{E}}(\wh{\phi})(\wh{w})).
$$
In diagram form this looks as follows
$$
\begin{CD}
\wh{w} @> (\wh{\phi},\wh{w})>> \mu_{\wh{E}}(\wh{\phi})(\wh{w})\\
@VVV @VVV\\
\Xi(\wh{w}) @>(\tau(\wh{\phi}),\Xi(\wh{w}))>>  \Xi(\mu_{\wh{E}}(\wh{\phi})(\wh{w})).
\end{CD}
$$
Since the following identity holds
\begin{eqnarray}
\Xi(\mu_{\wh{E}}(\wh{\phi})(\wh{w}))=\mu_E(\tau(\wh{\phi}))(\Xi(\wh{w}))
\end{eqnarray}
so that we can rewrite the diagram as
$$
\begin{CD}
\wh{w} @> (\wh{\phi},\wh{w})>> \mu_{\wh{E}}(\wh{\phi})(\wh{w})\\
@VVV @VVV\\
\Xi(\wh{w}) @>(\tau(\wh{\phi}),\Xi(\wh{w}))>>  \mu_E(\tau(\wh{\phi}))(\Xi(\wh{w})).
\end{CD}
$$
\begin{definition}\index{D- Strong bundle covering model}
We shall refer to the diagram ${\bf K}$ described in (\ref{EQNX19.13}) and with the properties 
just listed as a {\bf strong bundle covering model}.
\qed
\end{definition}

Next we can introduce the notion of a uniformizer for ${\bf E}$.

\begin{definition}\index{D- Uniformizer for a bundle covering square}
Given a bundle covering square ${\bf E}$ let $\beta$ be an object in $\mathscr{B}$.
A {\bf uniformizer} for ${\bf E}$ at the object $\beta$,
written as $\bm{\wh{\Gamma}}:{\bf K}\rightarrow {\bf E}$,  consists of a strong bundle covering model ${\bf K}$
and four uniformizers 
\begin{eqnarray}
&\Psi : E\rightarrow \mathscr{B}\ \ \ \text{and}\ \ \ \wh{\Psi}:\wh{E}\rightarrow \mathscr{A}&\\
&\Gamma:W\rightarrow \mathscr{E}_{\mathscr{B}}\ \ \text{and}\ \ \wh{\Gamma}:\wh{W}\rightarrow \mathscr{E}_{\mathscr{A}}&\nonumber
\end{eqnarray}
such that $\beta$ is in the image of $\Psi$ and the following holds:
\begin{itemize}
\item[(1)]\  $(\wh{\Psi},\Psi)$ is a uniformizer for the finite-to-one covering $\mathsf{P}:\mathscr{A}\rightarrow \mathscr{B}$.
\item[(2)] \ $(\Gamma,\Psi)$ is a bundle uniformizer for $P_{\mathscr{B}}:\mathscr{E}_{\mathscr{B}}\rightarrow \mathscr{B}$.
\item[(3)]\   $(\wh{\Gamma},\wh{\Psi})$ is a bundle uniformizer for $P_{\mathscr{A}}:\mathscr{E}_{\mathscr{A}}\rightarrow \mathscr{A}$.
\end{itemize}
Moreover,  the uniformizers described by $\bm{\wh{\Gamma}}$ fit into the commutative diagram
$$
\begin{tikzcd}[row sep=scriptsize, column sep=scriptsize]
& \wh{W} \arrow[dl] \arrow[rr] \arrow[dd] & & \mathscr{E}_{\mathscr{A}} \arrow[dl] \arrow[dd] \\
\wh{E} \arrow[rr, crossing over] \arrow[dd] & & \mathscr{A} \\
& W \arrow[dl] \arrow[rr] & & \mathscr{E}_{\mathscr{B}} \arrow[dl] \\
E \arrow[rr] & &\mathscr{B} \arrow[from=uu, crossing over]\\
\end{tikzcd}
$$
where the  horizontal arrows represent the uniformizers.  
\qed
\end{definition}
 We explain the form of the uniformizers next.
$\Psi$
 is defined on the special ep-groupoid $E=G\ltimes O$ and maps
an object $e\in O$ to an object $\Psi(e)$ in $\mathscr{B}$ and a morphism $\phi$ to $\Psi(\phi)$.
The functor $\Gamma:W\rightarrow \mathscr{E}_{\mathscr{B}}$ 
maps $w\in W$ to $\Gamma(w)$ which has the form
$$
\Gamma(w) = (\Psi(e), A_e(w)),
$$
where $e=Q(w)$ and $A_e:Q^{-1}(e)\rightarrow P_{\mathscr{B}}^{-1}(\Psi(e))$ is a topological linear 
isomorphism. A morphism in $W$ has the form $(\phi,w)$ with $Q(w)=s(\phi)$,
$$
(\phi,w):w\rightarrow \mu_E(\phi)(w),
$$
  and is mapped as
$$
\Gamma(\phi,w)=(\Psi(\phi), A_e(w)).
$$
The functor $\wh{\Psi}:\wh{E}\rightarrow \mathscr{A}$ maps an object $\wh{e}$ to $\wh{\Psi}(\wh{e})$
and a morphism $\wh{\phi}$ to $\wh{\Psi}(\wh{\phi})$. The functor $\wh{\Gamma}$ maps $\wh{w}\in \wh{W}$
to $\wh{\Gamma}(\wh{w})$. A morphism $(\wh{\phi},\wh{w}):\wh{w}\rightarrow \mu_{\wh{E}}(\wh{\phi})(\wh{w})$ 
is mapped to an element of the form
$$
\wh{\Gamma}(\wh{\phi},\wh{w}) =(\wh{\Psi}(\wh{\phi}), \wh{A}_{\wh{e}}(\wh{w})),
$$
where $\wh{A}_{\wh{e}}:\wh{Q}^{-1}(\wh{e})\rightarrow \mathsf{P}_{\mathscr{A}}^{-1}(\wh{\Psi}(\wh{\phi}))$ is a linear topological isomorphism, and
$\wh{Q}(\wh{w})=\wh{e}$.  The following commutative diagram lists the relevant maps on the object level
$$
\begin{tikzcd}[row sep=scriptsize, column sep=scriptsize]
& \wh{w} \arrow[dl] \arrow[rr] \arrow[dd] & & (\wh{\Psi}(\wh{e}),\wh{A}_{\wh{e}}(\wh{w})) \arrow[dl] \arrow[dd] \\
\wh{e}=\wh{Q}(\wh{w}) \arrow[rr, crossing over] \arrow[dd] & & \wh{\Psi}(\wh{e}) \\
& w=\Xi(\wh{w}) \arrow[dl] \arrow[rr] & & (\Psi(e), A_e(w)) =\wh{\mathsf{P}} (\wh{\Psi}(\wh{e}),\wh{A}_{\wh{e}}(\wh{w}))\arrow[dl] \\
e=\tau(\wh{e}) \arrow[rr] & &\Psi(e) \arrow[from=uu, crossing over].\\
\end{tikzcd}
$$
There are the following relationships 
\begin{eqnarray}
& \mathsf{P} \circ \wh{\Psi} (\wh{e}) = \Psi\circ \tau(\wh{e})\ \ \text{for}\ \ \wh{e}\in\wh{E}&\\
& \Theta_{\wh{\Psi}(\wh{e})}\circ \wh{A}_{\wh{e}}(\wh{w}) = A_e\circ \Xi(\wh{w})\ \  \text{for}\ \ \wh{w}\in \wh{W}.\nonumber&\\
& Q\circ \Xi (\wh{w}) = \tau\circ \wh{Q}(\wh{w})\ \ \text{and}\ \ \wh{w}\in \wh{W}.\nonumber
\end{eqnarray}
On the morphism level we obtain the diagram
$$
\begin{tikzcd}[row sep=scriptsize, column sep=scriptsize]
& (\wh{\phi},\wh{w}) \arrow[dl] \arrow[rr] \arrow[dd] & & (\wh{\Psi}(\wh{\phi}),\wh{A}_{\wh{e}}(\wh{w})) \arrow[dl] \arrow[dd] \\
\wh{\phi} \arrow[rr, crossing over] \arrow[dd] & & \wh{\Psi}(\wh{\phi}) \\
& (\phi,w) \arrow[dl] \arrow[rr] & & (\Psi(\phi), A_e(w)) \arrow[dl] \\
\phi=\tau(\wh{\phi}) \arrow[rr] & &\Psi(\phi) \arrow[from=uu, crossing over],\\
\end{tikzcd}
$$
and we have the following relationship
\begin{eqnarray}
&\mathsf{P}\circ \wh{\Psi}(\wh{\phi}) = \Psi\circ \tau(\wh{\phi})\ \ \text{for}\ \ \wh{\phi}\in \wh{E}.&
\end{eqnarray}
Next we consider the transition construction. We already assume that we have 
transition constructions, as described before, for the covering and strong bundle situations.
The main point here is to emphasize the necessary compatibility.
Associated to 
$$
{\bf K} \xrightarrow{\bm{\wh{\Gamma}}} {\bf E} \xleftarrow{\bm{\wh{\Gamma}'}} {\bf K'}
$$
we can take the weak fibered products which produces the diagram

$$
\begin{array}{cc}
\overline{\bm{M}}(\bm{\wh{\Gamma}},\bm{\wh{\Gamma}'})\colon \ \ \ &\ \ \ \ \ 
\begin{CD}
\bm{M}(\wh{\Gamma},\wh{\Gamma}')@>\bm{P} >> \bm{M}(\Gamma,\Gamma')\\
@V\bm{\wh{m}}VV @V\bm{m}VV\\
\bm{M}(\wh{\Psi},\wh{\Psi}')@>\bm{p}>> \bm{M}(\Psi,\Psi')
\end{CD}
\end{array}
$$
From the previous discussions we already understand (separately)
\begin{eqnarray*}
&\bm{M}(\wh{\Psi},\wh{\Psi}')\xrightarrow{\bm{p}}\bm{M}(\Psi,\Psi')\ \ \text{(covering, Section \ref{SECX19.7})}&\\
&\bm{M}(\wh{\Gamma},\wh{\Gamma}')\xrightarrow{\bm{\wh{m}}}\bm{M}(\wh{\Psi},\wh{\Psi}')\ \ \text{(strong bundle, Section \ref{SECX19.6})}&\\
&\bm{M}({\Gamma},{\Gamma}')\xrightarrow{\bm{m}}\bm{M}({\Psi},{\Psi}')\ \ \text{(strong bundle, Section \ref{SECX19.6})}.&
\end{eqnarray*}
The new ingredient is the top horizontal arrow. For each of the four ingredients we assume that we have transition constructions which are related
in a way we shall describe next.  We begin with the  bottom part of the diagram
$$
\bm{M}(\wh{\Psi},\wh{\Psi}')\xrightarrow{\bm{p}}\bm{M}(\Psi,\Psi')
$$
Recall that the morphisms in $\mathscr{A}$ are written as $\wh{\Phi}$ and in $\mathscr{B}$ as $\Phi$.
Given $\wh{h} =(\wh{o},\wh{\Phi},\wh{o}')\in 
\bm{M}(\wh{\Psi},\wh{\Psi}')$ let $h=(o,\Phi,o')=\bm{p}(\wh{h})=(\tau(\wh{o}),\tau(\wh{\Phi}),\tau'(\wh{o}'))$, 
and require that
$$
\wh{F}_{\wh{h}}:\mathscr{O}(\wh{E},s(\wh{h}))\rightarrow(\bm{M}(\wh{\Psi},\wh{\Psi}'),\wh{h}):\wh{u}\rightarrow \wh{F}_{\wh{h}}(\wh{u})
$$
is related to 
$$
F_h:\mathscr{O}(E,s(h))\rightarrow(\bm{M}(\Psi,\Psi'),h):u\rightarrow F_h(u)
$$
via
\begin{eqnarray}
\bm{p}\circ \wh{F}_{\wh{h}}(\wh{u}) = F_{h}\circ \tau(\wh{u})\ \ \text{for}\ \ \wh{u}\ \text{near}\ \ s(\wh{h}).
\end{eqnarray}
Since $\tau$ is a local sc-diffeomorphism $F_h$ is determined by $\wh{F}_{\wh{h}}$.  
We note that this requirement is precisely the one used in the proper covering construction discussed in a previous section.

Next we consider
the strong bundle
$$
\bm{M}(\Gamma,\Gamma')\xrightarrow {\bm{m}} \bm{M}(\Psi,\Psi').
$$
The transition map we shall denote by $\bm{F}_h$, which is a germ fitting into the commutative diagram
$$
\begin{CD}
W|\mathscr{O}(E,o)@> \bm{F}_h >> \bm{M}(\Gamma,\Gamma')\\
@V Q VV   @V\bm{m} VV\\
\mathscr{O}(E,o) @> F_h >>   (\bm{M}(\Psi,\Psi'),h),
\end{CD}
$$
and the following relationship holds
\begin{eqnarray}
\bm{m}\circ \bm{F}_h(w)=F_h\circ Q(w)\ \ \text{with}\ Q(w)\ \text{near}\ o.
\end{eqnarray}
Finally we consider 
$$
\bm{M}(\wh{\Gamma},\wh{\Gamma}')\xrightarrow {\bm{\wh{m}}} \bm{M}(\wh{\Psi},\wh{\Psi}').
$$
Then the relevant germ is $\bm{\wh{F}}_{\wh{h}}$ is
$$
\begin{CD}
\wh{W}|\mathscr{O}(\wh{E},\wh{o}) @> \bm{\wh{F}}_{\wh{h}} >>  \bm{M}(\wh{\Gamma},\wh{\Gamma}')\\
@V \wh{Q} VV    @V \bm{\wh{m}} VV\\
\mathscr{O}(\wh{E},\wh{o}) @> \wh{F}_{\wh{h}}>>  (\bm{M}(\wh{\Psi},\wh{\Psi}'),\wh{h}),
\end{CD}
$$
and the compatibility condition is given by
\begin{eqnarray}
\bm{\wh{m}}\circ \bm{\wh{F}}_{\wh{h}}(\wh{w})=\wh{F}_{\wh{h}}\circ \wh{Q}(\wh{w})\ \ \text{with}\ \wh{Q}(\wh{w})\ \text{near}\ \wh{o}.
\end{eqnarray}
So far the discussion is governed by the material from previous sections.

The additional requirement is a compatibility of $\bm{F}_h$ and $\bm{\wh{F}}_{\wh{h}}$.
Namely we have to require with $h=\bm{p}(\wh{h})$
\begin{eqnarray}
\bm{P}\circ \bm{\wh{F}}_{\wh{h}}(\wh{w}) = \bm{F}_h \circ \Xi (\wh{w})\ \ \text{for}\ \ \wh{Q}(\wh{w})\ \text{near}\ \wh{o}.
\end{eqnarray}
Since $\Xi$ is a local sc-diffeomorphism it follows that $\bm{F}_h$ is determined by $\bm{\wh{F}}_{\wh{h}}$.
From the previous discussion we already know that $\bm{F}_h$ determines $F_h$ and $\bm{\wh{F}}_{\wh{h}}$
determines $\wh{F}_{\wh{h}}$.

\begin{definition}
Let $\mathsf{P}:\mathscr{A}\rightarrow \mathscr{B}$ be a finite-to-one covering construction, and 
$$
\mu_{\mathscr{A}}:\mathscr{A}\rightarrow \text{Ban}\ \ \ \text{and}\ \ \ \mu_{\mathscr{B}}:\mathscr{B}\rightarrow \text{Ban}
$$
and $\Theta:\mu_{\mathscr{A}}\rightarrow \mu_{\mathscr{B}}\circ\mathsf{P}$ a natural transformation.  Denote the associated bundle covering square by ${\bf E}$.
A {\bf uniformizer construction} for ${\bf E}$ is given by a functor
$$
\overline{\bm{F}}: {\mathscr{B}}^-\rightarrow \text{SET},
$$
which associates to an object $\beta$ a set of uniformizers  
$$
\overline{\bm{F}}(\beta)=\left\{\bm{\wh{\Gamma}}\right\}
$$
at the object $\beta$, and 
 where $\bm{\wh{\Gamma}}:{\bf K}\rightarrow {\bf E}$ are defined on strong bundle covering models ${\bf K}$.
 \qed
\end{definition}

Finally, given a uniformizer construction $\bm{\bar{F}}$ for ${\bf E}$ we define a transition construction.
\begin{definition}
Let $\overline{\bm{F}}:\mathscr{B}\rightarrow \text{SET}$ be a uniformizer construction for the bundle covering square ${\bf E}$.
A transition construction $\overline{\bm{M}}$ for $\overline{\bm{F}}$ associates to the  four transition sets associated 
to $\overline{\bm{M}}(\bm{\wh{\Gamma}},\bm{\wh{\Gamma}'})$ the transition constructions $F_h$, $\wh{F}_{\wh{h}}$, $\bm{F}_h$, and
$\wh{\bm{F}}_{\wh{h}}$ with the compatibilities previously described.
\qed
\end{definition}
We obtain the following result.
\begin{theorem}
Assume that $\overline{\bm{F}}$ is a uniformizer construction and $\overline{\bm{M}}$ a transition construction for the bundle covering square ${\bf E}$.
Suppose further that the natural topology ${\mathcal T}$ for $\mathscr{B}$ satisfies the properness condition and it is paracompact.
Then he natural topologies for $|\mathscr{A}|$, $\vert\mathscr{B}\vert$, $\vert \mathscr{E}_{\mathscr{A}}\vert$, and $\vert\mathscr{E}_{\mathscr{B}}\vert$
are all metrizable. Further 
$$
\begin{CD}
\mathscr{E}_{\mathscr{A}}@> \wh{\mathsf{P}}>>\mathscr{E}_{\mathscr{B}}\\
@V P_{\mathscr{A}} VV   @V P_{\mathscr{B}} VV\\
\mathscr{A} @>\mathsf{P} >> \mathscr{B}
\end{CD}
$$
obtains the structure of proper strong bundle functor as described at the end of Section \ref{SECV17.7}.
\qed
\end{theorem}

\backmatter
\bibliographystyle{plain}

\appendix

\printindex


\end{document}